\documentclass[a4paper,10pt,english]{smfbook}
\usepackage{amsmath,amscd,amssymb}
\usepackage[dvips]{color}

\newcommand{\nbiga}{\mathcal{A}}
\newcommand{\nbigb}{\mathcal{B}}
\newcommand{\nbigc}{\mathcal{C}}
\newcommand{\nbigd}{\mathcal{D}}
\newcommand{\nbige}{\mathcal{E}}
\newcommand{\nbigf}{\mathcal{F}}
\newcommand{\nbigg}{\mathcal{G}}
\newcommand{\nbigh}{\mathcal{H}}
\newcommand{\nbigi}{\mathcal{I}}
\newcommand{\nbigj}{\mathcal{J}}
\newcommand{\nbigk}{\mathcal{K}}
\newcommand{\nbigl}{\mathcal{L}}
\newcommand{\nbigm}{\mathcal{M}}
\newcommand{\nbign}{\mathcal{N}}
\newcommand{\nbigo}{\mathcal{O}}
\newcommand{\nbigp}{\mathcal{P}}
\newcommand{\nbigq}{\mathcal{Q}}
\newcommand{\nbigr}{\mathcal{R}}
\newcommand{\nbigs}{\mathcal{S}}
\newcommand{\nbigt}{\mathcal{T}}
\newcommand{\nbigu}{\mathcal{U}}
\newcommand{\nbigv}{\mathcal{V}}
\newcommand{\nbigw}{\mathcal{W}}
\newcommand{\nbigx}{\mathcal{X}}
\newcommand{\nbigy}{\mathcal{Y}}
\newcommand{\nbigz}{\mathcal{Z}}

\newcommand{\proj}{\mathbb{P}}
\newcommand{\seisuu}{\mathbb{Z}}
\newcommand{\rnum}{{\boldsymbol Q}}

\newcommand{\cnum}{{\boldsymbol C}}
\newcommand{\real}{{\boldsymbol R}}
\newcommand{\hyperh}{\mathbb{H}}

\newcommand{\Tate}{\mathbb{T}}

\newcommand{\DD}{\mathbb{D}}
\newcommand{\EE}{\mathbb{E}}

\newcommand{\gbigc}{\mathfrak C}
\newcommand{\gbigd}{\mathfrak D}
\newcommand{\gbige}{\mathfrak E}
\newcommand{\gbigf}{\mathfrak F}

\newcommand{\gbigh}{\mathfrak H}

\newcommand{\gbigk}{\mathfrak K}
\newcommand{\gbigl}{\mathfrak L}
\newcommand{\gbigm}{\mathfrak M}

\newcommand{\gbigs}{\mathfrak S}
\newcommand{\gbigt}{\mathfrak T}
\newcommand{\gbigu}{\mathfrak U}
\newcommand{\gbigv}{\mathfrak V}

\newcommand{\gbigz}{\mathfrak Z}

\newcommand{\gminia}{\mathfrak a}
\newcommand{\gminib}{\mathfrak b}
\newcommand{\gminic}{\mathfrak c}
\newcommand{\gminid}{\mathfrak d}
\newcommand{\gminie}{\mathfrak e}

\newcommand{\gminih}{\mathfrak h}

\newcommand{\gminik}{\mathfrak k}

\newcommand{\gminim}{\mathfrak m}

\newcommand{\gminio}{\mathfrak o}
\newcommand{\gminip}{\mathfrak p}
\newcommand{\gminiq}{\mathfrak q}

\newcommand{\gminis}{\mathfrak s}

\newcommand{\vecxi}{{\boldsymbol \xi}}
\newcommand{\veceta}{{\boldsymbol \eta}}
\newcommand{\vecrho}{{\boldsymbol \rho}}
\newcommand{\vece}{{\boldsymbol e}}

\newcommand{\vecv}{{\boldsymbol v}}
\newcommand{\vecu}{{\boldsymbol u}}
\newcommand{\vecw}{{\boldsymbol w}}
\newcommand{\vecgamma}{{\boldsymbol \gamma}}
\newcommand{\vecl}{{\boldsymbol l}}
\newcommand{\veczero}{{\boldsymbol 0}}
\newcommand{\vecalpha}{{\boldsymbol \alpha}}
\newcommand{\veca}{{\boldsymbol a}}
\newcommand{\vecb}{{\boldsymbol b}}
\newcommand{\vecbeta}{{\boldsymbol \beta}}
\newcommand{\vecdelta}{{\boldsymbol \delta}}

\newcommand{\vecc}{{\boldsymbol c}}
\newcommand{\vecd}{{\boldsymbol d}}

\newcommand{\veck}{{\boldsymbol k}}
\newcommand{\vecm}{{\boldsymbol m}}

\newcommand{\vecI}{{\boldsymbol I}}
\newcommand{\vecomega}{{\boldsymbol \omega}}
\newcommand{\vecx}{{\boldsymbol x}}
\newcommand{\vecf}{{\boldsymbol f}}
\newcommand{\vecepsilon}{{\boldsymbol \epsilon}}
\newcommand{\vecF}{{\boldsymbol F}}
\newcommand{\vecE}{{\boldsymbol E}}
\newcommand{\vecW}{{\boldsymbol W}}
\newcommand{\vecn}{{\boldsymbol n}}
\newcommand{\vecp}{{\boldsymbol p}}

\newcommand{\vecz}{{\boldsymbol z}}
\newcommand{\vecC}{{\boldsymbol C}}

\newcommand{\vecV}{{\boldsymbol V}}
\newcommand{\vecA}{{\boldsymbol A}}
\newcommand{\vecB}{{\boldsymbol B}}
\newcommand{\vecS}{{\boldsymbol S}}
\newcommand{\vecG}{{\boldsymbol G}}

\newcommand{\llarr}{\longleftarrow}

\newcommand{\rarr}{\rightarrow}
\newcommand{\lrarr}{\longrightarrow}
\newcommand{\darr}{\downarrow}

\newcommand{\pf}{{\bf Proof}\hspace{.1in}}

\def\Hom{\mathop{\rm Hom}\nolimits}
\def\poly{\mathop{\rm poly}\nolimits}
\def\End{\mathop{\rm End}\nolimits}

\def\Cok{\mathop{\rm Cok}\nolimits}

\def\Image{\mathop{\rm Im}\nolimits}

\def\Re{\mathop{\rm Re}\nolimits}

\def\Gr{\mathop{\rm Gr}\nolimits}

\def\rank{\mathop{\rm rank}\nolimits}
\def\Spec{\mathop{\rm Spec}\nolimits}

\def\Ker{\mathop{\rm Ker}\nolimits}
\def\ker{\mathop{\rm Ker}\nolimits}

\def\Gr{\mathop{\rm Gr}\nolimits}
\def\Sym{\mathop{\rm Sym}\nolimits}

\def\ad{\mathop{\rm ad}\nolimits}
\def\Res{\mathop{\rm Res}\nolimits}

\def\ord{\mathop{\rm ord}\nolimits}
\def\degpar{\mathop{\rm par\textrm{-}deg}\nolimits}
\def\parch{\mathop{\rm par\textrm{-}ch}\nolimits}

\def\parchern{\mathop{\rm par\textrm{-}c}\nolimits}

\def\Ric{\mathop{\rm Ric}\nolimits}
\def\tr{\mathop{\rm tr}\nolimits}

\def\vol{\mathop{\rm dvol}\nolimits}
\def\dvol{\mathop{\rm dvol}\nolimits}
\def\Diff{\mathop{\rm Diff}\nolimits}

\def\can{\mathop{\rm can}\nolimits}
\def\var{\mathop{\rm var}\nolimits}
\def\id{\mathop{\rm id}\nolimits}

\def\codim{\mathop{\rm codim}\nolimits}
\def\gap{\mathop{\rm gap}\nolimits}

\def\Supp{\mathop{\rm Supp}\nolimits}

\def\VPTgen{\mathop{\rm VPT}\nolimits}
\def\VPTgenwild{\mathop{\rm VPT^{\wild}}\nolimits}
\def\MPT{\mathop{\rm MPT}\nolimits}

\def\MFL{\mathop{\rm MFL}\nolimits}
\def\RHB{\mathop{\rm RHB}\nolimits}
\def\MF{\mathop{\rm MF}\nolimits}

\def\Irr{\mathop{\rm Irr}\nolimits}

\newcommand{\del}{\partial}
\newcommand{\delbar}{\overline{\del}}

\newcommand{\pardeg}{\degpar}

\newcommand{\nhom}{{\mathcal Hom}}

\newcommand{\nbar}{\underline{n}}

\newcommand{\jbar}{\underline{j}}
\newcommand{\mbar}{\underline{m}}
\newcommand{\kbar}{\underline{k}}

\newcommand{\lbar}{\underline{l}}
\newcommand{\pbar}{\underline{p}}

\newcommand{\nibar}{\underline{2}}

\newcommand{\shikaku}{\sharp}
\newcommand{\sankaku}{\triangle}

\newcommand{\harmonicbundle}{(E,\delbar_E,\theta,h)}

\newcommand{\barz}{\overline{z}}
\newcommand{\zbar}{\barz}
\newcommand{\zetabar}{\overline{\zeta}}
\newcommand{\baralpha}{\overline{\alpha}}
\newcommand{\alphabar}{\baralpha}
\newcommand{\barlambda}{\overline{\lambda}}
\newcommand{\lambdabar}{\barlambda}
\newcommand{\fbar}{\overline{f}}

\newcommand{\etabar}{\overline{\eta}}

\newcommand{\tbar}{\overline{t}}
\newcommand{\xbar}{\overline{x}}
\newcommand{\sbar}{\overline{s}}

\newcommand{\Abar}{\overline{A}}

\newcommand{\Poin}{{\bf p}}
\newcommand{\poin}{\Poin}
\newcommand{\prolong}[1]{{}^{\diamond}{#1}}

\newcommand{\prolongg}[2]{{}_{#1}{#2}}

\newcommand{\DDlambda}{\DD^{\lambda}}
\newcommand{\DDlambdastar}{\DD^{\lambda\,\star}}
\newcommand{\DDlambdahat}{\widehat{\DD}^{\lambda}}
\newcommand{\DDhatlambda}{\DDlambdahat}

\newcommand{\DDlambdaast}{\DD^{\lambda\,\ast}}

\newcommand{\laplacian}{\Delta''}
\newcommand{\Deltalambda}{\Delta^{\lambda}}

\newcommand{\doublelangle}{\langle\langle}
\newcommand{\doublerangle}{\rangle\rangle}

\newcommand{\nbigelambda}{\nbige^{\lambda}}
\newcommand{\nbigelambdazero}{\nbige^{\lambda_0}}

\newcommand{\blowup}{\widetilde}

\newcommand{\nbigxlambda}{\nbigx^{\lambda}}

\newcommand{\nbigdlambda}{\nbigd^{\lambda}}

\newcommand{\KMS}{{\mathcal{KMS}}}

\newcommand{\Par}{{\mathcal Par}}
\newcommand{\Sp}{{\mathcal Sp}}

\newcommand{\kmsmap}{\gminik}
\newcommand{\paramap}{\gminip}
\newcommand{\eigenmap}{\gminie}

\newcommand{\lefttop}[1]{{}^{#1}}

\def\irr{\mathop{\rm irr}\nolimits}
\def\reg{\mathop{\rm reg}\nolimits}

\def\Harm{\mathop{\rm Harm}\nolimits}
\def\nil{\mathop{\rm nil}\nolimits}

\newcommand{\Fzero}{F^{(\lambda_0)}}
\newcommand{\EEzero}{\EE^{(\lambda_0)}}

\newcommand{\Vzero}{V^{(\lambda_0)}}

\newcommand{\psizero}{\psi^{(\lambda_0)}}
\newcommand{\tildepsizero}{\widetilde{\psi}^{(\lambda_0)}}
\newcommand{\tildepsi}{\widetilde{\psi}}
\newcommand{\psitilde}{\tildepsi}

\newcommand{\nbigfzero}{\nbigf^{(\lambda_0)}}
\newcommand{\nbigqzero}{\nbigq^{(\lambda_0)}}

\newcommand{\Uzero}{U^{(\lambda_0)}}

\newcommand{\deldel}{\eth}
\newcommand{\deldelbar}{\overline{\deldel}}

\newcommand{\lamda}{\lambda}
\newcommand{\lambdazero}{\lambda_{0}}

\newcommand{\distribution}{\gbigd\gminib}

\newcommand{\naiveprolong}[1]{\lefttop{\square}{#1}}

\newcommand{\supp}{\gminis}

\newcommand{\openclosed}[2]{]#1,#2]}
\newcommand{\openopen}[2]{]#1,#2[}

\newcommand{\omegatilde}{\widetilde{\omega}}
\newcommand{\vecEhat}{\widehat{\vecE}}

\newcommand{\Ehat}{\widehat{E}}

\newcommand{\vecvhat}{\widehat{\vecv}}

\newcommand{\Ahat}{\widehat{A}}

\newcommand{\Gtilde}{\widetilde{G}}

\newcommand{\Etilde}{\widetilde{E}}
\newcommand{\vecEtilde}{\widetilde{\vecE}}
\newcommand{\thetatilde}{\widetilde{\theta}}
\newcommand{\Vhat}{\widehat{V}}
\newcommand{\nablahat}{\widehat{\nabla}}
\newcommand{\Vtilde}{\widetilde{V}}
\newcommand{\nablatilde}{\widetilde{\nabla}}
\newcommand{\vecvtilde}{\widetilde{\vecv}}
\newcommand{\vecwtilde}{\widetilde{\vecw}}
\newcommand{\vecutilde}{\widetilde{\vecu}}

\newcommand{\vecvbar}{\overline{\vecv}}
\newcommand{\fhat}{\widehat{f}}

\newcommand{\wbar}{\overline{w}}

\newcommand{\nbigetilde}{\widetilde{\nbige}}
\newcommand{\DDlambdatilde}{\widetilde{\DD}^{\lambda}}
\newcommand{\htilde}{\widetilde{h}}
\newcommand{\nbigelambdatilde}{\widetilde{\nbige}^{\lambda}}
\newcommand{\ptilde}{\widetilde{p}}
\newcommand{\pitilde}{\widetilde{\pi}}
\newcommand{\Phitilde}{\widetilde{\Phi}}
\newcommand{\vtilde}{\widetilde{v}}
\newcommand{\phat}{\widehat{p}}

\newcommand{\ftilde}{\widetilde{f}}
\newcommand{\utilde}{\widetilde{u}}
\newcommand{\Ftilde}{\widetilde{F}}
\newcommand{\Phibar}{\overline{\Phi}}
\newcommand{\gminiabar}{\overline{\gminia}}
\newcommand{\atilde}{\widetilde{a}}

\newcommand{\mubar}{\overline{\mu}}

\newcommand{\nbigltilde}{\widetilde{\nbigl}}
\newcommand{\nbigdhat}{\widehat{\nbigd}}
\newcommand{\DDtilde}{\widetilde{\DD}}
\newcommand{\nbigxhat}{\widehat{\nbigx}}
\newcommand{\nbigkhat}{\widehat{\nbigk}}
\newcommand{\Rtilde}{\widetilde{R}}
\newcommand{\Dhat}{\widehat{D}}

\newcommand{\Itilde}{\widetilde{I}}

\newcommand{\vecwhat}{\widehat{\vecw}}

\newcommand{\pihat}{\widehat{\pi}}
\newcommand{\Omegabar}{\overline{\Omega}}
\newcommand{\Fhat}{\widehat{F}}
\newcommand{\nbigfhat}{\widehat{\nbigf}}
\newcommand{\Ghat}{\widehat{G}}

\newcommand{\puebar}{\overline{p}}
\newcommand{\Rhat}{\widehat{R}}
\newcommand{\vecuhat}{\widehat{\vecu}}
\newcommand{\Ohat}{\widehat{O}}
\newcommand{\Sbar}{\overline{S}}
\newcommand{\nbigehat}{\widehat{\nbige}}
\newcommand{\Phihat}{\widehat{\Phi}}

\newcommand{\nbigvtilde}{\widetilde{\nbigv}}
\newcommand{\nbigvhat}{\widehat{\nbigv}}
\newcommand{\ellsitabar}{\underline{\ell}}
\newcommand{\Utilde}{\widetilde{U}}
\newcommand{\Dtilde}{\widetilde{D}}
\newcommand{\Xtilde}{\widetilde{X}}
\newcommand{\Xbar}{\overline{X}}
\newcommand{\vecy}{\boldsymbol y}
\newcommand{\Zhat}{\widehat{Z}}
\newcommand{\What}{\widehat{W}}
\newcommand{\DDhat}{\widehat{\DD}}
\newcommand{\nbigrtilde}{\widetilde{\nbigr}}

\newcommand{\nbigmtilde}{\widetilde{\nbigm}}
\def\Der{\mathop{\rm Der}\nolimits}
\def\Stab{\mathop{\rm Stab}\nolimits}
\def\ord{\mathop{\rm ord}\nolimits}
\def\Herm{\mathop{\rm Herm}\nolimits}
\def\full{\mathop{\rm full}\nolimits}
\def\MT{\mathop{\rm MT}\nolimits}
\def\MTW{\mathop{\rm MTW}\nolimits}
\def\MTN{\mathop{\rm MTN}\nolimits}
\def\wild{\mathop{\rm wild}\nolimits}
\def\Gal{\mathop{\rm Gal}\nolimits}
\def\strict{\mathop{\rm strict}\nolimits}

\def\Can{\mathop{\rm Can}\nolimits}
\def\Var{\mathop{\rm Var}\nolimits}
\def\moderate{\mathop{\rm mod}\nolimits}
\def\op{\mathop{\rm op}\nolimits}
\def\Mod{\mathop{\rm Mod}\nolimits}
\def\ramification{\mathop{\rm rami}\nolimits}

\def\Hol{\mathop{\rm Hol}\nolimits}
\def\Sec{\mathop{\rm Sec}\nolimits}
\def\Sep{\mathop{\rm Sep}\nolimits}

\def\Gys{\mathop{\rm Gys}\nolimits}
\def\exchange{\mathop{\rm exchange}\nolimits}

\def\val{\mathop{\rm val}\nolimits}
\def\SDcat{\mathop{\rm SD}\nolimits}
\def\SDLcat{\mathop{\rm SDL}\nolimits}

\def\Des{\mathop{\rm Des}\nolimits}
\def\hol{\mathop{\rm hol}\nolimits}

\newcommand{\Dbar}{\overline{D}}

\newcommand{\Wtilde}{\widetilde{W}}
\newcommand{\Ptilde}{\widetilde{P}}
\newcommand{\psihat}{\widehat{\psi}}
\newcommand{\iotahat}{\widehat{\iota}}

\newcommand{\vecR}{{\boldsymbol R}}
\newcommand{\vbar}{\overline{v}}

\newcommand{\Tbar}{\overline{T}}
\newcommand{\gminibbar}{\overline{\gminib}}

\newcommand{\nbigehatlambda}{\widehat{\nbige}^{\lambda}}
\newcommand{\nbigelambdahat}{\nbigehatlambda}
\newcommand{\nbigutilde}{\widetilde{\nbigu}}
\newcommand{\nbigstilde}{\widetilde{\nbigs}}

\newcommand{\tauhat}{\widehat{\tau}}
\newcommand{\That}{\widehat{T}}
\newcommand{\Phat}{\widehat{P}}
\newcommand{\nbigmhat}{\widehat{\nbigm}}
\newcommand{\nbigzhat}{\widehat{\nbigz}}
\newcommand{\Nhat}{\widehat{N}}
\newcommand{\gbigehat}{\widehat{\gbige}}

\newcommand{\gtilde}{\widetilde{g}}
\newcommand{\Ztilde}{\widetilde{Z}}
\newcommand{\nbigttilde}{\widetilde{\nbigt}}
\newcommand{\gbigstilde}{\widetilde{\gbigs}}
\newcommand{\varphitilde}{\widetilde{\varphi}}
\newcommand{\Ctilde}{\widetilde{C}}

\newcommand{\gbigetilde}{\widetilde{\gbige}}
\newcommand{\gbigctilde}{\widetilde{\gbigc}}

\newcommand{\Tzero}{T^{(\lambda_0)}}
\newcommand{\Ttildezero}{\widetilde{T}^{(\lambda_0)}}

\newcommand{\tautilde}{\widetilde{\tau}}

\newcommand{\ybar}{\overline{y}}
\newcommand{\epsilonbar}{\overline{\epsilon}}
\newcommand{\Atilde}{\widetilde{A}}
\newcommand{\Otilde}{\widetilde{O}}

\newcommand{\Deltatilde}{\widetilde{\Delta}}
\newcommand{\nbigftilde}{\widetilde{\nbigf}}
\newcommand{\nbigtzero}{\nbigt^{(\lambda_0)}}
\newcommand{\nbigpzero}{\nbigp^{(\lambda_0)}}
\newcommand{\Mtilde}{\widetilde{M}}
\newcommand{\Ytilde}{\widetilde{Y}}

\newcommand{\wtilde}{\widetilde{w}}
\newcommand{\ztilde}{\widetilde{z}}
\newcommand{\nbiglbar}{\overline{\nbigl}}
\newcommand{\nbigkbar}{\overline{\nbigk}}

\newcommand{\nbigxtilde}{\widetilde{\nbigx}}
\newcommand{\gbigttilde}{\widetilde{\gbigt}}
\newcommand{\nbigvlambda}{\nbigv^{\lambda}}

\newcommand{\nbigdtilde}{\widetilde{\nbigd}}
\newcommand{\Qtilde}{\widetilde{Q}}
\newcommand{\Ttilde}{\widetilde{T}}

\newcommand{\what}{\widehat{w}}
\newcommand{\uhat}{\widehat{u}}

\newcommand{\betabar}{\overline{\beta}}
\newcommand{\vecatilde}{\widetilde{\veca}}

\newcommand{\forget}{for}
\newcommand{\nbigitilde}{\widetilde{\nbigi}}
\newcommand{\Secbar}{\overline{\Sec}}
\newcommand{\Multisector}{\nbigm\nbigs}
\newcommand{\SD}{\nbigs\nbigd}

\newcommand{\Irrbar}{\overline{\Irr}}

\newcommand{\nbigktilde}{\widetilde{\nbigk}}

\newcommand{\vecctilde}{\widetilde{\vecc}}
\newcommand{\vecnbigi}{{\boldsymbol \nbigi}}
\newcommand{\vecnbigf}{{\boldsymbol \nbigf}}
\newcommand{\vecnbigftilde}{\boldsymbol \nbigftilde}
\newcommand{\vecnbigg}{{\boldsymbol \nbigg}}

\newcommand{\DDprimelambda}{\DD^{\prime\lambda}}

\newcommand{\nbigmlambda}{\nbigm^{\lambda}}
\newcommand{\nrhom}{R{\mathcal Hom}}
\newcommand{\DDD}{\boldsymbol D}
\newcommand{\jtilde}{\widetilde{j}}
\newcommand{\vecH}{{\boldsymbol H}}
\newcommand{\nbiggzero}{\nbigg^{(\lambda_0)}}

\newcommand{\ctilde}{\widetilde{c}}
\newcommand{\taubar}{\overline{\tau}}
\newcommand{\Kbar}{\overline{K}}
\newcommand{\Ktilde}{\widetilde{K}}
\newcommand{\veckappa}{{\boldsymbol \kappa}}
\newcommand{\nbigdzero}{\nbigd^{(\lambda_0)}}
\newcommand{\nbigxzero}{\nbigx^{(\lambda_0)}}
\newcommand{\nbigyzero}{\nbigy^{(\lambda_0)}}
\newcommand{\nbigdhatzero}{\widehat{\nbigd}^{(\lambda_0)}}
\newcommand{\nbigditi}{\nbigd^{(\lambda_1)}}
\newcommand{\nbigxiti}{\nbigx^{(\lambda_1)}}
\newcommand{\nbigdhatiti}{\widehat{\nbigd}^{(\lambda_1)}}
\newcommand{\nbigdtildezero}{\nbigdtilde^{(\lambda_0)}}
\newcommand{\nbigxtildezero}{\nbigxtilde^{(\lambda_0)}}

\newcommand{\Psihat}{\widehat{\Psi}}
\newcommand{\btilde}{\widetilde{b}}
\newcommand{\vecbtilde}{\widetilde{\vecb}}
\newcommand{\vecetilde}{\widetilde{\vece}}
\newcommand{\etilde}{\widetilde{\vece}}
\newcommand{\Omegatilde}{\widetilde{\Omega}}
\newcommand{\TTtilde}{T\widetilde{T}}

\newcommand{\gminiatilde}{\widetilde{\gminia}}
\newcommand{\gminibtilde}{\widetilde{\gminib}}
\newcommand{\nbigfshitadhat}{\nbigf_{|\widehat{D}}}
\newcommand{\nbigfshitadhatN}{\nbigf_{|\widehat{D}^{(N)}}}

\newcommand{\jbikkuri}{j_!}
\newcommand{\FzeroVichiVni}{F_0Hom(V_1,V_2)}
\newcommand{\FminusVichiVni}{F_{<0}Hom(V_1,V_2)}
\newcommand{\nbigfminusEndESbar}{\nbigf_{<0}\End(E)_{|\Sbar}}
\newcommand{\nbigfzeroEndESbar}{\nbigf_{0}\End(E)_{|\Sbar}}

\newcommand{\nbigibar}{\overline{\nbigi}}

\newcommand{\Qhat}{\widehat{Q}}

\newcommand{\SDtilde}{\widetilde{\SD}}
\newcommand{\vecVtilde}{\widetilde{\vecV}}
\newtheorem{thm}{Theorem}[section]
\newtheorem{cor}[thm]{Corollary}

\newtheorem{rem}[thm]{Remark}
\newtheorem{lem}[thm]{Lemma}
\newtheorem{prop}[thm]{Proposition}
\newtheorem{df}[thm]{Definition}

\newtheorem{condition}[thm]{Condition}
\newtheorem{assumption}[thm]{Assumption}

\newtheorem{notation}[thm]{Notation}

\author{Takuro Mochizuki}
\address{
Research Institute for Mathematical Sciences,
Kyoto University, Kyoto 606-8502, Japan}
\email{takuro@kurims.kyoto-u.ac.jp}

\title[Wild harmonic bundles]
{Wild Harmonic Bundles and \\
Wild Pure Twistor $D$-modules}

\makeindex
\begin{document}
\frontmatter

\begin{abstract}

We study
(i) asymptotic behaviour of wild harmonic bundles,
(ii) the relation between semisimple
 meromorphic flat connections
 and wild harmonic bundles,
(iii)
 the relation between wild harmonic bundles
 and polarized wild pure twistor $D$-modules.
As an application, 
we show the hard Lefschetz theorem for
algebraic semisimple holonomic $D$-modules,
conjectured by M. Kashiwara.

We also study resolution of turning points
for algebraic meromorphic flat bundles,
and we show the existence of 
the good Deligne-Malgrange lattice after blow up,
which seems of foundational importance
for the study algebraic $D$-modules.

\end{abstract}



\maketitle
\tableofcontents

\mainmatter

\chapter{Introduction}
In our previous work
(\cite{mochi}, 
\cite{mochi4},
\cite{mochi2} and \cite{mochi5}),
we studied tame harmonic bundles
after C. Simpson
(\cite{s1}, \cite{s2}, \cite{s5}, \cite{s3}, for example),
and we established their foundational property.
(See also the works of O. Biquard \cite{b}
and J. Jost and K. Zuo \cite{JZ2}).
In cooperation with the work of C. Sabbah \cite{sabbah2}
on pure twistor $D$-modules (based on the theory
of pure Hodge modules
due to M. Saito \cite{saito1}, \cite{saito2}),
we obtained some deep results on 
algebraic regular holonomic $D$-modules.

In this monograph, we systematically study
{\em wild harmonic bundles} over
complex manifolds of arbitrary dimension,
and we obtain generalizations of our previous results 
for tame harmonic bundles 
in the case of wild harmonic bundles.
We also give applications to
algebraic meromorphic flat bundles
and algebraic holonomic $D$-modules.
Recently,
there has been a growing interest in 
wild harmonic bundles and 
(not necessarily regular) 
holonomic $D$-modules on curves.
(See \cite{biquard-boalch},
 \cite{boalch1},
 \cite{boalch2},
 \cite{Bridgeland-Laredo},
 \cite{hertling-sevenheck},
 \cite{Katzarkov-Kontsevich-Pantev},
 \cite{sabbah6}, 
 \cite{simpson_middle_convolution},
 \cite{szabo},
 \cite{witten}, for example.)
However,
in this monograph,
we will NOT consider
moduli spaces,
mirror symmetries,
geometric Langlands theory,
integrable systems,
non-abelian Hodge theory,
Painlev\'e equations,
and any other fashionable subjects
related with harmonic bundles and Higgs bundles.
It can be said that our goal is more modest and basic.
Nonetheless, 
the author has been deeply impressed with
what he saw in this study,
for example
an interaction between the theories of harmonic bundles
and $D$-modules.

\section{Contents of this monograph}
\label{subsection;08.10.24.1}
Briefly speaking,
this study consists of three main bodies
and preliminaries for them:
\begin{description}
\item[(A)]
 Asymptotic behaviour of wild harmonic bundles.
\item[(B)]
 Application to algebraic meromorphic flat bundles.
\item[(C)]
 Application to wild pure twistor $D$-modules
 and algebraic $D$-modules.
\end{description}

Let us briefly describe each part.
We will give some more detailed introductions later
(Sections
 \ref{subsection;08.10.24.25}--\ref{subsection;08.10.26.1}).

We have two main issues in Part (A).
Let $X$ be a complex manifold,
and let $D$ be a normal crossing hypersurface of $X$.
Let $\harmonicbundle$ be a good wild harmonic bundle
on $(X,D)$.
Although we will not explain the definition here
(see Section \ref{subsection;08.10.24.20},
 or Section \ref{subsection;07.12.4.1}
 for more precision),
it means that $\harmonicbundle$ is a harmonic bundle
on $X-D$ satisfying some conditions
around each point of $D$.
We would like to prolong it to something on $X$,
that is the first main issue,
and fundamental for us.
The role of such prolongation may be
compared with nilpotent orbit theorem 
in Hodge theory 
due to W. Schmid \cite{sch}.

Then, we would like to understand more detailed property.
It is achieved by showing that 
we obtain tame harmonic bundles 
from wild harmonic bundles
as Gr with respect to Stokes filtrations,
which is the second main issue in (A).
By this reduction, the study of the asymptotic behaviour
of wild harmonic bundles is reduced to the tame case 
investigated in \cite{mochi2}.

\vspace{.1in}

There are two main purposes in Part (B).
One is to characterize 
semisimplicity of algebraic meromorphic flat bundles
by the existence of pluri-harmonic metric
with some nice property.
The other is to show the existence of
resolutions of turning points for
algebraic meromorphic flat bundles.

Due to K. Corlette \cite{corlette},
a flat bundle on a smooth projective variety 
has a pluri-harmonic metric
if and only if it is semisimple, 
i.e., a direct sum of irreducible ones.
It was generalized to the case of meromorphic flat bundles
with regular singularity (\cite{JZ2} and \cite{mochi2}).
In this monograph,
we will establish such a characterization in the irregular case.
We also study the Kobayashi-Hitchin correspondence
for meromorphic flat bundles.

We have an interesting application of 
such an existence result of pluri-harmonic metric
in the resolution of turning points
for algebraic meromorphic flat bundles,
which is the other main result in Part (B).
It seems of foundational importance
in the study of algebraic holonomic $D$-modules,
and might be compared with a resolution of
singularities for algebraic varieties.

\begin{rem}
Recently, K. Kedlaya established it
in a more general situation
with a completely different way.
See {\rm\cite{kedlaya}} and {\rm\cite{kedlaya2}}.
\hfill\qed
\end{rem}

In Part (C),
we will establish the relation between 
wild harmonic bundles and 
polarized wild pure twistor $D$-modules.
Recently, Sabbah introduced the notion of
wild pure twistor $D$-modules \cite{sabbah5}.
Our result roughly says that
polarized wild pure twistor $D$-modules
are actually minimal extensions of
wild harmonic bundles.
Together with the result in Part (B),
we obtain the correspondence between
semisimple holonomic $D$-modules
and polarizable wild pure twistor $D$-modules 
on complex projective varieties.
As an application, 
we will show the Hard Lefschetz theorem for
algebraic semisimple  (not necessarily regular)
holonomic $D$-modules,
conjectured by M. Kashiwara.

We may also say that
the result in Part (C) makes us possible
to define the push-forward for wild harmonic bundles,
which will be useful to produce
new wild harmonic bundles,
or to enrich some operations for flat bundles 
by polarized twistor structures.

\vspace{.1in}

We need various preliminaries.
Let us mention some of major ones.
We will revisit asymptotic analysis
for meromorphic flat bundles
in our convenient way,
which was originally
studied by H. Majima \cite{majima}
and refined by Sabbah \cite{sabbah4}.
We will put a stress on 
canonically defined Stokes filtrations.
It is fundamental for us to consider
Gr with respect to Stokes filtrations, 
and deformations caused by
variation of irregular values.
(See Chapters
\ref{section;07.11.8.15}--\ref{chapter;10.5.26.11}.)

As another important preliminary,
we show that acceptable bundles are 
naturally prolonged to filtered bundles.
(See Section \ref{subsection;07.9.22.10}.)
After the work of M. Cornalba-P. Griffiths and Simpson 
(\cite{cg}, \cite{s1} and \cite{s2})
the author studied acceptable bundles,
and he obtained such prolongation
for acceptable bundles which
come from tame harmonic bundles
(\cite{mochi}, \cite{mochi4} and \cite{mochi2}).
To apply the theory in the wild case,
we will establish it for general acceptable bundles.
Although only small changes are required,
we will give rather detailed arguments
in view of its importance.
(See Chapter \ref{section;08.10.23.2}.)

\section{Asymptotic behaviour of wild harmonic bundles}
\label{subsection;08.10.24.25}
\subsection{Prolongation}
\label{subsection;08.10.24.20}

\subsubsection{Harmonic bundle}

Recall the definition of harmonic bundle
\cite{s5}.
\index{harmonic bundle}
Let $(E,\delbar_E,\theta)$ be a Higgs bundle on 
a complex manifold.
Let $h$ be a hermitian metric of $E$.
Then, we have the associated unitary connection
$\delbar_E+\del_E$
and the adjoint $\theta^{\dagger}$ of $\theta$
with respect to $h$.
The metric $h$ is called {\em pluri-harmonic},
if the connection
$\DD^1:=\delbar_E+\del_E+\theta+\theta^{\dagger}$
is flat.
In that case,
$(E,\delbar_E,\theta,h)$ is called a {\em harmonic bundle}.
We also have another equivalent definition.
Let $(V,\nabla)$ be a flat bundle on a complex manifold.
Let $h$ be a hermitian metric of $V$.
Then, we have the decomposition
$\nabla=\nabla^u+\Phi$,
where $\nabla^u$ is a unitary connection
and $\Phi$ is self-adjoint with respect to $h$.
We have the decompositions
$\nabla^u=\del_V+\delbar_V$
and $\Phi=\theta+\theta^{\dagger}$
into the $(1,0)$-part and the $(0,1)$-part.
We say that $h$ is a {\em pluri-harmonic metric}
of $(V,\nabla)$,
if $(V,\delbar_V,\theta)$ is a Higgs bundle.
In that case,
$(V,\nabla,h)$ is called a {\em harmonic bundle}.

Let $\harmonicbundle$ be a harmonic bundle
on a complex manifold $Y$.
For any complex number $\lambda$,
we have the flat $\lambda$-connection
$(\nbigelambda,\DDlambda)$,
i.e.,
 the holomorphic vector bundle
$\nbigelambda:=
 (E,\delbar_E+\lambda\theta^{\dagger})$
with the flat $\lambda$-connection
$\DDlambda:=\delbar_E+\lambda\theta^{\dagger}
 +\lambda\del_E+\theta$.
\index{$\lambda$-flat bundle}
We also have 
the family of $\lambda$-flat bundles $(\nbige,\DD)$
on $\cnum_{\lambda}\times Y$,
i.e.,
the holomorphic vector bundle
$\nbige:=\bigl(p_{\lambda}^{-1}(E),
 \delbar_E+\lambda\theta^{\dagger}
 +\delbar_{\lambda}\bigr)$
equipped with a family of flat $\lambda$-connections
$\DD:=\delbar_E+\lambda\theta^{\dagger}
+\lambda\del_E+\theta$,
where 
$p_{\lambda}:\cnum_{\lambda}\times Y\lrarr Y$
denotes the projection.
\index{family of $\lambda$-flat bundles}

\subsubsection{Good wild harmonic bundle}

Let $X$ be a complex manifold
with a simple normal crossing hypersurface $D$.
Let $M(X,D)$ be the set of meromorphic functions
whose poles are contained in $D$,
and $H(X)$ be the set of holomorphic functions on $X$.
\index{set $M(X,D)$}
\index{set $H(X)$}
As mentioned in Section \ref{subsection;08.10.24.1},
it is fundamental for us to prolong 
a {\em good wild harmonic bundle} on $X-D$
to something on $X$.
We would like to explain some more details.
To begin with, we explain
what is good wild harmonic bundle.

Let us consider the local case,
i.e., $X:=\Delta^n=\bigl\{
 (z_1,\ldots,z_n)\,\big|\,|z_i|<1
 \bigr\}$ and
$D:=\bigcup_{j=1}^{\ell}\{z_j=0\}$.
Let $\harmonicbundle$ be a harmonic bundle
on $X-D$.
It is called
{\em strongly unramifiedly good wild harmonic bundle}
on $(X,D)$,
if there exist a good set of irregular values
$\Irr(\theta)\subset M(X,D)/H(X)$,
a finite subset $\Sp(\theta)\subset\cnum^{\ell}$
and a decomposition
\begin{equation}
 \label{eq;08.10.24.10}
 (E,\theta)
=\bigoplus_{\gminia\in\Irr(\theta)}
  \bigoplus_{\vecalpha\in \Sp(\theta)}
 (E_{\gminia,\vecalpha},
 \theta_{\gminia,\vecalpha}),
\end{equation}
such that
the eigenvalues of $\theta_{\gminia,\vecalpha}$
are $d\gminia+\sum_{i=1}^{\ell} \alpha_i\cdot dz_i/z_i$
modulo multi-valued holomorphic one forms on $X$.
(See Definition \ref{df;07.6.1.15}
 for {\em good set of irregular values}.)
In another brief word,
for the expression
\[
 \theta_{\gminia,\vecalpha}
=d\gminia+\sum_{j=1}^{\ell} 
 (\alpha_j+f_{\gminia,\vecalpha,j}) 
 \cdot \frac{dz_j}{z_j}
+\sum_{j=\ell+1}^n 
 g_{\gminia,\vecalpha,j}\cdot dz_j,
\]
the characteristic polynomials
$\det(T-f_{\gminia,\vecalpha,j})$ $(j=1,\ldots,\ell)$
and $\det(T-g_{\gminia,\vecalpha,j})$
$(j=\ell+1,\ldots,n)$ are contained in
$H(X)[T]$,
and $\det(T-f_{\gminia,\vecalpha,j})_{|\{z_j=0\}}
=T^{\rank E_{\gminia,\vecalpha}}$.
\index{strongly unramifiedly good wild harmonic bundle}
We say that $(E,\delbar_E,\theta,h)$
is a {\em strongly good wild harmonic bundle} on $(X,D)$,
if there exists a ramified covering
$\varphi:(X',D')\lrarr (X,D)$
such that
$\varphi^{\ast}(E,\delbar_E,\theta,h)$
is a strongly unramifiedly good wild harmonic bundle,
i.e., $(E,\delbar_E,\theta,h)$ is the descent
of an unramifiedly good wild harmonic bundle.
\index{strongly good wild harmonic bundle}

The definitions can be easily globalized.
Let $X$ be a complex manifold,
and let $D$ be a normal crossing hypersurface.
A harmonic bundle $(E,\delbar_E,\theta,h)$
on $X-D$ is called a {\em good wild harmonic bundle}
on $(X,D)$,
if the following holds:
\begin{itemize}
\item For any point $P\in D$,
 there exists a coordinate neighbourhood $X_P$
 around $P$
 such that $\harmonicbundle_{|X_P\setminus D}$
 is a strongly good wild harmonic bundle on
 $(X_P,D\cap X_P)$.
\end{itemize}

We have one more additional notion of
{\em wild harmonic bundle}.
Let $Y$ be an irreducible complex 
(not necessarily smooth) variety.
Let $\harmonicbundle$ be a harmonic bundle
defined on the complement of 
a closed analytic subset $Z$ of $Y$.
We say that $\harmonicbundle$ is
wild on $(Y,Z)$ or $(Y,Y-Z)$,
if there exists a morphism $\varphi$
of a smooth complex variety $Y'$ to $Y$
such that 
(i) $\varphi$ is birational and projective,
(ii) $\varphi^{-1}(Z)$ is a normal crossing
 hypersurface of $Y'$,
(iii) $\varphi^{-1}\harmonicbundle$
 is a good wild harmonic bundle on
 $\bigl(Y',\varphi^{-1}(Z)\bigr)$.
\index{wild harmonic bundle}

Essential analysis is done for
unramifiedly good wild harmonic bundle.
The study of good wild harmonic bundle
can be easily reduced to the unramified case.
The notion of wild harmonic bundle
is rather auxiliary to give statements 
in the application to wild pure twistor $D$-modules.

\begin{rem}
In the one dimensional case,
any wild harmonic bundle is good,
and they can be defined in a much simpler way.
Namely, we have only to impose the condition
that the characteristic polynomial of 
the associated Higgs field
is meromorphic.
\hfill\qed
\end{rem}

\subsubsection{Sheaf of holomorphic sections 
with polynomial growth}

We explain various prolongation of
a good wild harmonic bundle 
$\harmonicbundle$ on $(X,D)$
to some meromorphic objects on $X$.
We consider the local case
$X:=\Delta^n$ and $D:=\bigcup_{i=1}^{\ell}\{z_i=0\}$.
The constructions and the results
can be easily globalized.

For any $\lambda$,
we would like to prolong
the $\lambda$-flat bundle
$(\nbigelambda,\DDlambda)$
on $X-D$ to a meromorphic $\lambda$-flat bundle
on $(X,D)$.
For any  $\veca\in\real^{\ell}$
and any open subset $U\subset X$,
we define
\begin{equation}
\label{eq;10.5.31.1}
 \nbigp_{\veca}\nbigelambda(U):=
 \Bigl\{
 f\in \nbigelambda(U\setminus D)\,\Big|\,
 |f|_h=O\Bigl(
 \prod_{j=1}^{\ell}|z_j|^{-a_j-\epsilon}
 \Bigr),\,\,\forall\epsilon>0
 \Bigr\}.
\end{equation}
By taking sheafification,
we obtain an increasing sequence of
$\nbigo_{X}$-modules
$\nbigp_{\ast}\nbigelambda:=
 \bigl(\nbigp_{\veca}\nbigelambda\,\big|\,
 \veca\in\real^{\ell}\bigr)$
and an $\nbigo_X(\ast D)$-module
$\nbigp\nbigelambda:=
 \bigcup_{\veca\in\real^{\ell}}
 \nbigp_{\veca}\nbigelambda$.
The following theorem is the starting point of 
the study on asymptotic behaviour
of good wild harmonic bundle.
\index{sheaf $\nbigp\nbigelambda$}
\index{sheaf $\nbigp_{\veca}\nbigelambda$}
\index{filtered sheaf $\nbigp_{\ast}\nbigelambda$}

\begin{thm}
[Theorem \ref{thm;07.12.2.55},
Theorem \ref{thm;07.11.18.20}]
\label{thm;07.12.4.20}
$(\nbigp_{\ast}\nbigelambda,\DDlambda)$
is a good filtered $\lambda$-flat bundle.
If $\harmonicbundle$ is unramified
with the decomposition in 
{\rm(\ref{eq;08.10.24.10})},
the set of irregular values is given by
$\Irr\bigl(\nbigp\nbigelambda,\DDlambda\bigr)
=\bigl\{
 (1+|\lambda|^2)\,\gminia\,\big|\,
 \gminia\in\Irr(\theta)\bigr\}$.
\hfill\qed
\end{thm}
We refer to Section \ref{section;10.5.20.1}
for details on 
the notion of {\em good filtered $\lambda$-flat bundle},
but we explain a brief meaning of the theorem
in the case that
$\harmonicbundle$ is {\em unramified}
with a decomposition (\ref{eq;08.10.24.10}).
\index{good filtered $\lambda$-flat bundle}
Let $O$ denote the origin $(0,\ldots,0)$.
Each $\nbigp_{\veca}\nbigelambda$
is a locally free $\nbigo_X$-module,
and we have the decomposition
for the completion at $O$
\begin{equation}
 \label{eq;08.10.24.6}
 (\nbigp_{\veca}\nbigelambda,\DDlambda)_{|\Ohat}
=\bigoplus_{\gminia\in \Irr(\theta)}
 \bigl(\nbigp_{\veca}\nbigelambdahat_{\gminia},
 \DDlambdahat_{\gminia}\bigr),
\end{equation}
where
$\DDlambdahat_{\gminia}
 -(1+|\lambda|^2)d\gminia$
have logarithmic singularities.
Moreover,
we have such decompositions
on the completions at any point of $D$.

\vspace{.1in}
Let us give a remark on the proof.
We use the essentially same arguments
as those in our previous paper \cite{mochi2}
to show that $\nbigp_{\veca}\nbigelambda$ are locally free.
Namely, we give an estimate on the Higgs field $\theta$
(Theorem \ref{thm;07.10.4.1} and
 Theorem \ref{thm;07.10.4.3})
as the wild version of {\em Simpson's Main estimate},
which briefly means that
the decomposition (\ref{eq;08.10.24.10})
is asymptotically orthogonal,
and that 
$\theta_{\gminia,\vecalpha}
-d\gminia-\sum \alpha_i\, dz_i/z_i$
are bounded with respect to
$h$ and the Poincar\'e metric of $X-D$.
Then, we obtain that $(\nbigelambda,h)$ is 
{\em acceptable},
i.e., the curvature of $(\nbigelambda,h)$
is bounded with respect to $h$
and the Poincar\'e metric
(Corollary \ref{cor;07.11.22.30}).
\index{acceptable bundle}
We obtain that $\nbigp_{\ast}\nbigelambda$
is a filtered bundle
by using a general theory of acceptable bundles
(Theorem \ref{thm;07.3.16.1}).
In the wild case,
we need some more argument to show that
$\DDlambda$ is meromorphic
and to obtain the decomposition (\ref{eq;08.10.24.6}).

\subsubsection{Family of 
meromorphic $\lambda$-flat bundles}

Next, we would like to consider 
the prolongation of the family of
$\lambda$-flat bundles $(\nbige,\DD)$ on 
$\cnum_{\lambda}\times(X-D)$
to a family of meromorphic $\lambda$-flat bundles
on $\cnum_{\lambda}\times (X,D)$.
A naive idea is to consider the family
$\bigcup_{\lambda\in\cnum_{\lambda}}
 \nbigp\nbigelambda$,
i.e., the sheaf of holomorphic sections
of $\nbige$ with polynomial growth.
In the tame case, it gives a nice {\em meromorphic} object.
(But, note that we need some more consideration
for {\em lattices}.)
However, 
$\bigcup_{\lambda\in\cnum} \nbigp\nbigelambda$
cannot be a good meromorphic prolongment
in the wild case,
because the irregular values $(1+|\lambda|^2)\gminia$
of $\DDlambda$
have non-holomorphic dependence on $\lambda$,
as mentioned in Theorem \ref{thm;07.12.4.20}.
Hence, we need the deformations of
$\nbigp\nbigelambda$ $(\lambda\neq 0)$
caused by variation of the irregular values
(Section \ref{subsection;10.5.17.51}).

\paragraph{}
We briefly and imprecisely explain such a deformation
in the case $\lambda=1$.
First, we can obtain such deformation
for a unramifiedly good meromorphic flat bundle,
by considering the deformation of Stokes structure.
It is easily extended to the (possibly) ramified case,
and then the global case.
Hence, we explain the essential part,
i.e., the local an unramified case.
Let $X:=\Delta^n$ and 
$D:=\bigcup_{i=1}^{\ell}\{z_i=0\}$.
Let $(\nbigv,\nabla)$ be 
a meromorphic flat bundle on $(X,D)$,
with a lattice $V$ which has decompositions
as in (\ref{eq;08.10.24.6})
on the completions at any point of $D$,
for example
\begin{equation}
 \label{eq;08.10.24.11}
 (V,\nabla)_{|\Ohat}
=\bigoplus_{\gminia\in\Irr(\nabla)}
 (\Vhat_{\gminia},\nablahat_{\gminia}),
\quad\quad
 \bigl(
\nablahat_{\gminia}-d\gminia
\bigr)\Vhat_{\gminia}
\subset
 \Vhat_{\gminia}\otimes
 \Omega^1_X(\log D).
\end{equation}
Such a lattice is called 
an {\em unramifiedly good lattice} 
of $(\nbigv,\nabla)$
(Section \ref{subsection;08.1.23.1}).
\index{unramifiedly good lattice}
We should remark that such a lattice 
may not exist in general.
For simplicity,
we also assume a non-resonance condition
to the lattice,
i.e., the difference of the distinct eigenvalues
of the residues are not integers.
Let $S$ be a small multi-sector in $X-D$
whose closure contains $O$.
Let $\Sbar$ denote its closure
in the real blow up $\Xtilde(D)$ of $X$ along $D$,
or more precisely, the fiber product of
the real blow up along the irreducible components 
of $D$, taken over $X$.
According to the asymptotic analysis
for meromorphic flat bundles
(see \cite{majima}, \cite{sabbah4},
 and Chapter \ref{section;07.12.26.6}
 in this monograph),
we can lift (\ref{eq;08.10.24.11})
to a flat decomposition of $V_{|\Sbar}$:
\begin{equation}
 \label{eq;08.10.24.12}
 V_{|\Sbar}
=\bigoplus_{\gminia\in\Irr(\nabla)}
 V_{\gminia,S}
\end{equation}
Namely, the completion of
(\ref{eq;08.10.24.12})
along $\Sbar\cap\pi^{-1}(O)$
is the same as the pull back of
(\ref{eq;08.10.24.11}).
We consider the order $\leq_S$
on $\Irr(\nabla)$
given by $\gminia\leq_S\gminib$
$\Longleftrightarrow$
$-\Re(\gminia)(Q)\leq -\Re(\gminib)(Q)$ on
any points $Q$ of $S$
which are sufficiently close to $O$.
Then, it can be shown that
\[
 \nbigftilde^S_{\gminia}(V_{|\Sbar}):=
 \bigoplus_{\gminib\leq_S\gminia}
 V_{\gminib,S}
\]
is independent of the choice of a lifting
(\ref{eq;08.10.24.12}).
Thus, we obtain the filtration $\nbigftilde^S$
of $V_{|\Sbar}$ or $V_{|S}$,
which is a $\nabla$-flat filtration
indexed by
$\bigl(\Irr(\nabla),\leq_S\bigr)$.
It is called the full Stokes filtration.
\index{full Stokes filtration}
For two small multi-sectors $S_i$ $(i=1,2)$,
the filtrations $\nbigftilde^{S_i}$
satisfy some compatibility on $S_1\cap S_2$.
(See Section \ref{subsection;08.9.28.200}.)
Such a system of filtrations 
$\bigl\{\nbigftilde^S(V_{|S})\,|\,S\subset X-D\bigr\}$
is called the Stokes structure of $\nbigv$.
\index{Stokes structure}
It can be shown that
the meromorphic flat bundle
$(\nbigv,\nabla)$ is recovered
from the flat bundle $\nbigv_{|X-D}$
with the Stokes structure.

For any $T>0$,
we set 
$\Irr\bigl(\nabla^{(T)}\bigr):=
 \bigl\{
 T\gminia\,\big|\,
 \gminia\in\Irr(\nabla) \bigr\}$.
The natural bijection
$\Irr(\nabla)\simeq
 \Irr(\nabla^{(T)})$
preserves the order $\leq_S$.
We naturally obtain 
a system of the filtrations
$\{\nbigftilde^S_{T\gminia}\,\big|\,
 \gminia\in\Irr(\nabla)\}$
for each small multi-sector
indexed by
$\bigl(
 \Irr(\nabla^{(T)}),\leq_S
 \bigr)$,
which gives a new Stokes structure.
We obtain the associated meromorphic flat bundle
$(\nbigv^{(T)},\nabla^{(T)})$ on $(X,D)$,
which is a deformation caused
by variation of irregular values.
\index{deformation caused by
variation of irregular values}
As mentioned,
it is not difficult to extend it 
to the global and (possibly) ramified case.

\begin{rem}
Such a deformation grew out of the discussion
with C. Sabbah.
We study it in a more general situation.

After finishing the first version
of this monograph,
the author found that
such a deformation in the curve case
appeared
in several works such as
{\rm \cite{boalch3}},
{\rm \cite{Concini-Kac-Procesi}}
and {\rm\cite{laredo}}.
The referee kindly informed that
the irregular values was used
as parameters for universal deformations
in {\rm\cite{malgrange-deformation}}.
\hfill\qed
\end{rem}

\paragraph{}
Applying this deformation procedure to
$(\nbigp\nbigelambda,\DDlambda)$
with $T=(1+|\lambda|^2)^{-1}$  for $\lambda\neq 0$,
we obtain the meromorphic $\lambda$-flat bundle
$(\nbigq\nbigelambda,\DDlambda)$.
\index{sheaf $\nbigq\nbigelambda$}
We put $\nbigq\nbige^0:=\nbigp\nbige^0$.
The next theorem says that
the family $\bigcup \nbigq\nbigelambda$
gives a nice meromorphic object
on $\cnum_{\lambda}\times (X,D)$.
For simplicity of the description,
let $\nbigxzero$ denote a neighbourhood
of $\{\lambda_0\}\times X$
in $\cnum_{\lambda}\times X$,
and $\nbigdzero:=
 \nbigxzero\cap\bigl(\cnum_{\lambda}\times D\bigr)$.
\begin{thm}
[Theorem \ref{thm;07.12.4.15}]
\label{thm;07.12.4.25}
We have a unique family of 
meromorphic $\lambda$-flat bundles
$(\nbigq\nbige,\DD)$ on
$\cnum_{\lambda}\times(X,D)$
such that the specialization of
$(\nbigq\nbige,\DD)$ to $\{\lambda\}\times X$
is isomorphic to $(\nbigq\nbigelambda,\DDlambda)$.
Moreover,
we have a family of
good filtered $\lambda$-flat bundles
$(\nbigqzero_{\ast}\nbige,\DD)$
on $(\nbigxzero,\nbigdzero)$ for each 
$\lambda_0\in\cnum_{\lambda}$
with the KMS-structure
such that
$\nbigq\nbige_{|\nbigxzero}
=\bigcup \nbigqzero_{\veca}\nbige$.
\index{sheaf $\nbigq\nbige$}
\hfill\qed
\end{thm}

We refer to Section {\rm\ref{subsection;07.10.14.15}}
for {\em KMS-structure}.
Although it is quite important in the study of
the asymptotic behaviour of wild harmonic bundles,
it has already appeared in the tame case \cite{mochi2},
where we have studied it in detail.
Hence, we omit to explain it in this introduction.
\index{KMS-structure}

\vspace{.1in}

Note that 
we have the uniqueness of the family
because 
$(\nbigq\nbige,\DD)_{\{\lambda\}\times X}
\simeq
 (\nbigq\nbigelambda,\DDlambda)$.
Hence, to show Theorem \ref{thm;07.12.4.25},
we have only to consider it
in the local and unramified case.
Then, we shall argue in two steps.
We first construct a family of meromorphic $\lambda$-flat
bundles $\nbigpzero\nbige$ on $(\nbigxzero,\nbigdzero)$
by using a deformed metric $\nbigpzero h$
given in (\ref{eq;07.7.21.4}).
\index{sheaf $\nbigpzero\nbige$}
Namely, for the metric $\nbigpzero h$,
we consider an increasing sequence of
$\nbigo_{\nbigxzero}$-modules
$\nbigpzero_{\ast}\nbige
:=\bigl(
 \nbigpzero_{\veca}\nbige\,\big|\,
 \veca\in\real^{\ell}
 \bigr)$ given as in (\ref{eq;10.5.31.1}),
and the $\nbigo_{\nbigxzero}(\ast\nbigdzero)$-module
$\nbigpzero\nbige:=
 \bigcup_{\veca\in\real^{\ell}}
 \nbigpzero_{\veca}\nbige$.
The restriction of $\nbigpzero\nbige$
to $\{\lambda\}\times X$ is denoted by
$\nbigpzero\nbigelambda$.
\index{sheaf $\nbigpzero\nbigelambda$}

\begin{thm}
[Theorem \ref{thm;07.11.20.10},
 Proposition \ref{prop;07.10.21.100}]
$(\nbigpzero_{\ast}\nbige,\DD)$
is a family of good filtered $\lambda$-flat bundles
with the KMS-structure.
Moreover,
$(\nbigpzero\nbigelambda,\DDlambda)$
is naturally isomorphic to
$\bigl(
 (\nbigp\nbigelambda)^{T_1(\lambda,\lambda_0)},
 \DDlambda
\bigr)$,
where
$T_1(\lambda,\lambda_0):=
 (1+\lambda\lambdabar_0)(1+|\lambda|^2)^{-1}$.
\hfill\qed
\end{thm}

In the second step,
we obtain a family of meromorphic 
$\lambda$-flat bundles $(\nbigqzero\nbige,\DD)$
as the deformation
$(\nbigpzero\nbige,\DD)^{(T_2(\lambda,\lambda_0))}$,
where $T_2(\lambda,\lambda_0):=
(1+\lambda\lambdabar_0)^{-1}$.
By varying $\lambda_0\in\cnum$
and gluing $\nbigqzero\nbige$,
we obtain the desired family $(\nbigq\nbige,\DD)$.

\begin{rem}
We should emphasize that
$\nbigq\nbige$ is given on $\cnum_{\lambda}\times X$.
Contrastively, the reduction in the next subsection
is given locally around any point of $D$.
\hfill\qed
\end{rem}

\subsection{Reduction from wild harmonic bundles
 to tame harmonic bundles}

We would like to analyze more closely
the behaviour of good wild harmonic bundle 
{\em around a given point of $D$}.
For that purpose, we consider Gr
with respect to Stokes structure.
The construction in this subsection is 
given only locally,
although the construction in the previous 
subsection can be easily globalized.
We set
$X:=\Delta^n$ and 
$D:=\bigcup_{i=1}^{\ell}\{z_i=0\}$,
and we will shrink $X$ around
the origin $O$ without mention it.

\subsubsection{Gr of meromorphic flat bundle
associated to the Stokes structure}
Before an explanation of 
the reduction of
(unramifiedly good) wild harmonic bundles,
let us explain the procedure
to take Gr with respect to Stokes structure
for a meromorphic flat bundle $(\nbigv,\nabla)$
with an unramifiedly good lattice $V$.
(See Chapter \ref{section;07.12.26.6}
 for more details.)
For each small multi-sector $S\subset X-D$,
we have the full Stokes filtration $\nbigftilde^S$
of $V_{|\Sbar}$ on $\Sbar$.
We obtain a graded bundle
$\Gr^{\nbigftilde}(V_{|\Sbar})=
 \bigoplus\Gr^{\nbigftilde}_{\gminia}(V_{|\Sbar})$ 
on $\Sbar$ associated to $\nbigftilde^S$.
Although the filtrations depend on multi-sectors $S$,
they satisfy some compatibility.
Hence, 
we can glue $\Gr^{\nbigftilde}_{\gminia}(V_{|\Sbar})$
and obtain the bundle
$\Gr^{\nbigftilde}_{\gminia}(V_{|\Xtilde(D)})$
with the induced meromorphic flat connection
$\nabla_{\gminia}$
on the real blow up $\Xtilde(D)$.
It can be shown that it is the pull back of
a meromorphic flat connection
$(\Gr^{\nbigftilde}_{\gminia}(V),\nabla_{\gminia})$
on $(X,D)$,
which is defined to be Gr of $(V,\nabla)$
with respect to the full Stokes structure.
(We can find that such a construction 
has already appeared in
\cite{Deligne-Malgrange-Ramis}
for meromorphic flat bundles
on curves.)

Although it is essentially the same as taking
direct summands of the decomposition 
(\ref{eq;08.10.24.11}),
there are some advantages.
The above construction fits
to our viewpoint that
a meromorphic flat bundles on $X$ is
prolongment of a flat bundle on $X-D$.
Moreover, it is suitable for the reduction of
a variation of pure twistor structure,
explained below.

\subsubsection{Gr of family of meromorphic 
$\lambda$-flat bundles associated to the Stokes structure}

Let $\harmonicbundle$ be an {\em unramifiedly}
good wild harmonic bundle on $(X,D)$
with a decomposition (\ref{eq;08.10.24.10}).
We use the notation in 
Section \ref{subsection;08.10.24.20}.
We set $W:=\nbigd\cup(\{0\}\times X)$.
Let $\nbigxtilde(W)$ be the real blow up
of $\nbigx$ along $W$.
Let $S$ be a small multi-sector in 
$\nbigx-W$.
As in the case of ordinary meromorphic flat bundles,
we have the full Stokes filtration
$\nbigftilde^S$ of $\nbigq\nbige_{|\Sbar}$.
By varying $S$ and gluing 
$\Gr^{\nbigftilde}_{\gminia}(\nbigq\nbige_{|\Sbar})$,
we obtain a family of $\lambda$-flat bundles
on $\nbigxtilde(W)$.
Moreover,
as the descent for
$\nbigxtilde(W)\lrarr\nbigx$,
we obtain a family of meromorphic $\lambda$-flat bundles
$\Gr^{\nbigftilde}_{\gminia}(\nbigq\nbige)$ on 
$\cnum_{\lambda}\times (X,D)$
for each $\gminia\in\Irr(\theta)$.
It has the unique irregular value $\gminia$.
They are called the full reduction
of $(\nbigq\nbige,\DD)$.

\subsubsection{Gr of 
 variation of pure twistor structure}

From the unramifiedly good wild harmonic bundle
$\harmonicbundle$ on $X-D$,
we obtain an unramifiedly good wild harmonic bundle
$(E,\del_E,\theta^{\dagger},h)$ on the conjugate
$(X^{\dagger},D^{\dagger})$.
We have the associated family of $\mu$-flat bundles
$(\nbige^{\dagger},\DD^{\dagger})$
on $\cnum_{\mu}\times (X^{\dagger}-D^{\dagger})$,
which is prolonged to
a family of meromorphic $\mu$-flat bundles
$\bigl(\nbigq\nbige^{\dagger},\DD^{\dagger}\bigr)$
on $\cnum_{\mu}\times (X^{\dagger},D^{\dagger})$.
We also obtain the full reductions
$\Gr^{\nbigftilde}_{\gminiabar}(\nbigq\nbige^{\dagger})$
on $\cnum_{\mu}\times(X^{\dagger},D^{\dagger})$
for any $\gminia\in\Irr(\theta)$.
Note $\Irr(\theta^{\dagger})=\bigl\{
 \gminiabar\,\big|\,\gminia\in\Irr(\theta)
 \bigr\}$.

Let $S$ be a small multi-sector in $X-D$.
Let $U(\lambda_0)$ be a small neighbourhood
of $\lambda_0\neq 0$ in $\cnum_{\lambda}$.
We have the full Stokes filtration
$\nbigftilde^S\bigl(
 \nbigq\nbige_{|U(\lambda_0)\times\Sbar}\bigr)$ of
$\nbigq\nbige_{|U(\lambda_0)\times\Sbar}$.
Let $U(\mu_0)$ be the neighbourhood of
$\mu_0=\lambda_0^{-1}$ in $\cnum_{\mu}$,
corresponding to $U(\lambda_0)$.
We also have the full Stokes filtration
$\nbigftilde^{S}\bigl(
 \nbigq\nbige^{\dagger}_{|U(\mu_0)\times
 \Sbar^{\dagger}}\bigr)$.
We can observe that
the filtrations are essentially the same on
$U(\lambda_0)\times S=U(\mu_0)\times S^{\dagger}$
(Proposition \ref{prop;07.11.24.29}).
Actually, they are characterized 
by the growth order of the norms of flat sections.
Hence, we have a natural isomorphism
$\Gr^{\nbigftilde}_{\gminia}
 (\nbigq\nbige_{|U(\lambda_0)\times S})
\simeq
 \Gr^{\nbigftilde}_{\gminiabar}
 (\nbigq\nbige^{\dagger}_{|U(\mu_0)\times S^{\dagger}})$.
By gluing them,
we obtain a natural identification:
\begin{equation}
 \label{eq;07.12.28.20}
 \Gr^{\nbigftilde}_{\gminia}
 (\nbigq\nbige)_{|\cnum_{\lambda}^{\ast}\times (X-D)}
\simeq
 \Gr^{\nbigftilde}_{\gminiabar}
 (\nbigq\nbige^{\dagger})_{
 |\cnum_{\mu}^{\ast}\times(X-D)}
\end{equation}

Recall that 
the gluing of $(\nbige,\DD)$ and 
$(\nbige^{\dagger},\DD^{\dagger})$
gives a variation of pure twistor structure
$(\nbige^{\sankaku},\DD^{\sankaku})$ of weight $0$
on $\proj^1\times (X-D)$
with a polarization
$\nbigs:\nbige^{\sankaku}
  \otimes\sigma^{\ast}\nbige^{\sankaku}
 \lrarr\Tate(0)$.
(See \cite{s3} or \cite{mochi2}.
 We will review it
 in Section \ref{subsection;08.10.25.20}.)
Because of the isomorphism (\ref{eq;07.12.28.20}),
we obtain a variation of twistor structure
$\bigl(
 \Gr_{\gminia}(\nbige^{\sankaku}),
 \DD^{\sankaku}_{\gminia}
\bigr)$
for each $\gminia\in\Irr(\theta)$,
on which we have the induced pairings 
$\nbigs_{\gminia}$.
The following theorem is one of the most important
results in this paper.

\begin{thm}[Theorem 
\ref{thm;07.10.11.120}]
 \label{thm;08.10.25.10}
$\bigl(\Gr_{\gminia}(\nbige^{\sankaku}),
 \DD^{\sankaku}_{\gminia}, \nbigs_{\gminia}\bigr)$ 
is a variation of pure polarized twistor structure
of weight $0$,
if we shrink $X$ appropriately.
It comes from a harmonic bundle,
which is the tensor product of
a tame harmonic bundle
$(E_{\gminia},\delbar_{\gminia},
 \theta_{\gminia},h_{\gminia})$
and $L(\gminia)$.

Here, $L(\gminia)$ denotes a harmonic bundle,
which consists of
a line bundle $\nbigo_{X-D}\cdot e$,
the Higgs field $\theta e=e\cdot d\gminia$
and the metric $h(e,e)=1$.
\hfill\qed
\end{thm}

In some sense,
Theorem \ref{thm;08.10.25.10}
reduces the study of the asymptotic behaviour
of wild harmonic bundles to the tame case.
For example,
the completion of $\nbigq\nbige$
along $\cnum_{\lambda}\times \{O\}$
is naturally isomorphic to 
the direct sum of the completion of
$\nbigl(\gminia)\otimes
 \nbigq\nbige_{\gminia}$.
With the detailed study on
$\nbigq\nbige_{\gminia}$ for tame harmonic bundles
in \cite{mochi2},
we can say that we have already understood
$\nbigl(\gminia)\otimes
 \nbigq\nbige_{\gminia}$ pretty well,
and hence $\nbigq\nbige$.
Such an observation is very useful 
when we apply the prolongment of
good wild harmonic bundles
to the theory of polarized wild pure twistor $D$-modules.
We can also derive the norm estimate.

\subsubsection{Uniqueness of prolongation}

Recall the uniqueness of prolongation
of a flat bundle on $X-D$
to a meromorphic flat bundle on $(X,D)$
with regular singularity,
by which we have a very easy characterization
of $\nbigq\nbige$ in the {\em tame} case.
Namely, assume that
we have some family of meromorphic
$\lambda$-flat bundles $\nbigv$ on 
$\cnum_{\lambda}\times(X,D)$
such that (i) the restriction to 
$\cnum_{\lambda}\times(X-D)$ is
$(\nbige,\DD)$,
(ii) each restriction $\nbigv^{\lambda}$ $(\lambda\neq 0)$
is regular.
Then, we have the natural isomorphism
$\nbigv\simeq\nbigq\nbige$,
if $\harmonicbundle$ is tame.
However, in the non-tame case,
we do not have such an obvious characterization,
which was one of the main obstacles
for the author in this study.
He has not yet known
whether there exists an easy characterization
for meromorphic prolongation of 
a family of $\lambda$-flat bundles
with good lattices.
However, we have a useful characterization
of meromorphic prolongation of a variation 
of polarized pure twistor structure.
Let $\bigl(\Vtilde_0,\Vtilde_{\infty}\bigr)$
be an unramifiedly good meromorphic prolongment of 
$(\nbige^{\sankaku},\DD^{\sankaku},\nbigs)$.
(See Section \ref{subsection;07.10.10.2}.)
We have the variation of
$\proj^1$-holomorphic bundles
 $\bigl(\Gr^{\Vtilde}_{\gminia}(\nbige^{\sankaku}),
 \DD^{\sankaku}_{\gminia}\bigr)$
with the pairing $\nbigs_{\gminia}$
for each $\gminia\in\Irr(\theta)$,
obtained as the full reduction
with respect to the Stokes structure 
of $(\Vtilde_0,\Vtilde_{\infty})$.

\begin{thm}[Theorem 
\ref{thm;07.10.11.10}]
\label{thm;08.10.25.15}
If $\bigl(\Gr^{\Vtilde}_{\gminia}(\nbige^{\sankaku}),
 \DD^{\sankaku}_{\gminia},\nbigs_{\gminia}
 \bigr)$
 are variations of 
polarized pure twistor structure of weight $0$ 
 for any $\gminia\in\Irr(\theta)$,
the prolongment is canonical,
i.e.,
$\Vtilde_0\simeq\nbigq\nbige$ and
$\Vtilde_{\infty}\simeq\nbigq\nbige^{\dagger}$.
\hfill\qed
\end{thm}

\section{Application to meromorphic flat bundles}

\subsection{Resolution of turning points
 of meromorphic flat bundles}

We recall the notion of Deligne-Malgrange lattice.

\subsubsection{Deligne lattice}

Let $X$ be a complex manifold,
and let $D$ be a normal crossing hypersurface
with the irreducible decomposition $D=\bigcup_{i\in I}D_i$.
Let $(\nbigv,\nabla)$ be a meromorphic flat connection
on $(X,D)$.
Namely, $\nbigv$ is an $\nbigo_X(\ast D)$-module
with a flat connection
$\nabla:\nbigv\lrarr\nbigv\otimes\Omega^1_X$.
If $(\nbigv,\nabla)$ has regular singularity along $D$,
there exists a lattice $V\subset\nbigv$ 
with the following property:
\begin{itemize}
\item
 $\nabla$ is logarithmic with respect to $V$
 in the sense 
 $\nabla(V)\subset V\otimes\Omega^{1}(\log D)$.
 Note that
 the residue $\Res_{D_i}(\nabla)$ is given
 in $\End(V_{|D_i})$
 for each irreducible component of $D$.
\item
 Any eigenvalue $\alpha$ of $\Res_{D_i}(\nabla)$ 
 satisfies $0\leq \Re(\alpha)<1$.
\end{itemize}
This lattice $V$ is called the {\em Deligne lattice} of
$(\nbigv,\nabla)$,
which plays an important role
in the study of meromorphic flat bundles
with regular singularities,
or more generally, 
regular holonomic $D$-modules.
\index{Deligne lattice}

\subsubsection{Deligne-Malgrange lattice 
 (one dimensional case)}

It is natural and important to ask the existence of 
such a lattice in the irregular case.
If the base manifold is a curve,
it is classically well known.
Let us consider the case $X:=\Delta$
and $D:=\{0\}$.
According to Hukuhara-Levelt-Turrittin theorem,
there exists 
an appropriate ramified covering
$\varphi:(X',D')\lrarr (X,D)$ such that
the formal structure of the pull back 
$\varphi^{\ast}(\nbigv,\nabla)$ is quite simple.
Namely, there exists a finite subset 
 $\Irr(\nabla)\subset M(X',D')/H(X')$
and a decomposition
\[
 \varphi^{\ast}(\nbigv,\nabla)_{|\Dhat'}
=
 \bigoplus_{\gminia\in \Irr(\nabla)}
 (\nbigvhat'_{\gminia},\nablahat'_{\gminia}),
\]
such that each 
$\nablahat_{\gminia}^{\reg}:=
 \nablahat'_{\gminia}-d\gminia$ 
has regular singularity.
We have the Deligne lattices
$\Vhat_{\gminia}'$ for
meromorphic flat bundles with regular singularity
$(\nbigvhat'_{\gminia},\nablahat^{\reg}_{\gminia})$,
and we obtain the formal lattice
\[
\bigoplus_{\gminia\in\Irr(\nabla)}
 \Vhat'_{\gminia}
 \subset
  \varphi^{\ast}(\nbigv)_{|\Dhat'}.
\]
It determines the lattice
$V'\subset\varphi^{\ast}\nbigv$
with a decomposition
\[
 (V',\varphi^{\ast}\nabla)_{|\Dhat'}
=\bigoplus_{\gminia\in\Irr(\nabla)}
 (\Vhat'_{\gminia},\nablahat'_{\gminia}),
\]
such that
(i) $\nablahat'_{\gminia}-d\gminia$
are logarithmic with respect to
$\Vhat_{\gminia}'$ for any $\gminia$,
(ii) any eigenvalues $\alpha$
of the residue satisfy $0\leq \Re(\alpha)<1$.
Since $V'$ is invariant under the action of
the Galois group of this ramified covering,
we obtain the lattice $V\subset\nbigv$
as the descent of $V'$.
This is the {\em Deligne-Malgrange lattice}
in the one dimensional case.
\index{Deligne-Malgrange lattice}

\subsubsection{Good Deligne-Malgrange lattice}

In the higher dimensional case,
the existence of such a lattice 
was proved by B. Malgrange \cite{malgrange}.
But, before recalling his result,
we explain what is an ideal generalization
in the higher dimensional case.
(See Section \ref{section;10.5.4.100}.)

Let $X$ be a complex manifold 
of arbitrary dimension
with a normal crossing hypersurface $D$.
Let $V$ be a lattice of 
a meromorphic flat bundle $(\nbigv,\nabla)$
on $(X,D)$.
We say that $V$ is 
an {\em unramifiedly good Deligne-Malgrange lattice}
if the following holds at any $P\in D$:
\index{unramifiedly good Deligne-Malgrange lattice}
\begin{itemize}
\item
 Let $X_P$ be a small neighbourhood of $P$ in $X$.
 Let $I(P):=\{i\,|\,P\in D_i\}$.
 We set $D_P:=X_P\cap D$
 and $D_{I(P)}:=X_P\cap \bigcap_{i\in I(P)}D_i$.
 Then, we have a finite subset
 $\Irr(\nabla,P)\subset M(X_P,D_P)/H(X_P)$
 and a decomposition
\[
 (V,\nabla)_{|\Dhat_{I(P)}}
=\bigoplus_{\gminia\in \Irr(\nabla,P)}
 (\Vhat_{\gminia},\nablahat_{\gminia})
\]
such that
(i) $\nablahat_{\gminia}-d\gminia$ are logarithmic,
(ii) the eigenvalues $\alpha$ of the residues
 satisfy $0\leq \Re(\alpha)<1$.
 Precisely, we impose the condition that
 $\Irr(\nabla,P)$ is a good set of irregular values
 (Section \ref{subsection;07.11.6.1}).
\end{itemize}
We say that $V$ is a {\em good Delinge-Malgrange lattice},
if the following holds for any $P\in D$:
\index{good Deligne-Malgrange lattice}
\begin{itemize}
\item
If we take an appropriate ramified covering
$\varphi:(X'_P,D'_P)\lrarr (X_P,D_P)$,
there exists an unramifiedly good Deligne-Malgrange
lattice $V'$ of $\varphi^{\ast}(\nbigv,\nabla)$,
and $V_{|X_P}$ is the descent of $V'$.
\end{itemize}
They are uniquely determined, if they exist.
In the one dimensional case,
a Deligne-Malgrange lattice is always good
in this sense.
(See \cite{mochi11} for a different but
equivalent definition.)

\subsubsection{Existence theorem of Malgrange}

However, in general,
a good Deligne-Malgrange lattice
may not exist for a meromorphic flat bundle.
Instead, Malgrange proved the following
in \cite{malgrange}.

\begin{prop}
\label{prop;08.10.24.30}
There exists a $\nbigo_X$-reflexive lattice 
$V\subset \nbigv$
and an analytic subset 
$Z\subset D$ with $\codim_X(Z)\geq 2$
such that 
$V_{|X-Z}$ is a good Deligne-Malgrange lattice
of $(\nbigv,\nabla)_{|X-Z}$.
\hfill\qed
\end{prop}
Although he called it the canonical lattice,
we would like to call it {\em Deligne-Malgrange lattice}.
We have the minimum $Z_0$ among the closed subset
$Z$ as above.
Any point of $Z_0$ is called a {\em turning point}
for $(\nbigv,\nabla)$.
\index{Deligne-Malgrange lattice}
\index{turning point}

\subsubsection{Resolution of turning points}

The Deligne-Malgrange lattice is a very nice clue 
for the study of meromorphic flat bundles.
For example, we will use it to obtain
a Mehta-Ramanathan type theorem
for simple meromorphic flat bundles,
i.e.,
a meromorphic flat bundle is simple
if and only if so is its restriction to sufficiently ample
generic ample hypersurface.
(Section \ref{subsection;08.9.29.100}).
However, the existence of turning points
is a serious obstacle for an asymptotic analysis
of meromorphic flat bundles,
as Sabbah observed in his study
of Stokes structure of meromorphic flat bundles
on complex surfaces.
(We have already mentioned
 asymptotic analysis for meromorphic flat bundles
 in Section \ref{subsection;08.10.24.20}.)
He proposed a conjecture to expect
the existence of a resolution of turning points.
We established it
for algebraic meromorphic flat bundles
on surfaces \cite{mochi6}.
We will establish the following theorem 
in the higher dimensional case.

\begin{thm}[Theorem
 \ref{thm;07.10.14.60},
 Corollary \ref{cor;08.12.18.5}]
\label{thm;07.12.5.1}
Let $X$ be a smooth proper complex algebraic variety
with a normal crossing hypersurface $D$.
Let $(\nbigv,\nabla)$ be a meromorphic flat connection
on $(X,D)$.
Then, there exists a birational projective morphism
$\varphi:X'\lrarr X$ such that
(i) $D':=\varphi^{-1}(D)$ is 
a simple normal crossing hypersurface,
(ii) $X'-D'\simeq X-D$,
(iii) $\varphi^{\ast}(\nbigv,\nabla)$
has no turning points.
\hfill\qed
\end{thm}

Theorem \ref{thm;07.12.5.1} seems
of foundational importance in the study
of algebraic meromorphic flat bundles
or more generally, algebraic holonomic $D$-modules.
Although it is argued in 
Chapter \ref{section;08.10.24.41},
it can be shown more shortly.
Actually,
it follows from some of the results in 
Chapters \ref{section;07.11.8.15},
\ref{section;07.6.2.1},
\ref{section;08.10.24.42}
and Part \ref{part;08.9.29.21}.
We will briefly discuss ideas of the proof
in Subsection \ref{subsection;08.12.24.10}.

\begin{rem}
Recently,
K. Kedlaya showed a generalization 
with a completely different method
in his excellent work
{\rm\cite{kedlaya}, \cite{kedlaya2}}.
In particular, he established it
for excellent schemes.
In the complex analytic situation,
he obtained a local result.
\hfill\qed
\end{rem}

\subsection{Characterization of semisimplicity
 of meromorphic flat bundles}
\label{subsection;08.10.24.40}

According to Corlette \cite{corlette},
a flat bundle on a {\em smooth projective variety}
has a pluri-harmonic metric,
if and only if 
it is a semisimple object
in the category of flat bundles on $X$.
Such a characterization was established
for meromorphic flat bundles
with regular singularity,
by the work of Jost-Zuo and us
(\cite{JZ2} and \cite{mochi2}).
We would like to generalize it in the irregular 
singular case.
We need a preparation to state the theorem.

\subsubsection{$\sqrt{-1}\real$-good wild harmonic bundles
 and good Deligne-Malgrange lattice}

We explain how good wild harmonic bundles
and good Deligne-Malgrange lattices are
related.
Although we explain the local and unramified
case,
it is easily generalized in the global
and (possibly) ramified case.

Let $X:=\Delta^n$ and
$D:=\bigcup_{j=1}^{\ell}\{z_j=0\}$.
Let $\harmonicbundle$ be 
an unramifiedly good wild harmonic bundle
on $X-D$ with a decomposition 
as in (\ref{eq;08.10.24.10}).
We say that 
$\harmonicbundle$ is an
{\em unramifiedly $\sqrt{-1}\real$-good 
 wild harmonic bundle},
if 
\[
 \Sp(\theta)\subset (\sqrt{-1}\real)^{\ell}.
\]
\index{unramifiedly $\sqrt{-1}\real$-good 
 wild harmonic bundle}
We say that it is a 
{\em $\sqrt{-1}\real$-good wild harmonic bundle},
if it is the descent of an unramifiedly 
$\sqrt{-1}\real$-good wild harmonic bundle.
\index{$\sqrt{-1}\real$-good wild harmonic bundle}
As mentioned,
these definitions can be extended
to the global and (possibly) ramified case.
Then, good Deligne-Malgrange lattices naturally appear
in the study of harmonic bundles 
by the following proposition,
which immediately follows from 
Proposition \ref{prop;07.7.19.31}.
\begin{prop}[Proposition
 \ref{prop;08.12.18.10}]
Let $\harmonicbundle$ be a $\sqrt{-1}\real$-good wild
harmonic bundle on $(X,D)$.
Then, $(\nbigp_0\nbige^1,\DD^1)$ is 
the good Deligne-Malgrange lattice of
the meromorphic flat bundle
$(\nbigp\nbige^1,\DD^1)$.
\hfill\qed
\end{prop}

\subsubsection{Characterization}

Let $X$ be a smooth projective variety,
and let $D$ be a normal crossing
hypersurface of $X$.
Let $(\nbigv,\nabla)$ be a meromorphic flat bundle
on $(X,D)$.
Recall that there exists 
a closed subset $Z\subset D$
with $\codim_X(Z)\geq 2$
such that $(\nbigv,\nabla)_{|X-Z}$
has a good Deligne-Malgrange lattice,
according to Proposition \ref{prop;08.10.24.30}.
The next theorem gives a nice characterization
of semisimplicity of $(\nbigv,\nabla)$.

\begin{thm}[Theorem
\ref{thm;07.10.15.1}]
\label{thm;08.9.23.2}
\label{thm;07.12.5.2}
The following conditions are equivalent.
\begin{itemize}
\item
 $(\nbigv,\nabla)$ is semisimple
 in the category of meromorphic flat bundles.
\item
 There exists a $\sqrt{-1}\real$-good wild harmonic bundle
 $\harmonicbundle$ on $(X-Z,D-Z)$
 such that $\nbigp \nbige^1_{|X-Z}\simeq
 \nbigv_{|X-Z}$.
\end{itemize}
Such a metric is unique
up to obvious ambiguity.
\hfill\qed
\end{thm}

\subsubsection{Outline of the proof of
 Theorems \ref{thm;07.12.5.1}
 and \ref{thm;07.12.5.2}}
\label{subsection;08.12.24.10}

At this moment,
the proof of the theorems
are given by the following
flow of the arguments,
which is an interesting interaction
between the theories of
harmonic bundles and meromorphic flat bundles.
(But, see Remarks \ref{rem;10.6.28.1}
and \ref{rem;10.6.28.2}.)
\begin{description}
\item[Step 0.]
 In the curve case,
 Theorem \ref{thm;07.12.5.1} is classical,
 and Theorem \ref{thm;07.12.5.2}
 was known by the work of
 Biquard-Boalch, Sabbah and Simpson.
\item[Step 1.]
 We established Theorem \ref{thm;07.12.5.2}
 in the case $\dim X=2$
 by using the mod $p$-reduction method in \cite{mochi6}.
\item[Step 2.
 (Theorem \ref{thm;08.9.23.2}
 in the case $\dim X=2$)]
This step is the motivation for the author
to study resolution of turning points.
We would like to find a pluri-harmonic metric
of $(\nbigv,\nabla)_{|X-D}$,
for which there is a standard framework in global analysis.
It is briefly and imprecisely as follows:
(i) take an initial metric,
(ii) deform it along the flow given by
 a heat equation,
(iii) the limit  of the heat flow
 should be a pluri-harmonic metric.
Simpson \cite{s1} essentially established
a nice general theory for (ii) and (iii),
once an appropriate initial metric is taken in (i).
To construct an initial metric,
we have to know the local normal form
of meromorphic flat bundles.
It requires a resolution of turning points
in Step 1.

We should remark that we cannot
directly use the above framework,
even if $(\nbigv,\nabla)$ has no turning points.
It will be achieved by the argument in 
\cite{mochi4} and \cite{mochi5}
prepared for Kobayashi-Hitchin correspondence
of meromorphic flat bundles with regular singularities.

\item[Step 3.
 (Theorem \ref{thm;07.12.5.1}
 in the case $\dim X=n$ $(n\geq 3)$)]
This is the easiest part.
We have the following Mehta-Ramanathan type
theorem.
\begin{prop}[Proposition 
 \ref{prop;06.8.12.15}]
$(\nbigv,\nabla)$ is simple
if and only if
$(\nbigv,\nabla)_{|Y}$ is simple
for an arbitrarily ample generic
hypersurface $Y$.
\hfill\qed
\end{prop}
The ``if'' part is clear,
and the other side is non-trivial.
This kind of claim is very standard
for classical stability conditions in algebraic geometry.
The Deligne-Malgrange lattice is equipped
with the natural parabolic structure,
and Sabbah essentially observed that
simplicity and parabolic stability are equivalent.
Then, applying the arguments
due to Mehta and Ramanathan
(\cite{mehta-ramanathan1} and 
 \cite{mehta-ramanathan2})
we will obtain the desired equivalence.

Then, the inductive argument is easy.
For any general and sufficiently ample 
hypersurface $Y\subset X$,
there exists a pluri-harmonic metric $h_Y$ for 
$(\nbigv,\nabla)_{|Y}$.
There exists a finite subset $Z\subset X$
such that $X-Z$ is covered by 
such hypersurfaces $Y$.
So, for $P\in X-Z$,
take $Y$ such as $P\in Y$,
and we would like to define
$h_{|P}:=h_{Y|P}$.
We have to check
$h_{Y_1|Y_1\cap Y_2}
=h_{Y_2|Y_1\cap Y_2}$,
but it follows from 
the uniqueness
because $\dim (Y_1\cap Y_2)\geq 1$.
Thus, we obtain the desired metric.

\item[Step 4.
 (Theorem \ref{thm;07.12.5.3}
 in the case $\dim X=n$ $(n\geq 3)$)]
It can be observed that
we have only to consider the case
in which $(\nbigv,\nabla)$ is simple
(Corollary \ref{cor;10.5.4.101}).
After Step 3,
we take a harmonic bundle
$\harmonicbundle$ on $(X-Z,D-Z)$
as in Theorem \ref{thm;07.12.5.2}.
If $\harmonicbundle$ is a 
$\sqrt{-1}\real$-good wild harmonic bundle
on $(X,D)$,
we have the good Deligne-Malgrange lattice
and $\nbigp_0\nbige^1$ of 
$(\nbigp \nbige^1,\DD^1)$ on $X$.
Because $\codim_X(Z)\geq 2$
and $\nbigp \nbige^1_{|X-Z}\simeq \nbigv_{|X-Z}$,
we have the isomorphism
$\nbigp \nbige^1\simeq \nbigv$ on $X$,
and $\nbigp_0 \nbige^1$ is 
the good Deligne-Malgrange lattice 
of $\nbigv$.
Hence, if $(E,\nabla,h)$ is 
$\sqrt{-1}\real$-good wild on $(X,D)$,
we have nothing to do.

Of course, in general,
$\harmonicbundle$ is not 
$\sqrt{-1}\real$-good wild on $(X,D)$.
However, we have replaced the problem
with the control of the eigenvalues of 
the Higgs field $\theta$,
for which we can use classical techniques
in algebraic or complex geometry.
(See Section \ref{subsection;08.9.28.150}.)
It is much easier than the control of
irregular values of meromorphic flat bundles,
and it can be done.
(See Section \ref{subsection;07.6.12.51}.)
\end{description}

\subsection{Kobayashi-Hitchin correspondence 
 for wild harmonic bundles}

We also have a subject
related to the characterization of stability of
good filtered flat bundles.
Let $X$ be a connected smooth projective variety
of dimension $n$ with an ample line bundle $L$,
and let $D$ be a simple normal crossing hypersurface.
If we are given a good wild harmonic bundle
$(E,\nabla,h)$ on $X-D$,
we obtain the filtered flat bundle
$(\nbigp_{\ast}\nbige^1,\DD^1)$
as in Theorem \ref{thm;07.12.2.55}.
We can show that it is $\mu_L$-polystable,
and each stable component has the trivial
characteristic numbers.
Conversely, we can show the following.

\begin{thm}
[Theorem \ref{thm;06.1.23.100}]
\label{thm;07.12.5.3}
Let $(\vecE_{\ast},\nabla)$ be a good filtered flat bundle
on $(X,D)$.
We put $(E,\nabla):=(\vecE_{\ast},\nabla)_{|X-D}$.
If $(\vecE_{\ast},\nabla)$ is
a $\mu_L$-stable good filtered flat bundle on $(X,D)$
with trivial characteristic numbers
$\pardeg_{L}(\vecE_{\ast})=
\int_X\parch_{2,L}(\vecE_{\ast})=0$,
there exists a pluri-harmonic metric $h$
of $(E,\nabla)$ with the following properties:
\begin{itemize}
\item $(E,\nabla,h)$ is a good wild harmonic bundle on $X-D$.
\item $h$ is adapted to the parabolic structure of
 $\vecE_{\ast}$.
\end{itemize}
Such a pluri-harmonic metric is unique up to obvious ambiguity.
\hfill\qed
\end{thm}

\begin{rem}
\label{rem;10.6.28.1}
If we know Theorem {\rm\ref{thm;07.12.5.1}},
it is not difficult to deduce
Theorem {\rm\ref{thm;07.12.5.2}}
from Theorem {\rm\ref{thm;07.12.5.3}}
as in the tame case,
which will be used in the surface case.
However, we use Theorem {\rm\ref{thm;07.12.5.2}}
to show Theorem {\rm\ref{thm;07.12.5.1}}
in the case $\dim X\geq 3$,
and so we need some different argument.
\hfill\qed
\end{rem}

\begin{rem}
\label{rem;10.6.28.2}
As already mentioned several times,
after the submission of this monograph,
Kedlaya {\rm\cite{kedlaya2}}
obtained the higher dimensional
version of the resolution of turning points
with a different argument.
But, we keep our original flow of the argument.
\hfill\qed
\end{rem}

\section[holonomic $D$-modules and
 wild pure twistor $D$-modules]{Application to 
 holonomic $D$-modules
 and wild pure twistor $D$-modules}
\label{subsection;08.10.26.1}

\subsection{A conjecture of Kashiwara
 on algebraic holonomic $D$-modules}
\label{subsection;08.10.25.2}

Let $\nbigm$ be an algebraic holonomic $D_X$-module 
on a smooth complex algebraic variety $X$.
Let us recall some operations on $D$-modules.
\begin{description}
\item[(Push-forward)]
Let $F:X\lrarr Y$ be a projective morphism
of smooth complex algebraic varieties:
\begin{itemize}
\item
We have the push-forward
$F_{\dagger}\nbigm$
in the derived category 
of cohomologically holonomic $D_Y$-modules.
The cohomology sheaves are denoted by
$F^i_{\dagger}\nbigm$,
which are algebraic holonomic $D_Y$-modules.
\item
We have the Lefschetz morphism
$c_1(L):F_{\dagger}^i\nbigm\lrarr
  F_{\dagger}^{i+2}\nbigm$
for any line bundle $L$ on $X$.
\end{itemize}
\item[(Nearby cycle and vanishing cycle)]
Let $g:X\lrarr \cnum$ be an algebraic function.
By applying the nearby cycle functor 
and the vanishing cycle functor,
we obtain algebraic holonomic $D_X$-modules
$\psi_g(\nbigm)$ and $\phi_g(\nbigm)$.
They are equipped with the induced
nilpotent maps $N$.
By taking Gr with respect to the weight filtrations $W(N)$,
we obtain algebraic holonomic $D_X$-modules
$\Gr^{W(N)}\psi_g(\nbigm)$
and $\Gr^{W(N)}\phi_g(\nbigm)$.

More generally,
according to P. Deligne,
for any $n\in\seisuu_{>0}$
and $\gminia\in\cnum[t_n^{-1}]$,
we obtain an algebraic holonomic $D_X$-module
$\psi_{g,\gminia}(\nbigm)$
by applying the nearby cycle functor
with ramification and the exponential twist by $\gminia$.
(See Section \ref{subsection;08.10.25.1}.
 The author learned this idea from Sabbah.)
\index{nearby cycle functor
 with ramification and exponential twist}
We also obtain a holonomic $D_X$-module
$\Gr^{W(N)}\psi_{g,\gminia}(\nbigm)$
by taking Gr with respect to
the weight filtration of the induced nilpotent map.
\end{description}

There are several works to find 
a subcategory $\nbigc$ of
the category of algebraic holonomic $D$-modules
with the following property:
\begin{itemize}
\item
 $\nbigo_X\in\nbigc$ for 
 any smooth quasi-projective variety $X$.
\item
 $\nbigm_1\oplus\nbigm_2\in\nbigc$
 if and only if $\nbigm_i\in\nbigc$. 
\item 
 $\nbigc$ is stable
 under push-forward 
 for any projective morphism $F:X\lrarr Y$.
 Moreover,
 Hard Lefschetz theorem holds for $\nbigc$
 in the sense that
 $c_1(L)^i:F_{\dagger}^{-i}\nbigm
 \lrarr F_{\dagger}^i\nbigm$ are isomorphisms
 for any $i\geq 0$,
 any projective morphism $F$,
 any relatively ample line bundle $L$,
 and any $\nbigm\in \nbigc$.
\item
 $\nbigc$ is stable
 under the functors $\Gr^W\psi_{g}$
 and $\Gr^W\phi_g$ for any function $g$.
\end{itemize}
For example,
the category of 
(regular) holonomic $D$-modules
is stable for the functors
$F^i_{\dagger}$, $\Gr^W\psi_g$
and $\Gr^W\phi_g$.
However, the Hard Lefschetz theorem does not hold
in general.
In their pioneering work \cite{bbd},
A. Beilinson, J. Bernstein, P. Deligne and O. Gabber
showed the existence of such a subcategory 
called {\em geometric origin}
by using the technique of the reduction
to positive characteristic,
which is the minimum among the subcategories
with the above property.
It is one of the main motivations for this study
to show the following theorem.
\begin{thm}[A conjecture of Kashiwara, 
 Theorem \ref{thm;08.3.6.2}]
\label{thm;08.3.6.1}
The category of algebraic semisimple 
holonomic $D$-modules
has the above property.
Namely, let $X$ be 
a smooth complex algebraic variety,
$\nbigm$ be an algebraic semisimple holonomic 
$D$-module.
Then, the following holds:
\begin{itemize}
\item
 Let $F:X\lrarr Y$ be a projective morphism 
 of smooth quasi-projective varieties.
 Then, $F_{\dagger}^j(\nbigm)$ are also semisimple 
 for any $j$,
 and the morphisms
 $c_1(L)^j:F_{\dagger}^{-j}\nbigm
 \lrarr F_{\dagger}^j\nbigm$
 are isomorphic for any $j\geq 0$
 and any relatively ample line bundle $L$.
 In particular,
 $F_{\dagger}\nbigm$ is isomorphic to
 $\bigoplus F^i_{\dagger}(\nbigm)[-i]$
 in the derived category of
 cohomologically holonomic $D_Y$-modules.
\item
 Let $g$ be an algebraic function on $X$,
 and let $\gminia\in\cnum[t_n^{-1}]$.
 Then,
 $\Gr^W\psi_{g,\gminia}(\nbigm)$
 and $\Gr^W\phi_{g}(\nbigm)$ 
 are also semisimple.
\hfill\qed
\end{itemize}
\end{thm}

The study of this kind of property of $D$-modules
or perverse sheaves was invented by
Beilinson-Bernstein-Deligne-Gabber,
which we have already mentioned above.
M. Saito proved the property for the category of
$D$-modules underlying polarizable pure Hodge modules
in his celebrated work {\rm\cite{saito1}}.
M. Kashiwara conjectured \cite{k5}
that the category of algebraic semisimple holonomic
$D$-modules has the above property.
Sabbah proved in {\rm\cite{sabbah2}} 
the property for regular holonomic $D$-modules 
underlying regular pure imaginary 
polarizable pure twistor $D$-modules.
Simpson {\rm\cite{s10}} also suggested 
such a line of the study.
In {\rm\cite{mochi2}}, we established
the correspondence between
algebraic semisimple regular holonomic $D$-modules
and regular pure imaginary polarizable pure twistor $D$-modules,
and hence the property was proved for
algebraic semisimple regular holonomic $D$-modules.
It was also established by
the works of V. Drinfeld, G. Boeckle-C. Khare and D. Gaitsgory
{\rm(\cite{drinfeld}, \cite{bk}, \cite{gaitsgory})}
via the method of arithmetic geometry
based on the work of L. Lafforgue.
And, M. de Cataldo-L. Migliorini {\rm \cite{cataldo-migliorini}} 
gave another proof of the original result
of Beilinson-Bernstein-Deligne-Gabber
by using their own Hodge theoretic method
but without Saito's method.

\begin{rem}
Contrast to the previous results,
regularity is not assumed
in Theorem {\rm\ref{thm;08.3.6.1}}.
\hfill\qed
\end{rem}

\subsection{Polarized wild pure twistor $D$-module}

Recall that {\em harmonic bundle} is suitable for
the study of {\em semisimplicity} of flat bundles
or $D$-modules from the beginning by Corlette's work.
(Recall also Theorem \ref{thm;07.12.5.2}.)
Then, a natural strategy to attack Theorem \ref{thm;08.3.6.1}
is the following,
which we call Sabbah's program:
\begin{itemize}
\item
 Introduce the category of 
 ``holonomic $D$-modules with pluri-harmonic metrics''
 which should have the property in 
 Section \ref{subsection;08.10.25.2}.
\item
 Show the functorial correspondence
 between ``holonomic $D$-modules 
with pluri-harmonic metrics''
 and algebraic semisimple holonomic $D$-modules.
\end{itemize}

Sabbah introduced the notion of
polarized wild pure twistor $D$-modules
as  ``holonomic $D$-modules with pluri-harmonic metrics''.
We refer to \cite{sabbah2} and 
\cite{sabbah5} for the precise definition
and the basic properties.
\index{polarized wild pure twistor $D$-module}
(We will review it
in Section \ref{subsection;08.9.8.40}
with a preparation given in Chapter \ref{section;08.10.20.40}.)
Needless to say,
it is not obvious at all 
how to think ``pluri-harmonic metrics''
for $D$-modules.

A very important hint was given by Simpson \cite{s3}.
From the beginning of his study \cite{s1},
he was motivated by similarity between
harmonic bundles
and variation of polarized pure Hodge structure.
In \cite{s3}, he introduced the notion of
{\em mixed twistor structure},
and he gave a new formulation of harmonic bundle
as {\em variation of polarized pure twistor structure},
which is formally parallel to the definition of
variation of polarized pure Hodge structure.
\index{mixed twistor structure}
\index{variation of polarized pure twistor structure}
It makes possible for us to formulate
``the harmonic bundle version'' 
(or ``twistor version'')
of  most objects in the theory of
variation of Hodge structure.
And, he proposed a principle,
called Simpson's Meta-Theorem,
that the theory of Hodge structure
should be generalized
to the theory of twistor structure.
\index{Simpson's Meta-Theorem}

In his highly original work (\cite{saito1} and \cite{saito2}),
Saito introduced the notion of
polarized pure Hodge modules
as a vast generalization of
variation of polarized pure Hodge structure,
and he showed that the category of
polarized pure Hodge modules has the nice property,
such as Hard Lefschetz theorem.
It is natural to expect that
we can define ``holonomic $D$-modules
with pluri-harmonic metrics''
as the twistor version of
polarized pure Hodge modules.

And, it was done by Sabbah.
Note that it was still a hard work.
We should emphasize that
he made various useful innovations 
and observations
such as sesqui-linear pairings,
their specialization by using Mellin transforms
(\cite{Barlet-Maire1} and \cite{Barlet-Maire2}),
the nearby cycle functor with ramification
and exponential twist for $\nbigr$-triples,
and so on.

\subsection{Correspondences}

One of the main purposes of our study
is to establish the relation between
algebraic semisimple holonomic $D$-modules
and polarizable wild pure twistor $D$-modules
through wild harmonic bundles:
\begin{equation}
 \label{eq;08.10.25.5}
 \fbox{\rm\begin{tabular}{c}
 semisimple algebraic \\
 holonomic $D$-module
 \end{tabular}}
\leftrightarrow
 \fbox{\rm \begin{tabular}{c}
 $\sqrt{-1}\real$-wild \\
 harmonic bundle
 \end{tabular} }
 \leftrightarrow
 \fbox{\rm
 \begin{tabular}{c}
 polarizable $\sqrt{-1}\real$-wild\\
 pure twistor $D$-module
 \end{tabular}
 }
\end{equation}

\subsubsection{Wild harmonic bundle
 and polarized wild pure twistor $D$-module}

We said that polarized wild pure twistor $D$-modules
were ``holonomic $D$-modules 
with pluri-harmonic metrics'',
as a heuristic explanation.
We make it rigorous by the next theorem.
For simplicity, we consider the case in which
$X$ is a smooth projective variety.
Let $Z$ be a closed irreducible subvariety of $X$.

\begin{thm}
[Theorem \ref{thm;07.10.28.30}]
\label{thm;07.12.6.30}
We have a natural equivalence 
of the categories of the following objects:
\begin{itemize}
\item
 Polarizable wild pure twistor $D$-modules
 whose strict supports are $Z$.
\item
 Wild harmonic bundles 
 defined on Zariski open subsets of $Z$.
\hfill\qed
\end{itemize}
\end{thm}

Theorem \ref{thm;07.12.6.30}
is not only one of the most important key points
in the proof of Theorem \ref{thm;08.3.6.1}
as mentioned in the next subsection,
but it also makes possible for us to consider
``push-forward'' of wild harmonic bundles
in some sense.
In other words, the push-forward for
holonomic $D$-modules is enriched
by polarized pure twistor structure.
It might be useful to investigate the property
of morphisms between moduli spaces
of flat bundles induced by push-forward.
For example,
the study of polarized wild pure twistor $D$-modules 
is related with Fourier transform 
(\cite{malgrange_book}, \cite{sabbah6},
 \cite{sabbah7})
or Nahm transforms 
(\cite{aker-szabo}, \cite{jardim}, \cite{szabo})
for meromorphic flat bundles or wild harmonic bundles
on $\proj^1$.
See also \cite{simpson_middle_convolution}.
In principle,
they should be the specialization
of the corresponding transforms of
polarizable wild pure twistor $D$-modules,
which could be useful for the study
of the corresponding morphisms of
the moduli spaces.

The proof of Theorem \ref{thm;07.12.6.30}
briefly consists of three parts.
\begin{itemize}
\item
 We have to prolong wild harmonic bundles
 $\harmonicbundle$ on a Zariski open subset $U$ of $Z$
 to polarized wild pure twistor $D$-modules on $Z$.
 The most essential case is that 
 $X=Z=\Delta^n$, 
 $D:=X-U=\bigcup_{i=1}^{\ell}\{z_i=0\}$,
 and $\harmonicbundle$ is 
 an unramifiedly good wild harmonic bundle
 on $(X,D)$.
 Since the family of meromorphic $\lambda$-flat bundles
 $(\nbigq\nbige,\DD)$ in Theorem \ref{thm;07.12.4.25}
 is too large,
 we replace it with the ``minimal extension''.
 We need the detailed study on the specialization
 along a function on $X$.
 In the tame case, it was done in \cite{mochi2}.
 The wild case is essentially reduced to the tame case 
 by using Theorem \ref{thm;08.10.25.10}.
\item
 For a given polarized wild pure twistor $D$-module $\nbigt$
 whose strict support is $Z$,
 it is not difficult to show the existence of
 a Zariski open subset $U\subset Z$
 such that $\nbigt_{|X-(Z-U)}$ comes from
 a harmonic bundle $(E,\delbar_E,\theta,h)$ on $U$.
 However, we have to show that
 $(E,\delbar_E,\theta,h)$ is a wild harmonic bundle.
 For that purpose,
 we need various preliminaries 
 such as
 resolution of turning points for meromorphic 
 Higgs field (Section \ref{subsection;07.6.12.51}),
 curve test for wild harmonic bundles
 (Section \ref{subsection;08.9.28.150})
 and so on.
\item
 We have to show the uniqueness of
 prolongation of wild harmonic bundles
 to polarized wild pure twistor $D$-modules.
 In the tame case, this is rather trivial.
 (Recall that a flat bundle is uniquely extended
 to a meromorphic flat bundles with regular singularity.)
 However, in the wild case, it is not obvious.
 It essentially follows from Theorem \ref{thm;08.10.25.15}.
\end{itemize}

\subsubsection{Semisimple holonomic $D$-modules
 and polarizable $\sqrt{-1}\real$-wild pure 
 twistor $D$-modules}

Before stating the next theorem,
we need some preparations.
\begin{itemize}
\item
Recall that there exists the subclass of
$\sqrt{-1}\real$-good wild harmonic bundles
(Section \ref{subsection;08.10.24.40}).
We have the corresponding subcategory of
{\em polarized $\sqrt{-1}\real$-wild 
pure twistor $D$-modules}.
\index{polarized $\sqrt{-1}\real$-wild pure
 twistor $D$-module}
\item
A polarized wild pure twistor $D$-module
is precisely a {\em wild pure twistor $D$-module
with a polarization}.
There is an obvious ambiguity
in the choice of a polarization,
as there exists an obvious ambiguity
in the choice of a pluri-harmonic metric
for a harmonic bundle.
A wild pure twistor $D$-module
is called a {\em polarizable wild pure twistor $D$-module},
if it has a polarization.
\index{polarizable wild pure twistor $D$-module}
\end{itemize}

For a polarizable $\sqrt{-1}\real$-wild pure twistor
$D$-module $\nbigt$,
let $\Xi_{DR}(\nbigt)$ denote
the underlying semisimple holonomic $D$-module.
The next theorem means the correspondence
(\ref{eq;08.10.25.5}).
It essentially follows from Theorem \ref{thm;07.12.6.30}
and Theorem \ref{thm;07.12.5.2}.

\begin{thm}[Theorem 
 \ref{thm;07.10.14.75}]
$\Xi_{DR}(\nbigt)$ is semisimple
for any polarizable $\sqrt{-1}\real$-wild
pure twistor $D$-module $\nbigt$.
Moreover,
$\Xi_{DR}$ gives an equivalence of 
the categories of
polarizable $\sqrt{-1}\real$-wild pure twistor
$D$-modules
and semisimple holonomic $D$-modules
on $X$.
\hfill\qed
\end{thm}

By transferring the operations for 
polarizable $\sqrt{-1}\real$-wild pure twistor $D$-modules,
we obtain that the category of
semisimple holonomic $D$-modules
has the desired property in Section 
\ref{subsection;08.10.25.2}.
It completes the second part of Sabbah's program.

\section{Acknowledgement}
Needless to say,
I regrettably cannot mention all the people
whom I owe gratitude.
Among all,
I am grateful to C. Sabbah who attracted
my attention to the study of
wild pure twistor $D$-modules and
irregular meromorphic flat bundles.
I received much benefit from 
the discussion with him
in the early stage,
and his previous works.
Actually, readers will find his names 
in many places of this monograph.
I also thank his patience and tolerance
to my lack of communication ability.
I express my special thanks to C. Simpson 
for some discussions and his kindness in Nice.
His works and principle (Simpson's Meta-Theorem)
still provide us with the most significant
foundation for our study of harmonic bundles.
During the study,
I realized again the depth of the study of M. Saito.
I also thank his comments to this work.
Comments due to C. Hertling on 
\cite{mochi8} was useful also 
for improving this monograph.
It is my great pleasure to express my gratitude
to K. Aker,
O. Biquard,
R. Budney,
R. Donagi,
H. Esnault,
B. Fantechi,
M. Fujii,
K. Fujiwara,
K. Fukaya,
L. G\"ottsche,
H. Hironaka,
L. Ionescu,
H. Iritani,
J. Jost,
H. Kajiura,
F. Kato,
T. Kato,
Y. Kida,
A. Kono,
M. Kontsevich,
J. Li,
K. Liang,
B. Malgrange,
P. Matsumi,
D. Nadler,
H. Nakajima,
H. Narita,
S. Ohta,
S. Okada,
T. Pantev,
U. A. Rozikov,
M-H. Saito,
C. Sevenheck,
Sz. Szabo,
H. Uehara,
A. Usnich,
D. Wei,
S. Yamaguchi,
K. Yoshioka,
M. Yoshinaga
for some discussions and their kindness.
I appreciate a stimulating communication
with M. Hien.
Special thanks go to K. Vilonen.
I thank Y. Tsuchimoto and A. Ishii
for their constant encouragement.
I would like to express my sincere gratitude
to the referee for his/her valuable comments
and reading this long manuscript.
I am grateful to M. Kashiwara
for his challenging conjecture
which was the strong motivation
for the development of this story.

\vspace{.1in}

This manuscript contains a part of my talks
given in the conferences
``From tQFT to tt$^{\ast}$ and integrability'' in Augsburg,
``New developments in Algebraic Geometry, 
 Integrable Systems and Mirror symmetry'' in Kyoto,
``Partial differential equations and 
differential Galois theory'' in Marseille.
I also improved the presentation for the conferences
``International Conference on Complex Geometry''
in Hanoi,
``Conference in Honor of Z. Mebkhout
 for his Sixtieth Birthday''
in Seville.
It also contains the contents of my lectures
in May 2008 in the University of Tokyo.
I would like to express my gratitude
to the organizers on this occasion.

\vspace{.1in}

I started the study on wild harmonic bundles
during my stay at Max-Planck Institute for Mathematics.
I spent a hard time of trial and error,
and got most of the essential ideas contained in this paper,
during my stay at Institut des Hautes \'Etudes Scientifiques.
Then, I wrote this manuscript at Kyoto University,
Department of Mathematics
and Research Institute for Mathematical Sciences.
I express my gratitude to
the institutions and the colleagues
for the excellent mathematical environments
and the supports.
I also thank the partial financial supports by
Sasakawa Foundation
and Ministry of Education, Culture, 
Sports, Science and Technology.

\part{Good Meromorphic $\varrho$-Flat Bundles}
\label{part;08.9.28.160}

\chapter{Good Formal Property of 
 Meromorphic $\varrho$-Flat Bundle}
\label{section;07.11.8.15}

In this chapter, we shall study
good formal property 
of meromorphic $\varrho$-flat bundle.

First two sections \ref{subsection;07.11.6.1}
and \ref{section;10.5.26.41}
are preliminary.
We recall the notion of
good set of irregular values
in Section \ref{subsection;07.11.6.1}.
We study unramifiedly good lattice
for meromorphic formal $\varrho$-flat bundle
in Section \ref{section;10.5.26.41}.
(See \cite{kedlaya} for deeper results.)
In Section \ref{subsection;08.1.23.1},
we introduce the notion of good lattice
of meromorphic $\varrho$-flat bundle.
It is defined as a lattice whose completions
at any points have nice property.
We hope that the completion along a divisor
has a nice property,
which is studied in Section \ref{subsection;10.5.4.102}.
In Section \ref{section;10.5.26.50},
we introduce the notion of
good filtered $\varrho$-flat bundle,
which will play an important role
in the study of wild harmonic bundles.
In Section \ref{subsection;08.9.28.120},
we introduce the notion of 
good lattice in the level $\vecm$.
It seems useful for our study
on unramifiedly good lattice
for which we use inductive arguments in level.
In Section \ref{section;10.5.4.100},
we restrict ourselves to
ordinary meromorphic flat bundles,
and we study good Deligne-Malgrange lattice.
In Section \ref{section;10.5.20.1},
we prepare some terminology
for family of filtered $\lambda$-flat bundles,
which is significant in our study on
wild harmonic bundles.

\section{Good set of irregular values and truncations}
\label{subsection;07.11.6.1}

\subsection{Definition}
\label{subsection;08.10.28.1}

\subsubsection{The partial order on $\seisuu^{n}$}

We use the partial order $\leq_{\seisuu^{n}}$
(or simply denoted by $\leq$) of $\seisuu^{n}$
given by the comparison of each component,
i.e.,
$\veca\leq_{\seisuu^n}\vecb
\Longleftrightarrow
 a_i\leq b_i,\,(\forall i)$.
\index{order $\leq_{\seisuu^{n}}$}
Let $\veczero$ denote the zero in $\seisuu^n$.
It is also denoted by $\veczero_n$
when we distinguish the dependence on $n$.
\index{vector $\vecdelta_j$}

\subsubsection{Order of poles of meromorphic functions}

Let $\Delta^{\ell}$ denote the multi-disc
\[
\bigl\{(z_1,\ldots,z_{\ell})\,\big|\,
 |z_i|<1,\,i=1,\ldots,\ell \bigr\}.
\]
Let $Y$ be a complex manifold.
Let $X:=\Delta^{\ell}\times Y$.
Let $D_i:=\{z_i=0\}$ and 
$D:=\bigcup_{i=1}^{\ell}D_i$
be hypersurfaces of $X$.
We also put
$D_{\ellsitabar}=\bigcap_{i=1}^{\ell}D_i$,
which is naturally identified with $Y$.
Let $M(X,D)$ (resp. $H(X)$)
denote the space of meromorphic (resp. holomorphic)
functions on $X$ whose poles are contained in $D$.
\index{set $M(X,D)$}
\index{set $H(X)$}
For $\vecm=(m_1,\ldots,m_{\ell})\in\seisuu^{\ell}$,
we put $\vecz^{\vecm}:=\prod_{i=1}^{\ell}z_i^{m_i}$.
For any $f\in M(X,D)$,
we have the Laurent expansion:
\[
 f=\sum_{\vecm\in\seisuu^{\ell}}
f_{\vecm}\bigl(y\bigr)\, \vecz^{\vecm}
\]
Here $f_{\vecm}$ are holomorphic functions
on $D_{\ellsitabar}$.
We will often use the following natural identification
without mention:
\[
 M(X,D)\big/\vecz^{\vecn}\, H(X)\simeq
 \Bigl\{f\in M(X,D)\,\Big|\,
 f_{\vecm}=0,\,\,\forall\vecm\geq\vecn\Bigr\}
\]
Namely, we will often regard an element of
$M(X,D)\big/\vecz^{\vecn}\, H(X)$
as an element of $M(X,D)$ via the above identification.

For $f\in M(X,D)$,
let $\ord(f)$ denote the minimum of the set 
\[
 \nbigs(f):=
 \bigl\{\vecm\in\seisuu^{\ell}\,
 \big|\,f_{\vecm}\neq 0\bigr\}
 \cup \{\veczero_{\ell}\},
\]
if it exists.
Note that we are interested in
the order of poles,
and that $\ord(f)$ is always contained in
$\seisuu_{\leq 0}^{\ell}$
according to this definition (if it exists).
\index{order $\ord(f)$}
We give some examples.
\begin{itemize}
\item
In the case
$f=0$,
we have $\nbigs(f)=
 \bigl\{\veczero\bigr\}$,
and hence 
$\ord(f)=\veczero$.
More generally,
for any holomorphic function $f$,
we have $\ord(f)=\veczero$.
\item
In the case
$f=z_1^{-1}z_2^{-1}
 +z_1^{-1}+z_2^{-1}$,
$\nbigs(f)$ is 
$\bigl\{
 (-1,-1), (-1,0), (0,-1), \veczero
 \bigr\}$, and hence
$\ord(f)=(-1,-1)$.
\item
In the case $f=z_1 z_2^{-1}$,
we have
$\nbigs(f)=\bigl\{
 \veczero,(1,-1)
 \bigr\}$,
and hence $\ord(f)$ does not exist.
\item
In the case
$f=z_1^{-1}+z_2^{-1}$,
we have 
$\nbigs(f)=\bigl\{
 (-1,0),(0,-1),\veczero
 \bigr\}$,
and hence $\ord(f)$ does not exist.
\end{itemize}

For any $\gminia\in M(X,D)/H(X)$,
we take any lift $\gminiatilde$ to $M(X,D)$,
and we set $\ord(\gminia):=\ord(\gminiatilde)$,
if the right hand side exists.
Note that it is independent of the choice of
a lift $\gminiatilde$.
If $\ord(\gminia)\neq 0$,
$\gminiatilde_{\ord(\gminia)}$ is independent
of the choice of a lift $\gminiatilde$,
which is denoted by $\gminia_{\ord(\gminia)}$.

\begin{rem}
Let $k$ be a ring.
The above notion of order makes sense
for the localization of $k[\![z_1,\ldots,z_n]\!]$
with respect to $z_i$ $(i=1,\ldots,\ell)$.
We will not this kind of remark
in the following.
\hfill\qed
\end{rem}

\subsubsection{Good set of irregular values}

We introduce the notion of
good set of irregular values,
which will be used as index sets of 
irregular decompositions and Stokes filtrations.

\begin{df}
\label{df;07.6.1.15}
A finite subset $\nbigi\subset M(X,D)\big/H(X)$ 
is called a good set of irregular values on $(X,D)$,
if the following conditions are satisfied:
\begin{itemize}
\item
$\ord(\gminia)$ exists 
for each element $\gminia\in \nbigi$.
If $\gminia\neq 0$ in $M(X,D)/H(X)$,
$\gminia_{\ord(\gminia)}$
is invertible on $D_{\ellsitabar}$.
\item
$\ord(\gminia-\gminib)$ exists
for any two distinct $\gminia,\gminib\in \nbigi$,
and 
$(\gminia-\gminib)_{\ord(\gminia-\gminib)}$
is invertible on $D_{\ellsitabar}$.
\item
The set 
$\bigl\{\ord(\gminia-\gminib)\,\big|\,
 \gminia,\gminib\in \nbigi
 \bigr\}$
is totally ordered
with respect to the partial order $\leq_{\seisuu^{\ell}}$
on $\seisuu^{\ell}$.
\hfill\qed
\end{itemize}
\end{df}
\index{good set of irregular values}
\index{irregular value}

The third condition is slightly stronger
than that considered in \cite{sabbah4},
which seems convenient for our inductive argument
in levels.

\begin{rem}
The condition in Definition {\rm\ref{df;07.6.1.15}}
does not depend on the choice of a holomorphic 
coordinate system
such that $D=\bigcup_{i=1}^{\ell}\{z_i=0\}$.
\hfill\qed
\end{rem}

\begin{rem}
\label{rem;07.11.15.10}
We will often use a coordinate system
such that $\ord(\gminia-\gminib)$ 
and $\ord(\gminia)$
are contained in the set 
$\coprod_{i=0}^{\ell}
 \seisuu_{<0}^i\times
 \veczero_{\ell-i}$
for any $\gminia,\gminib\in \nbigi$.
Such a coordinate system is called
admissible for $\nbigi$.
\hfill\qed
\end{rem}
\index{admissible coordinate}

\subsubsection{Examples}
\label{subsection;10.5.24.2}

The set $\nbigi_0:=\bigl\{
 z_1^{-1}z_2^{-1}+z_1^{-1}+z_2^{-1}
 \bigr\}$ is a good set of irregular values.
The order of pole is given by
$\ord(z_1^{-1}z_2^{-1}+z_1^{-1}+z_2^{-1})
 =(-1,-1)$.
Let us consider the following examples:
\[
 \nbigi_1:=\bigl\{
 z_1^{-1}z_2^{-1},
 z_1^{-1},0
 \bigr\},
\quad
 \nbigi_2:=\bigl\{
 z_1^{-1}z_2^{-1},
 z_1^{-1},z_2^{-1},0
 \bigr\}
\]
Then, $\nbigi_1$ is a good set of irregular values.
The orders of poles are given as follows:
\[
\ord(z_1^{-1}z_2^{-1})=(-1,-1),
\quad
\ord(z_1^{-1})=(-1,0),
\quad
\ord(0)=(0,0),
\]
\[
 \ord(z_1^{-1}z_2^{-1}-z_1^{-1})
=(-1,-1),
\quad
\ord(z_1^{-1}z_2^{-1}-0)=(-1,-1),
\quad
\ord(z_1^{-1}-0)=(-1,0)
\]
However $\nbigi_2$ is not.
Actually,
$\ord(z_1^{-1}-z_2^{-1})$ does not exist.

We consider the following examples:
\[
 \nbigi_3:=\bigl\{
 z_1^{-1}z_2^{-1}+z_1^{-1},\,
 z_1^{-1}z_2^{-1}
 \bigr\},
\quad
 \nbigi_4:=\bigl\{
 z_1^{-1}z_2^{-1}+z_1^{-1}+z_2^{-1},\,
 z_1^{-1}z_2^{-1}
 \bigr\}
\]
Then, $\nbigi_3$ is a good set of irregular values.
The orders of poles are given as follows:
\[
 \ord(z_1^{-1}z_2^{-1}+z_1^{-1})
=(-1,-1),
\quad
 \ord(z_1^{-1}z_2^{-1})
=(-1,-1),
\]
\[
  \ord\bigl(
 (z_1^{-1}z_2^{-1}+z_1^{-1})
-z_1^{-1}z_2^{-1}
 \bigr)=(-1,0)
\]
However, $\nbigi_4$ is not.
Actually,
$\ord\bigl(
 (z_1^{-1}z_2^{-1}+z_1^{-1}+z_2^{-1})
-z_1^{-1}z_2^{-1}
 \bigr)$ does not exist.

The property of good set of irregular values
is not preserved for some canonical constructions.
For example,
let us consider 
$\nbigi=\{\gminia_i\,|\,i=1,2,3,4\}$
given as follows:
\[
 \gminia_1=z_1^{-1},
\quad
 \gminia_2=2z_1^{-1},
\quad
 \gminia_3=3z_1^{-1}(1+z_2),
\quad
 \gminia_4=4z_1^{-1}(1+z_2)
\]
Then,
$\nbigi\otimes\nbigi_i^{\lor}:=\bigl\{
 \gminia_i-\gminia_j\,\big|\,
 i,j=1,2,3,4
  \bigr\}$
is not necessarily good.
Actually,
\[
 (\gminia_3-\gminia_4)
-(\gminia_1-\gminia_2)
=z_2z_1^{-1}
\]

\subsection{Auxiliary sequence}
\label{subsection;08.1.23.30}
\index{auxiliary sequence}

Let $\nbigi$ be a good set of irregular values on $(X,D)$.
Note that the set $\bigl\{
 \ord(\gminia)\,\big|\,\gminia\in \nbigi
 \bigr\}$ is totally ordered,
because
$\ord(\gminia)\not\leq\ord(\gminib)$
and $\ord(\gminia)\not\geq\ord(\gminib)$
imply that
$\ord(\gminia-\gminib)$ does not exist.
We set 
 $\vecm(0):=\min\bigl\{\ord(\gminia)\,\big|\,
 \gminia\in \nbigi\bigr\}$.
We have the set
$\nbigt(\nbigi):=\bigl\{
 \ord(\gminia-\gminib)\,\big|\,
 \gminia,\gminib\in \nbigi
 \bigr\}$ contained in $\seisuu_{\leq 0}^{\ell}$.
Note 
$\vecm(0)\leq_{\seisuu^{\ell}}\vecm$
for any $\vecm\in \nbigt(\nbigi)$,
because $\gminia_{\vecm}\neq 0$
for some $\gminia\in \nbigi$.
Since $\nbigt(\nbigi)$ is assumed to be
totally ordered with respect to
the partial order $\leq_{\seisuu^{\ell}}$,
we can take a sequence
$\nbigm:=\bigl(
 \vecm(0),\vecm(1),\vecm(2),\ldots,
 \vecm(L),\vecm(L+1)
 \bigr)$ in $\seisuu_{\leq 0}^{\ell}$
with the following property:
\begin{itemize}
\item $\nbigt(\nbigi)\subset \nbigm$
and $\vecm(L+1)=\veczero_{\ell}$.
\item
For each $i\leq L$,
there exists
$1\leq \gminih(i)\leq \ell$ such that
 $\vecm(i+1)=\vecm(i)+\vecdelta_{\gminih(i)}$,
where 
$\vecdelta_j:=
(\overbrace{0,\ldots,0}^{j-1},1,0,\ldots,0)$,
\end{itemize}
Such a sequence is called an auxiliary sequence
for $\nbigi$.
It is not uniquely determined for $\nbigi$.
We often omit to mention $\vecm(L+1)$,
because it is fixed to be $\veczero$.
It seems convenient for an inductive argument.

\begin{rem}
\label{rem;08.1.25.1}
In the case that $D$ is smooth,
i.e., $\ell=1$,
the auxiliary sequence is canonically determined.
We have
$m(0):=\min\bigl\{\ord(\gminia)\,\big|\,
 \gminia\in \nbigi
 \bigr\}$,
and $m(j):=m(0)+j$.
In this case,
we prefer to use the orders
$\ord(\gminia)$ directly.
\hfill\qed
\end{rem}

\subsubsection{Example}
In the example in Section \ref{subsection;08.10.28.1},
\[
 \nbigt(\nbigi_0)=\bigl\{\veczero \bigr\},
\quad
\nbigt(\nbigi_1)=\bigl\{
 (-1,-1),(-1,0),\veczero
 \bigr\},
\quad
\nbigt(\nbigi_3)=\bigl\{
 (-1,0),\veczero
 \bigr\}.
\]
Hence, 
$\nbigm=\bigl\{
 (-1,-1),(-1,0),(0,0)
 \bigr\}$
is an auxiliary sequence
for them.
Note
$\nbigm'=\bigl\{
 (-1,-1),(0,-1),(0,0)
 \bigr\}$
is also an auxiliary sequence for $\nbigi_0$,
but not for $\nbigi_i$ $(i=1,3)$.

\subsection{Truncation}
\label{subsection;07.11.12.20}
\index{truncation}

For any $\vecm\in\seisuu_{\leq 0}^{\ell}$,
let $\eta_{\vecm}:M(X,D)\lrarr M(X,D)$
be given as follows:
\begin{equation}
 \label{eq;07.11.6.2}
 \eta_{\vecm}(\gminia):=
 \sum_{\vecn\leq\vecm}
 \gminia_{\vecn}\,\vecz^{\vecn}
\end{equation}
\index{map $\eta_{\vecm}$}
Let $\xi_{\vecm}:M(X,D)\lrarr M(X,D)$
be given as follows:
\begin{equation}
 \label{eq;07.11.6.3}
 \xi_{\vecm}(\gminia):=
 \sum_{\vecn\not\geq \vecm}
 \gminia_{\vecn}\, \vecz^{\vecn}
\end{equation}
\index{map $\xi_{\vecm}$}
The induced maps $M(X,D)\big/H(X)\lrarr M(X,D)\big/H(X)$
are also denoted by 
$\eta_{\vecm}$ and $\xi_{\vecm}$.

Let $\nbigi$ be a good set of irregular values on $(X,D)$.
We take  an auxiliary sequence 
$\nbigm=\bigl(\vecm(0),\vecm(1),\ldots,\vecm(L+1)\bigr)$
for $\nbigi$.
The function
$\vecz^{-\vecm(i)}\,
\bigl(\xi_{\vecm(i+1)}(\gminia)-
\xi_{\vecm(i)}(\gminia)\bigr)$
is holomorphic on $X$,
and it is independent of the variable
$z_{\gminih(i)}$.
We define
\begin{equation}
 \label{eq;07.11.6.4}
 \zeta_{\vecm(i)}(\gminia):=
 \xi_{\vecm(i+1)}(\gminia)-
 \xi_{\vecm(i)}(\gminia).
\end{equation}
\index{map $\zeta_{\vecm(i)}$}
By construction,
we have
$\xi_{\vecm(i)}(\gminia)
=\sum_{j<i}\zeta_{\vecm(j)}(\gminia)$.

\begin{lem}
\label{lem;07.11.6.5}
For $\gminia_j\in \nbigi$ $(j=1,2)$,
the equality
$\eta_{\vecm(i)}(\gminia_1)
=\eta_{\vecm(i)}(\gminia_2)$
implies
$\xi_{\vecm(i+1)}(\gminia_1)
=\xi_{\vecm(i+1)}(\gminia_2)$,
and hence
$\zeta_{\vecm(i)}(\gminia_1)
=\zeta_{\vecm(i)}(\gminia_2)$.
In particular,
$\xi_{\vecm(i+1)}(\gminib)$
and $\zeta_{\vecm(i)}(\gminib)$ 
are well defined
for $\gminib\in \eta_{\vecm(i)}(\nbigi)$.
\end{lem}
\pf
Because
$\eta_{\vecm(i)}(\gminia_1)=
 \eta_{\vecm(i)}(\gminia_2)$,
we have $\ord(\gminia_1-\gminia_2)\geq \vecm(i+1)$.
Hence, we have $\gminia_{1\,\vecn}=\gminia_{2\,\vecn}$
for any $\vecn\not\geq\vecm(i+1)$,
which implies the claim of the lemma.
\hfill\qed

\vspace{.1in}

When we are given an auxiliary sequence,
it will often be convenient to use the following symbol
for $\gminia\in \nbigi$:
\index{map $\etabar_{\vecm(i)}$}
\begin{equation}
 \label{eq;08.9.13.2}
 \etabar_{\vecm(i)}(\gminia):=
 \xi_{\vecm(i+1)}(\gminia)
\end{equation}
Note $\etabar_{\vecm(L)}(\gminia)=\gminia$
in $M(X,D)\big/H(X)$ for any 
$\gminia\in \nbigi$.
We have the decomposition
$\etabar_{\vecm(i)}(\gminia)
=\sum_{j\leq i}\zeta_{\vecm(j)}(\gminia)$.
The set $\nbigi(\vecm(i)):=\etabar_{\vecm(i)}(\nbigi)$
is called the truncation of $\nbigi$ at 
the level $\vecm(i)$.
It is also a good set of irregular values.
\index{truncation $\nbigi(\vecm(i))$}
We should remark that
$\etabar_{\vecm(i)}$
and the set $\nbigi(\vecm(i))$
depend on the choice of an auxiliary sequence,
in general.

\begin{rem}
In the case that $D$ is smooth,
we have $\etabar_p=\eta_p$
for $p\in\seisuu_{\leq 0}$.
\hfill\qed
\end{rem}

\subsubsection{Example}

Let us consider $\nbigi_0$ in Section
\ref{subsection;08.10.28.1}.
If we take an auxiliary sequence 
$\vecm(0)=(-1,-1),
 \vecm(1)=(-1,0)$,
we have the following:
\[
 \etabar_{\vecm(0)}(z_1^{-1}z_2^{-1}
+z_1^{-1}+z_2^{-1})
=z_1^{-1}z_2^{-1}+z_2^{-1}
\]
If we take an auxiliary sequence 
$\vecm(0)=(-1,-1),
 \vecm(1)=(0,-1)$,
we have the following:
\[
 \etabar_{\vecm(0)}(z_1^{-1}z_2^{-1}
+z_1^{-1}+z_2^{-1})
=z_1^{-1}z_2^{-1}+z_1^{-1}
\]
Let us consider the example
$\nbigi_1$ with the auxiliary sequence
$\vecm(0)=(-1,-1)$,
$\vecm(1)=(-1,0)$.
We have the following:
\[
 \etabar_{\vecm(0)}(z_1^{-1}z_2^{-1})
=z_1^{-1}z_2^{-1},
\quad
 \etabar_{\vecm(0)}(z_1^{-1})=
 \etabar_{\vecm(0)}(0)=0
\]
Hence,
$\nbigi_1(\vecm(0))=\bigl\{
 z_1^{-1}z_2^{-1},0
 \bigr\}$.

Let us consider the example $\nbigi_3$
with an auxiliary sequence
$\vecm(0)=(-1,-1)$,
$\vecm(1)=(-1,0)$.
We have the following:
\[
 \etabar_{\vecm(0)}(z_1^{-1}z_2^{-1}+z_1^{-1})
=\etabar_{\vecm(0)}(z_1^{-1}z_2^{-1})
=z_1^{-1}z_2^{-1}
\]
Hence, $\nbigi_3(\vecm(0))=\{z_1^{-1}z_2^{-1}\}$.

\vspace{.2in}

We have the following picture in our mind
for truncation.

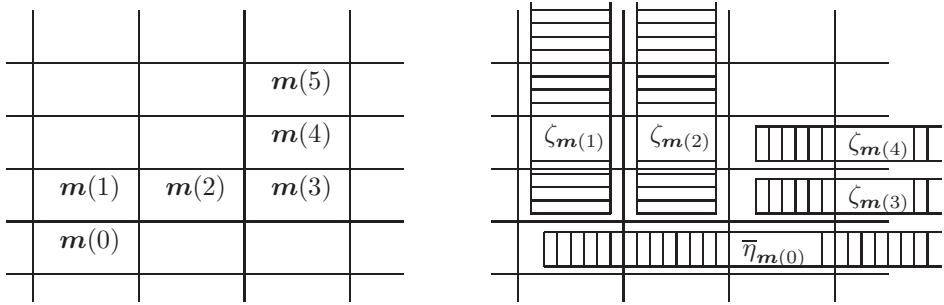
\begin{picture}(180,120)(-100,-80)
\put(0,0){$\vecm(5)$}
\put(0,-20){$\vecm(4)$}
\put(0,-40){$\vecm(3)$}
\put(-40,-40){$\vecm(2)$}
\put(-80,-40){$\vecm(1)$}
\put(-80,-60){$\vecm(0)$}

\put(-100,10){\line(1,0){150}}
\put(-100,-10){\line(1,0){150}}
\put(-100,-30){\line(1,0){150}}
\put(-100,-50){\line(1,0){150}}
\put(-100,-70){\line(1,0){150}}

\put(-90,-80){\line(0,1){110}}
\put(-50,-80){\line(0,1){110}}
\put(-10,-80){\line(0,1){110}}
\put(30,-80){\line(0,1){110}}

\end{picture}
\begin{picture}(150,120)(-100,-80)
\put(-100,10){\line(1,0){150}}
\put(-100,-10){\line(1,0){150}}
\put(-100,-30){\line(1,0){150}}
\put(-100,-50){\line(1,0){150}}
\put(-100,-70){\line(1,0){150}}

\put(-90,-80){\line(0,1){110}}
\put(-50,-80){\line(0,1){110}}
\put(-10,-80){\line(0,1){110}}
\put(30,-80){\line(0,1){110}}

\put(-80,-67){\line(1,0){150}}
\put(-80,-54){\line(1,0){150}}
\multiput(-80,-67)(5,0){15}{\line(0,1){13}}
\multiput(25,-67)(5,0){9}{\line(0,1){13}}
\put(-5,-63){$\etabar_{\vecm(0)}$}

\put(-85,-47){\line(0,1){80}}
\put(-55,-47){\line(0,1){80}}
\multiput(-85,-47)(0,5){5}{\line(1,0){30}}
\multiput(-85,-5)(0,5){8}{\line(1,0){30}}
\put(-80,-20){$\zeta_{\vecm(1)}$}

\put(-45,-47){\line(0,1){80}}
\put(-15,-47){\line(0,1){80}}
\multiput(-45,-47)(0,5){5}{\line(1,0){30}}
\multiput(-45,-5)(0,5){8}{\line(1,0){30}}
\put(-40,-20){$\zeta_{\vecm(2)}$}

\put(0,-47){\line(1,0){70}}
\put(0,-34){\line(1,0){70}}
\multiput(0,-47)(5,0){6}{\line(0,1){13}}
\multiput(60,-47)(5,0){2}{\line(0,1){13}}
\put(35,-43){$\zeta_{\vecm(3)}$}

\put(0,-27){\line(1,0){70}}
\put(0,-14){\line(1,0){70}}
\multiput(0,-27)(5,0){6}{\line(0,1){13}}
\multiput(60,-27)(5,0){2}{\line(0,1){13}}
\put(35,-23){$\zeta_{\vecm(4)}$}

\end{picture}

\begin{center}
$L=4$,
$\vecm(0)=(-2,-3)$,
$\vecm(1)=(-2,-2)$,
$\vecm(2)=(-1,-2)$, \\
$\vecm(3)=(0,-2)$,
$\vecm(4)=(0,-1)$,
$\vecm(5)=(0,0)$.
\end{center}

\section{Unramifiedly good lattice in the formal case}
\label{section;10.5.26.41}
\subsection{Definition}
\label{subsection;10.5.2.5}

We recall some definitions
related with formal connection
for our use.
See \cite{kedlaya} for more deep property
of formal connections.

Let $k$ be an integral domain over $\cnum$.
For some $1\leq\ell\leq n$,
we consider $R_0:=k[\![z_1,\ldots,z_n]\!]$
and its localization $R$ with respect to 
$z_i$ $(i=1,\ldots,\ell)$.
Let $\nbigx$ be the formal scheme
associated to $R_0$.
Let $\nbigd_i$ denote the formal subscheme of 
$\nbigx$
corresponding to $z_i=0$.
We put $\nbigd:=\bigcup_{i=1}^{\ell}\nbigd_i$.
For each $I\subset\ellsitabar$,
we set $\nbigd_I:=\bigcap_{i\in I}\nbigd_i$.
Let $\nbigk:=\Spec k$.
We use natural identifications
$R_0=\nbigo_{\nbigx}$ and 
$R=\nbigo_{\nbigx}(\ast \nbigd)$.

Let $\nbigm$ be an $\nbigo_{\nbigx}$-module.
Let $\varrho\in \nbigo_{\nbigk}$.
Recall that 
a $\varrho$-connection of $\nbigm$
relative to $\nbigk$ is a $k$-linear map
$\DD:\nbigm\lrarr
\nbigm\otimes\Omega^1_{\nbigx/\nbigk}$
such that
$\DD(g\cdot s)
=(\varrho dg)\cdot s+g\cdot \DD s$.
A pairing of 
$\DD(s)\in \nbigm\otimes\Omega^1_{\nbigx/\nbigk}$
and a vector field $v$ of $\nbigx$ is denoted by $\DD(v)s$.
It is called flat,
if the curvature
$\DD\circ\DD
:\nbigm\lrarr\nbigm\otimes\Omega^2_{\nbigx/\nbigk}$
is $0$.
A meromorphic $\varrho$-flat bundle on 
$(\nbigx,\nbigd)$ is 
a free $\nbigo_{\nbigx}(\ast \nbigd)$-module $\nbigm$
equipped with a flat $\varrho$-connection.
If $\varrho$ is nowhere vanishing,
we obtain a flat connection $\varrho^{-1}\DD$
relative to $\nbigk$.
It is often denoted by $\DD^f$.
\index{connection $\DD^f$}

\begin{rem}
We are mainly interested in the cases
(i) $k=\cnum$ and $\varrho=1$ (ordinary flat connection)
(ii) $k=\cnum$ and $\varrho=0$ (Higgs field)
(iii) $k=\cnum$ and $\varrho=\lambda\in\cnum$
 (flat $\lambda$-connection)
(iv)  $k$ is a $\cnum[\lambda]$-algebra,
 and $\varrho=\lambda$
 (family of flat $\lambda$-connections).

We will often omit to say
``relative to $\nbigk$'',
if there is no risk of confusion.
\hfill\qed
\end{rem}

Let $(\nbigm,\DD)$
be a meromorphic $\varrho$-flat bundle on 
$(\nbigx,\nbigd)$.
A coherent $\nbigo_{\nbigx}$-submodule
$\nbigl\subset\nbigm$
is called a lattice,
if $\nbigl\otimes \nbigo_{\nbigx}(\ast \nbigd)=\nbigm$.
The specialization
$\nbigl\otimes \nbigo_{\nbigd_I}$
is denoted by $\nbigl_{|\nbigd_I}$.

\begin{df}
A lattice $\nbigl$ of $\nbigm$ is called
$\gminia$-logarithmic for 
$\gminia\in\nbigo_{\nbigx}(\ast\nbigd)/\nbigo_{\nbigx}$,
if (i) $\nbigl$ is $\nbigo_{\nbigx}$-free,
(ii) $\DD-d\gminiatilde$ is logarithmic
for a lift $\gminiatilde$ of $\gminia$ to
$\nbigo_{\nbigx}(\ast\nbigd)$.
(We will often use the same symbol $\gminia$
to denote a lift to $\nbigo_{\nbigx}(\ast \nbigd)$ 
in the subsequent argument.)

If $\nbigm$ has an $\gminia$-logarithmic lattice,
it is called $\gminia$-regular.
\hfill\qed
\end{df}

\begin{df}
A lattice $\nbigl$ of $\nbigm$ is called
unramifiedly good,
if there exist a good set of irregular values
 $\Irr(\DD)\subset
\nbigo_{\nbigx}(\ast \nbigd)/\nbigo_{\nbigx}$
and a decomposition
\begin{equation}
\label{eq;10.5.2.2}
 (\nbigl,\DD)=
 \bigoplus_{\gminia\in\Irr(\DD)}
 (\nbigl_{\gminia},\DD_{\gminia})
\end{equation}
such that 
$\DD_{\gminia}$ are $\gminia$-logarithmic.

If $\nbigm$ has an unramifiedly good lattice,
we say that $\nbigm$ is unramifiedly good.
\hfill\qed
\end{df}
\index{unramifiedly good lattice}

The decomposition
(\ref{eq;10.5.2.2}) induces
\begin{equation}
 \label{eq;10.5.2.3}
 \nbigm=\bigoplus_{\gminia\in\Irr(\DD)}
 \nbigl_{\gminia}\otimes\nbigo_{\nbigx}(\ast \nbigd).
\end{equation}
The decompositions
(\ref{eq;10.5.2.2}) and (\ref{eq;10.5.2.3})
are called irregular decomposition of
$\nbigl$ and $\nbigm$, respectively.

\begin{lem}
Let $\nbigl$ and $\nbigl'$
be unramifiedly good lattices of $\nbigm$
with irregular decompositions
$\nbigl=
 \bigoplus_{\gminia\in\Irr(\DD)}
\nbigl_{\gminia}$
and
$\nbigl'=
 \bigoplus_{\gminia\in\Irr'(\DD)}
 \nbigl'_{\gminia}$.
Then, we have
\[
 \nbigl_{\gminia}\otimes\nbigo_{\nbigx}(\ast \nbigd)
=\nbigl'_{\gminia}\otimes\nbigo_{\nbigx}(\ast \nbigd)
\]
for any $\gminia\in\Irr(\DD)
 \cup\Irr'(\DD)$.
In particular,
the decomposition {\rm(\ref{eq;10.5.2.2})}
is uniquely determined for $\nbigl$,
and the decomposition {\rm(\ref{eq;10.5.2.3})}
is uniquely determined for $\nbigm$.
\end{lem}
\pf
Take $\gminia,\gminib\in
 \Irr(\DD)\cup\Irr'(\DD)$
such that 
$\gminia-\gminib\neq 0$ in $\nbigo_\nbigx(\ast \nbigd)/\nbigo_\nbigx$.
We would like to show that the induced morphism
$\varphi_{\gminib,\gminia}:
 \nbigl_{\gminia}\otimes\nbigo_\nbigx(\ast \nbigd)
\lrarr \nbigl_{\gminib}'\otimes
 \nbigo_{\nbigx}(\ast \nbigd)$
is $0$.
There exists $1\leq i\leq \ell$
such that the order of $\gminia-\gminib$
with respect to $i$ is strictly smaller than $0$.
We may assume $i=1$.
Let $\nbiga$ be the localization of
$k[\![z_2,\ldots,z_n]\!]$ with respect to
$\prod_{i=2}^{\ell}z_i$,
and let $\nbigr:=\nbiga(\!(z_1)\!)$.
By a standard result in the one variable case
(see Corollary \ref{cor;10.5.2.4} below),
we obtain that the induced morphism
$\nbigl_{\gminia}\otimes\nbigr
\lrarr
 \nbigl'_{\gminib}\otimes\nbigr$ is $0$,
and hence $\varphi_{\gminib,\gminia}=0$.
Then, the claim of the lemma
immediately follows.
\hfill\qed

\subsubsection{Residue and
induced $\DD$-connections}
\label{subsection;10.5.28.4}
\index{residue}
If we are given an unramifiedly good lattice $\nbigl$,
we obtain an endomorphism
$\Res_i(\DD)$
of $\nbigl_{|\nbigd_i}$ in a standard way.
Namely,
for any $f\in \nbigl_{\gminia|\nbigd_i}$,
we take $\ftilde\in \nbigl$ such that
$\ftilde_{|\nbigd_i}=f$, 
and put
$\Res_i(\DD_{\gminia})f:=\bigl(
 \DD_{\gminia}^{\reg}(z_i\del_i)\ftilde
 \bigr)_{|\nbigd_i}$,
where $\DD^{\reg}_{\gminia}:=
 \DD_{\gminia}-d\gminiatilde$
for a lift $\gminiatilde$ of $\gminia$.
We set
$\Res_i(\DD):=
 \bigoplus \Res_i(\DD_{\gminia})
 \in \End(\nbigl_{|\nbigd_i})$.
It is independent of the choice of 
lifts $\ftilde$ and $\gminiatilde$.
It is also independent of the choice of
the coordinate functions $z_i$.
For any $I\ni i$,
the induced endomorphism of $\nbigl_{|\nbigd_I}$
is also denoted by $\Res_i(\DD)$.

If $\gminia$ does not contain
the negative power of $z_i$,
we can define a meromorphic flat 
$\varrho$-connection
$\lefttop{i}\DD_{\gminia}$
on $\nbigl_{\gminia|\nbigd_i}$.
Let $\nbigd(i^c):=\bigcup_{j\neq i}\nbigd_i$.
The section $dz_i/z_i$ of
$\Omega_{\nbigx/\nbigk}^1(\log \nbigd_i)$
induces a splitting
$\Omega_{\nbigx/\nbigk}^1(\log \nbigd_i)_{|\nbigd_i}
\simeq
 \Omega_{\nbigd_i/\nbigk}^1
 \oplus\nbigo_{\nbigd_i}$.
Let $\pi$ denote the projection 
onto $\Omega_{\nbigd_i/\nbigk}^1$.
It induces
$\Omega_{\nbigx/\nbigk}^1(\log\nbigd_i)
 (\ast \nbigd(i^c))
\lrarr
 \Omega_{\nbigd_i/\nbigk}^1(\ast \nbigd(i^c))$.
Then, for $f\in \nbigl_{\gminia|\nbigd_i}$,
take $F\in \nbigl$ such that $F_{|\nbigd_i}=f$,
and put 
$\lefttop{i}\DD(f)
=\pi\bigl((\DD F)_{|\nbigd_i}\bigr)$.
It is independent of the choice of a lift $F$.
But, it depends on the choice of 
the function $z_i$.
If we replace $z_i$ with $\omega\,z_i$
for some invertible $\omega$,
the difference of
the induced $\varrho$-connections
is $(\omega^{-1}d\omega)_{|\nbigd_i}$.

If we are given lifts $\gminiatilde$ 
for any $\gminia\in \nbigi$,
we obtain a flat $\varrho$-connection
$\lefttop{i}\DD^{\reg}_{\gminia}$
of $\nbigl_{\gminia|\nbigd_i}$ by the same procedure,
and $\lefttop{i}\DD^{\reg}:=
 \bigoplus_{\gminia\in\nbigi}
 \lefttop{i}\DD^{\reg}_{\gminia}$.
It depends on the choice of
lifts $\gminiatilde$ and the function $z_i$.

\begin{lem}
\label{lem;10.5.8.1}
$\Res_i(\DD_{\gminia})$
is $\lefttop{i}\DD^{\reg}$-flat.
If $\rho\neq 0$,
the eigenvalues of $\Res_i(\DD)$
are algebraic over $k$.
\end{lem}
\pf
The first claim is clear from the above constructions.
The second claim follows from the first one.
\hfill\qed

\subsubsection{Good lattice}

\begin{df}
A lattice $\nbigl$ of $\nbigm$ is called good,
if there exists a ramified covering
$\varphi:(\nbigx',\nbigd')
\lrarr (\nbigx,\nbigd)$ and
an unramifiedly good lattice $\nbigl'$
of $\nbigm'=\varphi^{\ast}\nbigm$
such that $\nbigl$ is the descent of $\nbigl'$.
\hfill\qed
\end{df}

If we take an $e$-th root $\zeta_i$ of $z_i$
for appropriate $e$,
we have an extension of rings
$R_0\subset R'_0=k[\![\zeta_1,\ldots,\zeta_{\ell},
 z_{\ell+1},\ldots,z_n]\!]$.
Let $G$ be the Galois group of the extension.
We put $\nbigm':=\nbigm\otimes_{R_0}R_0'$.
Then, the above condition says that
$\nbigm'$ has a $G$-equivariant
unramifiedly good lattice $\nbigl'$,
and $\nbigl$ is the $G$-invariant part of $\nbigl'$.

\begin{lem}
\label{lem;10.6.10.1}
Let $\nbigl$ be a good lattice of $\nbigm$.
Put $e_1:=\rank\nbigm!$,
and let $\nbigx_1\lrarr\nbigx$
be a ramified covering such that
the ramification indices at $\nbigd_i$
are $e_1$.
Then, $\nbigl$ is the descent of 
an unramifiedly good lattice $\nbigl_1$
on $\nbigx_1$.
In other words, we have an a priori bound
on the minimal ramification indices.
\end{lem}
\pf
We take $e$, $\nbigx'$,
$\nbigl'$ and $\nbigm'$ as above.
We may assume that 
$e$ is divisible by $e_1$.
We have a factorization
$\nbigx'\lrarr\nbigx_1\lrarr\nbigx$.
Note that 
$\Irr(\DD)\subset
 \nbigo_{\nbigx'}(\ast \nbigd')\big/
 \nbigo_{\nbigx'}$
is contained in
$\nbigo_{\nbigx_1}(\ast \nbigd_1)\big/
 \nbigo_{\nbigx_1}$.
It is well known in the one dimensional case.
The higher dimensional case can be
reduced to the curve case easily.
Then, for the irregular decomposition of
$\nbigl'$,
each direct summand is stable
for the action of the Galois group of
$\nbigx'/\nbigx_1$.
Then, the descent of $\nbigl'$
to $\nbigx_1$ gives the desired lattice.
\hfill\qed

\subsection{A criterion for a lattice to be good}

\subsubsection{Statement}
\label{subsection;10.5.3.1}

Let $\nbigx\lrarr \nbigk$,
$\nbigd$ and $\varrho$ be as in
Section \ref{subsection;10.5.2.5}.
For simplicity, we assume that $k$ is a local ring.
Then, $\nbigx$ has a unique closed point $O$.
Let $(\nbigm,\DD)$ be a meromorphic
flat bundle on $(\nbigx,\nbigd)$.
Let $\nbigl$ be a lattice of $\nbigm$.
Assume that we are given the following:
\begin{itemize}
\item
 A good set of irregular values
 $\nbigi\subset \nbigo_\nbigx(\ast \nbigd)/\nbigo_\nbigx$.
\item
 A decomposition
 $\nbigl=\bigoplus_{\gminia\in \nbigi}\nbigl_{\gminia}$
 as an $\nbigo_\nbigx$-module,
 which is not necessarily compatible with
 $\DD$.
\item
 Let $p_{\gminia}$ denote the projection
 onto $\nbigl_{\gminia}$,
 and we  put
 $\Phi:=\sum_{\gminia\in \nbigi}
 d\gminia\cdot p_{\gminia}$.
 Then,
 $\DD^{(0)}:=\DD-\Phi$
 is logarithmic with respect to $\nbigl$.
 (It is not necessarily flat.)
\end{itemize}

\begin{prop}
\label{prop;10.5.2.21}
$\nbigl$ is an unramifiedly good lattice
of $\nbigm$,
and we have
$\Irr(\DD)=\nbigi$.
\end{prop}

\subsubsection{Preliminary}
\label{subsection;10.5.2.23}

Let $F$ be a free 
$\nbigo_{\nbigx}$-module
with a meromorphic flat $\varrho$-connection
$\DD:
F\lrarr F\otimes\Omega^1_{\nbigx/\nbigk}(\ast \nbigd)$.
Let $\vecv$ be a frame of $F$,
and let $A$ be 
the $\nbigo_\nbigx(\ast \nbigd)$-valued matrix
determined by
$\DD(z_1\del_1)\vecv=\vecv\, A$.
Assume that we have a decomposition
$\vecv=(\vecv_1,\vecv_2)$
such that the corresponding decomposition of
$A$ has the following form:
\[
  A=\left(\begin{array}{cc}
 \Omega_1 & 0 \\ 0 & \Omega_2
 \end{array}\right)
+\left(\begin{array}{cc}
 A_{11} & A_{12}\\
 A_{21} & A_{22}
 \end{array}\right)
\]
\begin{itemize}
\item
The entries of $A_{p,q}$ are regular,
i.e., sections of $\nbigo_{\nbigx}$.
\item
There exists $\vecm\in\seisuu_{<0}^{\ell_1}$
for some $1\leq\ell_1\leq \ell$,
such that the entries of
$\Omegabar_p:=
   \vecz^{-\vecm}\,\Omega_p$ $(p=1,2)$
are regular.
\item
 $\Omegabar_{1|O}$ and
 $\Omegabar_{2|O}$ 
 have no common eigenvalues.
\end{itemize}
Let us consider a change of the base
of the following form:
\[
 \vecv'=\vecv\, G,
\quad
 G=I+\left(\begin{array}{cc} 
 0 & T_2\\ T_1 & 0
 \end{array}
 \right) 
\]
Here, the entries of $T_p$ $(p=1,2)$
are sections of
$z_1\, \nbigo_{\nbigx}$.
We would like to take $G$ 
such that $\DD(z_1\del_1)$
is block-diagonalized with respect to
$\vecv'$ as follows:
\begin{equation}
 \label{eq;10.5.2.10}
 \DD(z_1\del_1)\vecv'=\vecv'\, B,
\quad
 B=\left(
 \begin{array}{cc}
 \Omega_1+Q_1 & 0 \\ 0 & \Omega_2+Q_2
 \end{array}
 \right)
\end{equation}
\begin{lem}
We have regular solutions $T_p$ and $Q_p$
$(p=1,2)$
such that {\rm(\ref{eq;10.5.2.10})} holds.
Moreover, $\vecz^{\vecm}T_p$
are also regular.
\end{lem}
\pf
The relation of $A$, $G$ and $B$ are given by
$AG+\varrho\,z_1\del_1G=GB$.
We obtain the equations:
\[
  A_{11}+A_{12}T_1+Q_{1}=0,
\quad
 \Omega_2T_1+A_{21}+A_{22}T_1
+\varrho\,z_1\del_1T_1
=T_1\Omega_1+T_1Q_1
\]
By eliminating $Q_1$, we obtain the following:
\[
 \Omega_2T_1-T_1\Omega_1+A_{21}+A_{22}T_1
+\varrho\,z_1\del_1T_1+T_1(A_{11}+A_{12}T_1)=0.
\]
We obtain the equation:
\begin{equation}
\label{eq;10.5.2.11}
 \Omegabar_2T_1-T_1\Omegabar_1
+\vecz^{-\vecm}(A_{21}+A_{22}T_1+
T_1A_{11}+\varrho\,z_1\del_1T_1+T_1A_{12}T_1)=0.
\end{equation}
We clearly have a regular solution $T_1$
of (\ref{eq;10.5.2.11}).
Moreover, 
because $\Omegabar_2T_1-T_1\Omegabar_1
\equiv 0$ modulo $\vecz^{-\vecm}$,
we obtain that $\vecz^{\vecm}T_1$ is also regular.
For such $T_1$, $Q_1$ is also regular.
Similarly, we have desired
$T_2$ and $Q_2$,
and $\vecz^{\vecm}T_2$ is also regular.
\hfill\qed

\vspace{.1in}
We put $V_i:=z_i\del_i$ $(i\leq \ell)$
and $V_i:=\del_i$ $(i>\ell)$.
Let $A^{(i)}$ be the $\nbigo_\nbigx(\ast \nbigd)$-valued
matrices determined by
$\DD(V_i)\vecv
=\vecv\,A^{(i)}$.
Assume that the decomposition of $A^{(i)}$
corresponding to $\vecv=(\vecv_1,\vecv_2)$
has the form
\[
 A^{(i)}=
 \left(\begin{array}{cc}
 \Omega^{(i)}_1 & 0 \\ 0 & \Omega^{(i)}_2
 \end{array}\right)
+\left(\begin{array}{cc}
 A^{(i)}_{11} & A^{(i)}_{12}\\
 A^{(i)}_{21} & A^{(i)}_{22}
 \end{array}\right),
\]
where the entries of $A^{(i)}_{p,q}$ are regular,
and the entries of $\vecz^{-\vecm}\Omega^{(i)}_p$
are regular.
Let $C^{(i)}$ be determined by
$\DD(V_i)\vecv'
=\vecv'\, C^{(i)}$.
Since $\vecz^{\vecm}T_p$ $(p=1,2)$ are regular,
 $C^{(i)}$ has the form
\[
 C^{(i)}=\left(
 \begin{array}{cc} 
 \Omega^{(i)}_1 & 0\\
 0 & \Omega^{(i)}_2
 \end{array}
 \right)
+\left(
 \begin{array}{cc} 
 C^{(i)}_{1,1} & C^{(i)}_{1,2}\\
 C^{(i)}_{2,1} & C^{(i)}_{2,2}
 \end{array}
 \right),
\]
where the entries of $C^{(i)}_{p,q}$ are regular.

\begin{lem}
\label{lem;10.5.2.20}
We have $C^{(i)}_{1,2}=0$
and $C^{(i)}_{2,1}=0$.
\end{lem}
\pf
Because
$\bigl[\DD(V_i),
 \DD(z_1\del_1)
 \bigr]=0$,
we have the relation
$\varrho\,V_iC^{(1)}+C^{(i)}C^{(1)}
=\varrho\,z_1\del_1C^{(i)}+C^{(1)}C^{(i)}$,
from which we obtain the following equality:
\[
 C_{12}^{(i)}(\Omega_2+Q_2)
=\varrho\,z_1\del_1C^{(i)}_{12}
+(\Omega_1+Q_1)C^{(i)}_{12}
\]
Then, it is easy to obtain
$C_{1,2}^{(i)}=0$.
Similarly, we can obtain
$C_{2,1}^{(i)}=0$.
\hfill\qed

\subsubsection{Proof of Proposition
 \ref{prop;10.5.2.21}}
\label{subsection;10.5.2.22}

Let us return to the setting in
Subsection \ref{subsection;10.5.3.1}.
We use an induction on the number $|\nbigi|$.
If $|\nbigi|=1$,
the claim of Proposition \ref{prop;10.5.2.21}
is obvious.
Assume that we have already proved
the claim of the proposition in the case $|\nbigi|<m_0$,
and let us show the case $|\nbigi|=m_0$.

We take an auxiliary sequence
$\vecm(0),\ldots,\vecm(L)$ for $\nbigi$.
Let $\nbigi(\vecm(0))$ denote the image of
$\nbigi$ via $\etabar_{\vecm(0)}$.
It is easy to observe that 
we have only to consider the case
$|\nbigi(\vecm(0))|>1$.

\begin{lem}
We have a flat decomposition
\[
(\nbigl,\DD)=
 \bigoplus_{\gminib\in \nbigi(\vecm(0))}
 (\nbigl^{\vecm(0)}_{\gminib},
 \DD^{\vecm(0)}_{\gminib}), 
\]
and an $\nbigo_{\nbigx}$-decomposition
\[
 \nbigl^{\vecm(0)}_{\gminib}
=\bigoplus_{\etabar_{\vecm(0)}(\gminia)=\gminib}
 (\nbigl^{\vecm(0)}_{\gminib})_{\gminia}
\]
with the following property:
\begin{itemize}
\item
 Let $p_{\gminia}'$ denote the projection 
 of $\nbigl^{\vecm(0)}_{\gminib}$
 onto $(\nbigl^{\vecm(0)}_{\gminib})_{\gminia}$,
and we put
 $\Psi^{\vecm(0)}_{\gminib}:=
 \sum_{\gminia\in\etabar_{\vecm(0)}^{-1}(\gminib)}
 d\gminia\, p'_{\gminia}$.
Then,
$\DD^{\vecm(0)}_{\gminib}
-\Psi_{\gminib}^{\vecm(0)}$
are logarithmic with respect to 
$\nbigl^{\vecm(0)}_{\gminib}$.
\end{itemize}
\end{lem}
\pf
Let $\vecv$ be a frame of $\nbigl$ compatible
with the decomposition
$\nbigl=\bigoplus_{\gminia\in \nbigi} 
 \nbigl_{\gminia}$.
Let $\gminia_i$ be determined by
$v_i\in \nbigl_{\gminia_i}$.
Let $\Omega$ be the diagonal matrix valued $1$-form
whose $(i,i)$-th entry is $d\gminia_i$.
By applying Lemma \ref{lem;10.5.2.20}
successively,
we obtain a frame $\vecw$ of $\nbigl$
such that 
$\DD\vecw=\vecw\,\bigl(
 \Omega+B \bigr)$,
where $B$ satisfies the following:
\begin{itemize}
\item
 $B$ is the matrix-valued logarithmic one form.
\item
 $B_{i,j}=0$ unless
 $\etabar_{\vecm(0)}(\gminia_i)
 =\etabar_{\vecm(0)}(\gminia_j)$.
\end{itemize}
For $\gminib\in \nbigi(\vecm(0))$,
let $\nbigl^{\vecm(0)}_{\gminib}$
be the subbundle generated by
$\bigl\{
 w_i\,\big|\,
 \etabar_{\vecm(0)}(\gminia_i)=\gminib
\bigr\}$.
For $\gminia\in \nbigi$ with
$\etabar_{\vecm(0)}(\gminia)=\gminib$,
let $(\nbigl^{\vecm(0)}_{\gminib})_{\gminia}$
be the subbundle generated by
$\bigl\{
w_i\,\big|\,\gminia_i=\gminia
 \bigr\}$.
Then, they have the desired property.
\hfill\qed

\vspace{.1in}
Because $\nbigl^{\vecm(0)}_{\gminib}$
satisfy the assumption in the proposition,
we may apply the hypothesis of the induction.
Thus, the proof of Proposition \ref{prop;10.5.2.21}
is finished.
\hfill\qed

\subsection{Unramifiedly good Deligne-Malgrange lattice}
\label{subsection;08.9.28.103}
\index{unramifiedly good Deligne-Malgrange lattice}

Let us use the setting in Subsection
\ref{subsection;10.5.2.5}
with $k=\cnum$ and $\varrho=1$.
We use the symbol $\nabla$
instead of $\DD$.
\begin{df}
A lattice $\nbigl$ of 
a meromorphic flat bundle $(\nbigm,\nabla)$
on $(\nbigx,\nbigd)$ is called 
an unramifiedly good Deligne-Malgrange lattice,
if (i) $\nbigl$ is an unramifiedly good lattice,
(ii) the eigenvalues $\alpha$
 of $\Res_i(\nabla)$ $(i=1,\ldots,\ell)$
satisfy $0\leq \Re\alpha<1$.
\hfill\qed
\end{df}
\index{unramifiedly good Deligne-Malgrange
lattice}

An unramifiedly good Deligne-Malgrange lattice
is uniquely determined,
if it exists.
If $\nbigl=\bigoplus \nbigl_{\gminia}$
is the unramifiedly good Deligne-Malgrange lattice
of $\nbigm$,
we have frames $\vecv_{\gminia}$ of 
$\nbigl_{\gminia}$ such that
$(\nabla_{\gminia}-d\gminia)\vecv_{\gminia}
=\vecv_{\gminia}\,
 \left(\sum_{j=1}^{\ell}A_j\, dz_j/z_j\right)$,
where $A_j$ are constant matrices.
They induce a frame
$\vecv=(\vecv_{\gminia})$ of $\nbigl$.
Such a frame is called normalizing frame.

\subsubsection{Extension}
\label{subsection;10.4.28.1}

Let $0\lrarr (\nbigm^{(1)},\nabla^{(1)})\lrarr
 (\nbigm^{(0)},\nabla^{(0)})\lrarr 
 (\nbigm^{(2)},\nabla^{(2)})\lrarr 0$
be an exact sequence of
meromorphic flat bundles on $(\nbigx,\nbigd)$.
Let us show the following proposition.

\begin{prop}
\label{prop;07.9.26.105}
Assume $(\nbigm^{(i)},\nabla^{(i)})$ ($i=1,2$)
have the unramifiedly good Deligne-Malgrange
lattices $\nbigl^{(i)}$,
and $\nbigi:=\Irr(\nbigm^{(1)})\cup\Irr(\nbigm^{(2)})$
is good.
Then, $(\nbigm^{(0)},\nabla^{(0)})$
also has the unramifiedly 
good Deligne-Malgrange lattice $\nbigl^{(0)}$.
We also have the exact sequence
$0\lrarr \nbigl^{(1)}\lrarr 
 \nbigl^{(0)}\lrarr\nbigl^{(2)}
 \lrarr 0$.
\end{prop}

Although this also seems to follow some deep results
in \cite{kedlaya},
we keep our elementary proof.
In the following argument,
we assume that the coordinate system is
admissible for $\nbigi$.
An element of $f\in \nbigo_{\nbigx}(\ast \nbigd)$ is
called $z_j$-regular,
if $f$ does not contain the negative power of $z_j$.

\paragraph{First step}

For $\gminia,\gminib\in \nbigi$,
we have the expansion
$\gminia-\gminib=\sum(\gminia-\gminib)_l\, z_p^l$.
We put $\ord_p(\gminia-\gminib):=
 \min\bigl\{l\,\big|\,(\gminia-\gminib)_l\neq 0\bigr\}$
and $\gminis(\gminia-\gminib):=
 \bigl\{p\,\big|\,\ord_p(\gminia-\gminib)<0\bigr\}$.
We put $V_i:=z_i\del_i$ $(i\leq \ell)$
and $V_i:=\del_i$ $(i>\ell)$.
We take normalizing frames $\vecv^{(i)}$ of $\nbigl^{(i)}$
compatible with the irregular decomposition.
We have the decomposition of the frames
$\vecv^{(i)}=\bigl(\vecv^{(i)}_{\gminia}\bigr)$.
Let $A^{(i)}_j$ be determined by
$\nabla^{(i)}(V_j)\vecv^{(i)}
=\vecv^{(i)}\,A^{(i)}_j$.
We have the decomposition $A^{(i)}_j=
 \bigoplus\bigl(V_j(\gminia)+\Abar^{(i)}_{j,\gminia}\bigr)$,
where $\Abar^{(i)}_{j,\gminia}$ 
are constant matrices.
We also have 
$\Abar^{(i)}_{j,\gminia}=0$ for $j>\ell$.

We take a lift $\vecvtilde^{(2)}$ of $\vecv^{(2)}$
to $\nbigm^{(0)}$.
The frame of $\nbigm^{(0)}$
given by $\vecv^{(1)}$ and $\vecvtilde^{(2)}$
is denoted by $(\vecv^{(1)},\vecvtilde^{(2)})$.
Then, we have the following:
\[
 \nabla^{(0)}(V_j)\bigl(\vecv^{(1)},\vecvtilde^{(2)}\bigr)
=\bigl(\vecv^{(1)},\vecvtilde^{(2)}\bigr)
 \,
 \left( \begin{array}{cc}
 A^{(1)}_j & U_j\\
 0 & A^{(2)}_j
 \end{array}
 \right)
\]
Here, the entries of $U_j$ are contained in 
$\nbigo_{\nbigx}(\ast \nbigd)$.
We have the decomposition
$U_j=(U_{j,\gminia,\gminib})$
corresponding to the decompositions
of the frames 
$\vecv^{(i)}=(\vecv_{\gminia}^{(i)})$.
We will consider transforms of frames
of the following form,
where the entries of $W$ are contained in 
$\nbigo_\nbigx(\ast \nbigd)$:
\begin{equation}
 \label{eq;07.12.29.1}
 \bigl(\vecv^{(1)},\vecvtilde^{(2)\prime}\bigr)=
  \bigl(\vecv^{(1)},\vecvtilde^{(2)}\bigr)\,
 \left(
 \begin{array}{cc}
 I & W\\ 0& I
 \end{array}
 \right),
\end{equation}
\[
  \nabla^{(0)}(V_j)
 \bigl(\vecv^{(1)},\vecvtilde^{(2)\prime}\bigr)
=\bigl(\vecv^{(1)},\vecvtilde^{(2)\prime}\bigr)
 \left(
 \begin{array}{cc}
 A^{(1)}_j & U' _j\\
 0& A^{(2)}_j
 \end{array}
 \right)
\]

\begin{lem}
\label{lem;07.9.26.101}
We can take a lift $\vecvtilde^{(2)}$
such that
$U_{j,\gminia,\gminib}$ are $z_p$-regular
for any $p\not\in \gminis(\gminia-\gminib)$
and for any $j$.
\end{lem}
\pf
We will inductively
take transforms as above,
such that the following claims hold:
\begin{description}
\item[$P(m)$]
 Let $p\leq m$.
 If $p\not\in \gminis(\gminia-\gminib)$,
 $U_{p,\gminia,\gminib}$ is $z_p$-regular.
\item[$Q(m)$]
  Let $p\leq m$.
 If $p\not\in \gminis(\gminia-\gminib)$,
 $U_{j,\gminia,\gminib}$ is $z_p$-regular
 for any $j$.
\end{description}
In the condition, we have only to consider
$p$ satisfying $p\leq \ell$.

\begin{lem}
\label{lem;07.9.26.100}
If $P(m)$ holds,
$Q(m)$ also holds.
\end{lem}
\pf
Due to the commutativity
$\bigl[\nabla^{(0)}(V_j),\nabla^{(0)}(V_p)\bigr]=0$,
we have the following relation:
\begin{equation}
 \label{eq;07.11.9.15}
  V_p(U_j)+U_p\, A_j^{(2)}+A_p^{(1)}\, U_j
=V_j(U_p)+U_j\, A_p^{(2)}+A_j^{(1)}\, U_p
\end{equation}
Hence, we obtain the following equality:
\begin{multline}
\label{eq;07.11.9.20}
 z_p\del_p U_{j,\gminia,\gminib}
+V_p(\gminia-\gminib)\, U_{j,\gminia,\gminib}
+\Abar^{(1)}_{p,\gminia}\, U_{j,\gminia,\gminib}
-U_{j,\gminia,\gminib}\, \Abar^{(2)}_{p,\gminib}\\
-V_j\bigl(U_{p,\gminia,\gminib}\bigr)
-V_j(\gminia-\gminib)\, U_{p,\gminia,\gminib}
-\Abar^{(1)}_{j,\gminia}\, U_{p,\gminia,\gminib}
+U_{p,\gminia,\gminib}\,\Abar^{(2)}_{j,\gminib}
=0
\end{multline}
Assume $U_{j,\gminia,\gminib}\neq 0$,
and let us consider the expansion
$U_{j,\gminia,\gminib}
=\sum_{l\geq N}U_{j,\gminia,\gminib,l}\, z_p^l$,
where $U_{j,\gminia,\gminib,N}\neq 0$.
Assume $N<0$, and we will derive a contradiction.
Because $\ord_p(\gminia-\gminib)\geq 0$,
we obtain the following relation:
\[
 N\, U_{j,\gminia,\gminib,N}
+\Abar^{(1)}_{p,\gminia}\, U_{j,\gminia,\gminib,N}
-U_{j,\gminia,\gminib,N}\,\Abar^{(2)}_{p,\gminib}
=0
\]
Because the eigenvalues $\alpha$ of 
$\Abar^{(1)}_{p,\gminia}$
and $\Abar^{(2)}_{p,\gminia}$ satisfy
$0\leq \Re(\alpha)<1$,
we obtain $U_{j,\gminia,\gminib,N}=0$,
which contradicts with our choice of $N$.
Hence, we obtain $N\geq 0$.
Thus, the proof of Lemma \ref{lem;07.9.26.100}
is finished.
\hfill\qed

\vspace{.1in}

Assume $Q(m-1)$ holds for
a lift $\vecvtilde^{(2)}$.
We would like to replace it with a lift
for which $P(m)$ holds,
by successive use of 
the transforms as in (\ref{eq;07.12.29.1}).
We use the following equality,
obtained from the relation
$A_j^{(1)}\, W+U_j+V_j(W)
=U_j'+W\, A_j^{(2)}$:
\begin{equation}
 \label{eq;07.11.9.10}
 V_j(\gminia-\gminib)\, W_{\gminia,\gminib}
+\Abar^{(1)}_{j,\gminia}\, W_{\gminia,\gminib}
-W_{\gminia,\gminib}\, \Abar^{(2)}_{j,\gminib}
+V_j\bigl(W_{\gminia,\gminib}\bigr)
+U_{j,\gminia,\gminib}
-U_{j,\gminia,\gminib}'=0
\end{equation}
If $m\in\gminis(\gminia-\gminib)$,
we have nothing to do,
and so we assume $m\not\in\gminis(\gminia-\gminib)$.
Let us consider the expansion
$U_{m,\gminia,\gminib}
=\sum_{l\geq N}U_{m,\gminia,\gminib,l}\, z_m^l$.
Assume $N<0$.
Let $W_{\gminia,\gminib,N}$ be the unique solution
of the following equation:
\[
  \Abar^{(1)}_{m,\gminia}\,
 W_{\gminia,\gminib,N}
-W_{\gminia,\gminib,N}\,
 \Abar^{(2)}_{m,\gminib}
+N\, W_{\gminia,\gminib,N}
+U_{m,\gminia,\gminib,N}=0.
\]
By the hypothesis $Q(m-1)$ of the induction,
$U_{m,\gminia,\gminib,N}$ is assumed to be
$z_p$-regular
for $p<m$ with $p\not\in\gminis(\gminia-\gminib)$.
Hence, $W_{\gminia,\gminib,N}$ is also $z_p$-regular.
We put
$W_{\gminia,\gminib}=W_{\gminia,\gminib,N}\, z_m^N$.
Then, because of (\ref{eq;07.11.9.10}) with $j=m$,
the obtained $U'_{m,\gminia,\gminib}$
has the expansion
$\sum_{l>N}U'_{m,\gminia,\gminib,l}\, z_m^l$.
Because of (\ref{eq;07.11.9.10}),
$U'_{j,\gminia,\gminib}$ is also $z_p$-regular
for any $j$ and for $p<m$ 
with $p\not\in\supp(\gminia-\gminib)$.
Hence, we can eliminate the negative power
in $U_{m,\gminia,\gminib}$ after the finite procedure,
preserving the condition $Q(m-1)$,
and we can arrive at a lift $\vecvtilde^{(2)}$
for which $P(m)$ holds.

Therefore, after a finite procedure,
we can arrive at a lift $\vecvtilde^{(2)}$
for which $Q(\ell)$ holds.
Thus, the claim of Lemma \ref{lem;07.9.26.101}
is proved.
\hfill\qed

\paragraph{End of the proof of
Proposition \ref{prop;07.9.26.105}}

Let $\vecvtilde^{(2)}$ be a lift
as in Lemma \ref{lem;07.9.26.101}.
We would like to replace it with a lift
for which the $U_j$-components are contained in 
$\nbigo_\nbigx$,
by successive use of (\ref{eq;07.12.29.1}).
We put $F:=z_1\del_1(\gminia-\gminib)$.
Note that
$F\, \vecz^{-\ord(\gminia-\gminib)}$ 
is invertible.
We put 
$W_{\gminia,\gminib}:=
 -F^{-1}\, U_{1,\gminia,\gminib}$.
Then, we have the following,
due to (\ref{eq;07.11.9.10}):
\[
 U'_{1,\gminia,\gminib}
=F^{-1}\,\bigl(
 U_{1,\gminia,\gminib}\,
 \Abar^{(2)}_{1,\gminib}
-\Abar^{(1)}_{1,\gminia}\, U_{1,\gminia,\gminib}
 \bigr)
-V_1\bigl(F^{-1}\, U_{1,\gminia,\gminib}\bigr)
\]
Let $k$ be determined by
$\ord(\gminia-\gminib)\in\seisuu_{<0}^k$,
i.e.,
$\gminis(\gminia-\gminib)=\{1,\ldots,k\}$.
We have the subset $\nbigs\subset\seisuu^k$
and the expansion:
\[
 U_{1,\gminia,\gminib}
=\sum_{\vecn\in\nbigs}
 U_{1,\gminia,\gminib,\vecn}(z_{k+1},\ldots,z_n)
 \,\vecz^{\vecn},
\quad U_{1,\gminia,\gminib,\vecn}\neq 0.
\]
Note that $\nbigs$ is bounded below
with respect to $\leq_{\seisuu^k}$.
Then, the expansion of
$U_{1,\gminia,\gminib}'$ is as follows:
\[
 U'_{1,\gminia,\gminib}
=\sum_{\vecn\in\nbigs_1}
 U'_{1,\gminia,\gminib,\vecn}(z_{k+1},\ldots,z_n)
 \,\vecz^{\vecn},
\quad U'_{1,\gminia,\gminib,\vecn}\neq 0
\]
Here
$\nbigs_1=
 \bigl\{\vecm-\ord(\gminia-\gminib)\,\big|\,
 \vecm\in\nbigs\bigr\}$.
Hence, we can make
$\ord_q(U_{1,\gminia,\gminib})$ sufficiently large
for any $q=1,\ldots,k$ after finite procedure.
So, we have arrived at a lift
$\vecvtilde^{(2)}$,
for which the entries of
$U_{1,\gminia,\gminib}$ are contained in $\nbigo_\nbigx$.

\vspace{.1in}
Let us show 
$\ord_q(U_{j,\gminia,\gminib})\geq 0$ 
for any $q=1,\ldots,k$
and for any $j$.
We have the subset $\nbigs\subset\seisuu^k$
bounded below with respect to $\leq_{\seisuu^k}$,
and the expansion as follows:
\[
 U_{j,\gminia,\gminib}
=\sum_{\vecn\in\nbigs}
 U_{j,\gminia,\gminib,\vecn}(z_{k+1},\ldots,z_n)
 \, \vecz^{\vecn},
\quad
 U_{j,\gminia,\gminib,\vecn}\neq 0
\]
Let $\vecn_0$ be a minimal element of $\nbigs$.
Assume $\vecn_0\not\in\seisuu_{\geq 0}^k$.
Let us look at  
the $\vecz^{\vecn_0+\ord(\gminia-\gminib)}$-term
of (\ref{eq;07.11.9.20}) with $p=1$.
Note that
$V_j(\gminia-\gminib)\, U_{1,\gminia,\gminib}$
does not have the 
$\vecz^{\vecn_0+\ord(\gminia-\gminib)}$-term,
because the entries of $U_{1,\gminia,\gminib}$
are contained in $\nbigo_\nbigx$.
Hence, we obtain
$U_{j,\gminia,\gminib,\vecn_0}\,
 (\gminia-\gminib)_{\ord(\gminia-\gminib)}=0$,
and thus
$U_{j,\gminia,\gminib,\vecn_0}=0$
which contradicts with our choice of $\nbigs$.
Hence we have $\nbigs\in\seisuu_{\geq\,0}^k$.

\vspace{.1in}
Therefore, we have arrived at a lift
$\vecvtilde^{(2)}$ for which
the entries of $U_j$ are contained in $\nbigo_\nbigx$.
Let $\nbigl^{(0)}$ be the submodule 
of $\nbigm^{(0)}$ generated by
$(\vecv^{(1)},\vecvtilde^{(2)})$
over $\nbigo_\nbigx$.
By Proposition \ref{prop;10.5.2.21},
$\nbigl^{(0)}$ is an unramifiedly good lattice.
We have the exact sequence
$0\lrarr\nbigl^{(1)}\lrarr\nbigl^{(0)}
 \lrarr \nbigl^{(2)}\lrarr 0$,
and we can easily deduce that
$\nbigl^{(0)}$ is also Deligne-Malgrange.
\hfill\qed

\subsection{Preliminary from the one variable case
(Appendix)}
\label{section;10.5.3.11}

Let $k$ be an integral domain over $\cnum$.
We consider $\nbigr_0:=k[\![t]\!]$
and $\nbigr:=k(\!(t)\!)$,
which are naturally equipped with
a derivation $\del_t$.
An $\nbigr$-module $\nbigm$ is called differential module,
if it is equipped with the action of $\del_t$
such that $\del_t(f\, s)=
\del_t(f)\, s+f\,\del_ts$
for $f\in \nbigr$ and $s\in \nbigm$.
We recall some basic facts
on differential $\nbigr$-modules from \cite{levelt}
for reference in our argument.

\subsubsection{Extension of decomposition}

Let $\nbigm$ be a finitely generated
differential $\nbigr$-free module
with an $\nbigr_0$-free lattice $\nbigl$
such that
$t^{M+1}\del_t\nbigl\subset\nbigl$
for some $M>0$.
Note that we have an induced endomorphism
$G$ of $\nbigl\otimes_{\nbigr_0}k$.
Assume that there exists decomposition
$(\nbigl\otimes_{\nbigr_0}k,G)
=(V_1,G_1)\oplus (V_2,G_2)$.
For $i\neq j$,
we have the endomorphism
$\Gtilde_{i,j}$ of $Hom(V_i,V_j)$ given by
$\Gtilde_{i,j}(f)=f\circ G_i-G_j\circ f$.
\begin{lem}
\label{lem;10.6.23.10}
If $\Gtilde_{i,j}$ are invertible
for $(i,j)=(1,2),(2,1)$,
then we have a decomposition
$\nbigl=\nbigl_1\oplus\nbigl_2$
such that
(i) $t^{M+1}\del_t\nbigl_i\subset\nbigl_i$,
(ii) $\nbigl_{i}\otimes k=V_i$.
\end{lem}
\pf
We give only a sketch of a proof,
by following \cite{levelt}.
Let $\vecv$ be a frame of $\nbigl$
with a decomposition
$\vecv=(\vecv_1,\vecv_2)$
such that $\vecv_{i|t=0}$ give frames
of $V_i$.
Let $A$ be the $\nbigr_0$-valued matrices
determined by
$t^{M+1}\del_t\vecv=\vecv\, A$.
Then, $A$ has the following decomposition
corresponding to $\vecv=(\vecv_1,\vecv_2)$:
\[
 A=\left(
 \begin{array}{cc}
 \Omega_1 & 0 \\
 0 & \Omega_2
 \end{array}
 \right)
+\left(
 \begin{array}{cc}
 A_{11} & A_{1,2}\\
 A_{2,1}& A_{2,2}
 \end{array}
 \right)
\]
Here, $\Omega_i$ are $k$-valued
matrices determined by
$G_i\vecv_i=\vecv_i\, \Omega_i$,
and $A_{i,j}$ are $t\nbigr_0$-valued
matrices.
We consider a change of the base
of the following form:
\[
 \vecv'=\vecv\, G,
\quad
 G=I+\left(\begin{array}{cc} 
 0 & X\\ Y & 0
 \end{array}
 \right) 
\]
Here, the entries of $X$ and $Y$
are contained in $t\nbigr_0$.
We would like to take $G$ 
such that 
\begin{equation}
 t^{M+1}\del_t\vecv'=\vecv'\, B,
\quad
 B=\left(
 \begin{array}{cc}
 \Omega_1+Q_1 & 0 \\ 0 & \Omega_2+Q_2
 \end{array}
 \right).
\end{equation}
The relation of $A$, $G$ and $B$ are given by
$AG+t^{M+1}\del_tG=GB$.
We obtain the equations
$A_{11}+A_{12}Y+Q_{1}=0$
and
$\Omega_2Y+A_{21}+A_{22}Y+t^{M+1}\del_tY
=Y\Omega_1+YQ_1$.
By eliminating $Q_1$, we obtain the equation
\begin{equation}
\label{eq;10.5.2.1}
\Omega_2Y-Y\Omega_1+A_{21}+A_{22}Y
+t^{M+1}\del_tY+Y(A_{11}+A_{12}Y)=0 
\end{equation}
By the assumption, we have 
the invertibility of the endomorphism
on the space of $k$-valued $(r_2,r_1)$-matrices,
given by $Z\longmapsto
 \Omega_2 Z-Z\Omega_1$,
where $r_i:=\rank\nbigl_i$ $(i=1,2)$.
By using a $t$-expansion,
we can find a solution of (\ref{eq;10.5.2.1})
in the space of
$t\nbigr_0$-valued matrices.
Similarly, we can find desired
$X$ and $Q_2$.
\hfill\qed

\subsubsection{Uniqueness}

\begin{lem}
Let $\nbigm$ be an $\nbigr$-free
differential module.
Assume that there exists
an $\nbigr_0$-free lattice $\nbigl\subset\nbigm$
and $\gminia\in \nbigr\setminus \nbigr_0$
such that 
$t\del_t-t\del_t\gminia$ preserves $\nbigl$.
Then, any flat section of $\nbigm$
is $0$.
\end{lem}
\pf
Take $f\in \nbigm$ such that $\del_tf=0$.
Assume $f\neq 0$,
and we will deduce a contradiction.
We can take $N\in\seisuu$ such that
$t^Nf\in\nbigl$ and the induced element of
$\nbigl/t\nbigl$ is non-zero.
By the assumption,
we have
\[
\nbigl\ni(t\del_t-t\del_t\gminia) (t^Nf)
=(N-t\del_t\gminia)t^Nf
\]
But, it is easy to see that
$(N-t\del_t\gminia)t^Nf\not \in\nbigl$,
and thus we have arrived at a contradiction.
\hfill\qed

\vspace{.1in}
Let $\nbigm_i$ $(i=1,2)$ be 
differential $\nbigr$-free modules
with $\nbigr_0$-free lattices $\nbigl_i$ such that
$t\del_t-t\del_t\gminia_i$ preserve $\nbigl_i$.

\begin{cor}
\label{cor;10.5.2.4}
Assume $\gminia_1-\gminia_2\neq 0$
in $\nbigr/\nbigr_0$.
Then, any flat morphism
$\nbigm_1\lrarr\nbigm_2$
is $0$.
\hfill\qed
\end{cor}

\subsubsection{}

Let $\nbigm$ be a differential $\nbigr$-module.
Let $E$ be an $\nbigr_0$-lattice of $\nbigm$
such that
$t^{m+1}\del_t E\subset E$ for some $m>0$.
We have the induced endomorphism 
$G$ of $E_{|t=0}$.

\begin{lem}
\label{lem;10.5.4.1}
Let $s\in \nbigm$.
If $G$ is invertible,
we have $\del_ts=0$ 
if and only if $s=0$.
\hfill\qed
\end{lem}

Let $E_i$ $(i=1,2)$ be lattices of $\nbigm$
such that
$t^{m_i+1}\del_t E_i\subset E_i$
for some $m_i>0$.
Let $G_i$ be the endomorphism of
$E_{i|t=0}$ induced by $t^{m_i+1}\del_t$.
\begin{lem}
\label{lem;10.5.4.2}
Assume that $G_i$ are semisimple and
non-zero.
Let $T_i$ be the set of eigenvalues of
$G_i$.
Then, we have
$m_1=m_2$ and $T_1=T_2$.
\end{lem}
\pf
By extending $k$,
we may assume that
the eigenvalues of $G_i$
are contained in $k$.
We have $\del_t$-decomposition
$E_i=\bigoplus_{\gminib\in T_i}E_{i,\gminib}$
such that
$E_{i,\gminib|t=0}$
is the eigen space of $G_i$
corresponding to $\gminib$.
We have the induced map
$\varphi_{\gminic,\gminib}:
 E_{1,\gminib}\otimes \nbigr
\lrarr E_{2,\gminic}\otimes \nbigr$.
If $m_1\neq m_2$
or if $m_1=m_2$ but $\gminib\neq\gminic$,
we have $\varphi_{\gminic,\gminib}=0$
by Lemma \ref{lem;10.5.4.1}.
Then, the claim of Lemma \ref{lem;10.5.4.2}
follows.
\hfill\qed

\section{Good lattice of 
 meromorphic $\varrho$-flat bundle}
\label{subsection;08.1.23.1}

\subsection{Definition}
\label{subsection;10.5.3.3}

Let $\nbigx\lrarr \nbigk$ be a smooth fibration of
complex manifolds.
Let $\nbigd$ be a simple normal
crossing hypersurface of $\nbigx$
such that any intersections of irreducible
components are smooth over $\nbigk$.
Let $\varrho$ be a holomorphic function on $\nbigk$.
For a point $P$ of $\nbigx$,
let $\Phat$ denote the completion of $\nbigx$ at $P$.
In the following,
for a given $\nbigo_{\nbigx}$-module $\nbigf$,
let $\nbigf_{|\Phat}$ denote the formal completion
$\nbigf\otimes_{\nbigo_{\nbigx}}\nbigo_{\Phat}$.

Let $(\nbige,\DD)$ be
a meromorphic $\varrho$-flat bundle
on $(\nbigx,\nbigd)$ relative to $\nbigk$,
i.e.,
$\nbige$ is a locally free 
$\nbigo_{\nbigx}(\ast \nbigd)$-coherent sheaf
with a flat $\varrho$-connection
$\DD:\nbige
 \lrarr\nbige\otimes\Omega^1_{\nbigx/\nbigk}$
relative to $\nbigk$.
(A flat $\varrho$-connection is defined
in a standard way as in the formal case.
See Subsection \ref{subsection;10.5.2.5}.
We will often omit ``relative to $\nbigk$''
if there is no risk of confusion.)

\begin{df}
A lattice $E$ of $\nbige$
is called unramifiedly good at $P\in \nbigd$,
if $E_{|\Phat}$ is an unramifiedly good lattice
of $(\nbige,\DD)_{|\Phat}$.
If $E$ is unramifiedly good at any point of $\nbigd$,
$E$ is called an unramifiedly good lattice of
$(\nbige,\DD)$.
\hfill\qed
\end{df}
\index{unramifiedly good lattice}

\begin{notation}
The set of irregular values of 
$(E,\DD)_{|\Phat}$
is often denoted by 
$\Irr(E,\DD,P)$, $\Irr(E,P)$ or $\Irr(\DD,P)$.
\hfill\qed
\end{notation}
\index{set $\Irr(E,\DD,P)$}
\index{set $\Irr(E,P)$}
\index{set $\Irr(DD,P)$}

For $P\in \nbigd$,
let $\nbigx_P$ denote a small neighbourhood of $P$
in $\nbigx$, and put $\nbigd_P:=\nbigx_P\cap \nbigd$.

\begin{df}
\mbox{{}}
\label{df;10.5.3.2}
\begin{itemize}
\item
$(E,\DD)$ 
is called good at $P$,
if there exist a small neighbourhood $\nbigx_P$
and a ramified covering 
$\varphi_P:(\nbigx_P',\nbigd_P')\lrarr (\nbigx_P,\nbigd_P)$
such that $E$ is the descent of 
an unramifiedly good lattice $E'$ of
$\varphi_P^{\ast}\nbige$.
\item
$E$ is called good,
if $E$ is good at any point of $\nbigd$.
\hfill\qed
\end{itemize}
\end{df}
\index{good lattice}
\index{good}

In the condition of Definition {\rm\ref{df;10.5.3.2}},
such $E'$ is not unique,
even if $\varphi_P$ is fixed.
We also remark that $\varphi_P^{\ast}E$
is not necessarily unramifiedly good.

\begin{rem}
We will often say that
$(E,\DD)$ is (unramifiedly) good
on $(\nbigx,\nbigd)$,
if $E$ is a (unramifiedly) good lattice of
a meromorphic $\varrho$-flat bundle
$\bigl(E(\ast \nbigd),\DD\bigr)$
on $(\nbigx,\nbigd)$.
\hfill\qed
\end{rem}

\begin{df}
A meromorphic $\varrho$-flat bundle
is called (unramifiedly) good,
if it locally has an (unramifiedly) good lattice.
\hfill\qed
\end{df}
\index{good meromorphic $\varrho$-flat bundle}
\index{unramifiedly good meromorphic
$\varrho$-flat bundle}

If $\nbigk$ is a point and $\varrho\neq 0$,
a good meromorphic $\varrho$-flat bundle
has a global good lattice.
Actually, it is given by a Deligne-Malgrange
lattice. (See Section {\rm\ref{section;10.5.4.100}}.)

\subsection{Some functoriality}

Let $E$ be a good lattice of
a meromorphic $\varrho$-flat bundle
$(\nbige,\DD)$ on $(\nbigx,\nbigd)$.
Let $E^{\lor}$ denote the dual of
$E$ in the category of $\nbigo_{\nbigx}$-modules,
and $\nbige^{\lor}$ denote the dual of
$\nbige$ in the category of
$\nbigo_{\nbigx}(\ast \nbigd)$-modules.
We have $\nbige^{\lor}=
E^{\lor}\otimes_{\nbigo_{\nbigx}}
 \nbigo_{\nbigx}(\ast \nbigd)$.
We have the naturally induced
flat $\varrho$-connection $\DD$
of $\nbige^{\lor}$.
We have the following functoriality
for dual.
\begin{itemize}
\item
$E^{\lor}$ is a good lattice of
$(\nbige^{\lor},\DD)$.
If $E$ is unramified,
$E^{\lor}$ is also unramified.
For each $P\in\nbigd$,
we have
$\Irr(E^{\lor},P)
=\bigl\{
 -\gminia\,\big|\,
 \gminia\in\Irr(E,P)
 \bigr\}$.
\end{itemize}

Let $E_i$ $(i=1,2)$
be unramifiedly good lattices of $(\nbige_i,\DD_i)$.
We have the following functoriality
for tensor product and direct sum.
\begin{itemize}
\item
If $\Irr(E_1,P)\otimes\Irr(E_2,P):=\bigl\{
 \gminia_1+\gminia_2\,\big|\,
 \gminia_i\in\Irr(E_i,P)
 \bigr\}$ is good for any $P\in\nbigd$,
then $E_1\otimes E_2$ is an unramifiedly
good lattices of $(\nbige_1\otimes\nbige_2,\DD)$
with $\Irr(E_1\otimes E_2,P)=
 \Irr(E_1,P)\otimes\Irr(E_2,P)$.
\item
If $\Irr(E_1,P)\oplus\Irr(E_2,P)=
 \Irr(E_1)\cup\Irr(E_2,P)$ is good 
for any $P\in\nbigd$,
then $E_1\oplus E_2$ is an unramifiedly
good lattices of $(\nbige_1\oplus\nbige_2,\DD)$
with $\Irr(E_1\oplus E_2,P)=
 \Irr(E_1,P)\oplus\Irr(E_2,P)$.
\end{itemize}

Let $\nbigx_1$ be a complex manifold
with a normal crossing hypersurface $\nbigd_1$.
Let $F:\nbigx_1\lrarr \nbigx$ be a morphism
such that
(i) $F^{-1}(\nbigd)\subset\nbigd_1$,
(ii) the induced morphism
 $\nbigx_1\lrarr \nbigk$ is a smooth fibration,
(iii) any intersection of some irreducible components
 of $\nbigd_1$ is smooth over $\nbigk$.
Let $E$ be a good lattice of $(\nbige,\DD)$
on $(\nbigx,\nbigd)$.
We have the following functoriality
for the pull back.
\begin{itemize}
\item
 $F^{\ast}E$ is a good lattice of
 $F^{\ast}(\nbige,\DD)
 =F^{-1}(\nbige,\DD)\otimes
 \nbigo_{\nbigx_1}(\ast \nbigd_1)$.
 If $E$ is unramifiedly good,
 $F^{\ast}E$ is also unramifiedly good,
 and we have
\[
 \Irr(F^{\ast}E,P)=\bigl\{
 F^{\ast}\gminia\,\big|\,
 \gminia\in \Irr(E,F(P)) \bigr\}. 
\]
\end{itemize}

\subsection{A criterion for a lattice to be good}

Let $\nbigx$, $\nbigd$ and $(\nbige,\DD)$
be as in Subsection \ref{subsection;10.5.3.3}.
Let $E$ be a lattice of $\nbige$.
Assume that we are given the following:
\begin{itemize}
\item
 A good set of irregular values
 $\nbigi\subset M(\nbigx,\nbigd)\big/H(\nbigx)$.
\item
 A holomorphic decomposition
 $E=\bigoplus_{\gminia\in \nbigi}E_{\gminia}$.
\item
 Let $p_{\gminia}$ denote the projection
 onto $E_{\gminia}$,
 and we put
 $\Phi:=\sum_{\gminia\in \nbigi}
 d\gminia\cdot p_{\gminia}$.
 Then, $\DD^{(0)}:=\DD-\Phi$
 is logarithmic with respect to $E$.
 (Note we do not assume $\DD^{(0)}$ is flat.)
\end{itemize}

We obtain the following proposition
as a corollary of Proposition \ref{prop;10.5.2.21}.
It will be useful in the proof of
Theorem \ref{thm;07.11.18.20}.
\begin{prop}
\label{prop;07.7.7.23}
$E$ is an unramifiedly good lattice 
of $(\nbige,\DD)$.
For any $P\in \nbigd$,
the set $\Irr(\DD,P)$
is equal to the image of
$\nbigi$ via
$M(\nbigx,\nbigd)\big/H(\nbigx)\lrarr
 \nbigo_{\Phat}(\ast \nbigd)\big/\nbigo_{\Phat}$.
\hfill\qed
\end{prop}

\subsection{Family of good lattice is good lattice}

Let $\nbigx$, $\nbigd$, $(\nbige,\DD)$
be as in Subsection \ref{subsection;10.5.3.3}.
For each $y\in \nbigk$,
we set $\nbigx^y:=\{y\}\times_{\nbigk}\nbigx$
and $\nbigd^y:=\{y\}\times_{\nbigk}\nbigx$.
We have the induced 
meromorphic $\varrho(y)$-flat bundle
$(\nbige^{y},\DD^{y})$
on $(\nbigx^y,\nbigd^y)$.
Let $\nbigi\subset M(\nbigx,\nbigd)/H(\nbigx)$
be a good set of irregular values.
The image of $\nbigi$ via
$M(\nbigx,\nbigd)/H(\nbigx)\lrarr
\nbigo_{\Phat}(\ast \nbigd)/\nbigo_{\Phat}$
is denoted by $\nbigi_{\Phat}$.
The image via 
$M(\nbigx,\nbigd)/H(\nbigx)\lrarr 
M(\nbigx^y,\nbigd^y)/H(\nbigx^y)$
is denoted by $\nbigi^y$.
If $P$ is contained in $\nbigx^y$,
let $\nbigo_{\Phat}^y$ be the completion
of the local ring $\nbigo_{\nbigx^y,P}$,
and let $\nbigi^y_{\Phat}$
denote the image of $\nbigi^y$
via $M(\nbigx^y,\nbigd^y)/H(\nbigx^y)\lrarr 
 \nbigo^y_{\Phat}(\ast \nbigd)/\nbigo^y_{\Phat}$.

\begin{prop}
\label{prop;10.5.5.32}
Let $E$ be a lattice of $\nbige$.
The following conditions are equivalent.
\begin{itemize}
\item
$E$ is unramifiedly good,
and $\Irr(\DD,P)=\nbigi_{\Phat}$
for any $P\in \nbigd$.
\item
For each $y\in \nbigk$,
the specialization 
$E^y=E\otimes \nbigo_{\nbigx^y}$
is an unramifiedly good lattice of
$(\nbige^y,\DD^{y})$,
and 
$\Irr(\DD^{y},P)=
 \nbigi^y_{\Phat}$ $(P\in \nbigd_y)$.
\end{itemize}
\end{prop}
\pf
It is easy to see that the first condition
implies the second one.
We would like to show the converse.
We have only to consider the case
$\nbigx=\Delta^n\times \nbigk$
and $\nbigd=\bigcup_{i=1}^{\ell}\{z_i=0\}$.
Let $O=(0,\ldots,0)\in\Delta^n$.
We have only to show that
the completion of $E$ along $O\times \nbigk$
has the irregular decomposition
with the set of irregular values $\nbigi$.

Let $H(\nbigk)$ denote the space of
holomorphic functions on $\nbigk$.
For each $y\in \nbigk$,
we have the specialization
$\val_y:H(\nbigk)\lrarr \cnum$
given by $\val_y(f)=f(y)$.
We put $R_0:=H(\nbigk)[\![z_1,\ldots,z_n]\!]$
and $k_0:=\cnum[\![z_1,\ldots,z_n]\!]$,
and let $R$ (resp. $k$)
be the localization of $R_0$ (resp. $k_0$)
with respect to 
$z_i$ $(i=1,\ldots,\ell)$.
The natural morphism
$\val_y:R_0\lrarr k_0$
induces a functor from
the category of $R_0$-modules 
to the category of $k_0$-modules.
The image of an $R_0$-module $E$
is denoted by $E^{y}$.
We use the symbol $\nbige^y$
for an $R$-module $\nbige$
in a similar meaning.
To show Proposition \ref{prop;10.5.5.32},
we have only to show 
the following lemma.
\begin{lem}
\label{lem;10.5.5.31}
Let $\nbigi\subset R/R_0$ be a good set of
irregular values.
For any $y\in \nbigk$,
let $\nbigi^{y}\subset k/k_0$
denote the specialization of $\nbigi$ at $y$.

Let $(\nbige,\DD)$ be 
a meromorphic $\varrho$-flat bundle over $R$.
Let $E$ be a free $R_0$-lattice of $\nbige$
such that 
for each $y\in \nbigk$
the restriction $E^{y}$ is 
an unramifiedly good lattice of
$(\nbige^y,\DD^{y})$
with $\Irr(\DD^{y})=\nbigi^y$,
i.e., we have a decomposition
$(E^y,\DD^y)=
 \bigoplus_{\gminia\in \nbigi}
 (E^y_{\gminia},\DD^y_{\gminia})$
such that $\DD^y_{\gminia}$ 
are $\gminia^y$-logarithmic.

Then, $E$ is an unramifiedly good lattice of 
$(\nbige,\DD)$
with $\Irr(\DD)=\nbigi$.
\end{lem}
\pf
We use an induction on $\bigl|\nbigt(\nbigi)\bigr|$.
(See Section \ref{subsection;08.1.23.30}
 for $\nbigt(\nbigi)$.)
By considering the tensor product with
a meromorphic $\varrho$-flat line bundle,
we may assume 
$\min\bigl\{\ord\gminia\,\big|\,
 \gminia\in \nbigt(\nbigi)\bigr\}=
 \min\bigl\{\ord\gminia\,\big|\,
 \gminia\in \nbigi\bigr\}$.
Let us take an auxiliary sequence
$\vecm(0),\ldots,\vecm(L)$ for $\nbigi$.
We have $\vecm(1)=\vecm(0)+\vecdelta_j$
for some $j$.
Let $F$ denote the endomorphism of $E_{|z_j=0}$
induced by
$\vecz^{-\vecm(0)}
 \DD(z_j\del_j)$.
The eigenvalues of $F$ are given 
by the set
$T=\bigl\{
 (\vecz^{-\vecm(0)}z_j\del_j\gminia)_{|z_j=0}
 \,\big|\,
 \gminia\in \nbigi \bigr\}$.
Hence, we have the eigen decomposition:
\begin{equation}
 \label{eq;07.10.8.20}
 E_{|z_j=0}=\bigoplus_{\gminib\in T} \EE_{\gminib}
\end{equation}

\begin{lem}
\label{lem;10.5.5.30}
We can take a $\DD$-flat decomposition
$E=\bigoplus_{\gminib\in T} E_{\gminib}$
such that $E_{\gminib|z_j=0}=\EE_{\gminib}$.
\end{lem}
\pf
We give only an outline of the proof.
Let $\vecv$ be a frame of $E$ 
whose restriction to $z_j=0$ is compatible with
(\ref{eq;07.10.8.20}).
We have the decomposition
$\vecv=(\vecv_{\gminib})$
corresponding to the decomposition (\ref{eq;07.10.8.20}).
We have the following:
\[
 \vecz^{-\vecm(0)}
 \DD(z_j\del_j)\vecv
=\vecv\,\left(
\bigoplus_{\gminib\in T} \Omega_{\gminib}
+z_j\, B
 \right)
\]
Here, the entries of $\Omega_{\gminib}$ 
and $B$ are regular,
and $\Omega_{\gminib|O}$ 
$(\gminib\in T)$ have 
no common eigenvalues.
Applying an argument in \cite{levelt}
(or the argument in Subsection 
\ref{subsection;10.5.2.23}),
we can take $\vecv$
for which $B$ is block diagonal,
i.e., $B=\bigoplus_{\gminib\in T} B_{\gminib}$.
Let $E_{\gminib}$ be the $R_0$-submodule
generated by $\vecv_{\gminib}$.
It can be shown that $E_{\gminib}$ is 
$\DD$-flat
using the argument
in the proof of Lemma \ref{lem;10.5.2.20}.
Thus Lemma \ref{lem;10.5.5.30} is proved.
\hfill\qed

\vspace{.1in}

Let us return to the proof of Lemma 
\ref{lem;10.5.5.31}.
For $\gminib\in T$,
let $\nbigi(\gminib)$ denote 
the inverse image of $\gminib$
by the natural map $\nbigi\lrarr T$.
Its specialization at $y$ is denoted by
$\nbigi(\gminib)^y$.
We can deduce 
$E_{\gminib}^y
=\bigoplus_{\gminia\in\nbigi(\gminib)^y}
 E_{\gminia}^y$.
Then, we can apply the hypothesis of the induction
on each $(E_{\gminib},\DD_{\gminib})$.
Therefore,
we obtain Lemma \ref{lem;10.5.5.31}
and thus Proposition \ref{prop;10.5.5.32}.
\hfill\qed

\section{Decompositions}
\label{subsection;10.5.4.102}

\subsection{Openness property}
\label{subsection;10.5.5.40}

Let $\nbigx$ and $\nbigd$ be as in 
Section \ref{subsection;10.5.3.3}.
Let $(\nbige,\DD)$ be a meromorphic 
flat bundle on $(\nbigx,\nbigd)$.
Let $E$ be a lattice of $\nbige$.
Assume that it
is unramifiedly good at a point $P\in \nbigd$,
i.e.,
there exist a good set of irregular values
$\Irr(\DD,P)\subset
 \nbigo_{\Phat}(\ast \nbigd)/\nbigo_{\Phat}$
and a decomposition
\[
 (E,\DD)_{|\Phat}
=\bigoplus_{\gminia\in\Irr(\DD,P)}
 (\lefttop{P}\Ehat_{\gminia},
 \lefttop{P}\DDhat_{\gminia})
\]
such that 
$\lefttop{P}\DDhat_{\gminia}
-d\gminia$ are logarithmic
with respect to $\Ehat_{\gminia}$.
For a small neighbourhood $\nbigx_P$ of $P$,
let $\nbigd_P:=\nbigd\cap \nbigx_P$.
We will prove the following proposition
in Subsection \ref{subsection;10.5.3.31}.

\begin{prop}
\label{prop;10.5.3.6}
If $\nbigx_P$ is sufficiently small,
the following claims hold:
\begin{itemize}
\item
$\Irr(\DD,P)\subset M(\nbigx_P,\nbigd_P)/H(\nbigx_P)$,
i.e.,
it is contained in the image of
$M(\nbigx_P,\nbigd_P)/H(\nbigx_P)\lrarr
 \nbigo_{\Phat}(\ast \nbigd_P)/\nbigo_{\Phat}$.
\item $E$ is unramifiedly good at any point of
 $P'\in \nbigd_P$.
\item
 The good set of irregular values
 $\Irr(\DD,P')$ of 
$(E,\DD)_{|\Phat'}$
 is the image of 
 $\Irr(\DD,P)$ by
 $M(\nbigx_P,\nbigd_P)\big/H(\nbigx)\lrarr
 \nbigo_{\Phat'}(\ast \nbigd)\big/\nbigo_{\Phat'}$.
\end{itemize}
\end{prop}

\subsubsection{Good system of irregular values}

Before proceeding,
we give a consequence
and prepare a terminology.
When we are given an unramifiedly good lattice
$(E,\DD)$ of a meromorphic $\varrho$-flat bundle,
we put 
$\Irr(\DD):=\bigl\{
 \Irr(\DD,P)\,\big|\,P\in\nbigd
 \bigr\}$.
Then, Proposition \ref{prop;10.5.3.6}
says that $\Irr(\DD)$
is a good system of irregular values
in the following sense.
\index{system $\Irr(\DD)$}

\begin{df}
A system $\vecnbigi$ of
finite subsets
$\nbigi_P\subset\nbigo_{\nbigx}(\ast \nbigd)_P\big/
 \nbigo_{\nbigx,P}$ $(P\in \nbigd)$ is called
a good system of irregular values,
if the following holds for each $P\in \nbigd$:
\begin{itemize}
\item
Take a neighbourhood $\nbigx_P$
of $P$
such that $\nbigi_P\subset 
 M(\nbigx_P,\nbigd_P)/H(\nbigx_P)$,
where $\nbigd_P=\nbigx_P\cap \nbigd$.
Then, for each $P'\in \nbigd_P$,
$\nbigi_{P'}$ is the image of
$\nbigi_P$ via
$M(\nbigx_P,\nbigd_P)/H(\nbigx_P)\lrarr
 \nbigo_{\nbigx}(\ast \nbigd)_{P'}\big/
 \nbigo_{\nbigx,P'}$.
\hfill\qed
\end{itemize}
\index{good system of irregular values}
\end{df}

\begin{rem}
A good set of irregular values
$\nbigi\subset\nbigo_{\nbigx}(\ast \nbigd)_P/
 \nbigo_{\nbigx,P}$
naturally induces a good system of irregular 
values on a neighbourhood of $P$.
In that case, we will not distinguish
the induced system and $\nbigi$.
\hfill\qed
\end{rem}

\subsection{Decompositions
along the intersection of irreducible components}
\label{subsection;10.5.3.7}
We will also prove a refinement 
of the second and third claims
of Proposition \ref{prop;10.5.3.6}.
For simplicity,
let us consider the case
$\nbigx=\Delta^n\times \nbigk$
and $\nbigd=\bigcup_{i=1}^{\ell}\{z_i=0\}$.
We put $\nbigd_i:=\{z_i=0\}$.
For a subset $I\subset\ellsitabar$,
we put $\nbigd_I:=\bigcap_{i\in I}\nbigd_i$
and $\nbigd(I):=\bigcup_{i\in I}\nbigd_i$.
The complement $\ellsitabar\setminus I$
is denoted by $I^c$.
Let $\nbigdhat_I$ and $\nbigdhat(I)$
denote the completion of $\nbigx$
along $\nbigd_I$ and $\nbigd(I)$,
respectively.
(See \cite{banica},
 \cite{bingener} and \cite{krasnov}.
 See also a brief review in
 Subsection \ref{subsection;07.11.6.10}.)
\index{completion $\nbigdhat_I$}
\index{completion $\nbigdhat(I)$}

We may assume $P\in \nbigd_{\ellsitabar}$.
For a given small neighbourhood 
$\nbigx_P$ of $P$ in $\nbigx$,
we put 
$\nbigd_{I,P}:=\nbigd_I\cap \nbigx_P$ and 
$\nbigd(I)_P:=\nbigd(I)\cap \nbigx_P$.

Let $(E,\DD)$ be unramifiedly good at $P$.
Once we know
$\Irr(\DD,P)$ is contained in
$M(\nbigx_P,\nbigd_P)/H(\nbigx_P)$,
let $\Irr(\DD,I)$ denote the image of
$\Irr(\DD,P)$
via $M(\nbigx_P,\nbigd_P)/H(\nbigx_P)\lrarr
 M(\nbigx_P,\nbigd(I^c)_P)$.
We will prove the following proposition
in Subsection \ref{subsection;10.5.3.31}.
\begin{prop}
\label{prop;10.5.3.5}
Let $\nbigx_P$ be a sufficiently
small neighbourhood of $P$
in $\nbigx$.
For any subset $I\subset\ellsitabar$,
we have a decomposition
\begin{equation}
 \label{eq;10.5.3.10}
 (E,\DD)_{|\nbigdhat_{I,P}}
=\bigoplus_{\gminib\in \Irr(\DD,I)}
 \bigl(\lefttop{I}\Ehat_{\gminib},
 \lefttop{I}\DDhat_{\gminib}\bigr)
\end{equation}
such that 
 $(\lefttop{I}\DDhat_{\gminib}-d\gminib)
 \bigl(\lefttop{I}\Ehat_{\gminib}\bigr)
\subset
 \lefttop{I}\Ehat_{\gminib}
 \otimes \Bigl(
 \Omega^1_{\nbigx/\nbigk}\bigl(\log \nbigd(I)\bigr)
+\Omega^1_{\nbigx/\nbigk}\bigl(\ast \nbigd(I^c)\bigr)
 \Bigr)_{|\nbigx_P}$,
where we take a lift of $\gminib$
to $M(\nbigx_P,\nbigd_P)$.
\end{prop}
The decomposition (\ref{eq;10.5.3.10})
is called the irregular decomposition
of $(E,\DD)_{|\nbigdhat_{I,P}}$.
It induces the irregular decomposition
at any point $P'\in \nbigd_{I,P}\setminus
 \bigcup_{j\not\in I}\nbigd_j$.
In that sense, Proposition \ref{prop;10.5.3.5}
refines the second and third claims
in Proposition \ref{prop;10.5.3.6}.

\begin{rem}
The property of Proposition {\rm\ref{prop;10.5.3.5}}
was adopted as definition of
``unramifiedly good at $P$''
in the older version of this monograph.
\hfill\qed
\end{rem}

\begin{rem}
\label{rem;10.5.3.30}
For $\gminib\in\Irr(\DD,I)$,
let $\Irr(\DD,P,\gminib)$
be the inverse image of 
$\gminib$ via the natural map
$\Irr(\DD,P)\lrarr
 \Irr(\DD,I)$.
If we are given the decomposition 
{\rm(\ref{eq;10.5.3.10})},
we have 
$\lefttop{I}\Ehat_{\gminib|\Phat}
=\bigoplus_{\gminia\in
 \Irr(\DD,P,\gminib)}
 \lefttop{P}\Ehat_{\gminia}$.
Hence, it is easy to deduce
\[
 \lefttop{I}\Ehat_{\gminic|\nbigdhat_{J,P}}
=\bigoplus_{\gminib\in
 \Irr(\DD,J,\gminic)}
 \lefttop{J}\Ehat_{\gminib} 
\]
for $I\subset J$
and $\gminic\in\Irr(\DD,I)$,
where $\Irr(\DD,J,\gminic)$ is 
the inverse image of
$\gminic$ via the natural map
$\Irr(\DD,J)\lrarr
 \Irr(\DD,I)$.
\hfill\qed
\end{rem}

\subsection{Decomposition
along the union of irreducible components}
\label{subsection;10.5.3.40}

We continue to use the setting in
Subsection \ref{subsection;10.5.3.7}.
By shrinking $\nbigx$,
we assume $\nbigx=\nbigx_P$.
We will often need a decomposition
on the completion along $\nbigd(I)$
for some $I\subset\ellsitabar$.
For simplicity,
let us take an admissible coordinate system
(Remark \ref{rem;07.11.15.10})
for the good set $\Irr(\DD,P)$,
and we consider decompositions
along $\nbigd(\jbar)$ for $1\leq j\leq \ell$,
where $\jbar:=\{1,\ldots,j\}$.

Let $\Irr(\DD,j)$ 
and $\Irr'(\DD,j)$ denote
the images of $\Irr(\DD,P)$
via the following natural maps:
\[
 M(\nbigx,\nbigd)/H(\nbigx)\lrarr M(\nbigx,\nbigd)\big/
 M\bigl(\nbigx,\nbigd(\ellsitabar\setminus\{j\})\bigr)
\]
\[
 M(\nbigx,\nbigd)/H(\nbigx)\lrarr 
 M(\nbigx,\nbigd)/M(\nbigx,\nbigd(\underline{j-1}))
\]
Note that the natural map
$\Irr'(\DD,j)\lrarr \Irr(\DD,j)$
is bijective by our choice of the coordinate system,
via which we identify them.
We have the naturally defined maps
$\Irr'(\DD,i)\lrarr\Irr'(\DD,j)$
for any $i\leq j$,
which induces
$\pi_{j,i}:
 \Irr(\DD,i)\lrarr\Irr(\DD,j)$.

\begin{lem}
\label{lem;10.5.3.25}
There is the following decomposition:
\begin{equation}
 \label{eq;07.10.8.3}
 (E,\DD)_{|\nbigdhat(\jbar)}
=\!\!\!\!\!\bigoplus_{\gminib\in \Irr(\DD,j)}
 \!\!\!\!\!
 \bigl(\Ehat_{\gminib,\nbigdhat(\jbar)},
 \DD_{\gminib}\bigr)
\quad
\mbox{\rm such that }\,\,
 \Ehat_{\gminib,\nbigdhat(\jbar)|\nbigdhat_{i}}
=\!\!\!\bigoplus_{
 \substack{\gminic\in \Irr(\DD,i)\\
 \pi_{j,i}(\gminic)=\gminib}} 
 \!\!\!
 \lefttop{i}\Ehat_{\gminic}
\end{equation}
\end{lem}
\pf
For $J\subset\ellsitabar$,
let $i(J)$ be the number determined by
$\underline{i(J)}\cap J=\emptyset$
and $i(J)+1\in J$.
We have the maps
\[
 M(\nbigx,\nbigd)\big/H(\nbigx)
 \stackrel{p_1}{\lrarr}
 M(\nbigx,\nbigd)\big/M(\nbigx,\nbigd(\underline{i(J)}))
 \stackrel{p_2}{\lrarr}
 M(\nbigx,\nbigd)\big/M(\nbigx,\nbigd(J^c))
\]
Let $\Irr'(\DD,J)$
denote the image of
$\Irr(\DD,P)$ by $p_1$.
Then, $p_2$ naturally gives a bijection
$\Irr'(\DD,J)\lrarr
 \Irr(\DD,J)$,
by which we identify them.
If $j\in J$,
we have the naturally defined map
$\Irr'(\DD,J)\lrarr
 \Irr(\DD,j)$,
which induces
$\pi_{j,J}:\Irr(\DD,J)\lrarr
 \Irr(\DD,j)$.
For $\gminib\in\Irr(\DD,j)$,
we put
$\lefttop{J}\Ehat_{\gminib}:=
 \bigoplus_{\gminia\in\pi_{j,J}^{-1}(\gminib)}
 \lefttop{J}\Ehat_{\gminia}$.
By Remark \ref{rem;10.5.3.30},
we have 
$\lefttop{I}\Ehat_{\gminib|\nbigdhat_J}
=\lefttop{J}\Ehat_{\gminib}$
for $I\subset J$ and $\gminib\in\Irr(\DD,j)$.
Then, we obtain the decomposition 
(\ref{eq;07.10.8.3})
by using a general lemma
(Lemma \ref{lem;10.5.3.21} below).
\hfill\qed

\subsection{Decomposition
in the level $\vecm(i)$}
\label{subsection;10.5.5.36}

We use the setting in
Subsection \ref{subsection;10.5.3.40}.
Take an auxiliary sequence 
$\vecm(0),\ldots,\vecm(L),\vecm(L+1)=\veczero$
for $\Irr(\DD)$.
(See Section \ref{subsection;08.1.23.30}.)
Let $\Irrbar(\DD,\vecm(i))$ denote the image of
$\Irr(\DD)$ via $\etabar_{\vecm(i)}$.
Let $k(i)$ denote the number determined by
$\vecm(i)\in\seisuu_{<0}^{k(i)}
 \times\veczero_{\ell-k(i)}$.
Let $\kbar(i):=\underline{k(i)}$.
We remark that
\[
\pi_j:
 \Irrbar(\DD,\vecm(i))\lrarr
 M(\nbigx,\nbigd)/M(\nbigx,\nbigd(\underline{j-1})) 
\]
is injective for $j\leq k(i)$.
We also have the map
\[
 \etabar_{\vecm(i),j}:
\Irr(\DD,j)\lrarr 
 M(\nbigx,\nbigd)/M(\nbigx,\nbigd(\underline{j-1})) 
\]
given as follows:
\begin{multline*}
 \Irr(\DD,j)
\simeq \Irr'(\DD,j)
\subset M(\nbigx,\nbigd)\big/
 M(\nbigx,\nbigd(\underline{j-1})) \\
\stackrel{b}{\lrarr} 
 M(\nbigx,\nbigd)\big/
 M(\nbigx,\nbigd(\underline{j-1}))
\end{multline*}
Here, $b$ is 
induced by $\etabar_{\vecm(i)}$.
As in Lemma \ref{lem;10.5.3.25},
we obtain the following decomposition:
\begin{equation}
\label{eq;07.11.12.10}
(E,\DD)_{|\nbigdhat(\kbar(i))}
=\!\!\!\!\!\!\bigoplus_{\gminib\in \Irrbar(\DD,\vecm(i))}
 \!\!\!\!\!\!\!\!\!
 \bigl(\Ehat^{\vecm(i)}_{\gminib},
 \DD_{\gminib}\bigr),
\,\,\,
\mbox{\rm where}\,\,
 \Ehat^{\vecm(i)}_{\gminib|\nbigdhat_{j}}=
 \!\!\!\!\!\!\!\!\!\!
 \bigoplus_{\substack{
 \gminic\in \Irr(\DD,j)\\
 \etabar_{\vecm(i),j}(\gminic)=
 \pi_{j}(\gminib)\\
 }}
 \!\!\!\!\!\!\!\!\!\!\!\!\!
\lefttop{j}\Ehat_{\gminic},
\,\,\, \bigl(j\leq k(i)\bigr)
\end{equation}
The decomposition (\ref{eq;07.11.12.10})
is called the irregular decomposition
in the level $\vecm(i)$.

\index{irregular decomposition
 in the level $\vecm(i)$}

\subsection{Zero of $\varrho$}

We use the setting in Section
\ref{subsection;10.5.5.40}.
It is easy to show the following lemma.

\begin{lem}
\label{lem;10.5.5.33}
Assume $\varrho$ is constantly $0$.
Then, for 
a sufficiently small neighbourhood $\nbigx_P$ of $P$,
we have 
$\Irr(\DD,P)\subset
 M(\nbigx_P,\nbigd_P)/H(\nbigx_P)$
and a decomposition
$(E,\DD)_{|\nbigx_P}
=\bigoplus_{\gminia\in\Irr(\DD,P)}
 (E_{\gminia},\DD_{\gminia})$
such that
$\DD_{\gminia}$ are $\gminia$-logarithmic.
\hfill\qed
\end{lem}

Let us consider the case in which
$\varrho$  is not constantly $0$.
For simplicity,
we assume that $d\varrho$ 
is nowhere vanishing on 
$\nbigk^0:=\varrho^{-1}(0)$.
We put $\nbigx^{0}:=\nbigx\times_{\nbigk}\nbigk^{0}$
and $\nbigd^{0}:=\nbigd\times_{\nbigk}\nbigk^{0}$.
We also assume that $P\in \nbigd^{0}$.
As remarked in Lemma \ref{lem;10.5.5.33},
by shrinking $\nbigx$ around $P$,
we have a decomposition
\begin{equation}
 \label{eq;10.5.5.35}
  (E,\DD)_{|\nbigx^{0}}
=\bigoplus_{\gminia\in\Irr(\DD,P)}
 (E_{\gminia,\nbigx^{0}},
 \DD_{\gminia})
\end{equation}
such that
$\DD_{\gminia}-d\gminia_{|\nbigx^{0}}$
are logarithmic.
Let $\nbigxhat^{0}$ be the completion of
$\nbigx$ along $\nbigx^{0}$.
The following lemma can be shown
by a standard argument.
\begin{lem}
We have a flat decomposition
\[
 (E,\DD)_{|\nbigxhat^{0}}
=\bigoplus_{\gminia\in\nbigi}
 (E_{\gminia,\nbigxhat^{0}},
 \DD_{\gminia})
\]
such that
(i) its restriction to $\nbigx^{0}$
is the same as {\rm(\ref{eq;10.5.5.35})},
(ii) its restriction to $\Phat$
is the same as the irregular decomposition
of $(E,\DD)_{|\Phat}$.
\hfill\qed
\end{lem}

We have refinements of Sections 
\ref{subsection;10.5.3.40}
and \ref{subsection;10.5.5.36}.
We use the setting there.
For $\gminib\in\Irr(\DD,j)$,
let $\Irr(\DD,P,\gminib)$
denote the inverse image of $\gminib$
by the natural map
$\Irr(\DD,P)\lrarr\Irr(\DD,j)$.
We put
$\Ehat_{\gminib,\nbigxhat^0}:=
 \bigoplus_{\gminia\in\Irr(\DD,P,\gminib)}
 \Ehat_{\gminia,\nbigxhat^0}$.
We put $W(\jbar):=\nbigd(\jbar)\cup \nbigx^0$.
As in Lemma \ref{lem;10.5.3.25},
we obtain the decomposition
\begin{equation}
 \label{eq;10.5.5.40}
 (E,\DD)_{|\What(\jbar)}
=\bigoplus_{\gminib\in\Irr(\DD,j)}
 \bigl(\Ehat^{\vecm(p)}_{\gminib,\What(\jbar)},
 \DD_{\gminib}\bigr)
\end{equation}
such that
(i) $\Ehat_{\gminib,\What(\jbar)|\nbigdhat(\jbar)}
=\Ehat_{\gminia,\nbigdhat(\jbar)}$,
(ii) $\Ehat_{\gminib,\What(\jbar)|\nbigxhat^0}
=\Ehat_{\gminia,\nbigxhat^0}$.
Similarly,
for $\gminib\in\Irr(\DD,\vecm(i))$,
let $\Irr(\DD,P,\gminib)$
be the inverse image of $\gminib$
by the natural map
$\Irr(\DD,P)\lrarr
\Irr(\DD,\vecm(i))$.
Then, we have the decomposition
\begin{equation}
 \label{eq;10.5.5.41}
 (E,\DD)_{|\What(\kbar(i))}
=\bigoplus_{\gminib\in\Irr(\DD,\vecm(i))}
 \bigl(\Ehat^{\vecm(i)}_{\gminib,\What(\kbar(i))},
 \DD_{\gminib}\bigr)
\end{equation}
such that 
(i) $\Ehat^{\vecm(i)}_{\gminib,
 \What(\kbar(i))|\nbigdhat(\kbar(i))}
=\Ehat^{\vecm(i)}_{\gminib,\nbigdhat(\kbar(i))}$,
(ii)$\Ehat^{\vecm(i)}_{\gminib,
 \What(\kbar(i))|\nbigxhat^0}
=\bigoplus_{\etabar_{\vecm(i)}(\gminia)=\gminib}
 \Ehat_{\gminia,\nbigxhat^0}$.

\subsection{Proof of
 Proposition \ref{prop;10.5.3.6} and
 Proposition \ref{prop;10.5.3.5}}
\label{subsection;10.5.3.31}

We have only to show the propositions
under the setting of 
Subsection \ref{subsection;10.5.3.7}.
In the following,
instead of considering a neighbourhood $\nbigx_P$,
we will replace $\nbigx$
by a small neighbourhood of $P$
without mention,
if it is necessary.

\subsubsection{Step 1}

We fix $I\subset\ellsitabar$ for a moment.
Let $E$ be a free $\nbigo_{\nbigdhat_I}$-module
with a meromorphic flat connection
$\DD:E\lrarr 
 E\otimes\Omega^{1}_{\nbigdhat_I/\nbigk}(\ast \nbigd)$.
Assume that we are given the following:
\begin{itemize}
\item
$\vecm\in\seisuu_{\leq 0}^{\ell}$
and $i\in I$ such that $m_i<0$.
We set $\vecm':=\vecm+\vecdelta_i$.
\item
$\nbigi\subset \nbigo_{\Phat}(\ast \nbigd)$
such that,  for any $\gminia\in\nbigi$,
(i) $z_i^{-m_i}\gminia$ is independent
of the variable $z_i$,
(ii) $\vecz^{-\vecm}\,\gminia
 \in\nbigo_{\Phat}$.
\item
A decomposition
$E_{|\Phat}=
 \bigoplus_{\gminia\in\nbigi}
 \lefttop{P}E_{\gminia}$
such that
$ \vecz^{-\vecm'}(\DD-d\gminia)
 \bigl(\lefttop{P}E_{\gminia}\bigr)
\subset
 \lefttop{P}E_{\gminia}
\otimes
 \Omega^1_{\Phat/\nbigk}(\log \nbigd)$.
\end{itemize}

We set $\nbigi_0:=
 \bigl\{
 (\vecz^{-\vecm}\, \gminia)(P)\,
 \big|\,
 \gminia\in\nbigi
 \bigr\}\subset \cnum$.
We have a naturally defined map
$\pi:\nbigi\lrarr\nbigi_0$.
We set $\lefttop{P}E_{\gminib}
 :=\bigoplus_{\pi(\gminia)=\gminib}
 \lefttop{P}E_{\gminia}$.

Let $H(\nbigd_I)$ denote the space of
holomorphic functions on $\nbigd_I$.
Let $R$ denote the localization of
$H(\nbigd_I)[\![z_i\,|\,i\in I]\!]$
with respect to $\prod_{i=1}^{\ell} z_i$.

\begin{lem}
\label{lem;10.5.3.20}
$\nbigi$ is contained in $R$,
and we have a flat decomposition
$E=\bigoplus_{\gminib\in \nbigi_0}
 E_{\gminib}$
such that
$E_{\gminib|\Phat}=
 \lefttop{P}E_{\gminib}$.
\end{lem}
\pf
First, we remark that
$\vecz^{-\vecm}\DD(z_i\del_i)
 \lefttop{P}E_{\gminia}
 \subset\lefttop{P}E_{\gminia}$,
and thus
$\vecz^{-\vecm}\DD(z_i\del_i)
 E\subset E$.
Let $F$ be the endomorphism
of $E_{|\nbigdhat_I\cap \nbigd_i}$ induced by
$\vecz^{-\vecm}
 \DD(z_i\del_{i})$.
The eigen decomposition
of $F_{|P}$ is given by
$E_{|P}=\bigoplus_{\gminib\in\nbigi_0}
 \lefttop{P}E_{\gminib|P}$.
We obtain the decomposition
$E_{|\nbigdhat_I\cap \nbigd_i}
=\bigoplus_{\gminib\in\nbigi_0}
 G_{\gminib}$
such that
(i) $F(G_{\gminib})\subset G_{\gminib}$
(ii) $G_{\gminib|P}=\lefttop{P}E_{\gminib|P}$.
By comparing $F$ and its completion at $P$,
we obtain that $\nbigi\subset R$.
By using a standard argument
(see Section \ref{section;10.5.3.11}),
we obtain the decomposition
$E=\bigoplus_{\gminib\in \nbigi_0}
 E_{\gminib}$ such that
 (i) $E_{\gminib|\nbigdhat_I\cap \nbigd_i}
=G_{\gminib}$,
(ii) it is preserved by
$\vecz^{-\vecm}\DD$.
By a standard argument
as in Corollary \ref{cor;10.5.2.4},
we can show that
$E_{\gminib|\Phat}=\lefttop{P}E_{\gminib}$.
\hfill\qed

\subsubsection{Step 2}

For $1\leq p\leq \ell$,
we put $\pbar:=\{1,\ldots,p\}$.
Let $E$ be a free $\nbigo_{\nbigdhat(\pbar)}$-module
with a meromorphic flat connection
$\DD:E\lrarr 
 E\otimes\Omega^1_{\nbigdhat(\pbar)/\nbigk}
 (\ast \nbigd)$.
Assume that we are given 
a good set of irregular values
$\Irr(\DD)\subset
 \nbigo_{\Phat}(\ast \nbigd(\pbar))/\nbigo_{\Phat}$
and a decomposition
\[
 (E,\DD)_{|\Phat}
=\bigoplus_{\gminia\in\Irr(\DD)}
 (\lefttop{P}E_{\gminia},
 \lefttop{P}\DD_{\gminia})
\]
such that 
$\lefttop{P}\DD_{\gminia}$
are $\gminia$-logarithmic.
For $I\subset \pbar$,
let $\Irr(\DD,I)$
denote the image of $\Irr(\DD)$
via the natural map
$p_I:
 \nbigo_{\Phat}(\ast \nbigd(\pbar))/
 \nbigo_{\Phat}
\lrarr
 \nbigo_{\Phat}(\ast \nbigd(\pbar))/
 \nbigo_{\Phat}(\ast\nbigd(I_1))$,
where $I_1:=\pbar\setminus I$.
For each $I$ and $\gminib\in \Irr(\DD,I)$,
we set
\[
 \lefttop{P}E_{\gminib}:=
 \bigoplus_{\substack{
 \gminia\in\Irr(\DD)\\
 p_I(\gminia)=\gminib }}
\lefttop{P}E_{\gminia}
\]

\begin{lem}
\label{lem;10.5.3.12}
If we shrink $\nbigx$ appropriately,
$\Irr(\DD)$ is contained in
the image of
$M(\nbigx,\nbigd(\pbar))/H(\nbigx)
\lrarr
 \nbigo_{\Phat}(\ast\nbigd(\pbar))
 \big/\nbigo_{\Phat}$.
For each $I\subset \pbar$,
we have a decomposition
$ E_{|\nbigdhat_I}=
 \bigoplus_{\gminib\in\Irr(\DD,I)}
 \lefttop{I}E_{\gminib}$
such that
$\lefttop{I}E_{\gminib|\Phat}
=\lefttop{P}E_{\gminib}$.
\end{lem}
\pf
We use an induction on 
the rank of $E$.
Assume that the coordinate system is 
admissible for $\Irr(\DD)$.
Take an auxiliary sequence
$\vecm(0),\ldots,\vecm(L)$
for $\Irr(\DD)$.
We put $T:=\bigl\{
 (\vecz^{-\vecm(0)}\gminia)(P)\,\big|\,
 \gminia\in\Irr(\nabla)
 \bigr\}$.
We have the naturally defined map
$q:\Irr(\nabla)\lrarr T$.
For each $\alpha\in T$,
we put 
$\lefttop{P}E_{\alpha}
=\bigoplus_{q(\gminia)=\alpha}
 \lefttop{P}E_{\gminia}$.
Then, $E_{|\Phat}=\bigoplus
 \lefttop{P}E_{\alpha}$ is a flat decomposition.
It is easy to observe that
if we are given a flat decomposition
$E_{|\Phat}=\bigoplus_{\alpha\in T}
 \lefttop{P}E'_{\alpha}$ such that
$\lefttop{P}E'_{\alpha|P}=
 \lefttop{P}E_{\alpha|P}$,
then we have
$\lefttop{P}E'_{\alpha}=
 \lefttop{P}E_{\alpha}$.

Due to Lemma \ref{lem;10.5.3.20}
with $I=\{\gminih(0)\}$,
$\etabar_{\vecm(0)}(\gminia)$
are meromorphic functions
for any $\gminia\in \Irr(\DD)$.
Hence, 
by considering the tensor product
with a meromorphic flat line bundle,
we have only to consider the case
in which $|T|\geq 2$.
Let $k$ be determined by
$\vecm(0)\in\seisuu_{<0}^k\times
 \veczero_{\ell-k}$.

Take $I\subset\pbar$.
If $I\cap\kbar=\emptyset$,
the trivial decomposition is desired one.
Let us consider the case $I\cap\kbar\neq\emptyset$.
We take $i\in I\cap\kbar$.
Let $m_i(0)$ be the $i$-th component of $\vecm(0)$.
We have
$\vecz^{-\vecm(0)}\nabla(z_i\del_i)E\subset E$.
Let $G$ be the induced endomorphism
of $E_{|P}$.
It is semisimple,
the eigenvalues are given by
$\bigl\{
 m_i(0)\alpha\,\big|\,\alpha\in T
 \bigr\}$,
and the eigen decomposition
is given by
$E_{|P}=\bigoplus\lefttop{P}E_{\alpha|P}$,
where $G_{|E_{\alpha|P}}$ is 
the multiplication of $m_i(0)\alpha$.
By applying Lemma \ref{lem;10.5.3.20},
we can extend it
to a flat decomposition of $E_{|\Dhat_I}$,
i.e.,
\[
 E_{|\Dhat_I}=
 \bigoplus_{\alpha\in T}
 \lefttop{I}E_{\alpha},
\quad\mbox{\rm such that\,\,}
\quad
 \lefttop{I}E_{\alpha|P}
=\lefttop{P}E_{\alpha|P}
\]
Then, we obtain
$\lefttop{I}E_{\alpha|\Phat}
=\lefttop{P}E_{\alpha}$.
For $I\subset J\subset\kbar$ as above,
we have $\lefttop{I}E_{\alpha|\Dhat_J}
=\lefttop{J}E_{\alpha}$.
Due to Lemma \ref{lem;10.5.3.21} below,
we obtain the flat decomposition
\[
 (E,\nabla)_{|\Dhat(\kbar)}
=\bigoplus_{\alpha\in T}
 (E_{\alpha},\nabla_{\alpha}),
\quad\mbox{\rm such that }
E_{\alpha|\Dhat_I}
=\lefttop{I}E_{\alpha}
\]
We may apply the hypothesis of the induction
to $(E_{\alpha},\nabla_{\alpha})$ 
on $D(\kbar)$,
and we obtain
Lemma \ref{lem;10.5.3.12}.
\hfill\qed

\subsubsection{Step 3}

We can complete the proof of
Proposition \ref{prop;10.5.3.6} and
Proposition \ref{prop;10.5.3.5}
by applying Lemma \ref{lem;10.5.3.12}
to $(E,\DD)_{|\nbigdhat}$.

\subsubsection{General lemma}

Recall the following general lemma.
\begin{lem}
\label{lem;10.5.3.21}
Let $\Vhat$ be a free $\nbigo_{\nbigdhat}$-module
on $\nbigx$.
Assume that we are given a decomposition
$\Vhat_{|\nbigdhat_I}=\bigoplus
 \lefttop{I}\Vhat_{\gminia}$ 
for each $I\subset\ellsitabar$,
such that
$\lefttop{I}\Vhat_{\gminia|\nbigdhat_J}
=\lefttop{J}\Vhat_{\gminia}$
for any $I\subset J$.
Then, we have a unique decomposition
$\Vhat=\bigoplus\Vhat_{\gminia}$
on $\nbigdhat$,
which induces the decompositions
on $\nbigdhat_I$.
\end{lem}
\pf
Let $\lefttop{I}\pi_{\gminia}$ be 
the projection of $\Vhat_{|\nbigdhat_I}$
onto $\lefttop{I}V_{\gminia}$.
Then, we have
$\lefttop{I}\pi_{\gminia|\nbigdhat_J}
=\lefttop{J}\pi_{\gminia}$.
Let $\vecv$ be a frame of $V$.
Let $\lefttop{I}\Pi_{\gminia}\in M_r(\nbigo_{\nbigdhat_I})$
be determined by
$\lefttop{I}\pi_{\gminia}(\vecv)
=\vecv\cdot \lefttop{I}\Pi_{\gminia}$,
where $r=\rank(V)$.
Because
$\lefttop{I}\Pi_{\gminia|\nbigdhat_J}
=\lefttop{J}\Pi_{\gminia}$,
we have $\Pi_{\gminia}\in M_r(\nbigo_{\nbigdhat})$
such that 
$\Pi_{\gminia|\nbigdhat_I}
=\lefttop{I}\Pi_{\gminia}$.
(Use the exact sequence
in the proof of Proposition 4.1 \cite{har2}, for example.)
Let $\pi_{\gminia}$ be the endomorphism
of $\Vhat$ given by 
$\pi_{\gminia}(\vecv_{|\nbigdhat})=
 \vecv_{|\nbigdhat}\cdot\Pi_{\gminia}$,
and let $\Vhat_{\gminia}$ be the image of
$\pi_{\gminia}$.
Then, $\Vhat=\bigoplus \Vhat_{\gminia}$
gives the desired decomposition.
\hfill\qed

\section{Good filtered $\varrho$-flat bundle}
\label{section;10.5.26.50}

\subsection{Good filtered $\varrho$-flat bundle}

Let $\nbigx\lrarr \nbigk$, $\nbigd$ and $\varrho$
be as in Subsection \ref{subsection;10.5.3.3}.
Let $\nbigd=\bigcup_{i\in \Lambda}\nbigd_i$
be the irreducible decomposition.
Recall that 
a filtered $\varrho$-flat sheaf on $(\nbigx,\nbigd)$ is
defined to be a filtered sheaf 
$\vecE_{\ast}=\bigl(
 \prolongg{\veca}{E}\,\big|\,
 \veca\in\real^{\Lambda}
 \bigr)$ on $(\nbigx,\nbigd)$
with a meromorphic flat $\varrho$-connection
$\DD$ of 
the $\nbigo_\nbigx(\ast \nbigd)$-module
$\vecE=
 \bigcup_{\veca\in\real^{\Lambda}}
 \prolongg{\veca}{E}$.
If $\vecE_{\ast}$ is a filtered bundle,
$(\vecE_{\ast},\DD)$ 
is called a filtered $\varrho$-flat bundle.
See Subsection {\rm\ref{subsection;10.5.28.1}} below
for a brief account on filtered sheaf and filtered bundle.
We shall use some notation and terminology 
given there.
\index{filtered $\varrho$-flat sheaf}
\index{filtered $\varrho$-flat bundle}

\begin{df}
\mbox{{}}\label{df;10.5.5.10}
Let $(\vecE_{\ast},\DD)$ be 
a filtered $\varrho$-flat bundle.
\begin{itemize}
\item
$(\vecE_{\ast},\DD)$
is called unramifiedly good,
if $\prolongg{\vecc}{E}$ are 
unramifiedly good lattices
for any $\vecc\in\real^{\Lambda}$.
\item
$(\vecE_{\ast},\DD)$ is 
called good at $P\in \nbigd$,
if there exists a ramified covering
$\varphi_e:(\nbigx'_P,\nbigd'_P)\lrarr (\nbigx_P,\nbigd_P)$
such that
$(\vecEtilde_{\ast},
 \varphi^{\ast}_e\DD)$
is unramifiedly good.
Here, 
$(\nbigx_P,\nbigd_P)$,
$(\nbigx_P',\nbigd_P')$ and $\varphi_P$
are as in Definition {\rm\ref{df;10.5.3.2}},
and $\vecEtilde_{\ast}$
is induced by $\varphi_P$ and 
$\vecE_{\ast}$ as 
in Section {\rm\ref{subsection;07.11.5.60}}
below.
\item
$(\vecE_{\ast},\DD)$ is called good,
if it is good at any point of $\nbigd$.
In other words,
$(\vecE_{\ast},\DD)$ is good,
if it is the descent of an unramifiedly good 
filtered $\varrho$-flat bundle
around any point of $\nbigd$.
\hfill\qed
\end{itemize}
\end{df}
\index{unramifiedly good filtered $\varrho$-flat bundle}
\index{good filtered $\varrho$-flat bundle}

\subsection{Residue}
\label{subsection;10.5.5.11}
\index{residue}

\subsubsection{Unramified case}

\index{residue}
\label{subsection;10.5.6.1}

Let $\nbigx\lrarr \nbigk$, $\nbigd$, $\varrho$ and 
$(\nbige,\DD)$ be as in
Section \ref{subsection;10.5.3.3}.
Let $E$ be an unramifiedly good lattice 
of $(\nbige,\DD)$.
Let $\nbigd_i$ be an irreducible component of $\nbigd$.
For each $P\in \nbigd_i$,
we have 
$\Res_{\nbigd_i}(\DD_{|\Phat})
 \in \End(E_{|\Phat\cap \nbigd_i})$.
(See Subsection \ref{subsection;10.5.2.5}.)

\begin{lem}
We have the residue endomorphism
$\Res_{\nbigd_i}(\DD)\in\End(E_{|\nbigd_i})$
such that $\Res_{\nbigd_i}(\DD)_{|\Phat}
=\Res_{\nbigd_i}(\DD_{|\Phat})$
for any $P\in \nbigd_i$.
The eigenvalues of $\Res_{\nbigd_i}(\DD)$
are the pull back of possibly multi-valued 
functions on $\nbigk$.
In particular,
their restriction to $\nbigd_i\times_\nbigk\{y\}$
are constant if $\varrho$ is not constantly $0$
around $y$.
\end{lem}
\pf
The first claim follows from the construction of
$\Res_{\nbigd_i}(\DD_{|\Phat})$
and Proposition \ref{prop;10.5.3.5}.
The second claim follows from
the first one and Lemma \ref{lem;10.5.8.1}.
\hfill\qed

\subsubsection{Ramified case}

Let $\nbigx\lrarr \nbigk$, $\nbigd$ and $\varrho$ be as above.
If $(\vecE_{\ast},\DD)$
is unramifiedly good on $(\nbigx,\nbigd)$,
we have the induced endomorphism
$\Res_i(\DD)$ on
$\prolongg{\vecc}{E}_{|\nbigd_i}$.
It preserves the induced filtration $\lefttop{i}F$
of $\prolongg{\vecc}{E}_{|\nbigd_i}$,
and hence we have the induced endomorphism
$\Gr^F_a\Res_i(\DD)$
of $\lefttop{i}\Gr^F_a(\prolongg{\vecc}{E})$.

\begin{prop}
\label{prop;10.5.5.12}
Even if a good filtered $\varrho$-flat bundle
$(\vecE_{\ast},\DD)$
is not necessarily unramified,
we have the induced endomorphism
$\Gr^F\Res_i(\DD)$
of $\lefttop{i}\Gr^F_a\bigl(\prolongg{\vecc}{E}\bigr)$
on $\nbigd_i$ for each $i\in \Lambda$.
It preserves the induced filtrations 
$\lefttop{j}F$ of 
$\lefttop{i}\Gr^F_a
 \bigl(\prolongg{\vecc}{E}\bigr)_{|\nbigd_i\cap \nbigd_j}$.

The eigenvalues of 
$\Gr^F\Res_i(\DD)$
are the pull back of possibly multi-valued 
functions on $\nbigk$.
In particular,
their restriction to $\nbigd_i\times_\nbigk\{y\}$
are constant if $\varrho$ is not constantly $0$
around $y$.
\end{prop}

Due to Proposition \ref{prop;10.5.5.12},
$\Gr^F\Res_i(\DD)$ ($i\in I$) induce
the endomorphisms of
$\lefttop{I}\Gr^F_{\veca}(\prolongg{\vecc}{E})$,
which are also denoted by
$\Gr^F\Res_i(\DD)$
or $\Res_i(\DD)$.
In the following,
$\Gr^F\Res_{i}(\DD)$
are often denoted by
$\Res_i(\DD)$
for simplicity of the description.

\subsubsection{Proof of 
Proposition \ref{prop;10.5.5.12}}
First, we consider the case
$\nbigx=\Delta^n$ and $\nbigd=\{z_1=0\}$.
We put $\nbigxtilde:=\nbigx$ and $\nbigdtilde=\nbigd$,
and we have a ramified covering
$\varphi_e:(\nbigxtilde,\nbigdtilde)\lrarr (\nbigx,\nbigd)$
given by $\varphi_e(z_1,\ldots,z_n)
=(z_1^{e},z_2,\ldots,z_n)$
such that 
the induced filtered $\varrho$-flat bundle
$(\vecEtilde_{\ast},\DDtilde)$
on $(\nbigxtilde,\nbigdtilde)$ is unramifiedly good.
We take $c\in\real$ and we put 
$\ctilde:=c\,e$.
We have the residue
$\Res(\DDtilde)$ of 
$\prolongg{\ctilde}{\Etilde}$ on $\nbigdtilde$.
Let $\mu_e:=\{\omega\in\cnum\,|\,\omega^e=1\}$,
which naturally acts on $\nbigxtilde$ 
by $\omega^{\ast}(z_1)=\omega\, z_1$,
and $(\vecEtilde_{\ast},\DDtilde)$ 
is $\mu_e$-equivariant.
The endomorphism $\Res(\DDtilde)$
is $\mu_e$-equivariant.

We can take a frame $\vecvtilde$ of
$\prolongg{\ctilde}{\Etilde}$
such that
(i) it is compatible with the induced filtration $F$ 
 of $\prolongg{\ctilde}{\Etilde}_{|\nbigdtilde}$,
(ii) for each $p$ we have $b_p\in\seisuu$
satisfying $0\leq b_p<e-1$ and
 $\omega^{\ast}\vtilde_p
=\omega^{-b_p}\,\vtilde_p$.
We put $a(\vtilde_p):=\deg^F(\vtilde_p)$.
They induce the frame of
$\Gr^F\bigl(\prolongg{\ctilde}{\Etilde}\bigr)
=\bigoplus_{\ctilde-1<\atilde\leq \ctilde}
 \Gr^F_{\atilde}\bigl(\prolongg{\ctilde}{\Etilde}\bigr)$.

We put
$v_p:=z^{b_p}\, \vtilde_p$,
which is a $\mu_e$-invariant section.
The tuple $\vecv=(v_p)$ naturally gives
a frame of $\prolongg{c}{E}$
compatible with the parabolic filtration.
Hence, they induce a frame
of $\Gr^F(\prolongg{c}{E})
=\bigoplus_{c-1<a\leq c}
 \Gr^F_a(\prolongg{c}{E})$.
The frames give an isomorphism:
\begin{equation}
 \label{eq;08.2.23.6}
 \Gr^F_{\atilde}(\prolongg{\ctilde}{\Etilde})
\simeq
 \bigoplus_{\atilde-e\,a\in\seisuu}
 \Gr^F_a(\prolongg{c}{E})
\end{equation}
The decomposition of 
$\Gr^F_{\atilde}(\prolongg{\ctilde}{\Etilde})$
corresponding to (\ref{eq;08.2.23.6})
is given by the eigen decomposition
with respect to the action of $\omega^{\ast}$.
Since $\Res(\DDtilde)$
is $\mu_m$-equivariant,
it induces endomorphisms
$G_a$ of $\Gr^{F}_a(\prolongg{c}{E})$.
It is easy to check that the isomorphism 
(\ref{eq;08.2.23.6}) is independent of 
the choice of a frame $\vecvtilde$.
It is also independent of the choice of 
a coordinate chart
up to constant multiplication on
each direct summand of the right hand side.
Hence, $G_a$ is independent
of the choice of frames and coordinate system.
For each $c-1<a\leq c$,
let $b(a)\in\seisuu$ be determined by
$\ctilde-1<e\,a+b(a)\leq \ctilde$.
In this case,
we define the endomorphism
$\Gr^F_a\Res(\DD)$
of $\Gr^F_a(\prolongg{c}{E})$
as follows:
\[
 \Gr^F_a\Res(\DD):=
e^{-1}\bigl(G_a+\varrho\,b(a)\bigr)
\]

\begin{lem}
\label{lem;08.2.23.10}
If $(\vecE_{\ast},\DD)$
is unramified,
it is the same as the endomorphism
induced by the residue $\Res(\DD)$.
\end{lem}
\pf
By considering the completion along $\nbigd$,
the problem can be reduced to
the regular case.
Then, the claim can be checked by a direct calculation.
\hfill\qed

\vspace{.1in}
By using Lemma \ref{lem;08.2.23.10},
we can check that 
$\Gr^F_a\Res(\DD)$ is independent
of the choice of a ramified covering
$(\nbigxtilde,\nbigdtilde)\lrarr (\nbigx,\nbigd)$.
Thus, we obtain the well defined endomorphism
$\Gr^F_a(\DD)$ 
of $\Gr^F_a(\prolongg{c}{E})$ 
in the case that $\nbigd$ is smooth.

\vspace{.1in}
Let $\nbigd_i^{\circ}:=\nbigd_i\setminus\bigcup_{j\neq i}\nbigd_j$.
We have obtained the endomorphism
$\Gr^F_a\Res_i(\DD)$ of
$\Gr^F_a(\prolongg{\vecc}{E})_{|\nbigd_i^{\circ}}$.
Let us show that it is extended
to an endomorphism of
$\Gr^F_a(\prolongg{\vecc}{E})$,
and that it preserves the parabolic filtrations
$\lefttop{l}F$ $(l\neq i)$.
Since it is a local property,
we have only to consider the case $\nbigx:=\Delta^n$
and $\nbigd=\bigcup_{i=1}^{\ell}\{z_i=0\}$.
We have a ramified covering
$\varphi_e:(\nbigxtilde,\nbigdtilde)\lrarr (\nbigx,\nbigd)$
given by $\varphi_e(z_1,\ldots,z_n)
=(z_1^{e},\ldots,z_{\ell}^e,
 z_{\ell+1},\ldots,z_{n})$,
such that the induced filtered $\varrho$-flat bundle
$(\vecEtilde_{\ast},\DDtilde)$
is unramifiedly good.
We put $\Gal(\nbigxtilde/\nbigx):=
 \bigl\{\vecomega=(\omega_1,\ldots,\omega_{\ell})\,\big|\,
 \omega_i\in\mu_e\bigr\}$.
We have the natural action of
$\Gal(\nbigxtilde/\nbigx)$ on $\nbigxtilde$
given by 
$\vecomega^{\ast}z_j=\omega_j\, z_j$
for $j=1,\ldots,\ell$.
It is lifted to the action on 
$(\vecEtilde_{\ast},\DDtilde)$,
and $(\vecE_{\ast},\DD)$
is the descent.

Let $\vecc\in\real^{\ell}$ and 
$\vecctilde:=e\,\vecc$.
Let $c_i$ and $\ctilde_i$ denote
the $i$-th components of 
$\vecc$ and $\vecctilde$,
respectively.
For any $\ctilde_i-1<\atilde\leq \ctilde_i$,
we have the endomorphism
$\Gr^F_{\atilde}
 \Res_{\nbigdtilde_i}(\DDtilde)$  of 
$\lefttop{i}\Gr^F_{\atilde}
 (\prolongg{\vecctilde}{\Etilde})$
on $\nbigdtilde_i$.
It is $\Gal(\nbigxtilde/\nbigx)$-equivariant,
and the restriction 
$\Gr^F_{\atilde}
 \Res_{\nbigdtilde_i}
 (\DDtilde)_{|\nbigdtilde_i\cap \nbigdtilde_j}$ 
preserves the induced filtration
$\lefttop{j}F$ of 
$\lefttop{i}\Gr^F_{\atilde}
 (\prolongg{\vecctilde}{\Etilde})
 _{|\nbigdtilde_i\cap\nbigdtilde_j}$.

We can take a frame
$\vecvtilde=(\vtilde_p)$ of
$\prolongg{\vecctilde}{\Etilde}$
such that
(i) it is compatible with the filtrations 
 $\lefttop{k}F$ $(k=1,\ldots,\ell)$,
(ii) there exist tuples of integers
 $\vecb_p=(b_{p,1},\ldots,b_{p,\ell})$
 satisfying $0\leq b_{p,k}<e-1$
 and 
 $\vecomega^{\ast}\vtilde_p
=\prod_{k=1}^{\ell}
   \omega_k^{-b_{p,k}}\, \vtilde_p$.
(See Section 2.3 of \cite{mochi2}.)
We put $a_k(\vtilde_p):=\lefttop{k}\deg^F(\vtilde_p)$.
Let $\Atilde^{(i)}$ be the matrix valued
holomorphic function on $\nbigdtilde_i$,
determined by
$\Res_i(\DDtilde)\vecvtilde_{|\nbigdtilde_i}
=\vecvtilde_{|\nbigdtilde_i}\, \Atilde^{(i)}$,
i.e.,
$\Res_i(\DDtilde)\vtilde_{q|\nbigdtilde_i}
=\sum_p \Atilde^{(i)}_{p,q}\,\vtilde_{p|\nbigdtilde_i}$.
We have $\Atilde^{(i)}_{p,q}=0$
unless $b_{p,i}=b_{q,i}$
and $a_i(\vtilde_p)\leq a_i(\vtilde_q)$.
Due to the $\Gal(\nbigxtilde/\nbigx)$-equivariance of
$\Res_i(\DDtilde)$,
the following functions are holomorphic on $\nbigdtilde_i$
and $\Gal(\nbigxtilde/\nbigx)$-invariant:
\begin{equation}
 A^{(i)}_{p,q}:=
\Atilde^{(i)}_{p,q}\,
 \prod_{k\neq i}
 z_k^{b_{p,k}-b_{q,k}}
\end{equation}
Moreover, we have
$A^{(i)}_{p,q|\nbigdtilde_l\cap\nbigdtilde_i}=0$ 
for $l\neq i$,
if either one of the following holds:
\begin{equation}
 \label{eq;08.2.23.11}
 \mbox{\rm (i) }
 b_{p,l}-b_{q,l}>0,
\quad\quad
 \mbox{\rm (ii) }
 b_{p,l}=b_{q,l},\,\,
 a_l(\vtilde_p)>a_l(\vtilde_q)
\end{equation}

We take $c_i-1<a\leq c_i$,
which determines $\ctilde_i-1<\atilde\leq \ctilde_i$
such that $b(a):=\atilde-e\, a\in\seisuu$.
We put $I(a,i):=\bigl\{p\,\big|\,
 a_i(\vtilde_p)=\atilde,\,\,
 b_{p,i}=b(a) \bigr\}$.
Let $\vecutilde_a$ be the tuple
$\bigl(
 \utilde_{a,p}:=
 \vtilde_{p}\,\big|\,
 p\in I(a,i)\bigr)$.
We put
\[
 u_{a,p}:=
 \prod_kz_k^{b_{p,k}}\, 
 \utilde_{a,p}.
\]
Then, $\vecu_a=(u_{a,p}\,|\,p\in I(a,i))$
naturally induces a frame of
$\lefttop{i}\Gr^F_a(\prolongg{\vecc}{E})$
compatible with the induced filtrations
$\lefttop{l}F$ $(l\neq i)$ on $\nbigd_l\cap \nbigd_i$.
Let $A^{(i)}_a$ be the matrix valued
holomorphic function on $\nbigd_i$
given by 
$\bigl(
 A^{(i)}_{p,q}\,\big|\,
 p,q\in I(a,i)\bigr)$.
By definition,
we have
\[
 \Gr^F_a\Res_i(\DD)\vecu_{a|\nbigd_i}
=\vecu_{a|\nbigd_i}\, 
e^{-1}\bigl(
 A^{(i)}_a+\varrho\, b(a)
\bigr).
\]
It implies that
$\Gr^F_a\Res_i(\DD)$
is extended to an endomorphism of
$\Gr^F_a(\prolongg{\vecc}{E})$
on $\nbigd_i$.
If $\lefttop{l}\deg^F(u_{p})>\lefttop{l}\deg^F(u_q)$,
one of (\ref{eq;08.2.23.11}) occurs,
and hence $A^{(i)}_{p,q|\nbigd_l\cap \nbigd_i}=0$.
It implies that
$\Gr^F_a\Res_i(\DD)_{|\nbigd_l\cap \nbigd_i}$
preserves the filtrations $\lefttop{l}F$ $(l\neq i)$.
\hfill\qed

\subsubsection{Some notation}
Let us consider the case  $\nbigk=\{y\}$.
If $\varrho\neq 0$,
the eigenvalues of $\Res_i(\DD)$
are constant. 
The endomorphisms
$\Res_i(\DD)$ $(i\in I)$ on
$\lefttop{I}\Gr^{F}_{\veca}(\prolongg{\vecc}{E})$
are commutative.
Hence, we have the generalized eigen decomposition
\[
 \lefttop{I}\Gr^F_{\veca}\bigl(
 \prolongg{\vecc}{E}\bigr)
=\bigoplus_{\vecalpha}
 \lefttop{I}\Gr^{F,\EE}_{(\veca,\vecalpha)}
 \bigl(\prolongg{\vecc}{E}\bigr),
\]
where the eigenvalues of 
$\Gr^F\Res_i(\DD)$ on 
$\lefttop{I}\Gr^{F,\EE}_{(\veca,\vecalpha)}
 \bigl(\prolongg{\vecc}{E}\bigr)$
are the $i$-th components of $\vecalpha$.
\index{bundle
 $\lefttop{I}\Gr^{F,\EE}_{(\veca,\vecalpha)}
 \bigl(\prolongg{\vecc}{E}\bigr)$}
Recall that
we often consider the following sets
in this situation:
\[
 \KMS(\prolongg{\vecc}{E},\DD,I):=
 \bigl\{
 (\veca,\vecalpha)\,\big|\,
 \lefttop{I}\Gr^{F,\EE}_{(\veca,\vecalpha)}
 (\prolongg{\vecc}{E})\neq 0
 \bigr\}
\]
\[
  \KMS(\vecE_{\ast},\DD,I):=
 \!\!\!\bigcup_{\vecc\in\real^{\Lambda}}\!\!\!
 \KMS(\prolongg{\vecc}{E},\DD,I)
\]
\index{set
$\KMS(\prolongg{\vecc}{E},\DD,I)$}
\index{set
$\KMS(\vecE_{\ast},\DD,I)$}
\[
\Sp(\prolongg{\vecc}{E},\DD,I):=
 \bigl\{
 \alpha\,\big|\,
 (\veca,\vecalpha)\in 
 \KMS\bigl(\prolongg{\vecc}{E},\DD,I\bigr)
 \bigr\}
\]
\[
  \Sp(\vecE_{\ast},\DD,I):=
 \bigcup_{\vecc\in\real^{\Lambda}}
 \Sp(\prolongg{\vecc}{E},\DD,I)
\]
\index{set $\Sp(\vecE_{\ast},\DD,I)$}
\index{set $\Sp(\prolongg{\vecc}{E},\DD,I)$}
Each element of $\KMS(\vecE_{\ast},\DD,I)$
is called the KMS-spectrum of 
$(\vecE_{\ast},\DD)$ at $\nbigd_I$.
\index{KMS-spectrum}

\begin{rem}
Even in the case $\varrho=0$,
a similar notion makes sense,
if the eigenvalues of $\Res_i(\DD^0)$
are assumed to be constant.
The condition will be satisfied 
when we consider wild harmonic bundles.
\hfill\qed
\end{rem}

\subsection{Filtered bundle (Appendix)}
\label{subsection;10.5.28.1}

Let $X$ be a complex manifold
with a simple normal crossing hypersurface 
$D=\bigcup_{i\in\Lambda}D_i$.
A filtered sheaf on $(X,D)$
is a datum $\vecE_{\ast}=
 \bigl(\vecE,\{\prolongg{\vecc}{E}\}\,
 \big|\,\vecc\in\real^{\Lambda}\bigr)$
as follows:
\begin{itemize}
\item
$\vecE$ is a torsion-free coherent
$\nbigo_X(\ast D)$-module.
\item
 $\{\prolongg{\vecc}{E}\}$
is an increasing filtration by
coherent $\nbigo_X$-submodules of
$\vecE$
indexed by $\real^{\Lambda}$
such that
(i) $\vecE_{|X-D}=
 \prolongg{\vecc}{E}_{|X-D}$
for any $\vecc$,
(ii) $\prolongg{\veca}{E}=
 \bigcap_{\veca<\vecb}
 \prolongg{\vecb}{E}$,
(iii)
 $\vecE=\bigcup_{\veca\in\real^{\Lambda}}
 \prolongg{\veca}{E}$.
Here, the order on $\real^{\Lambda}$
is given by 
$\veca\leq\vecb$
$\Longleftrightarrow$
$a_i\leq b_i$ $(\forall i)$.
\item
 $\prolongg{\veca'}{E}
=\prolongg{\veca}{E}
 \otimes\nbigo_X(-\sum n_j\, D_j)$
as submodules of $\vecE$,
where $\veca'=\veca-(n_j\,|\,j\in \Lambda)$.
\item
For each $\vecc\in\real^{\Lambda}$,
let $\lefttop{i}\nbigf$ be a filtration 
of $\prolongg{\vecc}{E}$
indexed by $\openclosed{c_i-1}{c_i}$
is given as follows:
\begin{equation}
 \label{eq;10.5.28.1}
 \lefttop{i}\nbigf_d(\prolongg{\vecc}{E})
:=\bigcup_{\substack{a_i\leq d\\ \veca\leq \vecc}}
 \prolongg{\veca}{E}.
\end{equation}
Then the tuple 
$\prolongg{\vecc}{E}_{\ast}:=
 \bigl(\prolongg{\vecc}{E},
 \{\lefttop{i}\nbigf\,|\,i\in S\}
 \bigr)$ is a $\vecc$-parabolic sheaf,
i.e.,
the sets
$\bigl\{a\,\big|\,
 \lefttop{i}\Gr^{\nbigf}_a(\prolongg{\vecc}{E})
 \neq 0\bigr\}$
are finite.
\end{itemize}
\index{filtered sheaf}
\index{filtered bundle}
See Subsection 3.2 of \cite{mochi4}
for some property of filtered sheaf.
Each $\prolongg{\vecc}{E}_{\ast}$
is called the $\vecc$-truncation of $\vecE_{\ast}$.
We can reconstruct
$\vecE_{\ast}$ from $\prolongg{\vecc}{E}_{\ast}$.
If each $\prolongg{\vecc}{E}$ is locally free,
$\vecE_{\ast}$ is called a filtered bundle.
(See Remark \ref{rem;10.5.28.2} below.)

\subsubsection{Induced filtrations}

Let $\vecE_{\ast}$ be a filtered bundle
on $(X,D)$.
For each $\vecc$-truncation
$\prolongg{\vecc}{E}$,
we have a filtration $\lefttop{i}\nbigf$
given as in (\ref{eq;10.5.28.1}).
Let $\lefttop{i}F_d(\prolongg{\vecc}{E}_{|D_i})$
denote the image of the induced map
$\lefttop{i}\nbigf_d(\prolongg{\vecc}{E}_{|D_i})
\lrarr
 \prolongg{\vecc}{E}_{|D_i}$.
It is called the parabolic filtration 
of $\prolongg{\vecc}{E}$.
\index{parabolic filtration $\lefttop{i}F$}
For $I\subset \Lambda$,
we have the induced filtrations
$\lefttop{i}F$ $(i\in I)$ of
$\prolongg{\vecc}{E}_{|D_I}$.
It is known 
(\cite{borne1}, \cite{borne2},
\cite{hertling-sevenheck2})
that they are compatible
in the sense that it has locally a splitting,
i.e.,
for each $P\in D_I$,
take a small neighbourhood $D_{I,P}$ of $P$
in $D_I$,
then we have a splitting
$\prolongg{\vecc}{E}_{|D_I,P}
=\bigoplus G_{\vecd}$
such that
$\lefttop{I}F_{\vecb}(\prolongg{\vecc}{E_{|D_I,P}})
:=\bigcap
 \lefttop{i}F_{d_i}(\prolongg{\vecc}{E_{|D_I,P}})
=\bigoplus_{\vecd\leq\vecb}
 G_{\vecd}$.
It also implies the locally abelian condition
in \cite{i-s},
i.e.,
for a small neighbourhood $X_P$ of $P$ in $X$,
we can take a decomposition
$\prolongg{\vecc}{E}_{|X_P}
=\bigoplus H_{\vecd}$
such that 
$\lefttop{i}F_b(\prolongg{\vecc}{E}_{|D_{i}\cap X_P})
=\bigoplus_{d_i\leq b}
 H_{\vecd|D_i}$.

\begin{rem}
\label{rem;10.5.28.2}
The above compatibility condition was imposed
in our old but equivalent definition
of filtered bundle (\cite{mochi4} and \cite{mochi2}).
It is used implicitly,
even if it does not appear in the definition.
\hfill\qed
\end{rem}

Let $I$ be a subset of $\Lambda$.
Let $D_I:=\bigcap_{i\in I}D_i$.
For $\veca\in\real^I$,
we will often consider
\[
 \lefttop{I}F_{\veca}\bigl(
 \prolongg{\vecc}{E}_{|D_I}
 \bigr):=
 \bigcap_{i\in I}
 \lefttop{i}F_{a_i}
 \bigl(\prolongg{\vecc}{E}_{|D_I}\bigr),
\quad\quad
 \lefttop{I}\Gr^F_{\veca}\bigl(
 \prolongg{\vecc}{E}
 \bigr):=
\frac{\lefttop{I}F_{\veca}\bigl(
 \prolongg{\vecc}{E}_{|D_I}\bigr)}
{\sum_{\vecb\lneq\veca}
 \lefttop{I}F_{\vecb}\bigl(
 \prolongg{\vecc}{E}_{|D_I}\bigr)}.
\]
\index{bundle
 $\lefttop{I}\Gr^F_{\veca}\bigl(
 \prolongg{\vecc}{E}
 \bigr)$}
Here, $\vecb\lneq\veca$
means ``$\vecb\leq\veca$ and $\vecb\neq\veca$''.
We often consider the following sets in this situation:
\[
 \Par\bigl(\prolongg{\vecc}{E},I
 \bigr):=
 \bigl\{ \veca\in\real^I\,\big|\,
 \lefttop{I}\Gr^F_{\veca}(\prolongg{\vecc}{E})\neq 0
 \bigr\}
\]
\[
  \Par(\vecE_{\ast},I):=
 \bigcup_{\vecc\in\real^{\Lambda}}
 \Par(\prolongg{\vecc}{E},I)
\]
\index{set $\Par(\vecE_{\ast},I)$}
\index{set $\Par(\prolongg{\vecc}{E},I)$}

\subsubsection{Compatible frame}

For $P\in X$,
let $X_P$ denote a small neighbourhood of $P$
in $X$, and we put $D_P:=D\cap X_P$,
and $D_{J,P}:=D_J\cap X_P$.
Let $\Lambda(P):=\{j\in\Lambda\,|\,P\in D_j\}$.
Let $\vecE_{\ast}$ be a filtered bundle 
on $(X,D)$.
We can take a frame $\vecv$ 
of $\prolongg{\vecc}{E}_{|X_P}$
with the following property:
\begin{itemize}
\item
 For each $v_p$,
 the tuple of numbers
 $\veca(v_p)\in \prod_{j\in\Lambda(P)}
 \openclosed{c_j-1}{c_j}$.
\item
 For $J\subset \Lambda(P)$,
 $\lefttop{J}F_{\vecb}(\prolongg{\vecc}{E}_{|D_{J,P}})$
  is generated by
 $v_p$ such that
 $a_j(v_p)\leq b_j$ $(\forall j\in J)$.
\end{itemize}
Such a frame is called compatible with
the parabolic structure of $\prolongg{\vecc}{E}$.
The numbers $a_j(v_p)$ is 
often written as $\lefttop{j}\deg^F(v_p)$.
\index{compatible frame}
\index{degree $\lefttop{j}\deg^F(v_p)$}

\subsubsection{Pull back of filtered bundles}

\label{subsection;07.11.5.60}

Let us recall the pull back of filtered bundles.
See \cite{i-s} for a more systematic treatment.
Let $X:=\Delta^n_z$ and 
$D:=\bigcup_{i=1}^{\ell}\{z_i=0\}$.
Let $\Xtilde:=\Delta^m_w$ and 
$\Dtilde:=\bigcup_{j=1}^k\{w_j=0\}$.
Let $\varphi:\Xtilde\lrarr X$ be a morphism
such that $\varphi^{-1}(D)\subset \Dtilde$.
Then, $\varphi^{\ast}(z_i)=
 \prod_{j=1}^{k}w_{j}^{\alpha_{j,i}}\, g_i$
for some invertible function $g_i$
$(i=1,\ldots,\ell)$.
Let $\varphi^{\ast}:\real^{\ell}\lrarr\real^k$
be given by
$\varphi^{\ast}_j(\vecb):=
 \sum_{i=1}^{\ell}\alpha_{j,i}\, b_i$.
For any $\vecb\in\real^{k}$,
we set
\[
 \nbigs(\vecb):=\bigl\{
 (\veca,\vecn)\in\real^{\ell}\times\seisuu_{\geq 0}^k
 \,\big|\,
 \varphi^{\ast}(\veca)+\vecn\leq\vecb
 \bigr\}
\]

Let $\vecE_{\ast}$ be a filtered sheaf on $(X,D)$.
We put
\[
 \prolongg{\vecb}{\Etilde}
:=\sum_{(\veca,\vecn)\in\nbigs(\vecb)}
 \vecw^{-\vecn}\,
 \varphi^{\ast}\bigl(\prolongg{\veca}{E}\bigr).
\]
Thus, we obtain a filtered sheaf
$\vecEtilde_{\ast}$ on $(\Xtilde,\Dtilde)$.
It is independent of the choice of the coordinate systems
$\vecz$ and $\vecw$.

\begin{lem}
If $\vecE_{\ast}$ is a filtered bundle,
$\vecEtilde_{\ast}$ is also a filtered bundle.
\end{lem}
\pf
Let $\vecv$ be a frame of
$\prolong{E}$ compatible with the parabolic filtrations.
We put $a_i(v_p):=\lefttop{i}\deg^F(v_p)$
and $\veca(v_p):=\bigl(a_i(v_p)\bigr)$.
Let $\vecc=(c_j)\in\real^k$.
Let $n_j(v_i)$ be the integers 
determined by the condition
$c_j-1<
 \varphi^{\ast}_j(\veca(v_p))+n_j(v_p)
 \leq c_j$.
We set
\[
 \prolongg{\vecc}{\vtilde_p}:=
 \prod_{j}w_j^{-n_j(v_p)}\,
 \varphi^{\ast}(v_p).
\]
Then, we can check
that $\prolongg{\vecc}{\vecvtilde}:=
 \bigl(\prolongg{\vecc}{\vtilde_p}\bigr)$
gives a frame of $\prolongg{\vecc}{\Etilde}$
compatible with the parabolic filtrations.
\hfill\qed

\vspace{.1in}

Let $X$ (resp. $\Xtilde$) be a complex manifold
with a simple normal crossing hypersurface $D$
(resp. $\Dtilde$).
Let $\varphi:(\Xtilde,\Dtilde)\lrarr (X,D)$ be a morphism.
Let $\vecE_{\ast}$ be a filtered bundle on $(X,D)$.
Applying the above procedure locally,
we obtain a filtered bundle $\vecEtilde_{\ast}$
on $(\Xtilde,\Dtilde)$ globally.

\subsubsection{Descent with respect to
 a ramified covering}
\label{subsection;08.9.29.25}

Let $X:=\Delta^n_z$ and 
$\Xtilde:=\Delta^n_w$.
Let $D:=\bigcup_{i=1}^{\ell}\{z_i=0\}$
and $\Dtilde:=\bigcup_{i=1}^{\ell}\{w_i=0\}$.
Let $\varphi:(\Xtilde,\Dtilde)\lrarr (X,D)$
be a ramified covering given by
$\varphi(w_1,\ldots,w_n)
=(w_1^{m_1},\ldots,w_{\ell}^{m_{\ell}},w_{\ell+1},
 \ldots,w_n)$.
Let $\varphi^{\ast}:\real^{\ell}\lrarr\real^{\ell}$
be given by $\varphi^{\ast}(a_1,\ldots,a_{\ell})
=(m_1\, a_1,\ldots,m_{\ell}\, a_{\ell})$.

Let $\vecEtilde_{\ast}=\bigl(
 \prolongg{\vecb}{\Etilde} \bigr)$ 
be a filtered sheaf on $(\Xtilde,\Dtilde)$,
which is equipped with the $\Gal(\Xtilde/X)$-action.
Let $\prolongg{\veca}{E}$ be the decent of 
$\prolongg{\varphi^{\ast}(\veca)}{\Etilde}$.
Thus, we obtain a filtered sheaf
$\vecE_{\ast}$ on $(X,D)$.
It is easy to see that
$\vecE_{\ast}$ is also a filtered bundle,
if $\vecEtilde_{\ast}$ is a filtered bundle.
The construction is independent of the choice of
a coordinate system.

For any general ramified covering
of complex manifolds,
we obtain the global decent 
by applying the above procedure locally.

\section{Good lattice in the level $\vecm$}
\label{subsection;08.9.28.120}

We introduce an auxiliary concept of
good lattice in the level $\vecm$.
It seems useful in the inductive study 
on Stokes structure.
Because we consider only the unramified case,
we omit to distinguish it.

\subsection{Order of the pole}
\label{subsection;10.7.2.2}

We introduce an auxiliary notion of  ``order'' of
the pole of a $1$-form
or a meromorphic flat $\varrho$-connection.
Let $\nbigx\lrarr \nbigk$, $\nbigd$, $\varrho$ be as before.
Let $\nbigd=\bigcup_{i\in\Lambda}\nbigd_i$
be the irreducible decomposition.
Let $\vecm\in\seisuu_{\leq 0}^{\Lambda}$.
We put $\nbigd^{(1)}:=\bigcup_{m_i<0}\nbigd_i$
and $\nbigd^{(2)}:=\bigcup_{m_i=0}\nbigd_i$.
\begin{df}
Let $\omega$ be a holomorphic section of
$F\otimes\Omega^1_{\nbigx}(\ast \nbigd)$,
where $F$ is a locally free $\nbigo_{\nbigx}$-module.
We say $\ord\omega\geq \vecm$, if 
it is contained in
\[
F\otimes\Bigl(
 \vecz^{\vecm}
 \Omega^1_{\nbigx/\nbigk}(\log \nbigd^{(1)})
+\Omega^{1}_{\nbigx/\nbigk}(\log \nbigd^{(2)})
 \Bigr).
\]
We have similar conditions for $1$-forms
on formal complex spaces.
\hfill\qed
\end{df}
\index{order}

Let $E$ be a locally free $\nbigo_\nbigx$-module
with a meromorphic $\varrho$-flat
connection $\DD$
of $E(\ast \nbigd)$.
\begin{df}
We say $\ord(\DD)\geq\vecm$,
if the following holds:
\begin{equation}
\label{eq;10.5.10.1}
 \DD E\subset 
 E\otimes\bigl(
 \vecz^{\vecm}\Omega^1_{\nbigx/\nbigk}(\log \nbigd^{(1)})
+\Omega^{1}_{\nbigx/\nbigk}(\log \nbigd^{(2)})
\bigr)
\end{equation}
We have a similar condition
for the lattice of a meromorphic
flat $\varrho$-connection
on formal complex spaces.
\hfill\qed
\end{df}
\index{order}
Let $\vecv$ be a frame of $E$.
Let $A$ be determined by
$\DD\vecv=\vecv\, A$.
We have $\ord\DD\geq \vecm$
if and only if
$\ord A\geq\vecm$.
\begin{rem}
For any $j$ such that $m_j=0$,
we have the induced endomorphism
$\Res_j(\DD)$ of
$E_{|\nbigd_j}$.
\hfill\qed
\end{rem}

\begin{rem}
The condition {\rm(\ref{eq;10.5.10.1})}
implies the following:
\begin{equation}
\label{eq;10.5.18.1}
\DD E\subset
 \vecz^{\vecm}E\otimes
 \Omega^1_{\nbigx/\nbigk}\bigl(\log \nbigd\bigr)
\end{equation}
It was adopted as the definition of order
in the older version of this monograph.
The difference is not essential for our purpose.
The condition {\rm(\ref{eq;10.5.10.1})}
might be more natural,
and {\rm(\ref{eq;10.5.18.1})} might be easier
to state.
\hfill\qed
\end{rem}

\subsection{Good set of irregular values
 in the level $\vecm$}
\label{subsection;07.12.26.1}

\index{good set of irregular values
in the level $\vecm$}
\index{good set of irregular values
in the level $(\vecm,i)$}

This subsection is a complement of
Section \ref{subsection;07.11.6.1}.
Let $Y$ be a complex manifold.
We put $X:=\Delta^{\ell}\times Y$,
and $D:=\bigcup_{i=1}^{\ell}\{z_i=0\}$.

\begin{df}
Let $\vecm\in\seisuu_{\leq 0}^{\ell}\setminus\{\veczero\}$.
A finite set of meromorphic functions
 $\nbigi=\bigl\{\gminia=
 \gminia_{\vecm}\vecz^{\vecm}\bigr\}
 \subset M(X,D)$
is called 
a weakly good set of irregular values on $(X,D)$
in the level $\vecm$,
if the following holds:
\begin{itemize}
\item
 $\gminia_{\vecm}-\gminib_{\vecm}$
 are invertible holomorphic functions on $X$
 for any two distinct $\gminia,\gminib\in \nbigi$.
\end{itemize}
If moreover the following condition holds
for an integer $i$ such that $m_{i}<0$,
$\nbigi$ is called a weakly good set of irregular values
on $(X,D)$ in the level $(\vecm,i)$.
\begin{itemize}
\item
 $\gminia_{\vecm}-\gminib_{\vecm}$ 
 are independent of the variable $z_{i}$.
\end{itemize}
A weakly good set of irregular values
in the level $(\vecm,i)$ is called
a good set of irregular values
in the level $(\vecm,i)$,
if the following condition is satisfied.
\begin{itemize}
\item
$\gminia_{\vecm}$ are 
 holomorphic functions on $X$,
 which are independent of $z_i$.
\hfill\qed
\end{itemize}
\end{df}

Let $\nbigi$ be a weakly good set of irregular values
in the level $(\vecm,i)$.
We choose any $\gminia^{(0)}\in\nbigi$.
Then, the set
$\nbigi':=\bigl\{
 \gminia-\gminia^{(0)}\,\big|\,
 \gminia\in \nbigi
 \bigr\}$
is a good set of irregular values
in the level $(\vecm,i)$.

For a weakly good set of irregular values $\nbigi$
in the level $(\vecm,i)$,
we put $\nbigi^{\lor}:=\bigl\{
 -\gminia\,\big|\,\gminia\in\nbigi
 \bigr\}$.
For $\nbigi_i$ $(i=1,2)$,
we put $\nbigi_1\otimes\nbigi_2:=\bigl\{
 \gminia_1+\gminia_2\,\big|\,\gminia_i\in\nbigi_I
 \bigr\}$
and $\nbigi_1\oplus\nbigi_2=\nbigi_1\cup\nbigi_2$,
which are not necessarily weakly good
in the level $(\vecm,i)$.
\index{set $\nbigi^{\lor}$}
\index{set $\nbigi_1\otimes\nbigi_2$}
\index{set $\nbigi\oplus\nbigi_2$}

\subsubsection{}
\label{subsection;10.6.11.1}
Let $\nbigj$ be a good set of irregular values on $(X,D)$.
Take an auxiliary sequence
$\vecm(0),\vecm(1),\ldots,\vecm(L)$.
Let us observe that 
we have the associated
good sets of irregular values on $(X,D)$
in the level $(\vecm(i),\gminih(i))$ 
for $i=0,1,\ldots,L$,
after shrinking $X$.
Recall that we have the truncations
$\nbigj(\vecm(i)):=\etabar_{\vecm(i)}(\nbigj)$.
Formally, we set $\nbigj(\vecm(-1)):=\{0\}$.
The set
$\nbigi^{\vecm(0)}_0:=\nbigj(\vecm(0))$
is a good set of
irregular values in the level 
$(\vecm(0),\gminih(0))$.
We have the naturally induced morphisms
$\etabar_{\vecm(i),\vecm(j)}:
 \nbigj(\vecm(j))\lrarr \nbigj(\vecm(i))$ for $j>i$.
For any $\gminia\in \nbigj(\vecm(i-1))$,
we define
\[
 \nbigi^{\vecm(i)}_{\gminia}:=
\etabar_{\vecm(i-1),\vecm(i)}^{-1}(\gminia),
\quad
 \nbigibar^{\vecm(i)}_{\gminia}:=
 \bigl\{
 \zeta_{\vecm(i)}(\gminib)\,\big|\,
 \gminib\in
 \nbigi^{\vecm(i)}_{\gminia}
 \bigr\}.
\]
Then, 
$\nbigi^{\vecm(i)}_{\gminia}$
are weakly good sets of irregular values
in the level $(\vecm(i),\gminih(i))$,
and $\nbigibar^{\vecm(i)}_{\gminia}$
are good sets of irregular values
in the level $(\vecm(i),\gminih(i))$.

\subsection{Good lattice in the level $\vecm$}
\label{subsection;10.5.10.2}
\label{subsection;08.9.4.40}

Let $Y$ be a complex manifold
with a simple normal crossing divisor $\nbigd_{Y}'$.
Let $\nbigk$ be a complex manifold
with a holomorphic function $\varrho$.
Let $\nbigx:=\Delta_z^k\times Y\times \nbigk$,
$\nbigd_{z,i}:=\{z_i=0\}$
and $\nbigd_z:=\bigcup_{i=1}^{k}\nbigd_{z,i}$.
We also put 
$\nbigd_{Y}:=
 \Delta_z^k\times \nbigd_{Y}'\times\nbigk$
and $\nbigd:=\nbigd_z\cup \nbigd_{Y}$.

Let $\vecm\in\seisuu_{<0}^k$,
and let $i(0)$ be an integer such that 
$1\leq i(0)\leq k$.
We put $\vecm(1):=\vecm+\vecdelta_{i(0)}$.
Let $E$ be a locally free $\nbigo_{\nbigx}$-module,
and let $\DD$ be 
a meromorphic flat $\varrho$-connection
of $E(\ast \nbigd)$.

\begin{df}
\label{df;10.5.6.1}
$(E,\DD)$ is called
a weakly good lattice
of a meromorphic $\varrho$-flat bundle
in the level $(\vecm,i(0))$,
if there exist
a weakly good set of irregular values $\nbigi$
in the level $(\vecm,i(0))$ on $(\nbigx,\nbigd)$,
and a decomposition
\begin{equation}
 \label{eq;07.12.15.1}
 (E,\DD)_{|\nbigdhat_z}
=\bigoplus_{\gminia\in\nbigi}
 (\Ehat_{\gminia},\DDhat_{\gminia})
\end{equation}
such that 
$\ord(\DDhat_{\gminia}
-d\gminia)\geq\vecm(1)$.

If $\nbigi$ is a good set of irregular values
in the level $(\vecm,i(0))$,
$(E,\DD,\nbigi)$
is called a good lattice in the level $(\vecm,i(0))$.
\hfill\qed
\end{df}
\index{good lattice in the level $\vecm$}

The decomposition (\ref{eq;07.12.15.1})
is called the irregular decomposition
in the level $(\vecm,i(0))$,
(or simply $\vecm$).
\index{irregular decomposition in the level $(\vecm,i(0))$}
In this situation,
we will often say that
$(E,\DD,\nbigi)$ is
a (weakly) good lattice
in the level $(\vecm,i(0))$.
 The rank of $\Ehat_{\gminia}$
 will be often denoted by $r(\gminia)$.
\index{rank $r(\gminia)$}
The following lemma is clear.
\begin{lem}
\label{lem;08.11.1.1}
Let $(E,\DD,\nbigi)$
be a good lattice
in the level $(\vecm,i(0))$.
For any $\gminia\in M(\nbigx,\nbigd)$,
we consider the line bundle
$\nbigl(\gminia)=\nbigo_{\nbigx}\, e$
with the meromorphic 
$\varrho$-connection
$\DD e=e\,(d\gminia)$.
We set
$\nbigi':=\bigl\{
 \gminib+\gminia\,\big|\,
 \gminib\in\nbigi
 \bigr\}$
and $(E',\DD^{\prime}):=
 (E,\DD)\otimes\nbigl(\gminia)$.
Then, $(E',\DD^{\prime},\nbigi')$
is a weakly good lattice in the level 
$(\vecm,i(0))$.

Conversely, let $(E,\DD,\nbigi)$ be
a weakly good lattice in the level $(\vecm,i(0))$.
Take any element $\gminia\in\nbigi$,
and we set $\nbigi':=\bigl\{
 \gminib-\gminia\,\big|\
 \gminib\in \nbigi
 \bigr\}$
and $(E',\DD^{\prime}):=
(E,\DD)\otimes\nbigl(-\gminia)$.
Then, $(E',\DD^{\prime},\nbigi')$
is a good lattice in the level $(\vecm,i(0))$.
\hfill\qed
\end{lem}

The following lemma can be shown
by the argument in the proof of
Proposition \ref{prop;10.5.3.5}.
\begin{lem}
The condition in Definition {\rm\ref{df;10.5.6.1}}
is equivalent to the following:
\begin{itemize}
\item For any $P\in \nbigd_z$,
 $(E,\DD)_{|\Phat}$
 has a decomposition
$(E,\DD)_{|\Phat}
=\bigoplus_{\gminia\in\nbigi}
 (\lefttop{P}E_{\gminia},\DDhat_{\gminia})$
such that
$\ord(\DDhat_{\gminia}
-d\gminia)\geq\vecm(1)$.
\hfill\qed
\end{itemize}
\end{lem}

We put $\nbigk^{0}:=\varrho^{-1}(0)$,
$\nbigx^0:=\nbigx\times_{\nbigk}\nbigk^0$, etc..
For simplicity, we assume $d\varrho$
is nowhere vanishing on $\nbigk^0$.
After shrinking $\nbigx$,
we have the irregular decomposition
$(E,\DD)_{|\nbigx^0}
=\bigoplus_{\gminia\in\nbigi}
 (E_{\gminia,\nbigx^0},\DD^{0}_{\gminia})$
such that
$\ord(\DD^0_{\gminia}-d\gminia)
 \geq \vecm(1)$.
It is uniquely extended to a decomposition
$(E,\DD)_{|\nbigxhat^0}
=\bigoplus_{\gminia\in\nbigi}
 (\Ehat_{\gminia,\nbigxhat^0},
 \DDhat_{\gminia})$
on the completion $\nbigxhat^0$.
We put $W:=\nbigx^0\cup \nbigd_z$.
By using Lemma \ref{lem;10.5.3.21},
we obtain a decomposition
\begin{equation}
\label{eq;07.12.15.2}
 (E,\DD)_{|\What}
=\bigoplus_{\gminia\in \nbigi}
 (\Ehat_{\gminia,\What},\DDhat_{\gminia})
\end{equation}
such that 
$\ord\bigl(
 \DDhat_{\gminia}-d\gminia
\bigr)\geq\vecm(1)$.
The decomposition (\ref{eq;07.12.15.2})
is also called the irregular decomposition
in the level $(\vecm,i(0))$.

\subsection{Residue}

Let $(E,\DD,\nbigi)$ be a good lattice
in the level $(\vecm,i(0))$.
Because
\[
 \DD E\subset
 E\otimes\Bigl(
 \vecz^{\vecm}
 \Omega^1_{\nbigx/\nbigk}(\log \nbigd_z)
+\Omega^1_{\nbigx/\nbigk}(\log \nbigd_Y)
 \Bigr)
\]
we obtain the residue
$\Res_{Y,j}(\DD)\in
 \End(E_{|\nbigd_{Y,j}})$
for each irreducible component 
$\nbigd_{Y,j}$ of $\nbigd_Y$
in a standard way.
We obtain the residue
even in the case that
$(E,\DD,\nbigi)$ is a weakly good lattice
in the level $(\vecm,i(0))$
by considering the tensor product
with a meromorphic $\varrho$-flat bundle
of rank one.

\begin{lem}
\label{lem;10.5.11.40}
Let $(\nbige,\DD)$ be a meromorphic
$\varrho$-flat bundle on $(\nbigx,\nbigd)$
with a good lattice $(E,\nbigi)$
in the level $(\vecm,i(0))$.
Let $P$ be any point of $\nbigd$
such that $\varrho(P)\neq 0$.
Then, we can find a good lattice
$(E',\nbigi)$ of $E(\ast \nbigd)$
in the level $(\vecm,i(0))$
on a small neighbourhood $\nbigx_P$
with the following non-resonance property:
\begin{itemize}
\item
Let $Q$ be any point of 
an irreducible component $\nbigd_{Y,j}\cap \nbigx_P$.
Then, 
distinct eigenvalues $\alpha$, $\beta$
of $\Res_{Y,j}(\varrho^{-1}\DD)_{|Q}$
satisfy $\alpha-\beta\not\in\seisuu$.
\end{itemize}
\end{lem}
\pf
It can be shown by the standard argument
in the proof of Proposition
\ref{prop;10.5.10.2}.
Because we give some more details there,
we omit it here.
\hfill\qed

\subsection{Some functoriality}
\label{subsection;08.1.2.1}

In general,
we use the symbol
$V_1^{\bot}$
to denote the subspace
$\bigl\{f\in V_2^{\lor}\,|\,
 f(V_1)=0\bigr\}\subset V_2^{\lor}$
for given vector spaces $V_1\subset V_2$,
where $V_2^{\lor}$ denotes the dual space
of $V_2$.
\index{vector space $V^{\bot}$}
\index{vector space $V^{\lor}$}
It is naturally extended 
in the case of vector bundles.
Let $(E,\DD,\nbigi)$ be 
a (weakly) good lattice 
in the level $(\vecm,i(0))$.
We set 
$\nbigi^{\lor}:=\bigl\{
 -\gminia\,\big|\,\gminia\in\nbigi
 \bigr\}$.
Then, the dual 
$(E^{\lor},\DD^{\lor},\nbigi^{\lor})$
is also a (weakly) good lattice
in the level $(\vecm,i(0))$.
The direct summands
$\Ehat^{\lor}_{\gminia}$  $(\gminia\in\nbigi^{\lor})$
in the irregular decomposition are given as follows:
\[
 \Ehat^{\lor}_{\gminia}=
\left(
\bigoplus_{\substack{
 \gminib\in\nbigi \\
 \gminib\neq-\gminia}}
\Ehat_{\gminib}
 \right)^{\bot}
\]

\vspace{.1in}

Let $(E_p,\DD_p,\nbigi_p)$ $(p=1,2)$ be 
(weakly) good lattices 
in the level $(\vecm,i(0))$.
We put $\nbigi_1\otimes\nbigi_2:=
 \bigl\{
 \gminia_1+\gminia_2\,\big|\,\gminia_p\in\nbigi_p
 \bigr\}$.
If $\nbigi_1\otimes\nbigi_2$
is a (weakly) good set of irregular values 
in the level $(\vecm,i(0))$,
then $E_1\otimes E_2$ is 
a (weakly) good lattice in the level $(\vecm,i(0))$.
The direct summands of the irregular decomposition
are given as follows:
\[
 \widehat{(E_1\otimes E_2)}_{\gminia}
=\bigoplus_{\substack{
 (\gminia_1,\gminia_2)\in
 \nbigi_1\times \nbigi_2\\
 \gminia_1+\gminia_2=\gminia}}
 \Ehat_{1,\gminia_1}
\otimes
 \Ehat_{2,\gminia_2}
\]

We put $\nbigi_1\oplus\nbigi_2:=\nbigi_1\cup\nbigi_2$.
If $\nbigi_1\oplus\nbigi_2$ is 
a (weakly) good set of irregular values
in the level $(\vecm,i(0))$,
the direct sum $E_1\oplus E_2$ is also
a (weakly) good lattice 
in the level $(\vecm,i(0))$.
The direct summands in the irregular decomposition
are given as follows:
\[
 (E_1\oplus E_2)_{\gminia}
=E_{1,\gminia}\oplus E_{2,\gminia}
\]

A morphism 
$f:(E_1,\DD_1)
 \lrarr (E_2,\DD_2)$
of (weakly) good lattices in the level 
$(\vecm,i(0))$
is defined to be just a flat morphism.
Note that the induced morphism
$\fhat:(E_1,\DD_1)_{|\What}
\lrarr (E_2,\DD_2)_{|\What}$
preserves the irregular decomposition.

\subsection{Remark on growth order
of a flat section}
\label{subsection;08.9.28.119}

Let $(E,\DD)$ be 
a good lattice of $(\nbige,\DD)$
in the level $(\vecm,i(0))$ on $(\nbigx,\nbigd)$.
Let $\vecv$ be a frame of $E$.
We have the matrix-valued functions $A_i$
determined by 
$\DD\vecv=\vecv\,\Bigl(
 \sum_{i=1}^nA_i\, dz_i \Bigr)$.
We have the following:
\begin{itemize}
\item
$A_i=O\bigl(|z_i^{-1}|\, |\vecz^{\vecm}|\bigr)$
for $i=1,\ldots,k$.
\item
$A_i=O\bigl(|\vecz^{\vecm}|\bigr)
 +O(|z_i|^{-1})$ for $i=k+1,\ldots,\ell$.
\item
 $A_i=O\bigl(|\vecz^{\vecm}|\bigr)$
 for $i=\ell+1,\ldots,n$.
\end{itemize}

Assume $\varrho$ is nowhere vanishing on $\nbigk$,
for simplicity.
Let $S$ be a small multi-sector of $\nbigx-\nbigd_{z}$,
and let $f$ be a $\DD$-flat
section of $E_{|S}$, which is nowhere vanishing.
We have the expression
$f=\sum f_i\, v_i$.
We obtain a $\cnum^r$-valued function $\vecf=(f_i)$
on $S$.
\begin{lem}
\label{lem;07.12.27.10}
The following holds for some $C>0$:
\[
\Bigl| \log|\vecf|\Bigr|
\leq
 C\,|\vecz^{\vecm}|
+C\sum_{i=k+1}^{\ell}\log|z_i|^{-1}
\]
\end{lem}
\pf
It follows from Lemma \ref{lem;07.12.27.5}.
\hfill\qed

\subsection{The induced good lattice
 in the level $\vecm$}

Let $\nbigx=\Delta^n\times \nbigk$,
$\nbigd_i=\{z_i=0\}$
and $\nbigd=\bigcup_{i=1}^{\ell}\nbigd_i$.
Let $(E,\DD)$ be 
an unramifiedly good lattice of
a meromorphic $\varrho$-flat bundle $(\nbigx,\nbigd)$.
For simplicity, we assume that 
the coordinate system is admissible
for the good set $\Irr(\DD)$.
We take an auxiliary sequence 
$\vecm(0),\vecm(1),\ldots,\vecm(L)$
for $\Irr(\DD)$.
Let $\Irrbar(\DD,\vecm(i))$
denote the image of
$\etabar_{\vecm(i)}:
 \Irr(\DD)\lrarr M(\nbigx,\nbigd)/H(\nbigx)$.

\begin{lem}
\label{lem;07.6.17.1}
$\bigl(E,\DD,
\Irrbar(\DD,\vecm(0))\bigr)$ is 
a good lattice 
in the level $(\vecm(0),\gminih(0))$.
The decomposition
is given by the irregular decomposition
in the level $(\vecm(0),\gminih(0))$:
\[
 (E,\DD)_{|\nbigdhat(\kbar(0))}
=\bigoplus_{\gminia\in
 \Irrbar(\DD,\vecm(0))}
 (\Ehat_{\gminia},\DDhat_{\gminia})
\]
Here $k(0)$ is determined by
$\vecm(0)\in\seisuu_{<0}^{k(0)}
 \times\veczero_{\ell-k(0)}$,
and $\kbar(0)=\{1,\ldots,k(0)\}$.
\end{lem}
\pf
Let $\nbigd_{\ellsitabar}:=
 \bigcap_{j=1}^{\ell}\nbigd_j$.
Let $\vecuhat$ be a frame of $E_{|\nbigdhat_{\ellsitabar}}$
which is compatible with the irregular decomposition
$(E,\DD)_{|\nbigdhat_{\ellsitabar}}
=\bigoplus_{\gminic\in\Irr(\DD)}
 (\lefttop{\ellsitabar}\Ehat_{\gminic},
 \DD_{\gminic})$.
The connection $1$-form of $\DD$
with respect to the frame $\vecuhat$
is decomposed as 
$\DD\vecuhat=\vecuhat
 \,\left(
 \bigoplus_{\gminic\in\Irr(\DD)}
 \bigl(d\gminic\, I_{\gminic}+T_{\gminic}\bigr)
 \right)$,
where $T_{\gminic}$ are logarithmic $1$-forms,
and $I_{\gminic}$ are the identity matrices.
For each $\gminia\in\Irrbar(\DD,\vecm(0))$,
we put $\lefttop{\ellsitabar}\Ehat_{\gminia}:=
 \bigoplus_{\gminic\in\etabar_{\vecm(0)}^{-1}(\gminia)}
 \lefttop{\ellsitabar}\Ehat_{\gminic}$.
Let $\vecuhat_{\gminia}$ denote the tuple
of $\uhat_{i}\in
 \lefttop{\ellsitabar}\Ehat_{\gminia}$,
which gives a frame of 
$\lefttop{\ellsitabar}\Ehat_{\gminia}$.
We put $\nbigd'=\bigcup_{k(0)<i\leq \ell}\nbigd_i$.
Then, we have
$\DD\vecuhat_{\gminia}
=\vecuhat_{\gminia}
 \,(d\gminia\, I_{\gminia}+T'_{\gminia})$,
where
$T'_{\gminia}\in
 \vecz^{\vecm(1)}\Omega^1_{\nbigdhat_{\ellsitabar}}
 (\log \nbigd(\kbar(0)))
+\Omega^1_{\nbigdhat_{\ellsitabar}}
 (\log \nbigd')$,
and $I_{\gminia}$ are the identity matrices.
Then, the claim of the lemma follows.
\hfill\qed

\vspace{.1in}

Let $\vecm(i)$ be the minimum of
$\nbigt(\Irr(\DD))$,
i.e., $\bigl|\Irr(\DD,\vecm(j))\bigr|=1$
for $j<i$,
and $\bigl|\Irr(\DD,\vecm(i))\bigr|\geq 2$.
Note that
$\Irrbar(\DD,\vecm(i))$
is a weakly good set
in the level $(\vecm(i),\gminih(i))$.

\begin{lem}
$\bigl(E,\DD,
 \Irrbar(\DD,\vecm(i))
 \bigr)$ is weakly good 
in the level $(\vecm(i),\gminih(i))$.
\end{lem}
\pf
Take any
$\gminia\in \Irr(\DD)$.
We consider a line bundle
$\nbigl(-\gminia):=\nbigo_{\nbigx}\, e$
with the meromorphic $\varrho$-flat connection
$\DD e=e\,(-d\gminia)$.
Then, $(E',\DD^{\prime}):=
(E,\DD)\otimes\nbigl(-\gminia)$
is an unramifiedly good lattice with
$\Irr(\DD^{\prime})=\bigl\{
 \gminib-\gminia\,\big|\,
 \gminib\in \Irr(\DD)
 \bigr\}$.
Note that
$\vecm(i),\vecm(i+1),\ldots,\vecm(L)$
gives an auxiliary sequence
for $\Irr(\DD^{\prime})$.
Applying Lemma \ref{lem;07.6.17.1}
to $(E',\DD^{\prime})$,
we obtain that 
$(E',\DD^{\prime})$ 
is a good lattice
in the level $(\vecm(i),\gminih(i))$.
Then, the claim of the lemma follows.
\hfill\qed

\section{Good Deligne-Malgrange lattice
and Deligne-Malgrange lattice}
\label{section;10.5.4.100}

In this section,
we consider ordinary flat connections,
except the remark in
Subsection \ref{subsection;10.5.8.20}.
We use the symbols
$X$, $D$ and $\nabla$
instead of $\nbigx$, $\nbigd$ and $\DD$,
respectively.

\subsection{Good Deligne-Malgrange lattice}
\label{subsection;08.2.22.5}

\index{Deligne-Malgrange}
\index{good Deligne-Malgrange lattice}
\index{unramifiedly good Deligne-Malgrange lattice}

Let $X$ be a complex manifold,
and let $D=\bigcup_{i\in \Lambda}D_i$
be a simple normal crossing hypersurface.

\begin{df}
Let $E$ be an unramifiedly good lattice 
of a meromorphic flat bundle
$(\nbige,\nabla)$ on $(X,D)$.
It is called unramifiedly good Deligne-Malgrange,
if any eigenvalues $\alpha$ of $\Res_{D_i}(\nabla)$
$(i\in\Lambda)$
satisfy $0\leq \Re(\alpha)<1$.
\hfill\qed
\end{df}

Let $E$ be an unramifiedly
good Deligne-Malgrange lattice
of a meromorphic flat bundle $(\nbige,\nabla)$ on $(X,D)$.
We have the generalized eigen decomposition
$E_{|D_i}=\bigoplus_{\alpha\in\cnum}
 \lefttop{i}\EE_{\alpha}$
with respect to $\Res_{D_i}(\nabla)$.
Then, we have the parabolic filtration 
$\lefttop{i}F$ of $E_{|D_i}$
for each irreducible component $D_i$ of $D$,
given in a standard manner:
\[
 \lefttop{i}F_a(E_{|D_i}):=
\bigoplus_{-\Re(\alpha)\leq a}
\EE_{\alpha},
\quad\quad
 (-1<a\leq 0)
\]
It is easy to observe that
the parabolic filtrations are compatible.
The associated filtered bundle
is denoted by $\vecE^{DM}_{\ast}$,
which is called the Deligne-Malgrange filtered flat bundle
associated to $(\nbige,\nabla)$.
We have $E=\prolong{E}^{DM}$,
i.e., $E$ is the $\veczero$-truncation of 
$\vecE^{DM}_{\ast}$.
\index{Deligne-Malgrange filtered flat bundle
$\vecE^{DM}_{\ast}$}

\begin{lem}
\label{lem;08.2.22.3}
Let $\Xtilde$ and $X$ be complex manifolds,
and let $\Dtilde$ and $D$ be simple normal crossing
hypersurfaces of $\Xtilde$ and $X$, respectively.
Let $\varphi:(\Xtilde,\Dtilde)\lrarr (X,D)$ be a ramified covering.
Let $(\nbige,\nabla)$ be a meromorphic flat bundle
on $(X,D)$, which has the unramifiedly good
Deligne-Malgrange lattice $E$.
Then, $\varphi^{\ast}(\nbige,\nabla)$
also has the unramifiedly good Deligne-Malgrange lattice
$\Etilde$,
and $E$ is the descent of $\Etilde$.
\end{lem}
\pf
We obtain the filtered bundle $\vecEtilde_{\ast}$
on $(\Xtilde,\Dtilde)$
induced by $\varphi$ and $\vecE_{\ast}^{DM}$,
as in Section \ref{subsection;07.11.5.60}.
We can observe that
$\prolong{\Etilde}$ is the unramifiedly
good Deligne-Malgrange lattice
of $\varphi^{\ast}\nbige$.
Then, the claims of the lemma are clear.
\hfill\qed

\begin{df}
Let $E$ be a good lattice of
a meromorphic flat bundle $(\nbige,\nabla)$
on $(X,D)$.
\begin{itemize}
\item
$E$ is called good Deligne-Malgrange at $P\in D$,
if there exists a ramified covering
$\varphi_e:(X_P',D_P')\lrarr (X_P,D_P)$
such that $E_{|X_P}$ is the descent of 
the unramifiedly good Deligne-Malgrange lattice
of $\varphi_e^{\ast}\nbige$.
Here, $(X_P,D_P)$, 
$(X_P',D_P')$ and $\varphi_e$
are as in Definition \ref{df;10.5.3.2}.
Note that we also obtain the natural filtered
flat bundle $(\vecE_{\ast}^{DM},\nabla)$
on $(X_P,D_P)$ as the decent of
the Deligne-Malgrange filtered flat bundle
associated to $\varphi_e^{\ast}\nbige$
in this case.
\item
$E$ is called good Deligne-Malgrange,
if it is good Deligne-Malgrange
at any point $P\in D$.
In this case,
we have the associated filtered flat bundle
$(\vecE^{DM}_{\ast},\nabla)$,
which is called the Deligne-Malgrange
filtered flat bundle associated to $(\nbige,\nabla)$.
\hfill\qed
\end{itemize}
\end{df}
\index{good Deligne-Malgrange lattice}
\index{Deligne-Malgrange filtered bundle}

We should remark that good Deligne-Malgrange lattice
does not necessarily exist 
for a given meromorphic flat bundle
over a higher dimensional variety,
contrast to the one dimensional case.
But, it is unique, if it exists,
which follows from the uniqueness of
(not necessarily good) Deligne-Malgrange lattice 
(see Section \ref{subsection;08.2.22.2} below)
or Lemma \ref{lem;08.2.22.3}.

\begin{prop}
\label{prop;10.5.8.10}
Let $(\nbige,\nabla)$ be a meromorphic flat bundle.
The following conditions are equivalent:
\begin{description}
\item[(a)]
$(\nbige,\nabla)$ has a good Deligne-Malgrange lattice.
\item[(b)]
$(\nbige,\nabla)$ has a good lattice.
\item[(c)]
For each point $P\in D$,
there exists a neighbourhood $U$ of $P$
such that 
$(\nbige,\nabla)_{|U}$ has a good lattice.
\item[(d)]
For each point $P\in D$,
there exists a neighbourhood $U$ of $P$
such that 
$(\nbige,\nabla)_{|U}$ has 
a good Deligne-Malgrange lattice.
\end{description}
\end{prop}
\pf
The implications
$(a)\Longrightarrow (b)\Longrightarrow (c)$
are obvious.
The implication $(d)\Longrightarrow(a)$
follows from the uniqueness of the Deligne-Malgrange lattice.
Let us show the implication $(c)\Longrightarrow (d)$,
which can be carried out
using a standard successive
use of elementary transform.
We give only an outline.
Due to the uniqueness of 
the good Deligne-Malgrange lattice,
we have only to consider the problem
in the unramified case.
Let $E$ be an unramifiedly good lattice.
Let $\Sp(i)$ denote the set of the eigenvalues
of $\Res_{D_i}(\nabla)$.
We have the generalized eigen decomposition
$E_{|D_i}=
 \bigoplus_{\alpha\in\Sp(i)}
 \EE_{\alpha}$
with respect to $\Res_{D_i}(\nabla)$.
We put $a_+:=\max\bigl\{\Re(\alpha)\,\big|\,
 \alpha\in \Sp(i)\bigr\}$
and
$a_-:=\min\bigl\{\Re(\alpha)\,\big|\,
 \alpha\in \Sp(i)\bigr\}$.
We have the subbundles 
$F_{\pm}:=
 \bigoplus_{\Re(\alpha)=a_{\pm}}\EE_{\alpha}$.

Assume $a_-<0$.
We regard $F_-$ as an $\nbigo_X$-module.
Let $E'$ denote the kernel of
the naturally defined morphism
of $\nbigo_X$-modules $E\lrarr F_-$.
It is easy to show that
$E'$ is also good,
the eigenvalues $\alpha$ of $\Res_{D_i}(\nabla')$
satisfy $a_-<\Re(\alpha)\leq \max(a_+,1+a_-)$.
Assume $a_+>1$.
Let $E''$ denote the kernel of
the naturally defined morphism
$E\otimes\nbigo(D_i)\lrarr F_{+}\otimes\nbigo(D_i)$.
It is easy to show that
$E''$ is also good,
and the eigenvalues of $\Res_{D_i}(\nabla'')$
satisfy $\min(a_-,-1+a_+)\leq \Re(\alpha)<a_+$.
Hence, by composition of the above procedure,
we can construct the unramifiedly good Deligne-Malgrange lattice
from an unramifiedly good lattice.
\hfill\qed

\subsubsection{Local existence 
of non-resonance lattice in the family case}
\label{subsection;10.5.8.20}

In this subsection,
we consider the case $\nbigk$ is not necessarily
a point.
We put $\nbigx:=X\times\nbigk$
and $\nbigd:=D\times\nbigk$.
Let $(\nbige,\DD)$ be an unramifiedly
good meromorphic 
$\varrho$-flat bundle on $(\nbigx,\nbigd)$
relative to $\nbigk$.
\begin{prop}
\label{prop;10.5.10.2}
If $\varrho(P)\neq 0$,
there exist a small neighbourhood 
$\nbigx_P$ of $P$
and a good lattice $F$ of 
$(\nbige,\DD)_{|\nbigd_P}$
with the following property:
\begin{itemize}
\item
 Let $\alpha$ and $\beta$ be
 distinct eigenvalues
 of $\Res_i(\varrho^{-1}\DD)_Q
 \in \End(F_{|\nbigd_{P,i}})$ 
 for some $Q\in \nbigd_{P,i}$.
 Then, $\alpha-\beta\not\in\seisuu$. 
\end{itemize}
\end{prop}
\pf
It can be shown by the argument
in the proof of Proposition \ref{prop;10.5.8.10}.
We give only an outline.
We set $\nabla:=\varrho^{-1}\DD$.
We may assume $\nbigx=\Delta^n\times\nbigk$,
$\nbigd=\bigcup_{i=1}^{\ell}\{z_i=0\}$
and $P=(O,y)$, where $O$ denotes the
origin of $\Delta^n$ and $y\in\nbigk$.
Let $E$ be an unramifiedly good lattice.
Let $\Sp(i)$ denote the set of the eigenvalues
of the endomorphism
$\Res_i(\nabla)$ on $E^y_{|\nbigd_i^y}$.
We have the generalized eigen decomposition
$E^y_{|\nbigd^y_i}=\bigoplus_{\alpha\in\Sp(i)}
 \EE_{\alpha}$
with respect to $\Res_i(\nabla)$.
We put $a_-(E,i):=\min\bigl\{
 \Re(\alpha)\,\big|\,\alpha\in\Sp(i)
 \bigr\}$
and $a_+(E,i):=\max\bigl\{
 \Re(\alpha)\,\big|\,\alpha\in\Sp(i)
 \bigr\}$.
We have the subbundle
$F_-:=\bigoplus_{\Re(\alpha)=a_-(E,i)}
 \EE_{\alpha}$.
If $\nbigx_P$ is sufficiently small,
we have the subbundle
$F'_-\subset E_{|\nbigd_i}$
such that
(i) $F'_{-|\nbigd_i^y}=F_-$,
(ii) $\Res_{\nbigd_i}(\nabla)(F'_-)\subset F'_-$.
Applying the procedure
in the proof of
Proposition \ref{prop;10.5.8.10}
to $E$ and $F'_-$,
we obtain a lattice $E'$
such that $a_-(E,i)<a_-(E',i)\leq a_-(E,i)+1$.
By successive use of this procedure,
we may assume $0\leq a_-(E,i)$.
Similarly, we may assume $a_+(E,i)< 1$.
Then, if $\nbigx_P$ is sufficiently small,
$\alpha-\beta\not\in\seisuu$
for distinct eigenvalues of
$\Res_i(\nabla)_{Q}$
$(Q\in D_{i,P})$.
\hfill\qed

\subsection{Deligne-Malgrange filtered flat sheaf}
\label{subsection;08.2.22.2}

\subsubsection{Deligne-Malgrange lattice}
Let $X$ be a complex manifold,
and let $D=\bigcup_{i\in \Lambda}D_i$ be 
a simple normal crossing divisor of $X$.
As already remarked,
a meromorphic flat connection $(\nbige,\nabla)$
does not necessarily have a good lattice.
However,
according to Malgrange's theorem,
we have a lattice which is generically
good Deligne-Malgrange lattice.

\begin{prop}[Malgrange \cite{malgrange}]
There always exists
a unique lattice $E\subset\nbige$
characterized by the following property:
\begin{itemize}
\item
 $E$ is a coherent reflexive $\nbigo_X$-module.
\item
 There exists a Zariski closed subset $Z$ of $D$
 with $\codim_X(Z)\geq 2$,
 such that $E_{|X\setminus Z}$
 is the good Deligne-Malgrange
 lattice of $(\nbige,\nabla)_{|X\setminus Z}$.
\end{itemize}
It is called the canonical lattice in \cite{malgrange},
but we called it
Deligne-Malgrange lattice.
\index{Deligne-Malgrange lattice}
\index{canonical lattice}
\end{prop}

\begin{rem}
In this monograph,
a subset $Z\subset X$ is called Zariski closed,
if it is a closed complex analytic subset.
And, a subset $U\subset X$ is called Zariski open,
if it is the complement of 
a Zariski closed set.
\hfill\qed
\end{rem}
\index{Zariski closed}
\index{Zariski open}

\pf
We give only a remark
that $Z$ can be Zariski closed in the above sense.
Let $N(E)$ denote the closed analytic subset of $D$
such that
$Q\in N(E)$ if and only if
$E_Q$ is not locally free.
Let $D^{[2]}$ denote the set of singular points
of $D$.
In \cite{malgrange},
it is shown that there exists
a closed subset $Z\subset D$
in the ordinary topology,
such that
(i) $D^{[2]}\cup N(E)\subset Z$,
(ii) $E$ is good Deligne-Malgrange around 
 any $Q\in D\setminus Z$,
(iii) for any $Q\in Z$,
 there exists a small neighbourhood $X_Q$ of $Q$
 in $X$ with a closed analytic subset
 $\Ztilde_Q$ of $D_Q:=X_Q\cap D$
 satisfying
 $\codim_{X_Q}(\Ztilde_Q)\geq 2$ and
 $X_Q\cap Z\subset\Ztilde_Q$.
If the closed subsets $Z_i$ $(i\in \Lambda)$ have
the above property,
$\bigcap_{i\in \Lambda}Z_i$ also has it.
Hence, we have the minimum among
the closed subsets with the above property,
which will be denoted by $Z$ in the following.
Then, $Q\in D$ is contained in $Z$
if and only if
one of the following holds:
(i) $Q\in D^{[2]}\cup N(E)$,
(ii) $Q\not\in D^{[2]}\cup N(E)$
and $E$ is not good Deligne-Malgrange 
at $Q$.

Let us show that $Z$ is closed analytic.
For any $P\in Z$,
we can take a small neighbourhood $X_P$
and a closed analytic subset
$\Ztilde_P\subset D_P:=D\cap X_P$
satisfying
$\codim_{X_P}(\Ztilde_P)\geq 2$
and
$Z_P=X_P\cap Z\subset\Ztilde_P$.
We put $A_P:=\Ztilde_P\setminus Z_P$,
which is an open subset of $\Ztilde_P$.
We put $D_P^{[2]}:=D^{[2]}\cap X_P$
and $N(E)_P:=N(E)\cap X_P$.
Let $\Ztilde_P=\bigcup_{i\in\Gamma}\Ztilde_{P,i}$
be the irreducible decomposition.
For each $i\in\Gamma$,
we put
$\Wtilde_{P,i}:=
 D^{[2]}_P\cup N(E)_P\cup
 \bigcup_{j\neq i}\Ztilde_{P,j}$.
If $\Ztilde_{P,i}\setminus A_P\subset
 \Wtilde_{P,i}$,
we have $Z_P\subset \Wtilde_{P,i}$.
Hence, we may and will assume 
$\Ztilde_{P,i}\setminus 
 (A_P\cup\Wtilde_{P,i})\neq\emptyset$
for each $i\in\Gamma$.
Then, let us show that $A_P=\emptyset$,
i.e., $Z_P=\Ztilde_{P,i}$,
which implies that 
$Z$ is closed analytic subset of $X$.
For that purpose,
we have only to show that
$\Ztilde_{P,i}\cap A_P=\emptyset$
for each $i\in\Gamma$.
Assume the contrary,
and we shall deduce a contradiction.

Let $\Ztilde_{P,i}^{\star}$ denote
the smooth part of
$\Ztilde_{P,i}\setminus \Wtilde_{P,i}$.
Because $\Ztilde_{P,i}$ is irreducible,
$\Ztilde_{P,i}^{\star}$ is connected 
and non-empty.
Because 
$A_P\cap\Ztilde_{P,i}\neq\emptyset$,
we have 
$A_P\cap\Ztilde_{P,i}^{\star}\neq\emptyset$.
We have the two cases:
(A) $Z_P\cap\Ztilde^{\star}_{P,i}\neq\emptyset$,
(B) $Z_P\cap \Ztilde^{\star}_{P,i}=\emptyset$.
In the case (B),
note that $Z_P\setminus \Wtilde_{P,i}$
is contained in a closed analytic subset
whose codimension in $X$ is larger than $3$.

We take a point $Q\in Z_P\cap\Ztilde^{\star}_{P,i}$
as follows.
In the case (A), 
$Q$ is a point in the intersection of
$Z_P\cap\Ztilde^{\star}_{P,i}$
and the closure of $A_P\cap\Ztilde^{\star}_{P,i}$
in $\Ztilde^{\star}_{P,i}$.
In the case (B),
$Q$ is any point of
$Z_P\setminus \Wtilde_{P,i}$.

We take a small coordinate neighbourhood
$(X_Q,z_1,\ldots,z_n)$ around $Q$
such that $D_Q=\{z_1=0\}$.
We put $Z_Q:=Z_P\cap X_Q$.
In the case (A), 
we may assume that 
$Z_Q$ is the complement
of a non-empty open subset in 
$\{z_1=z_2=0\}$.
In the case (B),
$Z_Q$ is contained in a closed analytic
subset $Z_Q'$ with
$\codim_{X_Q} Z_Q'\geq 3$.
By our choice,
$E$ is good Deligne-Malgrange
around any $Q'\in D_Q\setminus Z_Q$.
In this situation, we shall show that
$E$ is good Deligne-Malgrange around $Q$,
which contradicts with the choice of $Q$,
and we can conclude that 
$A_P\cap\Ztilde_{P,i}=\emptyset$.

We have only to consider the case that
$E$ is unramifiedly good
around any $Q'\in D_Q\setminus Z_Q$.
We recall that a holomorphic function on 
$D_Q\setminus Z_Q$ is naturally extended
to a holomorphic function on $D_Q$.
We recall that,
for a holomorphic function $f$ on $D_Q$,
if the zero of $f$ is contained in $Z_Q$,
then $f$ is actually nowhere vanishing.
We also remark that
the fundamental group of
$D_Q\setminus Z_Q$ is trivial.
Hence, we obtain that
there exists a good set of irregular values
$\nbigi\subset
 M(X_Q,D_Q)/H(X_Q)$
such that, for any $Q'\in D_Q\setminus Z_Q$,
the restriction of $\nbigi$ to a neighbourhood
of $Q'$ is $\Irr(\nabla,Q')$.
Then, we can show that 
$(E,\nabla)_{|\Dhat_Q}$ has 
the irregular decomposition
by using a standard argument in \cite{levelt}.
We give only an indication.

Let $F$ be a locally free
$\nbigo_{\Dhat_Q}$-module
with a meromorphic connection $\nabla$.
Let $\nbigi\subset M(X_Q,D_Q)$ be 
a good set of irregular values.
Assume that, for any $Q'\in D_Q\setminus Z_Q$,
we have a decomposition
$(F,\nabla)_{|\Dhat_{Q'}}
=\bigoplus_{\gminia\in\nbigi}
 (F_{\gminia,Q'},\nabla_{\gminia})$
such that $\nabla_{\gminia}$ are 
$\gminia$-logarithmic,
where $D_{Q'}:=D\cap X_{Q'}$
for a small neighbourhood $X_{Q'}$ of $Q'$.
Then, we shall show that
there exists a decomposition
$(F,\nabla)=
 \bigoplus_{\gminia\in\nbigi}
 (F_{\gminia},\nabla_{\gminia})$
such that $\nabla_{\gminia}$ are
$\gminia$-logarithmic.
(Then, we obtain the desired decomposition
of $(E,\nabla)_{|\Dhat_Q}$.)
We use an inductive argument.

Let $m(F,\nabla):=
\min\{\ord_{z_1}\gminia\,|\,\gminia\in\nbigi\}$.
If $m(F,\nabla)=0$, there is nothing to do.
Let us consider the case $m(F,\nabla)=m$.
We set
$T:=\bigl\{
 -m(z_1^{-m}\gminia)_{|D_Q}\,\big|\,
 \gminia\in \nbigi\bigr\}$.
If $|T|=1$, by considering
$\nabla+d\bigl(z_1^{m}\alpha\bigr)/m$
for $\alpha\in T$,
we can reduce the case $m(F,\nabla)=m+1$.
If $|T|\geq 2$,
we consider the endomorphism
$G$ of $F_{|D_Q}$ induced by
$z_1^{-m}\nabla(z_1\del_1)$.
The set of the eigenvalues of 
$G_{|D_Q\setminus Z_Q}$
is given by $T$,
and we have the eigen decomposition
$F_{|D_Q\setminus Z_Q}=
\bigoplus_{\alpha\in T} \EE_{\alpha}$.
Then, the set of the eigenvalues of
$G$ is $T$,
and we have the eigen decomposition
$F_{|D_Q}=
\bigoplus_{\alpha\in T} \EE_{\alpha}$
on $D_Q$.
By a standard argument explained
in Lemma \ref{lem;10.6.23.10},
it is extended to a decomposition
$F=\bigoplus_{\alpha\in T}
 F_{\alpha}$ such that 
(i) $z_1^{-m}\nabla(z_1\del_1)F_{\alpha}
\subset F_{\alpha}$,
(ii) $F_{\alpha|D_Q}=\EE_{\alpha}$.
If we restrict it to a small neighbourhood of
$Q'\in D_Q\setminus Z_Q$,
it is the same as the irregular decomposition
around $Q'$.
Hence, it is $\nabla$-flat.
We have the natural map
$\pi:\nbigi\lrarr T$,
and let $\nbigi(\alpha):=\pi^{-1}(\alpha)$.
For each $Q'\in D_Q\setminus Z_Q$,
we have a decomposition
$(F_{\alpha},\nabla)_{|\Dhat_{Q'}}
=\bigoplus_{\gminia\in\nbigi(\alpha)}
 (F_{\gminia,Q'},\nabla_{\gminia})$.
Hence, we can obtain the desired decomposition
of $(F,\nabla)$ by an easy inductive argument.
\hfill\qed

\subsubsection{Deligne-Malgrange filtered sheaf}

We have the Deligne-Malgrange filtered flat bundle
$(\vecE_{X\setminus Z\,\ast}^{DM},\nabla)$ 
on $(X\setminus Z,D\setminus Z)$
associated to $(\nbige,\nabla)_{|X\setminus Z}$.

\begin{lem}
It is extended to a filtered flat sheaf on $(X,D)$,
i.e., we have the filtered flat sheaf
$(\vecE^{DM}_{\ast},\nabla)$ on $(X,D)$
with the following property:
\begin{itemize}
\item
 $(\vecE^{DM}_{\ast},\nabla)_{|X\setminus Z}
=(\vecE^{DM}_{X\setminus Z\,\ast},\nabla)$.
\item
 $\prolongg{\veca}{E}^{DM}$
 are coherent reflexive $\nbigo_X$-modules
 for any $\veca\in\real^{\Lambda}$.
\end{itemize}
It is called the Deligne-Malgrange filtered sheaf
associated to $(\nbige,\nabla)$.
\end{lem}
\pf
Since we have only to shift the condition
on the eigenvalues of the residues,
the claim can be shown by 
repeating the argument of Malgrange.
Otherwise,
it can be reduced to the existence 
of Deligne-Malgrange lattice,
which is explained in the following.
Since the problem is local,
we may assume
$X=\Delta^n$ and 
$D=\bigcup_{i=1}^{\ell}\{z_i=0\}$.
For each $\veca=(a_i)\in\real^{\ell}$,
let us consider the line bundle
$L_{\veca}:=\nbigo_X\,e$
with the logarithmic flat connection
$\nabla_{\veca}$ such that
$\nabla_{\veca} e=e\,
 \sum_{i=1}^{\ell} a_i\,dz_i/z_i$.
We have the Deligne-Malgrange lattice
$\prolongg{\veca}{E'}$
of $(\nbige,\nabla)\otimes (L_{\veca},\nabla_{\veca})$,
and we put 
$\prolongg{\veca}{E}:=
 \prolongg{\veca}{E'}\otimes L_{-\veca}$.
Then, $\vecE_{\ast}=
 \bigl(\prolongg{\veca}{E}\,\big|\, 
\veca\in\real^{\Lambda}\bigr)$
has the desired property.
\hfill\qed

\vspace{.1in}
We mention an important property
of Deligne-Malgrange filtered sheaf
on a projective variety.
\begin{lem}
Assume that $X$ is projective,
provided with an ample line bundle $L$.
Then, $\mu_L(\vecE_{\ast}^{DM})=0$ always hold.
(See Section {\rm\ref{subsection;08.9.29.19}}
for $\mu_L$.)
\end{lem}
\pf
The claim can be reduced to the one dimensional case,
which was shown in \cite{sabbah3}, for example.
\hfill\qed

\subsection{Good formal structure and good lattice}

Let $X$ be a complex manifold
with a normal crossing hypersurface $D$.
Let $(\nbige,\nabla)$ be a meromorphic flat bundle
on $(X,D)$.

\begin{prop}
\label{prop;10.5.4.3}
If $(\nbige,\nabla)_{|\Phat}$
has an (unramifiedly)
good Deligne-Malgrange lattice
for each $P\in D$,
then $(\nbige,\nabla)$
has an (unramifiedly)
good Deligne-Malgrange lattice.
\end{prop}

The case $D$ is smooth will be argued in
Section \ref{subsection;10.5.4.60},
and the normal crossing case will be argued
in Section \ref{subsection;10.5.4.61}.
We remark that we have only to consider
the unramified case,
according to Lemma \ref{lem;10.6.10.1}.

Before going into the proof of the proposition,
we give a consequence,
which will be used
in Section \ref{subsection;07.10.14.111}
for the proof of Theorem \ref{thm;07.10.14.60}.

\begin{cor}
\label{cor;10.5.4.101}
Let 
$0\lrarr\bigl(\nbige^{(1)},\nabla^{(1)}\bigr)\lrarr
\bigl(\nbige^{(0)},\nabla^{(0)}\bigr) \lrarr
 \bigl(\nbige^{(2)},\nabla^{(2)}\bigr)\lrarr 0$
be an exact sequence of meromorphic
flat bundles on $(X,D)$.
If $(\nbige^{(i)},\nabla^{(i)})$ $(i=1,2)$
have good Deligne-Malgrange lattices,
then $(\nbige^{(0)},\nabla^{(0)})$
also has a good Deligne-Malgrange lattice.
\end{cor}
\pf
It immediately follows from
Proposition \ref{prop;07.9.26.105}
and Proposition \ref{prop;10.5.4.3}.
\hfill\qed

\begin{rem}
In the earlier version of this monograph,
the claim of Proposition {\rm\ref{prop;10.5.4.3}}
was proved in the case $\dim X=2$.
We also proved the claim of 
Corollary {\rm\ref{cor;10.5.4.101}}
with a slightly different argument.
Proposition {\rm\ref{prop;10.5.4.3}}
seems useful for simplification of the arguments.
Although it is also given 
in our survey paper \cite{mochi11},
we include it for the convenience of readers.
\hfill\qed
\end{rem}

\subsection{Proof of Proposition \ref{prop;10.5.4.3}
in the smooth divisor case}
\label{subsection;10.5.4.60}

Let $X:=\Delta^n$
and $D:=\{z_1=0\}$.
Let $\Dhat$ denote the completion of $X$
along $D$.
Let $\pi:X\lrarr D$ denote the natural projection.

\subsubsection{}

First, let us observe that
we can ignore the subsets
whose codimension in $X$ is larger than $3$.
The result in this section
is also used in Section \ref{subsection;08.9.3.1}.
We will use the following easy lemma implicitly.
\begin{lem}
Let $Z$ be a closed analytic subset of $D$
with $\codim_D(Z)\geq 2$.
Let $\nbigi$ be a good set of irregular values
on $(X\setminus Z,D\setminus Z)$.
Then, it is also a good set of irregular values
on $(X,D)$.
\end{lem}
\pf
Take $\gminia\in\nbigi$.
By Hartogs property,
we obtain that $\gminia\in M(X,D)/H(X)$.
Because $\gminia_{\ord(\gminia)}$
is nowhere vanishing on $D\setminus Z$,
we obtain that it is also nowhere vanishing
on $D$.
We can check the other property similarly.
\hfill\qed

\vspace{.1in}

Let $(\nbige,\nabla)$ be a meromorphic
flat connection on $(X,D)$,
i.e, $\nbige$ is 
a (not necessarily locally free)
coherent $\nbigo_{X}(\ast D)$-module
with a meromorphic flat connection
$\nabla:\nbige\lrarr\nbige\otimes\Omega_X^1$.
Assume the following:
\begin{itemize}
\item
 There exist 
 a closed analytic subset $Z\subset D$
 with $\codim_D(Z)\geq 2$,
 and a good set of irregular values 
 $\nbigi\subset M(X,D)/H(X)$,
 such that 
 the irregular values of
 the restrictions
 $(\nbige,\nabla)_{|\pi^{-1}(P)}$
 are given by 
 $\bigl\{\gminia_{|\pi^{-1}(P)}\,\big|\,
 \gminia\in \nbigi\bigr\}$
 for any $P\in D\setminus Z$.
\end{itemize}

\begin{lem}
 \label{lem;08.1.21.50}
If the above condition is satisfied,
the Deligne-Malgrange lattice $E$
of $(\nbige,\nabla)$ is a locally free $\nbigo_X$-module,
and unramifiedly good
with $\Irr(\nabla,Q)=\nbigi$
for any $Q\in D$.
\end{lem}
\pf
Since $E$ is reflexive,
by extending $Z$,
we may assume that
$E_{|X\setminus Z}$ is locally free.
By using a well established argument
(see \cite{levelt} or \cite{malgrange}, for example),
we can easily obtain the irregular decomposition
$E_{|\widehat{D\setminus Z}}
=\bigoplus_{\gminia\in \nbigi}
 F_{\gminia,D\setminus Z}$.
Let $\pi_{\gminia}$ denote the projection
onto $F_{\gminia,D\setminus Z}$,
which gives a section of
$\End(E)_{|\widehat{D\setminus Z}}$.

Let us see that $\pi_{\gminia}$
is extended to a section of
$\End(E)_{|\widehat{D}}$.
Although it follows from a general result,
we show it directly.
It is easy to show the following claim
by using Hartogs theorem:
\begin{itemize}
\item
Any section of
$\nbigo_{\widehat{D\setminus Z}}$
is extended to a section of
$\nbigo_{\Dhat}$.
\end{itemize}

Since $E$ is reflexive,
we can (locally) take an injection
$i:E\lrarr \nbigo_X^{\oplus N}$ for some large $N$
such that the cokernel $\Cok(i)$
is torsion-free.
We can also take a surjection
$\varphi:\nbigo_X^{\oplus M}\lrarr E$.
The morphisms
$i$, $\varphi$ and $\pi_{\gminia}$
induce a morphism
$F_{\gminia}:
 \nbigo^{\oplus M}_{X|\widehat{D\setminus Z}}
 \lrarr 
 \nbigo^{\oplus N}_{X|\widehat{D\setminus Z}}$.
It is extended to a morphism
$\Ftilde_{\gminia}:\nbigo^{\oplus M}_{\Dhat}
 \lrarr
 \nbigo^{\oplus N}_{\Dhat}$.
Since $\Cok(i)$ is torsion free,
$\Ftilde_{\gminia}$ factors through $E_{|\Dhat}$.
Let $\nbigk$ denote the kernel of
$\nbigo^{\oplus M}_{|\Dhat}\lrarr E_{|\Dhat}$.
The restriction of $\Ftilde_{\gminia}$ to $\nbigk$
on $\widehat{D\setminus Z}$ is $0$.
Then, we obtain $\Ftilde_{\gminia|\nbigk}=0$
because $\nbigo^{\oplus\,N}_{\Dhat}$ is torsion-free.
Thus, we obtain the induced maps
$\pi_{\gminia}:E_{|\Dhat}\lrarr E_{|\Dhat}$ 
for $\gminia\in \nbigi$,
which satisfy
$\pi_{\gminia}\circ\pi_{\gminia}=\pi_{\gminia}$,
$\pi_{\gminia}\circ\pi_{\gminib}=0$
and $\sum\pi_{\gminia}=\id$.
They give a decomposition
$E=\bigoplus_{\gminia\in \nbigi}\Ehat_{\gminia}$.
Let us show that $\Ehat_{\gminia}$
is $\gminia$-logarithmic.
We have only to consider the case $\gminia=0$.

Take a point $P\in D\setminus Z$.
We have the vector space $V:=\Ehat_{0|P}$.
We have the endomorphism $f$ of $V$
induced by the residue.
Let $E'_{0}:=V\otimes\nbigo_X$ 
and $\nabla'_{0}=d+f\cdot dz_1/z_1$.
We have the natural flat isomorphism
$(E'_{0},\nabla'_0)_{|\pi^{-1}(P)}
\simeq
 (\Ehat_0,\nabla_0)_{|\pi^{-1}(P)}$.
Since the codimension of $Z$ in $D$
is larger than $2$,
we obtain a flat isomorphism
$\Phi_{0,D\setminus Z}:
 (E'_0,\nabla'_0)_{|\widehat{D\setminus Z}}
\simeq
 (\Ehat_0,\nabla_0)_{|\widehat{D\setminus Z}}$.
Since $E_0'$ and $\Ehat_0$ are both reflexive,
the above argument shows
that $\Phi_{0,D\setminus Z}$
and $\Phi_{0,D\setminus Z}^{-1}$
are extended to morphisms on $\Dhat$.
Thus, we are done.
\hfill\qed

\subsubsection{}
\label{subsection;10.5.4.10}

Let $(\nbige,\nabla)$ be a meromorphic flat bundle
on $(\Dhat,D)$.
Assume the following:
\begin{description}
\item[(A)]
$(\nbige,\nabla)_{|\Phat}$
has an unramifiedly good Deligne-Malgrange lattice
$\lefttop{P}E$ for each $P\in D$.
\end{description}
We obtain the following lemma
from Proposition \ref{prop;10.5.3.5}.
\begin{lem}
\label{lem;10.5.4.5}
Let $E$ be an $\nbigo_{\Dhat}$-free lattice of $\nbige$
such that $E_{|\Phat}=\lefttop{P}E$
for any $P\in D$.
Then, the following holds:
\begin{itemize}
\item
 There exists $\nbigi\in z_1^{-1}H(D)[z_1^{-1}]$
 such that
 $\nbigi_{|\Phat}=\Irr(\nabla,P)$ 
 for any $P\in D$.
\item
 We have a flat decomposition
 $E=\bigoplus_{\gminia\in \nbigi}E_{\gminia}$
 whose restriction to $\Phat$ is the same as
 the irregular decomposition of
 $\lefttop{P}E$ for any $P\in D$.
\hfill\qed
\end{itemize}
\end{lem}

\subsubsection{}

We put $Z:=\{z_1=z_2=0\}$.
Let $(\nbige,\nabla)$ be a meromorphic flat bundle
on $(\Dhat,D)$ satisfying the condition
(A) in Subsection \ref{subsection;10.5.4.10}.
Assume there exists an $\nbigo_X$-free lattice $E$
of $\nbige$ such that $E_{|\Phat}=\lefttop{P}E$
for each $D\setminus Z$.
Take any point $P\in D\setminus Z$.
We have a naturally induced action of
the fundamental group
$\pi_1(D\setminus Z,P)$
on $\Irr(\nabla,P)$
by Lemma \ref{lem;10.5.4.5}.

\begin{lem}
\label{lem;10.5.4.21}
If the action of $\pi_1(D\setminus Z,P)$
on $\Irr(\nabla,P)$ is trivial, we have
$E_{|\Qhat}=\lefttop{Q}E$
for any $Q\in Z$.
In particular, the conclusion
of Lemma {\rm\ref{lem;10.5.4.5}} holds.
\end{lem}
\pf
Because the action of $\pi_1(D\setminus Z,P)$
on $\Irr(\nabla,P)$ is trivial,
we have 
$\nbigi\subset 
 z_1^{-1}H(D\setminus Z)[z_1^{-1}]$
such that $\nbigi_{|\Phat'}=\Irr(\nabla,P')$
for any $P'\in D\setminus Z$.
We set $m:=\min\bigl\{
\ord_{z_1}(\gminia)\,\big|\,\gminia\in\nbigi
\bigr\}$.
We use a descending induction on $m$.
If $m=0$,
we can deduce that $\nabla$ is logarithmic,
and hence the claim is obvious.
Let us consider $m+1\Longrightarrow m$.
We put
\[
 T:=\bigl\{
 (z_1^{-m}z_1\del_1\gminia)_{|D}\,\big|\,
 \gminia\in\nbigi
 \bigr\}
\subset
 H(D\setminus Z).
\]
Because 
$z_1^{-m}\nabla(z_1\del_1)(\lefttop{P}E)
\subset \lefttop{P}E$ for any $P\in D\setminus Z$,
we have 
$z_1^{-m}\nabla(z_1\del_1)E_{|D\setminus Z}
 \subset E_{|D\setminus Z}$,
and hence
$z_1^{-m}\nabla(z_1\del_1)E\subset E$.
Let $G$ be the endomorphism of
$E_{|D}$ induced by
$z_1^{-m}\nabla(z_1\del_1)$.
Because the elements of $T$
are the eigenvalues of
$G_{|D\setminus Z}$,
they are algebraic over $H(D)$.
Hence, we obtain $T\subset H(D)$.

Let $Q\in Z$.
We will shrink $X$ around $Q$
without mention.
Let $N$ be the $H(D)(\!(z_1)\!)$-module
corresponding to $\nbige$,
i.e., the space of the global sections
of $\nbige$.
Let $L$ be the $H(D)[\![z_1]\!]$-lattice
of $N$ corresponding to $E$.
We put
$N':=N\otimes M(D,Z)(\!(z_1)\!)$ and
$L':=L\otimes M(D,Z)[\![z_1]\!]$.
We have the eigen decomposition
of $L'/z_1L'$ with respect to $G$,
which is extended to a decomposition
$L'=\bigoplus_{\gminib\in T}L'_{\gminib}$
such that
$\bigl(
 z_1^{-m+1}\del_1-\gminib
 \bigr)L'_{\gminib}\subset L'_{\gminib}$
by Lemma \ref{lem;10.6.23.10}.
We put $m(Q):=\min\bigl\{
 \ord_{z_1}(\gminia)\,\big|\,
 \gminia\in\Irr(\nabla,Q)
 \bigr\}$ and
\[
 T(Q):=\bigl\{
 \bigl(
 z_1^{-m(Q)+1}\del_1\gminia
 \bigr)_{|D}\,\big|\,
 \gminia\in\Irr(\nabla,Q)
 \bigr\}
\]
\begin{lem}
\label{lem;10.5.4.20}
We have $m(Q)=m$ 
and $T(Q)=T$.
\end{lem}
\pf
We may assume $Q=(0,\ldots,0)$.
We put $\nbign:=N\otimes\nbigo_{\Qhat}$.
It is equipped with an unramifiedly good
Deligne-Malgrange lattice $\lefttop{Q}\nbigl$
with $z_1\del_1$-decomposition
$\lefttop{Q}\nbigl=\bigoplus_{\gminib\in T(Q)}
 \lefttop{Q}\nbigl_{\gminib}$
such that 
$(z_1^{-m(Q)+1}\del_1-\gminib)
 \lefttop{Q}\nbigl_{\gminib}
\subset
 \lefttop{Q}\nbigl_{\gminib}$
for any $\gminib\in T(Q)$.
By considering the extension
to the field
$\cnum(\!(z_n)\!)\cdots(\!(z_2)\!)(\!(z_1)\!)$,
and by using Lemma \ref{lem;10.5.4.2},
we obtain Lemma \ref{lem;10.5.4.20}.
\hfill\qed

\vspace{.1in}
Let us return to the proof of Lemma
\ref{lem;10.5.4.21}.
By Lemma \ref{lem;10.5.4.20},
we have the eigen decomposition
of $E_{|D}$ with respect to $G$.
By Lemma \ref{lem;10.6.23.10},
it is extended to a decomposition
$E=\bigoplus_{\gminib}E_{\gminib}$
such that
$(t^{-m+1}\del_t-\gminib)E_{\gminib}
\subset E_{\gminib}$.
Put $\nbige_{\gminib}=E_{\gminib}(\ast D)$.
We can apply the hypothesis of the induction
to $\nbige_{\gminib}\otimes 
 L\bigl(-z_1^{-m}\gminib/m\bigr)$,
and the proof of Lemma \ref{lem;10.5.4.21}
is finished.
\hfill\qed

\begin{lem}
\label{lem;10.5.4.22}
The action of $\pi_1(D\setminus Z)$
on $\Irr(\nabla,P)$ is trivial.
In particular,
the conclusion of Lemma {\rm\ref{lem;10.5.4.5}}
holds.
\end{lem}
\pf
Because $\Irr(\nabla,P)$ is finite,
we can find a ramified covering
$\varphi:X'\lrarr X$
given by
$\varphi(z_1,\zeta_2,z_3,\ldots,z_n)
=(z_1,\zeta_2^e,z_3,\ldots,z_n)$
such that we can apply Lemma \ref{lem;10.5.4.21}
to $\varphi^{\ast}(\nbige,\nabla)$ and 
$\varphi^{\ast}E$.
Then,
$\varphi^{\ast}\Irr(\nabla,P)
\subset
 z_1^{-1}H(X')[z_1^{-1}]$
and 
$\varphi^{\ast}\Irr(\nabla,P)_{|\Qhat}
=\varphi^{\ast}\Irr(\nabla,Q)$.
Hence, we can conclude that
the action of $\pi_1(D\setminus Z,P)$
is trivial.
\hfill\qed

\subsubsection{}

Let $(\nbige,\nabla)$ be a meromorphic flat sheaf
on $(X,D)$
satisfying the condition (A) in Subsection 
\ref{subsection;10.5.4.10}.
\begin{lem}
The Deligne-Malgrange lattice $E$ of
$(\nbige,\nabla)$ is 
unramifiedly good Deligne-Malgrange.
Namely, the claim of Proposition 
{\rm\ref{prop;10.5.4.3}} holds
if $D$ is smooth.
\end{lem}
\pf
There exists a closed analytic subset $Z\subset D$
with $\codim_D(Z)\geq 2$
such that $E_{|X-Z}$ is locally free.
By Lemma \ref{lem;10.5.4.21}
and Lemma \ref{lem;10.5.4.22},
we obtain that there exists a closed
analytic subset $Z'\subset D$
with $\codim_D(Z')\geq 2$
such that $E_{|X-Z'}$ is good Deligne-Malgrange.
Then, the claim of the lemma follows from
Lemma \ref{lem;08.1.21.50}.
\hfill\qed

\subsection{Proof of Proposition \ref{prop;10.5.4.3}
in the normal crossing case}
\label{subsection;10.5.4.61}

Let $X=\Delta^n$
and $D=\bigcup_{i=1}^{\ell}\{z_i=0\}$.
We put
$\del D_1:=D_1\cap
 \bigcup_{2\leq j\leq \ell}D_j$.
We put $D_1^{\circ}:=D_1\setminus\del D_1$.

\subsubsection{}

We regard $M(D_1,\del D_1)(\!(z_1)\!)$
be a differential ring equipped with 
the differential $\del_1$.
Let $\nbign$ be a differential
$M(D_1,\del D_1)(\!(z_1)\!)$-module
with a $M(D_1,\del D_1)[\![z_1]\!]$-free
lattice $\nbigl$.
We put 
$\nbigl':=\nbigl\otimes H(D_1^{\circ})[\![z_1]\!]$.
Assume that we have 
$\nbigi\in z_1^{-1}H(D_1^{\circ})[z_1^{-1}]$
and a decomposition
$\nbigl'=\bigoplus_{\gminia\in\nbigi}
 \nbigl'_{\gminia}$
such that
(i) $(z_1\del_1-z_1\del_1\gminia)\nbigl_{\gminia}'
 \subset\nbigl_{\gminia}'$,
(ii) the eigenvalues $\alpha$ of 
 the induced endomorphism
 of $\nbigl'_{\gminia}/z_1\nbigl'_{\gminia}$ satisfy
 $0\leq \Re(\alpha)<1$.

\begin{lem}
\label{lem;10.5.4.201}
We have 
$\nbigi\subset z_1^{-1}M(D_1,\del D_1)[z_1^{-1}]$
and a decomposition
$\nbigl=\bigoplus_{\gminia\in\nbigi}
 \nbigl_{\gminia}$
such that
(i) $(z_1\del_1-z_1\del_1\gminia)\nbigl_{\gminia}
 \subset\nbigl_{\gminia}$,
(ii) the eigenvalues $\alpha$ of the induced morphism
 of $\nbigl_{\gminia}/z_1\nbigl_{\gminia}$ satisfy
 $0\leq \Re(\alpha)<1$.
Moreover, we have
$\nbigl_{\gminia}\otimes
 H(D_1^{\circ})[\![z_1]\!]
=\nbigl'_{\gminia}$.
\end{lem}
\pf
We use a descending induction on
$m(\nbigl):=\min\bigl\{
 \ord_{z_1}(\gminia)\,\big|\,
 \gminia\in\nbigi
 \bigr\}$.
If $m(\nbigl)=0$,
there is nothing to prove.
Let us consider the case
$m(\nbigl)=m<0$.
We put
$T(\nbigl):=\bigl\{
 m(z_1^{-m}\gminia)_{|z_1=0}\,\big|\,
 \gminia\in\nbigi
 \bigr\}$.
Let us consider the endomorphism $G$
of $\nbigl/z_1\nbigl$
induced by 
$z_1^{-m}\nabla(z_1\del_1)$.
Because the elements of $T(\nbigl)$ are
the eigenvalues of $G$,
they are algebraic over $M(D_1,\del D_1)$.
Then, we can deduce 
$T(\nbigl)\subset M(D_1,\del D_1)$
from $T(\nbigl)\subset H(D_1^{\circ})$.
If $|T(\nbigl)|=1$, 
by considering the tensor product
with a meromorphic flat bundle of rank one,
we can reduce the issue to the case 
$m(\nbigl)=m+1$.
Let us consider the case $|T(\nbigl)|\geq 2$.
It is standard that
the eigen decomposition of
$\nbigl/z_1\nbigl$
is uniquely extended to 
a $\nabla$-flat decomposition
$\nbigl=\bigoplus_{\gminib\in T(\nbigl)}
 \nbigl_{\gminib}$.
It is easy to observe that
$m(\nbigl_{\gminib})\geq m$,
and $|T(\nbigl_{\gminib})|\leq 1$
if $m(\nbigl_{\gminib})=m$.
Thus, we are done.
\hfill\qed

\subsubsection{}
\label{subsection;10.5.4.210}

Let $(\nbige,\nabla)$ be a meromorphic flat bundle
on $(X,D)$.
Let $(\nbige,\nabla)$ be a meromorphic
flat bundle on $(X,D)$.
Assume the following:
\begin{itemize}
\item
 $(\nbige,\nabla)_{|\Phat}$
 is unramifiedly good
 for each $P\in D$.
\item
 The Deligne-Malgrange lattice $E$ of
 $(\nbige,\nabla)$ is $\nbigo_X$-locally free.
\end{itemize}
Let us show that
$E$ is unramifiedly good Deligne-Malgrange
under the above assumption.

We put
$D^{[2]}:=\bigcup_{i\neq j}(D_i\cap D_j)$.
We can take a ramified covering
$\varphi:(X,D)\lrarr (X,D)$
with the following property:
\begin{itemize}
\item
For each $P\in D_i\setminus D^{[2]}$,
the action of $\pi_1(D_i\setminus D^{[2]},P)$
on $\Irr(\varphi^{\ast}\nabla,P)$
is trivial.
\end{itemize}
By using the argument for Lemma \ref{lem;10.5.4.22},
we may and will assume 
that the above property
holds for $(\nbige,\nabla)$
from the beginning.
We have already known that
$E_{|X-D^{[2]}}$ is unramifiedly
good Deligne-Malgrange.
In particular,
we have
$\nbigi\subset z_1^{-1}H(D_1^{\circ})[z_1^{-1}]$
and a decomposition
$E_{|\Dhat_1^{\circ}}
=\bigoplus_{\gminia\in \nbigi}
 \Ehat_{\gminia}$
such that
$\bigl(
 \nabla(z_1\del_1)-z_1\del_1\gminia
\bigr)\Ehat_{\gminia}
\subset \Ehat_{\gminia}$.

Let $\nbigm$ be the differential
$M(D_1,\del D_1)(\!(z_1)\!)$-module
corresponding to $\nbige$,
and let $\nbigl$ be 
the $M(D_1,\del D_1)[\![z_1]\!]$-lattice
induced by $E$.
Applying Lemma \ref{lem;10.5.4.201},
we obtain
$\nbigi\subset z_1^{-1}M(D_1,\del D_1)[z_1^{-1}]$
and a decomposition
$\nbigl=\bigoplus_{\gminia\in\nbigi}
 \nbigl_{\gminia}$
such that
$(z_1\del_1-z_1\del_1\gminia)\nbigl_{\gminia}
\subset \nbigl_{\gminia}$.
Let $\gbigk:=\cnum(\!(z_n)\!)\cdots(\!(z_2)\!)$.
By the natural extension
$M(D_1,\del D_1)\subset\gbigk$,
$\nbigl\otimes\gbigk[\![z_1]\!]$
is the Deligne-Malgrange lattice
of the differential module
$\gbigm:=\bigl(
 \nbign\otimes\gbigk(\!(z_1)\!),\del_1
\bigr)$.

Let $\lefttop{O}E$ be the unramifiedly good
Deligne-Malgrange lattice of $\nbige_{|\Ohat}$
with the irregular decomposition
$\lefttop{O}E=
 \bigoplus_{\gminia\in\Irr(\nabla,O)}
 \lefttop{O}E_{\gminia}$.
Let $\Irr(\nabla,1)$ be the image of
$\Irr(\nabla,O)$ via the map
$\nbigo_{\Ohat}(\ast D)/\nbigo_{\Ohat}
\lrarr
 \nbigo_{\Ohat}(\ast D)/\nbigo_{\Ohat}(\ast D(\neq 1))$.
It is easy to see that
$\gbigk[\![z_1]\!]\otimes \lefttop{O}E$
is the good Deligne-Malgrange lattice of
$(\gbigm_1,\del_1)$,
and the set of the irregular values
is given by $\Irr(\nabla,1)$.
Hence, we obtain 
$\Irr(\nabla,1)=\nbigi_1$
in $z_1^{-1}\gbigk[z_1^{-1}]$
and
$\lefttop{O}E\otimes\gbigk[\![z_1]\!]
=\nbigl_1\otimes\gbigk[\![z_1]\!]$.
We can deduce 
a similar relation for each $i=2,\ldots,\ell$.
Then, we obtain that
$\Irr(\nabla)\subset M(X,D)/H(X)$.

We take a frame $\vecv$ of $\lefttop{O}E$.
Let $f$ be a section of $E$.
We have the expression
$f=\sum f_p\, v_p$.
We obtain 
$f_p\in \gbigk[\![z_1]\!]$,
and hence $f_p$ is $z_1$-regular,
i.e., $f_p$ does not contain
the negative power of $z_1$.
Similarly, we obtain that $f_p$ are $z_j$-regular
for $j=2,\ldots,\ell$.
Thus, we obtain
$E_{|\Ohat}\subset \lefttop{O}E$.
Similarly, we obtain
$\lefttop{O}E\subset E_{|\Ohat}$,
and hence $E_{|\Ohat}=\lefttop{O}E$.
Thus, we obtain that
$E$ is unramifiedly good Deligne-Malgrange lattice.

\subsubsection{}

Let us consider the case
in which 
we do not assume that $E$ is locally free.
We have a closed analytic subset $Z\subset D$
with $\codim_D(Z)\geq 2$
such that 
$E_{|X-Z}$ is locally free.
Then, it is an unramifiedly 
good Deligne-Malgrange lattice
of $(\nbige,\nabla)_{|X-Z}$,
according to the result in Subsection 
\ref{subsection;10.5.4.210}.
We put $D_1^{\ast}:=D_1\setminus Z$
and $\del D_1^{\ast}:=\del D_1\setminus Z$.
We have 
$\nbigi\subset z_1^{-1}M(D_1^{\ast},
 \del D_1^{\ast})[z_1^{-1}]$
and the irregular decomposition
$E_{|\widehat{D_1}^{\ast}}
=\bigoplus_{\gminia\in \nbigi}
 \Ehat_{\gminia,D_1^{\ast}}$.
By using the Hartogs property
and the argument in the proof of
Lemma \ref{lem;08.1.21.50},
we obtain 
$\nbigi\subset z_1^{-1}M(D_1,\del D_1)[z_1^{-1}]$
and a decomposition
$E_{|\Dhat_1}=\bigoplus_{\gminia\in\nbigi}
 \Ehat_{\gminia}$
such that
(i) $(\nabla(z_1\del_1)-z_1\del_1\gminia)
 \Ehat_{\gminia}
\subset
 \Ehat_{\gminia}$,
(ii) the eigenvalues $\alpha$ of
the induced endomorphism of
$\Ehat_{\gminia|D_1}$
satisfy $0\leq \Re(\alpha)<1$.
Then, as in Subsection \ref{subsection;10.5.4.210},
we obtain
$E\otimes\gbigk[\![z_1]\!]
=\lefttop{O}E\otimes\gbigk[\![z_1]\!]$.
Let $\vecv$ be a frame of $\lefttop{O}E$.
Let $f$ be a section of $E$.
We have the expression $f=\sum f_p\,v_p$.
Then, we obtain that $f_p$ is $z_1$-regular.
Similarly, we obtain that
$f_p$ are $z_j$-regular $(j=1,\ldots,\ell)$
and hence $E_{|\Ohat}\subset\lefttop{O}E$.

To show $\lefttop{O}E\subset E_{|\Ohat}$,
we consider the dual.
Put $\nbige^{\lor}:=
 \Hom_{\nbigo_X(\ast D)}(\nbige,\nbigo_{X}(\ast D))$,
which is equipped with a naturally induced flat
connection $\nabla$.
Put $E^{\lor}:=\Hom_{\nbigo_X}(E,\nbigo_X)$,
which is a lattice of $\nbige^{\lor}$.
It is generically unramifiedly good lattice,
and the eigenvalues $\alpha$ of the residue
satisfy $-1<\Re(\alpha)\leq 0$.
Put $\lefttop{O}E^{\lor}:=
 \Hom_{\nbigo_{\Ohat}}(\lefttop{O}E,\nbigo_{\Ohat})$
which is an unramifiedly good lattice of 
$(\nbige^{\lor},\nabla)_{|\Ohat}$.
The eigenvalues $\alpha$ of the residues
satisfy $-1<\Re(\alpha)\leq 0$.
Then, we obtain
$E^{\lor}_{|\Ohat}
\subset \lefttop{O}E^{\lor}$
by the above argument.
Note that
$E^{\lor}_{|\Ohat}
=\Hom_{\nbigo_X}(E,\nbigo_X)\otimes\nbigo_{\Ohat}
\simeq
 \Hom_{\nbigo_{\Ohat}}(E_{|\Ohat},\nbigo_{\Ohat})$.
Hence, we can conclude that
$\lefttop{O}{E}=E_{|\Ohat}$.
Thus, we obtain that
$E$ is an unramifiedly good Deligne-Malgrange
lattice of $(\nbige,\nabla)$ at $O$,
and the proof of Proposition \ref{prop;10.5.4.3}
is finished.
\hfill\qed

\section{Family of filtered $\lambda$-flat bundles
and KMS-structure}
\label{section;10.5.20.1}
\subsection{Notation}

We will be particularly interested in the case
$\nbigk\subset\cnum_{\lambda}$
and $\varrho(\lambda)=\lambda$.
Let $X$ be a complex manifold,
and let $D$ be a simple normal crossing 
hypersurface of $X$
with the irreducible decomposition
$D=\bigcup_{i\in \Lambda}D_i$.
Let $\nbigx$ be an open subset of
$\nbigk\times X$,
and $\nbigd:=\nbigx\cap(\nbigk\times D)$.
If $\nbigk$ is not a point,
a $\lambda$-flat connection on $(\nbigx,\nbigd)$
is called a family of $\lambda$-connections
to emphasize that it does not 
contain the differential in the $\lambda$-direction.
We use the words ``family of'' in similar meanings.
For example,
``good family of filtered $\lambda$-flat bundles''
is just a good filtered $\lambda$-flat bundle
on $(\nbigx,\nbigd)$ in the sense of 
Definition \ref{df;10.5.5.10}.
\index{unramifiedly good family of
 filtered $\lambda$-flat bundles}
\index{good family of
 filtered $\lambda$-flat bundles}

For a fixed $\lambda_0$,
a neighbourhood of $\{\lambda_0\}\times X$
in $\cnum_{\lambda}\times X$
will be often denoted by $\nbigxzero$.
In that case,
for any subset $Y$ of $X$,
we will put $\nbigyzero:=
 (\cnum_{\lambda}\times Y)\cap\nbigxzero$
and $\nbigy^{\lambda}:=
 (\{\lambda\}\times Y)\cap \nbigxzero$.
If we are given a family of flat $\lambda$-connections
on $\nbigx$,
its restriction to $\nbigxlambda$
is denoted by $\DDlambda$.
\index{$\lambda$-connection $\DDlambda$}
\index{set $\nbigxzero$}
\index{set $\nbigyzero$}
\index{set $\nbigxlambda$}
\index{set $\nbigy^{\lambda}$}

Let $(\vecE_{\ast},\DD)$ be a good family of
filtered $\lambda$-flat bundles
on $(\nbigxzero,\nbigdzero)$.
For a subset $I$ of $\Lambda$,
we put $D_I:=\bigcap_{i\in I}D_i$.
In this case, the induced filtrations of
$\prolongg{\vecc}{E}_{|\nbigd^{(\lambda_0)}_I}$
are denoted by $\lefttop{i}\Fzero$
for $i\in I$.
(See Subsection \ref{subsection;10.5.28.1}
for the induced filtration.)
For $\veca_I=(a_i)\in\real^I$,
we use the following notation:
\[
 \lefttop{I}\Fzero_{\veca_I}\bigl(
 \prolongg{\vecc}{E}_{|\nbigdzero_I}
 \bigr)
:=\bigcap_{i\in I}\lefttop{i}\Fzero_{a_i}
 \bigl(\prolongg{\vecc}{E}_{|\nbigdzero_I}\bigr),
\,\,\,
 \lefttop{I}\Gr^{\Fzero}_{\veca_I}\bigl(
 \prolongg{\vecc}{E}
 \bigr):=
 \frac{\lefttop{I}\Fzero_{\veca_I}
 (\prolongg{\vecc}{E}_{|\nbigdzero_I})}
 {\sum_{\vecb_I\lneq\veca_I}\lefttop{I}\Fzero
 \bigl(\prolongg{\vecc}{E}_{|
 \nbigdzero_I}\bigr)}
\]
\index{bundle
 $\lefttop{I}\Gr^{\Fzero}_{\veca_I}\bigl(
 \prolongg{\vecc}{E}\bigr)$}
We put
\[
 \Par\bigl(\prolongg{\vecc}{E},I\bigr):=
 \bigl\{
 \veca_I\in\real^I\,\big|\,
\lefttop{I}\Gr^{\Fzero}_{\veca_I}
 (\prolongg{\vecc}{E})\neq 0
 \bigr\},
\quad
 \Par\bigl(\vecE_{\ast},I\bigr):=
 \bigcup_{\vecc\in\real^{\Lambda}}
 \Par\bigl(\prolongg{\vecc}{E},I\bigr).
\]
\index{set 
 $\Par\bigl(\prolongg{\vecc}{E},I\bigr)$}
\index{set
 $\Par\bigl(\vecE_{\ast},I\bigr)$}
As in Proposition \ref{prop;10.5.5.12},
we have the induced endomorphism
$\Res_i(\DD)$ 
on $\lefttop{i}\Gr^{\Fzero}_a(\prolongg{\vecc}{E})$,
which preserves the induced filtrations
$\lefttop{k}\Fzero$ of
$\lefttop{i}\Gr^{\Fzero}_a
 (\prolongg{\vecc}{E})_{|\nbigdzero_i\cap \nbigdzero_k}$.
Hence, we have the well defined endomorphisms
$\Res_i(\DD)$ $(i\in I)$
on $\lefttop{I}\Gr^{\Fzero}_{\veca_I}
 \bigl(\prolongg{\vecc}{E}\bigr)$.

\subsection{KMS-structure}
\label{subsection;07.10.14.15}

For simplicity,
let us consider the case
in which $\nbigx^{(\lambda_0)}$
is the product of $X$
and some neighbourhood $U(\lambda_0)$
of $\lambda_0$ in $\cnum_{\lambda}$.
Let $\paramap(\lambda):\real\times\cnum\lrarr\real$
and $\eigenmap(\lambda):\real\times\cnum\lrarr\cnum$
be given as follows:
\[
 \paramap\bigl(\lambda,(a,\alpha)\bigr)
=a+2\Re(\lambda\,\alphabar),
\quad
 \eigenmap\bigl(\lambda,(a,\alpha)\bigr)
=\alpha-a\, \lambda-\alphabar\,\lambda^2
\]
The induced map
$\real\times\cnum\lrarr\real\times\cnum$
is denoted by $\kmsmap(\lambda)$.
\index{map
$\paramap(\lambda)$}
\index{map
 $\eigenmap(\lambda)$}
\index{map
 $\kmsmap(\lambda)$}
Let $(\vecE_{\ast},\DD)$ be a 
good family of filtered $\lambda$-flat bundle
on $(\nbigxzero,\nbigdzero)$.
\begin{df}
\label{df;07.11.23.5}
We say that 
$(\vecE_{\ast},\DD)$ has the KMS-structure
at $\lambda_0$ indexed by
$T(i)\subset\real\times\cnum$ $(i\in \Lambda)$,
if the following holds:
\begin{itemize}
\item
 $\Par(\vecE_{\ast},i)$ is the image of
 $T(i)$ via the map $\paramap(\lambda_0)$.
\item 
 For each $a\in\Par(\vecE_{\ast},i)$,
 we put $\nbigk(a,i):=\bigl\{
 u\in T(i)\,\big|\,\paramap(\lambda_0,u)=a
 \bigr\}$.
 Then, the set of the eigenvalues of
 $\Res_i(\DDlambda)$
 on $\lefttop{i}\Gr^{\Fzero}_{a}
 \bigl(\prolongg{\vecc}{E}\bigr)_{|\nbigdlambda_i}$
 is
 $\bigl\{
 \eigenmap(\lambda,u)\,\big|\,
 u\in \nbigk(a,i)
\bigr\}$.
\end{itemize}
In that case,
$0\neq \lambda\in U(\lambda_0)$ is called generic,
if $\eigenmap(\lambda):
 T(i)\lrarr \cnum$
are injective for any $i$.
\hfill\qed
\end{df}
\index{KMS-structure}
\index{generic}

Assume $(\vecE_{\ast},\DD)$
has the KMS-structure at $\lambda_0$.
We have the generalized eigen decomposition
with respect to the induced action of $\Res_i(\DD)$
\begin{equation}
 \label{eq;07.10.16.20}
 \lefttop{i}\Gr^{\Fzero}_{a}
 \bigl(\prolongg{\vecc}{E}\bigr)
=\bigoplus_{u\in\nbigk(a,i)}
 \lefttop{i}\nbigg^{(\lambda_0)}_{u}
\bigl(\prolongg{\vecc}{E}\bigr),
\end{equation}
where $\Res_i(\DD)-\eigenmap(\lambda,u)$
are nilpotent
on $\lefttop{i}\nbigg^{(\lambda_0)}_{u}
\bigl(\prolongg{\vecc}{E}\bigr)$.

\begin{rem}
\label{rem;08.9.4.10}
If $(\vecE_{\ast},\DD)$ has the KMS-structure
at $\lambda_0$,
we have the decomposition
$\prolongg{\vecc}{E}_{|\nbigdzero_i}
=\bigoplus
 \lefttop{i}\EEzero_{\alpha}(
 \prolongg{\vecc}{E}_{|\nbigdzero_i})$
indexed by the eigenvalues of
$\Res_i(\DD^{\lambda_0})$,
such that
(i) it is preserved by $\Res_i(\DD)$,
(ii) the restriction of the decomposition
to $\nbigd^{\lambda_0}_i$
is the generalized eigen decomposition
of $\Res_i(\DD^{\lambda_0})$.
Then, the decomposition $\lefttop{i}\EEzero$ and
the filtration $\lefttop{i}\Fzero$ are compatible,
and 
$\lefttop{i}\Gr^{\Fzero}_a\lefttop{i}\EEzero_{\alpha}
 \bigl(\prolongg{\vecc}{E}_{|\nbigdzero_i}\bigr)$
is naturally isomorphic to 
$\lefttop{i}\nbigg^{(\lambda_0)}_u
 \bigl(\prolongg{\vecc}{E}\bigr)$,
where $u$ is determined by
$\kmsmap(\lambda_0,u)=(a,\alpha)$.
\hfill\qed
\end{rem}

If $(\vecE_{\ast},\DD)$ is unramified,
the KMS-structure is compatible 
with the irregular decomposition
in the following sense.
For simplicity we consider the case
$X=\Delta^n$, $D_i=\{z_i=0\}$ and 
$D=\bigcup_{i=1}^{\ell}D_i$.
We have the irregular decomposition
$\prolongg{\vecc}{E}_{|U(\lambda_0)\times \Ohat}
=\bigoplus_{\gminia\in T}
 \prolongg{\vecc}{\Ehat}_{\gminia}$.
They give a family of filtered $\lambda$-flat bundles
$\vecEhat_{\ast}$ on $U(\lambda_0)\times \Ohat$,
with the irregular decomposition
$(\vecEhat_{\ast},\DD)
=\bigoplus_{\gminia\in T}
 (\vecEhat_{\gminia\,\ast},\DDhat_{\gminia})$.
Each $(\vecEhat_{\gminia\ast},\DDhat_{\gminia})$
has the KMS-structure at $\lambda_0$.
By the natural isomorphism
$\prolongg{\vecc}{E}_{|O}
\simeq
 \bigoplus_{\gminia}
 \prolongg{\vecc}{E}_{\gminia|O}$,
the filtrations $\Fzero$ and the decompositions
$\EEzero$ are the same.

\subsection{Uniqueness of the filtrations}

According to the following lemma,
it makes sense to say that
$(\vecE,\DD)$ has 
the KMS-structure at $\lambda_0$.

\begin{lem}
\label{lem;07.11.23.5}
Let $(\vecE_{i\ast},\DD_i)$ $(i=1,2)$
be good filtered $\lambda$-flat bundles
on $(\nbigxzero,\nbigdzero)$,
which have the $KMS$-structures at $\lambda_0$.
Assume that we are given an isomorphism
$\varphi:(\vecE_{1},\DD_1)\simeq (\vecE_{2},\DD_2)$
of families of meromorphic $\lambda$-flat bundles.
Then, it induces an isomorphism
$\varphi:(\vecE_{1\,\ast},\DD_1)\simeq 
 (\vecE_{2\,\ast},\DD_2)$
of families of filtered $\lambda$-flat bundles.
\end{lem}
\pf
It can easily be reduced to the case
in which $D$ is smooth.
Let $\nbigdhatzero$ denote 
the completion of $\nbigxzero$
along $\nbigdzero$.
We are given the induced isomorphism
$\vecE_{1|\nbigdhatzero}
 \simeq \vecE_{2|\nbigdhatzero}$.
We have only to show that
it induces $\prolongg{a}{E}_{1|\nbigdhatzero}
=\prolongg{a}{E}_{2|\nbigdhatzero}$
for each $a\in\real$.

Let us consider the case
in which $(\vecE_{i\ast},\DD_i)$ are unramified.
We have the irregular decompositions:
\[
 (\vecE_{i\ast},\DD_i)_{|\nbigdhatzero}
=\bigoplus_{\gminia\in \Irr(\DD_i)}
 (\vecEhat_{i,\gminia\,\ast},\DDhat_{i,\gminia})
\]
We have $\Irr(\DD_1)=\Irr(\DD_2)$
and $\vecEhat_{1,\gminia}=\vecEhat_{2,\gminia}$
for each $\gminia$.
Let $T_i$ $(i=1,2)$
denote the index sets of the KMS-structures
of $\vecE_{i\ast}$.
Note that there exists a discrete subset
 $Z$ in $U(\lambda_0)\setminus\{0\}$
such that 
$\eigenmap(\lambda):
 T_1\cup T_2\lrarr \cnum$ is injective
for any $\lambda\in U(\lambda_0)-Z$.
Take any $\lambda_1\in U(\lambda_0)-Z$,
and a neighbourhood 
$U(\lambda_1)$ in $U(\lambda_0)-Z$.
We set $\nbigx^{(\lambda_1)}:=
 U(\lambda_1)\times X$
and $\nbigd^{(\lambda_1)}:=
 U(\lambda_1)\times D$.

\begin{lem}
\label{lem;07.12.30.20}
We have
$\varphi\bigl(
 \prolongg{a}{E}_{1|\nbigxiti}
 \bigr)
\subset
 \prolongg{a}{E}_{2|\nbigxiti}$.
\end{lem}
\pf
Since this is a quite standard claim,
we give only an outline.
We put 
\[
 \nbigl(a,\lambda_0):=
 \bigl\{u\in T_1\cup T_2\,\big|\,
 a-1<\paramap(\lambda_0,u)\leq a
 \bigr\}.
\]
We have the generalized eigen decomposition
\[
 \prolongg{a}{E}_{i|\nbigditi}
=\bigoplus_{u\in \nbigl(a,\lambda_0)}
 \EE_{\eigenmap(\lambda,u)}
 \bigl(\prolongg{a}{E}_{i|\nbigditi}
 \bigr),
\]
where
$\Res(\DD)$ has a unique
eigen value $\eigenmap(\lambda,u)$
on $\EE_{\eigenmap(\lambda,u)}
 \bigl(\prolongg{a}{E}_{i|\nbigditi}
 \bigr)$.
It is compatible with the irregular decompositions,
i.e.,
\begin{equation}
 \label{eq;08.9.4.4}
 \prolongg{a}{\Ehat}_{i,\gminia|\nbigditi}
=\bigoplus_{u\in \nbigl(a,\lambda_0)}
 \EE_{\eigenmap(\lambda,u)}
 \bigl(\prolongg{a}{\Ehat}_{i,\gminia|\nbigditi}
 \bigr).
\end{equation}
It can be shown that
$\eigenmap(\lambda,u_1)
-\eigenmap(\lambda,u_2)$
are not contained in $\seisuu-\{0\}$
for any distinct $u_1,u_2\in \nbigl(\lambda_0,a)$
and for any $\lambda\in U(\lambda_1)$.
Hence, we obtain the flat decomposition
\[
 \prolongg{a}{\Ehat}_{i,\gminia|\nbigdhatiti}
=\bigoplus_{u\in\nbigl(\lambda_0,a)}
 \prolongg{a}{\Ehat}_{i,\gminia,u}
\]
whose restriction to $\nbigditi$ is the same as 
(\ref{eq;08.9.4.4}).
We put
$\prolongg{a}{\Ehat}_{i,u}:=\bigoplus_{\gminia}
 \prolongg{a}{\Ehat}_{i,\gminia,u}$.
For each $a\in\real$,
there exists $b>0$
such that we have the flat morphism
\[
 \varphi_{u_1,u_2}:
 \prolongg{a}{\Ehat}_{1,u_1}
\lrarr \prolongg{a+b}{\Ehat}_{2,u_2}
\]
induced by the flat isomorphism
$\vecE_{1|\nbigdhatiti}\simeq
 \vecE_{2|\nbigdhatiti}$.
The restriction
$\varphi_{u_1,u_2|\nbigditi}$
has to be compatible with the residues.
Hence, it is easy to observe that
$\varphi_{u_1,u_2}
 (\prolongg{a}{\Ehat}_{1,u_1})
\subset \prolongg{a+b-N}{\Ehat}_{2,u_2}$
for any $N\geq 0$, if $u_1\neq u_2$.
It implies 
$\varphi_{u_1,u_2}=0$ unless $u_1=u_2$.
Let us look at $\varphi_{u,u}$.
Again, by the comparison of the eigenvalues
of the residues,
we obtain 
$\varphi_{u,u}
 \bigl(\prolongg{a}{\Ehat}_{1,u}\bigr)
 \subset\prolongg{a}{\Ehat}_{2,u}$.
Thus, we obtain Lemma \ref{lem;07.12.30.20}.
\hfill\qed

\vspace{.1in}
Let us return to the proof of Lemma \ref{lem;07.11.23.5}.
We have the morphism
$\varphi:
 \prolongg{a}{E_{1}}
\lrarr
 \prolongg{a+b}{E_{2}}$
for some $b$.
By Lemma \ref{lem;07.12.30.20},
the induced map
$\prolongg{a}{E_1}\lrarr
 \Gr^{\Fzero}_{a+b}(E_{2|U(\lambda_0)\times D})$
is trivial if $b>0$.
Hence, we obtain
$\varphi(\prolongg{a}{E}_1)
\subset \prolongg{a}{E_2}$
in the unramified case.

Let us consider the case in which
$(\vecE_{i,\ast},\DD_i)$ are not necessarily unramified.
We set 
$\Xtilde:=X$ and $\Dtilde:=D$.
We use the notation
$\nbigxtildezero$ and $\nbigdtildezero$
in similar meanings.
We take an appropriate ramified covering
$\varphi:\Xtilde\lrarr X$ such that 
the induced filtered bundle
$(\vecEtilde_{i,\ast},\DDtilde_i)$ via $\varphi$
and $(\vecE_{i\ast},\DD_i)$
on $(\nbigxtildezero,\nbigdtildezero)$
is unramifiedly good.
It is easy to see that
$(\vecEtilde_{i,\ast},\DDtilde_i)$
have the KMS-structure at $\lambda_0$.
By applying the above result to them,
we obtain $\vecEtilde_{1\,\ast}=\vecEtilde_{2\,\ast}$.
Since $\vecE_{i\,\ast}$ are obtained
as the descent of $\vecEtilde_{i\,\ast}$,
we obtain $\vecE_{1\,\ast}=\vecE_{2\,\ast}$.
Thus, the proof of Lemma \ref{lem;07.11.23.5}
is finished.
\hfill\qed

\subsection{Openness and invariance}

Pick $\vecc\in \real^{\Lambda}$
such that $c_i\not\in\Par\bigl(\vecE_{\ast},i\bigr)$
for each $i\in \Lambda$.
Let $\pi_{i,a}$ be the projection:
\[
 \lefttop{i}\Fzero_a\bigl(
 \prolongg{\vecc}{E}_{|U(\lambda_0)\times D_i}\bigr)
\lrarr
 \lefttop{i}\Gr^{\Fzero}_a
 \bigl(\prolongg{\vecc}{E}\bigr)
\]
Let $\lambda_1\in U(\lambda_0)$ be sufficiently close to
$\lambda_0$,
and let $U(\lambda_1)\subset U(\lambda_0)$
be a neighbourhood of $\lambda_1$.
Let $c_i-1<b\leq c_i$.
If $b=\paramap(\lambda_1,v)$
for some $v\in \nbigk(a,i)$,
we put on $\nbigditi_i$
\[
 \lefttop{i}F^{(\lambda_1)}_{b}:=
 \bigoplus_{
 \substack{u\in\nbigk(a,i)\\
 \paramap(\lambda_1,u)\leq b}}
\pi_{i,a}^{-1}\bigl(
 \lefttop{i}\nbigg^{(\lambda_0)}_u
\bigr).
\]
Otherwise,
let $b_0:=\max\bigl\{
 \paramap(\lambda_1,v)<b\,\big|\,
 v\in \nbigk(a,i)
 \bigr\}$,
and 
we put $\lefttop{i}F^{(\lambda_1)}_b:=
 \lefttop{i}F^{(\lambda_1)}_{b_0}$.
Thus, we obtain the filtration 
$\lefttop{i}F^{(\lambda_1)}$ of
$\prolongg{\vecc}{E}_{|\nbigditi_i}$.
It induces a family of filtered $\lambda$-flat bundles
$(\vecE^{(\lambda_1)}_{\ast},\DD)$
on $(\nbigxiti,\nbigditi)$.
It is easy to observe that
$\Res_{i}(\DD)$
has a unique eigenvalue 
$\eigenmap(\lambda,u)$
on $\lefttop{i}\Gr^{F^{(\lambda_1)}}_{
 \paramap(\lambda_1,u)}
 \bigl(\prolongg{\vecc}{E}\bigr)_{|\nbigditi_1}$.
Hence, $(\vecE^{(\lambda_1)}_{\ast},\DD)$
has the KMS-structure at $\lambda_1$.
The index sets are equal to
those for $(\vecE^{(\lambda_0)}_{\ast},\DD)$.

\begin{rem}
More precisely,
there exists an open subset
$\lambda_0\in U'(\lambda_0)\subset U(\lambda_0)$
depending only on the sets
$\KMS(\vecE^{(\lambda_0)}_{\ast},i)$
such that the above construction can be applied
for any $\lambda_1\in U'(\lambda_0)$.
\hfill\qed
\end{rem}

For each $\lambda\in U(\lambda_0)$
sufficiently close to $\lambda_0$,
we put $\vecE_{\ast}^{\lambda}:=
 (\vecE_{\ast}^{(\lambda)})_{|
 \{\lambda\}\times(X,D)}$,
which is a good filtered $\lambda$-flat bundle.
The set
$\KMS(\vecE_{\ast}^{\lambda},i)$
is equal to
the image of $T(i)$ 
via the map $\kmsmap(\lambda)$.
Note $\KMS(\vecE_{\ast}^0,i)=T(i)$
in the case $0\in U(\lambda_0)$.
We often identify them.

\chapter[Stokes structure]{Stokes Structure
of Good $\varrho$-Meromorphic Flat Bundle}
\label{section;07.12.26.6}
In this chapter,
we shall study Stokes structure of
unramifiedly good $\varrho$-flat bundle.
In Section \ref{section;10.5.27.40},
we give preliminary to describe Stokes filtration.
In Section \ref{section;10.5.10.5},
we state some theorems and propositions
for full Stokes filtrations,
which will be proved in Section \ref{section;10.5.18.5}.
We state some theorems and propositions
for partial Stokes filtration
in Section \ref{section;10.5.11.1},
which will be proved in Sections
\ref{subsection;08.9.28.30}--\ref{section;10.5.18.5}.

\section{Preliminary}
\label{section;10.5.27.40}
\subsection{Filtration indexed by a finite ordered set}
\label{subsection;08.9.28.200}

\subsubsection{Compatibility}
\label{subsection;08.8.31.1}

\index{filtration indexed by an ordered set}

Let $(I,\leq)$ be a finite ordered set.
Let $V$ be a vector space.
In this section,
a filtration $F$ of $V$ indexed by $(I,\leq)$
means a family of subspaces
$F_{a}\subset V$ $(a\in I)$ with the following property:
\begin{itemize}
\item
$F_a\subset F_b$ if $a\leq b$.
\item
There exists a splitting
$V=\bigoplus V_{a}$
such that $F_{a}=\bigoplus_{b\leq a}V_b$.
\end{itemize}
We put $F_{<a}:=\sum_{b<a} F_{b}$
and $\Gr^F_a(V)=F_{a}/F_{<a}
\simeq V_a$.
\index{Gr $\Gr^F_a(V)$}
For a given subset $S\subset I$,
we set $F_{S}:=\sum_{a\in S}F_a$.

\begin{rem}
Note that we have assumed the existence
of splitting, which is unusual.
For example, we do not have
such a splitting for filtered flat bundle.
We consider the above type of filtration
just for Stokes filtration.
\hfill\qed
\end{rem}

Let $\varphi:(I,\leq )\lrarr (I',\leq')$ be a morphism
of ordered sets,
and let $F$ be a filtration of $V$
indexed by $(I,\leq)$.
Then, we have the induced filtration
$F^{\varphi}$ indexed by $(I',\leq')$ constructed inductively
as follows:
\[
 F^{\varphi}_b=
 F^{\varphi}_{<b}+\sum_{a\in \varphi^{-1}(b)} F_a
\]
We set $V^{\varphi}_b:=
 \bigoplus_{a\in \varphi^{-1}(b)}V_a$.
Then, $V=\bigoplus_{b\in I'} V^{\varphi}_b$
gives a splitting of $F^{\varphi}$.
We say that $F^{\varphi}$
is induced by $F$ and $\varphi$.

\begin{df}
 \label{df;07.12.13.20}
\index{compatibility of filtrations}
Let $F$ and $F'$ be filtrations of
$V$ indexed by $(I,\leq)$ and $(I',\leq')$,
respectively.
Let $\varphi:(I,\leq)\lrarr (I',\leq')$ be a morphism
of ordered sets.
We say that $F$ and $F'$ are compatible 
over $\varphi$,
if $F'$ is the same as $F^{\varphi}$ above.

In the case $I=I'$ 
(but possibly $(I,\leq)\neq (I',\leq')$)
and $\varphi'=\id$,
we just say $F$ and $F'$ are compatible.
\hfill\qed
\end{df}

In the case $I=I'$,
we have the natural isomorphism
$\Gr^F_a(V)\simeq \Gr^{F'}_{a}(V)$.
\begin{lem}
\label{lem;10.5.12.22}
Let $F$ be a filtration of $V$
indexed by $(I,\leq)$.
Let $\leq_i$ $(i\in \Lambda)$
be orders on $I$ such that
(i) the identity 
$\varphi_i:(I,\leq)\lrarr(I,\leq_i)$
 are order preserving,
(ii) $a\leq b$
if and only if 
 $a\leq_i b$ for any $i\in\Lambda$.
Then, $F$ can be reconstructed 
from $F^{\varphi_i}$ $(i\in\Lambda)$
in the sense
$F_{a}=\bigcap_{i\in\Lambda}F^{\varphi_i}_{a}$.
\end{lem}
\pf
We take a splitting
$V=\bigoplus_{a\in I}V_{a}$
of the filtration $F$.
Recall 
$F^{\varphi_i}_a=
 \bigoplus_{b\leq_ia}V_b$.
Then, the claim of the lemma is clear.
\hfill\qed

\begin{rem}
It can be generalized 
for vector bundles appropriately.
\hfill\qed
\end{rem}

\subsubsection{Dual}
\label{subsection;08.9.5.3}

Let $(I,\leq)$ be an ordered set,
and let $V$ be a finite dimensional vector space
equipped with a filtration $F$
indexed by $(I,\leq)$.
Let us give an induced filtration $F^{\lor}$
on the dual vector space $V^{\lor}$.
We set $I^{\lor}:=I$
and let $\leq^{\lor}$ be the order of $I^{\lor}$
defined by
$a\leq^{\lor}b\Longleftrightarrow a\geq b$.
For distinction,
we use the symbol $-a$
if we regard $a\in I$ as an element of $I^{\lor}$.
And, ``$-a\leq^{\lor}-b$''
is denoted by $-a\leq -b$.
We hope that there are no risk of confusion.
\index{ordered set $(I^{\lor},\leq)$}

We take a splitting
$V=\bigoplus_{a\in I} V_a$
of the filtration $F$.
In general,
for a vector subspace $U\subset V$,
let $U^{\bot}\subset V^{\lor}$
be $\bigl\{
 f\in V^{\lor}\,\big|\,
 f(v)=0\,\,\forall v\in U
 \bigr\}$.
\index{dual space $V^{\lor}$}
\index{subspace $U^{\bot}$}
For each $a\in I$,
let $S(a)$ denote the set of
$b\in I$ such that 
$b\not\geq a$.
We have the subspaces of $V^{\lor}$
given as follows:
\[
 V^{\lor}_{-a}:=
\Bigl(
\bigoplus_{b\neq a}V_b
\Bigr)^{\bot},
\quad
 F^{\lor}_{-a}(V^{\lor}):=
\Bigl(
 \bigoplus_{b\in S(a)}
 V_b
\Bigr)^{\bot}
\]
\begin{lem}
The subspaces
$\bigl\{
 F^{\lor}_{-a}(V^{\lor})\,\big|\,
 -a\in I^{\lor} \bigr\}$
are well defined,
and give a filtration of $V^{\lor}$
indexed by $(I^{\lor},\leq)$.
The decomposition
$V^{\lor}=\bigoplus_{-a\in I^{\lor}}
 V_{-a}^{\lor}$ gives a splitting
of the filtration $F^{\lor}$.
\index{dual filtration $F^{\lor}$}
\end{lem}
\pf
Let us show that
$F_{\not\geq a}(V):=\bigoplus_{b\in S(a)}V_b$
is independent of the choice of a splitting.
Let $V=\bigoplus_{a\in I} V'_{a}$ be another
splitting of $F$.
Note 
$ V_b'\subset F_b(V)
=\bigoplus_{c\leq b}V_c$.
We can observe that
$b\in S(a)$ and $c\leq b$
imply $c\in S(a)$.
Hence, 
$V'_b\subset F_{\not\geq a}(V)$
for any $b\in S(a)$,
which implies that $F_{\not\geq a}(V)$
is independent of the choice of a splitting.
Because $F^{\lor}_{-a}(V^{\lor})=
 F_{\not\geq a}(V)^{\bot}$,
we obtain that the subspaces
$F^{\lor}_{-a}(V^{\lor})$ are well defined.
If $-b\leq -a$,
we have $b\geq a$,
and $F_{\not\geq b}(V)\supset F_{\not\geq a}(V)$.
Hence, we obtain
$F^{\lor}_{-b}(V^{\lor})\subset 
 F_{-a}^{\lor}(V^{\lor})$.
It is easy to show
\[
 F^{\lor}_{-a}(V^{\lor})
=\bigoplus_{-b\leq -a}V_{-b}^{\lor}
\]
Thus, we obtain the claims of the lemma.
\hfill\qed

\vspace{.1in}
We give some property
of the induced filtration above.
\begin{lem}
We have the natural isomorphism
$\Gr^{F^{\lor}}_{-a}(V^{\lor})
\simeq
 \Gr^F_a(V)^{\lor}$.
\end{lem}
\pf
The perfect pairing of $V$ and $V^{\lor}$
induces a pairing $P$ of
$F_{a}(V)$ and $F^{\lor}_{-a}(V^{\lor})$.
By definition,
the restriction of $P$ to 
$F_{<a}(V)\otimes F^{\lor}_{-a}(V^{\lor})$
is $0$.
Let $F_{\not>a}(V):=
\bigoplus_{b\not>a} V_b$.
Then, we have $F^{\lor}_{<-a}(V^{\lor})
=F_{\not>a}(V)^{\bot}$.
Hence, the restriction of $P$
to $F_{a}(V)\otimes F^{\lor}_{<-a}(V^{\lor})$
is $0$.
Hence, we have the induced pairing
$P_a$ of $\Gr^F_a(V)$ and
$\Gr^{F^{\lor}}_{-a}(V^{\lor})$.
It is easy to check that $P_a$ is perfect
by using the induced isomorphisms
$\Gr^F_a(V)\simeq V_a $ and
$\Gr^{F^{\lor}}_{-a}(V^{\lor})
 \simeq V^{\lor}_{-a}$.
\hfill\qed

\begin{lem}
\label{lem;08.9.5.1}
Let $\varphi:(I_1,\leq_1)\lrarr (I_2,\leq_2)$ be 
a morphism of ordered sets.
Let $F_i$ $(i=1,2)$ be filtrations
of a finite dimensional vector space $V$
which are compatible over $\varphi$.
Then, the induced filtrations
$F_i^{\lor}$ $(i=1,2)$ are compatible over
the induced morphism 
$\varphi^{\lor}:(I_1^{\lor},\leq_1)
\lrarr (I_2^{\lor},\leq_2)$.
\end{lem}
\pf
It is easy to show the claim
by using the induced splitting
$V^{\lor}=\bigoplus_{a\in I_1} V_{-a}^{\lor}$.
\hfill\qed

\begin{lem}
\label{lem;08.9.5.2}
Let $V_i$ be finite dimensional vector spaces
with filtrations $F(V_i)$ indexed by
an ordered set $(I,\leq)$.
Let $f:V_1\lrarr V_2$ be a linear map
preserving filtrations.
Then, the dual
$f^{\lor}:V_2^{\lor}\lrarr V_1^{\lor}$
preserves the dual filtrations
$F^{\lor}(V^{\lor}_i)$.
\end{lem}
\pf
Because
$F^{\lor}_{-a}(V_i)
=F_{\not\geq a}(V_i)^{\bot}$,
the claim is clear.
\hfill\qed

\subsubsection{Hom-space}
\label{subsection;08.10.31.11}

We prepare a notation.
Let $(I,\leq)$ be an ordered set,
and $V_i$ $(i=1,2)$ be finite dimensional vector spaces
equipped with filtrations indexed by $(I,\leq)$.
Let $F_0Hom(V_1,V_2)$ 
(resp. $F_{<0}Hom(V_1,V_2)$)
be the vector subspace of $Hom(V_1,V_2)$
which consists of the linear maps
$f$ satisfying the following:
\[
 f\bigl(F_a(V_1)\bigr)\subset
 F_a(V_2),
\quad
 \mbox{\rm resp. }
 f\bigl(F_a(V_1)\bigr)\subset
 F_{<a}(V_2)
\]
If we take a splitting
$V_i=\bigoplus_{a\in I} V_{i,a}$ of the filtrations,
we have the following:
\begin{equation}
\label{eq;08.10.31.10}
 F_0Hom(V_1,V_2)=
 \bigoplus_{a\geq b}
 Hom(V_{1,a},V_{2,b}),
\quad
 F_{<0}Hom(V_1,V_2)=
 \bigoplus_{a>b}
 Hom(V_{1,a},V_{2,b})
\end{equation}
\index{subspace $\FzeroVichiVni$}
\index{subspace $\FminusVichiVni$}

\subsubsection{Induced filtration}

\label{subsection;07.12.12.2}

Let $(I_L,\leq)\lrarr (I_{L-1},\leq)
 \lrarr \cdots \lrarr (I_1,\leq)$
be a sequence of surjections of ordered sets.
The induced morphism
$I_j\lrarr I_k$ is denoted by $\varphi_{j,k}$.
Let $V$ be a vector space.
Let $F^{I_L}$ be a filtration indexed by $I_L$.
Then, we have the induced filtrations $F^{I_j}$
for $j=1,\ldots,L$
obtained by the procedure explained
in Section \ref{subsection;08.8.31.1}.
Moreover, we obtain the following inductive structure.
Let $\Gr^{I_{j}}(V)=\bigoplus_{\gminib\in I_j}
 \Gr^{I_j}_{\gminib}(V)$ 
 denote the graded vector space
associated to $F^{I_j}$.
For each $\gminib\in I_{j-1}$,
we have the induced filtration $F^{I_j}$ on
  $\Gr_{\gminib}^{I_{j-1}}(V)$
indexed by the ordered set
$\varphi_{j,j-1}^{-1}(\gminib)\subset I_{j}$,
and the associated graded space
$\Gr^{I_{j}}\Gr^{I_{j-1}}(V)
 =
 \bigoplus_{\gminib\in I_{j-1}}
 \bigoplus_{\gminia\in\varphi_{j,j-1}(\gminib)}
 \Gr^{I_j}_{\gminia}\Gr^{I_{j-1}}_{\gminib}(V)$
 is naturally isomorphic to $\Gr^{I_j}(V)$.

Conversely,
assume that we are given the following inductive data:
\begin{itemize}
\item
A filtration $F^{I_1}$ of $V$
indexed by $I_1$.
We put $V^{I_1}_{\gminia}:=
 \Gr^{F^{I_1}}_{\gminia}(V)$
or $\gminia\in I_1$.
\item
For each $\gminib\in I_{j-1}$,
a filtration $F^{I_j}$ of 
 $V^{I_{j-1}}_{\gminib}$
indexed by
$\varphi_{j,j-1}^{-1}(\gminib)$.
We put 
$V^{I_j}_{\gminia}:=
 \Gr^{F^{I_j}}_{\gminia}\bigl(
 V^{I_{j-1}}_{\gminib}
 \bigr)$
for $\gminia\in \varphi_{j,j-1}^{-1}(\gminib)$.
\end{itemize}
Then, we obtain the naturally induced filtration
$F^{I_j}$ of $V$
with the following property:
\begin{itemize}
\item
 $F^{I_j}$ and $F^{I_k}$ are compatible
over $\varphi_{j,k}$.
\item
Let $\Gr^{I_j}(V)=
 \bigoplus_{\gminia\in I_j}\Gr^{I_j}_{\gminia}(V)$
 denote the graded vector space associated to
 $F^{I_j}$.
 Then, $\Gr^{I_j}_{\gminia}(V)$ is naturally isomorphic
 to $V^{I_j}_{\gminia}$.
\end{itemize}
The construction is given as follows.
Assume that we have already obtained 
the desired filtration $F^{I_{j-1}}$.
Let $\gminia\in I_{j}$
and $\gminib:=\varphi_{j,j-1}(\gminia)$.
We have the natural morphism
$\pi_{\gminib}^{I_{j-1}}:
F^{I_{j-1}}_{\gminib}\lrarr
 \Gr^{I_{j-1}}_{\gminib}(V)
\simeq
 V^{I_{j-1}}_{\gminib}$.
Then, we put 
\[
 F^{I_j}_{\gminia}(V):=
 F^{I_j}_{<\gminia}(V)
+\bigl(\pi^{I_{j-1}}_{\gminib}\bigr)^{-1}
 \bigl(
 F^{I_j}_{\gminia}(V^{I_{j-1}}_{\gminib})
 \bigr)
\]
It is easy to check
$\Gr^{I_j}_{\gminia}(V)$ is isomorphic to
$\Gr^{I_j}_{\gminia}(V^{I_{j-1}}_{\gminib})
\simeq V^{I_j}_{\gminia}$.

\subsection{Orders on a good set of irregular values}
\label{subsection;10.5.28.3}

Let $\nbigx\lrarr \nbigk$ be a smooth
fibration of complex manifolds.
Let $\nbigd$ be a normal crossing hypersurface
of $\nbigx$ such that any intersections of
irreducible components are smooth over $\nbigk$.
Let $\varrho$ be a holomorphic function on $\nbigx$
such that $d\rho$ is nowhere vanishing on
$\nbigk^0:=\rho^{-1}(0)$.
We put $\nbigx^0:=\nbigx\times_\nbigk\nbigk^0$,
$\nbigd^0:=\nbigd\times_\nbigk\nbigk^0$
and $W:=\nbigx^0\cup \nbigd$.
Let $\pi:\nbigxtilde(W)\lrarr \nbigx$
denote the real blow up.

Let $\nbigi_P\subset\nbigo_{X,P}(\ast D)/\nbigo_{X,P}$
be a good set of irregular values,
where $P\in\nbigd$.
For each $Q\in \pi^{-1}(P)$,
we shall introduce an order $\leq_{Q}^{\varrho}$
(simply denoted by $\leq_Q$)
on the set $\nbigi_P$.
We can take a coordinate neighbourhood
$(\nbigx_P,z_1,\ldots,z_n)$ around $P$
such that $\nbigd_P=\nbigx_P\cap \nbigd$
is expressed as $\bigcup_{i=1}^{\ell}\{z_i=0\}$,
and that 
$\nbigi_P\subset M(\nbigx_P,\nbigd_P)/H(\nbigx_P)$.
We take a lift $\gminiatilde\in 
 M(\nbigx_P,\nbigd_P)$ for each 
$\gminia\in\nbigi_P$.
We put $W_P:=\nbigd_P\cup \nbigx_P^0$.
For each distinct $\gminia,\gminib\in\nbigi_P$,
we put
\begin{equation}
 \label{eq;10.6.11.1}
 F_{\gminia,\gminib}:=
 -\Re\bigl(\varrho^{-1}(\gminiatilde-\gminibtilde)\bigr)\,
 \bigl|\vecz^{-\ord(\gminia-\gminib)}\varrho\bigr|.
\end{equation}
It naturally induces a $C^{\infty}$-function
on $\nbigxtilde_P(W_P)$.
\begin{df}
Let $Q\in\pi^{-1}(P)$.
We say $\gminia<_Q^{\varrho}\gminib$
for distinct $\gminia,\gminib\in\nbigi_P$,
if $F_{\gminia,\gminib}(Q)<0$.
We say $\gminia\leq_Q^{\varrho}\gminib$
for $\gminia,\gminib\in\nbigi_P$,
if $\gminia<_Q^{\varrho}\gminib$ or
$\gminia=\gminib$.
\hfill\qed
\end{df}

It is easy to check that the condition
is independent of the choice of
a coordinate system and lifts $\gminiatilde$.
We will denote it by $\leq_{Q}$,
when $\varrho$ is fixed.

\begin{notation}
\label{notation;10.6.11.2}
Later, we will also use 
the relation $\leq^{\varrho}_{A}$
for a subset $A\subset\nbigxtilde_P(W_P)$,
defined as
$\gminia<_A^{\varrho}\gminib
\Longleftrightarrow
F_{\gminia,\gminib}<0$ on
$\Abar\cap\pi^{-1}(\nbigd_P)$
for $\gminia,\gminib\in\nbigi_P$,
where $\Abar$ denotes the closure of $A$
in $\nbigxtilde_P(W_P)$.
\hfill\qed
\end{notation}
\index{relation $\leq_A$}

\subsubsection{}
\label{subsection;10.6.11.3}

For any given $C^{\infty}$-manifold
(possibly with corners),
we mean by a ``compact region''
the closure of a relatively compact
connected open subset which is
$C^{\infty}$-manifold with corners.
\index{compact region}
Let $\gbigu\bigl(\nbigxtilde(W)\bigr)$ denote 
the set of compact regions in $\nbigxtilde(W)$.
For any point $Q\in\nbigxtilde(W)$,
let $\gbigu(Q,\nbigxtilde(W))$ denote the set of
$\nbigu\in\gbigu\bigl(\nbigxtilde(W)\bigr)$ such that
$Q$ is contained in the interior part of $\nbigu$.
\index{set $\gbigu\bigl(\nbigxtilde(W)\bigr)$}
\index{set $\gbigu\bigl(Q,\nbigxtilde(W)\bigr)$}

Let $\vecnbigi=(\nbigi_P\,|\,P\in\nbigd)$
be a good system of irregular values.
For $Q\in\pi^{-1}(P)$,
let $\gbigu(Q,\nbigxtilde(W),\vecnbigi)$
denote the set of 
$\nbigu\in\gbigu(Q,\nbigxtilde(W))$
such that,
for any $\gminia,\gminib\in\nbigi_P$,
we have
$\gminia\leq_Q^{\varrho}\gminib
\Longleftrightarrow
 F_{\gminia,\gminib}\leq 0$ on $\nbigu$.
The following claims are clear.
\begin{itemize}
\item
If $\nbigu\in\gbigu(Q,\nbigxtilde(W))$
is sufficiently small, it is contained in
$\gbigu(Q,\nbigxtilde(W),\vecnbigi)$.
\item
For a given 
$\nbigu\in\gbigu(Q,\nbigxtilde(W),\vecnbigi)$,
we have $\gminia\leq^{\varrho}_{\nbigu}\gminib$
if and only if
$F_{\gminia,\gminib}\leq 0$ on 
$\nbigu\setminus\pi^{-1}(W)$.
\item
Let $\nbigu\in\gbigu(Q,\nbigxtilde(W),\vecnbigi)$.
For $Q'\in\nbigu\cap\pi^{-1}(\nbigd)$,
the natural map
\[
 \bigl(
 \nbigi_{\pi(Q)},\leq_Q^{\varrho}
 \bigr)
\lrarr
 \bigl(\nbigi_{\pi(Q')},\leq_{Q'}^{\varrho}\bigr) 
\]
is order preserving.
\end{itemize}
We will often denote 
$\gbigu(Q,\nbigxtilde(W))$
and $\gbigu(Q,\nbigxtilde(W),\vecnbigi)$
by $\gbigu(Q)$
and $\gbigu(Q,\vecnbigi)$
if there is no risk of confusion.

\subsection{Holomorphic functions and
 vector bundles on real blow up}
\label{subsection;08.9.28.221}

\subsubsection{Functions and vector bundles}

Let us recall the framework of asymptotic analysis
in \cite{sabbah4}. 
See Chapter II.1 of \cite{sabbah4}
for more details and precision.
Let $X$ be a complex manifold,
and let $D$ be a normal crossing divisor of $X$.
In the following, $\Xtilde(D)$ denote 
the fiber product, taken over $X$,
of the real blow up of
$X$ along the irreducible components of $D$,
which is called the real blow up of $X$ along $D$
in this paper.
It is naturally a $C^{\infty}$-manifold with corners.
Hence, $C^{\infty}$-functions on open subsets of
$\Xtilde(D)$ make sense.
\index{real blow up $\Xtilde(D)$}
In the case $X=\Delta^n$ and 
$D=\bigcup_{i=1}^{\ell}\{z_i=0\}$,
we may use the local coordinate 
$(r_1,\theta_1,\ldots,r_{\ell},\theta_{\ell},
 z_{\ell+1},\ldots,z_n)$,
where $(r_i,\theta_i)$ are given by
$z_i=r_i\, e^{\sqrt{-1}\theta_i}$.

We have the $C^{\infty}$-vector bundle
$\Omega^{0,1}_{X}(\log D)$ on $X$
locally generated by
$d\zbar_i/\zbar_i$ $(i=1,\ldots,\ell)$
and $d\zbar_i$ $(i=\ell+1,\ldots,n)$.
The pull back of 
$\Omega^{0,1}_{X}(\log D)$
via $\pi$ is also denoted by the same symbol.
For any open subset $U\subset \Xtilde(D)$,
we have the well defined differential operator
$\delbar:
 C^{\infty}(U)\lrarr
 C^{\infty}\bigl(U,
 \Omega^{0,1}_{X}
 (\log D)\bigr)$.
We say $f\in C^{\infty}(U)$ is holomorphic,
if $\delbar f=0$.
It is easy to see that
$f\in C^{\infty}(U)$ is holomorphic if and only if
$f_{|U\cap(X-D)}$ is holomorphic in the standard sense.
The space of holomorphic functions on $U$
is denoted by $\nbiga(U)$
or $\nbigo(U)$.

Let $Z$ be a locally closed  subset of $\Xtilde(D)$.
Let $C^{\infty}(\Zhat)$ denote
the space of the $C^{\infty}$-functions on $Z$
in the sense of Whitney
(see Chapters I.2 and I.4 of \cite{malgrange2}).
Namely it is the space of the $\infty$-jets of $X$
on $Z$,
satisfying some compatibility condition
as in Theorem 2.2 of \cite{malgrange2}.
Let $U$ be an open subset of $\Xtilde(D)$
such that $Z\subset U$.
Then, $C^{\infty}(\Zhat)$ is the same as 
the limit $\varprojlim C^{\infty}(U)/I_Z^n$,
where $I_Z$ denotes the ideal of $C^{\infty}(U)$
which consists of the functions $f$ such that
$f_{|Z}=0$.
Equivalently, let $I^{(\infty)}_Z$ denote the ideal of
$C^{\infty}(U)$ which consists of the functions $f$
such that $(\del^J f)_{|Z}=0$ 
for any derivation $\del^J$ with any degree.
Then, the natural morphism
$C^{\infty}(U)/I_Z^{(\infty)}
\lrarr C^{\infty}(\Zhat)$ is isomorphic.
(See the proof of Theorem 4.1 of \cite{malgrange2}.)
We prefer to regard
elements of $C^{\infty}(\Zhat)$ 
as $C^{\infty}$-functions on the space $\Zhat$
which is the completion of $\Xtilde(D)$ along $Z$.
The image of a function $f\in C^{\infty}(U)$
to $C^{\infty}(\Zhat)$ is denoted by
$f_{|\Zhat}$.

The differential operator
$\delbar:C^{\infty}(U)
\lrarr C^{\infty}\bigl(U,
 \Omega_{X}^{0,1}(\log D)\bigr)$
induces 
$\delbar:C^{\infty}(\Zhat) \lrarr
C^{\infty}\bigl(\Zhat,
 \Omega_{X}^{0,1}(\log D)\bigr)$.
We say that $f\in C^{\infty}(\Zhat)$ is holomorphic 
if $\delbar f=0$.
The space of holomorphic functions on $\Zhat$
is denoted by $\nbiga(\Zhat)$ or $\nbigo(\Zhat)$.
We will use the following fact
without mention.
\begin{lem}
Let $X=\Delta^n$ and
$D=\bigcup_{j=1}^{\ell}\{z_j=0\}$.
Let $W$ be a closed region of $\Xtilde(D)$,
and let $Z:=U\cap \pi^{-1}(D)$.
Let $f$ be a holomorphic function on $W-Z$
such that $|f|=O\bigl(\prod_{j=1}^{\ell}|z_j|^N\bigr)$
for any $N$.
Then, $f$ gives a holomorphic function on $W$
such that $f_{|\Zhat}=0$.

Conversely, 
if $f$ is a holomorphic function on $W$
such that $f_{|\Zhat}=0$,
it satisfies $|f|=O\bigl(\prod_{j=1}^{\ell}|z_j|^N\bigr)$
for any $N$.
\end{lem}
\pf
Let $f$ be a function as in the first claim.
By using Cauchy's formula,
we obtain 
$\del^Jf=O\bigl(\prod_{j=1}^{\ell}|z_j|^N\bigr)$
for any $N$
and for any derivation $\del^J$
with any degree.
(See the proof of Theorem 8.8 of \cite{wasow}.)
Then, the first claim of the lemma follows.
The second claim is obvious.
\hfill\qed

\vspace{.1in}

We recall Borel-Ritt type theorem,
due to Majima \cite{majima} for his strongly
asymptotically developpable functions,
and due to  Sabbah  \cite{sabbah4}
in this framework.

\begin{prop}[Theorem 2.2 of 
\cite{majima},
Proposition 1.1.16 of \cite{sabbah4}]
\label{prop;07.11.12.2}
Let $Q\in \pi^{-1}(D)$.
Let $f\in \nbiga(\Zhat)$,
where $Z$ denotes a closed neighbourhood of $Q$
in $\pi^{-1}(D)$.
Then, we have a neighbourhood  $U$ of $Q$ in $\Xtilde(D)$
and a holomorphic function $F\in \nbiga(U)$
such that $F_{|\Zhat_U}=f_{|\Zhat_U}$,
where $Z_U:=Z\cap U$.

In particular,
$\nbiga(\Zhat)$ is isomorphic to
$\varprojlim \nbiga(U)/\nbigi_Z^n$,
where $\nbigi_Z$ denote the ideal of $\nbiga(U)$
which consists of the holomorphic functions $f$
such that $f_{|Z}=0$.
\hfill\qed
\end{prop}

\begin{rem}
Assume that 
$X=\Delta^n$ and $D=\{z_1=0\}$.
Let $S$ be a sector in $X-D$ and
let $\Sbar$ be the closure of $S$
in $\Xtilde(D)$.
We set $Z:=\Sbar\cap\pi^{-1}(D)$.
Let $f\in\nbiga(\Zhat)$.
If we shrink $S$ in the radius direction
(Definition {\rm\ref{df;08.12.23.1}}),
we can take $F\in\nbiga(\Sbar)$
such that $F_{|\Zhat}=f$,
which can be shown
by the argument 
in the proof of Theorem {\rm 9.3} of {\rm\cite{wasow}}.
\hfill\qed
\end{rem}

Let $U$ denote an open subset of $\Xtilde(D)$.
Since $\Xtilde(D)$ is a $C^{\infty}$-manifold
with corner,
we have the well defined concept of
$C^{\infty}$-vector bundle on $U$.
Let $E$ be a $C^{\infty}$-bundle on $U$.
The space of $C^{\infty}$-sections of $E$
on $U$ is denoted by $C^{\infty}(U,E)$.
Assume we are given a differential operator
$\delbar_E:C^{\infty}(U,E)\lrarr
 C^{\infty}\bigl(U,
 E\otimes
 \Omega^{0,1}_{X}(\log D)\bigr)$
such that
(i)  
$ \delbar_E(f\cdot s)
=\delbar(f)\cdot s+f\cdot \delbar_E(s)$ holds
for any $f\in C^{\infty}(U)$
and $s\in C^{\infty}(U,E)$,
(ii) $\delbar_E\circ\delbar_E=0$.
A section $s$ is called holomorphic,
if $\delbar_E(s)=0$ holds.
We always assume the existence of 
holomorphic local frames of $E$
around each point of $\Xtilde(D)$.

\subsubsection{Some equivalence}

Let $\Mod(\nbigo_X)$ denote
the category of $\nbigo_X$-modules.
Let $\Mod_{lf}(\nbigo_X)$ denote
the full subcategory of $\Mod(\nbigo_X)$
of locally free $\nbigo_X$-sheaves
of finite ranks.
To state some equivalence,
which would be useful
for understanding of our later construction,
let us introduce a category $\nbigc$.
\begin{itemize}
\item
 Objects of $\nbigc$ are tuples
 $(\nbigf_{\Xtilde(D)},\nbigf_{\Dhat},\iota)$:
\begin{itemize}
\item
 $\nbigf_{\Xtilde(D)}$ is 
 a locally free $\nbigo_{\Xtilde(D)}$-module.
\item
 $\nbigf_{\Dhat}$ is
 a locally free $\nbigo_{\Dhat}$-module.
\item
 $\iota$ is an isomorphism
 $\pi^{\ast}\nbigf_{\Dhat}$
 and the completion of $\nbigf_{\Xtilde(D)}$
 along $\pi^{-1}(D)$.
\end{itemize}
\item
 Morphisms
 $(\nbigf^{(1)}_{\Xtilde(D)},
 \nbigf^{(1)}_{\Dhat},\iota^{(1)})
\lrarr
 (\nbigf^{(2)}_{\Xtilde(D)},
 \nbigf^{(2)}_{\Dhat},\iota^{(2)})$
 are pairs of morphisms
 $f_{\Xtilde(D)}:
 \nbigf^{(1)}_{\Xtilde(D)}
\lrarr
 \nbigf^{(2)}_{\Xtilde(D)}$
and
  $f_{\Dhat}:
 \nbigf^{(1)}_{\Dhat}
\lrarr
 \nbigf^{(2)}_{\Dhat}$
which are compatible with
 $\iota^{(1)}$ and $\iota^{(2)}$.
\end{itemize}

We have the functors
$F_1:\Mod_{lf}(\nbigo_X)\lrarr\nbigc$
given by 
$F_1(\nbigf):=\bigl(
 \pi^{\ast}\nbigf,
 \nbigf_{|\Dhat},
 \iota
 \bigr)$,
where $\iota$ denotes the natural isomorphism
and $\pi^{\ast}\nbigf=
 \pi^{-1}\nbigf\otimes_{\pi^{-1}\nbigo_X}
 \nbigo_{\Xtilde(D)}$.
We also have the functor
$F_2:\nbigc\lrarr \Mod(\nbigo_X)$
given by
$F_2(\nbigf_{\Xtilde(D)},\nbigf_{\Dhat},\iota)
=\pi_{\ast}(\nbigf_{\Xtilde(D)})$,
where $\pi_{\ast}$ denotes the push-forward
of sheaves.

\begin{prop}
\label{prop;08.9.4.102}
For any 
$(\nbigf_{\Xtilde(D)},\nbigf_{\Dhat},\iota)
 \in\nbigc$,
$\pi_{\ast}\nbigf_{\Xtilde(D)}$ is 
$\nbigo_{X}$-locally free,
and hence the functor $F_2$ factors through 
$\Mod_{lf}(\nbigo_X)$.
Moreover,
we have an isomorphism
$\nbigf_{\Dhat}\simeq
 \pi_{\ast}(\nbigf_{\Xtilde(D)})_{|\Dhat}$,
which induces $\iota$.
(Note that it is unique.)
\end{prop}
\pf
Let us show the first claim.
Since it is a local property,
we consider the case
$X=\Delta^n$ and 
$D=\bigcup_{i=1}^{\ell}\{z_i=0\}$.
Moreover,
we will replace $X$ with 
a smaller neighbourhood of the origin
without mention,
if it is necessary.

Let $(\nbigf_{\Xtilde(D)},\nbigf_{\Dhat},\iota)
 \in\nbigc$.
Let $\vecvhat$ be a frame of $\nbigf_{\Dhat}$.
Let $Q$ be any point of $\pi^{-1}(D)$.
According to Proposition \ref{prop;07.11.12.2},
we can take a small multi-sector $S$ of $X-D$
and a frame $\vecv_S$
of $\nbigf_{\Xtilde(D)|\Sbar}$
such that
(i) $\vecv_{S|\Zhat}=\pi^{-1}\vecvhat$,
(ii) $Q$ is contained in the interior part of $\Sbar$,
where $\Sbar$ denotes the closure of $S$
in $\Xtilde(D)$,
and $Z:=\Sbar\cap\pi^{-1}(D)$.

We take a covering
$X-D=\bigcup_{p=1}^NS_p$ by multi-sectors
on which we have frames $\vecv_{S_p}$ as above.
We take a partition of unity
$\bigl(\chi_p\,\big|\,p=1,\ldots,N\bigr)$
subordinated to the covering 
such that 
$P\chi_p$
are polynomial order in
$|z_j^{-1}|$ $(j=1,\ldots,\ell)$
for each differential operator $P$ on $X$.
We obtain a $C^{\infty}$-frame
$\vecv_{C^{\infty}}$ of $\nbigf_{\Xtilde(D)}$
given by
$\vecv_{C^{\infty}}:=
 \sum \vecv_{S_p}\,\chi_p$,
or more precisely,
$v_{C^{\infty},i}=
 \sum \chi_p\, v_{S_p,i}$.

Let $S$ be a small multi-sector of $X-D$,
and let $\vecv_S$ be as above.
Let $C_S$ be a matrix-valued function on $\Sbar$
determined by 
$\vecv_{C^{\infty}}
=\vecv_S\, (I+C_S)$,
where $I$ is the identity matrix.
\begin{lem}
\label{lem;08.9.4.101}
$C_{S|\Zhat}=0$.
\end{lem}
\pf
Let $C_p$ be determined on $S\cap S_p$
by $\vecv_{S_p}=\vecv_S\, (I+C_p)$.
Then, $C_{p|\widehat{Z\cap Z_p}}=0$.
By our choice of the partition of unity $(\chi_p)$,
we obtain $C_{S|\Zhat}=0$.
\hfill\qed

\vspace{.1in}

Let $A$ be the matrix valued $(0,1)$-form
determined by
$\delbar\vecv_{C^{\infty}}
=\vecv_{C^{\infty}}\, A$.
\begin{lem}
\label{lem;08.9.4.100}
$A$ gives a matrix-valued $C^{\infty}$-function
on $X$, 
and $A_{|\Dhat}=0$.
\end{lem}
\pf
By Lemma \ref{lem;08.9.4.101},
for each small multi-sector $S$,
$\delbar\vecv_{C^{\infty}|\Sbar}
=\vecv_S\,\delbar C_S
=\vecv_{C^{\infty}|\Sbar}\,
 (I+C_S)^{-1}\delbar C_S$.
Hence, $A_{|\Zhat}=0$
for each $S$.
Then, we obtain the claim of the lemma.
\hfill\qed

\vspace{.1in}

Let $E$ be the $C^{\infty}$-vector bundle 
on $X$,
which is the extension of $\nbigf_{\Xtilde(D)|X-D}$
by the frame $\vecv_{C^{\infty}}$.
According to Lemma \ref{lem;08.9.4.100},
the $\delbar$-operator of
$E_{|X-D}$ is naturally extended
to a $\delbar$-operator of $E$
in $C^{\infty}$,
i.e.,
the holomorphic structure of $E_{|X-D}$
is prolonged to that of $E$.
Let $\nbige$ denote the sheaf of
holomorphic sections of $E$.
By construction,
we have the natural isomorphism
$\nbige\simeq \pi_{\ast}\nbigf_{\Xtilde(D)}$.
Thus, the first claim of Proposition
\ref{prop;08.9.4.102} is proved.

By Lemma \ref{lem;08.9.4.100},
$\vecv_{C^{\infty}|\Dhat}$
naturally gives a frame of $\nbige_{|\Dhat}$.
By the frames $\vecv_{C^{\infty}|\Dhat}$
and $\vecvhat$,
we obtain an isomorphism
$\nbigf_{\Dhat}\simeq\nbige_{|\Dhat}$,
which induces $\iota$.
\hfill\qed

\begin{cor}
\label{cor;08.9.5.1}
The functors $F_1$ and $F_2$
give equivalence of categories,
and they are mutually quasi-inverse.
\hfill\qed
\end{cor}

We give some complements.
\begin{lem}
$F_i$ $(i=1,2)$ preserve
dual, tensor product and direct sum.
\end{lem}
\pf
By construction,
the functor $F_1$ preserves
dual, tensor product and direct sum.
Then, by Proposition \ref{prop;08.9.4.102},
$F_2$ also preserves
dual, tensor product and direct sum.
\hfill\qed

\subsubsection{Decent of 
a family of meromorphic $\varrho$-connections}

We use the setting in Subsection
\ref{subsection;10.5.28.3}.
Let $\nbigf$ be a locally free
$\nbigo_{\nbigx}$-module of finite rank.
We set $\pi^{\ast}\nbigf:=
 \nbigo_{\nbigxtilde(W)}
 \otimes_{\pi^{-1}\nbigo_{\nbigx}}
 \pi^{-1}\nbigf$.

\begin{lem}
\label{lem;08.9.5.12}
Let $\DD_{\nbigxtilde(W)}:
 \pi^{\ast}\nbigf\lrarr 
 \pi^{\ast}\bigl(
 \nbigf\otimes 
 \Omega^{1}_{\nbigx/\nbigk}(\ast \nbigd)
  \bigr)$ be 
a meromorphic $\varrho$-connection
of $\pi^{\ast}\nbigf$.
Then, there exists
a meromorphic $\varrho$-connection
$\DD$ of $\nbigf$
such that
$\DD_{\nbigxtilde(W)}=\pi^{\ast}\DD$.

Assume moreover that we are given 
a meromorphic $\varrho$-connection
$\DD_{\What}$ of $\nbigf_{|\What}$
such that $\pi^{\ast}\DD_{\What}$
equals the completion of
$\DD_{\nbigxtilde(W)}$ along $\pi^{-1}(W)$.
Then, the completion of $\DD$
equals $\DD_{\What}$.
\end{lem}
\pf
Let $P$ be any point of $\nbigx$.
Let $f$ be any section of $\nbigf$
on a neighbourhood $\nbigx_P$ of $P$ in $\nbigx$.
If $\nbigx_P$ is shrinked appropriately,
there exists a number $N$ such that
$\DD_{\nbigxtilde(W)}(\pi^{\ast}f)$
gives a section of
$\pi^{\ast}\bigl(
 \nbigf\otimes
 \Omega^1_{\nbigx/\nbigk}(N \nbigd)\bigr)$
on $\pi^{-1}(\nbigx_P)$.
Because 
$\pi_{\ast}\pi^{\ast}\bigl(
 \nbigf\otimes
 \Omega^1_{\nbigx/\nbigk}(N\nbigd) \bigr)
 \simeq
 \nbigf\otimes
 \Omega^1_{\nbigx/\nbigk}(N\nbigd)$,
$\DD_{\nbigxtilde(\nbigd)}(\pi^{\ast}f)$
naturally gives a section of
$\nbigf\otimes
 \Omega^1_{\nbigx/\nbigk}(N D)$ on 
$\nbigx_P$.
Thus, we obtain a map of sheaves
$\DD:\nbigf\lrarr 
 \nbigf\otimes
 \Omega^1_{\nbigx/\nbigk}(\ast \nbigd)$.
We can observe that
it is a meromorphic $\varrho$-connection,
and satisfies
$\pi^{\ast}\DD=\DD_{\nbigxtilde(D)}$,
and thus the first claim is proved.

If we are given a 
meromorphic $\varrho$-connection
$\DD_{\What}$ of $\nbigf_{|\What}$
as in the second claim,
we have
$\pi^{\ast}(\DD_{\What})
=\DD_{\nbigxtilde(W)|\widehat{\pi^{-1}(W)}}
=\pi^{\ast}(\DD_{|\What})$,
and hence $\DD_{\What}=\DD_{|\What}$.
Thus we obtain the second claim.
\hfill\qed

\section{Good meromorphic $\varrho$-flat bundle}
\label{section;10.5.10.5}
In this section,
we state a theorem and some propositions.
They will be proved in Section \ref{section;10.5.18.5}.
(See Subsection \ref{subsection;10.5.18.6}.)
We use the following setting
unless otherwise specified.
Let $\nbigx\lrarr \nbigk$ be a smooth
fibration of complex manifolds.
Let $\nbigd$ be a normal crossing hypersurface
of $\nbigx$ such that any intersections of
irreducible components are smooth over $\nbigk$.
Let $\varrho$ be a holomorphic function on $\nbigx$
such that $d\rho$ is nowhere vanishing on
$\nbigk^0:=\rho^{-1}(0)$.
We put $\nbigx^0:=\nbigx\times_\nbigk\nbigk^0$,
$\nbigd^0:=\nbigd\times_\nbigk\nbigk^0$
and $W:=\nbigx^0\cup \nbigd$.
Let $\pi:\nbigxtilde(W)\lrarr \nbigx$
denote the real blow up.

\subsection{Full Stokes filtration}

Let $(\nbige,\DD)$
be an unramifiedly good meromorphic 
$\varrho$-flat bundle on $(\nbigx,\nbigd)$.
We put $\pi^{\ast}\nbige:=
 \pi^{-1}\nbige\otimes_{\pi^{-1}\nbigo_\nbigx}
 \nbigo_{\nbigxtilde(W)}$,
which is an 
$\nbigo_{\nbigxtilde(W)}(\ast \nbigd)$-module.
For each $Q\in \pi^{-1}(\nbigd)$,
let $\pi^{\ast}\nbige_Q$
denote the germ at $Q$,
and $\pi^{\ast}\nbige_{|\Qhat}$
denote the formal completion.
In the following,
$\Irr(\DD,\pi(Q))$ is also denoted by
$\Irr(\DD,Q)$ for simplicity of the description.

The irregular decomposition
of $(\nbige,\DD)_{|\widehat{\pi(Q)}}$
induces
$\pi^{\ast}\nbige_{|\Qhat}=
 \bigoplus_{\gminia\in\Irr(\DD,Q)}
 \lefttop{Q}\nbigehat_{\gminia}$.
We put
\[
 \nbigfhat^Q_{\gminia}\bigl(
 \pi^{\ast}\nbige_{|\Qhat}\bigr):=
 \bigoplus_{\gminib\leq^{\varrho}_Q\gminia}
 \lefttop{Q}\nbigehat_{\gminia}.
\]

\begin{thm}
\label{thm;10.5.12.110}
For any $Q\in \pi^{-1}(\nbigd)$,
there exists a $\DD$-flat filtration
$\nbigftilde^Q$ of $\pi^{\ast}\nbige_Q$
indexed by
$\bigl(\Irr(\DD,Q),
 \leq^{\varrho}_Q\bigr)$
with the following property:
\begin{itemize}
\item
 $\Gr^{\nbigftilde^Q}_{\gminia}
 (\pi^{\ast}\nbige_Q)$
 are free $\nbigo_{\nbigxtilde(W)}$-modules,
 and 
 $\nbigftilde^Q_{\gminia}(\pi^{\ast}\nbige_Q)
 _{|\Qhat}
=\nbigfhat^Q_{\gminia}(\pi^{\ast}\nbige_{|\Qhat})$.
\item
 Take $\nbigu_Q\in\gbigu(Q,\Irr(\DD))$
 such that
 $\nbigftilde^Q$ of 
$\pi^{\ast}\nbige_{|\nbigu_Q}$ is given.
 For any $Q'\in \nbigu_Q\cap\pi^{-1}(\nbigd)$,
 we have the induced filtration
 $\nbigftilde^{Q}$ of $\pi^{\ast}\nbige_{Q'}$.
 Then, 
 $\nbigftilde^Q$ and $\nbigftilde^{Q'}$
 are compatible over
 $\bigl(\Irr(\DD,Q),
 \leq_Q^{\varrho}\bigr)
\lrarr
 \bigl(\Irr(\DD,Q'),
 \leq_{Q'}^{\varrho}\bigr)$.
\end{itemize}
The conditions characterize the system of
filtrations 
$\bigl(\nbigftilde^Q\,
 \big|\,Q\in \pi^{-1}(\nbigd)\bigr)$.
If $\varrho(Q)\neq 0$,
the first property characterizes
the filtration $\nbigftilde^Q$.
\end{thm}

The filtrations $\nbigftilde^Q$ are called
the full Stokes filtration of $\nbige$ at $Q$.
It is also called the Stokes filtration,
if there is no risk of confusion.

\begin{rem}
Such a filtration appeared 
in the classical works
on meromorphic flat bundles
on curves,
for example 
{\rm\cite{malgrange3}} and
{\rm\cite{malgrange_book}}.
(See also 
{\rm\cite{Deligne-Malgrange-Ramis}}.)
T. Pantev informed that
it is called ``Deligne-Malgrange filtration''.
We keep our terminology ``Stokes filtration'',
partially because we use ``Deligne-Malgrange''
in a different meaning.
\hfill\qed
\end{rem}

Let $P\in\nbigd$
with a small  neighbourhood $\nbigx_P$.
Let $\nbigu\in\gbigu(\nbigxtilde_P(W_P))$.
Let $\leq_{\nbigu}^{\varrho}$ be the order
on $\Irr(\DD,P)$ defined as follows:
\[
 \gminia\leq_{\nbigu}^{\varrho}\gminib
\Longleftrightarrow
 \gminia_Q\leq_{Q}^{\varrho}\gminib_Q
 \,\,(\forall Q\in \nbigu\cap\pi^{-1}(\nbigd))
\]
Here, $\gminia_Q,\gminib_Q\in
 \Irr(\DD,\pi(Q)))$ be the induced elements.
The use is consistent with Notation 
{\rm\ref{notation;10.6.11.2}}.

If $\pi^{\ast}\nbige_{|\nbigu}$
has a $\DD$-flat filtration
$\nbigftilde^{\nbigu}$ 
indexed by $\bigl(\Irr(\DD,P),
 \leq_{\nbigu}^{\varrho}\bigr)$
with the following property,
which is called a full Stokes filtration
on $\nbigu$:
\begin{itemize}
\item
Let $Q$ be any point of
$\nbigu\cap\pi^{-1}(\nbigd)$.
We have the induced filtration
$\nbigftilde^{\nbigu}$ of
$\pi^{\ast}\nbige_Q$.
Then, $\nbigftilde^{\nbigu}$
and $\nbigftilde^Q$ are compatible over
$\bigl(
 \Irr(\DD,P),\leq_{\nbigu}
 \bigr)\lrarr
 \bigl(
 \Irr(\DD,\pi(Q)),\leq_Q
 \bigr)$.
\end{itemize}
Such a filtration is uniquely determined if it exists,
by Lemma {\rm\ref{lem;10.5.12.22}}.

For a multi-sector $S$ of $\nbigx\setminus W$,
if a full Stokes filtration for $\Sbar$ exists,
if is also called Stokes filtration for $S$.

\subsection{Functoriality}
\label{subsection;10.5.12.15}

Let $(\nbige_i,\DD_i)$ $(i=1,2)$
be unramifiedly good meromorphic
$\varrho$-flat bundles on $(\nbigx,\nbigd)$.

\begin{prop}
Let $F:\nbige_1\lrarr \nbige_2$ be a 
$\DD$-flat morphism.
For simplicity, we assume that
$\Irr(\DD_1)
\cup
 \Irr(\DD_2)$
is also good,
i.e.,
$\Irr(\DD_1,P)\cup\Irr(\DD_2,P)$
is good for each $P\in \nbigd$.
Then, for each $Q\in \pi^{-1}(\nbigd)$,
the induced morphism
$\pi^{\ast}\nbige_{1\,Q}\lrarr
 \pi^{\ast}\nbige_{2\,Q}$ is compatible with
the full Stokes filtrations.
\end{prop}

\begin{prop}
If $\Irr(\DD_1)\otimes\Irr(\DD_2)$ 
is good,
we have
\[
 \nbigftilde^{Q}_{\gminia}\bigl(
 \pi^{\ast}(\nbige_1\otimes\nbige_2)_Q
 \bigr)
=\sum_{\gminia_1+\gminia_2=\gminia}
 \nbigftilde^Q_{\gminia_1}
 (\pi^{\ast}\nbige_{1,Q})
\otimes
 \nbigftilde^Q_{\gminia_2}
 (\pi^{\ast}\nbige_{2,Q})
\]
If $\Irr(\DD_1)\oplus\Irr(\DD_2)$ is good,
we have
\[
 \nbigftilde^{Q}_{\gminia}\bigl(
 \pi^{\ast}(\nbige_1\oplus\nbige_2)_Q
 \bigr)
=\nbigftilde^Q_{\gminia}
 (\pi^{\ast}\nbige_{1,Q})
\oplus
 \nbigftilde^Q_{\gminia}
 (\pi^{\ast}\nbige_{2,Q})
\]
\end{prop}

\begin{prop}
Let $(\nbige,\DD)$ be an unramifiedly
good $\varrho$-meromorphic flat bundle
on $(\nbigx,\nbigd)$.
The Stokes filtration of
$(\nbige^{\lor},\DD^{\lor})$
at $Q\in\pi^{-1}(\nbigd)$
is given as follows:
\[
 \nbigftilde^Q_{\gminia}\bigl(
 \pi^{\ast}\nbige^{\lor}_{Q}
 \bigr)
=
 \Bigl(
 \sum_{\gminib\not\geq-\gminia}
 \nbigftilde^Q_{\gminib}(\pi^{\ast}\nbige_Q)
 \Bigr)^{\bot}
\]
\end{prop}

\subsection{Characterization by growth order
 (the case $\varrho(Q)\neq 0$)}

Let $(\nbige,\DD)$
be an unramifiedly good meromorphic
$\varrho$-flat bundle on $(\nbigx,\nbigd)$.
Let $Q\in \pi^{-1}(\nbigd)$.
Take a small neighbourhood 
$\nbigu\in\gbigu(Q,\Irr(\DD))$
and a frame $\vecv$ of $\nbige_{|\nbigu}$.
A $\DD$-flat section of 
$\nbige_{|\nbigu\setminus\pi^{-1}(W)}$
is expressed as $f=\sum f_j\, v_{j}$,
where $f_j\in\nbigo_{\nbigu\setminus\pi^{-1}(W)}$.
Let $\vecf$ denote the tuple $(f_j)$.
\begin{prop}
\label{prop;10.5.12.203}
Assume $\varrho(Q)\neq 0$.
We have
$f\in \nbigftilde^Q_{\gminia}
 \nbige_{|\nbigu\setminus\pi^{-1}(\nbigd)}$
if and only if
$\bigl|
 \exp(-\varrho^{-1}\gminia)\vecf
 \bigr|$ is of polynomial order
on $\nbigu\setminus \pi^{-1}(\nbigd)$.
\end{prop}

\begin{rem}
Let $(z_1,\ldots,z_n)$ be a coordinate system
around $\pi(Q)$ such that
$\nbigd$ is expressed as
$\bigcup_{i=1}^{\ell}\{z_i=0\}$ around $\pi(Q)$.
We say that a function $F$ on 
$\nbigu_Q\setminus \pi^{-1}(\nbigd)$ is
of polynomial order,
if $|F|=O\bigl(\prod_{i=1}^{\ell}|z_i|^{-N}\bigr)$
for some $N>0$.
\index{polynomial order}
\hfill\qed
\end{rem}

\subsection{The associated graded bundle}
\label{subsection;10.5.21.10}

Let $P\in \nbigd$.
For each $Q\in\pi^{-1}(P)$,
we take a small neighbourhood 
$\nbigu_Q$ of $Q$,
then and we have 
$\Gr^{\nbigftilde^Q}\bigl(
 \pi^{\ast}\nbige_{|\nbigu_Q},\DD\bigr)$
on $\nbigu_Q$.
Although the filtration $\nbigftilde^Q$
depends on $Q$,
the system 
$\bigl(\nbigftilde^Q\,\big|\,Q\in\pi^{-1}(P)\bigr)$
satisfies the compatibility condition.
Hence,
by varying $Q\in\pi^{-1}(P)$
and gluing
 $\Gr^{\nbigftilde^Q}\bigl(
 \pi^{\ast}\nbige_{|\nbigu_Q},\DD\bigr)$,
we obtain the associated graded 
locally free $\nbigo_{\nbigxtilde(W)}(\ast \nbigd)$-module
with a flat $\varrho$-connection
\[
 \Gr^{\nbigftilde}(\pi^{\ast}\nbige_{\pi^{-1}(P)},
\DD)
=\bigoplus_{\gminia\in\Irr(\DD,P)}
 \bigl(
 \Gr^{\nbigftilde}_{\gminia}
 (\pi^{\ast}\nbige_{\pi^{-1}(P)}),
 \DD_{\gminia}
\bigr)
\]
on a neighbourhood of $\pi^{-1}(P)$.
By taking the push-forward via $\pi$,
we obtain an $\nbigo_\nbigx$-module
with a flat $\varrho$-connection
\[
 \Gr^{\nbigftilde}(\nbige_P,\DD)
=\bigoplus_{\gminia\in\Irr(\DD,P)}
 \bigl(
 \Gr^{\nbigftilde}_{\gminia}(\nbige_P),
 \DD_{\gminia}
 \bigr)
\]
on a neighbourhood $\nbigx_P$ of $P$.

\begin{prop}
\label{prop;10.5.12.205}
$\Gr^{\nbigftilde}(\nbige_P,\DD)$
is a graded meromorphic $\varrho$-flat bundle
satisfying 
$\Gr^{\nbigftilde}(\pi^{\ast}\nbige_{\pi^{-1}(P)},
\DD)
\simeq
\pi^{\ast}
\Gr^{\nbigftilde}(\nbige_P,\DD)$.
We have a canonical isomorphism
$\Gr^{\nbigftilde}(\nbige_P,\DD)_{|\Phat}
\simeq
 (\nbige_P,\DD)_{|\Phat}$.
\end{prop}

In particular, each
$\bigl(
 \Gr^{\nbigftilde}_{\gminia}(\nbige_P),
 \DD_{\gminia}
 \bigr)
\otimes L(-\gminia)$ is a regular meromorphic
$\varrho$-flat bundle.

\subsubsection{Functoriality}

Let $(\nbige_1,\DD_1)
\lrarr(\nbige_2,\DD_2)$
be a morphism of unramifiedly good meromorphic
$\varrho$-flat bundles.
For simplicity, we assume
$\Irr(\DD_1,P)\cup
 \Irr(\DD_2,P)$ is good.
Then, we have the induced morphism
on a neighbourhood of $P$:
\[
 \Gr^{\nbigftilde}(\nbige_{1,P})
\lrarr
 \Gr^{\nbigftilde}(\nbige_{2,P})
\]
If $\Irr(\DD_1,P)\otimes
 \Irr(\DD_2,P)$ is good,
we have the following canonical isomorphism
on a neighbourhood of $P$:
\[
 \Gr^{\nbigftilde}\bigl(
 (\nbige_1\otimes\nbige_2)_P
 \bigr)
\simeq
 \Gr^{\nbigftilde}(\nbige_{1,P})\otimes
 \Gr^{\nbigftilde}(\nbige_{2,P})
\]
If $\Irr(\DD_1,P)\oplus
 \Irr(\DD_2,P)$ is good,
we have the following canonical isomorphism
on a neighbourhood of $P$:
\[
 \Gr^{\nbigftilde}\bigl(
 (\nbige_1\oplus\nbige_2)_P
 \bigr)
\simeq
 \Gr^{\nbigftilde}(\nbige_{1,P})\oplus
 \Gr^{\nbigftilde}(\nbige_{2,P})
\]
For an unramifiedly good meromorphic
$\varrho$-flat bundle 
$(\nbige,\DD)$,
we have the following canonical isomorphism
on a neighbourhood of $P$:
\[
 \Gr^{\nbigftilde}(\nbige^{\lor}_P)
\simeq
 \Gr^{\nbigftilde}(\nbige_P)^{\lor}
\]

\subsection{Lattice}
\label{subsection;10.5.16.20}

Let $E$ be an unramifiedly good lattice
of $(\nbige,\DD)$.
We set 
$\pi^{\ast}E:=
 \pi^{-1}(E)\otimes_{\pi^{-1}\nbigo_\nbigx}
 \nbigo_{\nbigxtilde(W)}$.
We have the induced filtration
$\nbigftilde^Q$ of 
$\pi^{\ast}E_Q\subset
 \pi^{\ast}\nbige_Q$.
\begin{prop}
\mbox{{}}\label{prop;10.5.12.210}
\begin{itemize}
\item
$\Gr^{\nbigftilde}(\pi^{\ast}E_Q)$
is free as an $\nbigo_{\nbigxtilde}$-module.
\item
We have the induced lattice
$\Gr^{\nbigftilde}(E_P)$
of $\Gr^{\nbigftilde}(\nbige_P)$.
It is functorial with respect to
morphism,
tensor product, direct sum
and dual.
\end{itemize}
\end{prop}

\subsection{Splitting}

We give some statements for
splitting of full Stokes filtrations.

\subsubsection{Flat splitting}

First, we consider
$\DD$-flat splittings.
Let $(\nbige,\DD)$
be an unramifiedly good meromorphic
$\varrho$-flat bundle on $(\nbigx,\nbigd)$.

\begin{prop}
\mbox{{}}\label{prop;10.5.12.211}
Assume $\varrho(Q)\neq 0$.
We can find a $\DD$-flat
splitting of $\nbigftilde^Q$,
i.e.,
we can find a $\DD$-flat decomposition
$\pi^{\ast}\nbige_Q=
 \bigoplus_{\gminia\in\Irr(\DD,Q)}
 \pi^{\ast}\nbige_{Q,\gminia}$
such that
$\nbigftilde^Q_{\gminia}(\pi^{\ast}\nbige_Q)
=\bigoplus_{\gminib\leq_Q^{\varrho}\gminia}
 \pi^{\ast}\nbige_{Q,\gminib}$.
\end{prop}

Let $E$ be an unramifiedly good lattice
of $(\nbige,\DD)$.
\begin{prop}
\label{prop;10.5.10.3}
Assume $\varrho(Q)\neq 0$
and that $E$ satisfies the non-resonance condition,
i.e., 
 for distinct eigenvalues $\alpha$ and $\beta$ of
 $\Res_{\nbigd_i}(\varrho^{-1}\DD)$
 we have $\alpha-\beta\not\in\seisuu$.
Then, we can find a $\DD$-flat splitting
$\pi^{\ast}E_{Q}=
 \bigoplus_{\gminia\in\Irr(\DD,Q)}
\pi^{\ast}E_{Q,\gminia}$  of
$(\pi^{\ast}E_{Q},\nbigftilde^Q)$.
\end{prop}
Note that we have 
the following natural isomorphisms
on a neighbourhood of $Q$
with $\varrho(Q)\neq 0$,
for the direct summands of 
the above decompositions:
\[
 \pi^{\ast}\nbige_{Q,\gminia}
\simeq
 \pi^{\ast}\Gr^{\nbigftilde}_{\gminia}(\nbige_{\pi(Q)}),
\quad
 \pi^{\ast}E_{Q,\gminia}
\simeq
 \pi^{\ast}\Gr^{\nbigftilde}_{\gminia}(E_{\pi(Q)})
\]
For example,
we can take a nice frame of $\pi^{\ast}E_Q$
by lifting a nice frame of
$\Gr^{\nbigftilde}(E_{\pi(Q)})$.

If $\nbigd$ is smooth, we do not have to impose
additional assumptions
such as ``non-resonance'' and
``$\varrho(Q)\neq 0$''.
\begin{prop}
\label{prop;10.5.12.212}
Assume that $\pi(Q)$ is a smooth point of $\nbigd$.
Then, we can find
a $\DD$-flat splitting of
$(\pi^{\ast}E_Q,\nbigftilde^Q)$
and hence 
$(\pi^{\ast}\nbige_Q,\nbigftilde^Q)$.
\end{prop}

\subsubsection{Partially flat splitting}

We consider rough splittings
in more general situations.
Because the claims are local,
we use the setting and the notation
in Subsection \ref{subsection;10.5.3.7}.
Let $(\nbige,\DD)$ be a good meromorphic 
$\varrho$-flat bundle on $(\nbigx,\nbigd)$
with the good set of irregular values
$\nbigi\subset M(\nbigx,\nbigd)/H(\nbigx)$,
which are lifted to $M(\nbigx,\nbigd)$.
Let $E$ be an unramifiedly good lattice.
For simplicity,
we assume that the coordinate system
$(z_1,\ldots,z_n)$ is admissible for $\nbigi$.
Let $p$ be determined by the following conditions:
\[
 \ord_p(\gminia)<0\,\,\,
 (\forall\gminia\in\nbigi,\,\forall i\leq p),
\quad
 \ord_{p+1}(\gminia)=0\,\,
 (\exists\gminia\in\nbigi)
\]
Let $\DD_{\leq p}$
denote 
the restriction of $\DD$
to the $(z_1,\ldots,z_p)$-direction.

\begin{prop}
\label{prop;10.5.12.220}
For any $Q\in \pi^{-1}(\nbigd_{\ellsitabar})$,
we can find 
a $\DD_{\leq p}$-flat splitting
of $(\pi^{\ast}E_Q,\nbigftilde^Q)$.
It induces a $\DD_{\leq p}$-flat
splitting of 
$\nbigftilde^{Q'}(\pi^{\ast}E_{Q'})$
for $Q'$ contained in a sufficiently small 
neighbourhood of $Q$.
\end{prop}

We give a more refined statement.
Let $\lefttop{i}\DD$ $(i\geq p+1)$
be the induced flat $\varrho$-connection
on $E_{|\nbigd_i}$ with respect to $z_i$.
(See Subsection \ref{subsection;10.5.28.4})
Assume that $E_{|\nbigd_i}$ $(i\geq p+1)$
are equipped with 
$\lefttop{i}\DD$-flat
filtrations $\lefttop{i}F$.

\begin{prop}
\label{prop;10.5.12.221}
For any $Q\in \pi^{-1}(\nbigd_{\ellsitabar})$,
we can find 
a $\DD_{\leq p}$-flat splitting
of $(\pi^{\ast}E,\nbigftilde^Q)$
on a neighbourhood $\nbigu_Q$
which is compatible with
$\lefttop{i}F$ and 
$\Res_i(\DD)$.
Namely 
we have a $\DD_{\leq p}$-flat
splitting $\pi^{\ast}E_{|\nbigu_Q}
=\bigoplus_{\gminia\in\nbigi}
 \pi^{\ast}E_{\gminia,\nbigu_Q}$
such that 
(i) $\Res_i(\DD)$
preserves
$\pi^{\ast}E_{\gminia,\nbigu_Q|\nbigu_Q\cap 
 \pi^{-1}(\nbigd_i)}$,
(ii) $\lefttop{i}F_a=
 \bigoplus_{\gminia\in\nbigi}
 \lefttop{i}F_a\cap 
 \pi^{\ast}E_{\gminia,\nbigu_Q}$.
\end{prop}

We give some complements.
\begin{prop}
\label{prop;10.5.12.240}
Assume $\varrho(Q)\neq 0$.
If $\nbigu_Q$ is 
a sufficiently small neighbourhood of $Q$,
any $\DD_{\leq p}$-flat splitting
of $(E,\nbigftilde^Q)$ on
$\nbigu_Q\setminus \pi^{-1}(\nbigd)$
is extended to 
a $\DD_{\leq p}$-flat splitting
on $\nbigu_Q$.
\end{prop}

The following proposition is
a refinement of Proposition \ref{prop;10.5.10.3}.

\begin{prop}
\label{prop;10.5.12.241}
 Assume $\varrho(Q)\neq 0$.
 If we have $\alpha-\beta\not\in\seisuu$
 for distinct eigenvalues $\alpha$ and $\beta$ of
 $\Res_{\nbigd_i}(\varrho^{-1}\DD)$
 $(i\geq p+1)$,
 then we can find
 a $\DD$-flat splitting of
 $\pi^{\ast}E_Q$ compatible with
 $\lefttop{i}F$ $(i\geq p+1)$.
\end{prop}

\section{Good meromorphic $\varrho$-flat bundle
in the level $\vecm$}
\label{section;10.5.11.1}
In this section,
we state some propositions
on Stokes filtrations in the level $\vecm$.
The proof will be given in Section
\ref{section;10.5.11.10}.
(See Subsection \ref{subsection;10.5.18.10}.)
We use the setting
in Subsection \ref{subsection;10.5.10.2}
unless otherwise specified.
Namely, let $Y$ be a complex manifold
with a simple normal crossing divisor $\nbigd_{Y}'$.
Let $\nbigk$ be a complex manifold
with a holomorphic function $\varrho$
such that $d\rho$ is nowhere vanishing
on $\nbigk^0:=\rho^{-1}(0)$.
Let $\nbigx:=\Delta_z^k\times Y\times \nbigk$,
$\nbigd_{z,i}:=\{z_i=0\}$
and $\nbigd_z:=\bigcup_{i=1}^{k}\nbigd_{z,i}$.
We also put 
$\nbigd_{Y}:=
 \Delta_z^k\times \nbigd_{Y}'\times\nbigk$
and $\nbigd:=\nbigd_z\cup \nbigd_{Y}$.
For any subset $I\subset\kbar$,
we put $\nbigd_{z,I}:=\bigcap_{i\in I}\nbigd_{z,i}$.
We set $\nbigx^0:=\nbigx\times_{\nbigk}\nbigk^0$
and $W:=\nbigx^0\cup \nbigd_z$.
Let $\pi:\nbigxtilde(W)\lrarr \nbigx$ 
denote the real blow up
of $\nbigx$ along $W$.

\subsection{Orders on weakly good sets of irregular 
values in the level $\vecm$}
\label{subsection;07.12.13.60}

Let $\vecm\in\seisuu_{< 0}^{k}$.
Let $\nbigi$ be 
a weakly good set of irregular values on
$(\nbigx,\nbigd)$ in the level $\vecm$.
We put 
$F_{\gminia,\gminib}:=
-\Re\bigl(\varrho^{-1}\,
 (\gminia-\gminib)
 \bigr)\,|\varrho\vecz^{-\vecm}|$
for any distinct $\gminia,\gminib\in \nbigi$.
They naturally induce $C^{\infty}$-functions
on $\nbigxtilde(W)$.

\begin{notation}
\label{notation;07.12.13.15}
Let $A$ be any subset of $\nbigxtilde(W)$.
For distinct $\gminia,\gminib\in\nbigi$,
we say 
$\gminia<_A^{\varrho}\gminib$
if $F_{\gminia,\gminib}(Q)<0$ for any 
$Q\in A$.
We say 
$\gminia\leq_A^{\varrho}\gminib$ for 
$(\gminia,\gminib)\in \nbigi^2$
if either $\gminia<_A^{\varrho}\gminib$
or $\gminia=\gminib$ holds.
The relation $\leq_A^{\varrho}$ gives
a partial order on $\nbigi$.
We use the symbol
$\leq_P^{\varrho}$ in the case $A=\{P\}$.
(We will use the symbol
$\leq_P$ instead of $\leq_P^{\varrho}$,
if there is no risk of confusion.)
\hfill\qed
\end{notation}
\index{order $\leq_A$, $\leq^{\varrho}_A$}
\index{order $\leq_P$, $\leq^{\varrho}_P$}

Let $\nbigi$ be a weakly good set of irregular values
in the level $\vecm$.
For any $f\in M(\nbigx,\nbigd_z)$,
$\nbigi':=\bigl\{
 \gminia-f\,\big|\,\gminia\in\nbigi
 \bigr\}$ is also a weakly good set of irregular values
in the level $\vecm$.
The natural bijection obviously preserves
the orders $\leq_A^{\varrho}$
for any subset $A$ of $\nbigxtilde(W)$.

For $Q\in\pi^{-1}(\nbigd)$,
let $\gbigu(Q,\nbigxtilde(W),\nbigi)$
or $\gbigu(Q,\nbigi)$ denote the set of
$\nbigu\in\gbigu(Q,\nbigxtilde(W))$
such that 
$\leq_Q^{\varrho}=\leq^{\varrho}_{\nbigu}$.
If $\nbigu\in\gbigu(Q,\nbigxtilde(W))$ is sufficiently small,
it is contained in $\gbigu(Q,\nbigxtilde(W),\nbigi)$.
For any point 
$Q'\in \pi^{-1}(\nbigd)\cap \nbigu$,
the map
$\bigl(\nbigi,\leq_Q^{\varrho}\bigr)
\lrarr \bigl(\nbigi,\leq_{Q'}^{\varrho}\bigr)$
preserves the orders.

\subsection{Stokes filtration in the level $\vecm$}
\label{subsection;10.5.10.10}

Let $(\nbige,\DD)$
be a meromorphic $\varrho$-flat bundle
on $(\nbigx,\nbigd)$
with a weakly good lattice 
$(E,\nbigi)$ in the level $(\vecm,i(0))$.
The restriction of $\DD$
to the $\Delta_{z}^k$-direction
is denoted by $\DD_z$.
Let $\pi:\nbigxtilde(W)\lrarr \nbigx$ be the real blow up.
We have the decomposition
induced by (\ref{eq;07.12.15.2}):
\[
 \pi^{\ast}E_{|\widehat{\pi^{-1}(W)}}
=\bigoplus_{\gminia\in\nbigi}
 E_{\widehat{\pi^{-1}(W)},\gminia}
\]

\begin{prop}
\label{prop;10.5.11.20}
For any $Q\in\pi^{-1}(\nbigd_{z,\kbar})$,
there uniquely exists a $\DD$-flat
filtration $\nbigf^Q$ of $\pi^{\ast}\nbige_Q$
indexed by
$\bigl(\nbigi,\leq_Q^{\varrho}\bigr)$
with the following property:
\begin{itemize}
\item
It comes from a filtration of $\pi^{\ast}E_Q$
such that $\Gr^{\nbigf^Q}(\pi^{\ast}E_Q)$
is a locally free $\nbigo_{\nbigxtilde(\nbigd)}$-module.
\item
Take $\nbigu\in\gbigu(Q,\nbigi)$ on which
$\nbigf^Q$ is given.
Then we have
\[
 \nbigf^Q_{\gminia}(\pi^{\ast}E)_{|
 \nbigu\cap \widehat{\pi^{-1}(W)}}
=\bigoplus_{\gminib\leq_Q^{\varrho}\gminia}
 E_{\widehat{\pi^{-1}(W)},\gminia|\nbigu}
\]
\end{itemize}
\end{prop}

The filtration $\nbigf^Q$ is called
a Stokes filtration in the level $(\vecm,i(0))$
(or simply in the level $\vecm$).

\begin{prop}
\label{prop;10.5.11.21}
Take $\nbigu\in \gbigu(Q,\nbigi)$ 
on which $\nbigf^Q$ is given.
For any $Q'\in\nbigu\cap\pi^{-1}(\nbigd_{z,\kbar})$,
the natural morphism
$\bigl(
 \pi^{\ast}\nbige_{Q'},\nbigf^{Q}
 \bigr)
\lrarr
 \bigl(\pi^{\ast}\nbige_{Q'},\nbigf^{Q'}\bigr)$
is compatible.
\end{prop}

\begin{rem}
\label{rem;10.5.18.11}
Let $\nbigu\in\gbigu(\nbigxtilde(W))$.
If $\pi^{\ast}\nbige_{|\nbigu}$
has a $\DD$-flat filtration
$\nbigf^{\nbigu}$ 
indexed by $\bigl(\nbigi,
 \leq_{\nbigu}^{\varrho}\bigr)$
with the following property,
which is also called a Stokes filtration
on $\nbigu$:
\begin{itemize}
\item
Let $Q$ be any point of
$\nbigu\cap\pi^{-1}(\nbigd_{z,\kbar})$.
We have the induced filtration
$\nbigf^{\nbigu}$ of
$\pi^{\ast}E_Q$.
Then, $\nbigf^{\nbigu}$
and $\nbigf^Q$ are compatible over
$\bigl(
 \nbigi,\leq_{\nbigu}^{\varrho}
 \bigr)\lrarr
 \bigl(
 \nbigi,\leq_Q^{\varrho}
 \bigr)$.
\end{itemize}
Such a filtration is uniquely determined if it exists,
by Lemma {\rm\ref{lem;10.5.12.22}}.
For a multi-sector $S$ of $\nbigx\setminus W$,
if the Stokes filtration for $\Sbar$ exists,
it is also called Stokes filtration for $S$.
\hfill\qed
\end{rem}

\subsection{Functoriality}
\label{subsection;10.5.11.22}

Let $(\nbige_i,\DD_i)$ $(i=1,2)$
be meromorphic
$\varrho$-flat bundles on $(\nbigx,\nbigd)$
with weakly good lattices $(E_i,\nbigi_i)$.

\begin{prop}
\label{prop;10.5.12.10}
Let $F:\nbige_1\lrarr\nbige_2$ be a 
$\DD$-flat morphism.
For simplicity,
we assume that
$\nbigi_1\cup\nbigi_2$ is weakly good
in the level $(\vecm,i(0))$.
Then, for each $Q\in \pi^{-1}(\nbigd_{z,\kbar})$,
the induced morphism
$\pi^{\ast}\nbige_{1\,Q}\lrarr
 \pi^{\ast}\nbige_{2\,Q}$ is compatible with
the Stokes filtrations in the level $\vecm$.
\end{prop}

\begin{prop}
If $\nbigi_1\otimes\nbigi_2$ is weakly good
in the level $(\vecm,i(0))$,
we have
\[
 \nbigf^{Q}_{\gminia}\bigl(
 \pi^{\ast}(\nbige_1\otimes\nbige_2)_Q
 \bigr)
=\sum_{\gminia_1+\gminia_2=\gminia}
 \nbigf^Q_{\gminia_1}
 (\pi^{\ast}\nbige_{1,Q})
\otimes
 \nbigf^Q_{\gminia_2}
 (\pi^{\ast}\nbige_{2,Q})
\]
If $\nbigi_1\oplus\nbigi_2$ is weakly good
in the level $(\vecm,i(0))$,
we have
\[
 \nbigf^{Q}_{\gminia}\bigl(
 \pi^{\ast}(\nbige_1\oplus\nbige_2)_Q
 \bigr)
=\nbigf^Q_{\gminia}(\pi^{\ast}\nbige_{1,Q})
\oplus
 \nbigf^Q_{\gminia}(\pi^{\ast}\nbige_{2,Q})
\]
\end{prop}

\begin{prop}
Let $(\nbige,\DD)$ be
a $\varrho$-meromorphic flat bundle
on $(\nbigx,\nbigd)$
with a weakly good lattice $(E,\nbigi)$
in the level $(\vecm,i(0))$.
The Stokes filtration of
$(\nbige^{\lor},\DD^{\lor})$
in the level $\vecm$
is given as follows:
\[
 \nbigf^Q_{\gminia}\bigl(
 \pi^{\ast}\nbige^{\lor}_{Q}
 \bigr)
=
 \Bigl(
 \sum_{\gminib\not\geq-\gminia}
 \nbigf^Q_{\gminib}(\pi^{\ast}\nbige_Q)
 \Bigr)^{\bot}
\]
\end{prop}

\subsection{Characterization by growth order}

Let $(\nbige,\DD)$ and $E$ be as above.
Let $\vecv$ be a frame of $E$.
Let $\nbigu_Q\in\gbigu(Q,\nbigi)$
be sufficiently small.
For any $\DD_{z}$-flat
section $f$ of 
$E_{|\nbigu_Q\setminus \pi^{-1}(W)}$,
we have the expression
$f=\sum f_j\,v_j$.
Put $\vecf:=(f_j)$.

\begin{prop}
\label{prop;10.5.11.23}
Assume $\varrho(Q)\neq 0$.
We have $f\in\nbigf^Q_{\gminia}$,
if and only if the following holds
for some $C>0$:
\[
 \bigl|
\vecf\exp(\varrho^{-1}\gminia)
\bigr|
=O\Bigl(
 \exp\bigl(C|\vecz^{\vecm(1)}|\bigr)
 \,|z_{i(0)}|^{-C}
 \Bigr)
\]
\end{prop}

\subsection{The associated graded bundle}

For any $P\in \nbigd_{z,\kbar}$,
we have the associated graded sheaf
with an induced $\DD$-flat connection
on a neighbourhood of $\pi^{-1}(P)$:
\[
 \Gr^{\nbigf}(\pi^{\ast}\nbige_{\pi^{-1}(P)},
\DD)
=\bigoplus_{\gminia\in\nbigi}
 \bigl(
 \Gr^{\nbigf}_{\gminia}
 (\pi^{\ast}\nbige_{\pi^{-1}(P)}),
 \DD_{\gminia}
\bigr)
\]
By taking the push-forward via $\pi$,
we obtain a graded $\nbigo_\nbigx(\ast \nbigd)$-module
with a flat $\DD$-connection
on a neighbourhood of $P$:
\[
 \Gr^{\nbigf}(\nbige_P,\DD)
=\bigoplus_{\gminia\in\nbigi}
 \bigl(
 \Gr^{\nbigf}_{\gminia}(\nbige_P),
 \DD_{\gminia}
 \bigr) 
\]
Similarly,
we have
$\Gr^{\nbigf}(\pi^{\ast}E_{\pi^{-1}(P)},
 \DD)$
and $\Gr^{\nbigf}(E_P,\DD)$
on neighbourhoods of
$\pi^{-1}(P)$ and $P$, respectively.

\begin{prop}
$\Gr^{\nbigf}(\nbige_P,\DD)$
is a graded meromorphic
$\varrho$-flat bundle 
on a neighbourhood of $P$ such that
$\Gr^{\nbigf}(\pi^{\ast}\nbige_{\pi^{-1}(P)},
\DD)
\simeq
\pi^{\ast}
\Gr^{\nbigf}(\nbige_P,\DD)$.
We also have
$\Gr^{\nbigf}(\pi^{\ast}E_{\pi^{-1}(P)},
\DD)
\simeq
\pi^{\ast}
\Gr^{\nbigf}(E_P,\DD)$.

For each $\gminia$,
we have a canonical isomorphism
$\Gr^{\nbigf}_{\gminia}
 \bigl(E,\DD
 \bigr)_{|\What}
\simeq
 (\Ehat_{\gminia\What},\DDhat^{\varrho}_{\gminia})$,
where the right hand side is the direct summand
in {\rm(\ref{eq;07.12.15.2})}.
In particular,
$\bigl(
 \Gr^{\nbigftilde}_{\gminia}(\nbige),
 \DD_{\gminia}
 \bigr)$ has a weakly good lattice
$\bigl(
 \Gr^{\nbigf}_{\gminia}(E),
 \{\gminia\}
 \bigr)$
in the level $(\vecm,i(0))$.
\end{prop}

\subsubsection{Functoriality}

If we are given a morphism of
meromorphic $\varrho$-flat bundles
$F:(\nbige_1,\DD_1)
 \lrarr(\nbige_2,\DD_2)$
with weakly good lattices
$(E_i,\nbigi_i)$ in the level $(\vecm,i(0))$,
we have the induced morphism
\[
 \Gr^{\nbigf}(F):
 \Gr^{\nbigf}(E_1,\DD_1)
\lrarr
 \Gr^{\nbigf}(E_2,\DD_2).
\]
In particular, we have
$\Gr^{\nbigf}(\nbige_1,\DD_1)
\lrarr
 \Gr^{\nbigf}(\nbige_2,\DD_2)$.

If $\nbigi_1\otimes\nbigi_2$ is weakly good
in the level $(\vecm,i(0))$,
we have the following canonical isomorphism
on a neighbourhood of $P$:
\[
 \Gr^{\nbigf}(E_1\otimes E_2)
\simeq
 \Gr^{\nbigf}(E_1)\otimes
 \Gr^{\nbigf}(E_2)
\]
If $\nbigi_1\oplus\nbigi_2$ is weakly good
in the level $(\vecm,i(0))$,
we have the following canonical isomorphism
on a neighbourhood of $P$:
\[
 \Gr^{\nbigf}(E_1\oplus E_2)
\simeq
 \Gr^{\nbigf}(E_1)\oplus
 \Gr^{\nbigf}(E_2).
\]
Let $(E,\DD,\nbigi)$
be a weakly good lattice in the level $(\vecm,i(0))$.
We have a canonical isomorphism
on a neighbourhood of $P$:
\[
 \Gr^{\nbigf}(E^{\lor})
\simeq
 \Gr^{\nbigf}(E)^{\lor}
\]

\subsection{Splitting}
\label{subsection;10.5.10.11}

Let $\nbigd_{Y}=
 \bigcup_{j\in\Lambda}\nbigd_{Y,j}$
be the irreducible decomposition.
For each $j\in\Lambda$,
we have the residue
$\Res_{Y,j}(\DD)$
on $E_{|\nbigd_{Y,j}}$.
Assume that we are given 
a filtration $\lefttop{j}F$
of $E_{|\nbigd_{Y,j}}$
which is flat with respect to
locally induced filtration
$\lefttop{j}\DD$.
\begin{prop}
\label{prop;10.5.11.30}
We can take a $\DD_z$-flat
splitting 
$\pi^{\ast}E_{|\nbigu_Q}=
 \bigoplus_{\gminia\in\nbigi}
 E_{\gminia,\nbigu_Q}$
which is compatible with
$\Res_{Y,j}(\DD)$
and $\lefttop{j}F$ $(j\in\Lambda)$.
\end{prop}

\begin{prop}
\label{prop;10.5.11.31}
Assume $\varrho(Q)\neq 0$.
If the non-resonance condition
is satisfied for each
$\Res_{Y,j}(\varrho^{-1}\DD)$,
we can take a $\DD$-flat splitting
compatible with
$\lefttop{j}F$ $(j\in \Lambda)$.
\end{prop}

\begin{prop}
\label{prop;10.5.11.32}
Assume $\varrho(Q)\neq 0$.
If $\nbigu_Q$ is 
a sufficiently small neighbourhood of $Q$,
any $\DD_{z}$-flat splitting
of $(\nbige,\nbigf^Q)$ on
$\nbigu_Q\setminus \pi^{-1}(\nbigd_z)$
is extended to
a $\DD_{\leq p}$-flat splitting
on $\nbigu_Q$.
\end{prop}

\begin{prop}
\label{prop;10.5.11.41}
If $\nbigd_z$ is smooth and if 
$\nbigd_{Y}=\emptyset$,
we can find a $\DD$-flat splitting
of $\nbigf^Q$.
\end{prop}

\section{Some notation}
\subsection{}
\index{sector}
We prepare some notation for the proof of
the claims in Sections \ref{section;10.5.10.5}
and \ref{section;10.5.11.1}.
Let $\Sec[\delta,\theta^{(0)},\theta^{(1)}]$ 
denote a closed sector
in $\Delta^{\ast}$:
\[
 \Sec[\delta,\theta^{(0)},\theta^{(1)}]:=
 \bigl\{
 z\in\Delta^{\ast}\,\big|\,
 0<|z|\leq \delta,\,\,
 \theta^{(0)}\leq\arg(z)\leq\theta^{(1)}
 \bigr\}
\]
\index{sector $\Sec[\delta,\theta^{(0)},\theta^{(1)}]$}
Let $\pi:\Deltatilde(0)\lrarr\Delta$ denote 
the real blow up along $0$.
Let 
$\Secbar\bigl[\delta,\theta^{(0)},\theta^{(1)}\bigr]$
denote the closure of 
$\Sec\bigl[\delta,\theta^{(0)},\theta^{(1)}\bigr]$
in the real blow up $\Deltatilde(0)$ along $0$:
\[
 \Secbar\bigl[\delta,\theta^{(0)},\theta^{(1)}\bigr]:=
 \bigl\{
 (t,\theta)\,\big|\,
 0\leq t\leq \delta,\,\,
 \theta^{(0)}\leq\theta\leq\theta^{(1)}
 \bigr\}
\]
\index{sector 
$\Secbar\bigl[\delta,\theta^{(0)},\theta^{(1)}\bigr]$}
We put 
 $\Sec\bigl[\theta^{(0)},\theta^{(1)}\bigr]:=
 \pi^{-1}(0)\cap 
 \Secbar\bigl[\delta,\theta^{(0)},\theta^{(1)}\bigr]$.
\index{set
 $\Sec\bigl[\theta^{(0)},\theta^{(1)}\bigr]$}

\vspace{.1in}

\index{multi-sector}
Let $X:=\Delta^{k}\times Y$
and $D:=\bigcup_{i=1}^{k}\{z_i=0\}$.
In this paper,
a multi-sector (or sector, for simplicity)
of $X\setminus D$
means a subset $S$ of the following form
\begin{equation}
\label{eq;08.9.2.2}
\prod_{j\in I}\Sec[\delta_j,\theta_j^{(0)},\theta_j^{(1)}]
\times U
\subset
 (\Delta^{\ast})^{I}
\times
\bigl(
 (\Delta^{\ast})^{I^c} \times Y
\bigr),
\end{equation}
where $I\subset\kbar$,
$I^c:=\kbar\setminus I$,
$\Sec[\delta_j,\theta_j^{(0)},\theta_j^{(1)}]
\subset\Deltatilde_{z_j}(0)$,
and $U$ denotes a compact region in 
$(\Delta^{\ast})^{I^c}\times Y$.
(We admit the case $I=\emptyset$.)
Let $\pi:\Xtilde(D)\lrarr X$ denote the real blow up
of $X$ along $D$.
The closure of $S$ in $\Xtilde(D)$ is denoted by $\Sbar$.
\index{closure $\Sbar$}

\begin{notation}
Let $\Multisector(X\setminus D)$ denote the set of
multi-sectors in $X\setminus D$.
For a point $Q\in \pi^{-1}(D)$,
let $\Multisector(Q,X\setminus D)$ denote the set of
multi-sectors $S$ in $X\setminus D$ such that
$Q$ is contained in the interior part of $\Sbar$.
\hfill\qed
\end{notation}
\index{set $\Multisector(X\setminus D)$}
\index{set $\Multisector(Q,X\setminus D)$}

\begin{df}
\label{df;08.12.23.1}
We say that we shrink
$S=
 \prod_{j\in I}
 \Sec[\delta_j,\theta^{(0)}_j,\theta^{(1)}_j]
 \times U$
in the radius direction,
when we replace it
with
$S'=\prod_{j\in I}
 \Sec[\delta_j',\theta^{(0)}_j,\theta^{(1)}_j]
 \times U$ 
for some $\delta_j'\leq\delta_j$.
\hfill\qed
\end{df}
\index{shrink in the radius direction}

\subsection{}
We use the setting in Subsection
\ref{section;10.5.11.1}.
If $\nbigk^0\neq \emptyset$,
we implicitly assume that
$\nbigk$ is a product
$\Delta_{\varrho}\times\nbigk'$
to consider a sector of
$\nbigk^{\ast}:=\nbigk\setminus\nbigk^0$.
Because we are interested in a local theory,
it is not essential.

Let $\nbigi$ be a weakly good set of irregular values
in the level $(\vecm,i(0))$.
Let $\Sep(\gminia,\gminib)$ denote 
the subset $(F_{\gminia,\gminib})^{-1}(0)
\subset \nbigxtilde(W)$,
and let $\Sep\bigl(\nbigi\bigr)$ denote
the union of $\Sep(\gminia,\gminib)$
for distinct pairs $(\gminia,\gminib)$ of $\nbigi$.

\begin{notation}
Let $\Multisector(\nbigx\setminus W,\nbigi)$ denote
the set of the multi-sectors $S$ of $\nbigx\setminus W$
which is the product of
\[
 \Sec[\delta_{i(0)},\theta_{i(0)}^{(0)},
 \theta_{i(0)}^{(1)}]
\subset\Delta_{z_{i(0)}}^{\ast},
\quad
U_1\subset \prod_{p\neq i(0)}\Delta_{z_p}^{\ast},
 \quad
U_2\subset Y\times\nbigk^{\ast},
\]
where $U_1$ (resp. $U_2$)
is a compact region or a multi-sector
of $\prod_{p\neq i(0)}\Delta_{z_p}^{\ast}$
(resp. $Y\times\nbigk^{\ast}$).
We assume $|m_{i(0)}|\cdot
 \bigl|\theta^{(0)}_{i(0)}-\theta^{(1)}_{i(0)}
 \bigr|<\pi$,
and moreover, 
 there exist 
 $\theta_{i(0)}^{\prime(0)}
<\theta_{i(0)}^{\prime(1)}$ in 
the open interval
 $\openopen{\theta_{i(0)}^{(0)}}{\theta_{i(0)}^{(1)}}$
 such that
$\Sbar\cap\Sep(\nbigi)\subset
 \Sec[\delta_{i(0)},
 \theta_{i(0)}^{\prime(0)},
 \theta_{i(0)}^{\prime(1)}]
\times U_1
\times U_2$.
\hfill\qed
\end{notation}
\index{set
 $\Multisector(\nbigx\setminus W,\nbigi)$}

The following lemma is clear.
\begin{lem}
\label{lem;10.5.10.50}
We put $\nbigd_{z,\kbar}^0:=
 \nbigd_{z,\kbar}\times_{\nbigk}\nbigk^0$.
\begin{itemize}
\item
Let $Z_0\subset
 \pi^{-1}(\nbigd_{z,\kbar}\setminus \nbigd_{z,\kbar}^0)$
be a product of
$\prod_{p=1}^{k}
 \Sec[\theta^{(0)}_p,\theta^{(1)}_p]$
and a compact region $U$ of 
$Y\times\nbigk^{\ast}$
 such that
\[
 Z_0\cap\Sep(\nbigi)
\subset
\prod_{p\neq i(0)}
 \Sec[\theta^{(0)}_p,\theta^{(1)}_p]\times
  \Sec[\theta^{\prime(0)}_{i(0)},
 \theta^{\prime(1)}_{i(0)}]
 \times U
\]
for some 
 $\theta^{\prime(0)}_{i(0)}
 <\theta^{\prime(1)}_{i(0)}$
 in $\openopen{\theta^{(0)}_{i(0)}}{\theta^{(1)}_{i(0)}}$,
 where 
 $\Sec[\theta^{(0)}_p,\theta^{(1)}_p]\subset
 \Deltatilde_{z_p}(0)$.
Take sufficiently small $\delta_p>0$ $(p=1,\ldots,k)$,
and put $S:=\prod_{p=1}^{k}
 \Sec[\delta_p,\theta^{(0)}_p,\theta^{(1)}_p]
 \times U$.
Then,
$S\in \Multisector(\nbigx\setminus W,\nbigi)$.
\item
Assume $\nbigk=\Delta_{\varrho}\times\nbigk'$.
Let $Z_0\subset\pi^{-1}(\nbigd_{z,\kbar}^0)$
be the product of
$\Sec[\theta^{(0)}_{\varrho},
 \theta^{(1)}_{\varrho}]
\times\prod_{p=1}^{k}
 \Sec[\theta^{(0)}_p,\theta^{(1)}_p]$
and a compact region $U$ of $Y\times\nbigk'$
such that
\[
Z_0\cap\Sep(\nbigi)
\subset
\Sec[\theta^{(0)}_{\varrho},
 \theta^{(1)}_{\varrho}]
 \times
 \prod_{p\neq i(0)}
 \Sec[\theta^{(0)}_p,\theta^{(1)}_p]\times
  \Sec[\theta^{\prime(0)}_{i(0)},
 \theta^{\prime(1)}_{i(0)}]
 \times U
\]
for some 
 $\theta^{\prime(0)}_{i(0)}
 <\theta^{\prime(1)}_{i(0)}$
 in $\openopen{\theta^{(0)}_{i(0)}}{\theta^{(1)}_{i(0)}}$,
 where 
 $\Sec[\theta^{(0)}_{\varrho},
 \theta^{(1)}_{\varrho}]\subset
 \Deltatilde_{\varrho}(0)$.
Take sufficiently small
$\delta_{\lambda}>0$
and $\delta_p>0$ $(p=1,\ldots,k)$.
We put $S:=
\Sec[\delta_{\lambda},
 \theta^{(0)}_{\lambda},\theta^{(1)}_{\lambda}]
 \times
 \prod_{p=1}^{\ell}
 \Sec[\delta_p,\theta^{(0)}_p,\theta^{(1)}_p]
\times U$.
Then, $S\in\Multisector(\nbigx\setminus W,\nbigi)$.
\hfill\qed
\end{itemize}
\end{lem}

\begin{notation}
\label{notation;07.12.9.10}
For $Q\in \pi^{-1}(\nbigd_{z,\kbar})$,
let $\Multisector(Q,\nbigx\setminus W,\nbigi)$ denote the set of
multi-sectors 
$S\in\Multisector(\nbigx\setminus W,\nbigi)\cap
 \Multisector(Q,\nbigx\setminus W)$
such that
 $\Sbar\in \gbigu(Q,\nbigxtilde(W),\nbigi)$.
\hfill\qed 
\end{notation}
\index{set $\Multisector(Q,\nbigx\setminus W,\nbigi)$}

We obtain the following lemma
from Lemma \ref{lem;10.5.10.50},
which will be used implicitly.
\begin{lem}
Let $Q$ be any point of $\pi^{-1}(\nbigd_{z,\kbar})$.
Let $\nbigu$ be any neighbourhood of $Q$
in $\nbigxtilde(W)$.
Then, there exists 
$S\in\Multisector(Q,\nbigx\setminus W,\nbigi)$
such that $\Sbar\subset \nbigu$.
\hfill\qed
\end{lem}

The following lemma is clear by
the condition,
which will be also used implicitly.
\begin{lem}
For any $S\in \Multisector(\nbigx\setminus W,\nbigi)$,
there exists $S'\subset S$ such that
(i) $S'\in\Multisector(\nbigx\setminus W,\nbigi)$,
(ii) $\Sbar'\cap \pi^{-1}(\nbigd_{z,j})=\emptyset$
 $(j\neq i(0))$
 and $\Sbar'\cap\pi^{-1}(\nbigx^0)=\emptyset$,
(iii) $\leq_{S}^{\varrho}=\leq_{S'}^{\varrho}$.
\hfill\qed
\end{lem}

\section{Preliminary in the smooth divisor case}
\label{subsection;08.9.28.30}
In the following sections,
$\pi^{\ast}E$ and $\pi^{\ast}\nbige$
are denoted by $E$ and $\nbige$,
if there is no risk of confusion.
We use the setting in Section
\ref{section;10.5.11.1}
with $k=1$ and $\nbigd_{Y}=\emptyset$,
i.e.,
$\nbigx=\Delta_{z_1}\times Y\times\nbigk$
and $\nbigd=\{0\}\times Y\times\nbigk$.
We also assume that $\varrho$
is nowhere vanishing on $\nbigk$.
Let $\pi:\nbigxtilde(\nbigd)\lrarr\nbigx$
denote the real blow up along $\nbigd$.
For a multi-sector $S$ in $\nbigx-\nbigd$,
let $\Sbar$ denote the closure of $S$
in $\nbigxtilde(\nbigd)$,
and the intersection $\Sbar\cap \pi^{-1}(\nbigd)$
is denoted by $Z$.
In the following,
$\leq_{S}^{\varrho}$
and $<_S^{\varrho}$
are denoted
by $\leq_S$ and $<_S$,
respectively.

\subsection{Existence of a decomposition}
\label{subsection;08.9.2.25}

Let $m\in\seisuu_{<0}$.
Let $\nbigi$ be a good set of irregular values
on $(\nbigx,\nbigd)$ in the level $m$.
Let $(E,\DD,\nbigi)$ be 
a good lattice of
a meromorphic $\varrho$-flat bundle
in the level $m$.
We have the irregular decomposition
in the level $m$:
\begin{equation}
\label{eq;08.9.4.31}
 (E,\DD)_{|\nbigdhat}
=\bigoplus_{\gminia\in\nbigi}
 (\Ehat_{\gminia},\DDhat_{\gminia})
\end{equation}

\begin{lem}
\label{lem;07.7.19.40}
For any 
$S\in \Multisector\bigl(\nbigx-\nbigd,\nbigi\bigr)$,
we have a decomposition of $E_{|\Sbar}$
\begin{equation}
 \label{eq;07.9.29.3}
  E_{|\overline{S}}
=\bigoplus_{\gminia\in \nbigi}
 E_{\gminia,S}
\end{equation}
such that (i) flat with respect to $\del_{z_1}$,
(ii) the restriction to $\Zhat$ is equal to
the pull back of the irregular decomposition
{\rm(\ref{eq;08.9.4.31})}
on $\nbigdhat$.
\end{lem}
\pf
We closely follow the standard argument
(Chapter 12 of \cite{wasow}, for example).
We give only an outline.
Let $\vecwhat=(\what_{i})$
be a frame of $E_{|\nbigdhat}$
which is compatible with the decomposition
(\ref{eq;08.9.4.31}).
When $\what_i\in \Ehat_{\gminia}$,
we put $\gminia(i):=\gminia$.
Let $M_r(\cnum)$ denote the space of
$r$-th square matrices,
where $r=\rank E$.
Let $M_{\gminia,\gminib}$ denote
the set of $A\in M_r(\cnum)$
such that  $A_{i,j}=0$
unless $\gminia(i)=\gminia$
and $\gminia(j)=\gminib$.
Thus, we obtain the decomposition:
\begin{equation}
 \label{eq;07.11.12.3}
 M_r(\cnum)=\bigoplus_{\gminia,\gminib}
 M_{\gminia,\gminib}
\end{equation}
We have the natural isomorphism
$M_{\gminia,\gminia}\simeq
 M_{\rank \Ehat_{\gminia}}(\cnum)$.
The element corresponding to the identity matrix
is denoted by $I_{\gminia}$.
Any element $A\in M_r(\cnum)$ has the decomposition
$A=\sum A_{\gminia,\gminib}$ corresponding to
(\ref{eq;07.11.12.3}).
Let $M_r(\cnum)^d$ (resp. $M_r(\cnum)^r$)
denote the subspace of $M_r(\cnum)$
which consists of $A$ such that
$A_{\gminia,\gminib}=0$ unless (resp. if)
$\gminia=\gminib$.
Any element $A\in M_r(\cnum)$
is decomposed into $A=A^d+A^r$,
where $A^d\in M_r(\cnum)^d$
and $A^r\in M_r(\cnum)^r$.
We have the corresponding decomposition
for $M_r(\cnum)$-valued functions.

By Proposition \ref{prop;07.11.12.2},
we can take a holomorphic frame $\vecw$
of $\pi^{\ast}E_{|\Sbar}$ such that
$\vecw_{|\Zhat}= 
  \pi^{-1}\bigl(\vecwhat\bigr)_{|\Zhat}$.
Let $A$ be determined by
$\DD^{f}(z_1\del_{z_1})\vecw=\vecw\cdot A$.
By construction,
$A^d$ is of the form 
$\sum (m\varrho^{-1}\gminia)\, I_{\gminia}
 +A_0^d$,
where $|A_0^d|=O\bigl(|z_1|^{m+1}\bigr)$.
We also have $A^r_{|\Zhat}=0$.

We consider the change of the frame of the form
$\vecw'=\vecw\,\bigl(I+B\bigr)$
such that 
$\DD^{f}(z_1\del_{z_1})\vecw'=
 \vecw'\bigl(A^d+R\bigr)$,
where $B=B^r$ and $R=R^d$
and $B_{|\Zhat}=R_{|\Zhat}=0$.
The condition is as follows:
\begin{equation}
 \label{eq;07.9.29.2}
 R=(A^r\, B)^d,
\quad
 A^r+A^d\, B+(A^r\, B)^r
+z_1\del_{z_1}B
=B\, (A^d+R)
\end{equation}
By eliminating $R$,
we obtain the following equation
for $M_{r}(\cnum)^r$-valued function $B$:
\begin{equation}
\label{eq;07.9.29.1}
 z_1\del_{z_1}B
=B\, A^d-A^dB
-(A^rB)^r-A^r+B\,(A^r B)^d
\end{equation}
Note $A^r=O(|z_1|^N)$ for any $N$.
By changing the variable $x=z_1^{-1}$,
we can apply Proposition \ref{prop;07.6.15.3}
to (\ref{eq;07.9.29.1}).
Recall $S$ is of the form
$\Sec[\delta,\theta^{(0)},\theta^{(1)}]\times U$,
where $U$ denotes a compact region in 
$Y\times\nbigk$.
We can find solutions $B$ and $R$
of (\ref{eq;07.9.29.2})
such that $B=O(|z_1|^N)$ and $R=O(|z_1|^N)$
for any $N$
on $\Sec[\delta',\theta^{(0)},\theta^{(1)}]\times U$.
Since $B$ and $R$ are holomorphic,
we also obtain $B_{|\Zhat}=0$
and $R_{|\Zhat}=0$.
Hence, we obtain a decomposition like (\ref{eq;07.9.29.3}).
We can extend it to a decomposition 
on $\Sec[\delta,\theta^{(0)},\theta^{(1)}]\times U$
by using $\DD^f$.
Thus, we are done.
\hfill\qed

\subsection{Stokes filtration in the level $m$}

\index{Stokes filtration in the level $m$}

Let $S\in\Multisector(\nbigx-\nbigd,\nbigi)$
and a decomposition $E_{|\Sbar}
=\bigoplus_{\gminia\in \nbigi} E_{\gminia,S}$
as in Lemma \ref{lem;07.7.19.40}.
We define
\[
 \nbigf^S_{\gminia}:=
 \bigoplus_{\gminib\leq_S\gminia}
 E_{\gminib,S},
\quad
\nbigf^S_{<\gminia}:=
 \sum_{\gminib<_S\gminia}
 \nbigf^S_{\gminib}
=
 \bigoplus_{\gminib<_S\gminia}
 E_{\gminib,S}.
\]

\begin{lem}
\label{lem;07.6.14.5}
$\nbigf^S_{\gminia}$ is independent of the choice of
a splitting {\rm (\ref{eq;07.9.29.3})},
and it is $\DD$-flat.
\end{lem}
\pf
Let $E_{|\Sbar}=
 \bigoplus _{\gminia\in \nbigi} 
 E'_{\gminia,S}$
be another splitting.
The inclusion
$i_{\gminia}:E_{\gminia,S}\subset E_{|\Sbar}$
and the projection
$p_{\gminib}':E_{|\Sbar}\lrarr E'_{\gminib,S}$
give $\del_{z_1}$-flat morphisms
$E_{\gminia,S}\lrarr E'_{\gminib,S}$
for any $\gminia,\gminib$.
Let $f_{\gminia,\gminib}$ denote the composite.
We have $f_{\gminia,\gminib|\Zhat}=0$
for $\gminia\neq\gminib$
by construction,
and $f_{\gminia,\gminia|\Zhat}=\id$.
In particular, $f_{\gminia,\gminib}$ are bounded.
Hence, we obtain $f_{\gminia,\gminib}=0$
unless $\gminib\leq_S\gminia$
due to Corollary \ref{cor;07.12.27.3}.

Let $V$ be a holomorphic vector field
in the $Y$-direction on $U$.
We have the frames
$\vecv_{\gminia}$ of $E_{\gminia,S}$ on $S$.
Let $p_{\gminia}$ denote the projection
$E_{|\Sbar}\lrarr E_{\gminia,S}$.
Let $A_{\gminia,\gminib}(V)$ $(\gminia\neq\gminib)$
be determined by
$\bigl(p_{\gminib}\circ\DD^{f}(V)
 \circ i_{\gminia}\bigr)
 \vecv_{\gminia}
=\vecv_{\gminib}\cdot A_{\gminia,\gminib}(V)$.
Then, $A_{\gminia,\gminib}(V)_{|\Zhat}=0$
and $\del_{z_1}$-flat.
Thus, we obtain $A_{\gminia,\gminib}(V)=0$
unless $\gminib\leq_S\gminia$
due to Corollary \ref{cor;07.12.27.3}.
\hfill\qed

\begin{lem}
\label{lem;07.9.29.20}
The filtration $\nbigf^S$ 
indexed by $\bigl(\nbigi,\leq_S\bigr)$
is characterized by the following property:
\begin{itemize}
\item
 $\nbigf^{S}$ is flat with respect to $\del_{z_1}$.
\item
 $\nbigf^{S}_{\gminia|\Zhat}
 =\bigoplus_{\gminib\leq_S\gminia}
 \Ehat_{\gminib}$.
\end{itemize}
We call $\nbigf^S$ Stokes filtration
in the level $m$.
\end{lem}
\pf
Let $\Gr_{\gminia}:=
 \nbigf^S_{\gminia}/\nbigf^S_{<\gminia}$,
and let $r(\gminia):=\rank \Gr_{\gminia}$.
For a frame $\vecvbar_{\gminia}$ 
of $\Gr_{\gminia}$ on $\Sbar$,
let $A_{\gminia}$ be determined by
$\DD_{\gminia}^f(\del_{z_1})\vecvbar_{\gminia}
=\vecvbar_{\gminia}\cdot A_{\gminia}$,
where $\DD_{\gminia}^f$ denotes
the induced family of the flat connections 
of $\Gr_{\gminia}$.
Then, $A_{\gminia}$ is of the form
$\varrho^{-1}\del_{z_1}\gminia
 +A_{\gminia}^{\circ}$
where $A_{\gminia}^{\circ}$ is
the holomorphic section
of $\varrho^{-1}\, M_{r(\gminia)}(\cnum)\otimes
 \nbigo_{\nbigx}\bigl(m\nbigd\bigr)$.

Let $\nbigf^{\prime\,S}$ be a filtration
of $E_{|\Sbar}$ on $\Sbar$,
which has the above property.
We set $\Gr'_{\gminia}:=
 \nbigf^{\prime\,S}_{\gminia}
 /\nbigf^{\prime\,S}_{<\gminia}$,
which is equipped with 
the family of connections
$\DD_{\gminia}^{\prime\,f}$
along the $\Delta_{z_1}$-direction.
For any frame $\vecvbar_{\gminia}'$
of $\Gr'_{\gminia}$,
let $A_{\gminia}'$ be determined by
$\DD_{\gminia}^{\prime\,f}(\del_{z_1})
 \vecvbar_{\gminia}'
=\vecvbar_{\gminia}'\cdot A_{\gminia}'$,
and then $A_{\gminia}'$ is of a form
similar to $A_{\gminia}$.

We use an induction on the order $\leq_S$.
We put
$\nbigf^{S}_{\nbigb}:=
 \sum_{\gminia\in\nbigb}\nbigf^S_{\gminia}$
for any subset $\nbigb\subset \nbigi$.
Assume that 
we have already known
$\nbigf^{S}_{\gminib}=\nbigf^{\prime\,S}_{\gminib}$
for any $\gminib<_S\gminia$,
and we will show 
$\nbigf^S_{\gminia}=\nbigf^{\prime\,S}_{\gminia}$.
Let $\nbigb$ be the subset of $\nbigi$
such that
$\nbigf^{\prime\,S}_{\gminia}
\subset
 \nbigf^{S}_{\nbigb}$
and $\nbigf^{\prime\,S}_{\gminia}
 \not\subset \nbigf^S_{\nbigb'}$
for any $\nbigb'\subsetneq\nbigb$.
Let $\gminic$ be any maximal element of $\nbigb$.
Then, we obtain the flat morphism
$\phi_{\gminia,\gminic}:
 \Gr'_{\gminia}
\lrarr
 \Gr_{\gminic}$.
Due to Corollary \ref{cor;07.12.27.3},
we obtain $\phi_{\gminia,\gminic}=0$
unless $\gminia\geq_S\gminic$.
Therefore, we obtain 
$\nbigf^{\prime\,S}_{\gminia}
\subset \nbigf^S_{\nbigb}
\subset\nbigf^S_{\gminia}$.
By comparison of the rank,
we obtain 
$\nbigf^S_{\gminia}=\nbigf^{\prime\,S}_{\gminia}$.
Thus, we are done.
\hfill\qed

\section{Proof of the statements in
Section \ref{section;10.5.11.1}}
\label{section;10.5.11.10}
We use the setting in Subsection
\ref{section;10.5.11.1}.
Since we are interested in a local theory,
we put $Y:=\Delta^n_{\zeta}$,
 $\nbigd_{Y,j}:=\{\zeta_j=0\}$
and $\nbigd_Y:=\bigcup_{j=1}^{\ell}\nbigd_{Y,j}$
for some $1\leq\ell\leq n$,
although it does not matter
until Subsection \ref{subsection;07.12.7.1}.
We put $\gbigz:=\pi^{-1}(\nbigd_{z,\kbar})$.
In this section,
$P$ denotes a point in the real blow up.
For any multi-sector $S$ in $\nbigx\setminus W$,
let $\Sbar$ denote the closure of $S$
in $\nbigxtilde(W)$,
and let $Z$ denote 
$\Sbar\cap \pi^{-1}(W)$.
The orders
$\leq_{S}^{\varrho}$
and $<_S^{\varrho}$
are denoted
by $\leq_S$ and $<_S$,
respectively.

\subsection{Existence of Stokes filtration}
\label{subsection;10.5.26.12}

Let $\vecm\in\seisuu_{<0}^k$.
Let $(E,\DD,\nbigi)$ be a good lattice
of a meromorphic $\varrho$-flat bundle
in the level $(\vecm,i(0))$ on $(\nbigx,\nbigd)$.
The irregular decomposition (\ref{eq;07.12.15.1})
induces the following on $\Zhat$:
\begin{equation}
\label{eq;07.6.4.3}
 (E,\DD)_{|\Zhat}=
 \bigoplus_{\gminia\in\nbigi}
 (\Ehat_{\gminia},\DDhat_{\gminia})_{|\Zhat}
\end{equation}
We put
$\nbigf^Z_{\gminia}:=
 \bigoplus_{\gminib\leq_S\gminia}
 \Ehat_{\gminib|\Zhat}$,
and thus we obtain a filtration
$\nbigf^Z$ indexed by
$\bigl(\nbigi,\leq_S\bigr)$.
Let $\DD_z$ denote the restriction of $\DD$
to the $\Delta_z^k$-direction.

\begin{prop}
\label{prop;10.5.28.10}
Take $S\in
 \Multisector(\nbigx\setminus W,\nbigi)$
such that $\Sbar\cap\gbigz\neq \emptyset$.
If we sufficiently shrink $S$ in the radius direction,
the following holds:
\begin{itemize}
\item
There uniquely exists 
a $\DD$-flat filtration $\nbigf^S$
of $E_{|\Sbar}$ indexed by $(\nbigi,\leq_S)$
such that
$\nbigf^S_{|\Zhat}=\nbigf^Z$.
Moreover,
if a $\DD_{z}$-flat filtration $\nbigf^{\prime\,S}$
of $E_{|\Sbar}$ indexed by $(\nbigi,\leq_S)$
satisfies $\nbigf^{\prime\,S}_{|\Zhat}=\nbigf^Z$,
then $\nbigf^{\prime\,S}=\nbigf^S$.
\item
There exists a $\DD_{z}$-flat splitting
of the filtration $\nbigf^S$ on $\Sbar$.
Note that the restriction of such a splitting 
to $\Zhat$ is
equal to {\rm(\ref{eq;07.6.4.3})}.
\item
If $\varrho$ is nowhere vanishing,
any $\DD_{z}$-flat splitting of $\nbigf^S$ on $S$
is extended to a splitting on $\Sbar$.
\end{itemize}
We call $\nbigf^S$ Stokes filtration of $(E,\DD)$
on $S$ in the level $\vecm$.
\index{Stokes filtration in the level $\vecm$}
\end{prop}
\pf
We may assume $i(0)=1$.  
We show the following lemma analogue to 
Lemma \ref{lem;07.9.29.20}.
\begin{lem}
\label{lem;07.9.29.15}
If we sufficiently shrink $S$ in the radius direction,
there exists a decomposition on $\Sbar$
\begin{equation}
\label{eq;08.12.1.1}
  E_{|\overline{S}}
=\bigoplus_{\gminia\in\nbigi}
 \lefttop{1}E_{\gminia,S}
\end{equation}
such that (i) 
it is flat with respect to $\del_{z_1}$,
(ii) its restriction to $\Zhat$ is the same as
{\rm(\ref{eq;07.6.4.3})},
where $Z:=\Sbar\cap\pi^{-1}(W)$.
\end{lem}
\pf
We closely follow the standard argument
as in the proof of Lemma \ref{lem;07.7.19.40}.
Hence, we give only an outline.
We will shrink $S$ in the radius direction
without mention.
We take a frame 
$\vecwhat=(\vecwhat_{\gminia})$
of $E_{|\What}$
compatible with the irregular decomposition 
(\ref{eq;07.12.15.1}) in the level $\vecm$.
We use the decomposition of matrices
as in the proof of Lemma \ref{lem;07.7.19.40}.

By Proposition \ref{prop;07.11.12.2},
we can take a holomorphic frame $\vecw$
of $\pi^{\ast}E_{|\Sbar}$ such that
$\vecw_{|\Zhat}=
  \pi^{-1}\bigl(\vecwhat\bigr)$.
Let $A$ be determined by
$\DD^f(z_1\del_1)\vecw=\vecw\, A$.
By construction,
$A^d$ is of the form
$\sum_{\gminia}
(\varrho^{-1}z_1\del_1\gminia)\,
 I_{\gminia} +A_0^d$,
where $|A_0^d|=O\bigl(|\vecz^{\vecm(1)}|\bigr)$.
We also have $A^r_{|\Zhat}=0$.

We consider a change of frames of the form
$\vecw'=\vecw\,\bigl(I+B\bigr)$
such that 
$\DD^f(z_1\del_1)\vecw'=
 \vecw'\bigl(A^d+R \bigr)$,
where $B=B^r$, $R=R^d$,
and $B_{|\Zhat}=R_{|\Zhat}=0$.
Then, we obtain the equation
(\ref{eq;07.9.29.2}) for $B$ and $R$,
and (\ref{eq;07.9.29.1})
for $B$ by eliminating $R$.

We use the change of variables
$x=z_1^{-1}$
and $y_i=z_{i+1}^{-1}$ $(i=1,\ldots,k-1)$,
and $y_{k}=\varrho^{-1}$
if $\nbigk^0\neq \emptyset$.
By applying Proposition \ref{prop;07.6.15.3}
to (\ref{eq;07.9.29.1}),
we can find solutions $B$ and $R$
of (\ref{eq;07.9.29.2})
such that 
$B=O\Bigl(\prod_{i=1}^k|z_i|^N
 |\varrho|^N \Bigr)$
and $R=O\Bigl(\prod_{i=1}^k|z_i|^N
 |\varrho|^N \Bigr)$
for any $N$.
Since $B$ and $R$ are holomorphic,
we also obtain $B_{|\Zhat}=0$
and $R_{|\Zhat}=0$.
Thus, the proof of Lemma \ref{lem;07.9.29.15}
is finished.
\hfill\qed

\vspace{.1in}

Let $S$ be as in Lemma \ref{lem;07.9.29.15}.
We put
$\nbigf^{S}_{\gminia}:=
 \bigoplus_{\gminib\leq_S\gminia}
 \lefttop{1}E_{\gminia,S}$.
\begin{lem}
\label{lem;10.5.11.2}
They are independent of the choice of
 the decomposition {\rm(\ref{eq;08.12.1.1})},
and they are $\DD_z$-flat.
\end{lem}
\pf
We can take a multi-sector $S'\subset S$
such that
(i) $S'\in \Multisector(\nbigx\setminus W,\nbigi)$,
(ii) $\Sbar'\cap \pi^{-1}(\nbigd_{z,j})=\emptyset$
 for $j\neq i(0)$
and $\Sbar'\cap\pi^{-1}(\nbigx^0)=\emptyset$,
(iii) $\leq_{S'}=\leq_{S}$.
Then, we can show the claim of the lemma
by considering the restriction to $\Sbar'$
and by using Lemma \ref{lem;07.6.14.5}.
\hfill\qed

\vspace{.1in}

Let us consider the second condition.
Let $\Gr_{\gminia}$ denote the induced flat bundle
on $\Sbar$.
The frame $\vecw'$ 
in the proof of Lemma \ref{lem;07.9.29.15}
induces a frame $\vecvbar_{\gminia}$ 
of $\Gr_{\gminia}$.
Let $B_{\gminia}$ be determined by
$\DD_{z,\gminia}\vecvbar_{\gminia}
=\vecvbar_{\gminia}
 \bigl(d_z\gminia+B_{\gminia}\bigr)$.
By shrinking $S$ in the radius direction,
we may assume that 
$|B_{\gminib}|$  and $|B_{\gminic}|$
are sufficiently smaller than
$\bigl|
 \Re\bigl(\varrho^{-1}(\gminib-\gminic)\bigr)
 \bigr|$
 on $S$
for any $\gminib>_S\gminic$.
Let us consider the following claim:
\begin{description}
\item[$A(\gminia)$]
 There exists a $\DD_{z}$-flat splitting
 $\nbigf^S_{\gminia}=
 \bigoplus_{\gminib\leq_S\gminia}
 E_{\gminib,S}$ on $\Sbar$.
\end{description}
We have a similar claim $A(<_S\gminia)$
on the existence of a $\DD_z$-flat splitting 
of $\nbigf^S_{<\gminia}$.

We show $A(\gminia)$ by the induction on the order $\leq_S$.
If $\gminia$ is minimal, the claim $A(\gminia)$ is clear.
The claim $A(<_S\!\gminia)$
follows from $A(\gminib)$ for any $\gminib<_S\gminia$.
Let us show $A(\gminia)$
by assuming $A(<_S\gminia)$.
Let $f_{\gminia}$ be the morphism
$\Gr_{\gminia}\simeq
 \lefttop{1}E_{\gminia,S}\subset\nbigf^S_{\gminia}$.
Then, we obtain the following morphism
\[
 \DD_{z}(f_{\gminia}):
 \Gr_{\gminia}\lrarr
 \nbigf^S_{<\gminia}\otimes\Omega^1_{z}
=\bigoplus_{\gminib<\gminia}
 E_{\gminib,S}\otimes\Omega^1_{z}.
\]
We have $\DD_z(f_{\gminia})_{|\Zhat}=0$
by construction of $\nbigf^S$.

Note that $S$ is the product of
a multi-sector
$S_{\vecz}\subset(\Delta^{\ast})^k$
and $U\subset Y\times\nbigk^{\ast}$
which is a sector or a compact region.
The closure of $S_{\vecz}$ in $\Delta^k$
contains $0$.
Take a point $Q$ of $S_{\vecz}$.   
We take 
$g_{\gminia,\gminib}:
 \Gr_{\gminia}\lrarr E_{\gminib,S}$
such that
$\DD_{z}(g_{\gminia,\gminib})=
 (\DD_{z} f_{\gminia})_{\gminib}$,
and
$g_{\gminia,\gminib|Q\times U}=0$.
By sufficiently shrinking $S$ in the radius direction,
we can apply Lemma \ref{lem;07.9.30.50}
in this situation,
and we obtain $g_{\gminia,\gminib|\Zhat}=0$.
We put
$\fbar_{\gminia}:=
 f_{\gminia}-\sum_{\gminib<\gminia}g_{\gminia,\gminib}$,
then $\fbar_{\gminia}:\Gr_{\gminia}\lrarr\nbigf^S_{\gminia}$
is $\DD_{z}$-flat
and $\fbar_{\gminia}(\vecvbar_{\gminia})_{|\Zhat}
= \vecvhat_{\gminia|\Zhat}$.
Thus, we obtain $A(\gminia)$,
and the second condition in
Proposition \ref{prop;10.5.28.10} is satisfied.

Assume $\varrho$ is nowhere vanishing.
If $\gminib>_S\gminic$,
any $\DD_{z}$-flat morphism
$\Gr_{\gminib}^{\nbigf}(E_{|S})
\lrarr
 \Gr_{\gminic}^{\nbigf}(E_{|S})$
has the order
$O\bigl(\exp(-\epsilon|\vecz^{\vecm}|)\bigr)$
due to Corollary \ref{cor;07.12.27.4}.
Then, the third condition in 
Proposition \ref{prop;10.5.28.10} is satisfied.
\hfill\qed

\vspace{.1in}
Let $\Multisector^{\ast}(\nbigx\setminus W,\nbigi)$
denote the set of 
$S\in\Multisector(\nbigx\setminus W,\nbigi)$
such that $E_{|\Sbar}$ has
a $\DD$-flat filtration $\nbigf^S$
as in Proposition {\rm\ref{prop;10.5.28.10}}.
If $\varrho$ is invertible and 
$\nbigd$ is smooth, we have
$\Multisector^{\ast}(\nbigx\setminus W,\nbigi)
=\Multisector(\nbigx\setminus W,\nbigi)$.
\index{set $\Multisector^{\ast}(\nbigx
 \setminus W,\nbigi)$}

\begin{cor}
\label{cor;07.6.4.11}
For any point $P\in\gbigz$,
there exists 
$\nbigu_P\in
 \gbigu\bigl(P,\nbigxtilde(W),\nbigi
 \bigr)$
such that, for any 
$S\in \Multisector(P,\nbigx\setminus W,\nbigi)$
with $\Sbar\subset\nbigu_P$,
there uniquely exists 
a $\DD$-flat filtration $\nbigf^S$
of $E_{|\Sbar}$ indexed by $(\nbigi,\leq_S)$
satisfying the conditions in
Proposition {\rm\ref{prop;10.5.28.10}}.

Let $\Multisector^{\ast}(P,\nbigx\setminus W,\nbigi)$
denote the set of such multi-sectors.
\hfill\qed
\end{cor}
\index{set $\Multisector^{\ast}(P,
\nbigx\setminus W,\nbigi)$}

\begin{rem}
Even if $(E,\DD)$ is a weakly good lattice
in the level $(\vecm,i(0))$,
$(\End(E),\DD)$ is not necessarily
a weakly good lattice in the level $(\vecm,i(0))$.
We may not have the Stokes filtration of $\End(E)_{|\Sbar}$
for $S\in \Multisector^{\ast}(\nbigx\setminus W,\nbigi)$.
However,
as remarked 
in Section {\rm\ref{subsection;08.10.31.11}},
the $\DD$-flat subbundles
$\nbigf_{0}\End(E)_{|\Sbar}$ and
$\nbigf_{<0}\End(E)_{|\Sbar}$  are well defined,
which will be implicitly used.
\index{subbundle $\nbigfminusEndESbar$}
\index{subbundle $\nbigfzeroEndESbar$}
\hfil\qed
\end{rem}

\subsection{Compatibility}
\label{subsection;08.9.28.22}

\begin{lem}
\label{lem;07.12.13.30}
Let $S,S'\in 
 \Multisector^{\ast}(\nbigx\setminus W,\nbigi)$
such that $S'\subset S$.
\begin{itemize}
\item
The filtrations
$\nbigf^{S}$ and $\nbigf^{S'}$
are compatible over
$\bigl(\nbigi,\leq_S\bigr)\lrarr 
 \bigl(\nbigi,\leq_{S'}\bigr)$ 
in the sense of Definition {\rm\ref{df;07.12.13.20}}.
In particular,
we have $\nbigf^{S}=\nbigf^{S'}$
if $\bigl(\nbigi,\leq_S\bigr)$ is isomorphic to
$\bigl(\nbigi,\leq_{S'}\bigr)$.
\item
The decomposition
$E_{|\Sbar'}
=\bigoplus_{\gminia\in\nbigi}
 E_{\gminia,S|\Sbar'}$
gives a $\DD_{z}$-flat
splitting of the filtration $\nbigf^{S'}$.
\end{itemize}
\end{lem}
\pf
It follows from the characterization of
the Stokes filtrations (Proposition \ref{prop;10.5.28.10}).
\hfill\qed

\subsection{Splitting with nice property}
\label{subsection;07.12.7.1}

We have the induced
$\varrho$-flat connection
$\lefttop{j}\DD$ of $E_{|\nbigd_{Y,j}}$.
Since $\nbigf^S$ is $\DD$-flat,
$\Res_{Y,j}(\DD)$ preserves the filtration
$\nbigf^S_{|\nbigd_{Y,j}}$
and $\lefttop{j}\DD$.
Assume that we are given filtrations
$\lefttop{j}F$ $(j=1,\ldots,\ell)$
of $E_{|\nbigd_{Y,j}}$
which are preserved by
$\Res_{Y,j}(\DD)$ and $\lefttop{j}\DD$.

\begin{prop}
\label{prop;07.9.30.2}
Let $S\in\Multisector^{\ast}(\nbigx\setminus W,\nbigi)$.
We have a $\DD_{z}$-flat splitting
of the filtration $\nbigf^S$,
whose restriction to $\Sbar\cap \nbigd_{Y,j}$
is compatible with the residues
$\Res_{Y,j}(\DD)$ 
and the filtrations $\lefttop{j}F$
for $j=1,\ldots,\ell$,
after we shrink $S$ in the radius direction
appropriately.
\end{prop}
\pf
Take a large number $N$.
Let $\What^{(N)}$ denote
the $N$-th infinitesimal neighbourhood
of $W$.
By Lemma \ref{lem;10.6.30.3} below,
we can take a decomposition
$E=\bigoplus_{\gminia\in \nbigi} E_{\gminia,N}$
such that
(i) it is the same as  the irregular decomposition
 in the level $\vecm$  on $\What^{(N)}$,
(ii) the restriction to $\nbigd_{Y,j}$
 is compatible with $\Res_j(\DD)$
 and $\lefttop{j}F$
 for each $j=1,\ldots,\ell$.

Recall that $S$ is the product of
a multi-sector $S_{\vecz}\subset(\Delta^{\ast})^k$
and $S_{0}\subset Y\times\nbigk^{\ast}$,
where $S_{0}$ is a sector
or a compact region.
The closure of $S_{\vecz}$ in $\Delta_z^k$
contains the origin $O_z$.
Let $\Sbar_{0}$ denote the closure
of $S_{0}$ in
the real blow up of
$Y\times\nbigk$ along $Y\times\nbigk^0$.
Let $Q$ be any point of $S_{\vecz}$.
We have the following morphisms:
\begin{equation}
 \label{eq;07.9.30.101}
 \nbigf^S_{\gminia|Q\times \Sbar_{0}}
\lrarr
 E_{|Q\times \Sbar_{0}}
\lrarr
 \bigoplus_{\gminib\leq_S\gminia}
 E_{\gminib,N|Q\times\Sbar_{0}}
\end{equation}
If $Q$ is sufficiently close to $O_z$,
the composite of the morphisms 
in (\ref{eq;07.9.30.101})
is an isomorphism.
Let $\nbigg_{\gminia}^Q\subset
 \nbigf^S_{\gminia|Q\times \Sbar_{0}}$
denote the inverse image of
$E_{\gminia,N|Q\times\Sbar_{0}}$.
We may assume
$\nbigf^S_{<\gminia|Q\times\Sbar_{0}}
\cap \nbigg_{\gminia}^Q=\{0\}$.
Then, $\nbigg_{\gminia}^Q$ can be extended
to a $\DD_{z}$-flat subbundle
$\nbigg_{\gminia}^S$ of $E_{|S}$.
If $S$ is shrinked in the radius direction appropriately,
they give subbundles of $E_{|\Sbar}$,
due to Corollary \ref{cor;07.12.17.10}.
By construction,
it gives a splitting with the desired property.
\hfill\qed

\subsubsection{Special case 1}

We consider a $\DD$-flat splitting
in the non-resonance case.
We assume that $\varrho$ is nowhere vanishing.
Let $O_{\zeta}$ denote the origin in 
$Y=\Delta_{\zeta}^n$.
Take $S \in \Multisector^{\ast}(\nbigx\setminus W,\nbigi)$
of the form $S_1\times U_1\times U_2$
where $S_1$ is a multi-sector in $(\Delta^{\ast})^k$
whose closure contains the origin,
$U_1$ is a neighbourhood of $O_{\zeta}$
in $Y$,
and $U_2$ is a small compact region in
$\nbigk$.

\begin{prop}
\label{prop;07.9.30.5}
Assume
$\alpha-\beta\not\in \seisuu$
for distinct eigenvalues $\alpha$, $\beta$
of $\Res_{Y,j}(\DD^f)_{|\nbigd^y_j}$
$(j=1,\ldots,\ell,\,\,y\in U_2)$.
Then, we have a $\DD$-flat splitting
of $\nbigf^S$,
whose restriction to $S\cap\nbigd_{Y,j}$ is compatible
with $\lefttop{j}F$ for each $j=1,\ldots,\ell$.
\end{prop}
\pf
Let $Q$ be as in the proof of 
Proposition \ref{prop;07.9.30.2},
where we construct the $\DD_z$-flat splitting of
the filtration $\nbigf^S$ on $Q\times U$.
In particular,
we have the splitting
\begin{equation}
 \label{eq;08.12.1.2}
E_{|Q\times O_{\zeta}\times U_2}=
 \bigoplus_{\gminia\in \nbigi}
 E_{\gminia}^{Q\times O_{\zeta}\times U_2},
\end{equation}
which is compatible with the endomorphisms
$\Res_j(\DD)$ and the filtrations $\lefttop{j}F$
$(j=1,\ldots,\ell)$.
By the assumption on the eigenvalues
of $\Res_j(\DD)$,
(\ref{eq;08.12.1.2}) is extended to a $\DD$-flat splitting
of $\nbigf^{S}$ on $Q\times U_1\times U_2$.
By extending it to a flat splitting of $\nbigf^S$ on $S$,
we obtain the desired splitting.
\hfill\qed

\subsubsection{Special case 2}

We consider the case $\nbigd$ is smooth,
i.e.,
$k=1$ and  $\nbigd_Y=\emptyset$.
\begin{lem}
\label{lem;08.1.25.51}
For any $P\in \gbigz$,
there exist
$S_P\in\Multisector^{\ast}(P,\nbigx\setminus W,\nbigi)$
and a $\DD$-flat splitting
of the filtration $\nbigf^{S_P}$ on $\Sbar_P$.
\end{lem}
\pf
According to Proposition \ref{prop;07.9.30.5},
we have only to consider the case
$\varrho(P)=0$.
We use the symbol $m$ instead of $\vecm$.
We would like to take
a $\DD$-flat morphism
$\phi_{\gminia}:
 \Gr_{\gminia}(E_{|\Sbar})
\lrarr E_{|\Sbar}$
for some $S\in\Multisector^{\ast}
 (P,\nbigx\setminus W,\nbigi)$.
We construct such a morphism
inductively with respect to the order $\leq_P$.
If $\gminia$ is minimal with respect to $\leq_P$,
we have nothing to do.
Assume that we have already taken such morphisms
for any $\gminib<_S\gminia$.

We have a $\DD_z$-flat splitting
$E_{|\Sbar}=\bigoplus E_{\gminic,S}$
of the filtration $\nbigf^S$.
By the hypothesis of the induction,
we may assume that
$\nbigf^S_{<\gminia}(E_{|\Sbar})
=\bigoplus_{\gminib<_P\gminia}E_{\gminib}$ 
is $\DD$-flat.
Let $f_{\gminia}$ be the morphism
$\Gr_{\gminia}(E_{|\Sbar})
 \simeq E_{\gminia,S}\lrarr 
E_{|\Sbar}\simeq 
 \bigoplus \Gr_{\gminib}(E_{|\Sbar})$.
We remark $dz_1$-component of 
$\DD(f_{\gminia})$ is $0$.
We have the decomposition
$\DD(f_{\gminia})=
 \sum_{\gminib<_S\gminia}
 \DD(f_{\gminia})_{\gminib}$
corresponding to the $\DD$-flat decomposition
$\nbigf^S_{<\gminia}
 =\bigoplus_{\gminib<_S\gminia}
 E_{\gminib,S}$:
\[
 \DD(f_{\gminia})_{\gminib}:
 \Gr_{\gminia}(E_{|\Sbar})
\lrarr
 \Gr_{\gminib}(E_{|\Sbar})
\otimes 
 \Omega^1_{\nbigx/\nbigk}
 \otimes\nbigo(m\,D)
\]
Let $\vecvbar_{\gminia}$ be a holomorphic frame
of $\Gr_{\gminia}(E_{|\Sbar})$
for each $\gminia\in\nbigi$.
Let $R_{\gminia}$ be determined by
$\DD_{\gminia}\vecvbar_{\gminia}
=\vecvbar_{\gminia}\,(d\gminia+R_{\gminia})$.
We have $R_{\gminia}=O\bigl(|z_1^{m+1}|\bigr)$.
Let $A$ be determined by
 $\DD(f_{\gminia})_{\gminib}\vecvbar_{\gminia}
 =\vecvbar_{\gminib}\, A$.
Since $\DD(f_{\gminia})_{\gminib}$ is $\DD_z$-flat,
we have the following estimate for some $C>0$:
\[
 A\,
 \exp\bigl(\varrho^{-1}(\gminib-\gminia)\bigr)=
\left\{\begin{array}{ll}
 O\bigl(
 \exp\bigl(C\,|\varrho^{-1}\, z_1^{m+1}|\bigr)
 \bigr) & (m<-1)\\
 \mbox{{}}\\
 O\bigl(\exp(C\, |\varrho^{-1}|\,
 \log|z_1^{-1}|)\bigr) & (m=-1)
 \end{array}
\right.
\]
By shrinking $S$,
we obtain the estimate
$A=O\bigl(
 \exp\bigl(-\epsilon|\varrho^{-1}\,z_1^{m}|\bigr)
 \bigr)$
for some $\epsilon>0$.
Let $P_1:=\pi(P)\in \nbigd^0$.
Recall that we have assumed
$\nbigk=\Delta^1_{\varrho}\times\nbigk'$,
which induces
$\nbigk\lrarr \nbigk^0$
and hence $q:\nbigx\lrarr\nbigd^0$.
If we shrink $S$,
we can take a section
$g_{\gminia,\gminib}:
 \Gr_{\gminia}(E_{|S})
 \lrarr
 \Gr_{\gminib}(E_{|S})$
satisfying the following conditions:
\begin{itemize}
\item
 $\DD(g_{\gminia,\gminib})=
 \DD(f_{\gminia})_{\gminib}$
 and $g_{\gminia,\gminib|q^{-1}(P_1)\cap S}=0$.
\item 
 $g_{\gminia,\gminib}=O\bigl(
 \exp\bigl(-\epsilon|\varrho^{-1}z_1^{m}|/2\bigr)
 \bigr)$.
\end{itemize}
We put
$\phi_{\gminia}:=
 f_{\gminia}-\sum g_{\gminia,\gminib}$.
Because of the estimate for $g_{\gminia,\gminib}$,
the $\DD$-flat morphism
$\phi_{\gminia}$ is extended on $\Sbar$.
Thus, the inductive argument can proceed.
\hfill\qed

\subsection{Functoriality}
\label{subsection;08.9.8.1}

\subsubsection{Dual}
Let $(E,\DD,\nbigi)$ be a weakly good lattice
in the level $(\vecm,i(0))$ on $(\nbigx,\nbigd)$.
Let $S\in\Multisector^{\ast}(\nbigx\setminus W,\nbigi)$.
Let $E_{|\Sbar}=\bigoplus_{\gminia\in\nbigi} E_{\gminia,S}$
be a $\DD_{z}$-flat splitting of
the filtration $\nbigf^S(E_{|\Sbar})$
whose restriction to $\nbigd_{Y,j}$
is compatible with
$\Res_{Y,j}(\DD)$ 
for each $j=1,\ldots,\ell$.
For any $\gminia\in\nbigi$,
we put
\begin{equation}
\label{eq;08.1.3.1}
  E^{\lor}_{-\gminia,S}:=
 \left(
\bigoplus_{\substack{
 \gminib\in \nbigi\\
 \gminib\neq \gminia}}
E_{\gminib,S}
 \right)^{\bot},
\quad\quad
\nbigf^S_{\not\geq \gminia}(E_{|\Sbar}):=
 \sum_{\substack{\gminib\in \nbigi\\
 \gminib\not\geq_S\gminia
 }}
 \nbigf^S_{\gminib}(E_{|\Sbar}).
\end{equation}

\begin{lem}
\mbox{{}}\label{lem;07.12.7.2}
\begin{itemize}
\item
The decomposition
$E^{\lor}_{|\Sbar}
=\bigoplus_{\gminia\in\nbigi}
 E^{\lor}_{-\gminia,S}$
gives a $\DD^{\lor}_{z}$-flat splitting 
of the Stokes filtration
$\nbigf^S(E^{\lor}_{|\Sbar})$,
whose restriction to $\nbigd_{Y,j}$
is compatible with
$\Res_{Y,j}(\DD^{\lor})$ 
for each $j=1,\ldots,\ell$.
\item
In particular,
$\nbigf^S_{-\gminia}(E^{\lor}_{|\Sbar})
=\nbigf^S_{\not\geq\gminia}(E_{|\Sbar})^{\bot}$
for any $\gminia\in\nbigi$.
\end{itemize}
\end{lem}
\pf
The second claim follows from the first claim.
We put
\[
 \nbigf^{\prime\,S}_{-\gminia}(E^{\lor}_{|\Sbar})
 :=
 \bigoplus_{\substack{
 \gminib\in\nbigi\\
 -\gminib\leq_S-\gminia
 }}
  E^{\lor}_{-\gminib,S}.
\] 
For the first claim,
we have only to show
$\nbigf^{\prime\,S}_{-\gminia}(E^{\lor}_{|\Sbar})
 =\nbigf^S_{-\gminia}(E^{\lor}_{|\Sbar})$.
It is easy to check
$\nbigf^S_{-\gminia}(E^{\lor}_{|\Sbar})_{|\Zhat}
=\nbigf^Z_{-\gminia}(E^{\lor}_{|\Zhat})$.
Then, the first claim follows from the uniqueness
in Proposition \ref{prop;10.5.28.10}.
\hfill\qed

\subsubsection{Tensor product and direct sum}
Let $(E_p,\DD_p,\nbigi_p)$ $(p=1,2)$ be 
weakly good lattices
in the level $(\vecm,i(0))$
on $(\nbigx,\nbigd)$.
Let us consider the case that
$\nbigi_1\otimes\nbigi_2$
is weakly good
in the level $(\vecm,i(0))$.
We put $(\Etilde,\DDtilde):=
 (E_1,\DD_1)\otimes (E_2,\DD_2)$.
For a multi-sector $S\in
 \bigcap_{p=1,2}
 \Multisector^{\ast}\bigl(\nbigx\setminus W,\nbigi_p\bigr)$,
we take $\DD_{z}$-flat splittings
$E_{p|\Sbar}=
 \bigoplus_{\gminia_p\in \nbigi_p} E_{p,\gminia_p,S}$
whose restrictions to $\nbigd_{Y,j}$
are compatible with $\Res_{Y,j}(\DD_p)$.
We put 
\begin{equation}
\label{eq;07.12.15.100}
 \Etilde_{\gminia,S}
:=\bigoplus_{\substack{
 (\gminia_1,\gminia_2)\in
 \nbigi_1\times \nbigi_2\\
 \gminia_1+\gminia_2=\gminia}}
 E_{1,\gminia_1,S}
\otimes
 E_{2,\gminia_2,S}.
\end{equation}
The following lemma can be shown by
an argument in the proof of
Lemma \ref{lem;07.12.7.2}.
\begin{lem}
\mbox{{}}\label{lem;07.12.7.3}
\begin{itemize}
\item
The decomposition
$\Etilde_{|\Sbar}=
 \bigoplus_{\gminia\in\nbigi_1\otimes\nbigi_2}
 \Etilde_{\gminia,S}$ gives a $\DD_{z}$-flat
 splitting of the Stokes filtration
 of $\Etilde_{|\Sbar}$,
 whose restriction to $\nbigd_{Y,j}$ is compatible with
 $\Res_{Y,j}(\DDtilde)$ for each $j=1,\ldots,\ell$.
\item
In particular,
$\nbigf^S_{\gminia}(\Etilde_{|\Sbar})$
is equal to
$\sum_{\gminia_1+\gminia_2\leq_S\gminia}
 \nbigf^S_{\gminia_1}(E_{1|\Sbar})
\otimes
 \nbigf^S_{\gminia_2}(E_{2|\Sbar})$.
\hfill\qed
\end{itemize}
\end{lem}

Let us consider the case in which
$\nbigi_1\oplus\nbigi_2:=\nbigi_1\cup\nbigi_2$ 
is weakly good in the level $(\vecm,i(0))$.
For a multi-sector $S\in
\bigcap_{p=1,2}
 \Multisector^{\ast}\bigl(
 \nbigx\setminus W,\nbigi_p\bigr)$,
we take $\DD_z$-flat splittings of
$E_{p|\Sbar}=\bigoplus E_{p,\gminia,S}$
whose restriction to $\nbigd_{Y,j}$ 
is compatible with $\Res_{Y,j}(\DD)$.
We put
\begin{equation}
\label{eq;07.12.15.101}
 (E_1\oplus E_2)_{\gminia,S}:=
 E_{1,\gminia,S}\oplus E_{2,\gminia,S}
\end{equation}
The following lemma is obvious.
\begin{lem}
\mbox{{}}
\begin{itemize}
\item
The decomposition
$(E_1\oplus E_2)_{|\Sbar}=
 \bigoplus_{\gminia\in\nbigi_1\oplus\nbigi_2}
 (E_1\oplus E_2)_{\gminia,S}$ gives 
a $\DD_{z}$-flat splitting of the Stokes filtration
 of $(E_1\oplus E_2)_{|\Sbar}$,
 whose restriction to $\nbigd_{Y,j}$ is 
 compatible with $\Res_{Y,j}(\DD)$.
\item
In particular,
$\nbigf^S_{\gminia}\bigl((E_1\oplus E_2)_{|\Sbar}\bigr)$
is equal to
$\nbigf^S_{\gminia}(E_{1|\Sbar})
\oplus
 \nbigf^S_{\gminia}(E_{2|\Sbar})$.
\hfill\qed
\end{itemize}
\end{lem}

\subsubsection{Morphism}

Let $(E_p,\DD_p,\nbigi_p)$ $(p=1,2)$
be as above.
Let $F:(E_1,\DD_1)\lrarr(E_2,\DD_2)$
be a flat morphism.
We assume that $\nbigi_1\cup\nbigi_2$ is 
weakly good in the level $(\vecm,i(0))$.
\begin{lem}
\label{lem;07.12.7.6}
Let $S\in\bigcap_{p=1,2}
  \Multisector^{\ast}(\nbigx\setminus W,\nbigi_p)$.
The restriction $F_{|\Sbar}$ preserves
the Stokes filtrations.
\end{lem}
\pf
Let $\vecvhat_p$ be frames of $E_{p|\What}$
compatible with the decompositions
$E_{p|\What}
=\bigoplus \Ehat_{p,\gminia|\What}$.
Let $\Ahat$ be determined by
$F(\vecvhat_1)=\vecvhat_2\,\Ahat$.
We have the decomposition
$\Ahat=\bigoplus \Ahat_{\gminia,\gminia}$.
Let $S\in\bigcap_{p=1,2}
 \Multisector^{\ast}(\nbigx\setminus W,\nbigi_p)$.
We have the Stokes filtrations
$\nbigf^S(E_{p|\Sbar})$ $(p=1,2)$.
We take $\DD_z$-flat splittings
$E_{p|\Sbar}=\bigoplus E_{p,\gminia,S}$
of the filtrations $\nbigf^S(E_{p|\Sbar})$,
and let $\vecv_{p,S}$ be lifts of
$\vecvhat_p$ to $E_{p|\Sbar}$
compatible with the splittings.
Let $F_{\gminib,\gminia}$ be 
the $\DD_{z}$-flat morphism
$E_{1,\gminia,S}\lrarr E_{2,\gminib,S}$
induced by $F_{|\Sbar}$.
Let $A_{\gminib,\gminia}$
be determined by
$F_{\gminib,\gminia}(\vecv_{1,\gminia,S})
=\vecv_{2,\gminib,S}\,A_{\gminib,\gminia}$.
Lemma \ref{lem;07.12.7.6}
can be reduced to the following lemma.

\begin{lem}
\label{lem;07.11.25.50}
\mbox{{}}
\begin{itemize}
\item $A_{\gminib,\gminia}=0$
 unless $\gminia\geq_S\gminib$.
\item
 In the case $\gminia>_S\gminib$,
 we have the following estimate for some $C>0$:
\[
 A_{\gminib,\gminia}
\exp\bigl(\varrho^{-1}(\gminib-\gminia)\bigr)
=O\Bigl(
 \exp\bigl(C\,|\varrho^{-1}\vecz^{\vecm(1)}|
 \,\log|z_{i(0)}^{-1}| \bigr)
 \Bigr)
\]
\item
$A_{\gminia,\gminia|\Zhat}=\Ahat_{\gminia,\gminia}$.
In particular,
$\bigl|A_{\gminia,\gminia}\bigr|$ is bounded.
\end{itemize}
\end{lem}
\pf
The first two claims follow from 
Corollary \ref{cor;07.12.27.3}
and Corollary \ref{cor;07.12.27.4}.
The third claim is clear.
Thus, we obtain 
Lemma \ref{lem;07.11.25.50}
and Lemma \ref{lem;07.12.7.6}.
\hfill\qed

\begin{cor}
\mbox{{}}\label{cor;07.12.7.7}
\begin{itemize}
\item
If the restriction of $F$ to $\nbigx\setminus\nbigd$
is an isomorphism,
we have $\nbigi_1=\nbigi_2$
and 
$\nbigf^S_{\gminia}(E_{1|S\setminus \nbigd})
=\nbigf^S_{\gminia}(E_{2|S\setminus \nbigd})$.
\item
In particular, the Stokes filtration $\nbigf^S$
in the level $(\vecm,i(0))$ depends only on 
the meromorphic $\varrho$-flat bundle
$\bigl(E(\ast \nbigd),\DD\bigr)$,
in the sense that
it is independent of
the choice of a weakly good lattice
$E\subset E(\ast \nbigd)$ in the level $(\vecm,i(0))$.
\end{itemize}
\end{cor}
\pf
$F$ induces an isomorphism
$E_1(\ast \nbigd)\simeq E_{2}(\ast \nbigd)$,
and hence
$E_{1}(\ast \nbigd)_{|\What_z}
\simeq
 E_2(\ast \nbigd)_{|\What_z}$.
Then, we obtain
an isomorphism
$\Ehat_{1,\gminia}(\ast \nbigd)
\simeq
 \Ehat_{2,\gminia}(\ast \nbigd)$
for $\gminia\in \nbigi_1\cup\nbigi_2$.
Hence, we have
$\nbigi_1=\nbigi_2$.
Since we have
the inclusion
$\nbigf^S_{\gminia}(E_{1|S})\subset
 \nbigf^S_{\gminia}(E_{2|S})$
by Lemma \ref{lem;07.12.7.6},
we obtain $\nbigf^S_{\gminia}(E_{1|S})=
 \nbigf^S_{\gminia}(E_{2|S})$
by the comparison of the rank.
Thus, the first claim is proved.
The second claim follows from the first claim.
\hfill\qed

\subsection{The associated graded bundle}
\label{subsection;07.12.12.1}

For any
$S\in\Multisector^{\ast}(\nbigx\setminus W,\nbigi)$
and $\gminia\in\nbigi$,
we obtain a bundle
$\Gr^{\vecm}_{\gminia}(E_{|\Sbar})$ 
with a meromorphic flat $\varrho$-connection
$\DD_{\gminia,\Sbar}$ on $\Sbar$,
by taking Gr with respect to the filtration $\nbigf^S$.
By definition of Stokes filtrations,
we have a natural isomorphism
\[
 \bigl(
 \Gr^{\vecm}_{\gminia}(E_{|\Sbar}),
 \DD_{\gminia,\Sbar}\bigr)
 _{|\Zhat}
\simeq
 (\Ehat_{\gminia},\DDhat_{\gminia})_{|\Zhat}.
\]
When $S$ is varied,
we can glue them and obtain a bundle
$\Gr^{\vecm}_{\gminia}(E_{|\nbigvtilde(W)})$
with a meromorphic flat $\varrho$-connection
$\DD_{\gminia,\nbigvtilde(W)}$
on $\nbigvtilde(W)$,
where $\nbigv$ denotes some neighbourhood 
of $\nbigd_{z,\kbar}$,
and $\nbigvtilde(W)$
denotes the real blow up of $\nbigv$
along $W':=\nbigv\cap W$.
By the construction,
we are given an isomorphism
\[
\bigl(
 \Gr^{\vecm}_{\gminia}
 (E_{|\nbigvtilde(W)}),
 \DD_{\gminia,\nbigvtilde(W)}\bigr)
_{|\widehat{\pi^{-1}(W')}}
\simeq
 (\Ehat_{\gminia},\DDhat_{\gminia})
 _{|\widehat{\pi^{-1}(W')}}.
\]
According to Proposition \ref{prop;08.9.4.102},
Corollary \ref{cor;08.9.5.1} and
Lemma \ref{lem;08.9.5.12},
there exists a holomorphic vector bundle
$\Gr^{\vecm}_{\gminia}(E)$
with a meromorphic flat $\varrho$-connections
$\DD_{\gminia}$ on $(\nbigv,\nbigd\cap\nbigv)$
such that
\begin{equation}
\label{eq;08.9.5.10}
 \pi^{\ast}
 \bigl(\Gr^{\vecm}_{\gminia}(E),
 \DD_{\gminia}
 \bigr)_{|\nbigvtilde(W)}
\simeq
 \bigl(\Gr^{\vecm}_{\gminia}(
 E_{|\nbigvtilde(W)}),
 \DD_{\gminia,\nbigvtilde(W)}
 \bigr),\,\,\,
 \bigl(\Gr^{\vecm}_{\gminia}(E),
 \DD_{\gminia}
 \bigr)_{|\What'}
\simeq
 (\Ehat_{\gminia},\DDhat_{\gminia})
 _{|\What'}.
\end{equation}
It is well defined on the germ of
neighbourhoods of $\nbigd_{z,\kbar}$
in $\nbigx$.

\index{Gr $\Gr^{\vecm}_{\gminia}(E)$}

\begin{cor}
If we shrink $X$,
we can take a frame 
$\vecvhat=(\vecvhat_{\gminia})$ 
of $E_{|\What}$
compatible with the irregular decomposition
in the level $\vecm$,
such that 
the power series
$\Rhat_{\gminia}$ is convergent,
where $\DD_{\gminia}\vecvhat_{\gminia}
 =\vecvhat_{\gminia}\, \Rhat_{\gminia}$.
\end{cor}
\pf
We take a holomorphic frame $\vecw_{\gminia}$
of $\Gr^{\vecm}_{\gminia}(E)$ on $\nbigv$.
It induces a frame $\vecw_{\gminia|\What}$ of 
$\Gr^{\vecm}_{\gminia}(E)
 _{|\What}\simeq \Ehat_{\gminia}$.
We have only to put 
$\vecvhat_{\gminia}:=\vecw_{\gminia|\What}$.
\hfill\qed

\vspace{.1in}

Let $\vecw_{\gminia}$ be a frame of
$\Gr^{\vecm}_{\gminia}(E)$.
Let $S\in\Multisector^{\ast}(\nbigx\setminus W,\nbigi)$,
and let $E_{|\Sbar}=\bigoplus E_{\gminia,S}$
be a $\DD_{z}$-flat splitting
of the Stokes filtration $\nbigf^S$.
By the natural isomorphism
$E_{\gminia,S}\simeq
 \Gr_{\gminia}^{\vecm}(E)_{|\Sbar}$,
we take a lift $\vecw_{\gminia,S}$ of $\vecw_{\gminia}$.
Thus, we obtain a frame
$\vecw_S=\bigl(\vecw_{\gminia,S}\bigr)$.
The following corollary is clear from the construction.
\begin{cor}
\label{cor;07.10.18.1}
Let $\vecv$ be a frame of $E$,
and let $G_S$ be the matrix 
by $\vecv=\vecw_S\, G_S$.
Then, $G_S$ and $G_S^{-1}$ are bounded on $S$.
\hfill\qed
\end{cor}

\begin{rem}
\label{rem;07.12.21.22}
If $k=1$ and if
$\varrho$ is nowhere vanishing,
we obtain 
$\bigl(\Gr^{\vecm}_{\gminia}(E),\,
 \DD_{\gminia}\bigr)$ 
on $(\nbigx,\nbigd)$,
not only on $(\nbigv,\nbigd\cap\nbigv)$.
Note that we can always extend
a good lattice of
a meromorphic $\varrho$-flat bundle
on $(\nbigv,\nbigd\cap\nbigv)$
to that on $(\nbigx,\nbigd)$.
\hfill\qed
\end{rem}

\subsubsection{Functoriality}
\label{subsection;08.9.8.3}

Let $(E,\DD,\nbigi)$ be a good lattice
in the level $(\vecm,i(0))$.
From $(E^{\lor},\DD^{\lor},\nbigi^{\lor})$,
we obtain the associated graded bundle
$\Gr^{\vecm}(E^{\lor})
=\bigoplus_{\gminia\in\nbigi^{\lor}}
 \Gr^{\vecm}_{\gminia}(E^{\lor})$
with an induced meromorphic
flat $\varrho$-connection.

\begin{lem}
 \label{lem;07.12.7.5}
We have a natural flat isomorphism
$\Gr^{\vecm}_{\gminia}(E^{\lor})
\simeq
 \Gr^{\vecm}_{-\gminia}(E)^{\lor}$.
\end{lem}
\pf
By Lemma \ref{lem;07.12.7.2},
we have the natural isomorphism
$\Gr^{\vecm}_{\gminia}(E^{\lor})_{
 |\nbigvtilde(W)}
\simeq
 \Gr^{\vecm}_{-\gminia}(E)^{\lor}_{
 |\nbigv(W)}$.
It induces the desired isomorphism
on $\nbigv$.
\hfill\qed

\vspace{.1in}

Let $(E_p,\nabla_p,\nbigi_p)$ $(p=1,2)$ be 
weakly good lattices in the level $(\vecm,i(0))$.
The following lemma can be shown
similarly.
\begin{lem}
\label{lem;07.12.7.25}
Assume $\nbigi_1\otimes\nbigi_2$ is 
weakly good in the level $(\vecm,i(0))$.
Let $(\Etilde,\DDtilde):=
 (E_p,\DD_1)\otimes (E_p,\DD_2)$.
Then, we have the following natural isomorphism
for each $\gminia\in\nbigi_1\otimes\nbigi_2$:
\begin{equation}
 \label{eq;07.12.7.4}
 \Gr_{\gminia}^{\vecm}(\Etilde)
\simeq
\bigoplus_{\substack{
 (\gminia_1,\gminia_2)\in
 \nbigi_1\times
 \nbigi_2\\
 \gminia_1+\gminia_2=\gminia}}
 \Gr^{\vecm}_{\gminia_1}(E_1)
\otimes
 \Gr^{\vecm}_{\gminia_2}(E_2)
\end{equation}
If $\nbigi_1\oplus\nbigi_2$ is 
weakly good in the level $(\vecm,i(0))$,
then we have 
$\Gr^{\vecm}_{\gminia}(E_1\oplus E_2)
\simeq
 \Gr^{\vecm}_{\gminia}(E_1)
\oplus
 \Gr^{\vecm}_{\gminia}(E_2)$.
\hfill\qed
\end{lem}

\begin{lem}
\label{lem;07.12.7.17}
Let $F:(E_1,\DD_1)\lrarr (E_2,\DD_2)$
be a flat morphism.
Assume that $\nbigi_1\cup\nbigi_2$ is good,
for simplicity.
Then, we have the naturally induced morphism
\[
 \Gr^{\vecm}_{\gminia}(F):
 \Gr^{\vecm}_{\gminia}(E_1)
\lrarr \Gr^{\vecm}_{\gminia}(E_2).
\]
\end{lem}
\pf
We have the induced morphism
$\Gr^{\vecm}_{\gminia}(E_1)
 _{|\nbigvtilde(W)}
\lrarr
 \Gr^{\vecm}_{\gminia}(E_2)
 _{|\nbigvtilde(W)}$
by Lemma \ref{lem;07.12.7.6}.
It induces the desired morphism on $\nbigv$.
\hfill\qed

\begin{cor}
\label{cor;07.12.16.3}
In the situation of Lemma {\rm\ref{lem;07.12.7.17}},
if $E_{1|\nbigx-\nbigd}\lrarr 
 E_{2|\nbigx-\nbigd}$ is an isomorphism,
we have induced isomorphisms:
\[
 \Gr^{\vecm}_{\gminia}(E_1)\otimes\nbigo(\ast \nbigd)
\simeq
 \Gr^{\vecm}_{\gminia}(E_2)\otimes\nbigo(\ast \nbigd)
\quad
(\gminia\in\nbigi)
\]
Hence, the graded meromorphic $\varrho$-flat bundle
\[
 \bigoplus_{\gminia}
 \bigl(\Gr^{\vecm}_{\gminia}(E)\otimes
 \nbigo(\ast \nbigd),\DD_{\gminia}\bigr)
\]
is well defined for
the meromorphic $\varrho$-flat bundle
$\bigl(E(\ast \nbigd),\DD\bigr)$.
\end{cor}
\pf
By Corollary \ref{cor;07.12.7.7},
the restriction 
$\Gr^{\vecm}_{\gminia}(E_1)_{|\nbigv\setminus D}
\simeq
 \Gr^{\vecm}_{\gminia}(E_2)_{|\nbigv\setminus D}$
is an isomorphism.
Hence, the induced morphism
$\Gr^{\vecm}_{\gminia}(E_1)\otimes\nbigo(\ast\nbigd)
\lrarr
 \Gr^{\vecm}_{\gminia}(E_2)\otimes\nbigo(\ast \nbigd)$
is an isomorphism.
\hfill\qed

\subsection{A characterization  by the growth order}
\label{subsection;08.9.28.23}
Assume that
$\varrho$ is nowhere vanishing.
Let $(E,\DD,\nbigi)$ be a good lattice
in the level $(\vecm,i(0))$.
Take any frame $\vecv$ of $E$.
Let $S\in\Multisector^{\ast}
 (\nbigx\setminus\nbigd_z,\nbigi)$.
Let $f$ be a $\DD_{z}$-flat section of $E_{|S}$.
We have the expression $f=\sum f_j\,v_j$,
and obtain $\vecf=(f_j)$.

\begin{lem}
\label{lem;10.5.11.24}
We have $f\in\nbigf^S_{\gminib}$
if and only if the following holds
for some $C>0$:
\[
 \bigl|
 \vecf\,\exp\bigl(\varrho^{-1}\,\gminib\bigr)
 \bigr|
=O\Bigl(
 \exp\bigl(C\,|\vecz^{\vecm(1)}|
 \bigr)
 |z_{i(0)}|^{-C}
 \Bigr)
\]
Recall
$\vecm(1):=\vecm+\vecdelta_{i(0)}$.
\end{lem}
\pf
We take a $\DD_{z}$-flat splitting
$E_{|\Sbar}=\bigoplus E_{\gminia,S}$ of $\nbigf^S$,
and take a frame $\vecv_S=(\vecv_{\gminia,S})$
of $E_{|\Sbar}$ compatible with the splitting.
Let $r(\gminia):=\rank E_{\gminia,S}$.
We have the expression
$f=\sum_{\gminia}\sum_{j=1}^{r(\gminia)}
 f_{\gminia,S,j}\,v_{\gminia,S,j}$,
and obtain 
$\vecf_{\gminia}=
(f_{\gminia,S,j}\,|\,j=1,\ldots,r(\gminia))$.
Note
$|\vecf|$ and $\sum |\vecf_{\gminia}|$
are mutually bounded.

Let $R_{\gminia}$ be determined by
$\DD^f_{z}\vecv_{\gminia\,S}=
 \vecv_{\gminia\,S}\,
 \Bigl(\bigoplus\bigl(d_{z}(\varrho^{-1}\,\gminia)
 +R_{\gminia}\bigr)\Bigr)$.
Then, 
it is a holomorphic section of
the following:
\[
\sum_{i=1}^k
 M_{r(\gminia)}(\cnum)\otimes
 \vecz^{\vecm(1)}\,\nbigo_{\Sbar}\,dz_i/z_i
\quad(\mbox{\rm if }m_{i(0)}<-1)
\]
\[
\sum_{i\neq i(0)}
M_{r(\gminia)}(\cnum)\otimes
 \vecz^{\vecm(1)}\,
 \nbigo_{\Sbar}\,dz_i/z_i
+M_{r(\gminia)}(\cnum)\otimes
 \nbigo_{\Sbar}\,dz_i/z_i
\quad
(\mbox{\rm if }m_{i(0)}=-1)
\]
Since each $f_{\gminia}$ is $\DD_{z}$-flat,
we obtain the following estimate
in the case $f_{\gminia}\neq 0$,
by using Lemma \ref{lem;07.12.27.5}:
\[
\Bigl|\log\bigl|
 \vecf_{\gminia}
 \exp\bigl(\varrho^{-1}\gminia\bigr)
 \bigr|\Bigr|
\leq
  C\,|\vecz^{\vecm(1)}|
+C\log\bigl|z_{i(0)}^{-1}\bigr|
\]
Then, the claim of the lemma follows.
\hfill\qed

\vspace{.1in}
Let us consider the case
$Y=\Delta_{\zeta}^n$
and $\nbigd_Y=
 \bigcup_{i=1}^{\ell}\{\zeta_i=0\}$.
Let $f$ be a $\DD$-flat section of 
$E_{|S\setminus\nbigd_Y}$.
We obtain the following lemma
from Lemma \ref{lem;07.12.27.10}
by the argument in the proof of
Lemma \ref{lem;10.5.11.24}.
\begin{lem}
\label{lem;10.5.12.51}
We have $f\in\nbigf^S_{\gminib}$
if and only if the following holds
for some $C>0$:
\[
 \bigl|
 \vecf\,\exp\bigl(\varrho^{-1}\,\gminib\bigr)
 \bigr|
=O\Bigl(
 \exp\bigl(C\,|\vecz^{\vecm(1)}|
 \bigr)
 |z_{i(0)}|^{-C}
 \prod_{j=1}^{\ell}
 |\zeta_j|^{-C}
 \Bigr)
\]
\hfill\qed
\end{lem}

\subsection{Proof of the claims in
Section \ref{section;10.5.11.1}}
\label{subsection;10.5.18.10}

Corollary \ref{cor;07.6.4.11}
implies Proposition \ref{prop;10.5.11.20}.
Lemma \ref{lem;07.12.13.30} implies
Proposition \ref{prop;10.5.11.21}.
The functoriality in Subsection 
\ref{subsection;10.5.11.22}
follows from those in
Subsection \ref{subsection;08.9.8.1}.
The growth estimate in Proposition
\ref{prop;10.5.11.23}
is implied by 
that in Lemma \ref{lem;10.5.11.24}.
The associated graded meromorphic
$\varrho$-flat bundle and its functoriality
are studied in Subsection \ref{subsection;07.12.12.1}.
Proposition \ref{prop;10.5.11.30}
is implied by Proposition \ref{prop;07.9.30.2}.
Proposition \ref{prop;10.5.11.31}
follows from Proposition \ref{prop;07.9.30.5}.
Proposition \ref{prop;10.5.11.32}
follows from Proposition \ref{prop;10.5.11.31}
and Lemma \ref{lem;10.5.11.40}.
Proposition \ref{prop;10.5.11.41}
is implied by Lemma \ref{lem;08.1.25.51}

\subsection{Appendix (Lifting of formal frames)}
\label{subsection;07.11.13.5}
We discuss liftings of frames.
Although we will use such concepts
in our later argument, readers can skip here.

\subsubsection{Holomorphic lift on small sectors}
\label{subsection;07.6.16.2}
\index{holomorphic lift}

We take a frame $\vecvhat=(\vecvhat_{\gminia})$
of $E_{|\What}$
compatible with the irregular decomposition.
Let $\Rhat_{\gminia}$ be determined by
$\DD_{\gminia}\vecvhat_{\gminia}
=\vecvhat_{\gminia}\,
 \bigl(d\gminia+\Rhat_{\gminia}\bigr)$.

Let $S\in\Multisector^{\ast}(\nbigx\setminus W,\nbigi)$.
We take a $\DD_z$-flat decomposition
$E_{|\Sbar}=\bigoplus_{\gminia} E_{\gminia,S}$
which gives a splitting of $\nbigf^{S}$
as in Proposition \ref{prop;10.5.28.10}.
We can take a frame $\vecv_{\gminia,S}$ of $E_{\gminia,S}$
such that
$\vecv_{\gminia,S|\Zhat}=
 \vecvhat_{\gminia|\Zhat}$,
and we put $\vecv_{S}:=(\vecv_{\gminia,S})$,
which is called a holomorphic lift of $\vecvhat$ on $\Sbar$.

Let $E_{|\Sbar}=
 \bigoplus_{\gminia}E'_{\gminia,S}$ be
another $\DD_z$-flat splitting
of $\nbigf^S$,
and let $\vecv'_{S}=\bigl(\vecv'_{\gminia,S}\bigr)$
be a holomorphic lift compatible with the splitting.
Let $C=(C_{\gminia',\gminia})$ be determined by
$\vecv_S=\vecv'_S\,(I+C)$,
where $I$ denotes the identity matrix.

\begin{lem}
\label{lem;07.6.15.110}
\mbox{{}}
\begin{itemize}
\item
We have $C_{\gminia',\gminia|\Zhat}=0$ and 
$C_{\gminia',\gminia}=0$
unless $\gminia'\leq_S\gminia$.
\item
If $\gminia'<_S\gminia$,
we have the following for some $C>0$:
\[
 C_{\gminia',\gminia}\,
 \exp\bigl(\varrho^{-1}(\gminia'-\gminia)
 \bigr)
=
 O\Bigl(\exp\bigl(
 C\,|\varrho^{-1}|\bigl(
 |\vecz^{\vecm(1)}|
+\log|z_{i(0)}^{-1}| \bigr)
 \bigr)\,
 \Bigr)
\]
\end{itemize}
\end{lem}
\pf
We have $C_{\gminia',\gminia|\Zhat}=0$
by construction.
Since $\vecv_S$ and $\vecv_S'$ are compatible with
the filtration $\nbigf^S$ by construction,
we have
$C_{\gminia',\gminia}=0$
unless $\gminia'\leq_S\gminia$.
The other follows from
the estimate of the norm of
$\DD_z$-flat sections
(Lemma \ref{lem;07.12.27.1}).
\hfill\qed

\subsubsection{$C^{\infty}$-lift on $X$}
\label{subsection;07.6.16.8}
\index{$C^{\infty}$-lift}

\begin{lem}
\label{lem;07.6.16.3}
We have a local $C^{\infty}$-frame
$\vecv_{C^{\infty}}=
 \bigl(\vecv_{\gminia,C^{\infty}}\bigr)$
of $E$
on some neighbourhood $\nbigv$
of $\nbigd_{z,\kbar}$
with the following property:
\begin{itemize}
\item
 $\vecv_{C^{\infty}|\nbigdhat_z}
=\vecvhat$.
\item
 Let $S\in\Multisector^{\ast}(\nbigx\setminus W,\nbigi)$,
 and let $\vecv_S$ be a holomorphic lift
 of $\vecvhat$ on $\Sbar$
 as in Section {\rm\ref{subsection;07.6.16.2}}.
 Let $B_S$ be determined by 
 $\vecv_{C^{\infty}}=\vecv_S\, (I+B_S)$
 on $\pi^{-1}(\nbigv)\cap \Sbar$.
 Then, the following holds:
\begin{itemize}
 \item
  $B_{S|\Zhat}=0$.
 \item
 Let $B_S=\bigl(B_{S,\gminia,\gminib}\bigr)$
 be the decomposition, corresponding to
 the decomposition of the frame
 $\vecvhat=(\vecvhat_{\gminia})$.
 Then, we have
  $B_{S,\gminia,\gminib}=0$ unless $\gminia\leq_S\gminib$.
 \item
 In the case $\gminia<_S\gminib$,
 we have the following estimate of
 the $C^q$-derivatives of
 $B_{S,\gminia,\gminib}\,
 \exp\bigl(\varrho^{-1}(\gminia-\gminib)\bigr)$
 for some $C>0$ 
 and $N(q)\geq 0$
$(q\in\seisuu_{\geq 0})$:
\[
 O\Bigl(
 \exp\Bigl(
 C\,|\varrho^{-1}|
 \bigl(
 |\vecz^{\vecm(1)}|
+\log\bigl|z_{i(0)}^{-1}\bigr|
 \bigr)
 \Bigr)
\prod_{j=1}^k|z_j|^{-N(q)}
 |\varrho^{-1}|^{-N(q)}
 \Bigr)
\]
 \end{itemize}
\end{itemize}
In particular,
the frame $\vecv_{C^{\infty}|S}$ is compatible with 
the Stokes filtration $\nbigf^S$
 for $S\in\Multisector^{\ast}(\nbigx-W,\nbigi)$
 in the sense that
 $\bigl(\vecv_{\gminib,C^{\infty}}\,\big|\,
 \gminib\leq_S\gminia \bigr)$
 gives a frame of $\nbigf^S_{\gminia}$
 for each $\gminia$. 

 Such a frame $\vecv_{C^{\infty}}$ is called 
 a $C^{\infty}$-lift of $\vecvhat$.
\end{lem}
\pf
In the following,
for a given multi-sector $S$,
let $S^{\circ}$ denote the interior part.
We take multi-sectors
$S^{(j)}\in
 \Multisector^{\ast}(\nbigx\setminus W,\nbigi)$ 
 $(j=1,\ldots,N)$
of $\nbigx\setminus W$ such that
$\bigcup S^{(j)\circ}=\nbigv\setminus W$,
where $\nbigv$ is some open neighbourhood of 
$\nbigd_{z,\kbar}$ in $\nbigx$. 
We take holomorphic lifts 
$\vecv_{S^{(j)}}
=\bigl(v_{S^{(j)},i}\bigr)$
 of $\vecvhat$ on $\Sbar^{(j)}$.
We have only to glue them
in $C^{\infty}$-sense as follows,
for example.
We take small sectors
$S^{(j)}_3\subset 
 S^{(j)}_2\subset
 S^{(j)}_1=S^{(j)}$
such that
(i) $\bigcup_{j=1}^N S^{(j)\circ}_3
 =\nbigv\setminus W$,
(ii) $S^{(j)}_a$ in $\nbigv\setminus W$
is contained in $S^{(j)\circ}_{a-1}$ for $a=2,3$.
We take $C^{\infty}$-functions $\chi_j$ 
on $\nbigv\setminus W$
such that
(i) $\chi_j\geq 0$,
(ii) $\chi_j>0$ on $S^{(j)}_3$,
and $\chi_j=0$ outside of $S^{(j)}_2$,
(iii) each 
 $(\del_{\varrho}^m\prod\del_i^{m_i})\chi_j$
is polynomial order 
in $|\varrho^{-1}|$
and $|z_i^{-1}|$ $(i=1,\ldots,k)$,
(iv) $\sum_{j=1}^N \chi_j=1$.
We put
$\vecv_{C^{\infty}}:=
 \sum_{j=1}^N \chi_j\,\vecv_{S^{(j)}}$,
or more precisely,
$v_{C^{\infty},i}:=
 \sum_{j=1}^N\chi_j\, v_{S^{(j)},i}$.
Then,
$\vecv_{C^{\infty}}:=(v_{C^{\infty},i})$
gives a $C^{\infty}$-frame on 
$\pi^{-1}(\nbigv)$.

Let $S\in\Multisector^{\ast}(\nbigx\setminus W,\nbigi)$.
Let $C^{(j)}$ be determined by
$\vecv_{S^{(j)}}=\vecv_S\,(I+C^{(j)})$
on $\Sbar\cap \Sbar^{(j)}$,
where $I$ denotes the identity matrix.
Let $Z(S,S^{(j)})=
 \Sbar^{(j)}\cap \Sbar\cap \pi^{-1}\bigl(\nbigd_z\bigr)$.
Due to Lemma \ref{lem;07.6.15.110},
we have
(i) $C^{(j)}_{|\Zhat(S,S^{(j)})}=0$,
(ii) $C^{(j)}_{\gminia,\gminib}=0$ unless
$\gminia\leq_{ S^{(j)}\cap S}\gminib$,
(iii) $C^{(j)}_{\gminia,\gminib}\,
 \exp\bigl(\varrho^{-1}(\gminia-\gminib)\bigr)
 =
 O\Bigl(\exp\Bigl(C'|\varrho^{-1}|\bigl(
  |\vecz^{\vecm(1)}|
 +\log |z_{i(0)}^{-1}|
 \bigr)\Bigr)\Bigr)$
if $\gminia<_{S^{(j)}\cap S}\gminib$.
By construction, $\vecv_{C^{\infty}}=
 \vecv_S\,\Bigl(I+
 \sum\chi_j\, C^{(j)} \Bigr)$ holds.
Hence, $\vecv_{C^{\infty}}$ satisfies
the desired estimate on $S$.
It also implies that
$\vecv_{C^{\infty}}$
gives a $C^{\infty}$-frame of $E$ on $\nbigv$.
\hfill\qed

\vspace{.1in}

Let us look at the connection form of $\DD$
with respect to $\vecv_{C^{\infty}}$.
Let $I_{\gminia}$ denote the identity matrix
whose size is $\rank \Ehat_{\gminia}$.
Then, we have the following:
\[
 (\DD+\delbar_{\lambda}) 
 \vecv_{C^{\infty}}
=\vecv_{C^{\infty}}
 \left(
 \Bigl(
 \bigoplus_{\gminia\in\nbigi} 
 d\gminia\, I_{\gminia}
\Bigr)+R
 \right)
\]
We put $\nbigd^{(1)}:=\nbigd_z$,
$\nbigd^{(2)}=\nbigd_Y$
if $m_{i(0)}<-1$,
or $\nbigd^{(1)}:=
 \bigcup_{j\neq i(0)}\nbigd_{z,j}$,
 $\nbigd^{(2)}:=\nbigd_{z,i(0)}\cup\nbigd_Y$
if $m_{i(0)}=-1$.
We can deduce the following
from the property of $\vecv_{C^{\infty}}$.

\begin{lem}
\label{lem;07.6.16.100}
$R$ is a $C^{\infty}$-section of
\[
 M_r(\cnum)\otimes\Bigl(
 \vecz^{\vecm(1)}
 \Omega_{\nbigx/\nbigk}^{1,0}(\log \nbigd^{(1)})
+\Omega_{\nbigx/\nbigk}^{1,0}(\log \nbigd^{(2)})
+\Omega_{\nbigx}^{0,1}
 \Bigr),
\]
and we have
$R_{\gminia,\gminia|\What}
=\Rhat_{\gminia}$
and $R_{\gminia,\gminib|\What}=0$
for $\gminia\neq\gminib$.
For each sufficiently small sector $S$,
the following holds:
\begin{enumerate}
\item
$R_{\gminia,\gminib|S}=0$
 unless $\gminia\leq_{S} \gminib$.
\item If $\gminia<_{S}\gminib$,
 the $C^q$-derivatives of
 $R_{\gminia,\gminib}\,
 \exp\bigl(\varrho^{-1}(\gminia-\gminib)\bigr)$ is
\[
 O\Bigl(
 \exp\Bigl(
 C\,|\varrho^{-1}|\bigl(
  |\vecz^{\vecm(1)}|
+\log|z_{i(0)}^{-1}|
 \bigr) \Bigr)
 |\varrho|^{-N(q)}
\prod_{j=1}^k|z_j|^{-N(q)}
\Bigr)
\]
 for some $C>0$ 
 and $N(q)\geq 0$
 $(q\in\seisuu_{\geq 0})$.
\hfill\qed
\end{enumerate}
\end{lem}

Let $\vecv'_{C^{\infty}}$ be another $C^{\infty}$-lift
of $\vecvhat$.
Let $B$ be determined by
$\vecv_{C^{\infty}}'=\vecv_{C^{\infty}}\, (I+B)$.
It is easy to deduce the following.
\begin{lem}
\label{lem;07.6.16.9}
We have $B_{|\What}=0$.
On each sufficiently small sector $S$,
we have  $B_{\gminia,\gminib|S}=0$
unless $\gminia\leq_S\gminib$.
If $\gminia<_S\gminib$,
 the $C^q$-derivatives of 
 $B_{\gminia,\gminib}\,
 \exp\bigl(\varrho^{-1}(\gminia-\gminib)\bigr)$
 is 
\[
 O\Bigl(
 \exp\Bigl(
 C|\varrho^{-1}|\bigl(
 |\vecz^{\vecm(1)}|
+\log|z_{i(0)}^{-1}| 
 \bigr)
 \Bigr)
 |\varrho|^{-N(q)}
\prod_{j=1}^k|z_j|^{-N(q)}
 \Bigr)
\]
 for some $\epsilon>0$ 
 and $N(q)\geq 0$
 $(q\in\seisuu_{\geq 0})$.
\hfill\qed
\end{lem}

\subsection{Approximation of formal decompositions
 (Appendix)}

Let $X:=\Delta^n$,
$D_i:=\{z_i=0\}$ and
$D:=\bigcup_{i=1}^{\ell} D_i$
for some $\ell\leq n$.
In the following argument,
$N$ will denote a large integer,
and we will shrink $X$
around the origin $(0,\ldots,0)$ without mention.
Let $\Dhat$ denote the completion of $X$
along $D$.
\index{completion $\Dhat$}
(See \cite{banica},
 \cite{bingener} and \cite{krasnov}.
 See also a brief review in
 Subsection \ref{subsection;07.11.6.10}.)
Let $\Dhat^{(N)}$ denote 
the $N$-th infinitesimal neighbourhood of $D$ in $X$.
\index{$N$-th infinitesimal 
neighbourhood $\Dhat^{(N)}$}
Let $\iota:\Dhat\lrarr X$ denote 
the natural morphism of ringed spaces.
For any $\nbigo_X$-module $\nbigf$,
let $\nbigf_{|\Dhat}$ denote 
$\iota^{-1}\nbigf
 \otimes_{\iota^{-1}\nbigo_X}
 \nbigo_{\Dhat}$.
\index{restriction $\nbigfshitadhat$}
We use the symbol
$\nbigf_{|\Dhat^{(N)}}$
in a similar meaning.
\index{restriction $\nbigfshitadhatN$}

Let $V$ be a free $\nbigo_X$-module of finite rank,
equipped with the following data:
\begin{itemize}
\item
Sections $f_k^{(i)}\in \End(V_{|D_k})$ 
for $k=\ell+1,\ldots,n$ and $i=1,\ldots,M$.
\item
Filtrations
$\lefttop{k}\nbigf$ of $V_{|D_k}$
for $k=\ell+1,\ldots,n$.
\end{itemize}
Here, $M$ denotes some positive integer.

\subsubsection{Approximation of
an endomorphism}
\label{subsection;10.6.30.1}

Let $\Fhat\in\End(V)_{|\Dhat}$
satisfy the following conditions:
\begin{itemize}
\item
$\bigl[f^{(i)}_{k|\Dhat\cap D_k},
 \Fhat_{|\Dhat\cap D_k}\bigr]=0$
for $k=\ell+1,\ldots,n$
and $i=1,\ldots,M$,
where $[\cdot,\cdot]$ denotes the commutator.
\item
$\Fhat_{|\Dhat\cap D_k}$ preserves
$\lefttop{k}\nbigf$
for $k=\ell+1,\ldots,n$.
\end{itemize}
Here, $\Dhat\cap D_k$ means
the completion of $D_k$
along $D\cap D_k$.

\begin{lem}
\label{lem;07.6.15.2}
For any large $N$,
we can take a section $F^{(N)}\in \End(V)$
such that 
$F^{(N)}_{|\Dhat^{(N)}}
=\Fhat_{|\Dhat^{(N)}}$,
with the following property:
\begin{itemize}
\item
$\bigl[f^{(i)}_k,F_{|D_k}\bigr]=0$
for $k=\ell+1,\ldots,n$ and $i=1,\ldots,M$.
\item
$F^{(N)}_{|D_k}$ preserves
the filtrations $\lefttop{k}\nbigf$
for $k=\ell+1,\ldots,n$.
\end{itemize}
\end{lem}
\pf
Let $\Cok\bigl(\ad(f^{(i)}_k)\bigr)$ denote 
the cokernel of 
the morphism 
\[
 \ad(f^{(i)}_k):\End(V)_{|D_k}\lrarr 
 \End(V)_{|D_k}
\]
given by $\ad(f^{(i)}_k)(g)=[f^{(i)}_k,g]$.
We put $\fhat_k^{(i)}:=f^{(i)}_{k|D_k\cap \Dhat}$.
Let $\Cok\bigl(\ad(\fhat^{(i)}_k)\bigr)$
denote the cokernel of the morphism
$\ad(\fhat^{(i)}_k):
 \End(V)_{|\Dhat\cap D_k}
 \lrarr
 \End(V)_{|\Dhat\cap D_k}$.

For $k=\ell+1,\ldots,n$,
let $\End'(V_{|D_k})$ denote the sheaf of
sections of $\End(V_{|D_k})$
preserving $\lefttop{k}\nbigf$.
Similarly, let $\End'(\Vhat_{|\Dhat\cap D_k})$
denote the sheaf of
sections of $\End(\Vhat_{|\Dhat\cap D_k})$
preserving $\lefttop{k}\nbigf$.
We put 
\[
 A_k:=\End(V_{|D_k})\big/\End'(V_{|D_k}),
\quad
\Ahat_k:=
 \End(\Vhat_{|D_k\cap\Dhat})\big/
 \End'(\Vhat_{|D_k\cap \Dhat}).
\]
Let $\iota_k:D_k\lrarr X$ and 
$\iotahat_k:D_k\cap \Dhat\lrarr \Dhat$
denote the inclusions.
We set
\[
 \nbigf:=\Ker\left(
\End(V)\lrarr
 \bigoplus_{k=\ell+1}^n
 \iota_{k\ast}A_k
\oplus
 \bigoplus_{k=\ell+1}^n
 \bigoplus_{i=1}^M
 \iota_{k\ast}
 \Cok\bigl(\ad(f_k^{(i)})\bigr)
 \right)
\]
\[
  \nbigfhat:=\Ker\left(
 \End(V)_{|\Dhat}
\lrarr
 \bigoplus_{k=\ell+1}^n
 \iota_{k\ast}\Ahat_k
\oplus
 \bigoplus_{k=\ell+1}^{n}
 \bigoplus_{i=1}^M
 \iotahat_{k\ast}
 \Cok\bigl(\ad(\fhat_k^{(i)})\bigr)
 \right)
\]
We have only to show that
$\nbigfhat$ is the completion of
$\nbigf$ along $D$,
which implies the claim of the lemma.
Note that $\Cok\bigl(\ad(\fhat^{(i)}_k)\bigr)$
and $\Ahat_k$ are the completions of
$\Cok\bigl(\ad(f^{(i)}_k)\bigr)$
and $A_k$, respectively.
Since $\nbigo_{\Dhat}$ is faithfully flat over $\nbigo_X$
(\cite{bingener}),
we obtain
$\nbigfhat\simeq \nbigf_{|\Dhat}$.
Thus Lemma \ref{lem;07.6.15.2}
is finished.
\hfill\qed

\subsubsection{Approximation of
a decomposition}

Assume that we are given 
a formal decomposition
$V_{|\Dhat}=
 \bigoplus_{\gminia\in\nbigi}\Vhat_{\gminib}$
satisfying the following condition:
\begin{itemize}
\item
 $\Vhat_{\gminia|\Dhat\cap D_k}$
 $(\gminia\in\nbigi)$
 are preserved by
 $f^{(i)}_{k|\Dhat\cap D_k}$
 for $k=\ell+1,\ldots,n$
 and $i=1,\ldots,M$.
\item
 $\lefttop{k}\nbigf=
 \bigoplus_{\gminia\in\nbigi}
 \lefttop{k}\nbigf\cap
 V_{\gminia|\Dhat\cap D_k}$.
\end{itemize}

\begin{lem}
\label{lem;10.6.30.3}
For any large $N$,
we can take a decomposition
$V=\bigoplus_{\gminia\in\nbigi}V_{\gminia}^{(N)}$
such that
$V^{(N)}_{\gminia|\Dhat^{(N)}}
=\Vhat_{\gminia|\Dhat^{(N)}}$,
with the following property:
\begin{itemize}
\item
 $V^{(N)}_{\gminia|D_k}$
 $(\gminia\in\nbigi)$
 are preserved by
 $f^{(i)}_{k}$
 for $k=\ell+1,\ldots,n$
 and $i=1,\ldots,M$.
\item
 $\lefttop{k}\nbigf=
 \bigoplus_{\gminia\in\nbigi}
 \lefttop{k}\nbigf\cap
 V_{\gminia|D_k}$.
\end{itemize}
\end{lem}
\pf
Let $\pihat_{\gminia}$ $(\gminia\in\nbigi)$
be the projection of $V_{|\Dhat}$
onto $\Vhat_{\gminia}$.
We take an injection
$\psi:\nbigi\lrarr\seisuu$,
and we put
$\Fhat:=
 \sum_{\gminia\in\nbigi}
 \psi(\gminia)\cdot\pihat_{\gminia}$.
We take $F^{(N)}$ for $\Fhat$
as in Lemma \ref{lem;07.6.15.2}.
After shrinking $X$,
we have the decomposition
$V=\bigoplus V_{\gminia,N}$
such that
(i) $F^{(N)}(V_{\gminia,N})\subset
 V_{\gminia,N}$,
(ii) the eigenvalues of
$F^{(N)}_{|V_{\gminia,N}}$
are close to $\psi(\gminia)$.
Then, it gives the desired decomposition.
\hfill\qed

\section{Proof of the statements
in Section \ref{section;10.5.10.5}}
\label{section;10.5.18.5}
\subsection{Preliminary}
\label{subsection;10.5.12.1}

We use the setting and the notation
in Subsection \ref{subsection;10.5.3.7}.
Let $\nbigj\subset M(\nbigx,\nbigd)/H(\nbigx)$
be a good set of irregular values.
We assume that the coordinate system
$(z_1,\ldots,z_n)$ is admissible
for $\nbigj$.
We take an auxiliary sequence
$\vecm(0),\ldots,\vecm(L),\vecm(L+1)=\veczero$
(Section \ref{subsection;07.11.6.1}).
Let $k(p)$ be determined by
$\vecm(p)\in\seisuu_{<0}^{k(p)}\times\veczero$.
Let $\nbigj(\vecm(p))$
and $\nbigi_{\gminic}^{\vecm(p)}$
be as in Subsection \ref{subsection;07.12.26.1}.
Let $\etabar_{\vecm(p)}:\nbigj\lrarr\nbigj(\vecm(p))$
be the induced map.
For each $\gminic\in\nbigj(\vecm(p))$,
we put $\nbigj_{\vecc}:=
 \etabar_{\vecm(p)}^{-1}(\gminic)$,
which is a good set of irregular values.
We also have the naturally induced map
$\etabar_{\vecm(p-1),\vecm(p)}:
 \nbigj(\vecm(p))\lrarr\nbigj(\vecm(p-1))$.

Let $W(\kbar):=
 \nbigx^0\cup\nbigd(\kbar)$.
Let $\pi_k:\nbigxtilde(W(\kbar))
 \lrarr \nbigx$ be the real blow up
of $\nbigx$ at $W(\kbar)$.
In particular, $\pi_{\ell}=:\pi$.
We have the naturally induced maps
$\varpi_{k,m}:\nbigxtilde(W(\mbar))
 \lrarr \nbigxtilde(W(\kbar))$
for $m\geq k$.
The map $\nbigxtilde(W)\lrarr
 \nbigxtilde(W(\kbar))$
is denoted by $\varpi_k$.
We will use the following obvious lemma
implicitly.
\begin{lem}
Let $P\in\nbigd_{\ellsitabar}$
and $Q\in\pi^{-1}(P)$.
For $\gminia,\gminib\in \nbigj$,
we have $\gminia\leq_Q^{\varrho}\gminib$,
if and only if 
$\etabar_{\vecm(p)}(\gminia)
 \leq_{\varpi_{\kbar(p)}(Q)}^{\varrho}
 \etabar_{\vecm(p)}(\gminib)$,
where $\vecm(p)=\ord(\gminia-\gminib)$.
\hfill\qed
\end{lem}

For any $P\in\nbigd$,
let $\nbigj_{P}$ denote the image of
$\nbigj$ by
$M(\nbigx,\nbigd)/H(\nbigx)\lrarr
 \nbigo_{\nbigx}(\ast \nbigd)_{P}/
 \nbigo_{\nbigx,P}$.
For any $i$,
we put
$\nbigd_i^{\ast}:=
 \nbigd_i\setminus
 \Bigl(
 \nbigx^0\cup
 \bigcup_{j\neq i}
 \nbigd_j
 \Bigr)$.
Note that 
the natural map
$\nbigj\lrarr\nbigj_{P}$ is bijective
for any $P\in \nbigd_1^{\ast}$.

\begin{lem}
\label{lem;10.5.12.21}
Let $Q\in\pi^{-1}(\nbigd_{\ellsitabar})$.
We take a small
$\nbigu_Q\in\gbigu(Q,\nbigj)$
and $\gminia,\gminib\in\nbigj$.
Then, we have
$\gminia\leq_{Q}^{\varrho}\gminib$,
if and only if
$\gminia\leq_{Q'}^{\varrho}
 \gminib$
for any $Q'\in
 \nbigu_Q\cap\pi^{-1}(\nbigd_1^{\ast})$.
More strongly, we have
$\gminia\leq_{Q}^{\varrho}\gminib$,
if and only if
there is a dense subset $B$ of
$\nbigu_Q\cap\pi^{-1}(\nbigd_1^{\ast})$
such that 
$\gminia\leq_{Q'}^{\varrho}
 \gminib$ for any $Q'\in B$.
\end{lem}
\pf
We have 
$\gminia<_{Q'}^{\varrho}\gminib$
if and only if
$F_{\gminia,\gminib}(Q')<0$.
For fixed $\gminia,\gminib$,
after an appropriate coordinate change,
we may assume 
$\gminia-\gminib=\vecz^{\vecm}$.
Then, the claim of the lemma is clear.
\hfill\qed

\vspace{.1in}
We state it in a slightly generalized form.
\begin{lem}
\label{lem;10.5.16.3}
Let $I\subset \ellsitabar$.
\begin{itemize}
\item
Take $P\in
 \nbigd_{I}\setminus 
 \bigcup_{j\not\in I}\nbigd_j$.
We have the naturally induced 
bijective map
$\nbigj_P\lrarr\nbigj_{P'}$
for any $P'\in \nbigd_{\min I}^{\ast}$.
\item
Let $Q\in\pi^{-1}(P)$.
We take a small $\nbigu_Q\in\gbigu(Q,\nbigj_P)$
and $\gminia,\gminib\in\nbigj_P$.
Then, we have
$\gminia\leq_Q^{\varrho}\gminib$
if and only if 
$\gminia\leq_{Q'}^{\varrho}\gminib$
for any $Q'\in\nbigu_P\cap
 \pi^{-1}(\nbigd_{\min I}^{\ast})$.
More strongly,
we have
$\gminia\leq_Q^{\varrho}\gminib$
if and only if 
there exists a dense subset
$B\subset \nbigu_P\cap
 \pi^{-1}(\nbigd_{\min I}^{\ast})$
such that
$\gminia\leq_{Q'}^{\varrho}\gminib$
for any $Q'\in B$.
\hfill\qed
\end{itemize}
\end{lem}

\subsection{Reduction}
\label{subsection;10.6.28.1}

Let $(\nbige,\DD)$ be an unramifiedly
good $\varrho$-meromorphic flat bundle
on $(\nbigx,\nbigd)$
with a good lattice $E$
and a good set of irregular values $\nbigj$.
We assume that the coordinate system is admissible
for $\nbigj$.
We use the notation in Subsection 
\ref{subsection;10.5.12.1}.
We shall construct the associated
graded meromorphic $\varrho$-flat bundle
$\Gr^{\vecm(p)}(\nbige,\DD)$
with an unramifiedly good lattice
$\Gr^{\vecm(p)}(E)$
for any $p$,
defined on a neighbourhood of 
$\nbigd_{\kbar(0)}$.
We remark $\nbigd_{\kbar(0)}\subset
\nbigd_{\kbar(p)}$,
which we will implicitly use.

\subsubsection{One step reduction}

Let us consider the case in which
$\nbigj(\vecm(p-1))$
consists of a unique element
$\gminia$.
Then, $(E,\DD)$ is a weakly good lattice
in the level $(\vecm(p),\gminih(p))$.
By the procedure in 
Subsection \ref{subsection;07.12.12.1},
after shrinking $\nbigx$ around $\nbigd_{\kbar(p)}$,
we obtain a graded holomorphic bundle
$\Gr^{\vecm(p)}(E)
=\bigoplus_{\gminia\in\nbigj(\vecm(p))}
 \Gr^{\vecm(p)}_{\gminia}(E)$
with a meromorphic flat $\varrho$-connection
$\DD^{\vecm(p)}
=\bigoplus \DD^{\vecm(p)}_{\gminia}$
on $(\nbigx,\nbigd)$.
Due to (\ref{eq;08.9.5.10}),
the completion of
$\Gr^{\vecm(p)}_{\gminia}(E,\DD):=
 (\Gr^{\vecm(p)}_{\gminia}(E),
\DD^{\vecm(p)}_{\gminia})$
along $W(\kbar(p))$ is
naturally isomorphic to
$\bigl(\Ehat^{\vecm(p)}_{\gminia,\What(\kbar(p))},
 \DD_{\gminia}\bigr)$ in (\ref{eq;10.5.5.41}):
\[
 (\Gr^{\vecm(p)}_{\gminia}(E),
\DD^{\vecm(p)}_{\gminia})_{|\What(\kbar(p))}
\simeq
\bigl(\Ehat^{\vecm(p)}_{\gminia,\What(\kbar(p))},
 \DD_{\gminia}\bigr)
\]
In particular,
 $\bigl(
 \Gr^{\vecm(p)}_{\gminia}(E),
 \DD^{\vecm(p)}_{\gminia}\bigr)$
is also an unramifiedly good lattice.
We have
$\Irr(\DD^{\vecm(p)}_{\gminia})
=\nbigj_{\gminia}$.
In particular,
its image by $\etabar_{\vecm(p)}$
consists of one element.

\subsubsection{Reduction in the level $\vecm(p)$}

By shrinking $\nbigx$ around $\nbigd(\kbar(0))$,
we shall inductively construct
the unramifiedly good lattices
\[
 \Gr^{\vecm(p)}_{\gminib}(E,\DD)
=\bigl(\Gr^{\vecm(p)}_{\gminib}(E),
 \DD^{\vecm(p)}_{\gminib}
 \bigr)
\]
on $(\nbigx,\nbigd)$
for $\gminib\in \nbigj(\vecm(j))$
$(j=0,\ldots,L)$ with the following property:
\begin{itemize}
\item
The completion of
$\Gr^{\vecm(p)}_{\gminib}(E,\DD)$
along $W(\kbar(p))$
are naturally isomorphic to
$(\Ehat^{\vecm(p)}_{\gminib},
 \DDhat_{\gminib})$
in (\ref{eq;10.5.5.41}):
\index{Gr $\Gr^{\vecm(p)}_{\gminib}(E)$}
\begin{equation}
\label{eq;08.11.1.2}
\Gr^{\vecm(p)}_{\gminib}(E,\DD)_{|\What(\kbar(p))}
\simeq
(\Ehat^{\vecm(p)}_{\gminib,\What(\kbar(p))},
 \DDhat_{\gminib})
\end{equation}
\end{itemize}
Namely, for any $\gminib\in \nbigj(\vecm(p))$,
we put $\gminia:=
\etabar_{\vecm(p-1),\vecm(p)}(\gminib)
 \in\nbigj(\vecm(p-1))$,
and we define
\[
 \Gr^{\vecm(p)}_{\gminib}(E,\DD):=
 \Gr^{\vecm(p)}_{\gminib}
 \Gr^{\vecm(p-1)}_{\gminia}(E,\DD).
\]
By (\ref{eq;08.11.1.2}),
the following holds:
\begin{itemize}
\item
We have 
$\Irr(\DD^{\vecm(p)}_{\gminib})=
 \etabar_{\vecm(p)}^{-1}(\gminib)$
for any $\gminib\in\nbigj(\vecm(p))$.
In particular,
its image by $\etabar_{\vecm(p)}$
consists of one element.
\item
We have the following natural isomorphism
for each $p$:
\[
  (E,\DD)_{|\What(\kbar(p))}
\simeq
 \bigoplus_{\gminia\in \Irrbar(\DD,\vecm(p))}
 \bigl(\Gr^{\vecm(p)}_{\gminia}(E),
 \DD^{\vecm(p)}_{\gminia}
 \bigr)_{|\What(\kbar(p))}
\]
\end{itemize}
In particular,
$\Gr^{\vecm(L)}_{\gminia}(E,\DD)$
are $\gminia$-logarithmic.

\begin{rem}
\label{rem;07.12.25.10}
In the following,
we often formally
put $\Irr(\DD,\vecm(-1)):=\{0\}$,
$\Gr_0^{\vecm(-1)}(E)=E$ 
and $\DD^{\vecm(-1)}_0:=\DD$.
\index{set $\Irr(\DD,\vecm(-1)):=\{0\}$}
\index{Gr $\Gr_0^{\vecm(-1)}(E)$}
We also often use the symbol
$\Gr^{\nbigftilde}_{\gminia}(E)$
instead of $\Gr^{\vecm(L)}_{\gminia}(E)$,
which is called the full reduction of $(E,\DD)$.
\index{Gr $\Gr^{\full}_{\gminia}(E)$}
\hfill\qed
\end{rem}

\subsubsection{Functoriality of the associated graded bundle}
\label{subsection;07.12.13.100}

Let $(E_r,\DD_r)$ $(r=1,2)$ be unramifiedly good
such that $\Irr(\DD_r)=\nbigj_r$.
By using Lemma \ref{lem;07.12.7.17} inductively,
we obtain the following lemma.
\begin{lem}
\label{lem;07.12.13.111}
Let $F:(E_1,\DD_1)\lrarr(E_2,\DD_2)$
be a morphism.
Assume $\nbigj_1\cup\nbigj_2$ is also good.
We have the naturally induced flat morphisms
$\Gr^{\vecm(p)}_{\gminia}(F):
 \Gr^{\vecm(p)}_{\gminia}(E_1)
\lrarr
 \Gr^{\vecm(p)}_{\gminia}(E_2)$
for any $\gminia$.
\hfill\qed
\end{lem}

\begin{cor}
\label{cor;07.12.7.18}
If the restriction of $F$ to $\nbigx-\nbigd$
is an isomorphism,
we obtain naturally induced isomorphisms
for any $\gminia\in\Irrbar(\DD,\vecm(p))$:
\[
 \Gr^{\vecm(p)}_{\gminia}(F):
 \Gr^{\vecm(p)}_{\gminia}(E_1)(\ast \nbigd)
\lrarr 
 \Gr^{\vecm(p)}_{\gminia}(E_2)(\ast \nbigd)
\]
In particular,
the graded meromorphic $\varrho$-flat bundle
$\Gr^{\vecm(p)}(E)(\ast \nbigd)$
is well defined for
$\bigl(E(\ast \nbigd),\DD\bigr)$,
in the sense that they
are independent of the choice of
an unramifiedly good lattice $E$.
\hfill\qed
\end{cor}

If $\nbigj_1\otimes\nbigj_2$ is good,
we obtain the following natural isomorphism
for any $\gminia\in\Irrbar(\DDtilde,\vecm(p))$
by using Lemma \ref{lem;07.12.7.25} inductively:
\[
 \Gr^{\vecm(p)}_{\gminia}
 \bigl(E_1\otimes E_2,
 \DDtilde\bigr)
\simeq
 \bigoplus_{
 \substack{\gminia_p\in\Irrbar(\DD_p,\vecm(p))
 \\
 \gminia_1+\gminia_2=\gminia}}
 \Gr^{\vecm(p)}_{\gminia_1}
 \bigl(E_1,\DD_1\bigr)
\otimes
 \Gr^{\vecm(p)}_{\gminia_2}
  \bigl(E_2,\DD_2\bigr)
\]
If $\nbigj_1\oplus\nbigj_2$ is good,
we have
\[
 \Gr_{\gminia}^{\vecm(p)}(E_1\oplus E_2)
\simeq
 \Gr_{\gminia}^{\vecm(p)}(E_1)
\oplus
 \Gr_{\gminia}^{\vecm(p)}(E_2).
\]

If $(E,\DD)$ is unramifiedly good,
the dual $(E^{\lor},\DD^{\lor})$ is also
unramifiedly good.
By using Lemma \ref{lem;07.12.7.5} inductively,
we obtain the following natural isomorphism
for any $\gminia\in\nbigj$:
\[
 \Gr^{\vecm(p)}_{-\gminia}(E^{\lor},\DD^{\lor})
\simeq
 \Gr^{\vecm(p)}_{\gminia}(E,\DD)^{\lor}
\]
\subsection{Full and partial Stokes filtrations}
\label{subsection;07.11.18.10}

Let $(\nbige,\DD)$ be an unramifiedly good
meromorphic $\varrho$-flat bundle
on $(\nbigx,\nbigd)$, with a good lattice $E$
and a good set of irregular values $\nbigj$.
We shall explain the construction
of full and partial Stokes filtrations
of the stalks $\nbige_{|Q}$ for
$Q\in\pi^{-1}(\nbigd)$.

As explained in 
Subsection \ref{subsection;10.6.28.1},
after shrinking $\nbigx$ around $\nbigd(\kbar(0))$,
we may have the graded meromorphic
$\varrho$-flat bundle
$\Gr^{\vecm(p)}(\nbige,\nabla)$
with the unramifiedly good lattice
$\Gr^{\vecm(p)}(E)$
on $(\nbigx,\nbigd)$.
For $k(0)\leq k\leq \ell$,
let us consider the real blow up
$\pi_k:\nbigxtilde(W(\kbar))\lrarr\nbigx$.
Let $Q$ be any point of $\pi^{-1}(\nbigd_{\kbar})$.
The image of $Q$ by
$\varpi_{k(p),k}:
 \nbigxtilde(W(\kbar))\lrarr
 \nbigxtilde(W(\kbar(p)))$
is denoted by $Q_p$.
We have a small neighbourhood
$\nbigu_{Q_p}$ of $Q_p$
in $\nbigxtilde(W(\kbar(p)))$
and the Stokes filtration $\nbigf^{Q_p}$
of $\Gr_{\gminia}^{\vecm(p-1)}(E)_{|\nbigu_{Q_p}}$
indexed by 
$\nbigj(\vecm(p),\gminia)$
with $\leq_Q$,
where $\gminia\in\nbigj(\vecm(p-1))$
and
 \[
 \nbigj(\vecm(p),\gminia):=
 \bigl\{
 \gminib\in\nbigj(\vecm(p)),
 \etabar_{\vecm(p-1),\vecm(p)}(\gminib)=\gminia
 \bigr\}.
 \]
We take a small neighbourhood $\nbigu_Q$
of $Q$ in $\nbigxtilde(W(\kbar))$
such that 
$\varpi_{k(p),k}(\nbigu_Q)\subset
 \nbigu_{Q_p}$.
We obtain the filtered bundle
$\bigl(
 \Gr^{\vecm(p-1)}_{\gminia}(E)_{\nbigu_Q},
 \nbigf^{Q_p}
 \bigr)$ for each $\gminia\in\nbigj(\vecm(p-1))$,
and the associated graded bundle
is naturally isomorphic to
\[
 \bigoplus_{\gminib\in\nbigj(\vecm(p),\gminia)}
 \Gr^{\vecm(p)}_{\gminib}(E)_{|\nbigu_Q},
\]
By applying the inductive procedure in 
Section \ref{subsection;07.12.12.2},
we obtain the $\DD$-flat filtration
$\nbigf^{Q\,\vecm(p)}$ of $E_{|\nbigu_Q}$
indexed by the ordered set
$\bigl(\nbigj(\vecm(p)), \leq^{\varrho}_Q\bigr)$.
It is called the partial Stokes filtration of 
$E_{|\nbigu_Q}$ in the level $\vecm(p)$.
\index{partial Stokes filtration $\nbigf^{\vecm(p)\,Q}$}
In particular,
$\nbigf^{\vecm(L)\,Q}$ is called 
the full Stokes filtration,
and denoted by $\nbigftilde^Q$.
\index{full Stokes filtration $\nbigftilde^Q$}
We have the induced filtrations
of the germs of $E$ and $\nbige$ at $Q$,
which are denoted by the same symbols.

The following compatibility is clear by construction.
\begin{lem}
\label{lem;08.2.3.10}
Let $Q\in \pi^{-1}(\nbigd_{\kbar})$.
We take neighbourhoods $\nbigu_Q$
as above.
Let $Q_1\in\pi^{-1}(\nbigd_{\kbar})\cap\nbigu_Q$.
We take $\nbigu_{Q_1}\subset\nbigu_Q$.
Then, 
the filtrations
$\nbigf^{Q\,\vecm(p)}$
and 
$\nbigf^{Q_1\,\vecm(p)}$
of $E_{|\nbigu_{Q_1}}$
are compatible over
$\bigl(\nbigj(\vecm(p)),\leq_{Q}^{\varrho}
 \bigr)
\lrarr
 \bigl(\nbigj(\vecm(p)),\leq^{\varrho}_{Q_1}\bigr)$.
\hfill\qed
\end{lem}

\subsubsection{Functoriality}
\label{subsection;10.5.12.201}

Let $(E_r,\DD_r)$ $(r=1,2)$ be unramifiedly good
such that $\Irr(\DD_r)=\nbigj_r$.
By using Proposition \ref{prop;10.5.12.10},
we obtain the following lemma.
\begin{lem}
Let $F:(E_1,\DD_1)\lrarr(E_2,\DD_2)$
be a morphism.
Assume $\nbigj_1\cup\nbigj_2$ is also good.
Let $Q\in\pi^{-1}(\nbigd_{\kbar})$,
where $k(0)\leq k\leq \ell$.
The morphisms of germs
$F_Q:E_{1,Q}\lrarr E_{2,Q}$
preserve the Stokes filtrations
in the level $\vecm(p)$.
In particular,
the filtrations $\nbigf^{Q\,\vecm(p)}$
of the germ $\nbige_Q$
are well defined for $(\nbige,\DD)$
in the sense that they
are independent of the choice of
an unramifiedly good lattice $E$.
\hfill\qed
\end{lem}

By using the propositions in
Subsection \ref{subsection;10.5.11.22}
inductively,
we also obtain the functoriality
of the full and partial Stokes filtrations
for dual, tensor product
and direct sum as in 
Subsection \ref{subsection;10.5.12.15}.
We have similar functoriality for
the partial Stokes filtrations.

\subsection{Compatibility and characterization}

By applying the procedure explained in 
Subsection \ref{subsection;07.11.18.10}
to the restriction of $(\nbige,\DD)$
to a small neighbourhood of any point of $\nbigd$,
we obtain the full Stokes filtration
of the stalk of $E$ at any point of
$\pi^{-1}(\nbigd)$.
We shall argue the comparison of the filtrations.

\subsubsection{Preliminary}

We consider the case $\nbigd=\nbigd_1$.
For simplicity,
we assume that
$\varrho$ is nowhere vanishing.
Let $(\nbige,\nabla)$ be an unramifiedly
good meromorphic $\varrho$-flat bundle
on $(\nbigx,\nbigd)$
with a good lattice $E$ and 
a good set of irregular values $\nbigj$.
We can use the order of the poles of
elements of $\nbigj$
as the auxiliary sequence.
For $Q\in \pi^{-1}(\nbigd)$,
we have the partial Stokes filtrations
$\nbigf^{(m)}$ of
$\nbige_{Q}$ and $E_{Q}$.

\begin{lem}
\label{lem;10.5.12.20}
Let $\overline{\nbigf}^{(m)}$
be a $\DD$-flat filtration of $E_{Q}$
indexed by
$(\nbigj(m),\leq_{Q}^{\varrho})$
such that
$\overline{\nbigf}^{(m)}_{|\Qhat}
=\nbigf^{(m)}_{|\Qhat}$.
Then, we obtain
$\overline{\nbigf}^{(m)}
=\nbigf^{(m)}$.
\end{lem}
\pf
If we take a sufficiently small neighbourhood
$\nbigu_Q$ of $Q$ in $\nbigxtilde(\nbigd)$,
we obtain 
$\overline{\nbigf}^{(m)}_{|
 \widehat{\pi^{-1}(\nbigd)}\cap\nbigu_Q}
=\nbigf^{(m)}_{|
 \widehat{\pi^{-1}(\nbigd)}\cap\nbigu_Q}$.
Then, 
we obtain 
$\overline{\nbigf}^{(m)}
=\nbigf^{(m)}$
by using the argument in the proof of
Lemma \ref{lem;07.9.29.20}.
\hfill\qed

\subsubsection{}

Let us return to the original setting.
Let $Q\in \pi^{-1}(\nbigd_{\ellsitabar})$.
We take a small
$\nbigu_Q\in\gbigu(Q,\nbigj)$.
We set $\nbigd_i^{\ast}:=
 \nbigd_i\setminus\bigl(
 W\cup\bigcup_{j\neq i}\nbigd_j
 \bigr)$.
\begin{lem}
\label{lem;10.5.12.30}
Take any
$Q'\in \nbigu_Q\cap\pi^{-1}(\nbigd_i^{\ast})$
for some $1\leq i\leq \ell$.
Then, the filtrations $\nbigftilde^Q$
and $\nbigftilde^{Q'}$
are compatible over
$\bigl(
 \nbigj,\leq_{Q}^{\varrho}
 \bigr)
\lrarr
 \bigl(
 \nbigj_{Q'},\leq_{Q'}^{\varrho}
 \bigr)$.
\end{lem}
\pf
We construct a filtration
$\overline{\nbigf}$
of $E_{|\nbigu_Q}$
from $\nbigftilde^Q$
and $\bigl(\nbigj,\leq_{Q}^{\varrho}\bigr)
\lrarr
 \bigl(\nbigj_{Q'},\leq_{Q'}^{\varrho}\bigr)$.
By construction of
$\nbigftilde^Q$,
we can easily check that
$\overline{\nbigf}_{|\Qhat'}
=\nbigftilde^{Q'}$.
Then, the claim follows from
Lemma \ref{lem;10.5.12.20}.
\hfill\qed

\vspace{.1in}

Let us observe that
$\nbigftilde^Q$ can be reconstructed
from the filtrations
$\nbigftilde^{Q'}$
$(Q'\in\nbigu_Q\cap\pi^{-1}(\nbigd_1^{\ast}))$
in the following sense.
\begin{lem}
\label{lem;10.5.12.31}
Let $\overline{\nbigf}$ be a filtration
of $E_{|\nbigu_Q}$ such that
$\overline{\nbigf}_{|\Qhat'}$
and $\nbigftilde^{Q'}_{|\Qhat'}$
are compatible 
over $\bigl(\nbigj,\leq_Q^{\varrho}\bigr)
 \lrarr
 \bigl(\nbigj_{\pi(Q')},\leq_{Q'}^{\varrho}\bigr)$
for any
$Q'\in \nbigu_Q\cap\pi^{-1}(\nbigd_1^{\ast})$.
Then, we have
$\nbigftilde^Q=\overline{\nbigf}$.
\end{lem}
\pf
It follows from Lemma \ref{lem;10.5.12.22}
and Lemma \ref{lem;10.5.12.21}.
\hfill\qed

\vspace{.1in}
We state it in a slightly generalized form.
\begin{lem}
\label{lem;10.5.12.33}
Let $Q\in \pi^{-1}(\nbigd)$.
Let $i(Q):=\min\bigl\{
 i\,\big|\,\pi(Q)\in\nbigd_i
 \bigr\}$.
Take a sufficiently small neighbourhood
$\nbigu_Q$ of $Q$.
Let $\overline{\nbigf}$ be a filtration
of $E_{|\nbigu_Q}$ indexed by
$\bigl(\nbigj_{\pi(Q)},\leq_Q^{\varrho}\bigr)$
with the following property:
\begin{itemize}
\item
 For any $Q'\in
 \nbigu_Q\cap \pi^{-1}(\nbigd_{i(Q)}^{\ast})$,
 the filtrations
 $\overline{\nbigf}_{|\Qhat'}$
 and $\nbigftilde^{Q'}_{|\Qhat'}$ are 
 compatible over
 $\bigl(\nbigj_{\pi(Q)},\leq_Q^{\varrho}\bigr)
 \lrarr
 \bigl(\nbigj_{\pi(Q')},\leq_{Q'}^{\varrho}\bigr)$.
\end{itemize}
Then, $\overline{\nbigf}=\nbigftilde^{Q}$.
\hfill\qed
\end{lem}

\begin{lem}
\label{lem;10.5.12.41}
Let $Q\in \pi^{-1}(\nbigd)$.
We take a sufficiently small neighbourhood
$\nbigu_Q$ of $Q$ in $\nbigxtilde(W)$.
Then, for any $Q'\in\nbigu_Q\cap\pi^{-1}(\nbigd)$,
the filtrations
$\nbigftilde^{Q}$ and $\nbigftilde^{Q'}$
are compatible over
$\bigl(\nbigj_{\pi(Q)},\leq_Q^{\varrho}\bigr)
\lrarr
 \bigl(\nbigj_{\pi(Q')},\leq_{Q'}^{\varrho}\bigr)$.
\end{lem}
\pf
We construct  $\overline{\nbigf}$
from $\nbigftilde^{Q}$
by $\bigl(\nbigj_{\pi(Q)},\leq_Q^{\varrho}\bigr)
\lrarr
 \bigl(\nbigj_{\pi(Q')},\leq_{Q'}^{\varrho}\bigr)$.
By Lemma \ref{lem;10.5.12.30},
$\overline{\nbigf}$ satisfies the condition
in Lemma \ref{lem;10.5.12.33}
for $\nbigftilde^{Q'}$.
Hence, we obtain
$\overline{\nbigf}=\nbigftilde^{Q'}$.
\hfill\qed

\vspace{.1in}
Let us compare the Stokes filtrations
of $(E,\DD)$ and $\Gr^{\vecm(p)}(E,\DD)$.
Let $Q\in\pi^{-1}(\nbigd_{\ellsitabar})$.
Let $\nbigu_Q$ be a small neighbourhood 
of $Q$ in $\nbigxtilde(W)$.
We have the full Stokes filtration
$\nbigftilde^Q$ 
and the partial Stokes filtration
$\nbigf^{Q\,\vecm(p)}$ of $E_{|\nbigu_Q}$.
By construction,
we have a natural isomorphism
\begin{equation}
\label{eq;10.5.12.40}
 \Gr^{\nbigf^{Q\,\vecm(p)}}(E_{|\nbigu_Q})
\simeq
 \Gr^{\vecm(p)}(E)_{|\nbigu_Q}.
\end{equation}
Let $R\in\nbigu_Q$.
Take a small neighbourhood 
$\nbigu_R\subset\nbigu_Q$.
We have the full Stokes filtration
$\nbigftilde^{R}(E_{\nbigu_{Q'}})$
which induces a filtration 
on the left hand side of (\ref{eq;10.5.12.40}).
We also have the full Stokes filtration
$\nbigftilde^R(\Gr^{\vecm(p)}(E)_{|\nbigu_{Q'}})$,
i.e., the right hand side of (\ref{eq;10.5.12.40}).

\begin{cor}
\label{cor;10.5.19.20}
{\rm(\ref{eq;10.5.12.40})}
is an isomorphism of filtered bundles.
\end{cor}
\pf
By Lemma \ref{lem;10.5.12.41},
the both filtrations are induced
by the full Stokes filtration $\nbigftilde$
of $E_{|\nbigu_Q}$.
\hfill\qed

\begin{cor}
Let $P\in\nbigd_i$ be a smooth point of $\nbigd$.
Assume that the $i$-th component of
$\vecm(p)$ is negative.
Then, we have a natural isomorphism
$\Gr^{(m)}(E_{P})\simeq
 \Gr^{(m)}\Gr^{\vecm(p)}(E_P)$
for any $m<0$.
\hfill\qed
\end{cor}

\subsubsection{}

Let us show a refinement of 
Lemma \ref{lem;10.5.12.20}
in the normal crossing case.
\begin{lem}
\label{lem;10.5.12.200}
Let $Q\in\pi^{-1}(\nbigd_{\ellsitabar})$
such that $\varrho(Q)\neq 0$.
Let $\overline{\nbigf}$ be
a $\DD$-flat filtration of $E_{Q}$
indexed by
$\bigl(\nbigj,\leq_{Q}^{\varrho}\bigr)$
such that
$\overline{\nbigf}_{|\Qhat}
=\nbigftilde^{Q}_{|\Qhat}$.
Then, we have 
$\overline{\nbigf}=\nbigftilde^{Q}$.
\end{lem}
\pf
We take a small $\nbigu_Q\in\gbigu(Q,\nbigj)$
on which $E_{|\nbigu_Q}$ has
the full Stokes filtration $\nbigftilde^{Q}$
and the filtration $\overline{\nbigf}$.
We can take a linear map
$\varphi:\Delta\lrarr \Delta^n$ such that
(i) the image of the induced map
$\varphi_{\nbigk}:\Delta\times\nbigk\lrarr\nbigx$
is not contained in $\nbigd$,
(ii) $Q$ is contained in the image
of the induced map
$\varphitilde_{\nbigk}:
 \widetilde{\Delta}(0)\times\nbigk
 \lrarr \nbigxtilde(\nbigd)$.
Let $R$ be the inverse image of $Q$
via $\varphitilde_{\nbigk}$.
We take a small neighbourhood 
$\nbigu_R$ of $R$ in
$\Deltatilde(0)\times\nbigk$.
By the $\DD$-flatness,
we have only to show
$\overline{\nbigf}_{|\varphi_{\nbigk}(\nbigu_R)}
=\nbigf^{Q}_{|\varphi_{\nbigk}(\nbigu_R)}$.

The pull back
$\varphi_{\nbigk}^{\ast}(\nbige,\DD)$
has the unramifiedly good lattice
$\varphi_{\nbigk}^{\ast}E$,
and the set of irregular values
is given by
$\nbigj_1:=\bigl\{
 \varphi^{\ast}\gminia\,\big|\,
 \gminia\in\nbigj
 \bigr\}$.
We remark that
the natural map
$\nbigj\lrarr\nbigj_1$ is bijective,
and the orders
$\leq_{Q}^{\varrho}$
and $\leq_R^{\varrho}$ are the same.
Then, by Lemma \ref{lem;10.5.12.20},
we obtain that the restriction of
$\nbigftilde^{Q}$ and
$\overline{\nbigf}$ to $\varphi(\nbigu_R)$
are equal to the full Stokes filtration
of $\varphi^{\ast}(E)_{|\nbigu_{R}}$.
\hfill\qed

\subsection{Splitting of the full Stokes filtration}

\label{subsection;08.1.5.1}

\subsubsection{Flat splitting}

Let $(E,\DD)$ be an unramifiedly good lattice
on $(\nbigx,\nbigd)$
with $\Irr(\DD)=\nbigj$.
First, we consider the non-resonance case.
\begin{condition}
\label{condition;07.10.2.30}
We have
$\alpha-\beta\not\in\seisuu$
for any distinct eigenvalues $\alpha,\beta$ of
$\Res_j(\varrho^{-1}\DD)_{|\nbigd_j}$
$(j=2,\ldots,\ell)$.
\hfill\qed
\end{condition}

If $(E,\DD)$ satisfies Condition
{\rm\ref{condition;07.10.2.30}},
the induced lattices
$\Gr^{\vecm(p)}_{\gminia}(E,\DD)$
also satisfy Condition 
{\rm\ref{condition;07.10.2.30}}
for any $\gminia\in\nbigj(\vecm(p))$,
which follows from (\ref{eq;08.11.1.2}).

\begin{prop}
\label{prop;07.10.2.40}
Assume that $(E,\DD)$ satisfies
Condition {\rm\ref{condition;07.10.2.30}}.
Let $k$ satisfy $k(0)\leq k\leq \ell$.
Take $Q\in\pi_k^{-1}(\nbigd)
\subset\nbigxtilde(W(\kbar))$
such that $\varrho(Q)\neq 0$.
Then, there exists a small neighbourhood
$\nbigu_{Q}$ on which 
we can take a $\DD$-flat splitting
$E_{|\nbigu_{Q}}
=\bigoplus_{\gminia\in\nbigj}
 E_{\gminia,\nbigu_{Q}}$
of the full Stokes filtration $\nbigftilde^S$.
\end{prop}
\pf
We have only to consider the case
$Q\in\pi_k^{-1}(\nbigd_{\ellsitabar})$.
We take $\nbigu_{Q_p}$ as in
Subsection \ref{subsection;07.11.18.10}.
By Proposition \ref{prop;07.9.30.5},
we can find a $\DD$-flat splitting 
of the Stokes filtration $\nbigf^{Q_p\,\vecm(p)}$,
i.e.,
$\Gr^{\vecm(p-1)}(E)_{|\nbigu_{Q_p}}
\simeq
 \Gr^{\vecm(p)}(E)_{|\nbigu_{Q_p}}$.
Then, we can construct a desired splitting
by lifting the splittings inductively.
\hfill\qed

\begin{prop}
\label{prop;10.5.12.50}
Assume that $\nbigd$ is smooth.
For any $Q\in\pi^{-1}(\nbigd)$
there exists a small neighbourhood $\nbigu_Q$
on which we can take a $\DD$-flat splitting
$E_{|\Sbar}
=\bigoplus_{\gminia\in\nbigj}
 E_{\gminia,S}$
of the full Stokes filtration $\nbigftilde^S$
of $E_{|\Sbar}$.
\end{prop}
\pf
We have only to apply 
Lemma \ref{lem;08.1.25.51}
inductively.
\hfill\qed

\subsubsection{Partially flat splitting}

Even in the general case,
we have partially flat splittings,
which can be shown 
by the argument
in the proof of Proposition \ref{prop;07.10.2.40}.
\begin{lem}
\label{lem;08.1.25.50}
Let $k(0)\leq k\leq \ell$,
and let $Q\in
 \pi_k^{-1}(\nbigd)\subset\nbigxtilde(W(\kbar))$.
There exists a neighbourhood $\nbigu_Q$
on which we can take
a $\DD_{\leq k(p)}$-flat splitting
of the partial Stokes filtration
$\nbigf^{Q\,\vecm(p)}$ in the level $\vecm(p)$.
\hfill\qed
\end{lem}

\subsection{Characterization of 
 holomorphic sections of $E$}

Let $k(0)\leq k\leq \ell$.
Let $\vecv$ be a frame of $E$,
and $\vecu$ be a frame of $\Gr^{\vecm(p)}(E)$.
Take $Q\in\pi_k^{-1}(\nbigd_{\ellsitabar})$
and a small neighbourhood $\nbigu_Q$
on which we have the Stokes filtrations
$\nbigf^{Q\,\vecm(p)}$
and its splitting
$E_{|\nbigu_Q}
\simeq
 \bigoplus_{\gminib\in\nbigj(\vecm(p))}
 \Gr^{\vecm(p)}_{\gminib}(E)_{|\nbigu_Q}$.
By using the splitting,
we obtain a frame $\vecu_{Q}$
of $E_{|\nbigu_Q}$.
\begin{lem}
\label{lem;10.5.12.100}
Let $G_{Q}$ be the matrix-valued function
determined by 
$\vecv_{|\nbigu_Q}=\vecu_Q\,G_Q$.
Then,
$G_Q$ and $G_Q^{-1}$ are bounded
on $\nbigu_Q$.
\hfill\qed
\end{lem}

We take $Q_{j}\in\pi^{-1}(\nbigd_{\ellsitabar})$
$(j=1,\ldots,N)$
such that $\bigcup_{i=1}^N\nbigu_{Q_i}$
contains
$\pi^{-1}(\nbigd_{\ellsitabar})$.
We take $\vecu_{Q_i}$ as above.
Let $f$ be a holomorphic section of
$E_{|\nbigx\setminus\nbigd}$.
We have the expressions
$f_{|\nbigu_{Q_i}\setminus\pi^{-1}(W)}
=\sum f_{Q_i,j}\,u_{Q_i,j}$,
where $f_{Q_i,j}$ are holomorphic functions
on $\nbigu_{Q_i}\setminus\pi^{-1}(W)$.
We obtain the following lemma
from Lemma \ref{lem;10.5.12.100}.
\begin{lem}
$f$ is a section of $E$
if and only if $f_{Q_i,j}$ are bounded.
\hfill\qed
\end{lem}

\subsection{Characterization by growth order}
\label{subsection;08.9.28.40}

Let $\vecv$ be a holomorphic frame of 
$E$ on $\nbigx$.
Let $Q\in\pi^{-1}(\nbigd_{\ellsitabar})$,
and let $\nbigu_Q$ be a small neighbourhood of $Q$.
Let $f$ be a flat section of 
$E_{|\nbigu_Q\setminus \pi^{-1}(W)}$.
We have the expression $f=\sum f_j\,v_j$,
and obtain $\vecf=(f_j)$

\begin{prop}
\label{prop;07.10.11.2}
$f$ is contained in 
$\nbigf^{Q,\vecm(p)}_{\gminia}
 (E_{|\nbigu_Q\setminus \pi^{-1}(W)})$
if and only if the following estimate holds for some
$C>0$ and $N>0$:
\[
 \Bigl|
 \vecf\, \exp\bigl(
 \varrho^{-1}\etabar_{\vecm(p)}(\gminia)\bigr)
 \Bigr|
=O\Bigl(
 \exp\bigl(C|\vecz^{\vecm(p+1)}|\bigr)
 \!\!\!\!\!
 \prod_{k(p+1)<j\leq \ell}
 \!\!\!\!
 |z_j|^{-N}
 \Bigr)
\]
\end{prop}
\pf
We may replace $E$
with a lattice satisfying Condition 
\ref{condition;07.10.2.30}
by a meromorphic transform.
Hence, we may and will assume that
the condition is satisfied from the beginning.
Due to Proposition \ref{prop;07.10.2.40},
we can take a flat splitting
$E_{|\nbigu_Q}=\bigoplus E_{\gminia,Q}$
of the partial Stokes filtration 
$\nbigf^{Q,\vecm(p)}$.
Let $\vecu$ be a frame of
$\Gr^{\vecm(p)}(E)$
compatible with the grading,
and let $\vecu_{Q}$ be the lift of $\vecu$
to $E_{\gminia,Q}$ on $\nbigu_Q$
via the splitting.

We have the expression
$f=\sum f_{Q,j}\cdot u_{Q,j}$.
Let $\vecf_Q:=(f_{Q,j})$.
Corresponding to the grading,
we have the decomposition
$\vecf_Q=(\vecf_{Q,\gminib}\,|\,\gminib
 \in\nbigj(\vecm(p)))$.
By Lemma \ref{lem;10.5.12.100},
$\sum|\vecf_{Q,\gminib}|$
and $|\vecf|$ are mutually bounded.
Then, the claim follows from
Lemma \ref{lem;10.5.12.51}.
\hfill\qed

\subsection{Proof of the claims in
Section \ref{section;10.5.10.5}}
\label{subsection;10.5.18.6}

The filtrations in Theorem \ref{thm;10.5.12.110}
was constructed in
Subsection \ref{subsection;07.11.18.10}.
It clearly satisfies the first claim
in the theorem.
The compatibility is given in
Lemma \ref{lem;10.5.12.41}.
By Lemma \ref{lem;10.5.12.20}
and Lemma \ref{lem;10.5.12.33},
the conditions characterize the filtrations.
If $\varrho(Q)\neq 0$,
the first property suffices for characterization
according to Lemma \ref{lem;10.5.12.200}.
Thus, Theorem \ref{thm;10.5.12.110} is proved.

As remarked in Subsection
\ref{subsection;10.5.12.201},
the functoriality of the full Stokes filtration
follows from the inductive construction
of the full Stokes filtration
and the functoriality in Stokes filtration
of weakly good lattices in 
Subsection \ref{subsection;10.5.11.22}.
Proposition \ref{prop;10.5.12.203}
is Proposition \ref{prop;07.10.11.2}.
Proposition \ref{prop;10.5.12.205}
and the functoriality of $\Gr^{\nbigftilde}$
is clear from our construction of 
the full Stokes filtration.
Proposition \ref{prop;10.5.12.210}
is also clear.

According to Proposition \ref{prop;10.5.10.2},
we can locally take a non-resonance lattice.
Hence, 
we obtain
Proposition \ref{prop;10.5.12.211}
from Proposition \ref{prop;07.10.2.40}.
Proposition \ref{prop;10.5.10.3}
also follows from Proposition \ref{prop;07.10.2.40}.
Proposition \ref{prop;10.5.12.212}
is Proposition \ref{prop;10.5.12.50}.
As remarked in Lemma \ref{lem;08.1.25.50},
we can obtain
Proposition \ref{prop;10.5.12.220}
and Proposition \ref{prop;10.5.12.221}
by an inductive use of 
Proposition \ref{prop;10.5.11.30}.
We obtain Proposition \ref{prop;10.5.12.240}
by successive use of Proposition
\ref{prop;10.5.28.10}.
We also obtain Proposition \ref{prop;10.5.12.241}
by successive use of 
Proposition \ref{prop;10.5.11.31}.

\chapter[Stokes data]{Full Stokes Data
 and Riemann-Hilbert-Birkhoff Correspondence}
\label{chapter;10.5.26.11}
In this chapter, we study the Stokes structure
in more detail,
assuming that $\varrho$ is nowhere vanishing.
It is our purpose to describe an irregular singularity
in terms of more (though not completely) topological data,
called full Stokes data.

In Section \ref{section;10.5.26.20},
we will introduce the notion of 
full pre-Stokes data,
which is a system of filtrations
of a $\varrho$-flat bundle on 
the real blow up.
If we are interested in only ordinary
meromorphic flat bundles,
we have only to consider pre-Stokes data. 
But, we are interested in families and lattices, too.
So we shall introduce the notion of
Stokes data in Section \ref{section;10.5.26.21},
which is a full pre Stokes data
with lattices.
Then, in Section \ref{section;10.5.26.30},
we will establish the correspondence
between Stokes data and unramifiedly good
lattice of meromorphic $\varrho$-flat bundles,
called Riemann-Hilbert-Birkhoff correspondence.

As an application, we study 
the extension of a $\varrho$-meromorphic flat bundle
in Section \ref{section;10.5.26.31}.
The special case given in Section \ref{section;10.5.26.32}
will play an important role 
in our study on wild harmonic bundles.

\section{Full pre-Stokes data}
\label{section;10.5.26.20}
\subsection{Definition}

Let $\nbigx\lrarr \nbigk$ be a smooth
fibration of complex manifolds.
Let $\nbigd$ be a normal crossing hypersurface
of $\nbigx$ such that any intersections of
irreducible components are smooth over $\nbigk$.
Let $\varrho$ be a nowhere vanishing
holomorphic function on $\nbigx$.
Let $\pi:\nbigxtilde(\nbigd)\lrarr\nbigx$
be the real blow up.
The pull back of $\nbigo_{\nbigk}$
via the projection $\nbigxtilde(\nbigd)\lrarr\nbigk$
is also denoted by $\nbigo_{\nbigk}$.
Its restriction to a subset of $\nbigxtilde(\nbigd)$
is also denoted by $\nbigo_{\nbigk}$.

\begin{df}
\label{df;10.5.16.10}
Let $\vecnbigi$ be a good system of
irregular values on $(\nbigx,\nbigd)$.
Let $\nbigu$ be a subset of $\nbigxtilde(\nbigd)$.
Let $\gbigv$ be a locally free
$\nbigo_{\nbigk}$-module on $\nbigu$.
A full pre-Stokes data of $\gbigv$ over $\vecnbigi$
is a system $\vecnbigftilde$ of filtrations
$\nbigftilde^Q$ of germs $\gbigv_Q$
$(Q\in\nbigu\cap\pi^{-1}(\nbigd))$
indexed by 
$\bigl(\nbigi_{\pi(Q)},\leq_Q^{\varrho}\bigr)$
satisfying the following compatibility condition:
\begin{itemize}
\item
Let $Q\in\nbigu\cap\pi^{-1}(\nbigd)$.
Take a small neighbourhood
$\nbigu_Q$ in $\nbigu$
on which the filtration $\nbigftilde^Q$ is given.
Note that when $\nbigu_Q$ is sufficiently small,
we have 
$\gminia\leq_{Q}\gminib$
if and only if 
$\gminia\leq_{Q'}\gminib$
for any $Q'\in \nbigu_Q\cap\pi^{-1}(\nbigd)$.
Then, for any $Q'\in \nbigu_Q\cap\pi^{-1}(\nbigd)$,
$\bigl(
 \gbigv_{Q'},\nbigftilde^Q
 \bigr)
 \lrarr
 \bigl(\gbigv_{Q'},\nbigftilde^{Q'}\bigr)$
is compatible over
$\bigl(\nbigi_{\pi(Q)},\leq_Q^{\varrho}\bigr)
\lrarr
 \bigl(\nbigi_{\pi(Q')},\leq_{Q'}^{\varrho}\bigr)$.
\hfill\qed
\end{itemize}
\end{df}
By the compatibility,
we have the associated graded 
$\nbigo_{\nbigk}$-module on
$\pi^{-1}(P)\cap \nbigu$
$(P\in \nbigd)$,
which is denoted by
$\Gr^{\nbigftilde}(\gbigv_{\pi^{-1}(P)\cap\nbigu})$.

\subsubsection{}

Let $\gbigv_i$ $(i=1,2)$ be 
$\nbigo_{\nbigk}$-modules on
$U\subset\nbigxtilde(\nbigd)$
equipped with full pre-Stokes data
$\vecnbigftilde_i$ over $\vecnbigi$.
A morphism
$F:(\gbigv_1,\vecnbigftilde_1)
 \lrarr
 (\gbigv_2,\vecnbigftilde_2)$
is defined to be a morphism
of $\nbigo_{\nbigk}$-modules
such that the induced morphisms
$\gbigv_{1\,Q}\lrarr \gbigv_{2\,Q}$
preserve filtrations
for any $Q\in\pi^{-1}(\nbigd)$.
If $\gbigv$ is equipped with
a full pre-Stokes data over $\vecnbigi$,
then the dual $\gbigv^{\lor}$
is equipped with an induced full pre-Stokes
data over $\vecnbigi^{\lor}$.
Let $\gbigv_i$ be equipped with
full pre-Stokes data over $\vecnbigi_i$.
If $\vecnbigi_1\otimes\vecnbigi_2$
is good,
then $\gbigv_1\otimes\gbigv_2$
is equipped with an induced 
full pre-Stokes data over 
$\vecnbigi_1\otimes\vecnbigi_2$.
If $\vecnbigi_1\oplus\vecnbigi_2$
is good,
then $\gbigv_1\oplus\gbigv_2$
is equipped with an induced 
full pre-Stokes data over 
$\vecnbigi_1\oplus\vecnbigi_2$.

\subsubsection{Uniqueness}

Let $U$ be an open subset of $\pi^{-1}(\nbigd)$.
Let $\gbigv$ be an $\nbigo_{\nbigk}$-module on $U$.
\begin{lem}
\label{lem;10.5.16.40}
Let $\vecnbigftilde_i$ $(i=1,2)$
be full pre-Stokes data of $\gbigv$.
\begin{itemize}
\item
If there exists a dense subset $U'\subset U$
such that
$\nbigftilde^{Q}_1=\nbigftilde^{Q}_2$
for $Q\in U'$.
Then, we have
$\vecnbigftilde_1=\vecnbigftilde_2$.
\item
Let $Z$ be any subset of $U$.
If $\nbigftilde^Q_1=\nbigftilde^Q_2$
for any $Q\in Z$,
there exists a neighbourhood $Z'$ of $Z$
such that
$\nbigftilde^Q_1=\nbigftilde^Q_2$
for any $Q\in Z'$.
\end{itemize}
\end{lem}
\pf
The first claim easily follows from
Lemma \ref{lem;10.5.16.3}
and Lemma \ref{lem;10.5.12.22}.
The second claim is clear
from the compatibility condition.
\hfill\qed

\vspace{.1in}
Let $\gbigv_i$ $(i=1,2)$
be $\nbigo_{\nbigk}$-modules on $U$
with a morphism
$F:\gbigv_1\lrarr\gbigv_2$.
It is easy to deduce the following corollary.
\begin{cor}
\label{cor;10.5.19.1}
Let $\vecnbigftilde_i$ $(i=1,2)$ be
full pre-Stokes structure of $\gbigv_i$.
\begin{itemize}
\item
If there exists a dense subset $U'\subset U$
such that
$F$ preserves $\nbigftilde^Q$ for $Q\in U'$.
Then, $F$ preserves $\vecnbigftilde$.
\item
Let $Z$ be any subset of $U$.
If $F$ preserves $\nbigftilde^Q$
for any $Q\in Z$,
there exists a neighbourhood $Z'$ of $Z$
such that
$F$ preserves $\nbigftilde^Q$
for any $Q\in Z'$.
\hfill\qed
\end{itemize}
\end{cor}

\subsection{Global filtration and splitting}

Let $\nbigx=\Delta^n\times\nbigk$
and $\nbigd=\bigcup_{i=1}^{\ell}\{z_i=0\}$.
Let $\nbigi\subset M(\nbigx,\nbigd)/H(\nbigx)$
be a good set of irregular values.
Let $P\in\nbigd_{\ellsitabar}$.
We identify $\pi^{-1}(P)\simeq (S^1)^{\ell}$,
and we use the polar coordinate
$(\theta_1,\ldots,\theta_{\ell})$.

\begin{condition}
\label{condition;10.5.17.30}
Let $U$ be a closed convex subset of $\pi^{-1}(P)$
satisfying the following:
\begin{itemize}
\item
 There exist 
 $(\theta_1^{(0)},\ldots,\theta_{\ell}^{(0)})$
 such that
 $U$ is contained in 
 $\bigl\{
 (\theta_1,\ldots,\theta_{\ell})\,\big|\,
 |\theta_i-\theta_i^{(0)}|<\pi/2
 \bigr\}$.
 In particular,
 we can identify $U$ with a closed region
 in $\real^{\ell}$.
\item
 For each distinct pair $(\gminia,\gminib)$ of $\nbigi$,
 if $U\cap \Sep(\gminia,\gminib)\neq\emptyset$,
 it divides $U$ into two closed regions.
\hfill\qed
\end{itemize}
\end{condition}

\begin{prop}
\label{prop;10.5.15.10}
Let $\gbigv$ be a free $\nbigo_K$-module
on $U$ with full pre-Stokes data
$\bigl(
 \nbigftilde^Q\,\big|\,Q\in U
 \bigr)$.
Then, there uniquely exists a global filtration
$\nbigftilde^U$ 
indexed by 
$\bigl(\nbigi_P,\leq_U^{\varrho}\bigr)$
such that
for any $Q\in U$,
$\nbigftilde^U$ and $\nbigftilde^Q$
are compatible over
$\bigl(
 \nbigi_P,\leq_U^{\varrho}
 \bigr)\lrarr
 \bigl(\nbigi_P,\leq_Q^{\varrho}\bigr)$.
In other words,
there exists a decomposition
$V=\bigoplus_{\gminia\in\nbigi}V_{\gminia}$,
which gives a splitting of
$\nbigftilde^Q$ for each $Q\in U$.
(See Notation {\rm\ref{notation;10.6.11.2}}
for $\leq_U$.)
\end{prop}
\pf
The uniqueness follows from Lemma 
\ref{lem;10.5.12.22}.
Put $H_{\gminia,\gminib}:=
F_{\gminia,\gminib}^{-1}(0)$
for distinct $\gminia,\gminib\in\nbigi$.
A connected component of
$U\setminus \bigcup H_{\gminia,\gminib}$
is called a chamber.
If $Q$ is contained in a chamber,
then $\leq_Q$ is totally ordered.
If $Q$ and $Q'$ are contained in
the same chamber,
we have
$\leq_{Q}^{\varrho}=\leq_{Q'}^{\varrho}$.

Let $Q$ and $Q'$ be contained in 
two chambers
divided by a wall given as $H_{\gminia,\gminib}$.
Then, 
$\gminia\leq_Q\gminib$
if and only if 
$\gminib\leq_{Q'}\gminia$.
If $F_{\gminia',\gminib'}\neq 0$ 
on the wall,
then $\gminia\leq_Q\gminib$
if and only if $\gminia\leq_{Q'}\gminib$.

We remark that
an element $\gminia\in\nbigi_P$
is minimal with respect to
$\leq_U^{\varrho}$,
if and only if there exists 
an interior point $Q_0\in U$
such that $\gminia$
is the minimum with respect to
$\leq_{Q_0}^{\varrho}$.

Let $V$ be the space of global sections
of $\gbigv$.
We have natural isomorphisms
$V\simeq \gbigv_Q$ for any $Q\in U$.
We regard that we are given filtrations
$\nbigftilde^Q$ $(Q\in U)$ on $V$.
Let $\gminia$ be minimal with respect to
$\leq_U^{\varrho}$, and take $Q_0$ as above.
Let us observe that $\nbigftilde_{\gminia}^{Q_0}$
is contained in
$\nbigftilde^{Q}_{\gminia}$ 
for any $Q\in U$.
We take an interval $I$ connecting
$Q$ and $Q_0$.
Let $R$ be an intersection of $I$
with a wall.
Let $R_-$ and $R_+$
be points in $I$ such that
(i) they are sufficiently close to $R$, 
(ii) $Q_0$ is closer to $R_-$ than $R_+$.
In general, we have
$\nbigftilde_{\gminia}^{R_-}\supset
 \nbigftilde_{\gminia}^{R}\subset
 \nbigftilde_{\gminia}^{R_+}$.
By our assumption,
we have $F_{\gminia,\gminib}(R_-)<0$
for any $\gminib$ such that
$F_{\gminia,\gminib}(R)=0$,
which implies
$\nbigftilde^{R}_{\gminia}
=\nbigftilde^{R_-}_{\gminia}$.
Therefore, we obtain
$\nbigftilde^{Q_0}_{\gminia}
\subset
 \nbigftilde^{Q}_{\gminia}$.
We can also deduce that
$\nbigftilde^{Q_0}_{\gminia}
\lrarr \Gr^{\nbigftilde^Q}_{\gminia}$
is an isomorphism
for any $Q$.
Hence, in particular,
if $\gminib\neq \gminia$ is minimal
with respect to $\leq_Q^{\varrho}$,
we have $\nbigftilde^{Q}_{\gminib}\cap
 \nbigftilde^{Q_0}_{\gminia}=0$.

We put 
$V^{(0)}:=V/\nbigftilde^{Q_0}_{\gminia}$.
For any $Q\in U$ and $\gminib\in\nbigi$,
let $\nbigftilde^Q_{\gminib}(V^{(0)})$
be the image of $\nbigftilde^Q_{\gminib}(V)$
to $V^{(0)}$.
Let $V=\bigoplus V_{\gminib,Q}$
be a splitting of $\nbigftilde^Q$.
We remark that we may assume
$V_{\gminia,Q}=\nbigftilde^{Q_0}_{\gminia}$.
Then, it is easy to see that the images
of $V_{\gminib,Q}$ to $V^{(0)}$
gives a splitting of
$\nbigftilde^Q_{\gminib}$ of 
$(V^{(0)},\nbigftilde^Q)$.
We can also easily observe that
the system of filtrations
$\bigl(\nbigftilde^Q(V^{(0)})\,\big|\,Q\in U
 \bigr)$ is a full pre-Stokes data of $V^{(0)}$.

Assume that 
we have filtrations $\nbigftilde^U$
with the desired property
for both $V^{(0)}$ and $V$.
Then, $\nbigftilde^U_{\gminib}(V^{(0)})$
is obtained as the image of
$\nbigftilde^U_{\gminib}(V)$.
Actually, 
let $V=\bigoplus_{\gminib\in\nbigi}V_{\gminib}$
be a splitting of $\nbigftilde^U(V)$.
We have
$V_{\gminia}=\nbigftilde^{Q_0}_{\gminia}$.
The decomposition also gives a splitting of
$\nbigftilde^Q(V)$ for each $Q\in U$.
We have the induced decomposition
$V^{(0)}=\bigoplus V^{(0)}_{\gminib}$,
which gives a splitting of
$\nbigftilde^Q(V^{(0)})$ for each $Q\in U$.
It implies that the decomposition
gives a splitting of
$\nbigftilde^U(V^{(0)})$ by the uniqueness,
and we can conclude that
$\nbigftilde^U(V^{(0)})$ is obtained 
as the image of $\nbigftilde^U(V)$.

Let us show the claim of the proposition
by using an induction on $|\nbigi|$.
The case $|\nbigi|=1$ is obvious.
Take a minimal $\gminia$ with respect to
$\leq_U^{\varrho}$,
and let $Q_0$ be as above.
We can apply the hypothesis of the induction
to $V^{(0)}=V/\nbigftilde^{Q_0}_{\gminia}$.
If $\gminia$ is the minimum with respect to
$\leq^{\varrho}_U$,
we can construct the desired filtration of $V$
as the pull back via $V\lrarr V^{(0)}$.
Assume that
there is another minimal element $\gminib$.
There exists an interior point $Q_1\in U$ such that
$\gminib$ is the minimum with respect to
$\leq_{Q_1}$.
We remark $\nbigftilde^{Q_1}_{\gminib}\cap
 \nbigftilde^{Q_0}_{\gminia}=0$.
We put $V^{(1)}:=V/\nbigftilde^{Q_1}_{\gminib}$
and $V^{(2)}:=V/\bigl(
 \nbigftilde^{Q_1}_{\gminib}
\oplus\nbigftilde^{Q_0}_{\gminia}
 \bigr)$.
As remarked above,
they are equipped with the induced
full pre-Stokes structure.
By construction, we have
$\nbigftilde^{Q}_{\gminic}(V)=
 \nbigftilde^Q_{\gminic}(V^{(1)})
 \times_{\nbigftilde^Q_{\gminic}(V^{(2)})}
 \nbigftilde^Q_{\gminic}(V^{(0)})$.

By the hypothesis of the induction,
$V^{(i)}$ are equipped with
the filtration $\nbigftilde^U$
with the desired property.
Note that $\nbigftilde^U(V^{(2)})$
is obtained as the image of
$\nbigftilde^U(V^{(i)})$ $(i=0,1)$.
We put
$\nbigftilde_{\gminic}^U(V):=
\nbigftilde^{U}_{\gminic}(V^{(0)})
\times_{\nbigftilde^U_{\gminic}(V^{(2)})}
\nbigftilde^{U}_{\gminic}(V^{(1)})$.
Let us check that
$\nbigftilde^U$ has the desired property.
Let $V^{(2)}=\bigoplus V^{(2)}_{\gminic}$
be a splitting of $\nbigftilde^U$.
Let $V^{(0)}_{\gminic}\subset V^{(0)}$ be a lift
of $V^{(2)}_{\gminic}$.
We put 
$V_{\gminia}:=\nbigftilde^{Q_0}_{\gminia}$
and $V_{\gminib}:=\nbigftilde^{Q_1}_{\gminib}$.
By using that
$\nbigftilde^U(V^{(2)})$ is obtained as the image of
$\nbigftilde^U(V^{(i)})$ $(i=0,1)$,
we can check that
$V=\bigoplus V_{\gminic}$ is a splitting of
$\nbigftilde^U$.
Similarly, we can check that it gives a splitting 
of each $\nbigftilde^Q(V)$.
Hence, $\nbigftilde^U$ is compatible with
$\nbigftilde^Q$ for each $Q\in U$.
\hfill\qed

\vspace{.1in}

It seems useful
to consider the following type of covering
of $\pi^{-1}(P)$ $(P\in\nbigd)$.
\begin{df}
\label{df;10.5.18.35}
A finite covering $\{U_i\,|\,i\in\Gamma\}$
of $\pi^{-1}(P)$ is called good for $\nbigi_P$,
if any intersection 
$U_I:=\bigcap_{i\in I}U_i$
$(I\subset \Gamma)$ satisfy 
Condition {\rm\ref{condition;10.5.17.30}}.
\hfill\qed
\end{df}

\begin{lem}
For a given $\nbigi_P$,
there exists a good cover of $\pi^{-1}(P)$.
\end{lem}
\pf
For example,
we can construct such a finite covering
by using polytopes surrounded by
hypersurfaces generic to
$H_{\gminia,\gminib}$ 
for any distinct $\gminia,\gminib\in\nbigi_P$.
(See the proof of Proposition \ref{prop;10.5.15.10}
for $H_{\gminia,\gminib}$.)
\hfill\qed

\subsection{Full pre-Stokes data
for $\varrho$-flat bundle}
\label{subsection;10.5.19.10}

Let $(V,\DD)$ be a $\varrho$-flat bundle
on $\nbigx\setminus\nbigd$.
By considering $\DD$-flat and $\nbigk$-holomorphic
sections of $V$,
we obtain an $\nbigo_{\nbigk}$-module
$\gbigv'$ on $\nbigx\setminus\nbigd$.
By taking the push-forward via
$\iota:\nbigx\setminus\nbigd\subset\nbigxtilde(\nbigd)$,
we obtain an $\nbigo_{\nbigk}$-module
$\gbigv$ on $\nbigxtilde(\nbigd)$,
that is what we are mainly concerned.
Let $A$ be any subset of $\nbigxtilde(\nbigd)$.
If $\gbigv$ is given in this way,
a full pre-Stokes data of $\gbigv_{|A}$
is equivalent to
a system of $\DD$-flat filtrations
$\nbigftilde^Q$ of germs  $\iota_{\ast}(V)_{Q}$
$(Q\in A\cap\pi^{-1}(\nbigd))$
such that 
(i) they satisfy the compatibility condition
as in Definition \ref{df;10.5.16.10},
(ii) each $\nbigftilde^Q$ has a $\DD$-flat
 splitting.

We immediately obtain the following lemma
from Proposition \ref{prop;10.5.15.10}.
\begin{prop}
\label{prop;10.5.17.30}
Let $P\in \nbigd$.
If $U_1\subset \pi^{-1}(P)$ is sufficiently
small so that there exists
$U_2\supset U_1$ satisfying the condition
in Proposition {\rm\ref{prop;10.5.15.10}},
then there uniquely exists
a filtration $\nbigftilde^{U_1}$ of 
$\iota_{\ast}(V)_{|U_1}$
such that 
$\nbigftilde^{U_1}$ and
$\nbigftilde^{Q}$ $(Q\in U_1)$
are compatible over
$\bigl(\nbigi_{P},\leq_{U_1}^{\varrho}\bigr)
\lrarr
 \bigl(\nbigi_P,\leq_{Q}^{\varrho}\bigr)$.
\hfill\qed
\end{prop}

\begin{lem}
\label{lem;10.5.16.90}
Let $(V,\DD)$ be as above.
Let $P\in\nbigd$.
If we are given a full pre-Stokes data of
$\iota_{\ast}(V)_{|\pi^{-1}(P)}$,
it is uniquely extended to
a full pre-Stokes data on a small
neighbourhood of $\pi^{-1}(P)$.
\end{lem}
\pf
Let $\{U_i\,|\,i\in\Gamma\}$
be a good cover of $\pi^{-1}(P)$
for $\nbigi_P$.
For each $I\subset\Gamma$,
we have the filtration
$\nbigftilde^{U_I}$ of
$\iota_{\ast}(V)_{|U_I}$.
If $I\subset J$,
the filtration $\nbigftilde^{U_J}$
is induced by the restriction of
$\nbigftilde^{U_I}$ to $U_J$
and $\bigl(\nbigi_P,\leq_{U_I}^{\varrho}\bigr)
\lrarr\bigl(\nbigi_P,\leq_{U_J}^{\varrho}\bigr)$.
We take a small neighbourhood $\nbigu_I$
of $U_I$ in $\nbigxtilde(\nbigd)$
such that (i) $\leq_{\nbigu_I}=\leq_{U_I}$,
(ii) $\nbigu_J\subset\nbigu_I$ for $I\subset J$.
We put 
$\nbigu:=\bigcup_{I\subset\Gamma}\nbigu_I$.
For each $Q\in\nbigu$,
we can find $\nbigu_I\ni Q$.
Let $\nbigftilde^Q$ be the filtration
of $\iota_{\ast}(V)_Q$ induced by 
$\nbigftilde^{U_I}$
and $\bigl(\nbigi_P,\leq_{U_I}^{\varrho}\bigr)
\lrarr\bigl(\nbigi_{\pi(Q)},\leq_{Q}^{\varrho}\bigr)$.
It is easy to check the well definedness
of $\nbigftilde^Q$,
and that the system of filtrations is
a full pre-Stokes data of
$\iota_{\ast}(V)_{|\nbigu}$.
\hfill\qed

\section{Full Stokes data}
\label{section;10.5.26.21}
\subsection{The associated graded bundles}
\label{subsection;10.5.14.2}

Let $(V,\DD)$ be a $\varrho$-flat bundle
on $\nbigx\setminus\nbigd$.
For each $P\in\nbigd,$
$\nbigx_P$ denotes a small neighbourhood of $P$.
We put $\nbigd_P:=\nbigx_P\cap\nbigd$.
We will shrink it without mention.
Let $\iota:\nbigx\setminus\nbigd
 \subset\nbigxtilde(\nbigd)$
be the inclusion.
Let $\vecnbigftilde$
be a full pre-Stokes data of $\iota_{\ast}V$
over a good system of irregular values $\vecnbigi$.
For any point $P\in \nbigd$,
we obtain a graded sheaf
$\Gr^{\vecnbigftilde}(\iota_{\ast}V_{|
 \pi^{-1}(\nbigx_P)})$
with a $\varrho$-flat connection $\DD$.
In particular,
we obtain a graded $\varrho$-flat bundle
on $\nbigx_P^{\ast}$,
where
$\nbigx_P^{\ast}=\nbigx_P\setminus \nbigd_P$:
\[
 V_{P,\gminia}:=
 \Gr^{\vecnbigftilde}_{\gminia}\bigl(
 V_{|\nbigx_P^{\ast}}\bigr)
\quad
 V_P:=\bigoplus_{\gminia\in\nbigi_P}
 V_{P,\gminia}
\]

\subsubsection{Variant}
\label{subsection;10.5.17.5}

Let $\nbigd=\bigcup_{i\in\Lambda}\nbigd_i$
be the irreducible decomposition.
Let $I(P):=\bigl\{i\in\Lambda\,\big|\,
 P\in\nbigd_i\bigr\}$.
For any $J\subset I(P)$,
we put $J^c:=I(P)\setminus J$.
Let $\nbigi_P^J$ be the image of 
$\etabar_J:\nbigi_P\lrarr
 \nbigo_{\nbigx}(\ast \nbigd)_P/
 \nbigo_{\nbigx}(\ast \nbigd(J^c))_P$.
We obtain a filtration
$\nbigf^{Q,J}$ of $\iota_{\ast}(V)_Q$
indexed by
$\bigl(\nbigi_P^J,\leq_{Q}^{\varrho}
 \bigr)$ for $Q\in\pi^{-1}(P)$,
and the induced filtration
$\nbigftilde^Q$ on 
$\Gr_{\gminib}^{\nbigf^{Q,J}}
 \bigl(\iota_{\ast}(V)_Q\bigr)$
indexed by 
$\bigl(\etabar_J^{-1}(\gminib),
 \leq^{\varrho}_Q\bigr)$.
In particular,
we obtain a system of filtrations
$\vecnbigf^J:=\bigl(
 \nbigf^{Q,J}\,\big|\,Q\in\pi^{-1}(P)
 \bigr)$.
It satisfies the compatibility condition
in the following sense.
\begin{itemize}
\item
Take a small $\nbigu_Q\in\gbigu(Q,\nbigi_P^J)$
on which $\nbigf^{Q,J}$ is given.
(Recall the notation in Section 
{\rm\ref{subsection;10.6.11.3}}.)
For $Q'\in \nbigu_Q\cap\pi^{-1}(P)$,
we have 
the induced filtration $\nbigf^{Q,J}$
and $\nbigf^{Q',J}$
of $\iota_{\ast}(V)_{Q'}$.
Then, they are compatible over
$\bigl(\nbigi_P^J,\leq_{Q}^{\varrho}
 \bigr)\lrarr
 \bigl(\nbigi_{P}^J,\leq_{Q'}^{\varrho}\bigr)$.
\end{itemize}
For $\gminib\in\nbigi_P^J$,
we obtain
$\Gr_{\gminib}^{\vecnbigf^J}
 (\iota_{\ast}(V)_{\pi^{-1}(\nbigx_P)})$
with an induced $\varrho$-flat connection
on $\pi^{-1}(\nbigx_P)$.
In particular, we obtain
a graded $\varrho$-flat bundle
on $\nbigx_P^{\ast}$:
\[
 V^J_{P,\gminib}:=
 \Gr_{\gminib}^{\vecnbigf^J}(V_{|\nbigx_P^{\ast}}),
 \quad
 V^J_{P}:=\bigoplus_{\gminib\in\nbigi_P^J}
 V^J_{P,\gminib}
\]
We have an induced full pre-Stokes data
of $\iota_{\ast}(V^J_{P,\gminib})_{|\pi^{-1}(P)}$
over the good set $\etabar_J^{-1}(\gminib)$.
By construction,
we have a natural isomorphism
of graded $\varrho$-flat bundle
\begin{equation}
\label{eq;10.5.17.2}
 \Gr^{\vecnbigftilde}V^J_P
\simeq
 V_P.
\end{equation}
Let $P'\in\nbigx_P\cap\nbigd_J$
with $I(P')=J$.
By the compatibility of the given $\vecnbigftilde$,
we have a natural isomorphism
\begin{equation}
 \label{eq;10.5.17.1}
 V^J_{P|\nbigx_{P'}^{\ast}}\simeq
 V_{P'} 
\end{equation}

\subsection{Full Stokes data}
\label{subsection;10.5.16.21}

We continue to use the notation
in Subsection \ref{subsection;10.5.14.2}.
A graded extension of
a full pre-Stokes data $\vecnbigftilde$
is a datum 
$\vecnbigg=
 \bigl(\nbigg_{P}\,\big|\,
 P\in\nbigd \bigr)$
as follows:
\begin{itemize}
\item
 $\nbigg_P=\bigoplus_{\gminia\in\nbigi_P}
 \nbigg_{P,\gminia}$ is a graded locally free 
 $\nbigo_{\nbigx_P}$-module
 with an isomorphism
 $\nbigg_{P,\gminia|\nbigx_P^{\ast}}
 \simeq V_{P,\gminia}$
 such that $\DD_{\gminia}$ is
 $\gminia$-logarithmic
 with respect to $\nbigg_{P,\gminia}$.
\end{itemize}
Note it induces a subsheaf
\[
 \pi_P^{\ast}\nbigg_P:=
 \pi_P^{-1}\nbigg_{P}
 \otimes\nbigo_{\nbigxtilde_P(\nbigd_P)}
\subset
 \iota_{P\ast}V_P,
\]
where
$\pi_P:\nbigxtilde_P(\nbigd_P)
\lrarr \nbigx_P$
and $\iota_P:\nbigx_P\setminus\nbigd_P\subset
 \nbigx_P$
are natural maps.
Namely, for each $P\in\nbigd$,
we are considering an
$\nbigo_{\nbigxtilde_P(\nbigd_P)}$-lattice
of $\iota_{P\ast}V_P$
whose push-forward to $\nbigx_P$
is a locally free $\nbigo_{\nbigx_P}$-module.
We remark that
$\nbigg_P$ is recovered as
$\pi_{P\ast}\pi_{P}^{\ast}\nbigg_P$,
which will be implicitly used.

We introduce a compatibility condition
for the graded extension $\vecnbigg$.
First, we impose the following compatibility.
\begin{itemize}
\item
If $P'\in \nbigd_{P,I(P)}$,
we impose 
$\nbigg_{P',\gminia}=
 \nbigg_{P,\gminia|\nbigx_{P'}}$.
\end{itemize}
Then, for each $J\subset I(P)$ and
$\gminib\in\nbigi_P^J$,
we obtain a locally free $\nbigo_X$-module
$\nbigg^J_{\gminib}$
on $\nbigv\setminus\nbigd_P(J^c)$,
where $\nbigv$ is a neighbourhood
of $\nbigd_{P,J}$.
Actually, it is obtained as the gluing of
$V^J_{P,\gminib}$ and $\nbigg_{P',\gminib}$
for $P'\in \nbigd_{P,J}\setminus
 \bigcup_{i\not\in J}\nbigd_{P,i}$.
We remark (\ref{eq;10.5.17.1}).
By shrinking $\nbigx_P$,
we may assume $\nbigx_P=\nbigv$.
Let $\pi_{P,J}$ be the restriction of
$\pi_P$ to 
$\nbigxtilde_P(\nbigd_P)\setminus
\pi_P^{-1}(\nbigd_P(J^c))$.
Let $\iota_{P,J}:\nbigxtilde_P(\nbigd_P)\setminus
 \pi_P^{-1}(\nbigd_P(J^c))
\subset\nbigxtilde_P(\nbigd_P)$.
Then, we obtain a subsheaf
$\iota_{P,J\ast}
 \pi_{P,J}^{\ast}
 \nbigg^J_{\gminib}
\subset
 \iota_{P\ast}V^J_{P,\gminib}$.
By construction,
it induces an isomorphism
\begin{equation}
 \label{eq;10.5.17.3}
 \iota_{P,J\ast}
 \pi_{P,J}^{\ast}
 \nbigg^J_{\gminib}
 \otimes
 \iota_{P\ast}\nbigo_{\nbigx_P^{\ast}}
\simeq
 \iota_{P\ast}V^J_{P,\gminib}.
\end{equation}
Then, the compatibility condition
is given as follows:
\begin{itemize}
\item
For any $Q\in \pi^{-1}(P)$,
the filtration $\nbigftilde^Q$ of
$\iota_{P\ast}V^J_{P,\gminib}$
is induced by a filtration of
$\iota_{P,J\ast}\pi_{P,J}^{\ast}
 \nbigg^J_{\gminib}$
and (\ref{eq;10.5.17.3}).
\item
The restrictions of
$\pi_P^{\ast}\nbigg_P$
and
$\Gr^{\vecnbigftilde} 
 \iota_{P,J\ast}\pi_{P,J}^{\ast}
 \nbigg^J_{\gminib}$
to $\nbigxtilde_P(\nbigd_P)\setminus
 \pi^{-1}(\nbigd_P(J^c))$
are isomorphic,
extending (\ref{eq;10.5.17.2}).
\end{itemize}
If the compatibility condition is satisfied,
$(\vecnbigftilde,\vecnbigg)$ is called
a full Stokes data.

\subsubsection{}
Let $(V_i,\DD_i)$ be  $\varrho$-flat bundles
on $\nbigx\setminus\nbigd$
equipped with full Stokes data 
$\SDtilde_i=(\vecnbigftilde_i,\vecnbigg_i)$.
A morphism 
$F:(V_1,\DD_1,\SDtilde_1)
 \lrarr
 (V_2,\DD_2,\SDtilde_2)$
is defined to be a morphism
$(V_1,\DD_1,\vecnbigftilde_1)\lrarr
 (V_2,\DD_2,\vecnbigftilde_2)$
such that, for any $P\in\nbigd$
and $\gminia\in\nbigi_P$,
the induced morphisms
$\Gr_{\gminia}^{\nbigftilde}(V_{1|\nbigx_P^{\ast}})
\lrarr \Gr_{\gminia}^{\nbigftilde}(V_{2|\nbigx_P^{\ast}})$
are extended to 
$\nbigg_{1,P,\gminia}\lrarr 
 \nbigg_{2,P,\gminia}$.
Full Stokes data has obvious functoriality
for dual, tensor product and direct sum.

\subsubsection{Meromorphic Stokes data}

A meromorphic graded extension 
of a full pre-Stokes data $\vecnbigftilde$
is a tuple
$\vecnbigg=
 \bigl(\nbigg_{P}\,\big|\,
 P\in\nbigd \bigr)$, where 
each
 $\nbigg_P=\bigoplus_{\gminia\in\nbigi_P}
 \nbigg_{P,\gminia}$ is a locally free 
 $\nbigo_{\nbigx_P}(\ast\nbigd_P)$-module
 with isomorphisms
 $\nbigg_{P,\gminia|\nbigx_P^{\ast}}
 \simeq
 \Gr^{\vecnbigftilde}_{\gminia}(V_{|\nbigx_P^{\ast}})$
 such that 
 $\DD_{\gminia}$ are $\gminia$-regular
 with respect to $\nbigg_{P,\gminia}$.
We can consider a  compatibility condition
similar to the above.
If the compatibility condition is satisfied,
$(\vecnbigftilde,\vecnbigg)$ is called
a meromorphic full Stokes data.

\subsubsection{Another formulation of compatibility}
\label{subsection;10.5.14.1}

We give another formulation
of compatibility condition
for $\vecnbigg$.
Let us consider
$\nbigxtilde_P(\nbigd_P)
\stackrel{\varpi}{\lrarr}
 \nbigxtilde_P(\nbigd_P(J^c))
\stackrel{\pi_1}{\lrarr}\nbigx_P$.
Take $Q_1\in\pi_1^{-1}(P)$.
We remark that
the filtrations $\nbigftilde^{Q}$
$(Q\in \varpi^{-1}(Q_1))$ on
$\iota_{\ast}V^J_{P,\gminib}$
are constant,
because the orders $\leq^{\varrho}_{Q}$
on $\etabar_J^{-1}(\gminib)$ 
are independent of
$Q\in\varpi^{-1}(Q_1)$.
Let $\iota_J:\nbigx_P^{\ast}\subset
 \nbigxtilde_P(\nbigd_P(J^c))$.
By the above consideration,
we obtain an induced $\DD$-flat filtration
$\nbigftilde^{Q_1}$ on
the stalk
$\iota_{J\,\ast}
 (V^J_{P,\gminib})_{Q_1}$.
The system of filtrations 
$\bigl(\nbigftilde^{Q_1}\,\big|\,
 Q_1\in\pi_1^{-1}(P)\bigr)$
satisfies the standard compatibility condition.

Let $\iota_J'$
denote the inclusions of
$\nbigxtilde_P\setminus\nbigd_P(J^c)$
into $\nbigxtilde_P(\nbigd_P(J^c))$.
We have the subsheaf
$\iota_{J\ast}'\bigl(
 \lefttop{J}\nbigg_{\gminia}\bigr)
\subset
 \iota_{J\ast}V^J_{P,\gminib}$.
Now, we can state the compatibility condition.
\begin{itemize}
\item
For any $Q_1\in\pi_1^{-1}(P)$,
the filtration $\nbigftilde^{Q_1}$ of
$\iota_{J\ast}(V^J_{P,\gminib})_{Q_1}$
comes from a filtration of
 $\iota'_{J\ast}\bigl(
 \lefttop{J}\nbigg_{\gminia}\bigr)_{Q_1}$.
 We obtain
 $\Gr_{\gminia}^{\vecnbigftilde}\bigl(
 \lefttop{J}\nbigg_{\gminib}\bigr)$
 $(\gminia\in \etabar_J^{-1}(\gminib))$
 on $\nbigx_P\setminus\nbigd_P(J^c)$
 in a standard way.
\item
 $\Gr_{\gminia}^{\nbigftilde}\bigl(
 \lefttop{J}\nbigg_{\gminib}\bigr)$
 is isomorphic to 
 $\nbigg_{P,\gminia|\nbigx_P\setminus\nbigd_P(J^c)}$,
 extending (\ref{eq;10.5.17.2}).
\end{itemize}

\subsection{Stokes data
associated to unramifiedly good lattice}
\label{subsection;10.5.16.30}

Let $(\nbige,\DD)$ be a 
meromorphic $\varrho$-flat bundle
with an unramifiedly good lattice $E$
on $(\nbigx,\nbigd)$.
Let $(V,\DD):=(\nbige,\DD)_{|\nbigx\setminus\nbigd}$.
According to Theorem \ref{thm;10.5.12.110},
we have an induced full pre-Stokes data
$\vecnbigftilde$ of $(V,\DD)$.
According to Proposition \ref{prop;10.5.12.210},
we have a graded extension
$\vecnbigg$ of the full pre-Stokes data.
\begin{lem}
$\vecnbigg$ satisfies
the compatibility condition
in Subsection {\rm\ref{subsection;10.5.16.21}}.
\end{lem}
\pf
The full Stokes filtrations $\nbigftilde^Q$
$(Q\in\pi^{-1}(P))$
are induced by filtrations of 
$\pi^{\ast}E$.
Hence, the induced filtration
$\nbigf^{Q,J}$ $(Q\in\pi^{-1}(P))$
are so.
We obtain $\Gr^{\vecnbigf^J}_{\gminib}(E)$
on $\nbigx_P$.
It is easy to see that
$\lefttop{J}\nbigg_{\gminib}$
is naturally isomorphic to
the restriction of
$\Gr^{\vecnbigf^J}_{\gminib}(E)$
to $\nbigx_P\setminus \nbigd_P(J^c)$.
We also have
$\Gr^{\vecnbigftilde}_{\gminia}
 \Gr^{\vecnbigf^J}_{\gminib}(E)
\simeq
 \nbigg_{P,\gminia}$.
Then, the claim of the lemma is clear.
\hfill\qed

\vspace{.1in}

Namely,
$(V,\DD)$ is equipped with a naturally induced Stokes data,
which is functorial,
according to the results
in Subsections \ref{subsection;10.5.12.15}
and \ref{subsection;10.5.16.20}.

\subsubsection{Complement
on splitting of Stokes filtrations}
\label{subsection;10.5.17.40}

We give a complement on a splitting of
Stokes filtration
of a good meromorphic $\varrho$-flat bundle.
We set $\nbigx:=\Delta^n\times\nbigk$
and $\nbigd=\bigcup_{i=1}^{\ell}\{z_i=0\}$.
We assume that the coordinate system is admissible.
We use the notation in 
Subsections \ref{subsection;10.5.3.7}
and \ref{subsection;10.5.12.1}.

Let $U\subset\pi^{-1}(P)$ be as in Condition
\ref{condition;10.5.17.30}.
Let $\nbigu$ be a small neighbourhood of $U$
in $\nbigxtilde(\nbigd)$,
which will be shrinked.
As in Proposition \ref{prop;10.5.17.30},
we have a $\DD$-flat filtration 
$\nbigftilde^{U}$ of $(\iota_{\ast}E)_{|U}$,
which is extended to a $\DD$-flat filtration
$\nbigftilde^{\nbigu}$ of 
$\iota_{\ast}E_{|\nbigu}$.
We can take a splitting
$E_{|\nbigu}=
 \bigoplus_{\gminia\in\nbigi_P}
 E_{\nbigu,\gminia}$
of $\nbigftilde^{\nbigu}$
with the following property:
\begin{itemize}
\item
For $\gminib\in\nbigi(\vecm(p))$,
we put
\[
 E^{\vecm(p)}_{\nbigu,\gminib}:=
 \bigoplus_{\gminia\in\etabar_{\vecm(p)}^{-1}(\gminib)}
 E_{\nbigu,\gminia}
\]
Then, the decomposition is $\DD_{\leq k(p)}$-flat.
\end{itemize}
Actually, we successively apply the third claim
of Proposition \ref{prop;10.5.28.10}.

Let $\vecv=(\vecv_{\gminia})$ be a frame of
$\Gr^{\vecnbigf}(E)$ compatible with the grading.
By a natural isomorphism
$\Gr^{\vecnbigf}(E)_{|\nbigu}
\simeq
 E_{|\nbigu}$ given by the above splitting,
we make a frame $\vecv_{\nbigu}$.

Let $E_{|\nbigu}=\bigoplus E'_{\nbigu,\gminia}$
be another decomposition
with the above property.
We obtain another frame $\vecv'_{\nbigu}$.
Let $C=(C_{\gminia,\gminib})$ be the matrix
determined by
$\vecv_{\nbigu}=\vecv'_{\nbigu}\, (I+C)$,
where $I$ denotes the identity matrix.
\begin{lem}
\label{lem;10.5.17.50}
We have $C_{\gminia,\gminib}=0$
unless $\gminia<_{\nbigu}\gminib$.
If $\gminia<_{\nbigu}\gminib$,
we have the estimate
\[
 C_{\gminia,\gminib}
 \exp\bigl(
 \varrho^{-1}(\gminia-\gminib)
 \bigr)
=O\Bigl(
 \prod_{i=1}^{k(\gminia,\gminib)}
 |z_i|^{-N}
 \Bigr)
\]
Here, $k(\gminia,\gminib)$
be determined by
$\ord(\gminia,\gminib)
 \in \seisuu_{<0}^{k(\gminia,\gminib)}
 \times \veczero$.
\end{lem}
\pf
It follows from that the induced morphism
$\Gr^{\vecnbigftilde}_{\gminib}(E)_{|\nbigu}
\lrarr \Gr^{\vecnbigftilde}_{\gminia}(E)_{|\nbigu}$
is $\DD_{\leq k(\gminia,\gminib)}$-flat.
\hfill\qed

\section{Riemann-Hilbert-Birkhoff correspondence}
\label{section;10.5.26.30}
\subsection{Statement}
\label{subsection;10.5.15.2}

Let $\MFL(\nbigx,\nbigd,\vecnbigi)$ be 
the category of 
unramifiedly good lattices $(E,\DD)$ of
meromorphic $\varrho$-flat bundles
on $(\nbigx,\nbigd)$
whose good system of irregular values
is contained in $\vecnbigi$,
i.e.,
$\Irr(\nabla,P)\subset\nbigi_P$
for any $P\in \nbigd$.
Let $\SDLcat(\nbigx,\nbigd,\vecnbigi)$
be the category of $\varrho$-flat bundle
with a full Stokes data over $\vecnbigi$.
As explained in 
Subsection \ref{subsection;10.5.16.30},
we have a functor
\[
\RHB:\MFL(\nbigx,\nbigd,\vecnbigi)
\lrarr\SDLcat(\nbigx,\nbigd,\vecnbigi).
\]
We will prove the following theorem in Subsections
\ref{subsection;10.5.18.20}--\ref{subsection;10.5.18.21}.
\begin{thm}
\label{thm;10.5.14.10}
The functor $\RHB$ is an equivalence.
\end{thm}

Let $\MF(\nbigx,\nbigd,\vecnbigi)$ be 
the category of 
unramifiedly good meromorphic $\varrho$-flat bundles
on $(\nbigx,\nbigd)$
whose good system of irregular values
is contained in  $\vecnbigi$.
Let $\SDcat(\nbigx,\nbigd,\vecnbigi)$
be the category of $\varrho$-flat bundle
with a full meromorphic Stokes data over $\vecnbigi$.

\begin{cor}
\label{cor;10.5.17.60}
The naturally defined functor
\[
 \RHB:\MF(\nbigx,\nbigd,\vecnbigi)
\lrarr\SDcat(\nbigx,\nbigd,\vecnbigi)
\]
is an equivalence.
\hfill\qed
\end{cor}

Let $G$ be a finite group
acting on $(\nbigx,\nbigd)$ over $\nbigk$.
Let $\vecnbigi$ be a good system of
irregular values
such that $\nbigi_P=g^{\ast}\nbigi_{g(P)}$
for any $g\in G$ and $P\in\nbigd$.
Let $(V,\DD)$ be a $\DD$-flat bundle
on $\nbigx\setminus\nbigd$ with
a $G$-action,
i.e.,
for each $g\in G$,
we are given an isomorphism
$g^{\ast}(V,\DD)\simeq(V,\DD)$
compatible with the group law.
Let $\vecnbigftilde$ be a full pre-Stokes
data of $(V,\DD)$.
For each $g\in G$,
we have the induced full pre-Stokes data
$g^{\ast}\vecnbigftilde$
of $(V,\DD)\simeq g^{\ast}(V,\DD)$.
The full pre-Stokes data is called
$G$-equivariant if
$g^{\ast}\vecnbigftilde
=\vecnbigftilde$ for any $g\in G$.
Similarly the $G$-equivariance of
a Stokes data is defined.
The category of $G$-equivariant
full Stokes data is denoted by
$\SDLcat(\nbigx,\nbigd,\vecnbigi)^G$.
Let $\MFL(\nbigx,\nbigd,\vecnbigi)^G$
denote the category of 
$G$-equivariant unramifiedly good lattices
of meromorphic $\varrho$-flat bundles
over $\vecnbigi$.
We use the symbols
$\SDcat(\nbigx,\nbigd,\vecnbigi)^G$
and $\MF(\nbigx,\nbigd,\vecnbigi)^G$
in similar meanings.
It is easy to deduce the following
as a corollary of Theorem \ref{thm;10.5.14.10}.
\begin{cor}
The functors $\RHB:
\MFL(\nbigx,\nbigd,\vecnbigi)^G\lrarr
\SDLcat(\nbigx,\nbigd,\vecnbigi)^G$
and 
$\RHB:
\MF(\nbigx,\nbigd,\vecnbigi)^G\lrarr
\SDcat(\nbigx,\nbigd,\vecnbigi)^G$
are equivalences.
\hfill\qed
\end{cor}

\subsubsection{Descent}

Let $\varphi:(\nbigx',\nbigd')\lrarr (\nbigx,\nbigd)$
be a ramified Galois covering with 
the Galois group $G$.
Let $\vecnbigi':=\varphi^{\ast}\vecnbigi$.
We have naturally defined descent functors
$\Des:\MFL(\nbigx',\nbigd',\vecnbigi')^G
\lrarr \MFL(\nbigx,\nbigd,\vecnbigi)$
and
$\Des:\SDLcat(\nbigx',\nbigd',\vecnbigi')^G
\lrarr
 \SDLcat(\nbigx,\nbigd,\vecnbigi)$.

\begin{prop}
We have a natural isomorphism
$\Des\circ\RHB\simeq \RHB\circ\Des$.
\end{prop}
\pf
Let $(E',\DD')\in\MFL(\nbigx',\nbigd',\vecnbigi')^G$.
We set $(E,\DD):=\Des(E',\DD')$,
$(V',\DD',(\vecnbigftilde',\vecnbigg')):=
 \RHB(E',\DD')$
and $(V,\DD,(\vecnbigftilde,\vecnbigg)):=
\RHB(E,\DD)$.
By using the characterization of
full Stokes filtration in Theorem \ref{thm;10.5.12.110},
we obtain that $\vecnbigftilde$
is the descent of $\vecnbigftilde'$.

Let $\pi:\nbigxtilde(\nbigd)\lrarr\nbigx$
and $\pi':\nbigxtilde'(\nbigd')\lrarr\nbigx'$
be blow up.
Let $P\in\nbigd$.
We take a small neighbourhood $\nbigx_P$ of $P$
in $\nbigx$.
We put $\nbigx_P':=\varphi^{-1}(\nbigx_P)$,
which is a neighbourhood of 
the finite set $\varphi^{-1}(P)$.
We have 
$\Gr^{\nbigftilde'}_{\gminia}\bigl(
 \pi^{\prime\ast}E'
 _{|\pi^{\prime\,-1}(\nbigx_P')}\bigr)$
on $\pi^{\prime\,-1}(\nbigx_P')$.
We can easily compare
$\nbigg_{P,\gminia}$
and the descent of
$\bigoplus_{P'\in\pi^{-1}(P)}\nbigg_{P',\gminia}$,
because both of them are induced by
$\Gr^{\nbigftilde'}_{\gminia}\bigl(\pi^{\prime\ast}E'
 _{|\pi^{\prime\,-1}(\nbigx_P')}\bigr)$.
\hfill\qed

\subsubsection{Classification of
unramifiedly good meromorphic flat bundles}

By setting 
(i) $\nbigk$ is a point, (ii) $\varrho=1$
in Corollary \ref{cor;10.5.17.60},
we obtain a classification of
unramifiedly good meromorphic flat bundles
in terms of full pre-Stokes data.
Note that we need only 
full pre-Stokes data,
because the unramifiedly good Deligne-Malgrange
lattices are canonically associated.

In the one dimensional case,
Malgrange {\rm\cite{malgrange3}}
showed the correspondence between
meromorphic flat bundles
and flat bundles with pre-Stokes data.
See also the work due to
Sibuya {\rm\cite{sibuya_book}}
on the classification of
meromorphic flat bundles on curves.

In the higher dimensional case,
such a classification was studied 
in {\rm\cite{majima}}
and {\rm\cite{sabbah4}},
from different viewpoints.
If we are interested in deformation
by a variation of irregular values,
the classification according to
pre-Stokes data is useful to the author.

\subsection{Fully faithfulness}
\label{subsection;10.5.18.20}

Let us show that $\RHB$ is fully faithful.
Let $(E_i,\DD_i)\in\MFL(\nbigx,\nbigd,\vecnbigi)$
$(i=1,2)$.
\begin{lem}
\label{lem;10.5.14.11}
The natural map 
\[
 \Hom\bigl((E_1,\DD_1),(E_2,\DD_2)\bigr)
\lrarr\Hom\bigl(
 \RHB(E_1,\DD_1),\RHB(E_2,\DD_2)
 \bigr) 
\]
is bijective.
\end{lem}
\pf
It is clearly injective.
Let us show the surjectivity.
Let 
\[
 F:\RHB(E_1,\DD_1)\lrarr\RHB(E_2,\DD_2)
\]
be a morphism.
We have only to show that
the underlying morphism
$F:(E_1,\DD_1)_{|\nbigx\setminus\nbigd}
\lrarr
 (E_2,\DD_2)_{|\nbigx\setminus\nbigd}$
is extended to a morphism $E_1\lrarr E_2$.
By Hartogs property,
we may assume that $\nbigd$ is smooth.
Let $Q\in\pi^{-1}(P)$.
We take a $\DD$-flat splitting
$E_{i|\nbigu_Q}\simeq
 \bigoplus E_{i,Q,\gminia}$
of the full Stokes filtration.
The flat morphism $F$ induces
$F_{\gminib,\gminia}:
 E_{i,\gminia|\nbigu_Q^{\ast}}\lrarr
 E_{i,\gminib|\nbigu_Q^{\ast}}$.
Because the full Stokes filtrations
are preserved by $F$,
we have $F_{\gminib,\gminia}=0$
unless $\gminia\geq^{\varrho}_Q\gminib$.
By construction $F_{\gminia,\gminia}$
is bounded.
Because $F_{\gminib,\gminia}$ is $\DD$-flat,
we have 
$F_{\gminib,\gminia}
 \exp\bigl(
 \varrho^{-1}(\gminia-\gminib)
 \bigr)$ is of polynomial order in $|z_1^{-1}|$.
Hence, we obtain that 
$F$ is extended to a morphism
$E_{1|\nbigu_Q}\lrarr E_{2|\nbigu_Q}$
by varying $Q$,
we obtain
$E_{1|\nbigxtilde(\nbigd)}\lrarr
 E_{2|\nbigxtilde(\nbigd)}$,
which induces
$E_1\lrarr E_2$.
\hfill\qed

\vspace{.1in}

It remains to show that $\RHB$
is essentially surjective.
Namely, we have to give a construction
of an unramifiedly good lattice
from a full Stokes data.
By Lemma \ref{lem;10.5.14.11},
we have only to give it locally.
We will give an inductive construction in levels.

\subsection{Pre-Stokes data and Stokes data
 in the level $\vecm$}
\label{subsection;07.12.11.30}

We use the setting in Section \ref{section;10.5.11.1}.
For a subset $\nbiga\subset\nbigx$,
we put $\nbiga^{\star}:=\nbiga\setminus\nbigd_z$.
Let $P\in\nbigd_{z,\kbar}$.
Let $\nbigx_P$ denote a small neighbourhood of $P$
in $\nbigx$.
It will be shrinked if it is necessary.

Let $V$ be a locally free 
$\nbigo_{\nbigx^{\star}}$-module
with a meromorphic $\varrho$-flat connection $\DD$
whose pole is contained in
$\nbigd^{\star}=\nbigd_Y^{\star}$.
Assume that $\DD$ is $\gminib$-logarithmic
for some $\gminib\in M(\nbigx,\nbigd)$.
Let $\iota:\nbigx^{\star}\subset
 \nbigxtilde(\nbigd_z)$
and $\pi:\nbigxtilde(\nbigd_z)\lrarr\nbigx$
be natural maps.
Let $\nbigi\subset M(\nbigx,\nbigd)$
be a weakly good set of irregular values
in the level $(\vecm,i(0))$
such that
$\bigl\{
 \gminia-\gminib\,\big|\,
 \gminia\in\nbigi
 \bigr\}$ is good
in the level $(\vecm,i(0))$.

\begin{df}
A pre-Stokes data $\vecnbigf$ of $(V,\DD)$
in the level $(\vecm,i(0))$ at $P$
is a system of $\DD$-flat filtrations
$\nbigf^Q$ $(Q\in \pi^{-1}(P))$
of germs $\iota_{\ast}(V)_Q$
indexed by $(\nbigi,\leq_Q^{\varrho})$,
satisfying the compatibility condition.
(We assume the existence of
a $\DD_z$-flat splitting,
instead of a $\DD$-flat splitting.)
We will often to distinguish 
``$P$'' and ``in the level $(\vecm,i(0))$''
if there is no risk of confusion.
\hfill\qed
\end{df}
\index{pre-Stokes data in the level $(\vecm,i(0))$}

Varying $Q\in\pi^{-1}(P)$,
for each $\gminia\in\nbigi$,
we obtain a locally free 
$\nbigo_{\nbigx_P^{\star}}$-module
$\Gr^{\vecnbigf}_{\gminia}(V_{|\nbigx_P^{\star}})$
with $\gminib$-logarithmic flat $\varrho$-connection
$\DD_{\gminia}$
on $(\nbigx_P^{\star},\nbigd_P^{\star})$.
(Note that $\gminia-\gminib$ is holomorphic
on $\nbigx_P^{\star}$.)

\begin{df}
\label{df;08.9.6.1}
A graded extension of $\vecnbigf$ is
a tuple 
of locally free $\nbigo_{\nbigx_P}$-modules
$E_{\gminia}$ $(\gminia\in \nbigi)$
such that 
(i) $E_{\gminia|\nbigx_P^{\star}}
 \simeq \Gr^{\vecnbigf}_{\gminia}
 (V_{\nbigx_P^{\star}})$,
(ii) $\ord(\DD_{\gminia}-d\gminia)\geq 
 \vecm(1):=\vecm+\vecdelta_{i(0)}$
for each $\gminia\in\nbigi$.
The tuple
$\SD=(\nbigi,\vecnbigf,\{E_{\gminia}\})$
is called a Stokes data in the level $(\vecm,i(0))$
at $P$.
\hfill\qed
\end{df}

Let $(V_p,\DD_p)$ $(p=1,2)$ be 
$\gminib$-logarithmic $\varrho$-flat bundles
$\DD_p$
with pre-Stokes data $\vecnbigf_p$ $(p=1,2)$
in the level $(\vecm,i(0))$ over $\nbigi_p$.
A morphism 
$\bigl(V_1,\DD_1,\vecnbigf_1\bigr)
\lrarr
 \bigl(V_2,\DD_2,\vecnbigf_2\bigr)$
is defined to be a flat morphism
which preserves Stokes filtrations 
at any $Q\in\pi^{-1}(P)$.
Note that we obtain a naturally induced morphism
of $\Gr^{\vecnbigf}_{\gminia}(V_1)
 \lrarr \Gr^{\vecnbigf}_{\gminia}(V_2)$.
If moreover they are equipped with graded extensions,
a morphism 
$\bigl(V_1,\DD_1,\SD_1\bigr)
\lrarr
 \bigl(V_2,\DD_2,\SD_2\bigr)$
is defined to be a morphism
$(V_1,\DD_1,\vecnbigf_1)\lrarr(V_2,\DD_2,\vecnbigf_2)$
such that the induced morphisms
$\Gr_{\gminia}^{\vecnbigf}(V_1)
 \lrarr \Gr^{\vecnbigf}_{\gminia}(V_2)$
are extended to 
$E_{1,\gminia}\lrarr E_{2,\gminia}$.

\subsection{Good lattice and Stokes data
in the level $\vecm$}
\label{subsection;07.12.13.13}

We continue to use the notation
in Subsection \ref{subsection;07.12.11.30}.
In terms of Stokes data,
we can summarize the results in 
Section \ref{section;10.5.11.1}
for weakly good lattice in the level $\vecm$.
\begin{prop}
Let $(E,\DD,\nbigi)$ be a weakly good lattice
of a meromorphic $\varrho$-flat bundle
in the level $(\vecm,i(0))$ on $(\nbigx,\nbigd)$.
We put $(V,\DD):=(E,\DD)_{|\nbigx-\nbigd_z}$.
Then, at each point $P\in\nbigd_{z,\kbar}$,
we have the Stokes data 
$\SD(E,\DD)$ in the level $(\vecm,i(0))$
for $(V,\DD)$ associated to $(E,\DD)$.
\index{Stokes data $\SD(E,\DD)$ in the level $(\vecm,i(0))$}
The correspondence is functorial, and it preserves
direct sum, tensor product and dual.
\hfill\qed
\end{prop}

We study the converse.
Namely,
for a given Stokes data in the level $(\vecm,i(0))$
at $P\in\nbigd_{z,\kbar}$,
we shall construct a good lattice 
in the level $(\vecm,i(0))$
on a neighbourhood of $P$.

\subsubsection{Construction}
\label{subsection;07.12.8.15}

\index{good lattice associated to Stokes data
in the level $(\vecm,i(0))$}

Let $(V,\DD)$ be a 
$\gminib$-logarithmic $\varrho$-flat bundle on
$(\nbigx^{\star},\nbigd^{\star})$
with a Stokes data at $P$.
We use a subscript ``$z$''
to indicate that we consider differentials relative to
$Y\times\nbigk$.
Because $\DD-d\gminib$ is logarithmic,
$\DD_z-d_z\gminib$ gives
a flat $\varrho$-connection
relative to $Y\times\nbigk$.

We take a holomorphic frame 
$\vecvbar_{\gminia}$ of $E_{\gminia}$
on $\nbigx_P$.
Let $R_{\gminia}$ be the connection one form
of $\DD_{\gminia}$ with respect to $\vecvbar_{\gminia}$,
i.e.,
$\DD_{\gminia}\vecvbar_{\gminia}
=\vecvbar_{\gminia}\,R_{\gminia}$.
Take $Q\in\pi^{-1}(P)$ 
and a small neighbourhood $\nbigu_Q$.
We put $\nbigu_Q^{\star}=\nbigu_Q\setminus
 \pi^{-1}(\nbigd_z)$.
We can take a $\DD_{z}$-flat splitting
$V_{|\nbigu_Q}=\bigoplus
 V_{\gminia,Q}$ of $\nbigf^Q$.
Let $\vecv_{\gminia,Q}$ be the lift of
$\vecvbar_{\gminia|\nbigu_Q^{\star}}$
to $V_{\gminia,Q}$.
Then, $\vecv_Q=
 (\vecv_{\gminia,Q})$ gives a frame of
$V_{|\nbigu_Q^{\star}}$.
Let $R_Q$ be determined by
$\DD\vecv_Q=
 \vecv_Q\,R_Q$.
We have the decomposition
$R_Q=(R_{\gminia,\gminic,Q})$
corresponding to $\vecv_Q=(\vecv_{\gminia,Q})$.

\begin{lem}
\mbox{{}}\label{lem;08.9.5.30}
We have the following:
\begin{itemize}
\item
$R_{\gminia,\gminic,Q}=0$
unless $\gminia\leq^{\varrho}_Q\gminic$.
\item
$R_{\gminia,\gminia,Q}=R_{\gminia}$.
\item
For $\gminia<_Q^{\varrho}\gminic$,
there exists $C>0$ such that
\[
  R_{\gminia,\gminic,Q}\,
 \exp\bigl(
 \varrho^{-1}(\gminia-\gminic)\bigr)=
 O\Bigl(
 \exp\bigl(C\, |\vecz^{\vecm(1)}|\bigr)
 \,|z_{i(0)}|^{-C}
 \Bigr).
\]
\end{itemize}
\end{lem}
\pf
Since the filtration $\nbigf^Q$ is $\DD$-flat,
we obtain the first claim.
The second claim is clear by construction.
Since the lift is taken for a $\DD_z$-flat splitting,
the $dz_i$-components
of $R_{\gminia,\gminic,Q}$ are $0$.
We have the expression
$ R_{\gminia,\gminic,Q}=
\sum_{j=1}^n R_{\gminia,\gminic,Q,j}
 \,d\zeta_j$.
Let $F_{\gminia,\gminic,j}:
 E_{\gminic|\nbigu_Q}\lrarr E_{\gminia|\nbigu_Q}$
be determined by
$F_{\gminia,\gminic,j}\vecvbar_{\gminic}
=\vecvbar_{\gminia}\,
 R_{\gminia,\gminic,Q,j}$.
They are $\DD_z$-flat.
Hence, we obtain
the desired estimate for
$R_{\gminia,\gminic,Q,j}$
by using Lemma \ref{lem;07.12.27.1}.
\hfill\qed

\begin{lem}
\label{lem;07.12.7.11}
Let $\vecv'_Q$ be a frame of $V_{|\nbigu_Q^{\star}}$
induced by another splitting 
$V_{|\nbigu_Q^{\star}}=
 \bigoplus V_{\gminia,Q}'$.
Let $C$ be determined by 
$\vecv_Q=\vecv_{Q}'\, (I+C)$.
We have
(i) $C_{|\Zhat}=0$,
(ii) $C_{\gminia,\gminic}=0$ unless 
$\gminia<_Q\gminic$,
(iii)
 $C_{\gminia,\gminic}\,
 \exp\bigl(\varrho^{-1}(\gminia-\gminic)\bigr)
=O\Bigl(
 \exp\bigl(C\,|\vecz^{\vecm(1)}|\bigr)
 \,|z_{i(0)}|^{-C}
 \Bigr)$ for some $C>0$.
\end{lem}
\pf
Two $\DD_{z}$-flat splittings induce
a $\DD_{z}$-flat map
$\Phi_{\gminia,\gminic}:
 E_{\gminic|\nbigu_Q}\lrarr E_{\gminia|\nbigu_Q}$
 for $\gminia<_Q\gminic$.
It can be shown that
$\Phi_{\gminia,\gminic}\,
 \exp\bigl(\varrho^{-1}(\gminia-\gminic)\bigr)
 =O\Bigl(
 \exp\bigl(C\, |\vecz^{\vecm(1)}|\bigr)
 |z_{i(0)}|^{-C}
  \Bigr)$ for some $C>0$
with respect to the frames
$\vecvbar_{\gminia}$ and $\vecvbar_{\gminic}$
by Lemma \ref{lem;07.12.27.1}.
Then, the claim of Lemma \ref{lem;07.12.7.11}
follows.
\hfill\qed

\vspace{.1in}

By using the natural isomorphisms of
holomorphic bundles
\[
 V_{|\nbigu_Q^{\star}}=
 \bigoplus_{\gminia}
 V_{\gminia,Q}
\simeq
 \bigoplus_{\gminia}
 E_{\gminia|\nbigu_Q^{\star}},
\]
we extend $V_{|\nbigu_Q^{\star}}$ to
a holomorphic vector bundle 
$\Vtilde_{\nbigu_Q}$ on $\nbigu_Q$.
By Lemma \ref{lem;08.9.5.30},
$\DD_{|\nbigu_Q^{\star}}$ is extended to 
a meromorphic flat $\varrho$-connection
$\DD_{Q}$ of $\Vtilde_{\nbigu_Q}$ on 
$\nbigu_{Q}$.
Moreover, we have the following isomorphism:
\begin{equation}
 \label{eq;08.9.5.31}
 (\Vtilde_{\nbigu_Q},\DD_{\nbigu_Q})_{|
 \widehat{\pi^{-1}(\nbigd_z)}\cap\nbigu_{Q}}
\simeq
 \bigoplus_{\gminia}
 (E_{\gminia},\DD_{\gminia})_{|
 \widehat{\pi^{-1}(\nbigd_z)}\cap\nbigu_Q}
\end{equation}
By Lemma \ref{lem;07.12.7.11},
$(\Vtilde_{\nbigu_Q},\DD_{\nbigu_Q})$ 
and the isomorphism (\ref{eq;08.9.5.31})
are independent of the choice of 
a $\DD_z$-flat splitting
$V=\bigoplus V_{\gminia,Q}$.
If $Q'\in\nbigu_{Q'}\subset\nbigu_Q$,
we have 
$(\Vtilde_{\nbigu_Q},\DD_{\nbigu_Q})_{|\nbigu_{Q'}}
=(\Vtilde_{\nbigu_{Q'}},\DD_{\nbigu_{Q'}})$.
By varying $Q$ and gluing them,
we obtain a holomorphic vector bundle
$\Vtilde_{\nbigxtilde_P(\nbigd_{P,z})}$
with a meromorphic flat $\varrho$-connection
$\DD_{\nbigxtilde_P(\nbigd_{P,z})}$.
Moreover, we have an isomorphism
\begin{equation}
 (\Vtilde_{\nbigxtilde_P(\nbigd_{P,z})},
 \DD_{\nbigxtilde_P(\nbigd_{P,z})})_{|
 \widehat{\pi^{-1}(\nbigd_{P,z})}}
\simeq
 \bigoplus_{\gminia}
 (E_{\gminia},\DD_{\gminia})_{|
 \widehat{\pi^{-1}(\nbigd_z')}}
\end{equation}
According to Proposition \ref{prop;08.9.4.102},
Corollary \ref{cor;08.9.5.1} and
Lemma \ref{lem;08.9.5.12},
there exists a holomorphic vector bundle $\Vtilde$
with a meromorphic flat $\varrho$-connection
$\DD$ on $(\nbigx_P,\nbigd_{P,z})$
such that
\begin{equation}
 \pi^{\ast}
 \bigl(\Vtilde,\DD \bigr)
\simeq
(\Vtilde_{\nbigxtilde_P(\nbigd_{P,z})},
 \DD_{\nbigxtilde_P(\nbigd_{P,z})}),
\quad
 \bigl(\Vtilde,\DD \bigr)_{|\nbigdhat_{P,z}}
\simeq
 \bigoplus_{\gminia}
 (E_{\gminia},\DD_{\gminia})
 _{|\nbigdhat_{P,z}}.
\end{equation}

\subsubsection{Functoriality}
\label{subsection;07.12.13.101}

Let $(V,\DD)$ be a $\gminib$-logarithmic
$\varrho$-flat bundle 
on $(\nbigx^{\star},\nbigd^{\star})$
with a Stokes data
in the level $(\vecm,i(0))$.
Then, $(V^{\lor},\DD^{\lor})$
is also equipped with an induced Stokes data
in the level $(\vecm,i(0))$.

\begin{lem}
\label{lem;07.12.7.13}
The associated extension of
$(V^{\lor},\DD^{\lor})$
is naturally isomorphic
to the dual of that of $(V,\DD)$.
\end{lem}
\pf
Let $V_{|\nbigu_Q^{\star}}=
 \bigoplus_{\gminia}V_{\gminia,Q}$
be a $\DD_z$-flat splitting.
It induces a $\DD_z$-flat splitting
$V^{\lor}_{|\nbigu_Q^{\star}}=
 \bigoplus V_{-\gminia,Q}^{\lor}$.
We extend $V^{\lor}_{|\nbigu_Q^{\star}}$
to $V^{\lor}_{\nbigu_Q}$ by using the splitting.
Then, we have a natural isomorphism
$(\Vtilde_{\nbigu_Q})^{\lor}
 \simeq
 \Vtilde^{\lor}_{\nbigu_Q}$.
Hence, we obtain
$(\Vtilde_{\nbigxtilde_P(\nbigd_{P,z})})^{\lor}
\simeq
 \Vtilde^{\lor}_{\nbigxtilde_P(\nbigd_{P,z})}$.
It induces the desired isomorphism.
\hfill\qed

\begin{lem}
Let $(V_i,\DD_i)$ $(i=1,2)$ be 
$\gminib$-logarithmic $\varrho$-flat bundles
on $(\nbigx^{\star},\nbigd^{\star})$
with Stokes data
in the level $(\vecm,i(0))$ at $P$.
\begin{itemize}
\item
If $\nbigi_1\otimes \nbigi_2$
is weakly good in the level $(\vecm,i(0))$,
we have the induced Stokes data of
$(V_1,\DD_1)\otimes(V_2,\DD_2)$
in the level $(\vecm,i(0))$ at $P$,
and the associated extension of
$(V_1,\DD_1)\otimes (V_2,\DD_2)$
is naturally isomorphic 
to $\Vtilde_1\otimes\Vtilde_2$.
\item
If $\nbigi_1\oplus\nbigi_2$
is weakly good in the level $(\vecm,i(0))$,
then we have the induced Stokes data of
$(V_1,\DD_1)\oplus(V_2,\DD_2)$
in the level $(\vecm,i(0))$ at $P$,
and the associated extension
is naturally isomorphic to
$\Vtilde_1\oplus \Vtilde_2$.
\hfill\qed
\end{itemize}
\end{lem}

Let $(V_i,\DD_i)$ $(i=1,2)$ be 
$\gminib$-logarithmic $\varrho$-flat bundles
on $\bigl(\nbigx^{\star},\nbigd^{\star}\bigr)$
equipped with Stokes data $\SD_i$
in the level $(\vecm,i(0))$ at $P$.
For simplicity,
we assume that $\nbigi_1\cup\nbigi_2$ is
also weakly good in the level $(\vecm,i(0))$.

\begin{lem}
\label{lem;07.12.13.110}
Let $F:(V_1,\DD_1,\SD_1)\lrarr (V_2,\DD_2,\SD_2)$
be a morphism.
We have the naturally induced morphism
$\Ftilde:\Vtilde_1\lrarr \Vtilde_2$
on $\nbigxtilde_P$.
\end{lem}
\pf
We take $\DD_z$-flat splittings
$V_{i|\nbigu_Q}=\bigoplus V_{i,\gminia,Q}$
of the filtrations $\nbigf^Q$.
Let $\vecvbar_{i,\gminia}$ be holomorphic frames
of $E_{i,\gminia}$.
Let $\vecv_{i,\gminia,Q}$ denote the lifts of
$\vecvbar_{i,\gminia,Q}$ to $V_{i,\gminia,Q}$.
They give frames $\vecv_{i,Q}$ of $V_{i|\nbigu_Q}$.

Let $A_Q$ be determined by 
$F(\vecv_{1,Q})=\vecv_{2,Q}\, A_Q$.
We have the decomposition
$A_Q=(A_{\gminia,\gminic,Q})$
corresponding to the decomposition
$\vecv_{i,Q}
=(\vecv_{i,\gminia,Q})$ $(i=1,2)$.
Since $F$ preserves the filtrations $\nbigf^Q$,
we have 
 $A_{\gminia,\gminic}=0$ unless 
 $\gminia\leq_Q\gminic$.
Since $A_{\gminia,\gminia,Q}$ satisfy
$F_{\gminia}\bigl(\vecvbar_{1,\gminia}\bigr)
=\vecvbar_{2,\gminia}\,A_{\gminia,\gminia,Q}$,
they are holomorphic on $\nbigu_Q$.
In the case $\gminia<_Q\gminic$,
 we obtain the estimate 
\[
 A_{\gminia,\gminic,Q}\,
 \exp\bigl(\varrho^{-1}(\gminia-\gminic)\bigr)
=O\Bigl(
 \exp\bigl(C|\vecz^{\vecm(1)}|\bigr)\,
 |z_{i(0)}|^{-C}
 \Bigr)
\]
for some positive constant $C$,
by using Lemma \ref{lem;07.12.27.1}.
Hence, $F$ is extended to a morphism
on $\pi^{-1}(\nbigx_P)$,
and the claim of Lemma \ref{lem;07.12.13.110}
follows.
\hfill\qed

\begin{cor}
If the restriction of $F_{|\nbigx-\nbigd}$ is
an isomorphism,
the induced morphism
$\Ftilde:\Vtilde_1(\ast \nbigd)\lrarr \Vtilde_2(\ast \nbigd)$
is an isomorphism.
\hfill\qed
\end{cor}

\subsection{Proof of Theorem \ref{thm;10.5.14.10}}
\label{subsection;10.5.18.21}

Let us show that $\RHB$ is essentially surjective.
We have only to argue it locally.
We put
$\nbigx:=\Delta^n\times\nbigk$
and $\nbigd:=\bigcup_{i=1}^{\ell}\{z_i=0\}$.
Let $\nbigi\subset M(\nbigx,\nbigd)/H(\nbigx)$
be a good set of irregular values.
We assume that the coordinate system is admissible
for $\nbigi$,
and we take an auxiliary sequence
$\vecm(p)$ for $\nbigi$.
We shall construct 
an unramifiedly good meromorphic flat bundle
around $P\in\nbigd_{\ellsitabar}$,
from a flat bundle with a full Stokes data.

\subsubsection{Graded bundles
associated to full pre-Stokes data}
\label{subsection;10.5.15.3}

We have a refinement of the construction
in Subsection \ref{subsection;10.5.17.5}.
We use the notation there.
Let $(V,\DD)$ be a $\varrho$-flat bundle
on $\nbigx\setminus\nbigd$
with a full pre-Stokes data $\vecnbigftilde$
over a good set of irregular values $\nbigi$.
For each $Q\in \pi^{-1}(P)$,
we obtain the induced filtration
$\nbigf^{Q\,\vecm(p)}$ of 
$\iota_{\ast}(V)_Q$
indexed by 
$\bigl(\nbigi(\vecm(p)),
 \leq^{\varrho}_Q \bigr)$.
(See Subsection \ref{subsection;10.6.11.1}
for $\nbigi(\vecm(p))$.)
On the associated graded sheaf
$\Gr_{\gminib}^{\vecm(p)}(\iota_{\ast}(V)_Q)$,
we have the induced filtration
$\nbigftilde^Q$
indexed by
$\bigl(\etabar_{\vecm(p)}^{-1}(\gminib),
 \leq^{\varrho}_Q \bigr)$.
By varying $Q\in\pi^{-1}(P)$,
we obtain 
$\Gr^{\vecm(p)}(
 \iota_{\ast}(V,\DD)_{|\pi^{-1}(\nbigx_P)})$,
and a $\DD$-flat bundle
$\Gr^{\vecm(p)}(V_{|\nbigx_P^{\ast}})$
on $\nbigx_P^{\ast}$.
Let $\pi_{\vecm(p)}:
 \nbigxtilde(\nbigd(\kbar(p)))
\lrarr \nbigx$ be the real blow up,
and let $\iota_{\vecm(p)}:
 \nbigx_P^{\ast}\subset
 \nbigxtilde(\nbigd(\kbar(p)))$.
For each $Q_1\in\pi_{\vecm(p)}^{-1}(P)$,
we have the $\DD$-flat induced filtration
$\nbigftilde^{Q_1}$ of
$\iota_{\vecm(p)\ast}
 \Gr^{\vecm(p)}(V_{|\nbigx_P^{\ast}})_{Q_1}$.

For a given $J\subset\ellsitabar$,
let $\vecm(p_J)$ be determined by
$m_i(p_J+1)=0$ for any $i\in J$
and 
$m_i(p_J)\neq 0$ for some $i\in J$.
Note that the the image of
$\nbigi(\vecm(p_J))$
by $M(\nbigx,\nbigd)/H(\nbigx)\lrarr
 M(\nbigx,\nbigd)/M(\nbigx,\nbigd(J^c))$
coincides with $\nbigi^J$.
We have $\vecnbigf^J=\vecnbigf^{\vecm(p_J)}$.

\subsubsection{Construction}

Let $(V,\DD,\vecnbigftilde,\vecnbigg)
 \in\SDcat(\nbigx,\nbigd,\vecnbigi)$.
Take $P\in\nbigd_{\ellsitabar}$.
In the following,
we will replace $\nbigx$
with a small neighbourhood of $P$
if it is necessary.
We have $\Gr^{\vecm(p)}_{\gminib}(V)$ 
on $\nbigx\setminus\nbigd$
for $\gminib\in\nbigi(\vecm(p))$.
We shall construct 
a $\gminib$-logarithmic extension
$E^{\vecm(p)}_{\gminib}$ of 
$\Gr^{\vecm(p)}_{\gminib}V$
on $\nbigx\setminus\nbigd(\kbar(p+1))$.

For $J(k):=\{k+1,\ldots,\ell\}$,
we take $\vecm(p_{J(k)})$
as in Subsection \ref{subsection;10.5.15.3}.
Then, for $\gminib\in\nbigi(\vecm(p_{J(k)}))$,
we have a locally free 
$\nbigo_{\nbigx\setminus\nbigd(\kbar)}$-module
$\nbigg_{\gminib}^{J(k)}$.
Let $\pi_k:\nbigxtilde(\nbigd(\kbar))\lrarr\nbigx$
and $\iota_k:\nbigx\setminus\nbigd(\kbar)\lrarr
 \nbigxtilde(\nbigd(\kbar))$
be natural maps.
For each $Q\in\pi_{k}^{-1}(P)$,
we have the $\DD$-flat filtration
$\nbigftilde^Q$ of
$\iota_{k\ast}(\nbigg^{J(k)}_{\gminib})_Q$
indexed by 
$\bigl(\etabar_{\vecm(p_{J(k)})}^{-1}(\gminib),
 \leq_Q^{\varrho}\bigr)$.
(See the argument in Subsection 
\ref{subsection;10.5.14.1}.)
We have the induced filtrations
$\nbigf^{Q\,\vecm(p)}$ 
for any $p\geq p_{J(k)}$.
Since they satisfy a compatibility condition,
we obtain 
$\Gr^{\vecm(p)}_{\gminia}\bigl(
 \nbigg^{J(k)}_{\gminib}\bigr)$
on $\nbigx\setminus\nbigd(\kbar)$,
for any $p\geq p_{J(k)}$
and $\gminia\in
 \etabar_{\vecm(p_{J(k)})}^{-1}(\gminib)$.
By using the compatibility condition
of lattices in Stokes data,
we obtain
\begin{equation}
\label{eq;10.5.16.41}
 \Gr^{\vecm(p_{J(k-1)})}_{\gminia}
 \bigl(
 \nbigg^{J(k)}_{\gminib}
 \bigr)
\simeq
 \nbigg^{J(k-1)}_{\gminia|
 \nbigx\setminus\nbigd(\kbar)} 
\end{equation}
for any $\gminia\in \nbigi(\vecm(p_{J(k-1)}))$.

For $p$,
we take $k$ such that
$p_{J(k)}\leq p<p_{J(k+1)}$.
For any $\gminic\in\nbigi(\vecm(p))$,
we put
$E_{\gminic}^p:=
 \Gr^{\vecm(p)}_{\gminic}
 \bigl(
 \nbigg^{J(k)}_{\gminib}
 \bigr)$,
which is equipped with the induced filtrations
$\nbigf^{Q\,\vecm(p')}$
for any $Q\in\pi^{-1}(P)$ and $p'\geq p$.
The system of the filtrations
satisfy a compatibility condition.
By (\ref{eq;10.5.16.41}),
we have a natural isomorphism
$\Gr^{\vecm(q)}_{\gminic}E^p_{\gminib}
\simeq
 E^q_{\gminic|\nbigx\setminus
 \nbigd(\kbar(p-1))}$.

Then, by a successive use of the construction in
Subsection \ref{subsection;07.12.8.15},
we can construct a locally free
$\nbigo_{\nbigx}$-module $E_P$
with an isomorphism
$E_{|\nbigx\setminus\nbigd}\simeq
 V_{|\nbigx\setminus\nbigd}$
and 
$(E,\DD)_{|\Phat}\simeq
 \bigoplus_{\gminia\in\nbigi_P}
 \nbigg_{P,\gminia|\Phat}$.
In particular,
$E$ is an unramifiedly good lattice.
By construction,
for each $Q\in \pi^{-1}(P)$,
the full Stokes filtration 
$\nbigftilde^Q$ of $E$ at $Q$
is the same as that in 
the given Stokes data $\SD$.
It implies that, 
for any $Q\in \pi^{-1}(\nbigd)$,
the full Stokes filtration of $E$ at $Q$
is the same as that in the given Stokes data,
according to Lemma \ref{lem;10.5.16.40}.
(Note that we have shrinked $\nbigx$
around $P$.)

\vspace{.1in}

Let $J\subset \ellsitabar$.
Take $P'\in\nbigd_J$ with $I(P')=J$.
It remains to show that
we have a natural isomorphism
\begin{equation}
\label{eq;10.6.30.2}
\Gr^{\vecnbigftilde}(E_{P'})\simeq
\bigoplus_{\gminia\in\nbigi_{P'}}
 \nbigg_{P',\gminia} 
\end{equation}
on a small neighbourhood $\nbigx_{P'}$ of $P'$.
We have
$\Gr^{\vecnbigftilde}(E_{P'})\simeq
\Gr^{\nbigf^{J}}(E)_{|\nbigx_{P'}}$.
If $J=J(k)$ for some $k$,
(\ref{eq;10.6.30.2}) is clear 
by our construction of $E$.
In the general case,
we put $k:=\min J$.
Because we naturally have
$\nbigg^{J(k)}_{\gminia|\nbigx_{P'}}
\simeq
 \nbigg_{P',\gminia}$
and 
$\Gr^{\nbigf^J}(E)
\simeq
 \Gr^{\nbigf^{J(k)}}(E)$,
we have the desired isomorphism.
Thus the proof of Theorem \ref{thm;10.5.14.10}
is finished.
\hfill\qed

\section{Extension of Stokes data}
\label{section;10.5.26.31}
\subsection{Statement}
\label{subsection;10.5.19.2}

Let $\nbigx\lrarr \nbigb\lrarr \nbigk$
be smooth fibrations of complex manifolds.
Let $\nbigd$ be a normal crossing hypersurface
of $\nbigx$ such that each intersection of
irreducible components is smooth over $\nbigb$.
For simplicity,
we assume the following:
\begin{itemize}
\item
$a:\nbigb\lrarr\nbigk$ is equipped
with a section $b:\nbigk\lrarr\nbigb$,
and each fiber of $a$ is simply connected.
\item
We put $\nbigx^b:=\nbigx\times_{\nbigb}b$
and $\nbigd^b:=\nbigd\times_{\nbigb}b$.
Then,
$(\nbigx,\nbigd)$ is topologically
a product of $(\nbigx^b,\nbigd^b)$
and $\nbigb$.
\end{itemize}
For example, we would like to consider
the case $\nbigb=\nbigk\times B$
and $(\nbigx,\nbigd)=(\nbigx^b,\nbigd^b)\times \nbigb$
as complex manifolds.

Let $\varrho$ be a nowhere vanishing 
holomorphic function on $\nbigk$.
Let $\vecnbigi$ be a good system of irregular values
on $(\nbigx,\nbigd)$.
Its restriction to $\nbigx^b$ is denoted by
$\vecnbigi^b$.

\begin{thm}
\label{thm;10.5.16.1}
The restriction
$\SDLcat(\nbigx,\nbigd,\vecnbigi)
\lrarr
 \SDLcat(\nbigx^b,\nbigd^b,\vecnbigi^b)$
is equivalent.
\end{thm}

Note that, under the assumption,
the restriction induces an equivalence
between $\varrho$-flat bundles over
$\nbigx\setminus\nbigd$
and $\nbigx^b\setminus\nbigd^b$.
For a $\varrho$-flat bundle $(V,\DD)$
on $\nbigx\setminus\nbigd$,
let $(V^b,\DD^b)$ denote its restriction to
$\nbigx^b\setminus\nbigd^b$.
Theorem \ref{thm;10.5.16.1} says
that a Stokes data of
$(V^b,\DD^b)$ over $\vecnbigi^b$
is uniquely extended to
a Stokes data of $(V,\DD)$
in a functorial way.
By using the uniqueness,
we easily obtain the following:
\begin{cor}
The above extension is functorial
with respect to dual, tensor product,
direct sum 
in an obvious sense.
\end{cor}

\begin{cor}
\label{cor;10.5.17.11}
The restriction
$ \MFL(\nbigx,\nbigd,\vecnbigi)
 \lrarr
 \MFL(\nbigx^b,\nbigd^b,\vecnbigi^b)$
is equivalent.
\hfill\qed
\end{cor}

\subsubsection{Variants}

We immediately obtain
the meromorphic variant.
\begin{cor}
The restrictions
\[
  \SDcat(\nbigx,\nbigd,\vecnbigi)
\lrarr
 \SDcat(\nbigx^b,\nbigd^b,\vecnbigi^b),
\quad
 \MF(\nbigx,\nbigd,\vecnbigi)
 \lrarr
 \MF(\nbigx^b,\nbigd^b,\vecnbigi^b)
\]
are equivalent.
\hfill\qed
\end{cor}

Let $G$ be a finite group acting on 
$(\nbigx,\nbigd)$ over $\nbigb$.
Let $\nbigc(\nbigx\setminus\nbigd)^G$
be the category of $G$-equivariant
$\varrho$-flat bundles on $\nbigx\setminus\nbigd$.
By using the uniqueness,
we easily obtain the following.
\begin{cor}
\label{cor;10.5.17.10}
The restrictions
\[
 \SDLcat(\nbigx,\nbigd,\vecnbigi)^G
\lrarr
 \SDLcat(\nbigx^b,\nbigd^b,\vecnbigi^b)^G,
\quad
 \MFL(\nbigx,\nbigd,\vecnbigi)^G
 \lrarr
 \MFL(\nbigx^b,\nbigd^b,\vecnbigi^b)^G
\]
are equivalent.
We have the meromorphic variant.
\hfill\qed
\end{cor}

\subsection{Extension of full pre-Stokes data}

\label{subsection;10.5.15.20}

Let us consider the claim for full pre-Stokes data.

\begin{prop}
\mbox{{}}\label{prop;10.5.15.30}
\begin{itemize}
\item
Let $(V,\DD)$ be a $\varrho$-flat bundle
on $\nbigx\setminus\nbigd$.
If we are given a full pre-Stokes data
of $(V^b,\DD^b)$ over 
$(\nbigx^b,\nbigd^b,\vecnbigi^b)$,
it is uniquely extended to
a full pre-Stokes data of $(V,\DD)$
over $(\nbigx,\nbigd,\vecnbigi)$.
\item
Let $(V_i,\DD_i)$ $(i=1,2)$ be $\varrho$-flat bundles
on $\nbigx\setminus\nbigd$
equipped with full pre-Stokes structures
$\vecnbigftilde_i$
over $(\nbigx,\nbigd,\vecnbigi)$.
Let $F:(V_1,\DD_1)\lrarr (V_2,\DD_2)$
be a morphism.
If its restriction $F^b$ preserves
full Stokes filtrations,
then $F$ does so.
\end{itemize}
\end{prop}

We only give the proof of the first claim.
The second claim can be shown similarly.

\subsubsection{}
\label{subsection;10.5.16.92}

We put
$\nbigx_0:=\Delta^n\times Y\times\nbigk$
and $\nbigd_0:=\bigcup_{i=1}^n\{z_i=0\}$.
We consider the case
$\nbigx:=B\times\nbigx_0$
and $\nbigd:=B\times\nbigd_0$
for some simply connected complex manifold $B$.
Let $\pi:\nbigxtilde(\nbigd)\lrarr\nbigx$
and $\pi_0:\nbigxtilde_0(\nbigd_0)\lrarr\nbigx_0$
be the real blow up.
Let $\iota:\nbigx\setminus\nbigd\subset
\nbigxtilde(\nbigd)$
and $\iota_0:\nbigx_0\setminus\nbigd_0
\subset\nbigxtilde_0(\nbigd_0)$
be natural maps.
Let $(V,\DD)$ be a $\varrho$-flat bundle
on $\nbigx\setminus\nbigd$.
Let $P$ be any point of $\nbigd_0$.
Note that we have natural identifications
$\nbigi_{(b,P)}\simeq\nbigi_{(b',P)}$
for $b,b'\in B$.
In the following,
for $Q\in \pi_0^{-1}(P)$,
let $U_Q$ denote a small neighbourhood
of $Q$ in $\pi_0^{-1}(P)$.

\subsubsection{}
Let $b_0\in B$.
Assume that we are given
a full pre-Stokes structure of
$\iota_{\ast}(V,\DD)_{|\pi^{-1}(b_0,P)}$.
According to Lemma \ref{lem;10.5.16.90},
there exists a neighbourhood $B_0$ of $b_0$
such that 
$\iota_{\ast}(V)_{|B_0\times\pi_0^{-1}(P)}$
has a full pre-Stokes structure
whose restriction to $b\times\pi_0^{-1}(P)$
is equal to the given one.
We also obtain that such a 
full pre-Stokes structure is uniquely determined.

\subsubsection{}

Let $b_1\in B$
with a neighbourhood $B_1$.
Assume that there is 
an open subset $B_1'\subset B_1$
such that 
(i) $b_1$ is contained in
the closure of $B'_1$ in $B_1$,
(ii) a full pre-Stokes structure of
$\iota_{\ast}(V)_{|B_1'\times \pi_0^{-1}(P)}$
is given.
Let $Q\in\pi_0^{-1}(P)$.
We can take a small neighbourhood
$U(b_1,Q)=B_2\times U_Q
 \subset B_1\times\pi_0^{-1}(P)$
such that
$\leq_{(b_1,Q)}=\leq_{U(b_1,Q)}$.

\begin{lem}
If $U(b_1,Q)$ is sufficiently small,
we have
$\leq_{(b_1,Q)}=
 \leq_{U'(b_1,Q)}$,
where
$U'(b_1,Q):=U_1(b_1,Q)\times_{B}B_1'$.
Moreover,
for any $b\in B_2$,
we have
$\leq_{(b_1,Q)}=
\leq_{U^b(b_1,Q)}$,
where $U^b(b_1,Q)=
 b\times_{B}U(b_1,Q)$.
\end{lem}
\pf
For any fixed pair
$\gminia,\gminib\in \nbigi_{(b,P)}$,
after appropriate coordinate change,
$F_{\gminia,\gminib}
=-\Re\bigl(|\vecz^{-\vecm}|\,\vecz^{\vecm}\bigr)$.
Hence,
$F_{\gminia,\gminib}^{-1}(0)\cap
 \pi_0^{-1}(B_2\times P)
 \lrarr B_2\times P$
is a smooth fibration,
if $B_2$ is sufficiently small.
Then, the claim is clear.
\hfill\qed

\vspace{.1in}

Take $b_2\in B_2\cap B_1'$.
By using Proposition \ref{prop;10.5.15.10},
we have a $\DD$-flat filtration of
$\iota_{\ast}(V)_{| b_2\times U_Q}$.
It is extended to a $\DD$-flat filtration
on $U(b_1,Q)$.
It is independent of the choice of $b_2$.
By varying $Q\in\pi_0^{-1}(P)$,
we obtain a neighbourhood $B_4$ of $b_1$
and full pre-Stokes structure of
$\iota_{\ast}(V)_{|B_4\times \pi_0^{-1}(P)}$.
Note that the uniqueness is also obtained.

\subsubsection{}

Assume that we are given
a full pre-Stokes structure of
$\iota_{\ast}(V)_{|
 b_3\times\pi_0^{-1}(\nbigd_0)}$.
Take any $b_4\in B$ and $P\in\nbigd_0$.
We take a path $\gamma$ connecting $b_3$ and $b_4$.
By a continuity method,
we can show that
there exists a neighbourhood $\nbigv$ of $\gamma$
and a unique full pre-Stokes structure
of $\iota_{\ast}(V)_{|\nbigv\times \pi_0^{-1}(P)}$
whose restriction to $b\times\pi_0^{-1}(P)$
is equal to the given one.
In particular, we obtain a unique full pre-Stokes structure
of $\iota_{\ast}(V)_{|b_4\times\pi_0^{-1}(P)}$.
It is easy to show that
the filtrations are independent of the choice of a path,
and that the compatibility condition is satisfied.
Thus, the first claim of Proposition \ref{prop;10.5.15.30}
is proved under the setting of 
Subsection \ref{subsection;10.5.16.92}.

\subsubsection{}
Let us return to the setting in
Subsection \ref{subsection;10.5.19.2}.
Let $\nbigd=\bigcup_{i\in\Lambda}\nbigd_i$
be the irreducible decomposition.
For $I\subset\Lambda$,
we put $\nbigd_I^{\ast}:=
 \bigcap_{i\in I}\nbigd_i\setminus
 \bigcup_{i\not\in I}\nbigd_i$.
Let $\nbigx\stackrel{c}\lrarr
\nbigb\stackrel{a}\lrarr\nbigk$.
Take $P\in\nbigd$.
Let $I(P):=\bigl\{i\in\Lambda\,\big|\,
 P\in\nbigd_i\bigr\}$.
We put $y:=a\circ c(P)$.
We can take a path $\gamma$
in $(a\circ c)^{-1}(y)\cap \nbigd_{I(P)}^{\ast}$
connecting $P$ and
$c^{-1}(b(y))\cap\nbigd_{I(P)}^{\ast}$.
Note that such a path is unique up to homotopy.
By a continuity method,
we can show that
there exists a neighbourhood $\nbigv$ of $\gamma$
and a unique full pre-Stokes structure
of $\iota_{\ast}(V)_{|\pi^{-1}(\gamma)}$
whose restriction to the intersection with 
$\pi^{-1}(\nbigd^b)$ is equal to the given one.
In particular, we obtain a unique full pre-Stokes structure
of $\iota_{\ast}(V)_{|\pi_0(P)}$.
It is easy to show that
the filtrations are independent of the choice of a path,
and that the compatibility condition is satisfied.
Thus, the first claim of Proposition \ref{prop;10.5.15.30}
is proved.

\subsection{Extension of $\gminib$-logarithmic bundle}

Let $(V,\DD)$ be a $\DD$-flat bundle
on $\nbigx\setminus\nbigd$.
Let $\gminib\in M(\nbigx,\nbigd)$.
\begin{lem}
\label{lem;10.5.16.2}
Assume that the restriction
$(V^b,\DD^b)$
is extended to a $\gminib^b$-logarithmic
$\varrho$-flat bundle $E_0$.
Then,
$(V,\DD)$ is uniquely extended to
a $\gminib$-logarithmic $\varrho$-flat bundle $E$
such that $E_{|\nbigx\times b}=E_0$.
\end{lem}
\pf
We have only to consider the case $\gminib=0$.
By the uniqueness,
the claim is a local property.
Hence, we may assume
$(\nbigx,\nbigd)=B\times(\nbigx_0,\nbigd_0)$.
Then the existence is clear,
because $(V,\DD)$ 
is isomorphic to the pull back of
$(V,\DD)_{|b\times(\nbigx_0\setminus\nbigd_0)}$.
Let us show the uniqueness.
Let $E$ be such an extension.
If we restrict $\DD$ to the $B$-direction,
there is no pole,
i.e.,
we obtain a $\varrho$-flat connection
relative to $\nbigx_0$
without any pole.
Then, the claim is clear.
\hfill\qed

\vspace{.1in}

We mention a consequence of this lemma
for the proof of Theorem \ref{thm;10.5.16.1}.
Let $(\vecnbigf^b,\vecnbigg^b)$
be a full Stokes data of $(V^b,\DD^b)$.
According to Proposition \ref{prop;10.5.15.30},
$\vecnbigf^b$ is uniquely extended
to a full pre-Stokes data of $(V,\DD)$.
By Lemma \ref{lem;10.5.16.2},
we obtain that
$\vecnbigg^b$
is also uniquely extended to
a graded extension of $\vecnbigf^b$.
It remains to show
the compatibility condition for
$\vecnbigg$.

\subsection{Compatibility}

We give a preparation.
We put 
$\nbigx_0:=\Delta^n\times Y\times\nbigk$
and $\nbigd_{0,z}:=\bigcup_{i=1}^n\{z_i=0\}$.
Let $\nbigd_{0,Y}$ be a hypersurface
obtained as the pull back of a normal crossing
hypersurface of $Y$.
We put $\nbigd_0:=\nbigd_{0,z}\cup\nbigd_{0,Y}$.
Let $B$ be a complex manifold.
We put $\nbigx:=B\times\nbigx_0$,
$\nbigd_z:=B\times\nbigd_{0,z}$,
and $\nbigd:=B\times\nbigd_{0}$.
Let $\nbigi$ be a good set of irregular values
on $(\nbigx,\nbigd_z)$.

We have the following blow up:
\[
 \pi_{0,z}:\nbigxtilde_0(\nbigd_{0,z})
\lrarr\nbigxtilde_0
\quad\quad
 \pi_{0}:\nbigxtilde_0(\nbigd_0)
\lrarr\nbigxtilde_0
\]
\[
  \pi_{z}:\nbigxtilde(\nbigd_{z})
\lrarr\nbigxtilde
\quad\quad
 \pi:\nbigxtilde(\nbigd)
\lrarr\nbigxtilde
\]
We have the following inclusions:
\[
 \iota_0:\nbigx_0\setminus\nbigd_0
\lrarr\nbigxtilde_0(\nbigd_{0,z})
\quad\quad
 \iota_0':\nbigx_0\setminus\nbigd_{0,z}
\lrarr\nbigxtilde_0(\nbigd_{0,z})
\]
\[
 \iota:\nbigx\setminus\nbigd
\lrarr\nbigxtilde(\nbigd_z)
\quad\quad
 \iota':\nbigx\setminus\nbigd_z
\lrarr\nbigxtilde(\nbigd_z) 
\]

Let $(E,\DD)$ be a logarithmic
$\varrho$-flat bundle
on $(\nbigx\setminus\nbigd_z,\nbigd\setminus
 \nbigd_z)$.
We put $(V,\DD):=(E,\DD)_{|\nbigx\setminus\nbigd}$.
We have a natural inclusion
$\iota'_{\ast}E\subset \iota_{\ast}V$.

Let $\vecnbigftilde$ be a
full pre-Stokes structure of $(V,\DD)$ 
over $\nbigi$.
Because $\nbigi$ is contained in
$M(\nbigx,\nbigd_z)/H(\nbigx)$,
we have the induced filtrations
$\nbigftilde^{b,Q}$ 
of the germs $\iota_{\ast}(V)_{b,Q}$
for any $(b,Q)\in \nbigxtilde(\nbigd_z)$.

Let $b_0\in B$.
We have the restriction
$(V^{b_0},\DD^{b_0})$ 
on $\{b_0\}\times(\nbigx_0\setminus\nbigd_0)$
and $(E^{b_0},\DD^{b_0})$
on $\{b_0\}\times(\nbigx_0\setminus\nbigd_{0,z})$.
We have the natural inclusion
$\iota'_{0\ast}E^{b_0}
\subset
 \iota_{0\ast}V^{b_0}$.
We have the filtrations
$\nbigftilde^Q$ of $(\iota_{0\ast}V)_Q$
for any $Q\in\pi_{0,z}^{-1}(\nbigd_{0,z})$.

\begin{lem}
\label{lem;10.5.16.100}
Assume the following:
\begin{itemize}
\item
Take any $Q\in
 \pi_{0,z}^{-1}\bigl(\nbigd_{0,z}\bigr)$,
then $\nbigftilde^Q$ of
$\iota_{0\ast}(V^{b_0})_{Q}$
is induced by a filtration of
$\iota'_{0\ast}(E^{b_0})_Q$.
\end{itemize}
Then, for any $b\in B$
and $Q\in\pi_0^{-1}(\nbigd_z)$,
$\nbigftilde^{(b,Q)}$
of $\iota_{\ast}(V)_{(b,Q)}$
is induced by a filtration of
$\iota'_{\ast}(E)_{(b,Q)}$.

Note that the assumption and the claim
are trivial if $Q$ is not contained
in the inverse image of 
$\nbigd_{z}\cap\nbigd_Y$.
\end{lem}
\pf
By the assumption,
we have the filtration
$\nbigftilde^Q_1$
of the germ  $(\iota'_{0\ast}E^{b_0})_Q$
for each
$Q\in \pi_{0,z}^{-1}(\nbigd_{0,z})$.
The system is denoted by
$\vecnbigftilde_1$.
It satisfies the compatibility condition.
Because $(E,\DD)$ is logarithmic,
by forgetting the differentials in the $Y$-direction,
we obtain a flat $\varrho$-connection $\DD_1$
relative to $Y\times\nbigk$.
The system $\vecnbigftilde_1$
gives a full pre-Stokes structure of 
$\iota'_{0\ast}(E^{b_0},\DD_1^{b_0})$.
According to Proposition \ref{prop;10.5.15.30},
it is uniquely extended to
a full pre-Stokes structure 
$\vecnbigftilde_2$ of 
$\iota'_{\ast}(E,\DD_1)$.

The restriction of $\vecnbigftilde$
to $\nbigxtilde(\nbigd_z)
 \setminus\pi_z^{-1}(\nbigd_Y)$
gives a full pre-Stokes structure
of $\iota_{\ast}(V,\DD_1)_{|\nbigxtilde(\nbigd_z)
 \setminus\pi_z^{-1}(\nbigd_Y)}$.
By the uniqueness,
we obtain the coincidence of
the restrictions of
$\vecnbigftilde$ and $\vecnbigftilde_2$
to $\nbigxtilde(\nbigd_z)\setminus
 \pi_z^{-1}(\nbigd_Y)$.
Then, we can deduce that
$\vecnbigftilde$ is induced by
$\vecnbigftilde_2$
on $\nbigxtilde(\nbigd_z)$.
\hfill\qed

\subsection{End of the proof of 
Theorem \ref{thm;10.5.16.1}}

Let us finish the proof of Theorem
\ref{thm;10.5.16.1}.
It remains to check the compatibility
condition for $\vecnbigg$.
We use the setting in Subsection
\ref{subsection;10.5.16.92}.
Let $P\in\nbigd_0$.
In the following,
$\nbigx_P$ denotes a small neighbourhood
of $B\times\{P\}$.
We have a natural bijection
$\nbigi_{(b,P)}\simeq
 \nbigi_{(b',P)}$ for any $b,b'\in B$.
We will identify them naturally,
and denoted by $\nbigitilde_P$.
From a full pre-Stokes structure of $V$,
we have 
$\Gr^{\vecnbigf^J}_{\gminib}(V)$
on $\nbigx_P^{\ast}$ for
$\gminib\in\nbigitilde_P^J$.
The filtrations $\nbigftilde^Q$ 
$(Q\in\pi^{-1}(B\times P))$
induce a full pre-Stokes structure 
$\vecnbigftilde$
of $\Gr^{\vecnbigf^J}_{\gminib}(V)$,
by Lemma \ref{lem;10.5.16.90}.
We have $\lefttop{J}\nbigg_{\gminia}$
on $\nbigx_P\setminus\nbigd_P(J^c)$.
Assume the following:
\begin{itemize}
\item
 The filtrations
 $\nbigftilde^Q$
 of $\iota_{J\ast}\Gr^{\vecnbigf}_{\gminib}(V^{b_0})_Q$
 $(Q\in \pi_0^{-1}(P))$
 are induced by filtrations of
 $\iota'_{J\ast}\bigl(
 \lefttop{J}\nbigg^{b_0}_{\gminib}\bigr)_{Q}$.
\item
 $\Gr^{\vecnbigftilde}_{\gminia}\bigl(
 \lefttop{J}\nbigg^{b_0}_{\gminib}\bigr)$
 is isomorphic to
 $\nbigg^{b_0}_{P,\gminib|\nbigx_P\setminus\nbigd_P(J^c)}$.
\end{itemize}
By Lemma \ref{lem;10.5.16.100},
we obtain that 
the filtrations
$\nbigftilde^{(b,Q)}$
of $\iota_{J\ast}\Gr^{\vecnbigf}_{\gminib}(V)_{(b,Q)}$
 $(Q\in \pi_0^{-1}(P),\,b\in B)$
are induced by filtrations of
$\iota'_{J\ast}\bigl(
\lefttop{J}\nbigg_{\gminib}\bigr)_{(b,Q)}$.
By using Lemma \ref{lem;10.5.16.2},
we obtain that
 $\Gr^{\vecnbigftilde}_{\gminia}\bigl(
 \lefttop{J}\nbigg_{\gminib}\bigr)$
 is isomorphic to
 $\nbigg_{P,\gminib|\nbigx_P\setminus\nbigd_P(J^c)}$.
Thus, we obtain Theorem \ref{thm;10.5.16.1}.
\hfill\qed

\section{Deformation}
\label{section;10.5.26.32}
\subsection{Deformation $E^{(\nbigt)}$}

\subsubsection{Unramified case}

Let $\nbigc$ be a simply connected
compact region in $\cnum^m$.
We put $(\nbigx^{\circ},\nbigd^{\circ})
 :=(\nbigx,\nbigd)\times\nbigc$.
Let $\nbigt$ be a holomorphic function 
on $\nbigk\times\nbigc$.
Let $\vecnbigi$ be a good system of irregular values
on $(\nbigx,\nbigd)$.
For each $(P,c)\in\nbigd^{\circ}$,
we put $\nbigi_{(P,c)}^{(\nbigt)}:=\bigl\{
 \nbigt\gminia\,\big|\,\gminia\in\nbigi_P
 \bigr\}$.
Thus, we obtain a good system of
irregular values $\vecnbigi^{(\nbigt)}$.

We have an obvious deformation of 
a $\varrho$-flat bundle
with a Stokes data on $(\nbigx,\nbigd)$
over $\vecnbigi$.
Let $(V,\DD,\SD)\in \SDcat(\nbigx,\nbigd,\vecnbigi)$.
Let $(V^{\circ},\DD^{\circ})$ be
the $\varrho$-flat bundle on 
$\nbigx^{\circ}\setminus\nbigd^{\circ}$
obtained as the pull back via the projection
to $\nbigx\setminus\nbigd$.
Applying Theorem \ref{thm;10.5.14.10},
we obtain
$(V^{\circ},\DD^{\circ},\SD^{(\nbigt)})
 \in \SDcat(\nbigx^{\circ},\nbigd^{\circ},
 \vecnbigi^{(\nbigt)})$.

We have the corresponding deformation
for unramifiedly good lattice
of a meromorphic $\varrho$-flat bundle.
Namely, 
for $(E,\DD)\in \MFL(\nbigx,\nbigd,\vecnbigi)$,
we have 
$(E^{(\nbigt)},\DD^{(\nbigt)})\in
 \MFL(\nbigx^{\circ},\nbigd^{\circ},\vecnbigi^{(\nbigt)})$,
corresponding to the obvious deformation
of the Stokes data as above.
It is unique up to canonical isomorphism.

\subsubsection{Remark on descent}

Let $\varphi:(\nbigx',\nbigd')\lrarr(\nbigx,\nbigd)$
be a ramified Galois covering over $\nbigk$
with the Galois group $G$.
We put $\vecnbigi':=\varphi^{\ast}\vecnbigi$.
Take $(E',\DD')\in
 \MFL\bigl(\nbigx',\nbigd',\vecnbigi'\bigr)^G$.
Let $(E,\DD)\in\MFL(\nbigx,\nbigd,\vecnbigi)$
be the descent of $(E',\DD')$.
According to Corollary \ref{cor;10.5.17.10},
$(E',\DD')^{(\nbigt)}$
is also $G$-equivariant.

\begin{lem}
\label{lem;10.5.17.15}
$(E,\DD)^{(\nbigt)}$
is the descent of $(E',\DD')^{(\nbigt)}$.
\end{lem}
\pf
Let $(E_1,\DD_1)$ be the descent of
$(E',\DD')^{(\nbigt)}$.
By construction,
the restrictions of
$(E,\DD)^{(\nbigt)}$
and $(E_1,\DD_1)$
to $\nbigx^b$ are naturally isomorphic.
By Corollary \ref{cor;10.5.17.11},
they are isomorphic on $\nbigx$.
\hfill\qed

\vspace{.1in}

Let us consider the case
$\vecnbigi'$ is not necessarily the pull back
of a good system of irregular values on 
$(\nbigx,\nbigd)$.
Let $(E'_i,\DD'_i)\in
 \MFL(\nbigx',\nbigd',\vecnbigi')$ $(i=1,2)$.
Their descent $(E_i,\DD_i)$ are not necessarily
unramified.
\begin{lem}
\label{lem;10.5.18.20}
If $(E_1,\DD_1)\simeq (E_2,\DD_2)$,
then 
$(E_1,\DD_1)^{(\nbigt)}\simeq
(E_2,\DD_2)^{(\nbigt)}$.
\end{lem}
\pf
It is easy to reduce the issue
to the case that
$\nbigd$ is smooth.
Moreover, we have only to 
consider the case $\dim\nbigx=1$.
Because $E_1'\cap E_2'$
is also an unramifiedly good lattice,
we may assume $E_1'\subset E_2'$.
We have
$E_1^{\prime(\nbigt)}\subset
 E_2^{\prime(\nbigt)}$,
and we have only to compare
their sections which are 
invariant with respect to the Galois actions.
Then, the claim is obvious.
\hfill\qed

\subsubsection{General case}

Let $(E,\DD)$ be a good lattice
of a meromorphic $\varrho$-flat bundle
on $(\nbigx,\nbigd)$,
which is not necessarily unramified.
For any $P\in \nbigd$,
we can take a small neighbourhood 
$\nbigx_P$ and a ramified covering
$\varphi_P:
 (\nbigx_P',\nbigd_P')\lrarr(\nbigx_P,\nbigd_P)$
such that $(E,\DD)$ is the descent
of an unramifiedly good lattice $(E',\DD')$
on $(\nbigx_P',\nbigd_P')$.
By applying the procedure in the unramified case,
we obtain the deformation
$\bigl(\varphi_P^{\ast}(E,\DD)\bigr)^{(\nbigt)}$
on $(\nbigx_P^{\prime\circ},\nbigd_P^{\prime\circ})$.
By taking the descent,
we obtain $(E,\DD)^{(\nbigt)}_P$
on $(\nbigx_P^{\circ},\nbigd_P^{\circ})$.
It is well defined up to canonical isomorphisms
as a germ of good lattice of a meromorphic
$\varrho$-flat bundle at $P$,
according to Lemma \ref{lem;10.5.17.15}
and Lemma \ref{lem;10.5.18.20}.
By gluing, we can globalize 
and obtain a good lattice
$(E,\DD)^{(\nbigt)}$
of a meromorphic $\varrho$-flat bundle
on $(\nbigx^{\circ},\nbigd^{\circ})$.

If we are given a good filtered $\varrho$-flat bundle
$(\vecE_{\ast},\DD)$ on $(\nbigx,\nbigd)$,
we obtain a good filtered $\varrho$-flat bundle
$(\vecE_{\ast}^{(\nbigt)},\DD^{(\nbigt)})$
by applying the above procedure.

\subsubsection{Functoriality}

The deformation is compatible 
with dual, tensor product and direct sum,
under the appropriate assumption
on the irregular values.
Namely, we have the following natural isomorphisms:
\[
 (E_1\oplus E_2)^{(\nbigt)}
\simeq
 E_1^{(\nbigt)}\oplus E_2^{(\nbigt)},
\quad
 (E_1\otimes E_2)^{(\nbigt)}
\simeq
 E_1^{(\nbigt)}\otimes E_2^{(\nbigt)},
\quad
\bigl(E^{\lor}\bigr)^{(\nbigt)}
\simeq
\bigl(
 E^{(\nbigt)}\bigr)^{\lor}
\]

Let $(E_p,\DD_p)\in \MFL(\nbigx,\nbigd,\vecnbigi)$
$(p=1,2)$
with a morphism $f:(E_1,\DD_1)\lrarr (E_2,\DD_2)$.
Then, we have the induced morphism
$f^{(\nbigt)}:
 (E_1^{(\nbigt)},\DD_1^{(\nbigt)})
\lrarr 
 (E_2^{(\nbigt)},\DD_2^{(\nbigt)})$.

\subsubsection{}

Let $\nbigx_1$ be a complex manifold
with a normal crossing hypersurface $\nbigd_1$.
Let $F:\nbigx_1\lrarr \nbigx$ be a morphism
such that
(i) $F^{-1}(\nbigd)\subset\nbigd_1$,
(ii) the induced morphism
 $\nbigx_1\lrarr \nbigk$ is a smooth fibration,
(iii) any intersection of some irreducible components
 of $\nbigd_1$ is smooth over $\nbigk$.
Let $E$ be a good lattice of $(\nbige,\DD)$
on $(\nbigx,\nbigd)$.
We obtain a good lattice
$E_1:=F^{\ast}E$ of 
$(\nbige_1,\DD_1):=
 F^{\ast}(\nbige,\DD)
\otimes\nbigo_{\nbigx_1}(\ast \nbigd_1)$.

\begin{lem}
\label{lem;10.5.25.30}
Let $F_{\nbigc}$ be the induced morphism
$\nbigx_1^{\circ}\lrarr\nbigx^{\circ}$.
We have natural isomorphisms
\[
 E_1^{(\nbigt)}
\simeq
 F^{\ast}_{\nbigc}E^{(\nbigt)},
\quad
 (\nbige_1,\DD)^{(\nbigt)}
\simeq
 F_{\nbigc}^{\ast}\bigl(\nbige,\DD\bigr)^{(\nbigt)}
 \otimes
 \nbigo_{\nbigx_1^{\circ}}(\ast\nbigd_1^{\circ})
\]
\end{lem}
\pf
We have only to consider the unramified case,
in which the claim follows from 
Theorem \ref{thm;10.5.16.1}.
\hfill\qed

\subsection{Deformation $E^{(T)}$}
\label{subsection;10.5.17.51}

Let $T$ be a nowhere vanishing holomorphic 
function on $\nbigk$
such that $|\arg(T)|<\pi/2$.
For a given good lattice $(E,\DD)$
on $(\nbigx,\nbigd)$,
we shall construct a good lattice
$(E^{(T)},\DD^{(T)})$ on $(\nbigx,\nbigd)$.
We take a compact region
$\nbigc\subset\cnum$
which contains $0$ and $1$,
and take a nowhere vanishing holomorphic function
$\nbigt:\nbigk\times\nbigc\lrarr \cnum$
such that (i) $\nbigt_{|\nbigk\times\{0\}}=1$,
(ii) $\nbigt_{|\nbigk\times\{1\}}=T$,
(iii) $|\arg(\nbigt)|<\pi/2$.
Then, we obtain the deformation
$(E^{(\nbigt)},\DD^{(\nbigt)})$
on $(\nbigx^{\circ},\nbigd^{\circ})$.
By taking the specialization at $c=1$,
we obtain the desired $(E^{(T)},\DD^{(T)})$.

\begin{lem}
\label{lem;10.5.17.21}
$(E^{(T)},\DD^{(T)})$
is independent of the choice of
$(\nbigc,\nbigt)$
up to canonical isomorphisms.
\end{lem}
\pf
We have only to consider the local and unramified case.
Let $\nbigt_i$ $(i=0,1)$ be functions as above.
We take a small neighbourhood
$\nbigc_2$ of $\{0\leq c_2\leq 1\}$ in $\cnum$.
We consider a holomorphic function
$\nbigt_2:=
 (1-c_2)\,\nbigt_0+c_2\,\nbigt_1$
on $\nbigk\times \nbigc\times\nbigc_2$.
We obtain
$(E,\DD)^{(\nbigt_2)}$
on $(\nbigx,\nbigd)\times(\nbigc\times\nbigc_2)$.
The specializations
at $(c,c_2)=(1,i)$ $(i=0,1)$
correspond to $(E,\DD)^{(\nbigt_i)}$.
Let $p:\nbigx\times\nbigc_2\lrarr\nbigx$
be the projection.
By Theorem \ref{thm;10.5.16.1},
we have a natural isomorphism
$(E,\DD)^{(\nbigt_2)}_{|
 \nbigx\times\{1\}\times\nbigc_2}
\simeq
 p^{\ast}(E,\DD)^{(\nbigt_0)}$.
Hence, 
$(E,\DD)^{(\nbigt_i)}$ $(i=0,1)$
are naturally isomorphic.
\hfill\qed

\vspace{.1in}

Let $\vecnbigi$ be a good system of
irregular values on $(\nbigx,\nbigd)$.
For each $P\in\nbigd$,
we put
$\nbigi^{(T)}_P:=\bigl\{
 T\,\gminia\,\big|\,
 \gminia\in\nbigi_P
 \bigr\}$,
and we obtain a good system of irregular values
$\vecnbigi^{(T)}$.
The above construction gives
$\MFL(\nbigx,\nbigd,\vecnbigi)
\lrarr
 \MFL(\nbigx,\nbigd,\vecnbigi^{(T)})$,
in the unramified case.

\begin{lem}
\label{lem;07.12.21.60}
Let $T_i$ $(i=1,2)$ be holomorphic functions on $\nbigk$
such that $\bigl|\arg(T_i)\bigr|<\pi/2$
and $\bigl|\arg(T_1\,T_2)\bigr|<\pi/2$.
We have a canonical isomorphism
$E^{(T_1\, T_2)}
\simeq
 (E^{(T_1)})^{(T_2)}$.
\end{lem}
\pf
We have only to check the claim
in the local and unramified case.
We take a small neighbourhood
$\nbigc_i$ of $\{0\leq c_i\leq 1\}$ in $\cnum$.
Let us consider the function
$\nbigt(c_1,c_2)
=(1-c_1+c_1\,T_1)(1-c_2+c_2\,T_2)$
on $\nbigk\times\nbigc_1\times\nbigc_2$.
We have the deformation 
$(E,\DD)^{(\nbigt)}$
on $\nbigx\times\nbigc_1\times\nbigc_2$.
We can easily show that
both $(E,\DD)^{(T_1T_2)}$
and $\bigl((E,\DD)^{(T_1)}\bigr)^{(T_2)}$
are naturally isomorphic to
the specialization of $(E,\DD)^{(\nbigt)}$
at $(1,1)$.
\hfill\qed

\begin{lem}
\mbox{{}}\label{lem;07.12.20.5}
Under the appropriate assumptions 
on the irregular values,
the following holds:
\begin{itemize}
\item
The deformation is compatible 
with dual, tensor product, and direct sum.
\item
Let $(E_1,\DD_1)
 \lrarr (E_2,\DD_2)$ be a flat morphism.
Then, we have an induced flat morphism
$(E_1^{(T)},\DD_1^{(T)})
\lrarr
 (E_2^{(T)},\DD_2^{(T)})$.
\hfill\qed
\end{itemize}
\end{lem}

Let $\nbigx_1$ be a complex manifold
with a normal crossing hypersurface $\nbigd_1$.
Let $F:\nbigx_1\lrarr \nbigx$ be a morphism
such that
(i) $F^{-1}(\nbigd)\subset\nbigd_1$,
(ii) the induced morphism
 $\nbigx_1\lrarr \nbigk$ is a smooth fibration,
(iii) any intersection of some irreducible components
 of $\nbigd_1$ is smooth over $\nbigk$.
Let $E$ be a good lattice of $(\nbige,\DD)$
on $(\nbigx,\nbigd)$.
We obtain a good lattice
$E_1:=F^{\ast}E$ of 
$(\nbige_1,\DD_1):=F^{\ast}(\nbige,\DD)$.
We obtain the following lemma
from Lemma \ref{lem;10.5.25.30}.
\begin{lem}
\label{lem;10.5.25.31}
We have a natural isomorphism
$E_1^{(T)}\simeq
 F^{\ast}E^{(T)}$.
\hfill\qed
\end{lem}

\subsection{Deformation in the case that 
 $|\arg(T)|$ is small}
\label{subsection;08.9.28.41}

We give a characterization of
holomorphic sections of 
$(E^{(T)},\DD^{(T)})$
when $|\arg(T)|$ is sufficiently small.
We explain it in the case
that $(E,\DD)$ is unramified.
We put $\nbigx=\Delta^n\times\nbigk$
and $\nbigd=\bigcup_{i=1}^{\ell}\{z_i=0\}$.
Let $\vecnbigi$ be a good system of
irregular values on $(\nbigx,\nbigd)$.

Let $P\in\nbigd_{\ellsitabar}$.
We take a covering
$\{U_i\,|\,i\in\Gamma\}$ of $\pi^{-1}(P)$,
which is good for $\nbigi_P$.
(See Definition \ref{df;10.5.18.35}.)
We take neighbourhoods $\nbigu_i$ of $U_i$
in $\nbigxtilde(\nbigd)$.
Assume the following for $T$:
\begin{itemize}
\item
 $\{U_i\}$ is good for
 $\nbigi_P^{(T_c)}$
 for $0\leq c\leq 1$,
 where $T_c:=1+c(T-1)$.
\end{itemize}
This is satisfied if $\arg(T)$ is sufficiently small
for fixed $\{U_i\}$.

Let $\vecv=(\vecv_{\gminia})$
be a frame of $\Gr^{\vecnbigftilde}(E)$
compatible with the grading.
Let $\nbigu_I:=\bigcap_{i\in I}\nbigu_i$
for some $I\subset\Gamma$.
As in Subsection \ref{subsection;10.5.17.40},
we take a splitting
\begin{equation}
\label{eq;10.5.17.55}
 E_{|\nbigu_I}=\bigoplus E_{\nbigu_I,\gminia}.
\end{equation}
We have the induced frame
$\vecv_{I}$ of $E_{|\nbigu_I}$.
We put
\[
 \vecv_{I,\gminia}^{(T)}:=
 \vecv_{I,\gminia}\,
 \exp\bigl((T-1)\varrho^{-1}\gminia\bigr),
\quad
 \vecv_{I}^{(T)}
:=\bigl(\vecv^{(T)}_{I,\gminia}\bigr).
\]
\begin{lem}
\label{lem;10.5.26.10}
Let $\vecw$ be a frame of $E^{(T)}$.
Let $G_{I}$ be determined by
\[
 \vecw_{|\nbigu_I\setminus\pi^{-1}(\nbigd)}
=\vecv^{(T)}_{I|\nbigu_I\setminus\pi^{-1}(\nbigd)}
 \,G_I.
\]
Then, the entries of
$G_I$ and $G_I^{-1}$ are
holomorphic on $\nbigu_I$.

The splitting
$E_{|\nbigu_I\setminus\pi^{-1}(\nbigd)}
=\bigoplus E_{\nbigu_I,\gminia|\nbigu_I\setminus
 \pi^{-1}(\nbigd)}$
gives a splitting of the full Stokes filtration
of $E^{(T)}$.
\end{lem}
\pf
Let $\nbigc$ be a small neighbourhood of
$\{0\leq c\leq 1\}\subset\cnum$.
We consider the function
$\nbigt=1+c(T-1)$.
We have an unramifiedly good lattice
$(E,\DD)^{(\nbigt)}$ 
over $\nbigc\times(\nbigx,\nbigd)$.
By the assumption on $\arg(T)$,
the natural map
$\nbigi_P^{(\nbigt)}\lrarr\nbigi_P$
induces isomorphisms of the ordered sets:
\[
 \bigl(
 \nbigi_P^{(\nbigt)},\leq_{\nbigc\times\nbigu_I}
 \bigr)
\lrarr
 \bigl(\nbigi_P,\leq_{\nbigu_I}
 \bigr)
\]
By using the flat $\varrho$-connection,
we obtain a filtration
$\nbigftilde^{\nbigc\times\nbigu_I}$
of $E^{(\nbigt)}_{|\nbigc\times\nbigu_I}$
from $\nbigftilde^{\nbigu_I}$
of $E_{|\nbigu_I}$.
We also obtain a splitting
of $\nbigftilde^{\nbigc\times\nbigu_I}$
from the splitting (\ref{eq;10.5.17.55})
of $\nbigftilde^{\nbigu_I}$.
We remark that,
for any $(c,Q)\in \nbigc\times\nbigu_I$
the filtrations $\nbigftilde^{\nbigu_I}$
and $\nbigftilde^{(c,Q)}$ are compatible
over 
$\bigl(\nbigi_P^{(\nbigt)},\leq_{\nbigc\times\nbigu_I}
 \bigr)
\lrarr
 \bigl(\nbigi_P^{(\nbigt)},\leq_{(c,Q)}\bigr)$,
which follows from the characterization of
the full Stokes filtration in 
Theorem \ref{thm;10.5.12.110}.
In particular, the restriction of the splitting
to $\{1\}\times\bigl(
 \nbigu_I\setminus\pi^{-1}(\nbigd)\bigr)$
gives a splitting of the filtration
$\nbigftilde^{\nbigu_I}$
of $E^{(T)}$.
Note that 
$\vecv^{(T)}_I$ naturally gives
a frame of
$\Gr^{\vecnbigftilde}(E)
\otimes\nbigl\bigl(
 (T-1)\gminia
 \bigr)$,
where $\nbigl\bigl((T-1)\gminia\bigr)$
denotes
$\nbigo_{\nbigx}\, e$
with a flat $\varrho$-connection
$\DD e=e\, (T-1)d\gminia$.
Then, we obtain the first claim of the lemma
from Lemma \ref{lem;10.5.12.100}.
The second claim follows from the first one.
\hfill\qed

\vspace{.1in}

Let $f$ be a holomorphic section 
of $E_{|\nbigx-\nbigd}$.
We have the corresponding decomposition
$f_{|\nbigu_I}=\sum f_{\gminia,I}$.
We have the expression
$f_{\gminia,I}=
 \sum f^{(T)}_{\gminia,I,j}
 \, v^{(T)}_{\gminia,I,j}$.
We put 
$\vecf_{\gminia,I}:=
 \bigl(f^{(T)}_{\gminia,\nbigu,j}\bigr)$.
\begin{cor}
\label{cor;10.5.26.1}
$f$ gives a section of $E^{(T)}$
if and only if
$\vecf^{(T)}_{\gminia,I}$
is bounded for any $\gminia$ and $I$.
\hfill\qed
\end{cor}

\chapter{$L^2$-Cohomology of
Filtered $\lambda$-Flat Bundle on Curves}
\label{section;07.10.7.1}
We would like to compare 
various cohomology groups
associated to a filtered $\lambda$-flat bundle
with harmonic metric on a curve,
which will be achieved in Section
\ref{subsection;07.10.23.30}.
This chapter is a preparation
for local comparisons.
In Sections 
\ref{subsection;07.10.20.1}--\ref{subsection;08.1.5.2},
we consider the case in which $\lambda$ is fixed.
In Sections 
\ref{subsection;07.10.20.6}--\ref{subsection;07.11.16.120},
we study the family version.
We separate them although they are essentially the same.
The statements are given in
Sections \ref{subsection;07.10.20.1}
and \ref{subsection;07.10.20.6}, respectively.

\section{Local quasi isomorphisms for fixed $\lambda$}
\label{subsection;07.10.20.1}
We put $X:=\Delta_z$.
For any subset $Y\subset X$,
we put $Y^{\ast}:=Y\setminus\{O\}$.
Let $(V_{\ast},\DDlambda)$ be
a good filtered $\lambda$-flat bundle
on $(X,O)$.

\subsection{Sheaves of $L^2$-sections and
  holomorphic $L^2$-sections}
\label{subsection;07.9.11.11}

\subsubsection{Preliminary for metric}
\label{subsection;10.5.21.30}

Let $\vecv$ be a frame of $\prolong{V}$ compatible with
the parabolic filtration $F$
and the weight filtration $W$ on $\Gr^F$.
We put $a(v_i):=\deg^F(v_i)$ and 
$k(v_i):=\deg^W(v_i)$.
Note $-1<a(v_i)\leq 0$ and $k(v_i)\in\seisuu$.
Let $h$ be the hermitian metric given as follows:
\[
 h(v_i,v_j):=\delta_{i,j}\,
 |z|^{-2a(v_i)}\,(-\log|z|)^{k(v_i)}
\]
If a metric $h'$ comes from
another choice of a frame $\vecv'$
compatible with $F$ and $W$,
the metrics $h$ and $h'$ are mutually bounded.
Let $g_{\poin}$ denote the Poincar\'e metric
of $X^{\ast}$.

We recall basic property of the metric $h$
as above.
Let $f=\sum f_jv_j$ be a holomorphic section
of $V_{|X^{\ast}}$.
It is $L^2$ with respect to $h$ and $g_{\poin}$,
if and only if the following holds for 
each $f_i$:
\begin{itemize}
\item
 $f_i$ is holomorphic, if
 (i) $-1<a(v_i)<0$,
 or (ii) $a(v_i)=0$ and $k(v_i)\leq 0$.
\item
 $f_i$ is holomorphic and
 $f_i(O)=0$, if 
 $a(v_i)=0$ and $k(v_i)>0$. 
\end{itemize}
Let $\omega=\sum \omega_i\,v_i\,dz$
be a holomorphic section of 
$V\otimes\Omega_X^{1,0}$ on $X^{\ast}$.
It is $L^2$ with respect to $h$ and $g_{\poin}$
if and only if the following holds:
\begin{itemize}
\item
 $f_i$ is holomorphic, if
 (i) $-1<a(v_i)<0$,
 or (ii) $a(v_i)=0$ and $k(v_i)\leq -2$.
\item
 $f_i$ is holomorphic and
 $f_i(O)=0$, if 
 $a(v_i)=0$ and $k(v_i)>-2$.
\end{itemize}
We can check these claims 
by direct computations.

\subsubsection{Sheaf of $L^2$-sections}

\label{subsection;10.5.21.1}

Note that
the $\lambda$-connection $\DDlambda$ of $V$
and the derivation $\lambda\del_X+\delbar_X$
induce a derivation of
$V\otimes\Omega^{\bullet,\bullet}_X$,
which is also denoted by $\DDlambda$.
For an open subset $U\subset X$,
let $\nbigl^p(V_{\ast},\DDlambda)(U)$ 
be the space of sections $\tau$ of
$V\otimes\Omega_X^p$ on $U^{\ast}$
with the following property:
\begin{itemize}
\item
$\tau$ and $\DDlambda\tau$
are $L^2$ locally on $U$,
with respect to $h$ and $g_{\poin}$.
\end{itemize}
Let $\nbigl^p_{\poly}(V_{\ast},\DDlambda)(U)$
be the space of
$C^{\infty}$-sections $\tau$ of
$V\otimes\Omega_X^p$ on $U^{\ast}$
with the following property:
\begin{itemize}
\item
$\tau$ and $\DDlambda\tau$ are
$L^2$ and of polynomial order in $|z^{-1}|$
locally on $U$
with respect to $h$ and $g_{\poin}$.
\end{itemize}

\begin{rem}
In the following,
we say just ``polynomial order''
instead of ``polynomial order in $|z^{-1}|$''.
\hfill\qed
\end{rem}

Thus, we obtain complexes of sheaves
$\nbigl^{\bullet}(V_{\ast},\DDlambda)$
and 
$\nbigl_{\poly}^{\bullet}(V_{\ast},\DDlambda)$.
Let $\nbigl_{\hol}^{p}(V_{\ast},\DDlambda)$
$(p=0,1)$ be the subsheaves
$\nbigl^{p}(V_{\ast},\DDlambda)$,
which consists of holomorphic $p$-forms.
By a general theory of holomorphic functions,
we have $\nbigl^p_{\hol}(V_{\ast},\DDlambda)
\subset \nbigl^p_{\poly}(V_{\ast},\DDlambda)$.

We shall prove the following proposition in Sections
\ref{subsection;07.9.12.10}--\ref{subsection;07.9.12.2}.
The arguments are minor modification of
those in \cite{sabbah2} and \cite{z}.
\begin{prop}
\label{prop;07.9.12.20}
The naturally defined morphisms
\[
\begin{CD}
 \nbigl_{\hol}^{\bullet}(V_{\ast},\DDlambda)
@>{\varphi_0}>>
 \nbigl_{\poly}^{\bullet}(V_{\ast},\DDlambda)
@>{\psi_1}>>
 \nbigl^{\bullet}(V_{\ast},\DDlambda)
\end{CD}
\]
are quasi isomorphisms.
\end{prop}

\subsubsection{Algebraically determined sheaf}

Let $X'=\Delta_{z'}$,
and let $\varphi_n$ denote the ramified covering
$X'\lrarr X$ given by $\varphi_n(z')=z^{\prime\,n}$.
Recall that we have the induced 
good filtered $\lambda$-flat bundle on $(X',O')$
as in Section \ref{subsection;07.11.5.60},
which is denoted by 
$(V'_{\ast},\DD^{\prime\lambda})$.
If we choose $n$ appropriately,
$(V'_{\ast},\DD^{\prime\lambda})$
is unramified,
and we have the irregular decomposition:
\begin{equation}
 \label{eq;08.9.9.10}
 (V'_{\ast},\DD^{\prime\lambda})_{|\Ohat'}=
 \bigoplus_{\gminia\in \Irr(\DD^{\prime\lambda})}
 \bigl(\Vhat'_{\gminia\,\ast},
   \DD^{\prime\lambda}_{\gminia}\bigr)
\end{equation}
Since $\Vhat'_{0\ast}$ and
$\bigoplus_{\gminia\neq 0}
 \Vhat'_{\gminia\,\ast}$ are $\Gal(X'/X)$-equivariant,
we have the descent to $\Ohat$
which are denoted by 
$\Vhat_{\reg\ast}$ and $\Vhat_{\irr\ast}$,
respectively.

Let $a\in\real$ and 
$\gminia\in\Irr(\DD^{\prime\lambda})$.
We have the weight filtration $W$ of
the nilpotent part of the residue 
$\Gr^F_a\Res(\DD^{\prime\lambda})$ on 
$\Gr^F_a\bigl(\Vhat'_{\gminia}\bigr)$.
Let $W_k\bigl(\prolongg{a}{\Vhat'_{\gminia}}\bigr)$
denote the pull back of 
$W_k\Gr^F_a\bigl(\Vhat'_{\gminia}\bigr)$
via the natural projection
$\prolongg{a}{\Vhat'_{\gminia}}
\lrarr \Gr^F_a\bigl(\Vhat'_{\gminia}\bigr)$.
For each $\gminia\neq 0$,
we put
\[
 \nbigs\bigl(\Vhat'_{\gminia\,\ast}
 \otimes\Omega^{0,0}_{X'}\bigr):=
 W_{-2}\bigl(\prolongg{\ord(\gminia)}
 {\Vhat'_{\gminia}}\bigr),
\quad
\nbigs(\Vhat'_{\gminia\,\ast}
 \otimes\Omega^{1,0}_{X'}):=
 W_{-2}\bigl(
 \prolong{\Vhat'_{\gminia}}\bigr)
 \,\frac{dz'}{z'}.
\]
We have the descent of
$\bigoplus_{\gminia\neq 0}
\nbigs\bigl(\Vhat'_{\gminia\ast}
 \otimes\Omega_X^{p,0}\bigr)$,
which is denoted by
$\nbigs(\Vhat_{\irr\ast}\otimes\Omega_X^{p,0})$.

We have the generalized eigen decomposition
with respect to the residue
\[
 \Gr^F_0(\Vhat_0)
=\bigoplus_{\alpha}
 \EE_{\alpha}\Gr^F_0(\Vhat_0),
\]
where the restriction of $\Res(\DDlambda)$ 
to $\EE_{\alpha}\Gr^F_0(\Vhat_0)$
has the unique eigenvalue $\alpha$.
We have the weight filtration $W$
of the nilpotent part of the residue
$\Res(\DDlambda)$ on each 
$\EE_{\alpha}\Gr^F_0(\Vhat_0)$
and $\Gr^F_0(\Vhat_0)$.
Let $\nbigs(\Vhat_{0\ast}\otimes\Omega_X^{0,0})$
denote the inverse image of
$\bigoplus_{\alpha\neq 0}
 W_{-2}\EE_{\alpha}\Gr^{F}_0(\prolong\Vhat_{0})
\oplus
 W_0\EE_0\Gr^F_0\bigl(\prolong{\Vhat_{0}}\bigr)$
via the projection
$\prolong{\Vhat_0}\lrarr
 \Gr_0^F(\Vhat_{0})$.
Let $\nbigs(\Vhat_{0\ast}\otimes\Omega_X^{1,0})$
denote the inverse image of
$W_{-2}\Gr^F_0(\Vhat_{0})$
via the projection
$\prolong{\Vhat_0}\, dz/z\lrarr
 \Gr^F_0(\Vhat_{0})$.
Thus, we obtain a lattice
\[
\nbigs(\Vhat_{\ast}\otimes\Omega_X^{p,0})
:=\nbigs(\Vhat_{0\ast}\otimes\Omega_X^{p,0})
 \oplus
 \nbigs(\Vhat_{\irr\ast}\otimes\Omega_X^{p,0}) 
\]
of $\bigl(V\otimes\Omega_X^{p,0}(\ast O)\bigr)_{|\Ohat}$.
They induce lattices
$\nbigs(V_{\ast}\otimes\Omega_X^{p,0})$
of $V\otimes\Omega_X^{p,0}(\ast O)$.
The $\lambda$-connection $\DDlambda$ on $V_{\ast}$
and the differential $\lambda\, d_X$ on 
$\Omega^{\bullet,0}_X$
induce
$\DDlambda:
 \nbigs(V_{\ast}\otimes\Omega_X^{0,0})
\lrarr 
 \nbigs(V_{\ast}\otimes\Omega_X^{1,0})$.
Thus, we obtain a complex of sheaves
$\nbigs(V_{\ast}\otimes\Omega^{0,0})
  \stackrel{\DDlambda}{\lrarr}
 \nbigs(V_{\ast}\otimes\Omega^{1,0})$.
\begin{lem}
\label{lem;10.5.22.10}
We have a natural inclusion
$\nbigs(V_{\ast}\otimes\Omega^{\bullet,0})
\lrarr
 \nbigl_{\hol}^{\bullet}(V_{\ast},\DDlambda)$,
which is an isomorphism.
\end{lem}
\pf
We have only to compare
the germs at $O$.
By the condition in Section 
\ref{subsection;10.5.21.30},
it is easy to check that
sections of $\nbigs(V_{\ast}\otimes\Omega^{p,0})$
are $L^2$.
It also implies that 
$\DD f$ is $L^2$ for a section $f$
of $\nbigs(V_{\ast})$.
Hence, we obtain 
$\nbigs(V_{\ast}\otimes\Omega^{p,0})
\subset
 \nbigl_{\hol}^p(V_{\ast},\DDlambda)$
naturally.
In the case $p=1$,
it is clearly an isomorphism.
Let $f\in \nbigl_{\hol}^0(V_{\ast},\DDlambda)$.
Because $f$ is $L^2$,
each $f_i$ is holomorphic.
Then, we obtain $f\in \nbigs(V_{\ast})$
from $\DDlambda f\in 
 \nbigs(V_{\ast}\otimes\Omega^{1,0})$.
\hfill\qed

\subsubsection{Remark}
This kind of theorems, such as 
Proposition \ref{prop;07.9.12.20}
and Lemma \ref{lem;10.5.22.10},
was first proved by
S. Zucker {\rm\cite{z}}
in his study on singular variation of
Hodge structure.
Namely, he used the quasi isomorphism
to obtain a Hodge structure on
the intersection cohomology of
a variation of polarized pure Hodge structure.
If $h$ comes from a variation of
polarized pure Hodge structure,
it is easy to see that
$\nbigs(V_{\ast}\otimes\Omega^{\bullet,0})$
with $\lambda=1$
is naturally quasi isomorphic to
the de Rham complex of the minimal
extension of $V_{|X^{\ast}}$ on $X$.
The quasi isomorphism with $\lambda=1$
plays the role
connecting the intersection cohomology
and the $L^2$-cohomology.
Moreover, he obtained a Hodge structure
on $L^2$-cohomology
by using the quasi isomorphism with $\lambda=0$
with some global harmonic analysis.

For regular filtered $\lambda$-flat bundles,
it was proved in \cite{sabbah2}
and \cite{mochi2}
for the study on tame harmonic bundles
on curves.
In \cite{mochi2},
we used Zucker's method
in a straightforward way.
Sabbah \cite{sabbah2}
introduced an improvement
to argue it in a unified way for various $\lambda$.
He also studied the irregular singular case
in {\rm\cite{sabbah3}}.
We will use a different argument
to deal with Stokes structure.

\subsection{Variants in the case $\lambda\neq 0$}
\label{subsection;07.11.15.1}

We shall introduce 
complexes of sheaves 
$\nbigltilde^{\bullet}_{\poly}
 (V_{\ast},\DDlambda)$
and
$\nbiglbar^{\bullet}_{\poly}
 (V_{\ast},\DDlambda)$
on $X$,
whose restrictions to $X^{\ast}$
are the same as 
$\nbigl^p_{\poly}(V_{\ast},\DDlambda)_{|X^{\ast}}$.
Let $S$ be a small sector in $X^{\ast}$,
and let $\Sbar$ denote its closure in
the real blow up $\Xtilde(O)$.
We can take a lift of $\Sbar$
in the real blow up $\Xtilde'(O')$
via the covering map
$\Xtilde'(O')\lrarr\Xtilde(O)$,
which is denoted by the same notation $\Sbar$.
If $\lambda\neq 0$,
we have the full Stokes filtration
$\nbigftilde^S$ of $V'_{|\Sbar}$. 
We can take a flat splitting
$\prolongg{a}{V'}_{|\Sbar}
=\bigoplus_{\gminia}
 \prolongg{a}{V'}_{\gminia,S}$.
We can naturally identify
$V'_{|S}$ and $V_{|S}$.
Each section $f$ of $V\otimes\Omega_X^p$ on $S$
has the corresponding decomposition
$f=\sum f_{\gminia,S}$.

\subsubsection{Variant 1}

For an open $U\subset X$ with $O\in U$,
let $\nbigltilde^{p}_{\poly}(V_{\ast},\DDlambda)(U)$
be the space of $C^{\infty}$-sections $\tau$ 
of $V\otimes\Omega_X^p$ on $U^{\ast}$,
such that the following estimate holds 
with respect to $h$ and $g_{\poin}$
on each small sector $S$:
\begin{description}
\item[(a1)]
 $\tau_{\gminia,S}$ and 
 $\DDlambda\tau_{\gminia,S}$ $(\gminia\neq 0)$
 are of polynomial order.
\item[(a2)]
 $\tau_{0,S}$ and $\DDlambda\tau_{0,S}$
 are $L^2$ and of polynomial order.
\end{description}
 The conditions are independent of the choice 
 of a flat splitting 
 and a lift of sector to $\Xtilde'(O')$.
Then, we obtain the complex of sheaves
$\nbigltilde^{\bullet}_{\poly}(V_{\ast},\DDlambda)$.
We will prove the following proposition
in Section \ref{subsection;07.9.12.11}.
\begin{prop}
\label{prop;07.9.12.21}
The naturally defined morphisms
\[
 \begin{CD}
 \nbigs(V_{\ast}\otimes\Omega^{\bullet,0}_X)
@>{\varphi_0}>>
 \nbigl_{\poly}^{\bullet}(V_{\ast},\DDlambda)
@>{\psi_2}>>
 \nbigltilde_{\poly}^{\bullet}(V_{\ast},\DDlambda)
 \end{CD}
\]
are quasi isomorphisms.
\end{prop}

\subsubsection{Variant 2}

For an open $U\subset X$ with $O\in U$,
let $\nbiglbar^p_{\poly}(V_{\ast},\DDlambda)(U)$
be the space of $C^{\infty}$-sections $\tau$
of $V\otimes\Omega_X^p$ on $U^{\ast}$,
such that the following estimate holds 
with respect to $h$ and $g_{\poin}$
on each small sector $S$:
\begin{description}
\item[(b1)]
 $\tau_{\gminia,S}$ and $\DDlambda\tau_{\gminia,S}$
 ($\gminia\neq 0$)
 are $O\bigl(|z|^N\bigr)$ for any $N>0$
 with respect to $h$ and $g_{\poin}$.
\item[(b2)]
 $\tau_{0,S}$ and
 $\DDlambda\tau_{0,S}$ satisfy (a2).
\end{description}
Then, we obtain the complex of sheaves
$\nbiglbar^{\bullet}_{\poly}(V_{\ast},\DDlambda)$.
We will prove the following proposition
in Section \ref{subsection;08.9.9.1}.

\begin{prop}
\label{prop;07.10.22.55}
The naturally defined morphisms
\[
 \begin{CD}
 \nbiglbar_{\poly}^{\bullet}(V_{\ast},\DDlambda)
@>{\varphi_1}>>
 \nbigl_{\poly}^{\bullet}(V_{\ast},\DDlambda)
@>{\psi_2}>>
 \nbigltilde_{\poly}^{\bullet}(V_{\ast},\DDlambda)
 \end{CD}
\]
are quasi isomorphisms.
\end{prop}

\subsection{Deformation of 
the Stokes structure $(\lambda\neq 0)$}
\label{subsection;07.10.20.2}

Let $T>0$.
We have the deformation 
$\bigl(V^{(T)}_{\ast},\DDlambda\bigr)$
as in Section \ref{subsection;10.5.17.51}.
To compare
$\nbiglbar_{\poly}^{\bullet}
 \bigl(V^{(T)}_{\ast},\DDlambda\bigr)$
and 
$\nbiglbar_{\poly}^{\bullet}
 \bigl(V_{\ast},\DDlambda\bigr)$,
we shall introduce a complex of sheaves
$\nbiglbar^{\bullet}_{\poly}
 \bigl(V,h,h_{C^{\infty}}^{(T)}\bigr)$.

\subsubsection{Preliminary for metrics}

On each small sector $S$ in $X^{\ast}$,
we take a flat splitting
$V_{|\Sbar}=\bigoplus V_{\gminia,S}$
as in Section \ref{subsection;07.11.15.1}.
Let $G^{(T)}_S$ be the endomorphism
given by
$\bigoplus_{\gminia}
 \exp\bigl((1-T)\,\lambda^{-1}\,\gminia\bigr)\, 
\id_{V_{\gminia,S}}$.
Let $h^{(T)}_S$ be the hermitian metric of
$V_{|S}$ given by
$h^{(T)}_S(u,v):=
 h\bigl(G^{(T)}_S(u),G^{(T)}_S(v)\bigr)$.
If we construct $h^{\prime\,(T)}_S$
from other $h'$ and $G_S^{\prime(T)}$,
$h^{\prime(T)}_S$ and $h^{(T)}_S$
are mutually bounded.
By varying $S$ and gluing $h^{(T)}_S$ in $C^{\infty}$,
we obtain a $C^{\infty}$-metric
$h^{(T)}_{C^{\infty}}$.

We can construct a metric $h^{(T)}$ 
for $V_{\ast}^{(T)}$
as in Section \ref{subsection;10.5.21.1},
by taking a frame $\vecv^{(T)}$ of 
$\prolong{V^{(T)}}$
compatible with the parabolic filtration $F$
and the weight filtration $W$ on $\Gr^F$.

\begin{lem}
\label{lem;08.1.19.30}
The metrics $h^{(T)}_{C^{\infty}}$
and $h^{(T)}$ are mutually bounded.
\end{lem}
\pf
We may assume that 
$V_{\ast}$ is unramified.
Let $\vecwhat=(\vecwhat_{\gminia})$
be a frame of $V_{|\Ohat}$,
such that 
(i) it is compatible with the irregular decomposition
and the parabolic filtration,
(ii) the induced frame of
$\Gr^F(\prolong{V})$ is compatible with
the weight filtration $W$.
We put $a(w_i):=\deg^F(w_i)$
and $k(w_i):=\deg^W(w_i)$.
We take a lift $\vecw_{S}=(\vecw_{\gminia,S})$
of $\vecwhat$ to $\prolong{V}_{|\Sbar}$,
compatible with a $\DDlambda$-flat splitting 
of the full Stokes filtration.
Let $\htilde_S$ be the hermitian metric of
$V_{|S}$ determined by
$\htilde_S(w_{i},w_j):=\delta_{i,j}\,
 |z|^{-2a(w_i)}\, \bigl(-\log|z|\bigr)^{k(w_i)}$.
Then, it is mutually bounded with $h_{|S}$.

We have the frame $\vecwhat^{(T)}$
of $\prolong{V}^{(T)}_{|\Ohat}$,
which is obtained from $\vecwhat$
by the natural (non-flat) isomorphism
$V_{\ast|\Ohat}\simeq
 V^{(T)}_{\ast|\Ohat}$.
We put
$\vecw^{(T)}_{\gminia,S}:=
 \exp\bigl((T-1)\,\lambda^{-1}\,\gminia\bigr)
 \,\vecw_{\gminia,S}$,
and $\vecw^{(T)}_S:=(\vecw^{(T)}_{\gminia,S})$.
Then $\vecw^{(T)}_{S}$ is a lift
of $\vecwhat^{(T)}$ to $V^{(T)}_{|\Sbar}$
compatible with a $\DDlambda$-flat splitting
of the full Stokes filtration.
Let $\htilde^{(T)}_S$ be the metric of
$V^{(T)}_{|S}$ given by
$\htilde^{(T)}_S(w^{(T)}_i,w^{(T)}_j)
=\delta_{i,j}\, |z|^{-2a(w^{(T)}_i)}
 \, \bigl(-\log|z|\bigr)^{k(w^{(T)}_i)}$.
Then, $\htilde^{(T)}_S$ and $(h^{(T)})_{|S}$
are mutually bounded.

By the construction,
we have 
$\htilde^{(T)}_S\bigl(u,v\bigr)
=\htilde_S\bigl(G_S^{(T)}(u),G_S^{(T)}(v)\bigr)$.
Then, the claim of the lemma follows.
\hfill\qed

\subsubsection{A complex of sheaves}

We shall introduce the complex of sheaves
 $\nbiglbar^{\bullet}_{\poly}
 \bigl(V,h,h_{C^{\infty}}^{(T)}\bigr)$
on $X$.
We set
$\nbiglbar^p_{\poly}
 \bigl(V,h,h_{C^{\infty}}^{(T)}\bigr)_{|X^{\ast}}
=
\nbiglbar^p_{\poly}
 \bigl(V_{\ast},\DDlambda\bigr)_{|X^{\ast}}$.
For an open set $U\subset X$ with $O\in U$,
let $\nbiglbar^p_{\poly}
 \bigl(V,h,h_{C^{\infty}}^{(T)}\bigr)(U)$ be 
the space of $C^{\infty}$-sections $\tau$ of
$V\otimes\Omega^{p}$ on $U^{\ast}$
such that the following estimate holds
on each small sector $S$:
\begin{description}
\item[(c1)]
$\tau_{\gminia,S}$ and $\DDlambda\tau_{\gminia,S}$
 $(\gminia\neq 0)$
 are $O\bigl(|z|^N\bigr)$ for any $N$
with respect to  both
$(h,g_{\poin})$ and 
$(h_{C^{\infty}}^{(T)},g_{\poin})$.
\item[(c2)]
 $\tau_{0,S}$ and 
 $\DDlambda \tau_{0,S}$
 are $L^2$ and of polynomial order
with respect to $(h,g_{\poin})$.
 In other words,
 they satisfy the condition (a2).
 Note that the restrictions of 
 $h$ and $h^{(T)}_{C^{\infty}}$
 to $V_{0,S}$ are mutually bounded.
\end{description}
Thus, we obtain the complex of sheaves
$\nbiglbar^{\bullet}_{\poly}
 \bigl(V,h,h_{C^{\infty}}^{(T)}\bigr)$.

\subsubsection{Statement}
By construction of the complexes
and Lemma \ref{lem;08.1.19.30},
we have the following 
natural morphisms:
\begin{equation}
 \label{eq;10.5.21.32}
 \nbiglbar^{\bullet}_{\poly}
 \bigl(V_{\ast},\DDlambda\bigr)
\llarr
\nbiglbar^{\bullet}_{\poly}
 \bigl(V,h,h_{C^{\infty}}^{(T)}\bigr)
\lrarr
 \nbiglbar^{\bullet}_{\poly}
 \bigl(V^{(T)}_{\ast},\DDlambda\bigr)
\end{equation}
We will prove the following proposition
in Section \ref{subsection;08.9.9.5}.

\begin{prop}
\label{prop;07.7.22.35}
The morphisms in {\rm(\ref{eq;10.5.21.32})}
are quasi isomorphisms.
\end{prop}

\subsubsection{Remark}

We give a consequence of Proposition
\ref{prop;07.7.22.35}
for holonomic $D$-modules
on {\em projective curves}.
Let $C$ be a smooth projective curve.
Let $V$ be a meromorphic flat bundle on $C$.
For a given $T_i>0$ $(i=1,2)$,
we have the deformation $V^{(T_i)}$.
Let $M^{(T_i)}$ be the minimal
extensions of $V$ and $V^{(T_i)}$,
respectively.
We can deduce a natural isomorphism
of the cohomology of $D$-modules
$H^{\ast}(C,M^{(T_1)})
\simeq H^{\ast}(C,M^{(T_2)})$
by using the above quasi isomorphisms.
Actually, let $V^{(T_i)}_{\ast}$
be the Deligne-Malgrange filtered bundle
associated to $V^{(T_i)}$.
(See Section \ref{section;10.5.4.100}.)
It is standard that
$\nbigs(V_{\ast}^{(T_i)}\otimes
 \Omega^{\bullet,0})$
is naturally quasi isomorphic to
the de Rham complex of $M^{(T_i)}$.
Hence, the quasi isomorphisms in
Propositions {\rm\ref{prop;07.9.12.21}},
{\rm\ref{prop;07.10.22.55}}
and {\rm\ref{prop;07.7.22.35}}
induce the desired isomorphism 
$H^{\ast}(C,M^{(T_1)})
\simeq H^{\ast}(C,M^{(T_2)})$.

We can also obtain such an isomorphism
directly from the construction of $M^{(T_i)}$.
Recall that it is obtained as the specialization
of a meromorphic flat bundle $M^{(\nbigt)}$
on $\nbigc\times C$ for some appropriate
complex manifold $\nbigc$ with a function $\nbigt$.
By taking push-forward to $\nbigc$,
we obtain a flat bundle whose fiber over
$c\in\nbigc$ is naturally quasi isomorphic
to $H^{\ast}(\{c\}\times C,M^{(\nbigt(c))})$.
Hence, the parallel transport 
induces the desired isomorphism.
It seems to be able to check that 
two isomorphisms are the same
by using the family version of
the quasi-isomorphisms,
which might simplify our argument.
We would like to give more details
somewhere.

\section{Proof for fixed $\lambda$}
\label{subsection;08.1.5.2}
\subsection{An estimate in \cite{z}}
\label{subsection;07.9.12.10}

We recall a result due to Zucker \cite{z}.
We use the Poincar\'e metric $g_{\poin}$
and the associated volume form $\dvol_{g_{\poin}}$
of $X^{\ast}$ around $O$.
We use the polar coordinate $z=r\, e^{\sqrt{-1}\theta}$.
Let $\nbigl$ be a holomorphic line bundle on $X^{\ast}$
with a holomorphic frame $\sigma$
and a metric $h$ such that 
$|\sigma|_h\sim r^{-a}\, |\log r|^{k/2}$,
where $-1<a\leq 1/2$ and $k\in\seisuu$.
Let $\|\omega\|_{h,g_{\poin}}$ denote 
the $L^2$-norm of a section $\omega$ of 
$\nbigl\otimes\Omega^p$ 
with respect to $h$ and $g_{\poin}$.

Let $0<R<1/2$.
We put
$X(R):=\bigl\{z\in\cnum\,\big|\,|z|\leq R\bigr\}$
and $X^{\ast}(R)=X(R)\setminus\{O\}$.
Let $\omega=g\, \sigma\, d\zbar/\zbar$
be a $C^{\infty}$-section of $\nbigl\otimes\Omega^{0,1}_X$ on 
$X^{\ast}(R)$
with compact support.
We have the Fourier expansion
$g=\sum_{m\in\seisuu}
 g_m(r)\, e^{\sqrt{-1}m\theta}$.
We put 
 $g^{(1)}:=\sum_{m\neq 0}
 g_m(r)\, e^{\sqrt{-1}m\theta}$,
and thus we have the decompositions
$g=g_0+g^{(1)}$
and $\omega=\omega_0+\omega^{(1)}$.
In the cases ($n<0$) or ($n=0,a=0,k>1$),
we put
\begin{equation}
 \label{eq;08.1.26.5}
 u_n:=2r^n\int_{0}^r\rho^{-n-1}\, g_n(\rho)\, d\rho.
\end{equation}
In the cases ($n>0$), 
($n=0,-1<a<0$) or ($n=0,a=0,k<1$),
we put
\begin{equation}
 \label{eq;08.1.26.6}
 u_n:=-2r^{n}\int_{r}^R
 \rho^{-n-1}\, g_n(\rho)\, d\rho.
\end{equation}
We will not consider the cases
$(n=0,a>0)$ and $(n=0,a=0,k=1)$.
Then, we set
\begin{equation}
 \Phi^{(1)}(\omega):=
 \sum_{n\in\seisuu,n\neq 0}
 u_n(r)\, e^{\sqrt{-1}n\theta}.
\end{equation}
When $a\leq 0$ and $(a,k)\neq(0,1)$ 
are satisfied,
we also put
\begin{equation}
 \Phi(\omega):=
 \sum_{n\in\seisuu}
 u_n(r)\, e^{\sqrt{-1}n\theta}.
\end{equation}
The following proposition is proved 
in the proof of Proposition 6.4 and Proposition 11.5
of \cite{z}.
\begin{prop}
\label{prop;08.1.26.100}
Assume $R>0$ is sufficiently small.
Let $\omega$ be as above.
\begin{itemize}
\item
 We have $\delbar\Phi^{(1)}(\omega)=\omega^{(1)}$.
There exists a positive constant $C_1$,
such that 
$\|\Phi^{(1)}(\omega)\|_{h,g_{\poin}}
\leq C_1\, \|\omega^{(1)}\|_{h,g_{\poin}}$.
If we fix a compact subset $K_1$
of $\{-1<t\leq 1/2\}\times\seisuu$,
the constant $C_1$ can be taken independently
from $(a,k)\in K_1$.
(But it may depend on $K_1$.)
\item
Assume $a\leq 0$ and $(a,k)\neq(0,1)$.
We have $\delbar\Phi(\omega)=\omega$.
We also have a constant $C_2$
such that
$\|\Phi(\omega)\|_{h,g_{\poin}}
\leq C_2\, \|\omega\|_{h,g_{\poin}}$.
If we fix a compact subset $K_2$
of $\{-1<t<0\}\times\seisuu$,
the constant $C_2$ can be taken independently
from $(a,k)\in K_2$.
\hfill\qed
\end{itemize}
\end{prop}

Let $L^2(X^{\ast}(R),\nbigl)$ 
(resp. $L^2(X^{\ast}(R),\nbigl\otimes\Omega^{0,1})$)
denote the space of
$L^2$-sections of $\nbigl$ 
(resp. $\nbigl\otimes\Omega^{0,1}$)
on $X^{\ast}(R)$,
which are $L^2$ with respect to 
$h$ and $g_{\poin}$.
We obtain the following corollary.
\begin{cor}
\label{cor;08.1.26.2}
\mbox{{}}
\begin{itemize}
\item
$\Phi^{(1)}$ induces a bounded linear map
$L^2\bigl(X^{\ast}(R),\nbigl\otimes\Omega^{0,1}\bigr)
\lrarr L^2\bigl(X^{\ast}(R),\nbigl\bigr)$.
The range of $\Phi^{(1)}$ is contained in
the domain of $\delbar$,
and $\delbar\circ\Phi^{(1)}(\omega)=\omega^{(1)}$
 holds for any $\omega\in 
 L^2\bigl(X^{\ast}(R),\nbigl\otimes\Omega^{0,1}\bigr)$.
If we fix a compact subset $K_1$
of $\{-1<t\leq 1/2\}\times\seisuu$,
the norm of $\Phi^{(1)}$ is uniformly bounded
for $(a,k)\in K_1$.
\item
If $a\leq 0$ and $(a,k)\neq(0,1)$,
$\Phi$ induces a bounded linear map
$ L^2\bigl(X^{\ast}(R),\nbigl\otimes\Omega^{0,1}\bigr)
\lrarr L^2\bigl(X^{\ast}(R),\nbigl\bigr)$.
The range of $\Phi$ is contained in
the domain of $\delbar$,
and $\delbar\circ\Phi$ is the identity
of $L^2\bigl(X^{\ast}(R),\nbigl\otimes\Omega^{0,1}\bigr)$.
If we fix a compact subset $K_2$
of $\{-1<t<0\}\times\seisuu$,
the norm of $\Phi$ is uniformly bounded
for $(a,k)\in K_2$.
\hfill\qed
\end{itemize}
\end{cor}

\subsection{An estimate in \cite{s2}}

Let $f(z)\, d\zbar$ be a $(0,1)$-form on $\Delta^{\ast}$
such that $|f(z)|\leq |z|^{a}(-\log|z|)^{k}$
and that the support of $f$ is compact in $\Delta$.
We often need a solution $g$ of the equation
$\delbar g=f\, d\zbar$
satisfying some growth estimate
around the origin.
For that purpose,
we put
\[
 H(f)(z):=\int \frac{f(w)}{z-w}
\,
 \frac{\sqrt{-1}}{2\pi}dw\, d\wbar
\]

\begin{lem}\mbox{{}}
\label{lem;07.2.6.11}
\begin{itemize}
\item
In the case $-2<a<-1$,
we have $|H(f)|=O\bigl(|z|^{a+1}(-\log|z|)^k\bigr)$.
\item
In the case $(a=-2$, $k<-1)$
or $(a=-1$, $k>-1)$,
we have $|H(f)|=O\bigl(|z|^{a+1}(-\log|z|)^{k+1}\bigr)$.
\item
In the case $a=-1$ and $k=-1$,
 $|H(f)|=O\bigl(|z|^{a+1}(-\log|z|)^{k+1}
 \log(-\log|z|)\bigr)$.
\end{itemize}
\end{lem}
\pf
See Page 759--760 of \cite{s2}.
\hfill\qed

\subsection{Preliminary}
\label{subsection;08.1.27.10}

Let us start the proof of the propositions
in Section \ref{subsection;07.10.20.1}.
By an easy argument to use the decent,
we can reduce the problem to the unramified case.
Therefore, we may and will assume that
$(V_{\ast},\DDlambda)$ is unramified.
We use the polar coordinate $z=r\,e^{\sqrt {-1}\theta}$.
We may assume the following for the frame $\vecv$,
moreover:
\begin{enumerate}
\item $\vecv$ is compatible with
the irregular decomposition in $N$-th order
for some large $N$,
i.e., 
$\vecv_{|\Ohat^{(N)}}$
is compatible with the decomposition of
$V_{\ast|\Ohat^{(N)}}$ induced by
(\ref{eq;08.9.9.10}),
where $\Ohat^{(N)}$ denotes
the $N$-th infinitesimal neighbourhood of $O$.
\item $\vecv_{|O}$ is compatible with
the generalized eigen decomposition of $\Res(\DDlambda)$.
\item 
 Let $N_{a,\alpha,\gminia}$
 denote the nilpotent part of the endomorphisms
 on $\Gr^{F}_a\EE_{\alpha}\bigl(V_{\gminia|O}\bigr)$
 induced by $\Res(\DDlambda)$.
 Then,
 $N_{a,\alpha,\gminia}$ are represented
 by Jordan matrices with respect to the induced frames.
\end{enumerate}
We have the irregular value $\gminia(v_i)$,
and the eigenvalue $\alpha(v_i)$ of $\Res(\DDlambda)$
corresponding to $v_i$.
We also put $a(v_i):=\deg^F(v_i)$
and $k(v_i):=\deg^W(v_i)$.
We define
\begin{multline*}
 \nbigb(k):=
\bigl\{
 v_i\,\big|\,
 a(v_i)= \gminia(v_i)= \alpha(v_i)=0,\,\,
 k(v_i)=k
\bigr\} \\
\cup
\bigl\{
 v_i\,\big|\,
 a(v_i)=0,\,\,
 \bigl(\gminia(v_i),\alpha(v_i)\bigr)\neq (0,0)
 \bigr\}
\end{multline*}

Let $A$ be determined by 
$\DDlambda\vecv=\vecv\, A$.
Let $\Gamma$ be the diagonal matrix
whose $(i,i)$-entries are
$\alpha(v_i)\, dz/z+d\gminia(v_i)$.
We put $A_0:=A-\Gamma$.
We use the symbol $F_A$
to denote the section of
$\End(V)\otimes\Omega^{1,0}$
determined by 
$F_A(\vecv)=\vecv\, A$.
We use the symbol $F_{A_0}$ in a similar meaning.
Then, $F_{A_0}$ is bounded
with respect to $h$ and $g_{\poin}$.
We have the following decomposition:
\[
 A_0=\bigoplus_{a,\alpha,\gminia} 
 J_{a,\alpha,\gminia}\, \frac{dz}{z}
+A_0'\, \frac{dz}{z}
\]
Here $A_0'$ is holomorphic
and $F_{A_0'|O}$ strictly decreases the parabolic filtration.
And $J_{a,\alpha,\gminia}$ are
constant Jordan matrices
and represent
$N_{a,\alpha,\gminia}$
with respect to the induced frames of
$\Gr^F_a\EE_{\alpha}(V_{\gminia|O})$.

The $(1,0)$-operator $\del$ is defined by 
$\del\bigl(\sum f_i\, v_i\bigr)=
 \sum \del f_i\cdot v_i$.
Then, we have $\DDlambda=\delbar+\lambda\del+F_A$.

Let $\dvol_{g_{\poin}}$ denote 
the volume form of the Poincar\'e metric.
Recall that a section 
$\sum f_i\,v_i$ is $L^2$ if and only if
the following holds:
\[
 \sum\int
 |f_i|^2\,|z|^{-2a(v_i)}(-\log|z|^2)^{k(v_i)}
 \dvol_{g_{\poin}}
<\infty
\]
A section 
$\sum f_i\,v_i\, dz/z
+\sum g_j\,v_j\,d\zbar/\zbar$ is $L^2$
if and only if the following holds:
\[
 \sum\int
 |f_i|^2\,|z|^{-2a(v_i)}(-\log|z|^2)^{k(v_i)+2}
 \dvol_{g_{\poin}}<\infty
\]
\[
 \sum\int |g_j|^2\,|z|^{-2a(v_j)}(-\log|z|^2)^{k(v_j)+2}
 \dvol_{g_{\poin}}
<\infty
\]
A section 
$\sum f_i\,v_i\,dz\, d\zbar/|z|^2$ is $L^2$,
if and only if the following holds:
\[
 \sum\int
 |f_i|^2\,|z|^{-2a(v_i)}(-\log|z|^2)^{k(v_i)+4}
 \dvol_{g_{\poin}}<\infty
\]

\subsection{Vanishing of
 $\nbigh^2$ of $\nbigl^{\bullet}(V_{\ast},\DDlambda)$
 and $\nbigl^{\bullet}_{\poly}(V_{\ast},\DDlambda)$}
\label{subsection;07.9.12.1}

Let us consider Proposition \ref{prop;07.9.12.20}.
Let $\omega$ be an $L^2$-section of $V\otimes\Omega^2$.
We have the expression:
\[
 \omega=f\, \frac{dz\,d\zbar}{|z|^2}
\quad
 f=\sum f_i\, v_i
\]
Each $f_i$ has the Fourier expansion
$f_i=\sum_{m\in\seisuu}
 f_{i,m}(r)\, e^{\sqrt{-1}m\theta}$.
We set
\[
  \nbiga^{(0)}(f):=
 \sum_{v_i\in\nbigb(-1)}
 f_{i,0}(r)\, v_i,
\quad
 \nbiga^{(1)}(f):=
 f-\nbiga^{(0)}(f).
\]
We have the decomposition
$f=\nbiga^{(0)}(f)+\nbiga^{(1)}(f)$.
We have the corresponding decomposition
$\omega=\nbiga^{(0)}(\omega)+\nbiga^{(1)}(\omega)$.
Recall we have the following equalities:
\begin{equation}
 \label{eq;08.1.27.5}
 \del\bigl(f_{i,0}(r)\bigr)=
 \frac{1}{2}r\frac{\del f_{i,0}}{\del r}\,\frac{dz}{z},
\quad
 \delbar\bigl(f_{i,0}(r)\bigr)=
 \frac{1}{2}r\frac{\del f_{i,0}}{\del r}\, \frac{d\zbar}{\zbar}
\end{equation}
We show the following lemma based on 
Sabbah's idea contained in \cite{sabbah2}.
\begin{lem}
 \label{lem;07.11.16.7}
\mbox{{}}
\begin{itemize}
\item
 We have an $L^2$-section $\tau^{(1)}$ 
of $V\otimes\Omega^{1,0}$
 such that 
 $\delbar\tau^{(1)}=\nbiga^{(1)}(\omega)$.
\item
 We have an $L^2$-section $\tau^{(0)}$ 
 of $V\otimes\Omega^1$
 such that (i) $\DDlambda\tau^{(0)}$ is also $L^2$,
 (ii) $\nbiga^{(0)}
 \bigl(\omega-\DDlambda\tau^{(0)}\bigr)=0$.
\item
 In particular,
 we can take an $L^2$-section $\tau$ 
of $V\otimes\Omega^1$
 such that $\DDlambda\tau=\omega$.
\end{itemize}
If $\omega$ is $C^{\infty}$ and of polynomial order,
$\tau^{(i)}$ $(i=0,1)$,
$\DDlambda\tau^{(0)}$ and
$\tau$ are also $C^{\infty}$ 
and of polynomial order.
\end{lem}
\pf
The first claim follows from Corollary 
\ref{cor;08.1.26.2}.
Let us show the second claim.
We give the proof which also works 
in the family case.
We give preliminary arguments.

\subsubsection{(A)}
In the case 
$a(v_i)=0$, $\gminia(v_i)=0$ and $\alpha(v_i)\neq 0$,
we put
\begin{equation}
 \label{eq;08.1.27.2}
 \tau_1:=
 f_{i,0}\, v_i\, \frac{d\zbar}{\zbar}
+\lambda f_{i,0}\, v_i\,\frac{dz}{z}.
\end{equation}
Due to (\ref{eq;08.1.27.5}),
we have 
$\DDlambda\tau_1=F_A\bigl(f_{i,0}\, v_i\bigr)$.
Hence, we have the following:
\begin{multline}
\label{eq;07.7.22.2}
 f_{i,0}\, v_i\, \frac{dz\, d\zbar}{|z|^2}
-\alpha(v_i)^{-1}\DDlambda\tau_1
=
 f_{i,0}\, v_i\frac{dz\, d\zbar}{|z|^2}
-\alpha(v_i)^{-1}
 F_A\bigl(f_{i,0}\, v_i\bigr) \frac{d\zbar}{\zbar}
 \\
=\alpha(v_i)^{-1}F_{A_0}\bigl(f_{i,0}\, v_i\bigr)
 \frac{d\zbar}{\zbar}
=:\sum B_j\, v_j
 \frac{dz\,d\zbar}{|z|^2}
\end{multline}
Because 
$f_{i,0}\, v_i\, (dz/z)\, (d\zbar/\zbar)$ is $L^2$,
the sections
\[
 F_A(f_{i,0}\, v_i)\, d\zbar/\zbar,
\quad\quad
f_{i,0}\, v_i\, dz/z,
\quad\quad
f_{i,0}\, v_i\, d\zbar/\zbar
\]
are also $L^2$.
In particular,
$\tau_1$ and $\DDlambda\tau_1$ 
are $L^2$.
Because $F_{A_0}$ is bounded,
the right hand side of (\ref{eq;07.7.22.2})
are also $L^2$.
Let us look at $B_j$ more closely.
Because $A_0$ is holomorphic,
we have $B_j=\sum_{m\geq 0}B_{j,m}(r)
 \,e^{\sqrt{-1}m\theta}$.
If $a(v_j)=0$,
we have $B_{j,0}(r)=0$ unless
$\bigl(\gminia(v_j),\alpha(v_j)\bigr)=
 \bigl(\gminia(v_i),\alpha(v_i)\bigr)$
and $N_{a,\alpha,\gminia}v_{i|O}=v_{j|O}$.
Note $\deg^W(v_j)<\deg^W(v_i)$
for such $v_j$.

\subsubsection{(B)}

Let us consider the case 
$a(v_i)=0$ and $\gminia(v_i)\neq 0$.
Let $k$ be determined by
$\gminia(v_i)=\sum_{j=1}^k\gminia_j(v_i)z^{-j}$
and $\gminia_k(v_i)\neq 0$.
Recall we have the following:
\begin{equation}
 \del\Bigl(
 z^{k}f_{i,0}\, v_i\, \frac{d\zbar}{\zbar}
 \Bigr)
=z^k\del\bigl(f_{i,0}\, v_i\bigr)\frac{d\zbar}{\zbar}
+k\, z^k\, f_{i,0}\,
 v_i\frac{dz\, d\zbar}{|z|^2},
\end{equation}
\[
  \delbar\Bigl(
 z^k\, f_{0,i}\, v_i\,\frac{dz}{z}
 \Bigr)
=z^k\delbar\bigl(f_{i,0}\, v_i\bigr)\frac{dz}{z}
\]
We consider the following:
\begin{equation}
 \label{eq;07.11.16.6}
 \tau_1:=z^k\, f_{i,0}\, v_i\,\frac{d\zbar}{\zbar}
+\lambda\, z^k\, f_{i,0}\, v_i\, \frac{dz}{z}
\end{equation}
It is $L^2$,
and we have the following:
\begin{multline}
 \DDlambda(\tau_1)
=F_A\Bigl(z^k\, f_{i,0}\,v_i\,\frac{d\zbar}{\zbar}\Bigr)
+\lambda\, k\,
 z^k\, f_{i,0}\,v_i\frac{dz\, d\zbar}{|z|^2}
 \\
=\Bigl(
 z\frac{\del\gminia(v_i)}{\del z}
+\alpha(v_i)+k\,\lambda
 \Bigr)\, z^k\, f_{i,0}\, v_{i}
 \frac{dz\, d\zbar}{|z|^2}
+z^k\, F_{A_0}\bigl(f_{i,0}\, v_i\bigr)\,
 \frac{d\zbar}{\zbar}
\end{multline}
Hence, $\DDlambda(\tau_1)$ is also $L^2$.
Let $B_j$ be determined by the following:
\[
 f_{i,0}\, v_i\, \frac{dz\, d\zbar}{|z|^2}
-\frac{1}{-k\,\gminia_k(v_i)} \DDlambda\tau_1
=:\sum B_j\, v_j\, \frac{dz\,d\zbar}{|z|^2}
\]
We have 
$B_j=\sum_{m>0}
 B_{j,m}(r)\, e^{\sqrt{-1}m\theta}$.
It means
\[
\nbiga^{(0)}\Bigl(
 f_{i,0}\, v_i\, \frac{dz\, d\zbar}{|z|^2}
-\frac{1}{-k\,\gminia_k(v_i)} \DDlambda\tau_1
\Bigr)
=0. 
\]

\subsubsection{(C)}
Let us consider the case 
$a(v_i)=0$, $\gminia(v_i)=0$, $\alpha(v_i)=0$
and $k(v_i)=-1$.
Let $i(1)$ be determined by 
$N_{0,0,0} v_{i(1)|O}=v_{i|O}$
in $\Gr^F_0\EE_0(V_{0|O})$.
We put 
\begin{equation}
 \label{eq;07.11.16.10}
 \tau_1:=f_{i,0}\, v_{i(1)}\, \frac{d\zbar}{\zbar}
+\lambda\, f_{i,0}\, v_{i(1)}\, \frac{dz}{z}.
\end{equation}
It is $L^2$, and we have the following:
\[
  \DDlambda(\tau_1)
=f_{i,0}\, v_i\,\frac{dz\,d\zbar}{|z|^2}
+F_{A_0'}(\tau_1)
\]
Hence, $\DDlambda(\tau_1)$ is also $L^2$.
Let $B_j$ be determined by the following:
\[
 f_{i,0}\, v_i\,\frac{dz\, d\zbar}{|z|^2}
-\DDlambda(\tau_1)
=\sum B_j\, v_j\, \frac{dz\, d\zbar}{|z|^2}
\]
If $a(v_j)=0$,
we have 
$B_j=\sum_{m>0}
 B_{j,m}(r)\, e^{\sqrt{-1}m\theta}$.
In particular,
we have the following:
\[
 \nbiga^{(0)}\Bigl(
 f_{i,0}\, v_i\,\frac{dz\, d\zbar}{|z|^2}
-\DDlambda(\tau_1)
 \Bigr)=0
\]

\subsubsection{}

Let us show the second claim of Lemma
\ref{lem;07.11.16.7}.
We have
\[
 \nbiga^{(0)}(\omega)
=\sum_{v_i\in\nbigb(-1)}
 f_{i,0}\,v_i\frac{dz\,d\zbar}{|z|^2}
\]
Applying the procedure in (B) and (C),
we may and will assume that 
$f_{i,0}=0$ unless
$a(v_i)=0$, $\gminia(v_i)=0$
and $\alpha(v_i)\neq 0$.
Let $k_0$ be determined by
$\max\bigl\{k(v_i)\,\big|\,
 a(v_i)=0,\gminia(v_i)=0,
 \alpha(v_i)\neq 0
 \bigr\}$.
Applying the procedure in (A),
we can kill the coefficients 
in $\nbiga^{(0)}(\omega)$,
of $v_i$ with $a(v_i)=0$,
$\gminia(v_i)=0$ and $k(v_i)=k_0$.
By using an easy descending induction,
we can kill $\nbiga^{(0)}(\omega)$.
Thus, we obtain the second claim
of Lemma \ref{lem;07.11.16.7}.

\vspace{.1in}

Let us finish the proof of
Lemma \ref{lem;07.11.16.7}.
Assume that 
$\omega$ is $C^{\infty}$
and of polynomial order.
By construction, $\tau^{(0)}$
is $C^{\infty}$ and of polynomial
order by construction.
We can check that
$\DDlambda\tau^{(0)}$ is also
$C^{\infty}$ and of polynomial order
from its explicit description.
Because $\omega$ is
$C^{\infty}$ and of polynomial order
by construction,
$\nbiga^{(1)}(\omega)$ is
$C^{\infty}$ and of polynomial order.
Let $\nbiga^{(1)}_j(\omega)$
and $\tau^{(1)}_j$
denote the coefficients of $v_j$
in $\nbiga^{(1)}$ and $\tau^{(1)}$,
respectively.
We obtain that $\tau^{(1)}_j$ is $C^{\infty}$
on $X^{\ast}$
by the equation $\delbar\tau^{(1)}_j=
\nbiga_j^{(1)}(\omega)$ 
and the elliptic regularity of $\delbar$.
If $M$ is sufficiently large,
we have
(i) $\tau^{(1)}_j$ is $L^2$ on $X$,
(ii) $z^M\nbiga^{(1)}_j(\omega)$ is
 $L^{\infty}$ on $X$,
(iii) $\delbar(z^M\tau^{(1)}_j)
=z^M\nbiga^{(1)}_j(\omega)$
as a distribution on $X$.
By Sobolev's embedding,
we obtain that $z^M\tau^{(1)}_j$
is $L^p$ for some $p>2$.
Then, by using Sobolev's embedding again,
we obtain that $z^M\tau^{(1)}_j$
is $L^{\infty}$ on $X$.
Namely $\tau^{(1)}_j$
is of polynomial order.
Thus, the proof of Lemma \ref{lem;07.11.16.7}
is finished.
\hfill\qed

\subsection{$\nbigh^j(\psi_1\circ\varphi_0)$
 and $\nbigh^j(\varphi_0)$  for $j=0,1$}
\label{subsection;07.9.12.2}
Let us prove Proposition \ref{prop;07.9.12.20}.
Let $\omega$ be an $L^2$-section of 
$V\otimes\Omega^1$ such that $\DDlambda\omega=0$.
We have the expression
$\omega=f^{1,0}\, dz/z+f^{0,1}\, d\zbar/\zbar$.
We set
\[
 \nbiga^{(0)}(f^{0,1})
=\sum_{v_i\in\nbigb(1)}
 f^{(0,1)}_{i,0}(r)\, v_i,
\quad
 \nbiga^{(1)}(f^{0,1}):=
 f^{0,1}-\nbiga^{(0)}(f^{0,1})
\]
We have the decomposition
$f^{0,1}=\nbiga^{(0)}(f^{0,1})
 +\nbiga^{(1)}(f^{0,1})$.
We have the corresponding decomposition
$\omega^{0,1}=
\nbiga^{(0)}(\omega^{0,1})
+\nbiga^{(1)}(\omega^{0,1})$.

\begin{lem}
\label{lem;07.9.11.10}
\mbox{{}}
\begin{itemize}
\item
 We have an $L^2$-section $\tau^{(1)}$ of $V$
 such that 
 $\delbar\tau^{(1)}=\nbiga^{(1)}(\omega^{0,1})$.
\item
 We have an $L^2$-section $\tau^{(0)}$ of $V$
 such that (i) $\delbar\tau^{(0)}$ is also $L^2$,
 (ii) $\nbiga^{(0)}(\omega^{0,1}-\delbar\tau^{(0)})=0$.
\end{itemize}
As a result,
we can take an $L^2$-section $\tau$ of $V$
such that $\delbar\tau=\omega^{0,1}$.
If $\omega^{0,1}$ is $C^{\infty}$ 
and of polynomial order,
$\tau^{(i)}$ $(i=0,1)$,
$\delbar\tau^{(0)}$
and $\tau$ are also $C^{\infty}$ and of 
polynomial order.
\end{lem}
\pf
The first claim follows from Corollary \ref{cor;08.1.26.2}.
Let $C_j$ be the functions determined by the following:
\begin{equation}
\label{eq;10.5.21.40}
 F_{A_0}\Bigl(f^{0,1}d\zbar/\zbar\Bigr)
=\sum C_j\, v_j\, \frac{dz\, d\zbar}{|z|^2}
\end{equation}
From $\DDlambda\omega=0$,
we obtain the following relation
by considering the $v_i$-component:
\begin{equation}
\label{eq;07.7.22.3}
 \lambda\del f_i^{0,1}\, v_i\frac{d\zbar}{\zbar}
+f_i^{0,1}\, 
\Bigl(d\gminia(v_i)+\alpha(v_i)\frac{dz}{z}\Bigr)
 \, v_i\, \frac{d\zbar}{\zbar}
+\delbar f_i^{1,0}\, v_i\, \frac{dz}{z}
+C_{i}\, v_i\, \frac{dz\, d\zbar}{|z|^2}=0
\end{equation}
We use the Fourier expansion
$C_j=\sum C_{j,m}\, e^{\sqrt{-1}m\theta}$.
We give some preliminary arguments.

\vspace{.1in}

\noindent
{\bf (A)} 
Let us consider the case $a(v_i)=0$ and $\gminia(v_i)\neq 0$.
Let $k$ be determined by 
$\gminia(v_i)=\sum_{j=1}^k
 \gminia_j(v_i)\, z^{-j}$
with $\gminia_k(v_i)\neq 0$.
Let us look at
the $e^{-\sqrt{-1}k\theta}$-component
of (\ref{eq;07.7.22.3}).
Multiplying it by $r^k$, we obtain the following,
where we omit to write $dz\, d\zbar/|z|^2$:
\begin{multline}
\label{eq;07.11.16.20}
 -k\, f^{0,1}_{i,0}\,\gminia_k(v_i)
-\sum_{\substack{0<m<k \\ m+j=k}}
 f^{0,1}_{i,-m}\, j\, \gminia_{j}(v_i)\, r^{k-j}
+r^k\alpha(v_i)\, f^{0,1}_{i,-k}
+\frac{1}{2}r\frac{\del}{\del r}\bigl(
 \lambda r^k \, f^{0,1}_{i,-k}\bigr) \\
-k\lambda f^{0,1}_{i,-k}\, r^k
-\frac{1}{2}r\frac{\del} {\del r}
 \bigl(r^k\, f^{1,0}_{i,-k}\bigr)
+r^k\, C_{i,-k}=0
\end{multline}
We consider the following:
\begin{equation}
 \label{eq;07.11.16.19}
 \rho:=\left(
 -\sum_{\substack{0<m<k\\ m+j=k}}
 f^{0,1}_{i,-m}\, j\,\gminia_{j}(v_i)\, r^{k-j}
+r^k\,(\alpha(v_i)-k\lambda)
 \, f^{0,1}_{i,-k}
+r^kC_{i,-k}
\right)\, v_i\, \frac{d\zbar}{\zbar}
\end{equation}
Then, we have
$\int|\rho|_h^2\, r^{-2\epsilon}\,\dvol_{g_{\poin}}<\infty$
for some $\epsilon>0$.
By Corollary \ref{cor;08.1.26.2},
we can take $\rho_1$ such that
$\int|\rho_1|_h^2\, r^{-2\epsilon}\dvol_{g_{\poin}}<\infty$
and $\delbar\rho_1=\rho$.
Note that we have the following:
\begin{equation}
 \label{eq;07.11.16.22}
 \delbar\Bigl(
 (r^k\lambda f^{0,1}_{i,-k}-r^k f^{1,0}_{i,-k})\, v_i
 \Bigr)
=\frac{1}{2}r\frac{\del}{\del r}
 \bigl(r^k\lambda f^{0,1}_{i,-k}- r^k f^{1,0}_{i,-k}\bigr)
 \, v_i\, \frac{d\zbar}{\zbar}
\end{equation}
Hence, we have an $L^2$-section $\tau_2$
such that $f^{0,1}_{i,0}v_i\, d\zbar/\zbar=\delbar\tau_2$.

\vspace{.1in}
\noindent
{\bf (B)}\,\,
Let us consider the case $a(v_i)=0$, $\gminia(v_i)=0$
and $\alpha(v_i)\neq 0$.
Let us look at 
the $e^{\sqrt{-1}0\theta}$-component 
of (\ref{eq;07.7.22.3}).
We have the following:
\[
 \frac{\lambda}{2}r\frac{\del f^{0,1}_{i,0}}{\del r}
 \frac{dz\, d\zbar}{|z|^2}
+\alpha(v_i)\, f^{0,1}_{i,0}\frac{dz\, d\zbar}{|z|^2}
-\frac{1}{2}r\frac{\del f^{1,0}_{i,0}}{\del r}\,\frac{dz\,
d\zbar}{|z|^2}
+C_{i,0}\, \frac{dz\, d\zbar}{|z|^2}=0
\]
Hence, we have the following:
\[
\frac{1}{2}r\frac{\del}{\del r} 
\bigl(
 \lambda f^{0,1}_{i,0}-f^{1,0}_{i,0}
\bigr)
+\alpha(v_i)\, f^{0,1}_{i,0}
=\left\{
 \begin{array}{ll}
 -f_{i(1),0}^{0,1}+R & (\mbox{\rm if } 
\exists v_{i(1)},\,N_{\gminia,\alpha,a}v_{i(1)}=v_i)\\
  R & \mbox{\rm otherwise}
 \end{array}
 \right.
\]
Here,
$\int |R\, v_i|_h^2\, r^{-\epsilon}
 \,\dvol_{g_{\poin}}<\infty$.
Then, by using an easy inductive argument,
we can show that there exists $\tau_2$ such that
$\delbar\tau_2=f_{i,0}^{0,1}\, v_i\, d\zbar/\zbar$.

\vspace{.1in}
\noindent
{\bf (C)}\,\,
Let us consider the case $a(v_i)=0$,
$\gminia(v_i)=0$, $\alpha(v_i)=0$ and $k(v_i)=1$.
Let $i(-1)$ be determined by $Nv_{i|O}=v_{i(-1)|O}$.
Let us look at the $e^{\sqrt{-1}0\theta}$-component
of (\ref{eq;07.7.22.3}) for $v_{i(-1)}$.
We have the following:
\[
 \lambda\frac{r}{2}\frac{\del f^{0,1}_{i(-1),0}}{\del r}
 \frac{dz\, d\zbar}{|z|^2}
-\frac{r}{2}\frac{\del}{\del r}f^{1,0}_{i(-1),0}
 \frac{dz\, d\zbar}{|z|^2}
+f^{0,1}_{i,0}\frac{dz\, d\zbar}{|z|^2}
+\bigl(C_{i(-1),0}-f^{0,1}_{i,0}\bigr)
 \frac{dz\, d\zbar}{|z|^2}=0
\]
Hence,
we have the following:
\[
 f^{0,1}_{i,0}
=\frac{1}{2}r\frac{\del}{\del r}
 \bigl(-\lambda f^{0,1}_{i(-1),0}+f^{1,0}_{i(-1),0}\bigr)
+R
\]
Here,
$\int |R\, v_{i}|^2\,
 |z|^{-\epsilon}\,\dvol_{g_{\poin}}<\infty$
for some $\epsilon>0$.
We also have the following:
\[
 \bigl|f^{0,1}_{i(-1),0}\, v_i\bigr|_h^2
\sim
 \bigl|f^{0,1}_{i(-1),0}\, v_{i(-1)}\, 
 d\zbar/\zbar
 \bigr|_{h,g_{\poin}}^2
\]
Hence, we obtain that
$f^{0,1}_{i(-1),0}\, v_i$ is $L^2$.
Similarly,
$f^{1,0}_{i(-1),0}\, v_i$ is also $L^2$.
Thus, there exists an $L^2$-section $\tau_2$
such that 
$\delbar\tau_2=
 f^{0,1}_{i,0}\, v_i\, d\zbar/\zbar$.

\vspace{.1in}
The second claim of 
Lemma \ref{lem;07.9.11.10}
follows from
the above considerations (A), (B), (C).
Assume that $\omega$ is $C^{\infty}$
and of polynomial order.
By the argument in the proof of
Lemma \ref{lem;07.11.16.7},
we can show that
$\tau^{(1)}$ is $C^{\infty}$
and of polynomial order.
In (A),
if $\omega$ is $C^{\infty}$ and
of polynomial order,
$\rho$, $f_{i,-k}^{0,1}$,
$f_{i,-k}^{1,0}$ are $C^{\infty}$
and of polynomial order,
and we can show that $\rho_1$ is
$C^{\infty}$ and of polynomial order
by the argument in the proof of
Lemma \ref{lem;07.11.16.7}.
Hence, $\tau_2$ in (A) is
$C^{\infty}$ and of polynomial order.
We can show that $\tau_2$ in (B), (C)
are $C^{\infty}$ and of polynomial order.
Then, $\tau^{(0)}$
and $\delbar\tau^{(0)}$
in Lemma \ref{lem;07.9.11.10}
are $C^{\infty}$ and of polynomial order.
\hfill\qed

\vspace{.1in}

We put $\rho:=\omega-\DDlambda\tau$
which is a holomorphic section of 
$V\otimes\Omega^{1,0}$
on $X^{\ast}$.
We have the decomposition $\rho=\sum \rho_i$,
where each $\rho_i$ is the product of
$v_i$ and a holomorphic $(1,0)$-form
on $X^{\ast}$.
\begin{lem}
Let $l(v_i)\in\seisuu_{\geq \,0}$ be determined
as follows:
\begin{itemize}
\item
We put $l(v_i):=-\ord(\gminia(v_i))+1$
in the case $\gminia(v_i)\neq 0$.
\item
We put $l(v_i):=1$ in the case
$\gminia(v_i)=0$ and 
$\alpha(v_i)+\lambda a(v_i)\neq 0$.
\item
We put $l(v_i):=0$ otherwise.
\end{itemize}
Then, $z^{l(v_i)}\,\rho_i$
is $L^2$
with respect to $h$ and $g_{\poin}$.
In the second case,
$(-\log|z|)^{-1}\rho_i$ is $L^2$,
more strongly.

In particular,
$\rho$ is $C^{\infty}$
and of polynomial order.
\end{lem}
\pf
Let $\delta'$ denote the $(1,0)$-operator
determined by $h$ and $\delbar$.
Let $B$ be determined by $\delta'\vecv=\vecv\, B$.
Then, $B$ is diagonal,
and the $(i,i)$-entries are as follows:
\[
 -a(v_i)\frac{dz}{z}+\frac{k(v_i)}{-2\log|z|}\frac{dz}{z}
\]
The curvature $R(h)$ of $\delbar+\delta'$ is bounded
with respect to $h$ and $g_{\poin}$.
Hence, $\delta'\tau$ is also $L^2$.
(See the argument 
in the proof of Lemma \ref{lem;07.7.18.60},
for example.)

Let $\DD^{\lambda(1,0)}$ denote 
the $(1,0)$-part of $\DDlambda$.
We put $G:=\DD^{\lambda(1,0)}-\lambda\delta'$,
which is a section of $\End(V)\otimes\Omega^{1,0}$.
Let $A_1$ be determined by
$G\,\vecv=\vecv\, A_1$.
Then, we have the decomposition $A_1=\Gamma'+C$,
where $F_C$ is bounded with respect to $h$ and $g_{\poin}$,
and $\Gamma'$ is the diagonal matrix whose $(i,i)$-entry is
as follows:
\[
 d\gminia(v_i)+\bigl(\alpha(v_i)+\lambda a(v_i)\bigr)
 \frac{dz}{z}
\]

We have the decomposition
$\rho=\omega^{1,0}-\lambda\delta'\tau
 -\lambda F_C(\tau)-\lambda F_{\Gamma'}(\tau)$.
Note 
$\omega^{1,0}-\lambda\delta'\tau-
\lambda F_C(\tau)$ is $L^2$.
Then, the claim of the lemma follows.
\hfill\qed

\vspace{.1in}

Let $A_0$ and $\Gamma$ be as in
Section \ref{subsection;08.1.27.10}.
We put $\DDlambda_0:=\DDlambda-F_{\Gamma}$.
We have
$\DDlambda_0\vecv=\vecv\, A_0$.
Recall $F_{A_0}$ is bounded
with respect to $h$ and $g_{\poin}$.

\begin{itemize}
\item
In the case $\gminia(v_i)\neq 0$,
we have the $L^2$-holomorphic section $\kappa_i$
such that
$\rho_i=\bigl(d\gminia(v_i)+
 \alpha(v_i)\, dz/z
 \bigr)\, \kappa_i$.
Note $\DDlambda_0(\kappa_i)$ is also $L^2$.
\item
In the case $\gminia(v_i)=0$ and $a(v_i)<0$,
we have $(-\log|z|)^{-1}\,\rho_i$ is $L^2$.
Then, we obtain the $L^2$-property of $\rho_i$.
See Section \ref{subsection;10.5.21.30}.
\item
In the case $\gminia(v_i)=0$, $a(v_i)=0$
and $\alpha(v_i)\neq 0$,
we have $z\,\rho_i$ is $L^2$.
Hence, we have the $L^2$-holomorphic section
$\kappa_i$ such that
$\alpha(v_i)\, \kappa_i\, dz/z=\rho_i$.
Note $\DDlambda_0(\kappa_i)$ is also $L^2$.
\item
In the case $\gminia(v_i)=0$, $a(v_i)=0$ and $\alpha(v_i)=0$,
we have $\rho_i$ is $L^2$.
\end{itemize}

Hence, we obtain the following lemma.
\begin{lem}
\label{lem;07.7.22.10}
There exists an $L^2$-section $\nu$ of $V$
such that (i) $\DDlambda\nu$ is also $L^2$,
(ii) $\omega-\DDlambda\nu$ is a 
holomorphic $(1,0)$-form.
If $\omega^{0,1}$ is $C^{\infty}$ 
and of polynomial order,
$\nu$ is also $C^{\infty}$ and 
of polynomial order.
\hfill\qed
\end{lem}

Let us finish the proof of 
Proposition \ref{prop;07.9.12.20}.
Lemma \ref{lem;07.9.11.10} and
Lemma \ref{lem;07.7.22.10}
imply that $\psi_1\circ\varphi_0$
is a quasi isomorphism.

Let $\omega\in
 \nbigl^{\bullet}_{\poly}(V_{\ast},\DDlambda)$.
We take $\nu$ as in Lemma \ref{lem;07.7.22.10}.
Because $\omega-\DDlambda\nu$
is $C^{\infty}$ and of polynomial order,
we obtain that
$\DDlambda\nu$ is also $C^{\infty}$
and of polynomial order.
Hence, 
$\nu\in\nbigl_{\poly}^0(V_{\ast},\DDlambda)$.
It implies that
$\nbigh^1(\varphi_0)$ is surjective.
The injectivity of $\nbigh^1(\varphi_0)$
follows from the injectivity of
$\nbigh^1(\psi_1\circ\varphi_0)$.
It is easy to see $\nbigh^0(\varphi_0)$
is an isomorphism
by Lemma \ref{lem;10.5.22.10}.
Thus Proposition \ref{prop;07.9.12.20}
is proved.
\hfill\qed

\subsection{Proof of
Proposition \ref{prop;07.9.12.21}}
\label{subsection;07.9.12.11}

By using an easy argument to use the descent
with respect to the ramified covering $X'\lrarr X$,
we may and will assume that $(V_{\ast},\DDlambda)$
is unramified.
We have the full reduction
$\bigl(
 \Gr^{\nbigftilde}_{\gminia}(\prolong{V}),
 \DDlambda_{\gminia}\bigr)$
for each $\gminia\in\Irr(\DDlambda)$.
(See Section \ref{subsection;10.5.21.10}.)
Let $\vecvbar_{\gminia}$ be
a holomorphic frame of
$\Gr^{\nbigftilde}_{\gminia}(\prolong{V})$
such that 
(i) compatible with the induced parabolic structure,
the generalized eigen decomposition
of $\Res(\DDlambda_{\gminia})$,
(ii) the induced frame of
$\Gr^F\Gr^{\nbigftilde}_{\gminia}(\prolong{V})$
is compatible with the weight filtration 
of the nilpotent part of $\Res(\DDlambda)$.
Let $R_{\gminia}$ be determined by
$\DDlambda_{\gminia}
\vecvbar_{\gminia}
=\vecvbar_{\gminia}\, (d\gminia+R_{\gminia})$.

Let $X^{\ast}=\bigcup_{j=1}^M S_j$ be a covering
by small sectors.
We have the full Stokes filtration $\nbigftilde^{S_j}$
and a flat splitting $V_{|\Sbar_j}
=\bigoplus_{\gminia}V_{\gminia,S_j}$.
We take the lift of $\vecvbar_{\gminia}$
to $V_{\gminia,S_j}$.
By gluing them as in Section \ref{subsection;07.6.16.8},
we obtain a $C^{\infty}$-frame
$\vecv_{C^{\infty}}=
 \bigl(\vecv_{\gminia,C^{\infty}}\bigr)$.
Let $V_{\gminia,C^{\infty}}$ 
denote the subbundle
generated by $\vecv_{\gminia,C^{\infty}}$.
We obtain a decomposition 
$V_{|X^{\ast}}=
 \bigoplus V_{\gminia,C^{\infty}}$.
Let $\tau$ be a local $C^{\infty}$-section
of $V\otimes \Omega^p$,
which has the corresponding decomposition
$\tau=\sum \tau_{\gminia}$.
Then, $\tau$ is a section of
$\nbigltilde^{\bullet}_{\poly}(V_{\ast},\DDlambda)$, 
if and only if the following estimate holds
for $h$ and $g_{\poin}$:
\begin{itemize}
\item
 $\tau$ and $\DDlambda\tau$ are
 of polynomial order.
\item
 $\tau_{0,C^{\infty}}$ and 
 $(\DDlambda\tau)_{0,C^{\infty}}$
 are $L^2$.
\end{itemize}

Let $C$ be determined by the following:
\[
 \DDlambda\vecv_{C^{\infty}}
=\vecv_{C^{\infty}}
 \,\left(
 \bigoplus \bigl(d\gminia+R_{\gminia}\bigr)+C
 \right)
\]
Then, we have $C=O\bigl(|z|^N\bigr)$ for any $N>0$.
Let $\DD^{\lambda\prime}$ be 
the flat $\lambda$-connection
determined by the following:
\[
 \DD^{\lambda\prime}\vecv_{C^{\infty}}
=\vecv_{C^{\infty}}
 \,\left(
 \bigoplus\bigl(d\gminia+R_{\gminia}\bigr)
 \right)
\]
The $(0,1)$-part of $\DD^{\lambda\prime}$ is denoted by
$\delbar'$.
We put $F:=\DDlambda-\DD^{\lambda\prime}$,
and then we have $|F|_h=O\bigl(|z|^N\bigr)$ for any $N>0$.
We obtain a complex
$\nbigltilde_{\poly}^{\bullet}
 (V_{\ast},\DD^{\lambda\prime})$
from $\DD^{\lambda\prime}$
as in Section \ref{subsection;07.9.11.11}.
As a sheaf, we have
$\nbigltilde_{\poly}^{p}(V_{\ast},\DD^{\lambda\prime})
=\nbigltilde_{\poly}^p(V_{\ast},\DDlambda)$
for $p=0,1,2$.

\begin{lem}
\label{lem;07.9.12.4}
For any $\omega\in
 \nbigltilde_{\poly}^2
 (V_{\ast},\DD^{\lambda\prime})
=\nbigltilde_{\poly}^2(V_{\ast},\DDlambda)$,
we can take $\tau\in
 \nbigltilde^1_{\poly}(V_{\ast},\DD^{\lambda})$
such that $\DD^{\lambda\prime}\tau=\omega$.
Note that 
$\nbigltilde_{\poly}^{p}(V_{\ast},\DD^{\lambda\prime})
=\nbigltilde_{\poly}^p(V_{\ast},\DDlambda)$
as remarked above.
\end{lem}
\pf
We have only to consider the case
where $\DD^{\lambda}$ 
has a unique irregular value
$\gminia$.
(Note that $\DD^{\lambda\prime}$
has the same irregular value.)
In the case $\gminia=0$,
we may apply the results
in Section \ref{subsection;07.9.12.1}.
In the case $\gminia\neq 0$,
we have only to use Lemma \ref{lem;07.2.6.11},
for example.
\hfill\qed

\begin{lem}
\label{lem;07.11.16.1}
For any $\omega\in\nbigltilde_{\poly}^2(V_{\ast},\DDlambda)$,
we can take $\tautilde\in\nbigltilde^1_{\poly}(V_{\ast},\DDlambda)$
such that $\DD^{\lambda}\tautilde=\omega$.

In particular,
we obtain the vanishing of $\nbigh^2$
of the complex of sheaves
$\nbigltilde_{\poly}^{\bullet}(V_{\ast},\DDlambda)$
at $O$.
\end{lem}
\pf
Take $\tau$ as in Lemma \ref{lem;07.9.12.4}.
We have 
$\omega-\DD^{\lambda}\tau=O(|z|^N)$ for any $N$.
Take some large $M$.
According to Lemma \ref{lem;07.2.6.11},
we can take
a section $\kappa$ of $V\otimes\Omega^{1,0}$
such that
(i) $\delbar\kappa=\omega-\DD^{\lambda}\tau$,
(ii) $|\kappa|=O(|z|^M)$.
We have
$\kappa\in\nbigltilde_{\poly}^1(V_{\ast},\DDlambda)$,
and 
$\DDlambda(\tau+\kappa)=\omega$.
Thus, we obtain Lemma \ref{lem;07.11.16.1}
\hfill\qed

\vspace{.1in}

Let $\omega\in\nbigltilde_{\poly}^1(V_{\ast},\DDlambda)$
such that $\DDlambda\omega=0$.
We have
$\DD^{\lambda\prime}\omega=
  -F\, \omega=O(|z|^N)$
for any $N$.
Hence, we can take a large $M>0$
and a $C^{\infty}$-section $\kappa$
of $V\otimes\Omega^{1,0}$
such that $\delbar'\kappa=\DD^{\lambda\prime}\omega$
and $\bigl|\kappa\bigr|=O\bigl(|z|^M\bigr)$.
We put $\omega':=\omega-\kappa$,
and then $\DD^{\lambda\prime}\omega'=0$.
Note that the $(0,1)$-part of $\omega'$
and $\omega$ are equal.
\begin{lem}
\label{lem;07.9.12.2}
There exists
a local section 
 $\tau\in\nbigltilde^0(V_{\ast},\DDlambda)$
around $O$
such that
$\delbar'\tau=\omega^{0,1}$.
\end{lem}
\pf
We may assume that
$\DD^{\lambda}$ 
(and hence $\DD^{\lambda\prime}$)
has a unique irregular value $\gminia$.
In the case $\gminia=0$,
we can apply the result in
Section \ref{subsection;07.9.12.2}.
In the case $\gminia\neq 0$,
we can take $\tau$ such that
(i) $\delbar'\tau=\omega^{0,1}$,
(ii) $|\tau|=O(|z|^{-M})$ for some large $M$.
Let us show that $\DD^{\lambda\prime}\tau$
is of polynomial order.
Let $h'$ be the $C^{\infty}$-hermitian metric
of $V_{|X-O}$ such that
$h'\bigl(v_{C^{\infty},i},v_{C^{\infty},j}\bigr)=\delta_{i,j}$.
Note that $h'$ and $h$ are mutually bounded
up to polynomial order.
Let $\delta_1'$ be the $(1,0)$-operator
determined by $\delbar'$ and $h'$.
Note that $z^{M_1}\,\omega^{0,1}$ 
and $z^{M_1}\tau$ are
bounded for some large $M_1$,
and that
$\delbar'(z^{M_1}\tau)=z^{M_1}\omega^{0,1}$.
Then, it can be shown that
$\delta'_1(z^{M_1}\tau)$ is $L^2$.
(See the argument in the proof of 
Lemma \ref{lem;07.7.18.60},
 for example.)
Thus, we obtain that
$z^{M_1}\, \delta_1'\tau$  is $L^2$.
Taking large $M_2$,
we obtain $z^{M_2}\,(\delta_1'\tau-\omega^{1,0})$
is also $L^2$.

Since $\omega^{1,0}-\delta_1'\tau$ is holomorphic
with respect to $\delbar'$,
we obtain
$\delta_1'\tau-\omega^{1,0}
=O(|z|^{-M_3})$ for some large $M_3$.
Then, we obtain the desired estimate for
$\delta_1'\tau$.
\hfill\qed

\begin{lem}
\label{lem;07.9.12.3}
We can take a section
$\tautilde\in\nbigltilde^0(V_{\ast},\DDlambda)$
such that
$\delbar\tautilde=\omega^{0,1}$.
\end{lem}
\pf
Let $\tau$ be as in Lemma \ref{lem;07.9.12.2}.
We have
$\DDlambda\tau-\DD^{\lambda\prime}\tau
=F\, \tau=O(|z|^N)$ for any $N>0$.
According to Lemma \ref{lem;07.2.6.11},
we can take a section $\nu$ of $V$
such that
(i) $|\nu|=O\bigl(|z|^M\bigr)$
 for some large $M>0$,
(ii) $\delbar\nu=F^{0,1}\,\tau$.
Let $\delta'$ be the $(1,0)$-operator determined by
$h$ and $\delbar$.
Then, $\delta'\nu$ is $L^2$.
If $M$ is sufficiently large,
$(\DD^{\lambda(1,0)}-\lambda\delta')\nu$ is
 $O(|z|^{M/2})$.
We put
$\rho=\sum\rho_{\gminia}
 :=\omega-\kappa-\DDlambda\tau+\DDlambda\nu$.
Then,
$\rho$ is a holomorphic section of $V\otimes\Omega^{1,0}$,
$z^{L}\rho_{\gminia}$ is $L^2$ for some $L$
 for any $\gminia$,
and $\rho_0$ is $L^2$.
Hence, 
we obtain
$\tautilde:=
 \tau+\nu\in
 \nbigltilde^0_{\poly}(V_{\ast},\DDlambda)$.
Thus, Lemma \ref{lem;07.9.12.3}
is proved.
\hfill\qed

\vspace{.1in}

Let $\nbigstilde(V_{\ast}\otimes\Omega^{p,0})$ be
the sheaf of meromorphic sections $\tau$
of $V\otimes\Omega^{p,0}$ with the following property:
\begin{itemize}
\item Let $\tau_{|\Ohat}=\tauhat_{\reg}+\tauhat_{\irr}$ 
 be the decomposition corresponding
 to the irregular decomposition.
 Then, $\tauhat_{\reg}$ is contained in
 $\nbigs\bigl(\Vhat_{\reg\ast}\otimes\Omega^{p,0}\bigr)$.
\end{itemize}

By using Lemma \ref{lem;07.11.16.1}
and Lemma \ref{lem;07.9.12.3},
it is easy to show that
the natural inclusion
$\nbigstilde\bigl(V_{\ast}\otimes\Omega^{\bullet,0}\bigr)
 \lrarr
 \nbigltilde_{\poly}^{\bullet}(V_{\ast},\DDlambda)$
is a quasi isomorphism.
It is also standard and easy to show that
the natural inclusion
$\nbigs(V_{\ast}\otimes\Omega^{\bullet,0})\lrarr
 \nbigstilde(V_{\ast}\otimes\Omega^{\bullet,0})$
 is a quasi isomorphism.
Hence, we obtain that
$\psi_2\circ\varphi_0$ is a quasi isomorphism.
Thus, the proof of Proposition \ref{prop;07.9.12.21}
is finished.
\hfill\qed

\subsection{Proof of Proposition \ref{prop;07.10.22.55}}
\label{subsection;08.9.9.1}

Let $\pi:\Xtilde(O)\lrarr X$ denote the projection.
For an open subset $U$ of $\Xtilde(O)$,
let $\nbigltilde^p_{\poly}
 (V_{\ast},\DDlambda)_{\Xtilde(O)}(U)$
be the space of $C^{\infty}$-sections $\tau$ 
of $V\otimes\Omega_X^p$ on 
$U\setminus \pi^{-1}(O)$,
such that $\tau_{|U\cap S}$ and 
$\DDlambda\tau_{|U\cap S}$
satisfy the conditions {\bf (a1)} and {\bf (a2)}
with respect to $h$ and $g_{\poin}$
on each small sector $S\cap U$.
By taking sheafification,
we obtain a complex of sheaves
$\nbigltilde^{\bullet}_{\poly}
 (V_{\ast},\DDlambda)_{\Xtilde(O)}$
on $\Xtilde(O)$.
Similarly, we obtain a complex of sheaves
$\nbiglbar^{\bullet}_{\poly}
 (V_{\ast},\DDlambda)_{\Xtilde(O)}$
on $\Xtilde(O)$
from {\bf (b1)} and {\bf (b2)}.

\begin{lem}
\label{lem;08.9.9.2}
The natural inclusion
$(\psi_2\circ\varphi_1)_{\Xtilde(O)}:
 \nbiglbar^{\bullet}_{\poly}
 (V_{\ast},\DDlambda)_{\Xtilde(O)}
\lrarr
 \nbigltilde^{\bullet}_{\poly}
 (V_{\ast},\DDlambda)_{\Xtilde(O)}$
is a quasi-isomorphism.
\end{lem}
\pf
By choosing a flat splitting,
we may assume that $(V_{\ast},\DDlambda)$
has a unique irregular value $\gminia$.
In the case $\gminia=0$,
the claim is trivial.
By using the results
in Section \ref{subsection;07.11.16.100},
we can show the vanishing of the higher cohomology sheaves
of both $\nbigltilde^{\bullet}_{\poly}(V_{\ast},h)$
and $\nbiglbar^{\bullet}_{\poly}(V_{\ast},h)$
on $\Xtilde(O)$.
The comparison of the $0$-th cohomology sheaves is easy.
\hfill\qed

\vspace{.1in}
Note that
$\nbiglbar^{\bullet}_{\poly}
 (V_{\ast},\DDlambda)_{\Xtilde(O)}$
and
$\nbigltilde^{\bullet}_{\poly}
 (V_{\ast},\DDlambda)_{\Xtilde(O)}$
are $c$-soft in the sense of
Definition 2.5.5 of \cite{ks},
and that
$\nbiglbar^{\bullet}_{\poly}
 (V_{\ast},\DDlambda)$
and
$\nbigltilde^{\bullet}_{\poly}
 (V_{\ast},\DDlambda)$
are obtained as their push-forward.
Hence, we obtain Proposition \ref{prop;07.10.22.55}
from Lemma \ref{lem;08.9.9.2}.
\hfill\qed

\subsection{Proof of Proposition
 \ref{prop;07.7.22.35}}
\label{subsection;08.9.9.5}

We have only to consider the first morphism.
We use the notation in 
Section \ref{subsection;08.9.9.1}.
We obtain a complex of sheaves
$\nbiglbar_{\poly}^{\bullet}
 (V_{\ast},h,h^{(T)})_{\Xtilde(O)}$
on $\Xtilde(O)$
from the conditions {\bf (c1)} and {\bf (c2)}.

\begin{lem}
\label{lem;08.9.9.6}
The inclusion 
$\iota:\nbiglbar_{\poly}^{\bullet}
 (V_{\ast},h,h^{(T)})_{\Xtilde(O)}
\lrarr
\nbiglbar^{\bullet}_{\poly}
 (V_{\ast},\DDlambda)_{\Xtilde(O)}$
is a quasi isomorphism.
\end{lem}
\pf
We may assume that $(V_{\ast},\DDlambda)$
has the unique irregular value $\gminia$.
In the case $\gminia=0$,
the claim is trivial.
By using Lemmas \ref{lem;07.7.22.20}
and \ref{lem;07.7.22.21},
we can show the vanishing of the higher cohomology sheaves
$\nbiglbar^{\bullet}_{\poly}(V_{\ast},h)$
on $\Xtilde(O)$.
By the same argument as that
in the proof of Lemmas
\ref{lem;07.7.22.20} and \ref{lem;07.7.22.21},
we can also show the vanishing of 
the higher cohomology sheaves
$\nbiglbar^{\bullet}_{\poly}(V_{\ast},h,h^{(T)})$
on $\Xtilde(O)$.
The comparison of 
the $0$-th cohomology sheaves is easy.
Thus, we obtain Lemma \ref{lem;08.9.9.6}.
\hfill\qed

\vspace{.1in}
Because
$\nbiglbar_{\poly}^{\bullet}
 (V_{\ast},h,h^{(T)})_{\Xtilde(O)}$
and
$\nbiglbar^{\bullet}_{\poly}
 (V_{\ast},\DDlambda)_{\Xtilde(O)}$
are $c$-soft,
we obtain Proposition \ref{prop;07.7.22.35}
as the push-forward of Lemma \ref{lem;08.9.9.6}.
\hfill\qed

\section{Local quasi isomorphisms in family}
\label{subsection;07.10.20.6}
Let $\lambda_0\in\cnum_{\lambda}$.
Let $\nbigk$ be a neighbourhood of $\lambda_0$
in $\cnum_{\lambda}$,
which will be shrinked if it is necessary.
We put 
$\nbigx:=\nbigk\times X$
and $\nbigd:=\nbigk\times D$.
Let $p_{\lambda}$ denote the projection
forgetting the $\nbigk$-component.
For any subset $\nbigu$ of $\nbigx$,
we use the symbol
$\nbigu^{\ast}$ to denote
$\nbigu\setminus \nbigd$.
Let $(V_{\ast},\DD)$ be
a good family of filtered $\lambda$-flat bundles
on $(\nbigx,\nbigd)$
with KMS-structure at $\lambda_0$
indexed by $T$.
We shall consider the family versions
of the complexes 
in Section \ref{subsection;07.9.11.11}.

We have the induced endomorphism
$\Gr^{\Fzero}_a\Res(\DD)$ on 
$\Gr^{\Fzero}_a(V)$.
We assume that the conjugacy classes
of the nilpotent part of 
$\Gr^{\Fzero}_a\Res(\DD)$ are 
independent of $\lambda\in \nbigk$.

\subsection{Preliminary for metrics}
\label{subsection;10.5.23.1}

We have the filtration $\Fzero$
and the decomposition $\EEzero$
of $\prolong{V}$ as in Section
\ref{subsection;07.10.14.15}.
Let $\vecv$ be a frame of $\prolong{V}$ compatible with
$\Fzero$, $\EEzero$
and the weight filtration $W$ on 
$\Gr^{\Fzero,\EEzero}(\prolong{V})$.
We have $u(v_i)\in\real\times\cnum$
such that $\kmsmap(\lambda_0,u(v_i))=
 \deg^{\Fzero,\EEzero}(v_i)$.
We put $k(v_i):=\deg^{W}(v_i)$.
Let $h$ be the hermitian metric given as follows:
\[
 h(v_i,v_j):=\delta_{i,j}\, 
 |z|^{-2\paramap(\lambda,u(v_i))}\, (-\log|z|)^{k(v_i)}
\]
Note $-1<\paramap(\lambda_0,u(v_i))\leq 0$
for each $i$.
Recall $\paramap(\lambda,u(v_i))=
 a(v_i)+\Re(\lambdabar\alpha(v_i))$,
where $u(v_i)=(a(v_i),\alpha(v_i))\in\real\times\cnum$.
Hence, 
if $\paramap(\lambda_0,u(v_i))<0$,
we may and will assume that
$\paramap(\lambda,u(v_i))<0$ on $\nbigk$.
If $\paramap(\lambda_0,u(v_i))=0$
and $u(v_i)\neq (0,0)$,
we have $\paramap(\lambda,u(v_i))>0$
on an open subset whose closure contains $\lambda_0$.
If $u(v_i)=(0,0)$,
$\paramap(\lambda,u(v_i))$ is constantly $0$.

Note if $h'$ comes from another choice of $\vecv'$,
the metrics $h$ and $h'$ are mutually bounded.
In the following,
we will use the metric
$\gtilde_{\poin}:=
 g_{\poin}+d\lambda\, d\lambdabar$
and the induced volume form
$\dvol_{\gtilde}$ for the base space
$\nbigx\setminus\nbigd$.

\subsubsection{Condition for
a holomorphic section to be $L^2$}
\label{subsection;10.5.22.2}

Let $f=\sum f_i\,v_i$ be a holomorphic section
on an open subset $U$ of $\nbigk\times X^{\ast}$.
Let us describe the condition for $f$
to be $L^2$ on a neighbourhood of 
$(\lambda,O)\in \nbigd$.
By the orthogonality,
we have only to consider each $f_i\,v_i$.
\begin{description}
\item[$\paramap(\lambda,u(v_i))<0$]
$f_i$ is holomorphic at $(\lambda,O)$.
\item[$u(v_i)=(0,0)$, $k(v_i)\leq 0$]
$f_i$ is holomorphic at $(\lambda,O)$.
\item[$u(v_i)=(0,0)$, $k(v_i)>0$]
$f_i$ is holomorphic at $(\lambda,O)$
and $f_i(\lambda,O)=0$.
\item[$\paramap(\lambda,u(v_i))=0$,
$u(v_i)\neq (0,0)$]
$f_i$ is holomorphic at $(\lambda,O)$,
and $f_i(\lambda,0)=0$.
\item[$\paramap(\lambda,u(v_i))>0$]
$f_i$ is holomorphic at $(\lambda,O)$,
and $f_i(\lambda,0)=0$.
\end{description}
We have similar conditions
for a holomorphic section $f_i\,v_i\,dz/z$
of $V\otimes\Omega^{1,0}$
to be $L^2$ with respect to 
$(h,g_{\poin})$.
\begin{description}
\item[$\paramap(\lambda,u(v_i))<0$]
$f_i$ is holomorphic at $(\lambda,O)$.
\item[$u(v_i)=(0,0)$, $k(v_i)\leq -2$]
$f_i$ is holomorphic at $(\lambda,O)$.
\item[$u(v_i)=(0,0)$, $k(v_i)>-2$]
$f_i$ is holomorphic at $(\lambda,O)$
and $f_i(\lambda,O)=0$.
\item[$\paramap(\lambda,u(v_i))=0$,
$u(v_i)\neq (0,0)$]
$f_i$ is holomorphic at $(\lambda,O)$,
and $f_i(\lambda,0)=0$.
\item[$\paramap(\lambda,u(v_i))>0$]
$f_i$ is holomorphic at $(\lambda,O)$,
and $f_i(\lambda,0)=0$.
\end{description}

\subsubsection{Remark}
\label{subsection;10.5.23.11}

Let $U$ be a subset of $\nbigk\times X^{\ast}$.
Let $f$ be a $C^{\infty}$-section of 
$V\otimes p_{\lambda}^{\ast}\Omega^p$
on $U$.
We say that $f$ is $\lambda$-holomorphic,
if $\delbar_{\lambda}f=0$.
Assume $U=\nbigk\times X^{\ast}(R)$,
for simplicity.
Assume that $f$ is $\lambda$-holomorphic.
We have the expression
$f=\sum f_i\, v_i$.
We have
$-\delbar_{\lambda}\del_{\lambda}
 \log|f_i\, v_i|^2_{h,\gtilde_{\poin}}=0$.
Hence, $|f|_{h,\gtilde_{\poin}}^2$
is subharmonic
with respect to $\lambda$.
It implies the following,
for any 
 $\lambda$-holomorphic $C^{\infty}$-section
 $f$ which is $L^2$ with respect to 
$h$ and $\gtilde_{\poin}$:
\begin{enumerate}
\item
 The restrictions
 $f_{|\{\lambda\}\times X^{\ast}(R)}$
 are also $L^2$.
 Moreover, 
 when we fix a compact subset 
$\nbigk_1$ of the interior of $\nbigk$,
 there exists a constant $C_1>0$
 such that 
 $\|f_{|\{\lambda\}\times X^{\ast}(R)}
 \|_{h,g_{\poin}}
 \leq
 C_1\, \|f\|_{h,\gtilde_{\poin}}$
for any $\lambda\in \nbigk_1$.
 The constant $C_1$ is independent of $R$.
\item
 For any $\epsilon>0$,
 we take a small $R'$ such that
 $\bigl\|
 f_{|\nbigk\times X^{\ast}(R')}
 \bigr\|_{h,\gtilde_{\poin}}<\epsilon$.
 Then, for any $\lambda\in \nbigk_1$,
 we obtain
 $\|f_{|\{\lambda\}\times X^{\ast}(R')}\|
 _{h,g_{\poin}}
 \leq C_1\, \epsilon$.
 Hence, for the expression
 $f=\sum f_i\, v_i$,
 we can show the continuity of
 the function from $\nbigk$
 to the space of $L^2$-forms on $X^{\ast}(R)$
 given by 
 $\lambda\longmapsto
 f_{i|\{\lambda\}\times X^{\ast}(R)}
 \, |v_i|_{h}$,
 for example.
 (The claim $2$ follows from $1$.)
\end{enumerate}
Note if a hermitian metric $h'$ is mutually bounded
with $h$,
then the claims $1$ and $2$ above
also hold for $h'$,
although $|f|^2_{h',\gtilde_{\poin}}$
is not necessarily
subharmonic with respect to $\lambda$.

\subsection{Sheaves of $L^2$-sections}

\subsubsection{$\lambda$-holomorphic
$L^2$-sections and holomorphic
$L^2$-sections}

The family of $\lambda$-connections $\DD$ of $V$
and the differential
$\lambda\del_X+\delbar_X$ of $\Omega^{\bullet}_X$
induce the differential of
the $\lambda$-holomorphic $C^{\infty}$-sections of 
$V\otimes p_{\lambda}^{\ast}\Omega^{\bullet}_X$
on open subsets of $\nbigx^{\ast}$,
which is also denoted by $\DD$.
We shall introduce a complex of sheaves
$\nbigl^{\bullet}_{\poly}(V_{\ast},\DD)(\nbigu)$
on $X$,
which is an extension of
the sheaf of $\lambda$-holomorphic
$C^{\infty}$-sections of
$V\otimes\Omega^{\bullet}$.

For any open subset $\nbigu$ of $\nbigk$,
let $\nbigl^p_{\poly}(V_{\ast},\DD)(\nbigu)$ be 
the space of $\lambda$-holomorphic
$C^{\infty}$-sections $\tau$ of $V\otimes\Omega^p$ on 
open subsets $\nbigu^{\ast}$
with the following property:
\begin{itemize}
\item $\tau$ and $\DD\tau$
are $L^2$ and of polynomial order
locally on $\nbigu$
with respect to $h$ and $\gtilde_{\poin}$.
\end{itemize}
Thus, we obtain a complex of sheaves
$\nbigl^{\bullet}_{\poly}(V_{\ast},\DD)$
on $\nbigx^{(\lambda_0)}$.

Let $\nbigl^{p}_{\hol}(V_{\ast},\DD)$ $(p=1,2)$
be the subsheaf of
$\nbigl^p_{\poly}(V_{\ast},\DD)$,
which consists of holomorphic $(p,0)$-forms.
We will show the following proposition in Section 
\ref{subsection;07.11.16.37}--\ref{subsection;07.11.16.36}.
\begin{prop}
\label{prop;07.11.16.35}
The natural inclusion
$\varphi_0:
 \nbigl^p_{\hol}(V_{\ast},\DD)
\lrarr
 \nbigl^p_{\poly}(V_{\ast},\DD)$
is a quasi isomorphism.
\end{prop}

\subsubsection{Algebraically determined sheaf}

We shall give a rather algebraic description
of the complexes
up to quasi isomorphisms.
Let $X'$ and $\varphi_n$ be 
as in Section \ref{subsection;07.9.11.11}.
We put $D':=\{O'\}$
and $\nbigx':=\nbigk\times X'$
and $\nbigd':=\nbigk\times D'$.
Recall that we have the induced 
good family of filtered $\lambda$-flat bundles
on $(\nbigx', \nbigd')$
as in Section \ref{subsection;07.11.5.60},
which is denoted by $(V'_{\ast},\DD')$.
It is easy to see that 
$(V'_{\ast},\DD')$ also has the KMS-structure
at $\lambda_0$
if we shrink $\nbigk$ appropriately.
If we choose $n$ appropriately,
we have the irregular decomposition:
\[
 (V'_{\ast},\DD')_{|\nbigdhat^{\prime}}=
 \bigoplus_{\gminia\in \Irr(\DD^{\prime})}
 \bigl(\Vhat'_{\gminia\ast},\DD'_{\gminia}\bigr)
\]
Since $\Vhat'_{0\ast}$ and
$\bigoplus_{\gminia\neq 0}
 \Vhat'_{\gminia\ast}$ are $\Gal(X'/X)$-equivariant,
we have the descent to $\nbigdhat$
which are denoted  by 
$\Vhat_{\reg}$ and $\Vhat_{\irr}$,
respectively.

Let $a\in\real$ and $\gminia\in\Irr(\DD^{\prime})$.
We use the notation 
in Section \ref{subsection;07.10.14.15}.
We have the generalized eigen decomposition
\[
 \Gr^{\Fzero}_a(\Vhat'_{\gminia})
=\bigoplus_{u\in \nbigk(a)}
 \EE_{\eigenmap(\lambda,u)}
 \Gr^{\Fzero}_a(\Vhat'_{\gminia}),
\]
where the restriction of 
$\Res(\DD)-\eigenmap(\lambda,u)$
to $\EE_{\eigenmap(\lambda,u)}
 \Gr^{\Fzero}_a(\Vhat'_{\gminia})$
are nilpotent.
We have the weight filtration $W$
of the nilpotent part of 
$\Gr^{\Fzero}_a\Res(\DD')$
on $\EE_{\eigenmap(\lambda,u)}
 \Gr^{\Fzero}_a
 \bigl(\Vhat'_{\gminia}\bigr)$
and $\Gr^{\Fzero}_a\bigl(\Vhat'_{\gminia}\bigr)$.

For any $\gminia\neq 0$,
let $\nbigs\bigl(
 \Vhat'_{\gminia\ast}\otimes\Omega^{0,0}_{X'}\bigr)$
denote the pull back of
$W_{-2}\EE_{-\ord(\gminia)\lambda}
 \Gr^{\Fzero}_{\ord(\gminia)}
 \bigl(\Vhat'_{\gminia}\bigr)$
via the projection
$\prolongg{\ord(\gminia)}{\Vhat'_{\gminia}}
\lrarr
 \Gr^{\Fzero}_{\ord(\gminia)}
 \bigl(\Vhat'_{\gminia}\bigr)$.
Let $\nbigs\bigl(
 \Vhat'_{\gminia\ast}\otimes
 \Omega^{1,0}_{X'}\bigr)$
denote the pull back of 
$W_{-2}\EE_0\Gr^{\Fzero}_{0}
  \bigl(\Vhat'_{\gminia}\bigr)$
via the projection
$\prolong{\Vhat'_{\gminia}} dz'/z'
\lrarr
 \Gr^{\Fzero}_0
 \bigl(\Vhat'_{\gminia}\bigr)$.
The descent of
$\bigoplus_{\gminia\neq 0}
\nbigs(\Vhat'_{\gminia\ast}
 \otimes\Omega_{X'}^{p,0})$
is denoted by
$\nbigs(\Vhat_{\irr\ast}\otimes\Omega_X^{p,0})$.

Let $\nbigs(\Vhat_{0\ast}\otimes\Omega_{X}^{0,0})$
denote the inverse image of
$W_0\EE_0\Gr^{\Fzero}_0\bigl(\Vhat_{0}\bigr)$
via the projection
$\prolong{\Vhat_0}\lrarr
 \Gr^{\Fzero}_0\bigl(\Vhat_0\bigr)$.
Let $\nbigs(\Vhat_{0\ast}\otimes\Omega_X^{1,0})$
denote the inverse image of
$W_{-2}\EE_0\Gr^{\Fzero}_0\Vhat_{0}$
via the projection
$\prolong{\Vhat_0} dz/z\lrarr
 \Gr^{\Fzero}_0\Vhat_{0}$.

Thus, we obtain lattices
$\nbigs(\Vhat_{\ast}\otimes\Omega_X^{p,0})
=\nbigs(\Vhat_{0\ast}\otimes\Omega_X^{p,0})
 \oplus
 \nbigs(\Vhat_{\irr\ast}\otimes\Omega_X^{p,0})$
of $\bigl(V\otimes 
 p_{\lambda}^{\ast}\Omega^{p,0}_X(\ast D)
 \bigr)_{|\nbigdhat}$.
They induce lattices
$\nbigs(V_{\ast}\otimes\Omega_X^{p,0})$
of $V\otimes p_{\lambda}^{\ast}\Omega^{p,0}_X(\ast D)$.
The family of $\lambda$-connections 
$\DD$ on $V_{\ast}$
and the differential $\lambda\, d_X$ on 
$\Omega^{\bullet,0}_X$
induce 
$\DD:\nbigs(V_{\ast}\otimes\Omega_X^{0,0})
\lrarr \nbigs(V_{\ast}\otimes\Omega_X^{1,0})$.
Thus, we obtain a complex of sheaves
$\nbigs(V_{\ast}\otimes\Omega^{0,0}_X)
  \stackrel{\DD}{\lrarr}
 \nbigs(V_{\ast}\otimes\Omega^{1,0}_X)$.
We will show the following lemma
in Section \ref{subsection;10.5.22.20}.
\label{lem;10.5.22.15}
\begin{lem}
We have a natural inclusion
$\nbigs(V_{\ast}\otimes\Omega^{\bullet,0}_X)
\lrarr
 \nbigl^{\bullet}_{\hol}(V_{\ast},\DD)$.
It is a quasi isomorphism.
If $\lambda_0$ is generic, moreover,
it is an isomorphism.
\end{lem}

\subsection{Variants in the case $\lambda_0\neq 0$}
\label{subsection;07.11.15.30}

We shall introduce complexes of sheaves
$\nbigltilde^p_{\poly}(V_{\ast},\DD)$
and 
$\nbiglbar^p_{\poly}(V_{\ast},\DD)$
on $\nbigx$,
which are extensions of
the sheaf of $\lambda$-holomorphic
$C^{\infty}$-sections of $V\otimes\Omega^{\bullet}$
on $\nbigx^{\ast}$.
Although they can be given
as in Section \ref{subsection;07.11.15.1},
we define them in a slightly different 
but equivalent way.
We assume $\lambda_0\neq 0$.

\subsubsection{Decomposition}

We have the full reduction
$\Gr^{\nbigftilde}_{\gminia}(\prolong{V'})$ 
for $\gminia\in\Irr(\DD)$.
We take a frame $\vecvbar_{\gminia}$
of $\Gr^{\nbigftilde}_{\gminia}(\prolong{V'})$
compatible with
the induced decomposition $\EEzero$,
the induced filtration $\Fzero$,
and the weight filtration $W$
on $\Gr^{\EEzero,\Fzero}
 \bigl(\prolong{V'}\bigr)$.
For each $g\in\Gal(X'/X)$,
we have a naturally induced isomorphism
$g:\Gr^{\nbigftilde}_{\gminia}(\prolong{V'})
 \lrarr \Gr^{\nbigftilde}_{g\,\gminia}(\prolong{V'})$.
We assume
$g^{\ast}\vecvbar_{\gminia}
=\vecvbar_{g^{\ast}\gminia}$
for any $g\in\Gal(X'/X)$.

Let $S$ be a small sector in 
$\nbigx^{\prime}\setminus\nbigd^{\prime}$,
and let $\Sbar$ denote
its closure in the real blow up 
$\nbigxtilde'(\nbigd')$.
When $S$ is small,
we have the full Stokes filtration
$\nbigftilde^{S}$ of $V'_{|\Sbar}$,
and we can take a $\DD$-flat splitting
$V'_{|\Sbar}
=\bigoplus_{\gminia}V'_{\gminia,S}$.
We have the holomorphic lift $\vecv_{\gminia,S}$
of $\vecvbar_{\gminia}$ to $V'_{\gminia,S}$.

By shrinking $\nbigk$,
we take a covering
$\nbigx^{\prime}\setminus\nbigd^{\prime}
 =\bigcup_{j=1}^MS_j$
by small sectors
with the following property:
\begin{itemize}
\item
 Each $S_j$ is the product of
 a small sector of $X^{\ast}(R)$
 and $\nbigk$.
\item
We have the full Stokes filtrations $\nbigftilde^{S_j}$,
flat splittings 
$V'_{|\Sbar_j}=
 \bigoplus V_{\gminia,S_j}$,
and lifts $\vecv_{\gminia,S_j}$
of $\vecvbar_{\gminia}$ to $V_{\gminia,S_j}$.
\item
We may assume that $g^{-1}(S_j)$
is the same as some $S_{i(j,g)}$
for any $g\in \Gal(X'/X)$ and $j$,
and
$g^{\ast}\, \vecv_{\gminia,S_j}
=\vecv_{g^{\ast}\gminia,S_{i(g,j)}}$.
\end{itemize}
By gluing them as in Section \ref{subsection;07.6.16.8},
we obtain a $C^{\infty}$-frame
$\vecv_{C^{\infty}}=
 \bigl(\vecv_{\gminia,C^{\infty}}\bigr)$
of $\prolong{V}'$ on $\nbigx$.
We may assume the following:
\begin{itemize}
\item
$g^{\ast}\vecv_{\gminia,C^{\infty}}
=\vecv_{g^{\ast}\gminia,C^{\infty}}$
for $g\in\Gal(\Xtilde/X)$.
In particular, $\vecv_{0,C^{\infty}}$
is $\Gal(\Xtilde/X)$-equivariant,
which induces a tuple of sections of $V$
on $\nbigx\setminus\nbigd$.
It is also denoted by $\vecv_{0,C^{\infty}}$.
\item
$\vecv_{\gminia,C^{\infty}}$ are $\lambda$-holomorphic.
\end{itemize}
Let $V'_{\irr,C^{\infty}}$ be 
the $C^{\infty}$-subbundle
of $V'_{|
 \nbigx^{\prime}
\setminus\nbigd^{\prime}}$
generated by 
$\vecv_{\gminia,C^{\infty}}$ $(\gminia\neq 0)$,
and let $V'_{\reg,C^{\infty}}$ be 
the $C^{\infty}$-subbundle
of $V'_{|\nbigx^{\prime}
 \setminus\nbigd^{\prime}}$
generated by $\vecv_{0,C^{\infty}}$.
Since they are $\Gal(X'/X)$-equivariant,
they induce a decomposition
$V_{|\nbigx\setminus\nbigd}
=V_{\reg,C^{\infty}}\oplus V_{\irr,C^{\infty}}$.
The decomposition is $\lambda$-holomorphic.
Note that it is not canonical.

\begin{lem}
\label{lem;07.11.16.10}
Let $S$ be any small sector of 
$\nbigx^{\prime}
 \setminus\nbigd^{\prime}$
such that we have the full Stokes filtration
$\nbigftilde^{S}$ of $V'_{|\Sbar}$.
Let $Z:=\Sbar\cap \pi^{-1}(\nbigd)$,
where $\pi$ denotes the projection
of $\nbigxtilde^{\prime}
 (\nbigd^{\prime})$
to $\nbigx^{\prime}$.
Let $V'=\bigoplus_{\gminia} V_{\gminia,S}$
be any splitting of $\nbigftilde^S$.
Let $\vecv_{\gminia,S}$ 
be any holomorphic lift of $\vecvbar_{\gminia}$
to $V_{\gminia,S}$,
and let $\vecv_S=(\vecv_{\gminia,S})$.

Let $C$ be determined by
$\vecv_S=\vecv_{C^{\infty}}\, (I+C)$.
Corresponding to the decomposition of the frames,
we have the decomposition
$C=\bigl(C_{\gminia,\gminib}\bigr)$.
Then, the following holds:
\begin{itemize}
\item
 $C_{\gminia,\gminib|\Zhat}=0$.
\item
 $C_{\gminia,\gminib}=0$
 unless $\gminia\leq_S\gminib$.
\item
 $C_{\gminia,\gminib}\,
 \exp\bigl(\lambda^{-1}(\gminia-\gminib)\bigr)$
 is $O\bigl(|z|^{-N}\bigr)$ for some $N$
 in the case $\gminia<_S\gminib$.
\end{itemize}
\end{lem}
\pf
It can be shown using the same argument
as that in the proof of Lemma \ref{lem;07.6.16.3}.
\hfill\qed

\vspace{.1in}

For any local $C^{\infty}$-section $\tau$
of $V\otimes\Omega^p$,
we have the corresponding decompositions
$\tau=\tau_{\reg,C^{\infty}}+\tau_{\irr,C^{\infty}}$
and $\DD\tau=(\DD\tau)_{\reg,C^{\infty}}
+(\DD\tau)_{\irr,C^{\infty}}$.

\subsubsection{Complexes}

For any open $\nbigu\subset\nbigx$,
let $\nbigltilde^p_{\poly}(V_{\ast},\DD)(\nbigu)$ denote 
the space of $\lambda$-holomorphic
$C^{\infty}$-sections $\tau$ of 
$V\otimes\Omega^p$ on $\nbigu^{\ast}$
such that the following estimates hold
for $\gtilde_{\poin}$ and $h$:
\begin{description}
\item[(a1)]
 $\tau_{\irr,C^{\infty}}$
 and $(\DD\tau)_{\irr,C^{\infty}}$ are
 of polynomial order 
 locally on $\nbigu$.
\item[(a2)]
 $\tau_{\reg,C^{\infty}}$ 
 and $(\DD\tau)_{\reg,C^{\infty}}$
 are $L^2$ and of polynomial order
 locally on $\nbigu$.
\end{description}
Let $\nbiglbar^p_{\poly}(V_{\ast},\DD)(\nbigu)$ denote
the space of 
$\lambda$-holomorphic 
$C^{\infty}$-sections $\tau$ of 
$V\otimes\Omega^p$ on $\nbigu^{\ast}$
such that the following estimates hold
for $h$ and $\gtilde_{\poin}$:
\begin{description}
\item[(b1)]
 $\tau_{\irr,C^{\infty}}$
 and 
 $(\DD\tau)_{\irr,C^{\infty}}$ are
 $O\bigl(|z|^N\bigr)$ for any $N$
 locally on $\nbigu$.
\item[(b2)]
 $\tau_{\reg,C^{\infty}}$ and
 $(\DD\tau)_{\reg,C^{\infty}}$ 
 are $L^2$ and of polynomial order
 locally on $\nbigu$.
\end{description}
Thus, we obtain the complexes of sheaves
$\nbigltilde^{\bullet}_{\poly}(V_{\ast},\DD)$
and $\nbiglbar^{\bullet}_{\poly}(V_{\ast},\DD)$.
The following lemma is clear from
Lemma \ref{lem;07.11.16.10}.
\begin{lem}
On the germ of neighbourhoods
of $(\lambda_0,O)$ in $\nbigx$,
the complexes of sheaves
$\nbigltilde^{\bullet}_{\poly}(V_{\ast},\DD)$
and $\nbiglbar^{\bullet}_{\poly}(V_{\ast},\DD)$
are well defined,
i.e., they are independent of the choices of
covering $\nbigx\setminus\nbigk=\bigcup S_i$,
flat splittings on $\Sbar_i$ of $\nbigf^{S_i}$
and $C^{\infty}$-gluings of them.
\hfill\qed
\end{lem}

We will show the following proposition
in Sections \ref{subsection;07.11.16.60}
and \ref{subsection;07.11.16.61}.

\begin{prop}
\label{prop;10.5.23.25}
The natural morphisms
\[
 \nbiglbar_{\poly}^{\bullet}(V_{\ast},\DD)
\stackrel{\varphi_2}{\lrarr}
 \nbigl_{\poly}^{\bullet}(V_{\ast},\DD)
\stackrel{\varphi_1}{\lrarr}
 \nbigltilde_{\poly}^{\bullet}(V_{\ast},\DD)
\]
are quasi isomorphisms.
\end{prop}

\subsection{Deformation of the Stokes structure}
\label{subsection;07.12.2.25}

Assume $\lambda_0\neq 0$.
We consider family versions of
the complexes in Section 
\ref{subsection;07.10.20.2}.

\subsubsection{Construction of a metric}
Let us take a finite covering
$\nbigx\setminus\nbigd=\bigcup_{j=1}^M S_j$
as in Section \ref{subsection;07.11.15.30}.
Let $T=T(\lambda):=1+|\lambda|^2$,
and $F^{(T)}_{S_j}:=
 \bigoplus_{\gminia}
 \exp\bigl((1-T)\lambda^{-1}\gminia\bigr)
\,\id_{V_{\gminia,S_j}}$.
We consider the metric
$h^{(T)}_{S_j}(u,v):=
 h\bigl(F^{(T)}_{S_j}(u),F^{(T)}_{S_j}(v)\bigr)$.
By gluing $h^{(T)}_{S_j}$ in $C^{\infty}$,
we obtain a $C^{\infty}$-metric $h^{(T)}$
of $V'_{|\nbigx^{\prime}\setminus\nbigd^{\prime}}$.
We may assume that it is $\Gal(X'/X)$-equivariant.
The induced metric of $V_{|\nbigx\setminus\nbigd}$
is also denoted by $h^{(T)}$.
Note that the restrictions of the metrics
$h^{(T)}$
and $h$ to $V_{\reg,C^{\infty}}$
are mutually bounded.

\begin{rem}
$h^{(T)}$ is not obtained as a metric
for a good family of filtered
$\lambda$-flat bundles as in 
Section {\rm\ref{subsection;10.5.23.1}}.
The specialization of $h^{(T)}$
to $\{\lambda\}\times X$
is obtained as a metric for
the  good filtered $\lambda$-flat bundle
$(V^{\lambda}_{\ast},\DDlambda)^{(T(\lambda))}$.
(See Section {\rm\ref{subsection;10.5.23.10}}
for $V^{\lambda}_{\ast}$.)
\hfill\qed
\end{rem}

\subsubsection{Complexes}

For any open $\nbigu\subset\nbigx$,
let $\nbiglbar^p_{\poly}(V,h^{(T)})(\nbigu)$ be 
the space of 
$\lambda$-holomorphic $C^{\infty}$-sections $\tau$ of
$V\otimes\Omega^{p}$ on $\nbigu^{\ast}$
with the following growth estimate
with respect to $h^{(T)}$
 and $\gtilde_{\poin}$ locally on $\nbigu$:
\begin{description}
\item[(c1)]
 $\tau_{\irr,C^{\infty}}$ and 
 $(\DD\tau)_{\irr,C^{\infty}}$
 are $O(|z|^N)$ for any $N>0$.
\item[(c2)]
 $\tau_{\reg,C^{\infty}}$ and 
 $(\DD\tau)_{\reg,C^{\infty}}$
 are $L^2$ and of polynomial order.
 (Equivalently, they are 
 $L^2$ and of polynomial order
 with respect to $h$ and $\gtilde_{\poin}$.)
\end{description}
Let $\nbiglbar^p_{\poly}(V,h,h^{(T)})(\nbigu)$ be 
the space of $\lambda$-holomorphic
$C^{\infty}$-sections $\tau$ of
$V\otimes\Omega^{p}$
with the following property:
\begin{description}
\item[(d1)]
 $\tau_{\irr,C^{\infty}}$ and 
 $(\DD\tau)_{\irr,C^{\infty}}$
 are $O(|z|^N)$ for any $N>0$
 locally on $\nbigu$,
 with respect to
 both $(h,\gtilde_{\poin})$ and 
 $(h^{(T)},\gtilde_{\poin})$. 
\item[(d2)]
 $\tau_{\reg,C^{\infty}}$ and 
$(\DD\tau)_{\reg,C^{\infty}}$
 satisfy the conditions (c2).
\end{description}
Thus, we obtain complexes of sheaves
$\nbiglbar^{\bullet}_{\poly}\bigl(V,h^{(T)}\bigr)$ and
$\nbiglbar^{\bullet}_{\poly}\bigl(V,h,h^{(T)}\bigr)$.
The following lemma is clear 
from Lemma \ref{lem;07.11.16.10}.
\begin{lem}
$\nbiglbar^{\bullet}_{\poly}\bigl(V,h^{(T)}\bigr)$ and
$\nbiglbar^{\bullet}_{\poly}\bigl(V,h,h^{(T)}\bigr)$
are well defined for $(V_{\ast},\DD)$
on the germ of neighbourhoods 
of $(\lambda_0,O)$ in $\nbigx$.
\hfill\qed
\end{lem}

\subsubsection{Statement}
We will prove the following lemma
in Section \ref{subsection;07.11.16.61}.
\begin{prop}
\label{prop;07.11.16.63}
The naturally defined morphisms
\[
\begin{CD}
\nbiglbar^{\bullet}_{\poly}(V_{\ast},\DD)
 @<<<
\nbiglbar^{\bullet}_{\poly}(V,h,h^{(T)})
 @>>>
 \nbiglbar^{\bullet}_{\poly}(V,h^{(T)})
\end{CD}
\]
are quasi isomorphisms.
\end{prop}

\subsection{Complement for varying $\lambda$}

Assume $\lambda_0\neq 0$.
Let $\lambda_1\in \nbigk$.
We take a small neighbourhood $\nbigk_1$
of $\lambda_1$  in $\nbigk$,
and we put 
$\nbigx^{(\lambda_1)}:=
 \nbigk_1\times X$.
We have the good family of filtered 
$\lambda$-flat bundles
$V_{\ast}^{(\lambda_1)}$ obtained from $V_{\ast}$,
as explained in Section \ref{subsection;07.10.14.15}.
If we construct $h_1$
for $V_{\ast}^{(\lambda_1)}$ as above,
$h_1$ and $h_{|\nbigx^{(\lambda_1)}}$
are mutually bounded.

\subsubsection{}
\label{subsection;10.5.23.30}
By the previous procedures,
we obtain complexes of sheaves
on $\nbigx^{(\lambda_1)}$:
\[
 \nbigs\bigl(V^{(\lambda_1)}_{\ast}
 \otimes\Omega^{\bullet,0}\bigr),
\quad
 \nbigl^{\bullet}_{\poly}(V^{(\lambda_1)}_{\ast},\DD),
\quad
 \nbigltilde^{\bullet}_{\poly}(V^{(\lambda_1)}_{\ast},\DD),
\quad
 \nbiglbar^{\bullet}_{\poly}(V^{(\lambda_1)}_{\ast},\DD)
\]
By construction,
$\nbigl^{\bullet}_{\poly}(V^{(\lambda_1)}_{\ast},\DD)$,
$\nbigltilde^{\bullet}_{\poly}(V^{(\lambda_1)}_{\ast},\DD)$
and
$\nbiglbar^{\bullet}_{\poly}(V^{(\lambda_1)}_{\ast},\DD)$
are the same as the restrictions of
$\nbigl^{\bullet}_{\poly}(V_{\ast},\DD)$,
$\nbigltilde^{\bullet}_{\poly}(V_{\ast},\DD)$
and
$\nbiglbar^{\bullet}_{\poly}(V_{\ast},\DD)$
to $\nbigx^{(\lambda_1)}$,
respectively.
The sheaves 
for $(V_{\ast}^{(\lambda_1)},\DD)$
as in Section \ref{subsection;07.12.2.25}
are also the same as the restriction
of the sheaves for $(V_{\ast},\DD)$.

We have the following commutative diagram:
\[
\begin{CD}
 \nbigs\bigl(V_{\ast}\otimes\Omega^{\bullet,0}\bigr)_{|
 \nbigx^{(\lambda_1)}}
 @>>>
 \nbigl^{\bullet}_{\poly}(V_{\ast},\DD)
 _{|\nbigx^{(\lambda_1)}}\\
@V{\rho_1}VV @V{=}VV \\
 \nbigs\bigl(V^{(\lambda_1)}_{\ast}
 \otimes\Omega^{\bullet,0}\bigr)
 @>>>
 \nbigl^{\bullet}_{\poly}(V^{(\lambda_1)}_{\ast},\DD)
\end{CD}
\]
It is easy to see that $\rho_1$ is a quasi isomorphism.
(See the proof of Lemma \ref{lem;10.5.22.15}.)

\subsubsection{}
\label{subsection;10.5.23.10}
Let $\nbigx^{\lambda_1}:=\{\lambda_1\}\times X$.
By considering the specialization at 
$\nbigx^{\lambda_1}$,
i.e.,
taking the cokernel of the multiplication of
$\lambda-\lambda_1$,
we obtain a good filtered $\lambda_1$-flat bundle
$(V^{\lambda_1}_{\ast},\DD^{\lambda_1}):=
(V^{(\lambda_1)}_{\ast},\DD)_{|\nbigx^{\lambda_1}}$.
We obtain the complex of sheaves,
as in Section \ref{subsection;07.10.20.1}:
\[
 \nbigs\bigl(
 V^{\lambda_1}_{\ast}\otimes\Omega^{\bullet,0}\bigr),
\quad
 \nbigl^{\bullet}_{\poly}(
 V^{\lambda_1}_{\ast},\DD^{\lambda_1}),
\quad
  \nbiglbar^{\bullet}_{\poly}(V^{\lambda_1}_{\ast},
 \DD^{\lambda_1})
\]
The restriction of $h$ to $\nbigx^{\lambda_1}$
is denoted by $h_{\lambda_1}$.
We also obtain the following complexes of sheaves:
\[
  \nbiglbar^{\bullet}_{\poly}
 (V^{\lambda_1},h_{\lambda_1},
 h_{\lambda_1\,C^{\infty}}^{(T(\lambda_1))}),
\quad
  \nbiglbar^{\bullet}_{\poly}
 (V^{\lambda_1},
 h_{\lambda_1,C^{\infty}}^{(T(\lambda_1))})
\]
By taking the specialization at $\lambda_1$,
we obtain the following commutative diagrams,
which will be used in 
Section \ref{subsection;07.10.23.30}.
\begin{equation}
 \label{eq;08.1.31.2}
\begin{CD}
 \nbigs\bigl(V_{\ast}\otimes\Omega^{\bullet,0}\bigr)
  _{|\nbigx^{\lambda_1}}
 @>>>
 \nbigl^{\bullet}_{\poly}(V_{\ast},\DD)
 _{|\nbigx^{\lambda_1}}
 @<<<
 \nbiglbar^{\bullet}_{\poly}(V_{\ast},\DD)
 _{|\nbigx^{\lambda_1}}
 \\
@VVV @VVV @ VVV\\
 \nbigs\bigl(V^{\lambda_1}_{\ast}
 \otimes\Omega^{\bullet,0}\bigr)
 @>>>
 \nbigl^{\bullet}_{\poly}(V^{\lambda_1}_{\ast},\DD)
 @<<<
 \nbiglbar^{\bullet}_{\poly}(V^{\lambda_1}_{\ast},\DD)
\end{CD}
\end{equation}
\begin{equation}
 \label{eq;08.1.31.3}
 \begin{CD}
  \nbiglbar^{\bullet}_{\poly}(V_{\ast},\DD)
 _{|\nbigx^{\lambda_1}}
 @<<<
 \nbiglbar^{\bullet}_{\poly}(V_{\ast},h,h^{(T)})
 _{|\nbigx^{\lambda_1}}
 @>>> 
\nbiglbar^{\bullet}_{\poly}(V_{\ast},h^{(T)})
 _{|\nbigx^{\lambda_1}}\\
 @VVV @VVV @VVV \\
 \nbiglbar^{\bullet}_{\poly}(V^{\lambda_1}_{\ast},\DD)
 @<<<
 \nbiglbar^{\bullet}_{\poly}(
 V^{\lambda_1}_{\ast},h_{\lambda_1},
 h_{\lambda_1,C^{\infty}}^{(T(\lambda_1))})
 @>>> 
\nbiglbar^{\bullet}_{\poly}(
 V^{\lambda_1}_{\ast},
 h_{\lambda_1,C^{\infty}}^{(T(\lambda_1))})
 \end{CD}
\end{equation}

\section{Proof in the family case}
\label{subsection;07.11.16.120}
Because we will use the almost same arguments
as those in Section \ref{subsection;08.1.5.2},
we will omit some details.

\subsection{Variant of an estimate
 in \cite{z}}

Let $u=(a,\alpha)$
such that $-1<\paramap(\lambda_0,u)\leq 0$,
and let $k\in\seisuu$.
Let $\nbigl$ be a holomorphic line bundle
on $\nbigk\times X^{\ast}$
with the holomorphic frame $\sigma$
and the metric $h$ such that
$h(\sigma,\sigma)
=|z|^{-2\paramap(\lambda,u)}
 \,\bigl|\log |z|\bigr|^k$.
Let $\omega=g\,\sigma\, d\zbar/\zbar$ be 
a $\lambda$-holomorphic $C^{\infty}$-section
of $\nbigl\otimes p_{\lambda}^{\ast}\Omega^{0,1}_{X}$
on $\nbigk\times X^{\ast}(R)$,
which is $L^2$ and of polynomial order
with respect to 
$h$ and $\gtilde_{\poin}$.
We have the Fourier decomposition
$g=\sum g_n(r,\lambda)\, e^{\sqrt{-1}n\theta}$.
As in Section \ref{subsection;07.9.12.10},
we put
 $g^{(1)}:=\sum_{m\neq 0}
 g_m(r)\, e^{\sqrt{-1}m\theta}$,
and thus we have the decompositions 
$g=g_0+g^{(1)}$
and $\omega=\omega_0+\omega^{(1)}$.
Let $\delbar_z$ and $\delbar_{\lambda}$
denote the natural $(0,1)$-operator
along $X$-direction and $\nbigk$-direction,
respectively.
We have a direct consequence of
the estimate in \cite{z},
which we formulate 
for the reference in the subsequent argument.

\begin{lem}
\label{lem;08.1.27.1}
 We have a $\lambda$-holomorphic
 $C^{\infty}$-section $\tau$,
 which is 
 $L^2$ and of polynomial order
 with respect to $h$ and $\gtilde_{\poin}$,
 with the following property:
\begin{itemize}
\item
We have $\delbar_z\tau=\omega$
in the case that one of 
$(\paramap(\lambda_0,u)<0)$
or $(u=(0,0),k\neq 1)$ holds.
\item
We have $\delbar_z\tau=\omega^{(1)}$
otherwise.
\end{itemize}
\end{lem}
\pf
Note a remark in Section \ref{subsection;10.5.23.11}.
Let $\omega_{\lambda}$ denote the restriction 
of $\omega$ to $\{\lambda\}\times X^{\ast}(R)$.
Let us consider the case in which one of
$(\paramap(\lambda_0,u)<0)$
or $(u=(0,0),k\neq 1)$ holds.
Applying $\Phi$  to each $\omega_{\lambda}$,
we obtain the $L^2$-section $\Phi(\omega)$ of $\nbigl$.
By construction,
it satisfies $\delbar_z\Phi(\omega)=\omega$
and $\delbar_{\lambda}\Phi(\omega)=0$
as a distribution.
Then, it follows that $\Phi(\omega)$ is
$C^{\infty}$ and $\lambda$-holomorphic
due to standard ellipticity.
By using Sobolev embedding,
we obtain $\Phi(\omega)$ is of polynomial order.
(See the last part of the proof of Lemma
\ref{lem;07.11.16.7}.)

If neither $(\paramap(\lambda_0,u)<0)$
or $(u=(0,0),k\neq 1)$ are satisfied,
we obtain the desired section
by applying $\Phi^{(1)}$.
\hfill\qed

\subsection{Preliminary}
\label{subsection;07.11.16.37}
Let us start the proof of the propositions
in Section \ref{subsection;07.10.20.6}.
By an easy argument to use a decent,
we can reduce the problem to the unramified case.
Therefore, we may and will assume that
$(V_{\ast},\DD)$ is unramified.
We use the polar coordinate 
$z=r\, e^{\sqrt {-1}\theta}$.
We may assume the following for the frame $\vecv$:
\begin{enumerate}
\item $\vecv$ is compatible with
the irregular decomposition in $N$-th order
for some large $N$.
(See Section \ref{subsection;08.1.27.10}.)
\item 
 Let $N_{a,\alpha,\gminia}$
 be the nilpotent part of the endomorphisms
 on $\Gr^{\Fzero}_a
 \EEzero_{\alpha}\bigl(V_{\gminia|\nbigd}\bigr)$
 induced by $\Res(\DD)$.
 Then,
 $N_{a,\alpha,\gminia}$ are represented
 by Jordan matrices with respect to the induced frames.
\end{enumerate}
We have the irregular value $\gminia(v_i)$.
We also put $a(v_i):=\deg^{\Fzero}(v_i)$,
$\alpha(v_i):=\deg^{\EEzero}(v_i)$
and $k(v_i):=\deg^W(v_i)$.
Let $u(v_i)\in\real\times\cnum$ be determined by
$\kmsmap(\lambda_0,u(v_i))=
 \bigl(a(v_i),\alpha(v_i)\bigr)$.
We define
\begin{multline*}
 \nbigb(k):=
\bigl\{v_i\,\big|\,
 \gminia(v_i)=0,\,\,u(v_i)=(0,0),\,\, k(v_i)=k
\bigr\} 
\cup \\
\bigl\{v_i\,\big|\,
 a(v_i)=0,\,\,
 \bigl(\gminia(v_i),\alpha(v_i)\bigr)\neq (0,0)
 \bigr\}.
\end{multline*}

Let $A$ be determined by 
$\DD\vecv=\vecv\, A$.
Let $\Gamma$ be the diagonal matrix
whose $(i,i)$-entries are
\[
 \eigenmap\bigl(\lambda,u(v_i)\bigr)\, dz/z
+d\gminia(v_i).
\]
We put $A_0:=A-\Gamma$.
We use the symbol $F_A$
to denote the section of
$\End(V)\otimes\Omega^{1,0}$
determined by 
$F_A(\vecv)=\vecv\, A$.
We use the symbol $F_{A_0}$ in a similar meaning.
Then, $F_{A_0}$ is bounded
with respect to $h$ and $\gtilde_{\poin}$.
We have the following decomposition:
\[
 A_0=\bigoplus_{a,\alpha,\gminia} 
 J_{a,\alpha,\gminia}\, \frac{dz}{z}
+A_0'\, \frac{dz}{z}
\]
Here $A_0'$ is holomorphic
and $F_{A_0'|\nbigd}$ strictly decreases 
the filtration $\Fzero$.
And 
$J_{a,\alpha,\gminia}$ are 
constant Jordan matrices
and represent
$N_{a,\alpha,\gminia}$
with respect to the induced frames of
$\Gr^{\Fzero\,\EEzero}_{a,\alpha}\bigl(
 V_{\gminia|\nbigd}\bigr)$.

The $(1,0)$-operator $\del$ is defined by 
$\del\bigl(\sum f_i\, v_i\bigr)=\sum \del_X f_i\, v_i$.
Then, we have $\DD=\delbar+\lambda\del+F_A$.

\subsection{Proof of Lemma \ref{lem;10.5.22.15}}
\label{subsection;10.5.22.20}

Let $f=\sum f_i\,v_i\,dz/z$ be a section of
$\nbigs(V_{\ast}\otimes\Omega^{1,0})$
on an open subset $U$ of $\nbigx$.
Let us show that it is $L^2$ 
around $(\lambda,O)\in U\cap\nbigd$.
By construction, $f_i$ is holomorphic.
If $\paramap(\lambda,u(v_i))> 0$,
we have $\paramap(\lambda_0,u(v_i))=0$
and $u(v_i)\neq (0,0)$.
Hence, we have $f_i(\lambda,O)=0$.
If $\paramap(\lambda,u(v_i))=0$ and 
$u(v_i)\neq (0,0)$,
we have $\paramap(\lambda_0,u(v_i))=0$,
and $f_i(\lambda,O)=0$ by construction
of $\nbigs(V_{\ast}\otimes\Omega^{1,0})$.
If $u(v_i)=(0,0)$ and $k(v_i)>-2$,
we have $f_i(\lambda,O)=0$ by construction.
Hence, $f$ is $L^2$ by the condition described
in Section \ref{subsection;10.5.22.2},
i.e.,
$\nbigs(V_{\ast}\otimes\Omega^{1,0})
\subset
 \nbigl_{\hol}^1(V_{\ast},\DD)$.
Similarly and more easily,
we can check that a section of
$\nbigs(V_{\ast})$ is $L^2$.
Because
$\DD f\in\nbigs(V_{\ast}\otimes\Omega^{1,0})
 \subset
 \nbigl^1_{\hol}(V_{\ast},\DD)$,
we obtain $f\in \nbigl^0_{\hol}(V_{\ast},\DD)$.
Hence, we obtain a natural inclusion
$\nbigs(V_{\ast}\otimes\Omega^{\bullet,0})
\lrarr
 \nbigl_{\hol}^{\bullet}(V_{\ast},\DDlambda)$.

Assume $\lambda_0$ satisfies the following condition
for the index set $T$ of the KMS-structure:
\begin{description}
\item[(A)]
 If $u\in T$ satisfies $\paramap(\lambda_0,u)=0$,
 then $u=(0,0)$.
\end{description}
If $\nbigk$ is small,
we have $\paramap(\lambda,u(v_i))\leq 0$
for any $\lambda\in\nbigk$,
and $\paramap(\lambda,u(v_i))=0
\Longleftrightarrow
 \paramap(\lambda_0,u(v_i))=0
\Longleftrightarrow  u(v_i)=(0,0)$.
Hence, we obtain that
$\nbigs(V_{\ast}\otimes\Omega^{1,0})
\subset
 \nbigl^1_{\hol}(V_{\ast},\DD)$
is an isomorphism.
Let $f\in\nbigl^0_{\hol}(V_{\ast},\DD)$.
We have a description $f=\sum f_i\,v_i$.
Because it is $L^2$,
each $f_i$ are holomorphic.
Then, we can deduce $f\in\nbigs(V_{\ast})$
from $\DD f\in
 \nbigl^1_{\hol}(V_{\ast},\DD)
=\nbigs(V_{\ast}\otimes\Omega^{1,0})$.
Even if $\lambda_0$ does not satisfy (A),
we obtain that the inclusion of the germs
at $(\lambda_0,O)$ is an isomorphism
by the same argument.

Let us show that the inclusion of the germs
$\nbigs(V_{\ast}\otimes
 \Omega^{\bullet,0})_{(\lambda_1,O)}
\lrarr
 \nbigl^{\bullet}_{\hol}(V_{\ast},\DD)_{(\lambda_1,O)}$
is a quasi isomorphism at each 
$(\lambda_1,O)\in\nbigd\setminus\{(\lambda_0,O)\}$.
We may assume that $(\lambda_1,O)$
satisfies (A).

Let $\nbigk_1\subset \nbigk$
be a small neighbourhood of $\lambda_1$.
We put $\nbigx^{(\lambda_1)}:=
\nbigk_1\times X$.
We have the good family of filtered 
$\lambda$-flat bundles
$V_{\ast}^{(\lambda_1)}$ obtained from $V_{\ast}$,
as explained in Section \ref{subsection;07.10.14.15}.
Then, we obtain complexes of sheaves
$\nbigs\bigl(V^{(\lambda_1)}_{\ast}
 \otimes\Omega^{\bullet,0}\bigr)$
and $\nbigl^{\bullet}_{\hol}(
 V^{(\lambda_1)}_{\ast},\DD)$
on $\nbigx^{(\lambda_1)}$,
and the following morphisms:
\[
  \nbigs\bigl(V_{\ast}\otimes\Omega^{\bullet,0}\bigr)_{|
 \nbigx^{(\lambda_1)}}
\stackrel{a}{\lrarr}
 \nbigs\bigl(V^{(\lambda_1)}_{\ast}
 \otimes\Omega^{\bullet,0}\bigr)
\stackrel{b}{\lrarr}
 \nbigl^{\bullet}_{\hol}(V^{(\lambda_1)}_{\ast},\DD)
=\nbigl^{\bullet}_{\hol}(V_{\ast},\DD)
 _{|\nbigx^{(\lambda_1)}}
\]
By the previous consideration,
we have already known that 
(b) is an isomorphism.
It is easy to show $a$ is a quasi isomorphism
by a direct computation.
Thus, the proof of Lemma \ref{lem;10.5.22.15}
is finished.
\hfill\qed

\subsection{Vanishing of
 $\nbigh^2\bigl(
 \nbigl^{\bullet}_{\poly}(V_{\ast},\DD)\bigr)$}
\label{subsection;07.11.16.50}

Let us show
$\nbigh^2\bigl(
\nbigl^{\bullet}_{\poly}(V_{\ast},\DD)\bigr)=0$.
We have only to show such a claim
for the germ at $(\lambda_0,O)$.
(See a remark
in Section \ref{subsection;10.5.23.30}.
We will omit similar remarks in the following.)
Let $\omega$ be a $\lambda$-holomorphic
$C^{\infty}$-section of 
$V\otimes\Omega^2$,
which is $L^2$ with respect to $h$ and $\gtilde_{\poin}$.
We have the expression:
\[
 \omega=f\, \frac{dz}{z}\, \frac{d\zbar}{\zbar},
\quad
 f=\sum f_i\, v_i
\]
Each $f_i$ has the Fourier expansion
$f_i=\sum_{m\in\seisuu}
 f_{i,m}(r)\, e^{\sqrt{-1}m\theta}$.
We set
\[
  \nbiga^{(0)}(f):=
 \sum_{v_i\in\nbigb(-1)}
 f_{i,0}(r)\, v_i,
\quad
 \nbiga^{(1)}(f):=
 f-\nbiga^{(0)}(f).
\]
We have the decomposition
$f=\nbiga^{(0)}(f)+\nbiga^{(1)}(f)$.
We have the corresponding decomposition
$\omega=\nbiga^{(0)}(\omega)+\nbiga^{(1)}(\omega)$.
We show the following lemma based on an idea of Sabbah
contained in \cite{sabbah2}.
\begin{lem}
\label{lem;07.11.16.15}
\mbox{{}}
\begin{itemize}
\item
 We have a section $\tau^{(1)}$ of 
 $\nbigl^{1}_{\poly}(V_{\ast},\DD)$
 such that $\DD\tau^{(1)}=\nbiga^{(1)}(\omega)$.
\item
 We have a section $\tau^{(0)}$ of
 $\nbigl^{1}_{\poly}(V_{\ast},\DD)$
 such that
 $\nbiga^{(0)}\bigl(\omega-\DD\tau^{(0)}\bigr)=0$.
\item
 In particular,
 we can take a section $\tau$ of
 $\nbigl^1_{\poly}(V_{\ast},\DD)$
 such that $\DD\tau=\omega$.
\end{itemize}
\end{lem}
\pf
The argument is almost the same
as the proof of Lemma \ref{lem;07.11.16.7}.
Briefly,
we have only to replace $\alpha(v_i)$
with $\eigenmap(\lambda,u(v_i))$.
We give only an indication.
The first claim follows from 
Lemma \ref{lem;08.1.27.1}.
Let us show the second claim.
We give preliminary arguments.

\noindent
{\bf (A)}\,\,
If $a(v_i)=0$, $\gminia(v_i)=0$ and $\alpha(v_i)\neq 0$ hold,
let $\tau_1$ be given by (\ref{eq;08.1.27.2}).
We have 
$\DD\tau_1=F_A(f_{i,0}\, v_i)\, d\zbar/\zbar$.
We obtain that 
$\tau_1$ and $\DD\tau_1$ are 
$L^2$ and of polynomial order,
by using the estimate for
$f_{i,0}\, v_i\, (dz\,d\zbar/|z|^2)$.
Moreover, we have the following formula:
\begin{multline}
\label{eq;07.11.16.5}
 f_{i,0}\, v_i\cdot \frac{dz\, d\zbar}{|z|^2}
-\eigenmap\bigl(\lambda,u(v_i)\bigr)^{-1}
 \,\DD\tau_1
=
 f_{i,0}\, v_i\frac{dz\, d\zbar}{|z|^2}
-\eigenmap\bigl(\lambda,u(v_i)\bigr)^{-1}\,
 F_A\bigl(f_{i,0}\, v_i\bigr) \frac{d\zbar}{\zbar}
 \\
=\eigenmap\bigl(\lambda,u(v_i)\bigr)^{-1}
 \,
 F_{A_0}\bigl(f_{i,0}\, v_i\bigr)
 \frac{d\zbar}{\zbar}
=:\sum B_j\, v_j
\end{multline}
Because $F_{A_0}$ is bounded,
the right hand side of 
(\ref{eq;07.11.16.5})
is also $L^2$ and of polynomial order.
Let us see $B_j$ more closely.
If $a(v_j)=0$,
we have the Fourier expansion
$B_j=\sum_{m\geq 0}
 B_{j,m}(r,\lambda)e^{\sqrt{-1}m\theta}$,
and $B_{j,0}(r,\lambda)=0$ unless 
$\bigl(\gminia(v_j),\alpha(v_j)\bigr)=
 \bigl(\gminia(v_i),\alpha(v_i)\bigr)$
and $N_{a,\alpha,\gminia}v_{i|\nbigd}=v_{j|\nbigd}$.
Note $\deg^W(v_j)<\deg^W(v_i)$
for such $v_j$.

\vspace{.1in}
\noindent
{\bf (B)}\,\,
Let us consider the case in which
$a(v_i)=0$ and $\gminia(v_i)\neq 0$ hold.
Let $k$ be determined by
$\gminia(v_i)=\sum_{j=1}^k\gminia_j(v_i)z^{-j}$
and $\gminia_k(v_i)\neq 0$.
Let $\tau_1$ be given as in (\ref{eq;07.11.16.6}),
which is $L^2$ and of polynomial order
with respect to $h$ and $\gtilde_{\poin}$.
And, we have the following equality,
as in the proof of Lemma \ref{lem;07.11.16.7}:
\begin{multline}
 \DD(\tau_1)
=F_A\Bigl(z^k\, f_{i,0}\, v_i\, \frac{d\zbar}{\zbar}\Bigr)
+\lambda\, k\, 
 z^k\, f_{i,0}\, v_i\frac{dz\, d\zbar}{|z|^2}
 \\
=\Bigl(
 z\frac{\del\gminia(v_i)}{\del z}
+\eigenmap\bigl(\lambda,u(v_i)\bigr)+k\,\lambda
 \Bigr)\, z^k\, f_{i,0}\, v_{i}
 \frac{dz\, d\zbar}{|z|^2}
+z^k\, F_{A_0}\bigl(f_{i,0}\, v_i\bigr)\,\frac{d\zbar}{\zbar}
\end{multline}
Hence, $\DD(\tau_1)$ is also $L^2$ and of polynomial order.
Let $B_j$ be determined by the following:
\[
 f_{i,0}\, v_i\, \frac{dz\, d\zbar}{|z|^2}
-\frac{1}{-k\, \gminia_k(v_i)} \DD\tau_1
=:\sum B_j\, v_j\, \frac{dz\, d\zbar}{|z|^2}
\]
If $a(v_j)=0$,
we have 
$B_j=\sum_{m>0}
 B_{j,m}(r)\, e^{\sqrt{-1}m\theta}$.

\vspace{.1in}
\noindent
{\bf (C)}\,\,
Let us consider
the case in which 
$a(v_i)=0$, $\gminia(v_i)=0$, $\alpha(v_i)=0$
and $k(v_i)=-1$ hold.
Let $i(1)$ be determined by 
$N_{0,0,0} v_{i(1)|\nbigd}=v_{i|\nbigd}$
in $\Gr^{\Fzero,\EEzero}_{(0,0)}(V)$.
Let $\tau_1$ be given by (\ref{eq;07.11.16.10}).
Then, it is $L^2$ and of polynomial order,
and we have the following:
\[
  \DD(\tau_1)
=f_{i,0}\, v_i\, \frac{dz\, d\zbar}{|z|^2}
+F_{A_0'}(\tau_1)
\]
Hence, $\DD(\tau_1)$ is also
$L^2$ and of polynomial order.
Let $B_j$ be determined by the following:
\[
 f_{i,0}\, v_i\,\frac{dz\, d\zbar}{|z|^2}
-\DD(\tau_1)
=\sum B_j\, v_j\, \frac{dz\, d\zbar}{|z|^2}
\]
If $a(v_j)=0$,
we have 
$B_j=\sum_{m>0}
 B_{j,m}(r)\, e^{\sqrt{-1}m\theta}$.

\vspace{.1in}
By using the above preliminary arguments (A), (B) and (C)
with an easy induction,
we can show the second claim of 
Lemma \ref{lem;07.11.16.15}.
\hfill\qed

\subsection{Morphisms
 $\nbigh^j(\varphi_0)$  for $j=0,1$}
\label{subsection;07.11.16.36}

Let us show that
$\nbigh^j(\varphi_0)$  ($j=0,1$)
are quasi-isomorphisms.
Let $\omega$ be a section of
$\nbigl^1_{\poly}(V_{\ast},\DD)$
such that $\DD\omega=0$.
We have the expression
$\omega=f^{1,0}\, dz/z+f^{0,1}\, d\zbar/\zbar$.
We set
\[
  \nbiga^{(0)}(f^{0,1})
=\sum_{v_i\in\nbigb(1)}
 f^{(0,1)}_{i,0}(r)\, v_i,
\quad
\nbiga^{(1)}(f^{0,1}):=
 f^{0,1}-\nbiga^{(0)}(f^{0,1}).
\]
We have the decomposition
$f^{0,1}=\nbiga^{(0)}(f^{0,1})
 +\nbiga^{(1)}(f^{0,1})$.
We have the corresponding decomposition
$\omega^{0,1}=
 \nbiga^{(0)}(\omega^{0,1})
+\nbiga^{(1)}(\omega^{0,1})$.

\begin{lem}
\mbox{{}}\label{lem;07.11.16.25}
\begin{itemize}
\item
 There exists a $\lambda$-holomorphic
 $C^{\infty}$-section $\tau^{(1)}$ of $V$
 such that 
 (i) $L^2$ and of polynomial order,
 (ii) $\delbar\tau^{(1)}=\nbiga^{(1)}(\omega^{0,1})$.
\item
 We have a $\lambda$-holomorphic
 $C^{\infty}$-section $\tau^{(0)}$ of $V$
 such that 
 (i) $L^2$ and of polynomial order,
 (ii) $\nbiga^{(0)}(\omega^{0,1}-\delbar\tau^{(0)})=0$.
\end{itemize}
As a result,
we can take a $\lambda$-holomorphic
$C^{\infty}$-section $\tau$ of $V$
such that 
(i) $L^2$ and of polynomial order,
(ii) $\delbar\tau=\omega^{0,1}$.
\end{lem}
\pf
The argument is almost the same as 
the proof of Lemma \ref{lem;07.9.11.10}.
We have only to replace
$\alpha(v_i)$ with $\eigenmap(\lambda,u(v_i))$.
The first claim follows from Lemma \ref{lem;08.1.27.1}.
Let $C_j$ be the functions determined by the following:
\[
 F_{A_0}\Bigl(f^{0,1}d\zbar/\zbar\Bigr)
=\sum C_j\, v_j\, \frac{dz\, d\zbar}{|z|^2}
\]
From $\DD\omega=0$,
we obtain the following relation
by taking the $v_i$-component:
\begin{equation}
\label{eq;07.11.16.18}
 \lambda\del f_i^{0,1}v_i\frac{d\zbar}{\zbar}
+f_i^{0,1}
 \Bigl(d\gminia(v_i)
+\eigenmap\bigl(\lambda,u(v_i)\bigr)\frac{dz}{z}
 \Bigr)v_i \frac{d\zbar}{\zbar}
+\delbar f_i^{1,0}v_i \frac{dz}{z}
+C_{i} v_i \frac{dz\, d\zbar}{|z|^2}=0
\end{equation}
We use the Fourier expansion
$C_j=\sum C_{j,m}\, e^{\sqrt{-1}m\theta}$.
We give some preliminary arguments.

\vspace{.1in}

\noindent
{\bf (A)} 
Let us consider the case 
that $a(v_i)=0$ and $\gminia(v_i)\neq 0$ hold.
Let $k$ be determined by 
$\gminia(v_i)=\sum_{j=1}^k\gminia_j(v_i)\, z^{-j}$
and $\gminia_k(v_i)\neq 0$.
By looking at the $e^{-\sqrt{-1}k\theta}$-component
of (\ref{eq;07.11.16.18}),
we obtain (\ref{eq;07.11.16.20}),
with $\eigenmap(\lambda,u(v_i))$
instead of $\alpha(v_i)$,
as in the proof of Lemma \ref{lem;07.9.11.10}.
Let $\rho$ be given by (\ref{eq;07.11.16.19}).
If $\nbigk$ is sufficiently small,
we have
$\int|\rho|_h^2 r^{-2\epsilon}
 \dvol_{\gtilde_{\poin}}<\infty$
for some $\epsilon>0$.
By Lemma \ref{lem;08.1.27.1},
we can take a $\lambda$-holomorphic
$C^{\infty}$-section $\rho_1$ such that
(i) $\rho_1\, |z|^{-\epsilon}$ is $L^2$
 and of polynomial order,
(ii) $\delbar\rho_1=\rho$. 
Note (\ref{eq;07.11.16.22}).
Hence, we have a $\lambda$-holomorphic
$C^{\infty}$-section $\tau_2$
such that 
(i) $\tau_2$ is $L^2$ and of polynomial order
 with respect to $h$ and $\gtilde_{\poin}$,
(ii) $f^{0,1}_{i,0}v_i\, d\zbar/\zbar=\delbar\tau_2$.

\vspace{.1in}
\noindent
{\bf (B)}\,\,
We can argue
the case in which $a(v_i)=0$, $\gminia(v_i)=0$
and $\alpha(v_i)\neq 0$ hold,
by using the method in the part (B)
of the proof of Lemma \ref{lem;07.9.11.10}.
We have only to replace
$\alpha(v_i)$ with $\eigenmap(\lambda,u(v_i))$.

\vspace{.1in}
\noindent
{\bf (C)}\,\,
We can argue
the case in which $a(v_i)=0$,
$\gminia(v_i)=0$, $\alpha(v_i)=0$ and $k(v_i)=1$ hold,
by using the method
in the proof of Lemma \ref{lem;07.9.11.10}.

\vspace{.1in}

It is easy to obtain Lemma \ref{lem;07.11.16.25}
by using the above considerations.
\hfill\qed

\vspace{.1in}

We put $\rho:=\omega-\DD\tau$ on 
$\nbigx\setminus\nbigd$,
which gives a holomorphic section of
$V\otimes\Omega^{1,0}$
on $\nbigx\setminus\nbigd$.
We have the decomposition
$\rho=\sum \rho_i$,
where each $\rho_i$ is the product of $v_i$
and a holomorphic $(1,0)$-form
 on $\nbigx\setminus\nbigd$.
\begin{lem}
\label{lem;08.1.27.20}
Let $l(v_i)\in\seisuu_{\geq \,0}$ be determined
as follows:
\begin{itemize}
\item
We put $l(v_i):=-\ord(\gminia(v_i))+1$
in the case $\gminia(v_i)\neq 0$.
\item
We put $l(v_i):=1$ in the case
that $\gminia(v_i)=0$ and $u(v_i)\neq (0,0)$ hold.
\item
We put $l(v_i):=0$ otherwise.
\end{itemize}
Then, $z^{l(v_i)}\rho_i$
is $L^2$
with respect to $h$ and $\gtilde_{\poin}$.
In the second case,
$(-\log|z|)^{-1}\rho_i$ is $L^2$,
more strongly.
\end{lem}
\pf
We consider differentials
only along the direction of $X$
in the following argument.
Let $\delta'$ denote the $(1,0)$-operator
determined by $h$ and $\delbar$.
Let $B$ be determined by $\delta'\vecv=\vecv\, B$.
Then, $B$ is diagonal,
and the $(i,i)$-entries are as follows:
\[
 -\paramap\bigl(\lambda,u(v_i)\bigr)\frac{dz}{z}
+\frac{k(v_i)}{-2\log|z|}\frac{dz}{z}
\]
The curvature $R(h)$ of $\delbar+\delta'$
is expressed by the diagonal matrix
with respect to the frame $\vecv$,
whose $(i,i)$-entries are
$-k(v_i)\, |z|^{-2}(-\log|z|^2)^{-2}
 dz\, d\zbar$.
Hence, $\delta'\tau$  is also $C^{\infty}$
and of polynomial order.
(See the argument in the proof of
Lemma \ref{lem;07.7.18.60}.)

Let $\DD^{(1,0)}$ denote the $(1,0)$-part of $\DD$.
We put $G:=\DD^{(1,0)}-\lambda\delta'$,
which is a section of $\End(V)\otimes\Omega^{1,0}$.
Let $A_1$ be determined by
$G\,\vecv=\vecv\, A_1$.
Then, we have the decomposition $A_1=\Gamma'+C$,
where $F_C$ is bounded with respect to $h$ and $\gtilde_{\poin}$,
and $\Gamma'$ is the diagonal matrix 
whose $(i,i)$-entries are
as follows:
\[
 d\gminia(v_i)
+\Bigl(\eigenmap\bigl(\lambda,u(v_i)\bigr)
+\lambda\, \paramap(\lambda,u(v_i))\Bigr)
 \frac{dz}{z}
\]

We have the decomposition
$\rho=\omega^{1,0}-\lambda\delta'\tau
 -\lambda F_C(\tau)-\lambda F_{\Gamma'}(\tau)$.
Note that
$\omega^{1,0}-\lambda\delta'\tau-\lambda F_C(\tau)$ 
is $L^2$ and of polynomial order
with respect to $h$ and $g_{\poin}$.
Then, the claim of Lemma \ref{lem;08.1.27.20} 
follows.
\hfill\qed

\vspace{.1in}

Let $\Gamma$ and $A_0$ be as in 
Section \ref{subsection;07.11.16.37}.
We put $\DD_0:=\DD-F_{\Gamma}$.
We have
$\DD_0\vecv=\vecv A_0$.
Recall $F_{A_0}$ is bounded
with respect to $h$ and $\gtilde_{\poin}$.

\begin{itemize}
\item
In the case $\gminia(v_i)\neq 0$,
we have a holomorphic section $\kappa_i$
such that
(i) $L^2$ and of polynomial order
 with respect to $h$ and $\gtilde_{\poin}$,
(ii) $\rho_i=
 \bigl(d\gminia(v_i)
 +\eigenmap(\lambda,u(v_i)) dz/z\bigr)
 \kappa_i$.
Note that $\DD_0(\kappa_i)$ is also
$L^2$ and of polynomial order
with respect to $h$ and $\gtilde_{\poin}$.
\item
If $\gminia(v_i)=0$ and $a(v_i)<0$ hold,
we have $(-\log|z|^2)^{-1}\,\rho_i$ is $L^2$.
Then, we obtain that $\rho_i$ is a section of
$\nbigs\bigl(V_{\ast}\otimes\Omega^{1,0},\DD\bigr)$
from the holomorphic property.
\item
If $\gminia(v_i)=0$, $a(v_i)=0$
and $\alpha(v_i)\neq 0$ hold,
we have $z\,\rho_i$ is $L^2$.
Hence, we have 
the $L^2$-holomorphic section $\kappa_i$
such that
$\eigenmap(\lambda,u(v_i))
 \kappa_i dz/z=\rho_i$.
Note $\DD_0(\kappa_i)$ is also $L^2$.
\item
If $\gminia(v_i)=0$, $a(v_i)=0$ and $\alpha(v_i)=0$ hold,
we have $\rho_i$ is contained in
$\nbigs\bigl(V_{\ast}\otimes\Omega^{1,0},\DD\bigr)$.
\end{itemize}

Hence, we obtain the following lemma.
\begin{lem}
\label{lem;07.11.16.30}
There exists a section $\nu$ of
$\nbigl^0_{\poly}(V_{\ast},\DD)$ such that 
$\omega-\DD\nu$ is a 
holomorphic $(1,0)$-form.
\hfill\qed
\end{lem}

From Lemma \ref{lem;07.11.16.30},
it is easy to obtain that
$\nbigh^1(\varphi_0)$ is an isomorphism.
It is easy to show that $\nbigh^0(\varphi_0)$
is an isomorphism.
Thus the proof of Proposition \ref{prop;07.11.16.35}
is finished.
\hfill\qed

\subsection{The morphism $\varphi_1$}
\label{subsection;07.11.16.60}

Let us show that $\varphi_1$ in Proposition
\ref{prop;10.5.23.25} is a quasi isomorphism.
We have only to show that
the induced morphism
$\nbigs(V_{\ast}\otimes\Omega^{\bullet,0})
\lrarr \nbigltilde_{\poly}(V_{\ast},\DD)$
is a quasi-isomorphism.
We have only to show it for the germ
at $(\lambda_0,O)$.
By an easy argument to use a descent
with respect to the ramified covering $X'\lrarr X$,
we may and will assume that $(V_{\ast},\DD)$
is unramified.
Let $\vecvbar_{\gminia}$
and $\vecv_{\gminia,C^{\infty}}$
be as in Section \ref{subsection;07.11.15.30}.
Let $R_{\gminia}$ be determined by
$\DD\vecvbar_{\gminia}
=\vecvbar_{\gminia}\, 
 (d\gminia+R_{\gminia})$.
Let $C$ be determined by the following:
\[
 \DD\vecv_{C^{\infty}}
=\vecv_{C^{\infty}}
 \,\left(
 \bigoplus \bigl(d\gminia+R_{\gminia}\bigr)+C
 \right)
\]
Then, $C$ is $\lambda$-holomorphic,
and it satisfies $C=O\bigl(|z|^N\bigr)$ for any $N>0$.
Let $\DD^{\prime}$ be the family of 
flat $\lambda$-connections
determined by the following:
\[
 \DD^{\prime}\vecv_{C^{\infty}}
=\vecv_{C^{\infty}}
 \,\left(
 \bigoplus\bigl(d\gminia+R_{\gminia}\bigr)
 \right)
\]
The $(0,1)$-part of $\DD^{\prime}$ is 
denoted by $\delbar'$.
We put $F:=\DD-\DD^{\prime}$,
and then we have $|F|_h=O\bigl(|z|^N\bigr)$ for any $N>0$.
We obtain a complex of sheaves
$\nbigltilde_{\poly}^{\bullet}(V_{\ast},\DD^{\prime})$
from $\DD^{\prime}$
as in Section \ref{subsection;07.11.15.30}.
As sheaves, we have
$\nbigltilde_{\poly}^{p}\bigl(V_{\ast},\DD^{\prime}\bigr)
=\nbigltilde_{\poly}^p\bigl(V_{\ast},\DD\bigr)$
for $p=0,1,2$.

\begin{lem}
\label{lem;07.11.16.40}
For any 
$\omega\in\nbigltilde_{\poly}^2(V_{\ast},\DD)$,
we can take 
$\tau\in\nbigltilde^1_{\poly}(V_{\ast},\DD)$
such that $\DD^{\prime}\tau=\omega$.
\end{lem}
\pf
We have only to consider the case
in which $\DD$ has a unique irregular value
$\gminia$.
In the case $\gminia=0$,
we may apply the results
in Section \ref{subsection;07.11.16.50}.
In the case $\gminia\neq 0$,
we have only to use Lemma \ref{lem;07.2.6.11},
for example.
\hfill\qed

\begin{lem}
\label{lem;07.11.16.41}
For any $\omega\in\nbigltilde_{\poly}^2(V_{\ast},\DD)$,
we can take $\tautilde\in\nbigltilde^1_{\poly}(V_{\ast},\DD)$
such that $\DD\tautilde=\omega$.

In particular,
we obtain the vanishing of $\nbigh^2$
of $\nbigltilde^{\bullet}_{\poly}(V_{\ast},\DD)$.
\end{lem}
\pf
Take $\tau$ as in Lemma \ref{lem;07.11.16.40}.
We have 
$\omega-\DD\tau=O(|z|^N)$ for any $N$.
Take some large $M$.
According to Lemma \ref{lem;07.2.6.11},
we can take
a section $\kappa$ of $V\otimes\Omega^{1,0}$
such that
(i) $\delbar\kappa=\omega-\DD\tau$,
(ii) $|\kappa|=O(|z|^M)$.
We have
$\kappa\in\nbigltilde_{\poly}^1(V_{\ast},\DD)$,
and 
$\DD(\tau+\kappa)=\omega$.
Thus, we obtain Lemma \ref{lem;07.11.16.41}
\hfill\qed

\vspace{.1in}

Let $\omega\in\nbigltilde_{\poly}^1(V_{\ast},\DD)$
such that $\DD\omega=0$.
We have
$\DD^{\prime}\omega=
  -F\, \omega=O(|z|^N)$
for any $N$.
Hence, we can take a large $M>0$
and a $\lambda$-holomorphic
$C^{\infty}$-section $\kappa$
of $V\otimes\Omega^{1,0}$
such that $\delbar'\kappa=\DD^{\prime}\omega$
and $|\kappa|=O(|z|^M)$.
We put $\omega':=\omega-\kappa$,
and then $\DD'\omega'=0$.
Note that the $(0,1)$-parts of 
$\omega$ and $\omega'$ are equal.
\begin{lem}
\label{lem;08.1.5.3}
There exists
a local section 
 $\tau\in\nbigltilde^0(V_{\ast},\DD)$
around $(\lambda_0,O)$
such that
$\delbar'\tau=\omega^{0,1}$.
\end{lem}
\pf
We may assume that
$\DD$ has 
a unique irregular value $\gminia$.
In the case $\gminia=0$,
we can apply the result in
Section \ref{subsection;07.11.16.36}.
Let us consider the case $\gminia\neq 0$.
We can take $\tau$ such that
(i) $\delbar'\tau=\omega^{0,1}$,
(ii) $|\tau|=O(|z|^{-M})$ for some large $M$.
Let us show $\DD^{\prime}\tau$
is of polynomial order.
Let $h'$ be the $C^{\infty}$-hermitian metric
of $V_{|\nbigx^{(\lambda_0)}\setminus
 \nbigd^{(\lambda_0)}}$ 
such that
$h'\bigl(v_{C^{\infty},i},v_{C^{\infty},j}\bigr)=
 \delta_{i,j}$.
Note that $h'$ and $h$ are mutually bounded
up to polynomial order.
Let $\delta_1'$ be the $(1,0)$-operator
determined by $\delbar'$ and $h'$.
We consider differentials
only along the direction of $X$.
Note that $z^{M_1}\,\omega^{0,1}$ 
and $z^{M_1}\,\tau$ are
bounded for some large $M_1$.
We also have
$\delbar'(z^{M_1}\,\tau)
=z^{M_1}\,\omega^{0,1}$.
Since the curvature of $R(h,\delbar')$ is $0$,
it can be shown that
$\delta'_1(z^{M_1}\,\tau)$ is $L^2$
uniformly for $\lambda$.
(See the proof of Lemma \ref{lem;07.7.18.60},
for example.)
Thus, we obtain that
$z^{M_1}\, \delta_1'\tau$  is $L^2$
uniformly for $\lambda$.
Taking large $M_2$,
we obtain $z^{M_2}\,(\delta_1'\tau-\omega^{1,0})$
is also $L^2$ uniformly for $\lambda$.

Since $\omega^{1,0}-\delta_1'\tau$ is holomorphic
with respect to $\delbar'$,
we obtain
$\delta_1'\tau-\omega^{1,0}
=O(|z|^{-M_3})$ for some large $M_3$.
Then, we obtain the desired estimate for
$\delta_1'\tau$.
\hfill\qed

\begin{lem}
\label{lem;07.11.16.55}
We can take a section
$\tautilde\in\nbigltilde_{\poly}^0(V_{\ast},\DD)$
such that
$\delbar\tautilde=\omega^{0,1}$
\end{lem}
\pf
Let $\tau$ be as in Lemma \ref{lem;08.1.5.3}.
We have
$\DD\tau-\DD^{\prime}\tau
=F\, \tau=O(|z|^N)$ for any $N>0$.
Because of Lemma \ref{lem;07.2.6.11},
we can take a section $\nu$ of $V$
such that
(i) $|\nu|=O\bigl(|z|^M\bigr)$
 for some large $M>0$,
(ii) $\delbar\nu=F^{0,1}\,\tau$.
Let $\delta'$ be the $(1,0)$-operator determined by
$h$ and $\delbar$.
We consider the differentials 
only in the direction of $X$.
Since the curvature $R(h,\delbar)$ is uniformly bounded
with respect to $h$ and $g_{\poin}$,
it can be shown that $\delta'\nu$ is $L^2$
uniformly for $\lambda$.
If $M$ is sufficiently large,
$\bigl(\DD^{(1,0)}-\lambda\delta'\bigr)\nu$ is
 $O(|z|^{M/2})$.
We put
$\rho=\sum\rho_{\gminia}
 :=\omega-\kappa-\DD\tau+\DD\nu$.
Then, we obtain that
(i) $\rho$ is a holomorphic section 
 of $V\otimes\Omega^{1,0}$,
(ii) $z^{M_1}\,\rho_{\gminia}$ 
 $(\gminia\in\Irr(\DD))$
 are $L^2$ for some large $M_1$
 with respect to $\gtilde_{\poin}$ and $h$
 uniformly for $\lambda$,
(iii) $\rho_0$ is $L^2$ 
 with respect to $h$ and $\gtilde_{\poin}$
 uniformly for $\lambda$.
Hence, 
we obtain
$\nu\in \nbigltilde^0_{\poly}(V_{\ast},\DD)$.
Thus, Lemma \ref{lem;07.11.16.55} is proved.
\hfill\qed

\vspace{.1in}

Let $\nbigstilde(V_{\ast}\otimes\Omega^{p,0})$ be
the sheaf of meromorphic sections $\tau$
of $V\otimes\Omega^{p,0}$ with the following property:
\begin{itemize}
\item Let $\tau_{|\nbigdhat}=\tauhat_{\reg}+\tauhat_{\irr}$ 
 be the decomposition corresponding
 to the irregular decomposition.
 Then, $\tauhat_{\reg}$ is contained in
 $\nbigs\bigl(\Vhat_{\reg\ast}\otimes\Omega^{p,0}\bigr)$.
\end{itemize}

By using Lemma \ref{lem;07.11.16.41}
and Lemma \ref{lem;07.11.16.55},
it is easy to show that the natural inclusion
$\nbigstilde\bigl(V_{\ast}\otimes\Omega^{\bullet,0}\bigr)
 \lrarr
 \nbigltilde_{\poly}(V_{\ast},\DD)$
is a quasi isomorphism.
It is also standard and easy to show that
the natural inclusion
$\nbigs(V_{\ast}\otimes\Omega^{\bullet,0})\lrarr
 \nbigstilde(V_{\ast}\otimes\Omega^{\bullet,0})$
 is a quasi isomorphism.
Hence, we obtain that
$\varphi_1$ is a quasi isomorphism.
Thus, a half of Proposition 
\ref{prop;10.5.23.25} is finished.

\subsection{Proof of 
 Propositions \ref{prop;10.5.23.25}
 and \ref{prop;07.11.16.63}}
\label{subsection;07.11.16.61}

Let $\pi:\nbigxtilde(\nbigd)\lrarr\nbigx$ 
denote the projection.
For any open subset $\nbigu$ of $\nbigxtilde(\nbigd)$,
let $\nbigltilde^{\bullet}_{\poly}
 (V_{\ast},\DD)_{\nbigxtilde(\nbigd)}(\nbigu)$,
denote the space of $\lambda$-holomorphic
$C^{\infty}$-sections $\tau$ of 
$V\otimes\Omega^p$ on
$\nbigu\setminus\pi^{-1}(\nbigd)$
such that the conditions {\bf (a1)}
and {\bf (a2)} are satisfied.
By taking sheafification,
we obtain a complex of sheaves
$\nbigltilde^{\bullet}_{\poly}
 (V_{\ast},\DD)_{\nbigxtilde(\nbigd)}$
on $\nbigxtilde(\nbigd)$.
Similarly, we obtain a complex of sheaves
$\nbiglbar^{\bullet}_{\poly}
 (V_{\ast},\DD)_{\nbigxtilde(\nbigd)}$,
$\nbiglbar^{\bullet}_{\poly}(V,h^{(T)})
 _{\nbigxtilde(\nbigd)}$
and 
$\nbiglbar^{\bullet}_{\poly}
 (V,h,h^{(T)})_{\nbigxtilde(\nbigd)}$
on the real blow up $\nbigxtilde(\nbigd)$,
corresponding to the sheaves
$\nbiglbar^{\bullet}_{\poly}(V_{\ast},\DD)$,
$\nbiglbar^{\bullet}_{\poly}(V,h^{(T)})$
and 
$\nbiglbar^{\bullet}_{\poly}(V,h,h^{(T)})$.

Let $S$ be a small sector in
$\nbigx\setminus\nbigd$
such that
we have the full Stokes filtration 
$\nbigftilde^S$ of $V_{|\Sbar}$.
Let $\Sbar$ denote the closure of $S$
in $\nbigxtilde(\nbigd)$.
We can take a flat splitting
$V_{|S}=\bigoplus V_{\gminia,S}$.
For a $\lambda$-holomorphic $C^{\infty}$-section 
$\tau$ of $V\otimes \Omega^p$
on $S$,
we have the decomposition
$\tau=\sum \tau_{\gminia,S}$
corresponding to $V_{|S}=\bigoplus V_{\gminia,S}$.

\begin{lem}
\mbox{{}}

$\tau$ is a section of
$\nbigltilde_{\poly}^p(V_{\ast},\DD)_{\nbigxtilde(\nbigd)}$
on $\Sbar$,
if and only if the following estimate holds
locally on $\Sbar$
with respect to $h$ and $\gtilde_{\poin}$:
\begin{description}
\item[(a1')]
 $\tau_{\gminia,S}$ and $\DD \tau_{\gminia,S}$
 are of polynomial order 
 for $\gminia\neq 0$.
\item[(a2')]
 $\tau_{0,S}$ and $\DD \tau_{0,S}$ are 
 $L^2$ and of polynomial order.
\end{description}
$\tau$ is a section of
$\nbiglbar_{\poly}^p(V_{\ast},\DD)
 _{\nbigxtilde(\nbigd)}$
on $\Sbar$,
if and only if the following estimate holds
locally on $\Sbar$
with respect to $h$ and $\gtilde_{\poin}$:
\begin{description}
\item[(b1')]
 $\tau_{\gminia,S}$ and $\DD \tau_{\gminia,S}$
 are $O(|z|^N)$ for any $N>0$.
\item[(a2')]
 $\tau_{0,S}$ and $\DD\tau_{0,S}$
 are $L^2$ and of polynomial order.
\end{description}
$\tau$ is a section of 
$\nbiglbar^p_{\poly}(V,h^{(T)})_{\nbigxtilde(\nbigd)}$,
if and only if the following holds
locally on $\Sbar$:
\begin{description}
\item[(c1')]
 $\tau_{\gminia,S}$ and $\DD\tau_{\gminia,S}$
 are $O(|z|^N)$ for any $N>0$
 with respect to $h^{(T)}$
 and $\gtilde_{\poin}$.
\item[(a2')]
 $\tau_{0,S}$ and $\DD\tau_{0,S}$
 are  $L^2$ and of polynomial order
 with respect to $h$ and $\gtilde_{\poin}$.
\end{description}
$\tau$ is a section of 
$\nbiglbar^p_{\poly}(V,h,h^{(T)})
 _{\nbigxtilde(\nbigd)}$,
if and only if the following holds
locally on $\Sbar$:
\begin{description}
\item[(d1')]
 $\tau_{\gminia,S}$ and $\DD\tau_{\gminia,S}$
 are $O(|z|^N)$ for any $N>0$
 with respect to 
 both $(h^{(T)},\gtilde_{\poin})$
 and $(h,\gtilde_{\poin})$.
\item[(a2')]
$\tau_{0,S}$ and $\DD\tau_{0,S}$
are $L^2$  and of polynomial order
with respect to $(h,\gtilde_{\poin})$.
\end{description}
\end{lem}
\pf
It follows from 
Lemma \ref{lem;07.11.16.10}.
\hfill\qed

\begin{lem}
\label{lem;08.9.9.30}
The following natural morphisms
are quasi isomorphisms.
\begin{multline*}
\nbigltilde^{\bullet}_{\poly}
 (V_{\ast},\DD)_{\nbigxtilde(\nbigd)}
\llarr
\nbiglbar^{\bullet}_{\poly}
 (V_{\ast},\DD)_{\nbigxtilde(\nbigd)} \\
\llarr
\nbiglbar^{\bullet}_{\poly}
 (V,h,h^{(T)})_{\nbigxtilde(\nbigd)}
\lrarr
\nbiglbar^{\bullet}_{\poly}(V,h^{(T)})
 _{\nbigxtilde(\nbigd)}
\end{multline*}
\end{lem}
\pf
We have only to consider the case
that $\DD$ has a unique irregular value
$\gminia$.
If $\gminia=0$,
the sheaves are the same.
Let us consider the case $\gminia\neq 0$.
Let $\vecvbar_{\gminia}$ be as in 
Section \ref{subsection;07.11.15.30}.
We have the lift $\vecv_{\gminia,S}$
of $\vecvbar_{\gminia}$ to $V_{\gminia,S}$.
We may take a $\DD$-flat frame $\vecu_{\gminia,S}$
of $V_{\gminia,S}$.
Let $G_{\gminia}$ be determined by
$\vecu_{\gminia,S}
=\vecv_{\gminia,S}\, G_{\gminia}$.
Then,
$G_{\gminia}$ and $G_{\gminia}^{-1}$
are bounded up to polynomial order,
uniformly for $\lambda$
(See Lemma \ref{lem;07.12.27.5},
for example.)
Then, we can show that the vanishing 
of the higher cohomology sheaves of
$\nbigltilde_{\poly}^p(V_{\ast},\DD) 
_{\nbigxtilde(\nbigd)}$,
$\nbiglbar_{\poly}^p(V_{\ast},\DD)
 _{\nbigxtilde(\nbigd)}$,
$\nbiglbar_{\poly}^p(V,h^{(T)})
 _{\nbigxtilde(\nbigd)}$
and $\nbiglbar_{\poly}^p(V,h,h^{(T)})
 _{\nbigxtilde(\nbigd)}$,
by using the results in Section
\ref{subsection;07.11.16.100}.
The comparison of the $0$-th cohomology sheaves
are easy.
Thus, we obtain Lemma \ref{lem;08.9.9.30}.
\hfill\qed

\vspace{.1in}

By applying the push-forward
to the quasi isomorphisms
in Lemma \ref{lem;08.9.9.30},
we obtain 
the rest of Proposition \ref{prop;10.5.23.25}
and Proposition \ref{prop;07.11.16.63}.
\hfill\qed

\chapter{Meromorphic Variation of Twistor Structure}
One of the main results in this monograph
is the reduction from
unramifiedly good wild harmonic bundle
to tame harmonic bundle
(Theorem \ref{thm;07.10.11.120}).
It is convenient to prepare the procedure 
for reduction with respect to Stokes filtrations
in a more general situation.
That is the main purpose 
in Section \ref{subsection;07.10.10.2}.
We introduce the notion of
meromorphic prolongment
of a variation of twistor structure
with a symmetric pairing
(Definition \ref{df;08.9.14.51}),
and we explain the procedure 
to take Gr with respect to Stokes structure
in Section \ref{subsection;08.9.14.52}.

In Section \ref{subsection;08.10.25.20},
we give a review on the notion of variation of polarized
pure twistor structure due to Simpson \cite{s3}
(see also 
 \cite{sabbah2}, \cite{mochi} and \cite{mochi2}).

\section{Variation of polarized pure twistor structure}
\label{subsection;08.10.25.20}
We recall the notion of twistor structure
introduced by Simpson in \cite{s3},
in our convenient way.
See also 
\cite{Hertling},
\cite{sabbah2},
\cite{mochi} and \cite{mochi2}.

\subsection{Some sheaves and differential operators
 on $\proj^1\times X$}

Let $\proj^1$ denote a one dimensional
complex projective space.
We regard it as the gluing of 
two complex lines
$\cnum_{\lambda}$ and $\cnum_{\mu}$
by $\lambda=\mu^{-1}$.
We set $\cnum_{\lambda}^{\ast}:=
 \cnum_{\lambda}-\{0\}$.

Let $X$ be a complex manifold.
We set $\nbigx:=\cnum_{\lambda}\times X$
and $\nbigx^0:=\{0\}\times X$.
Let $\Omegatilde^{1,0}_{\nbigx}$ be 
the $C^{\infty}$-bundle
associated to
$\Omega^{1,0}_{\nbigx}(\log \nbigx^0)
\otimes \nbigo_{\nbigx}(\nbigx^0)$.
We put $\Omegatilde^{0,1}_{\nbigx}:=
\Omega^{0,1}_{\nbigx}$,
and we define
\[
 \Omegatilde^1_{\nbigx}:=
 \Omegatilde^{1,0}_{\nbigx}
\oplus
 \Omegatilde^{0,1}_{\nbigx},
\quad
 \Omegatilde^{\bullet}_{\nbigx}:=
 \mbox{$\bigwedge^{\bullet}$}
 \Omegatilde^{1}_{\nbigx}
\]
The associated sheaves of $C^{\infty}$-sections
are denoted by the same symbols.
Let $\DDtilde_X^f: 
 \Omegatilde^{\bullet}_{\nbigx}
\lrarr
 \Omegatilde^{\bullet+1}_{\nbigx}$
denote the differential operator induced by
the exterior differential $d$ of $\nbigx$.

Let $X^{\dagger}$ denote the conjugate of $X$.
We set $\nbigx^{\dagger}:=
 \cnum_{\mu}\times X^{\dagger}$.
By the same procedure,
we obtain the $C^{\infty}$-bundles
$\Omegatilde^{\bullet}_{\nbigx^{\dagger}}$
with the differential operator
$\DDtilde^{\dagger\,f}_X$.

Their restrictions to 
$\cnum_{\lambda}^{\ast}\times X
=\cnum_{\mu}^{\ast}\times X^{\dagger}$
are naturally isomorphic:
\[
 \bigl(\Omegatilde^{\bullet}_{\nbigx},
 \DDtilde_X^f\bigr)
 _{|\cnum_{\lambda}^{\ast}\times X}
=
 \bigl(\Omega^{\bullet}
 _{\cnum_{\lambda}^{\ast}\times X},d
 \bigr)
=\bigl(\Omegatilde^{\bullet}_{\nbigx^{\dagger}},
 \DDtilde_X^{\dagger\,f}\bigr)
 _{|\cnum_{\mu}^{\ast}\times X^{\dagger}}
\]
By gluing them,
we obtain a graded $C^{\infty}$-bundle
$\Omegatilde^{\bullet}_{\proj^1\times X}$
with a differential operator
$\DDtilde^{\sankaku}_X$.

\begin{rem}
$\DDtilde^f_X$ and 
$\DDtilde^{\dagger\,f}_X$
are denoted also by $d$,
if there is no risk of confusion.
\hfill\qed
\end{rem}

We have the decomposition
$\Omegatilde^1_{\proj^1\times X}
=\xi\Omega^1_{X}
\oplus
 \Omegatilde^1_{\proj^1}$
into the $X$-direction and the $\proj^1$-direction.
The restriction of $\DDtilde^{\sankaku}_X$
to the $X$-direction is denoted by
$\DD^{\sankaku}_X$.
The restriction to the $\proj^1$-direction
is denoted by $d_{\proj^1}$.
We have the decomposition
\[
 \Omegatilde^{1}_{\proj^1}
=\pi^{\ast}
 \Omega^{1,0}_{\proj^1}(2\cdot\{0,\infty\})
\oplus 
 \pi^{\ast}\Omega^{0,1}_{\proj^1},
\]
into the $(1,0)$-part and the $(0,1)$-part,
where $\pi$ denotes the projection
$\proj^1\times X\lrarr \proj^1$.
We have the corresponding decomposition
$d_{\proj^1}=\del_{\proj^1}+\delbar_{\proj^1}$.

\vspace{.1in}

Let $\sigma:\proj^1\lrarr\proj^1$
be the anti-holomorphic involution
given by $\sigma([z_0:z_1])=[-\zbar_1:\zbar_0]$.
The induced diffeomorphism
$\proj^1\times X\lrarr\proj^1\times X$
is also denoted by $\sigma$.
The multiplication on 
$\sigma^{\ast}\Omegatilde^{\bullet}_{\proj^1\times X}$
is twisted as 
$g\cdot \sigma^{\ast}(\omega)
=\sigma^{\ast}\bigl(
 \overline{\sigma^{\ast}(g)}\cdot\omega
 \bigr)$
for a function $g$ and 
a section $\omega$ of 
$\Omegatilde^{\bullet}_{\proj^1\times X}$.
Then, we have the $C^{\infty}$-isomorphism
$\Phi_{\sigma}:
 \sigma^{\ast}\Omegatilde^{\bullet}_{\proj^1\times X}
 \simeq
 \Omegatilde^{\bullet}_{\proj^1\times X}$
given by the complex conjugate
and the ordinary pull back
\[
 \Phi_{\sigma}(\sigma^{\ast}\omega)
=\overline{\sigma^{\ast}(\omega)}.
\]
It is easy to check that 
$\Phi_{\sigma}\circ
 \sigma^{\ast}(\DDtilde_X^{\sankaku})
=\DDtilde_X^{\sankaku}\circ\Phi_{\sigma}$.
Similar relations hold for
$\DD^{\sankaku}_X$ and $d_{\proj^1}$.
If we are given an additional bundle $\nbigf$,
the induced isomorphism
$\nbigf\otimes \sigma^{\ast}\bigl(
 \Omegatilde^{\bullet}_{\proj^1\times X}
 \bigr)
\simeq
 \nbigf\otimes\Omegatilde^{\bullet}_{\proj^1\times X}$
is also denoted by $\Phi_{\sigma}$.

\subsection{Definitions}
\subsubsection{Variation of twistor structure
 (variation of $\proj^1$-holomorphic bundle)}

Let $V$ be a $C^{\infty}$-vector bundle
on $\proj^1\times X$.
We use the same symbol
to denote the associated sheaf of
$C^{\infty}$-sections.
A $\proj^1$-holomorphic structure of
$V$ is defined to be a differential operator
\[
 d''_{\proj^1,V}:
 V\lrarr V\otimes
 \pi^{\ast}\Omega^{0,1}_{\proj^1}
\]
satisfying
(i) $d''_{\proj^1,V}(f\cdot s)
=f\cdot d''_{\proj^1,V}(s)
+\delbar_{\proj^1}(f)\cdot s$ 
for a $C^{\infty}$-function $f$
and a section $s$ of $V$,
(ii) $d''_{\proj^1,V}\circ d''_{\proj^1,V}=0$.
Such a tuple $(V,d''_{\proj^1,V})$
is called a $\proj^1$-holomorphic vector bundle.
\index{$\proj^1$-holomorphic vector bundle}

A $\TTtilde$-structure of
$(V,d''_{\proj^1,V})$ is a differential operator
\[
 \DD^{\sankaku}_V:
 V\lrarr V\otimes\xi\Omega^1_X
\]
such that
(i) $\DD^{\sankaku}_{V}(f\cdot s)
=f\cdot \DD^{\sankaku}_{V}(s)
+\DD^{\sankaku}_X(f)\cdot s$ 
for a $C^{\infty}$-function $f$
and a section $s$ of $V$,
(ii) $(d''_{\proj^1,V}+\DD^{\sankaku}_{V})^2=0$.
Such a tuple 
$(V,d''_{\proj^1,V},\DD^{\sankaku}_V)$
is called {\em a $\TTtilde$-structure} in \cite{Hertling},
or {\em a variation of 
$\proj^1$-holomorphic vector bundle}
in \cite{mochi2}.
In this section, we prefer to call it
{\em variation of twistor structure.}
We will not distinguish them.
\index{variation of $\proj^1$-holomorphic bundle}
\index{variation of twistor structure}
\index{$\TTtilde$-structure}

If $X$ is a point,
it is just a holomorphic vector bundle
on $\proj^1$.

\begin{rem}
We will often omit to specify $d''_{\proj^1,V}$
when we consider 
$\proj^1$-holomorphic bundles
or variations of twistor structure
(variations of $\proj^1$-holomorphic bundle).
\hfill\qed
\end{rem}

A morphism of variation of twistor structure
\[
 F:(V_1,d''_{\proj^1,V_1}, \DD^{\sankaku}_{V_1})
\lrarr (V_2,d''_{\proj^1,V_2},
 \DD^{\sankaku}_{V_2})
\]
is defined to be 
a morphism of the associated sheaves
of $C^{\infty}$-sections,
compatible with the differential operators.
If $X$ is a point,
it is equivalent to an $\nbigo_{\proj^1}$-morphism.

\subsubsection{Some functoriality}

Let $(V,\DD^{\sankaku}_V)$
be a variation of twistor structure.
Let $f:Y\lrarr X$ be a holomorphic map
of complex manifolds.
Then, we have the naturally induced
variation of twistor structure
$f^{\ast}(V,\DD^{\sankaku}_V)$
as in the case of ordinary connections.

Let $\sigma:\proj^1\lrarr\proj^1$ be as above.
Then, 
$\sigma^{\ast}V$
is naturally equipped with
a $\proj^1$-holomorphic structure and 
a $\TTtilde$-structure
$\DD^{\sankaku}_{\sigma^{\ast}V}$
given as follows:
\[
 \bigl(
 \DD^{\sankaku}_{\sigma^{\ast}V}
 +d''_{\sigma^{\ast}V}
\bigr)
\bigl(
 \Phi_{\sigma}(\sigma^{\ast}s)
\bigr)
=\Phi_{\sigma}\Bigl(
 \sigma^{\ast}\bigl(
 (\DD^{\sankaku}_V+d''_V)s\bigr)
 \Bigr)
\]
Here, $s$ denotes a section of
$V\otimes \xi\Omega_X^{\bullet}$.
Thus, we obtain the pull back of
variation of twistor structure by $\sigma$.

Direct sum,
tensor product, and dual
for variation of twistor structure
are defined in obvious manners.

\subsubsection{Variation of pure twistor structure}

Let $(V,d''_{\proj^1,V})$ be 
a $\proj^1$-holomorphic vector bundle
on $\proj^1\times X$.
It is called pure of weight $w$
if the restrictions 
$V_P:=(V,d''_{\proj^1,V})_{|\proj^1\times\{P\}}$
are pure twistor structure of weight $w$
for any $P\in X$,
i.e., $V_P$ are isomorphic to
direct sums of $\nbigo_{\proj^1}(w)$.
A variation of twistor structure
is called pure of weight $w$,
if the underlying $\proj^1$-holomorphic
vector bundle is pure of weight $w$.
\index{variation of pure twistor structure}

\subsubsection{Example (Tate objects)}

\index{Tate object}
Let $\Tate(w)$ be a Tate object 
in the theory of twistor structure.
(See \cite{s3} and 
 Section 3.3.1 of \cite{mochi2}.)
It is isomorphic to $\nbigo_{\proj^1}(-2w)$,
and equipped with the distinguished frames
\[
 \Tate(w)_{|\cnum_{\lambda}}
=\nbigo_{\cnum_{\lambda}}\,
 t^{(w)}_0,
\quad
  \Tate(w)_{|\cnum_{\mu}}
=\nbigo_{\cnum_{\mu}}\,
 t^{(w)}_{\infty},
\quad
   \Tate(w)_{|\cnum_{\lambda}^{\ast}}
=\nbigo_{\cnum_{\lambda}^{\ast}}\,
 t^{(w)}_{1}.
\]
The transformation is given by
\[
  t_0^{(w)}=(\sqrt{-1}\lambda)^{w}\, t_1^{(w)},
\quad
 t_{\infty}^{(w)}=(-\sqrt{-1}\mu)^w\, t_1^{(w)}.
\]
In particular,
$(\sqrt{-1}\lambda)^{-2w}t_0^{(w)}=
 t_{\infty}^{(w)}$.

We may identify $\Tate(w)$ with
$\nbigo_{\proj^1}\bigl(-w\cdot 0-w\cdot\infty\bigr)$
by the correspondence $t_1^{(w)}\longleftrightarrow 1$,
up to constant multiplication.
In particular,
we will implicitly use the identification
of $\Tate(0)$ with $\nbigo_{\proj^1}$
by $t_1^{(0)}\longleftrightarrow 1$.
We will also implicitly use 
the identification
$\Tate(m)\otimes\Tate(n)
\simeq \Tate(m+n)$
given by 
$t_a^{(m)}\otimes
 t_a^{(n)}
\longleftrightarrow
 t_a^{(m+n)}$.

\subsubsection{Example}

In Section 3.3.2 of \cite{mochi2},
we considered 
a line bundle $\nbigo(p,q)$
on $\proj^1$
with a natural $\cnum^{\ast}$,
which is isomorphic to $\nbigo_{\proj^1}(p+q)$
and equipped with the distinguished frames:
\[
 \nbigo(p,q)_{|\cnum_{\lambda}}
=\nbigo_{\cnum_{\lambda}}
 \, f_0^{(p,q)},
\quad
 \nbigo(p,q)_{|\cnum_{\mu}}
=\nbigo_{\cnum_{\mu}}
 \, f_{\infty}^{(p,q)},
\quad
 \nbigo(p,q)_{|\cnum_{\lambda}^{\ast}}
=\nbigo_{\cnum_{\lambda}^{\ast}}
 \, f_{1}^{(p,q)}.
\]
The transformation is given by 
\[
 f_0^{(p,q)}
=(\sqrt{-1}\lambda)^{-p}\,
 f_1^{(p,q)},
 \quad
  f_{\infty}^{(p,q)}
=(-\sqrt{-1}\mu)^{-q}\,
 f_1^{(p,q)}.
\]
In particular, 
$(\sqrt{-1}\lambda)^{p+q}f_0^{(p,q)}
=f_{\infty}^{(p,q)}$.

We may identify
$\nbigo(p,q)$ with
$\nbigo_{\proj^1}(p\cdot 0+q\cdot\infty)$
by the correspondence
$f_1^{(p,q)}\longleftrightarrow 1$,
up to constant multiplication.
We will implicitly use the identification
$\nbigo(p,q)\otimes\nbigo(p',q')
\simeq
 \nbigo(p+p',q+q')$
given by
$f_a^{(p,q)}\otimes
 f_a^{(p',q')}
\longleftrightarrow
 f_a^{(p+p',q+q')}$.
We will also implicitly identify 
$\Tate(w)$ with $\nbigo(-w,-w)$
by $t_a^{(w)}=f_a^{(-w,-w)}$
for $a=0,1,\infty$.

If we forget the natural $\cnum^{\ast}$-actions,
$\nbigo(p,q)$ and $\nbigo(p+r,q-r)$
are identified by
$f_{\kappa}^{(p,q)}\longleftrightarrow
 f_{\kappa}^{(p+r,q-r)}$
for $\kappa=0,\infty$.
In that case,
$f_{\kappa}^{(p,q)}$ are denoted by
$f_{\kappa}^{(p+q)}$.

\vspace{.1in}

Let $X$ be a complex manifold.
We have the pull back of
$\Tate(w)$ and $\nbigo(p,q)$
via the map from $X$ to a point.
They are denoted by
$\Tate(w)_X$ and $\nbigo(p,q)_X$,
respectively.
We will often omit the subscript $X$,
if there is no risk of confusion.
For a variation of twistor structure
$(V,\DD^{\sankaku})$,
the tensor product
$(V,\DD^{\sankaku})\otimes\Tate(w)$
is called the $w$-th Tate twist
of $(V,\DD^{\sankaku})$.
\index{Tate twist}

\subsubsection{Polarization}

Recall that we have the isomorphism (\cite{mochi2})
\[
 \iota_{\Tate(w)}:
 \sigma^{\ast}\Tate(w)
 \simeq
 \Tate(w),
\]
given by the natural identification
$\sigma^{\ast}\nbigo\bigl(-w\cdot 0-w\cdot\infty\bigr)
\simeq
 \nbigo\bigl(-w\cdot 0-w\cdot\infty\bigr)$
via $\sigma^{\ast}(1)\longleftrightarrow 1$,
or equivalently,
\[
 \sigma^{\ast}t_1^{(w)}\longleftrightarrow
 t_1^{(w)},
\quad
 \sigma^{\ast}t_0^{(w)}\longleftrightarrow
 (-1)^w\, t_0^{(w)},
\quad
  \sigma^{\ast}t_{\infty}^{(w)}\longleftrightarrow
 (-1)^w\, t_{\infty}^{(w)}.
\]

For a variation of twistor structure
$(V,\DD^{\sankaku}_V)$ on $\proj^1\times X$,
a morphism
\[
 \nbigs:(V,\DD^{\sankaku}_V)
\otimes
 \sigma^{\ast}(V,\DD^{\sankaku}_V)
\lrarr \Tate(-w)_X
\]
is called a pairing of weight $w$,
if it is $(-1)^w$-symmetric
in the following sense:
\[
  \iota_{\Tate(-w)}\circ
 \sigma^{\ast}\nbigs
=(-1)^{w}
 \nbigs\circ
 \exchange:
 \sigma^{\ast}V\otimes V\lrarr \Tate(-w)_X
\]
Here, $\exchange$ denotes the natural morphism
$\sigma^{\ast}V\otimes V\lrarr
 V\otimes\sigma^{\ast}V$
induced by the exchange of the components.
It is also called $(-1)^w$-symmetric pairing,
if we would like to emphasize $(-1)^w$-symmetric
property.
\index{pairing}
\index{$(-1)^w$-symmetric pairing}

Let $(V,\DD^{\sankaku}_V)$ be 
a variation of {\em pure} twistor structure
of weight $w$ on $\proj^1\times X$.
Let $\nbigs:(V,\DD^{\sankaku}_V)\otimes
 \sigma^{\ast}(V,\DD^{\sankaku}_V)\lrarr\Tate(-w)_X$
be a pairing of weight $w$.
We say that $\nbigs$ is a polarization of 
$(V,\DD^{\sankaku}_V)$,
if $\nbigs_{P}:=\nbigs_{|\proj^1\times\{P\}}$
is a polarization of
$V_{P}:=
(V,d''_{\proj^1})_{|\proj^1\times \{P\}}$
for each $P\in X$.
Namely, the following holds:
\begin{itemize}
\item
 If $w=0$,
 the induced Hermitian pairing
 $H^0(\nbigs_P)$ of $H^0(\proj^1,V_{P})$
 is positive definite.
\item
 In the general case,
 the induced pairing
 $\nbigs_P\otimes\nbigs_{0,-w}$
 of $V_P\otimes\nbigo(0,-w)$
 is a polarization of
 the pure twistor structure. 
 (See Example $2$ below for
 $\nbigs_{0,-w}$.)
\end{itemize}
\index{polarization}

When $S$ is a polarization of 
a pure twistor structure $V$ with weight $n$,
the induced pairing 
$\sigma(S):
 \sigma^{\ast}(V)\otimes V\lrarr \Tate(-n)$
and
$S^{\lor}:V^{\lor}\otimes
 \sigma^{\ast}(V^{\lor})
\lrarr\Tate(n)$ are also polarizations.
(See Lemma 3.38 of \cite{mochi2}.)

\subsubsection{Example 1}

The identification $\iota_{\Tate(w)}$
induces a flat morphism
$\nbigs_{\Tate(w)}:
 \Tate(w)\otimes \sigma^{\ast}\Tate(w)
 \lrarr  \Tate(2w)$.
It is a polarization of $\Tate(w)$ of weight $-2w$.

\subsubsection{Example 2}
The flat isomorphism
$\iota_{(p,q)}:
 \sigma^{\ast}\nbigo(p,q)
\simeq
 \nbigo(q,p)$ in \cite{mochi2} is given by
$ \sigma^{\ast}f_0^{(p,q)}
 \longmapsto
 (\sqrt{-1})^{p+q}f_{\infty}^{(q,p)}$,
$ \sigma^{\ast}f_{\infty}^{(p,q)}
\longmapsto
 (-\sqrt{-1})^{p+q}
 f_0^{(q,p)}$,
and
$\sigma^{\ast}f_1^{(p,q)}
 \longmapsto
 (\sqrt{-1})^{q-p}f_1^{(q,p)}$.
Hence, we obtain the morphism
\[
 \nbigs_{p,q}:
 \nbigo(p,q)\otimes\sigma^{\ast}\nbigo(p,q)
\lrarr \Tate(-p-q).
\]
it is a polarization of weight $p+q$.

\begin{rem}
It is essential to fix an isomorphism
$\iota:\sigma^{\ast}\nbigo_{\proj^1}(1)
\simeq\nbigo_{\proj^1}(1)$
such that $\sigma^{\ast}\iota\circ\iota=-1$.
It is unique up to isomorphism.
There could be a choice of a frame 
to reduce signatures.
\hfill\qed
\end{rem}

\subsubsection{Relation with harmonic bundles}

Simpson observed the equivalence between
the notions of variation of polarized pure twistor structure
and harmonic bundles.
(See \cite{s3}. See also \cite{mochi2} and \cite{sabbah2}.)
Let $p:\proj^1\times X\lrarr X$ be the projection.
Let $\harmonicbundle$ be a harmonic bundle on $X$.
We set $\nbige^{\sankaku}:=p^{\ast}E$,
which is naturally a $\proj^1$-holomorphic bundle.
It is equipped with the differential operator
$\DD^{\sankaku}:=
 \delbar_E+\lambda\theta^{\dagger}
+\del_E+\lambda^{-1}\theta$,
and $(\nbige^{\sankaku},\DD^{\sankaku})$
is a variation of pure twistor structure of weight $0$.
The polarization $\nbigs$ is given by
$\nbigs(u\otimes\sigma^{\ast}v):=
 p^{\ast}(h)(u,\sigma^{\ast}v)$.

\subsection{Gluing construction}
\label{subsection;08.7.29.10}

Recall the gluing construction of
variation of pure twistor structure in \cite{s3}.
See also \cite{mochi2}.
We have the decomposition
$ \Omegatilde^1_{\nbigx}
=\xi\Omegatilde^1_{X|\nbigx}
\oplus
 \Omegatilde^1_{\cnum_{\lambda}}$
into the $X$-direction
and the $\cnum_{\lambda}$-direction.
Let $d_X$ denote the restriction of 
the exterior differential to the $X$-direction.
Similarly,
we have the decomposition
$ \Omegatilde^1_{\nbigx^{\dagger}}
=\xi\Omegatilde^1_{X|\nbigx^{\dagger}}
\oplus
 \Omegatilde^1_{\cnum_{\mu}}$,
and the restriction of $\DDtilde^{\dagger\,f}_X$
to the $X$-direction is denoted by $d_{X^{\dagger}}$.
The notions of 
$\cnum_{\lambda}$-holomorphic bundles
or $\cnum_{\mu}$-holomorphic bundles
are defined
as in the case of $\proj^1$-holomorphic bundles.

Let $(V_0,d''_{\cnum_{\lambda},V_0})$ 
be a $\cnum_{\lambda}$-holomorphic 
bundle on $\nbigx$.
A $T$-structure \cite{Hertling} of $V_0$
is a differential operator
\[
 \DD^f_{V_0}:V_0\lrarr 
 V_0\otimes
 \xi\Omega^1_{X|\nbigx}
\]
satisfying 
(i) $\DD^f_{V_0}(f\cdot s)
=d_Xf\cdot s+f\cdot \DD^f_{V_0}(s)$
for a function $f$ and a section $s$ of $V$,
(ii) 
 $\bigl(
 d''_{\cnum_{\lambda},V_0}+\DD^f_{V_0}
\bigr)^2=0$.

Let $(V_{\infty},d''_{\cnum_{\mu},V_{\infty}})$
be a $\cnum_{\mu}$-holomorphic vector bundle
on $\nbigx^{\dagger}$.
A $\Ttilde$-structure \cite{Hertling}
is defined to be a differential operator
\[
 \DD^{\dagger\,f}_{V_{\infty}}:
 V_{\infty} \lrarr 
 V_{\infty}\otimes
 \xi\Omega^1_{X|\nbigx^{\dagger}}
\]
satisfying conditions
similar to (i) and (ii) above.

Assume that we are given 
an isomorphism $\Phi$:
\begin{equation}
\label{eq;08.9.10.30}
 \Phi:
 (V_0,d''_{\cnum_{\lambda},V_0},
 \DD^f_{V_0})
 _{|\cnum_{\lambda}^{\ast}\times X}
\simeq
 (V_{\infty},d''_{\cnum_{\mu},V_{\infty}},
 \DD^{\dagger\,f}_{V_{\infty}})_{
 |\cnum_{\mu}^{\ast}\times X^{\dagger}}
\end{equation}
We obtain a $C^{\infty}$-vector bundle $V$
on $\proj^1\times X$
by gluing $V_0$ and $V_{\infty}$ via $\Phi$.
By the condition (\ref{eq;08.9.10.30}),
$d''_{\cnum_{\lambda},V_0}$
and $d''_{\cnum_{\mu},V_{\infty}}$
give $\proj^1$-holomorphic structure 
$d''_{\proj^1,V}$,
and 
$\DD^f_{V_0}$
and $\DD^{\dagger\,f}_{V_{\infty}}$
induce the $\TTtilde$-structure
$\DD^{\sankaku}_V$.
Thus, we obtain a variation of
twistor structure
$(V,d''_{\proj^1,V},\DD^{\sankaku}_V)$.

Conversely,
we naturally obtain such 
$ (V_0,d''_{\cnum_{\lambda},V_0},
 \DD^f_{V_0})$,
$(V_{\infty},d''_{\cnum_{\mu},V_{\infty}},
 \DD^{\dagger\,f}_{V_{\infty}})$
and 
$\Phi$
from a variation of twistor structure
$(V,d''_{\proj^1,V},\DD^{\sankaku}_V)$
as the restriction to
$\nbigx$ and $\nbigx^{\dagger}$,
respectively.

Under the natural isomorphism
\[
 \xi\Omega^1_{X|\nbigx}
=\lambda^{-1}\cdot
  \Omega^{1,0}_{\nbigx/\cnum_{\lambda}}
 \oplus
  \Omega^{0,1}_{\nbigx/\cnum_{\lambda}}
\simeq
 \Omega^{1,0}_{\nbigx/\cnum_{\lambda}}
 \oplus
 \Omega^{0,1}_{\nbigx/\cnum_{\lambda}}
=\Omega^1_{\nbigx/\cnum_{\lambda}}
\]
a $T$-structure $\DD^f_{V_0}$ induces
a holomorphic family of flat $\lambda$-connections
$\DD_{V_0}$.
Similarly,
a $\Ttilde$-structure of 
$\DD^{\dagger\,f}_{V_{\infty}}$
naturally induces 
a holomorphic family of flat $\mu$-connections
$\DD^{\dagger}_{V_{\infty}}$.
Hence, 
a variation of twistor structure
is regarded as the gluing of
families of $\lambda$-flat bundles
and $\mu$-flat bundles.

\subsubsection{}
Let $(V,\DD^{\sankaku}_V)$
be a variation of twistor structure
on $\proj^1\times(X-D)$.
Let $\DD_{\sigma^{\ast}V_{\infty}}$
(resp. $\DD^{\dagger}_{\sigma^{\ast}V_0}$)
denote the associated family of 
flat $\lambda$-connections
(resp. $\mu$-connections)
on $\sigma^{\ast}V_{\infty}$
(resp. $\sigma^{\ast}V_{0}$).
The following lemma can be checked
by an easy and direct calculation.
We remark the signature.
\begin{lem}
\label{lem;08.1.29.1}
Let $f$ be a local section of $V_{\infty}$,
and let $A_i$ and $B_i$ $(i=1,\ldots,n)$
be determined by
$\DD^{\dagger}f=
 \sum A_i\cdot dz_i+\sum B_i\cdot d\zbar_i$,
where $A_i$ and $B_i$ are local sections of
$V_{\infty}$.
Then, we have
\[
  \DD_{\sigma^{\ast}V_{\infty}}
 (\sigma^{\ast}f)
=\sum \sigma^{\ast}(A_i)\cdot d\zbar_i
-\sum \sigma^{\ast}(B_i)\cdot dz_i.
\]
Similarly, 
we have
\[
 \DD^{\dagger}_{\sigma^{\ast}V_0}(\sigma^{\ast}g)
=-\sum \sigma^{\ast}(A_i)\cdot d\zbar_i
 +\sum \sigma^{\ast}(B_i)\cdot dz_i
\]
for a local section $g$ of $V_0$
with $\DD g=
 \sum A_i\cdot dz_i+\sum B_i\cdot d\zbar_i$.
\hfill\qed
\end{lem}

\section[Meromorphic variation of
 twistor structure]{Good meromorphic prolongment
of variation of twistor structure}
\label{subsection;07.10.10.2}
We shall introduce the notion of
unramifiedly good meromorphic prolongment
for variation of twistor structure.
We shall also observe that
a graded variation of twistor structure
is obtained as the graduation
with respect to Stokes structure.

\subsection{Unramifiedly good meromorphic prolongment}
\label{subsection;07.12.19.20}

Let $X$ be a complex manifold
with a normal crossing hypersurface $D$.
We put $\nbigx:=\cnum_{\lambda}\times X$
and $\nbigx^{\shikaku}:=
 \cnum_{\lambda}^{\ast}\times X$.
Let $p_{\lambda}$ be the projection
forgetting the $\lambda$-component.
For any subset $\nbigk$ of 
$\cnum_{\lambda}$,
we put $\nbigx_{\nbigk}:=\nbigk\times X$.
We use the symbols $\nbigd$,
$\nbigd^{\shikaku}$,
and $\nbigd_{\nbigk}$
in similar meanings.

\begin{df}
\label{df;10.5.19.20}
Let $(V_0,\DD)$ be a family of
$\lambda$-flat bundles on $\nbigx\setminus\nbigd$.
A family of meromorphic $\lambda$-flat bundles
$(\Vtilde_0,\DD)$ on $(\nbigx_{\nbigk},\nbigd_{\nbigk})$
is called an unramifiedly 
good meromorphic prolongment
of $(V_0,\DD)$,
if the following holds:
\index{unramifiedly good meromorphic prolongment}
\begin{itemize}
\item
 The restriction of $(\Vtilde_0,\DD)$ to 
 $\nbigx_{\nbigk}\setminus\nbigd_{\nbigk}$
 is $(V_0,\DD)_{|\nbigx_{\nbigk}
 \setminus\nbigd_{\nbigk}}$.
\item
 $(\Vtilde_0,\DD)$ locally has 
 an unramifiedly good lattice,
 i.e., 
 for each $P\in\nbigd_{\nbigk}$,
 there exists a small neighbourhood 
 $\nbigx_{P}$
 such that $(\Vtilde_0,\DD)_{|\nbigx_{P}}$
 has an unramifiedly good lattice.
\item
$\Irr(\Vtilde_0,\DD,P)\subset
 \nbigo_{X}(\ast D)_{p_{\lambda}(P)}
 \big/\nbigo_{X,p_{\lambda}(P)}$,
i.e.,
the elements of $\Irr(\Vtilde_0,\DD,P)$ 
are independent of the variable $\lambda$.
\hfill\qed
\end{itemize}
\end{df}
Under the third condition,
we have $\Irr(\Vtilde_0,\DD,P)=\Irr(\Vtilde_0,\DD,P')$
for any $p_{\lambda}(P)=p_{\lambda}(P')$,
if $\nbigk$ is connected.
In that case, for $R\in D$,
we take $P\in\nbigd$ such that 
$p_{\lambda}(P)=R$,
and put $\Irr(\Vtilde_0,\DD,R):=
\Irr(\Vtilde_0,\DD,P)$.
They will be denoted also by
$\Irr(\Vtilde_0,R)$ or $\Irr(\DD,R)$.

If we are interested only on
family of good filtered $\lambda$-flat bundles,
it is too strong to impose 
the independence from $\lambda$
for irregular values.
However, it seems appropriate to impose it
when we consider meromorphic prolongment
of a variation of twistor structure.

\subsection{Meromorphic prolongment of 
a variation of twistor structure}

Let $X^{\dagger}$ denote the conjugate of $X$.
We put $\nbigx^{\dagger}:=
 \cnum_{\mu}\times X^{\dagger}$
and $\nbigx^{\dagger\shikaku}:=
 \cnum_{\mu}^{\ast}\times X^{\dagger}$.
Let $p_{\mu}$ denote the projection
forgetting the $\mu$-component.
For any subset $\nbigh$ of $\cnum_{\mu}$,
we put $\nbigx^{\dagger}_{\nbigh}:=
 \nbigh\times X^{\dagger}$.
We use the symbols 
$D^{\dagger}$,
$\nbigd^{\dagger}$,
$\nbigd^{\dagger\shikaku}$,
and $\nbigd^{\dagger}_{\nbigh}$
in similar meanings.
We use the $C^{\infty}$-identification
$\nbigx^{\shikaku}=\nbigx^{\dagger\shikaku}$
given by $\lambda=\mu^{-1}$,
which preserves the $\cnum_{\lambda}^{\ast}$-holomorphic
structure.

Let $(V,\DD^{\sankaku})$
be a variation of twistor structure
on $(X\setminus D)\times\proj^1$.
We have the associated family of $\lambda$-flat
bundles
$(V_0,\DD)$ on $\nbigx\setminus\nbigd$,
and the associated family of $\mu$-flat bundles
$(V_{\infty},\DD^{\dagger})$ on
$\nbigx^{\dagger}\setminus\nbigd^{\dagger}$.
Let $\DD^f$ and 
$\DD^{\dagger\,f}$ denote
the associated families of flat connections.
We have the isomorphism
$V_{0|\nbigx^{\shikaku}\setminus\nbigd^{\shikaku}}
 \simeq
 V_{\infty|
 \nbigx^{\dagger\shikaku}
 \setminus\nbigd^{\dagger\shikaku}}$
preserving 
the families of the flat connections
and the holomorphic structures along 
the $\cnum_{\lambda}^{\ast}$-direction.

Let $\nbigk\subset\cnum_{\lambda}$
and $\nbigh\subset\cnum_{\mu}$
be connected compact regions
such that the union of the interior points 
of $\nbigk$ and $\nbigh$ is $\proj^1$.
Assume that we are given the following:
\begin{itemize}
\item
An unramifiedly good meromorphic prolongment 
$(\Vtilde_0,\DD)$ on 
 $(\nbigx_{\nbigk},\nbigd_{\nbigk})$
of $(V_0,\DD)$.
\item
An unramifiedly good meromorphic prolongment 
 $(\Vtilde_{\infty},\DD^{\dagger})$ on 
$(\nbigx^{\dagger}_{\nbigh},
 \nbigd^{\dagger}_{\nbigh})$ of 
 $(V_{\infty},\DD^{\dagger})$.
\item
For any $R\in D=D^{\dagger}$,
the sets
$\Irr(\Vtilde_{\infty},\DD^{\dagger},R)$
and $\Irr(\Vtilde_0,\DD,R)$ are related as
\[
 \Irr(\Vtilde_{\infty},\DD^{\dagger},R)
=\bigl\{
 \gminiabar\,\big|\,\gminia\in\Irr(\Vtilde_0,\DD,R)
 \bigr\}.
\]
\end{itemize}

Let $\pi:\nbigxtilde^{\shikaku}(\nbigd^{\shikaku})
\lrarr \nbigx^{\shikaku}$ denote the real blow up
of $\nbigx^{\shikaku}$ along $\nbigd^{\shikaku}$.
We have the $\nbigo_{
 \cnum_{\lambda}^{\ast}}$-module
$\gbigv_0$ on 
$\nbigxtilde^{\shikaku}(\nbigd^{\shikaku})$
associated to 
$(V_0,\DD)_{|\nbigx^{\shikaku}\setminus
 \nbigd^{\shikaku}}$
as in Subsection \ref{subsection;10.5.19.10}.
For any $Q\in \pi^{-1}\bigl(
 \nbigd^{\shikaku}_{\nbigk}\bigr)$,
we have the full Stokes filtration 
$\nbigftilde^Q(\gbigv_{0,Q})$
of the stalk of $\gbigv_{0}$ at $Q$
for the meromorphic extension $\Vtilde_0$.
Similarly,
let $\pi:\nbigxtilde^{\dagger\,\shikaku}
 (\nbigd^{\dagger\,\shikaku})
\lrarr \nbigx^{\dagger}$ denote the real blow up
of $\nbigx^{\dagger\,\shikaku}$
along $\nbigd^{\dagger\,\shikaku}$,
and let $\gbigv_{\infty}$ denote the sheaf on 
$\nbigxtilde^{\dagger\,\shikaku}
 (\nbigd^{\dagger\,\shikaku})$
associated to $V_{\infty}$.
For any point 
$Q\in \pi^{-1}\bigl(
 \nbigd^{\dagger\,\shikaku}_{\nbigh}
 \bigr)$,
we have the full Stokes filtration
$\nbigftilde^Q(\gbigv_{\infty,Q})$.

Note the natural identifications
$\nbigxtilde^{\shikaku}(\nbigd^{\shikaku})
=\nbigxtilde^{\dagger\,\shikaku}
(\nbigd^{\dagger\,\shikaku})$
and $\gbigv_0=\gbigv_{\infty}$,
and hence $\gbigv_{0,P}=\gbigv_{\infty,P}$
for any 
$P\in \pi^{-1}(\nbigd_{\nbigk}
 \cap\nbigd^{\dagger}_{\nbigh})$.
We also remark
\[
 \Re(\mu^{-1}\gminiabar)
=|\mu|^{-2}\Re(\lambda^{-1}\gminia)
\]
for $\lambda=\mu^{-1}$.
Hence, the natural bijection
$\Irr(\Vtilde_0,P)\lrarr
 \Irr(\Vtilde_{\infty},P)$ induces an isomorphism
of ordered sets 
$(\Irr(\Vtilde_0,P),\leq^{\lambda}_Q)$ and 
$(\Irr(\Vtilde_{\infty},P),\leq^{\mu}_Q)$
for any $Q\in \pi^{-1}(P)$.

\begin{df}
\label{df;07.11.13.20}
\mbox{{}}
\begin{itemize}
\item
We say that the Stokes structure of
$\Vtilde_0$ and $\Vtilde_{\infty}$ are the same,
if the filtrations
 $\nbigftilde^Q(\gbigv_{0,Q})$ and
 $\nbigftilde^{Q}(\gbigv_{\infty,Q})$
are the same
for any 
 $Q\in\pi^{-1}(
 \nbigd_{\nbigk}\cap
 \nbigd^{\dagger}_{\nbigh})$,
 under the above identification of
 the index sets.
\index{Stokes structures are the same}
\item
If the Stokes structures of $\Vtilde_0$ and $\Vtilde_{\infty}$
are the same,
$(\Vtilde_0,\Vtilde_{\infty})$ is called
an unramifiedly good meromorphic prolongment
of variation of twistor structure $(V,\DD^{\sankaku})$.
\index{unramifiedly good meromorphic prolongment}
\hfill\qed
\end{itemize}
\end{df}

\subsection{Meromorphic prolongment
 of the conjugate}
\label{subsection;10.5.24.20}

Let $\sigma:\proj^1\lrarr\proj^1$
be the anti-holomorphic involution
given by 
$\sigma([z_0:z_{\infty}])=
 [-\zbar_{\infty}:\zbar_0]$.
If we regard 
$\proj^1=\cnum_{\lambda}\cup\cnum_{\mu}$
by $\lambda=z_0/z_{\infty}$ and 
$\mu=z_{\infty}/z_0$,
we have the induced map
$\sigma:\cnum_{\lambda}\lrarr\cnum_{\mu}$
given by $\sigma(\lambda)=-\lambdabar$.
We have the naturally induced maps such as
$\nbigx\lrarr\nbigx^{\dagger}$ and
$\nbigx^{\shikaku}\lrarr\nbigx^{\dagger\shikaku}$,
which are also denoted by $\sigma$.

Let $(V,\DD^{\sankaku}_V)$
be a variation of twistor structure
on $\proj^1\times(X\setminus D)$.
Let $\DD_{\sigma^{\ast}V_{\infty}}$
(resp. $\DD^{\dagger}_{\sigma^{\ast}V_0}$)
denote the associated family of 
flat $\lambda$-connections
(resp. $\mu$-connections)
on $\sigma^{\ast}V_{\infty}$
(resp. $\sigma^{\ast}V_{0}$).
Let $(\Vtilde_{\infty},\DD^{\dagger})$
be an unramifiedly good meromorphic prolongment of 
$(V_{\infty},\DD^{\dagger})$
on $(\nbigx^{\dagger}_{\sigma(\nbigk)},
 \nbigd_{\sigma(\nbigk)}^{\dagger})$.
Let $\sigma^{\ast}\Vtilde_{\infty}$ be
the sheaf on $\nbigx_{\nbigk}$,
given by
$\sigma^{\ast}(\Vtilde_{\infty})(\nbigu):=
 \Vtilde_{\infty}(\sigma(\nbigu))$
for any open subset $\nbigu$ of $\nbigx_{\nbigk}$.
We have the natural $\nbigo_{\nbigx}(\ast\nbigd)$-module
structure on $\sigma^{\ast}(\Vtilde_{\infty})(\nbigu)$
given by
$f\cdot\sigma^{\ast}(s):=
 \sigma^{\ast}\bigl(
 \overline{\sigma^{\ast}(f)}\cdot s\bigr)$.
The family of $\lambda$-flat connections
$\DD_{\sigma^{\ast}V_{\infty}}$
naturally gives a family of meromorphic
flat $\lambda$-connections 
on $\sigma^{\ast}\Vtilde_{\infty}$.
Although the following lemma is clear,
we remark the signature
(Lemma \ref{lem;08.1.29.1}).
\begin{lem}
\mbox{{}}
\begin{itemize}
\item
 Let 
$(\sigma^{\ast}\Vtilde_{\infty},
 \DD_{\sigma^{\ast}V_{\infty}})$ be as above.
Then, it gives an unramifiedly good 
meromorphic prolongment
of $(\sigma^{\ast}V_{\infty},
 \DD_{\sigma^{\ast}V_{\infty}})$.
For each $R\in D$,
we have
\[
 \Irr(\sigma^{\ast}V_{\infty},R)
=\bigl\{
 -\gminiabar\,\big|\,
 \gminia\in\Irr(V_{\infty},R)
 \bigr\}.
\]
\item
Similarly, 
let $(\Vtilde_0,\DD)$ be
an unramifiedly good meromorphic prolongment
of $(V_0,\DD)$.
Then, 
$\bigl(\sigma^{\ast}\Vtilde_0,
 \DD_{\sigma^{\ast}V_0}^{\dagger}\bigr)$
 on
$(\nbigx^{\dagger},\nbigd^{\dagger})$
gives an unramifiedly good meromorphic
prolongment of 
$(\sigma^{\ast}V_0,
 \DD^{\dagger}_{\sigma^{\ast}V_0})$.
For each $R\in D$, we have
\[
 \Irr(\sigma^{\ast}V_0,R)
=\bigl\{-\gminiabar\,\big|\,
 \gminia\in\Irr(V_0,R)
 \bigr\}.
\]
\item
If $(\Vtilde_0,\Vtilde_{\infty})$ is 
an unramifiedly good meromorphic prolongment
of $(V,\DD^{\sankaku})$,
then $(\sigma^{\ast}\Vtilde_0,
 \sigma^{\ast}\Vtilde_{\infty})$
is an unramifiedly good meromorphic prolongment
of $\sigma^{\ast}(V,\DD^{\sankaku})$.
\hfill\qed
\end{itemize}
\end{lem}

\subsection{Meromorphic prolongment of the pairing}

Let $\nbigs:(V,\DD^{\sankaku})
 \otimes\sigma^{\ast}(V,\DD^{\sankaku})
 \lrarr \Tate(0)$
be a pairing of variation of 
twistor structure of weight $0$.
It consists of morphisms
$\nbigs_{0}:
 V_0\otimes\sigma^{\ast}V_{\infty}\lrarr
 \nbigo_{\nbigx\setminus\nbigd}$
and 
$\nbigs_{\infty}:
 V_{\infty}\otimes
 \sigma^{\ast}V_{0}\lrarr 
 \nbigo_{\nbigx^{\dagger}\setminus\nbigd^{\dagger}}$
which are compatible with
(i) the families of $\lambda$-connections
 or $\mu$-connections,
(ii) the gluing
on $\nbigx-W=\nbigx^{\dagger}-W^{\dagger}$.

Let $(\Vtilde_0,\Vtilde_{\infty})$ be an unramifiedly good
meromorphic prolongment of $(V,\DD)$.
For simplicity, we assume that 
$\Vtilde_0$ and $\Vtilde_{\infty}$ are
given on $\nbigx_{\nbigk}$
and $\nbigx^{\dagger}_{\sigma(\nbigk)}$.

\begin{df}
A pairing 
$\nbigstilde_0:\Vtilde_{0}\otimes
 \sigma^{\ast} \Vtilde_{\infty}
\lrarr \nbigo_{\nbigx_{\nbigk}}(\ast\nbigd_{\nbigk})$
is called a meromorphic prolongment of $\nbigs_0$
if $\nbigstilde_{0|\nbigx_{\nbigk}\setminus\nbigd_{\nbigk}}
=\nbigs_{0|\nbigx_{\nbigk}\setminus\nbigd_{\nbigk}}$.
Similarly, a pairing
$\nbigstilde_{\infty}:\Vtilde_{\infty}\otimes
 \sigma^{\ast}\Vtilde_{0}\lrarr
 \nbigo_{\nbigx_{\sigma(\nbigk)}^{\dagger}}
 (\ast\nbigd_{\sigma(\nbigk)}^{\dagger})$
is called a meromorphic prolongment of
$\nbigs_{\infty}$,
if $\nbigstilde_{\infty|
 \nbigx^{\dagger}_{\sigma(\nbigk)}
\setminus\nbigd^{\dagger}_{\sigma(\nbigk)}}
=\nbigs_{\infty|
 \nbigx^{\dagger}_{\sigma(\nbigk)}
\setminus\nbigd^{\dagger}_{\sigma(\nbigk)}}$
\hfill\qed
\end{df}
\index{meromorphic prolongment}

The following lemma is clear.
\begin{lem}
\label{lem;07.12.19.10}
Let $\nbigs:V\otimes\sigma^{\ast}V\lrarr\Tate(0)$
be a pairing of weight $0$.
\begin{itemize}
\item
A meromorphic prolongment of $\nbigs_0$
is unique, if it exists.
Similarly, a meromorphic prolongment of $\nbigs_{\infty}$
is unique, if it exists.
\item
$\nbigs_0$ has a meromorphic prolongment,
if and only if $\nbigs_{\infty}$ has a meromorphic
prolongment.
\item
$\nbigs_0$ has a meromorphic prolongment,
if and only if 
the induced morphism
$V_{0|\nbigx_{\nbigk}\setminus\nbigd_{\nbigk}}^{\lor}
\simeq\sigma^{\ast}V_{\infty|
 \nbigx^{\dagger}_{\sigma(\nbigk)}
\setminus\nbigd^{\dagger}_{\sigma(\nbigk)}}$
is extended to
$\Vtilde_0^{\lor}\simeq
 \sigma^{\ast}(\Vtilde_{\infty})$.
\hfill\qed
\end{itemize}
\end{lem}

\begin{df}
\mbox{{}}\label{df;08.9.14.51}
Let $(V,\DD^{\sankaku})$ be a variation of
twistor structure
with a pairing $\nbigs$ of weight $0$.
Let $(\Vtilde_0,\Vtilde_{\infty})$ be 
an unramifiedly good meromorphic prolongment
of $(V,\DD^{\sankaku})$.
We say that
$(\Vtilde_0,\Vtilde_{\infty})$ is
an unramifiedly good meromorphic prolongment
of $(V,\DD^{\sankaku},\nbigs)$,
if $\nbigs_0$ has the meromorphic prolongment.
\hfill\qed
\end{df}
\index{unramifiedly good meromorphic prolongment}

\subsection{Reduction of
 meromorphic variation of twistor structure}
\label{subsection;08.9.14.52}

Let $(V,\DD^{\sankaku})$ be a variation of
twistor structure on $\proj^1\times(X\setminus D)$
with a pairing $\nbigs$ of weight $0$.
Let $(\Vtilde_0,\Vtilde_{\infty})$ be 
an unramifiedly good meromorphic prolongment of
$(V,\DD^{\sankaku},\nbigs)$.
For simplicity,
we assume that
$\Vtilde_0$ and $\Vtilde_{\infty}$
are given on $\nbigx_{\nbigk}$
and $\nbigx^{\dagger}_{\sigma(\nbigk)}$,
respectively.

\subsubsection{Full reduction}

Let $R\in D$.
By taking Gr with respect to full Stokes filtrations
on a small neighbourhood $\nbigx_R$
of $p_{\lambda}^{-1}(R)\cap \nbigd_{\nbigk}$,
we obtain a graded meromorphic family of
$\lambda$-flat bundles
\[
 \Gr^{\nbigftilde}(\Vtilde_{0,R})
=\bigoplus_{\gminia\in\Irr(\Vtilde_0,R)}
 \Gr^{\nbigftilde}_{\gminia}(\Vtilde_{0,R})
\]
on $\nbigx_R$.
Similarly,
we obtain a graded meromorphic family of
$\mu$-flat bundles on 
a small neighbourhood 
$\nbigx^{\dagger}_R$ 
of $p_{\mu}^{-1}(R)\cap
 \nbigd^{\dagger}_{\sigma(\nbigk)}$:
\[
 \Gr^{\nbigftilde}(\Vtilde_{\infty,R})
=\bigoplus_{\gminib\in\Irr(\Vtilde_{\infty},R)}
 \Gr^{\nbigftilde}_{\gminib}(\Vtilde_{\infty,R})
\]
We may assume that
$\nbigx_R$ and $\nbigx^{\dagger}_R$
are of the form $\nbigk\times X_R$
and $\sigma(\nbigk)\times X_{R}^{\dagger}$.
We set $D_R:=X_R\cap D$.
Because of the coincidence of the Stokes filtrations,
we have a natural isomorphism
for each $\gminia\in\Irr(\Vtilde_0,R)$:
\[
 \Gr^{\nbigftilde}_{\gminia}(\Vtilde_0)
 _{|(\nbigk\cap\sigma^{\ast}(\nbigk))\times
 (X_R\setminus D_R)}
\simeq
 \Gr^{\nbigftilde}_{\gminiabar}(\Vtilde_{\infty})
_{|(\nbigk\cap\sigma^{\ast}(\nbigk))\times
 (X_R\setminus D_R)}
\]
By gluing, we obtain 
a variation of twistor structure
on $\proj^1\times (X_R\setminus D_R)$,
which is denoted by
$\Gr^{\nbigftilde}_{\gminia}
 (V_R,\DD_R^{\sankaku})$.
It is equipped with
an induced unramifiedly 
good meromorphic prolongment
$\bigl(\Gr^{\nbigftilde}_{\gminia}(\Vtilde_{0,R}),
 \Gr^{\nbigftilde}_{\gminiabar}(\Vtilde_{\infty,R})\bigr)$.

From the isomorphism
$\Vtilde_0^{\lor}\simeq
 \sigma^{\ast}(\Vtilde_{\infty})$
on $\nbigx_{\nbigk}$ induced by $\nbigs$,
we obtain an isomorphism
$\Gr^{\nbigftilde}_{-\gminia}(\Vtilde_{0,R}^{\lor})
\simeq
 \sigma^{\ast}
 \Gr^{\nbigftilde}_{\gminiabar}(\Vtilde_{\infty,R})$
on $\nbigx_{R,\nbigk}$.
Namely, 
we have an induced pairing:
\[
 \Gr_{\gminia}^{\nbigftilde}(\nbigs_{0,R}):
 \Gr^{\nbigftilde}_{\gminia}(\Vtilde_{0,R})
\otimes
 \sigma^{\ast}
  \Gr^{\nbigftilde}_{\gminiabar}(\Vtilde_{\infty,R})
\lrarr
 \nbigo_{\nbigx_{R,\nbigk}}(\ast\nbigd_{R,\nbigk}) 
\]
Similarly,
we have an induced pairing:
\[
 \Gr^{\nbigftilde}_{\gminiabar}(\nbigs_{\infty,R}):
 \Gr^{\nbigftilde}_{\gminiabar}(\Vtilde_{\infty,R})
\otimes
 \sigma^{\ast}
  \Gr^{\nbigftilde}_{\gminia}(\Vtilde_{0,R})
\lrarr
 \nbigo_{\nbigx^{\dagger}_{R,\nbigk}}
 (\ast\nbigd^{\dagger}_{R,\nbigk})
\]
It is easy to observe that 
their restrictions
to $(\nbigk\cap\sigma(\nbigk))\times
 (X_R\setminus D_R)$
are the same.
Hence, we have an induced symmetric pairing
$\Gr^{\nbigftilde}_{\gminia}(\nbigs_R)$
of $\Gr^{\nbigftilde}_{\gminia}(V_R,\DD_R^{\sankaku})$
equipped with a meromorphic prolongment.
The tuple is denoted by
$\Gr^{\nbigftilde}_{\gminia}
(V,\DD^{\sankaku},\nbigs)$.
We obtain 
\[
 \Gr^{\nbigftilde}
 (V,\DD^{\sankaku},\nbigs)
 :=\bigoplus_{\gminia\in \Irr(\Vtilde_0,R)}
 \Gr^{\nbigftilde}_{\gminia}
 (V,\DD^{\sankaku},\nbigs).
\]
It is called the full reduction of
$(V,\DD^{\sankaku},\nbigs)$,
and it is equipped 
with an unramifiedly
good meromorphic prolongment
$\bigl(\Gr^{\nbigftilde}(\Vtilde_{0,R}),
 \Gr^{\nbigftilde}(\Vtilde_{\infty,R})
 \bigr)$.

\subsubsection{Refinement}
\label{subsection;10.5.24.22}
Let us consider the case 
$X=\Delta^n$,
$D=\bigcup_{i=1}^{\ell}\{z_i=0\}$
and $R\in \bigcap_{i=1}^{\ell}\{z_i=0\}$.
We take an auxiliary sequence 
$\vecm(0),\ldots,\vecm(L),\vecm(L+1)=\veczero$
for the good set 
$\nbigi:=\Irr(\Vtilde_0,P)$.
As in the case of full reduction,
we obtain graded variation of twistor structure
\[
 \Gr^{\vecm(p)}(V_R,\DD_R^{\sankaku},\nbigs)
=\bigoplus_{\gminib\in 
 \etabar_{\vecm(p)}(\nbigi)}
 \Gr^{\vecm(p)}_{\gminib}
 (V_R,\DD_R^{\sankaku},\nbigs).
\]
It is naturally equipped with
an unramifiedly good meromorphic prolongment
$\bigl(\Gr^{\vecm(p)}(\Vtilde_{0,R}),
 \Gr^{\vecm(p)}(\Vtilde_{\infty,R})\bigr)$.
In particular,
$\bigl(\Gr^{\vecm(0)}(\Vtilde_{0,R}),
 \Gr^{\vecm(0)}(\Vtilde_{\infty,R})\bigr)$
is called the one step reduction.

\begin{rem}
\label{rem;07.12.19.25}
Let us consider the case
$X:=\Delta^n$ and $D:=\{z_1=0\}$.
Assume that we are given
an unramifiedly good meromorphic prolongment 
$(\Vtilde_0,\Vtilde_{\infty})$ 
of $(V,\DD^{\sankaku})$
on $(\nbigx,\nbigx^{\dagger})$.
If $D$ is smooth,
the reductions $\Gr^{(j)}(V,\DD^{\sankaku},\nbigs)$
are equipped with unramifiedly good meromorphic 
prolongments
$\bigl(
 \Gr^{(j)}(\Vtilde_0),
 \Gr^{(j)}(\Vtilde_{\infty})
\bigr)$
given on $(\nbigx',\nbigx^{\prime\dagger})$.
(Namely, we do not have to consider
the restriction to a compact region
in $\cnum_{\lambda}$, in this case.)
\hfill\qed
\end{rem}

\subsubsection{Compatibility}

We give a remark on compatibility.
We assume that the coordinate system is
admissible for a good set $\nbigi$
(Remark \ref{rem;07.11.15.10}).
Let $k$ be determined by
$\vecm(0)\in 
 \seisuu_{<0}^k\times\veczero_{\ell-k}$.

Take $1\leq j\leq k$
and $R_1\in D_j\cap X_R$,
which is not a singular point of $D(\kbar)$.
Let $X_{R_1}$ be a small neighbourhood of $R_1$
in $X$.
We set $X_{R_1}^{\ast}:=X_{R_1}\setminus D$.
We put 
\[
 (V_1,\DD^{\sankaku}_1,\nbigs_1):=
 (V,\DD^{\sankaku},\nbigs)
 _{|\proj^1\times X_Q^{\ast}},
\quad
  (V_2,\DD^{\sankaku}_2,\nbigs_2):=
 \Gr^{\vecm(0)}
 (V_R,\DD_R^{\sankaku},\nbigs)
 _{|\proj^1\times X_Q^{\ast}}
\]
They are equipped with
the induced unramifiedly meromorphic 
prolongments.
\begin{lem}
\label{lem;08.9.17.105}
We have a natural isomorphism:
\[
 \Gr^{\nbigftilde}
 \bigl(V_1,\DD^{\sankaku}_1,\nbigs_1\bigr)
\simeq
 \Gr^{\nbigftilde}
 \bigl(V_2,\DD^{\sankaku}_2,\nbigs_2\bigr)
\]
\end{lem}
\pf
It follows from Corollary \ref{cor;10.5.19.20}.
\hfill\qed

\part{Prolongation of Wild Harmonic Bundle}
\label{part;08.9.28.161}

\chapter[Prolongments $\nbigp\nbigelambda$]
{Prolongments $\nbigp\nbigelambda$ for
 Unramifiedly Good Wild Harmonic Bundles}
\label{section;07.6.2.1}
We start to study wild harmonic bundle.
In Section \ref{subsection;07.12.4.1},
we state the definition of wild harmonic bundles
and some related conditions for Higgs fields.

In Section \ref{subsection;07.10.7.1},
we state some estimates
related with the Higgs field
of an unramifiedly good wild harmonic bundle
(the wild version of Simpson's main estimate),
which will be proved in Section 
\ref{subsection;08.9.13.1}.
These estimates are the most foundational.

In Section \ref{subsection;08.9.14.40},
we consider the sheaves of
holomorphic sections whose norms are 
of polynomial orders,
and we show that they form 
a good filtered $\lambda$-flat bundle
(Theorem \ref{thm;07.12.2.55}).
We also obtain a characterization
of the Stokes filtrations
in terms of the growth order of the norms
of flat sections
(Proposition \ref{prop;08.9.14.45}).

In Section \ref{subsection;08.9.14.41},
we study the comparison of
the irregular decompositions
for $(\nbigp\nbige^0,\DD^0)$
and $(\nbigp\nbigelambda,\DDlambda)$.
We note that the family version will be studied in 
Section \ref{subsection;08.9.14.13}.
In the proof,
we give an estimate for the connection form
of the unitary connection associated to
$(\nbigelambda,h)$
(Lemma \ref{lem;07.6.2.25}).
It will also be useful for other purposes.

We would like to compare
the deformations caused by variation of irregular values
(Section \ref{subsection;10.5.17.51})
and by modification of the hermitian metrics.
It will be achieved in Proposition 
\ref{prop;07.10.21.100}.
We make a preparation in Section 
\ref{subsection;07.10.7.21}.

In Section \ref{subsection;08.9.14.46},
we give a criterion for a holomorphic section
to be bounded with respect to $h$.

\section{Definition of wild harmonic bundle}
\label{subsection;07.12.4.1}
\subsection{Local condition for Higgs fields}

Let $(E,\delbar_E,\theta)$ be a Higgs bundle
on $X-D$,
where $X$ is a complex manifold,
and $D$ is a normal crossing divisor of $X$.
We would like to state some conditions
for the Higgs field $\theta$.
First, let us consider the case 
$X=\Delta^n=\bigl\{\vecz=(z_1,\ldots,z_n)\,\big|\,|z_i|<1\bigr\}$,
$D_i=\{z_i=0\}$
and $D=\bigcup_{i=1}^{\ell}D_i$.
In that case, we have the expression:
\[
 \theta=\sum_{j=1}^{\ell} F_j\, \frac{dz_j}{z_j}
+\sum_{j=\ell+1}^{n}G_j\, dz_j
\]
We have the characteristic polynomials
\[
\det\bigl(T-F_j(\vecz)\bigr)
=\sum A_{j,k}(\vecz)\, T^k,
\quad
\det\bigl(T-G_j(\vecz)\bigr)=
 \sum B_{j,k}(\vecz)\, T^k.
\]
The coefficients $A_{j,k}$ and $B_{j,k}$
are holomorphic on $X-D$.

\begin{df}
\label{df;07.12.17.15}
We say that $\theta$ is tame,
if $A_{j,k}$ and $B_{j,k}$ are holomorphic on $X$
for any $k$,
and moreover,
if the restriction of $A_{j,k}$ to $D_j$ are constant
for any $j=1,\ldots,\ell$
and any $k$.
\hfill\qed
\end{df}
\index{tame Higgs field}
\begin{rem}
If $\theta$ comes from a tame harmonic bundle,
$\theta$ is tame in the above sense.
We do not have to assume that
the $A_{j,k|D_j}$ are constant
for the definition of tame harmonic bundle.
It is automatically satisfied.
\hfill\qed
\end{rem}

Let $\nbiga$ be a $\rnum$-vector subspace of $\cnum$.

\begin{df}
\label{df;07.11.25.30}
We say that $\theta$ is $\nbiga$-tame,
if $\theta$ is tame and the roots of the polynomials
$\det(T-F_j(\vecz))_{|D_j}$ are contained in
$\nbiga$, for any $j$.
\hfill\qed
\end{df}
\index{$\nbiga$-tame Higgs field}

\begin{df}
\mbox{{}}\label{df;08.9.28.150}
\begin{itemize}
\item
We say that $\theta$ is 
strongly unramifiedly ($\nbiga$-)good
on $(X,D)$,
if we have the good set of irregular values
$\Irr(\theta)\subset M(X,D)\big/H(X)$
and the decomposition
\begin{equation}
 \label{eq;10.5.24.1}
 (E,\theta)=\bigoplus_{\gminia\in\Irr(\theta)}
 (E_{\gminia},\theta_{\gminia}),
\end{equation}
such that
$\theta_{\gminia}-d\gminia\,\id_{E_{\gminia}}$
are ($\nbiga$-)tame.
\index{strongly unramifiedly good Higgs field}
\index{strongly unramifiedly good $\nbiga$-good Higgs field}
\item
We say that $\theta$ is 
strongly ($\nbiga$-)good
on $(X,D)$,
if $\varphi_e^{\ast}(\theta)$ is unramifiedly good
for some $e\in\seisuu_{>0}$,
where $\varphi_e$ is the covering
given by
$ \varphi_e(z_1,\ldots,z_n)=
 (z_1^{e},\ldots,z_{\ell}^{e},z_{\ell+1},\ldots,z_n)$.
\index{strongly good Higgs field}
\index{strongly $\nbiga$-good Higgs field}
\hfill\qed
\end{itemize}
\end{df}

The condition is independent of the choice of
a coordinate system.
The adjective ``strongly'' means
the existence of the global decomposition 
(\ref{eq;10.5.24.1}).
(Compare with Definition \ref{df;07.12.28.1} below.)
But, we will often omit ``strongly'',
if there is no risk of confusion.

\subsection{Global condition for Higgs fields}

Definition \ref{df;08.9.28.150}
can be globalized easily.

\begin{df}
\label{df;07.12.28.1}
\mbox{{}}
Let $X$ be a general complex manifold.
Let $D$ be a normal crossing hypersurface of $X$,
and let $(E,\theta)$ be a Higgs bundle on $X-D$.
\begin{itemize}
\item
$\theta$ is called
(unramifiedly, $\nbiga$-)good at $P\in D$,
if it is strongly (unramifiedly, $\nbiga$-)good 
on a coordinate neighbourhood of $P$.
\index{unramifiedly good Higgs field}
\index{unramifiedly $\nbiga$-good Higgs field}
\item
$\theta$ is called
(unramifiedly, $\nbiga$-)good on $(X,D)$,
if it is (unramifiedly, $\nbiga$-)good 
at any point $P\in D$.
\index{good Higgs field}
\index{$\nbiga$-good Higgs field}
\hfill\qed
\end{itemize}
\end{df}

We also introduce the more general condition.

\begin{df}
\label{df;07.12.28.2}
Let $X$ be a complex analytic space,
and let $Z$ be a closed analytic subset of $X$
such that $X-Z$ is smooth.
Let $(E,\theta)$ be a Higgs bundle 
on $X-Z$.
The Higgs field $\theta$ is called ($\nbiga$-)wild
on $(X,Z)$,
if there exists a complex manifold $X'$
with a projective birational map
$\varphi:X'\lrarr X$
such that (i) $\varphi^{-1}(Z)$ is normal crossing,
(ii) $\varphi^{-1}\theta$ is ($\nbiga$-)good
on $(X',\varphi^{-1}(Z))$.
\index{wild Higgs field}
\index{$\nbiga$-wild Higgs field}
\hfill\qed
\end{df}

\subsection{Condition for harmonic bundles}

We introduce the notion of
(unramifiedly) good wild harmonic bundles,
which is one of the main subjects
in this monograph.

\begin{df}
Let $X$ be a complex manifold,
and let $D$ be a normal crossing hypersurface of $X$.
Let $\harmonicbundle$ be a harmonic bundle
on $X-D$.
\begin{itemize}
\item
 It is called $\nbiga$-tame on $(X,D)$,
 if $\theta$ is $\nbiga$-tame on $(X,D)$.
\index{$\nbiga$-tame harmonic bundle}
\item
 It is called (unramifiedly, $\nbiga$-)good 
 wild harmonic bundle on $(X,D)$,
 if $\theta$ is 
 (unramifiedly, $\nbiga$-)good on $(X,D)$.
\index{unramifiedly good wild harmonic bundle}
\index{unramifiedly $\nbiga$-good wild harmonic bundle}
\index{good wild harmonic bundle}
\index{$\nbiga$-good wild harmonic bundle}
\hfill\qed
\end{itemize}
\end{df}

In the case $X=\Delta^n$
and $D=\bigcup_{i=1}^{\ell}\{z_i=0\}$,
we will often implicitly assume that
$\theta$ is {\em strongly} (unramifiedly $\nbiga$-)good.
We will often say that
$\harmonicbundle$ is a (unramifiedly $\nbiga$-)good
wild harmonic bundle on $X-D$
instead of $(X,D)$.
In our previous paper {\rm\cite{mochi2}},
$\sqrt{-1}\real$-tame harmonic bundle
is called tame pure imaginary harmonic bundle.
\index{$\sqrt{-1}\real$-tame harmonic bundle}

We introduce a more general notion.
\begin{df}
\label{df;08.1.7.1}
Let $Z$ be a closed analytic subset of $X$.
A harmonic bundle $\harmonicbundle$ on $X-Z$
is called ($\nbiga$-)wild on $(X,Z)$.
if $\theta$ is ($\nbiga$-)wild on $(X,Z)$.
\index{wild harmonic bundle}
\index{$\nbiga$-wild harmonic bundle}
(We will also say that
$(E,\delbar_E,\theta,h)$ is 
($\nbiga$-)wild on $(X,X-Z)$,
instead of $(X,Z)$.)
\hfill\qed
\end{df}

Analysis will be mainly done 
for (unramifiedly) good wild harmonic bundles.
In the curve case,
any wild harmonic bundle is good.
We remark that
even if $\harmonicbundle$ is 
a wild harmonic bundle on $(X,D)$,
where $D$ is normal crossing,
it is not necessarily good.
In Chapter \ref{section;10.6.1.1},
we will study when a harmonic bundle is wild.

\section[Simpson's main estimate]
 {Simpson's main estimate and acceptability
 of the associated bundles}
\label{subsection;07.10.7.1}
\subsection{Setting}
\label{subsection;07.11.17.1}

Let $X=\Delta^n$,
$D_i=\{z_i=0\}$
and $D=\bigcup_{i=1}^{\ell}D_i$.
Let $\harmonicbundle$ be 
an unramifiedly good wild harmonic bundle on $(X,D)$.
By shrinking $X$,
we assume to have the irregular decomposition:
\begin{equation}
\label{eq;07.7.18.10}
 (E,\theta)=
 \bigoplus_{\gminia\in\Irr(\theta)}
 (E_{\gminia},\theta_{\gminia})
\end{equation}
\index{subbundle $E_{\gminia}$}
Here,
$\theta_{\gminia}
-d\gminia\,\id_{E_{\gminia}}$
are tame.
We have the expression:
\[
 \theta=\sum_{j=1}^nf_j\, dz_j
\]
We put
$f_j^{\reg}:=
 f_j-\sum \del_j\gminia\, \pi_{\gminia}$
for $j=1,\ldots,n$.
\index{endomorphism $f_j^{\reg}$}
By the assumption,
$\det\bigl(T-z_j f^{\reg}_{j}\bigr)_{|D_j}$
($j=1,\ldots,\ell$)
are polynomials in a formal variable $T$
with constant coefficients.
Let $\Sp(\theta,j)\subset\cnum$ denote
the set of the solutions of 
$\det\bigl(T-z_j f^{\reg}_{j}\bigr)_{|D_j}=0$.
\index{set $\Sp(\theta,j)$}
For each $j=1,\ldots,\ell$,
we assume to have the decomposition
\begin{equation}
 \label{eq;07.7.18.11}
 (E,f_j^{\reg})=\bigoplus_{
 \alpha\in\Sp(\theta,j)}
 (E_{j,\alpha},f^{\reg}_{j,\alpha})
\end{equation}
such that the eigenvalues $\beta$ of
$(z_j\, f^{\reg}_{j,\alpha})_{|Q}$
satisfies $|\beta-\alpha|\leq C_0\, |z_j(Q)|^{\epsilon_0}$
for some $C_0,\epsilon_0>0$.
\index{subbundle $E_{j,\alpha}$}

We put
$X(R):=\bigl\{(z_1,\ldots,z_n)\in X\,\big|\,
 |z_i|<R
 \bigr\}$ for $R<1$
and $X^{\ast}(R):=(X-D)\cap X(R)$.

\subsection{Main estimate for the irregular part}
\label{subsection;07.11.17.20}

We take an auxiliary sequence
$\nbigm:=\bigl(
 \vecm(0),\vecm(1),\vecm(2),\ldots,
 \vecm(L) \bigr)\subset
 \seisuu_{\leq 0}^{\ell}-\{\veczero_{\ell}\}$ 
for the good set of the irregular values
$\Irr(\theta)$,
i.e.,
(i) $\vecm(0)=\min\bigl\{\ord(\gminia)\,\big|\,
 \gminia\in \Irr(\theta)\bigr\}$,
(ii) we have 
$\gminih(i)$ such that
 $\vecm(i+1)=\vecm(i)+\vecdelta_{\gminih(i)}$
 for each $i\leq L-1$
and $\vecm(L)+\vecdelta_{\gminih(L)}=\veczero_{\ell}$,
(iii) $\nbigt(\theta):=
 \bigl\{\ord(\gminia-\gminib)\,\big|\,
 \gminia,\gminib\in\Irr(\theta)
 \bigr\}
 \subset \nbigm$.
\index{set $\nbigt(\theta)$}

\vspace{.1in}

Let $\etabar_{\vecm}:\Irr(\theta)\lrarr M(X,D)$ 
be given as in (\ref{eq;08.9.13.2}).
\index{map $\etabar_{\vecm}$}
Let $\Irrbar(\theta,\vecm)$ denote the image of
$\Irr(\theta)$ via $\etabar_{\vecm}$
for $\vecm\in\nbigm$.
\index{set $\Irr(\theta,\vecm)$}
We may and will fix auxiliary
total orders $\leq$ on 
the finite sets $\Irr(\theta)$
and $\Irrbar(\theta,\vecm)$
for  any $\vecm\in \nbigm$
such that $\etabar_{\vecm}$ 
are order preserving.

For each $\vecm\in \nbigm$ and
$\gminib\in\Irrbar(\theta,\vecm)$,
we define
\[
 E^{\vecm}_{\gminib}:=
 \bigoplus_{\etabar_{\vecm}(\gminia)=\gminib}
 E_{\gminia},
\quad
 F^{\vecm}_{\gminib}(E):=
 \bigoplus_{\substack{\gminic\in\Irrbar(\theta,\vecm)\\
 \gminic\leq \gminib}}
 E_{\gminic}^{\vecm},
\quad
 F^{\vecm}_{<\gminib}(E):=
 \bigoplus_{\substack{\gminic\in\Irrbar(\theta,\vecm)\\
 \gminic< \gminib}}
 E_{\gminic}^{\vecm}.
\]
\index{subbundle $E^{\vecm}_{\gminib}$}
Let $E^{\vecm\prime}_{\gminib}$
denote the orthogonal complement of
$F_{<\gminib}^{\vecm}(E)$ 
in $F_{\gminib}^{\vecm}(E)$.
\index{subbundle $E^{\vecm\prime}_{\gminib}$}
Let $\pi^{\vecm}_{\gminib}$
denote the projection of $E$
onto $E_{\gminib}^{\vecm}$
with respect to the decomposition
$E=\bigoplus_{\gminib\in\Irrbar(\theta,\vecm)}
 E_{\gminib}^{\vecm}$.
\index{projection $\pi^{\vecm}_{\gminib}$}
Let $\pi_{\gminib}^{\vecm\prime}$
denote the orthogonal projection onto
$E_{\gminib}^{\vecm\prime}$.
\index{projection $\pi_{\gminib}^{\vecm\prime}$}
We put 
$\nbigr^{\vecm}_{\gminia}:=
 \pi^{\vecm}_{\gminia}
-\pi^{\vecm\prime}_{\gminia}$.

We will prove the following theorem
in Sections
\ref{subsection;07.10.4.2}--\ref{subsection;07.7.18.2}.
\begin{thm}
\label{thm;07.10.4.1}
There exist positive constants
$R_1$, $\epsilon_1$ and $A_1$
such that the following holds on $X^{\ast}(R_1)$
for any $\gminib\in\Irrbar\bigl(\theta,\vecm\bigr)$:
\[
 \bigl|\nbigr^{\vecm}_{\gminib}\bigr|_h\leq
A_1\,\exp\bigl(-\epsilon_1|\vecz^{\vecm}|\bigr)
\]
 In particular, the decomposition
$E=\bigoplus_{\gminib\in\Irrbar(\theta,\vecm)}
 E^{\vecm}_{\gminib}$
 is $\exp\bigl(
 -\epsilon_1|\vecz^{\vecm}|\bigr)$-asymptotically
 orthogonal,
in the sense that there exists $A_1'>0$
such that the following holds
for any $u_i\in E_{\gminib_i|Q}$
with $\gminib_1\neq\gminib_2$:
\[
 |h(u_1,u_2)|
\leq 
 A_1'\,|u_1|_h\,|u_2|_h\,
 \exp\bigl(
 -\epsilon_1|\vecz^{\vecm}(Q)|
 \bigr) 
\]

The constants $A_1$, $A_1'$,
$\epsilon_1$
and $R_1$ may depend only on 
$\rank(E)$,
$C_0$, $\epsilon_0$, $\Irr(\theta)$ and
$\Sp(\theta,j)$ $(j=1,\ldots,\ell)$ in 
Section {\rm\ref{subsection;07.11.17.1}}.
\end{thm}

\begin{cor}
\label{cor;07.11.22.1}
We have the estimate
$\bigl|
 \pi^{\vecm}_{\gminia}
-\pi^{\vecm\dagger}_{\gminia}\bigr|_h
\leq 2A_1\,
 \exp\bigl(-\epsilon_1|\vecz^{\vecm}|\bigr)$,
where $\pi^{\vecm\dagger}_{\gminia}$ denotes
the adjoint of $\pi^{\vecm}_{\gminia}$
with respect to $h$.
\hfill\qed
\end{cor}

\begin{rem}
The main part of Theorem {\rm\ref{thm;07.10.4.1}}
is the claim for $\vecm\in\nbigt(\theta)$,
but it is convenient to take
an auxiliary sequence $\nbigm$
for inductive arguments
in both the proof and the use.
\hfill\qed
\end{rem}

\subsection{Main estimate for the regular part}
\label{subsection;07.11.20.11}

We take an order $\leq$ on $\Sp(\theta,j)$
for each $j$.
We put $F_{j,\alpha}(E):=
 \bigoplus_{\beta\leq\alpha}E_{j,\alpha}$.
Let $E_{j,\alpha}'$ denote the orthogonal complement
of $F_{j,<\alpha}$ in $F_{j,\alpha}$.
\index{subbundle $E_{j,\alpha}'$}
Let $\pi_{j,\alpha}$ denote
the projection of $E$ onto $E_{j,\alpha}$
with respect to the decomposition
$E=\bigoplus E_{j,\alpha}$.
\index{projection $\pi_{j,\alpha}$}
Let $\pi_{j,\alpha}'$ denote the orthogonal projection
onto $E'_{j,\alpha}$.
\index{projection $\pi_{j,\alpha}'$}
For any $j=1,\ldots,\ell$,
we set
\[
 q_j:=f_j^{\reg}
-\sum_{\alpha\in\Sp(\theta,j)}
 \frac{\alpha}{z_j}\,\pi_{j,\alpha}'
\]
We put
$\nbigr^{\reg}_{j,\alpha}:=
 \pi_{j,\alpha}-\pi'_{j,\alpha}$.
We will show the following theorem in Sections
\ref{subsection;07.10.4.4}--\ref{subsection;07.10.4.5}.
\begin{thm}
\label{thm;07.10.4.3}
We have the following estimates:
\begin{description}
\item[$P(\reg)$]
We have 
$\bigl|z_j\, f_j^{\reg}\bigr|_h\leq A_{2}$ $(j=1,\ldots,\ell)$
and $\bigl|f_j^{\reg}\bigr|_h\leq A_{2}$ 
 $(j=\ell+1,\ldots,n)$ on $X^{\ast}(R_{2})$.
\item[$Q(\reg)$]
$|q_j|_h\leq A_{2}\,
 |z_j|^{-1}\bigl(-\log|z_j|\bigr)^{-1}$
for $j=1,\ldots,\ell$
on $X^{\ast}(R_{2})$.
\item[$R(\reg)$]
$\bigl|
 \nbigr^{\reg}_{j,\alpha}\bigr|_h\leq
 A_{2}\,|z_j|^{\epsilon_2}$
on $X^{\ast}(R_{2})$ for $j=1,\ldots,\ell$.
In particular,
the decomposition
$E=\bigoplus E_{j,\alpha}$
is $O(|z_j|^{\epsilon_2})$-asymptotically orthogonal.
\end{description}
The positive constants $A_{2}$, $R_2$ and $\epsilon_2$
may depend only on 
$\rank(E)$,
$C_0$, $\epsilon_0$, $\Irr(\theta)$
and $\Sp(\theta,j)$ 
in Section {\rm\ref{subsection;07.11.17.1}}.
\end{thm}

\begin{cor}
\label{cor;07.11.22.5}
We have the estimate
$\bigl|
 \pi_{j,\alpha}-\pi_{j,\alpha}^{\dagger}
 \bigr|_h
 \leq 2A_2\,|z_j|^{\epsilon_2}$,
where $\pi_{j,\alpha}^{\dagger}$ denotes
the adjoint of $\pi_{j,\alpha}$ with respect to $h$.
\hfill\qed
\end{cor}

For $j=1,\ldots,\ell$,
we consider the following:
\[
 f^{\nil}_j:=
 f^{\reg}_j-\sum_{\alpha\in\Sp(\theta,j)}
 \frac{\alpha}{z_j}\,\pi_{j,\alpha}
\]
\index{endomorphism $f_j^{\nil}$}
\begin{cor}
We have 
$\bigl|f^{\nil}_j\bigr|_h
\leq A_2'\,|z_j|^{-1}\bigl(-\log|z_j|\bigr)^{-1}$.
\hfill\qed
\end{cor}

For any 
$\vecalpha,\vecbeta\in
 \Sp(\theta):=\prod_{j=1}^{\ell}\Sp(\theta,j)$,
we set 
$ \Diff(\vecalpha,\vecbeta):=
 \bigl\{j\,\big|\,\alpha_j\neq\beta_j\bigr\}$
and 
\[
\nbigq_{\epsilon}(\vecalpha,\vecbeta):=
 \prod_{j\in \Diff(\vecalpha,\vecbeta)}
 |z_j|^{\epsilon}.
\]
\index{set $\Diff(\vecalpha,\vecbeta)$}
\index{set $\Sp(\theta)$}
\index{function 
 $\nbigq_{\epsilon}(\vecalpha,\vecbeta)$}
We put 
$E_{\gminia,\vecalpha}:=
 E_{\gminia}\cap \bigcap_{j=1}^{\ell}E_{j,\alpha_j}$.
\index{subbundle $E_{\gminia,\vecalpha}$}
We have the following immediate corollary of
Theorem \ref{thm;07.10.4.1}
and Theorem \ref{thm;07.10.4.3}.

\begin{cor}
\label{cor;07.7.18.11}
In the case
$(\gminia,\vecalpha)
\neq (\gminib,\vecbeta)$,
the subbundles
$E_{\gminia,\vecalpha}$ and
$E_{\gminib,\vecbeta}$ are
$\exp\bigl(
 -\epsilon|\vecz^{\ord(\gminia-\gminib)}|\bigr)
 \,
 \nbigq_{\epsilon}(\vecalpha,\vecbeta)$-asymptotically
orthogonal for some $\epsilon>0$.
\hfill\qed
\end{cor}

\subsection{Complementary estimates
 for the Higgs field}

We give some refinements,
which immediately follow from the theorems.
The proof will be given in Section 
\ref{subsection;07.11.22.12}.
Any constants and any estimates
may depend only on 
$\rank(E)$,
$C_0$, $\epsilon_0$, $\Irr(\theta)$
and $\Sp(\theta,j)$ 
in Section {\rm\ref{subsection;07.11.17.1}}.
We have the decomposition:
\begin{equation}
 \label{eq;07.11.22.23}
 \End(E)=
 \bigoplus_{\gminia,\gminia'\in\Irr(\theta)}
 \bigoplus_{\vecalpha,\vecalpha'\in\Sp(\theta)}
 \Hom\bigl(E_{\gminia,\vecalpha},E_{\gminia',\vecalpha'}\bigr)
\end{equation}
For a section $F$ of $\End(E)$,
we have the corresponding decomposition:
\begin{equation}
 \label{eq;08.12.8.1}
 F= \sum_{\gminia,\gminia'\in\Irr(\theta)}
 \sum_{\vecalpha,\vecalpha'\in\Sp(\theta)}
 F_{(\gminia,\vecalpha),(\gminia',\vecalpha')}
\end{equation}
Let $\nbigq_{\epsilon}(\vecalpha,\vecalpha')$
be as in Section \ref{subsection;07.11.20.11}.

\begin{prop}
\label{prop;07.11.22.11}
We have the following estimates
on $X^{\ast}(R_3)$:
\[
\Bigl|
\bigl(
 \pi^{\vecm(p)\dagger}_{\gminib}
-\pi^{\vecm(p)}_{\gminib}
\bigr)_{(\gminia,\vecalpha),(\gminia',\vecalpha')}
\Bigr|_h
\leq
 A_{3}\,
 \exp\bigl(-\epsilon_3|\vecz^{\vecm(p)}|
-\epsilon_3|\vecz^{\ord(\gminia-\gminia')}|\bigr)
\,
 \nbigq_{\epsilon_3}(\vecalpha,\vecalpha')
\]
For $j=1,\ldots,\ell$,
\[
\Bigl|
 \bigl(
 \pi^{\dagger}_{j,\gamma}
-\pi_{j,\gamma}
\bigr)_{(\gminia,\vecalpha),(\gminia',\vecalpha')}
\Bigr|_h
\leq A_{3}
 |z_j|^{\epsilon_3}\,
 \exp\bigl(-\epsilon_3|\vecz^{\ord(\gminia-\gminia')}|\bigr)
\,
 \nbigq_{\epsilon_3}(\vecalpha,\vecalpha')
\]
\[
\Bigl|
 \bigl(
 f^{\nil\dagger}_j
\bigr)_{(\gminia,\vecalpha),(\gminia',\vecalpha')}
\Bigr|_h
\leq A_3\,
 |z_j|^{-1}(-\log|z_j|)^{-1}\,
 \exp\bigl(-\epsilon_3|\vecz^{\ord(\gminia-\gminia')}|
 \bigr)\,
 \nbigq_{\epsilon_3}(\vecalpha,\vecalpha')
\]
For $j=\ell+1,\ldots,n$,
\[
\Bigl|
 \bigl(
 f^{\reg\dagger}_j
\bigr)_{(\gminia,\vecalpha),(\gminia',\vecalpha')}
\Bigr|_h
\leq
 A_3\,
 \exp\bigl(-\epsilon_3|\vecz^{\ord(\gminia-\gminia')}|
 \bigr)\,
 \nbigq_{\epsilon_3}(\vecalpha,\vecalpha').
\]
\end{prop}

\subsection{Estimate of the curvature}

Let $g_{\poin}$ denote the Poincar\'e metric
of $X-D$.
For any section $F$ of $\End(E)\otimes\Omega^p$,
we have the decomposition
$F=\sum F_{(\gminia,\vecalpha),(\gminia',\vecalpha')}$
corresponding to the decomposition (\ref{eq;07.11.22.23}).
We will prove the following proposition
in Section \ref{subsection;08.1.7.2}.
\begin{prop}
\label{prop;07.11.22.25}
We have the following estimates on $X^{\ast}(R_4)$
with respect to $h$ and $g_{\poin}$:
\[
\Bigl|
 \bigl[\theta,\theta^{\dagger}\bigr]_{
 (\gminia,\vecalpha),(\gminia',\vecalpha')}
\Bigr|_{h,g_{\poin}}
\leq A_4\,
 \exp\bigl(-\epsilon_4|\vecz^{\ord(\gminia-\gminia')}|\bigr)
 \,\nbigq_{\epsilon_4}(\vecalpha,\vecalpha')
\]
In particular, 
$\bigl| [\theta,\theta^{\dagger}]\bigr|_{h,g_{\poin}}\leq A_5$ 
on $X^{\ast}(R_4)$.
Here the constants $R_4,A_4,A_5$
may depend only on 
$\rank(E)$,
$C_0$, $\epsilon_0$, $\Irr(\theta)$
and $\Sp(\theta,j)$ 
in Section {\rm\ref{subsection;07.11.17.1}}.
\end{prop}

Let $d_{\lambda}'':=\delbar_E+\lambda\theta^{\dagger}$.
The holomorphic bundle
$(E,d_{\lambda}'')$ is denoted by $\nbigelambda$.
The curvature of the unitary connection
associated to $h$ and $d_{\lambda}''$ 
is denoted by $R(h,d_{\lambda}'')$.
Recall the relation
$R(h,d_{\lambda}'')=
 -(1+|\lambda|^2)\, [\theta,\theta^{\dagger}]$.
Hence, we obtain the following direct corollary
of Proposition \ref{prop;07.11.22.25}.
\begin{cor}
\label{cor;07.11.22.30}
We have the following estimates on $X^{\ast}(R_4)$
with respect to $h$ and $g_{\poin}$:
\[
\Bigl|
 R(h,d_{\lambda}'')_{
 (\gminia,\vecalpha),(\gminia',\vecalpha')}
\Bigr|_{h,g_{\poin}}
\leq (1+|\lambda|^2)\, A_4\,
 \exp\bigl(-\epsilon_4|\vecz^{\ord(\gminia-\gminia')}|\bigr)
 \,\nbigq_{\epsilon_4}(\vecalpha,\vecalpha')
\]
In particular, 
$\bigl|R(h,d''_{\lambda}) \bigr|_{h,g_{\poin}}\leq 
 (1+|\lambda|^2)\, A_5$ 
on $X^{\ast}(R_4)$.
Therefore,
$(\nbigelambda,h)$ is acceptable.
In particular,
$\bigl(
 \End(\nbigelambda),h\bigr)$
is also acceptable.
\hfill\qed
\end{cor}

\section{Proof of the estimates}
\label{subsection;08.9.13.1}
\subsection{Inductive statement for the irregular part}
\label{subsection;07.10.4.2}

In the following argument,
any constants and any estimates
may depend only on 
$\rank(E)$,
$C_0$, $\epsilon_0$, $\Irr(\theta)$
and $\Sp(\theta,j)$ 
in Section {\rm\ref{subsection;07.11.17.1}}.
Let $f$ and $g$ denote
some functions on $X^{\ast}$.
We say $f=O(g)$
if $|f|\leq A\cdot |g|$ holds on $X^{\ast}(R)$,
and we say $f\sim g$ 
if $A^{-1}\cdot |g|\leq |f|\leq A\cdot|g|$ holds
on $X^{\ast}(R)$,
for some positive constants $A,R>0$.

For $\gminia\in\Irr(\theta)$,
let $\zeta_{\vecm(i)}(\gminia)$ be as in 
(\ref{eq;07.11.6.4}).
It is also well defined for 
$\gminia\in \Irrbar(\theta,\vecm(i))$.
Note 
$\etabar_{\vecm(i)}(\gminia)
=\sum_{j\leq i}\zeta_{\vecm(j)}(\gminia)$
for $\gminia\in\Irr(\theta)$.
We define $f_j^{\vecm(i)}$
and $\mu_j^{\vecm(i)}$ as follows:
\begin{multline}
 f_{j}^{\vecm(i)}:=f_j-\sum_{\gminia\in\Irr(\theta)}
 \del_j\etabar_{\vecm(i-1)}(\gminia)\cdot \pi_{\gminia}
=f_j^{\vecm(i-1)}
-\sum_{\gminia\in\Irr(\theta)}\del_j\zeta_{\vecm(i-1)}(\gminia)
 \cdot\pi_{\gminia} \\
=f_j^{\vecm(i-1)}
-\sum_{\gminib\in \Irrbar(\theta,\vecm(i-1))}
 \del_j\zeta_{\vecm(i-1)}(\gminib)\cdot
 \pi^{\vecm(i-1)}_{\gminib}
\end{multline}
\begin{equation}
 \mu_j^{\vecm(i)}:=f_j^{\vecm(i)}
-\sum_{\gminia\in\Irr(\theta)}
 \del_j\zeta_{\vecm(i)}(\gminia)
 \cdot \pi^{\prime}_{\gminia}
=
 f_j^{\vecm(i)}-
 \sum_{\gminib\in\Irrbar(\theta,\vecm(i))}
\del_j\zeta_{\vecm(i)}(\gminib)
\cdot \pi^{\vecm(i)\prime}_{\gminib}
\end{equation}

For the proof of Theorem \ref{thm;07.10.4.1},
we show the following claims inductively on $i$.
\begin{description}
\item[$P(i,\irr)$]
 $|f_j^{\vecm(i')}|_h
=O\Bigl(
 \bigl|\vecz^{\vecm(i')-\vecdelta_j}\bigr|
 \Bigr)$
for $j=1,\ldots,\ell$ and for any $i'\leq i$.
\item[$Q(i,\irr)$]
 $\bigl|
 \mu_{\gminih(i')}^{\vecm(i')}\bigr|_h
=O\Bigl(\bigl|\vecz^{\vecm(i')}\bigr|\Bigr)$
 for any $i'\leq i$.
\item[$R(i,\irr)$]
 $\bigl|\nbigr^{\vecm(i')}_{\gminib}\bigr|_h
=O\Bigl(
  \exp\bigl(-\epsilon(3,i')\cdot|\vecz^{\vecm(i')}|\bigr)
 \Bigr)$
 for any $\gminib\in\Irrbar\bigl(\theta,\vecm(i')\bigr)$
 and for any $i'\leq i$.
 In particular, the decomposition
$E=\bigoplus_{\gminib\in\Irrbar(\theta,\vecm(i'))}
 E^{\vecm(i')}_{\gminib}$
 is $\exp\bigl(-\epsilon(3,i')
 \cdot|\vecz^{\vecm(i')}|\bigr)$-asymptotically
 orthogonal for any $i'\leq i$.
\end{description}

\begin{rem}
\label{rem;07.6.1.5}
By considering the pull back of $\harmonicbundle$
via a ramified covering
$\varphi:X\lrarr X$ given by
$\varphi(z_1,\ldots,z_n)
=(z_1^d,\ldots,z_{\ell}^d,z_{\ell+1},\ldots,z_n)$,
and by shrinking $X$,
we may and will assume that 
there exists a constant $C_{10}$,
which is independent of the choice of $Q\in X-D$,
such that 
\[
 \left|
 \beta-\frac{\del\gminia}{\del z_j}(Q)
-\frac{\alpha}{z_j(Q)}
 \right|\leq C_{10}
\]
for any eigenvalue $\beta$ of $f_{j,(\gminia,\alpha)|Q}$,
where $f_{j,(\gminia,\alpha)}$ denotes
the restriction of $f_j$ to $E_{j,\alpha}\cap E_{\gminia}$.
\hfill\qed
\end{rem}

\begin{rem}
By tensoring rank one meromorphic connection,
we may and will assume,
for any $\vecm\in\nbigm$ ,
there exists $\gminia\in\Irr(\theta)$
such that $\gminia_{\vecm}\neq 0$.
\hfill\qed
\end{rem}

\subsection{$P(i-1,\irr)+Q(i-1,\irr)+R(i-1,\irr)$
$\Longrightarrow$ $P(i,\irr)$}
\label{subsection;07.10.4.30}

Let $\Delta_j$
denote the Laplacian with respect to the variable $z_j$,
i.e., $\Delta_j:=-\del^2/\del z_j\del\zbar_j$.
Because of the commutativity
$\bigl[f_j,f_j^{\vecm(i)}\bigr]=0$,
we have the following inequality:
\begin{equation}
\label{eq;07.5.18.3}
 \Delta_j\log\bigl|f_j^{\vecm(i)}\bigr|_h^2
\leq
 -\frac{\bigl|[f_j^{\dagger},f_j^{\vecm(i)}]\bigr|_h^2}
 {\bigl|f_j^{\vecm(i)}\bigr|_h^2}
\end{equation}
(See Page 731 of \cite{s2}.)
Recall
$f_j=f_j^{\vecm(i)}
 +\sum\del_j\etabar_{\vecm(i-1)}(\gminia)
 \cdot\pi_{\gminia}$.
We have the following equality:
\[
 \sum_{\gminia}
 \del_j\etabar_{\vecm(i-1)}(\gminia)
 \cdot\pi_{\gminia}
=\sum_{\gminia}\sum_{i'< i}
 \del_j\zeta_{\vecm(i')}(\gminia)\cdot\pi_{\gminia}
=\sum_{i'< i}
 \sum_{\gminib\in\Irrbar(\theta,\vecm(i'))}
 \del_j\zeta_{\vecm(i')}(\gminib)
 \cdot\pi^{\vecm(i')}_{\gminib}
\]
For any section $F$ of $\End(E)$,
we have the decomposition
$F=\nbigc^{(i')}(F)+\nbigd^{(i')}(F)$
as follows:
\[
\nbigc^{(i')}(F)\in
 \!\!\!\!\!
 \bigoplus_{
 \substack{\gminib,\gminib'\in\Irrbar(\theta,\vecm(i'))\\
 \gminib\neq\gminib'}}
 \!\!\!\!\!
 \Hom\bigl(E^{\vecm(i')\prime}_{\gminib},
 E^{\vecm(i')\prime}_{\gminib'}\bigr),
\quad
\nbigd^{(i')}(F)\in
 \!\!\!
 \bigoplus_{\gminib\in\Irrbar(\theta,\vecm(i'))}
 \!\!\!
 \End\bigl(E^{\vecm(i')\prime}_{\gminib}\bigr)
\]
Then, we have 
$ \bigl[
 \pi_{\gminib}^{\vecm(i')\dagger},
 f_j^{\vecm(i)}
\bigr]
=\bigl[
 \nbigc^{(i')}(\pi_{\gminib}^{\vecm(i')\dagger}),
 f_j^{\vecm(i)}
 \bigr]
+\bigl[
 \nbigd^{(i')}(\pi^{\vecm(i')\dagger}_{\gminib}),
 \nbigc^{(i')}(f_j^{\vecm(i)})
 \bigr]$
for any $i'<i$.
Both the first and second terms are
$O\bigl(
 \exp(-\epsilon(4,i')|\vecz^{\vecm(i')}|)
 \bigr)\cdot |f_j^{\vecm(i)}|$
because of $R(i-1,\irr)$.
Note 
$\bigl|\del_j\zeta_{\vecm(i')}(\gminib)\bigr|=
O\bigl(|\vecz^{\vecm(i')}|\bigr)$
if the $j$-th component of
$\vecm(i')$ is $0$.
Hence, we have the estimates
\[
 |\del_j\zeta_{\vecm(i')}(\gminib)|
 \cdot
 \exp\bigl(
 -\epsilon(5,i')|\vecz^{\vecm(i')}|
 \bigr)
=O\Bigl(
 \exp\bigl(-\epsilon(6,i')|\vecz^{\vecm(i')}|\bigr)
 \Bigr).
\]
Therefore, we obtain the following estimate:
\[
 \Bigl|
 \bigl[\bigl(\del_j\zeta_{\vecm(i')}(\gminib)\cdot
 \pi^{\vecm(i')}_{\gminib} \bigr)^{\dagger},\,\,
 f^{\vecm(i)}_{j}
 \bigr]
 \Bigr|_h
=O\Bigl(
\exp\bigl(-\epsilon(7,i')\cdot|\vecz^{\vecm(i')}|\bigr)\Bigr)
 \cdot |f_j^{\vecm(i)}|_h
\]
Hence, we obtain the following inequality
on $X^{\ast}(R_{4,i})$ from (\ref{eq;07.5.18.3})
for some constant $C_{20,i}$:
\begin{equation}
\label{eq;07.6.18.5}
 \Delta_j\log\bigl|f_j^{\vecm(i)}\bigr|_h^2\leq
 -\frac{\bigl|[f_j^{\vecm(i)\dagger},f_j^{\vecm(i)}]\bigr|_h^2}
 {\bigl|f_j^{\vecm(i)}\bigr|_h^2}
+C_{20,i}
\end{equation}
Any eigenvalues of 
 $\bigl(z_j\cdot f_j^{\vecm(i)}\bigr)_{|Q}$
 for $Q\in X^{\ast}(R_{5,i})$ are
dominated by
$C_{21,i}\cdot\bigl|\vecz^{\vecm(i)}(Q)\bigr|$.
Hence, we can show $P(i)$
by a standard argument.
We give only an outline of the argument.
Let $\pi_j:X\lrarr D_j$ 
 denote the natural projection.
Let $P$ be any point of
$D_j\setminus\bigcup_{1\leq k\leq \ell,k\neq j}D_k$.
Let $\hyperh$ denote the upper half plane.
We have the universal covering 
$\varphi:\hyperh\lrarr\pi_j^{-1}(P)-\{P\}$
given by $z_j=\exp\bigl(\sqrt{-1}\zeta\bigr)$.
Let $Y(R_{5,i}):=\varphi^{-1} \bigl(X(R_{5,i})\bigr)$.
We put 
$F:=\varphi^{-1}\bigl(z_j\cdot f_j^{\vecm(i)}\bigr)$.
Let $\Delta_{\zeta}$ denote 
the Laplacian $-\del^2/\del \zeta\del \zetabar$.
Because of (\ref{eq;07.6.18.5}),
we have the following inequality
on $Y(R_{5,i})$:
\[
 \Delta_{\zeta}\log|F|_h^2\leq
 -\frac{\bigl|[F,F^{\dagger}]\bigr|_h^2}
 {|F|_h^2}
+C_{22,i}
\]
By using the argument for the proof of
Proposition 2.10 in \cite{mochi4}
(based on \cite{a} and \cite{s2}),
we obtain that 
$|F|^2_h$ is dominated by
the sum of the square of the absolute values 
of the eigenvalues of $F$ on $Y(R_{5,i})$,
and the estimate may depend only on $C_{22}$.
Thus, we obtain $P(i,\irr)$.
(Actually,
we do not use $Q(i-1,\irr)$.)

\subsection{$P(i,\irr)+Q(i-1,\irr)+R(i-1,\irr)$
$\Longrightarrow$ $Q(i,\irr)$}
\label{subsection;07.10.4.10}

We put $j:=\gminih(i)$ for simplicity of the description.
We remark the $j$-th component of $\vecm(i)$ is negative,
which we will use implicitly.
For any point $Q\in X-D$,
let $G_j^{\vecm(i)}(Q):=\sum \gminim(\beta)\cdot
 |\beta|^2$,
where $\beta$ runs through the eigenvalues
of $f^{\vecm(i)}_{j|Q}$,
and $\gminim(\beta)$ denotes the multiplicity of $\beta$.
We put $H_j^{\vecm(i)}:=
\bigl|f^{\vecm(i)}_{j|Q}\bigr|^2
 -G_j^{\vecm(i)}(Q)$.
We have only to show
$H_j^{\vecm(i)}
 =O\bigl( \bigl|\vecz^{\vecm(i)}\bigr|^2 \bigr)$.
We consider the following:
\[
 \tau_j^{\vecm(i)}:=
 \sum_{\gminia\in\Irr(\theta)}
 \del_j\bigl(\gminia-\etabar_{\vecm(i-1)}(\gminia)\bigr)
 \cdot
 \pi_{\gminia}'
+\sum_{\alpha\in\Sp(\theta,j)}
 \frac{\alpha}{z_j}
 \cdot\pi_{j,\alpha}'
\]
We have
$\bigl|\tau^{\vecm(i)}_j\bigr|_h\sim
 |\vecz^{\vecm(i)-\vecdelta_j}|$,
and 
$\Delta_j\log\bigl|\tau_j^{\vecm(i)}\bigr|_h^2$
is $C^{\infty}$.
We set
\[
 k_j^{\vecm(i)}:=
 \log\left(
 \frac{\bigl|f_j^{\vecm(i)}\bigr|^2_h}
 {\bigl|\tau_j^{\vecm(i)}\bigr|^2_h}
\right)
=\log\left(
 1+\frac{
 H_j^{\vecm(i)}
+G_j^{\vecm(i)}
-\bigl|\tau^{\vecm(i)}_j\bigr|_h^2}
 {\bigl|\tau_j^{\vecm(i)}\bigr|_h^2}
 \right).
\]
\begin{lem}
$\bigl|G_j^{\vecm(i)}
 -|\tau^{\vecm(i)}_j|^2\bigr|
=
 O\Bigl(
 \vecz^{\vecm(i)-\vecdelta_j}
 \Bigr)$.
\end{lem}
\pf
Let $f^{\vecm(i)}_{j,(\gminia,\alpha)}$ 
denote the restriction
of $f_j^{\vecm(i)}$
to $E_{\gminia}\cap E_{j,\alpha}$,
and let $\beta$ be any eigenvalue of
$f^{\vecm(i)}_{j,(\gminia,\alpha)|Q}$.
We have the following estimate
(see Remark \ref{rem;07.6.1.5}):
\[
\left|
 \frac{\del\gminia}{\del z_j}(Q)
+\frac{\alpha}{z_j(Q)}
-\beta
 \right|
\leq C_{30}
\]
Hence, we have the following:
\[
 |\beta|^2
-\left|
 \frac{\del\gminia}{\del z_j}(Q)
+\frac{\alpha}{z_j(Q)}
 \right|^2
=O\Bigl(\bigl|
 \vecz^{\vecm(i)-\vecdelta_j}
\bigr|
 \Bigr)
\]
Then, the claim of the lemma follows.
\hfill\qed

\vspace{.1in}

By using a standard argument
using elementary linear algebra
(see Page 729 of \cite{s2}),
we can show that 
there exists a constant $C_{31}>0$ such that
$\bigl|[f_j^{\vecm(i)\dagger},f_j^{\vecm(i)}]\bigr|_h
\geq C_{31}\cdot H_j^{\vecm(i)}$.
Recall (\ref{eq;07.6.18.5}),
where we have used $R(i-1,\irr)$.
Hence, we have the following estimate
on $X^{\ast}(R_{30,i})$
for some constants $C_{i}$:
\[
 \Delta_jk_j^{\vecm(i)}\leq
 -C_{32}\cdot \frac{\bigl(H_j^{\vecm(i)}\bigr)^2}
 {|\vecz^{\vecm(i)-\vecdelta_j}|^2}
+C_{33}
\leq
 -C_{34}|\vecz^{\vecm(i)-\vecdelta_j}|^2
\left(
 \frac{H_j^{\vecm(i)}}
 {\bigl|\tau_j^{\vecm(i)}\bigr|_h^2}
 \right)^2
+C_{35}
\]
We have already known
$H_j^{\vecm(i)}=O\Bigl(
 \bigl|\vecz^{\vecm(i)-\vecdelta_j}\bigr|^2
 \Bigr)$ due to $P(i,\irr)$.
Hence, we have the following estimate:
\begin{equation}
\label{eq;07.11.17.10}
 k_j^{\vecm(i)}
\sim
\frac{\bigl|
 H_j^{\vecm(i)}
+G_j^{\vecm(i)}
 -|\tau_j^{\vecm(i)}|_h^2\bigr|}
 {\bigl|\tau_j^{\vecm(i)}\bigr|_h^2}
\end{equation}

Put
$\nbigz(L):=\bigl\{Q\in X^{\ast}(R_{30,i})\,\big|\,
 H_j^{\vecm(i)}\geq L\cdot
 \bigl|\vecz^{\vecm(i)}\bigr|^2
 \bigr\}$ for some large $L>0$.
Note $2\vecm(i)\leq \vecm(i)-\vecdelta_j$.
Hence, $H_j^{\vecm(i)}$ is sufficiently larger than
$\bigl|G^{\vecm(i)}_j-|\tau^{\vecm(i)}_j|^2\bigr|$
on $\nbigz(L)$,
and we obtain the following
on $\nbigz(L)\cap X^{\ast}(R_{31,i})$:
\begin{equation}
 \label{eq;07.7.18.1}
C_{36}^{-1}\frac{H_j^{\vecm(i)}}
 {\bigl|\tau_j^{\vecm(i)}\bigr|^2}
\leq
 k_j^{\vecm(i)}
\leq
 C_{36}\frac{H_j^{\vecm(i)}}
 {\bigl|\tau_j^{\vecm(i)}\bigr|^2}
\end{equation}
We have
$2\bigl(\vecm(i)-\vecdelta_j\bigr)+4\vecm(i)
 -4\bigl(\vecm(i)-\vecdelta_j\bigr)
=2\vecm(i)+2\vecdelta_j\leq \veczero$.
Hence, we have the following inequality on 
$\nbigz(L)\cap X^{\ast}(R_{31,i})$:
\begin{equation}
 \label{eq;07.5.18.6}
 \Delta_j k_j^{\vecm(i)}
\leq
 -C_{37}\cdot \bigl|
 \vecz^{\vecm(i)-\vecdelta_j}
 \bigr|^2\bigl(k_j^{\vecm(i)}\bigr)^2
\leq
 -C_{38}|z_j|^{-4}
 \bigl(k_j^{\vecm(i)}\bigr)^2
\end{equation}
We would like to compare the functions
$|z_j|^2$ and $k_j^{\vecm(i)}$.
\begin{lem}
Let $R_{32,i}<R_{31,i}$.
We can take $C_{39}$ 
with the following property:
\begin{itemize}
\item
$C_{39}\cdot |z_j|^2>k_j^{\vecm(i)}$
on $\bigl\{|z_j|=R_{32,i}\bigr\}$.
\item
If $k_j^{\vecm(i)}(Q)\geq C_{39}\cdot |z_j(Q)|^2$
for $Q\in X^{\ast}(R_{32,i})$, we have $Q\in \nbigz(L)$.
\end{itemize}
\end{lem}
\pf
Due to (\ref{eq;07.11.17.10})
and $H^{\vecm(i)}_j=O(|z^{\vecm(i)-\vecdelta_j}|^2)$,
we have the boundedness of $k_j^{\vecm(i)}$.
Hence, the first condition is satisfied
for sufficiently large $C_{39}$.
Due to (\ref{eq;07.7.18.1}),
the second condition is satisfied 
for sufficiently large $C_{39}$.
\hfill\qed

\vspace{.1in}

Take $C_{40}>\max\bigl\{C_{38}^{-1},C_{39}\bigr\}$.
Note the following inequality for any $\eta\geq 0$:
\[
\Delta_j\Bigl(C_{40}\cdot |z_j|^2-\eta\log|z_j|\Bigr)
>
 -C_{38}|z_j|^{-4}\bigl(C_{40}\cdot|z_j|^2-\eta\log|z_j|\bigr)^2
\]
Let $P$ be any point of 
$D_j-\bigcup_{1\leq k\leq \ell,k\neq j}D_k$.
Let us consider the following set:
\[
\nbigz_{P}(C_{40},\eta):=
 \bigl\{
 Q\in \pi_j^{-1}(P)\cap X^{\ast}(R_{32,i})\,\big|\,
  k_j^{\vecm(i)}(Q)>C_{40}|z_j(Q)|^2
-\eta\log|z_j(Q)|^2
 \bigr\}
\]
Since $k_j^{\vecm(i)}$ is bounded,
the closure of
$\nbigz_{P}(C_{40},\eta)$ in
$\pi^{-1}(P)-\{P\}$ is compact.
By the choice of $C_{39}$,
the closure of $\nbigz_{P}(C_{40},\eta)$
has no intersection with $\{|z_j|=R_{32,i}\}$.
Because $\nbigz_P(C_{40},\eta)
\subset \nbigz(L)\cap X^{\ast}(R_{31,i})$,
we have the following on $\nbigz_{P}(C_{40},\eta)$:
\[
 \Delta_j\bigl(
 k_j^{\vecm(i)}
-\bigl(C_{40}|z_j|^2-\eta\log|z_j|\bigr)
 \bigr)
<
-C_{38}|z_j|^{-4}\Bigl(
 \bigl(
 k_j^{\vecm(i)}\bigr)^2
-\bigl(C_{40}|z_j|^2-\eta\log|z_j| \bigr)^2
\Bigr)
\leq 0
\]
Therefore,
$k_j^{\vecm(i)}
-\bigl(C_{40}|z_j|^2-\eta\log|z_j|\bigr)$
takes the maximum at the boundary of
$\nbigz_{P}(C_{40},\eta)$, which has to be $0$.
Thus, we have arrived at the contradiction,
and we can conclude 
$\nbigz_{P}(C_{40},\eta)=\emptyset$.
By taking $\eta\to 0$,
we obtain 
$k_j^{\vecm(i)}\leq C_{40}|z_j|^2$
on $\pi_j^{-1}(P)\cap X^{\ast}(R_{32,i})$.
Due to (\ref{eq;07.7.18.1}),
we obtain
$H_j^{\vecm(i)}=
O\bigl(
|\vecz^{\vecm(i)}|^2\bigr)$
and thus $Q(i,\irr)$.
(Actually, we do not use $Q(i-1,\irr)$.)

\subsection{$P(i,\irr)+Q(i,\irr)+R(i-1,\irr)$
$\Longrightarrow$ $R(i,\irr)$}
\label{subsection;07.7.18.2}

We continue to put $j:=\gminih(i)$.
Because of the commutativity
$[f_j,\pi_{\gminib}^{\vecm(i)}]=0$,
we have the following inequality:
\begin{equation}
 \label{eq;07.11.17.11}
 \Delta_j\log\bigl|\pi_{\gminib}^{\vecm(i)}\bigr|_h^2
\leq
 -\frac{\bigl|[f_j^{\dagger},
 \pi^{\vecm(i)}_{\gminib}]\bigr|_h^2}
 {|\pi^{\vecm(i)}_{\gminib}|_h^2}
\end{equation}
We consider the following function:
\[
 k^{\vecm(i)}_{\gminib}
:=\log\left(\frac{|\pi^{\vecm(i)}_{\gminib}|_h^2}
 {|\pi^{\vecm(i)\prime}_{\gminib}|_h^2}
 \right)
=
\log\left(
 1+\frac{|\nbigr^{\vecm(i)}_{\gminib}|_h^2}
 {|\pi^{\vecm(i)\prime}_{\gminib}|_h^2}
 \right)
\]
Because $[f_j^{\vecm(i)},\pi^{\vecm(i)}_{\gminib}]=0$,
we have the following:
\begin{multline*}
 0=\bigl[f_j^{\vecm(i)},\nbigr^{\vecm(i)}_{\gminib}\bigr]
+\bigl[f_j^{\vecm(i)},\pi^{\vecm(i)\prime}_{\gminib}\bigr]
 \\
=\left[
 \sum\del_j\zeta_{\vecm(i)}(\gminib)\cdot
 \pi^{\vecm(i)\prime}_{\gminib}
+\mu^{\vecm(i)}_j,\,
 \nbigr^{\vecm(i)}_{\gminib}
 \right]
+\bigl[f_j^{\vecm(i)},
 \pi^{\vecm(i)\prime}_{\gminib}\bigr]
\end{multline*}
We have
$\bigl[\mu_j^{\vecm(i)},
 \pi^{\vecm(i)\prime}_{\gminib}\bigr]
=\bigl[f_j^{\vecm(i)},
 \pi^{\vecm(i)\prime}_{\gminib}\bigr]
=O\bigl(|\vecz^{\vecm(i)}|\bigr)$
due to $Q(i,\irr)$,
and hence
$\nbigr^{\vecm(i)}_{\gminib}=O\bigl(|z_j|\bigr)$.
In particular,
we have the following estimate:
\begin{equation}
 \label{eq;07.11.17.12}
 k^{\vecm(i)}_{\gminib}
\sim
 \bigl|\nbigr^{\vecm(i)}_{\gminib}\bigr|_h^2
\end{equation}
From (\ref{eq;07.11.17.11}),
we have the following:
\[
 \Delta_j k^{\vecm(i)}_{\gminib}
\leq
 -\frac{\bigl|[f_j^{\dagger},
 \pi^{\vecm(i)}_{\gminib}]\bigr|_h^2}
 {\bigl|\pi^{\vecm(i)}_{\gminib}\bigr|_h^2}
\]
Note the following:
\[
 f_j^{\dagger}
=\sum_{i'\leq i}
 \sum_{\gminic\in\Irrbar(\theta,\vecm(i'))}
 \!\!\!
 \overline{\del_j\zeta_{\vecm(i')}(\gminic)}\cdot
 \pi^{\vecm(i')\prime}_{\gminic}
+\sum_{i'<i}\sum_{\gminib\in \Irrbar(\theta,\vecm(i'))}
 \!\!\!
 \overline{\del_j\zeta_{\vecm(i')}(\gminic)}\cdot
 \nbigr^{\vecm(i')\dagger}_{\gminic}
+\mu_j^{\vecm(i)\dagger}
\]
Because of $Q(i,\irr)$ and $R(i-1,\irr)$,
the second and third terms in the right hand side
are assumed to be much smaller than the first term
on $X^{\ast}(R_{40,i})$.
Let $p^{\vecm(i)}$ denote the projection
of $\End(E)$ onto the direct summand
$\bigoplus_{\gminib>\gminib'}
 Hom\bigl(E_{\gminib}^{\vecm(i)},E^{\vecm(i)}_{\gminib'}\bigr)$.
From the equality
$ [f_j^{\dagger},\pi^{\vecm(i)}_{\gminib}]
=[f_j^{\dagger},\nbigr^{\vecm(i)}_{\gminib}]
+[f_j^{\dagger},\pi^{\vecm(i)\prime}_{\gminib}]$,
we have the following inequality on $X^{\ast}(R_{42,i})$:
\[
 \bigl|
 [f_j^{\dagger},\pi^{\vecm(i)}_{\gminib}]
 \bigr|_h^2
\geq
 \bigl|
 p^{\vecm(i)}[f_j^{\dagger},\nbigr^{\vecm(i)}_{\gminib}]
 \bigr|_h^2
\geq
 C_{41}\cdot |\vecz^{\vecm(i)-\vecdelta_j}|^2
 \cdot \bigl|\nbigr^{\vecm(i)}_{\gminib}\bigr|_h^2
\]
Hence, we obtain the following on $X^{\ast}(R_{42,i})$:
\[
 \Delta_j k^{\vecm(i)}_{\gminib}
\leq
 -C_{42}\cdot |\vecz^{\vecm(i)-\vecdelta_j}|^2
 \cdot
 \bigl|\nbigr^{\vecm(i)}_{\gminib}\bigr|^2
\leq
 -C_{43}\cdot|\vecz^{\vecm(i)-\vecdelta_j}|^2
 \cdot
 k^{\vecm(i)}_{\gminib}
\]
We take small $\epsilon(10,i)>0$ such that
the following holds:
\begin{equation}
 \label{eq;08.9.12.20}
  \left(\frac{m_j(i)}{2}\right)^2
 \cdot\epsilon(10,i)^2
\leq
 C_{43}
\end{equation}
We have the following inequality
for any $\eta\geq 0$:
\begin{multline}
\label{eq;08.1.7.10}
  \Delta_j
\Bigl(
 \exp\bigl(-\epsilon(10,i)|\vecz^{\vecm(i)}|\bigr)
-\eta\log|z_j|
\Bigr)
 \geq \\
 -\left(\frac{m_j(i)}{2}\right)^2
 \cdot\epsilon(10,i)^2\cdot
 \bigl|\vecz^{\vecm(i)-\vecdelta_j}\bigr|^2\cdot
\Bigl(
\exp\bigl(-\epsilon(10,i)|\vecz^{\vecm(i)}|\bigr)
-\eta\log|z_j|
\Bigr)
\end{multline}
For $\eta>0$ and 
$P\in D_j-\bigcup_{1\leq k\leq \ell,k\neq j}D_k$,
we put
\[
 \Psi_{\eta}(P,z_j):=
 \exp\Bigl(
 \epsilon(10,i)
 \prod_{p\neq j}
 \bigl|z_p^{m_p(i)}(P)\bigr|
\cdot
 \bigl(R_{42,i}^{m_j(i)}
 -\bigl|z_j^{m_j(i)}\bigr|\bigr)
 \Bigr)
-\eta\cdot\log|z_j|
\]
Because of (\ref{eq;08.1.7.10}),
we obtain the following:
\[
 \Delta_j\bigl(\Psi_{\eta}(P,z_j)\bigr)
\geq
 -\left(\frac{m_j(i)}{2}\right)^2\epsilon(10,i)^2
 \prod_{p\neq j}\bigl|z_p^{2m_p(i)}(P)\bigr|
 \cdot\bigl|z_j^{2m_j(i)-2}\bigr|
\times \Psi_{\eta}(P,z_j)
\]
We have $\Psi_{\eta}(P,z_j)=1-\eta\log R_{42,i}$
for $|z_j|=R_{42,i}$,
which is larger than $1/2$,
if $\eta$ is sufficiently small.
We have already known that
$k_{\gminib}^{\vecm(i)}$ is bounded.
Hence, we can take a constant $C_{44}$
such that 
$k^{\vecm(i)}_{\gminib|\pi_j^{-1}(P)}(z_j)
<C_{44}\cdot\Psi_{\eta}(P,z_j)$
on $\{|z_j|=R_{42,i}\}$.
Let us consider the following set:
\[
 \nbigz(P,\eta):=
\Bigl\{
 z_j\in\pi_j^{-1}(P)\cap X^{\ast}(R_{42,i})\,\big|\,
 k^{\vecm(i)}_{\gminib|\pi_j^{-1}(P)}(z_j)
>C_{44}\cdot\Psi_{\eta}(P,z_j)
\Bigr\}
\]
Since 
$k^{\vecm(i)}_{\gminib}$ is bounded,
the closure of
$\nbigz(P,\eta)$ in $\{0<|z_j|\leq 1\}$ is compact.
By our choice of $C_{44}$,
the closure of $\nbigz(P,\eta)$ has no intersection
with $\bigl\{|z_j|=R_{42,i}\bigr\}$.
On $\nbigz(P,\eta)$, we have the following inequality
by our choice of $\epsilon(10,i)$ 
as in (\ref{eq;08.9.12.20}):
\[
 \Delta_j\Bigl(
 k^{\vecm(i)}_{\gminib|\pi_j^{-1}(P)}
-C_{44}\cdot\Psi_{\eta}(P,z_j)
 \Bigr)
\leq 0
\]
Hence, the values of
$k^{\vecm(i)}_{\gminib|\pi_j^{-1}(P)}
-C_{44}\Psi_{\eta}(P,z_j)$ on 
$\nbigz(P,\eta)$
is not larger than the boundary values,
which is $0$.
Thus, we have arrived at the contradiction,
and we obtain $\nbigz(P,\eta)=\emptyset$.
By taking $\eta\to 0$,
we obtain the following inequality
on $\pi^{-1}(P)\cap X^{\ast}(R_{42,i})$:
\[
 k^{\vecm(i)}_{\gminib|\pi_j^{-1}(P)}\leq
 C_{44}\cdot
 \exp\Bigl(
 \epsilon(10,i)\cdot
 \prod_{p\neq j}\bigl|z_p^{m_p(i)}(P)\bigr|
\cdot
 \bigl(R^{m_j(i)}_{42,i}-\bigl|z_j^{m_j(i)}\bigr|\bigr)
 \Bigr)
\] 
Let $R_{43,i}:=R_{42,i}/2$.
Then, the following holds
on $\pi_j^{-1}(P)\cap X^{\ast}(R_{43,i})$:
\[
 k^{\vecm(i)}_{\gminib|\pi_j^{-1}(P)}\leq
 C_{45}\cdot
 \exp\Bigl(
 -\epsilon(11,i)\cdot
 \prod_{p\neq j}\bigl|z_p^{m_p(i)}(P)\bigr|
 \cdot
 \bigl|z_j^{m_j(i)}\bigr|
 \Bigr)
\] 
Hence, we obtain
$\bigl|\nbigr^{\vecm(i)}_{\gminib}\bigr|_h
=O\Bigl(
 \exp\bigl(-\epsilon(11,i)
 \cdot\bigl|\vecz^{\vecm(i)}\bigr|\bigr)
 \Bigr)$.

\vspace{.1in}

Thus the proof of 
Theorem \ref{thm;07.10.4.1}
is finished.
\hfill\qed

\subsection{$P(\reg)$}
\label{subsection;07.10.4.4}

Let us begin the proof of Theorem \ref{thm;07.10.4.3}.
We have the following:
\[
 f_j^{\dagger}=f_j^{\reg\dagger}
+\sum_i\sum_{\gminib\in\Irrbar(\theta,\vecm(i))}
 \!\!\!\!\!\!
 \overline{\del_j\zeta_{\vecm(i)}(\gminib)}\cdot
 \pi^{\vecm(i)}_{\gminib}
+\sum_i\sum_{\gminib\in\Irrbar(\theta,\vecm(i))}
 \!\!\!\!\!\!
 \overline{\del_j\zeta_{\vecm(i)}(\gminib)}\cdot
 \bigl(
 \pi^{\vecm(i)\dagger}_{\gminib}
-\pi^{\vecm(i)}_{\gminib}
 \bigr)
\]
We have the estimate
$\bigl|
\pi^{\vecm\dagger}_{\gminib}
-\pi^{\vecm}_{\gminib}\bigr|_h
=O\bigl(\exp(-\epsilon_1|\vecz^{\vecm}|)\bigr)$.
We note that 
$|\del_j\zeta_{\vecm(i)}(\gminib)|
=O(|\vecz^{\vecm(i)}|)$
in the case that the $j$-th component of
$\vecm(i)$ is $0$.
Therefore,
we have the following inequality
on $X^{\ast}(R_{60})$:
\begin{equation}
 \label{eq;07.5.18.8}
 \Delta_j\log|f_j^{\reg}|^2
\leq
 -\frac{\bigl|[f_j^{\dagger},f_j^{\reg}]\bigr|_h^2}
 {|f^{\reg}_j|^2}
\leq
 -\frac{\bigl|[f_j^{\reg\dagger},f_j^{\reg}]\bigr|_h^2}
 {|f_j^{\reg}|_h^2}+C_{60}
\end{equation}

Since the eigenvalues of $f_j^{\reg}$ 
(resp. $z_j\cdot f_j^{\reg}$)
are bounded on $X-D$
in the case $j>\ell$ (resp. $j\leq \ell$),
we obtain the desired estimate
by using the argument in Section
\ref{subsection;07.10.4.30}
and the inequality (\ref{eq;07.5.18.8}).

\subsection{$Q(\reg)$}

We put 
$\rho_j:=\sum_{\alpha\in\Sp(j)}
 \gminim(\alpha)\cdot|\alpha|^2$,
where $\gminim(\alpha)$ denotes 
the multiplicity of $\alpha$.
By considering the tensor products 
with the rank one Higgs bundle,
we may assume $\rho_j\neq 0$.
Similarly,
we also put
$G^{\reg}_j(Q):=\sum \gminim(\beta)\cdot |\beta|^2$,
where $\beta$ runs through the eigenvalues
of $f^{\reg}_{j|Q}$.
We put
$H_j^{\reg}:=
 \bigl|f^{\reg}_{j|Q}\bigr|_h^2-G^{\reg}_j(Q)$.
By Remark \ref{rem;07.6.1.5},
we have 
$|q_j|^2-H_j^{\reg}=O(1)$
and $\rho_j|z_j|^{-2}-G_j^{\reg}=O(|z_j|^{-1})$.
We have only to show
$H_j^{\reg}\leq C\bigl(|z_j|^{-2}(-\log|z_j|)^{-2}\bigr)$.
We put
\[
 k_j^{\reg}:=
 \log\left(\frac{|f_j^{\reg}|^2}{\rho_j|z_j|^{-2}}\right)
=
 \log\left(1+
 \frac{
 H_j^{\reg}+(G_j^{\reg}-\rho_j|z_j|^{-2})}
 {\rho_j|z_j|^{-2}}\right).
\]
Let us consider the set
\[
 \nbigz(L):=\bigl\{
 Q\in X^{\ast}(R_{60})\,\big|\,
 H_j^{\reg}(Q)
>L|z_j(Q)|^{-2}(-\log|z_j(Q)|)^{-2}\bigr\}
\]
for some large $L>0$.
On $\nbigz(L)$,
we have 
$k_j^{\reg}\sim  H_j^{\reg}\cdot|z_j|^{-2}$.
From (\ref{eq;07.5.18.8}),
we have the following inequality
on $\nbigz(L)$:
\[
 \Delta_j k_j^{\reg}\leq
-C_{61} \frac{(H_j^{\reg})^2}{|z_j|^{-2}}+C_{62}
\leq \frac{-C_{63}}{|z_j|^2}
 \left(\frac{H_j^{\reg}}{\rho_j|z_j|^{-2}}\right)^2
 \!\!\!
+C_{64}
\leq
 -C_{65} \frac{
 \bigl(k_j^{\reg}\bigr)^2}{|z_j|^2}
+C_{66}
\]
We have
$k_j^{\reg}\geq C_{70}(-\log|z_j|)^{-2}$
on $\nbigz(L)$.
Hence, we have the following on $\nbigz(L)$,
if $L$ is sufficiently large:
\begin{equation}
 \label{eq;08.1.7.11}
 \Delta_j(k^{\reg}_j)
\leq
-C_{67}\frac{\bigl(k^{\reg}_j\bigr)^2}{|z_j|^2}
\end{equation}

\begin{lem}
\label{lem;08.1.7.12}
We can take $C_{71}$ and $R_{61}$
with the following property:
\begin{itemize}
\item
 $\Delta_j\bigl(
 C_{71}(-\log|z_j|)^{-2}
 \bigr) \geq
 -C_{67}\cdot |z_j|^{-2}\bigl(
 C_{71}(-\log|z_j|)^{-2}
 \bigr)^{2}$.
\item
 $k_j^{\reg}<C_{71}(-\log|z_j|)^{-2}$
 on $\bigl\{|z_j|=R_{61}\bigr\}$.
\item
If $k_j^{\reg}(Q)>C_{71}(-\log|z_j(Q)|)^{-2}$
for $Q\in X^{\ast}(R_{61})$,
we have $Q\in \nbigz(L)$.
\end{itemize}
\end{lem}
\pf
The first condition can be checked by a direct calculation
as in \cite{s2} or \cite{mochi2}.
Since we have already known that $k^{\reg}_j$ is bounded,
the second condition can be satisfied.
Since we have 
$k_j^{\reg}\sim 
 \bigl(H_j^{\reg}+(G_j^{\reg}-\rho_j|z_j|^{-2})\bigr)
 \cdot |z_j|^{2}$,
the third condition can be satisfied.
\hfill\qed

\vspace{.1in}

Let $P$ be any point of 
$D_j-\bigcup_{1\leq k\leq \ell,k\neq j}D_k$,
and let us consider the following set:
\[
 \nbigz(\eta):=
 \bigl\{
 Q\in\pi_j^{-1}(P)\cap X^{\ast}(R_{60})\,\big|\,
 k^{\reg}_j(Q)>C_{71}(-\log|z_j|)^{-2}
 -\eta\cdot\log|z_j|
 \bigr\}
\]
Then, we can show $\nbigz(\eta)=\emptyset$
by using a standard argument as in 
Section \ref{subsection;07.10.4.10},
with (\ref{eq;08.1.7.11})
and Lemma \ref{lem;08.1.7.12}.
(See \cite{s2} or \cite{mochi2}.)
By taking $\eta\to 0$,
we obtain the estimate
$k_j^{\reg}\leq
 C_{71}\bigl(-\log|z_j|\bigr)^{-2}$,
which implies
$\bigl|H_j^{\reg}\bigr|=O\bigl(|z_j|^{-2}(-\log|z_j|)^{-2}\bigr)$.
Therefore, we obtain $Q(j,\reg)$.

\subsection{$R(\reg)$}
\label{subsection;07.10.4.5}

We have
$0=\bigl[f_j^{\reg},\pi_{j,\alpha}\bigr]
=\bigl[f_j^{\reg},\nbigr^{\reg}_{j,\alpha}\bigr]
+\bigl[f_j^{\reg},\pi_{j,\alpha}'\bigr]$.
We also have
$ \bigl[
 f_j^{\reg},\pi_{j,\alpha}'
\bigr]
=O\left(|z_j|^{-1}(-\log|z_j|)^{-1}\right)$
by $Q(\reg)$.
By using it,
we obtain the preliminary estimate
$\nbigr^{\reg}_{j,\alpha}
=O\Bigl(
 \bigl(-\log|z_j|\bigr)^{-1}\Bigr)$.
We have the following:
\begin{multline}
 f_j^{\dagger}=\sum_i\sum_{
 \gminib\in\Irrbar(\theta,\vecm(i))}
 \overline{\del_j\zeta_{\vecm(i)}(\gminib)}\cdot
 \pi^{\vecm(i)}_{\gminib}
+\sum_{\alpha\in\Sp(\theta,j)}
 \frac{\alphabar}{\zbar_j}
 \cdot \pi'_{j,\alpha} \\
+\sum_i\sum_{\gminib\in\Irrbar(\theta,\vecm(i))}
 \overline{\del_j\zeta_{\vecm(i)}(\gminib)}\cdot
 \bigl(
 \pi^{\vecm(i)\dagger}_{\gminib}
-\pi^{\vecm(i)}_{\gminib}
 \bigr)
+q_j^{\dagger}
\end{multline}
The last two terms are much smaller
than the first two terms
on $X^{\ast}(R_{80})$.
We have the following:
\begin{multline}
\bigl[f^{\dagger}_j,\pi_{j,\alpha}\bigr]
=\Bigl[
 \sum_i\sum_{\gminib}
 \overline{\del_j\zeta_{\vecm(i)}(\gminib)}
 (\pi^{\vecm(i)\dagger}_{\gminib}
 -\pi^{\vecm(i)}_{\gminib}),
\pi'_{j,\alpha}
 \Bigr]
+[q^{\dagger}_j,\,\pi'_{j,\alpha}]\\
+\Bigl[
 \sum_{\beta}\frac{\betabar}{\zbar_j}
 \cdot \pi'_{j,\beta},
 \nbigr^{\reg}_{j,\alpha}
 \Bigr]
+\Bigl[
 \sum_i\sum_{\gminib}
 \overline{\del_j\zeta(\gminib)}
 (\pi^{\vecm(i)\dagger}_{\gminib}
 -\pi^{\vecm(i)}_{\gminib}),
\nbigr^{\reg}_{j,\alpha}
 \Bigr]
+[q^{\dagger}_j,\,\nbigr^{\reg}_{j,\alpha}] 
\end{multline}
By using an argument in 
Section \ref{subsection;07.7.18.2},
we obtain the following on $X^{\ast}(R_{80})$:
\begin{multline}
  \bigl|[f_j^{\dagger},\pi_{j,\alpha}]\bigr|_h
 \geq
\Bigl|
 \Bigl[\sum \frac{\betabar}{\zbar_j}\pi'_{j,\beta},\,\,
 \nbigr^{\reg}_{j,\alpha}\Bigr]
 \Bigr|_h
-\left|
\Bigl[
 \sum_i\sum_{\gminib}
 \overline{\del_j\zeta(\gminib)}
 (\pi^{\vecm(i)\dagger}_{\gminib}
 -\pi^{\vecm(i)}_{\gminib}),
\nbigr^{\reg}_{j,\alpha}
 \Bigr]\right|_h
\\
-\left|[q^{\dagger},\nbigr^{\reg}_{j,\alpha}] 
\right|_h
-\left|
 \Bigl[
 \sum_i\sum_{\gminib}
 \overline{\del_j\zeta(\gminib)}
 (\pi^{\vecm(i)\dagger}_{\gminib}
 -\pi^{\vecm(i)}_{\gminib}),
\pi'_{j,\alpha}
 \Bigr]
 \right|_h
\geq
 C_{80}|z_j|^{-1}\bigl|\nbigr^{\reg}_{j,\alpha}\bigr|_h
-C_{81}
\end{multline}
We consider the following function:
\[
 k_{j,\alpha}:=
 \log\left(\frac{|\pi_{j,\alpha}|_h^2}
 {|\pi'_{j,\alpha}|_h^2}\right)
=\log\left(
 1+\frac{|\nbigr^{\reg}_{j,\alpha}|_h^2}
 {|\pi_{j,\alpha}'|_h^2}
 \right)
\]
Recall that we have already known
that $k_{j,\alpha}$ is bounded.
Hence, we have the following inequality
on $X^{\ast}(R_{80})$:
\[
 \Delta_j k_{j,\alpha}\leq
-\frac{\bigl|[f_j^{\dagger},\pi_{j,\alpha}]\bigr|_h^2}
 {|\pi_{j,\alpha}|^2}
\leq
-\frac{C_{82}}{|z_j|^2}
 \frac{\bigl|\nbigr^{\reg}_{j,\alpha}\bigr|_h^2}
 {\bigl|\pi_{j,\alpha}'\bigr|_h^2}
+C_{83}
 \leq
 -\frac{C_{84}}{|z_j|^2}\cdot k_{j,\alpha}
+C_{85}
\]
By using an argument as in Section
\ref{subsection;07.10.4.10}
(see also \cite{s2} or \cite{mochi2}),
we obtain
$k_{j,\alpha}=O\bigl(
 |z_j|^{\epsilon(20)}\bigr)$,
and hence
$\bigl|\nbigr^{\reg}_{j,\alpha}\bigr|_h
 =O\bigl(
 |z_j|^{\epsilon(20)}\bigr)$.
Thus, the proof of Theorem \ref{thm;07.10.4.3}
is accomplished.
\hfill\qed

\subsection{Proof of Proposition
\ref{prop;07.11.22.11}}
\label{subsection;07.11.22.12}

\begin{lem}
\label{lem;08.9.15.20}
We have the following estimates:
\begin{equation}
\label{eq;07.11.22.2}
 \Bigl|
 \bigl[ \pi^{\vecm(p)\dagger}_{\gminia},
 \pi^{\vecm(q)}_{\gminib} \bigr]
\Bigr|_h
=O\Bigl(
 \exp\bigl(
 -\epsilon_{10}|\vecz^{\vecm(p)}|
 -\epsilon_{10}|\vecz^{\vecm(q)}|\bigr)
 \Bigr)
\end{equation}
\begin{equation}
\label{eq;07.11.22.3}
 \Bigl|
\bigl[
 \pi^{\vecm(p)\dagger}_{\gminia},
 \pi_{j,\alpha}
\bigr]
\Bigr|_h
=\Bigl|
 \bigl[
  \pi^{\vecm(p)}_{\gminia},
 \pi_{j,\alpha}^{\dagger}
 \bigr]
 \Bigr|
=O\Bigl(
 \exp\bigl(
 -\epsilon_{10}\cdot |\vecz^{\vecm(p)}|
 \bigr)\cdot 
 |z_j|^{\epsilon_{10}}
 \Bigr)
\end{equation}
\begin{equation}
\label{eq;07.11.22.4}
 \Bigl|
\bigl[ \pi_{i,\alpha},\pi^{\dagger}_{j,\alpha}\bigr]
\Bigr|_h
=O\Bigl(
 |z_i|^{\epsilon_{10}}
\cdot
 |z_j|^{\epsilon_{10}}
 \Bigr)
\end{equation}
\end{lem}
\pf
Due to Corollary \ref{cor;07.11.22.1},
we obtain the following:
\[
 \bigl[\pi^{\vecm(p)\dagger}_{\gminia},\,
 \pi^{\vecm(q)}_{\gminib}\bigr]
=\bigl[
(\pi^{\vecm(p)\dagger}_{\gminia}-\pi^{\vecm(p)}_{\gminia}),\,
 \pi^{\vecm(q)}_{\gminib}\bigr]
=O\Bigl(\exp\bigl(-\epsilon_1|\vecz^{\vecm(p)}|\bigr)
 \Bigr)
\]
Similarly, we obtain 
$\bigl[
 \pi^{\vecm(p)\dagger}_{\gminia},
 \pi^{\vecm(q)}_{\gminib}
 \bigr]
=O\Bigl(\exp\bigl(-\epsilon_1|\vecz^{\vecm(q)}|\bigr)
 \Bigr)$.
Then, we obtain (\ref{eq;07.11.22.2}).
The estimates (\ref{eq;07.11.22.3}) and (\ref{eq;07.11.22.4})
can be shown using a similar argument
together with Corollary \ref{cor;07.11.22.1}
and Corollary \ref{cor;07.11.22.5}.
\hfill\qed

\begin{lem}
\label{lem;08.12.8.3}
We have the following estimates
for $j=1,\ldots,\ell$:
\begin{equation}
\label{eq;07.11.22.6}
\bigl[f^{\nil}_{j},
 \pi^{\vecm(p)\dagger}_{\gminia}
 \bigr]
=O\Bigl(
 |z_j|^{-1}\bigl(-\log|z_j|\bigr)^{-1}
 \cdot\exp\bigl(-\epsilon_1|\vecz^{\vecm(p)}|\bigr)
 \Bigr),
\end{equation}
\[
\bigl[f^{\nil}_{j},
 \pi^{\dagger}_{i,\alpha}
 \bigr]
=O\Bigl(
 |z_j|^{-1}\bigl(-\log|z_j|\bigr)^{-1}
\cdot|z_i|^{\epsilon_2}
 \Bigr)
\]
We have the following estimates
for $j=\ell+1,\ldots,n$:
\begin{equation}
\label{eq;07.11.22.7}
\bigl[f^{\reg}_{j},
 \pi^{\vecm(p)\dagger}_{\gminia}
 \bigr]
=O\Bigl(
 \exp\bigl(-\epsilon_1|\vecz^{\vecm(p)}|\bigr)
 \Bigr),
\quad
\bigl[f^{\reg}_{j},
 \pi^{\dagger}_{i,\alpha}
 \bigr]
=O\Bigl(
 |z_i|^{\epsilon_2}
 \Bigr)
\end{equation}
\end{lem}
\pf
We have the following equalities:
\[
 \bigl[f^{\nil}_j,\,\pi^{\vecm(p)\dagger}_{\gminia}\bigr]
=\bigl[f^{\nil}_j,\,
 \pi^{\vecm(p)\dagger}_{\gminia}
 -\pi^{\vecm(p)}_{\gminia}
 \bigr],
\quad
 \bigl[f^{\nil}_j,\,\pi^{\dagger}_{i,\alpha}\bigr]
=\bigl[f^{\nil}_j,\,
 \pi^{\dagger}_{i,\alpha}
 -\pi_{i,\alpha}
 \bigr]
\]
Then, the estimate (\ref{eq;07.11.22.6}) follows.
The estimate (\ref{eq;07.11.22.7}) can be shown similarly.
\hfill\qed

\vspace{.1in}

The naturally defined map
$\Irrbar\bigl(\theta,\vecm(p)\bigr)
\lrarr
 \Irrbar\bigl(\theta,\vecm(p-1)\bigr)$ is denoted by
$\etabar_{\vecm(p-1),\vecm(p)}$.
We set
\[
 \gbigu(p):=
 \Bigl\{
 (\gminic,\gminic')\in
 \Irrbar\bigl(\theta,\vecm(p)\bigr)^2
 \,\Big|\,
 \gminic\neq\gminic',\,\,
 \etabar_{\vecm(p-1),\vecm(p)}(\gminic)
=\etabar_{\vecm(p-1),\vecm(p)}(\gminic')
 \Bigr\}.
\]
We have the decomposition:
\[
 \End(E)=\bigoplus_p
 \bigoplus_{(\gminic,\gminic')\in\gbigu(p)}
 \Hom\bigl(E^{\vecm(p)}_{\gminic},
 E^{\vecm(p)}_{\gminic'}\bigr)
\oplus
 \bigoplus_{\gminia\in\Irr(\theta)}
 \End(E_{\gminia},E_{\gminia})
\]
For any section $F$ of $\End(E)$,
we have the corresponding decomposition:
\[
 F=\sum_{p}
 \sum_{(\gminic,\gminic')\in\gbigu(p)}
 F^{\vecm(p)}_{\gminic,\gminic'}
+\sum_{\gminia\in\Irr(\theta)}
 F_{\gminia}
\]

\begin{lem}
\label{lem;07.11.22.10}
We have the following estimates:
\begin{equation}
 \label{eq;07.11.22.8}
 \bigl(\pi^{\vecm(p)\dagger}_{\gminib}
 -\pi^{\vecm(p)}_{\gminib}
 \bigr)^{\vecm(q)}_{\gminia,\gminia'}
=O\Bigl(
 \exp\bigl(-\epsilon_{10}|\vecz^{\vecm(p)}|
 -\epsilon_{10}|\vecz^{\vecm(q)}|
 \bigr)
 \Bigr)
\end{equation}
\begin{equation}
\label{eq;07.11.22.9}
\bigl(
 \pi^{\dagger}_{j,\alpha}-\pi_{j,\alpha}
\bigr)^{\vecm(q)}_{\gminia,\gminia'}
=O\Bigl(
 \exp\bigl(-\epsilon_{10}|\vecz^{\vecm(q)}|\bigr)
\cdot|z_j|^{\epsilon_{10}}
 \Bigr),
\quad
 (j=1,\ldots,\ell)
\end{equation}
\begin{equation}
\bigl(f_j^{\nil\dagger}\bigr)^{\vecm(p)}_{\gminia,\gminib}
=O\Bigl(
 |z_j|^{-1}(-\log|z_j|)^{-1}
 \cdot\exp\bigl(-\epsilon_{1}|\vecz^{\vecm(p)}|\bigr)
 \Bigr),
\quad
 (j=1,\ldots,\ell)
\end{equation}
\begin{equation}
\bigl(f_j^{\reg\dagger}\bigr)^{\vecm(p)}_{\gminia,\gminib}
=O\Bigl(
 \exp\bigl(-\epsilon_{1}|\vecz^{\vecm(p)}|\bigr)
 \Bigr),
\quad
 (j=\ell+1,\ldots,n)
\end{equation}
\end{lem}
\pf
We have the following equalities
for $\gminia\neq\gminia'$:
\begin{multline}
  \bigl(\pi^{\vecm(p)\dagger}_{\gminib}
 -\pi^{\vecm(p)}_{\gminib}
 \bigr)^{\vecm(q)}_{\gminia,\gminia'}
=\pi^{\vecm(q)}_{\gminia'}\circ
 \bigl(\pi^{\vecm(p)\dagger}_{\gminib}
 -\pi^{\vecm(p)}_{\gminib} \bigr)
 \circ\pi^{\vecm(q)}_{\gminia}\\
=\pi^{\vecm(q)}_{\gminia'}
\circ\pi^{\vecm(p)\dagger}_{\gminib}
\circ\pi^{\vecm(q)}_{\gminia} 
=\pi^{\vecm(q)}_{\gminia'}\circ
 \Bigl(
 \pi^{\vecm(p)\dagger}_{\gminib}
  \circ\pi^{\vecm(q)}_{\gminia}
-\pi^{\vecm(q)}_{\gminia}
  \circ\pi^{\vecm(p)\dagger}_{\gminib}
 \Bigr).
\end{multline}
Then, the estimate (\ref{eq;07.11.22.8})
follows from (\ref{eq;07.11.22.2}).
The other estimates 
can be shown similarly.
\hfill\qed

\vspace{.1in}

For each $j=1,\ldots,\ell$,
we have the decomposition:
\[
 \End(E)=
 \bigoplus_{\alpha,\beta\in\Sp(\theta,j)}
 \Hom(E_{j,\alpha},E_{j,\beta})
\]
For a section $F$ of $\End(E)$,
we have the corresponding decomposition
$F=\sum F_{j,\alpha,\beta}$.

\begin{lem}
\label{lem;07.11.22.21}
We have the following estimates:
\[
  \bigl(
\pi^{\vecm(p)\dagger}_{\gminia}
-\pi^{\vecm(p)}_{\gminia}\bigr)_{j,\alpha,\beta}
=O\Bigl(
\exp\bigl(-\epsilon_{11}|\vecz^{\vecm(p)}|\bigr)\cdot
 |z_j|^{\epsilon_{11}}\Bigr)
\]
For $p=1,\ldots,\ell$
and $\gamma\in\Sp(\theta,\gamma)$,
\[
  \bigl(
\pi^{\dagger}_{p,\gamma}
-\pi_{p,\gamma}\bigr)_{j,\alpha,\beta}
=O\bigl(|z_p|^{\epsilon_{11}}
 \cdot|z_j|^{\epsilon_{11}}\bigr)
\]
For $p=1,\ldots,\ell$,
\[
 \bigl(
 f_p^{\nil\dagger}
\bigr)_{j,\alpha,\beta}
=O\bigl(|z_j|^{\epsilon_{11}}
 |z_p|^{-1}(-\log|z_p|)^{-1}\bigr)
\]
For $p=\ell+1,\ldots,n$,
\[
  \bigl(
 f_p^{\reg\dagger}
\bigr)_{j,\alpha,\beta}
=O\bigl(|z_j|^{\epsilon_{11}}\bigr)
\]
\end{lem}
\pf
We can show them using the argument
in the proof of Lemma \ref{lem;07.11.22.10}.
\hfill\qed

\vspace{.1in}

Proposition \ref{prop;07.11.22.11}
follows from Lemma \ref{lem;07.11.22.10}
and Lemma \ref{lem;07.11.22.21}.
\hfill\qed

\subsection{Proof of Proposition \ref{prop;07.11.22.25}}
\label{subsection;08.1.7.2}

Let us consider the following:
\begin{equation}
\label{eq;07.11.22.35}
 \Phibar:=
 \sum_{\gminia\in\Irr(\theta)} d\gminiabar\cdot\pi_{\gminia}
+\sum_{j=1}^{\ell}
 \sum_{\alpha\in\Sp(\theta,j)}
 \alpha\cdot\frac{d\zbar_j}{\zbar_j}\cdot
 \pi_{j,\alpha}
\end{equation}
We have the following:
\begin{multline}
 \theta^{\dagger}-\Phibar
=\sum_{p}
 \sum_{\gminia\in\Irrbar(\theta,\vecm(i))}
 \overline{d\zeta_{\vecm(i)}(\gminia)}
\cdot \bigl(
 \pi^{\vecm(i)\dagger}_{\gminia}
-\pi^{\vecm(i)}_{\gminia}
\bigr) \\
+\sum_{j=1}^{\ell}
 \sum_{\alpha\in\Sp(\theta,j)}
 \alpha\cdot\frac{d\zbar_j}{\zbar_j}\cdot
 \bigl(\pi^{\dagger}_{j,\alpha}-\pi_{j,\alpha}\bigr) 
+\sum_{j=1}^{\ell}
 f_j^{\nil\dagger}\cdot d\zbar_j
+\sum_{j=\ell+1}^n
 f_j^{\reg\dagger}\cdot d\zbar_j
\end{multline}
From Proposition \ref{prop;07.11.22.11},
we obtain the following estimate
with respect to $h$ and $g_{\poin}$:
\begin{equation}
 \label{eq;10.5.24.10}
 \bigl(\theta^{\dagger}-\Phibar\bigr)_{
 (\gminia,\vecalpha),(\gminia',\vecalpha')}
=O\Bigl(
 \exp\bigl(
 -\epsilon|\vecz^{\ord(\gminia-\gminia')}|
 \bigr)
\cdot \nbigq_{\epsilon}(\vecalpha,\vecalpha')
 \Bigr)
\end{equation}
Then, we obtain the following estimate
with respect to $h$ and $g_{\poin}$:
\begin{multline}
 \bigl[
 \theta,\,\theta^{\dagger}
 \bigr]_{(\gminia,\vecalpha),(\gminia',\vecalpha')}
=\bigl[
 \theta,\,\theta^{\dagger}-\Phibar
 \bigr]_{(\gminia,\vecalpha),(\gminia',\vecalpha')} 
\\
=\left(
 d(\gminia-\gminia')
+\sum_{j=1}^{\ell} 
 (\alpha_j-\alpha_j')\frac{Oz_j}{z_j}
+O(1)
\right)\cdot
 \bigl(\theta^{\dagger}-\Phibar\bigr)
 _{(\gminia,\vecalpha),(\gminia',\vecalpha')} 
 \\
=O\Bigl(
 \exp\bigl(-\epsilon|\vecz^{\ord(\gminia-\gminia')}|\bigr)
\cdot\nbigq_{\epsilon}(\vecalpha,\vecalpha')
 \Bigr)
\end{multline}
Thus, the proof of Proposition \ref{prop;07.11.22.25}
is accomplished.
\hfill\qed

\section{The associated
 good filtered $\lambda$-flat bundle}
\label{subsection;08.9.14.40}
\subsection{Statements and some notation}
\label{subsection;07.11.24.20}

Let $X$ be a complex manifold,
and let $D$ be a simple normal crossing hypersurface of $X$
with the irreducible decomposition 
$D=\bigcup_{i\in \Lambda} D_{i}$.
Let $\harmonicbundle$ be 
a good wild harmonic bundle on $X-D$.
Let $\nbigelambda$ denote the holomorphic vector bundle
$(E,\delbar_E+\lambda\theta^{\dagger})$
on $X-D$.
\begin{notation}
\label{notation;07.12.2.45}
Let $\veca=(a_i\,|\,i\in \Lambda)\in\real^{\Lambda}$.
Let $U$ be an open subset of $X$ with 
a holomorphic coordinate
$(z_1,\ldots,z_n)$ such that
$U\cap D=\bigcup_{j=1}^{\ell}\{z_j=0\}$.
For each $j=1,\ldots,\ell$,
we have some $i(j)\in \Lambda$ such that
$D_{i(j)}\cap U=\{z_j=0\}$.
Let $b_j:=a_{i(j)}$.
We define
\begin{equation}
\label{eq;07.12.2.40}
 \nbigp_{\veca}\nbigelambda(U):=
 \Bigl\{
 f\in \nbigelambda(U\setminus D)\,\Big|\,
 |f|_h=O\Bigl(
 \prod_{j=1}^{\ell}|z_j|^{-b_j-\epsilon}
 \Bigr),\,\,\forall\epsilon>0
 \Bigr\}.
\end{equation}
Taking the sheafification,
we obtain the $\nbigo_X$-module
$\nbigp_{\veca}\nbigelambda$.
We also obtain the $\nbigo_X(\ast D)$-module
$\nbigp\nbigelambda:=\bigcup_{\veca}
 \nbigp_{\veca}\nbigelambda$.
The filtered sheaf 
$\bigl(\nbigp_{\veca}\nbigelambda\,
 \big|\,\veca\in\real^{\ell}\bigr)$
is denoted by $\nbigp_{\ast}\nbigelambda$.
\hfill\qed
\end{notation}
\index{sheaf $\nbigp_{\veca}\nbigelambda$}
\index{sheaf $\nbigp\nbigelambda$}
\index{filtered bundle  $\nbigp_{\ast}\nbigelambda$}

\begin{rem}
In our previous papers 
(for example {\rm\cite{mochi2}}),
we used the symbol
$\prolongg{\veca}{\nbigelambda}$
to denote $\nbigp_{\veca}\nbigelambda$.
Since we will consider several kinds of prolongation
in the wild case,
we prefer the symbol $\nbigp_{\veca}\nbigelambda$
for distinction.
\hfill\qed
\end{rem}

\begin{thm}
\label{thm;07.12.2.55}
$\nbigp_{\ast}\nbigelambda$ is
a filtered bundle on $(X,D)$.
The weak norm estimate up to log order holds
in the sense of Theorem 
{\rm \ref{thm;07.10.9.1}}.
\end{thm}
\pf
Due to Corollary \ref{cor;07.11.22.30},
$(\nbigelambda,h)$ is acceptable.
Hence, the claim follows from
Theorem \ref{thm;07.3.16.1}.
\hfill\qed

\vspace{.1in}
We use the symbol $\lefttop{i}F$
to denote the induced filtration
of $\nbigp_{\veca}\nbigelambda_{|D_i}$
for $i\in \Lambda$.

\vspace{.1in}

We have the induced Higgs field
and the induced hermitian metric of $\End(E)$,
which are denoted by 
the same symbols $\theta$ and $h$,
respectively.
Note that 
the harmonic bundle
$\bigl(\End(E),\delbar_{\End(E)},h,\theta\bigr)$
is a wild harmonic bundle on $X-D$,
but not necessarily good.
(See an example in Section 
\ref{subsection;10.5.24.2}.)
We use the symbol 
$\bigl(\End(\nbigelambda),
 \DDlambda\bigr)$ to denote 
the associated $\lambda$-flat bundle.
Although $\bigl(\End(E),\delbar_{\End(E)},h,\theta\bigr)$
is not necessarily good,
$\bigl(\End(\nbigelambda),h\bigr)$ 
is acceptable,
as remarked in Corollary \ref{cor;07.11.22.30}.
The prolongment corresponding to $\veca$
is denoted by
$\nbigp_{\veca}\End(\nbigelambda)$.
\index{sheaf $\nbigp_{\veca}\End(\nbigelambda)$}
We will implicitly use the following proposition
in the argument below,
which immediately follows from 
Proposition \ref{prop;07.11.18.5}.
\begin{prop}
\label{prop;08.9.14.45}
$\nbigp_0\End(\nbigelambda)$ is naturally
isomorphic to the sheaf of
local endomorphisms $f$
of $\nbigp_{\veca}\nbigelambda$
such that $f_{|D_i}$ preserve the parabolic filtrations
$\lefttop{i}F$ for $i\in \Lambda$.
\hfill\qed
\end{prop}

We will prove the following theorem 
in Section \ref{subsection;07.10.11.1}.
\begin{thm}
\mbox{{}}\label{thm;07.11.18.20}
\begin{itemize}
\item
$\DDlambda$ is a meromorphic flat $\lambda$-connection
with respect to $\nbigp_{\veca}\nbigelambda$.
\item
$(\nbigp_{\ast}\nbigelambda,\DDlambda)$
is a good filtered $\lambda$-flat bundle.
If $\theta$ is unramifiedly good,
$(\nbigp_{\ast}\nbigelambda,\DDlambda)$
is an unramified good filtered $\lambda$-flat bundle.
\item
Under the setting in Section
{\rm\ref{subsection;07.11.17.1}},
the set of the irregular values is given as follows:
\[
 \Irr(\nbigp\nbigelambda,\DDlambda)
=\bigl\{
 (1+|\lambda|^2)\cdot\gminia\,\big|\,
 \gminia\in\Irr(\theta)
 \bigr\}
\]
\end{itemize}
\end{thm}

The claims are local,
and it can be easily reduced 
to the unramified case.
Hence, we may and will give the proof
under the setting in Section 
\ref{subsection;07.11.17.1}.
Then, the claims in the case $\lambda=0$ 
is a direct consequence of Simpson's Main estimate:
The decomposition
(\ref{eq;07.7.18.10})
is prolonged to 
$\bigl(\nbigp_{\veca}\nbige^0,\DD^0\bigr)
=\bigoplus \bigl(
 \nbigp_{\veca}\nbige^0_{\gminia},
 \DD^0_{\gminia} \bigr)$ on $X$
due to the asymptotic orthogonality
in Theorem \ref{thm;07.10.4.1},
and we obtain that 
$\DD^0_{\gminia}
 -d\gminia\cdot \id_{\nbigp_{\veca}\nbige_{\gminia}^0}$
are logarithmic
due to the estimate of the norm
of the Higgs field in Theorem \ref{thm;07.10.4.3}.

\subsubsection{Characterization of
the Stokes filtration}

Before going to the proof for the case $\lambda\neq 0$,
let us state the characterization of the Stokes
filtration of the meromorphic $\lambda$-flat bundle
$(\nbigp\nbigelambda,\DDlambda)$
$(\lambda\neq 0)$
in terms of the growth order of 
the norms of the flat sections with respect to $h$.
Since the property is local,
we give the statement
under the setting in Section
\ref{subsection;07.11.17.1}.
Let $S$ be a small multi-sector in $X-D$.
We have the partial Stokes filtration
$\nbigf^{S,\vecm(i)}$ of $\nbigp\nbigelambda_{|\Sbar}$
in the level $\vecm(i)$
(Section \ref{subsection;07.11.18.10})
indexed by the ordered set 
$\bigl(\Irrbar(\theta,\vecm(i)),\leq^{\lambda}_S\bigr)$.
It is flat with respect to $\DDlambda$,
and it is characterized by the growth order
of the flat sections with respect to $h$.

\begin{prop}
\label{prop;07.10.11.3}
Let $f$ be a flat section of $\nbigelambda_{|S}$.
We have $f\in\nbigf^{S,\vecm(i)}_{\gminib}$
if and only if the following estimate holds
for some $C>0$ and $M>0$:
\[
 \bigl|f\cdot
 \exp\bigl((\lambda^{-1}+\lambdabar)\cdot
 \etabar_{\vecm(i)}(\gminib)\bigr)\bigr|_h
=O\Bigl(
 \exp\bigl(C|\vecz^{\vecm(i+1)}|\bigr)
\cdot\prod_{k(i+1)<j\leq \ell}|z_j|^{-M}
 \Bigr)
\]
Here, $k(i)$ are determined by
$\vecm(i)\in\seisuu_{<0}^{k(i)}\times\veczero$.
\end{prop}
\pf
It follows from 
Proposition \ref{prop;07.10.11.2}
and the weak norm estimate for 
the acceptable bundles.
(See Theorems \ref{thm;07.12.2.55}
and \ref{thm;07.10.9.1}).
\hfill\qed

\subsubsection{Example 1}
\label{subsection;10.5.24.4}

Let $X:=\Delta^n$ and 
$D:=\bigcup_{i=1}^{\ell}\{z_i=0\}$.
Take $(\veca,\vecalpha)\in
 \real^{\ell}\times\cnum^{\ell}$
and $\gminia\in M(X,D)$.
We assume that 
$\vecz^{-\vecm}\,\gminia$
is nowhere vanishing holomorphic function
on $X$ for some 
$\vecm\in\seisuu_{<0}^{k}\times\veczero_{\ell-k}$.
Let $L(\veca,\vecalpha,\gminia)$
be the harmonic bundle on $X-D$
given as follows:
\[
 L(\veca,\vecalpha,\gminia)
=\nbigo_{X-D}\cdot e,
\quad
 \theta=d\gminia+
 \sum_{i=1}^{\ell}\alpha_i\frac{dz_i}{z_i},
\quad
 h(e,e)=\prod_{i=1}^{\ell}|z_i|^{-2a_i}
\]
The associated operators
$\del$ and $\theta^{\dagger}$
are as follows:
\[
 \del\,e=e\,\Bigl(
 -\sum_{i=1}^{\ell}a_i\frac{dz_i}{z_i}
 \Bigr),
\quad
 \theta^{\dagger}
=d\gminiabar+
 \sum_{i=1}^{\ell}
 \alphabar_i\frac{d\zbar_i}{\zbar_i}
\]
Let $\nbigl^{\lambda}(\veca,\vecalpha,\gminia)$
be the associated $\lambda$-flat bundle.
We have the holomorphic section
$u^{\lambda}$ of 
$\nbigl^{\lambda}(\veca,\vecalpha,\gminia)$
given as follows:
\[
 u^{\lambda}
=\exp\Bigl(
 -\lambda\gminiabar
+\lambdabar\gminia
-\sum_{i=1}^{\ell}
 \lambda\,\alphabar_i\log|z_i|^2
 \Bigr)\,e
\]
We can easily check the following:
\[
 |u^{\lambda}|_h=
 \prod_{i=1}^{\ell}
 |z_i|^{-\paramap(\lambda,a_i,\alpha_i)},
\quad
 \DDlambda u^{\lambda}
=u^{\lambda}\Bigl(
 (1+|\lambda|^2)\,d\gminia
+\sum_{i=1}^{\ell}
 \eigenmap(\lambda,a_i,\alpha_i)
 \frac{dz_i}{z_i}
 \Bigr)
\]
See Subsection \ref{subsection;07.10.14.15}
for the maps
$\paramap(\lambda)$ and $\eigenmap(\lambda)$.
We put $\paramap(\lambda,\veca,\vecalpha)=\bigl(
 \paramap(\lambda,a_i,\alpha_i) \bigr)
 \in\real^{\ell}$.
We have a natural isomorphism
$\nbigp_{\paramap(\lambda,\veca,\vecalpha)}
 \nbigl^{\lambda}
 (\veca,\vecalpha,\gminia)
 \simeq\nbigo_{X}\cdot u^{\lambda}$.
The residues
$\Res_i(\DDlambda)$ on
$\nbigp_{\paramap(\lambda,\vecu)}\nbigl^{\lambda}
 (\veca,\vecalpha,\gminia)_{|D_i}$
are given by the multiplication of
$\eigenmap(\lambda,a_i,\alpha_i)$.

Note that if $\gminia\neq 0$,
$u^{\lambda}$ depends on $\lambda$
in a non-holomorphic way.

\subsubsection{Example 2}
\label{subsection;10.5.24.3}

Let $X:=\Delta$ and $D:=\{O\}$.
Let $V$ be a finite dimensional
$\cnum$-vector space
with a nilpotent map $N$.
Recall that we have a tame harmonic bundle
$E(V,N)$ on $X-D$ with the following property:
\begin{itemize}
\item
Let $\nbigelambda(V,N)$ be the associated
$\lambda$-flat bundle.
The parabolic structure of
$\nbigp_0\nbigelambda(V,N)$ is trivial,
and we have an isomorphism:
\[
 (\nbigp_0\nbigelambda(V,N)_{|O},
 \Res(\DDlambda))
\simeq
 (V,N)
\]
\end{itemize}
For example, it can be constructed as follows.
We put $y:=-\log|z|^2$.
Recall that we have the harmonic bundle
$Mod(2):=(E,\delbar_E,\theta,h)$ given as follows:
\[
E:=\nbigo_{X-D}\,e_1\oplus\nbigo_{X-D}\,e_{-1},
\quad
 H(h,\vece)=
 \left(
 \begin{array}{cc}
 y & 0 \\ 0 & y^{-1}
 \end{array}
 \right),
\quad
 \theta\,\vece=\vece\,
 \left(
 \begin{array}{cc}
 0 & 0 \\ 1 & 0
 \end{array}
 \right)
\frac{dz}{z}
\]
Here, $H(h,\vece)$ is the matrix
whose $(i,j)$-entries are $h(e_i,e_j)$.
We take a frame $\vecv^{\lambda}$ of
the holomorphic bundle $\nbigelambda$
given as follows:
\begin{equation}\label{eq;b12.9.10}
 \vecv^{\lambda}=\vece\,
 \left(
 \begin{array}{cc}
 1 & -\lambda\,y^{-1}\\
 0 & 1
 \end{array}
 \right)
\end{equation}
Then $\vecv^{\lambda}$ gives 
a frame of $\nbigp_0\nbige$,
and the following holds:
\[
 \DDlambda\vecv^{\lambda}=
 \vecv^{\lambda}\,
 \left(
 \begin{array}{cc}
 0 & 0 \\ 1 & 0
 \end{array}
 \right)\frac{dz}{z},
\quad
 \Res(\DDlambda)\vecv^{\lambda}_{|O}
=\vecv^{\lambda}_{|O}\,
 \left(
 \begin{array}{cc}
 0 & 0 \\ 1 & 0
 \end{array}
 \right)
\]
Then, we can construct $E(V,N)$
as a direct sum of some
symmetric tensor products
of $\Mod(2)$.

\subsection{Construction of a decomposition}
\label{subsection;07.7.7.23}

We use the setting in Section
\ref{subsection;07.11.17.1}.
Let $\lambda\neq 0$.
For $\vecm(i)=(m_1(i),\ldots,m_{\ell}(i))\in\nbigm$,
we put 
\[
 \supp(i):=\bigl\{j\,\big|\,m_j(i)\neq 0\bigr\}.
\]

\begin{lem}
\label{lem;07.7.7.20}
Let $\epsilon$ be any small positive number,
and let $N$ be any large number.
We can take holomorphic sections
$p_{\gminib}^{\vecm(i)}$ of
$\nbigp_0\End(\nbigelambda)$
for $\vecm(i)\in\nbigm$ and
 $\gminib\in\Irrbar(\theta,\vecm(i))$,
such that the following holds:
\[
 \bigl|
 p^{\vecm(i)}_{\gminib}
 -\pi^{\vecm(i)}_{\gminib}
\bigr|_{h}
=O\left(
 \prod_{j\in\supp(i)}|z_j|^{2N}
\cdot
 \prod_{j\not\in\supp(i)}
 |z_j|^{-\epsilon}
 \right)
\]
\[
 \bigl(
 p^{\vecm(i)}_{\gminib}
\bigr)^2=p^{\vecm(i)}_{\gminib},
\quad
 [p^{\vecm(i)}_{\gminib_1},
 p^{\vecm(i)}_{\gminib_2}]=0,
\quad
 p^{\vecm(i-1)}_{\gminia}
=\sum_{\substack{
 \gminib\in\Irrbar(\theta,\vecm(i)) \\
 \etabar_{i-1,i}(\gminib)=\gminia
 }}
 p^{\vecm(i)}_{\gminib}
\]
Here, $\etabar_{i-1,i}:\Irrbar(\theta,\vecm(i))\lrarr 
 \Irrbar(\theta,\vecm(i-1))$ denotes the naturally defined map.
\end{lem}
\pf
Let $d_{\lambda}'':=\delbar_E+\lambda\theta^{\dagger}$.
We have $d_{\lambda}''\pi^{\vecm(i)}_{\gminib}
=\lambda\cdot\bigl[\theta^{\dagger},
 \pi^{\vecm(i)}_{\gminib}\bigr]
=O\bigl(\exp(-\epsilon|\vecz^{\vecm(i)}|)\bigr)$
with respect to $h$ and $g_{\poin}$,
due to Proposition \ref{prop;07.11.22.11}.
Let $\epsilon_{11}$ be 
sufficiently smaller than $\epsilon$,
and let $N_{11}$ and $N_{12}$ be 
sufficiently larger than $N$.
By Lemma \ref{lem;07.9.23.36},
we can take sections $s^{\vecm(i)}_{\gminib}$
of $\End(E)$ on $X-D$
such that the following holds:
\[
 d_{\lambda}''s^{\vecm(i)}_{\gminib}
=d_{\lambda}''\pi^{\vecm(i)}_{\gminib},
\]
\[
 \int \bigl|s_{\gminib}^{\vecm(i)}\bigr|_h^2
 \prod_{j\in\supp(i)}
 |z_j|^{-10N_{11}}
 \prod_{j\not\in\supp(i)}
 |z_j|^{\epsilon_{11}}
\left(-\sum_{j=1}^{\ell} \log|z_j|\right)^{N_{12}}
\dvol_{g_{\poin}}<\infty
\]
By Lemma \ref{lem;07.9.24.5},
we obtain the following estimate:
\[
 \bigl|s_{\gminib}^{\vecm(i)}\bigr|_h
=O\left(
 \prod_{j\in\supp(i)}|z_j|^{4N_{11}}
 \prod_{j\not\in\supp(i)}
 |z_j|^{-2\epsilon_{11}}
 \right)
\]
We put
$\puebar^{\vecm(i)}_{\gminib}:=
 \pi^{\vecm(i)}_{\gminib}-s^{\vecm(i)}_{\gminib}$.
If $\epsilon_{11}$ is sufficiently small,
$\puebar^{\vecm(i)}_{\gminib}$
gives a section of
$\nbigp_0\End(\nbigelambda)$.
Moreover, we have the following:
\[
\Bigl(
 \bigl(\puebar_{\gminib}^{\vecm(i)}\bigr)^2
-\puebar^{\vecm(i)}_{\gminib}
\Bigr)\nbigp_{\veca}{\nbigelambda}
\subset
 \nbigp_{\veca-3N_{11}\vecdelta(i)}\nbigelambda
\]
\[
  \bigl[
 \puebar^{\vecm(i)}_{\gminib_1},
 \puebar^{\vecm(j)}_{\gminib_2}
 \bigr]\bigl(\nbigp_{\veca}{\nbigelambda}\bigr)
\subset
 \nbigp_{\veca-3N_{11}\vecdelta(j)}\nbigelambda\,\,\,\,
 (i\leq j)
\]
\[
\Bigl(
 \puebar^{\vecm(i-1)}_{\gminia}
-\sum_{\substack{
 \gminib\in\Irrbar(\theta,\vecm(i)) \\
 \etabar_{i-1,i}(\gminib)=\gminia
 }}
 \puebar^{\vecm(i)}_{\gminib}
\Bigr)
 \nbigp_{\veca}\nbigelambda
\subset
 \nbigp_{\veca-3N_{11}\vecdelta(i)}\nbigelambda
\]
Here, $\vecdelta(i)$ denote the elements
of $\seisuu^{\ell}$ 
such that the $j$-th components
$\delta_j(i)$ are given by
\[
 \delta_j(i):=\left\{
 \begin{array}{ll}
 1 & (j\in\supp(i))\\
 0 & (j\not\in\supp(i))
 \end{array}
 \right.
\]

We would like to modify $\puebar^{\vecm(i)}_{\gminib}$
inductively on $i$ so that the desired conditions
are satisfied.
Consider the following state $P(i)$:
\begin{description}
\item[$P(i)$]
 We have $p^{\vecm(l)}_{\gminib}
 \in \nbigp_0\End(\nbigelambda)$ 
 for any $l<i$
 and for $\gminib\in\Irrbar(\theta,\vecm(l))$,
such that the following holds:
\[
 \bigl(
 p^{\vecm(l)}_{\gminib}
-\puebar^{\vecm(l)}_{\gminib}
\bigr)
 \nbigp_{\veca}\nbigelambda
\subset
 \nbigp_{\veca-2N_{11}\vecdelta(l)}
 \nbigelambda
\]
\[
 \bigl(
 p^{\vecm(l)}_{\gminib}
\bigr)^2=p^{\vecm(l)}_{\gminib},
\quad
 \bigl[p^{\vecm(l)}_{\gminib_1},
 p^{\vecm(l)}_{\gminib_2}\bigr]=0,
\quad
 p^{\vecm(l-1)}_{\gminia}
=\sum_{\substack{
 \gminib\in\Irrbar(\theta,\vecm(l))\\
 \etabar_{l-1,l}(\gminib)=\gminia
 }}
 p^{\vecm(l)}_{\gminib}
\]
\[
\bigl[
 p_{\gminib_1}^{\vecm(l)},
 \puebar_{\gminib_2}^{\vecm(j)}
\bigr]=0\,\,
(l< i\leq j)
\]
\end{description}
Let us give a procedure 
from $P(i-1)$ to $P(i)$.

First,
we give a procedure
$P(0)\Longrightarrow P(1)$.
Take an injection
$\varphi:\Irrbar(\theta,\vecm(0))\lrarr\seisuu$.
We consider the following:
\[
 \Phibar_0:=
 \sum_{\gminia\in\Irrbar(\theta,\vecm(0))}
 \varphi(\gminia)\cdot 
 \puebar^{\vecm(0)}_{\gminia}
\]
Note that $\Phibar_0$ gives an endomorphism
of $\nbigp_{\veca}\nbigelambda$
for each $\veca\in\real^{\ell}$,
which preserves the filtrations 
$\lefttop{j}F$ ($j=1,\ldots,\ell$).
We have the decomposition
$\nbigp_0\nbigelambda
=\bigoplus_{m\in\seisuu}V_m$
such that (i) $\Phibar_0(V_m)\subset V_m$,
(ii) the eigenvalues of $\Phibar_0$ on $V_m$
are close to $m$,
(iii) it is compatible with the filtrations
 $\lefttop{j}F$ $(j=1,\ldots,\ell)$.
Let $p^{\vecm(0)}_{\gminia}$ denote
the projection onto $V_{\varphi(\gminia)}$.
Then, we have
$\bigl(p^{\vecm(0)}_{\gminia}
 -\puebar^{\vecm(0)}_{\gminia}\bigr)
 \bigl(\nbigp_{\veca}\nbigelambda\bigr)
\subset
 \nbigp_{\veca-3N_{11}\vecdelta(0)}\nbigelambda$.
We also have
(i) $\bigl[p^{\vecm(0)}_{\gminia},
   p^{\vecm(0)}_{\gminib}\bigr]=0$,
(ii) $(p^{\vecm(0)}_{\gminia})^2=
 p^{\vecm(0)}_{\gminia}$,
(iii) $\sum p^{\vecm(0)}_{\gminia}=\id$.

For $i>0$
and $\gminib\in\Irr(\theta,\vecm(i))$,
we have the decomposition
\[
 \puebar^{\vecm(i)}_{\gminib}
=\sum (\puebar^{\vecm(i)}_{\gminib})_{m,n},
\quad
(\puebar^{\vecm(i)}_{\gminib})_{m,n}
 \in \Hom(V_n,V_m).
\]
We put
$\puebar^{\vecm(i)\prime}_{\gminib}:=
  \sum_n (\puebar^{\vecm(i)}_{\gminib})_{n,n}$.
\begin{lem}
\label{lem;07.10.5.30}
We have
$(\puebar^{\vecm(i)\prime}_{\gminib}
-\puebar^{\vecm(i)}_{\gminib})
 \bigl(\nbigp_{\veca}\nbigelambda\bigr)
\subset
 \nbigp_{\veca-3N_{11}\vecdelta(i)}\nbigelambda$.
\end{lem}
\pf
 We put
 $\Phi_0:=
 \sum\varphi(\gminia)p^{\vecm(0)}_{\gminia}$.
 We have
 $(\Phi_0-\Phibar_0)\nbigp_{\veca}\nbigelambda
 \subset
 \nbigp_{\veca-3N_{11}\vecdelta(i)}\nbigelambda$.
For any positive integer $M$,
let $\Phi_0^M$ denote the $M$-iteration
of $\Phi_0$.
Then, we have the following:
\[
 \bigl[\Phi_0^M,\puebar^{\vecm(i)}_{\gminib}\bigr]
 \nbigp_{\veca}\nbigelambda
 \subset
 \nbigp_{\veca-3N_{11}\vecdelta(i)}\nbigelambda,
\quad
 \bigl[\Phi_0^M,\puebar^{\vecm(i)}_{\gminib}\bigr]
 =\sum (m^M-n^M)
 \bigl(\puebar^{\vecm(i)}_{\gminib}\bigr)_{m,n}
\]
Then, we can easily derive the claim of
Lemma \ref{lem;07.10.5.30}.
\hfill\qed

\vspace{.1in}
We replace $\puebar^{\vecm(i)}_{\gminia}$
with $\puebar^{\vecm(i)\prime}_{\gminia}$,
and then we arrive at the state $P(1)$.
Assume we are in the state $P(i-1)$.
We apply the above argument for 
$P(0)\Longrightarrow P(1)$
to each $\Image p^{\vecm(i-1)}_{\gminia}$
with the endomorphisms
$\puebar^{\vecm(j)}_{\gminib}$ $(j\geq i)$.
Then, we can arrive at $P(i)$.
When we arrive at $P(L)$,
the proof of Lemma \ref{lem;07.7.7.20}
is finished.
\hfill\qed

\subsection{Proof of Theorem
\ref{thm;07.11.18.20}}
\label{subsection;07.10.11.1}

Let $\delta_{\lambda}'$ be the $(1,0)$-operator
determined by the condition
$d''_{\lambda}+\delta_{\lambda}'$ is unitary,
i.e.,
$\delta_{\lambda}'=\del-\lambdabar\theta$.
Let $R(h)$ denote the curvature of
$d''_{\lambda}+\delta_{\lambda}'$.
We have
$ \DDlambda=\delbar_E+\lambda\theta^{\dagger}
+\lambda\del_E+\theta
=\delbar_E+\lambda\theta^{\dagger}
+\lambda\delta'_{\lambda}
+(1+|\lambda|^2)\theta$.
Let $p_{\gminib}^{\vecm(i)}$ be as in
Lemma \ref{lem;07.7.7.20}.
We put 
\[
 \DDlambda_0:=
\delbar_E+\lambda\theta^{\dagger}
+\lambda\delta_{\lambda}'
+(1+|\lambda|^2)
 \left(\theta-\sum_{i=1}^L\sum_{\gminib}
 d\zeta_{\vecm(i)}(\gminib)
\cdot p^{\vecm(i)}_{\gminib}
 \right).
\]
It gives a holomorphic connection of
$\nbigelambda$ on $X-D$,
which is not necessarily flat.
We have
$ \DDlambda=
 \DDlambda_0
+\sum_i\sum_{\gminib}
 d\zeta_{\vecm(i)}(\gminib)\cdot
 p^{\vecm(i)}_{\gminib}$.

\begin{prop}
\label{prop;07.7.18.70}
$\DDlambda_0$ is logarithmic
with respect to 
$\nbigp_{\veca}\nbigelambda$.
\end{prop}
\pf
By considering the tensor product
with the rank one harmonic bundle,
we may and will assume $\veca=(0,\ldots,0)$.
First, we consider the case in which $D$ is smooth,
say $D=D_1$.
Let $\pi_j:X-D\lrarr D_j$ denote the natural projection
for $j=1,\ldots,n$.
Let $d_{\lambda,j}''$,  $\delta_{\lambda,j}'$
and $R_j(h)$ denote
the restriction of $d_{\lambda}''$,
$\delta'_{\lambda}$ and $R(h)$ 
to the curves $\pi_j^{-1}(Q)$
$(Q\in D_j)$,
respectively.
Let $f$ be a holomorphic section of
$\nbigp_0\nbigelambda$.
Because of the acceptability of
$(\nbigelambda,h)$,
we have 
$|f|_h\leq C_{20}\cdot 
 \bigl(-\log|z_1|\bigr)^{N_{20}}$
(Proposition \ref{prop;07.7.7.3}).
In the following estimate,
we do not have to be concerned with the signature.

\begin{lem}
\label{lem;07.7.18.50}
Let $\chi(z_j)$ be any test function
on $\bigl\{z_j\in\cnum\,\big|\,|z_j|<1\bigr\}$.
Let $j\neq 1$.
For any $\epsilon>0$,
there exists a constant $C_{\epsilon}$ such that
the following holds for any $Q\in D_j\setminus D_1$:
\begin{equation}
 \label{eq;07.7.18.20}
 \int_{\pi_j^{-1}(Q)}
 \bigl(\delta_{\lambda,j}'(\chi\cdot f),\,
 \delta_{\lambda,j}'(\chi\cdot f)\bigr)_h
< C_{\epsilon}\cdot |z_1(Q)|^{-2\epsilon}
\end{equation}
\end{lem}
\pf
The left hand side of (\ref{eq;07.7.18.20})
can be rewritten as follows:
\begin{multline}
\pm\int_{\pi_j^{-1}(Q)}
 \bigl(\chi\cdot f,\,
 d_{\lambda,j}''
 \delta_{\lambda,j}'(\chi \cdot f)
 \bigr)_h
 =
\\
\pm\int_{\pi_j^{-1}(Q)}
 \bigl(\chi\cdot f,\,R_j(h)(\chi\cdot f)
 \bigr)_h
\pm\int_{\pi_j^{-1}(Q)}
 \bigl(\chi\cdot f,\,
 \delta_{\lambda,j}'
 d_{\lambda,j}''(\chi \cdot f)\bigr)_h=
 \\
\pm\int_{\pi_j^{-1}(Q)}
 \bigl(\chi\cdot f,\,R_j(h)(\chi\cdot f)
 \bigr)_h
\pm\int_{\pi_j^{-1}(Q)}
 \bigl(d_{\lambda,j}''(\chi\cdot f), \,
d_{\lambda,j}''(\chi\cdot f)
 \bigr)_h
\end{multline}
Then, the claim of Lemma \ref{lem;07.7.18.50}
easily follows.
\hfill\qed

\begin{lem}
\label{lem;07.7.18.60}
Fix a small $\epsilon>0$.
There exists a constant $C$ such that the following holds
for any $Q\in D_1$:
\[
 \int_{\pi_1^{-1}(Q)}
 \bigl(\delta'_{\lambda,1}f,\,
 \delta'_{\lambda,1}f
 \bigr)_h\cdot |z_1|^{2\epsilon}
<C
\]
\end{lem}
\pf
Let $\rho$ be a non-negative 
 $C^{\infty}$-function on $\real$
such that $\rho(t)=1$ for $t\leq 1/2$ and 
$\rho(t)=0$ for $t\geq 2/3$.
We put $\chi_M(z_1):=
 \rho\bigl(-M^{-1}\log|z_1|\bigr)$
for any $M>1$.
We have only to show that 
there exists a constant $C>0$,
which is independent of 
$Q\in D_1$ and $M>1$,
such that the following holds:
\begin{equation}
 \label{eq;07.7.18.30}
 \int_{\pi_1^{-1}(Q)}
 \bigl(\delta_{\lambda,1}'(\chi_M f),\,
 \delta_{\lambda,1}'(\chi_M f) \bigr)_h
 |z_1|^{2\epsilon}
<C
\end{equation}
In the following argument,
we will ignore the contribution to the Stokes formula
from the integral over $\del\pi_1^{-1}(Q)$,
because they are uniformly dominated.
(Recall that $D$ is assumed to be smooth.)
The left hand side of (\ref{eq;07.7.18.30})
can be rewritten as follows,
up to the contribution from $\del\pi_1^{-1}(Q)$:
\begin{equation}
\label{eq;07.7.18.35}
 \pm\int_{\pi_1^{-1}(Q)}
 \bigl(\chi_M f,\,
 d_{\lambda,1}''\delta_{\lambda,1}'
 (\chi_M f)\bigr)_h\,
 |z_1|^{2\epsilon}
\pm\int_{\pi_1^{-1}(Q)}
 \bigl(\chi_M f,\,
 \delta_{\lambda,1}'(\chi_M f)\bigr)_h
 \,\epsilon\, |z_1|^{2\epsilon}
 \frac{dz_1}{z_1}
\end{equation}
The first term can be rewritten as follows,
up to the contribution of $\del\pi_1^{-1}(Q)$:
\begin{multline}
 \label{eq;07.10.5.35}
\pm\int_{\pi_1^{-1}(Q)}
 \bigl(\chi_M f,\,
 R_1(h)(\chi_M f)\bigr)\cdot |z_1|^{2\epsilon}
\pm\int_{\pi_1^{-1}(Q)}
 \bigl(\chi_M f,\,
 \delta_{\lambda,1}'d_{\lambda,1}''
 (\chi_M f)\bigr)\cdot
 |z_1|^{2\epsilon}\\
=\pm\int_{\pi_1^{-1}(Q)}
 \bigl(\chi_M f,\, R_1(h)(\chi_M f)\bigr)
 \cdot |z_1|^{2\epsilon}
\pm\int_{\pi_1^{-1}(Q)}
 \bigl(d_{\lambda,1}''(\chi_M f),\,
 d_{\lambda,1}''(\chi_M f) \bigr)
 \cdot |z_1|^{2\epsilon}\\
\pm\int_{\pi_1^{-1}(Q)}
 \bigl(\chi_M f,\,d_{\lambda,1}''(\chi_M f)\bigr)
 \cdot\epsilon\cdot |z_1|^{2\epsilon}\frac{dz_1}{z_1}
\end{multline}
The right hand side of (\ref{eq;07.10.5.35})
is uniformly dominated.
The second term in (\ref{eq;07.7.18.35}) is as follows,
up to the contribution of $\del\pi_1^{-1}(Q)$:
\[
\pm \int_{\pi_1^{-1}(Q)}
 \bigl(\chi_M f,\,\chi_M f\bigr)\,
 \epsilon^2\,|z_1|^{2\epsilon}
 \frac{dz_1\, d\zbar_1}{|z_1|^2}
\pm \int_{\pi_1^{-1}(Q)}
 \bigl(d_{\lambda,1}''(\chi_M f),\,
 \chi_M f\bigr)
 \,\epsilon\,|z_1|^{2\epsilon}
 \frac{dz_1}{z_1}
\]
It is bounded uniformly for $Q\in D_1$.
Then, 
it is easy to show the existence of a constant $C$,
which is independent of $N$ and $P$
such that (\ref{eq;07.7.18.30}) holds.
Thus, the proof of Lemma \ref{lem;07.7.18.60}
is finished.
\hfill\qed

\begin{lem}
\label{lem;07.7.18.55}
We put
$\Psi:=\theta-\sum_i\sum_{\gminib}
 d\zeta_{\vecm(i)}(\gminib)
 \cdot p^{\vecm(i)}_{\gminib}$.
Then we have the estimates
$\Psi(\del_j)=O(1)$ ($j\neq 1$)
and $\Psi(\del_1)=O(|z_1|^{-1})$ with respect to $h$.
\end{lem}
\pf
We have the boundedness of
$\del_j\zeta_{\vecm}(\gminib)\cdot
 \bigl(\pi^{\vecm}_{\gminib}
 -p^{\vecm}_{\gminib}\bigr)$
by construction of $p^{\vecm}_{\gminib}$.
According to Theorem \ref{thm;07.10.4.1}
and Theorem \ref{thm;07.10.4.3},
the $dz_j$-components $(j\neq 1)$ of
$\theta-\sum_i\sum_{\gminib}
 d\zeta_{\vecm(i)}(\gminib)\cdot
 \pi^{\vecm(i)}_{\gminib}$
are $O(1)$,
and the $dz_1$-component is 
$O\bigl(|z_1|^{-1}\bigr)$.
Then, the claim of Lemma \ref{lem;07.7.18.55} follows.
\hfill\qed

\vspace{.1in}

Let $j\neq 1$.
According to Lemma \ref{lem;07.7.18.50}
and Lemma \ref{lem;07.7.18.55},
there exists a constant $C_{10}>0$
such that the following holds
for any $Q\in D_j\setminus D_1$:
\[
 \int_{\substack{\pi_j^{-1}(Q)\\
 |z_j|<1/2}}
 \bigl|\DDlambda_0(\del_j)f_{|\pi_j^{-1}(Q)}\bigr|^2
 \cdot|dz_j\cdot d\zbar_j|
<C_{10} \bigl|z_1(Q)\bigr|^{-2\epsilon}
\]

Let $\vecv=(v_p)$ be a frame of $\nbigp_0{\nbigelambda}$
compatible with the parabolic structure.
We put $a(v_p):=\deg^F(v_p)$.
We have the expression:
\[
 \DDlambda_0(\del_j)f=\sum_p A^j_p\cdot v_p
\]
Then, $A^j_p$ are holomorphic functions on $X-D$.
We also have the following:
\[
 \int_{\substack{\pi_j^{-1}(Q)\\
 |z_j|<1/2}}
 \bigl|A^j_{p|\pi_i^{-1}(Q)}\bigr|^2\cdot 
 |z_1(Q)|^{2\epsilon-2a(v_p)}
<C_{11}
\]
By using Fubini's theorem,
we have 
\[
 \int_{\pi_1^{-1}(P)}
 \bigl|A^j_{p|\pi_1^{-1}(P)}\bigr|^2
 \cdot |z_1|^{2\epsilon'-2a(v_p)-2}
 \cdot|dz_1\cdot d\zbar_1|
 <\infty
\]
for almost all $P\in D_1$ such that $|z_j(P)|<1/2$
and for some $0<\epsilon'<\epsilon$.
Hence $A^j_p$ are holomorphic around 
the origin $O$,
and we obtain 
$\DDlambda_0(\del_j)f\in\nbigp_0{\nbigelambda}$.

According to Lemma \ref{lem;07.7.18.60}
and Lemma \ref{lem;07.7.18.55},
there exists a constant $C_{12}>0$ 
such that the following holds
for any $Q\in D_1$:
\[
 \int_{\pi_1^{-1}(Q)}
 \bigl|\DDlambda_0(z_1\del_1)f\bigr|_h^2
\cdot |z_1|^{2\epsilon-2}\cdot |dz_1\cdot d\zbar_1|
<C_{12}
\]
We can deduce
$\DDlambda_0(z_1\del_1)f\in \nbigp_0\nbigelambda$
as in the case $j\neq 1$.
Thus, we obtain that $\DDlambda_0$ is logarithmic,
when $D$ is smooth.

\vspace{.1in}

The general case can easily be reduced to the case
in which $D$ is smooth, due to the theorem of Hartogs.
Thus, the proof of Proposition \ref{prop;07.7.18.70}
is finished.
\hfill\qed

\vspace{.1in}

Now, the first claim of
Theorem \ref{thm;07.11.18.20}
immediately follows from Proposition \ref{prop;07.7.18.70}.
The second and the third claims of 
Theorem \ref{thm;07.11.18.20}
follows from Proposition \ref{prop;07.7.7.23}
and Proposition \ref{prop;07.7.18.70}.
\hfill\qed

\section{Comparison of the irregular decompositions
 in the case $\vecm<\veczero_{\ell}$}
\label{subsection;08.9.14.41}
\subsection{Statements}

We use the setting in 
Section \ref{subsection;07.11.17.1}.
Let $\lambda\neq 0$.
Assume there exists
$\vecm\in\nbigm$ such that
$\vecm<\veczero_{\ell}$.
We have the irregular decomposition on $\Dhat$
as in (\ref{eq;07.11.12.10}):
\begin{equation}
 \label{eq;07.7.20.4}
 \nbigp_{\veca}\nbigelambda_{|\Dhat}
=\bigoplus_{\gminib\in\Irrbar(\theta,\vecm)}
 \nbigp_{\veca}\nbigehatlambda_{\gminib}
\end{equation}

In the following,
$N$ will denote a large integer.
Let $\Dhat^{(N)}$ denote 
the $N$-th infinitesimal neighbourhood of $D$.
Let $\nbigp_{\veca}\nbigelambda=
 \bigoplus_{\gminib\in\Irrbar(\theta,\vecm)}
 \nbigp_{\veca}\nbige^{\lambda}_{\gminib,N}$
be a decomposition
whose restriction to $\Dhat^{(N)}$ is the same as
the restriction of (\ref{eq;07.7.20.4})
to $\Dhat^{(N)}$.
Let $q^{\vecm}_{\gminib,N}$ be the projection
onto $\nbigp_{\veca}\nbige^{\lambda}_{\gminib,N}$
with respect to the decomposition.
We will prove the following proposition
in Sections
\ref{subsection;07.10.5.51}--\ref{subsection;07.11.18.30}.

\begin{prop}
\label{prop;07.10.5.50}
We have the estimate
$\bigl|
 q^{\vecm}_{\gminib,N}-\pi^{\vecm}_{\gminib}
 \bigr|_h
=O\bigl(\prod_{i=1}^{\ell}|z_i|^{N}\bigr)$.
In particular, 
the decomposition
$\nbigp_{\veca}\nbigelambda
=\bigoplus
 \nbigp_{\veca}\nbigelambda_{\gminib,N}$
is $\prod_{i=1}^{\ell}|z_i|^{N}$-asymptotically
orthogonal with respect to $h$.
\end{prop}

\begin{rem}
See Section {\rm\ref{subsection;07.10.9.11}}
for the refinement in the general case.
\hfill\qed
\end{rem}

Before going into the proof,
we give some consequences.
Let $S$ be a small multi-sector in $X-D$,
and let $\Sbar$ denote the closure of $S$ in
the real blow up $\Xtilde(D)$.
As explained in Section 
\ref{subsection;07.11.18.10},
we have the partial Stokes filtration 
$\nbigf^{S,\vecm}$ of
$\nbigp_{\veca}\nbigelambda_{|\Sbar}$
 $(\lambda\neq 0)$
in the level $\vecm$,
indexed by the ordered set
$\bigl(\Irrbar(\theta,\vecm),\leq^{\lambda}_S\bigr)$.
Because $\vecm<\veczero_{\ell}$,
we can take a $\DDlambda$-flat splitting:
\begin{equation}
 \label{eq;07.7.20.10}
 \nbigp_{\veca}{\nbigelambda}_{|\Sbar}=
 \bigoplus_{\gminib\in\Irrbar(\theta,\vecm)}
 \nbigp_{\veca}\nbige^{\lambda}_{\gminib,S}
\end{equation}
Let $p^{\vecm}_{\gminib,S}$ be the projection
onto $\nbigp_{\veca}\nbige^{\lambda}_{\gminib,S}$
with respect to the decomposition (\ref{eq;07.7.20.10}).

\begin{cor}
\label{cor;07.7.7.42}
We have 
$\bigl|p^{\vecm}_{\gminib,S}
 -\pi^{\vecm}_{\gminib}\bigr|_h
\leq
 C_N\prod_{i=1}^{\ell}|z_i|^N$
for any $N>0$.
In particular, 
the decomposition {\rm(\ref{eq;07.7.20.10})}
is $\prod_{i=1}^{\ell}|z_i|^N$-asymptotically orthogonal
for any $N>0$.
\end{cor}
\pf
For any $N$,
we take $q^{\vecm}_{\gminib,N}$ as above.
We have
$\bigl|q^{\vecm}_{\gminib,N}
 -p^{\vecm}_{\gminib,S}\bigr|_{h}
\leq C_N''\prod_{i=1}^{\ell}|z_i|^{N}$.
Hence, the claim follows from 
Proposition \ref{prop;07.10.5.50}.
\hfill\qed

\vspace{.1in}
By varying $S$ and gluing
 $p^{\vecm}_{\gminib,S}$ in $C^{\infty}$
as in Section \ref{subsection;07.6.16.8},
we obtain $p^{\vecm}_{\gminib,C^{\infty}}$.

\begin{cor}
\label{cor;07.7.7.43}
We have 
$\bigl|
 p^{\vecm}_{\gminib,C^{\infty}}-\pi^{\vecm}_{\gminib}
 \bigr|_h\leq
 C_N\,\prod_{i=1}^{\ell}|z_i|^N$
for any $N>0$.
\end{cor}
\pf
It follows from Corollary \ref{cor;07.7.7.42}.
\hfill\qed

\subsection{Estimate of 
 $\del_j\pi^{\vecm}_{\gminia}$
 in the case $\vecm<\veczero_{\ell}$}
\label{subsection;07.10.5.51}

Recall $\nbige^0=E$ as holomorphic bundles.
We use the symbol $E$ in this section.
For a section $f$ of $E$ with 
$\del f=\sum_{j=1}^n A_j\, dz_j$,
we put $\del_jf:=A_j\, dz_j$.
We take a holomorphic frame $\vecv$ of
$\prolong{E}$ compatible with the parabolic structure
and the decompositions 
(\ref{eq;07.7.18.10}) and (\ref{eq;07.7.18.11}).
For $j=1,\ldots,n$,
let $F_j$ be the section of
$\End(E)\otimes\Omega^{1,0}_{X-D}$
determined by $F_j\vecv=\del_j\vecv$.
Then, we have the following estimate
for some $M>0$,
due to Lemma \ref{lem;07.6.2.20} below:
\begin{equation}
 \label{eq;07.6.2.21}
 |F_j|_h=O\left(
 \Bigl(
 \sum_{i=1}^{\ell}\bigl(-\log|z_i|\bigr)
 \Bigr)^M
 \right)
\end{equation}
We have the decomposition:
\[
 \End(E)=
\nbigd^{\vecm}(\End(E))\oplus
\nbigc^{\vecm}(\End(E))
\]
\[
 \nbigd^{\vecm}(\End(E)):=
\bigoplus_{\gminib\in\Irrbar(\theta,\vecm)}
 \End(E^{\vecm}_{\gminib})
\]
\[
 \nbigc^{\vecm}(\End(E)):=
\bigoplus_{
 \substack{\gminib,\gminib'\in\Irrbar(\theta,\vecm)\\
 \gminib\neq\gminib'}}
 \Hom(E^{\vecm}_{\gminib},E^{\vecm}_{\gminib'})
\]
According to Theorem \ref{thm;07.10.4.1},
$\nbigd^{\vecm}\bigl(\End(E)\bigr)$
and $\nbigc^{\vecm}\bigl(\End(E)\bigr)$
are $\exp\bigl(-\epsilon|\vecz^{\vecm}|\bigr)$-asymptotically
orthogonal.
For a section $g$ of $\End(E)\otimes\Omega^{1,0}$,
we have the corresponding decomposition
$g=\nbigd^{\vecm}(g)+\nbigc^{\vecm}(g)$.

\begin{lem}
\label{lem;07.6.2.25}
Around the origin $O\in X$,
we have the following estimates
for some $\epsilon>0$:
\[
 \del_{j}\pi^{\vecm}_{\gminia}
=O\bigl(\exp\bigl(-\epsilon|\vecz^{\vecm}|\bigr)\bigr)
\,\,\bigl(\gminia\in\Irrbar(\theta,\vecm)\bigr),
\quad
\nbigc^{\vecm}(F_j)=
O\Bigl(
 \exp\bigl(-\epsilon|\vecz^{\vecm}|\bigr)
 \Bigr)
\]
\end{lem}
\pf
We have
$ \del_{j}\pi^{\vecm}_{\gminia}
=\bigl[\pi^{\vecm}_{\gminia},F_j\bigr]
 \in \nbigc^{\vecm}(\End(E))\otimes\Omega^{1,0}$.
Due to the estimate (\ref{eq;07.6.2.21}),
$\del_j\pi^{\vecm}_{\gminia}$ is bounded
up to log order, with respect to $h$ and the Poincar\'e metric.
Since $\nbigd^{\vecm}(\End(E))$
and $\nbigc^{\vecm}(\End(E))$ are 
$\exp\bigl(-\epsilon|\vecz^{\vecm}|\bigr)$-asymptotically
orthogonal,
we obtain the following estimate:
\[
 (\pi^{\vecm}_{\gminia},
 \del_{j}\pi^{\vecm}_{\gminia})_h
=O\bigl(
 \exp(-\epsilon|\vecz^{\vecm}|)
 \bigr)
\]
Let $R_j(h)$ denote the $dz_j\,d\zbar_j$-component
of $R(h)$.
We have $\delbar_j\del_j\pi^{\vecm}_{\gminia}
=\bigl[\pi^{\vecm}_{\gminia},R_j(h)\bigr]
=O\bigl(\exp(-\epsilon|\vecz^{\vecm}|)\bigr)$
due to Corollary \ref{cor;07.11.22.30}.
Hence, we obtain the following:
\[
 \bigl(
 \pi^{\vecm}_{\gminia},
 \delbar_j\del_j\pi^{\vecm}_{\gminia}
 \bigr)_h
=O\Bigl(\exp\bigl(-\epsilon|\vecz^{\vecm}|\bigr)\Bigr)
\]
Let $\pi_j$ denote the projection $X-D\lrarr D_j$.
Let us consider the case $j\leq \ell$.
We put $S_j:=\{1,\ldots,\ell\}-\{j\}$.
Let $P$ be any point of $D_j^{\circ}$.
Note the following equality on $\pi_j^{-1}(P)$,
for any $\eta>0$:
\begin{multline}
 \bigl(\del_j\pi^{\vecm}_{\gminia},
 \del_j\pi^{\vecm}_{\gminia}\bigr)_h\,
\exp\Bigl(
 \eta\,|z_j|^{m_j}\,
 \prod_{i\in S_j}|z_i(P)|^{m_i}
 \Bigr)
=\\
\del_j\left(
 \bigl(\pi^{\vecm}_{\gminia},
 \del_j\pi^{\vecm}_{\gminia}\bigr)_h\,
\exp\Bigl(
 \eta\,|z_j|^{m_j}\,
 \prod_{i\in S_j}|z_i(P)|^{m_i}
 \Bigr)\right)
 \\
+\bigl(\pi^{\vecm}_{\gminia},
 \del_j\pi^{\vecm}_{\gminia}\bigr)_h\,
\exp\Bigl(
 \eta\,|z_j|^{m_j}\,
 \prod_{i\in S_j}|z_i(P)|^{m_i}
 \Bigr)
\,\frac{\eta \, m_j}{2}
 |z_j|^{m_j}\frac{dz_j}{z_j} \\
-\bigl(\pi^{\vecm}_{\gminia},
 \delbar_j\del_j\pi^{\vecm}_{\gminia}\bigr)_h\,
\exp\Bigl(
 \eta\,|z_j|^{m_j}\,
 \prod_{i\in S_j}|z_i(P)|^{m_i}
 \Bigr)
\end{multline}
Hence, we obtain the following finiteness
for some $\eta>0$ which is sufficiently smaller than $\epsilon$:
\[
 \int_{\pi_j^{-1}(P)}
 \bigl(\del_{j}\pi^{\vecm}_{\gminia},
 \del_{j}\pi^{\vecm}_{\gminia}\bigr)_h
 \,
 \exp\Bigl(\eta\,|z_j|^{m_j}
 \prod_{i\in S_j}|z_i(P)|^{m_i}
 \Bigr)<C
\]

Let $\varphi_P:\hyperh\lrarr \pi_j^{-1}(P)$
be the covering given by 
$\zeta\longmapsto \exp(2\pi\sqrt{-1}\zeta)$.
We put
\[
 K_n:=\bigl\{(\xi,\eta)\in\hyperh\,\big|\,
 -1<\xi<1,n-1<\eta<n+1\bigr\}.
\]
We have 
$  e^{2\pi m_j}\cdot e^{-2\pi m_j n}
<|z_j|^{m_j}
<e^{-2\pi m_j}\cdot e^{-2\pi m_j n}$
on $K_n$.
(Note $m_j<0$.)
Hence we have the following:
\begin{equation}
 \label{eq;07.6.2.25}
 \int_{K_n}
 \bigl(\del_{j}\pi^{\vecm}_{\gminia},
 \del_j\pi^{\vecm}_{\gminia}\bigr)_h
\leq
 C\,
 \exp\Bigl(
- \eta\, e^{2\pi m_j}\,
 e^{2\pi m_j n}
 \,
 \prod_{i\in S_j} |z_i(P)|^{m_i}
 \Bigr)
\end{equation}
Since $(E,\delbar_E,h)$ is acceptable,
there exists an $n_0$ such that 
we can apply Uhlenbeck's theorem \cite{u1}
to the pull back
$\varphi_P^{-1}(E,\delbar_E,h)_{|K_n}$
for any $n\geq n_0$ and any $P\in D_j$.
Hence,
we can take an orthonormal frame 
$\vece$
of $\varphi_P^{-1}\End(E)$ 
on $K_n$ $(n\geq n_0)$
such that $A$ is small,
where $A$ is determined by
$(\del_j+\delbar_j) \vece=\vece\, A$.
We have the expression
$ \del_j\pi^{\vecm}_{\gminia}
=\sum_i \rho_i \, e_i \, dz_j$.
We put $\vecrho:=(\rho_i)$.
The inequality (\ref{eq;07.6.2.25})
gives the estimate of the $L^2$-norm of
$\vecrho$.
From the estimate of
$\delbar_j\del_j\pi^{\vecm}_{\gminia}$,
we obtain the following:
\begin{equation}
 \label{eq;07.7.20.2}
 \bigl|
 \delbar_j\vecrho+A^{0,1}\vecrho
 \bigr|
\leq
 C_2\,\exp\Bigl(
 -\eta_2\, e^{2\pi n m_j}\,
\prod_{i\in S_j}|z_i(P)|^{m_i}
 \Bigr)
\end{equation}
By using a standard argument,
we obtain the estimate of the sup norm 
of $\vecrho$ on $K_n$ from (\ref{eq;07.7.20.2})
and (\ref{eq;07.6.2.25}).
Then we obtain the following
estimate on $K_n$:
\begin{equation}
 \label{eq;07.7.20.6}
 \bigl|
 \del_{j}\pi^{\vecm}_{\gminia}
 \bigr|_h
\leq 
 C_3\exp\Bigl(
 -\eta_3\, e^{2\pi nm_j}\,
\prod_{i\in S_j} |z_i(P)|^{m_i}
 \Bigr)
\leq
 C_4\exp\Bigl(
 -\eta_4\,|z_j|^{m_j}\,
\prod_{i\in S_j}|z_i(P)|^{m_i}
 \Bigr)
\end{equation}
Thus, we are done in the case $j\leq \ell$.

Let us consider the case $j>\ell$.
We have the following:
\[
  \int_{\pi_j^{-1}(Q)}
 \bigl(\del_j\pi^{\vecm}_{\gminia},
 \del_j\pi^{\vecm}_{\gminia} \bigr)_h
=
 \\
\int_{\substack{\pi_j^{-1}(Q)\\ |z_j|=1}}
 \bigl(\pi^{\vecm}_{\gminia},
 \del_j\pi^{\vecm}_{\gminia}\bigr)_h
-\int_{\pi_j^{-1}(Q)}
 \bigl(\pi^{\vecm}_{\gminia},
 \delbar_j\del_j\pi^{\vecm}_{\gminia}\bigr)_h
\]
Thus, we obtain 
$\|\del_j\pi^{\vecm}_{\gminia}\|_{L^2}
\leq
 C\exp\bigl(-\epsilon\prod_{p=1}^{\ell}|z_p(Q)|^{m_p}\bigr)$.
We can obtain the estimate for the sup norm of
$\del_j\pi^{\vecm}_{\gminia}$
from the $L^2$-estimate of $\del_j\pi^{\vecm}_{\gminia}$
and the estimate for the sup norm of
$\delbar_j\del_j\pi^{\vecm}_{\gminia}$.
(See the argument to obtain (\ref{eq;07.7.20.6}) above.)
Thus, the proof of Lemma \ref{lem;07.6.2.25} is finished.
\hfill\qed

\subsection{Estimate of
 $\DDlambda p^{\vecm}_{\gminia,N}$
 in the case $\vecm<\veczero_{\ell}$}
\label{subsection;07.10.5.52}

Let $p^{\vecm}_{\gminib,N}$ 
$(\gminib\in\Irrbar(\theta,\vecm))$
be as in Lemma \ref{lem;07.7.7.20}.
We have the decomposition:
\begin{equation}
\label{eq;07.7.20.3}
 \nbigp_{\veca}\nbigelambda
=\bigoplus_{\gminib\in\Irrbar(\theta,\vecm)}
 \Image p^{\vecm}_{\gminib,N}
\end{equation}

\begin{lem}
\label{lem;07.7.20.5}
$\DDlambda(p^{\vecm}_{\gminib,N})
=O\Bigl(\prod_{j=1}^{\ell}|z_j|^{N}\Bigr)$.
\end{lem}
\pf
We have
$\del_E\pi^{\vecm}_{\gminib}=
 O\Bigl(\exp\bigl(-\epsilon|\vecz^{\vecm}|\bigr)\Bigr)$
due to Lemma \ref{lem;07.6.2.25}.
We also have 
$[\theta,\pi^{\vecm}_{\gminib}]=0$
and $[\theta^{\dagger},\pi^{\vecm}_{\gminib}]
=O\Bigl(\exp\bigl(-\epsilon|\vecz^{\vecm}|\bigr)\Bigr)$.
Hence, $\DDlambda\pi^{\vecm}_{\gminib}=O\Bigl(
  \exp\bigl(-\epsilon|\vecz^{\vecm}|\bigr)\Bigr)$.

We put 
$s^{\vecm}_{\gminib}:=
\pi^{\vecm}_{\gminib}-p^{\vecm}_{\gminib}$.
We have
$\bigl|s^{\vecm}_{\gminib}\bigr|_h
=O\bigl(\prod_{j=1}^{\ell}|z_j|^{2N}\bigr)$.
We have the following:
\[
 \DDlambda s^{\vecm}_{\gminib}
=d_{\lambda}''\pi^{\vecm}_{\gminib}
+\lambda\delta'_{\lambda} s^{\vecm}_{\gminib}
+(1+|\lambda|^2)\theta s^{\vecm}_{\gminib}
\]
We have
$d_{\lambda}''\pi^{\vecm}_{\gminib}=O\Bigl(
 \exp\bigl(-\epsilon|\vecz^{\vecm}|\bigr)\Bigr)$.
We have 
$\theta s^{\vecm}_{\gminib}
=O(\prod_{i=1}^{\ell}|z_i|^{3N/2})$.
Let us look at 
$\delta'_{\lambda}s^{\vecm}_{\gminib}$.
Let $\pi_j$ denote the projection $X-D\lrarr D_j$.

\begin{lem}
\label{lem;07.10.5.55}
\mbox{{}}
We have the following estimates
independently from 
$j=1,\ldots,n$ and
$Q\in\pi_j^{-1}(D)$.
\begin{itemize}
\item
$\bigl\|
 |z_j|^{-3N/2}
 \delta'_{\lambda,j}
 s^{\vecm}_{\gminib|\pi_j^{-1}(Q)}
 \bigr\|_{L^2}\leq
 C\, \prod_{i\in S_j}\bigl|z_i(Q)\bigr|^{N}$
for $j=1,\ldots,\ell$,
where 
$S_j:=\{1,\ldots,\ell\}-\{j\}$.
\item
 $\bigl\|
 \delta'_{\lambda,j}
 s^{\vecm}_{\gminib|\pi_j^{-1}(Q)}
 \bigr\|_{L^2}\leq
 C\, \prod_{i=1}^{\ell}|z_i(Q)|^N$
 for $j=\ell+1,\ldots,n$.
\end{itemize}
\end{lem}
\pf
Let us consider the case $j\leq \ell$.
We put $L:=3N/2$, for simplicity of the notation.
Let $\rho_1$ be a non-negative $C^{\infty}$-function
on $\real$ such that
$\rho_1(t)=1$ for $t\leq 1/2$ and $\rho_1(t)=0$
for $t\geq 2/3$.
Let $\kappa$ be a non-negative $C^{\infty}$-function
on $\Delta$ such that
$\kappa(z)=1$ for $|z|\leq 1/2$
and $\kappa(z)=0$ for $|z|\geq 2/3$.
We put $\chi_M(z_j):=
 \rho_1\bigl(-M^{-1}\log|z_j|\bigr)\, \kappa(z_j)$.
In the following estimate,
we do not have to be careful
on the signature.
We have the following:
\begin{multline}
  \int_{\pi_j^{-1}(Q)}
 \bigl(\delta'_{\lambda,j}(\chi_M s^{\vecm}_{\gminib}),
 \delta'_{\lambda,j}(\chi_M s^{\vecm}_{\gminib})\bigr)_h
 \, |z_j|^{-2L} 
=\\
\pm\int_{\pi_j^{-1}(Q)}
 \bigl(
 \chi_Ms^{\vecm}_{\gminib},
 d''_{\lambda,j}\delta_{\lambda,j}'
 (\chi_M s^{\vecm}_{\gminib})
 \bigr)_h\, |z_j|^{-2L}
 \\
\pm\int_{\pi_j^{-1}(Q)}
 \bigl(\chi_Ms^{\vecm}_{\gminib},
 \delta_{\lambda,j}'(\chi_Ms^{\vecm}_{\gminib})\bigr)_h
 \, L\, |z_j|^{-2L}\frac{dz_j}{z_j}
\end{multline}
The first term of the right hand side
is as follows:
\begin{multline*}
\pm\int_{\pi_j^{-1}(Q)}
 \!\!\!
 \bigl(\chi s^{\vecm}_{\gminib},
 R(h,d_{\lambda}'')(\chi s^{\vecm}_{\gminib})\bigr)_h
 |z_j|^{-2L}
\pm\int_{\pi_j^{-1}(Q)}
 \!\!\!
 \bigl(\chi s^{\vecm}_{\gminib},
 \delta'_{\lambda,j}d''_{\lambda,j}
 (\chi s^{\vecm}_{\gminib})
 \bigr)_h |z_j|^{-2L} \\
=
\pm\int_{\pi_j^{-1}(Q)}\!\!\!
 \bigl(\chi s^{\vecm}_{\gminib},
 R(h,d_{\lambda}'')(\chi s^{\vecm}_{\gminib})\bigr)_h
 |z_j|^{-2L}
\pm\int_{\pi_j^{-1}(Q)}\!\!\!
 \bigl(d''_{\lambda,j}(\chi s^{\vecm}_{\gminib}),
 d''_{\lambda,j}(\chi s^{\vecm}_{\gminib}) \bigr)
 |z_j|^{-2L} \\
\pm\int\bigl(\chi s^{\vecm}_{\gminib},
 d''_{\lambda,j}(\chi s^{\vecm}_{\gminib}) \bigr)_h
\, L|z_j|^{-2L}\frac{d\zbar_j}{\zbar_j}
\end{multline*}
($M$ is omitted.) 
The second term is as follows:
\begin{multline*}
\pm\int_{\pi_j^{-1}(Q)}
 \bigl(d_{\lambda}''(\chi_Ms^{\vecm}_{\gminib}),
 \chi_M s^{\vecm}_{\gminib} \bigr)_h\, 
 L\, |z_j|^{-2L}\frac{dz_j}{z_j} \\
\pm\int_{\pi_j^{-1}(Q)}
 \bigl(\chi_M s^{\vecm}_{\gminib},
 \chi_M s^{\vecm}_{\gminib}\bigr)\, L^2
 |z_j|^{-2L}\frac{d\zbar_j\, dz_j}{|z_j|^2}
\end{multline*}
Then, the claim in the case $j\leq \ell$ follows
from the estimate for
$|s^{\vecm}_{\gminib}|_h$ and 
$|d_{\lambda}''s^{\vecm}_{\gminib}|_h$
in the limit $M\to\infty$.

Let us consider the case $j>\ell$.
We have the following:
\begin{multline*}
 \int_{\pi_j^{-1}(Q)}
 \bigl(\delta'_{\lambda,j}(\kappa s^{\vecm}_{\gminib}),
 \delta'_{\lambda,j}(\kappa s^{\vecm}_{\gminib}) \bigr)
= \\
\pm\int_{\pi_j^{-1}(Q)}
 \bigl(d''_{\lambda,j}(\kappa s^{\vecm}_{\gminib}),
 d''_{\lambda,j}(\kappa s^{\vecm}_{\gminib})
 \bigr)
\pm\int_{\pi_j^{-1}(Q)}
 \bigl(\kappa s^{\vecm}_{\gminib},
 R_j(d''_{\lambda},h)\kappa s^{\vecm}_{\gminib}
 \bigr)
\end{multline*}
Hence, the claim in the case $j>\ell$
also follows from the estimate of
$|s^{\vecm}_{\gminib}|_h$ and
$|d''_{\lambda}s^{\vecm}_{\gminib}|_h$.
Thus, the proof of Lemma \ref{lem;07.10.5.55}
is finished.
\hfill\qed

\vspace{.1in}
Let us finish the proof of Lemma \ref{lem;07.7.20.5}.
From the estimate of $\DDlambda\pi^{\vecm}_{\gminib}$
and $\DDlambda s^{\vecm}_{\gminib}$,
we obtain the following:
\[
 \bigl\|
 \bigl(
 |z_j|^{-L}\,
 \DDlambda_jp^{\vecm}_{\gminia,N}
\bigr)_{|\pi_j^{-1}(Q)}
 \bigr\|_{L^2}
\leq C \prod_{i\in S_j}|z_i(Q)|^{4N/3},
\quad
 (j=1,\ldots,\ell)
\]
\[
 \bigl\|
 \DDlambda_j p^{\vecm}_{\gminia,N|\pi_j^{-1}(Q)}
 \bigr\|_{L^2}
\leq
 C\prod_{i=1}^{\ell}|z_i(Q)|^{4N/3},
\quad (j=\ell+1,\ldots,n)
\]
Then, by the holomorphic property,
we obtain 
$\DDlambda p^{\vecm}_{\gminia,N}
=O\bigl(\prod_{j=1}^{\ell}|z_j|^N\bigr)$.
Thus the proof of Lemma \ref{lem;07.7.20.5} is accomplished.
\hfill\qed

\subsection{End of Proof of Proposition
 \ref{prop;07.10.5.50}}
\label{subsection;07.11.18.30}

\begin{lem}
\label{lem;08.9.13.21}
The restrictions of 
the decompositions {\rm(\ref{eq;07.7.20.4})} 
and {\rm(\ref{eq;07.7.20.3})}
to $\Dhat^{(N)}$
are the same.
\end{lem}
\pf
Let $\vecv_{N}=(\vecv_{\gminia,N})$
and $\vecw_N=(\vecw_{\gminia,N})$ 
be frames of 
$\nbigp_{\veca}\nbigelambda$
whose restrictions to $\Dhat^{(N)}$
are compatible with the decompositions
{\rm(\ref{eq;07.7.20.4})} and
{\rm(\ref{eq;07.7.20.3})},
respectively.
Namely,
$\vecv_{\gminia,N|\Dhat^{(N)}}$
and $\vecw_{\gminia,N|\Dhat^{(N)}}$
give frames of
$\nbigp_{\veca}
 \nbigehatlambda_{\gminia|\Dhat^{(N)}}$ and 
$\bigl(\Image 
 p^{\vecm}_{\gminia,N}\bigr)_{|\Dhat^{(N)}}$,
respectively.
Let $A$ and $B$ be determined by the following:
\[
 \vecz^{-\vecm(0)}\DDlambda\vecv_N
 =\vecv_N\, A,
\quad
 \vecz^{-\vecm(0)}\DDlambda\vecw_N
=\vecw_N\, B
\]
We have the decompositions
$A=(A_{\gminia,\gminib})$
and $B=(B_{\gminia,\gminib})$,
corresponding to the decompositions of the frames.
By our choice,
we have 
$A_{\gminia,\gminib}
\equiv
 B_{\gminia,\gminib}
\equiv 0$
for $\gminia\neq\gminib$
modulo 
$\vecz^{-\vecm(0)}\prod_{j=1}^{\ell}z_j^N$.
Let $C$ be determined by
$\vecv=\vecw\, C$,
which has the decomposition
$C=(C_{\gminia,\gminib})$.
We obtain the following
modulo
$\vecz^{-\vecm(0)}\prod_{j=1}^{\ell}z_j^N$:
\begin{equation}
 \label{eq;08.9.13.20}
B_{\gminia,\gminia}\, C_{\gminia,\gminib}
-C_{\gminia,\gminib}\, A_{\gminib,\gminib}
+\vecz^{-\vecm(0)}C_{\gminia,\gminib}\equiv 0
\end{equation}
Assume $C_{\gminia,\gminib}\not\equiv 0$
modulo $\prod_{j=1}^{\ell}z_j^{N}$.
We have the expansion
$C_{\gminia,\gminib}
=\sum_{\vecn\in\seisuu_{\geq\,0}^{\ell}}
 C_{\gminia,\gminib;\vecn}\vecz^{\vecn}$.
We set $\vecdelta:=(1,\ldots,1)$.
Let $\vecn_0\not\geq N\,\vecdelta$ 
be a minimal among $\vecn$ such that 
$C_{\gminia,\gminib;\vecn}\neq 0$.
Then, we obtain
$(\gminia-\gminib)_{\ord(\gminia-\gminib)}
\, C_{\gminia,\gminib;\vecn_0}=0$
from (\ref{eq;08.9.13.20}).
Note 
$-\vecm(0)+\ord(\gminia-\gminib)+\vecn_0
\not\geq -\vecm(0)+N\,\vecdelta$.
Hence, we obtain
$C_{\gminia,\gminib;\vecn_0}=0$
which contradicts with our choice of $\vecn_0$.
Hence, we obtain $C_{\gminia,\gminib}\equiv 0$
modulo $\prod_{j=1}^{\ell}z_j^N$.
It implies the claim of Lemma 
\ref{lem;08.9.13.21}.
\hfill\qed

\vspace{.1in}
Since we have
$q^{\vecm}_{\gminib,N}
-p^{\vecm}_{\gminib,N}\equiv 0$
modulo $\prod_{i=1}^{\ell}z_i^N$,
the claim of Proposition \ref{prop;07.10.5.50}
is obtained.
\hfill\qed

\subsection{Complement}
\label{subsection;07.10.9.11}

We use the setting in Section \ref{subsection;07.11.17.1}.
Moreover, we assume 
that the coordinate system is admissible
for the good set $\Irr(\theta)$,
for simplicity.
Let $k$ be the integer determined by
$\vecm(0)\in\seisuu_{<0}^k\times\veczero_{\ell-k}$.
Let $\vecm\in\nbigm$  such that
$\vecm\in\seisuu_{<0}^k\times\veczero_{\ell-k}$.
We put $D(\leq k):=\bigcup_{j=1}^kD_j$.
\index{hypersurface $D(\leq k)$}
We have the irregular decomposition 
as in (\ref{eq;07.11.12.10}):
\[
 \nbigp_{\veca}{\nbigelambda}_{|\Dhat(\leq k)}
=\bigoplus_{\gminib\in\Irrbar(\theta,\vecm)}
 (\nbigp_{\veca}{\nbigehatlambda})^{\vecm}_{\gminib}
\]
The projection onto 
$(\nbigp_{\veca}\nbigehatlambda)^{\vecm}_{\gminib}$
is denoted by $\phat^{\vecm}_{\gminib}$.

By Lemma \ref{lem;10.6.30.3},
we can take endomorphisms
$p^{\vecm}_{\gminib,N}\in
\nbigp_0\End(\nbigelambda)$
such that
(i)
$p^{\vecm}_{\gminib,N|\Dhat^{(N)}(\leq k)}
 =
 \phat^{\vecm}_{\gminib|\Dhat^{(N)}(\leq k)}$,
(ii)
 $\bigl[p^{\vecm}_{\gminib,N|D_i},
 \Res_i(\DDlambda)\bigr]=0$ for $i=k+1,\ldots,\ell$.
According to the norm estimate of tame harmonic bundles
\cite{mochi2},
$\bigl|p^{\vecm}_{\gminib,N}\bigr|_h$
is bounded on
$\prod_{i=1}^k\{1/2\leq |z_i|\leq 2/3\}
\times
 (\Delta^{\ast})^{\ell-k}
\times\Delta^{n-\ell}$.
(Otherwise, we can use 
Lemma \ref{lem;07.11.5.2} below.)

Recall that
we have the projection $\pi^{\vecm}_{\gminib}$
in the Higgs side,
as given in Section \ref{subsection;07.11.17.20}.

\begin{lem}
\label{lem;07.7.7.41}
$p^{\vecm}_{\gminib,N}-\pi^{\vecm}_{\gminib}
=O\bigl(\prod_{i=1}^k|z_i|^N\bigr)$.
\end{lem}
\pf
Let $q:X-D\lrarr (\Delta^{\ast})^{\ell-k}\times \Delta^{n-\ell}$
be given by
$q(z_1,\ldots,z_n)=(z_{k+1},\ldots,z_n)$.
Let $d_{\lambda}'':=\delbar_E+\lambda\theta^{\dagger}$.
The restriction to $q^{-1}(Q)$ is also denoted by 
the same notation,
where $Q\in(\Delta^{\ast})^{\ell-k}\times\Delta^{n-\ell}$.
By using Proposition \ref{prop;07.11.22.11},
we obtain the following estimate on $q^{-1}(Q)$,
which is uniform for 
$Q\in(\Delta^{\ast})^{\ell-k}\times\Delta^{n-\ell}$:
\[
 d''_{\lambda}
 \bigl(
 p^{\vecm}_{\gminib,N|q^{-1}(Q)}
-\pi^{\vecm}_{\gminib|q^{-1}(Q)}
 \bigr)
=-d''_{\lambda}
 \bigl(
 \pi^{\vecm}_{\gminib|q^{-1}(Q)}
 \bigr)
=O\Bigl(
 \exp\bigl(-\epsilon|\vecz^{\vecm}|\bigr)
 \Bigr)
\]
The restrictions
$\bigl(\nbigelambda,h\bigr)_{|q^{-1}(Q)}$ 
are acceptable,
and the curvatures are uniformly bounded
for $Q\in (\Delta^{\ast})^{\ell-k}\times\Delta^{n-\ell}$.
By using Lemma \ref{lem;07.9.23.36}
and Lemma \ref{lem;07.9.24.5},
we can take a section $t_Q$ 
of $\End(\nbigelambda)_{|q^{-1}(Q)}$
such that
(i) $d_{\lambda}''t_Q=
 d_{\lambda}''\bigl(p^{\vecm}_{\gminib,N|q^{-1}(Q)}
 -\pi^{\vecm}_{\gminib|q^{-1}(Q)}\bigr)$,
(ii) $|t_Q|_h\leq C_1\prod_{i=1}^k|z_i|^N$,
where the constant $C_1$ is independent of $Q$.

By construction,
we have
$d_{\lambda}''\bigl(
 p^{\vecm}_{\gminib,N|q^{-1}(Q)}
-\pi^{\vecm}_{\gminib|q^{-1}(Q)}-t_Q\bigr)=0$
and
\[
\bigl|
p^{\vecm}_{\gminib,N|q^{-1}(Q)}
-\pi^{\vecm}_{\gminib|q^{-1}(Q)}-t_Q
 \bigr|_h\leq C_2 
\]
on $\prod_{i=1}^k\{1/2\leq |z_i|\leq 2/3\}$
independently from $Q$.
From the estimate for $t_Q$
and Proposition \ref{prop;07.10.5.50},
we have the estimate
\[
\bigl|p^{\vecm}_{\gminib,N|q^{-1}(Q)}
-\pi^{\vecm}_{\gminib|q^{-1}(Q)}-t_Q
 \bigr|_h\leq C_{Q,N} \prod_{i=1}^k|z_i|^N 
\]
depending on $Q$.
Then, by Proposition \ref{prop;07.7.7.3},
we obtain the estimate
\[
\bigl|p^{\vecm}_{\gminib,N|q^{-1}(Q)}
-\pi^{\vecm}_{\gminib|q^{-1}(Q)}-t_Q\bigr|_h
\leq C_N\prod_{i=1}^k|z_i|^N 
\]
independently from $Q$.
Thus, we are done.
\hfill\qed

\section{Small deformation of $\nbigp\nbigelambda$
 in the smooth divisor case}
\label{subsection;07.10.7.21}
\subsection{Characterization of the deformation
 $(\nbigp_{a}\nbigelambda)^{(T)}$ via the metric}
\label{subsection;07.12.21.50}

We use the setting in Section
\ref{subsection;07.11.17.1}
with $\ell=1$,
i.e., $D$ is the smooth divisor $D_1$.
We use a slightly different notation.
We use the symbol $j$
instead of $\vecm(i)$ to denote 
an element of $\nbigm$.
Let $m(0)$  denote the minimum
of the numbers $j$ such that
$\Irr(\theta,j)\neq \{0\}$.
For any $\gminib\in\Irr(\theta,j)$,
we put
\[
E^{(j)}_{\gminib}
:=\bigoplus_{\substack{\gminia\in\Irr(\theta)\\
 \eta_j(\gminia)=\gminib }}
 E_{\gminia}
\]
Let $\pi^{(j)}_{\gminib}$ denote the projection
of $E$ onto $E^{(j)}_{\gminib}$.
Recall 
we have the following estimate 
with respect to $h$ for some $\epsilon_1>0$,
due to Theorems \ref{thm;07.10.4.1},
\ref{thm;07.10.4.3}
and Lemma \ref{lem;07.6.2.25}:
\begin{equation}
\label{eq;07.12.21.5}
 \DDlambda\pi^{(j)}_{\gminib}=
  O\bigl(\exp(-\epsilon_1|z_1|^{j})\bigr)
\end{equation}
For $w\in\cnum$,
we consider the following:
\[
 g_{\irr}(w):=
 \sum_{\gminia\in\Irr(\theta)}
 \exp(w\cdot \overline{\gminia})\cdot\pi_{\gminia}
\]
Let $g_{\irr}(w)^{\ast}h$ be the metric
given by $g_{\irr}(w)^{\ast}h(u,v)=
 h\bigl(g_{\irr}(w)u,g_{\irr}(w)v\bigr)$.
We put
\begin{equation}
 \label{eq;07.11.23.25}
T_1(w):=1-\frac{\lambda\cdot\wbar}{1+|\lambda|^2}
\end{equation}

For any $j$ and any $\gminia\in\Irr(\theta,j-1)$,
we put $\nbigi^{(j)}_{\gminia}:=
 \eta_{j,j-1}^{-1}(\gminia)$.
Formally,
we put $\nbigi^{(m(0))}_{0}:=
 \Irr(\theta,m(0))$.
Let $\pi:\Xtilde(D)\lrarr X$ denote the real blow up of $X$
along $D$.
Fix $\lambda\neq 0$.
We regard
$X$ and $D$ as $\{\lambda\}\times X$
and $\{\lambda\}\times D$.
For any $j$ and any distinct
$\gminib_1,\gminib_2\in 
 \eta_{j,j-1}^{-1}(\gminia)$,
we put 
$F_{\gminib_1,\gminib_2}:=
 -|z_1^{-j}|
 \Re\bigl(\lambda^{-1}(\gminib_1-\gminib_2)\bigr)$,
which are $C^{\infty}$-functions on $\Xtilde(D)$.
Let $S_1,\ldots,S_N$ be small multi-sectors
of $X-D$,
such that (i) the union of their interior parts
is $X-D$,
(ii) $S_i\in\Multisector(X-D,\nbigi^{(j)}_{\gminia})$
 for any $j$ and 
 any $\gminia\in\Irr(\theta,j-1)$.
Let $\Sbar_i$ denote the closure of $S_i$
in $\Xtilde(D)$,
and let $Z_i$ denote $\Sbar_i\cap \pi^{-1}(D)$.
We may assume the following
for each $S_i$, each $j$,
each $\gminia\in\Irr(\theta,j-1)$
and each distinct pair $(\gminib_1,\gminib_2)$
in $\nbigi^{(j)}_{\gminia}$:
\begin{description}
\item[(A1)]
If the intersection
$Z_i\cap\{F_{\gminib_1,\gminib_2}=0\}$
is not empty,
$F_{\gminib_1,\gminib_2}$
is monotone with respect to $\arg(z_1)$.
\end{description}
If we choose sufficiently small
$0<\epsilon_2\leq \epsilon_1$,
either one of the following holds 
for each $S_i$, each $j$,
each $\gminia\in\Irr(\theta,j-1)$
and each distinct pair $(\gminib_1,\gminib_2)$
in $\nbigi^{(j)}_{\gminia}$:
\begin{description}
\item[(A2)]
If the intersection
$Z_i\cap\{F_{\gminib_1,\gminib_2}=0\}$
is empty,
 $\bigl|F_{\gminib_1,\gminib_2}\bigr|
 \geq \epsilon_2/2$
holds on $Z_i$.
\end{description}

We will prove the following proposition
in Sections
\ref{subsection;08.3.5.5}--\ref{subsection;08.3.5.6}.

\begin{prop}
\label{prop;07.7.20.230}
Assume the following:
\begin{itemize}
\item $|w|<\epsilon_2/100$.
\item $\bigl|\arg(T_1(w))\bigr|$ is small
 such that 
 the natural bijection
\[
  \bigl\{
\bigl(
t+(1-t)T_1(w)
\bigr)
 \cdot\gminia\,\big|\,
 \gminia\in\Irr(\theta)
 \bigr\}\lrarr \Irr(\theta) 
\]
 preserves the orders $\leq_{S_i}^{\lambda}$
 $(i=1,\ldots,N)$
 for any $0\leq t\leq 1$.
\end{itemize}
Then, the following holds:
\begin{itemize}
\item
Let $f$ be a holomorphic section of
$\nbigelambda$ on $X-D$.
It gives a section of the sheaf
$(\nbigp_a\nbigelambda)^{(T_1(w))}$,
if and only if it satisfies the following growth condition:
\[
 |f|_{g_{\irr}(w)^{\ast}h}
 =O\bigl(|z_1|^{-a-\epsilon}\bigr)
\,\,\,(\forall\epsilon>0)
\]
\end{itemize}
\end{prop}

\subsection{Comparison of the irregular decompositions
 on small multi-sectors}
\label{subsection;08.3.5.5}

Let $\pi:\Xtilde(D)\lrarr X$ denote the real blow up of
$X$ along $D$.

\begin{prop}
\label{prop;07.12.21.1}
For any point $P\in\pi^{-1}(D)$,
there exist a multi-sector
$S_P\in\Multisector(P,X-D)$
and a $\DDlambda$-flat splitting
$\nbigp_a\nbigelambda_{|\Sbar_P}
=\bigoplus_{\gminia\in\Irr(\theta,j)} 
 \nbigp_a\nbige^{\lambda\natural}_{\gminia,S_P}$
of the Stokes filtration $\nbigf^{S_P\,(j)}$ on $\Sbar_P$
in the level $j$
with the following property:
\begin{itemize}
\item
Let $p^{(j)\natural}_{\gminia,S_P}$ 
denote the projection onto 
$\nbigp_a\nbige^{\lambda\natural}_{\gminia,S_P}$
with respect to the decomposition.
Then, 
$\pi^{(j)}_{\gminia}-p^{(j)\natural}_{\gminia,S_P}
=O\bigl(\exp(-\epsilon_1|z_1^{j}|/2)\bigr)$
with respect to $h$.
\end{itemize}
In particular,
the decomposition 
$\nbigp_a\nbigelambda_{|\Sbar_P}
=\bigoplus_{\gminia\in\Irr(\theta,j)} 
 \nbigp_a\nbige^{\lambda\natural}_{\gminia,S_P}$
is $O\bigl(\exp(-\epsilon_1|z_1^j|/2)\bigr)$-asymptotically
orthogonal.
\end{prop}
\pf
By using Lemma \ref{lem;07.7.19.102},
we can take $\nbigq^{(j)}_{\gminib,S}$
such that
$\DDlambda\nbigq^{(j)}_{\gminib,S}
=\DDlambda\pi^{(j)}_{\gminib}$
and 
$\nbigq^{(j)}_{\gminib,S}=O\bigl(
 \exp\bigl(-\epsilon_1|z_1|^{j}/2\bigr)\bigr)$.
We put
$p^{(j)\natural}_{\gminib,S}:=
 \pi^{(j)}_{\gminib}-\nbigq^{(j)}_{\gminib,S}$,
which is $\DDlambda$-flat.
By applying the modification as in the proof of
Lemma \ref{lem;07.7.7.20},
we may and will assume that
$\bigl[p^{(j)\natural}_{\gminib,S},\,
 p^{(j)\natural}_{\gminic,S}\bigr]=0$
and $p^{(j)\natural}_{\gminib,S}
 \circ p^{(j)\natural}_{\gminib,S}
 =p^{(j)\natural}_{\gminib,S}$.

We put
$\nbigf^{S\,(j)\natural}_{\gminib}
 :=\bigoplus_{\gminic\leq^{\lambda}_S \gminib}
 \Image p^{(j)\natural}_{\gminic,S}$.
Let us compare 
the filtrations $\nbigf^{S\,(j)\natural}$
and $\nbigf^{S\,(j)}$.
We take a splitting 
$\nbigp_{a}\nbigelambda_{|\Sbar}
=\bigoplus \nbigp_{a}\nbigelambda_{\gminia,S}$,
and let $p_{\gminia,S}^{(j)}$ denote 
the projection onto 
$\nbigp_a\nbigelambda_{\gminia,S}$
with respect to the decomposition.
Since we have already known
that $p_{\gminia,S}^{(j)}-\pi^{(j)}_{\gminia}
=O(|z_1|^N)$ for any $N$
(Corollary \ref{cor;07.7.7.42}),
we have 
$p^{(j)\natural}_{\gminia,S}-p^{(j)}_{\gminia,S}
=O\bigl(|z_1|^N\bigr)$ for any $N>0$.
Hence, both
$\nbigf^{S\,(j)}_{\gminib|\Zhat}$
and $\nbigf^{S\,(j)\natural}_{\gminib|\Zhat}$
are the same as  $\nbigf^{Z\,(j)}_{|\gminib}$.
Hence, we obtain 
$\nbigf^{S\,(j)}_{\gminib}=
\nbigf^{S\,(j)\,\natural}_{\gminib}$.
(Use the uniqueness in Proposition \ref{prop;10.5.28.10}
successively.)
In other words,
the decomposition
$\bigoplus \Image p^{(j)\natural}_{\gminib,S}$
gives a splitting of the filtration
$\nbigf^{S\,(j)}$.
Thus, Proposition \ref{prop;07.12.21.1}
is proved.
\hfill\qed

\vspace{.1in}
We have the following corollaries.
\begin{cor}
Let $S$ be a small multi-sector in $X-D$,
and let 
 $\nbigp_{a}\nbigelambda_{|\Sbar}
=\bigoplus_{\gminic\in\Irr(\theta,j)}
 \nbigp_a\nbigelambda_{\gminic,S}$
be any $\DDlambda$-flat splitting 
of the filtration $\nbigf^{S\,(j)}$.
Let $p^{(j)}_{\gminic,S}$ denote the projection
onto $\nbigp\nbigelambda_{\gminic,S}$.
Then, we have
$\pi^{(j)}_{\gminic}-p^{(j)}_{\gminic,S}
=O\bigl(\exp(-\epsilon'|z_1|^{j})\bigr)$
for some $\epsilon'>0$.
\end{cor}
\pf
Take a finite covering of $S$
by multi-sectors $S_P$
as in Proposition \ref{prop;07.12.21.1},
and compare
$p^{(j)}_{\gminic,S}$,
$p^{(j)\natural}_{\gminic,S_P}$
and $\pi^{(j)}_{\gminic}$
on each $S_P$.
\hfill\qed

\subsection{Comparison of the irregular decompositions
 on $S_i$}

Let $\nbigp\nbigelambda_{|\Sbar_i}
=\bigoplus_{\gminia\in\Irr(\theta)}
 \nbigp\nbigelambda_{\gminia,S_i}$
be a $\DD$-flat splitting of
the full Stokes filtration $\nbigftilde^{S_i}$.
For any $\gminib\in\Irr(\theta,j)$,
we put 
\[
 \nbigp\nbige^{\lambda\,(j)}_{\gminib,S_i}
=\bigoplus_{
 \substack{\gminia\in\Irr(\theta)\\
 \eta_j(\gminia)=\gminib}}
 \nbigp\nbige^{\lambda}_{\gminia,S_i}.
\]
Let $p^{(j)}_{\gminib,S_i}$ denote the projection
onto $\nbigp\nbige^{\lambda\,(j)}_{\gminib,S_i}$.

\begin{prop}
\label{prop;08.1.8.15}
We have the following estimate
with respect to $h$:
\[
 p^{(j)}_{\gminib,S_i}
-\pi^{(j)}_{\gminib}
=O\Bigl(
 \exp\bigl(-\epsilon_2|z_1^j|/10\bigr)
 \Bigr)
\]
\end{prop}
\pf
For simplicity of the description,
we set $S:=S_i$, $\Sbar:=\Sbar_i$ and $Z:=Z_i$.
We have the decomposition:
\[
 \End(\nbigp\nbigelambda_{|\Sbar})
=\bigoplus_{(\gminib_1,\gminib_2)\in\Irr(\theta,j)^2}
 \Hom\bigl(
 \nbigp\nbige^{\lambda\,(j)}_{\gminib_1,S},
 \nbigp\nbige^{\lambda\,(j)}_{\gminib_2,S}
 \bigr)
\]
For any point $P\in Z$,
we take a small multi-sector $S_P$ as in 
Proposition \ref{prop;07.12.21.1}.
Let $\Sbar_P$ denote the closure of $S_P$
in $\Xtilde(D)$,
and let $Z_P:=\Sbar_P\cap\pi^{-1}(D)$.
If $S_P$ is sufficiently small,
we may assume that 
one of the following holds on $Z_P$,
for each distinct
$\gminib_1,\gminib_2\in\nbigi^{(j)}_{\gminia}$
$(\gminia\in\Irr(\theta,j))$:
\begin{itemize}
\item
 $F_{\gminib_1,\gminib_2}\leq -\epsilon_2/8$.
\item
 $F_{\gminib_1,\gminib_2}\geq -\epsilon_2/4$.
\end{itemize}
We have the decomposition:
\[
 p^{(j)\natural}_{\gminib,S_P}
-p^{(j)}_{\gminib,S|S_P}
=\sum_{
 \substack{(\gminic_1,\gminic_2)\in\Irr(\theta,j)^2\\
 \gminic_1<_{S_P}\gminic_2}}
 \bigl(
 p^{(j)\natural}_{\gminib,S_P}
 \bigr)_{\gminic_1,\gminic_2}
\in
\bigoplus_{
 \substack{(\gminic_1,\gminic_2)\in\Irr(\theta,j)^2\\
 \gminic_1<_{S_P}\gminic_2}}
 \Hom\bigl(
 \nbigp\nbige^{\lambda\,(j)}_{\gminic_1,S},
 \nbigp\nbige^{\lambda\,(j)}_{\gminic_2,S}
 \bigr)_{|S_P}
\]

\begin{lem}
\label{lem;08.1.8.17}
In the case $j=m(0)$,
we have 
\[
 p^{(m(0))\natural}_{\gminib,S_P}
-p^{(m(0))}_{\gminib,S|S_P}
=O\Bigl(
 \exp\bigl(-\epsilon_2|z_1^{m(0)}|/10\bigr)
 \Bigr)
\]
with respect to $h$
for any $\gminib\in\Irr(\theta,m(0))$.
\end{lem}
\pf
We have only to have estimates of
$\bigl(
 p^{(m(0))\natural}_{\gminib,S_P}
 -p^{(m(0))}_{\gminib,S|S_P}
\bigr)_{\gminic_1,\gminic_2}$
for $\gminic_1<_{S_P}\gminic_2$.
In the case 
$F_{\gminic_1,\gminic_2}
 \leq-\epsilon_2/8$ on $S_P$,
we have the following by flatness:
\[
 \bigl(
 p^{(m(0))\natural}_{\gminib,S_P}
 -p^{(m(0))}_{\gminib,S|S_P}
\bigr)_{\gminic_1,\gminic_2}
=O\Bigl(
 \exp\bigl(-\epsilon_2|z_1^{m(0)}|/10\bigr)
 \Bigr).
\]
In the case $F_{\gminic_1,\gminic_2}\geq -\epsilon_2/4$
on $S_P$,
(A2) cannot happen
for the $\gminic_1,\gminic_2$ and $S=S_i$.
Hence,
we can take a sequence
$P_1=P,P_2,\ldots,P_t\in Z$ with the following property:
\begin{itemize}
\item
 $Z_{P_i}\cap Z_{P_{i+1}}\neq\emptyset$.
\item
 $Z_{P_t}$ intersects with
 $\{F_{\gminic_1,\gminic_2}>0\}$.
\item
 $F_{\gminic_1,\gminic_2}\geq -\epsilon_2/4$
  on each $Z_{P_i}$.
\end{itemize}
By the second condition,
we have
$\bigl(
  p^{(m(0))\natural}_{\gminib,S_{P_t}}
 -p^{(m(0))}_{\gminib,S|S_{P_t}}
\bigr)_{\gminic_1,\gminic_2}=0$ 
on $S_{P_t}$.
On $S_{P_i}\cap S_{P_{i+1}}$,
we have
\begin{multline}
 \label{eq;08.1.8.10}
 p^{(m(0))\natural}_{\gminib,S_{P_i}}
-p^{(m(0))\natural}_{\gminib,S_{P_{i+1}}}
=
\bigl(
 p^{(m(0))\natural}_{\gminib,S_{P_i}}
-\pi^{(m(0))}_{\gminib}
\bigr)
-\bigl(
 p^{(m(0))\natural}_{\gminib,S_{P_{i+1}}}
-\pi^{(m(0))}_{\gminib}
\bigr) \\
=O\Bigl(
 \exp\bigl(-\epsilon_1|z_1^{m(0)}|/2\bigr)
 \Bigr).
\end{multline}
If $p^{(m(0))\natural}_{\gminib,S_{P_i}}
-p^{(m(0))\natural}_{\gminib,S_{P_{i+1}}}$ is not $0$,
we have the following for some $B,C,N>0$,
because of Lemma \ref{lem;07.7.10.1}:
\[
\Bigl|
\bigl(
 p^{(m(0))\natural}_{\gminib,S_{P_i}}
-p^{(m(0))\natural}_{\gminib,S_{P_{i+1}}}
\bigr)_{\gminic_1,\gminic_2}
 \exp\bigl((\gminic_1-\gminic_2)/\lambda\bigr)
\Bigr|
\geq
 B\cdot
 \exp\bigl(-C|z_1^{m(0)+1}|\bigr)
 \cdot|z_1|^{N}
\]
Hence, we obtain the following for some $B'>0$:
\begin{equation}
 \label{eq;08.1.8.11}
 \bigl|\bigl(
 p^{(m(0))\natural}_{\gminib,S_{P_i}}
-p^{(m(0))\natural}_{\gminib,S_{P_{i+1}}}
\bigr)_{\gminic_1,\gminic_2}
\bigr|
\geq
 B'\cdot
 \exp\bigl(-\epsilon_2|z_1^{m(0)}|/3\bigr)
\end{equation}
From (\ref{eq;08.1.8.10})
and (\ref{eq;08.1.8.11}),
we obtain the contradiction.
Hence, we have
$\bigl(
 p^{(m(0))\natural}_{\gminib,S_{P_i}}
-p^{(m(0))\natural}_{\gminib,S_{P_{i+1}}}
\bigr)_{\gminic_1,\gminic_2}=0$
on $S_{P_i}\cap S_{P_{i+1}}$.
Then, we can show 
$\bigl(
 p^{(m(0))\natural}_{\gminib,S_{P_i}}
-p^{(m(0))\natural}_{\gminib,S}
\bigr)_{\gminic_1,\gminic_2}=0$
on $S_{P_i}$ by using an inductive argument.
In particular,
we obtain 
$\bigl(
 p^{(m(0))\natural}_{\gminib,S_{P}}
-p^{(m(0))}_{\gminib,S}
\bigr)_{\gminic_1,\gminic_2}=0$
on $S_P$.
Thus, we obtain the desired estimate,
and the proof of Lemma \ref{lem;08.1.8.17}
is finished.
\hfill\qed

\begin{lem}
\label{lem;08.1.8.16}
We have 
$p^{(j)\natural}_{\gminib,S_P}
-p^{(j)}_{\gminib,S|S_P}
=O\Bigl(
 \exp\bigl(-\epsilon_2|z_1^j|/10\bigr)
 \Bigr)$
with respect to $h$
for any $j$ and
any $\gminib\in\Irr(\theta,j)$.
\end{lem}
\pf
We omit to denote $|S_P$
for simplicity of the description.
We use an induction on $j$.
In the case $j=m(0)$,
we have already done
(Lemma \ref{lem;08.1.8.17}).
Let us look at
$\bigl(
 p^{(j)\natural}_{\gminib,S_P}
-p^{(j)}_{\gminib,S}
 \bigr)_{\gminic_1,\gminic_2}$
for $(\gminic_1,\gminic_2)\in\Irr(\theta,j)^2$.
Let us consider the case
$\ord(\gminic_1-\gminic_2)<j$.
We have the following:
\begin{equation}
 \label{eq;08.2.2.20}
 \bigl(
 p^{(j)\,\natural}_{\gminib,S_P}
-p^{(j)}_{\gminib,S}
 \bigr)_{\gminic_1,\gminic_2}
=p^{(j)}_{\gminic_2,S}\circ
  \bigl(
 p^{(j)\,\natural}_{\gminib,S_P}
-p^{(j)}_{\gminib,S}
 \bigr)\circ p^{(j)}_{\gminic_1,S}
=p^{(j)}_{\gminic_2,S}
\circ
 p^{(j)\,\natural}_{\gminib,S_P}
\circ
 p^{(j)}_{\gminic_1,S}
\end{equation}
Because 
$\ord(\gminic_1-\gminic_2)<j$,
one of $\ord(\gminic_1-\gminib)$
or $\ord(\gminic_2-\gminib)$ is strictly smaller than $j$.
In the case $q:=\ord(\gminic_1-\gminib)<j$,
we have the following:
\begin{equation}
 \label{eq;08.2.2.21}
 p^{(j)}_{\gminic_2,S}
\circ
 p^{(j)\,\natural}_{\gminib,S_P}
\circ
 p^{(j)}_{\gminic_1,S}
=
 p^{(j)}_{\gminic_2,S}
\circ
 p^{(j)\,\natural}_{\gminib,S_P}
\circ p^{(q)\natural}_{\eta_q(\gminib),S_P}
\circ p^{(q)}_{\eta_q(\gminic_1),S}
\circ p^{(j)}_{\gminic_1,S}
\end{equation}
By the hypothesis of the induction,
we have the following:
\begin{equation}
 \label{eq;08.2.2.22}
 p^{(q)\natural}_{\eta_q(\gminib),S_P}
\circ p^{(q)}_{\eta_q(\gminic_1),S}
=\bigl(
 p^{(q)\natural}_{\eta_q(\gminib),S_P}
-p^{(q)}_{\eta_q(\gminib),S}
\bigr)\circ
 p^{(q)}_{\eta_q(\gminic_1),S}
=O\Bigl(
 \exp\bigl(-\epsilon_2|z_1^q|/10\bigr)
 \Bigr)
\end{equation}
From (\ref{eq;08.2.2.20}),
(\ref{eq;08.2.2.21})
and (\ref{eq;08.2.2.22}),
we obtain the desired estimate for
$\bigl(
 p^{(j)\,\natural}_{\gminib,S_P}
-p^{(j)}_{\gminib,S}
 \bigr)_{\gminic_1,\gminic_2}$
in the case $\ord(\gminic_1-\gminib)<j$.
We can obtain a similar estimate
in the case
$\ord(\gminic_2-\gminib)<j$.

Let us consider the case
$\ord(\gminic_1-\gminic_2)=j$.
Note we have 
$\eta_{j-1}(\gminic_1)=\eta_{j-1}(\gminic_2)$.
In the case
$F_{\gminic_1,\gminic_2}\leq-\epsilon_2/8$
on $S_P$,
we have the following by flatness:
\[
 \bigl(
 p^{(j)\natural}_{\gminib,S_P}
 -p^{(j)}_{\gminib,S}
\bigr)_{\gminic_1,\gminic_2}
=O\Bigl(
 \exp\bigl(-\epsilon_2|z_1^{j}|/10\bigr)
 \Bigr)
\]
In the case $F_{\gminic_1,\gminic_2}\geq -\epsilon_2/4$
on $S_P$,
we can take a sequence
$P_1=P,P_2,\ldots,P_t\in Z$ with the following property:
\begin{itemize}
\item
 $Z_{P_i}\cap Z_{P_{i+1}}\neq\emptyset$.
\item
 $Z_{P_t}$ intersects with
 $\{F_{\gminic_1,\gminic_2}>0\}$.
\item
 $F_{\gminic_1,\gminic_2}\geq -\epsilon_2/4$
  on each $Z_{P_i}$.
\end{itemize}
By the second condition,
we have
$\bigl(
  p^{(j)\natural}_{\gminib,S_{P_t}}
 -p^{(j)}_{\gminib,S|S_{P_t}}
\bigr)_{\gminic_1,\gminic_2}=0$ 
on $S_{P_t}$.
On $S_{P_i}\cap S_{P_{i+1}}$,
we have
\begin{equation}
\label{eq;08.1.8.12}
 p^{(j)\natural}_{\gminib,S_{P_i}}
-p^{(j)\natural}_{\gminib,S_{P_{i+1}}}
=
\bigl(
 p^{(j)\natural}_{\gminib,S_{P_i}}
-\pi^{(j)}_{\gminib}
\bigr)
-\bigl(
 p^{(j)\natural}_{\gminib,S_{P_{i+1}}}
-\pi^{(j)}_{\gminib}
\bigr)
=O\Bigl(
 \exp\bigl(-\epsilon_1|z_1^{j}|/2\bigr)
 \Bigr)
\end{equation}
If $p^{(j)\natural}_{\gminib,S_{P_i}}
-p^{(j)\natural}_{\gminib,S_{P_{i+1}}}$ is not $0$,
we have the following for some $B,C,N>0$,
due to Lemma \ref{lem;07.7.10.1}:
\[
\bigl|
\bigl(
 p^{(j)\natural}_{\gminib,S_{P_i}}
-p^{(j)\natural}_{\gminib,S_{P_{i+1}}}
\bigr)_{\gminic_1,\gminic_2}
 \exp\bigl((\gminic_1-\gminic_2)/\lambda\bigr)
\bigr|
\geq
 B\cdot
 \exp\bigl(-C|z_1^{j+1}|\bigr)
 \cdot|z_1|^{N}
\]
Hence, we obtain the following
for some $B'>0$:
\begin{equation}
 \label{eq;08.1.8.13}
 \bigl|\bigl(
 p^{(j)\natural}_{\gminib,S_{P_i}}
-p^{(j)\natural}_{\gminib,S_{P_{i+1}}}
\bigr)_{\gminic_1,\gminic_2}
\bigr|
\geq
 B'\cdot
 \exp\bigl(-\epsilon_2|z_1^{j}|/3\bigr)
\end{equation}
From (\ref{eq;08.1.8.12}) and (\ref{eq;08.1.8.13}),
we obtain the contradiction.
Hence, we have
$\bigl(
 p^{(j)\natural}_{\gminib,S_{P_i}}
-p^{(j)\natural}_{\gminib,S_{P_{i+1}}}
\bigr)_{\gminic_1,\gminic_2}=0$
on $S_{P_i}\cap S_{P_{i+1}}$.
We obtain
$\bigl(
 p^{(j)\natural}_{\gminib,S_{P_i}}
-p^{(j)\natural}_{\gminib,S}
\bigr)_{\gminic_1,\gminic_2}=0$
on $S_{P_i}$ by an inductive argument.
In particular,
$\bigl(
 p^{(j)\natural}_{\gminib,S_{P}}
-p^{(j)}_{\gminib,S}
\bigr)_{\gminic_1,\gminic_2}=0$
on $S_P$.
Thus, we obtain the desired estimate,
and the proof of Lemma \ref{lem;08.1.8.16}
is finished.
\hfill\qed

\vspace{.1in}

Proposition \ref{prop;08.1.8.15}
immediately follows from 
Proposition \ref{prop;07.12.21.1}
and Lemma \ref{lem;08.1.8.16}.
\hfill\qed

\subsection{Comparison of some metrics}

We consider the following:
\begin{equation}
 \label{eq;07.11.24.1}
 F(w):=
 \exp\bigl(w\cdot \nbigb\bigr),
\quad
 \nbigb:=\sum_{\gminia\in\Irr(\theta)}
 \gminia\cdot\pi^{}_{\gminia}
=\sum_{m(0)\leq j\leq -1}
 \sum_{\gminib\in\Irr(\theta,j)}
 \zeta_j(\gminib)\cdot \pi^{(j)}_{\gminib}
\end{equation}
Let $F(w)^{\ast}h$ be the metric
given by $F(w)^{\ast}h(u,v)=h\bigl(F(w)u,F(w)v\bigr)$.

\begin{lem}
 \label{lem;07.10.7.50}
$F(\wbar)^{\ast}h$ and $g_{\irr}(w)^{\ast}h$
are mutually bounded.
\end{lem}
\pf
Because
$g_{\irr}(w)\circ F(\wbar)^{-1}
=\sum \exp\bigl(
 2\sqrt{-1}\Image(w\cdot \overline{\gminia})
 \bigr)\cdot \pi_{\gminia}$,
we have the boundedness of
$\bigl|g_{\irr}(w)\circ F(\wbar)^{-1}\bigr|_h$.
Similarly, 
$\bigl|F(\wbar)\circ g_{\irr}(w)^{-1}\bigr|_h$ is bounded.
Then, the claim of Lemma \ref{lem;07.10.7.50} follows.
\hfill\qed

\vspace{.1in}
For each $S:=S_i$,
we put
\begin{equation}
 \label{eq;07.11.23.26}
 \nbigb_S:=
 \sum_j\sum_{\gminib\in\Irr(\theta,j)}
 \zeta_j(\gminib)\cdot p^{(j)}_{\gminib,S},
\quad\quad
F_S(w):=\exp\bigl(w\cdot\nbigb_S\bigr).
\end{equation}

\begin{lem}
\label{lem;07.7.20.150}
For $|w|<\epsilon_2/100$,
we have the following estimate: 
\[
 \bigl|F(w)\circ F_S(w)^{-1}-1\bigr|_h
=O\Bigl( \exp\bigl(-\epsilon_2 |z_1|^{-1}/100\bigr)\Bigr)
\]
\[
  \bigl|F_S(w)\circ F(w)^{-1}-1\bigr|_h
=O\Bigl(\exp\bigl(-\epsilon_2|z_1|^{-1}/100\bigr) \Bigr)
\]
\end{lem}
\pf
We have $F(w)=\prod F^{(j)}(w)$
and $F_S(w)=\prod F^{(j)}_S(w)$,
where $F^{(j)}$ and $F_S^{(j)}$ are given as follows:
\[
 F^{(j)}(w):=\exp\Bigl(
 \sum_{\gminib\in\Irr(\theta,j)}
 w\cdot\zeta_j(\gminib)\cdot\pi^{(j)}_{\gminib}
 \Bigr),
\quad
 F^{(j)}_S(w):=\exp\Bigl(
 \sum_{\gminib\in\Irr(\theta,j)}
 w\cdot\zeta_j(\gminib)\cdot p^{(j)}_{\gminib,S}
\Bigr)
\]
We consider $G^{(j)}:=F^{(j)}-F_S^{(j)}$.
Because
$\pi^{(j)}_{\gminib}-p^{(j)}_{\gminib,S}
=O\bigl(\exp(-\epsilon_2|z_1|^{j}/10)\bigr)$,
we have $G^{(j)}=
 O\bigl(\exp(-\epsilon_2|z_1|^{j}/15)\bigr)$.
We set
\[
 \Gtilde^{(j)}:=
 \prod_{i>j}F_S^{(i)}\circ
 G^{(j)}\circ(F_S^{(j)})^{-1}\circ
 \prod_{i>j}\bigl(
 F_S^{(i)}\bigr)^{-1}.
\]
We have 
$|\Gtilde^{(j)}|_h=
O\bigl(\exp(-\epsilon_2|z_1|^{j}/20)\bigr)$.
We have the following equality:
\[
F\circ F_S^{-1}
=
(1+\Gtilde^{(-1)})\circ\cdots
\circ
(1+\Gtilde^{(m(0)+1)})\circ(1+\Gtilde^{(m(0))})
\]
Then, the claim of Lemma \ref{lem;07.7.20.150}
follows.
\hfill\qed

\vspace{.1in}
Let $F_S(w)^{\ast}h$ be the metric given by
$F_S(w)^{\ast}h(u,v)
=h\bigl(F_S(w)u,F_S(w)v\bigr)$.

\begin{lem}
\label{lem;07.7.20.200}
For $|w|<\epsilon_2/100$,
the metrics
$F(w)^{\ast}h$ and $F_S(w)^{\ast}h$ 
are mutually bounded.
\end{lem}
\pf
We have the following:
\[
\bigl|v\bigr|_{F(w)^{\ast}h}
\leq
\bigl|v\bigr|_{F_S(w)^{\ast}h}
\cdot \bigl|F_S(w)^{-1}\circ F(w)\bigr|_{F_S(w)^{\ast}h}
=\bigl|v\bigr|_{F_S(w)^{\ast}h}
\cdot\bigl|F(w)\circ F_S(w)^{-1}\bigr|_{h}
\]
Hence, we obtain
$\bigl|v\bigr|_{F(w)^{\ast}h}
\leq C\cdot \bigl|v\bigr|_{F_S(w)^{\ast}h}$
due to Lemma \ref{lem;07.7.20.150}.
We obtain 
$\bigl|v\bigr|_{F_S(w)^{\ast}h}
\leq C\cdot \bigl|v\bigr|_{F(w)^{\ast}h}$
by a similar argument.
Thus, the proof of Lemma \ref{lem;07.7.20.200}
is finished.
\hfill\qed

\subsection{End of proof of 
Proposition \ref{prop;07.7.20.230}}
\label{subsection;08.3.5.6}

According to Corollary \ref{cor;10.5.26.1},
$f$ gives a section of
$\bigl(\nbigp_a\nbigelambda\bigr)^{(T_1(w))}$
if and only if 
$\bigl|f\bigr|_{F_{S_i}(w)^{\ast}h}$
is bounded on $S_i$ for $i=1,\ldots,N$.
In Lemma \ref{lem;07.7.20.200},
we have obtained that 
$F(w)^{\ast}h$ and $F_{S_i}(w)^{\ast}h$
are mutually bounded on $S_i$ $(i=1,\ldots,N)$.
Hence, the claim of Proposition \ref{prop;07.7.20.230}
follows.
\hfill\qed

\section{Boundedness of some section}
\label{subsection;08.9.14.46}
We use the setting in 
Section \ref{subsection;07.11.17.1}.
For simplicity,
we assume that the coordinate system is
admissible for the good set $\Irr(\theta)$.
Let $k$ be determined by 
$\vecm(0)\in\seisuu_{<0}^k\times\veczero_{\ell-k}$.
We put $D(\leq k):=\bigcup_{j=1}^kD_j$.
Let $M$ be a sufficiently large integer,
say $M>100\cdot\ell\cdot |\vecm(0)|$,
where $|\vecm(0)|=\sum_{i=1}^k |m(0)_i|$.

\begin{lem}
\label{lem;07.11.5.2}
Let $f$ be a section of 
$\nbigp_0\End(\nbigelambda)$
with the following property:
\begin{itemize}
\item
 $\DDlambda f=0$ on $\Dhat^{(M)}(\leq k)$.
\item
 $\bigl[\Res_{i}(\DDlambda),f_{|D_i}\bigr]=0$ 
for $i=k+1,\ldots,\ell$.
\end{itemize}
Then, $|f|_h$ is bounded.
\end{lem}
\pf
We use an induction on $\ell$.
In the case $\ell=0$,
the claim is trivial.
In the following argument,
we will assume that the claim of the lemma 
holds in the case $\ell-1$.

Let us show the following claims for $m=1,\ldots,\ell$
by a descending induction on $m$:
\begin{description}
\item[A(m)]
Let $f$ be a section of $\nbigp_0\End(\nbigelambda)$
with the following property:
\begin{itemize}
\item
 $\DDlambda f=0$ on $\Dhat^{(M_m)}(\leq k)$,
 where $M_m=M-(\ell-m)\cdot|\vecm(0)|$.
\item
 $\bigl[\Res_{i}(\DDlambda),f_{|D_i}\bigr]=0$ for 
 $i=m+1,\ldots,\ell$.
\end{itemize}
Then, we have the following estimate
for some $C>0$ and $N>0$:
\[
 |f|_h\leq
 C\cdot \Bigl(
 -\sum_{i=1}^m\log|z_i|
 \Bigr)^N
\]
\end{description}

The claim $A(\ell)$ holds
due to the general result of the acceptable bundles
(Theorem \ref{thm;07.10.9.1}).
Let us show $A(m-1)$ by assuming $A(m)$.
Let $f\in\nbigp_0\End(\nbigelambda)$
such that
(i) $\DDlambda f_{|\Dhat^{(M_{m-1})}(\leq k)}=0$,
(ii)
$\bigl[
 \Res_{i}(\DDlambda), f_{|D_i}\bigr]=0$ 
for $i=m,\ldots,\ell$.
We would like to obtain the estimate:
\[
 |f|_h\leq C\cdot
\left(-\sum_{i=1}^{m-1}\log|z_i|\right)^N
\]
We put 
$g:=\DDlambda(\del_m)f
 \in\End(\nbigp_0\nbigelambda)$.
We have
(i) $\DDlambda g_{|\Dhat^{(M_m)}(\leq k)}=0$,
(ii) $\bigl[\Res_{i}(\DDlambda),\,
 g_{|D_i}\bigr]=0$ for $i=m+1,\ldots,\ell$.
Hence, we can apply $A(m)$ to $g$,
and we obtain the estimate
for some $C>0$ and $N>0$:
\[
 |g|_h\leq C\cdot\left(-\sum_{i=1}^{m}\log|z_i|\right)^N
\]
Let $\pi_m:X\lrarr D_m$ denote the projection.
We put $\pi_m^{-1}(Q)^{\ast}:=\pi_m^{-1}(Q)-\{Q\}$
for $Q\in D_m^{\circ}$.
Then, we obtain the following:
\[
 \int_{\pi_m^{-1}(Q)^{\ast}}
 \bigl|g_{|\pi_m^{-1}(Q)^{\ast}}\bigr|_h^2\cdot
 \bigl|dz_m\cdot d\zbar_m\bigr|
\leq
 C_1\cdot\left(-\sum_{j=1}^{m-1} \log|z_j(Q)|\right)^{N_1}
\]
Let $\Delta_m$ denote 
the Laplacian $-\del_{z_m}\del_{\zbar_m}$.
By Corollary \ref{cor;07.11.23.10} below,
we have the following inequality 
on $\pi_m^{-1}(Q)^{\ast}$:
\[
 \Delta_m\bigl(|f_{|\pi_m^{-1}(Q)^{\ast}}|_h^2\bigr)
\leq 
 \bigl|g_{|\pi_m^{-1}(Q)^{\ast}}\bigr|_h^2
\]
We can take $G_Q(z_m)$ satisfying the following:
\[
 \Delta_m G_Q=|g_{|\pi_m^{-1}(Q)}|^2,
\quad
 \sup\bigl|G_Q\bigr|\leq C'\cdot
\left(-\sum_{j=1}^{m-1}\log|z_j|(Q)\right)^N
\]
Note $\bigl|f_{|\pi_m^{-1}(Q)}\bigr|_h^2-G_Q$ is 
bounded on $\pi_m^{-1}(Q)^{\ast}$
for each $Q$,
which follows from the norm estimate
in the curve case
(Proposition \ref{prop;07.10.6.40}, below).
Therefore, we have
$\Delta_m
 \bigl(\bigl| f_{|\pi_m^{-1}(Q)}\bigr|_h^2-G_Q\bigr)
 \leq 0$ as distributions
on $\pi_m^{-1}(Q)$
(Lemma 2.2 of \cite{s2}).
Hence, we obtain the following:
\begin{multline}
 \sup\bigl|f_{|\pi_m^{-1}(Q)}\bigr|_h^2
\leq
 \max_{|z_m'|=1/2}
 \bigl|f_{|\pi_m^{-1}(Q)}(z_m')\bigr|_h^2
+C''\left(-\sum_{i=1}^{m-1}\log|z_i|^2(Q)\right)^N
 \\
\leq
 C'''\left(-\sum_{i=1}^{m-1}\log|z_i|^2(Q)\right)^N
\end{multline}
Here, we have used the hypothesis
of the induction on $\ell$
for the estimate of
$ \max_{|z_m'|=1/2}
 \bigl|f_{|\pi_m^{-1}(Q)}(z_m')\bigr|_h^2$.
Thus, we obtain $A(m-1)$,
and the descending induction on $m$ can proceed.
The claim $A(0)$ means Lemma  \ref{lem;07.11.5.2}.
\hfill\qed

\vspace{.1in}
We can show the following lemma
by the same argument.
\begin{lem}
Let $f$ be a section of $\nbigp_0\nbigelambda$
with the following property:
\begin{itemize}
\item
 $\DDlambda f=0$ on $\Dhat^{(M)}(\leq k)$.
\item
 $\Res_{i}(\DDlambda)(f)=0$ for $i=k+1,\ldots,\ell$.
\end{itemize}
Then, $|f|_h$ is bounded.
\hfill\qed
\end{lem}

Let $(E_i,\delbar_i,\theta_i,h_i)$ $(i=1,2)$
be unramifiedly good wild harmonic bundles
as in Section \ref{subsection;07.11.17.1},
with the good set of irregular values $\Irr(\theta_i)$.
For simplicity,
we assume that $\Irr(\theta_1)=\Irr(\theta_2)$,
and that the coordinate system is admissible
for $\Irr(\theta_1)$.
Let $\vecm(0)$ be the minimum of $\Irr(\theta_1)$,
and let $k$ be determined by
$\vecm(0)\in\seisuu_{<0}^k\times\veczero_{\ell-k}$.
Let $M$ be a sufficiently large integer,
say $M>100\cdot\ell\cdot |\vecm(0)|$,
where $|\vecm(0)|=\sum_{i=1}^k |m(0)_i|$.
Note that $\Hom(\nbigelambda_1,\nbigelambda_2)$
with the naturally induced metric is acceptable.
The following lemma can be shown
by the same argument.
\begin{lem}
\label{lem;10.6.13.50}
Let $f$ be a section of 
$\nbigp_0\Hom(\nbigelambda_1,\nbigelambda_2)$
with the following property:
\begin{itemize}
\item
 $\DDlambda f=0$ on $\Dhat^{(M)}(\leq k)$.
\item
 $\bigl[\Res_{i}(\DDlambda),f_{|D_i}\bigr]=0$ 
for $i=k+1,\ldots,\ell$.
\end{itemize}
Then, $|f|_h$ is bounded.
\hfill\qed
\end{lem}

\subsection{Weitzenb\"ock formula (Appendix)}
\label{subsection;08.9.29.13}

Let us recall a variant of Weitzenb\"ock formula
for harmonic bundles due to Simpson
with a slightly refined form.
It will be used in Section
\ref{subsection;08.9.15.43}.
The original one (Corollary \ref{cor;07.11.23.10})
was used in the proof of Lemma \ref{lem;07.11.5.2}.
It will be also useful in Sections
\ref{subsection;08.9.14.12}
and \ref{subsection;08.9.15.42}.

Let $z$ be a coordinate of $\cnum$,
and let $U$ be an open subset of $\cnum$.
We use the Euclidean metric $g=dz\cdot d\zbar$.
Let $\Delta''$ denote the Laplacian $-\del_z\del_{\zbar}$.

Let $\harmonicbundle$ be a harmonic bundle
on $U$.
For any complex number $\lambda$,
let $d_{\lambda}'':=\delbar_E+\lambda\theta^{\dagger}$,
$\delta_{\lambda}':=\del_E-\lambdabar\theta$,
and $\DDlambda:=
 d_{\lambda}''+\lambda\del_E+\theta$.

\begin{lem}
\label{lem;08.9.29.12}
Let $s$ be a $C^{\infty}$-section of $E$.
We have the following inequality:
\[
\Delta''|s|_h^2
\leq
 \bigl|\DDlambda s\bigr|^2_{h,g}
+2|s|_h\cdot \bigl|\delta_h'd_{\lambda}''s\bigr|_{h,g}
\]
\end{lem}
\pf
We use the pairing
$(E\otimes \Omega^{p,q})
\otimes (E\otimes\Omega^{r,s})
\lrarr \Omega^{p+s,q+r}$
induced by $h$,
which is denoted by $(\cdot,\cdot)_h$.
Let $R(h)$ denote the curvature 
of the unitary connection 
$d_{\lambda}''+\delta_{\lambda}'$.
We have the following equality:
\begin{multline}
 \del\delbar|s|_h^2
=\bigl(\delta_{\lambda}'d_{\lambda}''s,s\bigr)_h
-\bigl(d''_{\lambda}s,d''_{\lambda}s\bigr)_h
+\bigl(\delta_{\lambda}'s,\delta_{\lambda}'s\bigr)_h
+\bigl(s,d''_{\lambda}\delta_{\lambda}'s\bigr)_h\\
=\bigl(\delta_{\lambda}'d_{\lambda}''s,s\bigr)_h
-\bigl(s,\delta_{\lambda}'d_{\lambda}''s\bigr)_h
+\bigl(s,R(h)s\bigr)_h
-\bigl(d_{\lambda}''s,d_{\lambda}''s\bigr)_h
+\bigl(\delta_{\lambda}'s,\delta_{\lambda}'s\bigr)_h
\end{multline}
We have the following equality:
\[
 \bigl(s,R(h)s\bigr)_h
=-(1+|\lambda|^2)\cdot
 \bigl(\theta\cdot s,\theta\cdot s\bigr)_h
-(1+|\lambda|^2)\cdot
 \bigl(\theta^{\dagger}\cdot s,\theta^{\dagger}\cdot s\bigr)_h
\]
Let $\DD^{\lambda\prime}:=\lambda\del_E+\theta$.
Then, we have the following:
\[
 \bigl(\delta'_{\lambda}s,\delta_{\lambda}'s\bigr)_h
-(1+|\lambda|^2)\cdot \bigl(\theta\cdot s,\theta\cdot s\bigr)_h
=(1+|\lambda|^2)\cdot\bigl(\del_Es,\del_Es\bigr)_h
-\bigl(\DD^{\lambda\prime}s,\DD^{\lambda\prime}s\bigr)_h
\]
Hence, we obtain the following:
\begin{multline*}
 \del\delbar|s|_h^2=
 \bigl(\delta_{\lambda}'d_{\lambda}''s,\, s \bigr)_h
-\bigl(s,\delta_{\lambda}'d_{\lambda}''s\bigr)_h
-(1+|\lambda|^2)\cdot
 \bigl(\theta^{\dagger}s,\theta^{\dagger}s\bigr)_h
-\bigl(d''_{\lambda}s,d''_{\lambda}s\bigr)_h
 \\
+(1+|\lambda|^2)\cdot\bigl(\del_E s,\del_Es\bigr)_h
-\bigl(\DD^{\lambda\prime}s,
 \DD^{\lambda\prime}s\bigr)_h
\end{multline*}
Then, the claim of the lemma follows.
\hfill\qed

\begin{cor}[Lemma 4.18 
of \cite{mochi}] \label{cor;07.11.23.10}
Let $s$ be a holomorphic section of $\nbigelambda$.
Then, we have the inequality
$\Delta''|s|_h^2\leq
 \bigl|\DDlambda s\bigr|^2_{h,g}$.
(Lemma {\rm 4.18} of {\rm\cite{mochi}}
should be corrected.)
\hfill\qed
\end{cor}

\chapter{Some Basic Results in the Curve Case}
\label{section;08.10.24.42}
In this chapter,
we study the one dimensional case.
In Section \ref{subsection;08.9.14.30},
we show the norm estimate 
for holomorphic sections of $\nbigp\nbigelambda$.
In Section \ref{subsection;08.9.14.31},
we show the correspondence between 
the parabolic weights and the residues 
for $(\nbigp\nbige^0,\DD^0)$
and $(\nbigp\nbigelambda,\DDlambda)$.
These are natural generalizations of
Simpson's results in the tame case.
The arguments are also essentially the same.
In Section \ref{subsection;08.9.14.32},
we give a characterization 
of the lattice $\nbigp_a\nbigelambda$
by the eigenvalues of the residues
in the case that $\lambda$ is generic.
In Section \ref{subsection;08.9.14.33},
we argue some basic property
of harmonic forms for wild harmonic bundles
on quasiprojective curves,
which will be used
in the proof of Hard Lefschetz Theorem
for polarized wild pure twistor $D$-modules.
(See Section \ref{subsection;07.10.23.30}.)

\section{Norm estimate for holomorphic sections
 of $\nbigp_c\nbigelambda$}
\label{subsection;08.9.14.30}
\subsection{Statement}
\label{subsection;08.12.7.2}

We put $X=\Delta$ and $D=\{O\}$.
Let $\harmonicbundle$ be 
an unramifiedly good wild harmonic bundle
on $X-D$.
We assume to have the decomposition on $X$:
\begin{equation}
 \label{eq;07.7.18.100}
 (E,\theta)=
\bigoplus_{
 (\gminia,\alpha)\in\Irr(\theta)\times\Sp(\theta)}
 (E_{\gminia,\alpha},\theta_{\gminia,\alpha})
\end{equation}
We have the expression $\theta=f\cdot dz$.

As explained in Chapter \ref{section;07.6.2.1},
we have the prolongment
$\nbigp_{c}{\nbigelambda}$ 
for each $\lambda\in\cnum$
and $c\in\real$.
Let $F$ denote the parabolic filtration
of $\nbigp_{c}\nbigelambda_{|O}$.
We have the endomorphism
$\Res(\DDlambda)$ of $\nbigp_c\nbigelambda_{|O}$,
which preserves the parabolic filtration.
The induced endomorphism of
$\Gr^F(\nbigp_c\nbigelambda_{|O})$ is also denoted by
$\Res(\DDlambda)$.
Let $W$ denote the weight filtration of
$\Gr^F(\nbigp_c\nbigelambda_{|O})$ 
associated to the nilpotent part of $\Res(\DDlambda)$.

Let $\vecv$ be a holomorphic frame of $\nbigp_{c}\nbigelambda$
such that 
(i) it is compatible with the parabolic filtration $F$,
(ii) the induced frame on $\Gr^F(\nbigp_{c}\nbigelambda_{|O})$
 is compatible with the weight filtration $W$.
We put $a(v_i):=\deg^F(v_i)$ and $k(v_i):=\deg^W(v_i)$.
Let $h_0$ be the $C^{\infty}$-metric of $E$
given as follows:
\begin{equation}
\label{eq;07.7.18.101}
 h_0(v_i,v_j):=
 \delta_{i,j}\cdot |z|^{-2a(v_i)}
 \cdot (-\log|z|^2)^{k(v_i)}
\end{equation}

\begin{prop}
\label{prop;07.10.6.40}
The metrics $h$ and $h_0$ are mutually bounded.
In other words, the standard norm estimate holds.
\end{prop}
Since it can be shown 
by the argument in \cite{s2},
we give only an indication.

\subsection{The case $\lambda=0$}

First let us consider the case $\lambda=0$.
Let $\pi_{\gminia,\alpha}$ be the projection
onto $E_{\gminia,\alpha}$
in the decomposition (\ref{eq;07.7.18.100}).
According to  Theorems \ref{thm;07.10.4.1} 
and \ref{thm;07.10.4.3},
$\pi_{\gminia,\alpha}$ are bounded.
Hence, we have the following decomposition
as the prolongment of
the decomposition (\ref{eq;07.7.18.100}):
\begin{equation}
 \label{eq;07.10.6.10}
 (\nbigp_c\nbige^0,\DD^0)
=\bigoplus_{(\gminia,\alpha)\in
 \Irr(\theta)\times\Sp(\theta)}
 \bigl(
 \nbigp_c\nbige^0_{\gminia,\alpha},
 \DD^0_{\gminia,\alpha}
\bigr)
\end{equation}
The decomposition is compatible with the parabolic filtration
and the residue $\Res(\theta)$.
We may assume $\vecv$ is compatible with 
the decomposition (\ref{eq;07.10.6.10}).

For each $(a,\alpha,\gminia)$,
we have the endomorphism of 
$V_{a,\alpha,\gminia}:=
 \Gr^F_a(\nbigp_c\nbige^0_{\alpha,\gminia})$
induced by $\Res(\theta)$.
The nilpotent part is denoted by 
$N_{a,\alpha,\gminia}$.
We have the model harmonic bundle 
on $X-D$ obtained from
$\bigl(V_{a,\alpha,\gminia},N_{a,\alpha,\gminia}\bigr)$
denoted by
$E\bigl(V_{a,\alpha,\gminia},N_{a,\alpha,\gminia}\bigr)$.
(See Subsection \ref{subsection;10.5.24.3}).
We have the rank one harmonic bundle
$L(a,\alpha,\gminia)=\bigl(
 \nbigo\!\cdot\! e,
 \theta^{L}_{a,\alpha,\gminia},
 h^{L}_{a,\alpha,\gminia}\bigr)$
as in Subsection \ref{subsection;10.5.24.4}.
Then, we obtain the wild harmonic bundle
$(\Etilde,\delbar_{\Etilde},\htilde,\thetatilde):=
 \bigoplus_{(a,\alpha,\gminia)}
E\bigl(V_{a,\alpha,\gminia},N_{a,\alpha,\gminia}\bigr)
\otimes
 L(a,\alpha,\gminia)$.
According to Theorem \ref{thm;07.12.2.55},
we obtain $\nbigp_c\nbigetilde^0$ for each $c\in\real$,
and the decomposition:
\begin{equation}
\label{eq;07.10.6.11}
 (\nbigp_c\nbigetilde^0,\thetatilde)=
 \bigoplus_{(\alpha,\gminia)}
  \bigoplus_a
 \bigl(\nbigp_c\nbigetilde^{0}_{a,\alpha,\gminia}, 
 \thetatilde_{a,\alpha,\gminia}\bigr)
\end{equation}

We can take a holomorphic isomorphism
$\Phi:\nbigp_{c}{\nbigetilde^0}
 \lrarr\nbigp_c\nbige^0$
such that 
(i) it preserves the decompositions (\ref{eq;07.10.6.10})
 and (\ref{eq;07.10.6.11}),
(ii) it preserves the parabolic filtrations,
(iii) $\Gr^F(\Phi_{|O})$ is compatible with
 the residues.
It can be checked by a direct calculation
that $\Phi(\htilde)$ and $h_0$
are mutually bounded.
Recall that $h_0$ and $h$ are mutually bounded
up to log order,
and hence $\Phi(\htilde)$ and $h$ are 
mutually bounded up to log order.
Due to a general result of Simpson
(Corollary 4.3 of \cite{s2}),
we have the following:
\begin{itemize}
\item
Let $K$ be a hermitian metric of $E=\nbige^0$
with the following property:
\begin{description}
\item[(i)]
 It is adapted to the filtered bundle
 $\nbigp_{\ast}\nbige^0$.
 The induced metric of the dual $\nbige^{0\lor}$
 is also adapted to the induced parabolic filtration
 of $\nbigp_{\ast}\nbige^{0\lor}$.
\item[(ii)]
 Let $\theta^{\dagger}_K$ denote the adjoint of $\theta$
 with respect to $K$,
 and let $\del_K$ be the $(1,0)$-operator 
 determined by $\delbar_E$ and $K$.
 Let $F_K$ denote the curvature of
 the connection $\delbar_E+\del_K+\theta+\theta^{\dagger}_K$.
 Then, $F_K$ is $L^p$ with respect to
 $K$ and the Euclid metric $dz\cdot d\zbar$.
\end{description}
Then, the metrics $K$ and $h$ are mutually bounded.
\end{itemize}
It is easy to check that the conditions (i) and (ii) hold
for $\Phi(\htilde)$.
Thus, we can conclude that
$h$ and $h_0$ are mutually bounded 
in the case $\lambda=0$.

\subsection{The case $\lambda\neq 0$}

Let us consider the case $\lambda\neq 0$.
We have the irregular decomposition:
\begin{equation}
 \label{eq;07.10.6.12}
 \bigl(\nbigp_c\nbigelambda,
 \DDlambda\bigr)_{|\Dhat}=
 \bigoplus_{\gminia\in \Irr(\theta)}
 \bigl(
 \nbigp_c\nbigehatlambda_{\gminia},
 \DDlambda_{\gminia}
\bigr)
\end{equation}
We have the parabolic filtration $F$
and the endomorphism $\Res(\DDlambda)$
of $\nbigp_c\nbigehatlambda_{\gminia|O}$.
They are compatible.
Let $\EE_{\alpha}(\nbigp_c\nbigehatlambda_{\gminia|O})$
denote the generalized eigen space of $\Res(\DDlambda)$
corresponding to the eigenvalue $\alpha$.
We obtain a vector space
$V_{a,\alpha,\gminia}:=
 \Gr^F_{a}\EE_{\alpha}
 (\nbigp_c\nbigehatlambda_{\gminia|O})$.
Let $N_{a,\alpha,\gminia}$ denote the nilpotent part
of the endomorphism of $V_{a,\alpha,\gminia}$ induced by
$\Res(\DDlambda)$.
We have the model bundle 
$E\bigl(V_{a,\alpha,\gminia},N_{a,\alpha,\gminia}\bigr)$
obtained from
$\bigl(V_{a,\alpha,\gminia},N_{a,\alpha,\gminia}\bigr)$.
Let $(b,\beta)\in\real\times\cnum$ be given by the condition
$\kmsmap(\lambda,b,\beta)=(a,\alpha)$.
(See Subsection \ref{subsection;07.10.14.15}
for $\kmsmap(\lambda)$.)
Let $L(b,\beta,\gminia)$ be 
the harmonic bundle of rank one
as in Subsection \ref{subsection;10.5.24.4}.
We obtain a wild harmonic bundle
\[
 \bigl(\Etilde,\delbar_{\Etilde},\htilde,\thetatilde\bigr):=
 \bigoplus_{(a,\alpha,\gminia)}
E\bigl(V_{a,\alpha,\gminia},N_{a,\alpha,\gminia}\bigr)
\otimes
 L(b,\beta,\gminia).
\]
We have the associated meromorphic $\lambda$-flat
bundle $(\nbigp_c\nbigelambdatilde,\DDlambdatilde)$
which has the decomposition:
\begin{equation}
\label{eq;07.10.6.17}
 \nbigp_c\nbigetilde^{\lambda}
=\bigoplus_{\gminia\in\Irr(\theta)}
 \Bigl(
 \bigoplus_{a,\alpha}
 \nbigp_c\nbigetilde^{\lambda}_{a,\alpha,\gminia}
\Bigr)
\end{equation}
By construction,
we have an isomorphism
\[
 \Phi_{O,a,\alpha,\gminia}:
 (\Gr^F_a\EE_{\alpha}\nbigp_c
 \nbigelambdatilde_{a,\alpha,\gminia},
 \Res\DDlambdatilde)
\simeq
 (\Gr^F_a\EE_{\alpha}\nbigp_{c}
 \nbigehat^{\lambda}_{\gminia},
 \Res\DDlambda_{\gminia})
\]
We can take an isomorphism
$\Phi_{O,\gminia}:
 \bigoplus_{a,\alpha}
 \nbigp_c\nbigelambdatilde_{a,\alpha,\gminia|O}
\simeq
\nbigp_{c}
 \nbigehat^{\lambda}_{\gminia|O}$,
which preserves the parabolic filtration,
and induces 
$\bigoplus_{a,\alpha}\Phi_{O,a,\alpha,\gminia}$.
Let $\Dhat^{(N)}$ denote the $N$-th
infinitesimal neighbourhood of $D$
for a large integer $N$.
We can take an isomorphism
\[
\Phi_{\Dhat^{(N)},\gminia}:
\bigoplus_{a,\alpha}
 \nbigp_c
 \nbigelambdatilde_{a,\alpha,\gminia|\Dhat^{(N)}}
\simeq
 \nbigp_c
 \nbigehat^{\lambda}_{\gminia|\Dhat^{(N)}}
\]
such that 
$\Phi_{\Dhat^{(N)},\gminia|O}=\Phi_{O,\gminia}$.
By a general theory,
we can take a holomorphic
decomposition of $\nbigp_c\nbigelambda$
whose restriction to $\Dhat^{(N)}$
is the same as (\ref{eq;07.10.6.12}):
\begin{equation}
\label{eq;07.10.6.15}
 \nbigp_c\nbigelambda
=\bigoplus_{\gminia\in\Irr(\theta)}
 \nbigp_c\nbigelambda_{\gminia,N}
\end{equation}
Let $q_{\gminia,N}$ denote the projection
onto $\nbigp_c\nbigelambda_{\gminia,N}$.
We can take a holomorphic isomorphism
$\bigoplus_{a,\alpha}
 \nbigp_c\nbigelambdatilde_{a,\alpha,\gminia}
\simeq
\nbigp_c\nbigelambda_{\gminia,N}$
whose restriction to $\Dhat^{(N)}$
is equal to
$\Phi_{\Dhat^{(N)},\gminia}$.
In particular,
we obtain 
$\Phi:\nbigp_c\nbigelambdatilde\lrarr
\nbigp_c\nbigelambda$
such that
(i) it preserves the decompositions
 (\ref{eq;07.10.6.17}) and (\ref{eq;07.10.6.15}),
(ii) $\Phi_{|O}$ preserves the parabolic filtrations,
(iii) the induced map $\Gr^F(\Phi)$
 is compatible with the residues.
We identify $\nbige$ and $\nbigetilde$
in the following argument by $\Phi$.
By a direct calculation,
we can show that 
$\htilde$ and $h_0$ are mutually bounded.
Hence, we have only to show
$\htilde$ and $h$ are mutually bounded.
We remark that we have already known
that $\htilde$ and $h$ are mutually bounded
up to log order
by the weak norm estimate
(Theorem \ref{thm;07.12.2.55}),
which we will use implicitly.

Let $\gbigf$ be the endomorphism of 
$\nbigp_{c}\nbigelambda$
determined by
$\gbigf\cdot dz/z=\DDlambda-\DDlambdatilde$.
By construction, we have
$\gbigf(\nbigp_{c}{\nbigelambda})
\subset \nbigp_{c-\epsilon}{\nbigelambda}$
for some $\epsilon>0$.
We have the estimate
$\bigl[\gbigf,q_{\gminia,N}\bigr]
=\bigl[\DDlambda,q_{\gminia,N}\bigr]
=O(|z|^{N/2})$ by construction of
$q_{\gminia,N}$.
(Note that we have already known that
$\htilde$ and $h$ are mutually bounded,
this estimate does not depend on
the choice of a metric.
We will often omit this type of remark.)
Recall that $\DDlambda$ and $\htilde$
determine the operators
$\theta_{\htilde}\in\End(E)\otimes\Omega^{1,0}$,
$\theta_{\htilde}^{\dagger}\in\End(E)\otimes\Omega^{0,1}$
and the pseudo-curvature
$G(\DDlambda,\htilde) 
:=-\lambda^{-1}(1+|\lambda|^2)^2\cdot
 (\delbar_{\htilde}+\theta_{\htilde})^2
=-\lambda^{-1}(1+|\lambda|^2)^2\cdot
 \delbar_{\htilde}\theta_{\htilde}$.
(See \cite{mochi5}.)
As in the case $\lambda=0$,
we obtain that $h$ and $\htilde$ are mutually bounded,
once we show that $G(\DDlambda,\htilde)$ is
$L^p$ for some $p>1$,
thanks to a general result of Simpson.
(Corollary 4.3 of \cite{s2}.
See also Section 4.3 of \cite{mochi}.
Note, in \cite{mochi},
 the pseudo-curvature is considered as
 $G(\DDlambda,\htilde)
 =\delbar_{\htilde}\theta_{\htilde}$,
which does not make any essential difference
in the conclusion.)

Let $d''_{\lambda}:=\delbar_E+\lambda\theta^{\dagger}$,
which is the holomorphic structure of $\nbigelambda$.
We have the equalities:
\[
 0=(1+|\lambda|^2)^{-1}\cdot
 G(\DDlambdatilde,\htilde)
=R(d_{\lambda}'',\htilde)
+(1+|\lambda|^2)\bigl[\thetatilde,
 \thetatilde_{\htilde}^{\dagger}\bigr]
\]
\[
 (1+|\lambda|^2)^{-1}
G(\DDlambda,\htilde) 
=R(d_{\lambda}'',\htilde)
+(1+|\lambda|^2)
 \bigl[\theta_{\htilde},\theta^{\dagger}_{\htilde}\bigr]
=-(1+|\lambda|^2)\bigl(
 [\thetatilde,\thetatilde^{\dagger}_{\htilde}]
-[\theta_{\htilde},\theta^{\dagger}_{\htilde}]
 \bigr)
\]
Recall the following equalities
(Section 2.2 of \cite{mochi5}):
\[
 \theta_{\htilde}=
 \frac{\gbigf}{1+|\lambda|^2}+\thetatilde,
\quad
 \theta^{\dagger}_{\htilde}
 =\frac{\gbigf^{\dagger}_{\htilde}}{1+|\lambda|^2}
 +\thetatilde^{\dagger}_{\htilde}
\]
We have 
$[\gbigf,\gbigf^{\dagger}_{\htilde}]=O(|z|^{\epsilon})$.
Let us estimate
$\bigl[\gbigf,\thetatilde^{\dagger}_{\htilde}\bigr]$.
We have the decomposition:
\[
 \End(\nbigp_c\nbigelambda)
=\nbigc(\nbigp_c\nbigelambda)
\oplus
 \nbigd(\nbigp_c\nbigelambda),
\]
\[
  \nbigc(\nbigp_c\nbigelambda)
=\bigoplus_{\gminia\neq\gminib}
 \Hom\bigl(\nbigp_c\nbigelambda_{\gminia,N},
 \nbigp_c\nbigelambda_{\gminib,N}\bigr),
\quad
 \nbigd(\nbigp_c\nbigelambda)
=\bigoplus_{\gminia}
 \End\bigl(\nbigp_c\nbigelambda_{\gminia,N})
\]
We have the corresponding decomposition
$\gbigf=\nbigc(\gbigf)+\nbigd(\gbigf)$.
By construction,
we have $\nbigc(\gbigf)=O(|z|^N)$
with respect to $\htilde$,
and thus
$\bigl[\nbigc(\gbigf),
 \thetatilde^{\dagger}_{\htilde}\bigr]
=O\bigl(|z|^{N/2}\bigr)$.
We have the decomposition:
\[
 \thetatilde
=\bigoplus_{a,\alpha,\gminia}
 \Bigl(
 \bigl( d\gminia+\beta\cdot dz/z
 \bigr)\cdot \id_{\nbigp_c\nbigetilde_{a,\alpha,\gminia}}
+\thetatilde'_{a,\alpha,\gminia}
\Bigr)
\]
The terms $(\thetatilde'_{a,\alpha,\gminia})^{\dagger}_{\htilde}$
are $O\bigl(|z|^{-1}(-\log|z|)^{-1}\bigr)$
with respect to $\htilde$,
and the terms
$ \bigl( d\gminiabar+\alphabar\cdot d\zbar/\zbar^{-1}
 \bigr)\cdot \id_{\nbigp_c\nbigetilde_{a,\alpha,\gminia}}$
do not contribute to
$\bigl[\nbigd(\gbigf),\,
 \thetatilde^{\dagger}_{\htilde}\bigr]$.
Hence, we obtain
$\bigl[\nbigd(\gbigf),
 \thetatilde^{\dagger}_{\htilde}\bigr]
=O\bigl(|z|^{\epsilon-2}dzd\zbar\bigr)$,
and 
$G(\DDlambda,\htilde)$ is $L^p$ for some $p>1$.
Thus the proof of Proposition \ref{prop;07.10.6.40}
is finished.
\hfill\qed

\section{Comparison of the data at $O$}
\label{subsection;08.9.14.31}
\subsection{Statement}

Let $X$, $D$ and $\harmonicbundle$
be as in Subsection \ref{subsection;08.12.7.2}.
Let us compare the data at the origin $O$
of $\nbigp\nbigelambda$ and $\nbigp\nbige^0$.
For $\lambda\neq 0$,
we have the vector space
$\Gr^F_a\bigl(\nbigp\nbigehatlambda_{\gminia}\bigr)$
with the endomorphism $\Res(\DDhatlambda_{\gminia})$.
We have the generalized eigen decomposition:
\[
 \Gr^F_a\bigl(\nbigp\nbigehatlambda_{\gminia}\bigr)
=\bigoplus_{\alpha\in\cnum}
 \Gr^{F,\EE}_{a,\alpha}
 \bigl(\nbigp\nbigehatlambda_{\gminia}\bigr)
\]
The residue $\Res(\DDlambda)$ induces the endomorphism,
whose nilpotent part is denoted by
$N^{\lambda}_{\gminia,a,\alpha}$.
Similarly, 
we have the vector spaces
$\Gr^{F,\EE}_{a,\alpha}\bigl(
 \nbigp\nbige^0_{\gminia}
 \bigr)$
with the nilpotent endomorphism
$N^0_{\gminia,a,\alpha}$.
We consider the following sets:
\[
 \KMS\bigl(\nbigp\nbigehatlambda_{\gminia}\bigr):=
 \bigl\{
 (a,\alpha)\in\real\times\cnum\,\big|\,
\Gr^{F,\EE}_{a,\alpha}
\bigl(\nbigp\nbigehatlambda_{\gminia}\bigr)\neq 0
 \bigr\}
\quad (\lambda\neq 0)
\]
\[
 \KMS\bigl(\nbigp\nbige^0_{\gminia}\bigr):=
 \bigl\{
 (a,\alpha)\in\real\times\cnum\,\big|\,
 \Gr^{F,\EE}_{a,\alpha}\bigl(
 \nbigp\nbige^0_{\gminia}
 \bigr)\neq 0
 \bigr\}
\]
The dimension of the vector space corresponding to
$(a,\alpha)$ is called the multiplicity of $(a,\alpha)$,
and denoted by $\gminim(\lambda,a,\alpha)$.
The following proposition was observed by Biquard-Boalch
(\cite{biquard-boalch}),
at least if the wild harmonic bundle is given on a quasiprojective curve.
\begin{prop}
\label{prop;07.7.19.31}
$\kmsmap(\lambda)$ gives the bijective map
$\KMS\bigl(\nbigp\nbige^0_{\gminia}\bigr)
\lrarr
 \KMS\bigl(\nbigp\nbigehatlambda_{\gminia}\bigr)$
preserving the multiplicities.
The conjugacy classes of
$N^{0}_{\gminia,a,\alpha}$
and $N^{\lambda}_{\gminia,
 \kmsmap(\lambda,a,\alpha)}$ 
are the same.
\end{prop}
The claims can be shown using
the essentially same argument 
as that in \cite{s2}.
We indicate only an outline.

\subsection{An estimate in \cite{s2}}

Let $E$ be a holomorphic vector bundle 
on a punctured disc $\Delta^{\ast}$
with a frame $\vecv$.
Let $K$ be a hermitian metric of $E$
for which $\vecv$ is orthogonal
and $|v_i|_K=|z|^{-a_i}(-\log|z|)^{-k_i/2}$
for some $a_i\in\real$ and $k_i\in\seisuu$.

Let $F$ be a $C^{\infty}$-section of $E$ 
such that $|F|_K$ is bounded.
For the expression $F=\sum F_i\cdot v_i$,
we know $|F_i|\cdot |z|^{a_i}(-\log|z|)^{-k_i/2}$
are bounded.

\begin{prop}
\label{prop;07.2.6.31}
Let $M$ be the section of $E$ determined by
$\delbar F=M\cdot d\zbar$.
Assume the following conditions:
\[
|M|_K=O\bigl(|z|^{-1}(-\log|z|)^{-1}\bigr),
\quad
\int|M|_K^2\cdot(-\log|z|) \cdot
 \bigl|dz\cdot d\zbar\bigr|<\infty
\]
Then, the following holds:
\begin{itemize}
\item
 In the case $a_i\neq 0$,
 we have $|F_i|\cdot |z|^{a_i}(-\log|z|)^{-k_i/2}
 =O\bigl((-\log|z|)^{-1}\bigr)$.
\item
 In the case $a_i=0$ and $k_i\neq 0$,
 we have
 $\bigl\|F_i\cdot |z|^{a_i}(-\log|z|)^{-k_i/2}\bigr\|_W
 <\infty$,
 where $\bigl\|\cdot\bigr\|_W$ is given as follows:
\begin{equation}
 \label{eq;08.9.13.30}
 \|G\|_W^2:=\int |G|_K^2\frac{\bigl|dz\cdot d\zbar\bigr|
 }{|z|^2(-\log|z|)^2\log(-\log|z|)}
\end{equation}
\end{itemize}
\end{prop}
\pf
See the argument in Page 765--767
of \cite{s2}.
\hfill\qed

\subsection{Decompositions}
\label{subsection;08.9.29.4}

As the special case of Lemma \ref{lem;07.7.7.20},
we have the following lemma.
\begin{lem}
Take any large number $N$.
We can take holomorphic sections
$p_{\gminia,N}$ of $\End(\nbigp_c\nbigelambda)$
for $\gminia\in\Irr(\theta)$
such that the following holds:
\[
 p_{\gminia,N}-\pi_{\gminia}
=O\bigl(|z|^{2N} \bigr),
\quad
 \bigl(
 p_{\gminia,N}
\bigr)^2=p_{\gminia,N},
\quad
 \bigl[p_{\gminia_1,N},p_{\gminia_2,N}\bigr]=0,
\quad
 \sum_{\gminia\in\Irr(\theta)}
 p_{\gminia,N}=\id
\]
They preserve the parabolic structure.
\hfill\qed
\end{lem}

We take a refinement of the decomposition:
\begin{lem}
We can take holomorphic sections
$p_{\gminia,\alpha}$
of $\End(\prolong{\nbigelambda})$
for $(\gminia,\alpha)\in \Irr(\theta)\times\Sp(\theta)$
such that the following holds:
\[
 p_{\gminia,\alpha}-\pi_{\gminia,\alpha}
=O\bigl(|z|^{\epsilon}\bigr),
\quad
\bigl[
 p_{\gminia,\alpha},
 p_{\gminib,\beta}
\bigr]=0,
\quad
\bigl(p_{\gminia,\alpha}\bigr)^2
=p_{\gminia,\alpha},
\quad
 \sum_{\alpha}p_{\gminia,\alpha}
=p_{\gminia,N}
\]
\end{lem}
\pf
The argument is essentially
the same as that in the proof of Lemma \ref{lem;07.7.7.20}.
We give only an outline.
Let $d_{\lambda}'':=\delbar_E+\lambda\theta^{\dagger}$,
which is the holomorphic structure of $\nbigelambda$.
We have
$d_{\lambda}''\pi_{\gminia,\alpha}=
O\bigl(|z|^{\epsilon}\bigr)$
with respect to $h$ and the Poincar\'e metric.
By Lemma \ref{lem;07.9.23.36},
we can take sections $s_{\gminia,\alpha}$ of $\End(E)$
satisfying the following:
\[
 d_{\lambda}''s_{\gminia,\alpha}
=d_{\lambda''}\pi_{\gminia,\alpha},
\quad
 \int \bigl|s_{\gminia,\alpha}\bigr|_h^2
 \cdot |z|^{-2\epsilon'}
 \bigl(-\log|z|\bigr)^{N}\cdot\dvol_{g_{\poin}}
\]
By Lemma \ref{lem;07.9.24.5},
we obtain $s_{\gminia,\alpha}=O\bigl(|z|^{\epsilon''}\bigr)$
for some $\epsilon''>0$.

We put 
$\puebar_{\gminia,\alpha}:=
 \pi_{\gminia,\alpha}-s_{\gminia,\alpha}$.
Then, we have the following:
\[
 \bigl(
\puebar_{\gminia,\alpha}\bigr)^2
-\puebar_{\gminia,\alpha}
=O\bigl(|z|^{\epsilon}\bigr),
\quad
 \bigl[
 \puebar_{\gminia,\alpha},
 \puebar_{\gminib,\beta}
\bigr]=O\bigl(|z|^{\epsilon}\bigr),
\quad
 \sum_{\alpha}\puebar_{\gminia,\alpha}
-\puebar_{\gminia,\alpha}
=O\bigl(|z|^{\epsilon}\bigr)
\]
By modifying $\puebar_{\gminia,\alpha}$
with order $|z|^{\epsilon}$
as in the proof of Lemma \ref{lem;07.7.7.20},
we obtain the desired $p_{\gminia,\alpha}$.
\hfill\qed

\subsection{The eigenvalues of $\Res(\DDlambda)$}
\label{subsection;08.12.7.1}

We put
$\Phi:=\sum_{\gminia,\alpha}
 \bigl(d\gminia+\alpha \cdot dz/z\bigr)\cdot 
 p_{\gminia,\alpha}$.

\begin{lem}
\label{lem;07.2.6.4}
Let $v$ be a holomorphic section of $\nbigelambda$
such that 
\[
 |v|_h\sim |z|^{-a}(-\log|z|)^k.
\]
Then, the following estimate holds:
\[
 \Bigl|
 \DDlambda v-(1+|\lambda|^2)\Phi(v)
+\lambda \cdot a\cdot v\cdot \frac{dz}{z}
 \Bigr|_h
=O\bigl(|z|^{-a}(-\log|z|)^{k-1}\bigr)\frac{dz}{z}
\]
\end{lem}
\pf
We have only to consider the case $a=0$.
Let $\delta_{\lambda}'$ be determined by
$d''_{\lambda}$ and $h$.
Then the following holds:
\[
 \DDlambda v-(1+|\lambda|^2)\Phi(v)
=\lambda\delta_{\lambda}'v-(1+|\lambda|^2)
 (\theta-\Phi)(v)
\]
Because
$(\theta-\Phi)(v)=O((-\log|z|)^{k-1})dz/z$,
we have only to show 
\[
  \int|\delta_{\lambda}'v|^2
 (-\log|z|)^{1-\epsilon-2k}|dz\cdot d\zbar|<\infty
\]
for some $\epsilon>0$,
which can be shown using the argument
in Page 761--762 in \cite{s2}.
\hfill\qed

\vspace{.1in}

We take a holomorphic frame $\vecv$
of $\nbigp_a{\nbigelambda}$ such that
(i) it is compatible with the decomposition
$\nbigp_a\nbigelambda=
 \bigoplus \Image p_{\gminia,\alpha}$,
(ii) it is compatible with the parabolic filtration,
(iii) the induced frame of
$\Gr^F_a\bigl(\nbigp\nbigehatlambda_{\gminia}\bigr)$
is compatible with the weight filtration.
We put $a(v_i):=\deg^F(v_i)$.
Let $(\gminia(v_i),\alpha(v_i))$ be determined by
the condition that
$v_i\in \Image p_{\gminia(v_i),\alpha(v_i)}$.

We consider the following:
\[
 \DDlambda_0:=
 \DDlambda-(1+|\lambda|^2)
 \sum_{\gminia\in\Irr(\theta)} d\gminia\cdot p_{\gminia}
\]
\begin{lem}
\label{lem;07.10.6.50}
Let $A$ be the matrix determined by
$\Res(\DDlambda_0)\vecv_{|O}
=\vecv_{|O}\cdot A$.
If the order of $v_1,\ldots,v_r$ is compatible with
the parabolic filtration and the weight filtration,
$A$ is triangulated,
and the $i$-th diagonal entries are
\[
 (1+|\lambda|^2)\alpha(v_i)-\lambda \cdot a(v_i).
\]
\end{lem}
\pf
It follows from Lemma \ref{lem;07.2.6.4}.
\hfill\qed

\subsection{Comparison map}

We put $V_{a,\alpha,\gminia}:=
 \Gr^F_a(\nbigp_c\nbige^0_{\gminia,\alpha})$
on which we have the nilpotent endomorphism
 $N^0_{a,\alpha,\gminia}$.
We take model bundles:
\[
 (\Etilde_{a,\alpha,\gminia},
 \thetatilde_{a,\alpha,\gminia},
 \htilde_{a,\alpha,\gminia})
=\bigl(V_{a,\alpha,\gminia}\otimes\nbigo_{\Delta^{\ast}},
 N^0_{a,\alpha,\gminia}dz/z,
 h_{a,\alpha,\gminia}\bigr)
\otimes
 L\bigl(a,\alpha,\gminia\bigr)
\]
\[
 (\Etilde,\thetatilde,\htilde)
:=\bigoplus_{(a,\alpha,\gminia)}
 \bigl(\Etilde_{a,\alpha,\gminia},
 \thetatilde_{a,\alpha,\gminia},
 \htilde_{a,\alpha,\gminia}\bigr)
\]
Let $\Psi:\prolong{\Etilde}\lrarr \prolong{E}$
be a holomorphic isomorphism such that
(i) it preserves the decompositions,
(ii) it preserves the parabolic filtration,
(iii) $\Gr^F(\Res(G))=0$
where 
$G:=\Psi\circ\thetatilde-\theta\circ\Psi
 \in Hom(\prolong{\Etilde},\prolong{E})
 \otimes\Omega^{1,0}(\log D)$.
Note that $\Psi$ and $\Psi^{-1}$ are bounded,
due to Proposition \ref{prop;07.10.6.40}.
We identify $E$ and $\Etilde$
via $\Psi$ as $C^{\infty}$-bundles.

By construction, we have
$|\theta-\thetatilde|_{\htilde}
=O(|z|^{\epsilon})\cdot dz/z$.
Due to the asymptotic orthogonality
(Theorems \ref{thm;07.10.4.1} and \ref{thm;07.10.4.3}),
we have the following estimate:
\[
 \bigl|\thetatilde^{\dagger}_{\htilde}
-\theta^{\dagger}_{h}\bigr|_{\htilde}
=O\Bigl(\frac{dz}{|z|\bigl(-\log|z|\bigr)}\Bigr)
\]
The following lemma can be shown
by the same argument as 
that in the proof of Lemma 7.7 of \cite{s2}.
\begin{lem}
\label{lem;07.2.6.10}
We have the finiteness
$\int|M|_h^2(-\log |z|)\cdot\dvol_{g_{\poin}}<\infty$,
where $M$ is determined by
$\thetatilde^{\dagger}_{\htilde}
 -\theta^{\dagger}_h=M\cdot d\zbar$.
\hfill\qed
\end{lem}

\subsection{End of the proof of 
Proposition \ref{prop;07.7.19.31}}

We have the induced $\lambda$-connection
$(\nbigelambdatilde,\DDlambdatilde)
=\bigoplus
 (\nbigelambdatilde_{a,\alpha,\gminia},
 \DDlambdatilde_{a,\alpha,\gminia})$,
with the canonical frame
$\vecvtilde=(\vecvtilde_{a,\alpha,\gminia})$.
We put $\ptilde_{\gminia,\alpha}:=
 \pitilde_{\gminia,\alpha}
=\pi_{\gminia,\alpha}$.
We have the decomposition
$\nbigelambdatilde=\bigoplus
 \Image\ptilde_{\gminia,\alpha}$.

Let $\vecv$ be as in Subsection
\ref{subsection;08.12.7.1}.
Let $\deg^F(v_i)$ and $\deg^W(v_i)$
denote the degree of $v_i$ with respect to 
the parabolic filtration and the weight filtration.
Let $(\gminia(v_i),\alpha(v_i))$ be determined by the condition
that $v_i\in\Image p_{\gminia(v_i),\alpha(v_i)}$.
We use the symbols $\deg^F(\vtilde_i)$,
$\deg^W(\vtilde_i)$ and 
$(\gminia(\vtilde_i),\alpha(\vtilde_i))$,
in similar meanings.
Let $I=(I_{j,i})$ be determined by the relation
$\vtilde_i=\sum_{j} I_{j,i}\cdot v_j$.

\begin{lem}
\label{lem;07.10.6.51}
We put $\nbigb_{j,i}:=
 I_{j,i}|z|^{-\deg^F(v_j)+\deg^F(\vtilde_i)
 +(\deg^W(v_j)-\deg^{W}(\vtilde_i))/2}$.
\begin{itemize}
\item
In the case $(\gminia(v_j),\alpha(v_j))\neq 
 (\gminia(\vtilde_i),\alpha(\vtilde_i))$,
we have $\bigl|\nbigb_{j,i}\bigr|\leq C\cdot |z|^{\epsilon}$
for some $\epsilon>0$ and $C>0$.
\item
In the case $\deg^{F}(v_j)\neq \deg^{F}(\vtilde_i)$,
we have $\bigl|\nbigb_{j,i}\bigr|\leq C(-\log|z|)^{-1}$
for some $C>0$.
\item
In the case $\deg^W(v_j)\neq \deg^W(\vtilde_i)$,
we have $\|\nbigb_{j,i}(-\log|z|)\|_W<\infty$,
where $\|\cdot\|_W$ is given as in
{\rm(\ref{eq;08.9.13.30})}.
\end{itemize}
\end{lem}
\pf
We have the estimates
$\ptilde_{\gminia,\alpha}-p_{\gminia,\alpha}
=O\bigl(|z|^{\epsilon}\bigr)$ with respect to $h$.
Hence, the first claim follows.
The other claims can be shown using the argument
in \cite{s2},
using Lemma \ref{lem;07.2.6.10} and
Proposition \ref{prop;07.2.6.31}.
(See also \cite{mochi2}.)
\hfill\qed

\vspace{.1in}

Then, we can show the following equality
by the argument 
in the proof of Proposition 7.6 and Theorem 7 of \cite{s2},
using Lemma \ref{lem;07.10.6.50}
and Lemma \ref{lem;07.10.6.51}:
\begin{equation}
 \label{eq;07.10.6.52}
 \dim \Gr^W_k\EE_{\beta}
  \Gr^F_b(\nbigelambda_{\gminib})
=\dim \Gr^W_k\EE_{\beta}
  \Gr^F_b(\nbigelambdatilde_{\gminib})
\end{equation}
We can check the claim of Proposition \ref{prop;07.7.19.31}
for $(\Etilde,\delbar_{\Etilde},\thetatilde,\htilde)$
by a direct calculation.
Together with (\ref{eq;07.10.6.52}),
we obtain
the claims of Proposition \ref{prop;07.7.19.31}
for $(E,\delbar_E,\theta,h)$.
\hfill\qed

\section{A characterization of the lattices
 for generic $\lambda$}
\label{subsection;08.9.14.32}
Let $X$, $D$ and $\harmonicbundle$
be as in Subsection \ref{subsection;08.12.7.2}.
We have the set
$\KMS(E)\subset\real\times\cnum$ 
of the KMS-spectra at $\lambda=0$.
A complex number $\lambda$ is called generic
with respect to $\KMS(E)$,
if $\eigenmap(\lambda):\KMS(E)\lrarr\cnum$
is injective.
For generic $\lambda$,
we have the following characterization of
$\nbigp_a\nbigelambda$ ($a\in\real$).

\begin{prop}
\label{prop;08.1.22.2}
Let $\lambda$ be generic with respect to
$\KMS(E)$.
Let $V$ be a good lattice of
the meromorphic $\lambda$-flat bundle
$(\nbigp\nbigelambda,\DDlambda)$
with the following property:
\begin{itemize}
\item
 The set of the eigenvalues of
 $\Res(\DDlambda)$ on $V_{|D}$
 is the same as the following:
\begin{equation}
 \label{eq;08.1.22.1}
  \bigl\{
 \eigenmap(\lambda,u)\,\big|\,
 u\in \KMS(E),\,\,a-1<\paramap(\lambda,u)\leq a
 \bigr\}
\end{equation}
Then, $V=\nbigp_a\nbigelambda$.
\end{itemize}
\end{prop}
\pf
Let $S$ denote the set (\ref{eq;08.1.22.1}).
Note that $\lambda^{-1}(\alpha-\beta)\not\in\seisuu$
for any distinct $\alpha,\beta\in S$.
Hence, we have the flat decompositions:
\[
(\nbigp_a\nbigelambda,\DDlambda)_{|\Dhat}
=\bigoplus_{\gminia\in\Irr(\theta)}
 \bigoplus_{\alpha\in S}
 \bigl(\nbigp_a\nbigelambda_{\gminia,\alpha},\,
 \DDlambda_{\gminia,\alpha}
 \bigr)
\quad\quad
 (V,\DDlambda)_{|\Dhat}=
 \bigoplus_{\gminia\in \Irr(\theta)}
 \bigoplus_{\alpha\in S}
\bigl(
 \Vhat_{\gminia,\alpha},\,
 \DDlambda_{\gminia,\alpha}
\bigr)
\]
Here, 
 $\DDlambda_{\gminia,\alpha}-d\gminia
 -\alpha\cdot dz/z$ are logarithmic 
 with respect to 
 $\nbigp_{a}\nbigelambda_{\gminia,\alpha}$
 or $\Vhat_{\gminia,\alpha}$,
and the residues are nilpotent.
Then, we obtain
$\nbigp_a\nbigelambda_{\gminia,\alpha}
=\Vhat_{\gminia,\alpha}$
using a standard and classical argument:
(i) We can show
 $\bigoplus_{\alpha}
 \nbigp_a\nbigelambda_{\gminia,\alpha}
 \otimes\nbigo(\ast D)
=\bigoplus_{\alpha}
 \Vhat_{\gminia,\alpha}\otimes\nbigo(\ast D)$
 by using Corollary \ref{cor;10.5.2.4}.
(ii)
 By the assumption for $S$,
 we can also obtain 
$\nbigp_a\nbigelambda_{\gminia,\alpha}
 \otimes\nbigo(\ast D)
=\Vhat_{\gminia,\alpha}\otimes\nbigo(\ast D)$
 by a similar argument with minor modification.
(iii)
 In the regular case,
 this kind of claim is well known.
 (See Proposition II. 5.4 of \cite{d},
 for example.)
\hfill\qed

\begin{rem}
The proposition can easily be
generalized to the higher dimensional case.
\hfill\qed
\end{rem}

\section{Harmonic forms}
\label{subsection;08.9.14.33}
\subsection{The space of
harmonic forms of $\nbigelambda$}
\label{subsection;08.1.20.13}

Let $Y$ be a smooth projective curve,
and let $D$ be a finite subset of $Y$.
Let $\harmonicbundle$ be a wild harmonic bundle
on $Y-D$.
Let $\omega$ be a Kahler form of $Y-D$,
which is Poincar\'e like around $D$.
We reformulate the result in 
Section \ref{subsection;08.1.20.11} below.
Let $\DDlambdaast$ denote the formal adjoint
of $\DDlambda$ with respect $h$ and $\omega$.
We put 
$\Deltalambda:=
 \DDlambda\circ\DDlambdaast
+\DDlambdaast\circ\DDlambda$.
Recall
$\Delta^{\lambda}=(1+|\lambda|^2)\cdot\Delta^0$.

\begin{prop}
\label{prop;08.1.20.11}
Let $\phi$ be an $L^2$-section of
$E$ on $Y-D$.
We have the following equivalence:
\begin{multline*}
\mbox{\rm $\Deltalambda\phi=0$ for some $\lambda$}
\Longleftrightarrow
\mbox{\rm $\Deltalambda\phi=0$ for any $\lambda$}
 \\
\Longleftrightarrow
\mbox{\rm $\DDlambda\phi=0$ for some $\lambda$}
\Longleftrightarrow
\mbox{$\DDlambda\phi=0$ for any $\lambda$}
\end{multline*}
Let $\phi$ be an $L^2$-section of
$E\otimes\Omega^{1,1}$ on $Y-D$.
We have the following equivalence:
\begin{multline*}
\mbox{$\Deltalambda\phi=0$ for some $\lambda$}
\Longleftrightarrow
\mbox{$\Deltalambda\phi=0$ for any $\lambda$}
 \\
\Longleftrightarrow
\mbox{$\DDlambdaast\phi=0$ for some $\lambda$}
\Longleftrightarrow
\mbox{$\DDlambdaast\phi=0$ for any $\lambda$}
\end{multline*}
\end{prop}
\pf
 Let us show the claim for $0$-forms.
 According to Proposition \ref{prop;08.1.20.10} below,
 we have the equivalence
 $\Deltalambda\phi=0
 \Longleftrightarrow
 \DDlambda\phi=0$ for a fixed $\lambda$.
 Then, the desired equivalence
 follows from $\Deltalambda=(1+|\lambda|^2)\Delta^0$.
The claim for $2$-forms
 can be reduced to the case of $0$-forms.
\hfill\qed

\begin{prop}
\label{prop;08.1.20.12}
Let $\phi$ be an $L^2$-section of $E\otimes\Omega^1$
on $Y-D$.
The following conditions are equivalent.
\begin{description}
\item[(a)]
 $\Deltalambda\phi=0$ for some $\lambda$.
\item[(b)]
 $\Deltalambda\phi=0$ for any $\lambda$.
\item[(c)]
 $\DDlambda\phi=\DDlambdaast\phi=0$
 for some $\lambda$.
\item[(d)]
 $\DDlambda\phi=\DDlambdaast\phi=0$
 for any $\lambda$.
\item[(e)]
 $(\delbar+\theta)\phi
=(\del+\theta^{\dagger})\phi=0$.
\end{description}
\end{prop}
\pf
The equivalence of the conditions
(a)--(d) can be shown
using an argument in the proof of
Proposition \ref{prop;08.1.20.11}.
Note 
$\DD^{0\ast}=\sqrt{-1}\Lambda_{\omega}
 \bigl(\del+\theta^{\dagger}\bigr)$
on the one forms.
Hence, the condition (e) is equivalent to
(c) with $\lambda=0$.
\hfill\qed

\begin{rem}
The equivalence of 
{\rm (c)} and {\rm (d)} in
Proposition {\rm\ref{prop;08.1.20.12}}
can be shown by a direct calculation.
Note that $\DDlambdaast\phi=0$
is equivalent to $\DDlambdastar\phi=0$
for a $1$-form $\phi$.
(See {\rm\cite{mochi5}} for the notation.)
Let $A:=\delbar+\theta$
and $B:=\del+\theta^{\dagger}$.
Then, we have the equalities
$(1+|\lambda|^2)A=
 \DDlambda-\lambda\DDlambdastar$
and
$(1+|\lambda|^2)B=
 \DDlambdastar+\lambda\DDlambda$.
Hence, $\DD^{\lambda_1}$
and $\DD^{\lambda_1\star}$
can be expressed as 
the linear combinations of
$\DDlambda$ and $\DDlambdastar$
for any $\lambda_1$.
\hfill\qed
\end{rem}

The following notation
will be used in Section 
\ref{subsection;07.10.23.30}.
\begin{notation}
Let $\Harm^i$ denote the space
of the $L^2$-sections of $E\otimes\Omega^i$
satisfying the conditions in 
Propositions {\rm\ref{prop;08.1.20.11}}
and {\rm\ref{prop;08.1.20.12}}.
\hfill\qed
\end{notation}

Needless to say,
the equivalence in the propositions
does not hold if $Y$ is not projective.
Hence, we have to distinguish the conditions.

\subsection{Decay of $L^2$-harmonic one forms
 around the singularity}
\label{subsection;08.9.29.10}

Let us study the behaviour of
harmonic one forms around the singularity
in more details.
Let $X$, $D$ and $\harmonicbundle$
be as in Subsection \ref{subsection;08.12.7.2}.
We use the Poincar\'e metric of $X-D$.
Let $\tau$ be an $L^2$-section of $E\otimes\Omega^1$
on $X^{\ast}$ such that 
$(\delbar_E+\theta)\tau=
 (\del_E+\theta^{\dagger})\tau=0$.
We have the decomposition
$\tau=\sum_{\gminia\in\Irr(\theta)}\tau_{\gminia}$
corresponding to the decomposition
$E=\bigoplus_{\gminia\in\Irr(\theta)} E_{\gminia}$.
We will use the following proposition
in Section \ref{subsection;07.10.23.30}.

\begin{prop}
\label{prop;07.7.19.15}
$\tau$ is of polynomial order with respect to $h$
and the Poincar\'e metric.
For $\gminia\neq 0$,
we have the estimate
$\tau_{\gminia}=O\bigl(
 \exp(-\epsilon|z|^{\ord(\gminia)})
 \bigr)$.
\end{prop}
\pf
For simplicity,
we use the symbol $\del A/\del z$
to denote the coupling of $\del A$ and $\del_z$
for a section $A$ of $E$.
For the expression $\tau=A\cdot dz+B\cdot d\zbar$,
we have the following equalities:
\begin{equation}
 \label{eq;07.7.18.110}
 \frac{\del A}{\del \zbar}-f(B)=0,
\quad
 \frac{\del B}{\del z}-f^{\dagger}(A)=0.
\end{equation}
\begin{lem}
\label{lem;07.7.19.11}
$A$ and $B$ are of polynomial order
with respect to $h$.
\end{lem}
\pf
Let $\vecv$ be a frame of $\prolong{E}$
compatible with the decomposition (\ref{eq;07.7.18.100}),
the parabolic filtration and the weight filtration.
Let $\Theta$, $\Theta^{\dagger}$
and $C$ be the matrix-valued functions
determined by the conditions:
\[
  \theta \vecv=\vecv\cdot \Theta\cdot dz,
\quad
 \theta^{\dagger}\vecv=\vecv\cdot \Theta^{\dagger}d\zbar,
\quad
 \del\vecv=\vecv\cdot C\cdot dz.
\]
Then, $\Theta$, $\Theta^{\dagger}$ and $C$
are of polynomial order,
which follows from the estimate for
the Higgs field 
(Theorems \ref{thm;07.10.4.1}
 and \ref{thm;07.10.4.3})
and the acceptability of $(E,\delbar_E,h)$
(Theorem \ref{thm;07.10.9.1}
and Lemma \ref{lem;07.6.2.20}).
Let
$\vecA=(A_i)$ and $\vecB=(B_i)$
be $\cnum^{\rank E}$-valued functions determined by 
$ A=\sum A_i\cdot v_i$ and $B=\sum B_i\cdot v_i$.
Due to (\ref{eq;07.7.18.110}),
the following equalities hold on $\Delta^{\ast}$:
\[
 \frac{\del \vecA}{\del \zbar}-\Theta(\vecB)=0,
\quad
 \frac{\del \vecB}{\del z}+[C,\vecB]
-\Theta^{\dagger}(\vecA)=0
\]
Hence, we have the following 
for any non-negative integer $l\geq 0$:
\begin{equation}
\label{eq;07.7.18.111}
 \frac{\del (z^l\vecA)}{\del \zbar}
-\Theta(z^l\vecB)=0,
\quad
 \frac{\del (z^l\vecB)}{\del z}
-lz^{l-1}\vecB
+z^l[C,\vecB]
-z^l\Theta^{\dagger}(\vecA)=0
\end{equation}
If $l$ is sufficiently large,
(\ref{eq;07.7.18.111}) holds on $X$ as distributions.
For large $N$,
$z^N\vecA$ and $\Theta(z^N\vecB)$ are $L^p$
for some $p>0$.
Then $z^N\vecA$ is $L_1^p$ for some $p>2$
because of the first equality in (\ref{eq;07.7.18.111}),
and hence $z^N\vecA$ is bounded.
Thus, $A$ is shown to be of polynomial order.
By applying a similar argument
to $(E,\del_E,\theta^{\dagger},h)$,
it can be shown that
$B$ is also of polynomial order.
\hfill\qed

\vspace{.1in}

We give a refinement.
We put $\nbigs_1(j):=\bigl\{
 \gminia\in\Irr(\theta)\,\big|\,
 \ord(\gminia)\leq j\bigr\}$
and $\nbigs_{0}(j):=\bigl\{
 \gminia\in\Irr(\theta)\,\big|\,
 \ord(\gminia)>j \bigr\}$.
We put $E^{(j)}_a:=\bigoplus_{\gminia\in \nbigs_a(j)}E_{\gminia}$
for $a=0,1$.
Let $B=B_1^{(j)}+B_0^{(j)}$
and $A=A_1^{(j)}+A_0^{(j)}$
be the decomposition corresponding to
$E=E^{(j)}_1\oplus E^{(j)}_0$.
For a differential operator $\gbigd:E\lrarr E$,
we have the decomposition
$\gbigd=\sum_{a,b=1,2} \gbigd_{a,b}$,
where $\gbigd^{(j)}_{a,b}:E^{(j)}_a\lrarr E^{(j)}_b$.
We put $\nbigd^{(j)}(\gbigd):=\gbigd^{(j)}_{1,1}+\gbigd^{(j)}_{0,0}$
and $\nbigc^{(j)}(\gbigd):=\gbigd^{(j)}_{1,0}+\gbigd^{(j)}_{0,1}$.

\begin{lem}
\label{lem;07.7.19.10}
$A_1^{(j)}=O\bigl(\exp(-\epsilon|z|^{j})\bigr)$ 
for some $\epsilon>0$, with respect to $h$.
\end{lem}
\pf
Let us look at
the $E_1^{(j)}$-component of the equality
$\del B/\del z-f^{\dagger}(A)=0$.
Due to Lemma \ref{lem;07.6.2.25} and
the $\exp(-\epsilon|z|^{j})$-asymptotic orthogonality
of the decomposition
  $E=\bigoplus_{\gminia\in\Irr(\theta,j)}
     E^{(j)}_{\gminia}$ (Theorem \ref{thm;07.10.4.1}),
we obtain the following:
\begin{multline}
 0=\nbigd^{(j)}(\del_z)B_1^{(j)}
+\nbigc^{(j)}(\del_z) B_0^{(j)}
-\nbigd^{(j)}(f^{\dagger})A_1^{(j)}
-\nbigc^{(j)}(f^{\dagger})A_0^{(j)} \\
=\nbigd^{(j)}(\del_z) B_1^{(j)}
-\nbigd^{(j)}(f^{\dagger})A_1^{(j)}
+O\bigl(\exp(-C|z|^{j})\bigr)
\end{multline}
We have the unique $L^2$-section
$g_1^{(j)}$ of $E_1^{(j)}$ such that
$\theta \bigl(g_1^{(j)}\bigr)
=A_1^{(j)}\cdot dz$.
Then, we have
$ (\delbar+\theta)g_1^{(j)}
=A_1^{(j)}dz+B_1^{(j)}d\zbar$
on $X-D$,
i.e.,
$\del g_1^{(j)}/\del \zbar
=B_1^{(j)}$
and
$f(g_1^{(j)})=A_1^{(j)}$.
Since $A_1^{(j)}$ is of polynomial order,
we obtain that
$g_1^{(j)}$ is of polynomial order.
We have the following equality on $X-D$:
\[
 -\frac{\del^2}{\del z\del \zbar}
 |g_1^{(j)}|^2_h
=-\Bigl(
 \frac{\del^2 g_1^{(j)}}{\del z\del\zbar},g_1^{(j)}
 \Bigr)_h
-\Bigl|\frac{\del g_1^{(j)}}{\del z}\Bigr|_h^2
-\Bigl|\frac{\del g_1^{(j)}}{\del \zbar}\Bigr|_h^2
-\Bigl(
 g_1^{(j)},\frac{\del^2 g_1^{(j)}}{\del \zbar\del z}
 \Bigr)_h
\]
Since $A_1^{(j)}$ and $B_1^{(j)}$ are of polynomial order,
we obtain the following:
\begin{multline}
 \left(
 \frac{\del^2 g_1^{(j)}}{\del z\del\zbar},
 g_1^{(j)}\right)_h
=\left(
 \nbigd^{(j)}(\del_ z)\frac{\del g_1^{(j)}}{\del \zbar},
 g_1^{(j)}
 \right)_h
+O\bigl(\exp(-C|z|^{j})\bigr)
 \\
=\bigl(\nbigd^{(j)}(f^{\dagger})f(g_1^{(j)}),
 g_1^{(j)}\bigr)_h
+O\bigl(\exp(-C|z|^{j})\bigr)
=\bigl(f^{\dagger}f(g_1^{(j)}),g_1^{(j)}\bigr)_h
+O\bigl(\exp(-C|z|^{j})\bigr) \\ 
=\bigl(f(g_1^{(j)}),f(g_1^{(j)})\bigr)_h
+O\bigl(\exp(-C|z|^{j})\bigr)
\end{multline}
We also have the following:
\[
 \Bigl(
 \frac{\del^2g_1^{(j)}}{\del z\del \zbar}
-\frac{\del^2g_1^{(j)}}{\del \zbar\del z}
 \Bigr)dz \cdot d\zbar
=\bigl(\del\delbar+\delbar\del\bigr)g_1^{(j)}
=-\bigl(\theta\theta^{\dagger}
   +\theta^{\dagger}\theta\bigr)g_1^{(j)}
=-\bigl(ff^{\dagger}-f^{\dagger}f\bigr)
 g_1^{(j)}\cdot dz\cdot d\zbar
\]
Therefore, we obtain the following:
\begin{multline}
 -\frac{\del^2}{\del z\del\zbar}|g_1^{(j)}|_h^2= \\
-\Bigl(
 \frac{\del^2 g_1^{(j)}}{\del z\del\zbar},g_1^{(j)}
 \Bigr)_h
-\Bigl(g_1^{(j)},
 \frac{\del^2 g_1^{(j)}}{\del z\del\zbar}\Bigr)_h
-\bigl(g_1^{(j)},
 (ff^{\dagger}-f^{\dagger}f)g_1^{(j)}
 \bigr)_h
-\Bigl|\frac{\del g_1^{(j)}}{\del z} \Bigr|_h^2
-\Bigl|\frac{\del g_1^{(j)}}{\del \zbar}\Bigr|_h^2 \\
=-\bigl|f (g_1^{(j)})\bigr|_h^2
   -\bigl|f^{\dagger}(g_1^{(j)})\bigr|_h^2
-\Bigl|\frac{\del g_1^{(j)}}{\del z} \Bigr|_h^2
-\Bigl|\frac{\del g_1^{(j)}}{\del \zbar}\Bigr|_h^2 
+O\bigl(\exp(-C|z|^{j})\bigr) \\
\leq
-C_1|z|^{2(j-1)}\cdot\bigl|g_1^{(j)}\bigr|_h^2
+C_2\exp(-C_3|z|^{j})  
\end{multline}
\begin{lem}
\label{lem;07.7.19.3}
The following inequality holds on $X$
as distributions:
\begin{equation}
 \label{eq;07.7.19.2}
-\frac{\del^2}{\del z\del \zbar}|g_1^{(j)}|_h^2
\leq
-C_1\cdot |z|^{2(j-1)}\cdot |g_1^{(j)}|_h^2
+C_2\cdot \exp\bigl(-C_3|z|^{j}\bigr)
\end{equation}
\end{lem}
\pf
We have already known that the inequality holds on $X-D$.
Let $\varphi$ be a test function.
We have the following:
\begin{multline*}
 \int_{|z|\geq\delta}
 \bigl|g_1^{(j)}\bigr|_h^2\frac{\del^2\varphi}{\del z\del\zbar}
 dz\cdot d\zbar
=\\
\pm\int_{|z|=\delta}
 |g_1^{(j)}|_h^2\frac{\del \varphi}{\del\zbar}d\zbar
\pm\int_{|z|=\delta}
 \frac{\del |g_1^{(j)}|^2}{\del z}\varphi\cdot dz
+\int_{|z|\geq\delta}
 \frac{\del^2|g_1^{(j)}|_h^2}{\del z\del\zbar}\varphi
 \cdot dz\cdot d\zbar
\end{multline*}
Hence, we have only to show the existence of a sequence
$\{\delta_i\}$ with $\delta_i\lrarr 0$ 
such that the following holds:
\[
\lim \int_{|z|=\delta_i}
 |g_1^{(j)}|_h^2\cdot\frac{\del\varphi}{\del \zbar}\cdot d\zbar=0,
\quad
 \lim\int_{|z|=\delta_i}
 \frac{\del |g_1^{(j)}|_h^2}{\del z}\varphi\cdot dz=0.
\]
Let us show the second convergence. The first one can be shown
by a similar argument.
By construction, we have the following finiteness:
\[
 \int |g_1^{(j)}|_h^2\cdot |z|^{2(j-1)}\cdot
 |dz\cdot d\zbar|<\infty,
\quad
 \int\bigl(\delbar g_1^{(j)},\delbar g_1^{(j)}\bigr)_h<\infty
\]
Let $\rho$ be a non-negative $C^{\infty}$-function on $\real$
such that $\rho(t)=1$ for $t\leq 1/2$ and $\rho(t)=0$ for $t\geq 2/3$.
We put $\chi_N(z):=\rho\bigl(-N^{-1}\log|z|\bigr)$.
Note $\del \chi_N(z)$ and $\delbar \chi_N(z)$ are uniformly bounded
with respect to the Poincar\'e metric.
We have the following estimate, which is independent of $N$:
\begin{multline*}
 \int\bigl(\del(\chi_N \cdot g_1^{(j)}),
 \del(\chi_N\cdot g_1^{(j)})\bigr)
=\\
\pm\int\bigl(\chi_N\!\cdot\! g_1^{(j)},\,
 R(h)\!\cdot\!\chi_N\!\cdot\! g_1^{(j)}\bigr)
\pm\int\bigl(\delbar(\chi_N\cdot g_1^{(j)}),
 \delbar(\chi_N\cdot g_1^{(j)})\bigr)
<C
\end{multline*}
Thus, we obtain the finiteness
$\int\bigl(\del g_1^{(j)},\del g_1^{(j)}\bigr)<\infty$.
Then, we obtain the following finiteness:
\begin{multline}
 \label{eq;07.7.19.1}
 \int r^{j-1}dr
\left(
 \int \left|\frac{\del |g_1^{(j)}|_h^2}{\del z}\varphi\right|
\cdot r\cdot d\theta
\right)
\leq  \\
 C
\left(
 \int\bigl(\del g_1^{(j)},\del g_1^{(j)}\bigr)_h
\right)^{1/2}
 \left(
 \int|g_1^{(j)}|_h^2
 |z|^{2(j-1)}\cdot
 |dz\cdot d\zbar|
 \right)^{1/2}
<\infty
\end{multline}
The existence of the desired sequence $\{\delta_i\}$
follows from (\ref{eq;07.7.19.1}).
Thus, the proof of Lemma \ref{lem;07.7.19.3} is finished.
\hfill\qed

\vspace{.1in}

Let us return to the proof of Lemma \ref{lem;07.7.19.10}.
In general, we have the following inequality for $l>0$:
\[
 -\frac{\del^2}{\del z\del\zbar}
 \exp\bigl(-C|z|^{-l}\bigr)
=\exp\bigl(-C|z|^{-l}\bigr)\frac{l^2\cdot C}{4}|z|^{-l-2}
-\exp\bigl(-C|z|^{-l}\bigr)\frac{(l\cdot C)^2}{4}|z|^{-2l-2}
\]
Hence, we obtain the following inequality
for appropriate constant $G>0$:
\begin{equation}
 \label{eq;07.7.19.4}
 -\frac{\del^2}{\del z\del\zbar}
 \exp\bigl(-G|z|^{j}\bigr)
\geq
 -C_1\cdot |z|^{2(j-1)}\exp\bigl(-G|z|^{j}\bigr)
+C_2\cdot\exp(-C_3|z|^{j})
\end{equation}
We obtain
$|g_1^{(j)}|_h\leq C_4\exp(-G|z|^{j})$,
by using (\ref{eq;07.7.19.4}), 
Lemma \ref{lem;07.7.19.3}
and the standard argument 
as in \cite{a} and \cite{s2}.
(See also the proof of Theorem \ref{thm;07.10.4.1}
and Theorem \ref{thm;07.10.4.3}.)
Because
$\bigl|A_1^{(j)}\bigr|_h\leq
 C_6|z|^{j-1}\cdot |g_1^{(j)}|_h$,
the proof of Lemma \ref{lem;07.7.19.10} 
is accomplished.
\hfill\qed

\begin{lem}
\label{lem;07.7.19.12}
$B_1^{(j)}=O\bigl(\exp(-C|z|^{j})\bigr)$
with respect to $h$.
\end{lem}
\pf
We can apply a similar argument to
a harmonic bundle $(E,\del_E,h,\theta^{\dagger})$
on $(X-D)^{\dagger}$.
We have the decomposition similar to (\ref{eq;07.7.18.100}):
\[
 (E,f^{\dagger})=
 \bigoplus_{(\gminiabar,\alphabar)
   \in\Irr(\theta^{\dagger})\times\Sp(\theta^{\dagger})}
 \bigl(E^{\dagger}_{\gminiabar,\alphabar},
 f^{\dagger}_{\gminiabar,\alphabar}\bigr)
=\bigoplus_{\gminiabar
   \in\Irr(\theta^{\dagger})}
 \bigl(E^{\dagger}_{\gminiabar},
 f^{\dagger}_{\gminiabar}\bigr)
\]
We put
$\nbigs^{\dagger}_1(j):=\bigl\{
 \gminiabar\in\Irr(\theta^{\dagger})\,\big|\,
 \ord(\gminiabar)\leq j \bigr\}$
and $\nbigs^{\dagger}_0(j):=
 \bigl\{\gminiabar\in\Irr(\theta^{\dagger})\,\big|\,
 \ord(\gminiabar)>j \bigr\}$.
As in the previous argument,
we put $E^{\dagger(j)}_a:=
 \bigoplus_{\gminiabar\in\nbigs^{\dagger}_a(j)}
 E^{\dagger}_{\gminiabar}$.
By Lemma \ref{lem;07.7.19.10},
we have the corresponding decomposition
$B=B_0^{\dagger(j)}+B_1^{\dagger(j)}$,
and we obtain $B_1^{\dagger(j)}=O(\exp(-C|z|^{j}))$
with respect to $h$.

Let $\pi^{\dagger(j)}_a$ denote the projection onto
$E^{\dagger(j)}_a$ with respect to the decomposition 
$E=E^{\dagger(j)}_0\oplus E^{\dagger(j)}_1$.
Let $\pi^{(j)}_a$ denote the projection onto
$E^{(j)}_a$ with respect to the decomposition
$E=E^{(j)}_0\oplus E^{(j)}_1$.
According to  Theorem \ref{thm;07.10.4.1},
we have the estimate
$\pi^{(j)}_a-\pi^{(j)\dagger}_a
=O\Bigl(\exp\bigl(-\epsilon|z|^{j}\bigr)\Bigr)$
with respect to $h$.
We have the following:
\[
 B_1^{(j)}=\pi_1^{(j)}(B)
=B_1^{\dagger(j)}
+\bigl(\pi_1^{(j)}-\pi_1^{\dagger(j)}\bigr)(B)
\]
Since $B$ is of polynomial order,
we obtain the desired estimate for $B_1^{(j)}$.
Thus, the proof of Lemma \ref{lem;07.7.19.12} is finished
\hfill\qed

\vspace{.1in}
Now, the claim of Proposition \ref{prop;07.7.19.15}
immediately follows from Lemmas \ref{lem;07.7.19.11},
\ref{lem;07.7.19.10} and \ref{lem;07.7.19.12}.
\hfill\qed

\begin{rem}
Let $Y$, $D$, $\omega$ and $\harmonicbundle$
be as in Subsection {\rm\ref{subsection;08.1.20.13}}.
Due to Proposition {\rm\ref{prop;08.1.20.12}},
we can apply Proposition {\rm\ref{prop;07.7.19.15}}
to study the behaviour of
any $L^2$-harmonic one form 
of $(\nbigelambda,\DDlambda)$ around $D$,
after taking an appropriate ramified covering.
\hfill\qed
\end{rem}

\begin{rem}
Let $\tau$ be an $L^2$-section of $E$
such that $(\delbar+\theta)\tau=0$.
For the decomposition
$\tau=\sum \tau_{\gminia}$
corresponding to $E=\bigoplus E_{\gminia}$,
we obviously have the vanishings
$\tau_{\gminia}=0$
for $\gminia\neq 0$.
We also have a similar claim
for $2$-forms.
\hfill\qed
\end{rem}

\begin{rem}
In {\rm\cite{szabo}},
the exponential decay 
of harmonic forms is shown in a special case 
with a different method.
\hfill\qed
\end{rem}

\subsection{General remarks on
 $L^2$-harmonic forms (Appendix)}
\label{subsection;08.9.29.11}

\subsubsection{$L^2$-cohomology and harmonic forms}
\label{subsection;08.1.20.15}

We recall a general remark
on the $L^2$-cohomology,
following \cite{z}.
Let $(Y,g)$ be a complete Kahler manifold.
Let $\dvol_g$ denote the volume form associated to $g$.
Let $(V,\DDlambda)$ be a $\lambda$-flat bundle
on $Y$ with a hermitian metric $h$.
(See \cite{s4} and \cite{mochi5}.)
Let $\nbigl^j(V)$ be the space of
sections $f$ of $V\otimes\Omega^j$ 
on $Y$,
such that $f$ and $\DDlambda f$ are $L^2$
with respect to $h$ and $g$.
(Here, $\DDlambda f$ is taken
in the sense of distributions.
But, we do not have to be concerned with it,
because $(Y,g)$ is complete.
See \cite{av}.)
Thus, we obtain a complex
$(\nbigl^{\bullet}(V),\DDlambda)$.
Let $H^j\bigl(\nbigl^{\bullet}(V)\bigr)$
denote the $j$-th cohomology group
of the complex.

Let $\DDlambdaast$ denote the formal adjoint
of $\DDlambda$ with respect to $h$ and $g$.
Let $\vecH^j$ denote the space of
$L^2$-sections $f$ of
$V\otimes \Omega^j$ satisfying
$\DDlambda f=0$ and $\DDlambdaast f=0$.

\begin{lem}
\label{lem;08.1.19.111}
Assume that
$\bigoplus_j H^j\bigl(\nbigl^{\bullet}(V)\bigr)$
is finite dimensional.
Then, the natural map
$\bigoplus_j\vecH^j\lrarr 
 \bigoplus_jH^j\bigl(\nbigl^{\bullet}(V)\bigr)$
is an isomorphism.
\end{lem}
\pf
We use the inner product
$(f,g):=\int h(f,g)\cdot\dvol_g+
 \int h(\DDlambda f,\DDlambda g)\cdot\dvol_g$,
via which $\nbigl^j(V)$ is the Hilbert space.
Let $Z^j$ denote the kernel of
$\DDlambda:\nbigl^{j}(V)\lrarr
 \nbigl^{j+1}(V)$,
which is the closed subspace of
$\nbigl^{j}(V)$.
Let $Z^{j\bot}$ denote the orthogonal complement of
$Z^j$ in $\nbigl^{j}(V)$.
Let us see the continuous operator
$\DDlambda:\nbigl^{j-1}(V)
\lrarr Z^j$.
Because we have assumed
$\dim H^j\bigl(\nbigl^{\bullet}(V)\bigr)<\infty$,
the image $R^j$ of $\DDlambda$ is closed.
Let $\vecH_1^j$ denote the orthogonal complement
of $R^j$ in $Z^j$.
Then, we have the orthogonal decomposition
$\nbigl^j(V)
=R^j\oplus \vecH_1^j\oplus Z^j$,
and $\vecH_1$ is naturally isomorphic to
$H^j\bigl(\nbigl^{\bullet}(V)\bigr)$.

Let us show $\vecH_1^j=\vecH^j$.
Let $f\in \vecH_1^j$.
For any $C^{\infty}$-section $\varphi$
of $V\otimes\Omega^{j-1}$ with compact support,
we have 
$\int\bigl(\DDlambda\varphi,f\bigr)_{h,\omega}
 \dvol_{\omega}=0$.
Then, we obtain $\DDlambdaast f=0$
in the sense of distributions,
i.e., $f\in \vecH^j$.
For any $f\in \vecH_1^j$ and for any
$\varphi$ as above,
we have 
$\int\bigl(\DDlambda\varphi,f\bigr)_{h,\omega}=
 \int\bigl(\varphi,\DDlambdaast f\bigr)_{h,\omega}=0$,
and hence 
$f\in \vecH_1^j$.
Thus, we obtain Lemma \ref{lem;08.1.19.111}.
\hfill\qed

\subsubsection{Harmonic forms 
on a complete Kahler manifold}

Let $(Y,g)$ and $(V,\DDlambda)$
be as in Section \ref{subsection;08.1.20.15}.
Let $\Deltalambda:=
\DDlambda\circ\DDlambdaast
+\DDlambdaast\circ\DDlambda$.
We recall the following general remark.
\begin{lem}
\label{lem;08.1.20.16}
Let $\phi$ be an $L^2$-section of
$V\otimes\Omega^j$
such that $\DDlambda\phi$ and $\DDlambdaast\phi$
are $L^2$.
Then, $\Deltalambda\phi=0$
if and only if 
$\DDlambda\phi=\DDlambdaast\phi=0$.
\end{lem}
\pf
We have only to show the ``only if'' part.
We fix a base point $x_0\in X$.
Let $d(x,y)$ denote the distance of
$x,y\in X$ induced by $g$.
Let $B(R)$ denote 
$\bigl\{x\in X\,\big|\,d(x,x_0)\leq R\bigr\}$.
As in Page 90--91 of \cite{av},
we can take a sequence of
Lipschitz functions $\chi_n$ ($n=1,2,\ldots$)
such that 
(i) $0\leq \chi_n\leq 1$,
(ii) $\chi_n(x)=1$ on $B(n)$
 and $\chi_n(x)=0$ on $X-B(2n)$,
(iii) $|d\chi_n|\leq C$ for some fixed $C$
 and  for almost every $x\in X$.
Assume $\phi$ satisfies
$\Deltalambda\phi=0$.
We have the following:
\begin{multline*}
0=\int 
 \bigl(\Deltalambda\phi,\,\chi_n\phi\bigr)_{h,g}
 \dvol_g
-\int\bigl(\DDlambda\phi,\,
 \DDlambda(\chi_n\phi)\bigr)_{h,g}
 \dvol_{g} \\
-\int\bigl(\DDlambdaast\phi,\,\,
 \DDlambdaast(\chi_n\phi) \bigr)_{h,g}\dvol_{g}
\end{multline*}
Then, we obtain the vanishing
$\big\|\DDlambda\phi\bigr\|_{h,\omega}^2
+\bigl\|\DDlambdaast\phi\bigr\|_{h,\omega}^2
=0$
due to the theorem of Lesbegue.
\hfill\qed

\subsubsection{A general remark on $L^2$-conditions}

Let $X$ denote a closed disc
$\bigl\{z\in\cnum\,\big|\,|z|\leq 1\bigr\}$,
and we put $D=\{O\}$.
Let $\omega$ be a Poincar\'e like Kahler form 
of  $X-D$.
Let $\dvol_{\omega}$ denote
the volume form associated to $\omega$.
Let $(V,\DDlambda)$ be a $\lambda$-flat bundle
on $X-D$
with a hermitian metric $h$.
We recall a general remark.
\begin{lem}
\label{lem;08.1.19.110}
Let $\phi$ be a $C^{\infty}$-section of 
$V\otimes\Omega^j$ on $X-D$.
Assume that $\phi$
and $\Deltalambda\phi$ are $L^2$
with respect to $h$ and $\omega$.
Then, $\DDlambda\phi$ and 
$\DDlambdaast\phi$ are also $L^2$
with respect to $h$ and $\omega$.
\end{lem}
\pf
For any $C^{\infty}$-section $f$ of 
$V\otimes\Omega^j$ on $X-D$,
let $\|f\|_{h,\omega}$
denote the $L^2$-norm with respect to
$h$ and $\omega$.
If the support of $f$ is compact,
we have the following for some number  $B_0$,
which is independent of $f$:
\begin{equation}
 \label{eq;08.1.20.2}
  \bigl\|\DDlambda f
 \bigr\|_{h,\omega}^2
+\bigl\|\DDlambdaast f\bigr\|_{h,\omega}^2
=
 B_0\int_{X-D}
 \bigl(f,\,
 \Deltalambda(f)\bigr)_{h,\omega}
 \dvol_{\omega}.
\end{equation}

Let $\phi$ be a $C^{\infty}$-section of $V$
on $X-D$,
such that $\phi$
and $\Deltalambda\phi$ are $L^2$.
To show the $L^2$-property of
$\DDlambda\phi$,
we have only to be concerned with
the behaviour around $D$.
Hence, we may and will assume that
the support of $\phi$ is contained in
$\{|z|\leq 1/2\}$,
i.e, there are no contribution
of the boundary $\{|z|=1\}$
to the Stokes formula.
Let $\chi$ be any test function on $X-D$.
Let $\DDlambdastar_h$ be the operator
given as in Section 2.2 of \cite{mochi5}.
(We will use the symbol $\DDlambdastar$
 for simplicity.)
Recall
$\DDlambdaast=-\sqrt{-1}
 \bigl[\Lambda_{\omega},\,
 \DDlambdastar\bigr]$.
We have the following equality:
\begin{multline}
\label{eq;08.1.20.3}
 \Deltalambda(\chi\phi)
=\DDlambdaast\DDlambda(\chi\phi)
=\DDlambdaast\bigl(
 \chi\DDlambda\phi
 \bigr)
+\DDlambdaast\bigl(\DDlambda(\chi)\cdot\phi\bigr) \\
=\chi\Deltalambda(\phi)
-\sqrt{-1}\Lambda_{\omega}
 \bigl(\DDlambdastar(\chi)\cdot\DDlambda\phi\bigr)
+\DDlambdaast\bigl(\DDlambda(\chi)\cdot\phi\bigr)
\end{multline}
We obtain the following:
\begin{multline}
\label{eq;08.1.20.1}
 \int_{X-D}\bigl(
 \Deltalambda(\chi\phi),\,\,\chi\phi
 \bigr)_{h}\dvol_{\omega}
=
 \int_{X-D}\bigl(
 \chi\Deltalambda\phi,\,\chi\phi\bigr)_h
 \dvol_{\omega}
 \\
-\sqrt{-1}
 \int_{X-D}
 \Bigl(
\Lambda_{\omega}
 \bigl(\DDlambdastar(\chi)\cdot\DDlambda\phi\bigr),\,
 \chi\phi
\Bigr)_h\dvol_{\omega}
+\int_{X-D}\bigl(
 \DDlambda(\chi)\cdot\phi,\,\,
 \DDlambda(\chi\phi)
\bigr)_h\dvol_{\omega}
\end{multline}
The second term in the right hand side of
(\ref{eq;08.1.20.1}) can be rewritten as follows,
up to constant multiplication:
\begin{multline}
 \label{eq;08.1.20.4}
 \int_{X-D}
 \Bigl(
 \Lambda_{\omega}
 \bigl(\DDlambdastar(\chi)\cdot
 \chi\DDlambda\phi\bigr),\,
 \phi
 \Bigr)_{h}\dvol_{\omega}=\\
 \int_{X-D}
 \Bigl(
 \Lambda_{\omega}
 \bigl(\DDlambdastar(\chi)\cdot
 \DDlambda(\chi\phi)
-\DDlambdastar(\chi)\cdot\DDlambda(\chi)\cdot\phi
 \bigr),\,
 \phi
 \Bigr)_{h}\dvol_{\omega}
\end{multline}

Let $\rho$ be an $\real_{\geq 0}$-valued
$C^{\infty}$-function on $\real$ such that
$\rho(t)=1$ for $t\leq 1/2$
and $\rho(t)=0$ for $t\geq 1$.
Let $\chi_N:=\rho\bigl(-N^{-1}\log|z|\bigr)$.
Note $\del\chi_N$,
$\delbar\chi_N$ and $\delbar\del\chi_N$
are uniformly bounded with respect to $\omega$.
We obtain the following inequality
for some constants $B>0$,
from (\ref{eq;08.1.20.2}),
(\ref{eq;08.1.20.3}),
(\ref{eq;08.1.20.1})
and (\ref{eq;08.1.20.4}):
\[
 \bigl\|
 \DDlambda(\chi_N\cdot \phi)
 \bigr\|_{h,\omega}^2
\leq
 B\cdot\|\phi\|_{h,\omega}\cdot
  \Bigl(
 \|\phi\|_{h,\omega}
+\|\Deltalambda\phi\|_{h,\omega}
+\|\DDlambda(\chi_N\phi)\|_{h,\omega}
 \Bigr)
\]
We obtain the uniform boundedness of
$\bigl\|\DDlambda(\chi_N\cdot\phi)\bigr\|_{h,\omega}$,
and hence 
$\|\DDlambda(\phi)\|_{h,\omega}<\infty$.

\vspace{.1in}

Let $\phi$ be a $C^{\infty}$-section 
of $V\otimes\Omega^1$ on $X-D$
such that 
$\phi$ and $\Deltalambda\phi$ are $L^2$.
As in the case of $0$-form,
we may assume that the support of
$\phi$ is contained in $\{|z|\leq 1/2\}$.
Let $\chi$ be any test function on $X-D$.
We have the following equalities:
\begin{multline*}
\DDlambdaast\DDlambda(\chi\phi)
=\DDlambdaast\bigl(
 \DDlambda(\chi)\cdot\phi
+\chi\cdot\DDlambda(\phi)
 \bigr)= \\
\DDlambdaast\bigl(
 \DDlambda(\chi)\cdot\phi
 \bigr)
+\chi\cdot\DDlambdaast\DDlambda(\phi)
+\sqrt{-1}\DDlambdastar(\chi)
\cdot\Lambda_{\omega}
 (\DDlambda\phi)
\end{multline*}
Similarly, we have the following:
\[
 \DDlambda\DDlambdaast(\chi\phi)
=\DDlambda\bigl(-\sqrt{-1}\Lambda_{\omega}
 (\DDlambdastar(\chi)\cdot\phi)\bigr)
+\chi\cdot\DDlambda\DDlambdaast\phi
+\DDlambda\chi\cdot\DDlambdaast\phi
\]
Hence, we have the following:
\begin{multline}
 \label{eq;08.1.19.106}
 \int\bigl(\Deltalambda(\chi\phi),\,\,
 \chi\phi
 \bigr)_{h,\omega}\dvol_{\omega}= \\
\int_{X-D}\Bigl(
\DDlambdaast(\DDlambda(\chi)\cdot\phi) 
-\sqrt{-1}\DDlambda\bigl(
 \Lambda_{\omega}(\DDlambdastar\chi\cdot\phi)
 \bigr),\,\,\,
 \chi\phi
 \Bigr)_{h,\omega}\dvol_{\omega}\\
+\int_{X-D}
 \Bigl(\sqrt{-1}\DDlambdastar\chi
 \cdot\Lambda_{\omega}\DDlambda\phi 
+\DDlambda\chi\cdot\DDlambdaast\phi,\,\,\,
 \chi\phi
\Bigr)_{h,\omega}\dvol_{\omega} \\
+\int_{X-D}\bigl(\chi\Deltalambda(\phi),\,\,\,
 \chi\phi \bigr)_{h,\omega}
 \dvol_{\omega}
\end{multline}
The first term in the right hand side can be rewritten
as follows:
\begin{equation}
 \label{eq;08.1.19.105}
 -\int_{X-D}\bigl(
\DDlambda(\chi)\cdot\phi,\,\,
 \DDlambda(\chi\phi)
\bigr)_{h,\omega}\dvol_{\omega}
+\int_{X-D}\bigl(
\sqrt{-1}
 \Lambda_{\omega}(\DDlambdastar\chi\cdot\phi)
 ,\,\,\,
 \DDlambdaast(\chi\phi)
 \bigr)_{h,\omega}\dvol_{\omega}
\end{equation}
The second term in the right hand side of
(\ref{eq;08.1.19.106})
can be rewritten as follows:
\begin{multline}
\label{eq;08.1.19.107}
 \int_{X-D}
 \Bigl(
 \sqrt{-1}\DDlambdastar\chi
 \cdot\Lambda_{\omega}\DDlambda(\chi\phi)
+\DDlambda\chi\cdot\DDlambdaast(\chi \phi),\,\,\,
 \phi
\Bigr)_{h,\omega}\dvol_{\omega} \\
-\int_{X-D}
 \Bigl(
 \sqrt{-1}\DDlambdastar\chi
 \cdot\Lambda_{\omega}\bigl(
 \DDlambda(\chi)\cdot \phi\bigr)
+\sqrt{-1}\DDlambda\chi\cdot
 \Lambda_{\omega}
 \bigl(\DDlambdastar(\chi)\cdot \phi
 \bigr),\,\,\,
 \phi
 \Bigr)_{h,\omega}\dvol_{\omega}
\end{multline}
Let $\chi_N$ be as above.
Due to (\ref{eq;08.1.19.106}),
(\ref{eq;08.1.19.105})
and (\ref{eq;08.1.19.107}),
we obtain the following for some $B>0$:
\begin{multline}
 \bigl\|
 \DDlambda(\chi_N\cdot\phi)
 \bigr\|_{h,\omega}^2
+\bigl\|
 \DDlambdaast(\chi_N\cdot\phi)
 \bigr\|_{h,\omega}^2 
\leq \\
 B
 \|\phi\|_{h,\omega}
\Bigl(
 \|\phi\|_{h,\omega}
+\bigl\|
 \DDlambda(\chi_N\cdot\phi)
 \bigr\|_{h,\omega}
+\bigl\|
 \DDlambdaast(\chi_N\cdot\phi)
 \bigr\|_{h,\omega}
+\bigl\|\Delta^{\lambda}\phi\bigr\|_{h,\omega}
 \Bigr)
\end{multline}
Then, we obtain the uniform boundedness of
$\bigl\|\DDlambda(\chi_N\cdot\phi)
 \bigr\|_{h,\omega}
+\bigl\|\DDlambdaast(\chi_N\cdot\phi)
 \bigr\|_{h,\omega}$,
and hence
$\bigl\|\DDlambda\phi\bigr\|_{h,\omega}<\infty$
and 
$\bigl\|\DDlambdaast\phi\bigr\|_{h,\omega}<\infty$.

The claim for $2$-forms
can be reduced to that for $0$-forms.
Thus, the proof of Lemma 
\ref{lem;08.1.19.110} is finished.
\hfill\qed

\vspace{.1in}

A section $\phi$ of 
$V\otimes\Omega^{j}$
is called a harmonic $i$-form of $V$,
if $\Deltalambda \phi=0$.
We obtain the following corollary.
\begin{cor}
\label{cor;08.1.19.110}
Let $\phi$ be a harmonic $i$-form of $V$
which is $L^2$ with respect to $h$ and $\omega$.
Then, $\DDlambda\phi$
and $\DDlambdaast\phi$ are also $L^2$
with respect to $\omega$ and $h$.
\hfill\qed
\end{cor}

\begin{rem}
The claims of Lemma {\rm\ref{lem;08.1.19.110}}
and Corollary {\rm\ref{cor;08.1.19.110}}
should hold in the case of 
general complete Kahler manifold 
with an appropriate exhaustion function.
We omit the details.
\hfill\qed
\end{rem}

\subsubsection{Refinement of
Lemma \ref{lem;08.1.20.16}}
\label{subsection;08.1.20.11}

Let $C$ be a smooth projective curve,
and let $Z$ be a finite subset of $Z$.
Let $(V,\DDlambda)$ be a $\lambda$-flat bundle
on $C-Z$ with a hermitian metric $h$.
Let $\omega$ be a Kahler form of $C-D$
which is Poincar\'e like around $D$.

\begin{prop}
\label{prop;08.1.20.10}
Let $\phi$ be an $L^2$-section of
$V\otimes\Omega^i$ on $C-D$.
Then, $\Deltalambda\phi=0$
if and only if $\DDlambda\phi=\DDlambdaast\phi=0$.
\end{prop}
\pf
We have only to show the ``only if'' part.
Assume $\Deltalambda\phi=0$.
Because Corollary \ref{cor;08.1.19.110},
we obtain that $\DDlambda\phi$
and $\DDlambdaast\phi$ are $L^2$.
Then, we obtain $\DDlambda\phi=0$
and $\DDlambdaast\phi=0$
due to Lemma \ref{lem;08.1.20.16}.
\hfill\qed

\begin{rem}
Proposition {\rm\ref{prop;08.1.20.10}}
should hold for a $\lambda$-flat bundle
on a complete Kahler manifold
with appropriate exhaustion functions.
We omit the details.
\hfill\qed
\end{rem}

\chapter{Associated Family of 
Meromorphic $\lambda$-Flat Bundles}
\label{section;07.11.24.24}
Let $\harmonicbundle$ be a good wild harmonic bundle
on $X-D$,
where $X$ is a complex manifold
and $D$ is a simple normal crossing hypersurface.
We have the family of $\lambda$-flat bundles
$(\nbige,\DD)$ on $\cnum_{\lambda}\times (X-D)$
associated to $\harmonicbundle$.
We would like to extend it on $\cnum_{\lambda}\times X$
in a meromorphic way.
We have already obtained
a meromorphic prolongment
$\nbigp\nbigelambda$ of
$(\nbigelambda,\DDlambda)$
for each fixed $\lambda$
in Chapter \ref{section;07.6.2.1}.
However,
as was mentioned in {\em Introduction},
the family $\bigcup \nbigp\nbigelambda$
cannot be a nice meromorphic object
unless the harmonic bundle is tame.

In this chapter,
for a given complex number $\lambda_0$,
we study a preliminary prolongment 
$\nbigpzero\nbige$ on
a neighbourhood of $\{\lambda_0\}\times X$
obtained as the sheaf of holomorphic sections
whose norms are of polynomial growth
with respect to a modified metric
$\nbigpzero h$.
It will be deformed to
$\nbigqzero\nbige$ in 
Section \ref{subsection;07.11.5.1},
that is the desired family.

In Section \ref{subsection;08.9.14.10},
we construct a filtered bundle 
$\nbigpzero_{\ast}\nbige$.
In Section \ref{subsection;08.9.14.11},
we show that 
$(\nbigpzero_{\ast}\nbige,\DD)$ is
a good family of filtered $\lambda$-flat bundles.
We also show that the specializations
$\nbigpzero\nbigelambda$
are obtained as the deformation of
$\nbigp\nbigelambda$
caused by a variation of irregular values
(Section \ref{subsection;10.5.17.51}).

In Section \ref{subsection;08.9.14.12},
we give a remark on the growth order
of the norms of partially flat sections
uniformly for $\lambda$.
This is a preparation for the proof of Theorem 
\ref{thm;07.10.11.120}.
(See Section \ref{subsection;07.10.15.20}.)

We study a locally uniform comparison 
of irregular decompositions
of $(\nbigp\nbige^0,\DD^0)$
and $(\nbigpzero\nbige,\DD)$
in Section \ref{subsection;08.9.14.13}.
Because we will not use it 
in the other part of this monograph,
the reader can skip it.

\section{Filtered bundle
$\nbigpzero_{\ast}\nbige$}
\label{subsection;08.9.14.10}
\subsection{Local construction of
 $\nbigpzero_{\veca}\nbigelambda$
 and $\nbigpzero_{\veca}\nbige$
  in the unramified case}
\label{subsection;07.11.23.1}

We use the setting and the notation in Section
\ref{subsection;07.10.7.1}.
We put
$g(\lambda):=
g_{\irr}(\lambda)\cdot
g_{\reg}(\lambda)$,
where $g_{\irr}$ and $g_{\reg}$
are given as follows:
\[
g_{\irr}(\lambda):=
 \exp\left(
 \sum_{\gminia\in\Irr(\theta)}
 \lambda\gminiabar\cdot\pi_{\gminia}
 \right),
\,\,\,
g_{\reg}(\lambda):=
\prod_{j=1}^{\ell}
\exp\left(
 \sum_{\alpha\in\Sp(\theta,j)}
\lambda\cdot\alphabar\cdot\log|z_j|^2
 \cdot \pi_{j,\alpha}
 \right)
\]
Here, $\pi_{\gminia}$
denote the projections
onto $E_{\gminia}$ in the decomposition 
(\ref{eq;07.7.18.10}),
and $\pi_{j,\alpha}$ denote
the projections onto $E_{j,\alpha}$
in the decomposition (\ref{eq;07.7.18.11}).
\index{endomorphism
 $g_{\irr}(\lambda)$}
\index{endomorphism
 $g_{\reg}(\lambda)$}
\index{endomorphism
 $g(\lambda)$}

Let $U(\lambda_0)$ denote a neighbourhood of
$\lambda_0$ in $\cnum_{\lambda}$.
We set $\nbigxzero:=U(\lambda_0)\times X$.
We use the symbols
$\nbigdzero$ and $\nbigdzero_i$
in similar meanings.
We have the following hermitian metrics of
$\nbige_{|\nbigxzero-\nbigdzero}$:
\begin{equation}
 \label{eq;07.7.21.4}
 \nbigp^{(\lambda_0)}h(u,v):=
 h\bigl(
 g(\lambda-\lambda_0)u,
 g(\lambda-\lambda_0)v
 \bigr)
\end{equation}
\[
  \nbigp_{\irr}^{(\lambda_0)}h(u,v):=
 h\bigl(
 g_{\irr}(\lambda-\lambda_0)u,
 g_{\irr}(\lambda-\lambda_0)v
 \bigr)
\]
The naturally induced metric
of $\nbigelambda$ $(\lambda\in U(\lambda_0))$
is also denoted by the same symbols.
\index{metric $\nbigp^{(\lambda_0)}h$}
\index{metric $\nbigp_{\irr}^{(\lambda_0)}h$}

\begin{notation}
Let $\veca\in\real^{\ell}$.
Let $V$ be an open subset of $\nbigxzero$.
We set 
\[
 \nbigpzero_{\veca}\nbige(V):=
 \Bigl\{
 f\in\nbige(V^{\ast})\,\Big|\,
 |f|_{\nbigpzero h}=O\Bigl(
 \prod_{j=1}^{\ell}
 |z_j|^{-a_j-\epsilon}
 \Bigr),\,\,\forall\epsilon>0
 \Bigr\},
\]
where 
$V^{\ast}:=V\setminus \nbigdzero$.
By taking sheafification,
we obtain an $\nbigo_{\nbigxzero}$-module
$\nbigpzero_{\veca}\nbige$.
We put $\nbigpzero\nbige:=
 \bigcup_{\veca\in\real^{\ell}}
 \nbigpzero_{\veca}\nbige$.
The filtered sheaf on $(\nbigxzero,\nbigdzero)$
is denoted by $\nbigpzero_{\ast}\nbige$.

The specialization of 
$\nbigpzero_{\veca}\nbige$
and $\nbigpzero\nbige$
to $\{\lambda\}\times X$ are denoted by
$\nbigpzero_{\veca}\nbigelambda$
and $\nbigpzero\nbigelambda$, respectively.
The specialization of
$\nbigpzero_{\ast}\nbige$
to $\{\lambda\}\times (X,D)$
is denoted by $\nbigpzero_{\ast}\nbigelambda$.
\hfill\qed
\end{notation}
\index{sheaf $\nbigpzero_{\veca}\nbige$}
\index{sheaf $\nbigpzero\nbige$}
\index{filtered bundle $\nbigpzero_{\ast}\nbige$}
\index{sheaf $\nbigpzero_{\veca}\nbigelambda$}
\index{sheaf $\nbigpzero\nbigelambda$}
\index{filtered sheaf $\nbigpzero_{\ast}\nbigelambda$}

\subsection{Global construction of
 $\nbigpzero_{\ast}\nbigelambda$
 and $\nbigpzero_{\ast}\nbige$}

\label{subsection;07.12.2.100}

The construction of the filtered sheaf 
$\nbigpzero_{\ast}\nbige$
in Section \ref{subsection;07.11.23.1}
can be obviously globalized
and extended to the ramified case.
Let $X$ be a general complex manifold,
and let $D$ be a simple normal crossing hypersurface 
with the irreducible decomposition 
$D=\bigcup_{i\in\Lambda}D_i$.
Let $(E,\delbar_E,\theta,h)$ be
a good wild harmonic bundle on $X-D$,
which is not necessarily unramified.
Let $U$ be an open subset of $X$
with a holomorphic coordinate $(z_1,\ldots,z_n)$
such that $U\cap D=\bigcup_{j=1}^{\ell}\{z_j=0\}$.
We take a ramified covering
$\varphi:U'\lrarr U$
given by 
$\varphi(\zeta_1,\ldots,\zeta_n)
=(\zeta_1^{m_1},\ldots,\zeta_{\ell}^{m_{\ell}},
 \zeta_{\ell+1},\ldots,\zeta_n)$
such that $\varphi^{\ast}\harmonicbundle$
is unramified.
Then, we obtain the hermitian metrics 
as in (\ref{eq;07.7.21.4}).
It is equivariant with respect to $\Gal(U'/U)$,
we obtain the hermitian metrics
$\nbigpzero h_{U}$ and $\nbigpzero_{\irr}h_U$
of $\nbige_{|U(\lambda_0)\times (U\setminus D)}$.
By the same procedure,
we obtain the filtered sheaf
$\nbigpzero_{\ast}\nbige_{U}$
on $U(\lambda_0)\times (U,D\cap U)$.
When we are given two such open sets
$U_i$ $(i=1,2)$ of $X$,
the hermitian metrics 
$\nbigpzero h_{U_i}$ $(i=1,2)$
are mutually bounded.
Hence, the restrictions of the filtered sheaves 
$\nbigpzero_{\ast}\nbige_{U_i}$
to $U(\lambda_0)\times 
 \bigl(U_1\cap U_2,D\cap(U_1\cap U_2)\bigr)$
are the same.
By varying $U$ and gluing
$\nbigpzero_{\ast}\nbige_{U}$,
we obtain the filtered sheaf
$\nbigpzero_{\ast}\nbige$ on 
$(\nbigx^{(\lambda_0)},\nbigd^{(\lambda_0)})$,
where $\nbigx^{(\lambda_0)}$ denotes
a neighbourhood of $\{\lambda_0\}\times X$
in $\cnum_{\lambda}\times X$
and $\nbigd^{(\lambda_0)}:=\nbigx^{(\lambda_0)}\cap
 (\cnum_{\lambda}\times D)$.

\vspace{.1in}
We will show the following theorem
in Section \ref{subsection;07.11.20.50}.

\begin{thm}
\label{thm;07.11.20.10}
$\nbigpzero_{\ast}\nbige$ is a filtered vector bundle
on $(\nbigx^{(\lambda_0)},\nbigd^{(\lambda_0)})$,
if $\nbigx^{(\lambda_0)}$ is sufficiently small.
\end{thm}

We have the induced filtration
$\lefttop{i}\Fzero$ of
$\nbigpzero_{\veca}\nbige_{|\nbigd^{(\lambda_0)}_i}$.
The tuple $\bigl(\lefttop{i}\Fzero\,\big|\,i\in \Lambda\bigr)$
is denoted by $\vecF^{(\lambda_0)}$.
We have the weak norm estimate up to small polynomial
order in the following sense.
For simplicity, we consider the case
$X=\Delta^n$ and $D=\bigcup_{j=1}^{\ell}\{z_j=0\}$.
Assume that $U(\lambda_0)$ is sufficiently small.
Let $\vecv$ be a frame of
$\nbigpzero_{\veca}\nbige$
compatible with $\vecF^{(\lambda_0)}$.
Let $a_j(v_i):=\lefttop{j}\deg^{\Fzero}(v_i)$.
We put 
$v_i':=v_i\,
 \prod_{j=1}^{\ell}|z_j|^{a_j(v_i)}$.
Let $H(\nbigpzero h,\vecv')$ 
denote the Hermitian matrix-valued
function whose $(i,j)$-entries are given 
by $\nbigpzero h(v_i',v_j')$.
We have the weak norm estimate
up to small polynomial order,
which will also be proved in Section 
 \ref{subsection;07.11.20.50}.
\begin{prop}
\label{prop;07.11.21.10}
For any $\epsilon>0$,
there exist a positive constant $C_{\epsilon}$
such that 
\[
 C_{\epsilon}^{-1}\,
 \prod_{j=1}^{\ell}|z_j|^{\epsilon}\leq
 H(\nbigpzero h,\vecv')\leq
 C_{\epsilon}\,
 \prod_{j=1}^{\ell}|z_j|^{-\epsilon}.
\]
\end{prop}

\subsection{Prolongment
 $\nbigtzero_{\veca}\nbige$}
\label{subsection;07.11.20.50}

Let us return to the setting in Section 
\ref{subsection;07.11.23.1}.
We put
$\nbigt^{(\lambda_0)}d_{\lambda}'':=
 g(\lambda-\lambda_0)\circ\bigl(
 \delbar_E+\lambda\theta^{\dagger}
 \bigr)\circ g(\lambda-\lambda_0)^{-1}$.
Let $\nbigt^{(\lambda_0)}\nbige$ denote 
the following holomorphic bundle
on $\nbigxzero-\nbigdzero$:
\[
 \bigl(p_{\lambda}^{-1}E,
 \nbigt^{(\lambda_0)}d''_{\lambda}+\delbar_{\lambda}
 \bigr)
\]
Let $\veca\in\real^{\ell}$.
For any open subset $V\subset \nbigxzero$,
we define
\[
 \nbigtzero_{\veca}\nbige(V):=
 \Bigl\{
 f\in\nbigtzero\nbige(V^{\ast})\,\Big|\,
 |f|_h=O\Bigl(
 \prod_{i=1}^{\ell}|z_i|^{-a_i-\epsilon}
 \Bigr),\,\,
 \forall\epsilon>0
 \Bigr\},
\]
where $V^{\ast}:=V\setminus \nbigdzero$.
By taking sheafification,
we obtain a filtered sheaf 
$\nbigtzero_{\ast}\nbige$
on $(\nbigxzero,\nbigdzero)$.
The following lemma is clear from the construction.
\begin{lem}
The multiplication of 
$g(\lambda-\lamda_0)$
induces the holomorphic isomorphisms
$\nbige\simeq \nbigtzero\nbige$
and $\nbigpzero_{\veca}\nbige\simeq
 \nbigtzero_{\veca}\nbige$.
\hfill\qed
\end{lem}

Therefore, 
we obtain Theorem \ref{thm;07.11.20.10}
and Proposition \ref{prop;07.11.21.10}
from the following proposition.
\begin{prop}
\label{prop;07.11.21.2}
If $\nbigxzero$ is sufficiently small,
$\nbigtzero_{\ast}\nbige$ is a filtered bundle
on $(\nbigxzero,\nbigdzero)$.
The weak norm estimate 
up to small polynomial order
holds for $(\nbigtzero_{\ast}\nbige,h)$.
(See Proposition {\rm\ref{prop;07.11.21.10}}
 for weak norm estimate up to small polynomial order.)
\end{prop}
\pf
We put
$\Lambda(\lambda-\lambda_0)
:=\nbigt^{(\lambda_0)}d''_{\lambda}
-d''_{\lambda_0}$.
According to Theorem \ref{thm;08.10.18.11},
Proposition \ref{prop;07.11.21.2}
follows from the following lemma.

\begin{lem}
\label{lem;07.6.2.5}
If $U(\lambda_0)$ is sufficiently small,
we have
$\bigl|\Lambda(\lambda-\lambda_0)\bigr|_{h,g_{\poin}}
\leq C|\lambda-\lambda_0|$
on $\nbigxzero-\nbigdzero$,
where $g_{\poin}$ denotes the Poincar\'e metric
for $X-D$.
\end{lem}
\pf
Let $\Phibar$ be given by (\ref{eq;07.11.22.35}).
We have the following equality:
\begin{multline}
 \label{eq;07.11.22.36}
 \Lambda(\lambda-\lambda_0)=
\lambda_0\cdot\Bigl(
 g(\lambda-\lambda_0)\circ
\bigl(\theta^{\dagger}-\Phibar\bigr) 
 \circ g(\lambda-\lambda_0)^{-1}
-\bigl(\theta^{\dagger}-\Phibar\bigr)
 \Bigr) 
 \\
+(\lambda-\lambda_0)\cdot
 g(\lambda-\lambda_0)
\circ\Bigl(
 \theta^{\dagger}
-\Phibar
 \Bigr)
\circ
 g(\lambda-\lambda_0)^{-1}
\end{multline}
Let us look at the first term of (\ref{eq;07.11.22.36}).
We use the decomposition as in 
(\ref{eq;08.12.8.1}).
If $(\gminia,\vecalpha)\neq(\gminia',\vecalpha')$,
we obtain the following
from (\ref{eq;10.5.24.10}):
\begin{multline}
 \Bigl(
 g(\lambda-\lambda_0)\circ
\bigl(\theta^{\dagger}-\Phibar\bigr) 
 \circ g(\lambda-\lambda_0)^{-1}
-\bigl(\theta^{\dagger}-\Phibar\bigr)
 \Bigr)_{(\gminia,\vecalpha),(\gminia',\vecalpha')} \\
=\Bigl(
\exp\bigl((\lambda-\lambda_0) (\gminiabar'-\gminiabar)
 \bigr)
\cdot\prod_{j=1}^{\ell}
 |z_j|^{2(\lambda-\lambda_0)(\alphabar_j'-\alphabar_j)}
-1\Bigr)
\cdot
\bigl(\theta^{\dagger}-\Phibar\bigr)_{
 (\gminia,\vecalpha),(\gminia',\vecalpha')}
 \\
=|\lambda-\lambda_0|\cdot
O\Bigl(
 \exp\bigl((\lambda-\lambda_0) (\gminiabar'-\gminiabar)
 \bigr)
\cdot\prod_{j=1}^{\ell}
 |z_j|^{2(\lambda-\lambda_0)(\alphabar_j'-\alphabar_j)}
\Bigr)
 \times \\
O\Bigl(
 \exp\bigl(-\epsilon|\vecz^{\ord(\gminia-\gminia')}|\bigr)
 \cdot\nbigq_{\epsilon}(\vecalpha,\vecalpha')
 \Bigr)
\end{multline}
If $|\lambda-\lamda_0|$ is sufficiently small,
we have
\[
 \Bigl|
(\lambda-\lambda_0)(\gminiabar'-\gminiabar)
\Bigr|
-\epsilon|\vecz^{\ord(\gminia-\gminia')}|
\leq
 -\epsilon|\vecz^{\ord(\gminia-\gminia')}|\big/2
\]
\[
\prod_{j=1}^{\ell}
 |z_j|^{2(\lambda-\lambda_0)
 (\alphabar_j'-\alphabar_j)}
\times\nbigq_{\epsilon}(\vecalpha,\vecalpha')
=O\Bigl(
\nbigq_{\epsilon/2}(\vecalpha,\vecalpha')
 \Bigr)
\]
We also have
$\Bigl(
 g(\lambda-\lambda_0)\circ
\bigl(\theta^{\dagger}-\Phibar\bigr) 
 \circ g(\lambda-\lambda_0)^{-1}
-\bigl(\theta^{\dagger}-\Phibar\bigr)
 \Bigr)_{(\gminia,\vecalpha),(\gminia,\vecalpha)}=0$.
Hence, we obtain the desired estimate
for the first term.

For the second term,
we have the following in the case
$(\gminia,\vecalpha)\neq (\gminia',\vecalpha')$:
\begin{multline}
 \Bigl(
 g(\lambda-\lambda_0)\circ
\bigl(\theta^{\dagger}-\Phibar\bigr) 
 \circ g(\lambda-\lambda_0)^{-1}
 \Bigr)_{(\gminia,\vecalpha),(\gminia',\vecalpha')} \\
=
\exp\bigl((\lambda-\lambda_0) (\gminiabar'-\gminiabar)
 \bigr)
\cdot\prod_{j=1}^{\ell}
 |z_j|^{2(\lambda-\lambda_0)(\alphabar_j'-\alphabar_j)}
\cdot
\bigl(\theta^{\dagger}-\Phibar\bigr)_{
 (\gminia,\vecalpha),(\gminia',\vecalpha')}
 \\
= O\Bigl(
 \exp\bigl((\lambda-\lambda_0) (\gminiabar'-\gminiabar)
 \bigr)
 \prod_{j=1}^{\ell}
 |z_j|^{2(\lambda-\lambda_0)(\alphabar_j'-\alphabar_j)}
\Bigr)
\cdot
O\Bigl(
 \exp\bigl(-\epsilon|\vecz^{\ord(\gminia-\gminia')}|\bigr)
 \nbigq_{\epsilon}(\vecalpha,\vecalpha')
 \Bigr)
\end{multline}
We also have  the following:
\[
 \Bigl(
 g(\lambda-\lambda_0)\circ
\bigl(\theta^{\dagger}-\Phibar\bigr) 
 \circ g(\lambda-\lambda_0)^{-1}
 \Bigr)_{(\gminia,\vecalpha),(\gminia,\vecalpha)}
=
\Bigl( \sum_{j=1}^{\ell}f_j^{\nil}\cdot dz_j
+\sum_{j=\ell+1}^{n}f_j^{\reg}\cdot dz_j
\Bigr)
\]
Therefore, we obtain the desired estimate
in Lemma \ref{lem;07.6.2.5}.
The proof of the propositions is also finished.
\hfill\qed

\section{Family of meromorphic 
flat $\lambda$-connections on $\nbigpzero_{\ast}\nbige$}
\label{subsection;08.9.14.11}
\subsection{Comparison of
 $\nbigpzero\nbigelambda$
 and $\nbigp\nbigelambda$}
\label{subsection;08.12.8.2}

We continue to use the setting in Section
\ref{subsection;07.12.2.100}.
We would like to give another interpretation
of the specialization $\nbigpzero\nbigelambda$,
as the deformation of meromorphic
$\lambda$-flat bundles 
caused by variation of the irregular values,
explained in Section \ref{subsection;10.5.17.51}.
For any complex number $\lambda$,
we set 
\begin{equation}
 \label{eq;07.10.14.20}
 T(\lambda):=
 \frac{1+\lambda\lambdabar_0}{1+|\lambda|^2}.
\end{equation}
We take $\veca=(a_i)\in\real^{\Lambda}$ such that
$a_i\not\in\Par(\nbigp\nbigelambdazero,i)$ for each $i$.
\begin{prop}
\label{prop;07.10.21.100}
Let $P$ be any point of $X$.
There exist a neighbourhood 
$X_P$ of $P$ in $X$
and a neighbourhood $U_P(\lambda_0)$ of $\lambda_0$
in $\cnum_{\lambda}$,
such that the following holds
for any $\lambda\in U_P(\lambda_0)$:
\begin{itemize}
\item
We have the natural isomorphism
$\nbigpzero_{\veca}\nbigelambda_{|X_P}
\simeq
 (\nbigp_{\veca}\nbigelambda_{|X_P})^{(T(\lambda))}$
of $\nbigo_{X_P}$-modules,
which is the extension of
the identity on $X_P\setminus D$.
(We put $(\nbigp_{\veca}\nbige^0)^{(1)}
 =\nbigp_{\veca}\nbige^0$
formally, 
in the case $\lambda_0=0$.)
\item
In particular,
we have the isomorphism
$\nbigpzero\nbigelambda_{|X_P}
\simeq
 (\nbigp\nbigelambda_{|X_P})^{(T(\lambda))}$
of $\nbigo_{X_P}(\ast D)$-modules,
which is the extension of
the identity on $X_P\setminus D$.
\end{itemize}
\end{prop}
\pf
The claims are local property.
Therefore, we may and will use the setting in
Section \ref{subsection;07.11.23.1}.
We have only to show the first claim.
Let $U(\lambda_0)$ denote a small neighbourhood 
of $\lambda_0$ in $\cnum_{\lambda}$.
Let $\veca\in\real^{\ell}$ be as above.
We consider the prolongment for the metric
$\nbigpzero_{\irr}h$,
i.e.,
for any open subset $V\subset\nbigxzero$,
we set 
\[
 \nbigpzero_{\irr\,\veca}\nbige(V):=
 \Bigl\{
 f\in\nbige(V^{\ast})\,\Big|\,
 |f|_{\nbigpzero_{\irr} h}=
 O\Bigl(
 \prod_{j=1}^{\ell}
 |z_j|^{-a_j-\epsilon}
 \Bigr),\,\,\forall\epsilon>0
 \Bigr\},
\]
where 
$V^{\ast}:=V\setminus \nbigdzero$.
Thus, we obtain an $\nbigo_{\nbigxzero}$-module
$\nbigpzero_{\irr\,\veca}\nbige$.

\begin{lem}
\label{lem;07.11.23.2}
We have 
$\nbigpzero_{\irr\veca}\nbige
=\nbigpzero_{\veca}\nbige$
for $\veca$ as in this subsection,
if $U(\lambda_0)$ is sufficiently small.
\end{lem}
\pf
Note that $\nbigpzero_{\irr}h$
and $\nbigpzero h$ are mutually bounded
up to $|z|^{-\eta|\lambda-\lambda_0|}$-order
for some $\eta>0$.
Then, the claim of Lemma \ref{lem;07.11.23.2}
follows from Proposition \ref{prop;07.11.21.10}.
\hfill\qed

\vspace{.1in}
We remark 
$\nbigpzero_{\irr}h=
 g_{\irr}(\lambda-\lambda_0)^{\ast}h$.
(See Section \ref{subsection;07.12.21.50}
 for $g_{\irr}(w)^{\ast}h$.)

\begin{lem}
\label{lem;07.12.21.51}
If the divisor $D$ is smooth,
the claim of Proposition {\rm\ref{prop;07.10.21.100}}
holds.
\end{lem}
\pf
We would like to apply Proposition
\ref{prop;07.7.20.230}.
We use the notation in Section
\ref{subsection;07.12.21.50}.
Note $T(\lambda)=T_{1}(\lambda-\lambda_0)$.

Let us consider the case $\lambda_0=0$.
Let $\nbigx^{(0)}:=U(0)\times X$
and $W:=(U(0)\times D)\cup (\{0\}\times X)$.
We take a finite covering
$\nbigx^{(0)}-W=\bigcup S_i$
satisfying the following:
\[
 S_i\in\bigcap_j
 \bigcap_{\gminia\in\Irr(\theta,j)}
 \Multisector\bigl(
 \nbigx^{(0)}-W,\nbigi^{(j)}_{\gminia}
 \bigr)
\]
For each $\lambda\in U(0)-\{0\}$,
let $S^{\lambda}_i:=S_i\cap(\{\lambda\}\times X)$.
Then, we have
\[
 S_i^{\lambda}\in\bigcap_j
 \bigcap_{\gminia\in\Irr(\theta,j)}
 \Multisector\bigl(
 \{\lambda\}\times (X-D),
 \nbigi^{(j)}_{\gminia}
 \bigr)
\]
If $U(0)$ is sufficiently small,
we may assume $\epsilon_2=\epsilon_1$
in the condition (A2).
We have 
$T(\lambda)=(1+|\lambda|^2)^{-1}>0$
for any $\lambda$,
and hence the second assumption in 
Proposition \ref{prop;07.7.20.230}
is trivial.
Therefore, the claim of Lemma \ref{lem;07.12.21.51}
immediately follows from Proposition \ref{prop;07.7.20.230}
in the case $\lambda_0=0$.

Let us consider the case $\lambda_0\neq 0$.
We can take a finite covering
$X-D=\bigcup_{i=1}^N S_i$,
where $S_i$ are multi-sectors satisfying the following:
\[
 \{\lambda_0\}\times S_i
 \in\bigcap_{j}\bigcap_{\gminia\in\Irr(\theta,j)}
 \Multisector\bigl(\{\lambda_0\}\times(X-D),\,
 \nbigi_{\gminia}^{(j)}\bigr)
\]
If $U(\lambda_0)$ is sufficiently small,
we have the following for any 
$\lambda\in U(\lambda_0)$:
\[
  \{\lambda\}\times S_i
 \in\bigcap_{j}\bigcap_{\gminia\in\Irr(\theta,j)}
 \Multisector\bigl(\{\lambda\}\times(X-D),\,
 \nbigi_{\gminia}^{(j)}\bigr)
\]
We may assume 
that the conditions (A1) and (A2) are also satisfied 
for each $\lambda\in U(\lambda_0)$.
We put 
$T_t(\lambda):=t+(1-t)T(\lambda)$.
Let $(\nbigi_{\gminia}^{(j)})^{(T_t(\lambda))}:=
 \bigl\{
 T_t(\lambda)\cdot \gminic\,\big|\,
 \gminic\in \nbigi_{\gminia}^{(j)}
 \bigr\}$.
If $U(\lambda_0)$ is sufficiently small,
we may also have the following:
\[
  \{\lambda\}\times S_i
 \in\bigcap_{j}\bigcap_{\gminia\in\Irr(\theta,j)}
 \Multisector\bigl(\{\lambda\}\times(X-D),\,
 (\nbigi_{\gminia}^{(j)})^{(T_t(\lambda))} \bigr)
\]
Then, the second assumption
of Proposition \ref{prop;07.7.20.230}
is satisfied,
and the claim of Lemma \ref{lem;07.12.21.51}
follows.
\hfill\qed

\vspace{.1in}

Due to Theorem \ref{thm;07.11.20.10},
$\nbigpzero_{\veca}\nbigelambda$ is locally free.
We also know that
$\bigl(\nbigp_{\veca}\nbigelambda
 \bigr)^{(T(\lambda))}$
is also locally free.
Let $D_j(<\epsilon):=
 \bigl\{\vecz\in X\,\,|\,
 z_j=0,\,\,\|\vecz\|<\epsilon
 \bigr\}$.
By using Lemma \ref{lem;07.12.21.51}
and Lemma \ref{lem;10.5.25.31},
we obtain a natural isomorphism
of $\nbigpzero_{\veca}\nbigelambda$
and 
$\bigl(\nbigp_{\veca}\nbigelambda
 \bigr)^{(T(\lambda))}$ on 
$X-\bigcup_{j=1}^{\ell}D_j(<\epsilon)$.
By Hartogs theorem,
it is extended to the isomorphism on $X$.
Thus, the proof of Proposition \ref{prop;07.10.21.100}
is finished.
\hfill\qed

\vspace{.1in}
We have immediate consequences of
Proposition \ref{prop;07.10.21.100}.
\begin{cor}
\label{cor;07.10.14.10}
Let $P$ be any point of $X$.
Let $X_P$ and $U_P(\lambda_0)$
be as in Proposition {\rm\ref{prop;07.10.21.100}}.
\begin{itemize}
\item
For any $\lambda\in U_P(\lambda_0)$,
the flat $\lambda$-connection $\DDlambda$ of
$\nbigp^{(\lambda_0)}\nbigelambda_{|X_P}$ 
is meromorphic.
\item
$\bigl(
\nbigp^{(\lambda_0)}_{\ast}
 \nbigelambda,\DDlambda\bigr)_{|X_P}$
is good.
\item
The set of the irregular values of 
$\bigl(\nbigpzero\nbigelambda,\DDlambda\bigr)$
is given by
\[
 \Irr(\nbigpzero\nbigelambda,\DDlambda):=
 \bigl\{
 (1+\lambda\lambdabar_0)\cdot\gminia\,\big|\,
 \gminia\in\Irr(\theta)
 \bigr\} 
\]
under the setting in Section 
{\rm \ref{subsection;07.10.7.1}}.
\hfill\qed
\end{itemize}
\end{cor}

\subsection{The family of 
 flat $\lambda$-connections of $\nbigpzero\nbige$
 and the KMS-structure at $\lambda_0$}

We continue to use the setting in 
Section \ref{subsection;07.12.2.100}.
\begin{cor}
\label{cor;08.9.28.117}
Let $P$ be any point of $X$.
Let $X_P$ and $U_P(\lambda_0)$
be as in Proposition {\rm\ref{prop;07.10.21.100}}.
\begin{itemize}
\item
The family of flat $\lambda$-connections $\DD$ of
$\nbigp^{(\lambda_0)}\nbige_{|U_P(\lambda_0)\times X_P}$
is meromorphic.
\item
The family of filtered $\lambda$-flat bundles
$(\nbigpzero_{\ast}\nbige,\DD)$ is good
on $U(\lambda_0)\times \bigl(X_P,X_P\cap D\bigr)$.
\item
The irregular values of 
$\bigl(\nbigpzero\nbige,\DD\bigr)$
is given by 
\[
  \Irr(\nbigpzero\nbige,\DD):=
 \bigl\{
 (1+\lambda\lambdabar_0)\cdot\gminia\,\big|\,
 \gminia\in\Irr(\theta)
 \bigr\}
\]
under the setting in Section 
{\rm \ref{subsection;07.10.7.1}}.
\end{itemize}
\end{cor}
\pf
It follows from Corollary \ref{cor;07.10.14.10}
and Proposition \ref{prop;10.5.5.32}.
\hfill\qed

\begin{prop}
\label{prop;07.10.14.21}
$(\nbigpzero_{\ast}\nbige,\DD)$
has the KMS-structure at $\lambda_0$
in the sense of Definition {\rm\ref{df;07.11.23.5}}.
\end{prop}
\pf
We may assume that $D$ is smooth.
We consider 
$\nbigpzero_{a}\nbige/\nbigpzero_{<a}\nbige$
on $D$.
Take $\epsilon>0$.
If $U(\lambda_0)$ is sufficiently small,
we have
$\bigl(
 \nbigpzero_a\nbige\big/\nbigpzero_{<a}\nbige
 \bigr)_{|D\times\{\lambda\}}
\simeq
 \bigl(
 \nbigp_{a+\epsilon}\nbigelambda
\bigr)^{T(\lambda)}
 \big/
\bigl(
 \nbigp_{a-\epsilon}\nbigelambda
\bigr)^{T(\lambda)}$ for any $\lambda\in U(\lambda_0)$.
Hence,
we can conclude that the set of
the eigenvalues of $\Res(\DDlambda)$
on $\bigl(
 \nbigpzero_a\nbige\big/\nbigpzero_{<a}\nbige
 \bigr)_{|D_i\times\{\lambda\}}$
is given by the following,
according to Proposition \ref{prop;07.7.19.31}:
\[
 \bigl\{
 \eigenmap(\lambda,u)\,\big|\,
 u\in \KMS(\nbige^0),\,\,
 \paramap(\lambda_0,u)=a
 \bigr\}
\]
Thus, we are done.
\hfill\qed

\vspace{.1in}
From the proof of Proposition \ref{prop;07.10.14.21},
we also obtain the following corollary.

\begin{cor}
\label{cor;07.11.23.20}
Let $P$ be any point of $X$,
and let $X_P$ and $U_P(\lambda_0)$
be as in Proposition {\rm\ref{prop;07.10.21.100}}.
For each $\lambda\in U(\lambda_0)$,
we have the induced filtered bundle
$(\nbigpzero\nbige_{|X_P})_{\ast}^{\lambda}$
as given in Section {\rm\ref{subsection;07.10.14.15}}.
If $U_P(\lambda_0)$ is sufficiently small,
we have the natural isomorphism
$(\nbigpzero\nbige_{|X_P})_{\ast}^{\lambda}
\simeq
 (\nbigp_{\ast}\nbigelambda_{|X_P})^{(T(\lambda))}$
of the family of $\lambda$-flat bundles
for each $\lambda\in U_P(\lambda_0)$,
where $T(\lambda)$ is given as in {\rm(\ref{eq;07.10.14.20})}.
In particular,
we have the isomorphism of
meromorphic $\lambda$-flat bundles
$(\nbigpzero\nbige_{|X_P})^{\lambda}
\simeq
 (\nbigp\nbigelambda_{|X_P})^{(T(\lambda))}$
for each $\lambda\in U_P(\lambda_0)$.
\hfill\qed
\end{cor}

\subsection{Comparison of 
 $\nbigp^{(\lambda_0)}\nbige$
 and $\nbigp^{(\lambda_1)}\nbige$}
\label{subsection;07.12.22.110}

We continue to use the setting in 
Section \ref{subsection;07.12.2.100}.
Let $\veca\in\real^{\Lambda}$ be
as in Subsection \ref{subsection;08.12.8.2}.
Let $P$ be any point of $X$,
and let $X_P$ and $U_P(\lambda_0)$
be as in Proposition \ref{prop;07.10.21.100}.
We take $\lambda_1\in U_P(\lambda_0)$.
By shrinking $X_P$,
we may also assume to have
$(\nbigp^{(\lambda_1)}\nbige,\DD)$
on $U_P(\lambda_1)\times (X_P,D\cap X_P)$
for some 
neighbourhood
$U_P(\lambda_1)\subset U_P(\lambda_0)$
of $\lambda_1$.
We assume $0\not\in U_P(\lambda_1)$,
although $\lambda_0$ may be $0$.

\begin{prop}
\label{prop;07.10.8.100}
\mbox{{}}
\begin{itemize}
\item
The family of meromorphic $\lambda$-flat bundles
$\bigl(\nbigp^{(\lambda_0)}
 \nbige, \DD\bigr)_{|U_P(\lambda_1)\times X_P}$
is naturally isomorphic to the deformation
$\bigl(\nbigp^{(\lambda_1)}\nbige,\DD\bigr)^{(T')}
 _{|U_P(\lambda_1)\times X_P}$
with 
\[
 T'=(1+\lambda\lambdabar_0)
 \cdot(1+\lambda\lambdabar_1)^{-1}.
\]
\item
The family of meromorphic $\lambda$-flat bundles
$\bigl(\nbigp_{\veca}^{(\lambda_0)}
 \nbige,\DD\bigr)_{|U_P(\lambda_1)\times X_P}$
is naturally isomorphic to the deformation
$\bigl(\nbigp_{\veca}^{(\lambda_1)}\nbige,
 \DD\bigr)^{(T')}_{|U_P(\lambda_1)\times X_P}$.
\end{itemize}
\end{prop}
\pf
We have only to show the second claim.
According to Lemma \ref{lem;07.12.21.60}
and Proposition \ref{prop;07.10.21.100},
both the restrictions
$\nbigp_{\veca}^{(\lambda_0)}
 \nbige_{|\{\lambda\}\times X_P}$
and
$\bigl(\nbigp_{\veca}^{(\lambda_1)}\nbige\bigr)^{(T')}
 _{|\{\lambda\}\times X_P}$
are naturally isomorphic to 
$(\nbigp_{\veca}\nbigelambda)^{(T_1)}_{|X_P}$,
where $T_1=(1+\lambda\lambdabar_0)
 (1+|\lambda|^2)^{-1}$.
Then, it is easy to see that
the natural isomorphism of
$\nbigp_{\veca}^{(\lambda_0)}
 \nbige_{|U(\lambda_0)\times (X_P\setminus D)}$
and
$\bigl(\nbigp_{\veca}^{(\lambda_1)}\nbige\bigr)^{(T')}
 _{|U(\lambda_0)\times (X_P\setminus D)}$
are extended to the isomorphism
on $X_P$.
\hfill\qed

\vspace{.1in}
From the family of $\lambda$-flat bundles
$(\nbigpzero_{\ast}\nbige,\DD)$
on $U_P(\lambda_0)\times (X_P,X_P\cap D)$ 
with the KMS-structure at $\lambda_0$,
we obtain the family of $\lambda$-flat bundles
$\bigl((\nbigp^{(\lambda_0)}\nbige)^{(\lambda_1)}_{\ast},\DD\bigr)$
on $U_P(\lambda_1)\times (X_P,X_P\cap D)$ 
as in Section \ref{subsection;07.10.14.15}.
\begin{prop}
We have the natural isomorphism
\[
 \bigl(
 (\nbigp^{(\lambda_0)}\nbige)^{(\lambda_1)}_{\ast},
 \DD\bigr)
 \simeq
 \bigl(\nbigp^{(\lambda_1)}_{\ast}\nbige,
 \DD\bigr)^{(T')}
\]
where $T'$ is given in Proposition {\rm\ref{prop;07.10.8.100}}.
\end{prop}
\pf
We have the isomorphism
$\bigl(
 (\nbigp^{(\lambda_0)}\nbige)^{(\lambda_1)},
 \DD\bigr)
 \simeq
 \bigl(\nbigp^{(\lambda_1)}\nbige,\DD\bigr)^{(T')}$
due to Proposition \ref{prop;07.10.8.100}.
Then, the claim follows from Lemma \ref{lem;07.11.23.5}.
\hfill\qed

\section{Estimate of the norms of 
partially flat sections}
\label{subsection;08.9.14.12}
\label{subsection;07.12.22.1}

We use the setting 
in Section \ref{subsection;07.11.17.1}.
For simplicity, we assume that the coordinate is 
admissible for the good set $\Irr(\theta)$.
Let $k$ be determined by the condition
$\vecm(0)\in\seisuu_{<0}^k\times\veczero_{\ell-k}$.
Let $\lambda_0\in\cnum_{\lambda}$.
Let $\nbigx^{(\lambda_0)}$ denote
a neighbourhood of $\{\lambda_0\}\times X$
in $\cnum_{\lambda}\times X$.
We use the symbols like
$\nbigd^{(\lambda_0)}(\leq k)$
in similar meanings.
In the case $\lambda_0=0$,
we put 
$W(\leq k):=
 \nbigd^{(\lambda_0)}(\leq k)
\cup\bigl(\{0\}\times X\bigr)$.
Otherwise,
we put $W(\leq k):=\nbigd^{(\lambda_0)}(\leq k)$.
Let $\pi: 
\nbigxtilde^{(\lambda_0)}(W(\leq k))
 \lrarr\nbigx^{(\lambda_0)}$
denote the real blow up of
$\nbigx^{(\lambda_0)}$ along $W(\leq k)$.
Let $\nbigi^{(\lambda_0)}$
denote the image of 
$\Irr\bigl(\nbigpzero\nbige,\DD\bigr)$
by $\etabar_{\vecm(0)}$.
Let $\DD_{\leq k}$ denote the restriction of
$\DD$ to derivations along
the $(z_1,\ldots,z_k)$-direction.

Let $S$ be a multi-sector in
$\nbigx^{(\lambda_0)}-W(\leq k)$,
and let $\Sbar$ denote the closure of $S$
in $\nbigxtilde^{(\lambda_0)}(W(\leq k))$.
If $S$ is sufficiently small,
we have the partial Stokes filtration $\nbigf^S$ 
of $\nbigp^{(\lambda_0)}\nbige$ on $\Sbar$
in the level $\vecm(0)$
indexed by the ordered set 
$(\nbigi^{(\lambda_0)},\leq_S)$
due to Proposition \ref{prop;10.5.11.20}.

\begin{lem}
\label{lem;07.6.2.60}
Let $f$ be a $\DD_{\leq k}$-flat section of
$\End(\nbigpzero\nbige)_{|S}$
such that
\begin{itemize}
\item
$\bigl[\Res_i(\DD),
 f_{|\nbigd^{(\lambda_0)}_i\cap S}\bigr]=0$
for $i=k+1,\ldots,\ell$.
\item
$f_{|\nbigd^{(\lambda_0)}_i\cap S}$
preserves the filtrations $\lefttop{i}\Fzero$
for $i=k+1,\ldots,\ell$.
\item
 $f\bigl(\nbigf^S_{\gminia}\nbigpzero_0\nbige\bigr)
 \subset \nbigf^S_{<\gminia}\nbigpzero_0\nbige$
for any $\gminia\in\Irr(\nbigpzero\nbige,\DD)$.
\end{itemize}
If we shrink $S$ in the radius direction,
we have the estimate
\[
 |f|_h\leq C\cdot
 \exp\bigl(-\epsilon|\lambda^{-1}\vecz^{\vecm(0)}|\bigr)
\]
for some $C>0$ and $\epsilon>0$.
\end{lem}
\pf
We use an induction on $\ell$.
We assume that the claim of the lemma holds
for any unramifiedly good wild harmonic bundles
on $X'-D'$,
where $X'=\Delta^{n}$,
$D'_i=\{z_i=0\}$
and $D'=\bigcup_{i=1}^{\ell'}D_i'$
for $\ell'<\ell$.

By shrinking $S$ in the radius direction,
we can take a $\DD_{\leq k}$-flat splitting of $\nbigf^S$:
\begin{equation}
 \label{eq;07.6.2.51}
 \nbigp_0^{(\lambda_0)}
 \nbige_{|\Sbar}
=\bigoplus_{\gminia\in
 \nbigi^{(\lambda_0)}}
 \nbigp_0^{(\lambda_0)}\nbige_{\gminia,S}
\end{equation}
Let $\vecu_S=(\vecu_{\gminia,S})$ be a frame of
$\nbigp^{(\lambda_0)}_0
 \nbige_{|\Sbar}$
compatible with the decomposition (\ref{eq;07.6.2.51}).
Let $B_{\gminia,S}$ be the matrix valued function
determined by
$\DD_{\leq k}\vecu_{\gminia,S}=
 \vecu_{\gminia,S}\cdot\bigl(
 d_{\leq k}\gminia+B_{\gminia,S}
 \bigr)$,
where $d_{\leq k}$ denote the restriction
of the exterior differential
along the $(z_1,\ldots,z_k)$-direction.
Shrinking $S$ in the radius direction,
we may assume that
$\bigl|\lambda^{-1}\cdot B_{\gminia,S}\bigr|$ 
are sufficiently smaller than 
$\bigl|\Re(\lambda^{-1}(\gminia-\gminib))\bigr|$ 
on $S$
for any
$\gminia,\gminib\in\nbigi^{(\lambda_0)}$
such that $\gminia<_S\gminib$.
By shrinking $X$,
we may assume that
$S$ is of the form
$\prod_{i=1}^k\Sec[1,\theta_i^{(0)},\theta_i^{(1)}]
 \times \Delta^{n-k}\times U(\lambda_0)$
in the following argument,
where $U(\lambda_0)$ is a neighbourhood of
$\lamda_0$ in $\cnum_{\lambda}$
in the case $\lambda_0\neq 0$,
or a small sector 
$\Sec[\delta_{\lambda},
 \theta^{(0)}_{\lambda},
 \theta^{(1)}_{\lambda}]
\subset\Delta_{\lambda}^{\ast}$
in the case $\lambda_0=0$.

Let us show the following statement
by a descending induction on $m\geq k$:
\begin{description}
\item[A(m)]
 Let $f$ be a $\DD_{\leq k}$-flat section of
 $\nbigp_0^{(\lambda_0)}\End\nbige_{|S}$
such that
\begin{itemize}
\item
 $\Res_i(\DD)f_{|\nbigd^{(\lambda_0)}_i\cap S}=0$
 for $i=m+1,\ldots,\ell$,
\item
 $f\bigl(\nbigf^S_{\gminia}\nbigpzero_0\nbige\bigr)
 \subset \nbigf^S_{<\gminia}\nbigpzero_0\nbige$
for any $\gminia$.
\end{itemize}
If we shrink $S$ in the radius direction,
we have the following estimate
for some $C>0$, $\epsilon>0$ and $N>0$:
\[
 |f|_h\leq C\cdot
 \exp\bigl(-\epsilon|\lambda^{-1}\vecz^{\vecm(0)}|\bigr)
 \cdot\left(
-\sum_{i=k+1}^{m}\log|z_i|
 \right)^N
\]
\end{description}

First, let us show $A(\ell)$.
We have the expression
$f=\sum_{\gminia<_S\gminib}
 f_{\gminia,\gminib,i,j}\cdot u_{\gminia,i}
 \otimes u_{\gminib,j}^{\lor}$.
Because $f$ is $\DD_{\leq k}$-flat,
we have
$\bigl|f_{\gminia,\gminib,i,j}\bigr|=
O\bigl(
 \exp(-\epsilon_1|\lambda^{-1}\vecz^{\vecm(0)}|)
 \bigr)$
due to Corollary \ref{cor;07.12.27.4}
and Corollary \ref{cor;07.12.17.10}.
Let $\vecv$ be a frame of
$\nbigp^{(\lambda_0)}_0\nbige$
for which we have the expression
$f=\sum \ftilde_{i,j}\cdot v_i\otimes v_j^{\lor}$.
Let $B$ be determined by
$\vecu_S=\vecv\cdot B$,
and then $B$ and $B^{-1}$ are bounded.
Hence, we have
$\bigl|\ftilde_{i,j}\bigr|=
 O\bigl(
 \exp(-\epsilon_2|\lambda^{-1}\vecz^{\vecm(0)}|)
 \bigr)$.
We put
\[
\nbigz:=\bigl\{
(\lambda,z_1,\ldots,z_n)\in S\,\big|\,
 1/2\leq |z_i|\leq 1,\,\,i=k+1,\ldots,\ell
 \bigr\}. 
\]

\begin{lem}
\label{lem;07.7.7.10}
We have the estimate
$|f|_h=O\Bigl(
 \exp\bigl(
 -\epsilon_3|\lambda^{-1}\vecz^{\vecm(0)}|
 \bigr)
 \Bigr)$
on $\nbigz$.
\end{lem}
\pf
Let $\nbigp^{(\lambda_0)}h$ be the metric
as in Section \ref{subsection;07.11.23.1}.
By construction of $\nbigpzero\nbige$,
we have
$|v_i|_{\nbigp^{(\lambda_0)}h}\leq 
 C\cdot \prod_{i=1}^k|z_i|^{-\delta}$
 on $\nbigz$ for some $\delta>0$.
Hence, we have
$|f|_{\nbigp^{(\lambda_0)}h}
=O\bigl(
 \exp(-\epsilon_4|\lambda^{-1}\vecz^{\vecm(0)}|)
 \bigr)$
on $\nbigz$ for some $\epsilon_4>0$.
The metrics
$\nbigp^{(\lambda_0)}h$
and $h$ are mutually bounded up to
$\exp\bigl(
 C'|\lambda-\lambda_0| \cdot |\vecz^{\vecm(0)}|
 \bigr)$-order.
Hence,
we obtain the estimate with respect to $h$
by shrinking $S$ appropriately.
\hfill\qed

\vspace{.1in}

Let us return to the proof of
Lemma \ref{lem;07.6.2.60}.
Let $\nbigz_1:=
 \bigcap_{i=k+1}^{\ell}(\nbigd_i\cap S)$.
Let $\pi_{\nbigz_1}$ denote the projection
$S\lrarr\nbigz_1$.
Let us consider the restriction to
$\pi_{\nbigz_1}^{-1}(\lambda,Q)$ for
$(\lambda,Q)\in \nbigz_1$.
The metrized holomorphic bundles
$(\End(\nbige),h)_{|\pi_{\nbigz_1}^{-1}(\lambda,Q)}$
are acceptable,
whose curvatures are dominated uniformly for 
$(\lambda,Q)\in \nbigz_1$.
Hence, we obtain the following estimate
due to Proposition \ref{prop;07.7.7.3}:
\[
 \bigl|f_{|\pi^{-1}(\lambda,Q)}
 (z_{k+1},\ldots,z_{\ell})\bigr|_h
\leq C\!\!
\max_{\substack{|z_i'|=1/2\\i=k+1,\ldots,\ell}}\!\!
 \bigl|
 f_{|\pi^{-1}(\lambda,Q)}(z'_{k+1},\ldots,z'_{\ell})
 \bigr|_h\!
 \left(-\sum_{i=k+1}^{\ell}
 \log|z_i|\right)^N
\]
Here, the constant $C$ can be independent of
$\lambda$ and $Q$.
Thus, we obtain $A(\ell)$.

\vspace{.1in}
Let us show $A(m-1)$,
by assuming $A(m)$.
We put $g:=\DD(\del_m)f$.
Since $g$ satisfies the assumption of $A(m)$,
we have the following estimate:
\[
 |g|_h\leq
 C\cdot\exp\bigl(
 -\epsilon_5\cdot|\lambda^{-1}\vecz^{\vecm(0)}|
 \bigr)
 \cdot
\Bigl(
 -\sum_{i=k+1}^m\log|z_i|
\Bigr)^N
\]
Let $\pi_m:S\lrarr  \nbigd_m^{(\lambda_0)}$
denote the projection.
We put $\pi_m^{-1}(\lambda,Q)^{\ast}
 :=\pi_m^{-1}(\lambda,Q)-\{(\lambda,Q)\}$
for $(\lambda,Q)\in \pi_m(S)$.
Let $\Delta_m$ denote the Laplacian
$-\del_{z_m}\cdot\del_{\zbar_m}$.
We have
$\Delta_m\bigl|f_{|\pi_m^{-1}(\lambda,Q)}\bigr|_h^2
\leq
 \bigl|g_{|\pi_m^{-1}(\lambda,Q)}\bigr|_h^2$
by Corollary \ref{cor;07.11.23.10}
on $\pi_m^{-1}(\lambda,Q)^{\ast}$.
Because $\Res_m(\DD)(f_{|(\lambda,Q)})=0$,
we obtain the boundedness of
the section $f_{|\pi_m^{-1}(\lambda,Q)}$
by the norm estimate in the curve case.
(See Proposition \ref{prop;07.10.6.40}.)
We can take $G_{\lambda,Q}$
with the following property:
\[
 \Delta_m G_{\lambda,Q}
=\bigl|g_{|\pi_m^{-1}(\lambda,Q)}\bigr|_h^2,
\quad\quad
 |G_{\lambda,Q}|
\leq
 C\cdot
 \exp\bigl(-\epsilon_6|\lambda^{-1}\vecz^{\vecm(0)}|\bigr)
\left(-\sum_{i=k+1}^{m-1}\log|z_i|\right)^N
\]
Here, the constant $C$ is independent of
the choice of $Q$ and $\lambda$.
Then, we have
$\Delta_m\bigl(
 \bigl|f_{|\pi_m^{-1}(\lambda,Q)}\bigr|_h^2
 -G_{\lambda,Q}
 \bigr)\leq 0$,
and hence we obtain the following:
\begin{multline*}
 \bigl|f_{|\pi_m^{-1}(\lambda,Q)}(z_m)\bigr|_h^2
\leq
 \max_{|z_m'|=1/2}
 \bigl|f_{|\pi_m^{-1}(\lambda,Q)}(z_m')\bigr|_h^2
 \\
+2C\exp(-\epsilon_7|\lambda^{-1}\vecz^{\vecm(0)}|)
 \left(-\sum_{i=k+1}^{m-1}\log|z_i|\right)^N
\end{multline*}
By the hypothesis of the induction on $\ell$,
we may assume to have the desired estimate
for the restriction of $f$ to 
$\bigl\{\vecz\in X-D\,\big|\,1/3<|z_m|<2/3\bigr\}
 \times U(\lambda_0)$.
Thus, we obtain $A(m-1)$,
and $A(k)$ means the claim of
Lemma \ref{lem;07.6.2.60}.
\hfill\qed

\vspace{.1in}

We can show the following lemmas
by using an argument
in the proof of Lemma \ref{lem;07.6.2.60}.

\begin{lem}
Let $S$ be a small multi-sector in
$\nbigxzero-W(\leq k)$.
Let $f$ be a $\DD_{\leq k}$-flat section of
$\nbigf^S_{<0}
 \nbigp_0^{(\lambda_0)}\nbige_{|S}$
such that
$\Res_i(\DD)
 f_{|\nbigd^{(\lambda_0)}_i\cap S}=0$
for $i=k+1,\ldots,\ell$.
When we shrink $S$ in the radius direction,
we have the estimate
$|f|_h\leq C
 \exp\bigl(-\epsilon|\lambda^{-1}\vecz^{\vecm(0)}|\bigr)$
for some $C>0$ and $\epsilon>0$.
\hfill\qed
\end{lem}

\begin{lem}
Let $\lambda\neq 0$.
Let $S$ be a small multi-sector of $X-D(\leq k)$,
such that 
we have the partial Stokes filtration
$\nbigf^S$ of $\nbigp\nbigelambda_{|\Sbar}$
in the level $\vecm(0)$.
Let $f$ be a $\DDlambda_{\leq k}$-flat section of
$\nbigf^S_{<0}\nbigelambda_{|\Sbar}$ such that
$\Res_i(\DDlambda)f_{|D_i\cap \Sbar}=0$
for $i=k+1,\ldots,\ell$.
When we shrink $S$ in the radius direction,
we have 
$ |f|_h\leq C
 \exp\bigl(-\epsilon|\vecz^{\vecm(0)}|\bigr)$
for some $C>0$ and $\epsilon>0$.
\hfill\qed
\end{lem}

\section{Locally uniform comparison
 of the irregular decompositions}
\label{subsection;08.9.14.13}
We will compare the irregular decompositions
of $(\nbigp\nbige^0,\DD^0)$
and $(\nbigpzero\nbige,\DD)$.
The main results of this section are 
Lemma \ref{lem;07.7.7.57},
Corollary \ref{cor;08.9.14.20},
Lemma \ref{lem;07.7.8.5} and
Corollary \ref{cor;08.9.14.21}.
Because we will not use them
in the rest of this paper,
the reader can skip here.

\subsection{Around $\lambda_0\neq 0$}
\label{subsection;08.9.28.6}

We continue to use the setting in Section 
\ref{subsection;07.12.22.1}.
Let us consider the case $\lambda_0\neq 0$.

\begin{lem}
\label{lem;07.7.7.55}
\mbox{{}}
Let $S$ be a small multi-sector in
$\nbigx^{(\lambda_0)}
-\nbigd^{(\lambda_0)}(\leq k)$.
\begin{itemize}
\item
We can take a $\DD_{\leq k}$-flat splitting 
$\nbigp^{(\lambda_0)}
 \nbige_{|\Sbar}
=\bigoplus
 \nbigp^{(\lambda_0)}\nbige_{\gminia,S}$
of the Stokes filtration $\nbigf^S$
in the level $\vecm(0)$,
whose restriction to
$S\cap\nbigd^{(\lambda_0)}_i$
is compatible with the residues
$\Res_{i}(\DD)$
and the filtrations $\lefttop{i}\Fzero$
for $i=k+1,\ldots,\ell$.
\item
If $\lambda_0$ is generic,
then we can take a $\DD$-flat splitting
with the above property.
\end{itemize}
\end{lem}
\pf
It follows from Proposition \ref{prop;07.9.30.2}
and Proposition \ref{prop;07.9.30.5}.
\hfill\qed

\vspace{.1in}
Let $\nbigx^{\lambda}:=\{\lambda\}\times X$.
We use the symbols like $\nbigd^{\lambda}(\leq k)$
in similar meanings.
We also use the symbol $S^{\lambda}$
to denote $S\cap \nbigx^{\lambda}$
for a multi-sector $S$ in 
$\nbigx^{(\lambda_0)}-\nbigd^{(\lambda_0)}(\leq k)$.

\begin{lem}
\label{lem;07.12.22.11}
Let $S$ be a small multi-sector in
$\nbigx^{(\lambda_0)}
-\nbigd^{(\lambda_0)}(\leq k)$.
Let $\nbigp^{(\lambda_0)}\nbige_{|\Sbar}
=\bigoplus \nbigpzero\nbige_{\gminia,S}$
be a $\DD_{\leq k}$-flat splitting as in
Lemma {\rm\ref{lem;07.7.7.55}}.
If $|\lambda-\lambda_0|$ is sufficiently small,
the restriction of the splitting to 
$\nbigx^{\lambda}$
gives a splitting
of the Stokes filtration
$\nbigf^{S^{\lambda}}$
of $\nbigp\nbigelambda_{|S^{\lambda}}$.
\end{lem}
\pf
The restriction
$\nbigp^{(\lambda_0)}
 \nbigelambda_{|\Sbar^{\lambda}}
=\bigoplus 
 \nbigpzero\nbige_{\gminia,S|\Sbar^{\lambda}}$
gives a splitting of
the Stokes filtrations 
$\nbigf^{S^{\lambda}}
 \bigl(\nbigpzero\nbigelambda_{|\Sbar^{\lambda}}\bigr)$
of $\nbigpzero\nbigelambda_{|\Sbar^{\lambda}}$.
Then, the claim follows from
Lemma \ref{lem;10.5.26.10}.
\hfill\qed

\vspace{.1in}

Let $p^{\vecm(0)}_{\gminia,S}$ denote the projection
onto $\nbigp^{(\lambda_0)}\nbige_{\gminia,S}$.
Let $p_{\gminia,S}^{\vecm(0)\prime}$
come from another decomposition
with the first property
in Lemma \ref{lem;07.7.7.55}.

\begin{lem}
\label{lem;07.7.7.56}
When we shrink $S$ in the radius direction,
we have 
$p^{\vecm(0)}_{\gminia,S}
-p^{\vecm(0)\prime}_{\gminia,S}=
 O\bigl(\exp(-\epsilon|\vecz^{\vecm(0)}|)\bigr)$
with respect to $h$ on 
$S\setminus \nbigd^{(\lambda_0)}$
for some $\epsilon>0$.
\end{lem}
\pf
We have
(i) $\bigl(p^{\vecm(0)}_{\gminia,S}
 -p^{\vecm(0)\prime}_{\gminia,S}\bigr)
\nbigf^S_{\gminia}\subset\nbigf^S_{<\gminia}$,
(ii) 
$\bigl[
 p^{\vecm(0)}_{\gminia,S}
-p^{\vecm(0)\prime}_{\gminia,S},
 \Res_i(\DD)\bigr]=0$
on $S\cap \nbigd^{(\lambda_0)}_i$
for each $i=k+1,\ldots,\ell$,
(iii) the restriction of
 $p^{\vecm(0)}_{\gminia,S}
  -p^{\vecm(0)\prime}_{\gminia,S}$
to $S\cap \nbigd^{(\lambda_0)}_i$
 preserves the filtration
for each $i=k+1,\ldots,\ell$.
Hence, the claim follows from Lemma \ref{lem;07.6.2.60}.
\hfill\qed

\vspace{.1in}

We take small multi-sectors
$S_j$ $(j=1,\ldots,N)$
such that the union of the interior point of
$S_j$ is $\nbigv_0\setminus\nbigd^{(\lambda_0)}(\leq k)$,
where $\nbigv_0$ denotes a neighbourhood
of $\nbigd^{(\lambda_0)}_{\kbar}$.
By gluing $p^{\vecm(0)}_{\gminia,S_j}$
in $C^{\infty}$ as in Section \ref{subsection;07.6.16.8},
we construct the $C^{\infty}$-map
$p^{\vecm(0)}_{\gminia,C^{\infty}}$.
Due to Lemma \ref{lem;07.7.7.56},
we have
$(\delbar_E+\lambda\theta^{\dagger})
 p^{\vecm(0)}_{\gminia,C^{\infty}}
=O\Bigl(
 \exp\bigl(-\epsilon|\vecz^{\vecm(0)}|\bigr)
 \Bigr)$
with respect to $h$ and the Poincar\'e metric $g_{\poin}$
on $\nbigv\setminus \nbigd^{(\lambda_0)}$,
where $\nbigv$ denotes some neighbourhood
of $\nbigd^{(\lambda_0)}_{\kbar}$.

\begin{lem}
\label{lem;07.7.7.57}
We have the estimate
$\bigl|\pi^{\vecm(0)}_{\gminia}
 -p^{\vecm(0)}_{\gminia,C^{\infty}}\bigr|_{h}
\leq C_N\cdot \prod_{i=1}^k|z_i|^N$
for any $N>0$
on $\nbigv_1\setminus \nbigd^{(\lambda_0)}$,
where $\nbigv_1$ denotes some neighbourhood
of $\nbigd^{(\lambda_0)}_{\kbar}$.
\end{lem}
\pf
By shrinking $X$,
we may assume 
$(\delbar_E+\lambda\theta^{\dagger})
 p^{\vecm(0)}_{\gminia,C^{\infty}}
=O\Bigl(
 \exp\bigl(-\epsilon|\vecz^{\vecm(0)}|\bigr)
 \Bigr)$ with respect to 
$h$ and $g_{\poin}$
on $\nbigx^{(\lambda_0)}-\nbigd^{(\lambda_0)}$.
Let $\pi:\nbigx^{(\lambda_0)}-\nbigd^{(\lambda_0)}
 \lrarr 
 \nbigd^{(\lambda_0)}_{\kbar}$ 
denote the natural projection.
Then, the restrictions
$(\nbige,h)_{|\pi^{-1}(\lambda,Q)}$ are acceptable,
and the curvatures are dominated uniformly for
$(\lambda,Q)\in 
\pi(\nbigx^{(\lambda_0)}-\nbigd^{(\lambda_0)})$.
We also have
$\bigl(
(\delbar_E+\lambda\theta^{\dagger})
 \pi^{\vecm(0)}_{\gminia}\bigr)_{|\pi^{-1}(\lambda,Q)}
=O\bigl(\exp(-\epsilon|\vecz^{\vecm(0)}|)\bigr)$
with respect to $(h,g_{\poin})$,
which is uniform for $(\lambda,Q)$.
Thus, we obtain the following estimate
uniformly for $(\lambda,Q)$:
\[
 \Bigl(
 \bigl(\delbar_E+\lambda\theta^{\dagger}\bigr)
 \bigl(p^{\vecm(0)}_{\gminia,C^{\infty}}
 -\pi^{\vecm(0)}_{\gminia}\bigr)
\Bigr)_{|\pi^{-1}(\lambda,Q)}
=O\bigl(
 \exp(-\epsilon|\vecz^{\vecm(0)}|)\bigr)
\]
Let $\nbigz:=\bigl\{
 (\lambda,z_1,\ldots,z_n)\,\big|\,
 1/2<|z_i|<1\,\,(i=1,\ldots,k)
 \bigr\}$.
We obviously have the boundedness of
$\pi^{\vecm(0)}_{\gminia}$
on $\nbigz$.
Due to Lemma \ref{lem;07.6.2.60}
with $k=0$,
we also have the boundedness of
$p^{\vecm(0)}_{\gminia,S}$
on $S\cap \nbigz$,
and hence we obtain 
the boundedness of 
$p^{\vecm(0)}_{\gminia,C^{\infty}}$
on $\nbigz$.
Moreover,
we have
$\bigl|p^{\vecm(0)}_{\gminia,C^{\infty}}
 -\pi^{\vecm(0)}_{\gminia}\bigr|_{h|
 \pi^{-1}(\lambda,Q)}
\leq
 C_{(\lambda,Q),N}\cdot \prod_{i=1}^k|z_i|^N$
for any $(\lambda,Q)\in\nbigd^{(\lambda_0)}_{\kbar}$
and any $N>0$,
due to Corollary \ref{cor;07.7.7.43}
and Lemma \ref{lem;07.12.22.11}.
Then, the claim of the lemma follows from
Lemma \ref{lem;07.7.7.40}.
\hfill\qed

\begin{cor}
\label{cor;08.9.14.20}
Let $S$ be a small sector in 
$\nbigx^{(\lambda_0)}-\nbigd^{(\lambda_0)}(\leq k)$.
Let $\nbigp^{(\lambda_0)}
 \nbige_{|\Sbar}
=\bigoplus
 \nbigp^{(\lambda_0)}\nbige_{\gminia,S}$
be a splitting of the Stokes filtration $\nbigf^S$
in the level $\vecm(0)$ as in Lemma 
{\rm\ref{lem;07.7.7.55}}.
We have the estimate
$\bigl|\pi^{\vecm(0)}_{\gminia}
 -p^{\vecm(0)}_{\gminia,S}\bigr|_{h}
\leq C_N\cdot \prod_{i=1}^k|z_i|^N$
for any $N>0$
on $(\nbigv_1\cap S)\setminus \nbigd^{(\lambda_0)}$.
\hfill\qed
\end{cor}

\subsection{Around $\lambda_0=0$}

Let us consider the case $\lambda_0=0$.

\begin{lem}
\label{lem;07.7.7.65}
Let $S=S_{\lambda}\times S_{\vecz}$ 
be a sufficiently small multi-sector 
in $\nbigx^{(0)}-W(\leq k)$,
where $S_{\lambda}$ and $S_{\vecz}$ denote
sectors in $U(0)-\{0\}$ and
$X-D$, respectively.
\begin{itemize}
\item
We have a $\DD_{\leq k}$-flat splitting
$\nbigp^{(0)}\nbige_{|\Sbar}
=\bigoplus
 \nbigp^{(0)}\nbige_{\gminia,S}$
of $\nbigf^S$ whose restriction to 
$S\cap\nbigd^{(0)}_i$
is compatible with the residues
$\Res_{i}(\DD)$
and the filtrations $\lefttop{i}\Fzero$
for $i=k+1,\ldots,\ell$.
\item
For any $\lambda\in S_{\lambda}$,
the restriction
$\nbigp^{(0)}\nbigelambda_{|S^{\lambda}}
=\bigoplus
 \nbigp^{(0)}\nbige_{\gminia,S|S^{\lambda}}$ 
gives a splitting of the Stokes filtration
of $\nbigp\nbigelambda_{|S^{\lambda}}$
with the above property.
\end{itemize}
\end{lem}
\pf
The first claim follows from
Proposition \ref{prop;07.9.30.2}.
The restriction
$\nbigp^{(0)}\nbigelambda_{|\Sbar^{\lambda}}
=\bigoplus
 \nbigp^{(0)}\nbige_{\gminia,S|\Sbar^{\lambda}}$ 
gives a splitting of the Stokes filtration
$\nbigf^{S^{\lambda}}\bigl(
 \nbigp^{(0)}\nbigelambda_{|\Sbar^{\lambda}}
 \bigr)$
of $\nbigp^{(0)}\nbigelambda_{|\Sbar^{\lambda}}$.
The filtration 
$\nbigf^{S^{\lambda}}\bigl(
 \nbigp\nbigelambda_{|S^{\lambda}}
 \bigr)$ given by
\[
  \nbigf^{S^{\lambda}}_{(1+|\lambda|^2)\gminia}
 \bigl(
 \nbigp\nbigelambda_{|S^{\lambda}}
 \bigr):=
 \nbigf^{S^{\lambda}}_{\gminia}
 \bigl(
 \nbigp^{(0)}\nbigelambda_{|S^{\lambda}}
 \bigr)
\]
is the same as the Stokes filtration of
$\nbigp\nbigelambda_{|S^{\lambda}}$.
Hence, the second claim follows.
\hfill\qed

\vspace{.1in}

Let $p^{\vecm(0)}_{\gminia,S}$ denote the projection
onto $\nbigp^{(0)}\nbige_{\gminia,S}$.
Let $p_{\gminia,S}^{\vecm(0)\prime}$
come from another decomposition with the property
in Lemma \ref{lem;07.7.7.65}.

\begin{lem}
\label{lem;07.7.7.66}
When we shrink $S$ in the radius direction,
we have 
$p^{\vecm(0)}_{\gminia,S}
-p^{\vecm(0)\prime}_{\gminia,S}=
 O\Bigl(
\exp\bigl(-\epsilon|\lambda^{-1}\vecz^{\vecm(0)}|\bigr)
 \Bigr)$
with respect to $h$
on $S\setminus\nbigd^{(0)}$.
\end{lem}
\pf
We have
(i) $\bigl(p^{\vecm(0)}_{\gminia,S}
 -p^{\vecm(0)\prime}_{\gminia,S}\bigr)
\nbigf^S_{\gminia}\subset\nbigf^S_{<\gminia}$,
(ii) 
$\bigl[
 p^{\vecm(0)}_{\gminia,S}-p^{\vecm(0)\prime}_{\gminia,S},
 \Res_i(\DD)\bigr]=0$
on $S\cap \nbigd^{(0)}_i$ for $i=k+1,\ldots,\ell$,
(iii)
the restriction of 
$p^{\vecm(0)}_{\gminia,S}-p^{\vecm(0)\prime}_{\gminia,S}$
to $S\cap \nbigd^{(0)}_i$
preserves the filtration $\lefttop{i}\Fzero$
for each $i=k+1,\ldots,\ell$.
Then, the claim follows from Lemma \ref{lem;07.6.2.60}.
\hfill\qed

\vspace{.1in}

We take small multi-sectors
$S_j$ $(j=1,\ldots,N)$
such that the union of the interior point of
$S_j$ is $\nbigv_0\setminus W(\leq k)$,
where $\nbigv_0$ denotes a neighbourhood
of $\{0\}\times D_{\kbar}$.
By gluing $p^{\vecm(0)}_{\gminia,S_j}$
in $C^{\infty}$ as in Section \ref{subsection;07.6.16.8},
we construct the $C^{\infty}$-map
$p^{\vecm(0)}_{\gminia,C^{\infty}}$
on $\nbigv_0\setminus W(\leq k)$,
which is extended to the $C^{\infty}$-map
on $\nbigv_0\setminus \nbigd^{(0)}$.
Shrinking $\nbigv_0$,
we have the estimate
$ (\delbar_E+\lambda\theta^{\dagger}+\delbar_{\lambda})
 p^{\vecm(0)}_{\gminia,C^{\infty}}
=O\Bigl(
 \exp\bigl(-\epsilon_1|\lambda^{-1}\vecz^{\vecm(0)}|\bigr)
 \Bigr)$
with respect to $h$ and $g_{\poin}$,
due to Lemma \ref{lem;07.7.7.66}.

\begin{lem}
\label{lem;07.7.8.5}
We have the estimate
$\bigl|\pi^{\vecm(0)}_{\gminia}
 -p^{\vecm(0)}_{\gminia,C^{\infty}}\bigr|_{h}
\leq C_N\cdot |\lambda|\cdot\prod_{i=1}^k|z_i|^N$
for any $N>0$
on $\nbigv_1\setminus \nbigd^{(0)}(\leq k)$,
where $\nbigv_1$ denotes some neighbourhood
of $\{0\}\times D_{\kbar}$.
\end{lem}
\pf
Take small $0<\delta_1<\delta_2$.
Due to Lemma \ref{lem;07.7.7.57},
we may assume to have
$\bigl|\pi^{\vecm(0)}_{\gminia}
 -p^{\vecm(0)}_{\gminia,C^{\infty}}\bigr|_{h}
\leq C_N'\cdot
 \prod_{i=1}^k|z_i|^N$
with respect to $h$ on
$\bigl\{\lambda\,\big|\,
 \delta_1\leq |\lambda|\leq \delta_2\bigr\}
  \times (X-D)$ for any $N>0$
by shrinking $X$.
We may assume to have the following estimate
on $\{0<|\lambda|<\delta_2\}\times (X-D)$,
with respect to $h$ and the Euclidean metric 
$d\lambda\cdot d\lambdabar$:
\[
g(\vecz,\lambda):=
 \delbar_{\lambda}
 \bigl(
 \pi^{\vecm(0)}_{\gminia}
-p^{\vecm(0)}_{\gminia,C^{\infty}}
 \bigr)
=-\delbar_{\lambda}
 p^{\vecm(0)}_{\gminia,C^{\infty}}
=O\Bigl(
 \exp\bigl(
 -\epsilon_2|\lambda^{-1}\vecz^{\vecm(0)}|\bigr)
 \Bigr)
\]

Let $\vece$ be an orthonormal frame of $\End(E)$ on $X-D$
with respect to the metric induced by $h$.
We have the expression
$g=\sum g_i\cdot e_i$.
We have the estimate
$ \bigl|g_i\bigr|\leq
 C\cdot\exp\bigl(
 -\epsilon_3|\lambda^{-1}\vecz^{\vecm(0)}|
 \bigr)$.
Let $\chi$ be a positive valued $C^{\infty}$-function
on $\cnum_{\lambda}$
such that $\chi(\lambda)=1$ for 
$|\lambda|\leq(\delta_1+\delta_2)/2$
and $\chi(\lambda)=0$ for $|\lambda|\geq \delta_2$.
We put 
\[
 G_i(\vecz,\lambda):=
 \lambda
 \int_{\cnum}
 \frac{ g_i(\vecz,\mu)\cdot\mu^{-1}\cdot \chi(\mu)}
 {\mu-\lambda}
\frac{\sqrt{-1}}{2\pi}
 d\mu\cdot d\mubar
\]
We put $G:=\sum G_i\cdot e_i$.
Then, we have $\delbar_{\lambda}G=g$
on $\bigl\{|\lambda|<(\delta_1+\delta_2)/2\bigr\}
 \times (X-D)$.
We also have the following estimate:
\[
 \bigl|G(\vecz,\lambda)\bigr|_h\leq
 C\cdot 
 |\lambda|\cdot
 \exp\bigl(-\epsilon_4\cdot|\vecz^{\vecm(0)}|\bigr)
\]
Let us see
$H:=\pi^{\vecm(0)}_{\gminia}
 -p^{\vecm(0)}_{\gminia,C^{\infty}}-G$.
We have
$\delbar_{\lambda}H=0$ and $H(\vecz,0)=0$.
By a lemma of Schwarz,
we obtain
$ \bigl|H(\vecz,\lambda)\bigr|_h
\leq
 C_N\cdot
 |\lambda|\cdot\prod_{i=1}^k|z_i|^N$
for any $N$.
Then, we obtain the desired estimate.
\hfill\qed

\begin{cor}
\label{cor;08.9.14.21}
We have the estimate
$\bigl|\pi^{\vecm(0)}_{\gminia}
 -p^{\vecm(0)}_{\gminia,S}\bigr|_{h}
\leq C_N\cdot |\lambda|\cdot\prod_{i=1}^k|z_i|^N$
for any $N>0$.
\hfill\qed
\end{cor}

\chapter{Smooth Divisor Case}
\label{section;08.10.23.1}
Let $\harmonicbundle$ be
an unramifiedly good wild harmonic bundle 
on $X-D$,
where $X$ is a complex manifold
and $D$ is a simple normal crossing hypersurface.
In Chapter \ref{section;07.11.24.24},
we studied the prolongment
$(\nbigpzero\nbige,\DD)$ of $(\nbige,\DD)$,
i.e.,
the sheaf of holomorphic sections
whose norms are of polynomial growth
with respect to the hermitian metric
$\nbigpzero h$.
In this chapter,
we restrict ourselves 
to the case that $D$ is smooth,
and we will do more refined analysis.
The results in this chapter are rather technical,
and preliminary for the later sections.

\vspace{.1in}
In Section \ref{subsection;08.9.15.40},
we show that
$(\nbige,\nbigpzero h)$ is acceptable
(Proposition \ref{prop;07.11.22.45}),
which is new even in the tame case.
For the proof,
we obtain a complementary estimate
of the connection form
(Lemma \ref{lem;07.7.18.150}
 and Corollary \ref{cor;07.7.19.100})
which is also useful for other purposes.

Recall we have studied
in Section \ref{subsection;08.9.14.13}
the comparison of irregular decompositions
for $(\nbigp\nbige^0,\DD^0)$
and $(\nbigpzero\nbige,\DD)$
in the level $\vecm(0)$.
In Section \ref{subsection;08.9.15.41},
we compare irregular decompositions in any level
under the assumption that $D$ is smooth
(Proposition \ref{prop;08.9.17.2}
and Proposition \ref{prop;08.9.17.3}).
This comparison will be used in
the proof of Theorem \ref{thm;07.10.11.120}
in one way.
We also use it for comparison
of the hermitian metrics with some twisting
(Corollary \ref{cor;07.11.24.3}),
which will be used for the family version
of norm estimate 
(Section \ref{subsection;08.9.15.43}).

In Section \ref{subsection;08.9.15.42},
we show a standard norm estimate 
(Proposition \ref{prop;07.11.22.55})
for $(\nbigp\nbigelambda,\DDlambda)$
under the assumption that $D$ is smooth,
which is a generalization of 
Proposition \ref{prop;07.10.6.40},
and preliminary for the family version
of norm estimate 
(Section \ref{subsection;08.9.15.43}).
We also show an estimate of a connection form
(Proposition \ref{prop;08.9.15.44}),
which can be skipped.

In Section \ref{subsection;08.9.15.43},
we show the family version of norm estimate
(Proposition \ref{prop;07.10.22.30}).
This is preliminary for Step 1 in the proof of
Theorem \ref{thm;07.10.23.20}
(Section \ref{subsection;07.10.23.30}).

\vspace{.1in}

{\em During this chapter,
we use the setting in Section
{\rm\ref{subsection;07.11.17.1}} with $\ell=1$,
i.e., $D$ is assumed to be smooth.}

\section{Acceptability of $\nbigpzero h$}
\label{subsection;08.9.15.40}
\subsection{Statement}
Let $U(\lambda_0)$ denote a small neighbourhood
of $\lambda_0$.
We will prove the following proposition 
in Section \ref{subsection;08.1.10.1}.
\begin{prop}
\label{prop;07.11.22.45}
$(\nbige,\nbigpzero h)$ is acceptable
on $(X-D)\times U(\lambda_0)$,
if $U(\lambda_0)$ is sufficiently small.
\end{prop}

In particular,
$\nbigpzero_{\ast}\nbige
 =\bigl(
 \nbigpzero_a\nbige\,\big|\,
 a\in\real
 \bigr)$ is a filtered bundle.
We have the induced filtration $\Fzero$
of $\nbigpzero_a\nbige_{|U(\lambda_0)\times D}$.
Let $\vecv$ be a frame of 
$\nbigpzero_{\veca}\nbige$
which is compatible with 
the filtrations $\Fzero$.
We put
\[
v_i':=v_i\cdot |z_1|^{\deg^{\Fzero}(v_i)}.
\] 
Let $H(\nbigpzero h,\vecv')$ denote
the Hermitian matrix-valued function
whose $(i,j)$-th entries are
$\nbigpzero h(v_i',v_j')$.

\begin{cor}
\label{cor;08.9.15.13}
We have the weak norm estimate,
i.e.,
$C_1\cdot(-\log|z_1|)^{-N}
\leq
 H(\nbigpzero h,\vecv')
\leq
 C_2\cdot \bigl(-\log|z_1|\bigr)^N$
for some $C_1,C_2,N>0$.
\end{cor}
\pf
It follows from the general result
for the acceptable bundles
(Theorem \ref{thm;07.10.9.1}).
\hfill\qed

\subsection{Estimate of the connection form}
\label{subsection;07.10.11.22}

Let
$\eta_m:M(X,D)\lrarr M(X,D)$  be given by
$\eta_m(\gminia)=
 \sum_{m(0)\leq j\leq m}\gminia_j\cdot z_1^{j}$.
The image $\eta_m\bigl(\Irr(\theta)\bigr)$ is denoted by
$\Irr(\theta,m)$.
We take total orders 
$\leq$ on $\Irr(\theta)$ and $\Irr(\theta,m)$
such that $\eta_m$ are order preserving.
We put 
\[
 E^{(m)}_{\gminib}=
 \bigoplus_{\substack{\gminia\in\Irr(\theta)\\
 \eta_m(\gminia)=\gminib }}
 E_{\gminia},
\quad
 F^{(m)}_{\gminib}:=
 \bigoplus_{\substack{\gminic\in\Irr(\theta,m)\\
 \gminic\leq \gminib}} E^{(m)}_{\gminib},
\quad
 F^{(m)}_{<\gminib}:=\bigoplus_{\substack{\gminic\in\Irr(\theta,m)\\
 \gminic<\gminib }} E^{(m)}_{\gminib}
\]
Let $\pi^{(m)}_{\gminib}$ denote the projection
onto $E^{(m)}_{\gminib}$
with respect to the decomposition
$E=\bigoplus_{\gminib\in\Irr(\theta,m)} E^{(m)}_{\gminib}$.
We also put
$E^{(0)}_{(\gminia,\alpha)}:=E_{(\gminia,\alpha)}$.
Let $\pi^{(0)}_{\gminia,\alpha}$ denote the projection
onto $E^{(0)}_{(\gminia,\alpha)}$
with respect to the decomposition (\ref{eq;07.7.18.100}).

For simplicity,
we use the symbol $\vecgamma$
to denote an element of $\Irr(\theta)\times\Sp(\theta)$.
We have the decomposition
$\End(E)=\nbigc^{(0)}\bigl(\End(E)\bigr)
\oplus \nbigd^{(0)}\bigl(\End(E)\bigr)$:
\[
 \nbigd^{(0)}\bigl(\End(E)\bigr):=
 \bigoplus_{\vecgamma\in\Irr(\theta)\times\Sp(\theta)}
 \End\bigl(E^{(0)}_{\vecgamma}\bigr),
\]
\[
 \nbigc^{(0)}\bigl(\End(E)\bigr):=
 \bigoplus_{\substack{
 \vecgamma,\vecgamma'\in\Irr(\theta)\times\Sp(\theta)\\
 \vecgamma\neq\vecgamma' }}
 \Hom\bigl(E^{(0)}_{\vecgamma},
 E^{(0)}_{\vecgamma'}\bigr)
\]
We use the symbols $\nbigd^{(m)}\bigl(\End(E)\bigr)$ and
$\nbigc^{(m)}\bigl(\End(E)\bigr)$ 
$(m\leq -1)$ in similar meanings.
For any section $G\in \End(E)\otimes\Omega^{p,q}$,
we have the corresponding decomposition
$G=\nbigd^{(m)}(G)+\nbigc^{(m)}(G)$.

Let $\del_E$ denote the $(1,0)$-operator
associated to $\delbar_E$ and $h$.
Let $\vecv$ be a frame of $\prolong{E}$
compatible with the decomposition (\ref{eq;07.7.18.100})
and the parabolic filtration.
Let $F\in\End(E)\otimes\Omega^{1,0}$
determined by $F(\vecv)=\del_E\vecv$.

\begin{lem}\mbox{{}}
\label{lem;07.7.18.150}
We have the following estimates on $X^{\ast}$
for some $\epsilon>0$:\\
\noindent
{\bf (The case $j<0$)}
\begin{equation}
 \label{eq;07.2.2.3}
\nbigc^{(j)}(F)=O\bigl(\exp(-\epsilon|z_1|^{j})\bigr),
\quad
\del_E\pi^{(j)}_{\gminib}
=O\bigl(\exp(-\epsilon|z_1|^{j})\bigr)
\end{equation}
{\bf (The case $j=0$)}
\begin{equation}
 \label{eq;07.2.2.4}
 \nbigc^{(0)}(F)=
 O\Bigl(|z_1|^{\epsilon}\frac{dz_1}{z_1}\Bigr),
\quad
\del_E\pi^{(0)}_{\gminia,\alpha}=
 O\bigl(|z_1|^{\epsilon}\bigr)
 \frac{dz_1}{z_1}
\end{equation}
\end{lem}
\pf
The estimate (\ref{eq;07.2.2.3}) is the special case of
Lemma \ref{lem;07.6.2.25}.
Let us consider the case $j=0$.
Since the argument is essentially the same,
we give only an outline.
In this proof,
$\epsilon$ will denote a positive constant,
and we will make it smaller without mention.
We have $\del_E\pi^{(0)}_{\vecgamma}
=[F,\pi^{(0)}_{\vecgamma}]\in 
 \nbigc^{(0)}(\End(E))\otimes\Omega^{1,0}$
and the estimate
\[
 \del_E\pi^{(0)}_{\vecgamma}
=O\Bigl((-\log|z_1|)^N\Bigr)
\]
for some $N$
with respect to $h$ and $g_{\poin}$.
Because 
$\nbigc^{(0)}(\End(E))$ and 
$\nbigd^{(0)}(\End(E))$
are $|z_1|^{\epsilon}$-asymptotically
orthogonal
(Theorem \ref{thm;07.10.4.1}
and Theorem \ref{thm;07.10.4.3}),
we obtain the following estimate 
with respect to $h$ and $g_{\poin}$:
\[
 \bigl(\pi^{(0)}_{\vecgamma},
 \del_E\pi^{(0)}_{\vecgamma}\bigr)
=O(|z_1|^{\epsilon}).
\]
We also have the following estimate
with respect to $h$ and $g_{\poin}$,
due to Corollary \ref{cor;07.11.22.30}:
\begin{equation}
\label{eq;07.11.22.40}
 \delbar_E\del_E\pi^{(0)}_{\vecgamma}
=\bigl[R(h),\pi^{(0)}_{\vecgamma}\bigr]
=O\bigl(|z_1|^{\epsilon}\bigr)
\end{equation}
We have the following equality:
\begin{multline*}
 \bigl(\del_E\pi^{(0)}_{\vecgamma},
 \del_E\pi^{(0)}_{\vecgamma}\bigr)_h
 |z_1|^{-\epsilon}
=\del\Bigl(
 \bigl(\pi^{(0)}_{\vecgamma},
 \del_E\pi^{(0)}_{\vecgamma}\bigr)_h
 |z_1|^{-\epsilon}
 \Bigr) \\
+\bigl(\pi^{(0)}_{\vecgamma},
 \del_E\pi^{(0)}_{\vecgamma}\bigr)_h
 \cdot\frac{\epsilon}{2}\cdot
 |z_1|^{-\epsilon}\frac{dz_1}{z_1}
-\bigl(\pi^{(0)}_{\vecgamma},
 \delbar_E\del_E\pi^{(0)}_{\vecgamma}\bigr)_h
 |z_1|^{-\epsilon}
\end{multline*}
Hence,
we obtain the following finiteness:
\begin{equation}
 \label{eq;07.11.22.41}
 \int \bigl(\del_E\pi_{\vecgamma}^{(0)},
 \del_E\pi^{(0)}_{\vecgamma}\bigr)_h
 \cdot|z_1|^{-\epsilon}<\infty
\end{equation}
From (\ref{eq;07.11.22.40})
and (\ref{eq;07.11.22.41}),
we obtain the desired estimate.
(See the argument in the proof of
Lemma \ref{lem;07.6.2.25}.)
\hfill\qed

\begin{cor}
\label{cor;07.7.19.100}
For $p< j$,
we have
$\nbigc^{(p)}(\del_E\pi^{(j)}_{\gminia})
=O\bigl(\exp(-\epsilon|z_1|^{p})\bigr)$.
\end{cor}
\pf
We have the following:
\[
 [F,\pi^{(j)}_{\gminia}]
=\bigl[
 \nbigc^{(p)}(F)+\nbigd^{(p)}(F),\,
 \nbigd^{(p)}(\pi^{(j)}_{\gminia})
 \bigr]
=O\bigl(
 \exp(-\epsilon|z_1|^{p})
 \bigr)
+\bigl[
 \nbigd^{(p)}(F),\nbigd^{(p)}(\pi^{(j)}_{\gminia})
 \bigr]
\]
Hence, $\nbigc^{(p)}\bigl([F,\pi^{(j)}_{\gminia}]\bigr)
=O\bigl(\exp(-\epsilon|z_1|^{p})\bigr)$,
and the claim of the corollary follows.
\hfill\qed

\subsection{Proof of Proposition \ref{prop;07.11.22.45}}
\label{subsection;08.1.10.1}

We use the notation in Section 
\ref{subsection;07.11.20.50}.
We set $\gtilde_{\poin}:=
 g_{\poin}+d\lambda\cdot d\lambdabar$.
We have only to show 
the curvature of $(\nbigtzero\nbige,h)$ is bounded
with respect to $h$ and $\gtilde_{\poin}$.
Let $\delta_{\lambda_0}'=\del_E-\lambdabar_0\theta$ 
denote the $(1,0)$-operator determined by
$d_{\lambda_0}''$ and $h$.
In the following,
$g(\lambda-\lambda_0)$
is also denoted by $g$.

\begin{lem}
\label{lem;08.9.15.11}
If $U(\lambda_0)$ is sufficiently small,
$\delta_{\lambda_0}'\Lambda(\lambda-\lambda_0)$
is  bounded 
with respect to $h$ and $g_{\poin}$
uniformly for $\lambda\in U(\lambda_0)$.
\end{lem}
\pf
In the proof,
we will make $U(\lambda_0)$ smaller
without mention.
Let $\Phibar$ be given by (\ref{eq;07.11.22.35}).
We have 
\[
 \Lambda(\lambda-\lambda_0)=
-(\lambda-\lambda_0)\cdot\Phibar
+\lambda\cdot g(\lambda-\lambda_0)\cdot
 \theta^{\dagger} \cdot g(\lambda-\lambda_0)^{-1}
-\lambda_0\cdot\theta^{\dagger}
\]
We have 
$\delta'_{\lambda_0}\theta^{\dagger}
=-\lambdabar_0\bigl[
 \theta,\theta^{\dagger}
 \bigr]$,
which is bounded 
(Proposition \ref{prop;07.11.22.25}).
Due to Lemma \ref{lem;07.7.18.150},
we have the following estimate
of $\delta'_{\lambda_0}\Phibar$
with respect to $h$ and $g_{\poin}$:
\[
  \bigl[\del_E-\lambdabar_0\theta,\Phibar\bigr]
=\!\!\!
 \sum_{m(0)\leq j\leq -1}
 \!\!\sum_{\gminib}
 d\overline{\zeta_{j}(\gminib)}
 \cdot\del_{E}\pi^{(j)}_{\gminib}
+\!\!\!\!
 \sum_{(\gminia,\alpha)\in\Irr(\theta)\times\Sp(\theta)}
 \!\!\!\!
 \alphabar\frac{d\zbar_1}{\zbar_1}\cdot
 \del_E\pi^{(0)}_{\gminia,\alpha}
=O\bigl(|z_1|^{\epsilon}\bigr)
\]

It remains to estimate 
$\delta'_{\lambda_0}\bigl(
 g(\lambda-\lambda_0)\cdot
 \theta^{\dagger}\cdot
 g(\lambda-\lambda_0)^{-1}
 \bigr)$,
which can be rewritten as follows:
\begin{equation}
\label{eq;08.9.15.5}
 \Bigl[\del_E-\lambdabar_0\theta,\,\,\,
 g\cdot\theta^{\dagger}\cdot g^{-1}
 \Bigr]
= g\cdot
 \Bigl[  g^{-1} \del_E g, \,\,
  \theta^{\dagger} \Bigr]\cdot g^{-1}
-\lambdabar_0\cdot
  g\cdot
\bigl[\theta,\theta^{\dagger}\bigr]\cdot
 g^{-1}
\end{equation}
Due to Proposition \ref{prop;07.11.22.25},
the second term in the right hand side of
(\ref{eq;08.9.15.5}) is bounded
with respect to $h$ and $g_{\poin}$
uniformly for $\lambda\in U(\lambda_0)$.
Let $F$ be as in Section
\ref{subsection;07.10.11.22}.
We have the following equality:
\[
 g^{-1}
 \del_E g
=g^{-1}
 \bigl[F,\,g\bigr]
+\sum_{\gminia,\alpha}
  \alphabar\cdot\frac{dz_1}{z_1}
 \pi^{(0)}_{\gminia,\alpha}
\]
Hence, we have the following estimate
with respect to $h$ and $g_{\poin}$:
\begin{equation}
\label{eq;08.9.15.10}
 \delta'_{\lambda_0}\bigl(
 g\cdot \theta^{\dagger}\cdot g^{-1}
 \bigr) 
=
g \bigl[ g^{-1}\cdot F\cdot g-F,\,\,
 \theta^{\dagger} \bigr]\cdot g^{-1}
 \\
+\sum_{\gminia,\alpha}
 g \cdot \alphabar\cdot \frac{dz_1}{z_1}
 \bigl[
 \pi^{(0)}_{\gminia,\alpha},\,\,
 \theta^{\dagger}
 \bigr]\cdot g^{-1}
+O(1)
\end{equation}
The second term of the right hand side of
(\ref{eq;08.9.15.10}) is bounded
uniformly for $\lambda$,
due to Lemma \ref{lem;08.9.15.20}
and Lemma \ref{lem;08.12.8.3}.
Let  us look at the first term in (\ref{eq;08.9.15.10}).
For the decomposition as in (\ref{eq;08.12.8.1}),
we obtain the following estimate
from Lemma \ref{lem;07.7.18.150}:
\[
 \bigl( g^{-1}\cdot F\cdot g-F\bigr)
_{(\gminia,\alpha),(\gminia',\alpha')}=
\left\{
 \begin{array}{ll}
 O\bigl(
 \exp(-\epsilon|z_1|^{\ord(\gminia-\gminia')})
 \bigr) & (\gminia\neq\gminia')\\
 O\bigl(
 |z_1|^{\epsilon}
 \bigr) & (\gminia=\gminia',\alpha\neq\alpha')
 \end{array}
\right.
\]
Then, we obtain the estimate for the first term
from Proposition \ref{prop;07.11.22.11}.
Thus, we obtain Lemma \ref{lem;08.9.15.11}.
\hfill\qed

\vspace{.1in}

Let 
$R\bigl(\nbigtzero d_{\lambda}'',h\bigr)$ 
be the curvature of
the unitary connection associated to
$\nbigtzero d_{\lambda}''$ and $h$.

\begin{lem}
\label{lem;08.9.15.12}
If $U(\lambda_0)$ is sufficiently small,
$R\bigl(\nbigtzero d_{\lambda}'',h\bigr)$ 
is bounded 
with respect to $h$ and $g_{\poin}$
uniformly for $\lambda\in U(\lambda_0)$.
\end{lem}
\pf
Let $\Lambda^{\dagger}(\lambda-\lambda_0)$
denote the adjoint of
$\Lambda(\lambda-\lambda_0)$
with respect to $h$.
Then, $R(\nbigtzero d''_{\lambda},h)$
can be rewritten as follows:
\begin{equation}
 \label{eq;08.9.15.1}
 \bigl[d_{\lambda_0}'',\delta_{\lambda_0}'\bigr]
+\delta_{\lambda_0}'\Lambda(\lambda-\lambda_0)
-d_{\lambda_0}''\Lambda^{\dagger}(\lambda-\lambda_0)
-\bigl[\Lambda(\lambda-\lambda_0),
 \Lambda^{\dagger}(\lambda-\lambda_0)\bigr]
\end{equation}
If $U(\lambda_0)$ is sufficiently small,
the first and fourth terms are bounded
uniformly for $\lambda\in U(\lambda_0)$
with respect to $h$ and $g_{\poin}$,
due to Corollary \ref{cor;07.11.22.30}
and Lemma \ref{lem;07.6.2.5}.
We have already obtained the estimate
for the second term in Lemma 
\ref{lem;08.9.15.11}.
Since the third term is adjoint to the second term
up to signature,
we obtain Lemma \ref{lem;08.9.15.12}.
\hfill\qed

\vspace{.1in}

The curvature of
$(\nbigtzero\nbige,h)$
is rewritten as follows:
\[
\bigl[\delbar_{\lambda}+ d_{\lambda_0}''+\Lambda,
 \del_{\lambda}+\delta_{\lambda_0}'-\Lambda^{\dagger} 
 \bigr]
=R(\nbigtzero d_{\lambda}'',h)
-\delbar_{\lambda}\Lambda^{\dagger}
+\del_{\lambda}\Lambda.
\]
We have only to show the uniform boundedness 
of $\del_{\lambda}\Lambda$
with respect to $h$ and $\gtilde_{\poin}$,
which can be rewritten as follows:
\begin{equation}
 \label{eq;07.11.22.52}
 \del_{\lambda}\Lambda
=\lambda_0\cdot g\cdot
 \bigl[g^{-1}\del_{\lambda}g,\,\,
 \theta^{\dagger}
 \bigr] \cdot g^{-1}
+
d\lambda\cdot g\cdot
 (\theta^{\dagger}-\Phibar)\cdot g^{-1}
+\lambda \cdot g\cdot
 \bigl[g^{-1} \del_{\lambda}g,
 \theta^{\dagger}-\Phibar\bigr] g^{-1}
\end{equation}
We have the following equality:
\[
 g^{-1}\del_{\lambda}g
=\sum_{\gminia,\alpha}
 \bigl(\gminiabar+\alphabar\log|z_1|^2\bigr)
 \cdot d\lambda\cdot\pi^{(0)}_{\gminia,\alpha}
\]
Then, it is easy to show that the right hand side 
of (\ref{eq;07.11.22.52}) is bounded
uniformly for $\lambda$ with respect to
$h$ and $\gtilde_{\poin}$,
by using Proposition \ref{prop;07.11.22.11}.
Thus, Proposition \ref{prop;07.11.22.45}
is proved.
\hfill\qed

\section{Locally uniform comparison
 of the irregular decompositions}
\label{subsection;08.9.15.41}
\subsection{Statements}
\label{subsection;08.12.9.1}

Let $U(\lambda_0)$ denote a neighbourhood
of $\lambda_0$ in $\cnum_{\lambda}$.
Let $\nbigxzero:=U(\lambda_0)\times X$
and $\nbigd^{(\lambda_0)}:=U(\lambda_0)\times D$.
We have two metrics of $\nbige_{|\nbigxzero-\nbigdzero}$.
One is $h$, and the other is
$\nbigpzero h$ given in {\rm(\ref{eq;07.7.21.4})}.

We put $W:=\nbigd^{(\lambda_0)}$
in the case $\lambda_0\neq 0$,
and $W:=\nbigd^{(0)}\cup\bigl(\{0\}\times X\bigr)$
in the case $\lambda_0=0$.
Let $\pi:\nbigxtilde^{(\lambda_0)}(W)
 \lrarr \nbigx^{(\lambda_0)}$
denote the real blow up along $W$.
Let $S$ be a small multi-sector
in $\nbigx^{(\lambda_0)}-W$.
Let $\Sbar$ denote 
the closure of $S$ in the real blow up
$\nbigxtilde^{(\lambda_0)}(W)$.
We have the Stokes filtrations $\nbigf^{(j)}$ of
$\nbigpzero\nbige_{|\Sbar}$ in the level $j$.
We will prove the estimate with respect to $h$
in the following proposition
in Section \ref{subsection;08.9.17.21}.

\begin{prop}
\label{prop;08.9.17.1}
Let $f$ be a $\DD$-flat section of 
$\End(\nbige)_{|
 \Sbar\setminus \pi^{-1}\nbigdzero}$
such that
$f\bigl(
 \nbigf^{(j)\,S}_{\gminia}\nbigpzero_0\nbige\bigr)
 \subset
 \nbigf^{(j)\,S}_{<\gminia}\nbigpzero_0\nbige$
for any $\gminia\in\Irr(\theta)$.
If we shrink $S$ in the radius direction, 
we have the following estimates
for some $\epsilon>0$:
\begin{equation}
 \label{eq;08.10.31.15}
 \bigl|f\bigr|_{\nbigpzero h}
=O\Bigl(
 \exp\bigl(-\epsilon|\lambda^{-1}z_1^{j}|\bigr)
 \Bigr),
\quad
 \bigl|f\bigr|_{h}
=O\Bigl(
 \exp\bigl(-\epsilon|\lambda^{-1}z_1^{j}|\bigr)
 \Bigr).
\end{equation}

Similarly,
let $f$ be a $\DD$-flat section of 
$\nbigf^{(j)\,S}_{<0}\nbigpzero\nbige_{|\Sbar
 \setminus \pi^{-1}(\nbigdzero)}$.
If we shrink $S$ in the radius direction, 
we have the following estimates
for some $\epsilon>0$:
\begin{equation}
\label{eq;08.10.31.20}
 \bigl|f\bigr|_{\nbigpzero h}
=O\Bigl(
 \exp\bigl(-\epsilon|\lambda^{-1}z_1^{j}|\bigr)
 \Bigr),
\quad
 \bigl|f\bigr|_{h}
=O\Bigl(
 \exp\bigl(-\epsilon|\lambda^{-1}z_1^{j}|\bigr)
 \Bigr).
\end{equation}
\end{prop}
Note that 
the estimate for $\nbigpzero h$ 
clearly follows from
Corollary \ref{cor;07.12.27.4}
and Corollary \ref{cor;08.9.15.13},
which we will use implicitly.

\vspace{.1in}

Let 
 $\nbigpzero_{a}\nbige_{|\Sbar}
=\bigoplus_{\gminia\in\Irr(\theta,j)}
 \nbigpzero_a\nbige^{(j)}_{\gminia,S}$
be any $\DD$-flat splitting 
of the filtration $\nbigf^{S\,(j)}$,
and let $p^{(j)}_{\gminia,S}$ denote 
the projection onto 
$\nbigpzero\nbige^{(j)}_{\gminia,S}$.
We will prove the following proposition
in Sections 
\ref{subsection;08.9.17.20}--\ref{subsection;08.9.17.22}.
\begin{prop}
\label{prop;08.9.17.2}
If we shrink $S$ in the radius direction,
we have the estimate
\[
 \pi^{(j)}_{\gminia}-p^{(j)}_{\gminia,S}
=O\bigl(|\lambda|\cdot\exp(-\epsilon|z_1^j|)\bigr)
\]
for some $\epsilon>0$ with respect to 
both the metrics $\nbigpzero h$ and $h$.
\end{prop}

We take a small multi-sectors $S_i$ $(i=1,\ldots,N)$
of $\nbigxzero-W$
such that the union of the interior points of
$S_i$ is $\nbigv-W$,
where $\nbigv$ denotes a neighbourhood
of $\nbigdzero$.
If $S_i$ are sufficiently small,
we have the partial Stokes filtrations
$\nbigf^{(j)\,S_i}$ of $\nbigpzero_a\nbige_{|\Sbar_i}$.
We can take $\DD$-flat splittings
$\nbigpzero_a\nbige_{|\Sbar_i}
=\bigoplus\nbigpzero_a\nbige^{(j)}_{\gminia,S_i}$.
Let $p^{(j)}_{\gminia,S_i}$ denote the projection
onto $\nbigpzero_a\nbige^{(j)}_{\gminia,S_i}$.
By gluing them in $C^{\infty}$
as in Section \ref{subsection;07.6.16.8},
we construct the $C^{\infty}$-section
$p^{(j)}_{\gminia,C^{\infty}}$
of $\End(\nbige)$
on $\nbigv\setminus W$,
which is extended to the $C^{\infty}$-section
on $\nbigv\setminus\nbigdzero$.
We will prove the following proposition
in Sections 
\ref{subsection;08.9.17.20}--\ref{subsection;08.9.17.22}.
\begin{prop}
\label{prop;08.9.17.3}
If we shrink $\nbigxzero$,
we have the estimate
\[
\pi^{(j)}_{\gminia}
-p^{(j)}_{\gminia,C^{\infty}}
=O\Bigl(|\lambda|\cdot
 \exp\bigl(-\epsilon|z_1^j|\bigr)\Bigr) 
\]
for some $\epsilon>0$
with respect to both 
the metrics $\nbigpzero h$ and $h$.
\end{prop}

Although 
the estimate with respect to $\nbigpzero h$
in Propositions \ref{prop;08.9.17.2} and 
\ref{prop;08.9.17.3} are equivalent,
we state them separately for our convenience
in the proof.

We have the following corollary
as a consequence of 
Proposition \ref{prop;08.9.17.1}.

\begin{cor}
Let $S$ be a small multi-sector in 
$\nbigxzero-W$,
and let 
\[
 \nbigpzero_{a}\nbige_{|\Sbar}
=\bigoplus_{\gminia\in\Irr(\theta,j)}
 \nbigpzero_a\nbige^{(j)}_{\gminia,S}
=\bigoplus_{\gminia\in\Irr(\theta,j)}
 \nbigpzero_a\nbige^{(j)\prime}_{\gminia,S}
\]
be two $\DD$-flat splittings
of the filtration $\nbigf^{S\,(j)}$.
Let $p^{(j)}_{\gminia,S}$ 
and $p^{(j)\prime}_{\gminia,S}$
denote the projections
onto $\nbigpzero\nbige_{\gminia,S}$
and $\nbigpzero\nbige_{\gminia,S}'$,
respectively.
Then, 
$p^{(j)\prime}_{\gminia,S}-p^{(j)}_{\gminia,S}
=O\bigl(\exp(-\epsilon|\lambda|^{-1}|z_1|^{j})\bigr)$
for some $\epsilon>0$
with respect to both the metrics
$\nbigpzero h$ and $h$.
\hfill\qed
\end{cor}

\subsection{Locally uniform comparison 
 of the twisted metrics}
\label{subsection;07.11.24.7}

Let $S$ be a small multi-sector in
$\nbigx^{(\lambda_0)}-W$.
We take a $\DD$-flat splitting
$\nbigpzero_{a}\nbige_{|\Sbar}
=\bigoplus_{\gminia\in\Irr(\theta)}
 \nbigpzero_{a}\nbige_{\gminia,S}$
of the full Stokes filtration $\nbigftilde^S$.
We have the induced splittings
$\nbigpzero_{a}\nbige_{|\Sbar}
=\bigoplus_{\gminib\in\Irr(\theta,j)}
 \nbigpzero_{a}\nbige^{(j)}_{\gminib,S}$
of the partial Stokes filtration $\nbigf^{S\,(j)}$.
We consider the following for $w\in\cnum$,
as in (\ref{eq;07.11.23.26}):
\begin{equation}
 \label{eq;08.1.20.30}
 F_S(w):=\exp\bigl(w\cdot\nbigb_S\bigr),
\quad
  \nbigb_S:=
 \sum_{m(0)\leq j\leq -1}\sum_{\gminib\in\Irr(\theta,j)}
 \zeta_j(\gminib)\cdot p^{(j)}_{\gminib,S}
\end{equation}
We also consider the following
as in (\ref{eq;07.11.24.1}):
\[
 F(w):=
 \exp\bigl(w\cdot \nbigb\bigr),
\quad
 \nbigb:=\sum_{\gminia\in\Irr(\theta)}
 \gminia\cdot\pi^{}_{\gminia}
=\sum_{m(0)\leq j\leq -1}
 \sum_{\gminib\in\Irr(\theta,j)}
 \zeta_j(\gminib)\cdot \pi^{(j)}_{\gminib}
\]

\begin{cor}
\label{cor;07.11.24.2}
Take $\eta>0$.
If we shrink $U(\lambda_0)$,
there exist some constants $\epsilon>0$ and $C>0$
such that the following estimate holds 
for any $|w|<\eta$
and for both the metrics
$h$ and $\nbigpzero_{\irr}h$:
\[
 F_S(w)\circ F(w)^{-1}-1
=O\Bigl(
 \exp\bigl(-\epsilon|z_1^{-1}|\bigr)
 \Bigr),
\quad
 F(w)\circ F_S(w)^{-1}-1
=O\Bigl(
 \exp\bigl(-\epsilon|z_1^{-1}|\bigr)
\Bigr)
\]
\end{cor}
\pf
It can be shown by using the argument
in the proof of Lemma \ref{lem;07.7.20.150},
together with Proposition \ref{prop;08.9.17.2}.
\hfill\qed

\begin{cor}
\label{cor;07.11.24.3}
If we shrink $U(\lambda_0)$, the following holds:
\begin{itemize}
\item
$F_S(\lambdabar-\lambdabar_0)^{\ast}h$
and $\nbigpzero_{\irr} h$ are mutually 
bounded, uniformly for $\lambda$.
\item
$F_S(-\lambdabar+\lambdabar_0)^{\ast}\nbigpzero_{\irr} h$
and $h$ are mutually bounded,
uniformly for $\lambda$.
\end{itemize}
\end{cor}
\pf
We have
$\nbigpzero_{\irr}h=
 F(\lambdabar-\lambdabar_0)^{\ast}h$
and
$h=F(-\lambdabar+\lambdabar_0)^{\ast}
 \nbigpzero_{\irr}h$.
Hence, the claim of Corollary \ref{cor;07.11.24.3}
can be shown using the same argument
as that in the proof of Lemma \ref{lem;07.7.20.200},
together with Corollary \ref{cor;07.11.24.2}.
\hfill\qed

\subsection{Preliminary}
\label{subsection;08.9.17.20}

We have the decomposition:
\begin{equation}
 \label{eq;07.11.23.80}
 \End(E)=\bigoplus_{\gminia,\gminib\in\Irr(\theta)}
\Hom(E_{\gminia},E_{\gminib})
\end{equation}
For any section $G$ of $\End(E)$,
we have the corresponding decomposition
$G=\sum G_{\gminia,\gminib}$.
We will use the following lemma implicitly,
which is obvious from the construction
of $\nbigpzero h$.
\begin{lem}
\label{lem;08.10.29.1}
Let $Y$ be a subset of $\nbigxzero-W$.
Let $f$ be a section of $\End(E)$ on $Y$.
Assume that the following holds for some $\epsilon>0$:
\begin{equation}
 \label{eq;08.9.17.10}
 |f_{\gminia,\gminib}|_{h}
=O\bigl(
 \exp\bigl(
 -\epsilon|\lambda^{-1}z_1^q|
 -\epsilon|z_1^{\ord(\gminia-\gminib)}|\bigr)
 \bigr)
\end{equation}
Then, if we shrink $U(\lambda_0)$,
the following holds for some $\epsilon'>0$:
\begin{equation}
\label{eq;08.9.17.11}
  |f_{\gminia,\gminib}|_{\nbigpzero h}
=O\Bigl(
 \exp\bigl(-\epsilon'|\lambda^{-1}z_1^q|
 -\epsilon'|z_1^{\ord(\gminia-\gminib)}|\bigr)
 \Bigr)
\end{equation}
Conversely,
if {\rm(\ref{eq;08.9.17.11})} holds
for some $\epsilon'>0$,
we obtain {\rm(\ref{eq;08.9.17.10})}
for some $\epsilon>0$
by shrinking $U(\lambda_0)$.
\hfill\qed
\end{lem}

\subsection{Comparison around $\lambda_0\neq 0$}

We use the notation in Section 
\ref{subsection;07.10.11.22}.
We have the following estimates
for some $\epsilon_1>0$
with respect to $h$ and $g_{\poin}$,
due to 
Theorem \ref{thm;07.10.4.1},
Theorem \ref{thm;07.10.4.3},
Lemma \ref{lem;07.7.18.150}
and 
Corollary \ref{cor;07.7.19.100}:
\[
\DD\pi^{(j)}_{\gminib}=
  O\bigl(\exp(-\epsilon_1|z_1|^{j})\bigr),
\quad
\nbigc^{(p)}\bigl(
 \DD\pi^{(j)}_{\gminib}
\bigr)= O\bigl(\exp(-\epsilon_1|z_1|^{p})\bigr),
\,\,\,(p<j)
\]
According to Lemma \ref{lem;08.10.29.1},
we may assume to
have the following estimate 
on $\nbigx^{(\lambda_0)}-\nbigd^{(\lambda_0)}$
with respect to $\nbigpzero h$:
\[
 \DD\pi^{(j)}_{\gminib}=
  O\bigl(\exp(-\epsilon_1|z_1|^{j}/2)\bigr)
\]

Let $\pi:
 \nbigxtilde^{(\lambda_0)}(\nbigd^{(\lambda_0)})
 \lrarr \nbigx^{(\lambda_0)}$ denote the real blow up
of $\nbigx^{(\lambda_0)}$ along $\nbigd^{(\lambda_0)}$.

\begin{lem}
\label{lem;07.12.22.30}
For any point $P\in\pi^{-1}(\nbigd^{(\lambda_0)})$,
there exist a multi-sector
$S\in\Multisector(P, 
 \nbigx^{(\lambda_0)}-\nbigd^{(\lambda_0)})$
and a $\DD$-flat splitting
$\nbigpzero_a\nbige_{|\Sbar}
=\bigoplus_{\gminia\in\Irr(\theta,j)} 
 \nbigpzero_a\nbige^{(j)\natural}_{\gminia,S}$
of the Stokes filtration $\nbigf^{S\,(j)}$ on $\Sbar$
in the level $j$
with the following property:
\begin{itemize}
\item
Let $p^{(j)\natural}_{\gminia,S}$ 
denote the projection
onto $\nbigpzero_a\nbige^{(j)\natural}_{\gminia,S}$
with respect to the decomposition.
Then,
\[
 \pi^{(j)}_{\gminia}-p^{(j)\natural}_{\gminia,S}
=O\bigl(\exp(-\epsilon_1|z_1|^{j}/10)\bigr)
\]
with respect to $\nbigpzero h$
\end{itemize}
\end{lem}
\pf
The argument is essentially the same as
the proof of Proposition \ref{prop;07.12.21.1}.
By using Lemma \ref{lem;07.7.19.102},
we can take $\nbigq^{(j)}_{\gminib,S}$
such that
$\DD\nbigq^{(j)}_{\gminib,S}
=\DD\pi^{(j)}_{\gminib}$
and 
$\nbigq^{(j)}_{\gminib,S}=O\bigl(
 \exp\bigl(-\epsilon_1|z_1|^{j}\bigr)\bigr)$
with respect to $\nbigpzero h$.
We put
$p^{(j)\natural}_{\gminib,S}:=
 \pi^{(j)}_{\gminib}-\nbigq^{(j)}_{\gminib,S}$,
which is $\DD$-flat.
By applying the modification as in the proof of
Lemma \ref{lem;07.7.7.20},
we may and will assume 
$\bigl[p^{(j)\natural}_{\gminib,S},\,
 p^{(j)\natural}_{\gminic,S}\bigr]=0$
and $p^{(j)\natural}_{\gminib,S}
 \circ p^{(j)\natural}_{\gminib,S}
 =p^{(j)\natural}_{\gminib,S}$.

We put
$\nbigf^{S\,(j)\natural}_{\gminib}
 :=\bigoplus_{\gminic\leq_S \gminib}
 \Image p^{(j)\natural}_{\gminic,S}$.
Let us compare 
the filtrations $\nbigf^{S\,(j)\natural}$
and $\nbigf^{S\,(j)}$.
Let $S^{\lambda}$ denote $S\cap\nbigx^{\lambda}$.
Let us consider the filtrations of
$\nbigelambda_{|S^{\lambda}}$
given as follows:
\[
 \overline{\nbigf}^{S^{\lambda}\,(j)}
 _{(1+|\lambda|^2)\gminia}
 \bigl(\nbigelambda_{|S^{\lambda}}
 \bigr):=
 \nbigf^{S^{\lambda}\,(j)}
 _{(1+\lambda\lambdabar_0)\gminia}
 \bigl(\nbige_{|\Sbar}\bigr)_{|S^{\lambda}},
\quad
 \overline{\nbigf}^{S^{\lambda}\,(j)\,\natural}
 _{(1+|\lambda|^2)\gminia}
 \bigl(\nbigelambda_{|S^{\lambda}}
 \bigr):=
 \nbigf^{S^{\lambda}\,(j)\,\natural}
 _{(1+\lambda\lambdabar_0)\gminia}
 \bigl(\nbige_{|\Sbar}\bigr)_{|S^{\lambda}}
\]
If $|\lambda-\lambda_0|$ is sufficiently small,
$\overline{\nbigf}^{S^{\lambda}\,(j)}
 (\nbigelambda_{|S^{\lambda}})$
is the same as 
the Stokes filtration 
$\nbigf^{S^{\lambda}\,(j)}\bigl(
 \nbigp\nbigelambda_{|S^{\lambda}}
 \bigr)$ of $\nbigp\nbigelambda_{|S^{\lambda}}$.
By construction,
$\overline{\nbigf}^{S^{\lambda}\,(j)\,\natural}
 \bigl(\nbigelambda_{|S^{\lambda}}\bigr)$
is the same as 
the filtration $\nbigf^{S^{\lambda}\,(j)\,\natural}$
of $\nbigp\nbigelambda_{|S^{\lambda}}$
in the proof of Proposition \ref{prop;07.12.21.1}.
Hence, we obtain that
the specializations of
$\nbigf^{S\,(j)\natural}$
and $\nbigf^{S\,(j)}$ to $S^{\lambda}$ are the same.
As a result,
we obtain that
the filtrations $\nbigf^{S\,(j)\natural}$
and $\nbigf^{S\,(j)}$ are the same.
Thus, we obtain Lemma \ref{lem;07.12.22.30}.
\hfill\qed

\vspace{.1in}
Let us show Proposition \ref{prop;08.9.17.2}
in the case $\lambda_0\neq 0$.
We take a finite covering
$S\subset \bigcup S_{P_i}$,
where $S_{P_i}$ are as in 
Lemma \ref{lem;07.12.22.30}.
We have 
$p^{(j)}_{\gminia,S}
-p^{(j)\natural}_{\gminia,S_{P_i}}
=O\bigl(
 \exp(-\epsilon|z_1|^j)
 \bigr)$ 
with respect to $\nbigpzero h$
on $S_{P_i}$ for some $\epsilon>0$
by using the estimate for $\nbigpzero h$
in Proposition \ref{prop;08.9.17.1}.
Then, we obtain the estimate of
$\bigl|
 p^{(j)}_{\gminia,S}-\pi^{(j)}_{\gminia}
 \bigr|_{\nbigpzero h}$
by Lemma \ref{lem;07.12.22.30}.

Let us show the estimate for $h$.
We may assume to have
$[p^{(j)}_{\gminia\,S},p^{(i)}_{\gminib\,S}]=0$.
We also have
$[\pi^{(j)}_{\gminia},
 \pi^{(i)}_{\gminib}]=0$.
Then, 
we obtain
\[
 [p^{(j)}_{\gminia\,S},\pi^{(i)}_{\gminib}]
=O\Bigl(
 \exp\bigl(
 -\epsilon_1|z_1|^{i}-\epsilon_1|z_1|^j
 \bigr)\Bigr)
\]
with respect to $\nbigpzero h$
from the previous estimate.
(See the proof of Lemma \ref{lem;08.9.15.20},
for example.)
By the same argument as the proof of 
Lemma \ref{lem;07.11.22.10},
we obtain the following:
\begin{equation}
\label{eq;07.12.22.50}
\bigl|
 \bigl(
 p^{(j)}_{\gminib,S}-\pi^{(j)}_{\gminib}
\bigr)_{\gminia,\gminia'}
\bigr|_{\nbigpzero h}
=O\Bigl(
\exp\bigl(
 -\epsilon_2|z_1|^j
 -\epsilon_2|z_1|^{\ord(\gminia-\gminia')}
 \bigr)
\Bigr)
\end{equation}
By Lemma \ref{lem;08.10.29.1},
we obtain the desired
estimate for
$\pi^{(j)}_{\gminia}-p^{(j)}_{\gminia,S}
=O\bigl(\exp(-\epsilon|z_1|^{j})\bigr)$
with respect to $h$.
Thus, the proof of Proposition 
\ref{prop;08.9.17.2} in the case $\lambda_0\neq 0$
is finished.
Proposition \ref{prop;08.9.17.3}
in the case $\lambda_0\neq 0$
immediately follows.

\subsection{Comparison around $\lambda_0=0$}
\label{subsection;08.9.17.22}

Let us show Proposition \ref{prop;08.9.17.3}
in the case $\lambda_0=0$.
Let $S$ be a small multi-sector in $\nbigx^{(0)}-W$.
Let $p_{\gminia,S}^{(j)}$ and 
$p_{\gminia,S}^{(j)\prime}$
come from $\DD$-flat splittings 
of the filtration $\nbigf^{S(j)}$.
If we shrink $S$ in the radius direction,
we have the following estimate
for some $\epsilon_1>0$
with respect to the metric $\nbigp^{(0)} h$:
\begin{equation}
\label{eq;07.12.22.55}
p^{(j)}_{\gminia,S}-p^{(j)\prime}_{\gminia,S}=
 O\Bigl(
 \exp\bigl(-\epsilon_1|\lambda^{-1}z_1^{j}|\bigr)
 \Bigr)
\end{equation}
Let $p^{(j)}_{\gminia,C^{\infty}}$
be as in Subsection \ref{subsection;08.12.9.1}.
By shrinking $\nbigv$,
we obtain the following estimate on 
$\nbigv\setminus \nbigd^{(0)}$
with respect to $\nbigp^{(0)}h$ 
by using (\ref{eq;07.12.22.55})
and $\delbar_{\lambda}\pi^{(j)}_{\gminia}=0$:
\begin{equation}
\label{eq;07.12.22.70}
g(\lambda,\vecz)\cdot d\lambdabar:=
\delbar_{\lambda}
\bigl(
 p^{(j)}_{\gminia,C^{\infty}}
-\pi^{(j)}_{\gminia}
\bigr)
=
O\Bigl(
 \exp\bigl(-\epsilon_2|\lambda^{-1}z_1^{j}|\bigr)
 \cdot d\lambdabar
\Bigr)
\end{equation}

Take small $0<\delta_1<\delta_2$.
Let us estimate
$\pi^{(j)}_{\gminia}
-p^{(j)}_{\gminia,C^{\infty}}$
on $\nbigz:=\{\lambda\,|\,\delta_1\leq|\lambda|\leq \delta_2\}
\times (X-D)$.
\begin{lem}
\label{lem;07.12.22.80}
Let $S$ be a small multi-sector 
in $\nbigz-(\nbigz\cap\nbigd^{(0)})$.
Let $\nbigp^{(0)}\nbige_{|\Sbar}
=\bigoplus \nbigp^{(0)}\nbige_{\gminia,S}^{(j)}$
be a $\DD$-flat splitting of
the Stokes filtration $\nbigf^{S\,(j)}$,
and let $p^{(j)}_{\gminia,S}$ denote the projection
onto $\nbigp^{(0)}\nbige_{\gminia,S}^{(j)}$.
Then, 
$\pi^{(j)}_{\gminia}-
 p^{(j)}_{\gminia,S}
=O\bigl(\exp(-\epsilon_3|z_1|^j)\bigr)$
with respect to $\nbigp^{(0)}h$
for some $\epsilon_3>0$.
As a result,
$\pi^{(j)}_{\gminia}
-p^{(j)}_{\gminia,C^{\infty}}
=O\bigl(
 \exp(-\epsilon_4|z_1^j|) \bigr)$
for some $\epsilon_4>0$
with respect to $\nbigp^{(0)}h$
on $\nbigz-(\nbigz\cap \nbigd^{(0)})$.
\end{lem}
\pf
It can be shown by the arguments
in the proof of Lemma \ref{lem;07.12.22.30}
and Proposition \ref{prop;08.9.17.2}
in the case $\lambda_0\neq 0$.
We remark that
the Stokes filtrations of
$\nbigp^{(0)}\nbigelambda$
and $\nbigp\nbigelambda$
are essentially the same,
i.e.,
$\nbigf^{S^{\lambda}}_{\gminia}
 \bigl(\nbigp^{(0)}\nbigelambda_{|S^{\lambda}}\bigr)
=\nbigf^{S^{\lambda}}_{(1+|\lambda|^2)\gminia}
 \bigl(\nbigp\nbigelambda_{|S^{\lambda}}\bigr)$.
\hfill\qed

\vspace{.1in}

Let us continue the proof of 
Proposition \ref{prop;08.9.17.3}.
Let $\vecv$ be a holomorphic frame of
$\nbigp^{(0)}_{0}\End\nbige$ on $\nbigx^{(0)}$.
We have the expression
$g=\sum g_i\cdot v_i$.
We have the estimate
$|g_i|\leq
 C\cdot 
 \exp\bigl(-\epsilon_5|\lambda^{-1}z_1^j|\bigr)$.
Let $\chi$ be a non-negative
valued $C^{\infty}$-function
on $\cnum_{\lambda}$
such that $\chi(\lambda)=1$ for 
$|\lambda|\leq(\delta_1+\delta_2)/2$
and $\chi(\lambda)=0$ for $|\lambda|\geq \delta_2$.
We consider the following:
\[
 G_i(\vecz,\lambda):=
 \lambda\cdot
 \int\frac{ g_i(\vecz,\mu)\cdot\mu^{-1}\cdot \chi(\mu)}
 {\mu-\lambda}
 \cdot 
\frac{\sqrt{-1}}{2\pi}\cdot
 d\mu\cdot d\mubar
\]
We put $G:=\sum G_i\cdot v_i$.
Then, we have $\delbar_{\lambda}G=g$
on $\bigl\{|\lambda|<(\delta_1+\delta_2)/2\bigr\}
 \times (X-D)$.
We also have the following estimate:
\[
 \bigl|G(\vecz,\lambda)\bigr|_{\nbigp^{(0)}h}\leq
 C\cdot 
 |\lambda|\cdot
 \exp\bigl(-\epsilon_6\cdot|z_1^j|\bigr)
\]
Let us look at
$H:=\pi^{(j)}_{\gminia}
 -p^{(j)}_{\gminia,C^{\infty}}-G$.
We have
$\delbar_{\lambda}H=0$ and $H(\vecz,0)=0$.
By a lemma of Schwarz 
and Lemma \ref{lem;07.12.22.80},
we obtain
$ \bigl|H(\vecz,\lambda)\bigr|_{\nbigp^{(0)}h}
= O\Bigl(
|\lambda|\cdot
 \exp\bigl(-\epsilon_7\cdot |z_1|^j\bigr) \Bigr)$.
Then, we obtain 
\begin{equation}
 \label{eq;08.12.9.2}
 \bigl|
 \pi^{(j)}_{\gminia}
-p^{(j)}_{\gminia,C^{\infty}}
 \bigr|_{\nbigpzero h}
=O\Bigl(
|\lambda|\cdot
 \exp\bigl(-\epsilon_8\cdot |z_1|^j\bigr) \Bigr),
\end{equation}
i.e., the estimate with respect to
$\nbigp^{(0)}h$ in Proposition \ref{prop;08.9.17.3}.

Let us show the estimate with respect to $h$.
By construction, 
we have the following estimate:
\begin{equation}
\label{eq;07.11.23.131}
 \Bigl| \bigl[
 p^{(j)}_{\gminia,C^{\infty}},
 p^{(i)}_{\gminib,C^{\infty}}
 \bigr] \Bigr|_{\nbigpzero h}
\leq
 C_{21}\cdot\exp\Bigl(
 \bigl(
 -\epsilon_{21}|\lambda^{-1}z_1^i|-\epsilon_{21}
 |\lambda^{-1}z_1^j|\bigr)
 \Bigr)
\end{equation}
We also have
$[\pi^{(j)}_{\gminia},
 \pi^{(i)}_{\gminib}]=0$.
From (\ref{eq;08.12.9.2})
and (\ref{eq;07.11.23.131}),
we obtain the following estimate:
\begin{equation}
 \label{eq;07.11.23.132}
\Bigl|\bigl[
 p^{(j)}_{\gminia,C^{\infty}},
\pi^{(i)}_{\gminib}
\bigr]\Bigr|_{\nbigpzero h}
\leq
C_{22}\cdot |\lambda|\cdot \exp\bigl(
 -\epsilon_{22}|z_1|^j
 -\epsilon_{22}|z_1|^i
 \bigr)
\end{equation}
We have the decomposition
$ p^{(j)}_{\gminia,C^{\infty}}
-\pi^{(j)}_{\gminia}
=\sum_{\gminib,\gminib'}
\bigl(
p^{(j)}_{\gminia,C^{\infty}}
-\pi^{(j)}_{\gminia}\bigr)_{\gminib,\gminib'}$
corresponding to (\ref{eq;07.11.23.80}).
By the same argument as that in the proof of 
Lemma \ref{lem;07.11.22.10},
we obtain the following estimate
from (\ref{eq;07.11.23.132}):
\begin{equation}
 \label{eq;07.11.23.133}
\Bigl| \bigl(p^{(j)}_{\gminib,C^{\infty}}
-\pi^{(j)}_{\gminib}\bigr)_{\gminia,\gminia'}
\Bigr|_{\nbigpzero h}
\leq
 C_{23}\cdot|\lambda|\cdot
 \exp\bigl(-\epsilon_{23}|z_1|^j
 -\epsilon_{23}|z_1|^{\ord(\gminia-\gminia')}
 \bigr)
\end{equation}
Then the estimate with respect to $h$
in Proposition \ref{prop;08.9.17.3} follows
from Lemma \ref{lem;08.10.29.1}.

\vspace{.1in}
The estimate with respect to $\nbigp^{(0)} h$
in Proposition \ref{prop;08.9.17.2}
immediately follows from Proposition
\ref{prop;08.9.17.3}.
The estimate with respect to $h$
follows from the estimate with respect to
$\nbigp^{(0)}h$  as above.

\subsection{Proof of Proposition \ref{prop;08.9.17.1}}
\label{subsection;08.9.17.21}

Let us consider the first claim.
We take a $\DD$-flat splitting
$ \nbigpzero\nbige_{|\Sbar}
=\bigoplus_{\gminia\in\Irr(\theta)}
 \nbigpzero\nbige_{\gminia,S}$
of the full Stokes filtration $\nbigftilde^{S}$.
We have only to consider the case
$f\in 
 \Hom\bigl(\nbigpzero\nbige_{\gminia,S},
 \nbigpzero\nbige_{\gminib,S}
 \bigr)$
for some $\gminia>_S\gminib$.
We have the following:
\[
 \bigl|f\bigr|_{\nbigpzero h}
=O\Bigl(
 \exp\bigl(-\epsilon|\lambda^{-1}
 z_1^{\ord(\gminia-\gminib)}|\bigr)
 \Bigr)
\]
We would like to obtain the estimate
with respect to $h$,
by using Lemma \ref{lem;08.10.29.1}.
For any $l$ and $\gminic\in\Irr(\theta,l)$
such that $\gminic\neq \eta_l(\gminia)$,
we have the following estimate
with respect to $\nbigpzero h$:
\begin{equation}
 \label{eq;08.9.16.1000}
 f\circ\pi^{(l)}_{\gminic}
=f\circ (\pi^{(l)}_{\gminic}-p^{(l)}_{\gminic})
=O\Bigl(
 \exp\bigl(-\epsilon|z_1^l|
 -\epsilon|\lambda^{-1}
 z_1^{\ord(\gminia-\gminib)}|
 \bigr)
 \Bigr)
\end{equation}
We have similar estimate for
$\pi^{(l)}_{\gminid}\circ f$
for $\gminid\neq\eta_l(\gminib)$.
Then, we obtain the desired estimate
(\ref{eq;08.10.31.15})
by using Lemma \ref{lem;08.10.29.1}.

Let us consider the second claim.
We have only to consider the case
$f\in\nbigpzero\nbige_{\gminia,S}$
for some $\gminia<_S0$.
We have the decomposition
$f=\sum f_{\gminib}$
corresponding to $E=\bigoplus E_{\gminib}$.
We have the following estimate for 
$\gminic\neq\eta_l(\gminia)$
with respect to $\nbigpzero h$:
\[
 f_{\gminic}
=\pi_{\gminic}^{(l)}(f)
=\bigl(\pi_{\gminic}^{(l)}-p_{\gminic}^{(l)}\bigr)(f)
=O\Bigl(
 \exp\bigl(-\epsilon|z_1^l|
 -\epsilon|\lambda^{-1}z_1^{\ord(\gminia)}|
 \bigr)
 \Bigr)
\]
Then, it is easy to obtain (\ref{eq;08.10.31.20}).
\hfill\qed

\section{Norm estimate for a fixed $\lambda$}
\label{subsection;08.9.15.42}
\subsection{Norm estimate for $\nbigp\nbigelambda$}
\label{subsection;07.11.23.35}

Let $F$ denote the parabolic filtration of 
$\nbigp_a\nbigelambda_{|D}$.
The residue $\Res(\DDlambda)_{|Q}$ $(Q\in D)$
induces the endomorphism of
$\Gr^F_{a}(\nbigp\nbigelambda_{|Q})$.
Let $N^{\lambda}_{a|Q}$ denote the nilpotent part.

\begin{lem}
\label{lem;07.11.23.21}
The conjugacy classes of $N^{\lambda}_{a|Q}$
are independent of the choice of $Q\in D$.
In particular,
the weight filtration $W$ 
gives the filtration in the category of vector bundles.
\end{lem}
\pf
Let us consider the case
in which $\lambda$ is generic
(Definition \ref{df;07.11.23.5}).
We have the irregular decomposition
\[
 (\nbigp_a\nbigelambda,\DDlambda)_{|\Dhat}
=\bigoplus
 \bigl(\nbigp_a\nbigelambdahat_{\gminia},
 \DDlambdahat_{\gminia}\bigr),
\]
which is compatible with the parabolic filtrations $F$.
We have the flat $\lambda$-connection
$\DDlambda_{\gminia,a}$
on $\Gr^F_a\bigl(
 \nbigp\nbigelambdahat_{\gminia}
 \bigr)$ given by
$\DDlambda_{\gminia,a}(f)
=p_a\DDlambdahat_{\gminia}(F)$,
where $p_a$ is the naturally induced map
$\nbigp_a\nbigelambdahat_{\gminia}
\lrarr \Gr_a^F(\nbigp\nbigelambdahat_{\gminia})$,
and $F$ is a section of
$\nbigp_a\nbigelambdahat_{\gminia}$
such that $p_a(F)=f$.
Since the residue is flat with respect to
$\DDlambda_{\gminia,a}$,
the conjugacy classes are independent of $Q\in D$.

Let us consider the general case.
Let $\pi:X-D\lrarr D$ denote the projection.
Let $(E_Q,\delbar_{E_Q},\theta_Q,h_Q)$ denote
the restriction of $\harmonicbundle$ to
$\pi^{-1}(Q)$.
We have the $\lambda$-connections
$(\nbigp\nbigelambda_Q,\DDlambda_Q)$.
We have the correspondence 
of the conjugacy classes
of the nilpotent part of the residues
of $\Res(\DDlambda_Q)$ and 
$\Res(\DD^{\lambda_1}_Q)$
for any $\lambda_1$
in Proposition \ref{prop;07.7.19.31}.
Hence, the claim holds in the general case.
\hfill\qed

\vspace{.1in}

Let $\vecv$ be a holomorphic frame of
$\nbigp_c\nbigelambda$
such that 
(i) it is compatible with the parabolic filtration,
(ii) the induced frame on $\Gr^F(\nbigp_c\nbigelambda)$
 is compatible with the weight filtration $W$.
We put $a(v_i):=\deg^F(v_i)$ and $k(v_i):=\deg^W(v_i)$.
Let $h_0$ be the $C^{\infty}$-metric of $E$
given by
$h_0(v_i,v_j):=
 \delta_{i,j}\cdot 
 |z_1|^{-2a(v_i)}\cdot (-\log|z_1|^2)^{k(v_i)}$.

\begin{prop}
\label{prop;07.11.22.55}
The metrics $h$ and $h_0$ are mutually bounded.
\end{prop}
\pf
We give only an outline.
We have the irregular decomposition:
\begin{equation}
\label{eq;07.11.22.56}
 (\nbigp_{\ast}\nbigelambda,\DDlambda)_{|\Dhat}
=\bigoplus_{\gminia\in\Irr(\theta)}
 (\nbigp_{\ast}\nbigehatlambda_{\gminia},\DDlambda_{\gminia})
\end{equation}
Here, $\DDlambda_{\gminia}-(1+|\lambda|^2)d\gminia$
are logarithmic with respect to
 $\nbigp_{\ast}\nbigelambdahat_{\gminia}$.
By taking Gr
associated to the generalized eigen decomposition $\EE$ 
and the parabolic filtration $F$,
we obtain vector spaces
\[
 \Gr_{a,\alpha}^{F,\EE}\bigl(
 \nbigp\nbigelambdahat_{\gminia|O}
\bigr).
\]
For each $(a,\alpha,\gminia)$,
we have the model bundle on $\Delta^{\ast}$ 
obtained from the vector space
$\Gr^{F,\EE}_{a,\alpha}
 (\nbigp\nbigelambda_{\gminia|O})$
with the nilpotent part
$N_{a,\alpha\gminia|O}$ of the residue
(Section \ref{subsection;10.5.24.3}),
which is denoted by
$(E^1_{a,\alpha,\gminia},
 \delbar_{E^1_{a,\alpha,\gminia}},
 h^1_{a,\alpha,\gminia},
 \theta^1_{a,\alpha,\gminia})$.
Let $u=(b,\beta)\in\real\times\cnum$
 be determined by
$\kmsmap(\lambda,u)=(a,\alpha)$.
We have the rank one harmonic bundle
$L(u,\gminia)=\bigl(
 \nbigo\cdot e,\theta^2_{u,\gminia},
 h^2_{u,\gminia}\bigr)$ on $X-D$ given as follows:
\[
\theta^2_{u,\gminia}=d\gminia
+\beta\cdot\frac{dz_1}{z_1},
\quad
h^2_{u,\gminia}(e,e)=|z_1|^{-2b}
\]
Let $q_1:X-D\lrarr \Delta^{\ast}$ be the projection
$q_1(z_1,\ldots,z_n)=z_1$.
We set
\[
\bigl(E_{a,\alpha,\gminia},
 \delbar_{E_{a,\alpha,\gminia}},
 \theta_{a,\alpha,\gminia},
 h_{a,\alpha,\gminia}\bigr):=
q_1^{\ast}\bigl(E^1_{a,\alpha,\gminia},
 \delbar_{E^1_{a,\alpha,\gminia}},
 \theta^1_{a,\alpha,\gminia},
 h^1_{a,\alpha,\gminia}
 \bigr)
\otimes
L(u,\gminia).
\]
Let $(E_1,\delbar_{E_1},h_1,\theta_1)$
denote their direct sum.

For a large $N$,
let $\Dhat^{(N)}$ denote
the $N$-th infinitesimal neighbourhood
of $D$.
We take a holomorphic decomposition
\begin{equation}
 \label{eq;08.12.9.5}
 \nbigp_c\nbigelambda
=\bigoplus_{\gminia\in\Irr(\theta)}
 \nbigp_c\nbigelambda_{\gminia,N},
\end{equation}
whose restriction to $\Dhat^{(N)}$
is the same as the restriction of
(\ref{eq;07.11.22.56}).
We have the natural decomposition
\begin{equation}
\label{eq;08.12.9.6}
 \nbigp_c\nbigelambda_1
=\bigoplus_{\gminia\in\Irr(\theta)}
 \nbigp_c\nbigelambda_{1,\gminia},
\end{equation}
where $\nbigp_c\nbigelambda_{1,\gminia}$
are induced by
$\bigoplus_{(a,\alpha)}
\bigl(E_{a,\alpha,\gminia},
 \delbar_{E_{a,\alpha,\gminia}},
 \theta_{a,\alpha,\gminia},
 h_{a,\alpha,\gminia}\bigr)$.

We can take a holomorphic isomorphism
$\Phi:\nbigp_c\nbigelambda_1\lrarr\nbigp_c\nbigelambda$
such that 
(i) it preserves the decompositions
(\ref{eq;08.12.9.5})
and (\ref{eq;08.12.9.6}),
(ii) it preserves the parabolic filtration,
(iii) $\Gr^F(\Phi_{|D})$ is compatible with
 the residues.
It is easy to show that
$\Phi(h_1)$ is mutually bounded with $h_0$.

We have the harmonic bundle
$E_2:=\Hom(E_1,E)$
with the induced Higgs field $\theta_2$
and the induced pluri-harmonic metric $h_2$.
We can regard $\Phi$ 
as a section of $\nbigp_0\nbigelambda_2$.
For each $Q\in D$,
we have $g_Q:=\DDlambda_Q(\Phi_{|\pi^{-1}(Q)})\in
 \nbigp_{-\epsilon+1}{\nbigelambda_{2|Q}}$
for some $\epsilon>0$.
Therefore, it is $L^p$ for some $p>2$,
and the $L^p$-norm is bounded uniformly for $Q$.
We have
$\Delta\bigl|\Phi_{|\pi_1^{-1}(Q)}\bigr|^2_{h_2}
 \leq |g_Q|^2$
on $\pi^{-1}(Q)\setminus \{Q\}$,
due to Corollary \ref{cor;07.11.23.10}.
We have already known 
$\bigl|\Phi_{|\pi_1^{-1}(Q)}\bigr|_{h_2}$ is bounded
(Proposition \ref{prop;07.10.6.40}).
Hence, the inequality holds on $\pi_1^{-1}(Q)$
as distributions.
Therefore, $\Phi_{|\pi^{-1}(Q)}$ can be estimated
by the values on $\del\pi^{-1}(Q)$
and the $L^{p}$-norm of $|g_Q|$.
Thus, we obtain the claim of Proposition \ref{prop;07.11.22.55}.
\hfill\qed

\subsection{Estimate of the connection form 
(Appendix)}
\label{subsection;08.9.29.9}

We have a complement on the estimate
of the connection form in 
Lemma \ref{lem;07.6.2.20}
and Lemma \ref{lem;07.7.18.150}.
Since we will not use it below,
the reader can skip this subsection.
Let $\vecv$ be a frame of $\prolong{E}$
compatible with the decomposition 
$E=\bigoplus_{
 (\gminia,\alpha)\in\Irr(\theta)\times\Sp(\theta)}
 E_{\gminia,\alpha}$,
the parabolic filtration,
and the weight filtration $W$.
Let $F\in \End(E)\otimes\Omega^{1,0}$ be determined
by $F(\vecv)=\del\vecv$.
Let $F_0\in\End(E)\otimes\Omega^{1,0}$ be determined by
$F_0(v_i)=-a(v_i)\cdot v_i\cdot dz_1/z_1$.

\begin{prop}
\label{prop;08.9.15.44}
We have the following estimate:
\[
F_1=F_0
+O\left(\frac{\log(-\log|z_1|)}{-\log|z_1|}\right)
   \frac{dz_1}{z_1}
+\sum_{j=2}^nO(1)\cdot dz_j
\]
\end{prop}
\pf
We give only an outline.
Let $h_0$ be the metric as in 
Subsection \ref{subsection;07.11.23.35},
and let $s$ be determined by $h=h_0\cdot s$.
The connection form for $h_0$ satisfies the above estimate.
We have
$\del_h=\del_{h_0}+s^{-1}\del_{h_0}s$.
We have only to estimate $s^{-1}\del_{h_0}s$.
By the estimate of the curvatures $R(h)$ and $R(h_0)$,
we obtain the following estimate with respect to $h_0$:
\[
 \delbar(s^{-1}\del_{h_0}s)
=O\left(\frac{dz\cdot d\zbar}
 {|z|^2\cdot\bigl(-\log|z|^2\bigr)^2}
 \right)
\]
Let $\vecv^{\lor}=(v_i^{\lor})$ denote
the dual frame of $\vecv$.
We have the expressions
$s^{-1}\del_{h_0}s=
 \sum A_{i,j}^{(p)}\cdot 
 v_i\cdot v_j^{\lor}\cdot dz_p$
and
$\delbar(s^{-1}\del_{h_0}s)=
 \sum B_{i,j}^{(p,q)}v_i\cdot v_j^{\lor}
 \cdot dz_p\cdot d\zbar_q$.
We have 
$\del A_{i,j}^{(p)}/\del\zbar_p=B_{i,j}^{(p,p)}$.
Let $\del_{h_0,p}$ denote the restriction of
$\del_{h_0}$ to the $z_p$-direction.
Let $\pi_p:X-D\lrarr D_p$ denote the projection.

By using Lemma \ref{lem;07.2.6.11},
we obtain the following expression:
\begin{equation}
 \label{eq;07.1.5.1}
 s^{-1}\del_{h_0,1}s=
 O\left(\frac{\log(-\log|z|)}{-\log|z|}\right)
 \frac{dz_1}{z_1}
+\Phi\cdot dz_1
\end{equation}
Here, $\delbar_1\Phi=0$.
We also have the finiteness
$\int_{\pi_1^{-1}(Q)}
 |s^{-1}\del_{h_0,1}s|_{h_0}^2<\infty$
for any $Q\in D_1$.
(See Lemma 5.5 of \cite{mochi4}, for example.)
Hence, $\Phi_{|\pi_1^{-1}(Q)}$ is also $L^2$,
with respect to $h_0$.
Then, it is easy to obtain the desired estimate
for the $dz_1$-component
by using the maximum principle
for holomorphic functions.

Let us look at the $dz_p$-component $(p>1)$.
Let $\chi(z_p)$ be a non-negative $C^{\infty}$-function
on $\Delta^{\ast}$ such that
$\chi(z_p)=1$ for $|z_p|\leq 1/2$
and $\chi(z_p)=0$ for $|z_p|\geq 2/3$.
We can show 
$\int_{\pi_p^{-1}(Q)}
 \bigl(s^{-1}\del_{h_0,p}(\chi s),
 \del_{h_0,p}(\chi s)\bigr)$
is estimated uniformly for $Q\in D_p\setminus D_1$.
(See the argument in
 the proof of Lemma \ref{lem;07.6.2.20}.)
We also have the uniform estimate
of the sup norm of
$\delbar_p(s^{-1}\del_{h_0,p} s)_{|\pi_p^{-1}(Q)}$.
Hence, we obtain the uniform estimate of
the sup norm of $s^{-1}\del_p s_{|\pi_p^{-1}(Q)}$.
\hfill\qed

\section{Norm estimate in family}
\label{subsection;08.9.15.43}
\subsection{Statement}

We have the good family of 
filtered $\lambda$-flat bundles
$\nbigp^{(\lambda_0)}_{\ast}\nbige$.
We have the filtration $\Fzero$
and the decomposition
$\EEzero$ of 
$\nbigp^{(\lambda_0)}_a
 \nbige_{|U(\lambda_0)\times D}$.
(See Remark \ref{rem;08.9.4.10}
for the decomposition $\EEzero$.)
\begin{lem}
The conjugacy class of the nilpotent part of
$\Res(\DD)_{|(\lambda,P)}$
is independent of the choice of 
$(\lambda,P)\in U(\lambda_0)\times D$.
\end{lem}
\pf
It follows from Corollary \ref{cor;07.11.23.20}
and Lemma \ref{lem;07.11.23.21}.
\hfill\qed

\vspace{.1in}

Let $\vecv$ be a frame of 
$\nbigp^{(\lambda_0)}_a \nbige$
compatible with $\bigl(\EEzero,\Fzero,W\bigr)$.
Let $u(v_i)\in\real\times\cnum$ be determined by
$\kmsmap(\lambda_0,u(v_i))
=\bigl(\deg^{\Fzero}(v_i),\deg^{\EEzero}(v_i)\bigr)$.
We put $k(v_i):=\deg^{W}(v_i)$.
Let $h_0$ be the metric given as follows:
\[
 h_0(v_i,v_j):=
 \delta_{i,j}\cdot |z_1|^{-2\paramap(\lambda,u(v_i))} 
\cdot\bigl(-\log|z_1|\bigr)^{k(v_i)}
\]

We will prove the following proposition in Sections
\ref{subsection;07.11.24.5}--\ref{subsection;07.11.24.6}.
\begin{prop}
\label{prop;07.10.22.30}
If we shrink $U(\lambda_0)$,
$\nbigp_{\irr}^{(\lambda_0)}h$ and $h_0$
are mutually bounded
on $U(\lambda_0)\times (X-D)$.
(See Section {\rm\ref{subsection;07.11.23.1}}
for the hermitian metric
 $\nbigpzero_{\irr}h$.)
\end{prop}

Before going into the proof of the proposition,
we give a corollary.
Let $S$ and $F_S$ be 
as in Section \ref{subsection;07.11.24.7}.
\begin{cor}
\label{cor;07.12.2.30}
If we shrink $U(\lambda_0)$,
the metrics $h$ and 
$F_S(-\lambdabar+\lambdabar_0)^{\ast}h_0$
are mutually bounded on $S$.
\end{cor}
\pf
It follows from Proposition \ref{prop;07.10.22.30}
and Corollary \ref{cor;07.11.24.3}.
\hfill\qed

\subsection{Preliminary comparison}

As a preparation,
we give a weak comparison
up to small polynomial order.

\begin{lem}
\label{lem;08.1.8.20}
The metrics
$\nbigpzero_{\irr} h$
and $h_0$ are uniformly mutually bounded
up to $|z_1|^{-\eta|\lambda-\lambda_0|} (-\log|z_1|)^N$-order
for some $\eta>0$.
\end{lem}
\pf
Let $a(v_i):=\deg^{\Fzero}(v_i)$,
and let $h_0'$ be the hermitian metric
given by 
\[
 h_0'(v_i,v_j):=\delta_{i,j}
 \cdot |z_1|^{-2a(v_i)}.
\]
We have already known that
$\nbigpzero h$ and $h_0'$
are uniformly mutually bounded 
up to log order because of the acceptability of
$\nbigp^{(\lambda_0)}h$
(Proposition \ref{prop;07.11.22.45}).
Then, the claim of Lemma 
\ref{lem;08.1.8.20} follows.
\hfill\qed

\subsection{Comparison map for
the associated graded bundle}
\label{subsection;07.11.24.5}

We have the decomposition:
\[
 (\nbigpzero_{\ast}\nbige,\DD)_{|U(\lambda_0)\times \Dhat}
=\bigoplus _{\gminia\in\Irr(\theta)}
 \bigl(
 \nbigpzero_{\ast}\nbigehat_{\gminia},\DDhat_{\gminia}
\bigr)_{|U(\lambda_0)\times \Dhat}
\]
Here, $\DDhat_{\gminia}-(1+\lambda\lambdabar_0)\cdot d\gminia$
are logarithmic with respect to
$\nbigpzero_{\ast}\nbigehat_{\gminia}$.
Take any point $P\in D$.
We have the vector space
$V_{\gminia,u}:=
 \Gr^{\Fzero,\EEzero}_{\kmsmap(\lambda_0,u)}\bigl(
 \nbigpzero\nbigehat_{\gminia}\bigr)_{|(\lambda_0,P)}$.
We have the endomorphism $N_{\gminia,u}$,
which is the nilpotent part of
$\Res(\DD)_{|(\lambda_0,P)}$.
Let $\bigl(
E'_{1,\gminia,u},\delbar'_{1,\gminia,u},\theta'_{1,\gminia,u},
 h'_{1,\gminia,u} \bigr)$
be the model bundle
associated to $(V_{\gminia,u},N_{\gminia,u})$
(Section \ref{subsection;10.5.24.3}).
Let $L(u,\gminia)$ be the harmonic bundle of rank one
as in the proof of Proposition \ref{prop;07.11.22.55}.
We put 
\[
 (E_{1,\gminia},\delbar_{1,\gminia},
 \theta_{1,\gminia},h_{1,\gminia})
:=\bigoplus_{u}
 \bigl(E'_{1,\gminia,u},\delbar'_{1,\gminia,u},
 \theta'_{1,\gminia,u},h'_{1,\gminia,u} \bigr)
\otimes 
 L(u,\gminia),
\]
\[
 (E_{1},\delbar_1,\theta_1,h_1)
:=\bigoplus_{\gminia}
 (E_{1,\gminia},\delbar_{1,\gminia},
 \theta_{1,\gminia},h_{1,\gminia}).
\]
Let 
$(\nbigpzero_{\ast}\nbige_1,\DD_1)
=\bigoplus\bigl(
 \nbigpzero_{\ast}\nbige_{1,\gminia},
 \DD_{1,\gminia}\bigr)$
denote the associated family of
filtered $\lambda$-flat bundles,
equipped with the hermitian metric 
$\nbigp_{\irr}^{(\lambda_0)}h_1
 =\bigoplus \nbigpzero_{\irr}h_{1,\gminia}$.

Taking the associated graded bundle
of $(\nbigpzero_{\ast}\nbige,\DD)$
with respect to the Stokes structure,
we obtain a family of filtered $\lambda$-flat bundles
$\Gr^{\full}_{\gminia}
 \bigl(\nbigp_{\ast}^{(\lambda_0)}\nbige\bigr)$
with $\DD_{\gminia}$.
Note that 
$\DD_{\gminia}
  -(1+\lambda\lambdabar_0)\cdot d\gminia$
and $\DD_{1,\gminia}
 -(1+\lambda\lambdabar_0)\cdot d\gminia$ 
are logarithmic.
We can take a holomorphic isomorphism
$\Phi_{\gminia}:
 \nbigp^{(\lambda_0)}_{\ast}
 \nbige_{1,\gminia}
\lrarr
 \Gr^{\full}_{\gminia}
 (\nbigp_{\ast}^{(\lambda_0)}\nbige)$
with the following property:
\begin{itemize}
\item
 It preserves the filtration $\Fzero$,
\item
$\Gr^{F}(\Phi_{\gminia|U(\lambda_0)\times D})$
 is compatible with $\Res(\DD_{\gminia})$ and 
 $\Res(\DD_{1,\gminia})$.
\end{itemize}

\subsection{Comparison map 
around $\lambda_0\neq 0$}
\label{subsection;08.9.15.30}

Let us consider the case $\lambda_0\neq 0$.
We take a finite covering
$U(\lambda_0)\times (X-D)=\bigcup S_j$ by small multi-sectors.
We would like to construct a $C^{\infty}$-map
$\Phi:\nbigp^{(\lambda_0)}\nbige_1\lrarr
 \nbigp^{(\lambda_0)}\nbige$
such that 
(i) $\Phi_{|S_j}$ are compatible 
 with the full Stokes filtrations $\nbigftilde^{S_j}$ 
 on $S_j$,
(ii) $\Gr^{\nbigftilde^{S_j}}_{\gminia}
 (\Phi_{|S_j})
 =\Phi_{\gminia|S_j}$.
We take a flat splitting 
$\nbigp^{(\lambda_0)}\nbige_{|\Sbar_j}
=\bigoplus
 \nbigp^{(\lambda_0)}\nbige_{\gminia,S_j}$
of the full Stokes filtration $\nbigftilde^{S_j}$.
Then, we obtain a holomorphic isomorphism
$\Phi_{S_j}:\nbigp^{(\lambda_0)}
 \nbige_{1|\Sbar_j}
\simeq\nbigp^{(\lambda_0)}
 \nbige_{|\Sbar_j}$
preserving the decomposition.
By gluing them in $C^{\infty}$,
we obtain the desired map.

\begin{lem}
\label{lem;07.10.22.15}
$\Phi_{|\lambda\times (X-D)}$
is bounded with respect to $h$ and $h_1$
for each $\lambda$.
\end{lem}
\pf
$\Phi_{|\lambda\times (X-D)}$
can be extended to a $C^{\infty}$-isomorphism
$\nbigp\nbigelambda_1\lrarr\nbigp\nbigelambda$
such that
(i) it preserves the parabolic structure,
(ii) $\Gr^F(\Phi_{|\{\lambda\}\times D})$ is compatible
 with the residues.
Then, the claim of the lemma follows 
 from the norm estimate
 (Proposition \ref{prop;07.11.22.55}).
\hfill\qed

\vspace{.1in}
On a small sector $S=S_j$,
we have the decomposition:
\begin{equation}
\label{eq;07.11.24.12}
 \Hom\bigl(\nbigpzero\nbige_{1|\Sbar},
 \nbigpzero\nbige_{|\Sbar}\bigr)
=\bigoplus_{\gminia,\gminib}
 \Hom\bigl(
 \nbigpzero\nbige_{1,\gminia,S},
 \nbigpzero\nbige_{\gminib,S}
 \bigr)
\end{equation}
We give a general lemma.
\begin{lem}
\label{lem;07.11.24.11}
Let $s$ be a section of
$\Hom\bigl(\nbigpzero\nbige_{1,\gminia,S},
 \nbigpzero\nbige_{\gminib,S}\bigr)$
with the following estimate
with respect to $\nbigpzero_{\irr} h$
and $\nbigpzero_{\irr}h_1$:
\[
s=\left\{\begin{array}{ll}
O\bigl(\exp(-\epsilon_1|z_1|^{\ord(\gminia-\gminib)})\bigr)
 & (\gminia\neq\gminib)\\
 \mbox{{}}\\
 O\bigl(|z_1|^{c}\bigr) &(\gminia=\gminib)
\end{array}
\right.
\]
Here, $c\in\real$ and $\epsilon_1>0$.
If we shrink $U(\lambda_0)$, 
we also have the following estimate
with respect to $h$ and $h_1$
for some $\epsilon_2>0$:
\begin{equation}
\label{eq;07.11.24.13}
s=\left\{\begin{array}{ll}
O\bigl(\exp(-\epsilon_2|z_1|^{\ord(\gminia-\gminib)})\bigr)
 & (\gminia\neq\gminib)\\
 \mbox{{}}\\
 O\bigl(|z_1|^{c}\bigr) &(\gminia=\gminib)
\end{array}
\right.
\end{equation}
\end{lem}
\pf
Let $F_S$ be as in Section \ref{subsection;07.11.24.7}.
By the assumption, 
the estimate for $s$ like (\ref{eq;07.11.24.13})
holds with respect to the metrics
$F_S(-\lambdabar+\lambdabar_0)^{\ast}
 \nbigpzero_{\irr}h$
and $F_S(-\lambdabar+\lambdabar_0)^{\ast}
 \nbigpzero_{\irr}h_1$.
We may assume to have
$F_S(-\lambdabar+\lambdabar_0)^{\ast}
 \nbigpzero_{\irr}h_1
=h_1$.
By Corollary \ref{cor;07.11.24.3},
the metrics
$F_S(-\lambdabar+\lambdabar_0)^{\ast}
 \nbigpzero_{\irr}h$
and $h$ are mutually bounded.
Hence, we obtain Lemma \ref{lem;07.11.24.11}.
\hfill\qed

\vspace{.1in}
We have the decomposition
$\Phi=\sum \Phi_{\gminia,\gminib,S}$
corresponding to (\ref{eq;07.11.24.12}).

\begin{lem}
\label{lem;07.10.22.10}
If $\gminia\neq\gminib$, we have the following estimate
with respect to $h$, $h_1$
and the Euclid metric $g_X$ of $X$,
uniformly for $\lambda$:
\begin{equation}
\label{eq;07.10.22.1}
 \Phi_{\gminia,\gminib,S}=O\Bigl(
 \exp\bigl(-\epsilon|z_1|^{\ord(\gminia-\gminib)}\bigr)
 \Bigr)
\end{equation}
\begin{equation}
\label{eq;07.10.22.8}
 \DD(\Phi)_{\gminia,\gminib,S}
=O\Bigl(
 \exp\bigl(-\epsilon|z_1|^{\ord(\gminia-\gminib)}\bigr)
 \Bigr)
\end{equation}
\begin{equation}
\label{eq;07.10.22.3}
 \DD^{1,0}d''\Phi_{\gminia,\gminib,S}
=O\Bigl(
 \exp\bigl(-\epsilon|z_1|^{\ord(\gminia-\gminib)}\bigr)
 \Bigr)
\end{equation}
\begin{equation}
\label{eq;07.10.22.2}
 d''\Phi_{\gminia,\gminib,S}
=O\Bigl(
 \exp\bigl(-\epsilon|z_1|^{\ord(\gminia-\gminib)}\bigr)
 \Bigr)
\end{equation}
Here, 
$d''=\delbar_E+\lambda\theta^{\dagger}$
denotes the holomorphic structure of
$\nbige$ along the $X$-direction,
and 
$\DD^{1,0}$ denotes the $(1,0)$-part of
$\DD$, i.e.,
$\DD^{1,0}:=\lambda\del_E+\theta$.
\end{lem}
\pf
It is clear that
$\Phi_{\gminia,\gminib}=O\Bigl(
 \exp\bigl(-\epsilon|z_1|^{\ord(\gminia-\gminib)}\bigr)
 \Bigr)$ for $h_0$ and $\nbigp_{\irr}^{(\lambda_0)}h_1$.
Since $h_0$ and $\nbigp_{\irr}^{(\lambda_0)}h$
are uniformly mutually bounded up to small polynomial order,
we obtain
$\Phi_{\gminia,\gminib}=O\Bigl(
 \exp\bigl(-\epsilon|z_1|^{\ord(\gminia-\gminib)}\bigr)
 \Bigr)$ for $\nbigp_{\irr}^{(\lambda_0)}h$
and $\nbigp_{\irr}^{(\lambda_0)}h_1$.
Then, we obtain (\ref{eq;07.10.22.1}) due to 
Lemma \ref{lem;07.11.24.11}.
The others can be shown similarly.
\hfill\qed

\vspace{.1in}
We have $d''\Phi_{\gminia,\gminia,S}=0$
by construction.
Hence, we obtain the following uniform
estimate with respect to
$h$, $h_1$ and $g_X$ from (\ref{eq;07.10.22.2}):
\[
 d''\Phi=O\bigl(\exp(-\epsilon|z_1|^{-1})\bigr)
\]
Then, we have the following uniform estimates
with respect to $h$, $h_1$ and $g_X$:
\begin{equation}
\label{eq;07.11.23.150}
 \bigl(\theta\circ d''\Phi
-d''\Phi\circ \theta_1\bigr)
=O\bigl(\exp(-\epsilon|z_1|^{-1})\bigr)
\end{equation}
Let $\delta_{\lambda}':=\del_E-\lambdabar\theta$
and $\delta_{1,\lambda}':=\del_{E_1}-\lambdabar\theta_1$.
We have the relation
$\DD^{(1,0)}=
 \lambda\delta_{\lambda}'+(1+|\lambda|^2)\theta$.
and $\DD_1^{(1,0)}
=\lambda\delta_{1,\lambda}'+(1+|\lambda|^2)\theta_1$.
From (\ref{eq;07.10.22.3}) and (\ref{eq;07.11.23.150}),
we obtain the following uniform estimate 
with respect to $h$, $h_1$ and $g_X$:
\begin{equation}
\label{eq;07.10.22.9}
 \bigl(\delta'_{\lambda}\circ d''\Phi
-d''\Phi\circ \delta_{1,\lambda}'\bigr)
=O\Bigl(\exp(-\epsilon|z_1|^{-1})\Bigr)
\end{equation}

We prepare other estimates.

\begin{lem}
\label{lem;07.11.23.160}
$\Phi$ is uniformly bounded
up to small polynomial order 
with respect to $h$ and $h_1$.
\end{lem}
\pf
$\Phi$ is uniformly bounded with respect to
$h_0$ and $\nbigp^{(\lambda_0)}_{\irr}h_1$.
Hence, it is uniformly bounded up to small polynomial order
with respect to $\nbigp_{\irr}^{(\lambda_0)}h$ and
$\nbigp_{\irr}^{(\lambda_0)}h_1$
(Lemma \ref{lem;08.1.8.20}).
Since the components 
$\Phi_{\gminia,\gminib}$ $(\gminia\neq\gminib)$
are estimated as in (\ref{eq;07.10.22.1}),
we obtain Lemma \ref{lem;07.11.23.160}
by using Lemma \ref{lem;07.11.24.11}.
\hfill\qed

\begin{lem}
\label{lem;07.11.23.200}
We have the uniform estimate 
$\DD\Phi_{\gminia,\gminia,S}=
O\bigl(|z_1|^{-1+\epsilon}\bigr)$
for some $\epsilon>0$,
with respect to $h$, $h_1$ and $g_X$.
\end{lem}
\pf
By construction,
we have the uniform estimate
$\DD\Phi_{\gminia,\gminia,S}
 =O\bigl(|z_1|^{-1+\epsilon_1}\bigr)$
for some $\epsilon_1>0$,
with respect to 
$h_0$, $\nbigpzero_{\irr} h_1$ and $g_X$.
Since
$h_0$ and $\nbigpzero_{\irr}h$
are uniformly mutually bounded
up to $|z_1|^{-\eta|\lambda-\lambda_0|}
 (-\log|z_1|)^N$-order
(Lemma \ref{lem;08.1.8.20}),
we obtain the uniform estimate
$\DD\Phi_{\gminia,\gminia,S}
 =O\bigl(|z_1|^{-1+\epsilon_2}\bigr)$
for some $\epsilon_2>0$
with respect to $\nbigpzero_{\irr}h$,
$\nbigpzero_{\irr}h_1$ and $g_X$.
Then, we obtain the claim of 
Lemma \ref{lem;07.11.23.200}
by using Lemma \ref{lem;07.11.24.11}.
\hfill\qed

\vspace{.1in}

From Lemma \ref{lem;07.11.23.200}
and (\ref{eq;07.10.22.8}),
we obtain the following estimates
with respect to $h$, $h_1$ and $g_X$
uniformly for $\lambda$:
\begin{equation}
\label{eq;07.10.22.50}
\DD\Phi=O\bigl(|z_1|^{-1+\epsilon}\bigr)
\end{equation}

\subsection{Comparison map around $\lambda_0=0$}

We almost repeat the construction in
Section \ref{subsection;08.9.15.30}.
We take a finite covering
$(U(0)-\{0\})\times (X-D)
=\bigcup S_j$ by a small multi-sectors.
We would like to construct
a $C^{\infty}$-map 
$\Phi:\nbigp^{(0)}\nbige_1\lrarr
 \nbigp^{(0)}\nbige$
such that 
(i) $\Phi_{|S_j}$ are compatible 
 with the full Stokes filtrations $\nbigftilde^{S_j}$ 
 on $S_j$,
(ii) $\Gr^{\nbigftilde^{S_j}}_{\gminia}
 (\Phi_{|S_j}) =\Phi_{\gminia|S_j}$.
We take a flat splitting
$\nbigp^{(0)}\nbige_{|\Sbar_j}
=\bigoplus \nbigp^{(0)}\nbige_{\gminia,S_j}$
of the full Stokes filtration $\nbigftilde^{S_j}$.
Then, we obtain the isomorphism
$\Phi_{S_j}:\nbigp^{(0)}\nbige_{1|\Sbar_j}
=\nbigp^{(0)}\nbige_{|\Sbar_j}$
preserving the decomposition.
By gluing them in $C^{\infty}$,
we obtain the desired map.
The following lemma can be shown 
by the argument in the proof of 
Lemma \ref{lem;07.10.22.15}.
\begin{lem}
$\Phi_{|\lambda\times (X-D)}$
is bounded with respect to $h$ and $h_1$
for each $\lambda$.
\hfill\qed
\end{lem}

On a small sector $S=S_j$,
we have the decomposition:
\begin{equation}
\label{eq;07.11.24.15}
 \Hom\bigl(\nbigp^{(0)}\nbige_{1|\Sbar},
 \nbigp^{(0)}\nbige_{|\Sbar}\bigr)
=\bigoplus_{\gminia,\gminib}
 \Hom\bigl(
 \nbigp^{(0)}\nbige_{1,\gminia,S},
 \nbigp^{(0)}\nbige_{\gminib,S}
 \bigr)
\end{equation}

\begin{lem}
\label{lem;07.11.24.16}
Let $s$ be a section of
$\Hom\bigl(\nbigp^{(0)}\nbige_{1,\gminia,S},
 \nbigp^{(0)}\nbige_{\gminib,S}\bigr)$
with the following estimate
with respect to $\nbigp^{(0)}_{\irr} h$
and $\nbigp^{(0)}_{\irr}h_1$:
\[
s=\left\{\begin{array}{ll}
O\bigl(\exp(-\epsilon_1\cdot|\lambda|^{-1}\cdot
 |z_1|^{\ord(\gminia-\gminib)})\bigr)
 & (\gminia\neq\gminib)\\
 \mbox{{}}\\
 O\bigl(|z_1|^{c}\bigr) &(\gminia=\gminib)
\end{array}
\right.
\]
Here, $c\in\real$ and $\epsilon_1>0$.
If we shrink $U(0)$, 
we also have the following estimate
with respect to $h$ and $h_1$
for some $\epsilon_2>0$:
\[
 s=\left\{\begin{array}{ll}
O\bigl(\exp(-\epsilon_2\cdot|\lambda|^{-1}\cdot
 |z_1|^{\ord(\gminia-\gminib)})\bigr)
 & (\gminia\neq\gminib)\\
 \mbox{{}}\\
 O\bigl(|z_1|^{c}\bigr) &(\gminia=\gminib)
\end{array}
\right.
\]
\end{lem}
\pf
It can be shown using the same argument
as that in the proof of Lemma \ref{lem;07.11.24.11}.
\hfill\qed

\vspace{.1in}

We have the decomposition
$\Phi=\sum \Phi_{\gminia,\gminib,S}$
corresponding to (\ref{eq;07.11.24.15}).

\begin{lem}
If $\gminia\neq\gminib$, 
we have the following estimate
with respect to $h$ and $h_1$,
uniformly for $\lambda$:
\begin{equation}
 \label{eq;07.10.22.21}
 \Phi_{\gminia,\gminib,S}=O\Bigl(
 \exp\bigl(-\epsilon|\lambda|^{-1}
 |z_1|^{\ord(\gminia-\gminib)}\bigr)
 \Bigr)
\end{equation}
\begin{equation}
\label{eq;07.10.22.46}
 \DD(\Phi)_{\gminia,\gminib,S}
=O\Bigl(
 \exp\bigl(-\epsilon|\lambda|^{-1}|z_1|^{\ord(\gminia-\gminib)}\bigr)
 \Bigr)
\end{equation}
\begin{equation}
\label{eq;07.10.22.41}
 \DD^{1,0}d''\Phi_{\gminia,\gminib,S}
=O\Bigl(
 \exp\bigl(-\epsilon|\lambda|^{-1}|z_1|^{\ord(\gminia-\gminib)}\bigr)
 \Bigr)
\end{equation}
\begin{equation}
 \label{eq;07.10.22.22}
 d''\Phi_{\gminia,\gminib,S}
=O\Bigl(
 \exp\bigl(-\epsilon|\lambda|^{-1}|z_1|^{\ord(\gminia-\gminib)}\bigr)
 \Bigr)
\end{equation}
\end{lem}
\pf
We have 
$\Phi_{\gminia,\gminib,S}=O\Bigl(
 \exp\bigl(-\epsilon|\lambda|^{-1}\cdot
 |z_1|^{\ord(\gminia-\gminib)}\bigr)
 \Bigr)$
with respect to $h_0$ and $\nbigp^{(0)}_{\irr}h_1$.
Since $\nbigp_{\irr}^{(0)}h$ and $h_0$ are 
uniformly mutually bounded
up to small polynomial order,
we obtain $\Phi_{\gminia,\gminib,S}=O\Bigl(
 \exp\bigl(-\epsilon|\lambda|^{-1}\cdot
 |z_1|^{\ord(\gminia-\gminib)}\bigr)
 \Bigr)$
with respect to 
$\nbigp_{\irr}^{(0)}h$ and $\nbigp_{\irr}^{(0)}h_1$.
Then, we obtain (\ref{eq;07.10.22.21}) with respect to
$h$ and $h_1$ by Lemma \ref{lem;07.11.24.16}.
The others can be shown similarly.
\hfill\qed

\vspace{.1in}
By construction, we have
$d''\Phi_{\gminia,\gminia}=0$.
Hence, we obtain the following estimate
with respect to $h$ and $h_1$, from (\ref{eq;07.10.22.22}):
\[
 d''\Phi=O\bigl(\exp(-\epsilon|\lambda|^{-1}|z_1|^{-1})\bigr)
\]
Then,
we obtain the following estimate 
with respect to $h$, $h_1$ and $g_X$:
\begin{equation}
 \label{eq;07.10.22.40}
 \bigl(\theta\circ d''\Phi
-d''\Phi\circ \theta_1\bigr)
=O\Bigl(
 \exp\bigl(-\epsilon|\lambda|^{-1}|z_1|^{-1}\bigr)
 \Bigr)
\end{equation}
Thus, we obtain the following estimate
from (\ref{eq;07.10.22.41}) and (\ref{eq;07.10.22.40}):
\begin{equation}
\label{eq;07.10.22.51}
  \bigl(\delta_{\lambda}'\circ d''\Phi
-d''\Phi\circ \delta_{1,\lambda}'\bigr)
=O\Bigl(
 \exp\bigl(-\epsilon|\lambda|^{-1}|z_1|^{-1}\bigr)
 \Bigr)
\end{equation}

We have the counterparts of
Lemma \ref{lem;07.11.23.160}
and Lemma \ref{lem;07.11.23.200},
which can be shown in similar ways.
\begin{lem}
\label{lem;07.11.23.203}
$\Phi$ is uniformly bounded
up to small polynomial order,
with respect to $h$, $h_1$ and $g_X$.
\hfill\qed
\end{lem}

\begin{lem}
\label{lem;07.11.23.201}
We have the estimate 
$\DD\Phi_{\gminia,\gminia,S}=O(|z_1|^{-1+\epsilon})$
for some $\epsilon>0$
with respect to $h$, $h_1$ and $g_X$.
\hfill\qed
\end{lem}

From Lemma \ref{lem;07.11.23.201}
and (\ref{eq;07.10.22.46}),
we obtain the following estimate
with respect to $h$, $h_1$ and $g_X$:
\begin{equation}
\label{eq;07.10.22.47}
 \DD\Phi=O\bigl(|z_1|^{-1+\epsilon}\bigr)
\end{equation}

\subsection{End of proof of 
 Proposition \ref{prop;07.10.22.30}}
\label{subsection;07.11.24.6}

By construction,
it is easy to see that
$\nbigpzero_{\irr}h_1$ and 
$h_0$ are mutually bounded.
Hence,
we have only to show that
$\nbigpzero_{\irr}h$ and 
$\nbigpzero_{\irr}h_1$ are mutually bounded
uniformly for $\lambda$.

\begin{lem}
\label{lem;07.11.24.17}
Once we show 
that $h$ and $h_1$ are mutually bounded,
we obtain $\nbigpzero_{\irr}h$
and $\nbigpzero_{\irr}h_1$ are mutually bounded.
\end{lem}
\pf
Let $S$ and $F_S$
be as in Section \ref{subsection;07.11.24.7}.
If $h$ and $h_1$ are mutually bounded,
$F_S(\lambdabar-\lambdabar_0)^{\ast}h$
and $F_S(\lambdabar-\lambdabar_0)^{\ast}h_1$
are mutually bounded.
We may assume to have
$F_S(\lambdabar-\lambdabar_0)^{\ast}h_1
=\nbigpzero_{\irr}h_1$.
By Corollary \ref{cor;07.11.24.3},
$F_S(\lambdabar-\lambdabar_0)^{\ast}h$
and $\nbigpzero_{\irr}h$ are mutually bounded.
Thus, Lemma \ref{lem;07.11.24.17} is proved.
\hfill\qed

\vspace{.1in}

Let us show that $h$ and $h_1$ are mutually bounded.
We have only to show that
$\Phi$ and $\Phi^{-1}$ are bounded
uniformly for $\lambda$
with respect to $h$ and $h_1$.

Let $\Phi_{\lambda}:=\Phi_{|\{\lambda\}\times (X-D)}$.
We regard it as the $C^{\infty}$-section of
the harmonic bundle
$Hom(\nbigelambda_1,\nbigelambda)$
with the induced $\lambda$-connection $\DDlambda_2$
and the pluri-harmonic metric $h_2$.
Let $\pi:X\lrarr D$ be the projection.
Let $Q$ be any point of $D$.
We put $\pi^{-1}(Q)^{\ast}:=\pi^{-1}(Q)-\{Q\}$.
Applying Lemma \ref{lem;08.9.29.12},
we obtain the following inequality
on each $\pi^{-1}(Q)^{\ast}$:
\begin{equation}
 \label{eq;07.10.22.42}
 -\frac{\del^2}{\del z_1\del\zbar_1}
 \bigl|\Phi_{\lambda}\bigr|_{h_2}^2
\leq 
 \bigl|\DDlambda_2\Phi_{\lambda}\bigr|_{h_2}^2
+2\bigl|\Phi_{\lambda}\bigr|_{h_2}\cdot
 \bigl|
 \delta_{2,\lambda}'\circ d''_{2,\lambda}(\Phi_{\lambda})\bigr|_{h_2}
\end{equation}
From (\ref{eq;07.10.22.42}),
(\ref{eq;07.10.22.9}), (\ref{eq;07.10.22.50}),
(\ref{eq;07.10.22.51}), (\ref{eq;07.10.22.47}),
Lemma \ref{lem;07.11.23.160}
and Lemma \ref{lem;07.11.23.203},
we obtain the inequality 
on $\pi^{-1}(Q)^{\ast}$
\[
 -\frac{\del^2}{\del z_1\del\zbar_1}|\Phi_{\lambda}|^2
 \leq C|z_1|^{-2+\epsilon},
\]
where the constant 
$C$ is independent of $\lamda$ and $Q$.
Since we have already known the boundedness
of $|\Phi_{\lambda}|_{h_2}^2$ on 
$\pi^{-1}(Q)^{\ast}$
(Lemma \ref{lem;07.10.22.15}),
the inequality holds on $\pi^{-1}(Q)$
as distributions.
Hence, the values 
$\bigl|\Phi\bigr|_{h_2}$ on $\pi^{-1}(Q)^{\ast}$
is dominated by the values on the boundary
$\del\pi^{-1}(Q)$.
Thus, we obtain the uniform boundedness of $\Phi$.
Similarly, we obtain the uniform boundedness
for $\Phi^{-1}$.
Thus, we obtain that
$h$ and $h_1$ are mutually bounded,
and the proof of Proposition \ref{prop;07.10.22.30}
is finished.
\hfill\qed

\chapter[Prolongation and reduction]
 {Prolongation and Reduction of Variation of
 Polarized Pure Twistor Structure}
\label{section;08.9.14.100}
In this chapter,
we study the reduction with respect to
the Stokes structure of 
unramifiedly good wild harmonic bundles
(Theorem \ref{thm;08.9.16.100}).
This is one of the main results 
in the study of wild harmonic bundles.
Recall that we obtain a polarized mixed twistor 
structure from a tame harmonic bundle
by taking Gr with respect to KMS-structure,
which is one of the most important achievements
in the study of tame harmonic bundles.
We also have the reduction from
polarized mixed twistor structure
to polarized mixed Hodge structure
by taking Gr with respect to 
the weight filtration.
(See \cite{mochi2} for these reductions.
 See also a survey in \cite{mochi8}.)
By these reductions,
the study of the asymptotic behaviour of
wild harmonic bundles
is reduced to that of
variation of Hodge structures.

\vspace{.1in}
In Section \ref{subsection;07.11.5.1},
we define the canonical prolongment
of a variation of pure twistor structures
whose underlying harmonic bundle
is good wild.
In Section \ref{subsection;07.11.24.25},
we first construct a family of meromorphic
$\lambda$-flat bundles
$(\nbigq\nbige,\DDlambda)$
on $\cnum_{\lambda}\times (X,D)$
associated to a good wild harmonic bundle
$\harmonicbundle$ on $(X,D)$.
This family of meromorphic flat bundles
will be also useful in Chapter \ref{section;07.10.12.10}.
This subsection is a continuation of
Section \ref{subsection;08.9.14.11}.
Applying the construction to
the harmonic bundle
$(E,\del_E,\theta^{\dagger},h)$
on the conjugate complex manifolds
$(X^{\dagger},D^{\dagger})$,
we obtain the prolongment 
$(\nbigq\nbige^{\dagger},\DD^{\dagger})$
on $\cnum_{\mu}\times (X^{\dagger},D^{\dagger})$
in Section \ref{subsection;08.10.26.10}.
Then, we show in Section 
\ref{subsection;08.9.16.1}
that $(\nbigq\nbige,\nbigq\nbige^{\dagger})$
gives a  meromorphic prolongment 
(Section \ref{subsection;07.10.10.2})
of the variation of polarized pure twistor structure
$(\nbige^{\sankaku},\DD^{\sankaku},\nbigs)$
associated to $\harmonicbundle$,
which is called the canonical prolongment.

In Section \ref{subsection;07.11.24.50},
we give the statement of the reduction theorem
(Theorem \ref{thm;08.9.16.100}).
Note that it also implies a characterization
of the canonical prolongment 
$(\nbigq\nbige,\nbigq\nbige^{\dagger})$,
which is non-trivial
because of the absence of the uniqueness 
of an unramifiedly good meromorphic prolongment 
with irregular singularity.
We explain a plan of the proof
in Section \ref{subsection;08.9.25.1}.

In Section \ref{subsection;08.9.29.14},
we give a sufficient condition
for a holomorphic vector bundle
on $\proj^1$ with a symmetric pairing
to be a polarized pure twistor structure.
Intuitively,
the twistor property and
the positive definiteness
should be open conditions.
We can expect that
a holomorphic vector bundle
with a symmetric pairing on $\proj^1$
is a polarized pure twistor structure,
if it is ``close to''
a polarized pure twistor structure.
We give a condition
to describe ``close to'',
which is convenient for our purpose.
This is one of the main tools in the proof.

We study the one step reduction
in Section \ref{subsection;07.10.15.20}.
Briefly speaking,
we show that 
the one step reduction is a variation of
polarized pure twistor structure
if and only if the original one is so.
Main tools are 
Lemma \ref{lem;07.6.2.90}
and Lemma \ref{lem;07.6.2.60}.
However, we cannot conclude anything
about prolongments in this stage,
which will be postponed until Section
\ref{subsection;08.9.25.10}.

We study
in Section \ref{subsection;07.11.24.30}
the full reduction in the case that $D$ is smooth.
As in Section \ref{subsection;07.10.15.20},
we show that 
the full reduction is a variation of
polarized pure twistor structure
if and only if the original one is so.
Moreover, we can obtain a characterization
of the canonical prolongment.
The argument is essentially the same
as that in Section \ref{subsection;07.10.15.20},
except that we use Proposition \ref{prop;08.9.17.1}
instead of Lemma \ref{lem;07.6.2.60}.

The proof of Theorem  \ref{thm;08.9.16.100}
is finished in Section 
\ref{subsection;08.9.25.15}.

\vspace{.1in}
In Section \ref{subsection;08.9.25.20},
we argue the norm estimate for
holomorphic sections of
$\nbigp\nbigelambda$.
Briefly speaking,
our result says that
the pluri-harmonic metric $h$ is 
determined up to boundedness
by the parabolic structure and the residues,
as in the case of tame harmonic bundles.
In Sections
\ref{subsection;08.9.25.21}--\ref{subsection;08.9.25.22},
we construct comparison maps
between the reductions and the original one,
and we show that they are bounded
(Theorem \ref{thm;07.11.5.41}).
Since we can apply the norm estimate 
for a tame harmonic bundle
to a wild harmonic bundle
with a unique irregular value,
Theorem \ref{thm;07.11.5.41}
implies the norm estimate
for unramifiedly good wild harmonic bundles.
The norm estimate for good wild harmonic bundles
can be easily reduced to the unramified case.

In Sections 
\ref{subsection;08.9.25.30}--\ref{subsection;07.11.5.101},
we give a rather detailed description
of the norm estimate in the surface case,
which is essentially the same as 
those in Section 2.5 of \cite{mochi4}.
The result will be used in
Proposition \ref{prop;07.12.6.12}
to show the vanishing of some characteristic
numbers of the filtered bundle
associated to a good wild harmonic bundle.

Section \ref{subsection;08.9.14.50} is an appendix
in which we consider the reduction with respect to
KMS structure
for a {\em regular} meromorphic variation of
pure twistor structure with a pairing 
on $\proj^1\times\Delta^{\ast}$,
and we shall give a characterization of 
purity and polarizability
in terms of the limit twistor structure
(Lemma \ref{lem;07.7.11.12}).
For that purpose, in Subsection
\ref{subsection;10.6.1.2}--\ref{subsection;08.1.28.10},
we shall review the construction
of a limit vector bundle on $\proj^1$
from a variation of twistor structure
on a punctured disc
with regular meromorphic extension,
called the limit twistor structure.
When it comes from a tame harmonic bundle,
the construction is explained in \cite{mochi2}.
We need only minor change.
Then, we show the characterization.
Although a similar result was obtained in \cite{sabbah2}
for $\nbigr$-triples with a different method,
we would like to understand it
from our own viewpoint.
Although we restrict ourselves to
the one dimensional case,
the arguments can be generalized
in the higher dimensional case.

\section{Canonical prolongation}
\label{subsection;07.11.5.1}
\subsection{Prolongment $\nbigq\nbige$}
\label{subsection;07.11.24.25}
\index{filtered bundle $\nbigqzero_{\ast}\nbige$}
\index{sheaf $\nbigqzero_{\veca}\nbige$}
\index{sheaf $\nbigqzero\nbige$}

Let $X$ be a complex manifold
with a simple normal crossing hypersurface 
$D=\bigcup_{i\in\Lambda}D_i$.
Let $\harmonicbundle$ be a good wild harmonic bundle
on $X-D$.
Let $\nbigx^{(\lambda_0)}$ denote 
a neighbourhood of $\{\lambda_0\}\times X$
in $\cnum_{\lambda}\times X$.
We use the symbol $\nbigd^{(\lambda_0)}$
in a similar meaning.
We have the associated family of 
$\lambda$-flat bundles
$\bigl(\nbigp_{\ast}^{(\lambda_0)}\nbige,\DD\bigr)$
on $(\nbigx^{(\lambda_0)},\nbigd^{(\lambda_0)})$.
(See Section \ref{subsection;08.9.14.10}.)
In the case $\lambda_0\neq 0$,
we have the following deformation
in Section \ref{subsection;10.5.17.51}:
\[
 \nbigq_{\veca}^{(\lambda_0)}\nbige:=
 \bigl( 
 \nbigp_{\veca}^{(\lambda_0)}\nbige
 \bigr)^{(T(\lambda))},
\quad
 T(\lambda)=\frac{1}{1+\lambda\lambdabar_0}
\]
In the case $\lambda_0=0$,
we define 
$\nbigq_{\veca}^{(0)}\nbige:=
 \nbigp_{\veca}^{(0)}\nbige$.
Thus, we obtain a global filtered bundles
\[
 \nbigq_{\ast}^{(\lambda_0)}\nbige
 =\bigl(
 \nbigqzero_{\veca}\nbige\,\big|\,
 \veca\in\real^{\Lambda}
 \bigr) 
\]
on $(\nbigx^{(\lambda_0)},
 \nbigd^{(\lambda_0)})$
with the family of meromorphic flat $\lambda$-connections
$\DD$.
We put 
$\nbigqzero\nbige:=
 \bigcup_{\veca\in\real^{\Lambda}}
 \nbigqzero_{\veca}\nbige$,
which is an $\nbigo_{\nbigx^{(\lambda_0)}}
  \bigl(\ast\nbigd^{(\lambda_0)}\bigr)$-locally free sheaf.

\begin{lem}
Let $P$ be any point of $X$.
Let $X_P$, $U_P(\lambda_0)$ and $U_P(\lambda_1)$
be as in Section {\rm\ref{subsection;07.12.22.110}}.
Then, we have a natural isomorphism
$\nbigq^{(\lambda_0)}\nbige_{|U_P(\lambda_1)\times X_P}
\simeq
 \nbigq^{(\lambda_1)}\nbige_{|U_P(\lambda_1)\times X_P}$.
 We mean by ``natural'' that
 the restriction to $X-D$ is the identity.
\end{lem}
\pf
It follows from Proposition \ref{prop;07.10.8.100}
and Lemma \ref{lem;07.12.21.60}.
\hfill\qed

\vspace{.1in}

Thus, 
we obtain a family of meromorphic $\lambda$-flat bundles
$(\nbigq\nbige,\DD)$ on
$\cnum_{\lambda}\times (X,D)$.
The restriction of $(\nbigq\nbige,\DD)$ to $\nbigxlambda$
is denoted by $(\nbigq\nbigelambda,\DDlambda)$.
\index{sheaf $\nbigq\nbigelambda$}
The following theorem is clear from the construction.
\begin{thm}
 \mbox{{}}\label{thm;07.12.4.15}
\begin{itemize}
\item
 We have a natural isomorphism of
 the meromorphic $\lambda$-flat bundles
 $(\nbigq\nbigelambda,\DDlambda)\simeq
 (\nbigp\nbigelambda,\DDlambda)^{(T)}$
 where $T:=(1+|\lambda|^2)^{-1}$.
 We mean by ``natural'' that
 the restriction to $X-D$ is the identity.
 Moreover, we have the isomorphism of
 good filtered $\lambda$-flat bundles
$(\nbigq_{\ast}^{(\lambda)}\nbigelambda,\DDlambda)
 \simeq(\nbigp_{\ast}\nbigelambda,\DDlambda)^{(T)}$.
\item
 Under the setting in Section 
 {\rm\ref{subsection;07.11.17.1}},
 the set of
 the irregular values 
 $\Irr(\nbigq_{\veca}\nbigelambda,\DDlambda)$
 is $\Irr(\theta)$.
\item
 $(\nbigqzero_{\ast}\nbige,\DD)$
 is a good filtered $\lambda$-flat bundles
 with the KMS-structure at $\lambda_0$.
 If $\theta$ is unramified,
 $(\nbigqzero_{\ast}\nbige,\DD)$ is also unramified.
\item
 Under the setting of 
 Section {\rm\ref{subsection;07.11.17.1}},
 the set of the irregular values 
 $\Irr(\nbigqzero_{\veca}\nbige,\DD)$
 is $\Irr(\theta)$.
\hfill\qed
\end{itemize}
\end{thm}

Let us look at the Stokes filtrations
under the setting 
in Section \ref{subsection;07.11.17.1}.
We take an auxiliary sequence
for $\Irr(\theta)$.
Let $S$ be a small multi-sector in $X-D$.
If $\lambda\neq 0$,
we have the partial Stokes filtration
$\nbigf^{S,\vecm(i)}(\nbigq\nbigelambda_{|\Sbar})$ of
$\nbigq\nbigelambda_{|\Sbar}$
in the level $\vecm(i)$.
(See Section \ref{subsection;07.11.18.10}.)

\begin{prop}
\label{prop;07.11.24.27}
Let $f$ be a flat section of $\nbigelambda_{|S}$.
We have 
 $f\in\nbigf^{S,\vecm(i)}_{\gminib}
 (\nbigq\nbigelambda)_{|S}$,
if and only if the following estimate holds
for some $C>0$ and $N>0$:
\[
 \bigl|f\, \exp\bigl((\lambda^{-1}+\lambdabar)\,
 \etabar_{\vecm(i)}(\gminib)\bigr)\bigr|_h
=O\Bigl(
 \exp\bigl(C|\vecz^{\vecm(i+1)}|\bigr)
\,\prod_{k(i+1)<j\leq \ell}|z_j|^{-N}
 \Bigr)
\]
\end{prop}
\pf
It follows from Proposition \ref{prop;07.10.11.3}
and Theorem \ref{thm;07.12.4.15}.
\hfill\qed

\subsection{Prolongments 
 $\nbigp\nbige^{\dagger\mu}$,
 $\nbigp^{(\mu_0)}_{\ast}\nbige^{\dagger}$
 and  $\nbigq\nbige^{\dagger}$}
\index{sheaf $\nbigp\nbige^{\dagger\mu}$}
\index{filtered sheaf
 $\nbigp^{(\mu_0)}_{\ast}\nbige^{\dagger}$}
\index{sheaf $\nbigq\nbige^{\dagger}$}
\label{subsection;08.10.26.10}

For simplicity, we restrict ourselves
to the setting in Section \ref{subsection;07.11.17.1}.
We obtain an unramifiedly good wild harmonic bundle
$(E,\del_E,\theta^{\dagger},h)$ 
on $X^{\dagger}-D^{\dagger}$.
We have the family of $\mu$-flat bundles
$(\nbige^{\dagger},\DD^{\dagger})$
on $\cnum_{\mu}\times (X^{\dagger}-D^{\dagger})$,
and the restriction 
$(\nbige^{\dagger\mu},\DD^{\dagger\mu})$ to
$\{\mu\}\times(X^{\dagger}-D^{\dagger})$.
As in Section \ref{subsection;07.11.24.20},
we obtain an unramifiedly good filtered $\mu$-flat bundle
$\bigl(\nbigp_{\ast}\nbige^{\dagger\mu},
 \DD^{\dagger\mu}\bigr)$ on $(X^{\dagger},D^{\dagger})$
for each $\mu$.
The set of the irregular values of $\DD^{\dagger\mu}$
is given as follows:
\[
 \Irr\bigl(\nbigp\nbige^{\dagger\mu},
 \DD^{\dagger\mu}\bigr)
=\bigl\{
 (1+|\mu|^2)\,\gminiabar\,\big|\,
 \gminia\in\Irr(\theta)
 \bigr\}
\]

Let $U(\mu_0)$ denote 
a small neighbourhood of $\mu_0$.
Applying the construction in Chapter
\ref{section;07.11.24.24} to 
$(E,\del_E,\theta^{\dagger},h)$,
we obtain a family of filtered $\mu$-flat bundles
$(\nbigp_{\ast}^{(\mu_0)}\nbige^{\dagger},
 \DD^{\dagger})$ with the KMS-structure at $\mu_0$.
Applying the deformation as in
Section \ref{subsection;07.11.24.25},
we obtain a family of filtered $\mu$-flat bundles
$(\nbigq^{(\mu_0)}_{\ast}\nbige^{\dagger},
\DD^{\dagger})$ on 
$U(\mu_0)\times (X^{\dagger},D^{\dagger})$
and $\nbigo_{U(\mu_0)\times X^{\dagger}}
\bigl(\ast (U(\mu_0)\times D^{\dagger})\bigr)$-free module
$\nbigq^{(\mu_0)}\nbige^{\dagger}$.
By gluing them,
we obtain the family of 
meromorphic $\mu$-flat bundles
$(\nbigq\nbige^{\dagger},\DD^{\dagger})$
on $\cnum_{\mu}\times (X^{\dagger},D^{\dagger})$.
The restriction of $\nbigq\nbige^{\dagger}$
to $\{\mu\}\times X$
is denoted by
$\nbigq\nbige^{\dagger\mu}$.

We have a characterization of the Stokes filtration
of the meromorphic $\mu$-flat bundle
$\bigl(\nbigq\nbige^{\dagger\mu},
 \DD^{\dagger\mu}\bigr)$ in the case $\mu\neq 0$
in terms of the growth order with respect to $h$.
Let $S^{\dagger}$ be a small multi-sector in 
$X^{\dagger}-D^{\dagger}$.
We have the partial Stokes filtration
$\nbigf^{S^{\dagger},\vecm(i)}
 (\nbigq\nbige^{\dagger\mu}_{|\Sbar^{\dagger}})$ of 
$\nbigq\nbige^{\dagger\mu}_{|\Sbar^{\dagger}}$
in the level $\vecm(i)$
indexed by the ordered set 
$\bigl(\Irr(\theta^{\dagger},\vecm(i)),
 \leq^{\mu}_{S^{\dagger}}\bigr)$.

\begin{prop}
\label{prop;07.11.24.28}
Let $f$ be a flat section of 
$\nbige^{\dagger\mu}_{|S^{\dagger}}$.
We have 
$f\in\nbigf^{S^{\dagger},\vecm(i)}_{\gminibbar}$,
if and only if the following estimate holds
for some $C>0$ and $N>0$:
\[
 \Bigl|f\,
 \exp\bigl((\mu^{-1}+\mubar)\,
 \etabar_{\vecm(i)}(\gminibbar)\bigr)\Bigr|_h
=O\Bigl(
 \exp\bigl(C|\vecz^{\vecm(i+1)}|\bigr)
 \prod_{k(i+1)<j\leq \ell}|z_j|^{-N}
 \Bigr)
\]
\hfill\qed
\end{prop}

\subsection{Canonical prolongation
of a variation of polarized pure twistor structures}
\label{subsection;08.9.16.1}

We continue to use the setting in Section
\ref{subsection;07.11.17.1}.
We have the variation of
polarized pure twistor structures 
$(\nbige^{\sankaku},\DD^{\sankaku},\nbigs)$ of weight $0$
associated to $\harmonicbundle$.
(See \cite{s3}. See also \cite{mochi2}.)
In particular, it gives a variation of twistor structure
with a symmetric pairing.

\begin{prop}
\label{prop;07.11.24.29}
The pair $(\nbigq\nbige,\nbigq\nbige^{\dagger})$
gives an unramifiedly good meromorphic prolongment of
$\bigl(\nbige^{\sankaku},\DD^{\sankaku},\nbigs\bigr)$
(Section {\rm\ref{subsection;07.10.10.2}}).
\end{prop}
\pf
Let $\lambda=\mu^{-1}$.
Let us show that
the Stokes filtration $\nbigf^{S,\vecm(i)}$ of
$(\nbigq\nbigelambda,\DD^{\lambda,f})$
and $\nbigf^{\dagger\,S,\vecm(i)}$
of $(\nbigq\nbige^{\dagger\,\mu},\DD^{\dagger\,\mu\,f})$
are the same.
Note 
\[
 \Re\bigl((\lambda^{-1}+\lambdabar)\,\gminib\bigr)
=\Re\bigl((\mu^{-1}+\mubar)\,\gminibbar\bigr)
\]
for $\gminib\in\Irrbar(\theta,\vecm(i))$
and $\lambda=\mu^{-1}$.
Since both the flat filtrations
$\nbigf^{S,\vecm(i)}$ and 
 $\nbigf^{\dagger\,S,\vecm(i)}$
are characterized by the same condition
on the growth order of the norms of flat sections
with respect to $h$
(Proposition \ref{prop;07.11.24.27}
and Proposition \ref{prop;07.11.24.28}),
they are the same.

Let us show that the pairing
$\nbigs_0:
 \nbige\otimes\sigma^{\ast}\nbige^{\dagger}
 \lrarr \nbigo_{\nbigx-\nbigd}$
can be extended to 
$\nbigq\nbigs_0:
 \nbigq\nbige\otimes\sigma^{\ast}\nbigq\nbige^{\dagger}
\lrarr\nbigo_{\nbigx}(\ast\nbigd)$.
Take $\lambda_0\in \cnum_{\lambda}$,
and let $U(\lambda_0)$ denote a small neighbourhood
of $\lambda_0$.
Let $\nbigx^{(\lambda_0)}:=
 U(\lambda_0)\times X$
and $\nbigd^{(\lambda_0)}:=
 U(\lambda_0)\times D$.
\begin{lem}
\label{lem;07.11.25.60}
$\nbigs_0$ is prolonged to a pairing
\[
 \nbigpzero\nbigs_0:
\nbigpzero\nbige
 \otimes
 \sigma^{\ast}\nbigp^{(\sigma(\lambda_0))}
 \nbige^{\dagger}
\lrarr \nbigo_{\nbigx^{(\lambda_0)}}
 (\ast\nbigd^{(\lambda_0)}).
\]
\end{lem}
\pf
Let $f_1$ be a section of 
$\nbigpzero(\nbige)$
on $\nbigxzero$.
Let $\nbigpzero_{\irr}h$ be 
the hermitian metric given in
(\ref{eq;07.7.21.4}).
By construction,
we have the estimates 
$|f_1|_{\nbigpzero_{\irr}h}=
 O\Bigl(\prod_{j=1}^{\ell}|z_j|^{-N}\Bigr)$
for some $N$.
Let $f_2$ be a section of
$\nbigp^{(\sigma(\lambda_0))}\nbige$
on $\sigma\bigl(\nbigxzero\bigr)$.
Let $\nbigp^{\dagger(\mu_0)}_{\irr}h$
of $p_{\mu}^{-1}E$ on $U(\mu_0)\times (X-D)$
given as follows:
\[
 \nbigp^{\dagger(\mu_0)}_{\irr}h(u,v)
=h\Bigl(
 g^{\dagger}_{\irr}(\mu-\mu_0)u,
 g^{\dagger}_{\irr}(\mu-\mu_0)v
 \Bigr),
\quad
 g^{\dagger}_{\irr}(\mu):=
 \exp\left(
 \sum_{\gminia\in\Irr(\theta)}
\mu\,\gminia\,\pi^{\dagger}_{\gminia}
 \right)
\]
We have the estimate
$|f_2|_{\nbigp^{\dagger(\sigma(\lambda_0))}_{\irr}h}
=O\Bigl(\prod_{j=1}^{\ell}|z_j|^{-N}\Bigr)$
for some $N$.
Recall $\sigma(\lambda)=-\lambdabar$,
if we regard $\sigma$
as a morphism $\cnum_{\lambda}\lrarr \cnum_{\mu}$.
(See Subsection \ref{subsection;10.5.24.20}.)
Then, we have
\[
 g_{\irr}^{\dagger}
 \bigl(\sigma(\lambda)-\sigma(\lambda_0)\bigr)
=\bigl(g_{\irr}(\lambda-\lambda_0)^{-1}\bigr)^{\dagger}.
\]
We obtain the following:
\[
 \nbigs_0(f_1,\sigma^{\ast}f_2)
=h\bigl(f_1,\sigma^{\ast}f_2\bigr)
=h\Bigl(
 g_{\irr}(\lambda-\lambda_0)f_1,
\sigma^{\ast}\bigl(
 g_{\irr}^{\dagger}
 \bigl(\mu-\mu_0\bigr) f_2
\bigr) \Bigr)
=O\Bigl(\prod_{j=1}^{\ell}|z_j|^{-2N}
 \Bigr)
\]
Hence, we obtain Lemma \ref{lem;07.11.25.60}.
\hfill\qed

\vspace{.1in}

By the functoriality of the deformation
in Lemma \ref{lem;07.12.20.5},
we obtain a pairing
\[
 \nbigqzero\nbigs_0:
 \nbigq^{(\lambda_0)}\nbige
\otimes
 \sigma^{\ast}
 \nbigq^{(\sigma(\lambda_0))}
 \nbige^{\dagger}
\lrarr\nbigo_{\nbigxzero}
 \bigl(\ast \nbigd(\lambda_0)\bigr)
\]
for $\lambda_0\neq 0$,
which is a prolongment of $\nbigs_0$.
Note that such a prolongment is unique
if it exists.
Hence, we have 
$\nbigqzero\nbigs_{0|\nbigx^{(\lambda_1)}}
=\nbigq^{(\lambda_1)}\nbigs_0$
if $\nbigx^{(\lambda_1)}\subset 
 \nbigx^{(\lambda_0)}$.
By varying $\lambda_0$ and gluing them,
we obtain the desired prolongment
$\nbigq\nbigs_0$ of $\nbigs_0$.
Thus we obtain Proposition \ref{prop;07.11.24.29}.
\hfill\qed

\begin{df}
\label{df;08.9.16.30}
The prolongment
given by $\bigl(\nbigq\nbige,\nbigq\nbige^{\dagger}\bigr)$
is called the canonical prolongment of
$(\nbige^{\sankaku},\DD^{\sankaku},\nbigs)$.
\index{canonical prolongment}
\hfill\qed
\end{df}

\section{Reduction and uniqueness}
\label{subsection;07.11.24.50}
\subsection{Statement}

We set $X:=\Delta^n$,
$D_i:=\{z_i=0\}$
and $D:=\bigcup_{i=1}^{\ell}D_i$.
For any subset $\nbigk\subset\cnum_{\lambda}$,
we put $\nbigx_{\nbigk}:=
 \nbigk\times X$
and $\nbigx^{\dagger}_{\sigma(\nbigk)}:=
 \sigma(\nbigk)\times X^{\dagger}$.
For $0<R<1$,
we put
\[
 X(R):=\bigl\{(z_1,\ldots,z_n)\in X\,\big|\,
 |z_i|<R
 \bigr\},
\quad
X^{\ast}(R):=(X-D)\cap X(R).
\]

Let $(V,\DD^{\sankaku})$
be a variation of twistor structure
of weight $0$ on $\proj^1\times (X-D)$
equipped with a symmetric pairing 
$\nbigs:
 (V,\DD^{\sankaku})\otimes
 \sigma^{\ast}(V,\DD^{\sankaku})
 \lrarr \Tate(0)$.
Assume the following:
\begin{itemize}
\item
We are given
an unramifiedly good meromorphic prolongment
$\vecVtilde:=(\Vtilde_0,\Vtilde_{\infty})$
of $(V,\DD^{\sankaku},\nbigs)$
on $\bigl(\nbigx_{\nbigk},
 \nbigx^{\dagger}_{\sigma(\nbigk)}\bigr)$.
(See Section \ref{subsection;07.10.10.2}.)
Let $T$ denote the good set of irregular values
of the prolongment $(\Vtilde_0,\DD)$.
\item
$(\Vtilde_0,\DD)$ has the KMS-structure at 
each $\lambda_0\in\nbigk$,
i.e.,
there exists an increasing sequence of lattices
$\nbigqzero_{\veca} \Vtilde_0
\subset \Vtilde_{0}$ 
on a neighbourhood of $\{\lambda_0\}\times X$
indexed by $\veca\in\real^{\ell}$
such that
$(\nbigqzero_{\ast}\Vtilde_0,\DD)$
is a family of good filtered $\lambda$-flat bundles
with $KMS$-structure.
\item
Similarly, $(\Vtilde_{\infty},\DDtilde)$
has the KMS-structure
at each $\mu_0\in\sigma(\nbigk)$.
The corresponding family of filtered $\mu$-flat 
bundle is denoted by
$\nbigq^{(\mu_0)}_{\ast}\Vtilde_{\infty}$.
\end{itemize}
Recall that we have the reduction
$\Gr^{\vecVtilde}(V,\DD^{\sankaku},\nbigs)$
associated to the full Stokes structure.
(See Subsection \ref{subsection;08.9.14.52}.
Although it is denoted by 
$\Gr^{\nbigftilde}(V,\DD^{\sankaku},\nbigs)$
there,
we use the symbol
 $\Gr^{\vecVtilde}(V,\DD^{\sankaku},\nbigs)$
to emphasize the dependence on
the choice of $\vecVtilde$.)

\begin{rem}
If $(V,\DD^{\sankaku},\nbigs)$ as above is
a variation of polarized pure twistor structure
of weight $0$,
the underlying harmonic bundle
is unramifiedly good wild
by the above assumption,
because the Higgs field is given by
the restriction of $\DD^{\sankaku}$
to $\{0\}\times (X-D)$.
We will use it implicitly
in the subsequent argument.
\hfill\qed
\end{rem}

We will prove the following theorem in Sections
\ref{subsection;07.10.15.20}--\ref{subsection;08.9.25.15}.
\begin{thm}
\label{thm;08.9.16.100}
\label{thm;07.10.11.120}
\label{thm;07.10.11.10}
The following conditions are equivalent:
\begin{description}
\item[(P1)]
 $(V,\DD^{\sankaku},\nbigs)
 _{|\proj^1\times X^{\ast}(R)}$
 is a variation of polarized pure twistor structure
 of weight $0$
 for some $0<R<1$,
 and the prolongment 
 $\vecVtilde$ is canonical
 (Definition {\rm\ref{df;08.9.16.30}}).
\item[(P2)]
 The full reduction
 $\Gr^{\vecVtilde}(V,\DD^{\sankaku},\nbigs)
 _{|\proj^1\times X^{\ast}(R)}$
 is a variation of polarized pure twistor structure
 of weight $0$
 for some $0<R<1$.
\end{description}
\end{thm}

We give a consequence.
We take an auxiliary sequence
$\vecm(0),\ldots,\vecm(L)$
for the good set of irregular values $T$.
Recall that we have the reduction
$\Gr^{\vecm(p)}(V,\DD^{\sankaku},\nbigs)$
associated to the partial Stokes structure
in the level $\vecm(p)$.

\begin{cor}
\label{cor;08.9.16.101}
The following conditions are equivalent:
\begin{itemize}
\item
  $(V,\DD^{\sankaku},\nbigs)$
 is a variation of polarized pure twistor structure
 of weight $0$
 for some $0<R<1$,
 and the prolongment 
 $\vecVtilde$ is canonical.
\item
 The one step reduction
 $\Gr^{\vecm(0)}(V,\DD^{\sankaku},\nbigs)
 _{|\proj^1\times X^{\ast}(R)}$
 with respect to $\vecVtilde$
 is a variation of polarized pure twistor structure
 of weight $0$ 
 for some $0<R<1$,
 and the induced prolongment 
 $\Gr^{\vecm(0)}\vecVtilde:=
 (\Gr^{\vecm(0)}(\Vtilde_0),
 \Gr^{\vecm(0)}(\Vtilde_{\infty}))$ is canonical.
\end{itemize}
\end{cor}
\pf
According to Theorem \ref{thm;08.9.16.100},
both the conditions are equivalent to
the second condition in
Theorem \ref{thm;08.9.16.100}.
\hfill\qed

\subsection{Plan of the proof}
\label{subsection;08.9.25.1}

Before going into the proof,
we give a rough sketch of the argument.
In Section \ref{subsection;07.10.15.20},
we study the one step reduction,
and 
we will prove the following proposition,
which is a part of the claim of Corollary
\ref{cor;08.9.16.101}.
\begin{prop}
\mbox{{}}\label{prop;08.9.16.102}
The following holds:
\begin{description}
\item[(A)]
 Assume that $(V,\DD^{\sankaku},\nbigs)$ is
 a variation of polarized pure 
 twistor structure of weight $0$,
 and that $(\Vtilde_0,\Vtilde_{\infty})$ is canonical.
 Then, there exists $0<R<1$ such that
 $\Gr^{\vecm(0)}(V,\DD^{\sankaku},
 \nbigs)_{|\proj^1\times X^{\ast}(R)}$
 is a variation of polarized pure twistor structure
 of weight $0$.
\item[(B)]
 If $\Gr^{\vecm(0)}(V,\DD^{\sankaku},
 \nbigs)$
 is a variation of polarized pure twistor structure
 of weight $0$,
 there exists $0<R<1$ such that
 $(V,\DD^{\sankaku},\nbigs)
 _{|\proj^1\times X^{\ast}(R)}$ is 
 a variation of polarized pure twistor structure 
 of weight $0$.
\end{description}
\end{prop}

We remark that it is not so easy to
deduce a conclusion on prolongments.
For example, in the claim (A),
$\Gr^{\vecm(0)}(V,\DD^{\sankaku},\nbigs)$
comes from a harmonic bundle,
say $(E^{(1)},\delbar^{(1)},\theta^{(1)},h^{(1)})$
which is easily shown to be unramifiedly good wild.
In this stage,
$\Gr^{\vecm(0)}(V,\DD^{\sankaku},\nbigs)$
has two prolongments.
One is the induced prolongment
$\bigl(
 \Gr^{\vecm(0)}(\Vtilde_0),
 \Gr^{\vecm(0)}(\Vtilde_{\infty})
 \bigr)$.
The other is the canonical prolongment
$\bigl(\nbigq\nbige^{(1)},
 \nbigq\nbige^{(1)\dagger}
 \bigr)$
associated to 
$(E^{(1)},\delbar^{(1)},\theta^{(1)},h^{(1)})$.
For an inductive argument,
we would like to show that
they are the same.
We will show it eventually,
but it does not seem easy to show it directly.
Hence, we will postpone it.

\vspace{.1in}
In Section \ref{subsection;07.11.24.30},
we will study the full reduction
in the case that $D$ is smooth,
and we will show the following proposition.
\begin{prop}
\label{prop;07.9.10.3}
If $D$ is smooth,
the  claim of Theorem 
{\rm \ref{thm;08.9.16.100}} holds.
\end{prop}

Note that the prolongment of the full reduction
is easily controlled.
As in the previous case,
we have two prolongments.
One is the induced prolongment,
and the other is canonical.
But, the full reduction is graded
by the good set $T$,
and each graded piece has the unique
irregular value for both the prolongments.
Hence, the two prolongments have to be
the same by the uniqueness of
the Deligne extension of flat bundles.

\vspace{.1in}

Once we have established 
Theorem \ref{thm;08.9.16.100}
in the smooth divisor case,
it is not difficult to compare
the prolongments
$\nbigq\nbige^{(1)}$
and $\Gr^{\vecm(0)}(\Vtilde_0)$ for one step reduction
in the normal crossing divisor case.
Then, we obtain Theorem \ref{thm;08.9.16.100}
by an inductive argument,
which will be done in 
Section \ref{subsection;08.9.25.15}.

\section{Preliminary for convergence}
\label{subsection;08.9.29.14}
We will give a preliminary for convergence
of a sequence of twistor structures.

\subsection{Preliminary}
\label{subsection;07.12.23.20}

Let $V$ be an $n$-dimensional vector space
with a base $e_1,\ldots,e_n$.
Let $h$ be a hermitian pairing of $V$.
There exists a positive number $\epsilon_0$,
depending only on $n$,
with the following property:
\begin{itemize}
\item
If $\bigl|h(e_i,e_j)-\delta_{i,j}\bigr|<\epsilon_0$,
then $h$ is positive definite.
Here $\delta_{i,j}$  denotes $1$ $(i=j)$
or $0$ $(i\neq j)$.
\end{itemize}

Let $\sigma:\proj^1\lrarr\proj^1$ be given by
$\sigma\bigl([z_0:z_1]\bigr)
 =[-\zbar_1:\zbar_0]$.
\begin{lem}
\label{lem;07.12.23.22}
Let $V$ be a pure twistor structure of weight $0$
with a symmetric pairing
$S:V\otimes\sigma^{\ast}V\lrarr\Tate(0)$.
Let $\lambda$ be any point of $\proj^1$.
Assume the following :
\begin{itemize}
\item
 There exists a holomorphic frame 
 $w_1,\ldots,w_n$ of $V$ on $\proj^1$ such that 
\[
 \bigl|S_{|\lambda}\bigl(w_{i|\lambda},
 \sigma^{\ast}w_{j|\sigma(\lambda)}\bigr)
 -\delta_{i,j}\bigr|<\epsilon_0.
\]
\end{itemize}
Then, $S$ gives a polarization of $V$.
\end{lem}
\pf
Since any global section of $\Tate(0)=\nbigo_{\proj^1}$
is constant,
we have
\[
 S_{|\lambda}\bigl(w_{i|\lambda},
 \sigma^{\ast}w_{j|\sigma(\lambda)}\bigr)
=S(w_{i},\sigma^{\ast}w_j).
\]
Hence, the claim is clear
by the choice of $\epsilon_0$.
\hfill\qed

\subsection{Estimate on a bundle}

In this section,
we use the standard Fubini-Study metric of $\proj^1$.
Let $E$ be a holomorphic vector bundle on $\proj^1$
with a holomorphic base 
$\vecv=(v_1,\ldots,v_r)$,
i.e.,
$E\simeq\bigoplus_{i=1}^r
 \nbigo_{\proj^1}\, v_i$.
Let $h$ be the hermitian metric
given by $h(v_i,v_j)=\delta_{i,j}$.
Let $L^2_{1}(E)$ be the space of the sections $f$ of $E$
such that $f$ and 
$\delbar_E f$ are $L^2$ on $\proj^1$
with respect to $h$ and the Fubini-Study metric.
It is the $L^2$-space
with the hermitian product
$(f,g)_{L^2_1}:=\int h(f,g)+\int h(\delbar_E f,\delbar_E g)$.
Let $L^2(E\otimes\Omega^{0,1})$ denote the space
of $L^2$-sections of $E\otimes \Omega^{0,1}$.
We have the surjectivity
of $\delbar_E:L_{1}^2(E)\lrarr
 L^2(E\otimes\Omega^{0,1})$.
Let $H$ denote the orthogonal complement of
$\Ker(\delbar_E)$ in $L^2_{1}(E)$.
Then,
$\delbar_E:H\lrarr L^2(E\otimes\Omega^{0,1})$
is homeomorphic.
There exists a constant $C_1$ such that
$ C_1^{-1}\|u\|_{L^2}\leq \|\delbar_E u\|_{L^2}
 \leq C_1\|u\|_{L^2}$
for $u\in H$.

Let $A$ be a $C^{\infty}$-section of
$\End(E)\otimes\Omega^{0,1}$
such that $\sup|A|_h\leq \delta$.
Let $E_A$ denote the holomorphic bundle
$(E,\delbar_E+A)$.

\begin{lem}
\label{lem;07.7.7.70}
If $\delta$ is sufficiently small,
we have
$2^{-1}\|\delbar_E u\|_{L^2}
\leq \|(\delbar_E+A)u\|_{L^2}
\leq 2\|\delbar_E u\|_{L^2}$
for any $u\in H$.
Hence,
$\delbar_E+A:H\lrarr L^2(E\otimes\Omega^{0,1})$
is a homeomorphism.
In particular, $E_A$ is a pure twistor structure
of weight $0$.
\end{lem}
\pf
If $\delta$ is sufficiently small,
we have the following for any $u\in H$:
\[
 \|(\delbar_E+A)u\|_{L^2}\leq
 \|\delbar_E u\|_{L^2}+\|Au\|_{L^2}
 \leq 2\|\delbar_E u\|_{L^2}
\]
\[
 \|(\delbar_E+A)u\|_{L^2}\geq
 \|\delbar_E u\|_{L^2}-\|Au\|_{L^2}\geq 
 \frac{1}{2}\|\delbar_E u\|_{L^2}
\]
Then, the claim of the lemma immediately follows.
\hfill\qed

\begin{lem}
\label{lem;07.12.23.23}
There exist constants $\delta>0$ and $C_{10}>0$
with the following property:
\begin{itemize}
\item
 If $\sup|A|_h<\delta$,
 we can take $u_i\in H$ such that
 $(\delbar_E+A)u_i=(\delbar_E+A)v_i$
 and $\sup|u_i|_h<C_{10}\,\sup|A|_h$.
 Note $|v_i|_h=1$.
\end{itemize}
\end{lem}
\pf
It follows from Lemma \ref{lem;07.7.7.70}
and a standard bootstrapping argument.
\hfill\qed

\vspace{.1in}

Let $S_h:E\otimes\sigma^{\ast}E\lrarr\Tate(0)$ 
be the polarization of the pure twistor structure $E$
corresponding to $h$.
Let $d_{\Herm}$ denote the natural distance
of the space of the hermitian metrics of a vector space.

\begin{lem}
\label{lem;07.12.23.21}
Let $S:E_A\otimes\sigma^{\ast}E_A\lrarr\Tate(0)$
be a symmetric pairing.
Fix a point $\lambda$ of $\proj^1$,
and assume moreover 
the following for some $\eta>0$:
\begin{itemize}
\item
 $\bigl|S_{|\lambda}(u,\sigma^{\ast}v)-
 S_{h|\lambda}(u,\sigma^{\ast}v)\bigr|
<\eta\, |u|_h\, |v|_h$
for any $u\in E_{|\lambda}$ and $v\in E_{|\sigma(\lambda)}$.
\end{itemize}
If $\delta$ and $\eta$ are sufficiently small,
$(E_A,S)$ is a polarized pure twistor structure
of weight $0$.
Moreover, the following holds
for the metric $h_S$ corresponding to $S$:
\[
 \sup_{\mu\in\proj^1}
 d_{\Herm}(h_{S|\mu},h_{|\mu})
\leq
 C_{12}\, (\delta+\eta)
\]
Here, $C_{12}$ is a constant
depending only on $\rank E$.
\end{lem}
\pf
We may assume that $E_A$ is pure twistor structure
of weight $0$,
due to Lemma \ref{lem;07.7.7.70}.
Let $\epsilon_0$ be as 
in Section \ref{subsection;07.12.23.20}.
If $\delta$ is sufficiently small,
we can take sections $u_i\in H$
such that $(\delbar_E+A)u_i=(\delbar_E+A)v_i$
and $\sup|u_i|_h\leq C_{10}\,\delta<\epsilon_0/100$,
as remarked in Lemma \ref{lem;07.12.23.23}.
We put $w_i:=v_i-u_i$.
We may assume that 
$\vecw=(w_1\ldots,w_n)$
gives a holomorphic frame of $E_A$.
We have the following:
\begin{multline}
 S_{|\lambda}
 \bigl(w_{i|\lambda},
 \sigma^{\ast}w_{j|\sigma(\lambda)}\bigr)
-\delta_{i,j}
=S_{|\lambda}
 \bigl(w_{i|\lambda},
 \sigma^{\ast}w_{j|\sigma(\lambda)}\bigr)
-S_{h|\lambda}(v_{i|\lambda},
 \sigma^{\ast}v_{j|\sigma(\lambda)})\\
=\Bigl(
S_{|\lambda}
 \bigl(w_{i|\lambda},
 \sigma^{\ast}w_{j|\sigma(\lambda)}\bigr)
-S_{h|\lambda}\bigl(
 w_{i|\lambda},\sigma^{\ast}w_{j|\sigma(\lambda)}
 \bigr)
\Bigr) \\
+\Bigl(
S_{h|\lambda}\bigl(
 w_{i|\lambda},\sigma^{\ast}w_{j|\sigma(\lambda)}
 \bigr)
-S_{h|\lambda}(v_{i|\lambda},
 \sigma^{\ast}v_{j|\sigma(\lambda)})
\Bigr)
\end{multline}
We have the following:
\[
\Bigl|
 S_{|\lambda}
 \bigl(w_{i|\lambda},
 \sigma^{\ast}w_{j|\sigma(\lambda)}\bigr)
-S_{h|\lambda}\bigl(
 w_{i|\lambda},\sigma^{\ast}w_{j|\sigma(\lambda)}
 \bigr)
\Bigr|
\leq 4\eta
\]
We also have the following:
\[
\Bigl|
 S_{h|\lambda}\bigl(
 w_{i|\lambda},\sigma^{\ast}w_{j|\sigma(\lambda)}
 \bigr)
-S_{h|\lambda}\bigl(v_{i|\lambda},
 \sigma^{\ast}v_{j|\sigma(\lambda)}\bigl)
\Bigr|
\leq 4C_{10}\,\delta
\]
Hence, 
we obtain 
$\Bigl|
 S_{|\lambda}
 \bigl(w_{i|\lambda},
 \sigma^{\ast}w_{j|\sigma(\lambda)}\bigr)
-\delta_{i,j}\Bigr|\leq
 \epsilon_0$,
and
the first claim of Lemma \ref{lem;07.12.23.21}
follows from Lemma \ref{lem;07.12.23.22}.

Let us show the second claim.
We have 
$\Bigl|
 S_{|\lambda}
 \bigl(w_{i|\lambda},
 \sigma^{\ast}w_{j|\sigma(\lambda)}\bigr)
-\delta_{i,j}\Bigr|\leq C_{20}\,(\delta+\eta)$
due to the above argument.
Hence, there exists a constant matrix $B_1$
with the following property:
\begin{itemize}
\item
The frame $\vecwtilde=\vecw(I+B_1)$
is orthonormal with respect to $h_S$.
\item
$|B_1|<C_{21}\,(\delta+\eta)$.
\end{itemize}
Let $B_0$ be determined by
$\vecw=\vecv\, (I+B_0)$.
By construction,
we have $\sup|B_0|\leq C_{20}\,\delta$.
Therefore,
we have $\sup|B_2|\leq C_{22}\, (\delta+\eta)$,
where $B_2$ is determined by
$\vecwtilde=\vecv\, (I+B_2)$.
Since $\vecwtilde$ and $\vecv$ are
orthonormal frames with respect to 
$h_S$ and $h$ respectively,
we obtain the second claim of 
Lemma \ref{lem;07.12.23.21}.
\hfill\qed

\subsection{Estimate through a $C^{\infty}$-map}

Let $(V^{(1)},S^{(1)})$
be a polarized pure twistor structure.
Let $h^{(1)}$ be the corresponding metric
of $H^0(\proj^1,V^{(1)})$.
It induces a hermitian metric of $V^{(1)}$,
which is also denoted by $h^{(1)}$.

Let $V^{(2)}$ be a holomorphic vector bundle
on $\proj^1$.
Let $F:V^{(1)}\lrarr V^{(2)}$ be a $C^{\infty}$-isomorphism.
Assume the following for some $\delta_1>0$:
\[
\sup_{\proj^1} \bigl|F^{-1}\circ\delbar_{V^{(2)}}\circ F
-\delbar_{V^{(1)}}\bigr|_{h^{(1)}}
\leq \delta_1
\]
\begin{lem}
\label{lem;07.12.23.40}
If $\delta_1$ is sufficiently small,
$V^{(2)}$ is also a pure twistor structure.
\end{lem}
\pf
It follows from Lemma \ref{lem;07.7.7.70}.
\hfill\qed

\vspace{.1in}

Let $S^{(2)}$ be a symmetric pairing of $V^{(2)}$,
and let $F$ be as above.
Fix a point $\lambda$ of $\proj^1$.
Assume moreover that there exists a $\delta_2>0$
such that the following holds
for any $u\in V^{(1)}_{|\lambda}$
and $v\in V^{(1)}_{|\sigma(\lambda)}$:
\[
\bigl|
 S^{(2)}\bigl(F(u),\sigma^{\ast}(F(v))\bigr)
-S^{(1)}(u,\sigma^{\ast}v)
\bigr|
\leq
 \delta_2\,|u|_{h^{(1)}}\, 
 |v|_{h^{(1)}}
\]
\begin{lem}
\mbox{{}} \label{lem;07.6.2.90}
If $\delta_i$ $(i=1,2)$ are sufficiently small,
the following holds:
\begin{itemize}
\item
$(V^{(2)},S^{(2)})$ is a polarized pure twistor structure.
Let $h^{(2)}$ denote the corresponding hermitian metric.
\item
We have
$\sup_{\proj^1}
 d_{\Herm}\bigl(h^{(1)},F^{\ast}h^{(2)}\bigr)
 \leq C\, (\delta_1+\delta_2)$,
where $C$ is a constant depending only on the rank.
\end{itemize}
\end{lem}
\pf
It follows from Lemma \ref{lem;07.12.23.21}.
\hfill\qed

\section{One step reduction}
\label{subsection;07.10.15.20}
\subsection{Preliminary}

For simplicity,
we assume that the coordinate system is admissible
for the good set $T$.
Let $k$ be the number determined by
$\vecm(0)\in\seisuu_{<0}^k\times\veczero_{\ell-k}$.
By shrinking $X$,
we may assume that 
$\Gr^{\vecm(0)}(\Vtilde_0)$
and $\Gr^{\vecm(0)}(\Vtilde_{\infty})$
are given on 
$\nbigx_{\nbigk}$ and 
$\nbigx^{\dagger}_{\sigma(\nbigk)}$,
respectively.
Let $d''_{\proj^1,V}$
and $d''_{\proj^1,\Gr^{\vecm(0)}(V)}$
denote the $\proj^1$-holomorphic structures
of $V$ and $\Gr^{\vecm(0)}(V)$,
respectively.

In the following, we assume one of the following:
\begin{description}
\item[(Case A)]
$(V,\DD^{\sankaku},\nbigs)$ 
is a variation of polarized pure twistor structure
of weight $0$,
and the prolongment $\vecVtilde$ is canonical.
The underlying harmonic bundle
is denoted by $(E,\delbar_E,\theta,h)$.
\item[(Case B)]
$\Gr^{\vecm(0)}(V,\DD^{\sankaku},\nbigs)$
is a variation of polarized pure twistor structure
of weight $0$.
The underlying harmonic bundle
is denoted by
 $(E^{(1)},\delbar^{(1)},\theta^{(1)},h^{(1)})$.
It is graded as
\[
 (E^{(1)},\delbar^{(1)},\theta^{(1)},h^{(1)})
=\bigoplus
 _{\gminia\in T(\vecm(0))}
 (E^{(1)}_{\gminia},
 \delbar_{\gminia}^{(1)},
 \theta_{\gminia}^{(1)},
 h_{\gminia}^{(1)}),
\]
corresponding to the decomposition
$\Gr^{\vecm(0)}(V,\DD^{\sankaku},\nbigs)
=\bigoplus
 \Gr^{\vecm(0)}_{\gminia}(V,\DD^{\sankaku},\nbigs)$.
\end{description}
We have some remarks in Case (B),
which will be implicitly used
in the subsequent argument.
We have the canonical prolongment
$(\nbigq\nbige^{(1)},\DD^{(1)})
=\bigoplus
 \bigl(\nbigq\nbige^{(1)}_{\gminia},
 \DD^{(1)}_{\gminia}
 \bigr)$.
Note that we should distinguish
$\Gr^{\vecm(0)}(\Vtilde_0)$ and
$\nbigq\nbige^{(1)}$.
However, we have the natural isomorphism
\[
 \Gr^{\vecm(0)}(\Vtilde_0)_{|\cnum_{\lambda}
 \times (X-D(\leq k))}
\simeq 
 \nbigq\nbige^{(1)}_{|\cnum_{\lambda}
 \times (X-D(\leq k))},
\]
because both of them are regular
along $D_i$ for $i>k$.
The induced KMS-structures at any $\lambda_0$
are also the same due to the uniqueness 
of KMS-structure.
We also remark that
the grading naturally gives
a splitting of the Stokes filtration
in the level $\vecm(0)$
for each small multi-sector.

\subsection{Preliminary estimate}

Let $\lambda_0\in \nbigk$,
and let $U_{\lambda}(\lambda_0)$ 
denote a small neighbourhood of $\lambda_0$.
If $\lambda_0\neq 0$,
we assume $0\not\in U_{\lambda}(\lambda_0)$.
We will shrink it in the following argument without mention.
We set $W:=U_{\lambda}(\lambda_0)\times D(\leq k)$
if $\lambda_0\neq 0$,
and $W:=(U_{\lambda}(\lambda_0)\times D)
 \cup (\{0\}\times X)$ if $\lambda_0=0$.
Let $S$ be a multi-sector in 
$(U_{\lambda}(\lambda_0)\times X)-W$,
and let $\Sbar$ denote the closure of $S$
in the real blow up of $U_{\lambda}(\lambda_0)\times X$
along $W$.
If $S$ is sufficiently small,
we can take a $\DD_{\leq k}$-flat splitting 
$\nbigqzero_0\Vtilde_{0|\Sbar}=\bigoplus_{\gminia}
 \nbigqzero_0\Vtilde_{0,\gminia,S}$
of the partial Stokes filtration $\nbigf^S$ 
in the level $\vecm(0)$,
whose restrictions to 
$S\cap \bigl(U_{\lambda}(\lambda_0)\times D_i\bigr)$
is compatible with the endomorphism
$\Res_{i}(\DD)$
and the filtration $\lefttop{i}\Fzero$
for each $i=k+1,\ldots,\ell$.
(See Proposition \ref{prop;07.9.30.2}.)
The splitting induces
a $\DD_{\leq k}$-flat isomorphism on $S$:
\begin{equation}
 \label{eq;08.9.16.120}
 g_{S}:
 \Gr^{\vecm(0)}(\nbigqzero_0\Vtilde_0)_{|S}
\simeq
 \nbigqzero_0\Vtilde_{0|S}
\end{equation}
Its inverse is denoted by $f_S$.

\begin{lem}
\mbox{{}}\label{lem;07.7.8.10}
Let $g_S'$ and $f_S'$ be obtained from
another $\DD_{\leq k}$-flat splitting 
of $\nbigf^S$.
If $S$ is shrinked in the radius direction,
the following holds:
\begin{description}
\item[(Case A)]
 $\id-f_S^{-1}\circ f_S'=O\Bigl(
 \exp\bigl(
 -\epsilon|\lambda^{-1}\vecz^{\vecm(0)}|\bigr)
 \Bigr)$ for some $\epsilon>0$
 with respect to $h$.
\item[(Case B)]
 $\id-g_S^{-1}\circ g_S'=O\Bigl(
 \exp\bigl(-\epsilon|\lambda^{-1}\vecz^{\vecm(0)}|\bigr)
 \Bigr)$ for some $\epsilon>0$
 with respect to $h^{(1)}$.
\end{description}
\end{lem}
\pf
Let us consider Case (A).
Note $\nbigqzero_0\Vtilde_0=\nbigqzero_0\nbige$ in this case.
The $\DD_{\leq k}$-flat endomorphism 
$\gbigf:=\id-f_S^{-1}\circ f_S'$
strictly decreases the Stokes filtration of
$\nbigqzero_0\nbige_{|\Sbar}$ 
in the level $\vecm(0)$.
It is compatible with
the filtrations $\lefttop{i}\Fzero$ 
and the endomorphisms $\Res_i(\DD)$ on
$S\cap \bigl(U_{\lambda}(\lambda_0)\times D_i\bigr)$
$(i=k+1,\ldots,\ell)$.
If $U_{\lambda}(\lambda_0)$ is shrinked,
the Stokes filtrations of
$\nbigqzero_0\nbige_{|\Sbar}$
and $\nbigpzero_0\nbige_{|\Sbar}$
are the same.
By shrinking $S$ in the radius direction
we obtain the desired estimate
due to Lemma \ref{lem;07.6.2.60}.
We can argue Case (B)
in a similar way.
\hfill\qed

\vspace{.1in}

Let $g_S^{(p)}$ $(p=1,\ldots,L)$ be
$\DD_{\leq k}$-flat isomorphisms
as in (\ref{eq;08.9.16.120}).
Let $f_S^{(p)}$ denote their inverses.
Let $\chi_p$ $(p=1,\ldots,L)$ be non-negative valued
$C^{\infty}$-functions on $S$
such that
(i) $\sum \chi_p=1$,
(ii) $\delbar \chi_p$
 are of polynomial order
 in $|\lambda^{-1}|$
 and $|z_i^{-1}|$ for $i=1,\ldots,k$.
We set 
$g:=\sum \chi_p\, g_p$ and
$f:=\sum \chi_p\, f_p$.

\begin{lem}
\label{lem;08.9.16.132}
In Case (A),
we have the following estimate
with respect to $h$:
\begin{equation}
 \label{eq;08.9.16.130}
  f^{-1}\circ d''_{\proj^1,\Gr^{\vecm(0)}(V)}\circ f
-d''_{\proj^1,V}
=O\Bigl(
 \exp\bigl(-\epsilon|\lambda^{-1}\vecz^{\vecm(0)}|\bigr)
 \Bigr)
\end{equation}
In Case (B),
we have the following estimate
with respect to $h^{(1)}$:
\begin{equation}
 \label{eq;08.9.16.131}
 g^{-1}\circ d''_{\proj^1,V}\circ g
-d''_{\proj^1,\Gr^{\vecm(0)}(V)}
=O\Bigl(
 \exp\bigl( 
 -\epsilon|\lambda^{-1}\vecz^{\vecm(0)}|
 \bigr)
 \Bigr)
\end{equation}
\end{lem}
\pf
Let us show (\ref{eq;08.9.16.130})
in Case (A).
We remark the following estimate
with respect to $h$,
due to Lemma \ref{lem;07.7.8.10}:
\[
 (f_S^{(p)})^{-1}\circ f-\id
=\sum \chi_q\, 
\bigl((f_S^{(p)})^{-1}\circ f^{(q)}_S-\id\bigr)
 =O\Bigl(
 \exp\bigl(
 -\epsilon|\lambda^{-1}\vecz^{\vecm(0)}|
 \bigr)
 \Bigr)
\]
The left hand side of (\ref{eq;08.9.16.130})
can be rewritten as follows:
\[
 f^{-1}\circ\bigl(
 \delbar_{\proj^1,\Gr^{\vecm(0)}(V)}\circ f
-f\circ\delbar_{\proj^1,V}
 \bigr)
=f^{-1}\circ\Bigl(
 \sum \delbar\chi_p\,
 f_{S}^{(p)}
 \Bigr)
=\sum \delbar\chi_p\,
 \bigl(
 f^{-1}\, f_S^{(p)}-\id
 \bigr)
\]
Hence, we obtain the desired estimate
(\ref{eq;08.9.16.130}).
We can show (\ref{eq;08.9.16.131})
in Case (B)
by the same argument.
\hfill\qed

\subsubsection{The other side}

Let $U_{\mu}(\mu_0)$ denote 
a small neighbourhood of $\mu_0\in\sigma(\nbigk)$.
We set
 $W^{\dagger}:=
 \bigl(\{0\}\times X^{\dagger}\bigr)
 \cup \bigl(U_{\mu}(0)\times 
 D^{\dagger}(\leq k)\bigr)$.
If we are given a small multi-sector $S$ of
$\bigl(U_{\mu}(\mu_0)\times 
 X^{\dagger}\bigr)-W^{\dagger}$,
we can take a $\DD^{\dagger}_{\leq k}$-flat isomorphism
\begin{equation}
 \label{eq;08.9.16.121}
 g^{\dagger}_{S}:
 \Gr^{\vecm(0)}(
 \nbigq^{(\mu_0)}_0\Vtilde_{\infty})_{|S}
\simeq
 \nbigq^{(\mu_0)}_0\Vtilde_{\infty|S}
\end{equation}
in a similar way,
by taking a $\DD^{\dagger}_{\leq k}$-flat splitting
of the Stokes filtration in the level $\vecm(0)$.
Its inverse is denoted by $f^{\dagger}_S$.
We have estimates for them
similar to Lemma \ref{lem;07.7.8.10}
and Lemma \ref{lem;08.9.16.132}.

\subsubsection{Remark for gluing}

Let $\lambda_0\neq 0$
and $\mu_0=\lambda_0^{-1}$.
Let $U_{\lambda}(\lambda_0)$ be a small
neighbourhood of $\lambda_0$ in $\cnum_{\lambda}$,
and let $U_{\mu}(\mu_0)$ be 
the corresponding neighbourhood of
$\mu_0$ in $\cnum_{\mu}$.
Assume that $\lambda_0$ is generic
with respect to the KMS-structure of $\Vtilde_0$.
If $S$ is a sufficiently small multi-sector 
in $U_{\lambda}(\lambda_0)\times
 \bigl(X-D(\leq k)\bigr)$,
we can take a $\DD$-flat splitting
$\nbigqzero_0\Vtilde_{0|\Sbar}
=\bigoplus_{\gminia}
 \nbigqzero_0\Vtilde_{0,\gminia,S}$.
(See Proposition \ref{prop;07.9.30.5}.)
Note that $S$ naturally gives
a small multi-sector
of $U_{\mu}(\mu_0)\times 
 \bigl(X^{\dagger}-D^{\dagger}(\leq k)\bigr)$,
and that the splitting naturally induces
a $\DD^{\dagger}$-splitting
$\nbigq^{(\mu_0)}_0\Vtilde_{\infty|\Sbar}
=\bigoplus_{\gminia}
 \nbigq^{(\mu_0)}_0\Vtilde_{\infty,\gminia,S}$.

\subsubsection{Estimate for pairing}

Let $S$ be a sufficiently small multi-sector of
$\bigl(U_{\lambda}(\lambda_0)\times X\bigr)-W$.
Let $g_S$ be as in (\ref{eq;08.9.16.120}),
and let $f_S$ denote its inverse.
We obtain the small multi-sector $\sigma(S)$
of $\sigma\bigl(U_{\lambda}(\lambda_0)\bigr)
 \times \bigl(X^{\dagger}-D^{\dagger}(\leq k)\bigr)$.
Let $g^{\dagger}_{\sigma(S)}$
be as in (\ref{eq;08.9.16.121}) for $\sigma(S)$,
and let $f^{\dagger}_{\sigma(S)}$ denote its inverse.
\begin{lem}
\label{lem;08.9.16.200}
In Case (A),
we have the following estimate
with respect to $h$:
\begin{equation}
 \label{eq;08.9.16.150}
 \nbigs-\Gr^{\vecm(0)}\nbigs\circ
 \bigl(f_S\otimes 
 \sigma^{\ast}f_{\sigma(S)}^{\dagger}
 \bigr)
=O\Bigl(
 \exp\bigl(-\epsilon|\lambda^{-1}\vecz^{\vecm(0)}|\bigr)
 \Bigr)
\end{equation}
In Case (B),
we have the following estimate
with respect to $h^{(1)}$:
\begin{equation}
 \label{eq;08.9.16.151}
 \Gr^{\vecm(0)}\nbigs-\nbigs\circ
 \bigl(g_S\otimes 
 \sigma^{\ast}g^{\dagger}_{\sigma(S)}
 \bigr)
=O\Bigl(
 \exp\bigl(-\epsilon|\lambda^{-1}\vecz^{\vecm(0)}|\bigr)
 \Bigr)
\end{equation}
\end{lem}
\pf
By the perfect pairings,
we have the natural isomorphisms
\[
 \Vtilde_{\infty}\simeq
 \sigma^{\ast}(\Vtilde_0^{\lor}),
\quad
 \Gr^{\vecm(0)}(\Vtilde_{\infty})
\simeq
 \sigma^{\ast}\bigl(
 \Gr^{\vecm(0)}(\Vtilde_0)^{\lor}
 \bigr).
\]
Hence, $g_S$ induces 
$\sigma^{\ast}(g_S^{\lor})$
as in (\ref{eq;08.9.16.121})
for $\sigma(S)$.
If $g^{\dagger}_{\sigma(S)}=
 \sigma^{\ast}\big(g_S^{\lor})$,
both (\ref{eq;08.9.16.150})
and (\ref{eq;08.9.16.151}) vanish.

Let us show (\ref{eq;08.9.16.150}) 
in Case (A).
We set 
$ \gbigh:=\Gr^{\vecm(0)}\nbigs\circ\bigl(
 f_S\otimes (\sigma^{\ast}f^{\dagger}_{\sigma(S)}
 -f_S^{\lor})
 \bigr)$.
Note that $\nbigs$ gives the identification
$\sigma^{\ast}\nbigq\nbige^{\dagger}
 \simeq
 \nbigq\nbige^{\lor}$,
and we can regard $\gbigh$
as a $\DD_{\leq k}$-flat section of
$\nbigpzero\End(\nbigelambda)$
such that
(i) it is compatible with the filtrations
$\lefttop{i}F$ and the endomorphisms $\Res_i(\DD)$ on
$S\cap \bigl(U_{\lambda}(\lambda_0)\times D_i\bigr)$
$(i=k+1,\ldots,\ell)$,
(ii) $\gbigh$ strictly decreases the Stokes filtration
of $\nbigpzero\nbige_{|\Sbar}$.
Hence, we obtain the desired estimate 
from Lemma \ref{lem;07.6.2.60}.
We can show (\ref{eq;08.9.16.151})
in Case (B) by a similar argument.
\hfill\qed

\vspace{.1in}

Let $g_S^{(p)}$ $(p=1,\ldots,L)$ be 
as in (\ref{eq;08.9.16.120}),
and  $f_S^{(p)}$ denote their inverses.
Let $\chi_{S,p}$ $(p=1,\ldots,L)$ be non-negative valued
functions such that
(i) $\sum \chi_{S,p}=1$,
(ii) $\delbar \chi_{S,p}$ are of polynomial order
 in $|\lambda^{-1}|$ and $|z_i^{-1}|$ for $i=1,\ldots,k$.
We set
$g:=\sum \chi_p\, g^{(p)}_S$
and $f:=\sum \chi_p\, f^{(p)}_S$.
Let $g^{\dagger\,(q)}_{\sigma(S)}$ 
$(q=1,\ldots,M)$ be as in
(\ref{eq;08.9.16.121}) for $\sigma(S)$,
and $f^{\dagger\,(q)}_{\sigma(S)}$ 
denote their inverses.
Let $\chi_{\sigma(S),q}$ $(q=1,\ldots,M)$ 
be non-negative valued
functions such that
(i) $\sum \chi_{\sigma(S),q}=1$,
(ii) $\delbar \chi_{\sigma(S),q}$ are 
 of polynomial order
 in $|\mu^{-1}|$ and $|z_i^{-1}|$ for $i=1,\ldots,k$.
We set 
$g^{\dagger}:=\sum \chi_{\sigma(S),q}\, 
 g^{\dagger(q)}_{\sigma(S)}$
and
$f^{\dagger}:=\sum \chi_{\sigma(S),q}\, 
 f^{\dagger(q)}_{\sigma(S)}$.
\begin{lem}
\label{lem;08.9.17.51}
In Case (A),
we have the following estimate
with respect to $h$:
\begin{equation}
 \label{eq;08.9.16.140}
 \nbigs-\Gr^{\vecm(0)}\nbigs
 \circ(f\otimes \sigma^{\ast}f^{\dagger})
=O\Bigl(
 \exp\bigl(-\epsilon|\lambda^{-1}\vecz^{\vecm(0)}|\bigr)
 \Bigr)
\end{equation}
In Case (B),
we have the following estimate
with respect to $h^{(1)}$:
\begin{equation}
 \label{eq;08.9.16.141}
 \Gr^{\vecm(0)}\nbigs
 -\nbigs\circ(g\otimes \sigma^{\ast}g^{\dagger})
=O\Bigl(
 \exp\bigl(-\epsilon|\lambda^{-1}\vecz^{\vecm(0)}|\bigr)
 \Bigr)
\end{equation}
\end{lem}
\pf
It follows from Lemma
\ref{lem;08.9.16.200}.
\hfill\qed

\subsection{Proof of Proposition \ref{prop;08.9.16.102}}
\label{subsection;07.9.11.2}

By shrinking $\nbigk$,
we may assume that
any $\lambda\in\sigma(\nbigk)\cap \nbigk$
is generic
with respect to the KMS-structure of
$\Vtilde_0$ and $\Vtilde_{\infty}$.
We set 
$W_{\nbigk}:=\bigl(\nbigk\times D(\leq k)\bigr)
 \cup\bigl(\{0\}\times X\bigr)$.
We take a finite covering
$\nbigx_{\nbigk}-W_{\nbigk}
\subset \bigcup_{p=1}^L S_p$
by small multi-sectors $S_p$
satisfying the following:
\begin{itemize}
\item
On each $S_p$, we can take
a $\DD_{\leq k}$-flat morphism
 $g_{S_p}$ as in  (\ref{eq;08.9.16.120}).
Its inverse is denoted by $f_{S_p}$.
\item
On each $\sigma(S_p)$,
we can take a $\DD_{\leq k}^{\dagger}$-flat
 morphism $g^{\dagger}_{\sigma(S_p)}$
as in (\ref{eq;08.9.16.121}).
Its inverse is denoted by 
 $f_{\sigma(S_p)}^{\dagger}$.
\item
If $S_p\cap\nbigk\cap\sigma(\nbigk)\neq\emptyset$,
we assume that
$g_{S_p}$ comes from a $\DD$-flat splitting.
We assume a similar condition
for $g^{\dagger}_{\sigma(S_p)}$.
\end{itemize}
We take a partition of unity
$\bigl(
 \chi_{S_p},\chi_{\sigma(S_p)}\,\big|\,
 p=1,\ldots,L
 \bigr)$ subordinated to
the covering
$(S_p,\sigma(S_p)\,|\,p=1,\ldots,L)$
such that
 $\del \chi_{S_p}$ 
 (resp. $\del \chi_{\sigma(S_p)}$)
 are of polynomial order
 in $|\lambda|^{-1}$ (resp. $|\lambda|$)
 and
 $|z_i|^{-1}$ for $i=1,\ldots,k$.
We set 
\[
 g:=\sum_p \chi_{S_p}\, g_{S_p}
+\sum_p\chi_{\sigma(S_p)}
 \, g^{\dagger}_{\sigma(S_p)},
\quad\quad
  f:=\sum_p \chi_{S_p}\, f_{S_p}
+\sum_p\chi_{\sigma(S_p)}
 \, f^{\dagger}_{\sigma(S_p)}.
\]
Let us consider the case (A).
By Lemma \ref{lem;08.9.16.132}
and its analogue for the $\mu$-side,
we obtain the following with respect to $h$:
\[
  f^{-1}\circ d''_{\proj^1,\Gr^{\vecm(0)}(V)}\circ f
-d''_{\proj^1,V}
=O\Bigl(
 \exp\bigl(-\epsilon\,
 (|\lambda^{-1}|+|\lambda|)\,
 |\vecz^{\vecm(0)}|\bigr)
 \Bigr)
\]
By Lemma \ref{lem;08.9.17.51},
we have the following estimate
with respect to $h$:
\[
  \nbigs-\Gr^{\vecm(0)}\nbigs\circ
 (f\otimes \sigma^{\ast}f^{\dagger})
=O\Bigl(
 \exp\bigl(-\epsilon\,
 (|\lambda^{-1}|+|\lambda|)\,
 |\vecz^{\vecm(0)}|\bigr)
 \Bigr)
\]
Then, there exists $0<R<1$ such that
$\Gr^{\vecm(0)}(V,\DD^{\sankaku},\nbigs)$
is a variation of polarized pure twistor structure,
due to Lemma \ref{lem;07.12.23.40}
and Lemma \ref{lem;07.6.2.90}.
Hence, the claim (A) of Proposition
\ref{prop;08.9.16.102} is proved.
The claim (B) can be shown 
in a similar way.
\hfill\qed

\section{Full reduction in the smooth divisor case}
\label{subsection;07.11.24.30}
In this section, we assume $D=D_1$,
and we will prove Proposition \ref{prop;07.9.10.3}.
The argument is
essentially the same as that in the proof of 
Proposition  \ref{prop;08.9.16.102}.
We almost repeat it
by changing Lemma \ref{lem;07.6.2.60}
with Proposition \ref{prop;08.9.17.1}.

Let $d''_{\proj^1,V}$
and $d''_{\proj^1,\Gr^{\vecVtilde}(V)}$
denote the $\proj^1$-holomorphic structure
of $V$ and $\Gr^{\vecVtilde}(V)$.
In the following, we assume one of the following:
\begin{description}
\item[(Case A)]
$(V,\DD^{\sankaku},\nbigs)$ 
is a variation of polarized pure twistor structure 
of weight $0$,
and the prolongment $\vecVtilde$ is canonical.
The underlying harmonic bundle
is denoted by $(E,\delbar_E,\theta,h)$.
\item[(Case C)]
$\Gr^{\vecVtilde}(V,\DD^{\sankaku},\nbigs)$
is a variation of polarized pure twistor structure
of weight $0$.
The underlying harmonic bundles
are denoted by
 $(E^{(0)},\delbar^{(0)},\theta^{(0)},h^{(0)})$.
It is graded as
\[
 (E^{(0)},\delbar^{(0)},\theta^{(0)},h^{(0)})
=\bigoplus
 _{\gminia\in T}
 (E^{(0)}_{\gminia},
 \delbar_{\gminia}^{(0)},
 \theta_{\gminia}^{(0)},
 h_{\gminia}^{(0)}),
\]
corresponding to
$\Gr^{\vecVtilde}(V,\DD^{\sankaku},\nbigs)
=\bigoplus
 \Gr^{\vecVtilde}_{\gminia}(V,\DD^{\sankaku},\nbigs)$.
\end{description}
We have some remarks in Case (C),
which we will implicitly use in the subsequent argument.
We have the canonical prolongment
$(\nbigq\nbige^{(0)},\DD^{(0)})
=\bigoplus
 \bigl(\nbigq\nbige^{(0)}_{\gminia},
 \DD^{(0)}_{\gminia}
 \bigr)$.
Since 
$\Gr^{\vecVtilde}_{\gminia}(\Vtilde_0)$
and 
$\nbigq\nbige^{(0)}_{\gminia}$
have the unique irregular value $\gminia$,
the natural isomorphism 
on $\cnum_{\lambda}\times (X-D)$
is extended to the isomorphism
on $\cnum_{\lambda}\times X$
by the uniqueness of the Deligne extension
of flat bundles.
We also remark that the grading
gives a splitting of the full Stokes filtration
for each small multi-sector.

\subsection{Preliminary estimate}

Let $\lambda_0\in \nbigk$,
and let $U_{\lambda}(\lambda_0)$ 
denote a small neighbourhood of $\lambda_0$.
We will shrink it in the following argument without mention.
We set $W:=U_{\lambda}(\lambda_0)\times D$
if $\lambda_0\neq 0$,
and $W:=
 (\{0\}\times X)\cup
 \bigl(U_{\lambda}(0)\times D\bigr)$
if $\lambda_0=0$.
Let $S$ be a multi-sector in 
$\bigl(U_{\lambda}(\lambda_0)\times X\bigr)-W$,
and let $\Sbar$ denote the closure of $S$
in the real blow up of $U_{\lambda}(\lambda_0)\times X$
along $W$.
If $S$ is sufficiently small,
we can take a $\DD$-flat splitting 
$\nbigqzero_0\Vtilde_{0|\Sbar}=\bigoplus_{\gminia\in T}
 \nbigqzero_0\Vtilde_{0,\gminia,S}$
of the full Stokes filtration $\nbigftilde^S$.
The splitting induces a $\DD$-flat isomorphism
\begin{equation}
 \label{eq;08.9.16.320}
 g_{S}:
 \Gr^{\nbigftilde}(\nbigqzero_0\Vtilde_0,\DD)_{|S}
\simeq
 (\nbigqzero_0\Vtilde_0,\DD)_{|S}.
\end{equation}
Its inverse is denoted by $f_S$.
We can show the following lemma
by using Proposition \ref{prop;08.9.17.1}
and the argument in the proof of
Lemma \ref{lem;07.7.8.10}.

\begin{lem}
\mbox{{}}\label{lem;07.10.11.100}
Let $g_S'$ and $f_S'$ be obtained from
another $\DD$-flat splitting of $\nbigf^S$.
If $S$ is shrinked in the radius direction,
the following holds:
\begin{description}
\item[(Case A)]
 $\id-f_S^{-1}\circ f_S'=O\Bigl(
 \exp\bigl(-\epsilon|\lambda^{-1}z_1^{-1}|\bigr)
 \Bigr)$ for some $\epsilon>0$
 with respect to $h$.
\item[(Case C)]
 $\id-g_S^{-1}\circ g_S'=O\Bigl(
 \exp\bigl(-\epsilon|\lambda^{-1}z_1^{-1}|\bigr)
 \Bigr)$ for some $\epsilon>0$
 with respect to $h^{(0)}$.
\hfill\qed
\end{description}
\end{lem}

Let $g_S^{(p)}$ $(p=1,\ldots,L)$ be
$\DD$-flat isomorphisms
as in (\ref{eq;08.9.16.320}).
Let $f_S^{(p)}$ denote their inverses.
Let $\chi_p$ $(p=1,\ldots,L)$ be non-negative valued
$C^{\infty}$-functions on $S$
such that
(i) $\sum \chi_p=1$,
(ii) $\delbar \chi_p$
 are of polynomial order in $|\lambda^{-1}|$
 and $|z_1^{-1}|$.
We set 
$g:=\sum \chi_p\, g_p$ and
$f:=\sum \chi_p\, f_p$.

\begin{lem}
\label{lem;08.9.16.332}
In Case (A),
we have the following estimate
with respect to $h$:
\begin{equation}
 \label{eq;08.9.16.330}
  f^{-1}\circ d''_{\proj^1,\Gr^{\vecVtilde}(V)}\circ f
-d''_{\proj^1,V}
=O\Bigl(
 \exp\bigl(-\epsilon|z_1^{-1}|\bigr)
 \Bigr)
\end{equation}
In Case (C),
we have the following estimate
with respect to $h^{(0)}$:
\begin{equation}
 \label{eq;08.9.16.331}
 g^{-1}\circ d''_{\proj^1,V}\circ g
-d''_{\proj^1,\Gr^{\vecVtilde}(V)}
=O\Bigl(
 \exp\bigl(-\epsilon|z_1^{-1}|\bigr)
 \Bigr)
\end{equation}
\end{lem}
\pf
It follows from 
Lemma \ref{lem;07.10.11.100}.
\hfill\qed

\subsubsection{The other side}

Let $\mu_0\in\sigma(\nbigk)-\{0\}$,
and let $U_{\mu}(\mu_0)$ denote 
a small neighbourhood of $\mu_0$.
If we are given a small multi-sector $S$ of
$U_{\mu}(\mu_0)\times 
 \bigl(X^{\dagger}-D^{\dagger}\bigr)$,
we can take a $\DD^{\dagger}$-flat isomorphism
\begin{equation}
 \label{eq;08.9.16.321}
 g^{\dagger}_{S}:
 \Gr^{\nbigftilde}(
 \nbigq^{(\mu_0)}_0\Vtilde_{\infty},
 \DD^{\dagger})_{|S}
\simeq
 (\nbigq^{(\mu_0)}_0\Vtilde_{\infty},
 \DD^{\dagger})_{|S}
\end{equation}
in a similar way,
by taking a $\DD^{\dagger}$-flat splitting
of the Stokes filtration.
Its inverse is denoted by $f^{\dagger}_S$.
We have estimates for them
similar to Lemma \ref{lem;07.10.11.100}
and Lemma \ref{lem;08.9.16.332}.

\subsubsection{Remark for gluing}

Let $\lambda_0\neq 0$
and $\mu_0=\lambda_0^{-1}$.
Let $U_{\lambda}(\lambda_0)$ be a small
neighbourhood of $\lambda_0$ in $\cnum_{\lambda}$,
and let $U_{\mu}(\mu_0)$ be 
the corresponding neighbourhood of
$\mu_0$ in $\cnum_{\mu}$.
Let $S$ be a sufficiently small sector 
in $U_{\lambda}(\lambda_0)\times (X-D)$,
and let $\nbigqzero_0\Vtilde_{0|\Sbar}
=\bigoplus_{\gminia}
 \nbigqzero_0\Vtilde_{0,\gminia,S}$
be a $\DD$-flat splitting.
Note that $S$ naturally gives
a small multi-sector
of $U_{\mu}(\mu_0)\times 
 \bigl(X^{\dagger}-D^{\dagger}\bigr)$,
and that the splitting naturally induces
a $\DD^{\dagger}$-splitting
$\nbigq^{(\mu_0)}_0\Vtilde_{\infty|\Sbar}
=\bigoplus_{\gminia}
 \nbigq^{(\mu_0)}_0\Vtilde_{\infty,\gminia,S}$.

\subsubsection{Estimate for pairing}

Let $S$ be a sufficiently small multi-sector of
$U_{\lambda}(\lambda_0)\times \bigl(X-D\bigr)$.
Let $g_S$ be as in (\ref{eq;08.9.16.320}),
and let $f_S$ denote its inverse.
We obtain the small multi-sector $\sigma(S)$
of $\sigma\bigl(U_{\lambda}(\lambda_0)\bigr)
 \times \bigl(X^{\dagger}-D^{\dagger}\bigr)$.
Let $g^{\dagger}_{\sigma(S)}$
be as in (\ref{eq;08.9.16.321}) for $\sigma(S)$,
and let $f^{\dagger}_{\sigma(S)}$ 
denote its inverse.
We can show the following lemma
by using Proposition \ref{prop;08.9.17.1}
and the argument in the proof of 
Lemma \ref{lem;08.9.16.200}.
\begin{lem}
\label{lem;08.9.16.400}
In Case (A),
we have the following estimate
with respect to $h$:
\begin{equation}
 \label{eq;08.9.16.450}
 \nbigs-\Gr^{\vecVtilde}\nbigs\circ
 (f_S\otimes \sigma^{\ast}f_{\sigma(S)}^{\dagger})
=O\Bigl(
 \exp\bigl(-\epsilon|z_1^{-1}|\bigr)
 \Bigr)
\end{equation}
In Case (C),
we have the following estimate
with respect to $h^{(0)}$:
\begin{equation}
 \label{eq;08.9.16.451}
 \Gr^{\vecVtilde}\nbigs-\nbigs\circ
 (g_S\otimes \sigma^{\ast}g^{\dagger}_{\sigma(S)})
=O\Bigl(
 \exp\bigl(-\epsilon|z_1^{-1}|\bigr)
 \Bigr)
\end{equation}
\hfill\qed
\end{lem}

Let $g_S^{(p)}$ $(p=1,\ldots,L)$ be 
as in (\ref{eq;08.9.16.320}),
and  $f_S^{(p)}$ denote their inverses.
Let $\chi_{S,p}$ $(p=1,\ldots,L)$ be non-negative valued
functions such that
(i) $\sum \chi_{S,p}=1$,
(ii) $\delbar \chi_{S,p}$ are of polynomial order
 in $|\lambda^{-1}|$ and $|z_1^{-1}|$.
We set
$g=\sum \chi_p\, g^{(p)}_S$
and $f=\sum \chi_p\, f^{(p)}_S$.
Let $g^{\dagger\,(q)}_{\sigma(S)}$ 
$(q=1,\ldots,M)$ be as in
(\ref{eq;08.9.16.321}) for $\sigma(S)$,
and $f^{\dagger\,(q)}_{\sigma(S)}$ 
denote their inverses.
Let $\chi_{\sigma(S),q}$ $(q=1,\ldots,M)$ 
be non-negative valued
functions such that
(i) $\sum \chi_{\sigma(S),q}=1$,
(ii) $\delbar \chi_{\sigma(S),q}$ are 
 of polynomial order in $|\mu^{-1}|$
and $|z_1^{-1}|$.
We set 
$g^{\dagger}=\sum \chi_{\sigma(S),q}\, 
 g^{\dagger(q)}_{\sigma(S)}$
and
$f^{\dagger}=\sum \chi_{\sigma(S),q}\, 
 f^{\dagger(q)}_{\sigma(S)}$.

\begin{lem}
\label{lem;08.9.17.50}
In Case (A),
we have the following estimate
with respect to $h$:
\begin{equation}
 \label{eq;08.9.16.440}
 \nbigs-\Gr^{\vecVtilde}\nbigs
 \circ(f\otimes \sigma^{\ast}f^{\dagger})
=O\Bigl(
 \exp\bigl(-\epsilon|\vecz^{-1}|\bigr)
 \Bigr)
\end{equation}
In Case (C),
we have the following estimate
with respect to $h^{(0)}$:
\begin{equation}
 \label{eq;08.9.16.441}
 \Gr^{\vecVtilde}\nbigs
-\nbigs\circ(g\otimes \sigma^{\ast}g^{\dagger})
=O\Bigl(
 \exp\bigl(-\epsilon|\vecz^{-1}|\bigr)
 \Bigr)
\end{equation}
\end{lem}
\pf
It follows from Lemma
\ref{lem;08.9.16.400}.
\hfill\qed

\subsection{Proof of Proposition
\ref{prop;07.9.10.3}}

We take a compact region 
$\nbigk\subset\cnum_{\lambda}$
such that $\nbigk\cup\sigma(\nbigk)=\proj^1$.
We set 
$W_{\nbigk}:=\bigl(\nbigk\times D\bigr)
 \cup\bigl(\{0\}\times X\bigr)$.
We take a finite covering
$\nbigx_{\nbigk}-W_{\nbigk}
\subset \bigcup_{p=1}^L S_p$
by small multi-sectors $S_p$
satisfying the following:
\begin{itemize}
\item
On each $S_p$, we can take
a $\DD$-flat morphism
 $g_{S_p}$ as in  (\ref{eq;08.9.16.320}).
Its inverse is denoted by $f_{S_p}$.
\item
On each $\sigma(S_p)$,
we can take a $\DD^{\dagger}$-flat
 morphism $g^{\dagger}_{\sigma(S_p)}$
as in (\ref{eq;08.9.16.321}).
Its inverse is denoted by 
 $f_{\sigma(S_p)}^{\dagger}$.
\end{itemize}
We take a partition of unity
$\bigl(
 \chi_{S_p},\chi_{\sigma(S_p)}\,\big|\,
 p=1,\ldots,L
 \bigr)$ subordinated to the covering
$(S_p,\sigma(S_p)\,|\,p=1,\ldots,L)$
such that
 $\del \chi_{S_p}$ 
 (resp. $\del\chi_{\sigma(S_p)}$)
 are of polynomial order
 in $|\lambda|^{-1}$ 
 (resp. $|\lambda|$)
 and  $|z_1|^{-1}$.
We set 
\[
 g:=\sum_p \chi_{S_p}\, g_{S_p}
+\sum_p\chi_{\sigma(S_p)}
 \, g^{\dagger}_{\sigma(S_p)},
\quad\quad
  f:=\sum_p \chi_{S_p}\, f_{S_p}
+\sum_p\chi_{\sigma(S_p)}
 \, f^{\dagger}_{\sigma(S_p)}.
\]
Let us consider the case (A).
By Lemma \ref{lem;08.9.16.332}
and its analogue in the $\mu$-side,
we obtain the following with respect to $h$:
\[
  f^{-1}\circ d''_{\proj^1,\Gr^{\vecVtilde}(V)}\circ f
-d''_{\proj^1,V}
=O\Bigl(
 \exp\bigl(-\epsilon\, (|\lambda^{-1}|+|\lambda|)
\, |z_1^{-1}|\bigr)
 \Bigr)
\]
By Lemma \ref{lem;08.9.17.50},
we have the following estimate
with respect to $h$:
\[
  \nbigs-\Gr^{\vecVtilde}\nbigs\circ
 (f_S\otimes \sigma^{\ast}f_{\sigma(S)}^{\dagger})
=O\Bigl(
 \exp\bigl(-\epsilon\,
 (|\lambda^{-1}|+|\lambda|)\,
 |z_1^{-1}|\bigr)
 \Bigr)
\]
Then, there exists $0<R<1$ such that
$\Gr^{\vecVtilde}(V,\DD^{\sankaku},\nbigs)
 _{|\proj^1\times X^{\ast}(R)}$
is a variation of polarized pure twistor structure,
due to Lemma \ref{lem;07.12.23.40}
and Lemma \ref{lem;07.6.2.90}.

Let us consider the case (C).
By a similar argument,
we obtain that 
$(V,\DD^{\sankaku},\nbigs)$
is a variation of polarized pure twistor structure.
Namely,
by Lemma \ref{lem;08.9.16.332},
we obtain the following with respect to $h^{(0)}$:
\[
  g^{-1}\circ d''_{\proj^1,V}\circ g
-d''_{\proj^1,\Gr^{\vecVtilde}(V)}
=O\Bigl(
 \exp\bigl(-\epsilon\,
 (|\lambda^{-1}|+|\lambda|)\,
 |z_1^{-1}|\bigr)
 \Bigr)
\]
By Lemma \ref{lem;08.9.17.50},
we have the following estimate
with respect to $h$:
\[
  \Gr^{\vecVtilde}\nbigs-\nbigs\circ
 (g\otimes \sigma^{\ast}g^{\dagger})
=O\Bigl(
 \exp\bigl(-\epsilon\,
 (|\lambda^{-1}|+|\lambda|)\,
 |z_1^{-1}|\bigr)
 \Bigr)
\]
Then, there exists $0<R<1$ such that
$(V,\DD^{\sankaku},\nbigs)
 _{|\proj^1\times X^{\ast}(R)}$
is a variation of polarized pure twistor structure,
due to Lemma \ref{lem;07.12.23.40}
and Lemma \ref{lem;07.6.2.90}.
Let $(E,\delbar_E,\theta,h)$ be 
the underlying harmonic bundle.
Let $d_{\Herm}$ denote the natural distance
of the symmetric space of the hermitian metrics.
Due to Lemma \ref{lem;07.6.2.90},
we have the following:
\[
 d_{\Herm}\bigl(
 g^{\ast}h, h^{(0)}
 \bigr)
=O\Bigl(
 \exp\bigl(-\epsilon|z_1^{-1}|\bigr)
 \Bigr)
\]
In particular,
the hermitian metrics
$g^{\ast}h$ and $h^{(0)}$
are mutually bounded.
It means that the full Stokes filtrations $\nbigftilde^S$ 
of the family $\Vtilde_0$ of meromorphic $\lambda$-flat bundles
can be characterized by the growth orders of 
the norms of the flat sections with respect to $h$,
as in Proposition \ref{prop;07.11.24.27}.
Therefore, we have the natural isomorphisms of 
the meromorphic prolongments
$\Vtilde_0\simeq\nbigq\nbige$.
Similarly, we obtain
$\Vtilde_{\infty}\simeq\nbigq\nbige^{\dagger}$.
Hence $(\Vtilde_0,\Vtilde_{\infty})$
is canonical.
\hfill\qed

\section{End of Proof of
 Theorem \ref{thm;08.9.16.100}}
\label{subsection;08.9.25.15}
\subsection{From (P1) to (P2)}
\label{subsection;08.9.25.10}

Assume that $(V,\DD^{\sankaku},\nbigs)$
is a variation of polarized pure twistor structure
of weight $0$,
and that the prolongment $\vecVtilde$ is canonical.
If we shrink $X$,
we have the harmonic bundle
$\bigl(E^{(1)},
 \delbar^{(1)},
 \theta^{(1)},
 h^{(1)}\bigr)$ on $X-D$
which induces
$\Gr^{\vecm(0)}(V^{\sankaku},
 \DD^{\sankaku},\nbigs)$
by Proposition \ref{prop;08.9.16.102}.
\begin{lem}
\label{lem;08.9.17.110}
$\Gr^{\vecm(0)}\Vtilde:=
 \bigl(
 \Gr^{\vecm(0)}(\Vtilde_0),
 \Gr^{\vecm(0)}(\Vtilde_{\infty})
 \bigr)$
is canonical.
\end{lem}
\pf
Let $k$ be determined by
$\vecm(0)\in\seisuu_{<0}^{k}\times\veczero_{\ell-k}$.
We take $1\leq j\leq k$.
Let $Q\in D_j$ be sufficiently close to
$D_{\kbar}$,
which is not singular point of $D(\leq k)$.
We take a small neighbourhood $X_Q$
of $Q$ in $X$.
We set $X_Q^{\ast}:=X_Q\setminus D$.
We put
\[
 (V_1,\DD_1^{\sankaku},\nbigs_1):=
 (V,\DD^{\sankaku},\nbigs)
 _{|\proj^1\times X_Q^{\ast}},
\quad
 (V_2,\DD_2^{\sankaku},\nbigs_2):=
 \Bigl(
 \Gr^{\vecm(0)}(V,\DD^{\sankaku},\nbigs)
 \Bigr)
 _{|\proj^1\times X_Q^{\ast}}.
\]
They are equipped with the meromorphic
prolongments 
induced by $\vecVtilde$
and $\Gr^{\vecm(0)}\vecVtilde$,
respectively.
They are denoted by
$\vecVtilde_1$ and $\vecVtilde_2$,
respectively.
As remarked in Lemma \ref{lem;08.9.17.105},
we have the natural isomorphism:
\[
 \Gr^{\vecVtilde_1}(V_1,\DD_1^{\sankaku},\nbigs_1)
\simeq
 \Gr^{\vecVtilde_2}(V_2,\DD_2^{\sankaku},\nbigs_2)
\]
By Proposition \ref{prop;07.9.10.3},
$\Gr^{\vecVtilde_1}(V_1,\DD_1^{\sankaku},\nbigs_1)$
is a variation of polarized pure twistor structure.
Again, according to Proposition \ref{prop;07.9.10.3},
$\vecVtilde_2$ is the canonical prolongment of
$(V_2,\DD^{\sankaku}_2,\nbigs_2)$.
Then, it is easy to conclude that
$\Gr^{\vecm(0)}\vecVtilde$ is
the canonical prolongment of
$\Gr^{\vecm(0)}(V,\DD^{\sankaku},\nbigs)$.
\hfill\qed

\vspace{.1in}
Now the claim $(P1)\Longrightarrow (P2)$
can be shown by an easy induction.

\subsection{From (P2) to (P1)}

In the following,
we will shrink $X$ without mention.
By using Proposition \ref{prop;08.9.16.102}
in a descending inductive way,
we obtain that
$\Gr^{\vecm(i)}(V,\DD^{\sankaku},\nbigs)$
are variation of polarized pure twistor structure
for any $i$.
In particular,
$(V,\DD^{\sankaku},\nbigs)$
is a variation of polarized pure twistor structure.
Let $(E,\delbar_E,\theta,h)$ be 
the underlying harmonic bundle,
which is unramifiedly good wild
with the good set $T$.
For each $j=1,\ldots,\ell$,
let $T(j)$ denote the image of
$T$ via the natural map
\[
  M(X,D)/H(X)\lrarr M(X,D)/M(X,D(\neq j)),
\]
where $D(\neq j):=\bigcup_{i\neq j,1\leq i\leq \ell}D_i$.
We put $D_j^{sm}:=D_j\setminus D(\neq j)$.
We take an auxiliary sequence
$\vecm(0),\ldots,\vecm(L)$ for $T$.
There exists $q(j)$ such that
$m_j(q(j))=-1$ and $m_j(q(j)+1)=0$,
where $m_j(i)$ denotes the $j$-th component of
$\vecm(i)$.
Note we have the natural bijection
$\etabar_{\vecm(q(j))}(T)
 \simeq T(j)$.

Let $Q\in D_j^{sm}$.
Let $X_Q$ be a small neighbourhood of $Q$.
We put $D_Q:=X_Q\cap D$
and $X_Q^{\ast}:=X_Q-D_Q$.
We set
$(V_Q,\DD^{\sankaku}_Q,\nbigs_Q):=
 (V,\DD^{\sankaku},\nbigs)_{|\proj^1\times X_Q^{\ast}}$.
It is equipped with a prolongment
$\vecVtilde_Q$ induced by $\vecVtilde$.
By the choice of $q(j)$,
\[
 \Gr^{\vecVtilde_Q}(V_Q,\DD^{\sankaku}_Q,\nbigs_Q)
\simeq
  \Gr^{\vecm(q(j))}(V,\DD^{\sankaku},\nbigs)_{
 |\proj^1\times X_Q^{\ast}}.
\]
Since $\Gr^{\vecVtilde_Q}
 (V_Q,\DD^{\sankaku}_Q,\nbigs_Q)$
is a variation of polarized pure twistor structure,
the prolongment $\vecVtilde_Q$
is canonical due to Proposition \ref{prop;07.9.10.3}.
Hence, we obtain that $\vecVtilde$ is canonical,
and thus the proof of Theorem \ref{thm;08.9.16.100}
is finished.
\hfill\qed

\section{Norm estimate}
\label{subsection;08.9.25.20}
\subsection{One step reduction}
\label{subsection;08.9.25.21}

We use the setting 
in Section \ref{subsection;07.11.17.1}.
We have the variation of 
polarized pure twistor structure
$(\nbige^{\sankaku},\DD^{\sankaku},\nbigs)$
associated to $\harmonicbundle$.
For simplicity,
we assume that the coordinate system is 
admissible for the good set $\Irr(\theta)$.
In the following argument,
we will shrink $X$ without mention,
if it is necessary.

We take an auxiliary sequence
$\vecm(0),\vecm(1),\ldots,\vecm(L)$
for $\Irr(\theta)$.
Let $(E_0,\delbar_0,
 \theta_0,h_0)$
be the unramifiedly good wild harmonic bundle
underlying 
$\Gr^{\vecm(0)}(\nbige^{\sankaku},
 \DD^{\sankaku},\nbigs)$.
We have the associated meromorphic 
$\lambda$-flat bundles
$(\nbigp_{\veca}\nbigelambda,\DDlambda)$,
and $(\nbigp_{\veca}\nbigelambda_{0},
 \DD^{\lambda}_{0})$ on $(X,D)$.
Let $k$ be determined by
$\vecm(0)\in\seisuu_{<0}^k\times\veczero_{\ell-k}$.
By construction,
we have the following natural isomorphism:
\[
\Phihat:
 \bigl(
 \nbigp_{\veca}\nbigelambda_0,\DDlambda_0
\bigr)_{|\Dhat(\leq k)}
\simeq
\bigl(\nbigp_{\veca}\nbigelambda,\DDlambda\bigr)
 _{|\Dhat(\leq k)}
\]
We have the harmonic bundles
$E_1:=\Hom(E_0,E)$ and $E_2:=\Hom(E,E_0)$
with the naturally induced Higgs fields $\theta_i$
and pluri-harmonic metrics $h_i$ $(i=1,2)$.
Note that
$(\nbigelambda_i,h_i)$ $(i=1,2)$ are
acceptable.
Let $\nbigp\nbigelambda_i$ denote
the associated meromorphic $\lambda$-bundle.
We can regard $\Phihat$ and $\Phihat^{-1}$
as sections of
$\nbigp_0\nbigelambda_{1|\Dhat(\leq k)}$
and $\nbigp_0\nbigelambda_{2|\Dhat(\leq k)}$,
respectively.

Let $N$ be a sufficiently large number.
According to Lemma \ref{lem;07.6.15.2},
we can take a section $\Phi_N$
of $\nbigp_0\nbigelambda_{1|\Dhat(\leq k)}$
with the following property:
\begin{itemize}
\item
 $\Phihat_{|\Dhat^{(N)}(\leq k)}
=\Phi_{N|\Dhat^{(N)}(\leq k)}$,
where $\Dhat^{(N)}(\leq k)$ denotes
the $N$-th infinitesimal neighbourhood of
$D(\leq k)$.
\item
$\Res_{i}(\DDlambda_1)(\Phi_{N|D_i})=0$
for $i=k+1,\ldots,\ell$.
\end{itemize}

\vspace{.1in}
We can regard $\Phi_N$ as an isomorphism
of $\nbigp_{\veca}\nbigelambda_0$ and 
$\nbigp_{\veca}\nbigelambda$
preserving the parabolic filtration.
If we shrink $X$ appropriately,
$\Phi_N$ is an isomorphism.

\begin{prop}
\label{prop;07.11.5.30}
$\Phi_{N|X-D}$ and $\Phi_{N|X-D}^{-1}$
are bounded with respect to $h$ and $h_0$.
\end{prop}
\pf
If we regard $\Phi_N$ as a section of
$\nbigp_0\nbigelambda_1$,
we obtain the boundedness with respect to $h_1$
due to Lemma \ref{lem;10.6.13.50}.
Thus, $\Phi_N$ is bounded with respect to 
$h$ and $h_0$.
We can regard $\Phi_N^{-1}$
as a section of 
$\nbigp_{0}\nbigelambda_2$
satisfying
$\Phihat^{-1}_{|\Dhat^{(N)}(\leq k)}
=\Phi^{-1}_{N|\Dhat^{(N)}(\leq k)}$
and
 $\Res_{i}(\DDlambda_2)(\Phi^{-1}_{N|D_i})=0$
for $i=k+1,\ldots,\ell$.
Then, $\Phi_N^{-1}$ is bounded
with respect to $h_2$
due to Lemma \ref{lem;10.6.13.50}.
Thus, $\Phi_N^{-1}$ is also bounded
with respect to $h$ and $h_0$.
\hfill\qed

\subsection{Full reduction}
\label{subsection;08.9.25.22}

By using an inductive argument,
we can reduce
the norm estimate for 
unramifiedly good wild harmonic bundle
to that for tame harmonic bundles
studied in \cite{mochi2}.
Let $(E_3,\delbar_{3},\theta_3,h_3)$
be the harmonic bundle
underlying 
$\Gr^{\nbigftilde}(\nbige^{\sankaku},
 \DD^{\sankaku},\nbigs)$.
It is graded
\begin{equation}
 \label{eq;07.11.5.40}
\bigl(E_3,\delbar_{E_3},\theta_3,h_3\bigr):=
\bigoplus_{\gminia\in\Irr(\theta)}
(E_{\gminia},\delbar_{E_{\gminia}},
 \theta_{\gminia},h_{\gminia}),
\end{equation}
corresponding to 
the decomposition
$\Gr^{\nbigftilde}
 (\nbige^{\sankaku},\DD^{\sankaku},\nbigs)
=\bigoplus
 \Gr^{\nbigftilde}_{\gminia}
 (\nbige^{\sankaku},\DD^{\sankaku},
 \nbigs)$.
Each $\theta_{\gminia}-d\gminia$ is tame.
We have the associated 
meromorphic $\lambda$-flat bundles
$(\nbigp_{\veca}\nbigelambda,\DDlambda)$
and 
$(\nbigp_{\veca}\nbigelambda_3,\DDlambda_3)$.

\begin{thm}
\label{thm;07.11.5.41}
There exists a holomorphic isomorphism
$\Phi:\nbigp_{\veca}\nbigelambda_3\lrarr
 \nbigp_{\veca}\nbigelambda$
with the following property:
\begin{itemize}
\item
 It preserves the parabolic structures.
\item
 $\Res_i(\DDlambda)\circ\Phi_{|D_i}
 -\Phi_{|D_i}\circ\Res_i(\DDlambda_3)=0$.
\end{itemize}
Moreover,
 $\Phi$ and $\Phi^{-1}$ are bounded
with respect to $h$ and $h_3$.
\end{thm}
\pf
We have only to use Proposition 
\ref{prop;07.11.5.30} inductively.
\hfill\qed

\begin{rem}
We obtained
the norm estimate for holomorphic sections
of $\lambda$-flat bundles
associated to tame harmonic bundles,
in terms of the parabolic filtration
and the weight filtration.
(See Theorem {\rm 13.29} of {\rm\cite{mochi2}}.)
Since each
$(E_{\gminia},\delbar_{\gminia},
 \theta_{\gminia}-d\gminia,h_{\gminia})$
is tame,
it can be applied to
$(\nbigp\nbigelambda_3,\DDlambda_3,h_3)$,
and transferred to 
$(\nbigp\nbigelambda,\DDlambda,h)$
via the morphism $\Phi$
in Theorem {\rm\ref{thm;07.11.5.41}}.
Hence, it gives a satisfactory norm estimate
for holomorphic sections of good wild harmonic bundles.
\hfill\qed
\end{rem}

\subsection{Surface case}
\label{subsection;08.9.25.30}

We give a rather detailed description of
the norm estimate for holomorphic sections
of $\nbigp_{\vecc}\nbigelambda$
in the case $\dim X=2$ and $D=D_1\cup D_2$.
In this subsection,
$\harmonicbundle$ is assumed to be good wild,
but not necessarily unramified.

We have the parabolic filtrations
$\lefttop{i}F$ of
$\nbigp_{\vecc}\nbigelambda_{|D_i}$ $(i=1,2)$.
Let $\lefttop{i}\Gr^F_a
 \bigl(\nbigp_{\vecc}\nbigelambda\bigr):=
 \lefttop{i}F_{a}\bigl(\nbigp_{\vecc}\nbigelambda\bigr)
 \big/\lefttop{i}F_{<a}\bigl(\nbigp_{\vecc}\nbigelambda\bigr)$,
which are equipped with
the induced endomorphisms
$\lefttop{i}\Gr_a(\Res_i(\DDlambda))$.
The eigenvalues of
$\lefttop{i}\Gr_a(\Res_i(\DDlambda))_{|Q}$
are independent of the choice of $Q\in D_i$.
(It is clear in the case $\lambda\neq 0$.
 It follows from Proposition \ref{prop;07.7.19.31}
 or the definition of good wild harmonic bundles,
 in the case $\lambda=0$.)
Hence, we have the well defined nilpotent part
$N_{i,a}$ of $\Gr_a(\Res_i(\DDlambda))$.
Let $\Par(\nbigp_{\vecc}\nbigelambda,i)$
denote the set of $a\in\real$
such that 
$\lefttop{i}\Gr^F_a\bigl(
 \nbigp_{\vecc}\nbigelambda
 \bigr)\neq 0$.

We put 
$\lefttop{\nibar}F_{\veca}
 \bigl(\nbigp_{\vecc}\nbigelambda_{|O}
 \bigr):=
 \lefttop{1}F_{a_1}
 \bigl(\nbigp_{\vecc}\nbigelambda_{|O}\bigr)
 \cap
  \lefttop{2}F_{a_2}
 \bigl(\nbigp_{\vecc}\nbigelambda_{|O}\bigr)$
for $\veca=(a_1,a_2)\in\real^2$.
Then, we put
\[
 \lefttop{\nibar}\Gr^{F}_{\veca}
 \bigl(\nbigp_{\vecc}\nbigelambda\bigr):=
 \frac{\lefttop{\nibar}F_{\veca}
 \bigl(\nbigp_{\vecc}\nbigelambda_{|O}\bigr)
 }
 {\sum_{\vecb\lneq\veca}
  \lefttop{\nibar}F_{\vecb}
 \bigl(\nbigp_{\vecc}\nbigelambda_{|O}\bigr)}.
\]
Let $\Par\bigl(\nbigp_{\vecc}\nbigelambda,O\bigr)$
denote the set of $\veca\in\real^2$
such that 
$ \lefttop{\nibar}\Gr^{F}_{\veca}
 \bigl(\nbigp_{\vecc}\nbigelambda\bigr)\neq 0$.
On $\lefttop{\nibar}\Gr^F_{\veca}(\nbigp_{\vecc}\nbigelambda)$,
we have the induced endomorphisms
$\lefttop{\nibar}\Gr_{\veca}(\Res_i(\DDlambda))$
$(i=1,2)$.
The nilpotent parts are denoted by $N_{i,\veca}$.

\begin{lem}
\mbox{{}}
\begin{itemize}
\item
The conjugacy classes of $N_{i,a|Q}$
are independent of the choice of $Q\in D_i$.
\item
Let $q_i:\Par(\nbigp_{\vecc}\nbigelambda,O)
\lrarr\Par(\nbigp_{\vecc}\nbigelambda,i)$
denote the projection.
The conjugacy classes of
$N_{i,a|O}$ and 
$\bigoplus_{\veca\in q_i^{-1}(a)}N_{i,\veca}$
are the same.
\end{itemize}
\end{lem}
\pf
We take a ramified covering 
$\varphi:\Xtilde\lrarr X$ given by
$\varphi(\zeta_1,\zeta_2)
=(\zeta_1^e,\zeta_2^e)$
such that 
$(\Etilde,\delbar_{\Etilde},\thetatilde,\htilde)
:=\varphi^{-1}\harmonicbundle$
is unramified.
Let 
 $(\nbigp_{\ast}\nbigelambdatilde,
 \DDlambdatilde)$ 
denote the associated 
filtered $\lambda$-flat bundle.
Take 
$\Qtilde\in \Dtilde_i-\{\Otilde\}$
and $Q=\varphi(\Qtilde)\in D_i-\{O\}$.
We have the following isomorphisms
for any $-1<b\leq 0$
which preserve the conjugacy classes
of the nilpotent parts of the residues:
\[
 \lefttop{i}\Gr^F_{b}
 (\nbigp_0\nbigelambdatilde)_{|\Qtilde}
\simeq
\bigoplus_{
 \substack{a\in\Par(\nbigp_{\vecc}\nbigelambda,i)\\
 ae-b\in\seisuu}}
 \lefttop{i}\Gr^F_a(\nbigp_{\vecc}\nbigelambda)_{|Q}
\]
For $\vecb\in \openclosed{-1}{0}^2$,
we have the following isomorphisms
preserving the conjugacy classes of
the nilpotent parts of the residues:
\[
 \lefttop{\nibar}\Gr^F_{\vecb}
 (\nbigp_0\nbigelambdatilde)_{|\Otilde}
\simeq
\bigoplus_{\substack{
 \veca\in\Par(\nbigp_{\vecc}\nbigelambda,O)\\
 e\veca-\vecb\in\seisuu^2
 }}
\lefttop{\nibar}\Gr^F_{\veca}
 (\nbigp_{\vecc}\nbigelambda)_{|O}
\]
Hence, we have only to show the claim for 
$(\nbigp_{\ast}\nbigelambdatilde,\DDlambdatilde)$,
i.e.,
in the unramified case.

Let $(\Etilde_0,\delbar_{\Etilde_0},\thetatilde_0,\htilde_0)$
be obtained as the full reduction from 
$(\Etilde,\delbar_{\Etilde},\thetatilde,\htilde)$.
We can take an isomorphism 
$\Phitilde:\nbigp_0\nbigelambdatilde_0
 \lrarr\nbigp_0\nbigelambdatilde$
which preserves the parabolic filtrations
and the residues,
as in Theorem \ref{thm;07.11.5.41}.
Hence, we have only to show the claim
for tame harmonic bundles.

In the tame case, the claim is proven 
in \cite{mochi2}.
We indicate an outline.
If $\lambda$ is generic,
i.e., the maps
$\eigenmap(\lambda):
 \KMS(\nbige^0,i)\lrarr \cnum$
are injective for any $i$,
the claim is clear,
because the conjugacy classes of 
the nilpotent parts of
$\lambda^{-1}\Res_i(\DDlambda)_{|Q}$ 
for any $Q\in D_i$
are equal to those of the logarithm
of the unipotent part
of the monodromy around $D_i$.
By using Corollary 12.43 of \cite{mochi2},
we obtain that
the conjugacy classes of the nilpotent parts of
$\lefttop{\nibar}\Gr^F\Res_i(\DDlambda)_{|O}$
and $\lefttop{i}\Gr^F\Res_i(\DDlambda)_{|Q}$
for $Q\in D_i-\{O\}$ are independent
of the choice of $\lambda$.
Thus, we are done.
\hfill\qed

\vspace{.1in}

We put $N_{1}:=\bigoplus N_{1,a}$.
Let $W(N_{1})$ denote the weight filtration
of $N_{1}$ on 
$\lefttop{1}\Gr^F(\nbigp_{\vecc}\nbigelambda)$.
In particular,
we obtain the filtration
$W(N_{1})_{|O}$
of $\lefttop{1}\Gr^F(\nbigp_{\vecc}\nbigelambda)_{|O}$.
It induces the filtration
on $\lefttop{\nibar}\Gr^F
 (\nbigp_{\vecc}\nbigelambda)$,
which is denoted by
$W(N_{1})^{(1)}$.
We also put $N^{(1)}_{1}:=\bigoplus N_{1,\veca}$,
which is an endomorphism
of $\lefttop{\nibar}\Gr^{F}(\nbigp_{\vecc}\nbigelambda)$.
We have the weight filtration $W(N_1^{(1)})$
of $\lefttop{\nibar}
\Gr^F(\nbigp_{\vecc}\nbigelambda)$.

\begin{lem}
\label{lem;07.11.5.45}
$W(N_{1})^{(1)}$
and $W(N_{1}^{(1)})$ are the same.
\end{lem}
\pf
Although this claim is proved in \cite{mochi2},
we give an outline.
We have the induced filtration $\lefttop{2}F$
on $\lefttop{1}\Gr^F
 \bigl(\nbigp_{\vecc}\nbigelambda\bigr)_{|O}$.
Let $\nbigv$ be the vector bundle
on $\Spec\cnum[t]$ 
obtained as the Rees bundle associated to
$\lefttop{2}F$.
The endomorphism $N_{1|O}$ naturally induces
$\nbign_1$ on $\nbigv$.
The restriction to $t=0$ is $N^{(1)}_{1|O}$.
Since the degeneration of the conjugacy classes
does not happen,
the weight filtration of $W(\nbign_1)$ is
the filtration in the category of the vector bundles
on $\Spec\cnum[t]$.
The specialization at $t=0$ is equal to
$W(N^{(1)}_1)$,
and the specialization at $t\neq 0$
is equal to $W(N)$.

Let $\nbigw$ denote the filtration
naturally induced by $W(N_1)$.
The specialization of $\nbigw$ at $t=0$
is equal to $W(N_1)^{(1)}$.
We also have $\nbigw=W(\nbign_1)$.
Hence, we obtain $W(N_1^{(1)})=W(N_1)^{(1)}$.
\hfill\qed

\vspace{.1in}

We have the nilpotent endomorphism
$N_{\veca}(\nibar)=N_{1,\veca}+N_{2,\veca}$.

\begin{lem}
\label{lem;07.11.5.50}
There exists a decomposition
\[
 \nbigp_{\vecc}\nbigelambda
=\bigoplus_{(\veca,\veck)\in 
 \Par(\nbigp_{\vecc}\nbigelambda,O)\times\seisuu^2}
 U_{\veca,\veck}
\]
with the following property:
\begin{itemize}
\item
It gives a splitting of the parabolic filtrations:
\[
 \bigoplus_{q_i(\vecb)\leq a}
 \bigoplus_{\veck}
 U_{\vecb,\veck|D_i}
=\lefttop{i}F_a\bigl(\nbigp_{\vecc}\nbigelambda_{|D_i}\bigr)
\]
\[
 \bigoplus_{\vecb\leq\veca}
 \bigoplus_{\veck}
 U_{\vecb,\veck|O}
=\lefttop{\nibar}F_{\veca}
 \bigl(\nbigp_{\vecc}\nbigelambda_{|O}\bigr)
\]
\item
Under the isomorphism
$\lefttop{1}\Gr^F_a(\nbigp_{\vecc}\nbigelambda)\simeq
 \bigoplus_{\veca\in q_1^{-1}(a)}
 \bigoplus_{\veck}U_{\veca,\veck|D_1}$,
the following holds:
\[
  \bigoplus_{\veca\in q_1^{-1}(a)}
 \bigoplus_{k_1\leq l_1}
 U_{\veca,\veck|D_1}
=W_{l_1}(N_1)\bigl(
 \lefttop{1}\Gr^F_a
 \bigl(\nbigp_{\vecc}\nbigelambda\bigr)
\bigr)
\]
Under the isomorphism
$\lefttop{\nibar}\Gr^F_{\veca}\bigl(
 \nbigp_{\vecc}\nbigelambda \bigr)
\simeq
 \bigoplus_{\veck}U_{\veca,\veck|O}$,
the following holds:
\[
 \bigoplus_{\veck\leq\vecl}
 U_{\veca,\veck|O}
=W_{l_1}(N_1)\cap
 W_{l_2}(N(\nibar))
\Bigl(
 \lefttop{\nibar}\Gr^F_{\veca}
 \bigl(\nbigp_{\vecc}\nbigelambda\bigr)
\Bigr)
\]
\end{itemize}
\end{lem}
\pf
It follows from Lemma \ref{lem;07.11.5.45}.
(We can use a more general result,
 for example Corollary 4.47
 of \cite{mochi2}.)
\hfill\qed

\vspace{.1in}

We take a decomposition as in Lemma \ref{lem;07.11.5.50}.
Let $\vecv$ be a frame of 
$\nbigp_{\vecc}\nbigelambda$
compatible with the decomposition.
When $v_i\in U_{\veca,\veck}$,
we put $\veca(v_i):=\veca$
and $\veck(v_i):=\veck$.

Let $h_1$ be the hermitian metric given as follows:
\[
 h_1(v_i,v_j):=\delta_{i,j}\, 
 |z_1|^{-2a_1(v_i)}\,
 |z_2|^{-2a_2(v_i)}\,
 \bigl(-\log|z_1|\bigr)^{k_1(v_i)}\,
 \bigl(-\log|z_2|\bigr)^{k_2(v_i)-k_1(v_i)}
\]
Let $Z:=\bigl\{(z_1,z_2)\,\big|\,|z_1|<C|z_2|\bigr\}$.

\begin{prop}
$h$ and $h_1$ are mutually bounded
on $Z$.
\end{prop}
\pf
The problem can be reduced to the unramified case.
It follows from 
Theorem \ref{thm;07.11.5.41}
and the norm estimate in the tame case 
(Theorem 13.29 of \cite{mochi2}).
\hfill\qed

\subsection{Blow up}
\label{subsection;07.11.5.101}

Let $\Xtilde:=\Delta^2=\{(\zeta_1,\zeta_2)\}$.
Let $\pi:\Xtilde\lrarr X$ be given by
$\pi(\zeta_1,\zeta_2)=(\zeta_1\!\,\!\zeta_2,\,\zeta_2)$.
Let $(\Etilde,\delbar_{\Etilde},\thetatilde,\htilde):=
 \pi^{-1}\harmonicbundle$.
We have the associated filtered bundle
$\nbigp_{\ast}\nbigelambdatilde$.
\begin{lem}
$\nbigp_{\ast}\nbigelambdatilde$
is obtained from $\nbigp_{\ast}\nbigelambda$
by the procedure in Section
{\rm \ref{subsection;07.11.5.60}}.
\end{lem}
\pf
It follows from the weak norm estimate 
for the acceptable bundles.
\hfill\qed

\vspace{.1in}
For simplicity, we consider the case $\vecc=(0,0)$.
In the case $a_1(v_i)+a_2(v_i)>-1$,
we put 
$\vtilde_i:= \pi^{\ast}v_i$.
In the case $a_1(v_i)+a_2(v_i)\leq -1$,
we put $\vtilde_i:= \pi^{\ast}v_i\, \zeta_2^{-1}$.
Then, $\vecvtilde=(\vtilde_i)$
gives a frame of $\nbigp_0\nbigelambdatilde$
compatible with the parabolic structure.
We put $a_j(\vtilde_i):=\lefttop{j}\deg^F(\vtilde_i)$.
We also put $k_j(\vtilde_i):=k_j(v_i)$.
Let $\htilde_0$
be the hermitian metric given as follows:
\[
 \htilde_0(\vtilde_i,\vtilde_j)
=\delta_{i,j}\,
 |\zeta_1|^{-2a_1(\vtilde_i)}
\,
 |\zeta_2|^{-2a_2(\vtilde_i)}
\]
Let $\chi$ be a non-negative valued function on $\real$
such that $\chi(t)=1$ $(t\leq 1/2)$
and $\chi(t)=0$ $(t\geq 2/3)$.
Let $\rho(\zeta):\cnum^{\ast}\lrarr\real$
 be the function given by
$\rho(\zeta)=-\chi(|\zeta|)\, \log|\zeta|^2$.
We set
\[
 \htilde_1(\vtilde_i,\vtilde_j)
:=\htilde_0(\vtilde_i,\vtilde_j)
\,
  (1+\rho(\zeta_1)+\rho(\zeta_2))^{k_1(\vtilde_i)}
\,
 (1+\rho(\zeta_2))^{k_2(\vtilde_i)-k_1(\vtilde_i)}.
\]

\begin{lem}
\label{lem;07.11.5.110}
$\pi^{\ast}h$ and $\htilde_1$
are mutually bounded.
The curvature $R(\htilde_0)$ is $0$.
Moreover,
$R(\htilde_1)$ and $\del_{\htilde_1}-\del_{\htilde_0}$
are bounded with respect to
both $(\omegatilde,\htilde_i)$ $(i=0,1)$,
where $\omegatilde$ denotes the Poincar\'e metric
of $\Xtilde-\Dtilde$.
\end{lem}
\pf
The first claim follows from the norm estimate.
The other claim can be shown by direct calculations.
\hfill\qed

\section[Regular meromorphic variation
 of twistor structure]{Regular
 meromorphic variation 
 of twistor structure on a disc (Appendix)}
\label{subsection;08.9.14.50}
Let $X$ be the disc $\{z\in\cnum\,|\,|z|<1\}$,
and let $D=\{0\}$.
We consider a regular meromorphic extension
of a variation of pure twistor structure
on $\proj^1\times(X-D)$.
In \cite{mochi2},
we showed that
if the variation of pure twistor structure
is pure and polarized,
then the limit twistor structure
is mixed polarized.
We shall study the converse
Lemma \ref{lem;07.7.11.12}.
Although a similar result is given in \cite{sabbah2},
we would like to understand it from our viewpoint.

We shall review the construction of
the limit twistor structure
(Subsections
\ref{subsection;10.6.1.2}--\ref{subsection;08.1.28.10}).
Then, we study the characterization of
pure and polarized property.

\subsection{Preliminary}
\label{subsection;10.6.1.2}

Let $\lambda_0\in\cnum$,
and let $\nbigk$ be a neighbourhood of $\lambda_0$
in $\cnum$.
We put $\nbigx:=\nbigk\times X$
and $\nbigd:=\nbigk\times D$.
In this case,
$\nbigd$ can be identified with $\nbigk$ naturally.

Let $(\vecV_{\ast},\DD)$ be 
a good family of filtered $\lambda$-flat bundles,
which is regular in the sense
$\DD(\prolongg{a}V)\subset 
 \prolongg{a}V\otimes\Omega^{1,0}(\log D)$.
The restriction to $\nbigx\setminus\nbigd$
is denoted by $V$.
Assume that $\vecV_{\ast}$ 
has the KMS-structure at $\lambda_0$
indexed by $T\subset\real\times\cnum$,
i.e.,
$\KMS\bigl(\vecV_{\ast}\bigr)
=\bigl\{\kmsmap(\lambda_0,u)\,\big|\,
 u\in T\bigr\}$.
(See Section \ref{section;10.5.20.1}.)
We have the natural $\seisuu$-action on $T$
given by $n\cdot u=u+(n,0)$.
We recall some objects
induced by the KMS structure.

\subsubsection{The induced bundle $\nbigg_u(V)$}

Recall that we have the induced filtration
$\Fzero$ of $\prolongg{a}V$
and 
the generalized eigen decomposition
\[
 \Gr^{\Fzero}_a(V)=
 \bigoplus_{\substack{u\in T,\\
 \kmsmap(\lambda_0,u)=a}}
 \nbiggzero_u(V)
\]
on $\nbigk$,
as in Section \ref{subsection;07.10.14.15}.
Let $\nbign_u$ denote the nilpotent part
of the residue on $\nbiggzero_u(V)$.

\subsubsection{The KMS-structure
 of the space of multi-valued flat sections}

Assume $\lambda_0\neq 0$.
The restriction of $(V,\DD)$ to 
$\{\lambda\}\times (X\setminus D)$
is denoted by $(V^{\lambda},\DDlambda)$.
Let $\nbigh(V)$ denote the holomorphic vector bundle
on $\nbigk$,
whose fiber over $\lambda$ is the space 
$H(V^{\lambda})$ of the multi-valued
flat sections of $(V^{\lambda},\DD^{\lambda,f})$.
We have the monodromy automorphism $M$
along the loop with the counter-clockwise direction.
The restriction to $\lambda$ is denoted by 
$M^{\lambda}$.
The set of the eigenvalues of $M^{\lambda}$
is given by
 $\Sp^f(V^{\lambda}):=
 \bigl\{
 \eigenmap^f(\lambda,u)\,\big|\,
 u\in T/\seisuu \bigr\}$.
(Recall $\eigenmap^f(\lambda,u)
 =\exp\bigl(-2\pi\sqrt{-1}\lambda^{-1}
 \eigenmap(\lambda,u)
 \bigr)$.)
We have the unique monodromy invariant decomposition
\begin{equation}
 \label{eq;08.1.28.1}
 \nbigh(V)=
 \bigoplus_{\omega\in \Sp^f(V^{\lambda_0})}
 \EEzero_{\omega}\nbigh(V)
\end{equation}
whose restriction to $\lambda_0$ is the same as
the generalized eigen decomposition
of $M^{\lambda_0}$.

Let $\nbigk^{\ast}:=\nbigk\setminus\{\lambda_0\}$.
We may assume that
any $\lambda\in\nbigk^{\ast}$ is generic.
We have the generalized eigen decomposition:
\[
  \nbigh(V)_{|\nbigk^{\ast}}
=\bigoplus_{ u\in T/\seisuu}
 \EE_{\eigenmap^f(\lambda,u)}
 \nbigh(V)_{|\nbigk^{\ast}}
\]
Here, the fiber of
$\EE_{\eigenmap^f(\lambda,u)}
 \nbigh(V)_{|\nbigk^{\ast}}$ over $\lambda$
is the generalized eigen space of $M^{\lambda}$
corresponding to $\eigenmap^f(\lambda,u)$.
We put
\[
 \nbigfzero_b\nbigh(V)_{|\nbigk^{\ast}}
:=\bigoplus_{\substack{u\in T/\seisuu\\
 \paramap^f(\lambda_0,u)\leq b}}
\EE_{\eigenmap(\lambda,u)}
 \nbigh(V)_{|\nbigk^{\ast}},
\]
where $\paramap^f(\lambda,u):=
 \Re\bigl(\lambda\alphabar+\lambda^{-1}\alpha\bigr)
=\paramap(\lambda,u)+
 \Re\bigl(\lambda^{-1}\,\eigenmap(\lambda,u)\bigr)$.
Thus, we obtain a filtration
indexed by 
$\bigl\{\paramap^f(\lambda,u)\,\big|\,
 u\in\Tbar\bigr\}\subset\real$.
It is extended to a filtration
of $\nbigh(V)$ on $\nbigk$,
denoted by $\nbigfzero$.
It is monodromy invariant,
and compatible with the decomposition $\EEzero$.

We put $\nbiggzero_u\nbigh(V):=
 \Gr^{\nbigfzero,\EEzero}_{\kmsmap^f(\lambda_0,u)}
 \nbigh(V)
=\Gr^{\nbigfzero}_{\paramap^f(\lambda_0,u)}
 \EEzero_{\eigenmap^f(\lambda_0,u)}
 \nbigh(V)$ on $\nbigk$.
We have the automorphism $M_u$ induced by $M$,
whose unipotent part is denoted by $M_u^{\rm uni}$.
We put $N_u:=(-2\pi\sqrt{-1})^{-1}\log M_u^{\rm uni}$.

\begin{rem}
$\nbigfzero$ is not related with
Stoke filtrations.
Because we consider this kind of filtration
only in the regular case,
there is no risk of confusion.
\hfill\qed
\end{rem}

\subsubsection{The decomposition and the filtration
 of $\prolongg{a}{V}$}

The decomposition (\ref{eq;08.1.28.1}) induces
a $\DD$-flat decomposition
$\vecV_{\ast}
=\bigoplus_{\omega\in \Sp^f(V^{\lambda_0})} 
 \vecV_{\omega,\ast}$.
The filtration $\nbigfzero$ of $\nbigh(V)$
naturally induces a $\DD$-flat filtration 
of $V$ on $\nbigx\setminus\nbigd$,
which is also denoted by $\nbigfzero$.
By construction, the subbundles
$\nbigfzero_b$ of $V$ are naturally extended
to those of
$\prolongg{a}{V}_{|\nbigx-\{(\lambda_0,O)\}}$.
Moreover, we have the following lemma.
\begin{lem}
$\nbigfzero_b$ are naturally extended
to subbundles of 
$\prolongg{a}{V}$ on $\nbigx$.
Namely,
we obtain an induced filtration
$\nbigfzero$ of $\prolongg{a}{V}$ 
in the category of vector bundles
on $\nbigx$.
\end{lem}
\pf
It is easy to reduce the problem
to the case $\Sp^f(V^{\lambda_0})=\{1\}$.
First, let us consider the case in which
(i) $b$ is the minimum among the numbers $c$
such that $\Gr^{\nbigfzero}_c\nbigh(V)\neq 0$,
(ii) $\rank \nbigfzero_b\nbigh(V)=1$.
Let $s$ be a frame of
$\nbigfzero_b\nbigh(V)$.
We have the element $u_0\in T$
such that 
$\paramap^f(\lambda_0,u_0)=b$
and $a-1<\paramap(\lambda_0,u_0)\leq a$.
Then,
$t:=s\cdot 
 \exp\bigl(\lambda^{-1}
 \eigenmap(\lambda,u_0)\, \log z\bigr)$
naturally gives a single-valued holomorphic
section of 
$\nbigfzero_b\bigl(
 \prolongg{a}{V}_{|\nbigx\setminus\{(\lambda_0,O)\}}
 \bigr)$.
It is extended to a section of $\prolongg{a}{V}$.
Let us show
$t_{|(\lambda_0,O)}\neq 0$
in $\prolongg{a}{V^{\lambda_0}}_{|O}$.
Let $t^{\lambda_0}$ denote the restriction  of $t$
to $\{\lambda_0\}\times X$.
We have the relation
$\DD^{\lambda_0}t^{\lambda_0}
=t^{\lambda_0}\cdot
  \eigenmap(\lambda_0,u_0)\, dz/z$.
If $t^{\lambda_0}_{|O}=0$,
there exists 
$\bigl(c,\eigenmap(\lambda_0,u_0)\bigr)
 \in \KMS\bigl(\vecV_{\ast}^{\lambda_0}\bigr)$
such that $c<a$.
However, 
it is easy to show
$\paramap(\lambda_0,u)\geq 
 \paramap(\lambda_0,u_0)=a$
for any $u\in T$ such that
$\eigenmap(\lambda,u)
=\eigenmap(\lambda_0,u_0)$,
by using the relation
$\paramap(\lambda_0,u)-\paramap(\lambda_0,u_0)
=\paramap^f(\lambda_0,u)-\paramap^f(\lambda_0,u_0)$.
Hence, we obtain $t^{\lambda_0}_{|(\lambda_0,O)}\neq 0$
in $\prolongg{a}{V}_{|(\lambda,O)}$.
Thus, $\nbigfzero_{b}$ gives a subbundle
of $\prolongg{a}{V}$ on $\nbigx$ in this case.

We can reduce the general case
to the above case,
by using the exterior product.
\hfill\qed

\vspace{.1in}

By construction,
the filtration $\nbigfzero$ is compatible
with the decomposition
$\prolongg{a}{V}
=\bigoplus \prolongg{a}{V}_{\omega}$.
Let us look at the restriction
of the decomposition and the filtration
to $\nbigd$.
Clearly, we have the following:
\[
 \prolongg{a}{V}_{\omega|\nbigd}
=\bigoplus_{
 \exp(-2\pi\sqrt{-1}(\alpha/\lambda))=\omega
 }
 \EEzero_{\alpha}(\prolongg{a}V_{|\nbigd})
\]
(See Remark \ref{rem;08.9.4.10}
 for $\EEzero_{\alpha}$.)
We take $u\in T$ such that
$a-1<\paramap(\lambda_0,u)\leq a$
and $\eigenmap^f(\lambda_0,u)=\omega$.
It is easy to show the following:
\begin{equation}
 \label{eq;08.1.28.2}
\nbigfzero_{\paramap(\lambda_0,u)}\bigl(
 \prolongg{a}{V_{\omega}}
\bigr)_{|\nbigd}
=\Fzero_{\paramap(\lambda_0,u)}
 \EEzero_{\eigenmap(\lambda_0,u)}
 (\prolongg{a}{V_{|\nbigd}})
\oplus
 \bigoplus_{\alpha\in S_1}
 \EEzero_{\alpha}(\prolongg{a}{V}_{|\nbigd})
\end{equation}
Here,
$ S_1:=\bigl\{
 \exp\bigl(-2\pi\sqrt{-1}(\alpha/\lambda_0)\bigr)=\omega,
 \,\,
 \lambda_0^{-1}\bigl(
 \alpha-\eigenmap(\lambda_0,u)
 \bigr)\in\seisuu_{<0}
 \bigr\}$.

\subsubsection{The induced bundle
 $\nbigg_u\nbigv$
 and the isomorphism}
\label{subsection;08.1.28.110}

For any $u\in T$,
we set 
\[
 \nbiggzero_u\nbigv:=
 \frac{\nbigfzero_{\paramap^f(\lambda_0,u)}\Bigl(
 \prolongg{\paramap(\lambda_0,u)}
 V_{\eigenmap^f(\lambda_0,u)} \Bigr)}
{\nbigfzero_{<\paramap^f(\lambda_0,u)}\Bigl(
 \prolongg{\paramap(\lambda_0,u)}
 V_{\eigenmap^f(\lambda_0,u)} \Bigr)}.
\]
We have an induced family of flat $\lambda$-connections
$\DD_u$ on $\nbiggzero_u\nbigv$,
which is logarithmic.
The eigenvalue of the residue
$\Res(\DD_u)$ is $\eigenmap(\lambda,u)$.
The nilpotent part is denoted by $\nbign_u$.

By (\ref{eq;08.1.28.2}),
we have a natural isomorphism
$\bigl(
 \nbiggzero_u\nbigv_{|\nbigd},\,
 \nbign_u
\bigr)
\simeq
\bigl(
 \nbiggzero_u(V),\,\nbign_u
\bigr)$.
The space of the multi-valued flat sections
$\nbigh(\nbiggzero_u\nbigv)$
with the monodromy 
is naturally isomorphic to
$\nbiggzero_u\nbigh(V)$ with $M_u$.
In this situation,
we have an isomorphism
$\Phi^{\can}:(\nbiggzero_u\nbigh(V),\lambda\, N_u)
\lrarr (\nbiggzero_u(V),\nbign_u)$
which is given as follows.
(See Section 10.4.1 of \cite{mochi2},
 for example.)
Let $F$ be a section of $\nbiggzero_u\nbigh(V)$.
We can regard it as a multi-valued flat
section of $\nbigg_u\nbigv$ on $\nbigx\setminus\nbigd$.
Then, we put $\Phi^{\can}_u(F):=F_0$
for the expansion
\[
 F=\sum_{j=0}^m F_j\,
 \exp\bigl(-\lambda^{-1}
 \eigenmap(\lambda,u)\log z\bigr)
 \,(\log z)^j, 
\]
where $F_j$ are holomorphic sections of
$\nbigg_u\nbigv$.

\subsection{Globalization}
Let $\nbigx:=\cnum_{\lambda}\times X$
and $\nbigd:=\cnum_{\lambda}\times D$.
If $(V,\DD)$ is a family of
meromorphic $\lambda$-flat bundles
on $(\nbigx,\nbigd)$,
which has the KMS-structure 
at each $\lambda_0\in\cnum_{\lambda}$
indexed by $T$.
For each $\lambda_0$,
we have a neighbourhood $U(\lambda_0)$ of $\lambda_0$
in $\cnum_{\lambda}$,
and we obtain
$(\nbiggzero_u(V),\,\nbign_u)$ on $U(\lambda_0)$.
In the case $\lambda_0\neq 0$,
we also have $(\nbiggzero_u\nbigh(V),\,N_u)$
and the isomorphism
$\Phi^{\can}:(\nbiggzero_u\nbigh(V),\,\lambda\, N_u)
\lrarr (\nbiggzero_u(V),\nbign_u)$.
By using the uniqueness of KMS-structure
(see Lemma {\rm\ref{lem;07.11.23.5}},
 for example),
we can glue them for various $\lambda_0$,
and so we obtain the following:
\begin{itemize}
\item
A bundle
$\nbigg_u(V)$ with the nilpotent endomorphism
$\nbign_u$ on $\cnum_{\lambda}$.
\item
A bundle $\nbigg_u\nbigh(V)$ 
with the nilpotent endomorphism $N_u$
on $\cnum_{\lambda}^{\ast}$.
\item
An isomorphism
$\Phi^{\can}:(\nbigg_u\nbigh(V),\lambda\, N_u)
\lrarr (\nbigg_u(V),\nbign_u)_{|\cnum_{\lambda}^{\ast}}$.
\end{itemize}
We also obtain a vector bundle
$\nbigg_u\nbigv$ on $\cnum_{\lambda}^{\ast}\times X$
with a family of logarithmic flat $\lambda$-connections
$\DD_u$,
which is locally constructed as above.

\begin{rem}
The construction is functorial 
and compatible with
dual, tensor product, and direct sum.
\hfill\qed
\end{rem}

\subsection{Gluing}
\label{subsection;08.1.28.10}

Let $\nbigx:=\cnum_{\lambda}\times X$
and $\nbigx^{\dagger}:=\cnum_{\mu}\times X^{\dagger}$.
We use the symbols
$\nbigd$ and $\nbigd^{\dagger}$
in similar meanings.
For $\lambda_0\in\cnum_{\lambda}$,
let $\nbigx^{(\lambda_0)}$ denote 
a product of $X$ and a neighbourhood $U(\lambda_0)$
of $\lambda_0$.
We use the symbols
$\nbigd^{(\lambda_0)}$,
$\nbigx^{\dagger(\mu_0)}$,
and $\nbigd^{\dagger(\mu_0)}$
in similar meanings.

Let $(V,\DD^{\sankaku})$ be a variation of
twistor structure on
$\proj^1\times (X\setminus D)$
with an unramifiedly good meromorphic prolongment 
$(\Vtilde_0,\Vtilde_{\infty})$.
Moreover, we assume the following:
\begin{itemize}
\item
 $\Vtilde_0$ has the KMS-structure at each $\lambda_0$.
 Namely, 
 we have locally free 
 $\nbigo_{\nbigx^{(\lambda_0)}}$-lattices
 $\nbigpzero_a(V)\subset \Vtilde_0$ $(a\in\real)$
 such that 
 $(\nbigpzero_{\ast}(V_0),\DD_{V_0})$
 is a regular family of filtered $\lambda$-flat bundles
 with KMS-structure indexed by $T_0$.
\item
 Similarly,
 for each $\mu_0$,
 we have 
 locally free 
 $\nbigo_{\nbigx^{\dagger(\mu_0)}}$-lattices
 $\nbigp^{(\mu_0)}_a(V_{\infty})
 \subset \Vtilde_{\infty}$ $(a\in\real)$
 such that 
 $(\nbigp^{(\mu_0)}_{\ast}(V_{\infty}),
 \DD_{V_{\infty}}^{\dagger})$
 is a regular family of filtered $\mu$-flat bundles
 with KMS-structure
 indexed by $T_{\infty}$.
\end{itemize}
In that case,
we say that $(\Vtilde_0,\Vtilde_{\infty})$
gives a regular meromorphic extension of
$(V,\DD^{\sankaku})$.
\index{regular meromorphic extension}

\begin{lem}
\mbox{{}}
\begin{itemize}
\item
The map 
$u=(a,\alpha)\longmapsto u^{\dagger}=(-a,\alphabar)$
induces a bijection
$T_0\simeq T_{\infty}$.
\item
We have a natural isomorphism
$\bigl(\nbigg_{u}\nbigh(V_0),N_u\bigr)
\simeq
 \bigl(
 \nbigg_{u^{\dagger}}\nbigh(V_{\infty}),
 -N_{u^{\dagger}}^{\dagger}
 \bigr)$.
\end{itemize}
\end{lem}
\pf
Let $(V^{\lambda}_0,\DDlambda)$ 
denote the restriction
of $(V_0,\DD_{V_0})$ to $\{\lambda\}\times X$.
Let $\Sp\bigl(\nbigp_aV_0^{\lambda}\bigr)$
denote the set of the eigenvalues of
$\Res(\DDlambda)$,
and $\Sp(V^{\lambda}_0):=
 \bigcup_{a\in\real}\Sp\bigl(\nbigp_aV_0^{\lambda}\bigr)$.
Let $\Sp^f\bigl(V_0^{\lambda}\bigr)$ denote
the set of the eigenvalues of
the monodromy on the space of
the multi-valued flat sections of $V_0^{\lambda}$,
where the monodromy is taken along the loop
with the counter-clockwise direction.
Then, we have the bijective correspondence
between
$\Sp(V^{\lambda}_0)/\seisuu$ 
and $\Sp^f(V_0^{\lambda})$
given by
$\alpha\longleftrightarrow
 \exp\bigl(-2\pi\sqrt{-1}(\alpha/\lambda)\bigr)$.
(The action of $\seisuu$ on
$\Sp(V^{\lambda})$ is given by
$n\cdot \alpha=\alpha+n\lambda$.)
Hence,
the set $\Sp(V_0^{\lambda})$ is determined
by $\Sp^f(V_0^{\lambda})$.

Note that 
$\Sp(V^{\lambda}_0)$ is the image of $T_0$
via the map $\eigenmap(\lambda)$.
There exists a discrete subset $Z_0$
of $\cnum_{\lambda}^{\ast}$
such that $\eigenmap(\lambda):T_0\lrarr\cnum$
is injective for any 
$\lambda\in\cnum_{\lambda}^{\ast}-Z_0$.
Hence, $T_0$ is determined by the 
family of the sets
$\bigl\{
 \Sp^f(V^{\lambda}_0)\,\big|\,
\lambda\in \cnum_{\lambda}^{\ast}-Z_0
 \bigr\}$.

Similarly,
let $V^{\mu}_{\infty}$ denote the restriction
of $V_{\infty}$ to 
$\{\mu\}\times (X^{\dagger}-D^{\dagger})$,
and let $\Sp^f(V^{\mu}_{\infty})$ $(\mu\neq 0)$
denote the set of the eigenvalues 
of the monodromy on the space of
the multi-valued flat sections of  
$V^{\mu}_{\infty}$,
where the monodromy is taken along the loop
with the clockwise direction.
There exists a discrete subset 
$Z_{\infty}\subset \cnum_{\mu}^{\ast}$
such that 
$\eigenmap(\mu):T_{\infty}\lrarr \cnum$
is injective for any $\mu\in\cnum_{\mu}-Z_{\infty}$,
and the set $T_{\infty}$ is determined
by the family of the sets
$\bigl\{
 \Sp^f(V^{\mu}_{\infty})\,\big|\
 \mu\in \cnum_{\mu}^{\ast}-Z_{\infty}
 \bigr\}$.

For $\lambda=\mu^{-1}$,
we have the natural bijection
$\Sp^f(V^{\lambda}_0)
\simeq
 \Sp^f(V^{\mu}_{\infty})$
given by $\omega\longleftrightarrow \omega^{-1}$.
Then, we obtain the desired bijection
$T_0\simeq T_{\infty}$ by the formal calculation.
Note 
$\eigenmap(\lambda,u)
=\eigenmap(\lambda^{-1},u^{\dagger})$.
Hence, we obtain the first claim.

We have the natural identification
$\nbigh(V_0)\simeq \nbigh(V_{\infty})$
on $\cnum_{\lambda}^{\ast}=\cnum_{\mu}^{\ast}$.
By using
$\paramap^f(\lambda,u)
=\paramap^f(\lambda^{-1},u^{\dagger})$
and 
$\eigenmap^f(\lambda,u)
=\eigenmap^f(\lambda^{-1},u^{\dagger})^{-1}$,
we can show the second claim.
\hfill\qed

\vspace{.1in}

Then, we obtain the vector bundle 
$S_u(V)$ on $\proj^1$
by gluing
$\nbigg_u(V_0)$  and 
$\nbigg_{u^{\dagger}}(V_{\infty})$.
\index{vector bundle $S_u(V)$}
We also obtain the nilpotent map
$\nbign^{\sankaku}:
 S_u(V)\lrarr S_u(V)\otimes \Tate(-1)$,
where 
$\nbign^{\sankaku}_{|\cnum_{\lambda}}
:=\nbign_u\cdot t_0^{(-1)}$
and $\nbign^{\sankaku}_{|\cnum_{\mu}}
 =\nbign_{u^{\dagger}}^{\dagger}
 \cdot t_{\infty}^{(-1)}$.
\index{nilpotent map $\nbign^{\sankaku}$}

\begin{rem}
The construction is functorial
and compatible with dual,
tensor product, and direct sum.
\hfill\qed
\end{rem}

\subsubsection{Pairing}
\label{subsection;08.1.28.11}

Let $\nbigs:V\otimes\sigma^{\ast}(V)\lrarr\Tate(0)$
be a symmetric perfect pairing.
By the regularity,
$\nbigs_0$ is extended to a pairing
$\Vtilde_0\otimes\sigma^{\ast}\Vtilde_{\infty}
\lrarr \nbigo_{\nbigx}(\ast\nbigd)$.
Similarly,
$\nbigs_{\infty}$ is extended to a pairing
$\Vtilde_{\infty}\otimes\sigma^{\ast}\Vtilde_0
\lrarr \nbigo_{\nbigx}(\ast\nbigd)$,
i.e.,
$(\Vtilde_0,\Vtilde_{\infty})$
gives a meromorphic prolongment
of $(V,\DD^{\sankaku},\nbigs)$.
By using the uniqueness of the KMS-structure,
we obtain an induced pairings:
\begin{itemize}
\item
$\nbigs_u:
 \nbigg_u(V_0)\otimes
 \sigma^{\ast}\nbigg_{u^{\dagger}}(V_{\infty})
\lrarr \nbigo_{\cnum_{\lambda}}$
satisfying
$\nbigs_u(\nbign_u\otimes \id)
-\nbigs_u(\id\otimes
 \sigma^{\ast}\nbign^{\dagger}_{u^{\dagger}})
=0$.
(See Lemma \ref{lem;08.1.29.1}
 for the signature.)
\item
$\nbigs_u:
 \nbigg_u\nbigh(V_0)\otimes
 \sigma^{\ast}\nbigg_{u^{\dagger}}\
 \nbigh(V_{\infty})
 \lrarr \nbigo_{\cnum_{\lambda}^{\ast}}$.
\item
$\nbigs_u:
 \nbigg_u\nbigv_0\otimes
 \sigma^{\ast}
 \nbigg_{u^{\dagger}}\nbigv_{\infty}
 \lrarr \nbigo_{\cnum_{\lambda}^{\ast}\times X}$.
\end{itemize}
Then, we can show that the gluing
of $\nbigg_u(V_0)$ and 
$\nbigg_{u^{\dagger}}(V_{\infty})$
is compatible with the pairing.
Thus, we obtain 
the symmetric pairing
$\nbigs_u:S_u(V)\otimes\sigma^{\ast}S_u(V)
\lrarr \Tate(0)$,
which satisfies
$\nbigs_u\bigl(\nbign_u^{\sankaku}\otimes\id\bigr)
+\nbigs_u\bigl(
 \id\otimes\sigma^{\ast}(\nbign_u^{\sankaku})\bigr)
=0$.

\subsection{A characterization of purity
and polarizability}
\label{subsection;08.9.29.18}

A result similar to the following lemma was shown
in \cite{sabbah2} with a different argument.

\begin{lem}
\label{lem;07.7.11.12}
If $\bigl(S_u(V),\nbign_u^{\sankaku},\nbigs_u\bigr)$
is a polarized mixed twistor structure
of weight $0$ for each $u\in T$,
then $(V,\DD^{\sankaku},\nbigs)$ is 
a variation of polarized pure twistor structure
of weight $0$
after $X$ is shrinked around $D$.
\end{lem}
\pf
We have the polarized variation of
pure twistor structure 
$\bigl(V_u^{(1)},\DD_u^{\sankaku(1)},\nbigs_u^{(1)}\bigr)$
induced by
the polarized mixed twistor structure
$\bigl(S_u(V),\nbign_u^{\sankaku},\nbigs_u\bigr)$.
(See Sections 3.5.3, 3.6 and 3.7.5 of \cite{mochi2}.)
We have the corresponding tame harmonic bundle
$\bigl(E_{0,u}^{(1)},\delbar_{0,u},
 \theta_{0,u}^{(1)},h_{0,u}^{(1)}\bigr)$.
We put 
\[
 (E^{(1)},\delbar_{E^{(1)}},\theta^{(1)},h^{(1)}):=
\bigoplus_{u}
 \bigl(E_{u}^{(1)},\delbar_u^{(1)},    
 \theta_{u}^{(1)},h_{u}^{(1)}\bigr),
\]
\[
 \bigl(E_{u}^{(1)},\delbar_u^{(1)},     
 \theta_{u}^{(1)},h_{u}^{(1)}\bigr)
=
 \bigl(E_{0,u}^{(1)},\delbar_{0,u},
 \theta_{0,u}^{(1)},h_{0,u}^{(1)}\bigr)
 \otimes L(u).
\]
Let $\bigl(V^{(1)},\DD^{(1)\sankaku},
 \nbigs^{(1)\sankaku}\bigr)$
be the corresponding
variation of polarized pure twistor structure.
We have the regular meromorphic extension
$(\Vtilde^{(1)}_0,\Vtilde^{(1)}_{\infty})$.
We recover
$\bigoplus
\bigl(S_u(V),\nbign_u^{\sankaku},\nbigs_u^{\sankaku}\bigr)$
by applying the construction in 
Sections 
\ref{subsection;08.1.28.10}--\ref{subsection;08.1.28.11}
to $\bigl(V^{(1)},\DD^{(1)\sankaku},\nbigs^{(1)\sankaku}\bigr)$.
Namely, we have the canonical isomorphisms
$G_u:S_u(V^{(1)})\simeq S_u(V)$
compatible with the nilpotent maps
and the pairings.

\vspace{.1in}

Let us construct a $C^{\infty}$-map
$f:V^{(1)}\lrarr V$.
In the following, $U(\lambda_0)$
denotes a small neighbourhood of 
$\lambda_0\in\cnum_{\lambda}$.
 Note we have the isomorphisms:
\begin{equation}
 \label{eq;08.1.28.15}
 \Gr^{\Fzero}_a(V^{(1)}_0)\simeq
 \bigoplus_{\paramap(\lambda_0,u)=a}
 S_u(V^{(1)}_0)_{|U(\lambda_0)},
\quad
 \Gr^{\Fzero}_a(V_0)\simeq
 \bigoplus_{\paramap(\lambda_0,u)=a}
 S_u(V_0)_{|U(\lambda_0)}
\end{equation}
If $\lambda_0$ is generic,
we have the unique {\em flat} isomorphism
$f_{U(\lambda_0)}:\nbigp^{(\lambda_0)}_{\ast}(V^{(1)}_0)
\lrarr
 \nbigp^{(\lambda_0)}_{\ast}(V_0)$
with the following property:
\begin{itemize}
\item
Under the isomorphisms (\ref{eq;08.1.28.15}),
 the induced map 
\[
 \Gr^{\Fzero}_a(f_{U(\lambda_0)}):
 \Gr^{\Fzero}_a(V^{(1)}_0)\simeq
 \Gr^{\Fzero}_a(V_0)
\]
is equal to the restriction of
 $\bigoplus_{\paramap(\lambda_0,u)=a} G_u$.
\end{itemize}
Even if $\lambda_0$ is not generic,
we can take a holomorphic (not necessarily flat)
isomorphism 
$f_{U(\lambda_0)}:
 \nbigp^{(\lambda_0)}_{\ast}(V^{(1)}_0)
\lrarr
 \nbigp^{(\lambda_0)}_{\ast}(V_0)$
of the filtered bundles on $U(\lambda_0)\times (X,D)$
with the above property.
If $f'_{U(\lambda_0)}$ is another isomorphism
satisfying the above condition,
we have
\begin{equation}
 \label{eq;08.10.28.10}
\bigl|\id-f_{U(\lambda_0)}^{-1}\circ f_{U(\lambda_0)}'
 \bigr|_{h^{(1)}}=O\bigl(|z|^{\epsilon}\bigr)
\end{equation}

Let $U(\mu_0)$ denote a neighbourhood of
$\mu_0\in\cnum_{\mu}$.
We take isomorphisms
$f^{\dagger}_{U(\mu_0)}:
 \nbigp^{(\mu_0)}_{\ast}(V^{(1)}_{\infty})
\lrarr
 \nbigp^{(\mu_0)}_{\ast}(V_{\infty})$
of filtered flat bundles on 
$U(\mu_0)\times (X^{\dagger},D^{\dagger})$
with a similar property for each $\mu_0\in\cnum_{\mu}$.
Let us look at the morphism:
\[
J:=
 \nbigs^{(1)}_0-
 \nbigs_0\circ\bigl(
 f_{U(\lambda_0)}\otimes 
 \sigma^{\ast}f^{\dagger}_{U(\sigma(\lambda_0))}
 \bigr)
 :\nbigp^{(\lambda_0)}_{\ast}(V_0)
\otimes
 \sigma^{\ast}
 \nbigp^{(\sigma(\lambda_0))}_{\ast}(V_{\infty})
\lrarr \nbigo_{\nbigx^{(\lambda_0)}}(\ast\nbigd^{(\lambda_0)})
\]
Because of the conditions for $f_{U(\lambda_0)}$
and $f^{\dagger}_{U(\sigma(\lambda_0))}$,
we have
$|J|=O\bigl(|z|^{\epsilon}\bigr)$ for some $\epsilon>0$
with respect to the metric $h^{(1)}$.

We take compact regions $W_1\subset\cnum_{\lambda}$
and $W_2\subset\cnum_{\mu}$
such that (i) $W_1\cup W_2=\proj^1$,
(ii) any $\lambda\in W_1\cap W_2$ is generic.
We can take a finite covering
$W_1\subset \bigcup U(\lambda_0)$.
By gluing $f_{U(\lambda_0)}$ in $C^{\infty}$,
we obtain $f_{W_1}$.
Similarly, we obtain $f^{\dagger}_{W_2}$.
Note
$f_{W_1|W_1\cap W_2}
=f^{\dagger}_{W_2|W_1\cap W_2}$.
By gluing them in $C^{\infty}$, we obtain $f$.

Let $f_Q:=f_{|\proj^1\times Q}$.
Let $V_{Q}:=V_{|\proj^1\times Q}$
and $V^{(1)}_Q:=V^{(1)}_{|\proj^1\times Q}$.
Then, we have
\[
 f_Q^{-1}\circ\delbar_{V_Q}\circ f_Q
-\delbar_{V^{(1)}_Q}
 =O\bigl(|z(Q)|^{\epsilon}\bigr)
\]
with respect to $h^{(1)}$
due to (\ref{eq;08.10.28.10}).
Hence, $(V,\DD^{\sankaku})$ is a variation of
pure twistor structure
if $Q$ is sufficiently close to $O$,
due to Lemma \ref{lem;07.12.23.40}.
Because of the estimate of $J$ above,
there exists a constant $C>0$
such that
the following holds
for any $u\in V^{(1)}_{|(\lambda,Q)}$
and $v\in V^{(1)}_{(\sigma(\lambda),Q)}$:
\[
\bigl|
 \nbigs^{(1)}(u,\sigma^{\ast}v)
 -\nbigs\bigl(f_Q(u)\otimes\sigma^{\ast}f_Q(v)\bigr)
 \bigr|
\leq C\cdot |z(Q)|^{\epsilon}
 \cdot |u|_{h^{(1)}}\cdot |v|_{h^{(1)}}
\]
Hence, $S$ gives the polarization
due to Lemma \ref{lem;07.6.2.90}.
Thus, the proof of Lemma \ref{lem;07.7.11.12}
is finished.
\hfill\qed

\begin{rem}
Lemma \ref{lem;07.7.11.12}
can be generalized in the higher dimensional case.
We omit the details here.
\hfill\qed
\end{rem}

\begin{rem}
The converse of Lemma {\rm\ref{lem;07.7.11.12}}
was proved in {\cite{mochi2}}.
Namely, 
if $(V,\DD^{\sankaku,\nbigs})$
is a variation of polarized pure twistor structure
of weight $0$,
then $\bigl(S_u(V),\nbign_u^{\sankaku},\nbigs_u\bigr)$
is a polarized mixed twistor structure
of weight $0$ for each $u\in T$.
\hfill\qed
\end{rem}

\chapter{Prolongation as $\nbigr$-Triple}
\label{section;07.10.12.10}
In this chapter,
we consider the $\nbigr_X$-modules $\gbige$
and the $\nbigr$-triples $\gbigt(E)$
associated to unramifiedly good 
wild harmonic bundles
$\harmonicbundle$.
In particular, we study their specialization
along any monomial functions.
These results are mainly preliminary for
Theorem \ref{thm;07.10.28.30}.
(We refer to \cite{sabbah2}
for details on $\nbigr$-modules and 
$\nbigr$-triples.
We will give a review on
$\nbigr$-modules,
$\nbigr$-triples and variants
in Chapter \ref{section;08.10.20.40}.)

In Section \ref{subsection;08.9.28.1},
we construct $\nbigr$-modules $\gbige$
associated to unramifiedly good wild 
harmonic bundles
$\harmonicbundle$.
It is a natural generalization
of the construction in the tame case
which was studied in \cite{mochi2}.

We investigate the basic property of such $\nbigr$-modules
in Sections
\ref{subsection;08.9.28.2}--\ref{subsection;08.2.2.1}.
It is our basic strategy to reduce the study
to the tame case.
Hence, we review 
in Section \ref{subsection;08.9.28.2}
our previous result on $\gbige$
for the tame case \cite{mochi2}.
Then, we study in Section 
\ref{subsection;08.9.28.3}
the $\nbigr$-module
associated to the tensor product
of a tame harmonic bundle
and a rank one wild harmonic bundle.
Applying these preliminary result,
we show the basic property of $\gbige$
in Section \ref{subsection;08.2.2.1}.
In particular, we show 
in Section \ref{subsection;07.9.7.35}
the strict $S$-decomposability of $\gbige$
along any monomial functions.
And, we show 
in Section \ref{subsection;08.9.28.4}
the strict $S$-decomposability
of $P\Gr^{W(N)}_p\psitilde_{g,u}(\gbige)$
along any coordinate function.
Sections 
\ref{subsection;08.9.28.5}--\ref{subsection;08.2.2.5}
are preliminary for Section 
\ref{subsection;08.9.28.10}.
We study the specialization and the reduction.

In Section \ref{subsection;08.9.28.11},
we construct a Hermitian sesqui-linear
pairing $\gbigc$ of $\gbige$,
and thus we obtain an $\nbigr$-triple
$\gbigt(E)=(\gbige,\gbige,\gbigc)$.

In Section \ref{subsection;08.9.28.12},
we give a characterization of 
$\gbigt(E)$ in the case $\dim X=1$.
This is one of the main differences
between the tame case and the wild case.
In the tame case,
such a characterization is given
in a much simpler way,
because of the uniqueness of
meromorphic prolongment with regular singularity.
However, we do not have such a nice
uniqueness in the irregular case.
So we need a more consideration.

In Section \ref{subsection;08.9.28.10},
we study the specialization of
$\nbigr$-triples $\gbigt(E)$.
By Proposition \ref{prop;07.9.8.5},
we give a comparison between
the specializations of
the original $\nbigr$-triple $\gbigt(E)$
and the reduced one.
By using this proposition,
we can conclude that
the direct summands of
$P\Gr^{W(N)}_{p}\psitilde_{g,u}(\gbigt(E))$
{\em generically} come from unramifiedly good
wild harmonic bundles.
Then, by Proposition \ref{prop;07.9.8.26},
we show that they come from
unramifiedly good wild harmonic bundles.

\section[Prolongation as 
$\nbigr$-module]{$\nbigr$-module
associated to unramifiedly good wild harmonic bundle}
\label{subsection;08.9.28.1}
\label{subsection;08.1.11.2}

Let $X:=\Delta^n$,
$D_i:=\{z_i=0\}$ and $D:=\bigcup_{i=1}^nD_i$.
Let $\harmonicbundle$ be
an unramifiedly good wild harmonic bundle on $X-D$.
We assume that the coordinate system
is admissible for the good set $\Irr(\theta)$.
We use the notation $\nbigx:=\cnum_{\lambda}\times X$
and $\nbigd:=\cnum_{\lambda}\times D$.
We obtained
the family of meromorphic $\lambda$-flat bundles
$(\nbigq\nbige,\DD)$ on $(\nbigx,\nbigd)$
in Section \ref{subsection;07.11.5.1}.
We can naturally regard $\nbigq\nbige$ 
as a left $\nbigr_X$-module.
Let $U(\lambda_0)$ be a sufficiently small neighbourhood
of $\lambda_0$.
Let $\vecdelta=(1,\ldots,1)\in\real^n$.
For any $\veca\in\real^n$,
let $\nbigq^{(\lambda_0)}_{<\veca}\nbige$
be the union of
$\nbigq^{(\lambda_0)}_{\vecb}\nbige$
for $\vecb\in\real^n$ such that
$b_i<a_i$ $(i=1,\ldots,n)$.
It is equal to
$\nbigq^{(\lambda_0)}_{\veca-\epsilon\vecdelta}
 \nbige$ for some $\epsilon>0$,
and hence it is a locally free
$\nbigo_{U(\lambda_0)\times X}$-lattice of
$\nbigq^{(\lambda_0)}\nbige:=
 \nbigq\nbige_{|U(\lambda_0)\times X}$.

In particular, we have the lattice
$\nbigq^{(\lambda_0)}_{<\vecdelta}\nbige$.
Let $\gbige^{(\lambda_0)}$ denote 
the $\nbigr_X$-submodule
of $\nbigqzero\nbige$ generated by
$\nbigq^{(\lambda_0)}_{<\vecdelta}\nbige$
over $\nbigr_X$:
\[
 \gbige^{(\lambda_0)}:=
 \nbigr_X\cdot \nbigq^{(\lambda_0)}_{<\vecdelta}\nbige
\]
\begin{lem}
\label{lem;07.9.7.12}
$\gbige^{(\lambda_0)}$ is a coherent
$\nbigr_X$-module.
\end{lem}
\pf
Since $\nbigq\nbige$ is a locally free
$\nbigo_{\nbigx}(\ast \nbigd)$-module,
$\gbige^{(\lambda_0)}$ is a pseudo-coherent 
$\nbigo_{\nbigx}$-module.
Since $\gbige^{(\lambda_0)}$ is finitely generated
over $\nbigr_X$,
the claim of the lemma follows.
See Proposition \ref{prop;08.10.20.31}
below, for example.
\hfill\qed

\vspace{.1in}

Let us take $\lambda_1\in U(\lambda_1)\subset U(\lambda_0)$.
\begin{lem}
\label{lem;07.9.7.10}
We have
$\gbige^{(\lambda_0)}_{|U(\lambda_1)\times X}=
 \gbige^{(\lambda_1)}$.
As a result,
we obtain the global $\nbigr_X$-module $\gbige$
on $\nbigx$.
\index{sheaf $\gbige$}
\end{lem}
\pf
We have 
$\nbigq^{(\lambda_0)}_{<\vecdelta}
\nbige_{|U(\lambda_1)\times X}
\subset
 \nbigq^{(\lambda_1)}_{<\vecdelta}\nbige$.
Hence, $\gbige^{(\lambda_0)}_{|U(\lambda_1)\times X}
 \subset \gbige^{(\lambda_1)}$.
We would like to show the reverse implication.
The composite of the following morphisms
are denoted by $\pi_i$
for $i=1,\ldots,n$:
\[
 \nbigqzero_{\vecdelta}\nbige
\lrarr
 \nbigqzero_{\vecdelta}
 \nbige_{|U(\lambda_0)\times D_i}
\lrarr
 \lefttop{i}\Gr^{\Fzero}_1\bigl(
 \nbigqzero_{\vecdelta}
 \nbige_{|U(\lambda_0)\times D_i}
 \bigr)
\]
Here $\lefttop{i}\Gr^{\Fzero}$ is taken with respect to 
the naturally induced filtration
$\lefttop{i}\Fzero$ of 
$\nbigqzero_{\vecdelta}\nbige_{|U(\lambda_0)\times D_i}$.
Let $\nbigk(1,\lambda_0):=
 \bigl\{
 u\in\KMS(\nbige^0,i)\,\big|\,
 \paramap(\lambda_0,u)=1 \bigr\}$.
We have the generalized eigen decomposition
with respect to the induced endomorphism
$\Res_i(\DD)$:
\[
 \lefttop{i}\Gr^{\Fzero}_1\bigl(
 \nbigqzero_{\vecdelta}
 \nbige_{|U(\lambda_0)\times D_i}
 \bigr)
=\bigoplus_{u\in\nbigk(1,\lambda_0)}
 \lefttop{i}\EE_{u}
\]
Here, $\Res_i(\DD)-\eigenmap(\lambda,u)$ is 
nilpotent on $\lefttop{i}\EE_{u}$.
We put
\[
 \nbigk_m:=\bigcap_{i=1}^m\pi_i^{-1}\Bigl(
 \bigoplus_{u\neq (1,0)}\lefttop{i}\EE_u
 \Bigr)
\cap
 \bigcap_{i=m+1}^{n}\pi_i^{-1}(0).
\]
We have $\nbigk_0=
   \nbigq^{(\lambda_0)}_{<\vecdelta}\nbige$,
and hence $\nbigk_0\subset\gbige^{(\lambda_0)}$.
It is easy to see
$\nbigq^{(\lambda_1)}_{<\vecdelta}\nbige
\subset  \nbigk_{n|U(\lambda_1)\times X}$.
Hence, 
we have only to show 
$\nbigk_m\subset\gbige^{(\lambda_0)}$.
We put 
$\nbigkhat_{m}:=
 \nbigk_{m|\widehat{U(\lambda_0)\times O}}$
and 
$\gbigehat^{(\lambda_0)}:=
 \gbige^{(\lambda_0)}_{|
 \widehat{U(\lambda_0)\times O}}$.
We have only to show
$\nbigkhat_m\subset
 \gbigehat^{(\lambda_0)}$.
We use an induction on $m$.

Let us show $m-1\Longrightarrow m$.
Since $\nbigqzero_{<\vecdelta}\nbige$
is an unramifiedly good lattice,
we have the irregular decomposition
\begin{equation}
 \label{eq;10.6.14.1}
 \nbigqzero_{<\vecdelta}\nbigehat:=
 \nbigqzero_{<\vecdelta}\nbige_{|
 \widehat{U(\lambda_0)\times O}}
=\bigoplus_{\gminia\in \Irr(\theta)}
 \nbigqzero_{<\vecdelta}\nbigehat_{\gminia}.
\end{equation}
Let $\DD_0$ denote 
the family of logarithmic $\lambda$-connections
of $\nbigq_{<\vecdelta}\nbigehat$
given as follows:
\[
 \DD_0:=\DD-
 \bigoplus_{\gminia\in\Irr(\theta)} d\gminia\cdot 
 \id_{\nbigqzero_{<\vecdelta}\nbigehat_{\gminia}}
\]
Note $\Res_i(\DD_0)=\Res_i(\DD)$.
Let $\vecv$ be a frame of
$\nbigqzero_{<\vecdelta}\nbigehat$,
which is compatible with
the above irregular decomposition (\ref{eq;10.6.14.1}),
the parabolic filtrations $\lefttop{i}\Fzero$,
and the decomposition $\lefttop{i}\EEzero$.
We also assume that
the induced frame of
$\lefttop{m}
 \Gr^{\Fzero}(\nbigqzero_{<\vecdelta}\nbigehat)$
is compatible with the weight filtration
of the nilpotent part $\nbign_m$ of 
the endomorphism
on $\lefttop{m}\Gr^{\Fzero}
 (\nbigqzero_{<\vecdelta}\nbigehat)$
induced by $\Res_m(\DD)$.
(Note that the conjugacy classes of
$\nbign_{m|(\lambda,Q)}$ are independent of
$(\lambda,Q)\in U(\lambda_1)\times D_m$.
 It can be reduced to the tame case,
 by using the map in Theorem \ref{thm;07.11.5.41},
 or the completion at $(\lambda,O)$.
 The tame case was argued in \cite{mochi2}.
 See Lemma 12.47 of \cite{mochi2},
 for example.)
Let $\gminia_j$ be determined by
$v_j\in \nbigqzero_{<\vecdelta}\nbigehat_{\gminia_j}$.
Let $I(v_j)$ be the set of $i$
such that $\lefttop{i}\deg^{\EEzero}(v_j)\neq 0$
and $\lefttop{i}\deg^{\Fzero}(v_j)=0$.
For $l=m-1,m$, we put
\[
 \vtilde^{(l)}_j:=
 \prod_{\substack{
 i\leq l \\ i\in I(v_j)}}
 z_i^{-1}\, v_j.
\]
Then, $\vecvtilde^{(m)}$ is a frame of $\nbigkhat_m$.
We have only to show
$\vtilde_j^{(m)}\in \gbigehat^{(\lambda_0)}$.
If $m\not\in I(v_j)$,
we have $\vtilde_j^{(m)}\in\nbigkhat_{m-1}$,
and thus there is nothing to prove.
Let us consider the case $m\in I(v_j)$,
i.e.,
$\vtilde^{(m)}_j=z_m^{-1}\,\vtilde^{(m-1)}_j$.
Because $\DD_0\,\nbigkhat_{m-1}\subset
 \nbigkhat_{m-1}\otimes\Omega^{1,0}(\log D)$,
we have 
$ z_1\deldel_1\vtilde_j^{(m-1)}
-(z_1\del_1\gminia_j)\,\vtilde_j^{(m-1)}
 \in \nbigkhat_{m-1}$.
Therefore,
$\vecz^{\ord(\gminia_j)}\, \vtilde^{(m-1)}_j
 \in\gbigehat^{(\lambda_0)}$.
If the $m$-th component of 
$\ord(\gminia_j)$ is not $0$,
we are done.
Let us consider the case 
in which the $m$-th component of
$\ord(\gminia_j)$ is $0$.
We have the following:
\[
 \DD_0(\del_m)\vtilde^{(m-1)}_j
=\deldel_m(\vtilde^{(m-1)}_j)
-(\del_m\gminia_j)\, \vtilde^{(m-1)}_j
\]
Since we have already known
$\vecz^{\ord(\gminia_j)}\,\vtilde^{(m-1)}_j
 \in \gbigehat^{(\lambda_0)}$,
we have
$(\del_m\gminia_j)\,\vtilde^{(m-1)}_j
 \in \gbigehat^{(\lambda_0)}$.
We also have
$\deldel_m(\vtilde^{(m-1)}_j)
 \in\gbigehat^{(\lambda_0)}$,
and thus
 $\DD_0(\del_m)\vtilde^{(m-1)}_j
 \in\gbigehat^{(\lambda_0)}$.

Let $[v_j]$ denote the induced section of
$\lefttop{m}\Gr^{\Fzero}_0\bigl(
 \nbigqzero_{<\vecdelta}\nbigehat\bigr)$.
If $\nbign_m[v_j]=0$,
we have the following:
\[
 \Bigl(
 \DD_0(\del_m)\vtilde_j^{(m-1)}
-\eigenmap\bigl(\lambda,u_m(v_j)\bigr)\,
 z_m^{-1}\,\vtilde^{(m-1)}_j
 \Bigr)_{|z_m=0}
 \in \lefttop{m}\Fzero_{<1}
 \bigl(\nbigqzero_{\vecdelta}\nbigehat_{|z_m=0}\bigr)
\]
Hence,
$\DD_0(\del_m)\,\vtilde_j^{(m-1)}
-\eigenmap(\lambda,u_m(v_j))\,
 z_m^{-1}\,\vtilde^{(m-1)}_j$
is contained in 
$\nbigkhat_{m-1}\subset\gbigehat^{(\lambda_0)}$,
and we obtain
\[
\vtilde^{(m)}_j=z_m^{-1}\,\vtilde^{(m-1)}_j
 \in\gbigehat^{(\lambda_0)}.
\]
In general,
we have the number $p(v_j)$
determined by the following condition:
\[
\bigl(\Res_m(\DD)
-\eigenmap(\lambda,u_m(v_j))\bigr)^{p(v_j)}
 \bigl(\vtilde^{(m-1)}_j\bigr)=0
\]
\[
  \bigl(\Res_m(\DD)
-\eigenmap(\lambda,u_m(v_j))\bigr)^{p(v_j)-1}
 \bigl(\vtilde^{(m-1)}_j\bigr)\neq 0
\]
The induced section 
$\bigl[
 \DD_0(\del_m)\vtilde_{j}^{(m-1)}
-\eigenmap(\lambda,u_m(v_j))\,
 z_m^{-1}\,\vtilde_j^{(m-1)}
 \bigr]$ 
of $\lefttop{m}\Gr^{\Fzero}_1\bigl(
 \nbigqzero_{\vecdelta}\nbigehat \bigr)$
is contained in the subbundle
generated by
$[\vtilde_q^{(m)}]$ such that
$p(v_q)<p(v_j)$,
by our choice of the frame $\vecv$.
Then, we can show 
$\vtilde^{(m)}_j
 \in \gbigehat^{(\lambda_0)}$
by an induction on $p(v_j)$.
Thus, the induction on $m$ can proceed,
and the proof of Lemma \ref{lem;07.9.7.10}
is finished.
\hfill\qed

\begin{rem}
We have considered 
the $\nbigr$-module associated to
the good wild harmonic bundle
on $X-D$, where $D=\bigcup_{i=1}^n\{z_i=0\}$.
Let us consider the case that
a harmonic bundle  $(E',\delbar_{E'},\theta',h')$
is given on $X-D'$,
where $D'=\bigcup_{i=1}^{\ell}\{z_i=0\}$.
We put $\harmonicbundle:=
 (E',\delbar_{E'},\theta',h')_{|X-D}$.
We have the $\nbigr$-module
$\gbige$ on $\nbigx$
associated to $\harmonicbundle$.

We can construct an $\nbigr$-module
$\gbige'$ on $\nbigx$
from $(E',\delbar_{E'},\theta',h')$
in the same way.
Namely,
let $\vecdelta=(1,\ldots,1)\in\seisuu^{\ell}$
and let $\lambda_0\in\cnum_{\lambda}$.
We consider the $\nbigr$-submodule
$\gbige^{\prime\,(\lambda_0)}$
of $\nbigq\nbige'_{|\nbigx^{(\lambda_0)}}$
generated by
$\nbigqzero_{<\vecdelta}\nbige'$.
By varying $\lambda_0$ and gluing them,
we obtain the $\nbigr$-module $\gbige'$.

Note that we have a natural isomorphism
$\gbige'\lrarr\gbige$.
Indeed, we have natural isomorphisms
$\nbigq\nbige=\nbigq\nbige'(\ast\nbigd)$
and 
$\nbigqzero_{<\vecdelta}\nbige
=\nbigqzero_{<\vecdelta}\nbige'$,
and hence
$\gbige^{(\lambda_0)}
=\gbige^{\prime\,(\lambda_0)}$,
which induces the desired isomorphism.
We will use it implicitly.
\hfill\qed
\end{rem}

\section{Review of some results
 in the tame case}
\label{subsection;08.9.28.2}
\subsection{Some filtrations of $\gbige$}
\label{subsection;08.1.12.25}

Let $X:=\Delta^n$, $D_i:=\{z_i=0\}$ 
and $D:=\bigcup_{i=1}^n D_i$.
We put $\nbar:=\{1,\ldots,n\}$.
Let $\harmonicbundle$ be a tame harmonic bundle on $X-D$.
We have the associated $\nbigr_X$-module 
$\gbige$ on $\nbigx$.
For $\lambda_0\in\cnum$,
let $\nbigx^{(\lambda_0)}$ denote 
a small neighbourhood of
$\{\lambda_0\}\times X$ in $\nbigx$.
Note that 
$\lefttop{\nbar}\Vzero_{\veca}
 \bigl(\naiveprolong{\nbige}\bigr)$ 
and $\naiveprolong{\nbige}$ in \cite{mochi2}
are equal to
$\nbigqzero_{\veca+\vecdelta}(\nbige)$ 
and $\nbigq\nbige$ in this paper, respectively.
\index{sheaf $\naiveprolong{\nbige}$}
\index{sheaf $\lefttop{\nbar}\Vzero_{\veca}
 \bigl(\naiveprolong{\nbige}\bigr)$}
Here $\vecdelta=(1,\ldots,1)$.
We also have
$\lefttop{\nbar}\Vzero_{<0}(\gbige)
=\lefttop{\nbar}\Vzero_{<0}(\naiveprolong{\nbige})
=\nbigqzero_{<\vecdelta}(\nbige)$.
We will not distinguish them in the following argument.
\index{sheaf $\lefttop{\nbar}\Vzero_{<0}(\gbige)$}
\index{sheaf $\lefttop{\nbar}\Vzero_{<0}
 (\naiveprolong{\nbige})$}
\index{sheaf $\nbigqzero_{<\vecdelta}(\nbige)$}

Let $I$ be any subset of $\nbar$.
Let $q_I$ denote the projection of $\real^n$
to $\real^I$.
\index{projection $q_l$}
We have the filtration
$\bigl\{
 \lefttop{I}\Vzero_{\vecb}(\naiveprolong{\nbige})\,\big|\,
 \vecb\in\real^I
 \bigr\}$ of $\naiveprolong{\nbige}_{|\nbigx^{(\lambda_0)}}$,
which is equal to the following
in the notation of this paper:
\[
\lefttop{I}\Vzero_{\vecb}\bigl(\naiveprolong{\nbige}\bigr):=
\bigcup_{\substack{\vecb'\in\real^n\\q_I(\vecb')=\vecb}}
 \nbigqzero_{\vecb'+\vecdelta}\nbige
\]
\index{filtration 
 $\lefttop{I}\Vzero\bigl(\naiveprolong{\nbige}\bigr)$}
\index{sheaf
 $\lefttop{I}\Vzero_{\vecb}
 \bigl(\naiveprolong{\nbige}\bigr)$}
It is also equal to
$\nbigqzero_{\vecb'+\vecdelta}\nbige
\otimes\nbigo\bigl(\ast \sum_{j\not\in I}\nbigd_j\bigr)$
for any $\vecb'\in\real^n$
such that $q_I(\vecb')=\vecb$.
We also have the filtration
$\bigl(
\lefttop{I}\Vzero_{\vecb}\gbige\,\big|\,
\vecb\in\real^I\bigr)$
of $\gbige_{|\nbigx^{(\lambda_0)}}$
given by
\[
 \lefttop{I}\Vzero_{\vecb}(\gbige)
:=
\gbige\cap 
 \lefttop{I}\Vzero_{\vecb}(\naiveprolong{\nbige}).
\]
\index{filtration $\lefttop{I}\Vzero(\gbige)$}
\index{sheaf $\lefttop{I}\Vzero_{\vecb}(\gbige)$}
(See Section 15.2.4 and Corollary 15.63 
of \cite{mochi2}.)
We recall the following lemma.
\begin{lem}
\label{lem;08.1.31.11}
Let $I\subset \nbar$ and $J:=\nbar-I$.
Let $\vecb\in\real_{<0}^I$
and $\vecc\in\real^J$.
We set
\[
  \lefttop{I}
 \Tzero(\vecc,\vecb):=
 \frac{\lefttop{\nbar}\Vzero_{\vecb+\vecc}(\gbige)}
 {\sum_{ \vecd\lneq\vecc}
 \lefttop{\nbar}\Vzero_{\vecb+\vecd}(\gbige)}.
\]
Here,
$\vecd\lneq\vecc$ means
``$\vecd\leq\vecc$ and $\vecd\neq\vecc$''.
We put 
$c'_i:=
 \max\bigl\{c_i-n<0\,\big|\,
 n\in\seisuu_{\geq 0}\bigr\}$ for any $i\in J$,
$\vecc':=(c_i'\,|\,i\in J)$,
and
\[
 \lefttop{I}\Ttildezero(\vecc,\vecb):=
 \Image\Bigl(
 \prod_{i\in J}
 \deldel_i^{c_i-c_i'}:
 \lefttop{J}\Gr_{\vecc'}^{\Vzero}
 \lefttop{I}\Vzero_{\vecb}
 \bigl(\naiveprolong{\nbige}\bigr)
\lrarr
 \lefttop{J}\Gr_{\vecc}^{\Vzero}
 \lefttop{I}\Vzero_{\vecb}
 \bigl(\naiveprolong{\nbige}\bigr)
 \Bigr).
\]
Then, the following holds:
\begin{itemize}
\item
The multiplication of $\deldel_i$ induces
the surjection
$\lefttop{I}\Ttildezero(\vecc-\vecdelta_i,\vecb)
\lrarr
 \lefttop{I}\Ttildezero(\vecc,\vecb)$ if $c_i\geq 0$,
where
$\vecdelta_i$ denotes the element of $\real^{J}$
whose $j$-th element is $0$ $(j\neq i)$
or $1$ $(j=i)$.
\item
We have the natural isomorphism
$\lefttop{I}\Tzero(\vecc,\vecb)
\simeq
 \lefttop{I}\Ttildezero(\vecc,\vecb)$.
\end{itemize}
\end{lem}
\pf
The first claim is clear from the construction.
The second claim is
Lemma 15.46 of \cite{mochi2}.
\hfill\qed

\vspace{.1in}
For a subset $I\subset\nbar$,
let $\nbigr_{X,I}$ denote 
the sheaf of subalgebras of $\nbigr_X$
generated by $\nbigo_{\nbigx}$
and $\deldel_j$ $(j\in I)$.
\index{sheaf $\nbigr_{X,I}$}
We remark the following lemmas,
which implicitly appeared in \cite{mochi2}.
\begin{lem}
\label{lem;07.11.25.5}
Let $I\subset \nbar$
and $\vecb\in\real^I_{<0}$.
Let $\epsilon$ be any sufficiently small
positive number.
Then, we have 
the following equality on $\nbigx^{(\lambda_0)}$:
\begin{equation}
\label{eq;07.9.7.2}
 \lefttop{I}\Vzero_{\vecb}(\gbige)
=
 \nbigr_{X,\nbar-I}\cdot
 \lefttop{\nbar}
 \Vzero_{\vecb
              -\epsilon\vecdelta_{\nbar-I}}
 (\gbige)
\end{equation}
\end{lem}
\pf
It is clear that the right hand side of
(\ref{eq;07.9.7.2}) is contained in the left hand side.
Hence, we have only to show that
$\lefttop{\nbar}\Vzero_{\vecb+\vecc}(\gbige)$
is contained in the right hand side of (\ref{eq;07.9.7.2})
for any $\vecc\in\real^{\nbar-I}$,
where we regard 
$\vecb+\vecc\in\real^I\times\real^{\nbar-I}=\real^n$.

Let $\vecc\in\real^{\nbar-I}$,
and let $f$ be a section of
$\lefttop{\nbar}\Vzero_{\vecb+\vecc}(\gbige)$.
Let us show that $f$ is contained 
in the right hand side of (\ref{eq;07.9.7.2}).
We put $\gminiq(\vecc):=
 \#\bigl\{j\in I\,\big|\,c_j\geq 0\bigr\}$.
We use an induction on $\gminiq(\vecc)$.
If $\gminiq(\vecc)=0$,
there is nothing to prove.
Assume $\gminiq(\vecc)>0$,
and let $c_i>0$.
Due to Lemma \ref{lem;08.1.31.11},
we have a section 
$g\in \lefttop{\nbar}
 \Vzero_{\vecb+\vecc-\vecdelta_i}
 (\gbige)$ such that
\[
 f-\deldel_i g=
 \sum_{\vecc'\lneq \vecc} h_{\vecc'},
\quad
 h_{\vecc'}\in
 \lefttop{\nbar}\Vzero_{\vecb+\vecc'}(\gbige).
\]
We remark that
the set of the parabolic weights 
is discrete.
Hence, we can reduce the number $\gminiq$
after finite steps.
Thus, we obtain (\ref{eq;07.9.7.2}).
\hfill\qed

\begin{lem}
\label{lem;08.1.31.55}
Let $I\subset\nbar$ and $\vecb\in\real_{<0}^I$.
Let $i\in\nbar-I$ and $c\in\real$.
Then, we have the following:
\[
 \lefttop{i}\Vzero_c
 \lefttop{I}
 \Vzero_{\vecb}(\gbige)
=\sum_{(c',m)\in\nbigu}
\deldel_i^m\bigl(
 \lefttop{i}\Vzero_{c'}
  \lefttop{I}
 \Vzero_{\vecb}(\gbige)
\bigr),
\]
\[
\nbigu:=\bigl\{
 (c',m)\in\real_{<0}\times\seisuu_{\geq\,0}\,\big|\,
 c'+m\leq c
 \bigr\}.
\]
\end{lem}
\pf
If $c<0$, there is nothing to prove.
Assume $c\geq 0$.
We have only to show the following equality:
\begin{equation}
 \label{eq;08.1.31.15}
 \lefttop{i}\Vzero_c
 \lefttop{I}
 \Vzero_{\vecb}(\gbige)
=
 \lefttop{i}\Vzero_{<c}
   \lefttop{I}
 \Vzero_{\vecb}(\gbige)
+ \deldel_i\bigl(
 \lefttop{i}\Vzero_{c-1}
   \lefttop{I}
 \Vzero_{\vecb}(\gbige)
\bigr)
\end{equation}
Clearly, the right hand side is contained
in the left hand side.
We put $K:=\nbar-(\{i\}\cup I)$.
Let $\vecd\in\real^{K}$,
and let us show that
$\lefttop{\nbar}\Vzero
 _{\vecb+c\vecdelta_i+\vecd}(\gbige)$
is contained in the right hand side.
Let $\gminiq(\vecd):=
 \#\bigl\{j\,\big|\,d_j\geq 0 \bigr\}$.
We use an induction on $\gminiq(\vecd)$.
In the case
$\gminiq(\vecd)=0$,
i.e.,
$\vecd\in\real_{<0}^{K}$,
the claim easily follows 
from Lemma \ref{lem;08.1.31.11}.
Let $f$ be a section of
$\lefttop{\nbar}
 \Vzero_{\vecb+c\vecdelta_i+\vecd}(\gbige)$.
Due to Lemma \ref{lem;08.1.31.11},
there exists a section $g$ of
$\lefttop{\nbar}
 \Vzero_{\vecb+(c-1)\vecdelta_i+\vecd}(\gbige)$
such that the following holds:
\[
 f-\deldel_ig
=h_{1}
+\sum_{\vecd'\lneq\vecd}
 h_{\vecd'},
\quad\quad
 h_1\in 
 \lefttop{i}\Vzero_{<c}
\lefttop{I}
 \Vzero_{\vecb}(\gbige),
\quad\quad
 h_{\vecd'}
 \in 
 \lefttop{\nbar}
 \Vzero_{\vecb+c\vecdelta_i+\vecd'}
 (\gbige)
\]
Since the set of the parabolic weights is discrete,
we can reduce the number $\gminiq$
after finite steps.
\hfill\qed

\vspace{.1in}
For any subset 
$I\subset \nbar:=\{1,\ldots,n\}$,
we put $D_I:=\bigcap_{i\in I}D_i$,
and let $N_{D_I}^{\ast}X$ 
denote the conormal bundle
of $D_I$ in $X$.

\begin{lem}
\label{lem;08.10.20.10}
$\gbige$ is holonomic,
and the characteristic variety of $\gbige$
is contained in 
$\nbigs:=\bigcup_{I\subset\nbar} 
 (\cnum_{\lambda}\times N_{D_I}^{\ast}X)$.
\end{lem}
\pf
It is shown in 
Proposition 15.68 of \cite{mochi2} and its proof.
We give an outline.
For any $\veca=(a_i)\in\real^n_{\geq \,0}$,
we set $|\veca|=\sum a_i$.
For any non-negative integer $m$,
we set 
\[
 F_m(\gbige):=
 \sum_{\substack{
 \veca\in\real^n_{\geq\,0}\\
 |\veca|\leq m}}
 \lefttop{\nbar}\Vzero_{<\veca}(\gbige)
\]
on $\nbigxzero$.
It can be shown that
$F_m(\nbigr_X)\cdot F_0(\gbige)
=F_m(\gbige)$
by using Lemma \ref{lem;08.1.31.11}.
It implies that $\bigoplus_m G_m(\gbige)$
is finitely generated over 
the Rees algebra $\bigoplus_m F_m(\nbigr_X)$.
Hence, $F(\gbige)$ is a coherent filtration
of $\gbige$.
It is easy to see that
the support of
$\Gr^{F}(\gbige)$ as 
an $\nbigo_{\cnum_{\lambda}\times T^{\ast}X}$-module
is contained in $\nbigs$.
\hfill\qed

\subsection{Push-forward
 $i_{g\dagger}\gbige$ and 
 the $V$-filtration $\Uzero$}

Let $g$ be any function of the form
$\vecz^{\vecp}=\prod_{j\in I} z_j^{p_j}$ 
for $p_j>0$ $(j\in I)$.
Let $i_g:X\lrarr X\times \cnum_t$ denote the graph.
It is very important to investigate
the push-forward $i_{g\dagger}\gbige$
on $\nbigx\times\cnum_t$.
\index{sheaf $i_{g\dagger}\gbige$}
The support of $i_{g\dagger}\gbige$
is the graph of $g$,
which is naturally isomorphic to $\nbigx$.
And,
$i_{g\dagger}\gbige$ is identified with
$(i_{g\ast}\gbige)[\deldel_t]$ on $\nbigx\times\cnum_t$
(or simply denoted by $\gbige[\deldel_t]$),
where $i_{g\dagger}$ (resp. $i_{g\ast}$)
denotes the push-forward
of $\nbigr$-modules (resp. sheaves).
The action of $\nbigr_{X\times\cnum_t}$
is given by general formulas:
\begin{equation}
\label{eq;08.9.25.102}
\begin{array}{ll}
 a\cdot (\deldel_t^j\otimes u)
=\deldel_t^j\otimes (a\cdot u)
  & (a\in\nbigo_X)
 \\
 \deldel_i\cdot (\deldel_t^j\otimes u)
= \deldel_t^j\otimes (\deldel_i u)
-\deldel_t^{j+1}\otimes (\del_i g)\cdot u  \\
 t\cdot (\deldel_t^j\otimes u)
=\deldel_t^j\otimes (g\cdot u)
-j\lambda\deldel_t^{j-1}\otimes u\\
\deldel_t(\deldel_t^j\otimes u)
=\deldel_t^{j+1}\otimes u
\end{array}
\end{equation}
We will implicitly use the following formula
for $i\in I$:
\begin{equation}
\label{eq;08.9.25.101}
 \bigl(
 p_i\deldel_tt+\deldel_{i}z_i
\bigr)(\deldel_t^j\otimes u)
=-p_i\cdot j\cdot \lambda
 \deldel_t^j\otimes u
+\deldel_t^j\otimes\deldel_i (z_i u)
\end{equation}

Let $\pi$ denote the projection
$X\times\cnum_t\lrarr X$.
Let $V_0\nbigr_{X\times\cnum_t}$
denote the sheaf of subalgebras of 
$\nbigr_{X\times\cnum_t}$
generated by
$\pi^{\ast}\nbigr_X$
and $\deldel_tt$.
Recall that we have the $V$-filtration
$\Uzero$ of $i_{g\dagger}\gbige$ on
$\nbigx^{(\lambda_0)}$
given as follows
(Section 16.1 of \cite{mochi2}):
\[
 \Uzero_b\bigl(i_{g\dagger}\gbige\bigr):=
\pi^{\ast}\nbigr_{X}\cdot
\Bigl(
 \lefttop{I}\Vzero_{b\vecp}(\gbige)\otimes 1
\Bigr)\,\,\,
(b<0)
\]
\[
 \Uzero_b\bigl(i_{g\dagger}\gbige\bigr)
=\sum_{\substack{
 c<0,n\in\seisuu_{\geq\,0}\\
 c+n\leq b }}
 \deldel_t^n\cdot
 \Uzero_c\bigl(i_{g\dagger}\gbige\bigr)
 \,\,\,(b\geq 0)
\]
They are locally finitely generated
$V_0\nbigr_{X\times\cnum_t}$-modules,
which immediately follows from 
Lemma \ref{lem;07.11.25.5}.
Since they are contained in
a pseudo-coherent 
$\nbigo_{\nbigx\times\cnum_t}$-module
$i_{g\ast}\nbigq\nbige[\deldel_t]$,
they are $\nbigr_{X\times\cnum_t}$-coherent.
(See Proposition \ref{prop;08.9.26.10}.)
The following property can be checked
by a direct calculation:
\begin{itemize}
\item
 $t\cdot \Uzero_b\subset \Uzero_{b-1}$
 for any $b\in\real$,
and
 $t\cdot \Uzero_b=\Uzero_{b-1}$
 for $b<0$.
\item
 $\deldel_t\Uzero_b\subset \Uzero_{b+1}$
 for any $b\in\real$,
 and 
 the induced morphisms
 $\deldel_t:\Gr^{\Uzero}_{b-1}\lrarr
 \Gr^{\Uzero}_{b}$
 are surjective  for any $b>0$.
\end{itemize}
As in (16.3) of \cite{mochi2},
we set
\begin{equation}
 \label{eq;08.9.26.11}
 \KMS\bigl(
 i_{g\dagger}\gbige^0\bigr)
:=\bigcup_{i\in I}
 \bigl\{u\in\real\times\cnum\,\big|\,
 p_i\cdot u\in \KMS(\nbige^0,i)
 \bigr\}
\end{equation}
\index{set $\KMS\bigl(
 i_{g\dagger}\gbige^0\bigr)$}
\begin{equation}
 \label{eq;08.9.26.12}
 \nbigk\bigl(
 i_{g\dagger}\gbige,
 \lambda_0,b\bigr)
:=\bigl\{
u\in \KMS\bigl(i_{g\dagger}\gbige^0\bigr)
\,\big|\,
 \paramap(\lambda_0,u)=b
 \bigr\}
\end{equation}
\index{set $ \nbigk\bigl(
 i_{g\dagger}\gbige,
 \lambda_0,b\bigr)$}
Here $\gbige^0$
means the specialization of $\gbige$
at $\lambda=0$.
The following lemma is proved 
in \cite{mochi2}.
\begin{lem}
\mbox{{}}\label{lem;08.9.26.30}
\begin{itemize}
\item
The following endomorphism 
is nilpotent on $\Gr^{\Uzero}_b
\bigl(i_{g\dagger}\gbige\bigr)$:
\[
 \prod_{u\in\nbigk(i_{g\dagger}\gbige,
 \lambda_0,b)}
\bigl(-\deldel_tt+\eigenmap(\lambda,u)\bigr)
\]
(See Corollary {\rm 16.13} of 
{\rm \cite{mochi2}}.)
\item
We also obtain that
$\Gr^{\Uzero}_b\bigl(
 i_{g\dagger}\gbige
 \bigr)$ is strict,
i.e., the multiplication of $\lambda-\lambda_1$
is injective for any $\lambda_1$.
(See Proposition {\rm 16.47} of 
 {\rm \cite{mochi2}}.)
\hfill\qed
\end{itemize}
\end{lem}

More strongly,
we have the following proposition.
\begin{prop}[Proposition 16.49 
of \cite{mochi2}]
For any $\lambda_0$,
$i_{g\dagger}\gbige$ is strictly $S$-decomposable
along $t$ at $\lambda_0$ 
with the above $V$-filtration $\Uzero$.
(See Section {\rm\ref{subsection;08.10.18.101}}
for the notion of strict $S$-decomposability.)
\hfill\qed
\end{prop}

\begin{rem}
When we would like to know
some property of
$\psitilde_{g,u}(\gbige)$,
we have only to look at
$\Uzero_{b}$ for $b<0$.
(See Section {\rm\ref{subsection;08.10.18.101}}
for the functor $\psitilde_{g,u}$.)
\hfill\qed
\end{rem}

\subsection{Filtrations $\lefttop{i}\Vzero$}
\label{subsection;08.9.28.20}

Let $\lefttop{J}V_0\nbigr_X$
be as in Section \ref{subsection;08.9.26.1}.
Let $\lefttop{J,t}V_0\nbigr_{X\times\cnum_t}$
denote the sheaf of subalgebras of
$\nbigr_{X\times\cnum_t}$ generated by
$\pi^{\ast}(\lefttop{J}V_0\nbigr_X)$
and $\deldel_tt$.

Let $b<0$.
As in Section 16.1.4 of \cite{mochi2},
for any element
$\vecc\in
 \real^{I}_{\leq 0}\times\real^{\nbar-I}$,
we put
\[
  \lefttop{\nbar}\Vzero_{\vecc}
 \Uzero_b\bigl(i_{g\dagger}\gbige\bigr)
:=\pi^{\ast}(\lefttop{\nbar}V_0\nbigr_X)
 \cdot
 \bigl(\lefttop{\nbar}\Vzero_{\vecc+b\vecp}
 (\gbige)\otimes 1\bigr).
\]
For any $\vecc\in\real^n$,
we put
\[
 \lefttop{\nbar}\Vzero_{\vecc}
 \Uzero_b(i_{g\dagger}\gbige):=
 \sum_{(\vecn,\veca)\in S}
 \deldel^{\vecn}\Bigl(
 \lefttop{\nbar}\Vzero_{\veca}
 \bigl(\Uzero_b(i_{g\dagger}\gbige)\bigr)
\Bigr),
\]
where
$S:=\bigl\{
 (\vecn,\veca)\in\seisuu^I_{\geq\,0}\times
 (\real^I_{\leq 0}\times\real^{\nbar-I})\,\big|\,
 \vecn+\veca\leq\vecc \bigr\}$
and $\deldel^{\vecn}:=\prod \deldel_j^{n_j}$.
Note that they are 
$\lefttop{\nbar,t}V_0\nbigr_{X\times\cnum_t}$-modules,
which can be checked by using
(\ref{eq;08.9.25.102}) and 
(\ref{eq;08.9.25.101}).
For $1\leq i\leq n$ and $c\in\real$,
we put
\[
 \lefttop{i}\Vzero_c\Uzero_b(i_{g\dagger}\gbige)
:=\sum_{\substack{\vecc\in\real^n\\ q_i(\vecc)=c}}
 \lefttop{\nbar}\Vzero_{\vecc}
 \Uzero_b(i_{g\dagger}\gbige),
\]
where $q_i$ denotes the projection of $\real^n$
onto the $i$-th component.
It is easy to check that they are 
finitely generated over 
$\lefttop{i,t}V\nbigr_{X\times\cnum_t}$,
and pseudo-coherent
$\nbigo_{\nbigx\times\cnum_t}$-modules,
and hence coherent
$\lefttop{i,t}V\nbigr_{X\times\cnum_t}$-modules,
as remarked in Proposition \ref{prop;08.9.26.10}
below.
Thus, we obtain the filtration
$\lefttop{i}\Vzero$
on $\Uzero_b(i_{g\dagger}\gbige)$
for $b<0$ by coherent 
$\lefttop{i,t}V
 \nbigr_{X\times\cnum_t}$-modules.

We have the induced filtrations
on $\Gr^{\Uzero}_b(i_{g\dagger}\gbige)$
which are also denoted by $\lefttop{i}\Vzero$.
Since 
the filtrations $\lefttop{i}\Vzero$ are preserved by
the action of $-\deldel_tt$,
they are compatible with the decomposition
$\Gr^{\Uzero}_b(i_{g\dagger}\gbige)
=\bigoplus
 _{\paramap(\lambda_0,u)=b}
 \psizero_{t,u}(i_{g\dagger}\gbige)$.
Hence, the filtrations of
$\psitilde_{t,u}(i_{g\dagger}\gbige)$
on $\nbigx^{(\lambda_0)}$ are induced,
which are also denoted by $\lefttop{i}\Vzero$.
Note that
$\lefttop{i}\Vzero_a\psitilde_{t,u}
 (i_{g\dagger}\gbige)$
are coherent
$\lefttop{i}V_0\nbigr_X[\deldel_tt]$-modules.
Since $-\deldel_tt+\eigenmap(\lambda,u)$
is nilpotent on $\lefttop{i}\Vzero_a\psitilde_{t,u}
 (i_{g\dagger}\gbige)$,
they are also coherent $\lefttop{i}V_0\nbigr_X$-modules.

Let $N$ denote the nilpotent part
of the action of $-\deldel_tt$
on $\psitilde_{t,u}(i_{g\dagger}\gbige)$,
and let $W(N)$ denote the weight filtration.
We obtain the induced filtration
$\lefttop{i}\Vzero$
of $\Gr^{W(N)}_{h}
 \psitilde_{t,u}(i_{g\dagger}\gbige)$
by coherent $\lefttop{i}V_0\nbigr_X$-modules
on $\nbigx^{(\lambda_0)}$.

\begin{prop}
\mbox{{}}\label{prop;08.1.12.50}
\begin{itemize}
\item
 The filtrations $\lefttop{i}\Vzero$
 $(i=1,\ldots,n)$
 are compatible with 
 the primitive decomposition of
 $\Gr^{W(N)}_{h}
 \psitilde_{t,u}(i_{g\dagger}\gbige)$
 (Lemma {\rm 17.37} of {\rm\cite{mochi2}}).
\item
 The primitive part 
 $P\Gr^{W(N)}_h\psitilde_{t,u}(i_{g\dagger}\gbige)$
 is strictly $S$-decomposable along $z_i=0$
 with the filtration  $\lefttop{i}\Vzero$
 for any $i$
 (Corollary {\rm 17.45} and 
 Corollary {\rm 17.55} of {\rm\cite{mochi2}}).
\hfill\qed
\end{itemize}
\end{prop}

We give another description of the filtration
$\lefttop{i}\Vzero$
on $\Uzero_b(i_{g\dagger}\gbige)$
$(b<0)$.
Let $K_1:=I\cup\{i\}$
and $K_2:=\nbar-K_1$.
Let $\vecdelta_j$ denote the element of
$\real^n$ such that
(i) the $j$-th component is $1$, 
(ii) the other component is $0$.
\index{vector $\vecdelta_j$}
For any subset $J\subset \nbar$,
let $\vecdelta_J=\sum_{j\in J}\vecdelta_j$.
\index{vector $\vecdelta_J$}

\begin{lem}
\label{lem;08.1.12.40}
Let $b<0$.
Let $c\leq 0$ in the case $i\in I$,
or $c<0$ in the case $i\not\in I$.
We have the following 
for a sufficiently small $\epsilon>0$:
\begin{equation}
 \label{eq;08.1.31.20}
 \pi^{\ast}\bigl(
 \lefttop{i}V_0\nbigr_{X}\bigr)
 \cdot\bigl(
 \lefttop{\nbar}\Vzero_{b\vecp+c\vecdelta_i
 -\epsilon\vecdelta_{K_2}}(\gbige)\otimes 1
 \bigr)
=\lefttop{i}\Vzero_c
 \Uzero_b\bigl(i_{g\dagger}\gbige\bigr)
\end{equation}
\end{lem}
\pf
Let us consider the case $i\in I$.
Let $\vecc_1\in\real_{\leq 0}^{I}
 =\real_{\leq 0}^{K_1}$.
We have the following:
\begin{multline}
\label{eq;08.1.31.40}
 \sum_{\vecc_2\in\real^{K_2}}
 \lefttop{\nbar}\Vzero_{\vecc_1+\vecc_2}
 \Uzero_b(i_{g\dagger}\gbige)
=\sum_{\vecc_2\in\real^{K_2}}
 \pi^{\ast}\bigl(
 \lefttop{\nbar}V_0\nbigr_X\bigr)\cdot
 \bigl(\lefttop{\nbar}
 \Vzero_{\vecc_1+\vecc_2+b\vecp}
 (\gbige)\otimes 1
 \bigr) \\
=\pi^{\ast}\bigl(
 \lefttop{\nbar}V_0\nbigr_X\bigr)\cdot
 \bigl(
 \lefttop{K_1}\Vzero_{\vecc_1+b\vecp}
 (\gbige)\otimes 1
 \bigr) 
=\pi^{\ast}\bigl(
 \lefttop{\nbar}V_0\nbigr_X\cdot
  \nbigr_{X,K_2}\bigr)\cdot
 \bigl(
 \lefttop{\nbar}\Vzero_{\vecc_1+b\vecp
 -\epsilon\vecdelta_{K_2}}(\gbige)
\otimes 1
 \bigr)
\end{multline}
Here, we have used Lemma \ref{lem;07.11.25.5}.
Even in the case $i\not\in I$,
we have the formally same equality
as (\ref{eq;08.1.31.40})
for $\vecc_1\in
 \real^{I}_{\leq 0}\times\real_{<0}
 \subset \real^{K_1}$.
Then, it is easy to derive (\ref{eq;08.1.31.20}).
\hfill\qed

\begin{lem}
\label{lem;08.1.31.102}
Let $b<0$.
In the case $c>0$ and $i\in I$,
we have the following description:
\begin{equation}
\lefttop{i}\Vzero_c
 \Uzero_b\bigl(i_{g\dagger}\gbige\bigr)
=\sum_{(c',m)\in\nbigu}
 \deldel_i^m\bigl(
 \lefttop{i}\Vzero_{c'}
 \Uzero_b\bigl(i_{g\dagger}\gbige\bigr)
 \bigr)
\end{equation}
\[
  \nbigu=\bigl\{
 (c',m)\in\real_{\leq 0}\times\seisuu_{\geq 0}
 \,\big|\,
 c'+m\leq c
 \bigr\}
\]
In the case $c\geq 0$ and $i\not\in I$,
we have the following description:
\begin{equation}
 \label{eq;08.1.31.56}
\lefttop{i}\Vzero_c
 \Uzero_b\bigl(i_{g\dagger}\gbige\bigr)
=\sum_{(c',m)\in\nbigu}
 \deldel_i^m\bigl(
 \lefttop{i}\Vzero_{c'}
 \Uzero_b\bigl(i_{g\dagger}\gbige\bigr)
 \bigr)
\end{equation}
\[
  \nbigu=\bigl\{
 (c',m)\in\real_{<0}\times\seisuu_{\geq 0}
 \,\big|\,
 c'+m\leq c
 \bigr\}
\]
\end{lem}
\pf
The case $i\in I$ is easy
by the construction of $\lefttop{i}\Vzero$
on $\Uzero_b(i_{g\dagger}\gbige)$.
Let us consider the case $i\not\in I$.
For $\vecc_2\in\real^{K_2}$,
we have the following:
\begin{multline}
\sum_{\vecc_2}
 \lefttop{\nbar}
 \Vzero_{c\vecdelta_i+\vecc_2}
 \Uzero_b\bigl(i_{g\dagger}\gbige\bigr)
=\sum_{\vecc_2}
 \lefttop{\nbar}V_0\nbigr_X
 \cdot
 \bigl(\lefttop{\nbar}
 \Vzero_{c\vecdelta_i+\vecc_2+b\vecp}
 (\gbige)\otimes 1
 \bigr) \\
=\lefttop{\nbar}
 V_0\nbigr_X\cdot
 \bigl(
 \lefttop{K_1}\Vzero_{c\vecdelta_i+b\vecp}
 (\gbige)\otimes 1
 \bigr) 
=\sum_{(c',m)\in\nbigu}
 \deldel_i^m\Bigl(
 \lefttop{\nbar}V_0\nbigr_X
\cdot\bigl(
 \lefttop{K_1}
 \Vzero_{c'\vecdelta_i+b\vecp}
 (\gbige)\otimes 1
 \bigr)
 \Bigr) \\
=\sum_{(c',m)\in\nbigu}
 \deldel_i^m
 \Bigl(
 \lefttop{\nbar}V_0\nbigr_X
\cdot \nbigr_{X,K_2}
 \bigl(
 \lefttop{\nbar}\Vzero_{c'\vecdelta_i+b\vecp
  -\epsilon\vecdelta_{K_2}}(\gbige)\otimes 1
 \bigr)
 \Bigr)
\end{multline}
Here, we have used Lemma \ref{lem;08.1.31.55},
and we omit to denote $\pi^{\ast}$.
By using the description in Lemma \ref{lem;08.1.12.40},
we obtain (\ref{eq;08.1.31.56}).
\hfill\qed

\section{The case in which $\Irr(\theta)$ 
 consists of only one element}
\label{subsection;08.9.28.3}
\subsection{The splitting of
 the associated $\nbigr_X$-module
into the tensor product}
\label{subsection;08.9.29.203}

Let $X$, $D$,
and $\harmonicbundle$ be as in 
Section \ref{subsection;08.1.12.25}.
Let $\vecm\in\seisuu_{\leq 0}^{n}-\{\veczero\}$.
We put $\supp(\vecm):=\{j\,|\,m_j<0\}$.
Let $\gminia$ be a meromorphic function
of the form $\prod z_j^{m_j}\cdot \gminia'$
where $\gminia'$ is holomorphic and 
nowhere vanishing.
Note $\ord(\gminia)=\vecm$.
Let $L(\gminia)$ be 
the unramifiedly good wild harmonic bundle 
on $X-D$,
which consists of the line bundle $\nbigo_{X-D}$
with the Higgs field $d\gminia$
and the trivial hermitian metric.
We have the unramifiedly good wild harmonic bundle
$(E',\delbar_{E'},\theta',h'):=
 (E,\delbar_E,\theta,h)\otimes L(\gminia)$.
The associated 
coherent $\nbigr_X$-module is denoted by
$\gbige'$.

Let $\nbigl(\gminia)$ be 
the coherent $\nbigr_{X}$-module
$\nbigo_{\nbigx}
 \bigl(\ast \prod_{i\in \supp(\vecm)}z_i\bigr)\cdot e$
with $\deldel_j e=(\del_j\gminia)\cdot e$,
which is the $\nbigr_X$-module associated to 
the unramifiedly
good wild harmonic bundle $L(\gminia)$.
\index{sheaf $\nbigl(\gminia)$}
\index{harmonic bundle $L(\gminia)$}
We obtain the $\nbigr_X$-module
$\gbige\otimes_{\nbigo_{\nbigx}}\nbigl(\gminia)$.

\begin{lem}
 \mbox{{}}\label{lem;07.11.25.2}
\begin{itemize}
\item
 On $\nbigx^{(\lambda_0)}$,
 $\gbige\otimes\nbigl(\gminia)$ is generated
 by $\lefttop{\nbar}\Vzero_{<0}(\gbige)\otimes e$
 over $\pi^{\ast}\nbigr_X$.
\item
 We have the natural isomorphism
 $\gbige'\simeq\gbige\otimes\nbigl(\gminia)$.
\end{itemize}
\end{lem}
\pf
Since we have the natural isomorphisms
$\nbigq\nbige'\simeq
\nbigq\nbige\otimes\nbigl(\gminia)$ and
$\nbigqzero_{<\vecdelta}\nbige'
\simeq 
 \lefttop{\nbar}\Vzero_{<0}(\gbige)\otimes e$,
the second claim follows from
the first claim and the definition of $\gbige'$.

Let us consider the morphism
\[
\Phi:
\pi^{\ast}\nbigr_{X}\otimes
 \bigl(\lefttop{\nbar}\Vzero_{<0}(\gbige)\otimes e \bigr)
\lrarr \gbige\otimes\nbigl(\gminia)
\]
induced by the inclusion
$\lefttop{\nbar}\Vzero_{<0}(\gbige)\otimes e
\subset \gbige\otimes\nbigl(\gminia)$
and the $\nbigr_X$-action.
Let $F_m(\nbigr_X)$ denote the submodule of
$\nbigr_X$ which consists of the differential operators
of at most order $m$.
On $\nbigx^{(\lambda_0)}$,
we define
\[
 F_m(\gbige):=\Image\bigl(
 F_m(\nbigr_{X})\otimes
 \lefttop{\nbar}\Vzero_{<0}\gbige
\lrarr \gbige
 \bigr).
\]
Since $\bigcup_m F_m(\gbige)=\gbige$
on $\nbigx^{(\lambda_0)}$
by construction of $\gbige$,
we have only to show
$F_m(\gbige)\otimes\nbigl(\gminia)\subset \Image\Phi$
for any $m$.
Let us consider the following claims:
\begin{description}
\item[$(a_m)$]
 $F_m\bigl(\gbige\bigr)\otimes e\subset \Image\Phi$.
\item[$(b_m)$]
 $F_m\bigl(\gbige\bigr)\otimes\nbigl(\gminia)
 \subset \Image\Phi$.
\end{description}
The claim $(a_0)$ clearly holds.
Let us show $(a_{m})\Longrightarrow (b_m)$.
We have only to show
$\vecz^{N\vecm}
 \bigl(F_m(\gbige)\otimes e\bigr)
\subset\Image(\Phi)$ for any $N$,
which we show by an induction on $N$.
In the case $N=0$,
the claim directly follows from $(a_m)$.
Assume we have already obtained
$\vecz^{N\vecm}
 \bigl(F_m(\gbige)\otimes e\bigr)
\subset\Image(\Phi)$.
Take $i\in \supp(\vecm)$. 
Let $f\in F_m(\gbige)$.
We have the following:
\[
 \Image(\Phi)\ni
 z_i\deldel_i\bigl(
 \vecz^{N\vecm}\cdot f\otimes e
 \bigr)
=z_i\deldel_i(\vecz^{N\vecm}\cdot f)\otimes e
+(\vecz^{N\vecm}\cdot f)
 \otimes (z_i\del_i\gminia)\cdot e
\]
Because 
 $z_i\deldel_i(\vecz^{N\vecm}\cdot f)
 \in F_m(\gbige)$,
we have 
$z_i\deldel_i(\vecz^{N\vecm}\cdot f)\otimes e
  \in \Image(\Phi)$.
Hence, we obtain
$ \vecz^{N\vecm}f \otimes 
 (z_i\del_i\gminia)\cdot e\in\Image(\Phi)$,
which implies 
$\vecz^{(N+1)\vecm}f\otimes e\in\Image(\Phi)$.
Therefore, we obtain $(b_m)$.

Let us show $(b_{m-1})\Longrightarrow(a_m)$.
For $f\in F_m(\gbige)$,
we have
$\deldel_i(f\otimes e)=\deldel_i f\otimes e+
 (\del_i\gminia)f\otimes e$.
By the assumption,
$\deldel_i(f\otimes e)$ and
$(\del_i\gminia)f\otimes e$ are contained in
$\Image(\Phi)$.
Hence, $\deldel_if\otimes e$ is also contained in
$\Image(\Phi)$.
Then, $(a_m)$ follows.
Thus, the proof of Lemma \ref{lem;07.11.25.2}
is finished.
\hfill\qed

\vspace{.1in}
For any subset 
$I\subset \nbar:=\{1,\ldots,n\}$,
we put $D_I:=\bigcap_{i\in I}D_i$,
and let $N_{D_I}^{\ast}X$ 
denote the conormal bundle
of $D_I$ in $X$.

\begin{cor}
\label{cor;08.1.21.20}
$\gbige'$ is holonomic.
The characteristic variety of 
$\gbige'$ is contained in
$\nbigs:=\bigcup_{I\subset\nbar} 
 (\cnum_{\lambda}\times N_{D_I}^{\ast}X)$.
\end{cor}
\pf
We use the notation in the proof of
Lemma \ref{lem;08.10.20.10}.
We set $F_m(\nbigl(\gminia)):=
 F_m(\nbigr_X)\cdot e$.
Let $F_m\bigl(
 \gbige\otimes\nbigl(\gminia)
 \bigr)$ be the image of the following
naturally defined map:
\[
 \bigoplus_{p+q\leq m}
 F_p(\gbige)\otimes 
 F_q\bigl(\nbigl(\gminia)\bigr)
\lrarr
 \gbige\otimes\nbigl(\gminia)
\]
Let us show that
$F_m\bigl(\gbige\otimes\nbigl(\gminia)\bigr)$
is finitely generated.
Let $\nbigh$ denote the image of
\[
F_1(\nbigr_X)\cdot
 F_{m-1}\bigl(\gbige\otimes\nbigl(\gminia)\bigr)
\lrarr F_m\bigl(\gbige\otimes\nbigl(\gminia)\bigr)
\]
We have only to show 
$\nbigh=F_m\bigl(\gbige\otimes\nbigl(\gminia)\bigr)$.
A section $f$ of
$F_m\bigl(\gbige\otimes\nbigl(\gminia)\bigr)$
has an expression 
$\sum_{p+q\leq m} a_p\otimes b_q$,
where $a_p\in F_p(\gbige)$
and $b_q\in F_q(\nbigl(\gminia))$.
There exist sections
$b_m'\in F_{m-1}(\nbigl(\gminia))$
and $v\in F_1(\nbigr_X)$
such that $v\cdot b_m'=b_m$.
Then, $f$ is equivalent to
$\sum_{q\leq m-2}a_p\otimes b_q
+\bigl(
 a_1\otimes b_{m-1}+(va_0)\otimes b_m'
 \bigr)$ modulo $\nbigh$.
By an easy descending induction,
it is shown that $f$ is equivalent to
a section of
$F_m(\gbige)\otimes F_0\bigl(\nbigl(-\gminia)\bigr)$
modulo $\nbigh$.

Take $i\in\supp(\vecm)$.
Any section of 
$F_m(\gbige)\otimes F_0\bigl(\nbigl(-\gminia)\bigr)$
has an expression
\[
 a_m\otimes
 \bigl(\vecz^{-\vecm}\cdot z_i\del_i\gminia
 \cdot e\bigr).
\]
It is equivalent to
$-z_i\del_i a_m\otimes (\vecz^{-\vecm}e)
+a_m\otimes (m_i\vecz^{-\vecm}e)$
modulo $\nbigh$.
Hence, we have only to show 
$\vecz^{-\vecm}\cdot
 F_m\bigl(\gbige\bigr)\otimes e$
is contained in $\nbigh$.
It can be shown by using Lemma \ref{lem;08.1.31.11}
as in the proof of Lemma \ref{lem;07.11.25.5}.

Therefore, $F_m$ gives a coherent filtration of
$\gbige\otimes\nbigl(\gminia)$,
due to Proposition \ref{prop;08.10.20.21}.
It is easy to check that
the support of
$\Gr^F\bigl(\gbige\otimes\nbigl(\gminia)\bigr)$
as an $\nbigo_{\cnum_{\lambda}\times T^{\ast}X}$-module
is contained in $\nbigs$.
\hfill\qed

\subsection{Strictly $S$-decomposability of
 the associated $\nbigr_X$-modules
 along a monomial function}
\label{subsection;08.1.11.1}

Let $1\leq l\leq n$.
We put $\lbar:=\{1,\ldots,l\}$.
Let $g=\vecz^{\vecp}$
for some $\vecp\in\seisuu_{> 0}^{l}$.
We have the natural isomorphism
$i_{g\dagger}\bigl(\gbige\otimes\nbigl(\gminia)\bigr)
\simeq
 i_{g\dagger}\gbige\otimes
 \pi^{\ast}\nbigl(\gminia)$,
where $\pi:X\times\cnum_t\lrarr X$ be 
the natural projection.
Take $\lambda_0\in\cnum_{\lambda}$.
Let us consider the following for any $b\in\real$
on $\nbigx^{(\lambda_0)}\times\cnum_t$:
\begin{equation}
 \label{eq;10.6.5.1}
  \Uzero_b\bigl(
 i_{g\dagger}(\gbige\otimes\nbigl(\gminia))
 \bigr)
:=\Uzero_b\bigl(i_{g\dagger}\gbige \bigr)
 \otimes\pi^{\ast}\nbigl(\gminia)
\end{equation}
Note that we have the isomorphisms:
\begin{equation}
 \label{eq;08.1.12.15}
 \Gr^{\Uzero}_b\bigl(
 i_{g\dagger}(\gbige\otimes\nbigl(\gminia))
 \bigr)
\simeq
 \Gr_b^{\Uzero}\bigl(i_{g\dagger}\gbige\bigr)
\otimes\nbigl(\gminia)
\end{equation}

\begin{prop}
\mbox{{}}\label{prop;08.1.12.20}
\begin{itemize}
\item
$i_{g\dagger}\bigl(\gbige\otimes\nbigl(\gminia)\bigr)$
is strictly $S$-decomposable along $t$ 
at any $\lambda_0$,
and $\Uzero\bigl(i_{g\dagger}\gbige
 \otimes\nbigl(\gminia)\bigr)$
is the $V$-filtration.
As a result,
$\gbige\otimes\nbigl(\gminia)$ is
strictly $S$-decomposable along $g$.
\item
We have natural isomorphisms:
\[
 \psitilde_{t,u}
 \bigl(i_{g\dagger}\bigl(
 \gbige\otimes\nbigl(\gminia)
 \bigr)\bigr)
\simeq
 \psitilde_{t,u}
 \bigl(i_{g\dagger}\gbige\bigr)\otimes
 \nbigl(\gminia)
\]
Under the isomorphisms,
the nilpotent parts $N$ of $-\deldel_tt$ are the same.
In particular,
we have natural isomorphisms
for the primitive parts:
\[
 P\Gr_p^{W(N)}\psitilde_{t,u}
 \bigl(i_{g\dagger}\bigl(
 \gbige\otimes\nbigl(\gminia)
 \bigr)\bigr)
\simeq
 P\Gr_p^{W(N)}\psitilde_{t,u}
 \bigl(i_{g\dagger}\gbige\bigr)\otimes
 \nbigl(\gminia)
\]
\end{itemize}
\end{prop}
\pf
The second claim follows from
the first claim and the isomorphism 
(\ref{eq;08.1.12.15}).
By construction of the filtration,
we have the following:
\begin{multline}
 t\cdot \Uzero_a\bigl(
 i_{g\dagger}\bigl(\gbige\otimes\nbigl(\gminia)\bigr)
 \bigr)
=t\cdot\Uzero_a\bigl(i_{g\dagger}\gbige\bigr)
 \otimes\pi^{\ast}\nbigl(\gminia)
\subset \\
 \Uzero_{a-1}\bigl(i_{g\dagger}\gbige\bigr)
 \otimes\pi^{\ast}\nbigl(\gminia)
=\Uzero_{a-1}\bigl(
 i_{g\dagger}\bigl(\gbige\otimes\nbigl(\gminia)\bigr)
 \bigr)
\end{multline}
\begin{multline}
\deldel_t\cdot \Uzero_a\bigl(
 i_{g\dagger}\bigl(\gbige\otimes\nbigl(\gminia)\bigr)
 \bigr)
=
 \deldel_t\cdot\Uzero_{a}\bigl(i_{g\dagger}\gbige\bigr)
 \otimes\pi^{\ast}\nbigl(\gminia)
\subset \\
 \Uzero_{a+1}\bigl(i_{g\dagger}\gbige\bigr)
 \otimes\pi^{\ast}\nbigl(\gminia)
=\Uzero_{a+1}\bigl(
 i_{g\dagger}\bigl(\gbige\otimes\nbigl(\gminia)\bigr)
 \bigr)
\end{multline}
Moreover, we have
$t\cdot \Uzero_a=\Uzero_{a-1}$ $(a<0)$
and 
 $\deldel_t:\Gr^{\Uzero}_a\simeq
 \Gr^{\Uzero}_{a+1}$ $(a>-1)$.

From Lemma \ref{lem;08.9.26.30}
and the isomorphism (\ref{eq;08.1.12.15}),
we obtain that 
the following endomorphism 
is nilpotent on $\Gr^{\Uzero}_b
\bigl(i_{g\dagger}\gbige\otimes\nbigl(\gminia)\bigr)$:
\[
 \prod_{u\in\nbigk(i_{g\dagger}\gbige,\lambda_0,b)}
\bigl(-\deldel_tt+\eigenmap(\lambda,u)\bigr)
\]
We also obtain that
$\Gr^{\Uzero}_b\bigl(
 i_{g\dagger}\gbige\otimes\nbigl(\gminia)
 \bigr)$ is strict.

Let us prepare
to show that 
$\Uzero_b\bigl(i_{g\dagger}
 \gbige\otimes\nbigl(\gminia)\bigr)$
are $V_0\nbigr_{X\times\cnum_t}$-coherent.
We have only to consider the case $b<0$.
We have the natural inclusion
$\lefttop{\nbar}\Vzero_{
 b\vecp-\epsilon\vecdelta_{\nbar-\lbar}}
 (\gbige)\otimes 1
\subset \Uzero_b(i_{g\dagger}\gbige)$.
It induces the inclusion
\[
 \bigl(
 \lefttop{\nbar}\Vzero_{
 b\vecp-\epsilon\vecdelta_{\nbar-\lbar}}
 (\gbige)\otimes 1
\bigr)
 \otimes e
\subset 
 \Uzero_b\bigl(
 i_{g\dagger}\gbige\otimes \nbigl(\gminia)\bigr).
\]
Then, the $V_0\nbigr_{X\times\cnum_t}$-action
induces the following morphism:
\[
 \Phi_b:V_0\nbigr_{X\times\cnum_t}\otimes
\Bigl(
 \lefttop{\nbar}\Vzero_{
   b\vecp-\epsilon\vecdelta_{\nbar-\lbar}}(\gbige)
\otimes 1
\otimes e
\Bigr)
\lrarr
 \Uzero_{b}
 \bigl(i_{g\dagger}\gbige\otimes\nbigl(\gminia)\bigr)
\]
\begin{lem}
\label{lem;08.1.12.16}
$\Phi_b$ is surjective.
\end{lem}
\pf
For $J=(j_1,\ldots,j_n)$, we put
$\deldel^J=\prod_{i=1}^n\deldel_i^{j_i}$
and $|J|=\sum j_i$.
Let $F_m(V_0\nbigr_{X\times \cnum_t})$
be as follows:
\[
 F_m(V_0\nbigr_{X\times \cnum_t}):=
 \Bigl\{
 \sum_{|J|\leq m}a_J(\vecz,t,\deldel_tt)\cdot\deldel^{J}
 \Bigr\}
\]
The inclusion 
$\lefttop{\nbar}
 \Vzero_{b\vecp-\epsilon\vecdelta_{\nbar-\lbar}}
 (\gbige)\otimes 1\subset \Uzero_b(i_{g\dagger}\gbige)$
and the $\nbigr_{X\times\cnum_t}$-action
induce the following map:
\[
F_m\bigl(V_0\nbigr_{X\times\cnum_t}\bigr)
\otimes
\bigl(
 \lefttop{\nbar}
 \Vzero_{b\vecp-\epsilon\vecdelta_{\nbar-\lbar}}
(\gbige)\otimes 1
 \bigr)
\lrarr \Uzero_b(i_{g\dagger}\gbige)
\]
The image is denoted by
$F_m\Uzero_b(i_{g\dagger}\gbige)$.
\begin{lem}
 \label{lem;07.9.7.3}
\mbox{{}}
\begin{itemize}
\item
We have 
$\deldel_i\cdot F_m\Uzero(i_{g\dagger}\gbige)
\subset F_{m+1}\Uzero(i_{g\dagger}\gbige)$.
\item
We have
$\deldel_iz_i\cdot F_m\Uzero(i_{g\dagger}\gbige)
 \subset F_{m}\Uzero(i_{g\dagger}\gbige)$.
\end{itemize}
\end{lem}
\pf
The first claim is clear.
For $f\in \lefttop{\nbar}
 \Vzero_{b\vecp-\epsilon\vecdelta_{\nbar-\lbar}}(\gbige)$,
we have 
\[
 \deldel_iz_i\cdot P(\vecz,t,\deldel_tt)\deldel^J(f\otimes 1)
- P(\vecz,t,\deldel_tt)\deldel^J\cdot
 \deldel_iz_i(f\otimes 1)
\in
 F_m\Uzero_b(i_{g\dagger}\gbige)
\]
We also have
$\deldel_iz_i(f\otimes 1)
=(\deldel_iz_if)\otimes 1
-p_i\deldel_t t (f\otimes 1)$,
where $p_i:=0$ if $i\not\in\lbar$.
Thus the second claim of 
Lemma \ref{lem;07.9.7.3} follows.
\hfill\qed

\vspace{.1in}

Let us return to the proof of Lemma \ref{lem;08.1.12.16}.
According to Lemma \ref{lem;07.11.25.5},
we have the following:
\[
 \lefttop{\lbar}\Vzero_{b\vecp}(\gbige)\otimes 1
\subset
 \bigcup_m F_m\Uzero_b(i_{g\dagger}\gbige)
\]
Since $\Uzero_b(i_{g\dagger}\gbige)$
is generated by
$\lefttop{\lbar}\Vzero_{b\vecp}(\gbige)\otimes 1$
over $\nbigr_X$,
we obtain
\[
 \bigcup_mF_m\Uzero_b(i_{g\dagger}\gbige)
=\Uzero_b(i_{g\dagger}\gbige).
\]
Therefore,
to show the surjectivity of $\Phi_b$,
we have only to show
$F_m\Uzero_b(i_{g\dagger}\gbige)
\otimes \pi^{\ast}\nbigl(\gminia)
\subset\Image(\Phi_b)$
for any $m$.

Let us consider the following claims:
\begin{description}
\item[$(P_m)$]
 $F_m\Uzero_b\bigl(i_{g\dagger}\gbige\bigr)\otimes e
 \subset \Image(\Phi_b)$.
\item[$(Q_m)$]
 $F_m\Uzero_b\bigl(i_{g\dagger}\gbige\bigr)\otimes \nbigl(\gminia)
 \subset \Image(\Phi_b)$.
\end{description}
The claim $(P_0)$ is clear.
(Note the action of $\deldel_t$ on $\nbigl(\gminia)$
 is trivial.)
Let us show $(P_m)\Longrightarrow (Q_m)$.
We show 
$\vecz^{N\vecm}
\cdot F_m\Uzero_b(i_{g\dagger}\gbige)\otimes e
\subset \Image(\Phi_b)$
by an induction on $N$.
The case $N=0$ is $(P_m)$.
Let us show $N\Longrightarrow N+1$.
Take $i\in \supp(\vecm)$.
Let $f\in F_m\Uzero_b(i_{g\dagger}\gbige)$.
We have the following:
\[
  \Image(\Phi_b)
\ni
 z_i\deldel_i(\vecz^{N\vecm}\, f\otimes e) \\
=Nm_i\,
 \vecz^{N\vecm}\, f\otimes e
+\vecz^{N\vecm}(z_i\deldel_i f)\otimes e
+\vecz^{N\vecm}\, f\otimes (z_i\del_i\gminia)\, e
\]
We have 
$\vecz^{N\vecm}\cdot f\otimes e
 \in \Image(\Phi_b)$.
By Lemma \ref{lem;07.9.7.3},
we have $z_i\deldel_i f\in 
 F_m\Uzero_b(i_{g\dagger}\gbige)$,
and hence
$ \vecz^{N\vecm}(z_i\deldel_i f)\otimes e
 \in \Image(\Phi_b)$.
We obtain
$\vecz^{N\vecm}f\otimes 
 (z_i\del_i\gminia)\cdot e
 \in \Image(\Phi_b)$,
which implies
\[
 \vecz^{(N+1)\vecm}f\otimes e
 \in \Image(\Phi_b).
\]
Thus, we obtain $(Q_m)$.

Let us show $(Q_m)\Longrightarrow (P_{m+1})$.
Let $f\in F_m\Uzero_b(i_{g\dagger}\gbige)$.
We have 
$ \deldel_i(f\otimes e)
=(\deldel_if)\otimes e
+f\otimes\deldel_ie$.
By the assumption,
we have $\deldel_i(f\otimes e)\in\Image(\Phi_b)$
and $f\otimes\deldel_ie\in\Image(\Phi_b)$,
and thus
$(\deldel_if)\otimes e\in\Image(\Phi_b)$.
Then, the claim $(P_{m+1})$ follows.
Thus, the induction can proceed,
and the proof of Lemma \ref{lem;08.1.12.16}
is finished.
\hfill\qed

\vspace{.1in}

Let us check that
$\Uzero_b\bigl(i_{g\dagger}
 \gbige\otimes\nbigl(\gminia)\bigr)$
is $V_0\nbigr_{X\times\cnum_t}$-coherent.
We have the inclusion
\[
 \Uzero_b\bigl(i_{g\dagger}
 \gbige\otimes\nbigl(\gminia)\bigr)
\subset 
 i_{g\,\dagger}\bigl(
 \nbigq\nbige
 \bigr)\otimes\pi^{\ast}\nbigl(\gminia).
\]
Hence, it is easy to check that
$\Uzero_a\bigl(i_{g\dagger}\gbige\otimes\nbigl(\gminia)\bigr)$
is a pseudo-coherent $\nbigo_{\nbigx\times\cnum_t}$-module.
According to Lemma \ref{lem;08.1.12.16},
$\Uzero_a\bigl(i_{g\dagger}\gbige
 \otimes\nbigl(\gminia)\bigr)$
is finitely generated over $V_0\nbigr_{X\times\cnum_t}$.
Thus, the desired coherence follows
from Proposition \ref{prop;08.9.26.10}.
\hfill\qed

\begin{rem}
Let $\nbigm$ be an unramifiedly
good meromorphic flat bundle
with a unique irregular value.
By applying the above argument,
we obtain that the $V$-filtration of $\nbigm$
along $g$
is described as {\rm(\ref{eq;10.6.5.1})}.
\hfill\qed
\end{rem}

\subsection{Decomposition by the support}
\label{subsection;08.1.12.51}

Let $\nbar=\{1,\ldots,n\}$.
For each $J\subset \nbar$,
we put $D_J:=\bigcap_{j\in J}D_j$
and $D_J^{\circ}:=
 D_J\setminus \bigcup_{i\not \in J}(D_{J}\cap D_i)$.
Recall the decomposition
in Proposition 17.56 in \cite{mochi2}:
\[
 P\Gr^{W(N)}_p\psitilde_{t,u}(i_{g\dagger}\gbige)
=\bigoplus_{J\subset\nbar} \nbigm_{p,u,J}
\]
Here, $\nbigm_{p,u,J}$ are the $\nbigr_{X}$-modules
which come from 
tame harmonic bundles
$(E_{p,u,J},\delbar_{p,u,J},
 \theta_{p,u,J},h_{p,u,J})$
on $D_J^{\circ}$.
The strict supports of $\nbigm_{p,u,J}$ are 
$\nbigd_J:=\cnum_{\lambda}\times D_J$.
According to Proposition \ref{prop;08.1.12.20},
we obtain the following decomposition:
\begin{equation}
 \label{eq;07.9.7.37}
 P\Gr^{W(N)}_p\psitilde_{t,u}
 \bigl(i_{g\dagger}\gbige\otimes\nbigl(\gminia)\bigr)
=\bigoplus_{J\subset \nbar-I}
 \nbigm_{p,u,J}\otimes\nbigl(\gminia)
\end{equation}

We use the notation in
Section \ref{subsection;08.9.28.20}.
For any $b<0$ and $c\in\real$,
we consider
the $\lefttop{i,t}V_0\nbigr_{X\times\cnum_t}$-module:
\[
 \lefttop{i}\Vzero_c
  \Uzero_b\bigl(
 i_{g\dagger}(\gbige\otimes\nbigl(\gminia))
 \bigr)
:=\lefttop{i}\Vzero_c
 \Uzero_b\bigl(i_{g\dagger}\gbige \bigr)
 \otimes\pi^{\ast}\nbigl(\gminia)
\]

We would like to show that
$\nbigm_{p,u,J}\otimes\nbigl(\gminia)$ 
are strictly $S$-decomposable
along $z_i$ at $\lambda_0$
with the induced filtration $\lefttop{i}\Vzero$.

\begin{lem}
\label{lem;08.9.26.50}
The induced subsheaves
$\lefttop{i}\Vzero_c$ 
of $\Gr^{\Uzero}_b(i_{g\dagger}\gbige
 \otimes\nbigl(\gminia))$
are $\lefttop{i}V_0\nbigr_X$-coherent.
\end{lem}
\pf
Let $K:=\nbar-(\lbar\cup\{i\})$.
Let $b<0$.
Let $c\leq 0$ in the case $i\leq l$,
or $c<0$ in the case $l\leq i\leq n$.
We have the following naturally defined morphism:
\[
 \Phi:
 \lefttop{i,t}V_0\nbigr_{X\times\cnum_t}
 \otimes
 \bigl(
 \lefttop{\nbar}
 \Vzero_{b\vecp+c\vecdelta_i-\epsilon\vecdelta_K}
 (\gbige)\otimes 1\otimes e
 \bigr)
\lrarr
 \lefttop{i}\Vzero_c\Uzero_b\bigl(
 i_{g\dagger}\gbige\otimes\pi^{\ast}\nbigl(\gminia)
 \bigr)
\]
\begin{lem}
\label{lem;08.1.12.45}
The map $\Phi$ is surjective.
\end{lem}
\pf
Since the argument is similar to 
that in the proof of
Lemma \ref{lem;08.1.12.16},
we give only an outline.
Let 
 $F_m\lefttop{i,t}V_0\nbigr_{X\times\cnum_t}$
denote the intersection of
$F_mV_0\nbigr_{X\times\cnum_t}$ and
$\lefttop{i,t}V_0\nbigr_{X\times\cnum_t}$.
For $b$ and $c$ as above,
we have the following naturally defined map:
\[
 F_m\bigl(
 \lefttop{i,t}V_0\nbigr_{X\times\cnum_t}
 \bigr)
\otimes\bigl(
 \lefttop{\nbar}
 \Vzero_{b\vecp+c\vecdelta_i-\epsilon\vecdelta_K}
 (\gbige)\otimes 1
 \bigr)
\lrarr
 \lefttop{i}\Vzero_c
 \Uzero_b\bigl(i_{g\dagger}\gbige\bigr)
\]
The image is denoted by
$F_m\lefttop{i}\Vzero_c\Uzero_b\bigl(
 i_{g\dagger}\gbige\bigr)$.
As in Lemma \ref{lem;07.9.7.3},
the following holds:
\begin{itemize}
\item
 For $j\neq i$,
 we have 
 $\deldel_j F_m\lefttop{i}\Vzero_c
 \Uzero_b(i_{g\dagger}\gbige)
 \subset
 F_{m+1}\lefttop{i}\Vzero_{c}
 \Uzero_b(i_{g\dagger}\gbige)$.
\item
 For any $j$,
 we have 
 $\deldel_jz_j F_m\lefttop{i}\Vzero_c
 \Uzero_b(i_{g\dagger}\gbige)
 \subset
 F_m\lefttop{i}\Vzero_c\Uzero_b
 (i_{g\dagger}\gbige)$.
\end{itemize}

By the description in Lemma \ref{lem;08.1.12.40},
$\lefttop{i}\Vzero_c\Uzero_b(i_{g\dagger}\gbige)$
is 
$\bigcup F_m\bigl(
 \lefttop{i}\Vzero_c\Uzero_b(i_{g\dagger}\gbige)\bigr)$.
Hence, we have only to show 
$F_m\bigl(\lefttop{i}
 \Vzero_c\Uzero_b(i_{g\dagger}\gbige)
 \bigr)\otimes\nbigl(\gminia)\subset  \Image(\Phi)$
for any $m$,
which can be shown 
using an inductive argument
as in the proof of Lemma \ref{lem;08.1.12.16}.
Thus, we obtain Lemma \ref{lem;08.1.12.45}.
\hfill\qed

\begin{lem}
\label{lem;08.9.26.40}
If $i\in \supp(\vecm)$,
$\lefttop{i}\Vzero_a\Uzero_b(i_{g\dagger}\gbige
 \otimes\nbigl(\gminia))
=\Uzero_b(i_{g\dagger}\gbige
 \otimes\nbigl(\gminia))$.
\end{lem}
\pf
If $a<0$, we have
$z_i\cdot \lefttop{i}\Vzero_a
 \Uzero_b(i_{g\dagger}\gbige)
=\lefttop{i}\Vzero_{a-1}
 \Uzero_b(i_{g\dagger}\gbige)$.
Hence, they are the same
after tensoring $\nbigl(\gminia)$.
\hfill\qed

\begin{lem}
\label{lem;08.1.31.103}
In the case ($c>0$ and $i\leq l$),
we have the following description:
\begin{equation}
 \label{eq;08.1.31.100}
\lefttop{i}\Vzero_c
 \Uzero_b\bigl(i_{g\dagger}\gbige
 \otimes\nbigl(\gminia)\bigr)=
 \!\!
 \sum_{(c',m)\in\nbigu}
 \deldel_i^m\bigl(
 \lefttop{i}\Vzero_{c'}
 \Uzero_b\bigl(i_{g\dagger}\gbige
 \otimes\nbigl(\gminia)\bigr)
 \bigr),
\end{equation}
\[
  \nbigu=\bigl\{
 (c',m)\in\real_{\leq 0}\times\seisuu_{\geq 0}
 \,\big|\,
 c'+m\leq c
 \bigr\}
\]
In the case ($c\geq 0$ and $l<i\leq n$),
we have the following description:
\begin{equation}
 \label{eq;08.1.31.101}
\lefttop{i}\Vzero_c
 \Uzero_b\bigl(i_{g\dagger}\gbige
 \otimes\nbigl(\gminia)\bigr)=
 \!\!
 \sum_{(c',m)\in\nbigu}
 \deldel_i^m\bigl(
 \lefttop{i}\Vzero_{c'}
 \Uzero_b\bigl(i_{g\dagger}\gbige
 \otimes\nbigl(\gminia)\bigr)
 \bigr),
\end{equation}
\[
 \nbigu=\bigl\{
 (c',m)\in\real_{<0}\times\seisuu_{\geq 0}
 \,\big|\,
 c'+m\leq c
 \bigr\}
\]
\end{lem}
\pf
By Lemma \ref{lem;08.9.26.40},
we have only to consider the case
$i\not\in\supp(\vecm)$.
In the case ($c>0$ and $i\leq l$),
or ($c\geq 0$ and $l<i$),
we obtain the following equality
from Lemma \ref{lem;08.1.31.102}:
\[
 \lefttop{i}
 \Vzero_c\Uzero_b(i_{g\dagger}\gbige)
\otimes
 \nbigl(\gminia)
=\deldel_i\Bigl(
 \lefttop{i}\Vzero_{c-1}\Uzero_b
 (i_{g\dagger}\gbige)\otimes\nbigl(\gminia)
 \Bigr) 
+\lefttop{i}\Vzero_{<c}\Uzero_b
 (i_{g\dagger}\gbige)\otimes\nbigl(\gminia)
\]
Then, we obtain (\ref{eq;08.1.31.100})
and (\ref{eq;08.1.31.101})
by an easy induction.
\hfill\qed

\vspace{.1in}

Let us return to the proof of 
Lemma \ref{lem;08.9.26.50}.
According to Lemma \ref{lem;08.1.12.45}
and Lemma \ref{lem;08.1.31.103},
the induced subsheaves
$\lefttop{i}\Vzero_c$ 
of $\Gr^{\Uzero}_b(i_{g\dagger}\gbige
 \otimes\nbigl(\gminia))$
are finitely generated over 
$\lefttop{i}V_0\nbigr_X$.
It is easy to check 
the pseudo-coherence of $\lefttop{i}\Vzero$
as $\nbigo_{\nbigx}$-modules.
Hence, they are 
$\lefttop{i}V_0\nbigr_X$-coherent
by Proposition \ref{prop;08.9.26.10}.
Thus, the proof of Lemma \ref{lem;08.9.26.50}
is finished.
\hfill\qed

\vspace{.1in}
As a consequence,
the induced subsheaves
$\lefttop{i}V_c$ of 
\[
\psitilde_{t,u}(i_{g\dagger}\gbige
 \otimes\nbigl(\gminia))
\quad
\mbox{\rm and }
\quad
\Gr^{W(N)}\psitilde_{t,u}(i_{g\dagger}\gbige
 \otimes\nbigl(\gminia)) 
\]
are also $\lefttop{i}V_0\nbigr_X$-coherent modules.
Then, the following claim immediately
follows from Proposition \ref{prop;08.1.12.50}
and Proposition \ref{prop;08.1.12.20}.
\begin{prop}
\mbox{{}}\label{prop;08.1.12.61}
\begin{itemize}
\item
 The filtrations $\lefttop{i}\Vzero$ $(i=1,\ldots,n)$
 are compatible with 
 the primitive decomposition of
 $\Gr^{W(N)}_{p}
 \psitilde_{t,u}(i_{g\dagger}\gbige
 \otimes\nbigl(\gminia))$.
\item
 For each $i$,
 the primitive part 
 $P\Gr^{W(N)}_p
 \psitilde_{t,u}(i_{g\dagger}\gbige
 \otimes\nbigl(\gminia))$
 is strictly $S$-decomposable along $z_i=0$
 at any $\lambda_0$
 with the induced filtration  $\lefttop{i}\Vzero$.
\item
 In particular,
 $\nbigm_{p,u,J}
 \otimes\nbigl(\gminia)$ are strictly $S$-decomposable
along $z_i$ at any $\lambda_0$
with the induced filtration $\lefttop{i}\Vzero$.
\hfill\qed
\end{itemize}
\end{prop}

\section{Specialization of
 the associated $\nbigr_X$-module}
\label{subsection;08.2.2.1}
\subsection{Completion of 
 $\gbige$ along $\nbigd_{\nbar}$}

We use the setting in Section \ref{subsection;08.1.11.2}.
Let $\nbigdhat_{\nbar}$ denote the completion of $\nbigx$
along $\nbigd_{\nbar}$.
Let $\iota:\nbigdhat_{\nbar}\lrarr\nbigx$ 
denote the canonical morphism.
Recall that we have the irregular decomposition:
\begin{equation}
 \label{eq;07.9.7.11}
\iota^{\ast} \nbigqzero\nbige
=\bigoplus_{\gminia\in\Irr(\theta)}
 \nbigqzero\nbigehat_{\gminia},
\quad
 \iota^{\ast}\nbigqzero_{<\vecdelta}\nbige
=\bigoplus_{\gminia\in\Irr(\theta)}
 \nbigqzero_{<\vecdelta}\nbigehat_{\gminia}. 
\end{equation}
It is easy to see that
$\iota^{\ast}\gbige^{(\lambda_0)}$
is generated by 
$\iota^{\ast}\nbigqzero_{<\vecdelta}\nbige$
in $\iota^{\ast}\nbigqzero\nbige$ over
$\nbigr_{\Dhat_{\nbar}}$.
Hence, we have the decomposition
$\iota^{\ast}\gbige^{(\lambda_0)}
=\bigoplus_{\gminia}
 \gbigehat^{(\lambda_0)}_{\gminia}$
corresponding to {\rm(\ref{eq;07.9.7.11})},
where the $\nbigr_{\Dhat_{\nbar}}$-submodule
$\gbigehat_{\gminia}\subset
 \nbigqzero\nbigehat_{\gminia}$
is generated by
$\nbigqzero_{<\vecdelta}\nbigehat_{\gminia}$
on $U(\lambda_0)\times \Dhat_n$.
By varying $\lambda_0$ and gluing them,
we obtain
\begin{equation}
\label{eq;08.1.21.20}
 \iota^{\ast}\gbige
=\bigoplus_{\gminia}\gbigehat_{\gminia}.
\end{equation}

For each $\gminia\in\Irr(\theta)$,
we have the wild harmonic bundle
$(E_{\gminia},\delbar_{\gminia},
 \theta_{\gminia},h_{\gminia})$
such that $\Irr(\theta_{\gminia})=\{\gminia\}$,
which is obtained from $\harmonicbundle$ 
as the full reduction (Theorem \ref{thm;07.10.11.120}).

\begin{lem}
\label{lem;08.1.21.21}
We have the natural isomorphism
$\gbigehat_{\gminia}
\simeq
 \iota^{\ast}\gbige_{\gminia}$.
\end{lem}
\pf
By construction,
we have the natural isomorphisms
$ \nbigq\nbigehat_{\gminia}\simeq
 \iota^{\ast}\nbigq\nbige_{\gminia}$.
The KMS-structure at $\lambda_0$
is unique, if it exists,
which can be shown
by using an argument similar to 
that in the proof of
Lemma \ref{lem;07.11.23.5}. 
Hence, we have the equality
of the filtered bundles
$ \nbigqzero_{\ast}\nbigehat_{\gminia}=
 \iota^{\ast}\nbigqzero_{\ast}\nbige_{\gminia}$
for each $\lambda_0$
under the above isomorphism.
Then the claim of the lemma is clear
by construction in 
Section \ref{subsection;08.9.28.1}.
\hfill\qed

\subsection{Holonomicity}

For any subset 
$I\subset \nbar:=\{1,\ldots,n\}$,
we put $D_I:=\bigcap_{i\in I}D_i$,
and let $N_{D_I}^{\ast}X$
denote the conormal bundle of $D_I$ in $X$.
\index{conormal bundle $N_{D_I}^{\ast}X$}
Let $\nbigs$ denote the union
$\cnum_{\lambda}\times
 \bigcup_{I\subset\nbar}N_{D_I}^{\ast}X$.

\begin{prop}
The characteristic variety of $\gbige$
is contained in $\nbigs$.
In particular, $\gbige$ is holonomic.
\end{prop}
\pf
For $J=(j_1,\ldots,j_n)\in\seisuu_{\geq 0}^n$,
we put $\deldel^J=\prod_{i=1}^n\deldel_i^{j_i}$
and $|J|=\sum j_i$.
Let $F_m\nbigr_X=\bigl\{
 \sum_{|J|\leq m} a_J\cdot\deldel^J
 \bigr\}$.
We put $F_m(\gbige):=
 F_m\nbigr_X\cdot 
 \nbigqzero_{<\vecdelta}(\nbige)$
around $\{\lambda_0\}\times X$.
Thus, we obtain the coherent filtration of $\gbige$
around $\{\lambda_0\}\times X$.

Let $\Gr^F(\gbige)$ denote the associated
$\nbigo_{\cnum_{\lambda}
 \times T^{\ast}X}$-module.
Let us show that 
the support of $\Gr^F(\gbige)$ is 
contained in $\nbigs$,
when we shrink $X$.
We have the coherent filtrations of
$\bigl\{F_m(\gbige_{\gminia})\,\big|\,
 m=1,2,\ldots\bigr\}$ 
of $\gbige_{\gminia}$ constructed
in the same way.
Then, we have
$\iota^{\ast}F_m(\gbige)
\simeq
 \bigoplus \iota^{\ast}F_m(\gbige_{\gminia})$
under the isomorphism (\ref{eq;08.1.21.20})
and Lemma \ref{lem;08.1.21.21}.
Let $\nbigy$ denote 
$\cnum_{\lambda}
 \times (T^{\ast}X\times_X\{O\})$.
Let $Ch(\gbige)$ and $Ch(\gbige_{\gminia})$
denote the characteristic varieties 
of $\gbige$ and $\gbige_{\gminia}$.
Because
$\iota^{\ast}\Gr^F(\gbige)
\simeq
 \bigoplus \iota^{\ast}\Gr^{F}(\gbige_{\gminia})$,
we obtain that the completion of
$Ch(\gbige)$ along $\nbigy$
is the union of the completions of
$Ch(\gbige_{\gminia})$ along $\nbigy$.
Hence, the above claim follows from 
Corollary \ref{cor;08.1.21.20}.
\hfill\qed

\subsection{Strict $S$-decomposability
 along any monomial function}
\label{subsection;07.9.7.35}

Let us consider a monomial function
$g=\vecz^{\vecp}$
for some $\vecp\in\seisuu_{\geq 0}^n$.
\begin{prop}
\label{prop;07.10.28.50}
$\gbige$ is strictly $S$-decomposable along $g$.
We have natural isomorphisms:
\begin{equation}
 \label{eq;07.9.7.24}
 \iota^{\ast}\psitilde_{g,u}(\gbige)
\simeq
 \bigoplus_{\gminia\in\Irr(\theta)}
 \iota^{\ast}\psitilde_{g,u}(\gbige_{\gminia})
\end{equation}
It is compatible with the nilpotent part $N$
of $-\deldel_tt$.
\end{prop}
\pf
We put $\supp(\vecp):=\{i\,|\,p_i\neq 0\}$.
For any $b<0$,
we define
\begin{equation}
 \label{eq;07.9.7.15}
 \Uzero_{b}\bigl(i_{g\dagger}\gbige\bigr):=
V_0\nbigr_{X\times\cnum_t}\cdot
 \bigl(
 \lefttop{\nbar}
 \Vzero_{b\vecp-\epsilon\vecdelta_{\nbar-\supp(\vecp)}}
 (\nbigq\nbige)
\otimes 1
 \bigr).
\end{equation}
For any $b\geq 0$,
we define
\begin{equation}
 \label{eq;07.9.7.16}
 \Uzero_b(i_{g\dagger}\gbige):=
 \sum_{\substack{c<0,m\in\seisuu_{\geq \,0}
 \\ c+m\leq b}}
 \deldel_t^m\Uzero_c(i_{g\dagger}\gbige).
\end{equation}
They are $V_0\nbigr_{X\times\cnum_t}$-coherent modules.
By construction,
they satisfy the following conditions:
\begin{itemize}
\item $\bigcup_{a\in\real}
 \Uzero_a(i_{g\dagger}\gbige)=i_{g\dagger}\gbige$.
For any $b\in\real$,
there exists $\epsilon>0$ such that
$\Uzero_{b}(i_{g\dagger}\gbige)
 =\Uzero_{b+\epsilon}(i_{g\dagger}\gbige)$.
\item $t\cdot \Uzero_a(i_{g\dagger}\gbige)\subset
 \Uzero_{a-1}(i_{g\dagger}\gbige)$ for any $a\in\real$,
and $t\cdot \Uzero_a(i_{g\dagger}\gbige)=
 \Uzero_{a-1}(i_{g\dagger}\gbige)$ for any $a<0$.
\item
 $\deldel_t:\Uzero_a(i_{g\dagger}\gbige)
 \subset\Uzero_{a+1}(i_{g\dagger}\gbige)$ 
 for any $a\in\real$,
and the induced morphisms
 $\deldel_t:
 \Gr^{\Uzero}_a(i_{g\dagger}\gbige)\lrarr
 \Gr^{\Uzero}_{a+1}(i_{g\dagger}\gbige)$ 
 are surjective for any $a>-1$.
\end{itemize}

The induced morphism
$\nbigd_{\nbar}\times\cnum_t\lrarr
 X\times\cnum_t$ is also denoted by $\iota$.
Let us look at
$\iota^{\ast}\bigl(i_{g\dagger}\gbige\bigr)$
and the induced filtration
$\iota^{\ast}\Uzero\bigl(i_{g\dagger}\gbige\bigr)$.
We have the decomposition:
\begin{equation}
 \label{eq;07.9.7.17}
 \iota^{\ast}\bigl(i_{g\dagger}\gbige\bigr)
=\bigoplus
 i_{g\dagger}\gbigehat_{\gminia},
\quad
  \lefttop{\nbar}
 \Vzero_{b\vecp-\epsilon\vecdelta_{\nbar-\supp(\vecp)}}
 (\nbigq\nbige)
=\bigoplus_{\gminia}
 \lefttop{\nbar}
 \Vzero_{b\vecp-\epsilon\vecdelta_{\nbar-\supp(\vecp)}}
 (\nbigq\nbigehat_{\gminia})
\end{equation}
We have 
$\Uzero_b(i_{g\dagger}\gbigehat_{\gminia})$
defined as in (\ref{eq;07.9.7.15}) and (\ref{eq;07.9.7.16})
for $\gbigehat_{\gminia}$ instead of $\gbige$.
Then, we have the decomposition for any $b$
corresponding to (\ref{eq;07.9.7.17}):
\[
 \iota^{\ast}\Uzero_{b}(i_{g\dagger}\gbige)
=\bigoplus_{\gminia\in \Irr(\theta)}
 \Uzero_b(i_{g\dagger}\gbigehat_{\gminia})
\]

We have natural isomorphisms:
\[
 \gbigehat_{\gminia}
\simeq
 \iota^{\ast}\gbige_{\gminia},
\quad\quad
 \lefttop{\nbar}
 \Vzero_{b\vecp-\epsilon\vecdelta_{\nbar-\supp(\vecp)}}
 (\nbigq\nbigehat_{\gminia})
\simeq
 \iota^{\ast}
  \lefttop{\nbar}
 \Vzero_{b\vecp-\epsilon\vecdelta_{\nbar-\supp(\vecp)}}
 (\nbigq\nbige_{\gminia})
\]
Note the surjectivity of $\Phi_b$
in Lemma \ref{lem;08.1.12.16}.
Hence, we have the natural isomorphism:
\[
 \Uzero_b(i_{g\dagger}\gbigehat_{\gminia})
\simeq\iota^{\ast}\Uzero_{b}(i_{g\dagger}\gbige_{\gminia})
\]
Therefore,
we have the following equality of the filtrations
under the isomorphism
$\iota^{\ast}i_{g\dagger}\gbige
\simeq\bigoplus_{\gminia\in\Irr(\theta)}
 \iota^{\ast}i_{g\dagger}\gbige_{\gminia}$:
\begin{equation}
 \label{eq;07.9.7.22}
 \iota^{\ast}\Uzero_b\bigl(
 i_{g\dagger}\gbige \bigr)
=\bigoplus_{\gminia\in\Irr(\theta)}
 \iota^{\ast}\Uzero_b\bigl(
 i_{g\dagger}\gbige_{\gminia}\bigr)
\end{equation}

Recall that
$\gbige_{\gminia}$ is strictly $S$-decomposable along $g$
at any $\lambda_0$ with the $V$-filtration
$\Uzero(i_{g\dagger}\gbige_{\gminia})$
(Proposition \ref{prop;08.1.12.20}).
Hence, we can conclude that
$\gbige$ is strictly $S$-decomposable along $g$
at $\lambda_0$ with the $V$-filtration
$\Uzero(\gbige)$,
due to Proposition \ref{prop;08.1.23.16} below.
We obtain the isomorphism (\ref{eq;07.9.7.24})
from (\ref{eq;07.9.7.22}).
\hfill\qed

\vspace{.1in}
We have the tame harmonic bundles
$(E'_{\gminia},\delbar'_{\gminia},
 \theta'_{\gminia},h'_{\gminia})$
such that 
$(E_{\gminia},\delbar_{\gminia},
 \theta_{\gminia},h_{\gminia})
\simeq
 (E'_{\gminia},\delbar'_{\gminia},
 \theta'_{\gminia},h'_{\gminia})
 \otimes L(\gminia)$.
The associated $\nbigr$-modules
are denoted by $\gbige_{\gminia}'$.

\begin{cor}
\label{cor;07.9.7.36}
We have isomorphisms:
\begin{equation}
 \label{eq;08.1.12.62}
\iota^{\ast}
 \Gr^{W(N)}_p\psitilde_{g,u}(\gbige)
\simeq
 \bigoplus_{\gminia}
 \iota^{\ast}
 \Gr^{W(N)}_p\psitilde_{g,u}(\gbige'_{\gminia})
 \otimes
 \nbigl(\gminia)
\end{equation}
\begin{equation}
 \label{eq;07.11.25.25}
\iota^{\ast}
 P\Gr^{W(N)}_p
  \psitilde_{g,u}(\gbige)
\simeq
 \bigoplus_{\gminia}
\iota^{\ast}
 P\Gr^{W(N)}_p
 \psitilde_{g,u}(\gbige'_{\gminia})
\otimes\nbigl(\gminia)
\end{equation}
\end{cor}
\pf
It follows from Proposition \ref{prop;07.10.28.50}.
\hfill\qed

\vspace{.1in}
Recall the notion of
$\nbiga$-good wild harmonic bundle
in Definition \ref{df;07.11.25.30}.

\begin{cor}
\label{cor;07.11.2.40}
If $\harmonicbundle$ is 
an unramifiedly $\nbiga$-good wild
harmonic bundle,
then we have the vanishing
$P\Gr^{W(N)}_p
 \psitilde_{g,u}(\gbige)=0$
unless $u\in\real\times\nbiga$.
\end{cor}
\pf
Under the assumption,
$(E'_{\gminia},\delbar'_{\gminia},\theta'_{\gminia},
 h'_{\gminia})$ are $\nbiga$-tame harmonic bundles.
Hence, we obtain the vanishing of
$\iota^{\ast}P\Gr^{W(N)}_p
 \psizero_u(i_{g\dagger}\gbige'_{\gminia})=0$
unless $u\in\nbiga$
by the argument given in Section 19.5 of \cite{mochi2}.
Although we considered only 
the case $\nbiga=\sqrt{-1}\real$
in \cite{mochi2},
the other cases can be shown 
with the same argument.
\hfill\qed

\subsection{Strict $S$-decomposability of
 $P\Gr^{W(N)}_p
  \psitilde_{g,u}(\gbige)$ 
 along a coordinate function}
\label{subsection;08.9.28.4}

Let $b<0$.
We use the notation in 
Section \ref{subsection;08.1.12.51}.
For $1\leq i\leq n$,
we put $K:=\nbar-(\supp(\vecp)\cup\{i\})$.
In the case ($c\leq 0$ and $i\in\supp(\vecp)$)
or ($c<0$ and $i\not\in\supp(\vecp)$),
we define
\[
\lefttop{i}\Vzero_c
  \Uzero_b\bigl(
 i_{g\dagger}\gbige
 \bigr):=
\Image\Bigl(
 \lefttop{i,t}V_0\nbigr_{X\times\cnum_t}
\otimes
  \bigl(
 \lefttop{\nbar}
 \Vzero_{b\vecp-c\vecdelta_i-\epsilon\vecdelta_K}
 (\nbigq\nbige)\otimes 1
 \bigr)
\lrarr
\Uzero_b\bigl(
 i_{g\dagger}\gbige
 \bigr)
\Bigr).
\]
(See Section \ref{subsection;08.9.28.20}
for $\lefttop{i,t}V_0\nbigr_{X\times\cnum_t}$.)
In the case ($c>0$ and $i\in\supp(\vecp)$),
we put 
\[
\lefttop{i}\Vzero_c
  \Uzero_b\bigl(
 i_{g\dagger}\gbige
 \bigr):=
\sum_{(c',m)\in\nbigu}
 \deldel_i^m\cdot
 \lefttop{i}\Vzero_{c'}
  \Uzero_b\bigl(
 i_{g\dagger}\gbige
 \bigr),
\]
\[
 \nbigu:=\bigl\{
(c',m)\in\real_{\leq 0}\times\seisuu_{\geq 0}
 \,\big|\,
 c'+m\leq c
\bigr\}.
\]
In the case ($c\geq 0$ and $i\not\in\supp(\vecp)$),
we set
\[
\lefttop{i}\Vzero_c
  \Uzero_b\bigl(
 i_{g\dagger}\gbige
 \bigr):=
\sum_{(c',m)\in\nbigu}
 \deldel_i^m\cdot
 \lefttop{i}\Vzero_{c'}
  \Uzero_b\bigl(
 i_{g\dagger}\gbige
 \bigr),
\]
\[
 \nbigu:=\bigl\{
(c',m)\in\real_{<0}\times\seisuu_{\geq 0}
 \,\big|\,
 c'+m\leq c
 \bigr\}.
\]
Thus, we obtain a filtration
$\lefttop{i}\Vzero$ of
$\Uzero_b\bigl(i_{g\dagger}\gbige\bigr)$ for each $b<0$.
(Note that we have only to consider the case $b<0$
 to investigate the property of
 $\psitilde_{t,u}(i_{g\dagger}\gbige)$.)
The induced filtrations of
$\Gr^{\Uzero}_b\bigl(i_{g\dagger}\gbige \bigr)$,
$\psitilde_{t,u}(i_{g\dagger}\gbige)$ and
$\Gr^{W(N)}\psitilde_{t,u}(i_{g\dagger}\gbige)$
are also denoted by $\lefttop{i}\Vzero$.
They are $\lefttop{i}V_0\nbigr_X$-coherent.

\begin{prop}
\mbox{{}} \label{prop;08.1.12.60}
\begin{itemize}
\item
 The filtrations $\lefttop{i}\Vzero$ $(i=1,\ldots,n)$
 are compatible with 
 the primitive decomposition of
 $\Gr^{W(N)}_{p}
 \psitilde_{t,u}(i_{g\dagger}\gbige)$.
\item
 For each $i$,
 the primitive part 
 $P\Gr^{W(N)}_p
 \psitilde_{t,u}(i_{g\dagger}\gbige)$
 is strictly $S$-decomposable along $z_i=0$
 at any $\lambda_0$
 with the filtration  $\lefttop{i}\Vzero$.
\end{itemize}
\end{prop}
\pf
We also have the induced filtration
$\lefttop{i}\Vzero$ of
$\Gr^{W(N)}_p\psitilde_{g,u}(\gbige_{\gminia})$.
Under the isomorphism (\ref{eq;08.1.12.62}),
we have the equality of the filtrations
$\lefttop{i}\Vzero$,
by construction and Lemma \ref{lem;08.1.12.45}.
Then, the claim of Proposition \ref{prop;08.1.12.60}
follows from Proposition \ref{prop;08.1.12.61}
and Proposition \ref{prop;08.1.23.16} below.
\hfill\qed

\subsection{Comparison of the specializations}
\label{subsection;08.9.28.5}

According to Proposition \ref{prop;08.1.12.60},
we have the decomposition:
\begin{equation}
 \label{eq;08.9.26.100}
 P\Gr^{W(N)}_p
 \psitilde_{g,u}(\gbige)
=\bigoplus_I\nbigm_{u,p,I}
\end{equation}
\index{sheaf $\nbigm_{u,p,I}$}
Here, 
the support of $\nbigm_{u,p,I}$ is $\nbigd_I$,
and the support of 
any non-trivial coherent
$\nbigr$-submodule of $\nbigm_{u,p,I}$
is not contained in $\nbigd_J$ ($J\supset I$).
Since $\nbigm_{u,p,I}$ are strictly $S$-decomposable
along any $z_i$,
it is the push-forward of
the $\nbigr_{D_I}$-module $\nbigm'_{u,p,I}$.
Let $I^c:=\nbar-I$ and $z_{I^c}:=\prod_{j\in I^c}z_j$.
The localization
$\nbigm'_{u,p,I}\otimes\nbigo(\ast z_{I^c})$
is a family of meromorphic $\lambda$-flat bundles.

We take an auxiliary sequence
$\vecm(0),\ldots,\vecm(L)$  for $\Irr(\theta)$.
By shrinking $X$,
we obtain harmonic bundles
$\bigl(
 E^{\vecm(i)}_{\gminia},
 \delbar^{\vecm(i)}_{\gminia},
 \theta^{\vecm(i)}_{\gminia},
 h^{\vecm(i)}_{\gminia}
 \bigr)$
on $X-D$
for $\gminia\in\Irrbar(\theta,\vecm(i))$
as the reductions in the level $\vecm(i)$,
which induce
$\nbigr_X$-modules
$\gbige^{\vecm(i)}_{\gminia}$.
\index{sheaf $\gbige^{\vecm(i)}_{\gminia}$}
We have the decomposition:
\[
 P\Gr^{W(N)}_p\psitilde_{g,u}
 \bigl(\gbige^{\vecm(i)}_{\gminia}\bigr)
=\bigoplus_I\nbigm_{u,p,I}(E^{\vecm(i)}_{\gminia})
\]

\begin{lem}
\label{lem;08.2.2.10}
Let $k$ be determined by
$\vecm(0)\in\seisuu_{<0}^k\times\veczero_{n-k}$.
If $D_I\subset D(\leq k)$,
we have the natural isomorphism:
\begin{equation}
 \Psi:
 \nbigm_{u,p,I}(E^{\vecm(0)}_{0})
\simeq
 \nbigm_{u,p,I}
\end{equation}
\end{lem}
\pf
We may assume that $D_I\subset g^{-1}(0)$.
Let $\iota_I:\Dhat_I\lrarr X$ be the natural morphism.
We have the following natural isomorphisms:
\begin{multline}
  \iota_I^{\ast}P\Gr^{W(N)}_p
 \psitilde_{g,u}(\gbige)
\simeq
 P\Gr^{W(N)}_p
 \psitilde_{g,u}(\iota_I^{\ast}\gbige)
\simeq\bigoplus_{\gminia\in\Irr(\theta,\vecm(0))}
 P\Gr^{W(N)}_p
  \psitilde_{g,u}(\iota_I^{\ast}\gbige^{\vecm(0)}_{\gminia}) \\
\simeq
 P\Gr^{W(N)}_p\psitilde_{g,u}
 (\iota_I^{\ast}\gbige^{\vecm(0)}_{0})
\simeq
\iota_I^{\ast}
 P\Gr^{W(N)}_p\psitilde_{g,u}(\gbige^{\vecm(0)}_{0})
\end{multline}
Since the support of $\nbigm_{u,p,I}$ 
and $\nbigm_{u,p,I}(E^{\vecm(0)}_{0})$ are contained in
$\nbigd_I$,
we obtain the desired isomorphism.
\hfill\qed

\begin{cor}
\label{cor;08.2.2.11}
We have natural isomorphisms
$\Psi:\nbigm_{u,p,\nbar}\simeq
 \nbigm_{u,p,\nbar}(E^{\vecm(i)}_{0})$
for $i=0,1,\ldots,L$.
\hfill\qed
\end{cor}

\subsection{$C^{\infty}$-lift of a section}
\label{subsection;08.2.2.5}
\index{$C^{\infty}$-lift}

This subsection is a preparation for
Propositions \ref{prop;07.9.8.5} and 
\ref{prop;07.9.8.26}.
Let $\lambda_0$ be generic with respect to 
\[
 \bigcup_{i=1}^n
\bigl\{v\in\real\times\cnum\,\big|\,
 m_i\cdot v\in \KMS(\nbigp\nbige^0,i)\bigr\}.
\]
Let $\nbigk$ be a small neighbourhood of $\lambda_0$
in $\cnum_{\lambda}$.

We take an auxiliary sequence
$\vecm(0),\ldots,\vecm(L)$ for $\Irr(\theta)$.
Let $f$ be a section of
$\nbigm_{u,p,I}(E^{\vecm(i)}_{\gminia})$ 
on $\nbigk\times X$.
We would like to consider
a $C^{\infty}$-lift of $f$
to $\nbigm_{u,p,I}$.
Let $k$ be determined by
$\vecm(0)\in\seisuu_{<0}^k\times\veczero_{n-k}$.
For simplicity,
we consider the case in which
$\vecm(i)\in \seisuu_{<0}^k\times\veczero_{n-k}$.
We may assume $b:=\paramap(\lambda_0,u)<0$.
We will shrink $X$ and $\nbigk$
without mention.

We have a section 
$\ftilde\in 
 \Uzero_{b}(i_{g\dagger}
 \gbige^{\vecm(i)}_{\gminia})$
on $\nbigk\times X$
with the following property:
\begin{itemize}
\item
 The induced section
 $\ftilde^{(1)}$ of 
 $\Gr^{\Uzero}_b
 \bigl(i_{g\dagger}
 \gbige^{\vecm(i)}_{\gminia}\bigr)$
 is contained in 
 $W_p(N)\psitilde_{g,u}
 (\gbige^{\vecm(i)}_{\gminia})$.
 And, the induced section 
 $\ftilde^{(2)}$ of
 $\Gr^{W(N)}_p\psitilde_{g,u}\bigl(
 \gbige^{\vecm(i)}_{\gminia} \bigr)$
 is equal to $f$
 which is contained in 
 $P\Gr^{W(N)}_p\psitilde_{g,u}\bigl(
 \gbige^{\vecm(i)}_{\gminia} \bigr)$.
\end{itemize}
Recall $g=\vecz^{\vecp}$
for some $\vecp\in\seisuu_{\geq 0}^n$.
We set $K:=\bigl\{
 1\leq i\leq n\,\big|\,p_i=0
 \bigr\}$.
By the definition of the filtration $\Uzero$,
we have the expression of $f$
as a finite sum
\[
 f=\sum P_B\cdot \deldel^B
 (a_B\otimes 1),
\]
where $\deldel^A=\prod_{i=1}^n\deldel_i^{b_i}$
for $B=(b_1,\ldots,b_n)$,
$P_B\in \cnum[t\deldel_t]$,
and $a_B\in 
 \lefttop{\nbar}\Vzero_{b\vecp-\epsilon\vecdelta_K}
 (\nbigq\nbige^{\vecm(i)}_{\gminia})$
for some $\epsilon>0$.

Let $S$ be a small multi-sector in 
$\nbigk\times(X-D(\leq k))$.
Since $\lambda_0$ is assumed to be generic,
we have a $\DD$-flat splitting
of the Stokes filtration in the level $\vecm(i)$:
\[
 \lefttop{\nbar}\Vzero_{b\vecp-\epsilon\vecdelta_K}
 \bigl(\nbigq\nbige\bigr)_{|\Sbar}
=
 \bigoplus_{\gminib\in\Irr(\theta,\vecm(i))}
 \lefttop{\nbar}\Vzero_{b\vecp-\epsilon\vecdelta_K}
 \bigl(\nbigq\nbige\bigr)_{\gminib,S},
\quad\quad
\nbigq\nbige_{|\Sbar}
=\bigoplus_{\gminib\in\Irr(\theta,\vecm(i))}
 \nbigq\nbige_{\gminib,S}
\]
They induce a decomposition of the $\nbigr$-module
$ \gbige_{|S}
=\bigoplus \gbige^{\vecm(i)}_{\gminib,S}$,
and
each $\gbige^{\vecm(i)}_{\gminib,S}$
is naturally isomorphic to
$(\gbige^{\vecm(i)}_{\gminib})_{|S}$.
Let $a_{B,S}\in 
 \lefttop{\nbar}\Vzero_{b\vecp-\epsilon\vecdelta_K}
 \bigl(\nbigq\nbige\bigr)_{\gminia,S}$
be the lift of $a_B$.
We have the lift $\ftilde_S$ of the restriction $\ftilde_{|S}$
to $i_{g\dagger}\gbige^{\vecm(i)}_{\gminia,S}$
on $S\times\cnum_t$:
\[
 \ftilde_S=\sum P_B\cdot \deldel^B\cdot
 \bigl(a_{B,S}\otimes 1\bigr)
\]
If we are given another $\DD$-flat splitting
\[
 \lefttop{\nbar}\Vzero_{b\vecp-\epsilon\vecdelta_K}
 (\nbigq\nbige)_{|\Sbar}
=\bigoplus_{\gminia\in\Irr(\theta,\vecm(i))}
 \lefttop{\nbar}\Vzero_{b\vecp-\epsilon\vecdelta_K}
 (\nbigq\nbige)_{\gminib,S}',
\]
we obtain another lift:
\[
 \ftilde_S'=\sum P_B\cdot \deldel^B\cdot
 \bigl(a_{B,S}'\otimes 1\bigr)
\]
The following lemma is obvious by construction,
which is stated for the reference
in the subsequent argument.
\begin{lem}
\label{lem;08.2.2.2}
We take a frame $\vecw$ of
$\lefttop{\nbar}\Vzero_{b\vecp-\epsilon\vecdelta_K}
 \bigl(\nbigq\nbige\bigr)$,
and let $h_1$ be the hermitian metric 
for which $\vecw$ is orthonormal.
Then, we have 
$a_{B,S}-a_{B,S}'=O\Bigl(
 \exp\bigl(-\eta|\vecz^{\vecm(i)}|\bigr)
 \Bigr)$
for some $\eta>0$ with respect to $h_1$.
\hfill\qed
\end{lem}

Shrinking $X$ and $\nbigk$ appropriately,
we take a finite covering
$\nbigk\times (X-D)
=\bigcup S_i$ by small multi-sectors.
We take a partition of unity
$\{\chi_{i}\}$ subordinated to
the covering $\{S_i\}$
such that 
$\chi_i$ depend only on
$\arg(z_j)$ $(j=1,\ldots,k)$.
For each $S_i$,
we take a lift $\ftilde_{S_i}$ as above.
Let ``$\otimes C^{\infty}$'' denote 
the operation to take the tensor products
with the sheaf of $C^{\infty}$-functions.
(Note that it is faithfully flat,
 according to \cite{malgrange2}.)
Then, we obtain the following section:
\[
 \ftilde_{C^{\infty}}
=\sum \chi_{i}\cdot \ftilde_{S_i}
\in \bigl(
 i_{g\dagger}\gbige\otimes C^{\infty}
\bigr)
 _{|\nbigk\times (X-D(\leq k))\times \cnum_t}
\]
\begin{lem}
\label{lem;08.2.2.3}
$\ftilde_{C^{\infty}}$ gives
a section of
$\Uzero_b\bigl(
 i_{g\dagger}\gbige
 \bigr)\otimes C^{\infty}$
on $\nbigk\times X\times\cnum_t$.
\end{lem}
\pf
We have the following:
\[
 \chi_{i}\cdot \ftilde_{S_i}
=\sum_B P_B\cdot \chi_{i}\cdot\deldel^B
 (a_{B,S_i}\otimes 1)
=\sum_B P_B\cdot 
 \sum_{L}
 \deldel^L(c_{B,L,S_i}\otimes 1)
\]
Here, $c_{B,L,S_i}$ are the product 
of $a_{B,S_i}$ and 
$R_{B,L}(\chi_{i})$,
where $R_{B,L}\in\nbigr_X$
are independent of $S_i$.
By using $\sum \chi_i=1$
and Lemma \ref{lem;08.2.2.2},
we obtain that
$\sum_i c_{A,L,S_i}$ is a $C^{\infty}$-section
of $\lefttop{\nbar}
 \Vzero_{b\vecp-\vecdelta_K}(\nbigq)$.
Then, the claim of the lemma follows.
\hfill\qed

\vspace{.1in}

Because of Lemma \ref{lem;08.2.2.3},
we obtain an induced section
$\ftilde^{(1)}_{C^{\infty}}$ of
$\Gr^{\Uzero}_b(i_{g\dagger}\gbige)
\otimes C^{\infty}$
on $\nbigk\times X$.

\begin{lem}
\mbox{{}}\label{lem;08.2.2.8}
\begin{itemize}
\item
 $\ftilde^{(1)}_{C^{\infty}}$
 is contained in 
 $W_p(N)\psitilde_{g,u}(\gbige)
   \otimes C^{\infty}$.
 Hence, it induces a section of
 $\ftilde^{(2)}_{C^{\infty}}$
of $\Gr^{W(N)}_p\psitilde_{g,u}(\gbige)
 \otimes C^{\infty}$.
\item
$\ftilde^{(2)}_{C^{\infty}}$
is contained in
$P\Gr^{W(N)}_p\psitilde_{g,u}(\gbige)
 \otimes C^{\infty}$.
Moreover,
it is contained in 
$\nbigm_{u,p,I}\otimes C^{\infty}$
for the decomposition
$P\Gr^{W(N)}_p\psitilde_{g,u}(\gbige)
 \otimes C^{\infty}
=\bigoplus\nbigm_{u,p,J}\otimes C^{\infty}$.
\item
$(\ftilde^{(2)}_{C^{\infty}})_{| 
\nbigk\times \Dhat(\leq k)}$
is equal to $f_{|\nbigk\times\Dhat(\leq k)}$
under the isomorphism
\[
 \bigl(
 \nbigm_{u,p,I}\otimes C^{\infty}
\bigr)_{|\nbigk\times \Dhat(\leq k)}
\simeq
 \bigoplus_{\gminib\in\Irr(\theta,\vecm(i))}
  \bigl(
 \nbigm_{u,p,I}(E^{\vecm(i)}_{\gminib})
 \otimes C^{\infty}
\bigr)_{|\nbigk\times \Dhat(\leq k)}.
\]
\item
In particular, if $D_I\subset D(\leq k)$,
we have $\ftilde^{(2)}_{C^{\infty}}
=\Psi(f)$.
(See Lemma {\rm\ref{lem;08.2.2.10}}
for $\Psi$.)
\end{itemize}
\end{lem}
\pf
Let us show the first claim.
Let $v\neq u$,
and let $\ftilde^{(1)}_{C^{\infty},v}$
denote the 
$\psitilde_{g,v}(\gbige)\otimes C^{\infty}$-component
of $\ftilde^{(1)}_{C^{\infty}}$
for the decomposition
$\Gr^{\Uzero}_b(i_{g\dagger}\gbige)
 \otimes C^{\infty}$.
Assume that it is not $0$.
Then, there is an integer $l$ such that
$\ftilde^{(1)}_{C^{\infty},v}$
is contained in 
$W_{l}(N)\psitilde_{g,v}(\gbige)
 \otimes C^{\infty}$
but not in 
$W_{l-1}(N)\psitilde_{g,v}(\gbige)
 \otimes C^{\infty}$.
Let $[\ftilde^{(1)}_{C^{\infty},v}]$
denote the induced section of
$\Gr^{W(N)}_l\psitilde_{g,v}(\gbige)
 \otimes C^{\infty}$.
According to Proposition \ref{prop;08.1.12.60},
we have the following decomposition
as in (\ref{eq;08.9.26.100}):
\[
 \Gr^{W(N)}_l\psitilde_{g,v}(\gbige)
=\bigoplus_{J} M_{J}
\]
Let $[\ftilde^{(1)}_{C^{\infty},v}]
=\sum [\ftilde^{(1)}_{C^{\infty},v}]_J$
denote the corresponding decomposition.
To show the first claim,
we have only to show 
$[\ftilde^{(1)}_{C^{\infty},v}]_J=0$
for any $J$.
If $D_J\subset D(\leq k)$,
we can calculate
$[\ftilde^{(1)}_{C^{\infty},v}]_J$
after taking the completion along $D(\leq k)$.
Hence, it is easy to check 
$[\ftilde^{(1)}_{C^{\infty},v}]_J=0$
in this case.

Let us consider the case $D_J\not\subset D(\leq k)$.
Let $\ftilde^{(1)}_{S}$
denote the section of
$\Gr^{\Uzero}_b(i_{g\dagger}\gbige_{|S})$
induced by $\ftilde_S$.
By our choice of $\ftilde_S$,
the restriction of $\ftilde^{(1)}_{S}$ to
$S\cap\bigl(\nbigk\times(X-D(\leq k))\bigr)$
is contained in
$W_{p-1}(N)\psitilde_{g,u}(\gbige)_{|S}$.
Hence, the restriction of
$\ftilde^{(1)}_{C^{\infty}}$
to $\nbigk\times (X-D(\leq k))$
is contained in 
$W_{p-1}(N)\psitilde_{g,u}(\gbige)
\otimes C^{\infty}$.
As a result,
the restriction of
$[\ftilde^{(1)}_{C^{\infty},v}]$ to
$\nbigk\times \bigl(X-D(\leq k)\bigr)$ is $0$.

We have the inclusion
$M_J
\lrarr
 \Mtilde_J:=
 M_J\otimes\nbigo(\ast D(\leq k))$.
It induces the inclusion
$M_J\otimes C^{\infty}
\lrarr
 \Mtilde\otimes C^{\infty}$.
Recall that 
$\Mtilde_J$ is the push-forward
of a family of meromorphic $\lambda$-flat bundles
on $\nbigd_J$ to $\nbigx$.
Since the restriction of
$[\ftilde^{(1)}_{C^{\infty},v}]_J
 \in \Mtilde_J\otimes C^{\infty}$ 
to $\nbigk\times (X-D(\leq k))$ is $0$,
we can conclude that
$[\ftilde^{(1)}_{C^{\infty},v}]_J$ is $0$.

The second and third claims can be shown 
by the same argument.
The fourth claim follows from
the third one.
\hfill\qed

\section{Construction of the $\nbigr$-triple
 $(\gbige,\gbige,\gbigc)$}
\label{subsection;08.9.28.11}
\subsection{Statement}

We continue to use the setting in
Section \ref{subsection;08.1.11.2}.
Let $\sigma:\cnum^{\ast}_{\lambda}\lrarr
 \cnum_{\lambda}^{\ast}$
be given by
$\sigma(\lambda)=-\lambdabar^{-1}$.
\index{map $\sigma$}
We set 
$\vecS:= \bigl\{\lambda\in\cnum\,\big|\, 
|\lambda|=1\bigr\}$.
\index{set $\vecS$}
Recall that we have the naturally induced 
Hermitian sesqui-linear pairing
\begin{equation}
 \label{eq;10.5.24.30}
 C:\nbige_{|\vecS\times (X-D)}\otimes
 \sigma^{\ast}\nbige_{|\vecS\times (X-D)}\lrarr
 C^{\infty}(\vecS\times (X-D))
\end{equation}
given by 
$C(f_1\otimes\sigma^{\ast}f_2):=
 h(f_1,\sigma^{\ast}f_2)$.
We show the following proposition
in Sections 
\ref{subsection;08.9.27.2}--\ref{subsection;07.10.19.20}.
\begin{prop}
\label{prop;08.9.27.3}
We have a unique sesqui-linear pairing
\[
 \gbigc:\gbige_{|\vecS\times X}\otimes
\sigma^{\ast}\gbige_{|\vecS\times X}
\lrarr \distribution_{\vecS\times X/\vecS},
\]
whose restriction to $\vecS\times (X-D)$
is equal to $C$.
\index{pairing $\gbigc$}

As a result, we obtain an $\nbigr_X$-triple
$\gbigt(E):=(\gbige,\gbige,\gbigc)$
associated to the unramifiedly 
good wild harmonic bundle
$\harmonicbundle$.
\index{$\nbigr$-triple $\gbigt(E)$}
\end{prop}

\subsection{The induced pairing
 of $\nbigpzero_{<\vecdelta}\nbige$ and
 $\sigma^{\ast}\nbigp^{(-\lambda_0)}_{<\vecdelta}\nbige$}
\label{subsection;08.9.27.2}

Let $\lambda_0\in\vecS$.
Note
$\sigma(\lambda_0)=
 -\lambdabar_0^{-1}=-\lambda_0$.
Let $U(\lambda_0)$ be a sufficiently 
small neighbourhood of
$\lambda_0$.
Let $\vecI(\lambda_0):=
U(\lambda_0)\cap\vecS$.
\index{set $\vecI(\lambda_0)$}
We put 
$U(-\lambda_0):=\sigma(U(\lambda_0))$
and
$\vecI(-\lambda_0):=\sigma\bigl(\vecI(\lambda_0)\bigr)$.

\begin{lem}
\label{lem;07.10.18.15}
The restriction of {\rm(\ref{eq;10.5.24.30})}
to $\vecI(\lambda_0)\times(X-D)$
is extended to the following pairing:
\begin{multline}
 \label{eq;07.10.18.10}
  \nbigp^{(\lambda_0)}_{<\vecdelta}
 \nbige_{|\vecI(\lambda_0)\times X}
\otimes
 \sigma^{\ast} 
  \nbigp^{(-\lambda_0)}_{<\vecdelta}
 \nbige_{|\vecI(-\lambda_0)\times X}
\lrarr \\
 \Bigl\{
 F\in C^{\infty}\bigl(\vecI(\lambda_0)\times (X-D)\bigr)
 \,\Big|\,
 |F|=O\Bigl(\prod_{i=1}^n|z_i|^{-2+\epsilon}\Bigr)\,\,
 \mbox{\rm for some } \epsilon>0
 \Bigr\}
\end{multline}
\end{lem}
\pf
We use the notation in 
Section \ref{subsection;07.11.23.1}.
As a preparation,
we show that
$g_{\irr}(-\lambda+\lambda_0)^{\dagger}
\circ
 g_{\irr}(-\lambda+\lambda_0)^{-1}$
is bounded with respect to $h$
on $\vecI(\lambda_0)\times (X-D)$,
if $U(\lambda_0)$ is sufficiently small.
We put
\[
 g'_{\irr}(\lambda):=
 \exp\left(
 \sum_{\gminia\in\Irr(\theta)}
 \lambda\gminiabar\cdot \pi_{\gminia}^{\dagger}
 \right)
=\prod_{i=0}^{L}
 \exp\Bigl(
 \sum_{\gminib\in\Irrbar(\theta,\vecm(i))}
 \lambda\cdot
 \overline{\zeta_{\vecm(i)}(\gminib)}
 \cdot\pi^{\vecm(i)\dagger}_{\gminib}
 \Bigr)
\]
We have the boundedness of
\[
 g_{\irr}(-\lambda+\lambda_0)^{\dagger}
\circ 
 g'_{\irr}(-\lambda+\lambda_0)^{-1}
=\exp\Bigl(
 \sum_{\gminia}
 2\sqrt{-1}\Image\bigl(
 \overline{(\lambda-\lambda_0)}\gminia\bigr)
\cdot\pi^{\dagger}_{\gminia}
 \Bigr)
\]
with respect to $h$.
By the argument in the proof of
Lemma \ref{lem;07.7.20.150},
we can show the boundedness of 
$ g'_{\irr}(-\lambda+\lambda_0)
\circ
 g_{\irr}(-\lambda+\lambda_0)^{-1}$
with respect to $h$
if $|\lambda-\lambda_0|$ is sufficiently small.
Thus, we obtain the boundedness of
$g_{\irr}(-\lambda+\lambda_0)^{\dagger}
\circ
 g_{\irr}(-\lambda+\lambda_0)^{-1}$
with respect to $h$.

Then, we obtain the following estimate
for $f_1\in \nbigpzero_{<\vecdelta}
 \nbige_{|\vecI(\lambda_0)\times X}$
and $f_2\in \nbigp^{(-\lambda_0)}_{<\vecdelta}
 \nbige_{|\vecI(-\lambda_0)\times X}$:
\begin{multline}
\bigl| C(f_1,\sigma^{\ast}f_2)\bigr|=
\bigl|h\bigl(f_1(\lambda,\vecz),
 f_2(-\lambda,\vecz)\bigr)\bigr| \\
=\bigl|
 h\bigl(g_{\irr}(\lambda-\lambda_0)
 f_1(-\lambda,\vecz),
 g_{\irr}(-\lambda+\lambda_0)^{\dagger}
  f_2(-\lambda,\vecz)
 \bigr)
 \bigr| \\
\leq
 \bigl|f_1\bigr|_{\nbigp^{(\lambda_0)}_{\irr}h}
 \cdot \bigl|f_2\bigr|_{\nbigp_{\irr}^{(-\lambda_0)}h}
 \cdot
 \bigl|
 g_{\irr}(-\lambda+\lambda_0)^{\dagger}
 \circ
 g_{\irr}(-\lambda+\lambda_0)^{-1}
 \bigr|_h
\leq
 C
\bigl|f_1\bigr|_{\nbigp^{(\lambda_0)}_{\irr}h}
 \cdot \bigl|f_2\bigr|_{\nbigp_{\irr}^{(-\lambda_0)}h}
\end{multline}
Thus, we obtain Lemma \ref{lem;07.10.18.15}.
\hfill\qed

\subsection{Some estimates}
\label{subsection;07.11.26.12}

Recall that the coordinate system is 
assumed to be admissible for $\Irr(\theta)$.
We take an auxiliary sequence
$\vecm(0),\ldots,\vecm(L)$ for $\Irr(\theta)$.
Let $k$ be determined by
$\vecm(0)\in\seisuu_{<0}^{k}\times\veczero_{n-k}$.
Let $\DD_1$ denote the restriction of $\DD$
to the $z_1$-direction.
Let $S$ be a small multi-sector in
$U(\lambda_0)\times \bigl(X-D(\leq k)\bigr)$.
Note $\sigma(S)$ gives a small multi-sector
in $U(-\lambda_0)\times\bigl(X-D(\leq k)\bigr)$.
We take $\DD_1$-flat splittings
of the full Stokes filtrations
by the procedure in Section
\ref{subsection;08.1.5.1}:
\[
 \nbigp^{(\lambda_0)}_{<\vecdelta}
 \nbige_{|\Sbar}
=\bigoplus_{\gminia\in \Irr(\theta)}
 \nbigp^{(\lambda_0)}_{<\vecdelta}\nbige_{\gminia,S},
\quad\quad
 \nbigp^{(-\lambda_0)}_{<\vecdelta}
 \nbige_{|\sigma(\Sbar)}
=\bigoplus_{\gminia\in \Irr(\theta)}
 \nbigp^{(-\lambda_0)}_{<\vecdelta}\nbige_{\gminia,S}
\]
Let $\vecv_{\gminia}$ be frames of
$\Gr^{\nbigftilde}_{\gminia}\bigl(
 \nbigp^{(\lambda_0)}_{<\vecdelta}\nbige
 \bigr)$.
We have the natural lift 
$\vecv_{\gminia,S}$
of $\vecv_{\gminia|\Sbar}$
to $\nbigp^{(\lambda_0)}_{<\vecdelta}
   \nbige_{\gminia,S}$.
Similarly, let 
$\vecw_{\gminia}$ be frames of
$\Gr_{\gminia}^{\nbigftilde}\bigl(
 \nbigp^{(-\lambda_0)}_{<\vecdelta}\nbige
 \bigr)$,
and let $\vecw_{\gminia,S}$ 
denote the lift
to $\nbigp^{(-\lambda_0)}_{<\vecdelta}
  \nbige_{\gminia,S}$.

\begin{lem}
\mbox{{}}\label{lem;07.11.26.6}
Let $S_{\vecI(\lambda_0)}:=
 S\cap \bigl(\vecI(\lambda_0)\times X\bigr)$.
\begin{itemize}
\item
We have the vanishing
$ h\bigl(v_{\gminia,S,p},
 \sigma^{\ast}(w_{\gminib,S,q})
 \bigr)=0$
unless
$\gminia-\gminib\geq_S0$.
\item
In the case $\gminia-\gminib>_S0$,
we have the following estimate on $S_{\vecI}$
for some positive constants
$C$, $N$ and $\epsilon$:
\[
  h\bigl(v_{\gminia,S,p},
 \sigma^{\ast}w_{\gminib,S,q}
 \bigr)
 \exp\bigl(
-\bigl(\lambda^{-1}+\lambdabar_0\bigr)
 (\gminia-\gminib)
 \bigr) \\
=O\Bigl(
 |z_1|^{-N}
\cdot\prod_{j=2}^{n}
 |z_j|^{-2+\epsilon}
 \Bigr)
\]
\item
In the case $\gminia=\gminib$,
we have 
$h\bigl(v_{\gminia,S,p},
 \sigma^{\ast}w_{\gminia,S,q}
 \bigr)=
O\Bigl(\prod_{j=1}^{n}|z_j|^{-2+\epsilon}\Bigr)$
for some $\epsilon>0$.
\end{itemize}
\end{lem}
\pf
For fixed $\gminia,\gminib\in\Irr(\theta)$,
we put $C_{p,q}:=
 h\bigl(v_{\gminia,S,p},
 \sigma^{\ast}w_{\gminib,S,q}\bigr)$,
and then we obtain the matrix valued function
$\vecC=(C_{p,q})$.

\begin{lem}
\label{lem;07.11.26.5}
We have the estimate
$|\vecC|=O\Bigl(
 \prod_{j=1}^{n}|z_j|^{-2+\epsilon}
 \Bigr)$ for some $\epsilon>0$
on $\vecI(\lambda_0)\times S$.
In particular,
the third claim of
Lemma {\rm{\ref{lem;07.11.26.6}}}
holds.
\end{lem}
\pf
Let $\vecv$ and $\vecw$ be frames of
$\nbigp^{(\lambda_0)}_{<\vecdelta}
 \nbige_{|U(\lambda_0)\times X}$
and $\nbigp^{(-\lambda_0)}_{<\vecdelta}
 \nbige_{|U(-\lambda_0)\times X}$,
respectively.
Let $\vecv_{S}$ and $\vecw_S$
be the frames of 
$\nbigp^{(\lambda_0)}_{<\vecdelta}
 \nbige_{|\Sbar}$
and $\nbigp^{(-\lambda_0)}_{<\vecdelta}
 \nbige_{|\sigma(\Sbar)}$
given by $\vecv_{\gminia,S}$
and $\vecw_{\gminia,S}$
 $(\gminia\in \Irr(\theta))$.
Let $B_1$ and $B_2$
be determined by
$\vecv_S=\vecv\cdot B_1$
and $\vecw_S=\vecw\cdot B_2$.
Then, $B_i$ and $B_i^{-1}$ $(i=1,2)$
are bounded.
Let $\vecC'$ be the matrix valued function
whose $(p,q)$-entry is given by 
$h\bigl(v_{p},\sigma^{\ast}w_q\bigr)$.
By Lemma \ref{lem;07.10.18.15},
$|\vecC'|=
 O\Bigl(\prod_{j=1}^{n}|z_j|^{-2+\epsilon}\Bigr)$.
Then, the claim of Lemma \ref{lem;07.11.26.5}
follows.
\hfill\qed

\vspace{.1in}
Let us return to the proof of Lemma 
\ref{lem;07.11.26.6}.
Let $A_{\gminia}$ be the matrix-valued
holomorphic function
determined by the following:
\begin{equation}
 \label{eq;07.11.26.7}
 \DD^f_{1}\vecv_{\gminia}
=\vecv_{\gminia}\cdot 
\Bigl(
 (\lambda^{-1}+\lambdabar_0)\cdot
 \del_1\gminia\cdot dz_1
+A_{\gminia}\cdot dz_1/z_1\Bigr)
\end{equation}
Similarly, let $B_{\gminib}$ be 
the matrix-valued holomorphic function
determined by the following:
\begin{equation}
 \label{eq;07.11.26.8}
 \DD^f_1\vecw_{\gminib}
=\vecw_{\gminib}\cdot 
\bigl(
 (\lambda^{-1}-\lambdabar_0)\cdot
 \del_1\gminib\cdot dz_1
+B_{\gminib}\cdot dz_1/z_1\bigr)
\end{equation}
Put $r_1:=|z_1|$.
From (\ref{eq;07.11.26.7}) and (\ref{eq;07.11.26.8}),
we obtain the following relation on
$S_{\vecI(\lambda_0)}$:
\[
 r_1\frac{d\vecC}{dr_1}
=\Bigl(
 (\lambda^{-1}+\lambdabar_0)
 z_1\del_1\gminia
-\overline{(\lambda^{-1}+\lambdabar_0)
 z_1\del_1\gminib}
 \Bigr)\cdot\vecC
+\lefttop{t}A_{\gminia}\cdot\vecC
+\vecC\,
 \sigma^{\ast}
 \overline{B}_{\gminib}
\]
Then, the first and second claims of
Lemma \ref{lem;07.11.26.6}
follows from Lemma \ref{lem;07.11.26.5}
and Lemma \ref{lem;07.11.26.10} below.
\hfill\qed

\begin{cor}
\label{cor;08.9.27.20}
The pairing of
$\nbigftilde^{S}_{\gminia}
 \bigl(\nbigp^{(\lambda_0)}_{<\vecdelta}
 \nbige_{|S}\bigr)$
and 
$\sigma^{\ast}\nbigftilde^{\sigma(S)}_{\gminib}
 \bigl(\nbigp^{(-\lambda_0)}_{<\vecdelta}
 \nbige_{|\sigma(S)}\bigr)$
is $0$ unless
$\gminia-\gminib\geq_S0$.
\hfill\qed
\end{cor}
\subsection{The induced pairing of
 $\nbigq^{(\lambda_0)}_{<\vecdelta}\nbige$ and
 $\sigma^{\ast}
 \nbigq^{(-\lambda_0)}_{<\vecdelta}\nbige$}

\begin{lem}
\label{lem;07.11.26.14}
The pairing 
$C$ of $\nbige_{|\vecI(\lambda_0)\times (X-D)}$
and $\sigma^{\ast}\nbige_{|\vecI(-\lambda_0)\times (X-D)}$
is extended to the following:
\begin{multline}
  \label{eq;07.8.2.2}
  \nbigq^{(\lambda_0)}_{<\vecdelta}
 \nbige_{|\vecI(\lambda_0)\times X}
\otimes
 \sigma^{\ast}
 \nbigq^{(-\lambda_0)}_{<\vecdelta}
 \nbige_{|\vecI(-\lambda_0)\times X}
\lrarr \\
 \Bigl\{
 F\in C^{\infty}\bigl(\vecI(\lambda_0)\times (X-D)\bigr)
 \,\Big|\,|F|=O\Bigl(\prod_{i=1}^n|z_i|^{-2+\epsilon}\Bigr)\,\,
\mbox{\rm for some } \epsilon>0
 \Bigr\}
\end{multline}
\end{lem}
\pf
Let $S$, $\vecv_{\gminia,S}$
and $\vecw_{\gminia,S}$ be as in
Section \ref{subsection;07.11.26.12}.
We put 
\[
 \vecvtilde_{\gminia,S}:=
 \vecv_{\gminia,S}
 \cdot\exp(-\lambdabar_0\gminia),
\quad
 \vecwtilde_{\gminia,S}:=
 \vecw_{\gminia,S}
 \cdot\exp(\lambdabar_0\gminia).
\]
They give frames 
$\vecvtilde_{S}$ and $\vecwtilde_S$
of $\nbigq^{(\lambda_0)}_{<\vecdelta}
 \nbige_{|S}$
and $\nbigq^{(-\lambda_0)}_{<\vecdelta}
 \nbige_{|\sigma(S)}$,
respectively.
If $U(\lambda_0)$ is sufficiently small,
we have the following estimate
for some $\epsilon>0$
on $S_{\vecI(\lambda_0)}$,
according to Corollary \ref{lem;07.11.26.6}:
\[
 h(\vtilde_{S,p},\sigma^{\ast}\wtilde_{S,q})
=O\Bigl(\prod_{j}^{n}|z_j|^{-2+\epsilon}\Bigr)
\]
Let $\vecvtilde$ and $\vecwtilde$
be frames of $\nbigq^{(\lambda_0)}_{<\vecdelta}
 \nbige_{|U(\lambda_0)\times X}$
and $\nbigq^{(-\lambda_0)}_{<\vecdelta}
 \nbige_{|U(-\lambda_0)\times X}$,
respectively.
Let $B_1$ and $B_2$ be determined by
$\vecvtilde_S=\vecvtilde\cdot B_1$
and $\vecwtilde_S=\vecwtilde\cdot B_2$.
Then, as remarked in Lemma \ref{lem;10.5.26.10},
$B_1$ and $B_1^{-1}$ are 
bounded on $S_{\vecI(\lambda_0)}$.
We have a similar estimate
for $B_2$.
Hence, we obtain
$h\bigl(\vtilde_p,\sigma^{\ast}\wtilde_q\bigr)
=O\bigl(\prod_{j=1}^{n}|z_j|^{-2+\epsilon}\bigr)$
for some $\epsilon>0$
on $S_{\vecI(\lambda_0)}$.
By varying $S$,
we obtain Lemma \ref{lem;07.11.26.14}.
\hfill\qed

\subsection{Construction of $\gbigc$}
\label{subsection;07.10.19.20}

We can regard (\ref{eq;07.8.2.2}) as the pairing
\[
\gbigc: \nbigq^{(\lambda_0)}_{<\vecdelta}
 \nbige_{|\vecI(\lambda_0)\times X}
\otimes
 \sigma^{\ast}
 \nbigq^{(-\lambda_0)}_{<\vecdelta}
 \nbige_{|\vecI(-\lambda_0)\times X}
\lrarr
 \distribution_{\vecI(\lambda_0)
 \times X/\vecI(\lambda_0)}.
\]
We would like to extend it
to the $\nbigr_{|\vecS\times X}\otimes
     \sigma^{\ast}\nbigr_{|\vecS\times X}$-homomorphism
\[
\gbigc:
 \gbige^{(\lambda_0)}_{|\vecS\times X}\otimes
 \sigma^{\ast}\gbige^{(-\lambda_0)}_{|\vecS\times X}
\lrarr
 \distribution_{\vecS\times X/\vecS}.
\]
The argument is essentially the same 
as that in Section 18.1 of \cite{mochi2}.

\begin{lem}
\label{lem;07.9.7.30}
Assume $\lambda\in \vecI(\lambda_0)$ is generic.
If $\sum_{i=1}^m P_i\cdot u_i=0$
in $\nbigq\nbigelambda$
for $u_i\in \nbigq^{(\lambda_0)}_{<\vecdelta}\nbigelambda$
and $P_i\in\nbigr_X$ $(i=1,\ldots,m)$,
then we have 
$ \sum_{i=1}^m
 P_i\cdot C\bigl(u_i,\sigma^{\ast}v\bigr)=0$
in $\distribution_X$
for any 
$v\in \nbigq^{(-\lambda_0)}_{<\vecdelta}
 \nbige^{-\lambda}$.

Similarly,
if $\sum_i^m Q_i\cdot u_i=0$ 
in $\nbigq\nbige^{-\lambda}$
for $u_i\in
 \nbigq^{(-\lambda_0)}_{<\vecdelta}\nbige^{-\lambda}$
 and $Q_i\in\nbigr_X$ 
 $(i=1,\ldots,m)$,
then we have
 $\sum_{i=1}^m
 \sigma^{\ast}(Q_i)\cdot C(v,\sigma^{\ast}u_i)=0$
in $\distribution_X$
for any 
$v\in\nbigq_{<\vecdelta}^{(\lambda_0)}
 \nbigelambda$.
\end{lem}
\pf
Let $\vecv$ be a frame of 
$\nbigq_{<\vecdelta}^{(-\lambda_0)}\nbigelambda$.
We have 
$\sum P_i C(u_i,\sigma^{\ast}v_j)=0$
on $X-D$.
We set $\vecz^{N\vecdelta}:=\prod_{j=1}^{n}z_j^N$.
There exists a large number $N$ such that
\[
 \sum P_i\cdot C(u_i,
 \sigma^{\ast}(\vecz^{N\vecdelta}v_j))=0
\]
for any $j$.
Hence,
there exists a large $N$ such that
we have $\sum P_i\cdot C(u_i,\sigma^{\ast}v)=0$
for any $v\in 
 \nbigq_{-N\vecdelta}^{(\lambda_0)}\nbigelambda$.

Since we have assumed $\lambda$ is generic,
any $v\in
 \nbigq^{(-\lambda_0)}_{<\vecdelta}\nbigelambda$
have the expression
$v=\sum \deldel^J\cdot v_J$,
where $v_J\in 
  \nbigq^{(-\lambda_0)}_{-N\vecdelta}\nbigelambda$.
Then, the claim of the lemma can be shown
easily.
\hfill\qed

\begin{lem}
\label{lem;07.11.26.15}
If $\sum_{i=1}^m P_i\cdot u_i=0$
in $\gbige_{|\vecI(\lambda_0)\times X}$
for 
$u_i\in \nbigq^{(\lambda_0)}_{<\vecdelta}
   \nbige_{|\vecI(\lambda_0)\times X}$
and $P_i\in\nbigr_X$ $(i=1,\ldots,m)$,
then we have 
$ \sum_{i=1}^m
 P_i\cdot C\bigl(u_i,\sigma^{\ast}v\bigr)=0$
in $\distribution_{\vecI(\lambda_0)\times X/
 \vecI(\lambda_0)}$
for any 
$v\in \nbigq^{(-\lambda_0)}_{<\vecdelta} 
  \nbige_{|\vecI(-\lambda_0)\times X}$.

Similarly,
if $\sum_{i=1}^m Q_i\cdot u_i=0$ 
in $\gbige_{|\vecI(-\lambda_0)\times X}$
for $u_i\in \nbigq^{(-\lambda_0)}_{<\vecdelta}
 \nbige_{|\vecI(-\lambda_0)\times X}$
and $Q_i\in\nbigr_X$ $(i=1,\ldots,m)$,
then we have
 $\sum_{i=1}^m
 \sigma^{\ast}(Q_i)\cdot C(v,\sigma^{\ast}u_i)=0$
in $\distribution_{\vecI(\lambda_0)\times X/
 \vecI(\lambda_0)}$
for any 
$v\in\nbigq_{<\vecdelta}^{(\lambda_0)}
 \nbige_{|\vecI(\lambda_0)\times X}$.
\end{lem}
\pf
It follows from Lemma \ref{lem;07.9.7.30}
and the continuity.
\hfill\qed

\vspace{.1in}
Let us finish the proof of 
Proposition \ref{prop;08.9.27.3}.
Let $\vecv$ and $\vecw$ be frames of
$\nbigq^{(\lambda_0)}_{<\vecdelta}
 \nbige_{|U(\lambda_0)\times X}$
and $\nbigq^{(-\lambda_0)}_{<\vecdelta}
 \nbige_{|U(-\lambda_0)\times X}$,
respectively.
For any $f\in\gbige_{|\vecI(\lambda_0)\times X}$
and $g\in \gbige_{|\vecI(-\lambda_0)\times X}$,
we have the expressions
$f=\sum P_i\cdot v_i$
and $g=\sum Q_j\cdot w_j$
for some $P_i,Q_j\in\nbigr_X$.
We put 
$\gbigc(f,\sigma^{\ast}g):=
 \sum P_i\cdot \sigma^{\ast}Q_j
 \cdot \gbigc(v_i,\sigma^{\ast}w_j)$.
We can check the well-definedness
by using Lemma \ref{lem;07.11.26.15}.
It is easy to check the morphism is
$\nbigr_X\otimes\sigma^{\ast}\nbigr_{X}$-linear.
Thus, we obtain a sesqui-linear pairing
\begin{equation}
 \label{eq;08.9.27.4}
 \gbigc:\gbige_{|\vecI(\lambda_0)\times X}\otimes
\sigma^{\ast}\gbige_{|\vecI(-\lambda_0)\times X}
\lrarr \distribution_{\vecI(\lambda_0)\times X
 /\vecI(\lambda_0)}
\end{equation}
whose restriction to $\vecI(\lambda_0)\times (X-D)$ 
is equal to $C$.

Let $\gbigc'$ be another such 
$\nbigr_X\otimes
 \sigma^{\ast}\nbigr_{X}$-homomorphism.
Then, we have $\gbigc'=\gbigc$
because of the strict $S$-decomposability
of $\gbige$ along the function $\prod_{i=1}^n z_i$
(Proposition \ref{prop;07.10.28.50}).
See \cite{sabbah2},
\cite{mochi2}
or Proposition \ref{prop;a11.23.10} below.

Hence, 
by varying $\lambda_0$
and gluing (\ref{eq;08.9.27.4}),
we obtain the global sesqui-linear pairing
\[
 \gbigc:\gbige_{|\vecS\times X}\otimes
\sigma^{\ast}\gbige_{|\vecS\times X}
\lrarr \distribution_{\vecS\times X/\vecS}
\]
with the desired property.
Again,
the uniqueness follows from
the strict $S$-decomposability
of $\gbige$ along the function $\prod_{i=1}^n z_i$.
Thus, we obtain Proposition
\ref{prop;08.9.27.3}.
\hfill\qed

\section[A characterization in the one 
 dimensional case]{A characterization of 
 the prolongment in the one dimensional case}
\label{subsection;08.9.28.12}
\subsection{Statement}
\label{subsection;10.5.20.2}

We give a remark on the uniqueness 
of $\nbigr$-triples extending 
a variation of polarized pure twistor structure
in the one dimensional case.

Let $X:=\Delta$ and $D:=\{O\}$.
Let $\nbigt=(\nbigm,\nbigm,C)$ be 
a strictly specializable $\nbigr(\ast z)$-triple
on $X$ such that
$\nbigs=(\id,\id):\nbigt\lrarr\nbigt^{\ast}$
is a Hermitian sesqui-linear duality
of weight $0$.
Recall that the underlying 
family of meromorphic $\lambda$-flat bundles
$\nbigm$ has the KMS-structure
at each $\lambda_0\in\cnum$
due to a lemma of Sabbah \cite{sabbah5}.
(It is reviewed in Lemma \ref{lem;08.1.28.125}.)
For simplicity, we assume that it is unramified.
Let $\nbigi$ denote the set of irregular values.

To simplify the claim, we assume that
the restriction of $\nbigt$ to $X\setminus D$
comes from a harmonic bundle
$\harmonicbundle$.
By the assumption,
we obtain that $\harmonicbundle$ is
unramifiedly good wild.
According to Proposition \ref{prop;08.9.27.3},
we have the $\nbigr(\ast z)$-triple
$\gbigt(E)(\ast z):=
 \bigl(\nbigq\nbige,\nbigq\nbige,\gbigc(\ast z)\bigr)$ 
associated to
$\harmonicbundle$.
We would like to compare 
$\gbigt(E)(\ast z)$ and $\nbigt$.

\begin{thm}
\label{thm;10.5.20.4}
Assume that 
$(P\Gr^W_{\ell}\psitilde_{z,u,\gminia}(\nbigt),
 \nbigs_{\gminia,u,\ell})$
 are polarized pure twistor structure of weight $\ell$
for any $\gminia\in\nbigi$,
 $u\in\real\times\cnum$
and $\ell\in\seisuu_{\geq 0}$.
Then, we have a natural isomorphism
$\gbigt(E)(\ast z)\simeq\nbigt$.
\end{thm}

\begin{rem}
Conversely,
if $\nbigt\simeq\gbigt(E)(\ast z)$,
$(P\Gr^W_{\ell}\psitilde_{z,u,\gminia}(\nbigt),
 \nbigs_{\gminia,u,\ell})$
 are polarized pure twistor structure of weight $\ell$
for any $\gminia\in\nbigi$,
 $u\in\real\times\cnum$
and $\ell\in\seisuu_{\geq 0}$.
See Corollary {\rm\ref{cor;07.10.19.35}}.
\hfill\qed
\end{rem}

\begin{rem}
According to a result by Sabbah
in \cite{sabbah5},
we do not have to assume that
$(\nbigt,\nbigs)_{|X\setminus D}$ comes from
a harmonic bundle.
Namely, if
$(\nbigt,\nbigs)_{|X\setminus D}$ comes from
a variation of twistor structure
with a pairing of weight $0$,
and if the assumption of the theorem is satisfied,
the variation of twistor structure
comes from a harmonic bundle.
(It also follows from our argument below.)
\hfill\qed
\end{rem}

\subsection{Reduction of $\nbigr(\ast z)$-triples}
\label{subsection;07.12.20.7}

Let $\nbigt=(\nbigm,\nbigm,C)$ 
be an $\nbigr_X(\ast z)$-triple on $X$
such that $\nbigs=(\id,\id)$ gives
a Hermitian sesqui-linear duality of weight $0$.
We set $\nbigx:=\cnum_{\lambda}\times X$
and $\nbigd:=\cnum_{\lambda}\times D$.
Assume the following:
\begin{itemize}
\item
 $\nbigm$ is an unramifiedly good
 meromorphic prolongment of
 $\nbigm_{|\nbigx\setminus\nbigd}$
 in the sense of Definition \ref{df;10.5.19.20}.
 The set of irregular values is denoted by
 $\nbigi$.
\item $(\nbigt,\nbigs)_{|X\setminus D}$
 comes from a variation of twistor structure
 with a pairing of weight $0$.
 Namely, 
 $\nbigm_{|\nbigx\setminus \nbigd}$ comes from
 a family of $\lambda$-flat bundles,
 and $C$ is the restriction of
 a $\lambda$-holomorphic sesqui-linear pairing.
\end{itemize}

For any $\gminia\in\nbigi$,
we have the full reduction
$\nbigm_{\gminia}:=
 \Gr^{\nbigftilde}_{\gminia}(\nbigm)$
on the product of $\cnum_{\lambda}$
and a neighbourhood $X'$ of $D$.
For simplicity of description,
we replace $X$ with $X'$.

\subsubsection{}
Let $\lambda_0\in\vecS$.
Let $\nbigx^{(\lambda_0)}$ denote
a small neighbourhood of
$\{\lambda_0\}\times X$ in $\nbigx$.
We use the symbol $\nbigd^{(\lambda_0)}$
in a similar meaning.
Let $S$ be a small sector in 
$\nbigx^{(\lambda_0)}\setminus\nbigd^{(\lambda_0)}$.
Then, $\sigma(S)$
is a small sector in
$\nbigx^{(-\lambda_0)}\setminus
 \nbigd^{(-\lambda_0)}$.
Note that 
$\gminia\leq_S\gminib$ if and only if
$\gminia\geq_{\sigma(S)}\gminib$
for any $j$ and any $\gminia\in \nbigi$.
By Lemma \ref{lem;07.7.11.11} below,
we obtain an induced Hermitian sesqui-linear
pairing 
\[
 C_{0,\gminia}:
 \nbigm^{\circ}_{\gminia|\vecS\times(X\setminus D)}
 \otimes
 \sigma^{\ast}\nbigm^{\circ}_{\gminia|\vecS\times 
(X\setminus D)}
\lrarr \distribution_{\vecS\times(X\setminus D)/\vecS}.
\]
In particular,
we obtain an $\nbigr$-triple
$\Gr_{\gminia}(\nbigt)^{\circ}
:=(\nbigm_{\gminia}^{\circ},
 \nbigm_{\gminia}^{\circ},
 C_{0,\gminia})$ on $X\setminus D$
for $\gminia\in \nbigi$,
which is equipped with
a Hermitian sesqui-linear duality
$\Gr_{\gminia}(\nbigs)^{\circ}
 :=(\id,\id)$.
We will prove the following lemma in 
Subsection \ref{subsection;10.5.20.4}.
\begin{lem}
\label{lem;07.12.20.1}
$C_{0,\gminia}$
is naturally extended to 
a Hermitian sesqui-linear pairing:
\[
 C_{\gminia}:
 \nbigm_{\gminia|\vecS\times X}
\otimes
 \sigma^{\ast}\nbigm_{\gminia|\vecS\times X}
\lrarr \distribution_{\vecS\times X/\vecS}(\ast z)
\]
In particular, we obtain an $\nbigr(\ast z)$-triple
$\Gr_{\gminia}(\nbigt):=
 (\nbigm_{\gminia},\nbigm_{\gminia},C_{\gminia})$
on $X$ for $\gminia\in \nbigi$
with a Hermitian sesqui-linear duality
$\Gr_{\gminia}(\nbigs)=(\id,\id)$.
\end{lem}

\subsubsection{}
We assume that
$\nbigm$ is strictly specializable
and unramified.
(See a remark in Subsection
 \ref{subsection;10.5.20.2}.)
Since we have the natural isomorphism
$\nbigm_{|\nbigdhat}\simeq
\bigoplus_{\gminia\in \nbigi}
 \nbigm_{\gminia|\nbigdhat}$,
the induced $\nbigr$-triples
$\Gr_{\gminia}(\nbigt)$ are also
strictly specializable along $z=0$.
Moreover, we have natural isomorphisms
\begin{equation}
 \label{eq;10.5.20.3}
 \psitilde_{z,\gminia,u}(\nbigm)
 \simeq\psitilde_{z,\gminia,u}(\Gr_{\gminia}\nbigm)
\end{equation}
We will prove the following lemma
in Subsection \ref{subsection;10.5.20.5}.
\begin{lem}
\label{lem;07.12.20.9}
Under {\rm(\ref{eq;10.5.20.3})},
we have
$\psitilde_{z,\gminia,u}C
=\psitilde_{z,\gminia,u}C_{\gminia}$.
Namely,
we have a natural isomorphism
$\psitilde_{z,\gminia,u}\nbigt
\simeq 
 \psitilde_{z,\gminia,u}
 \bigl(\Gr_{\gminia}(\nbigt)\bigr)$
for any $\gminia\in \nbigi$.
\end{lem}

\subsubsection{}
Let $\nbigm^{\lor}$ denote the dual of
$\nbigm$ as a family of meromorphic 
$\lambda$-flat bundles.
Let $\Upsilon(\nbigt^{\circ},\nbigs^{\circ})$
be the variation of twistor structure 
with the pairing of weight $0$,
obtained as the gluing of
$\nbigm_{|\nbigx\setminus\nbigd}$
and 
$\sigma^{\ast}\nbigm^{\lor}_{|\nbigx\setminus\nbigd}$
by the isomorphism induced by $C$.

\begin{lem}
\label{lem;08.9.14.110}
The Stokes structures of
$\nbigm$ and $\sigma^{\ast}\nbigm^{\lor}$
are the same in the sense of 
Definition {\rm\ref{df;07.11.13.20}}.
Namely,
$(\nbigm,\sigma^{\ast}\nbigm^{\lor})$
is an unramifiedly good meromorphic prolongment
of $\Upsilon(\nbigt^{\circ},\nbigs^{\circ})$.
\end{lem}
\pf
Due to the existence of $C$
and Lemma \ref{lem;07.7.11.11} below,
the Stokes filtrations of
$\nbigm$ and $\sigma^{\ast}\nbigm$
are mutually dual on $\vecS\times\Xtilde(D)$.
Hence
their Stokes filtrations 
of $\nbigm$ and $\sigma^{\ast}\nbigm^{\lor}$
are the same
on $\vecS\times \Xtilde(D)$.
Then, it is easy to observe
that they are the same on
$\cnum_{\lambda}^{\ast}\times\Xtilde(D)$.
\hfill\qed

\vspace{.1in}

We have a variation of twistor structure
with an induced symmetric pairing
$\Gr_{\gminia}\bigl(
 \Upsilon(\nbigt^{\circ},\nbigs^{\circ})
\bigr)$
on $\proj^1\times(X\setminus D)$
as explained in Section \ref{subsection;07.10.10.2}.
By construction,
it is naturally isomorphic to
$\Upsilon\bigl(
 \Gr_{\gminia}(
 \nbigt^{\circ},\nbigs^{\circ})
 \bigr)$,
and $(\nbigm_{\gminia},
 \sigma^{\ast}\nbigm_{\gminia}^{\lor})$
gives a meromorphic prolongment of
$\Gr_{\gminia}\bigl(
 \Upsilon(\nbigt^{\circ},\nbigs^{\circ})
\bigr)$.
Note that $\nbigm_{\gminia}$ are
$\gminia$-regular.

\subsubsection{Proof of Lemma \ref{lem;07.12.20.1}}
\label{subsection;10.5.20.4}
Let $\lambda_0\in\vecS$.
Let $f$ be a section of
$\nbigm_{\gminia|\nbigx^{(\lambda_0)}}$,
and $g$ be a section of
$\nbigm_{\gminia|\nbigx^{(-\lambda_0)}}$.
We shall show that 
$C_{0,\gminia}\bigl(f,\sigma^{\ast}g \bigr)$
gives a section of
$\distribution_{\vecS\times X/\vecS}(\ast z)$
on $\vecI(\lambda_0)\times X$,
where $\vecI(\lambda_0):=U(\lambda_0)\cap \vecS$.

Let $E^{(\lambda_0)}$ be a good lattice of 
$\nbigm_{|\nbigx^{(\lambda_0)}}$.
We may assume that
$f$ is a section of
$\Gr_{\gminia}(E^{(\lambda_0)})$.
Let us take a covering 
$\nbigx^{(\lambda_0)}\setminus\nbigd^{(\lambda_0)}
=\bigcup_{i=1}^N S_i$ by small sectors
on which we have flat splittings
$E^{(\lambda_0)}_{|\Sbar_i}=
 \bigoplus_{\gminia\in \nbigi}
 E^{(\lambda_0)}_{\gminia,S_i}$
of the full Stokes filtration
$\nbigf^{S_i}$ of $E^{(\lambda_0)}_{|\Sbar_i}$.
We have the natural isomorphism
$E^{(\lambda_0)}_{\gminia,S_i}
\simeq
 \Gr_{\gminia}(E^{(\lambda_0)})_{|\Sbar_i}$.
Hence, we obtain a lift of $f_{|\Sbar_i}$
to $E^{(\lambda_0)}_{\gminia,S_i}$.
By gluing them in $C^{\infty}$
as in Section \ref{subsection;07.6.16.8},
we obtain a $C^{\infty}$-section
$\ftilde$ of $E^{(\lambda_0)}$,
which is called a $C^{\infty}$-lift of $f$
to $E^{(\lambda_0)}$.
Similarly,
we take a good lattice $E^{(-\lambda_0)}$
of $\nbigm_{|\nbigx^{(-\lambda_0)}}$
such that $g$ is a section
of $\Gr_{\gminia}(E^{(-\lambda_0)})$,
and we take a $C^{\infty}$-lift $\gtilde$
of $g$ to $E^{(-\lambda_0)}$.

Let $E^{(\lambda_0)}_{C^{\infty}}$
and $E^{(-\lambda_0)}_{C^{\infty}}$
denote the sheaf of $C^{\infty}$-sections of
$E^{(\lambda_0)}$ and $E^{(-\lambda_0)}$,
respectively.
The pairing 
$C:E^{(\lambda_0)}\otimes
 \sigma^{\ast}E^{(-\lambda_0)}\lrarr
 \distribution_{\vecS\times X/\vecS}(\ast z)$
can naturally be extended 
to $C:E^{(\lambda_0)}_{C^{\infty}}\otimes
 \sigma^{\ast}E^{(-\lambda_0)}_{C^{\infty}}\lrarr
 \distribution_{\vecS\times X/\vecS}(\ast z)$.
Hence, we have the section
$C(\ftilde,\sigma^{\ast}\gtilde)$
of $\distribution_{\vecS\times X/\vecS}(\ast z)$.
Because 
$C(\ftilde,\sigma^{\ast}\gtilde)
=C_{\gminia}(f,\sigma^{\ast}g)$
by construction of $\ftilde$ and $\gtilde$,
we are done.
\hfill\qed

\subsubsection{Proof of Lemma \ref{lem;07.12.20.9}}
\label{subsection;10.5.20.5}

By considering tensor product
with $\nbigl(-\gminia)$,
we may reduce the problem to the case $\gminia=0$.
We obtain the isomorphism of
the underlying $\nbigr$-modules
$\psitilde_{z,u}(\nbigm)
\simeq
 \psitilde_{z,u}(\nbigm_{0})$
from the isomorphism
$\nbigm_{|\nbigdhat}
\simeq
 \bigoplus_{\gminia\in T}
 \nbigm_{\gminia|\nbigdhat}$.
Let us compare the specializations of
sesqui-linear pairings
$C$ and $C_{0}$.
Let $\lambda_0\in\vecS$.
Let $[f_1]$ be a section of
$\psitilde_{z,u}
\bigl(\nbigm_{0}\bigr)_{|
\nbigd^{(\lambda_0)}}$.
Let $[f_2]$ be a section of
$\psitilde_{z,u}\bigl(
 \nbigm_{0}\bigr)_{|
 \nbigd^{(-\lambda_0)}}$.
We would like to show
$\psitilde_{z,u}C\bigl([f_1],
 \sigma^{\ast}[f_2]\bigr)
=\psitilde_{z,u}C_{0}\bigl([f_1],
 \sigma^{\ast}[f_2]\bigr)$.

We take $f_1\in 
 \Vzero_{\paramap(\lambda_0,u)}
 (\nbigm_{0|\nbigx^{(\lambda_0)}})$
which is a lift of $[f_1]$
in the sense that
$f_1$ induces the element
$[f_1]\in \psitilde_{z,u}(\nbigm_0)
\subset 
 \Gr^{\Vzero}_{\paramap(\lambda_0,u)}
 (\nbigm_0)$.
We also take $f_2\in\Vzero_{\paramap(-\lambda_0,u)}
 (\nbigm_{0|\nbigx^{(-\lambda_0)}})$
which is a lift of $[f_2]$.
Let $\chi$ be a test function on $X$
which is constantly $1$ around $0$.
We put $\omega_0:=
 \sqrt{-1}\, dz\, d\zbar/2\pi$.
By definition, we have the following:
\[
\bigl\langle
 \psitilde_{z,u}C_{0}([f_1],\sigma^{\ast}[f_2]),\,
 \rho
\bigr\rangle
=\underset{s+\eigenmap(\lambda,u)}{\Res}\,
 \Bigl\langle
 C_{0}(f_1,\sigma^{\ast}f_2),\,\,
 |z|^{2s}\,\chi(z)\, \rho\, \omega_0
 \Bigr\rangle
\]

We take a lattice $E^{(\lambda_0)}$ of 
$\nbigm_{|\nbigx^{(\lambda_0)}}$
such that
(i) 
$f_1$ is a section of 
$\Gr^{\nbigftilde}_{0}(E^{(\lambda_0)})$,
(ii)
$E^{(\lambda_0)}$ is contained in
$\Vzero_{\paramap(\lambda_0,u)}\nbigm
 _{|\nbigx^{(\lambda_0)}}$.
We take a $C^{\infty}$-lift $\ftilde_1$
of $f_1$ to $E^{(\lambda_0)}$
as in the proof of Lemma \ref{lem;07.12.20.1}.
Similarly,
we take a lattice $E^{(-\lambda_0)}$ of
$\nbigm_{|\nbigx^{(-\lambda_0)}}$
such that
(i) $f_2$ is a section of 
$\Gr^{\nbigftilde}_{0}(E^{(-\lambda_0)})$,
(ii) $E^{(-\lambda_0)}$ is contained in
 $V^{(-\lambda_0)}_{\paramap(-\lambda_0,u)}
 \bigl(\nbigm_{|\nbigx^{(-\lambda_0)}}\bigr)$.
We take a $C^{\infty}$-lift $\ftilde_2$ of $f_2$
to $E^{(-\lambda_0)}$.
Then, we have 
$C_{0}(f_1,\sigma^{\ast}f_2)
=C(\ftilde_1,\sigma^{\ast}\ftilde_2)$.

We have 
the section $[\ftilde_1]$ of
$\psitilde_{z,u}(\nbigm)_{|\nbigd^{(\lambda_0)}}$
induced by $[f_1]$
and the isomorphism
$\psitilde_{z,u}(\nbigm)
\simeq
 \psitilde_{z,u}(\nbigm_0)$.
Similarly,
$[f_2]$ naturally induces
a section $[\ftilde_2]$ of
$\psitilde_{z,u}(\nbigm)_{|\nbigd^{(-\lambda_0)}}$.
Let us show the following equality:
\begin{equation}
 \label{eq;07.12.20.2}
\underset{s+\eigenmap(\lambda,u)}{\Res}\,
 \bigl\langle
 C(\ftilde_1,\sigma^{\ast}\ftilde_2),\,\,
 |z|^{2s}\,\chi\, \rho\,\omega_0
 \bigr\rangle
=\bigl\langle
 \psitilde_{z,u}C([\ftilde_1],\sigma^{\ast}[\ftilde_2]),\,
 \rho
 \bigr\rangle
\end{equation}
We can take a lift $\ftilde'_1$ of $[\ftilde_1]$
to $E^{(\lambda_0)}$
such that 
$\ftilde'_1-\ftilde_1=
O(|z|^{N})$ for some sufficiently large $N$.
Similarly, we can take a lift
$\ftilde'_2$ of $[\ftilde_2]$ to $E^{(-\lambda_0)}$
such that
$\ftilde'_2-\ftilde_2=
O(|z|^{N})$ for some sufficiently large $N$.
Then, it can be shown that
both sides of (\ref{eq;07.12.20.2})
are equal to the following:
\[
 \underset{s+\eigenmap(\lambda,u)}{\Res}\,
 \bigl\langle
 C(\ftilde_1',\sigma^{\ast}\ftilde_2'),\,
 |z|^{2s}\,\chi\,\rho\,\omega_0
 \bigr\rangle
\]
Thus, we obtain
$\psitilde_{z,u}C=\psitilde_{z,u}C_{0}$.
\hfill\qed

\subsection{Proof of Theorem
\ref{thm;10.5.20.4}}

Let us return to the setting in 
Subsection \ref{subsection;10.5.20.2}.
According to Lemma \ref{lem;07.12.20.9},
$P\Gr^{W(N)}_j\psitilde_{z,\gminia,u}
 (\nbigt_{\gminia})$ 
and 
$P\Gr^{W(N)}_j\psitilde_{z,\gminia,u}
 (\nbigt)$
are naturally isomorphic,
compatible with the naturally induced
Hermitian sesqui-linear dualities.
Hence, we obtain that
$\Upsilon(\Gr_{\gminia}(\nbigt^{\circ},\nbigs^{\circ}))
=\Gr_{\gminia}(\Upsilon(\nbigt^{\circ},\nbigs^{\circ}))$
are variations of 
polarized pure twistor structure
from the assumption in the theorem,
according to a result by Sabbah \cite{sabbah2},
(It also follows from 
Lemma \ref{lem;07.7.11.12}
and Corollary \ref{cor;08.1.29.2} below.)
Then, the claim of Theorem \ref{thm;10.5.20.4}
follows from Theorem \ref{thm;07.10.11.10}.
\hfill\qed

\subsection{Strictly specializable 
holonomic $\nbigr(\ast z)$-module on a disc
(Appendix)}
\label{subsection;08.9.28.122}

We recall a lemma due to Sabbah.
It relates strict specializability of
an $\nbigr_X$-module
with KMS-structure of
the corresponding  family of
meromorphic $\lambda$-flat bundles.
Let $X:=\Delta$ and $D:=\{0\}$.
Let $\nbigk$ be a neighbourhood of 
$\lambda_0\in\cnum_{\lambda}$.
We put $\nbigx:=\nbigk\times X$
and $\nbigd:=\nbigk\times D$.

Let $\nbigm$ be a strict coherent holonomic
$\nbigr_X(\ast z)$-module on $\nbigx$
whose characteristic variety 
is contained in 
$\nbigk\times\bigl(
 N_X^{\ast}X\cup N_D^{\ast}X\bigr)$,
where $N_X^{\ast}X$ denotes the $0$-section
and $N_D^{\ast}X$ denotes the conormal
bundle of $D$ in $X$.
Then, $\nbigm$ gives a family of
meromorphic $\lambda$-flat bundles.
Moreover, we assume that
it is strictly specializable along $z$
with ramification and exponential twist.
(See \cite{sabbah2} or Sections 
\ref{subsection;08.9.4.11}--\ref{subsection;08.9.4.12}
 below.)
We would like to show the following lemma.
Although
it was already shown by Sabbah in \cite{sabbah5},
we would like to give a partially different proof
based on the arguments in \cite{malgrange},
\cite{sabbah4} and \cite{sabbah5},
just for our understanding.

\begin{lem}
\label{lem;08.1.28.125}
$\nbigm$ has the KMS-structure at $\lambda_0$,
i.e.,
we have the locally free $\nbigo_{\nbigx}$-submodules
$\nbigpzero_a(\nbigm)\subset\nbigm$
$(a\in\real)$
such that
(i) $\bigcup_{a}\nbigpzero_a\nbigm=\nbigm$,
(ii) 
$\nbigpzero_{\ast}\nbigm
=\bigl(\nbigpzero_a\nbigm\,\big|\,
 a\in\real \bigr)$ is 
a good family of filtered $\lambda$-flat bundles
with the KMS-structure at $\lambda_0$.
\end{lem}
\pf
For any $\lambda_1$,
let $\nbigm^{\lambda_1}$ denote the restriction
$\nbigm/(\lambda-\lambda_1)\nbigm$.
For simplicity, we assume that
$\nbigm^{\lambda_0}$ is unramified.
(The general case can be reduced to this case easily.)
We have the irregular decomposition
$\nbigm^{\lambda_0}_{|\Dhat}
=\bigoplus_{\gminia\in \nbigi}
 \nbigmhat^{\lambda_0}_{\gminia}$
with the good set of irregular values
$\nbigi\subset M(X,D)/H(X)$.

\begin{lem}
\mbox{{}}
\begin{itemize}
\item
$\nbigm^{\lambda}$ is also unramified
for each $\lambda$.
\item
The set of the irregular values of 
$\nbigm^{\lambda}$ is given by $\nbigi$.
\item
Let $\nbigm^{\lambda}_{|\Dhat}
=\bigoplus_{\gminia\in \nbigi}
 \nbigmhat^{\lambda}_{\gminia}$
be the irregular decomposition.
Then, $\rank(\nbigmhat_{\gminia}^{\lambda})
 =\rank(\nbigmhat_{\gminia}^{\lambda_0})$
\end{itemize}
\end{lem}
\pf
We give only an outline.
See \cite{sabbah5} for more details.
The following equality can be shown:
\[
 \rank(\nbigmhat^{\lambda}_{\gminia})
=\sum_{u\in \real\times\cnum/\seisuu\times\{0\}}
 \rank\psitilde_{z,\gminia,u}(\nbigm)_{|\lambda}
\]
Then, the claim of the lemma follows.
\hfill\qed

\vspace{.1in}
Let us return to the proof of Lemma \ref{lem;08.1.28.125}.
We put 
\[
 \Vzero_{\gminia,a}(\nbigm):=
 \Vzero_{a}\bigl(\nbigm\otimes\nbigl(-\gminia)\bigr)
 \otimes\nbigl(\gminia),
\quad
 \nbigpzero_{a}(\nbigm):=
 \bigcap_{\gminia\in T} 
 \Vzero_{\gminia,a-1}(\nbigm).
\]
We shall show
$\nbigpzero_a(\nbigm)$ is a good lattice of $\nbigm$.
Using the arguments in \cite{malgrange}
and \cite{sabbah4},
we can show the following general result:
\begin{itemize}
\item
 If we shrink $\nbigk$ appropriately,
 we have a ramified covering
 $\phi_{e}:\nbigk'\times X'\lrarr \nbigk\times X$ given by
 $\phi_e(\lambda,z)=
 \bigl(\lambda_0+(\lambda-\lambda_0)^e,
 z^e\bigr)$
 such that the restriction of
 $\phi_e^{\ast}\nbigm$ to
 $\nbigk^{\prime\ast}\times X'$
 has an unramifiedly good lattice,
 where 
 $\nbigk^{\prime\ast}:=\nbigk'\setminus\{\lambda_0\}$.
\end{itemize}
We have already known that
each $\nbigm^{\lambda}$ has a good lattice,
and the set of the irregular values is $T$.
Hence, it can be shown that
we do not have to take a ramified covering,
i.e.,
$\nbigm_{|\nbigk^{\ast}\times X}$
has an unramifiedly good lattice,
where $\nbigk^{\ast}:=
 \nbigk\setminus\{\lambda_0\}$.
Then, it is easy to see that
$\nbigpzero_{a}(\nbigm)_{|\nbigk^{\ast}\times X}$ 
is a good lattice of $\nbigm_{|\nbigk^{\ast}\times X}$.

Take small numbers
$0<\epsilon_1<\epsilon_2$
such that 
$\nbigk_0:=\bigl\{
 \lambda\,\big|\,
 \epsilon_1\leq|\lambda-\lambda_0|
 \leq \epsilon_2 \bigr\}$
is contained in $\nbigk$.
We can take a locally free $\nbigo_X$-submodule
$\nbigl\subset \nbigm$
such that 
$\nbigpzero_a(\nbigm)_{|\nbigk_0\times X}
\subset \nbigl_{|\nbigk_0\times X}$.
Because of the $\nbigo_{\nbigx}$-local freeness of 
$\nbigl$,
we obtain 
$\nbigpzero_{a}(\nbigm)
 \subset\nbigl$
by using Hartogs theorem.
Because $\nbigpzero_{a}(\nbigm)$ is
the limit of $\nbigo_X$-coherent submodules of
$\nbigpzero_a(\nbigm)$,
we can conclude that
$\nbigpzero_a(\nbigm)$ is $\nbigo_X$-coherent.

Let $\nbigpzero_a(\nbigm)^{\lor\lor}$
denote the $\nbigo_X$-submodule of
$\nbigm$ generated by $\nbigpzero_a(\nbigm)$.
It is $\nbigo_X$-locally free,
and
$\nbigpzero_a(\nbigm)^{\lor\lor}_{|\nbigx-\{(\lambda_0,O)\}}
=\nbigpzero_a(\nbigm)_{|\nbigx-\{(\lambda_0,O)\}}$.
We have already known that
$\nbigpzero_a(\nbigm)^{\lor\lor}_{|\nbigk^{\ast}\times X}$
is a good lattice.
Then,
by using the non-degeneration of the irregular values
and a well established argument in \cite{levelt},
it is easy to obtain that
$\nbigpzero_a(\nbigm)^{\lor\lor}$ 
is a good lattice of $\nbigm$.
Once we know that $\nbigm$ has good lattices,
it is easy to show
$\nbigpzero_a(\nbigm)^{\lor\lor}
=\nbigpzero_a(\nbigm)$.
By construction,
the family $\nbigpzero_{\ast}(\nbigm)$ of
good filtered $\lambda$-flat bundles
has the KMS-structure at $\lambda_0$.
Thus, we obtain Lemma \ref{lem;08.1.28.125}.
\hfill\qed

\subsection{Preliminary for sesqui-linear pairing
(Appendix)}
\label{subsection;08.9.29.120}

\subsubsection{Preliminary for 
the compatibility with Stokes structure}
\label{subsection;08.9.29.130}

Let $S$ be a small sector
$\bigl\{
 r\, e^{\sqrt{-1}\theta}\,\big|\,
 \theta_1<\theta<\theta_2,\,\,0<r<r_0
 \bigr\}$
in $\Delta^{\ast}$.
Let $\gminia,\gminib\in\cnum[z^{-1}]$ such that
$\Re(\gminia+\gminibbar)>0$ on $S$.

\begin{lem}
\label{lem;07.7.11.10}
There does not exist a distribution $F$ on $\Delta$
such that
$F_{|S}=
 \exp\bigl(\gminia+\gminibbar
 +\alpha\log z+\beta \log\zbar\bigr)$.
\end{lem}
\pf
Let $\epsilon>0$ be a sufficiently small positive number.
Let $\theta_1',\theta_2'$ be
$\theta_1<\theta_1'<\theta_2'<\theta_2$.
Let $\gbigf\subset C^{\infty}(\Delta)$ be the subspace
generated by the functions of the following form:
\[
 \exp\Bigl(-\frac{1}{r^{2\epsilon}}\Bigr)
 \cdot
 \rho(\theta)\cdot
 \frac{1}{r^A}
 Q(z,\zbar)
\]
Here,
(i) $\rho(\theta)$ is a $C^{\infty}$-function 
 whose support is contained in
 $[\theta_1',\theta_2']$,
(ii)
 $Q$ is a polynomial,
(iii) $A\in\real$.
The space $\gbigf$ is preserved by
$\del_z$ and $\del_{\zbar}$.

For $n\in\seisuu_{\leq 0}$,
let $L_{n}^2$ denote the space of
distributions which are contained
in the dual of the Sobolev space $L_{-n}^2$ on $\Delta$.
Let $F$ be a distribution on $\Delta$ such that
$F_{|S}=\exp\bigl(
 \gminia+\gminibbar+\alpha\log z+\beta\log\zbar\bigr)$.
We have some $n\in\seisuu_{<0}$
such that $F\in L_{n}^2$.
For any $g\in \gbigf$,
we have $g\cdot F\in L_n^2$.
We have the following:
\[
 \del_z(g\cdot F)_{|\Delta^{\ast}}
=(\del_z g)\cdot F_{|\Delta^{\ast}}
+g\cdot(\del_z F)_{|\Delta^{\ast}}
=\del g\cdot F_{|\Delta^{\ast}}
+\bigl(g\cdot(\del_z \gminia+\alpha\cdot z^{-1})\bigr)
\cdot F_{|\Delta^{\ast}}
\]
Note that $\del_zg$ and $g\cdot \del_z\gminia$
are contained in $\gbigf$.
By a similar calculation,
there exist $\alpha\in\gbigf$
such that 
$ \delbar\del\bigl(g\cdot F\bigr)_{|\Delta^{\ast}}
=\alpha\cdot F_{|\Delta^{\ast}}$
on $\Delta^{\ast}$.
Hence, 
there exist large integer $N$ and 
$\beta\in\gbigf$ 
such that 
$ \delbar\del\bigl(z^{N}\cdot g\cdot F\bigr)=\beta\cdot F$
on $\Delta$.
Using an inductive argument,
we can show that
there exist large number $N(\ell)$
and $\beta^{(\ell)}\in\gbigf$ 
for each $\ell\in\seisuu_{>0}$
such that
$(\delbar\del)^{\ell}(z^{N(\ell)}g\cdot F)
=\beta^{(\ell)}\cdot F
\in L_{n}^2$.
Then, we obtain
$z^{N(\ell)}g\cdot F$ is locally $L_{n+2\ell}^2$.
Hence, if $\ell$ is sufficiently large,
$z^{N(\ell)}g\cdot F$ has to be locally $L^2$.
But, we can directly show that
$z^N\cdot g\cdot F$ is not $L^2$ for any $N$
if $g\neq 0$.
Thus, we have arrived at the contradiction.
\hfill\qed

\vspace{.1in}

The lemma can also be shown by a more direct
argument to take a sequence of test functions
$\rho_n$ satisfying
(i) the supports of $\rho_n$ are contained in $S$,
(ii) $\rho_n\to 0$,
(iii) $\langle F,\rho_n\rangle\to\infty$.

\subsubsection{Compatibility of 
sesqui-linear pairing with Stokes filtrations}
\label{subsection;08.9.29.131}

Let $(\nbigv_{i},\nabla_i)$ $(i=1,2)$
be meromorphic flat bundles on $(\Delta,O)$,
i.e.,
$\nbigv_i$ are locally free $\nbigo(\ast O)$-modules
with connections
$\nabla_i:\nbigv_i\lrarr
 \nbigv_i\otimes\Omega^1_{\Delta}$.
For simplicity, we assume that they are unramified.
Let $\distribution_{\Delta}(\ast O)$ denote the sheaf of
distributions on $\Delta^{\ast}$ 
 with moderate growth at $O$
(see \cite{sabbah4} or 
 a review in Section \ref{subsection;08.1.15.20}
 below),
and let $\nbigd_{\Delta}(\ast O)$ denote the sheaf
of meromorphic differential operators on $\Delta$
whose poles are contained in $O$.
Let $\Delta^{\dagger}$
denote the conjugate of $\Delta$.
Let $C:\nbigv_1\otimes \overline{\nbigv_2}\lrarr 
 \distribution_{\Delta}(\ast O)$ be a sesqui-linear pairing,
i.e., $\bigl(\nbigd_{\Delta}(\ast O)\otimes
 \nbigd_{\Delta^{\dagger}}(\ast O)\bigr)$-homomorphism.
Let $S$ be a small sector in $\Delta^{\ast}$,
and let $\Sbar$ denote the closure of $S$
in the real blow up of $\Delta$ along $O$.
We take a $\nabla_i$-flat splitting of 
the full Stokes filtration $\nbigf^S$:
\begin{equation}
 \label{eq;07.9.10.1}
 \nbigv_{i|\Sbar}
=\bigoplus_{\gminia\in\Irr(\nabla_i)}
 \nbigv_{i,\gminia,S}
\end{equation}
Then, we obtain an induced
$\nbigd_S\otimes\nbigd_{S^{\dagger}}$-homomorphism
$C_{\gminia,\gminib}:
 \nbigv_{1,\gminia,S}\otimes
 \overline{\nbigv}_{2,\gminib,S}
\lrarr \distribution_S$ on $S$.

\begin{lem}
\label{lem;07.7.11.11}
We have $C_{\gminia,\gminib}=0$
unless $\Re(\gminia+\gminib)\leq 0$ on $S$.
\end{lem}
\pf
Let $\Ohat$ denote the completion of $\Delta$
along $O$.
We take a meromorphic frame of 
$\vecvhat^{(i)}=\bigl(\vecvhat^{(i)}_{\gminia}\,\big|\,
 \gminia\in\Irr(\nabla_i)\bigr)$
of $\nbigv_{i|\Ohat}$ which is compatible with
the irregular decomposition:
\[
 \nabla_i\vecvhat^{(i)}=\vecvhat^{(i)}
 \left(
 \bigoplus_{\gminia\in\Irr(\nabla_i)}
 \Bigl(d\gminia+A_{\gminia}^{(i)}\frac{dz}{z}\Bigr)
 \right)
\]
Here, $A_{\gminia}^{(i)}$ are constant matrices 
with Jordan forms.
For each small sector $S'$ in $\Delta^{\ast}$,
let $\vecv^{(i)}_{S'}$ be the lift of $\vecvhat^{(i)}$
to $\nbigv^{(i)}_{|S'}$,
which is compatible
with the decomposition (\ref{eq;07.9.10.1}).
By varying $S'$ and gluing them in $C^{\infty}$
as in Section \ref{subsection;07.6.16.8} below,
we obtain a $C^{\infty}$-frame 
$\vecv^{(i)}_{C^{\infty}}=\bigl(
 \vecv^{(i)}_{\gminia,C^{\infty}}
 \bigr)$.
We may assume that
$C\bigl(v^{(1)}_{\gminia,p,C^{\infty}},\,
 v^{(2)}_{\gminib,q,C^{\infty}}\bigr)$
also give distributions on $\Delta$.
On the sector $S$,
we may also assume that
$v^{(i)}_{\gminia,p,C^{\infty}|S}$
are holomorphic and compatible with the decompositions
(\ref{eq;07.9.10.1}).
Let us consider the distribution
$F_{p,q}:=C\bigl(v^{(1)}_{\gminia,p,C^{\infty}},
 v^{(2)}_{\gminib,q,C^{\infty}}\bigr)$
on $\Delta$.

First, we consider the case in which
 $v^{(1)}_{\gminia,p,C^{\infty}}$
and $v^{(2)}_{\gminib,q,C^{\infty}}$
are contained in the  the bottom part of
the weight filtration $W$
associated to the nilpotent parts of 
$A_{\gminia}^{(1)}$
and $A_{\gminib}^{(2)}$,
respectively.
Then, we may assume to have 
the following equality on $S$:
\[
 \nabla_1v^{(1)}_{\gminia,p,C^{\infty}|S}
=\Bigl(d\gminia+\alpha_p\frac{dz}{z}\Bigr)
 v^{(1)}_{\gminia,p,C^{\infty}|S},
\quad
 \nabla_2v^{(2)}_{\gminia,q,C^{\infty}|S}
=\Bigl(d\gminib+\beta_q\frac{dz}{z}\Bigr)
 v^{(2)}_{\gminib,q,C^{\infty}|S}
\]
We have
$ F_{p,q|S}=B\cdot
 \exp\bigl(\gminia+\alpha_p\log z
 +\gminibbar+\betabar_q\log \zbar\bigr)$
for some constant $B$.
Due to Lemma \ref{lem;07.7.11.10},
we have $B=0$,
i.e.,
$C_{\gminia,\gminib}(
 v^{(1)}_{\gminia,p,C^{\infty}|S}, 
 v^{(2)}_{\gminib,q,C^{\infty}|S})=0$.

We show $F_{p,q|S}=0$ in general
by using  an induction on
the degree with respect to the filtration $W$.
We put 
$d^{(1)}(p):=
 \deg^W(v^{(1)}_{\gminia,p,C^{\infty}|S})$
and $d^{(2)}(q):=
 \deg^W(v^{(2)}_{\gminib,q,C^{\infty}|S})$.
We will show $F_{p,q|S}=0$
assuming
$F_{p',q'|S}=0$ for any 
$d^{(1)}(p')\leq d^{(1)}(p)$
and $d^{(2)}(q')< d^{(2)}(q)$.
In that case,
$F_{p,q|S}$ satisfies the differential equation:
\[
 dF_{p,q|S}= 
 \Bigl(
 d\gminia+d\gminibbar
+\alpha_p\frac{dz}{z}
+\beta_q\frac{d\zbar}{\zbar}
\Bigr)\cdot F_{p,q|S}
\]
Hence, we obtain
$F_{p,q|S}=B\cdot 
 \exp\bigl(\gminia+\gminibbar
 +\alpha_p\log z+\betabar_q\log\zbar\bigr)$
for some constant $B$.
And, by Lemma \ref{lem;07.7.11.10},
we have $B=0$.
Similarly,
we can show $F_{p,q|S}=0$
assuming
$F_{p',q'|S}=0$ for any 
$d^{(1)}(p')< d^{(1)}(p)$
and $d^{(2)}(q')\leq d^{(2)}(q)$.
Thus, we obtain $C_{\gminia,\gminib}=0$.
\hfill\qed

\section{Specialization of the associated $\nbigr$-triples}
\label{subsection;08.9.28.10}
\subsection{Reduction and specialization}

Let $X:=\Delta^n$,
$D_i:=\{z_i=0\}$
and $D=\bigcup_{i=1}^nD_i$.
Let $\harmonicbundle$ be 
an unramifiedly $\nbiga$-good 
wild harmonic bundle
on $X-D$.
According to Proposition \ref{prop;08.9.27.3},
we have the $\nbigr_X$-triple
$\gbigt(E):=(\gbige,\gbige,\gbigc)$.
Let us study its specialization 
along a monomial function $g=\vecz^{\vecp}$,
where $\vecp=(p_i)\in\seisuu_{\geq 0}^n$.
We use the notation in Section \ref{subsection;08.2.2.1}.
According to Proposition \ref{prop;08.1.12.60},
we have the decomposition:
\begin{equation}
\label{eq;08.9.27.10}
 P\Gr^{W(N)}_p
 \psitilde_{g,u}(\gbige)
=\bigoplus_I\nbigm_{u,p,I}
\end{equation}
\index{sheaf $\nbigm_{u,p,I}$}
Here, 
the support of $\nbigm_{u,p,I}$ is $\nbigd_I$,
and the support of
any non-trivial $\nbigr$-submodule of $\nbigm_{u,p,I}$
is not contained in $\nbigd_J$ ($J\supset I$).
Since $\nbigm_{u,p,I}$ is strictly $S$-decomposable
along any $z_i$ $(i\in I)$,
it is the push-forward of
the $\nbigr_{D_I}$-modules $\nbigm'_{u,p,I}$.
We have the sesqui-linear pairing
\[
 C_{u,p,I}:
 \nbigm_{u,p,I|\vecS\times X}\otimes
 \sigma^{\ast}\nbigm_{u,p,I|\vecS\times X}
\lrarr
 \distribution_{\vecS\times X/\vecS},
\]
induced by
$\psitilde_{g,u}(\gbigc)\circ 
\bigl((\sqrt{-1}N)^p\otimes 1\bigr)$.
\index{pairing $C_{u,p,I}$}
We would like to know a more detailed property
of the $\nbigr$-triple
\[
 \nbigt_{u,p,I}(E):=
 (\nbigm_{u,p,I},
 \nbigm_{u,p,I},C_{u,p,I}),
\]
with the Hermitian sesqui-linear duality
$\nbigs_{u,p,I}=\bigl((-1)^p,1\bigr)$
of weight $p$.
\index{$\nbigr$-triple $\nbigt_{u,p,I}(E)$}

For a subset $J\subset\nbar$,
let $\Irr(\theta,J)$ denote the image of
$\Irr(\theta)$ by the natural map
\[
 M(X,D)/H(X)\lrarr 
 M(X,D)/M(X,D(\nbar-J)),
\]
where $D(\nbar-J):=\bigcup_{i\in\nbar-J}D_i$.
We have the unramifiedly 
$\nbiga$-good wild harmonic bundles
$(E_{J,\gminia},\delbar_{J,\gminia},
 \theta_{J,\gminia},h_{J,\gminia})$
for $\gminia\in\Irr(\theta,J)$ on $X-D$,
which is obtained as the reduction 
of $(E,\delbar_E,\theta,h)$.
\index{harmonic bundle
 $(E_{J,\gminia},\delbar_{J,\gminia},
 \theta_{J,\gminia},h_{J,\gminia})$}
(See Section \ref{subsection;07.11.24.50}.)
We have the associated $\nbigr_X$-triple
$(\gbige_{J,\gminia},\gbige_{J,\gminia},
 \gbigc_{J,\gminia})$ on $X$.
\index{$\nbigr$-triple
 $(\gbige_{J,\gminia},\gbige_{J,\gminia},
 \gbigc_{J,\gminia})$}
We have the decomposition
as in (\ref{eq;08.9.27.10}):
\[
 P\Gr^{W(N)}_p\psitilde_{g,u}(\gbige_{J,\gminia})
=\bigoplus_I
 \nbigm_{u,p,I}(E_{J,\gminia})
\]
\index{sheaf $\nbigm_{u,p,I}(E_{J,\gminia})$}
We have the induced sesqui-linear pairing
\[
 C_{u,p,I}(E_{J,\gminia}):
 \nbigm_{u,p,I}
 (E_{J,\gminia})_{|\vecS\times X}
\otimes
 \sigma^{\ast}\nbigm_{u,p,I}
 (E_{J,\gminia})_{|\vecS\times X}
\lrarr
 \distribution_{\vecS\times X/\vecS}.
\]
\index{pairing $C_{u,p,I}(E_{J,\gminia})$}
Thus, we obtain the $\nbigr_X$-triples
$\nbigt_{u,p,I}(E_{J,\gminia})$.
\index{$\nbigr$-triple
 $\nbigt_{u,p,I}(E_{J,\gminia})$}

\begin{prop}
\label{prop;07.9.8.5}
We have the natural isomorphism
$\nbigt_{u,p,I}(E)\simeq
 \nbigt_{u,p,I}(E_{I,0})$.
\end{prop}
\pf
We may assume  $D_I\subset g^{-1}(0)$.
Let $\iota_I:\Dhat_I\lrarr X$ be the natural morphism.
We have the following natural isomorphisms:
\begin{multline}
  \iota_I^{\ast}P\Gr^{W(N)}_p
 \psitilde_{g,u}(\gbige)
\simeq
 P\Gr^{W(N)}_p
 \psitilde_{g,u}(\iota_I^{\ast}\gbige)
\simeq\bigoplus_{\gminia\in\Irr(\theta,I)}
 P\Gr^{W(N)}_p
  \psitilde_{g,u}(\iota_I^{\ast}\gbige_{I,\gminia}) \\
\simeq
 P\Gr^{W(N)}_p\psitilde_{g,u}(\iota_I^{\ast}\gbige_{I,0})
\simeq
\iota_I^{\ast}
 P\Gr^{W(N)}_p\psitilde_{g,u}(\gbige_{I,0})
\end{multline}
Since the support of $\nbigm_{u,p,I}$ 
and $\nbigm_{u,p,I}(E_{I,0})$ are contained in
$\nbigd_I$,
we obtain the natural isomorphism
$\nbigm_{u,p,I}\simeq 
 \nbigm_{u,p,I}(E_{I,0})$.
In the following, we may and will 
implicitly identify them.
Let us compare the sesqui-linear pairings.
For simplicity of description, we put
$C'_{u,p,I}:=C_{u,p,I}(E_{I,0})$.
Because of the strict $S$-decomposability
of $\nbigm_{u,p,I}$ and $\nbigm_{u,p,I}(E_{I,0})$
along any $z_j$
(Proposition \ref{prop;08.1.12.60}),
we have only to compare $C_{u,p,I}$ 
and $C_{u,p,I}'$
on $X-\bigcup_{j\not\in I}D_j$.
The restrictions of $\nbigm_{u,p,I}$
and $\nbigm_{u,p,I}(E_{I,0})$
come from the family of $\lambda$-flat bundles
$\nbigm'_{u,p,I}$
and $\nbigm'_{u,p,I}(E_{I,0})$
on $\cnum_{\lambda}\times D_I^{\circ}$,
where 
$D_I^{\circ}:=D_I\setminus \bigcup_{I\subsetneq J}D_J$.
We have only to compare 
the restrictions of $C_{u,p,I}$ and $C_{u,p,I}'$
to $\nbigm'_{u,p,I}$
and $\nbigm'_{u,p,I}(E_{I,0})$.
(See Proposition \ref{prop;a11.23.10}.)
Let $\pi:X\lrarr D_I$ be the natural projection.
By considering the restriction to
$\pi^{-1}(P)$ for $P\in D_I^{\circ}$,
we may and will assume $I=\nbar$.

\vspace{.1in}
We take an auxiliary sequence
$\vecm(0),\ldots,\vecm(L)$ for $\Irr(\theta)$.
We have the reductions
$\bigl(E^{\vecm(i)}_{\gminia},
 \delbar^{\vecm(i)}_{\gminia},
 \theta^{\vecm(i)}_{\gminia}
 h^{\vecm(i)}_{\gminia}
 \bigr)$
for $\gminia\in\Irrbar(\theta,\vecm(i))$.
We have the associated $\nbigr$-triples
$\bigl(\gbige^{\vecm(i)}_{\gminia},
 \gbige^{\vecm(i)}_{\gminia},
 \gbigc^{\vecm(i)}_{\gminia}\bigr)$ on $X$.
We have the decomposition:
\[
 P\Gr^{W(N)}_p\psitilde_{g,u}
 \bigl(\gbige^{\vecm(i)}_{\gminia}\bigr)
=\bigoplus_J\nbigm_{u,p,J}(E^{\vecm(i)}_{\gminia})
\]
We have the induced sesqui-linear pairings
$C_{u,p,J}(E^{\vecm(i)}_{\gminia})$.
We have the natural isomorphisms
$\nbigm_{u,p,\nbar}(E^{\vecm(i)}_0)
\simeq
 \nbigm_{u,p,\nbar}(E^{\vecm(i+1)}_0)$,
as remarked in Corollary \ref{cor;08.2.2.11}.
We have only to show
$C_{u,p,\nbar}(E^{\vecm(i)}_{0})
=C_{u,p,\nbar}(E^{\vecm(i+1)}_0)$
for $i=-1,\ldots,L$
under the isomorphisms.
(We put $E^{\vecm(-1)}_0:=E$, formally.)
By an easy inductive argument,
we have only to compare
$C_{u,p,\nbar}$
and $C_{u,p,\nbar}(E^{\vecm(0)}_0)$.
To make a description simpler,
we put
$C^{\vecm(0)}_{u,p,\nbar}:=
 C_{u,p,\nbar}(E^{\vecm(0)}_0)$.

\vspace{.1in}

Take a generic $\lambda_0\in\vecS$.
We take a small neighbourhood
$U_1$ of $\lambda_0$,
and we put $U_2:=\sigma(U_1)$
which is a neighbourhood of $-\lambda_0$.
We put $U_i\cap \vecS=\vecI_i$ $(i=1,2)$.
Take 
$f_i\in \nbigm'_{u,p,\nbar|U_i\times D_{\nbar}}
=\nbigm'_{u,p,\nbar}
 (E^{\vecm(0)}_0)_{|U_i\times D_{\nbar}}$.
We would like to compare the distributions
$C_{u,p,\nbar}(f_1,\sigma^{\ast}f_2)$
and $C^{\vecm(0)}_{u,p,\nbar}(f_1,\sigma^{\ast}f_2)$.
We may assume 
$\paramap(\lambda_0,u)<0$
and $\paramap(-\lambda_0,u)<0$.

We take a lift
$\ftilde_1\in \bigl(
 i_{g\dagger}\gbige^{\vecm(0)}_0\bigr)_{
|U_1\times (X\times\cnum_t)}$
of $f_1$ as in Section \ref{subsection;08.2.2.5},
i.e.,
(i) 
$\ftilde_1\in
\Uzero_{\paramap(\lambda_0,u)}
 (i_{g\dagger}\gbige^{\vecm(0)}_{0})$,
(ii)
the induced section $\ftilde_1^{(1)}$ of
$\Gr^{\Uzero}_{\paramap(\lambda_0,u)}\bigl(
 i_{g\dagger}\gbige^{\vecm(0)}_{0}\bigr)$
is contained in 
$W_p(N)\psizero_{g,u}\gbige^{\vecm(0)}_{0}$,
(iii)
the induced section $\ftilde_1^{(2)}$ of
$\Gr^{W(N)}_p\psizero_{g,u}\gbige^{\vecm(0)}_{0}$
is equal to $f_1$.
Similarly, we take a lift
$\ftilde_2\in\bigl(
 i_{g\dagger}\gbige^{\vecm(0)}_0\bigr)_{|
 U_2\times (X\times\cnum_t)}$
of $f_2$.
Let $\chi$ be a test function
on $\cnum_t$ which is constantly $1$ around $t=0$.
Let $\omega_0:=(2\pi)^{-1}\sqrt{-1}dt\cdot d\tbar$.
By definition, we have the following:
\begin{multline}
\label{eq;07.9.8.3}
\bigl\langle
 C^{\vecm(0)}_{u,p,\nbar}
 \bigl(f_1,\sigma^{\ast}(f_2)\bigr),\,
 \rho
\bigr\rangle
=
\bigl\langle
 \psitilde_{g,u}\gbigc^{\vecm(0)}_{0}
 \bigl((\sqrt{-1}N)^p\ftilde^{(1)}_1,\,
 \sigma^{\ast}\ftilde^{(1)}_2\bigr),
 \rho
\bigr\rangle\\
=(\sqrt{-1})^p
 \underset{s+\eigenmap(\lambda,u)}{\Res}\,
 \Bigl\langle
 \gbigc^{\vecm(0)}_{0}\bigl(
(-\deldel_tt+\eigenmap(\lambda,u))^p 
 \ftilde_1,
 \sigma^{\ast}\ftilde_2
 \bigr),\,
 |t|^{2s}\!\cdot\!
 \chi\!\cdot\!\rho\!\cdot\!
 \omega_0
 \Bigr\rangle
\end{multline}

We take a lift $\ftilde_{1,C^{\infty}}$
of $\ftilde_1$ to 
$\Uzero\bigl(i_{g\dagger}\gbige\bigr)
\otimes C^{\infty}$,
as in Section \ref{subsection;08.2.2.5}.
Similarly,
we take a lift $\ftilde_{2,C^{\infty}}$
of $\ftilde_2$ to 
$U^{(-\lambda_0)}\bigl(i_{g\dagger}\gbige\bigr)
\otimes C^{\infty}$.
By construction, 
$\ftilde_{b,C^{\infty}}$ $(b=1,2)$
are contained in 
$\nbigf^{S\vecm(0)}_0$ 
for each small sector $S$.
By using Corollary \ref{cor;08.9.27.20},
we obtain the following equality
on $\vecI_1\times (X-D)$:
\begin{equation}
\label{eq;07.9.8.1}
 \gbigc\bigl(
 (-\deldel_tt+\eigenmap(\lambda,u))^p 
 \ftilde_{1,C^{\infty}},\,
  \sigma^{\ast}\ftilde_{2,C^{\infty}}\bigr)
=\gbigc^{\vecm(0)}_{0}
 \bigl((-\deldel_tt+\eigenmap(\lambda,u))^p 
 \ftilde_1,\,
 \sigma^{\ast}\ftilde_2\bigr)
\end{equation}
Due to Lemma \ref{lem;08.2.2.8},
the following holds:
\begin{multline}
\label{eq;07.9.8.2}
(\sqrt{-1})^p
\underset{s+\eigenmap(\lambda,u)}\Res
 \Bigl\langle
 \gbigc\bigl(
(-\deldel_tt+\eigenmap(\lambda,u))^p 
 \ftilde_{1,C^{\infty}},\,
 \sigma^{\ast}\ftilde_{2,C^{\infty}}\bigr),\,
 |t|^{2s}\!\cdot\!
 \chi\!\cdot\! \rho\!\cdot\!\omega_0
 \Bigr\rangle\\
=\bigl\langle
 \psi_{g,u}(\gbigc)\bigl(
 (\sqrt{-1}N)^p\ftilde^{(1)}_{1,C^{\infty}},\,
 \sigma^{\ast}\ftilde^{(1)}_{2,C^{\infty}}\bigr),\,\,
 \rho 
  \bigr\rangle 
=C_{u,p,\nbar}(f_1,\sigma^{\ast}f_2)
\end{multline}
From (\ref{eq;07.9.8.3}),
(\ref{eq;07.9.8.1})
and (\ref{eq;07.9.8.2}),
we obtain 
$C_{u,p,\nbar}=C^{\vecm(0)}_{u,p,\nbar}$
around generic $\lambda_0$.
By continuity,
we obtain $C_{u,p,\nbar}=C^{\vecm(0)}_{u,p,\nbar}$.
Thus, the proof of Proposition \ref{prop;07.9.8.5}
is finished.
\hfill\qed

\begin{cor}
\label{cor;07.10.19.35}
The restriction of $\nbigt_{u,p,I}$ to
$X\setminus \bigcup_{j\not\in I}D_j$
comes from an unramifiedly 
$\nbiga$-good wild harmonic bundle
$\bigl(E_{u,p,I},\delbar_{u,p,I},
 \theta_{u,p,I},h_{u,p,I}\bigr)$
on $D_I^{\circ}$
with the twist by $\Tate^S(-p/2)$.
The set of the irregular values of $\theta_{u,p,I}$
is contained in
$\bigl\{\gminia\in\Irr(\theta)\,\big|\,
 \gminia=0\,\,\mbox{\rm in }
 \Irr(\theta,I)
 \bigr\}$.
\end{cor}
\pf
Note that 
$(E_{I,0},\delbar_{E,0},\theta_{I,0},h_{I,0})$
is tame with respect to the divisors
$D_l$ $(l\in I)$.
Hence, 
$\nbigt_{u,p,I}(E_{I,0})_{|D_I^{\circ}}$
comes from the harmonic bundle
with the twist by $\Tate^S(-p/2)$.
(See Section 18.4.9 of \cite{mochi2}.)
By considering the restriction to $\{0\}\times X$,
we can easily check the other claims.
\hfill\qed

\subsection{The specializations come from
 wild harmonic bundles}

Let $\gbigt(E_{u,p,I})$ denote
the $\nbigr_X$-triples
$(\gbige_{u,p,I},\gbige_{u,p,I},\gbigc_{u,p,I})$
associated to 
the reductions $\bigl(E_{u,p,I},\delbar_{u,p,I},
  \theta_{u,p,I},h_{u,p,I}\bigr)$
(Corollary \ref{cor;07.10.19.35})
with the natural Hermitian sesqui-linear duality.
\index{$\nbigr$-triple
 $\gbigt(E_{u,p,I}):=
(\gbige_{u,p,I},\gbige_{u,p,I},\gbigc_{u,p,I})$}

\begin{prop}
\label{prop;07.9.8.26}
There exists an isomorphism
$\gbigt(E_{u,p,I})
\simeq
 \nbigt_{u,p,I}\otimes\Tate^{S}(p/2)$.
\end{prop}
\pf
Their restrictions to
$X\setminus \bigcup_{j\not\in I}D_j$ are the same.
We have only to show
$\gbige_{u,p,I}\simeq\nbigm_{u,p,I}$.
(Note Proposition \ref{prop;a11.23.10}.)
Since both of them are strictly $S$-decomposable
along $G:=\prod_{j\not\in I}z_j$,
we have only to compare their meromorphic
structures,
i.e.,
we have only to show
$\gbige_{u,p,I}\otimes\nbigo(\ast G)
=\nbigm_{u,p,I}\otimes\nbigo(\ast G)$.
(See Lemma \ref{lem;07.11.2.15}
 below.)
Note that they come from
families of meromorphic $\lambda$-flat bundles
on $\cnum_{\lambda}\times D_I$,
which we have to compare.
For that purpose,
we have only to compare their restrictions
to 
$\pi_j^{-1}(P)$ in $D_I$
for $j\not\in I$ and 
$P\in (D_j\cap D_I)\setminus
 \bigcup_{k\not\in \{j\}\cup I}
 \bigl(D_j\cap D_k\cap D_I\bigr)$,
where $\pi_j:D_I\lrarr D_j\cap D_I$ 
denotes the projection.
Therefore, we have only to consider the case
$\dim D_I=1$.
Because of Proposition \ref{prop;07.9.8.5},
we may assume $I=\{2,\ldots,n\}$
and $\Irr(\theta)\subset\cnum[z_1^{-1}]$.
We have the $\nbigr$-triples
$\Gr_{\gminia}\bigl(\nbigt_{u,p,I}
 \otimes\Tate^S(p/2)\bigr)^{\circ}$
for $\gminia\in\Irr(\theta)$ on $D_I^{\circ}$
obtained from 
$\nbigt_{u,p,I}\otimes\Tate^S(p/2)$
as the reductions 
(Section \ref{subsection;07.12.20.7}).
Because of the uniqueness
(Theorem \ref{thm;08.9.16.100}),
we have only to show that
$\Gr_{\gminia}\bigl(\nbigt_{u,p,I}\otimes
 \Tate^S(p/2)\bigr)^{\circ}$
are variations of 
polarized pure twistor structures of weight $0$.
For that purpose,
let us show that
$\Gr_{\gminia}\bigl(
 \nbigt_{u,p,I}\bigr)^{\circ}$
are isomorphic to the restriction of
$\nbigt_{u,p,I}(E_{1,\gminia})$
to $D_I^{\circ}$.

We have the isomorphism:
\[
 P\Gr^{W(N)}_p\psitilde_{g,u}
 (\gbige(\ast \nbigd_1))
 _{|\nbigdhat_1}
\simeq
 \bigoplus_{\gminia\in\Irr(\theta)}
 P\Gr^{W(N)}_p\psitilde_{g,u}
 (\gbige_{1,\gminia}(\ast \nbigd_1))
 _{|\nbigdhat_{1}}
\]
Thus, we have the isomorphism:
\begin{equation}
 \label{eq;08.1.12.1}
 \nbigm_{u,p,I}(\ast \nbigd_1)
 _{|\nbigdhat_{1}}
\simeq
 \bigoplus_{\gminia\in\Irr(\theta)}
 \nbigm_{u,p,I}(E_{1,\gminia})
 (\ast \nbigd_1)
 _{|\nbigdhat_{1}}
\end{equation}
Note that
$\nbigm_{u,p,I}(\ast \nbigd_1)$
and $\nbigm_{u,p,I}(E_{1,\gminia})(\ast \nbigd_1)$
come from families of
meromorphic $\lambda$-flat bundles
on $\nbigd_{I}$,
which are denoted by
$\nbigv$ and $\nbigv_{\gminia}$.
We regard them as
$\nbigo_{\nbigx}(\ast \nbigd_1)$-modules.
By (\ref{eq;08.1.12.1}),
we have the isomorphism
\begin{equation}
\label{eq;08.9.27.100}
\Psihat:
 \bigoplus
 \nbigv_{\gminia|\nbigdhat_{1}}
\simeq
 \bigoplus_{\gminia}\Gr_{\gminia}(\nbigv)_{|\nbigdhat_1}
=
\nbigv_{|\nbigdhat_{1}}
\end{equation}

\begin{lem}
\label{lem;08.9.27.101}
The reduction
$\Gr_{\gminia}(\nbigv)$
is naturally isomorphic to
$\nbigv_{\gminia}$
for each $\gminia\in \Irr(\theta)$.
\end{lem}
\pf
By considering the tensor product 
with rank one harmonic bundle,
we may assume that $\gminia=0$.
We are given the isomorphism:
\[
 \Phihat:
 \Gr_0(\nbigv)_{|\widehat{z_1=0}}
\simeq
 \nbigv_{0|\widehat{z_1=0}}
\]
Let $\lambda\neq 0$.
The restriction of $\Phihat$
to $\{\lambda=0\}$ is denoted by
$\Phihat^{\lambda}$.
Since the restrictions
$\Gr_0\nbigv_{|\{\lambda\}\times D_I}$
and $\nbigv_{0|\{\lambda\}\times D_I}$
are $\lambda$-flat meromorphic bundles
with regular singularity,
$\Phihat^{\lambda}$ is convergent,
and 
it gives the holomorphic isomorphism
$\Gr_{0}\nbigv_{|\{\lambda\}\times D_I}
\simeq
 \nbigv_{0|\{\lambda\}\times D_I}$.
Then, it is easy to see that
$\Phihat_{|\{\lambda\neq 0\}}$
comes from the flat isomorphism:
\begin{equation}
 \label{eq;08.1.12.5}
 \Gr_0\nbigv
 _{|\cnum_{\lambda}^{\ast}\times D_I}
\lrarr
 \nbigv_{0|\cnum_{\lambda}^{\ast}\times D_I}
\end{equation}
Let us show that
the isomorphism (\ref{eq;08.1.12.5})
is extended on $\cnum_{\lambda}\times D_I$.
We take local frames
$\vecv$ and $\vecw$  of
$\Gr_0\nbigv$ and $\nbigv_{0}$
on some neighbourhood $\nbigu$
around $(0,O)$ in $\nbigd_I$.
The morphism (\ref{eq;08.1.12.5})
is expressed by the matrix $A$
with respect to $\vecv$ and $\vecw$.
We may assume that 
the entries $a_{i,j}$ of $A$ are holomorphic 
on $\nbigu\setminus \{\lambda=0\}$.
Moreover,
$a_{i,j}$ has the 
formal expansion
$\sum a_{i,j,k}(\lambda)\cdot z_1^k$,
where $a_{i,j,k}(\lambda)$ are holomorphic
with respect to $\lambda$.
(Recall that we are given $\Phihat$
on $\cnum_{\lambda}\times\widehat{\{z_1=0\}}$.)
The power series are absolute convergent
on $\{|z_1|\leq R_1\}\times
 \{\delta_1\leq |\lambda|\leq \delta_2\}$.
Hence, we have the following
for any $l\geq 0$ and $|z_1|\leq R_1$:
\[
 \int_{|\lambda|=\delta_2}
 a_{i,j}\cdot \lambda^l\cdot d\lambda
=\sum_k\int_{|\lambda|=\delta_2}
 a_{i,j,k}(\lambda)\cdot \lambda^l
 \cdot z^k\cdot d\lambda
=0
\]
Then, we can conclude that
$a_{i,j}$ are holomorphic on
$\{|z_1|\leq R_1\}\times
 \{|\lambda|\leq \delta_2\}$,
which implies that
the morphism (\ref{eq;08.1.12.5})
is extended on 
$\cnum_{\lambda}\times D_I$.
Since the completion along $\{z_1=0\}$
is isomorphic,
it is isomorphic if we shrink $X$ appropriately.
\hfill\qed

\vspace{.1in}

Take a generic $\lambda_0$,
and its small neighbourhood $U(\lambda_0)$.
We have another direct construction
of the flat isomorphism 
$\Gr_{\gminia}
 \nbigv_{|U(\lambda_0)\times X}
\simeq
 \nbigv_{\gminia|U(\lambda_0)\times X}$.
Let $b=\paramap(\lambda_0,u)$.
Let $S$ be a small sector in 
$U(\lambda_0)\times (X-D_1)$.
Because we have assumed that $\lambda_0$ is generic,
we can take a $\DD$-flat splitting
of the full Stokes filtration
(Proposition \ref{prop;07.10.2.40}):
\begin{equation}
 \label{eq;08.9.27.110}
 \lefttop{\nbar}\Vzero_{
 b\vecp-\epsilon\vecdelta_{\nbar-\supp(\vecp)}}
 (\nbigq\nbige)_{|\Sbar}
=\bigoplus_{\gminia\in\Irr(\theta)}
 \lefttop{\nbar}\Vzero_{
 b\vecp-\epsilon\vecdelta_{\nbar-\supp(\vecp)}}
 (\nbigq\nbige)_{\gminia,\Sbar}
\end{equation}
It induces the flat decompositions
$\nbigq\nbige_{|S}
=\bigoplus_{\gminia} \nbigq\nbige_{\gminia,S}$
and
\begin{equation}
 \label{eq;07.9.8.16}
 \Uzero_b
 \bigl(i_{g\dagger}\gbige(\ast \nbigd_1)\bigr)_{|
 S\times\cnum_t}
=\bigoplus_{\gminia\in\Irr(\theta)}
 U^{(\lambda_0)}_b
 \bigl(i_{g\dagger}\gbige(\ast \nbigd_1)\bigr)_{\gminia,S}
\end{equation}
We have the corresponding decomposition
of $\Gr^{\Uzero}$.
Since the decomposition is compatible
with the action of $-\deldel_tt$,
we also obtain the decompositions:
\begin{equation}
 \label{eq;07.9.8.15}
 \psizero_{g,u}
 \bigl(\gbige(\ast \nbigd_1)\bigr)_{|S}
\!\!=\!\!\!\!\bigoplus_{\gminia\in\Irr(\theta)}
 \!\!\psizero_{g,u}
 \bigl(\gbige(\ast \nbigd_1)\bigr)_{\gminia,S}
\end{equation}
\begin{equation}
\label{eq;07.11.29.15}
 P\Gr^{W(N)}_p\psizero_{g,u}
 \bigl(\gbige(\ast \nbigd_1)\bigr)_{|S}
=\!\!\!\!\bigoplus_{\gminia\in\Irr(\theta)}
 \!\!P\Gr^{W(N)}_p
 \psizero_{g,u}
 \bigl(\gbige(\ast \nbigd_1)\bigr)_{\gminia,S}
\end{equation}
\begin{equation}
\label{eq;08.9.27.50}
\nbigv_{|S}
=\bigoplus_{\gminia\in\Irr(\theta)}
 \nbigv_{\gminia,S}
\end{equation}
By the order $\leq_S$ on the set $\Irr(\theta)$
associated to the sector $S$,
we obtain the filtrations $\nbigftilde^S$:
\[
 \nbigftilde^S_{\gminia}
 \Uzero_b
 \bigl(i_{g\dagger}\gbige(\ast \nbigd_1)\bigr)_{|S}
=\bigoplus_{\gminib\leq_S\gminia}
 U^{(\lambda_0)}_b
 \bigl(i_{g\dagger}\gbige(\ast \nbigd_1)\bigr)_{\gminib,S}
\]
\[
  \nbigftilde^S_{\gminia}
  \psizero_{g,u}
 \bigl(\gbige(\ast \nbigd_1)\bigr)_{|S}
=\bigoplus_{\gminib\leq_S\gminia}
 \psizero_{g,u}
 \bigl(\gbige(\ast \nbigd_1)\bigr)_{\gminib,S}
\]
\[
 \nbigftilde^S_{\gminia}
 P\Gr^{W(N)}_p\psizero_{g,u}
 \bigl(\gbige(\ast \nbigd_1)\bigr)_{|S}
=\bigoplus_{\gminib\leq_S\gminia}
 P\Gr^{W(N)}_p
 \psizero_{g,u}
 \bigl(\gbige(\ast \nbigd_1)\bigr)_{\gminib,S}
\]
\[
 \nbigftilde^S_{\gminia}
 \nbigv_{|S}
=\bigoplus _{\gminib\leq_S\gminia}
 \nbigv_{\gminib,S}
\]
Although the decompositions (\ref{eq;07.9.8.16}),
(\ref{eq;07.9.8.15}), (\ref{eq;07.11.29.15})
and (\ref{eq;08.9.27.50})
depend on the choice of
a splitting of $\nbigq\nbige$,
the filtrations $\nbigftilde^S$ 
are well defined.

By construction,
we have the natural isomorphism
$F_{\gminia,S}:
 \nbigv_{\gminia|S}\simeq
\nbigv_{\gminia,S}$.

\begin{lem}
If we shrink $S$,
$F_S:=\bigoplus F_{\gminia,S}$
is extended
to the isomorphism
$\bigoplus \nbigv_{\gminia|\Sbar}
\simeq
 \nbigv_{|\Sbar}$,
and $F_{S|\Zhat}=\pi^{-1}\circ\Psihat$.
The filtration
$\nbigftilde^{S}(\nbigv_{|S})$
is equal to the full Stokes filtration
of $\nbigv_{|S}$.
\end{lem}
\pf
Let $f$ be a section of
$\nbigv_{\gminia}$.
We take a lift $\ftilde$ to
$\Uzero_b\bigl(
 i_{g\dagger}\gbige_{1,\gminia}\bigr)$.
Namely,
the induced section
$\ftilde^{(1)}$ of
$\Gr^{\Uzero}_{b}
 (i_{g\dagger}\gbige_{1,\gminia})$
is contained in
$W_p(N)\psitilde_{u,g}(\gbige_{1,\gminia})$,
and it induces
$f\in \nbigv_{\gminia}
\subset
 \Gr^{W(N)}_p\psitilde_{u,g}
 (\gbige_{1,\gminia})$.
Note that
$\Uzero_b(i_{g\dagger}
 \gbige_{1,\gminia})(\ast D_1)$ $(b<0)$
is the $V_0\nbigr_{X\times\cnum_t}(\ast D_1)$-submodule
of $i_{g\dagger}\nbigqzero\nbige_{1,\gminia}$
generated by 
$\lefttop{\nbar}\Vzero_{
 b\vecp-\epsilon\vecdelta_{\nbar-\supp(\vecp)}}
 (\nbigq\nbige_{1,\gminia})\otimes 1$
over $\nbigr_X(\ast \nbigd_1)$.
Hence, we have the expression
\[
 \ftilde=
 \sum P_l\cdot (a_l\otimes 1),
\]
where $P_l\in \nbigr_X(\ast D_1)$
and $a_l\in 
 \lefttop{\nbar}\Vzero_{b\vecp
 -\epsilon\vecdelta_{\nbar-\supp(\vecp)}}
 (\nbigq\nbige_{1,\gminia})$.

Let $a_{l,S}$ be the lift of $a_{l}$ to
$ \lefttop{\nbar}\Vzero_{b\vecp
 -\epsilon\vecdelta_{\nbar-\supp(\vecp)}}
 (\nbigq\nbige)_{\gminia,\Sbar}$.
We obtain the section $\ftilde_S$
of $U^{(\lambda_0)}_b
 \bigl(i_{g\dagger}\gbige(\ast \nbigd_1)
 \bigr)_{\gminia,S}$
given as follows:
\[
\ftilde_S:=
 \sum P_l\cdot (a_{l,S}\otimes 1)
\]
It induces a section $f_S$ of $\nbigv_{\gminia,S}$,
which is equal to $F_{\gminia,S}(f_{|S})$.
To show the first claim of the lemma,
we have only to check that
$f_S$ induces the section of $\nbigv_{|\Sbar}$,
and that $f_{S|\Zhat}=\pi^{-1}(\Psihat(f_{|\nbigdhat_1}))$,
after shrinking $S$.

We take a finite covering
$U(\lambda_0)\times (X-D_1)
=\bigcup S_j$ by small multi-sectors
such that $S_1=S$.
We take a partition of unity $(\chi_i)$
subordinated to the covering $(S_i)$
such that $\chi_i$ depend only on $\arg(z_1)$.
On each $S_j$, let $a_{l,S_j}$
be the lift of $a_{l}$ to
$ \lefttop{\nbar}\Vzero_{b\vecp
 -\epsilon\vecdelta_{\nbar-\supp(\vecp)}}
 (\nbigq\nbige)_{\gminia,\Sbar_j}$.
As in Section \ref{subsection;08.2.2.5},
we obtain a section
$\ftilde_{C^{\infty}}$
of $U^{(\lambda_0)}_b
 \bigl(i_{g\dagger}\gbige(\ast \nbigd_1)
 \bigr)\otimes C^{\infty}$
given as follows:
\[
\ftilde_{C^{\infty}}:=
 \sum \chi_i\cdot P_l\cdot (a_{l,S_i}\otimes 1)
\]
As in Lemma \ref{lem;08.2.2.8},
the induced section of
$\Gr^{\Uzero}_b(i_{g\dagger}\gbige
 \otimes C^{\infty})$
is contained in
$W_p(N)\psitilde_{g,u}(\gbige)
 \otimes C^{\infty}$,
and moreover it induces
the section of
$\nbigv\otimes C^{\infty}
 \subset
 \Gr^{W(N)}_p\psitilde_{g,u}(\gbige)
 \otimes C^{\infty}$
whose restriction to
$\nbigdhat_1$ is equal to
$\Psihat(f_{|\nbigdhat_1})$.
Note that we have $f_{C^{\infty}|S}=f_S$
after $S$ is shrinked.
Hence, we can conclude that
$f_{S}$ gives a section of
$\nbigv_{|\Sbar}$,
and $f_{S|\Zhat}
=\pi^{-1}(\Psihat(f_{|\nbigdhat_1}))$.
\hfill\qed

\vspace{.1in}

Then, we obtain the isomorphism
$\bigoplus\nbigv_{\gminia|\Sbar}
\simeq \Gr(\nbigv)_{|\Sbar}$.
It is independent of the choice of
a $\DD$-flat spitting.
By gluing them, we obtain the isomorphism
$\bigoplus\nbigv_{\gminia|U(\lambda_0)\times X}
\simeq
 \Gr(\nbigv)_{|U(\lambda_0)\times X}$.
Since its restriction to $\nbigdhat_1$
is equal to $\Psihat$,
it is equal to the isomorphism
in Lemma \ref{lem;08.9.27.101}.

\vspace{.1in}

Note that the underlying $\nbigr_{D_I}$-modules
of $\Gr_{\gminia}(\nbigt_{u,p,I})$
and $\nbigt_{u,p,I}(E_{1,\gminia})$ are
$\Gr_{\gminia}(\nbigv)_
 {|\cnum_{\lambda}\times D_I^{\circ}}$
and $\nbigv_{\gminia|\cnum_{\lambda}\times D_I^{\circ}}$,
respectively.
To show that
 $\Gr_{\gminia}\bigl(
 \nbigt_{u,p,I}\bigr)$
and $\nbigt_{u,p,I}(E_{1,\gminia})$
are isomorphic,
we have only to compare the sesqui-linear pairings.
By the continuity,
we have only to compare them
on $\vecI(\lambda_0)\times S$,
where $\vecI(\lambda_0)$
denote the intersection of
$\vecS$ and a small neighbourhood of 
some generic $\lambda_0$,
and $S$ denotes a small sector in $X-D$.
Such a comparison can be done easily
by using the splittings 
(\ref{eq;07.9.8.16}),
(\ref{eq;07.9.8.15}),
(\ref{eq;07.11.29.15})
and (\ref{eq;08.9.27.50}).
Thus, the proof of Proposition \ref{prop;07.9.8.26}
is finished.
\hfill\qed

\section{Prolongation of 
ramified wild harmonic bundle on curve}
\label{subsection;07.10.19.30}
\subsection{The family of 
meromorphic $\lambda$-connections
 $\nbigq\nbige$}

Let $X:=\Delta$ and $D:=\{O\}$.
Let $\harmonicbundle$ be a wild harmonic bundle
defined on $X-D$, which is not necessarily unramified.
We take an appropriate ramified covering
$\varphi:(X',D')\lrarr (X,D)$ such that
$(E',\delbar_{E'},\theta',h')
:=\varphi^{-1}\harmonicbundle$ is unramified.
The ramification index is denoted by $e$.
We have the associated family of
meromorphic $\lambda$-flat bundles
$(\nbigq\nbige',\DD')$ on 
$\cnum_{\lambda}\times(X',D')$.
Since it is $\Gal(X'/X)$-equivariant,
we obtain the family of meromorphic $\lambda$-flat
bundles $(\nbigq\nbige,\DD)$
on $\cnum_{\lambda}\times (X,D)$.

We have the irregular decomposition:
\[
 (\nbigq\nbige',\DD')_{|
 U(\lambda_0)\times\Dhat'}
=\bigoplus_{\gminia\in\Irr(\theta)}
 \bigl(\nbigq\nbigehat'_{\gminia},\DDhat'_{\gminia}\bigr)
\]
We put
$\nbigq\nbigehat'_{\irr}:=
 \bigoplus_{\gminia\neq 0}\nbigq\nbigehat'_{\gminia}$.
Since $\nbigq\nbigehat'_{\irr}$ and
$\nbigq\nbigehat'_{0}$ are 
$\Gal(X'/X)$-equivariant,
we obtain the decomposition
\[
  \nbigq\nbigehat
=\nbigq\nbigehat_{0}\oplus
   \nbigq\nbigehat_{\irr}.
\]

\subsection{The good lattices}

For each $\lambda_0$,
we have the good lattices
$\nbigq^{(\lambda_0)}_{a}\nbige'$
of $\nbigq\nbige'$ on $U(\lambda_0)\times X'$.
They are $\Gal(X'/X)$-equivariant.
The descents are denoted by 
$\nbigq^{(\lambda_0)}_{a/e}\nbige'$.
It has the KMS-structure at $\lambda_0$
with the index set $\KMS(\nbigq\nbige^0)$.
We have the decomposition
$\nbigq^{(\lambda_0)}_a
 \nbige'_{|U(\lambda_0)\times\Dhat'}
=\bigoplus_{\gminia}
 \nbigq^{(\lambda_0)}_a\nbigehat'_{\gminia}$,
which induces a decomposition
\[
 \nbigq^{(\lambda_0)}_{c}
 \nbige_{|U(\lambda_0)\times\Dhat}
=\nbigq^{(\lambda_0)}_c\nbigehat_0
\oplus
 \nbigq^{(\lambda_0)}_c\nbigehat_{\irr}.
\]

\subsection{The $\nbigr$-module $\gbige$}

We naturally regard $\nbigq\nbige$
as an $\nbigr_X$-module.
Let $\gbige^{(\lambda_0)}$ be the $\nbigr_X$-submodule
of $\nbigq\nbige$ on $U(\lambda_0)\times X$,
generated by $\nbigq^{(\lambda_0)}_{<1}\nbige$.
The following lemma is clear from the construction.
\begin{lem}
We have the decomposition
$\gbige^{(\lambda_0)}_{|U(\lambda_0)\times\Dhat}
=\nbigq\nbigehat_{\irr}
\oplus \nbigq^{(\lambda_0)}_{\min}\nbigehat_0$,
where $\nbigq^{(\lambda_0)}_{\min}\nbigehat_0$ denotes
the $\nbigr_{\Dhat}$-submodule of
$\nbigq\nbigehat_0$ generated by
$\nbigq^{(\lambda_0)}_{<1}\nbigehat_0$.
\hfill\qed
\end{lem}

We have the wild harmonic bundles
$(E'_{\gminia},
 \delbar_{E'_{\gminia}},\theta'_{\gminia},h'_{\gminia})$
for $\gminia\in\Irr(\theta')$
on $X'-D'$, obtained as the full reduction of
$(E',\delbar_{E'},\theta',h')$.
Note that $(E'_0,\delbar_{E'_0},\theta'_0,h'_0)$ is 
equivariant,
and we have the descent $(E_0,\delbar_{E_0},\theta_0,h_0)$.
We have the family of $\lambda$-flat bundles
$(\nbigp\nbige_0,\DD_0)$
and the $\nbigr_X$-module $\gbige_0$,
associated to $(E_0,\delbar_{E_0},\theta_0,h_0)$.
Note we do not have to consider deformations
because $(E_0,\delbar_{E_0},\theta_0,h_0)$
is tame.

\begin{lem}
We have the natural isomorphisms
\[
 \nbigp\nbige_{0|\Dhat}\simeq
 \nbigq\nbigehat_0,
\quad
 \gbige_{0|\Dhat}\simeq
 \nbigq^{(\lambda_0)}_{\min}\nbigehat_0
\]
\end{lem}
\pf
We have the natural isomorphisms
$\nbigq^{(\lambda_0)}_c\nbigehat_0'
\simeq
 \nbigp^{(\lambda_0)}_c\nbige'_{0|\Dhat'}$.
Then, the first isomorphism is obtained.
Since $\gbige_{0|\Dhat}$ and
$\nbigq^{(\lambda_0)}_{\min}\nbigehat_0$ are
generated by 
$\nbigq^{(\lambda_0)}_{<1}\nbigehat_0$,
the second isomorphism is obtained.
\hfill\qed

\begin{lem}
For $\lambda_1\in U(\lambda_1)\subset U(\lambda_0)$,
we have 
$\gbige^{(\lambda_1)}=
 \gbige^{(\lambda_0)}_{|U(\lambda_1)\times X}$
on $U(\lambda_1)\times X$.
\end{lem}
\pf
We have only to compare the completion of them
along $U(\lambda_1)\times D$.
Both of them are isomorphic to
the direct sum of
$\nbigq_{\irr}\nbigehat$
and the completion of $\gbige_0$.
Hence, the claim of the lemma follows.
\hfill\qed

\vspace{.1in}
Therefore, we obtain the global $\nbigr_X$-module $\gbige$.

\begin{lem}
\label{lem;08.12.15.1}
$\gbige$ is strictly $S$-decomposable along 
the function $z^n$.
\end{lem}
\pf
We have only to show that 
the completion $\gbige_{|\Dhat}$ is strictly $S$-decomposable
along the function $z^n$.
It follows from the strict $S$-decomposability of
$\gbige_0$
(\cite{sabbah2} and \cite{mochi2}).
\hfill\qed

\subsection{The sesqui-linear pairing $\gbigc$}

Recall that we have the pairing for each $\lambda_0\in\vecS$:
\begin{multline*}
 C':
 \nbigq^{(\lambda_0)}_{<1}
 \nbige'_{|\vecI(\lambda_0)\times X'}
 \otimes
 \sigma^{\ast}
 \nbigq^{(-\lambda_0)}_{<1}
 \nbige'_{|\vecI(-\lambda_0)\times X'}
\lrarr \\
 \bigl\{f\in
 C^{\infty}\bigl(\vecI(\lambda_0)\times X'\bigr)\,\big|\,
 |f|=O\bigl(|z'|^{-2+\epsilon}\bigr) \,\,\exists\epsilon>0
 \bigr\}
\end{multline*}
It induces the pairing:
\begin{multline*}
 C:
 \nbigq^{(\lambda_0)}_{<1}
 \nbige_{|\vecI(\lambda_0)\times X}
 \otimes
 \sigma^{\ast}
 \nbigq^{(-\lambda_0)}_{<1}
 \nbige_{|\vecI(-\lambda_0)\times X}
\lrarr \\
 \bigl\{f\in
 C^{\infty}\bigl(\vecI(\lambda_0)\times X\bigr)\,\big|\,
 |f|=O\bigl(|z|^{-2+\epsilon}\bigr) \,\,\exists\epsilon>0
 \bigr\}
\end{multline*}
It induces 
$\gbigc:\gbige_{|\vecS\times X}
 \otimes\sigma^{\ast}\gbige_{|\vecS\times X}
\lrarr \distribution_{\vecS\times X/\vecS}$,
as in Section \ref{subsection;07.10.19.20}.
Thus, we obtain the $\nbigr$-triple
$\gbigt(E):=
 (\gbige,\gbige,\gbigc)$
even if $(E,\delbar_E,\theta,h)$ is ramified,
in the case $\dim X=1$.

\part{Kobayashi-Hitchin Correspondence}
\label{part;08.9.29.21}

\chapter{Preliminary}
\label{section;08.3.9.2}
This chapter is a collection of
miscellaneous preliminary for our study
in Part \ref{part;08.9.29.21}.

In Section \ref{subsection;08.9.29.19},
we recall some generality 
of filtered flat sheaf in \cite{mochi5}
related to $\mu_L$-stability,
with minor generalization.

In Section \ref{subsection;08.9.29.100},
we study Mehta-Ramanathan type theorem
for filtered flat sheaves.
This kind of result is always 
fundamental for 
the study on $\mu_L$-stability condition.
Note that we do not have to assume 
that the filtered sheaf is good.
This result will be useful in Chapter
\ref{section;08.10.24.41}.

In Section \ref{subsection;08.9.29.101},
we construct hermitian metrics
for good filtered $\lambda$-flat bundles
in a standard manner,
which have some nice property,
although not pluri-harmonic.
This is preliminary for 
Sections \ref{subsection;08.9.29.102}
and \ref{subsection;08.9.28.150}.

In Section \ref{subsection;08.9.29.102},
we collect some results related to
wild harmonic bundles on projective curves.
We review in Section \ref{subsection;08.1.14.5}
the Kobayashi-Hitchin correspondence 
for wild harmonic bundles on curves,
due to O. Biquard and P. Boalch 
\cite{biquard-boalch}.
Sections 
\ref{subsection;08.1.14.35}--\ref{subsection;08.9.29.105}
are preliminary for the proof of
Theorem \ref{thm;06.1.23.100}.
We show in Section \ref{subsection;08.1.14.35}
a convergence of 
the sequence of harmonic metrics for
$\epsilon$-perturbations.
Then, we study 
in Section \ref{subsection;08.1.14.36}
the convergence
of a sequence of hermitian metrics,
whose pseudo-curvatures 
converge to $0$ in some sense.
We also argue the continuity
of harmonic metrics 
for a holomorphic family of
stable good filtered flat bundles.

In Section \ref{subsection;08.9.28.150},
we give a sufficient condition
for a harmonic bundle to be good wild.
Proposition \ref{prop;07.10.13.11}
can be regarded as a kind of ``curve test''.

We explain in Section \ref{subsection;08.9.29.110}
some basic property
of the good filtered flat bundles
associated to good wild harmonic bundles.
Since the arguments for the proof
are essentially the same as those
in \cite{mochi4} and \cite{mochi5},
we give only outlines.

In Section \ref{subsection;07.10.15.80},
we explain a method of perturbation.
The point is to kill the nilpotent part
of the action of the residues on graded pieces.
The content is almost the same
as that in Section 2.1.6 of \cite{mochi5}.

\section{Preliminary for 
$\mu_L$-polystable filtered flat sheaf}
\label{subsection;08.9.29.19}
\subsection{Filtered flat sheaf and parabolic flat sheaf}

Let $X$ be a complex manifold
with a simple normal crossing divisor $D$
with the irreducible decomposition
$D=\bigcup_{i\in \Lambda}D_{i}$.
A filtered flat sheaf on $(X,D)$ is defined to be 
a pair of filtered sheaf 
$\vecE_{\ast}=\bigl(
 \prolongg{\veca}{E}\,\big|\,
 \veca\in\real^{\Lambda}
 \bigr)$ on $(X,D)$
and a meromorphic flat connection
$\nabla$ of the $\nbigo_X(\ast D)$-module
$\vecE=
 \bigcup_{\veca\in\real^{\Lambda}}
 \prolongg{\veca}{E}$.
(See Section 3.2 of \cite{mochi4}
 or Subsection \ref{subsection;10.5.28.1}
 in this paper for filtered sheaf.
 We use the symbol 
 $\prolong{E}$ instead of 
 $\prolongg{\veczero}{E}$.)
If $\vecE_{\ast}$ is a filtered bundle,
it is called a filtered flat bundle.
Recall that $(\vecE_{\ast},\nabla)$ is called regular,
if $\nabla\bigl(\prolongg{\veca}{E}\bigr)
\subset\prolongg{\veca}{E}\otimes\Omega^{1}(\log D)$
is satisfied.
\index{filtered flat sheaf}
\index{filtered flat bundle}
\index{sheaf $\prolong{E}$}

In this paper, a $\vecc$-parabolic sheaf
with a meromorphic connection 
is called a $\vecc$-parabolic flat sheaf.
By the operation to take the $\vecc$-truncation,
a filtered flat sheaf is equivalent to
a $\vecc$-parabolic flat sheaf.
\index{parabolic flat sheaf, $\vecc$-parabolic flat sheaf}

\begin{rem}
Usually, the logarithmic property is
contained in the definition of parabolic flat sheaf
or parabolic Higgs sheaf.
For example, in our previous paper {\rm\cite{mochi4}},
a parabolic Higgs sheaf means
a pair of parabolic sheaf $V_{\ast}$
with a Higgs field $\theta$
satisfying $\theta\bigl(\prolongg{\veca}{V}\bigr)
 \subset 
 \prolongg{\veca}{V}
 \otimes\Omega^{1,0}_X(\log D)$.
It would be more appropriate
that such an object is called 
a regular parabolic Higgs sheaf.
\hfill\qed
\end{rem}

\subsection{$\mu_L$-Stability}

\index{$\mu_L$-stable}
\index{$\mu_L$-semistable}
\index{$\mu_L$-polystable}

Let $X$ be a smooth irreducible 
$n$-dimensional projective variety
with a normal crossing hypersurface $D$
and an ample line bundle $L$.
The $\mu_L$-(semi)stability condition for
filtered flat sheaf $(\vecE_{\ast},\nabla)$ on $(X,D)$ 
is defined in a standard manner.
Namely,
we say $(\vecE_{\ast},\nabla)$ is $\mu_L$-stable
(resp. $\mu_L$-semistable)
if we have $\mu_L(\vecF_{\ast})<\mu_L(\vecE_{\ast})$
(resp. $\mu_L(\vecF_{\ast})\leq \mu_L(\vecE_{\ast})$)
for any sub-object
$(\vecF_{\ast},\nabla)\subset(\vecE_{\ast},\nabla)$
such that $0<\rank F<\rank E$,
where $\mu_L(\vecE_{\ast}):=
 (\rank E)^{-1}\,
 \int_Xc_1(L)^{n-1}\parchern_1(\vecE_{\ast})$.
We say that $(\vecE_{\ast},\nabla)$ 
is $\mu_L$-polystable,
if it is a direct sum of $\mu_L$-stable ones
$\bigoplus (\vecE_{i\ast},\nabla_i)$
such that $\mu_L(\vecE_{\ast})=\mu_L(\vecE_{i\ast})$.
By the correspondence of filtered sheaves
and parabolic sheaves,
we define
the $\mu_L$-stability,
$\mu_L$-semistability and $\mu_L$-polystability
conditions for parabolic flat sheaves.

As in \cite{mochi4},
we say that $(\vecE_{\ast},\nabla)$
is $\mu_L$-polystable 
with trivial characteristic numbers,
if it is $\mu_L$-polystable
and if each $\mu_L$-stable component 
$(\vecE_{i\ast},\nabla_i)$
satisfies 
$\mu_L(\vecE_{i\,\ast})=
 \int_{X}\parchern_{2,L}(\vecE_{i\,\ast})=0$.
\index{trivial characteristic number}

\subsection{Canonical decomposition}
\label{subsection;07.12.6.10}

Let $\bigl(\nbige^{(i)}_{\ast},\nabla^{(i)}\bigr)$ $(i=1,2)$
be $\mu_L$-semistable $\vecc$-parabolic flat sheaves
such that 
$\mu_L(\nbige^{(1)}_{\ast})
=\mu_L(\nbige^{(2)}_{\ast})$.
Let $f:(\nbige^{(1)}_{\ast},\nabla^{(1)})
 \lrarr (\nbige^{(2)}_{\ast},\nabla^{(2)})$ 
be a non-trivial morphism.
Let $(\nbigk_{\ast},\nabla_{\nbigk})$ 
denote the kernel of $f$,
which is equipped 
with the naturally induced 
parabolic structure and the flat connection.
Let $\nbigi$ denote the image of $f$,
and $\nbigitilde$ denote 
the saturated subsheaf of $\nbige^{(2)}$
generated by $\nbigi$,
i.e.,
$\nbigitilde/\nbigi$ is torsion,
and $\nbige^{(2)}/\nbigitilde$ is torsion-free.
The parabolic structures of $\nbige^{(1)}_{\ast}$ 
and $\nbige^{(2)}_{\ast}$
induce parabolic structures of
$\nbigi$ and $\nbigitilde$, respectively.
We denote the induced parabolic flat sheaves by
$(\nbigi_{\ast},\nabla_{\nbigi})$
and $(\nbigitilde_{\ast},\nabla_{\nbigitilde})$.
The following lemma can be shown
using the argument in the proof of Lemma 3.9 of \cite{mochi4}.
\begin{lem}
\label{lem;06.8.4.3}
$(\nbigk_{\ast},\nabla_{\nbigk})$,
$(\nbigi_{\ast},\nabla_{\nbigi})$
and $(\nbigitilde_{\ast},\nabla_{\nbigitilde})$
are also $\mu_L$-semistable
such that
$\mu_L(\nbigk_{\ast})=\mu_L(\nbigi_{\ast})
=\mu_L(\nbigitilde_{\ast})
=\mu_L(\nbige^{(i)}_{\ast})$.
Moreover,
 $\nbigi_{\ast}$ and $\nbigitilde_{\ast}$
are isomorphic in codimension one.
\hfill\qed
\end{lem}

\begin{lem}
\label{lem;06.8.22.20}
Let $\bigl(\nbige^{(i)}_{\ast},\nabla^{(i)}\bigr)$ $(i=1,2)$
be $\mu_L$-semistable reflexive saturated
parabolic flat sheaves
such that 
$\mu_L(\nbige^{(1)}_{\ast})
=\mu_L(\nbige^{(2)}_{\ast})$.
Assume either one of the following:
\begin{enumerate}
\item
One of $(\nbige^{(i)}_{\ast},\nabla^{(i)})$ is $\mu_L$-stable,
and $\rank(\nbige^{(1)})=\rank(\nbige^{(2)})$ holds.
\item
Both $(\nbige^{(i)}_{\ast},\nabla^{(i)})$ $(i=1,2)$
 are $\mu_L$-stable.
\end{enumerate}
If there is a non-trivial map 
$f:(\nbige^{(1)}_{\ast},\nabla^{(1)})
 \lrarr
 (\nbige^{(2)}_{\ast},\nabla^{(2)})$,
then $f$ is isomorphic.
\end{lem}
\pf
If $(\nbige^{(1)}_{\ast},\nabla^{(1)})$ is $\mu_L$-stable,
the kernel of $f$ is trivial due to Lemma \ref{lem;06.8.4.3}.
If $(\nbige^{(2)}_{\ast},\nabla^{(2)})$ is $\mu_L$-stable,
the image of $f$ and $\nbige^{(2)}$ are 
the same at the generic point of $X$.
Thus, we obtain that $f$ is generically isomorphic
in any case.
Then, we obtain that $f$ is isomorphic
in codimension one, due to Lemma 3.7 of \cite{mochi4}.
Since both $\nbige^{(i)}_{\ast}$ $(i=1,2)$
are reflexive and saturated,
we obtain that $f$ is isomorphic.
\hfill\qed

\begin{cor}
 \label{cor;05.9.14.5}
Let $(\nbige_{\ast},\nabla)$ be a $\mu_L$-polystable 
reflexive saturated parabolic flat sheaf.
Then, we have the unique decomposition:
\[
 (\nbige_{\ast},\nabla)
=\bigoplus_{j} 
 (\nbige^{(j)}_{\ast},\nabla^{(j)})\otimes\cnum^{m(j)}.
\]
Here, $(\nbige^{(j)}_{\ast},\nabla^{(j)})$ are $\mu_L$-stable
with $\mu_L(\nbige^{(j)}_{\ast})=\mu(\nbige_{\ast})$,
and they are mutually non-isomorphic.
It is called the canonical decomposition.
\index{canonical decomposition}
\hfill\qed
\end{cor}

\section{Mehta-Ramanathan type theorem}
\label{subsection;08.9.29.100}
\subsection{Statement}

Let $X$ be an $n$-dimensional smooth 
irreducible projective variety
with an ample line bundle $L$ over 
a field $k$ of characteristic $0$,
and let $D$ be a simple normal crossing divisor of $X$.

\begin{prop}
\label{prop;06.8.12.15}
Let $V_{\ast}$ be a parabolic sheaf on $(X,D)$
with a meromorphic flat connection $\nabla$.
Then, it is $\mu_L$-(semi)stable,
if and only if the following holds:
\begin{itemize}
\item
For any $m_1>0$,
there exists $m>m_1$ such that
$(V_{\ast},\nabla)_{|Y}$
is $\mu_L$-(semi)stable,
where $Y$ denotes a $1$-dimensional
complete intersection
of generic hypersurfaces of $L^m$.
\end{itemize}
\end{prop}

The proof will be given in Sections
\ref{subsection;07.6.11.2}--\ref{subsection;07.6.11.3}.
It is essentially the same as that in \cite{mochi4},
where we closely follow the argument
of V. Mehta, A. Ramanathan 
(\cite{mehta-ramanathan1}, \cite{mehta-ramanathan2})
and Simpson (\cite{s5}).
So, we will indicate only how to change in some points.

\vspace{.1in}

Before going into the proof,
we give a remark on
a characterization of semisimplicity
of meromorphic flat connection.
Let $(\nbige,\nabla)$ be a meromorphic flat connection
on $(X,D)$.
We recall the following lemma,
which we use implicitly in many places.
\begin{lem}
\label{lem;08.9.30.10}
$(\nbige,\nabla)$ is simple (resp. semisimple)
if and only if 
the associated Deligne-Malgrange filtered flat sheaf 
$(\vecE^{DM}_{\ast},\nabla)$ is $\mu_L$-stable
(resp. $\mu_L$-polystable).
\end{lem}
\pf
We have only to show the ``if'' part.
It was essentially observed by C. Sabbah
\cite{sabbah2}.
Namely,
we obtain
$(\vecF^{DM}_{\ast},\nabla)
 \subset(\vecE^{DM}_{\ast},\nabla)$
for any flat  subsheaf $\nbigf\subset\nbige$.
Because $\mu(\vecF^{DM}_{\ast})=
 \mu(\vecE^{DM}_{\ast})=0$,
$\vecF^{DM}_{\ast}$ breaks the stability of
$(\vecE_{\ast}^{DM},\nabla)$
if $\nbigf$ is non-trivial.
\hfill\qed

\begin{cor}
\label{cor;08.9.30.11}
$(\nbige,\nabla)$ is simple,
if and only if 
$(\nbige,\nabla)_{|Y}$ is simple,
where $Y$ denotes a complete intersection
of generic hypersurfaces of $L^m$ 
for arbitrary large $m$.
\end{cor}
\pf
The ``only if'' part is trivial.
The ``if'' part follows from
Proposition \ref{prop;06.8.12.15}
and Lemma \ref{lem;08.9.30.10}.
\hfill\qed

\subsection{Preliminary}
\label{subsection;07.6.11.2}

The following lemma can be shown
by using the argument in the proof of
Proposition 3.2 of \cite{mehta-ramanathan1}.

\begin{lem}
\label{lem;07.11.30.20}
Let $\nbigs_0$ be any bounded family of 
torsion-free sheaves on $X$.
There exists a large integer $m_0$
with the following property:
\begin{itemize}
\item
 Let $Y_1,\ldots,Y_r$ be generic hypersurfaces
 of $L^{m_i}$ for $m_i\geq m_0$.
 Let $Y=Y_1\cap\cdots \cap Y_r$.
 Let $p\geq m_0$,
and let $\nbigf$ be any member of $\nbigs_0$.
 Then, $H^0\bigl(Y,\nbigf\otimes L^{-p} \bigr)=0$.
\hfill\qed
\end{itemize}
\end{lem}

\begin{lem}
\label{lem;07.6.11.1}
Let $\nbigs_1$ be a bounded family
of torsion-free sheaves on $X$,
and let $\nbigs_2$ be a bounded family of
line bundles on $X$.
There exists a large integer $m_0$ 
with the following property:
\begin{itemize}
\item
 Let $Y_1,\ldots,Y_r$ be generic hypersurfaces
 of $L^{m_i}$ for $m_i\geq m_0$.
 Let $Y=Y_1\cap\cdots \cap Y_r$.
 Let $\nbigf\in\nbigs_1$ and $\nbigl\in\nbigs_2$.
 Let $\phi$ be any morphism
 $\nbigl\lrarr \nbigf\otimes\Omega^1_X$.
 Then, $\phi$ vanishes if and only if 
 the induced morphism
 $\phi':\nbigl_{|Y}\lrarr 
 \nbigf_{|Y}\otimes \Omega^1_{Y}$ 
 vanishes.
\end{itemize}
\end{lem}
\pf
We put $Y^{(i)}:=Y_1\cap \cdots \cap Y_i$.
If $m_0$ is sufficiently large,
we have
$H^0\bigl(Y^{(i)},\nbigl^{-1}\otimes\nbigf
 \otimes \Omega_X^1\otimes L^{-p}
 \bigr)=0$
and 
$H^0\bigl(Y^{(i)},\nbigl^{-1}\otimes\nbigf
 \otimes L^{-p}\bigr)=0$
for any $i$ and $p\geq m_0$,
according to Lemma \ref{lem;07.11.30.20}.
We have the exact sequence:
\[
 0\lrarr 
 L^{-m_{i+1}}\otimes
 \nbigl^{-1}\otimes\nbigf\otimes\Omega^1_{X|Y^{(i)}}
\lrarr
 \nbigl^{-1}\otimes\nbigf\otimes\Omega^1_{X|Y^{(i)}}
\lrarr
 \nbigl^{-1}\otimes\nbigf\otimes\Omega^1_{X|Y^{(i+1)}}
\lrarr 0
\]
Thus, the maps
$H^0\bigl(
 Y^{(i)},\nbigl^{-1}\otimes\nbigf\otimes\Omega^1_X
 \bigr)
\lrarr H^0\bigl(Y^{(i+1)},
 \nbigl^{-1}\otimes\nbigf\otimes\Omega^1_X\bigr)$
are injective,
if $m_0$ is sufficiently large.
Hence, we have only to show that
$\phi'=0$  implies that
$\phi_{|Y}:\nbigl_{|Y}\lrarr \nbigf\otimes\Omega^1_{X|Y}$
is trivial.
We have the exact sequence:
\[
0\lrarr \nbigl^{-1}\otimes\nbigf\otimes L^{-m_j}_{|Y}
\lrarr \nbigl^{-1}\otimes\nbigf\otimes \Omega^1_{Y^{(j-1)}|Y}
\lrarr \nbigl^{-1}\otimes\nbigf\otimes\Omega^1_{Y^{(j)}|Y}\lrarr 0.
\]
Hence, the maps
$H^0\bigl(Y,
 \nbigl^{-1}\otimes\nbigf\otimes\Omega^1_{Y^{(j)}}\bigr)
\lrarr
 H^0\bigl(Y,
 \nbigl^{-1}\otimes\nbigf
 \otimes\Omega^1_{Y^{(j+1)}}\bigr)$
are injective,
if $m_0$ is sufficiently large.
Thus we are done.
\hfill\qed

\vspace{.1in}

Let $E$ be a torsion-free sheaf on $X$
with a meromorphic flat connection
$\nabla:E\lrarr E\otimes\Omega^1_X(\ast D)$.
We remark the following obvious lemma
as a reference for the subsequent argument.
\begin{lem}
\label{lem;08.9.29.30}
Let $G$ be a subsheaf of $E$,
and let $\Gtilde$ be the saturated subsheaf
of $E$ generated by $G$,
i.e., $\Gtilde/G$ is torsion,
and $E/\Gtilde$ is torsion-free.
Then, $\nabla(G)\subset
 \Gtilde\otimes\Omega^1_X(\ast D)$
if and only if
$\nabla(\Gtilde)\subset
 \Gtilde\otimes\Omega^1_X(\ast D)$.
In this case, we say that
$G\subset E$ generates a flat subsheaf.
\hfill\qed
\end{lem}

Fix a constant $C$.
Let $\nbigs$ denote the family of
the couples $(\nbigl,\phi)$
of line bundles $\nbigl$ with $\deg(\nbigl)\geq C$
and a non-trivial morphism
$\phi:\nbigl\lrarr E$.
\begin{lem}
\label{lem;08.9.29.32}
$\nbigs$ is bounded.
We also have the boundedness of
the family of $\Cok(\phi)_{tf}$,
where $\Cok(\phi)_{tf}$ denotes the quotient
of $\Cok(\phi)$ by the torsion-part.
\end{lem}
\pf
Let $\nbigs'$ be the family of
torsion-free quotient sheaves $\nbigg$ of $E$
which are of the form $\Cok(\phi)_{tf}$
for some $(\nbigl,\phi)\in \nbigs$.
Note that $\deg(\Cok(\phi)_{tf})$ is bounded above.
Due to a result of Grothendieck 
(Lemma 2.5 of \cite{grothendieck}),
the family $\nbigs'$ is bounded.
Hence, 
we have the boundedness of
the family $\nbigs''$
of saturated subsheaves $\nbigk$ of $E$
generated by $\phi(\nbigl)$ 
for some $(\nbigl,\phi)\in\nbigs$.
By construction,
for any $(\nbigl,\phi)\in\nbigs$,
we have a member $\nbigk$ 
of the bounded family $\nbigs''$
such that $\phi(\nbigl)\subset\nbigk$.
We have 
$0\leq \deg(\nbigk)-\deg(\nbigl)\leq C'$
for some constant $C'$.
Hence, we obtain that 
the family $\nbigs$ is bounded.
\hfill\qed

\begin{lem}
\label{lem;07.6.11.5}
There exists an integer $m_0$,
depending only on $(E,\nabla)$ and 
the constant $C$,
with the following property:
\begin{itemize}
\item
 Let $m_i\geq m_0$.
 Let $Y=\bigcap_{i=1}^rY_i$,
 where $Y_i$ denotes general hypersurfaces
 of $L^{m_i}$.
 Let $(\nbigl,\phi)\in\nbigs$.
 Then, $\phi(\nbigl)$ generates a flat subsheaf of 
 $E\otimes\nbigo(\ast D)$,
 if and only if
 $\phi(\nbigl_{|Y})$ generates 
 a flat subsheaf of
 $E\otimes\nbigo(\ast D)_{|Y}$.
\end{itemize}
\end{lem}
\pf
There exists a large integer $N$
such that 
$\nabla(E)\subset
 E\otimes\Omega^1_X(N\, D)$.
Then,
$\nabla\circ \phi$ induces the $\nbigo_X$-homomorphism
$F_{\phi}:\nbigl\lrarr 
 \Cok(\phi)_{tf}\otimes\Omega^1_X(N\, D)$.
According to Lemma \ref{lem;08.9.29.30},
$F_{\phi}$ vanishes
if and only if
$\phi(\nbigl)$ generates a flat subsheaf of $E$.
Due to Lemmas \ref{lem;07.6.11.1}
and \ref{lem;08.9.29.32},
if $m_0$ is sufficiently large,
$F_{\phi}$ vanishes
if and only if the induced map
$F_{\phi,Y}:\nbigl_{|Y}\lrarr 
 \Cok(\phi)_{tf|Y}\otimes\Omega^1_Y(N\, D)$
vanishes.
The latter condition is equivalent to 
that $\phi_{|Y}(\nbigl_{|Y})$ generates 
a flat subsheaf of $E\otimes\nbigo(\ast D)_{|Y}$.
Thus, we are done.
\hfill\qed

\subsection{Family of degenerating curves}

We recall the setting in 
\cite{mehta-ramanathan1},
\cite{mehta-ramanathan2}
and Section 3.4 of \cite{mochi4}.
For simplicity, we assume $H^i(X,L^m)=0$
for any $m\geq 1$ and $i>0$.
We put $S_m:=H^0(X,L^m)$ for $m\in\seisuu_{\geq\,1}$.
For $\vecm=(m_1,\ldots,m_{n-1})\in\seisuu_{\geq\,1}^{n-1}$,
we put $S_{\vecm}:=\prod_{i=1}^{n-1} S_{m_i}$.
Let $Z_{\vecm}$ denote the correspondence variety,
i.e.,
$Z_{\vecm}=\bigl\{(x,s_1,\ldots,s_{n-1})
 \in X\times S_{\vecm}\,\big|\,
 s_i(x)=0,\,\,1\leq i\leq n-1 \bigr\}$.
The natural morphisms
$Z_{\vecm}\lrarr S_{\vecm}$
and $Z_{\vecm}\lrarr X$ are denoted by $q_{\vecm}$
and $p_{\vecm}$, respectively.
We put $Z^D_{\vecm}:=Z_{\vecm}\times_{X}D$
and $Z^{D_j}_{\vecm}:=Z_{\vecm}\times_XD_j$.
Let $K_{\vecm}$ denote the function field
of $S_{\vecm}$.
We put
$Y_{\vecm}:=
 Z_{\vecm}\times_{S_{\vecm}}K_{\vecm}$,
$Y^{D_j}_{\vecm}:=Z^{D_j}_{\vecm}\times_{S_{\vecm}}K_{\vecm}$
and $Y^D_{\vecm}:=Z^D_{\vecm}\times_{ S_{\vecm}}K_{\vecm}$.
The irreducible decomposition
of $Z^D_{\vecm}\times_{ S_{\vecm}}K_{\vecm}$
is given by $\bigcup_j Z^{D_j}_{\vecm}\times_{S_{\vecm}}K_{\vecm}$.

\vspace{.1in}

We fix a sequence of integers
$(\alpha_1,\ldots,\alpha_{n-1})$
with $\alpha_i\geq 2$.
We put $\alpha:=\prod\alpha_i$.
For a positive integer $m$,
let $(m)$ denote $(\alpha_1^{m},\ldots,\alpha_{n-1}^m)$.
Let $V_{\ast}$ be a coherent parabolic sheaf on $(X,D)$.
For each $m$,
we can take an open subset $U_m\subset S_{(m)}$
such that
(i) $q_{(m)}^{-1}(s)$ are smooth $(s\in U_m)$,
(ii) $q^{-1}_{(m)}(s)$ intersects with the smooth part 
of $D$ transversally,
(iii) $V_{\ast}$ is a parabolic bundle
on an appropriate neighbourhood of each $q^{-1}_{(m)}(s)\subset X$.
In the following, we will shrink $U_m$, if necessary.

Mehta and Ramanathan constructed a family of
degenerating curves.
Take integers $l>m>0$.
Let $A$ be a discrete valuation ring over $k$
with the quotient field $K$.
Then there exists a curve $C$ over $\Spec A$
with a morphism $\varphi:C\lrarr X\times\Spec A$
over $\Spec A$
with the properties:
(i) $C$ is smooth over $k$,
(ii) the generic fiber $C_K$ gives
 a sufficiently general $K$-valued point
 in $U_l$,
(iii) the special fiber $C_k$ is reduced
 with smooth irreducible components
 $C_k^i$ $(i=1,\ldots,\alpha^{l-m})$
 which are
 sufficiently general $k$-valued points in $U_m$.
 We use the symbol $D_C$
 to denote $C\times_X D$.

Then, we obtain the parabolic bundle
 $\varphi^{\ast}(V_{\ast})$ on $(C,D_C)$,
 which is denoted by $V_{\ast|C}$.
 The restriction to $C_K$ and $C_k^i$
 are denoted similarly.
 Let $W_{\ast}$ be a parabolic subbundle
 of $V_{\ast\,|C_K}$.
 Recall that 
 $W$ can be extended to
 the subsheaf $\Wtilde\subset V_{|C}$,
 flat over $\Spec A$  with the properties:
 (i) $\Wtilde$ is a vector bundle over $C$,
 (ii) $\Wtilde_{|C_k^i}\lrarr V_{|C_k^i}$ are injective.

We put $\Omega^1_{C/A}:=
 \Omega^1_{C}(\log C_k)
 \big/\Omega^1_{\Spec A}(\log t)$,
where $t$ denotes the closed point of $\Spec A$.
We have the induced meromorphic flat connection
of $V_{\ast}$ relative to $A$:
\[
 \nabla_C:V_{\ast}\lrarr V_{\ast}\otimes\Omega^1_{C/A}(\ast D_C)
\]
The restriction to $C_{k}^i$ are equal to 
the connection induced by the inclusion $C_k^i\subset X$.
If $W$ is a flat subbundle with respect to  $\nabla_{C_K}$,
then, $\nabla_C$ preserves $\Wtilde$,
and hence $\Wtilde_{|C_k^i}$ are also preserved by
$\nabla_{C_k^i}$.

\subsection{Proof of Proposition \ref{prop;06.8.12.15}}
\label{subsection;07.6.11.3}

\begin{lem}
\label{lem;06.8.13.1}
$(V_{\ast},\nabla)$ is $\mu_L$-semistable,
if and only if
there exists a positive integer $m_0$
such that $(V_{\ast},\nabla)_{|Y_{(m)}}$
is $\mu_L$-semistable
for any $m\geq m_0$.
\end{lem}
\pf
We have only to show
the ``only if'' part.
If $(V_{\ast},\nabla)_{|Y_{(m)}}$
is $\mu_L$-semistable for some $m$,
then $(V_{\ast},\nabla)_{|Y_{(l)}}$
is $\mu_L$-semistable for any $l> m$,
which we can show
by an argument in \cite{mehta-ramanathan1}.
(See also the first part of the proof of 
 Lemma 3.31 of \cite{mochi4}.)

We will show that $V_{\ast}$ is not semistable
if $V_{\ast|Y_{(m)}}$
are not semistable for any $m$.
By shrinking $U_m$ appropriately,
we may have the subsheaf
$W_{m\,\ast}$
of $p_{(m)}^{\ast}V_{\ast|q_{(m)}^{-1}U_m}$
such that 
(i) it is preserved by the induced relative connection
of $p_{(m)}^{\ast}V_{\ast|q_{(m)}^{-1}U_m}$,
(ii) $W_{m\,\ast|\,q_{(m)}^{-1}(s)}$
is the $\beta$-subobject of
$(V_{\ast},\nabla)_{|q_{(m)}^{-1}(s)}$
for any $s\in U_m$.
By an argument in \cite{mehta-ramanathan1}
(see also the proof of Lemma 3.31 of \cite{mochi4}),
we can show the existence of 
subsheaf $\Wtilde$ of $V$ such that
$\Wtilde_{|q_{(m)}^{-1}(s)}
=W_{m|q_{(m)}^{-1}(s)}$
for a sufficiently large $m$
and for some $s\in S_{(m)}$.
We can make $m$ arbitrarily large.
Then, $\Wtilde$ is preserved by $\nabla$
according to Lemma \ref{lem;07.6.11.5}.
Thus, $\Wtilde$ contradicts with
the $\mu_L$-semistability assumption of
$(V_{\ast},\nabla)$.
\hfill\qed

\vspace{.1in}

Now, Proposition \ref{prop;06.8.12.15}
follows from the next lemma.

\begin{lem}
\label{lem;06.8.12.10}
$(V_{\ast},\nabla)$ is $\mu_L$-stable,
if and only if
there exists a positive integer $m_0$
such that $(V_{\ast},\nabla)_{|Y_{(m)}}$
is $\mu_L$-stable
for any $m\geq m_0$.
\end{lem}
\pf
First, let us remark that
$(V_{\ast},\nabla)_{|q_{(m)}^{-1}(s)}$ 
has only obvious automorphisms
for any sufficiently large $m$
and general $s$,
if $(V_{\ast},\nabla)$ is $\mu_L$-stable.
To show it, we have only to consider the case
in which $V_{\ast}$ is reflexive and saturated
in the sense of Definition 3.17 of \cite{mochi4}.
Then, the sheaf $\nhom(V_{\ast},V_{\ast})$ is also reflexive.
For any $\nabla_{|q_{(m)}^{-1}(s)}$-flat
$f\in \nhom(V_{\ast},V_{\ast})_{|q_{(m)}^{-1}(s)}$,
we can take a lift $F\in\nhom(V_{\ast},V_{\ast})$
such that $F_{|q_{(m)}^{-1}(s)}=f$
if $m$ is sufficiently large.
We would like to show $\nabla(F)=0$,
which is a section of
$\nhom(V_{\ast},V_{\ast})\otimes
 \Omega_X^1(N\, D)$.
Because
$\nabla(F)_{|q_{(m)}^{-1}(s)}=0$,
the claim follows from Lemma \ref{lem;07.6.11.1}.
Hence, $F$ is an automorphism of
$(V_{\ast},\nabla)$,
which is a multiplication  of some constant.
Thus, $f$ is also a multiplication of some constant.

Assume that $(V_{\ast},\nabla)_{|Y^{(m)}}$ is not
stable for any $m$.
Then, by an argument in \cite{mehta-ramanathan2}
using the socle,
we can show the existence of a subsheaf
$0\neq \Wtilde_{m\ast}\subset V_{\ast}$
for an arbitrarily large $m$,
such that (i) $\Wtilde_{m|q_{(m)}^{-1}(s)}$
is preserved by $\nabla_{|q_{(m)}^{-1}(s)}$
for some general $s\in U_m$,
(ii) $\mu_L(\Wtilde_{m\ast})=\mu_L(V_{\ast})$.
(See also the proof of Lemma 3.32
 of \cite{mochi4}.)
Due to Lemma \ref{lem;07.6.11.5},
$\Wtilde$ is preserved by $\nabla$.
\hfill\qed

\section{Auxiliary metrics}
\label{subsection;08.9.29.101}
In this section,
we fix a non-zero $\lambda$.

\subsection{Regular case}

We recall a construction of a hermitian metric
for a regular filtered $\lambda$-flat bundles,
which are not harmonic
but satisfy some finiteness conditions.
We put $X:=\Delta$ and $D:=\{O\}$.
For any $\epsilon>0$,
let $g_{\epsilon}$ denote the metric of $X-D$
given by
$|z|^{-2+2\epsilon} dz\, d\zbar$.
We put $\mu_m:=\{\omega\in\cnum\,|\,\omega^m=1\}$,
which acts on $X$ by the multiplication
$(\omega,z)\longmapsto \omega\, z$.
Let $(\vecE_{\ast},\DDlambda)$ be 
a $\mu_m$-equivariant 
regular filtered $\lambda$-flat bundle.
The restriction to $X-D$ is denoted by $(E,\DDlambda)$.
We have the parabolic filtration
$F$ of $\prolong{E}_{|O}$.
For each $a\in\Par(\prolong{E})$,
we have the generalized eigen decomposition
$\Gr^{F}_a(\prolong{E}_{|O})
=\bigoplus_{\alpha\in\cnum}
 \Gr^{F,\EE}_{(a,\alpha)}(\prolong{E}_{|O})$.
Recall
$\KMS\bigl(\prolong{E} \bigr):=
 \bigl\{ (a,\alpha)\,\big|\,
 \Gr^{F,\EE}_{(a,\alpha)}(\prolong{E}_{|O})\neq 0
 \bigr\}$.
For any $(a,\alpha)\in\KMS\bigl(\prolong{E}\bigr)$,
we have the model bundle
$(E_{0,a,\alpha},
 \delbar_{0,a,\alpha},
 \theta_{0,a,\alpha},
 h_{0,a,\alpha})$
associated to
$\Gr^{F,\EE}_{(a,\alpha)}(\prolong{E}_{|O})$
with the nilpotent part $N_{a,\alpha}$
of $\Res(\DDlambda)$.
(See Section 6.2 of \cite{mochi2}.)
Namely, we set
\[
 E_{0,a,\alpha}
:=\Gr^{F,\EE}_{(a,\alpha)}(\prolong{E}_{|O})
 \otimes \nbigo_{X-D},
\quad
 \theta_{0,a,\alpha}
=N_{a,\alpha}\, \frac{dz}{z},
\]
and then there exists a model metric
$h_{0,a,\alpha}$
for $(E_{0,a,\alpha},\theta_{0,a,\alpha})$
such that the parabolic structure is trivial.

We take $u(a,\alpha)=(b,\beta)
 \in \real\times\cnum$
such that $\kmsmap(\lambda,u(a,\alpha))=(a,\alpha)$.
We have the model bundle $L\bigl(u(a,\alpha)\bigr)$.
(See Section 6.1 of \cite{mochi2}.)
Namely,
$L\bigl(u(a,\alpha)\bigr)$
is a line bundle $\nbigo_{X-D}\cdot e$
with the Higgs field
$\beta\, dz/z$
and the metric  $|e|=|z|^{-b}$.

Then, we obtain the harmonic bundle:
\[
 (E_0,\delbar_0,\theta_0,h_0):=
 \bigoplus_{(a,\alpha)}
(E_{0,a,\alpha},
 \delbar_{0,a,\alpha},
 \theta_{0,a,\alpha},
 h_{0,a,\alpha})
\otimes
 L\bigl(u(a,\alpha)\bigr)
\]
It is naturally equipped with 
the $\mu_m$-action.
Let $(\nbigp_{\ast}\nbigelambda_0,\DDlambda_0)$
denote the associated filtered $\lambda$-flat bundle.
By construction,
we can take a $\mu_m$-equivariant 
holomorphic isomorphism
$\Phi:\nbigp_0\nbigelambda_0
\lrarr \prolong{E}$ with the following property:
\begin{itemize}
\item
 $\Phi$ preserves the parabolic filtration.
\item
 The induced map
 $\Gr^F(\Phi):\Gr^F(\nbigp_0\nbigelambda_{0|O})
\lrarr
 \Gr^{F}(\prolong{E}_{|O})$
 is compatible with the residues.
\end{itemize}
We obtain the induced $\mu_m$-equivariant metric of $E$,
which is also denoted by $h_0$.

\begin{lem}
\mbox{{}}\label{lem;08.1.14.1}
\begin{itemize}
\item
 $(E,h_0)$ is acceptable.
\item
 Norm estimate holds for
 $(\vecE_{\ast},\DDlambda,h_0)$,
 i.e.,
 let $\vecv$ be a frame of
 $\prolong{E}$ such that
 (i) it is compatible with the parabolic filtration $F$,
 (ii) the induced frame of $\Gr^F(\prolong{E}_{|O})$
 is compatible with
 the weight filtration $W$ of
 the nilpotent part of $\Res(\DDlambda)$.
 We put $a(v_i):=\deg^F(v_i)$ and 
 $k(v_i):=\deg^W(v_i)$.
 Let $h_1$ be a hermitian metric of $E$
 given by
 $h_1(v_i,v_j):=\delta_{i,j}\,
 |z|^{-2a(v_i)}\,\bigl(-\log|z|\bigr)^{k(v_i)}$.
 Then, $h_0$ and $h_1$ are mutually bounded.
\item
 $\DDlambda-\DDlambda_0$
 is bounded with respect to $h_0$
 and $g_{\epsilon}$ for some $\epsilon>0$,
 under the identification of the bundles given by $\Phi$.
\item
 Let $\theta$ denote the section of
 $\End(E)\otimes\Omega^{1,0}$
 associated to $h_0$ and $\DDlambda$.
 (See {\rm\cite{s2}}, {\rm\cite{s5}} 
 or Section {\rm 2.2} of {\rm\cite{mochi5}}.)
 Then, 
 $\theta-\theta_0$ is bounded with respect to 
 $h_0$ and $g_{\epsilon}$
 with respect to $h_0$,
 under the above identification of the bundles.
\item
 $G(\DDlambda,h_0)$ is bounded
 with respect to $h_0$ and $g_{\epsilon}$
 for some $\epsilon>0$.
\end{itemize}
\end{lem}
\pf
The first and second claims 
follow from the property of tame harmonic bundles
\cite{mochi2}.
The third claim is clear by construction.
The fourth claim follows from
the third one and the relation
$\theta-\theta_0=
(1+|\lambda|^2)^{-1}
(\DDlambda-\DDlambda_0)$.
(See Subsection 2.2.1 of \cite{mochi5}.)
We have 
$\theta_0=O(1)\, dz/z$
with respect to $h_0$.
Hence $\theta=O(1)\, dz/z$.

Let $R(h_0)$ denote the curvature
of the holomorphic bundle $E$ with
the hermitian metric $h_0$.
Because $\dim X=1$,
we have the relation
\[
 G(\DDlambda,h_0)
=(1+|\lambda|^2)R(h)
+(1+|\lambda|^2)^2\bigl[
 \theta,\theta^{\dagger}
 \bigr].
\]
(See Subsection 2.2.3 of \cite{mochi5}.)
We have similar relation for
$G(\DDlambda_0,h_0)$.
Hence, we obtain the equality
\[
 G(\DDlambda,h_0)
=G(\DDlambda,h_0)
-G(\DDlambda_0,h_0)
=(1+|\lambda|^2)^2\bigl(
 [\theta,\theta^{\dagger}]
-[\theta_0,\theta_0^{\dagger}]
 \bigr).
\]
Then, the last claim follows.
\hfill\qed

\begin{cor}
If there exists 
a tame harmonic metric $h$
of $(E,\DDlambda)$ adapted to
$\vecE_{\ast}$,
the metrics $h$ and $h_0$ are mutually bounded.
\hfill\qed
\end{cor}

We put $\Xtilde:=\Delta^n$ and $\Dtilde:=\{z_1=0\}$.
We have the $\mu_m$-action on $\Xtilde$
by the multiplication on $z_1$,
which preserves $\Dtilde$.
Let $(\vecEtilde_{\ast},\DDlambdatilde)$ be 
a $\mu_m$-equivariant 
regular filtered $\lambda$-flat bundle
on $(\Xtilde,\Dtilde)$.
The restriction to $\Xtilde-\Dtilde$ is denoted by
$(\Etilde,\DDlambdatilde)$.
For any $\epsilon>0$,
let $\gtilde_{\epsilon}$ denote the metric
of $\Xtilde-\Dtilde$ given by
$|z_1|^{-2+2\epsilon}\, dz_1\, d\zbar_1
+\sum_{j=2}^n dz_j\, d\zbar_j$.

\begin{lem}
\label{lem;08.1.13.50}
There exists a $\mu_m$-equivariant hermitian metric 
 $\htilde_0$ of $\Etilde$
with the following property:
\begin{itemize}
\item
 $(\Etilde,\htilde_0)$ is acceptable.
\item
 Norm estimate holds for
 $(\vecEtilde_{\ast},\DDlambdatilde,\htilde_0)$,
  i.e.,
 let $\vecv$ be a frame of $\prolong{E}$
 such that
 (i) it is compatible with the parabolic filtration $F$,
 (ii) the induced frame of $\Gr^F(\prolong{E})$
 is compatible with the weight filtration $W$
 of the nilpotent part of $\Res(\DDlambdatilde)$.
 We put $a(v_i):=\deg^F(v_i)$ and 
 $k(v_i):=\deg^W(v_i)$.
 Let $h_1$ be a hermitian metric of $E$
 given by
 $h_1(v_i,v_j):=\delta_{i,j}\,
 |z_1|^{-2a(v_i)}\,\bigl(-\log|z_1|\bigr)^{k(v_i)}$.
 Then, $h_0$ and $h_1$ are mutually bounded.
\item
 $G(\DDlambdatilde,\htilde_0)$ is bounded
 with respect to $\htilde_0$
 and the metric $\gtilde_{\epsilon}$
 for some $\epsilon>0$.
\item
 Let $\thetatilde$ denote the section of
 $\End(\Etilde)\otimes\Omega^{1,0}$
 associated to $\DDlambdatilde$ and $\htilde_0$.
 Then, $\thetatilde=O(1)\, dz_1/z_1$.
\end{itemize}
\end{lem}
\pf
There exists 
a $\mu_m$-equivariant
regular filtered $\lambda$-flat bundle
$(\vecE_{\ast},\DDlambda)$
on $(X,D)$ such that
$(\vecEtilde_{\ast},\DDlambdatilde)$
is the pull back of $(\vecE_{\ast},\DDlambda)$
by the projection $(\Xtilde,\Dtilde)\lrarr (X,D)$.
Hence, the claim follows from Lemma \ref{lem;08.1.14.1}.
\hfill\qed

\subsection{Good filtered $\lambda$-flat bundle}
\label{subsection;08.1.14.2}

Let $X:=\Delta^n$
and $D:=\{z_1=0\}$.
Let $g_{\epsilon}$
denote the metric of $X-D$
given by
$|z_1|^{-2+2\epsilon}\, dz_1\,d\zbar_1
+\sum_{j=2}^ndz_j\, d\zbar_j$.
Let $(\vecE_{\ast},\DDlambda)$ 
be a good filtered $\lambda$-flat bundle on $(X,D)$.
The restriction to $X-D$ is denoted by $E$.
We take a ramified covering
$\varphi_e:(X',D')\lrarr (X,D)$ given by
$\varphi_e(z_1',z_2,\ldots,z_n)
 =(z_1^{\prime\,e},z_2,\ldots,z_n)$
such that
$\varphi_e^{\ast}(\vecE_{\ast},\DDlambda)$
is unramified.
We have the irregular decomposition:
\begin{equation}
 \label{eq;08.1.13.30}
 \varphi_e^{\ast}
 (\vecE_{\ast},\DDlambda)_{|\Dhat'}
=\bigoplus_{\gminia\in\Irr(\DDlambda)}
 (\vecEhat'_{\gminia\ast},
 \DDhat^{\prime\lambda}_{\gminia})
\end{equation}
Here, the filtration of 
$\varphi_e^{\ast}\vecE_{\ast}$
is given as in Section \ref{subsection;07.11.5.60}.
We have the natural $\Gal(X'/X)$-action
on $\Irr(\DDlambda)$.
For $\gminia\in\Irr(\DDlambda)$,
let $\Stab(\gminia)$ denote the stabilizer
of the $\Gal(X'/X)$-action.
Since $\DDhat^{\prime\lambda}_{\gminia}-d\gminia$
is logarithmic,
we have a $\Stab(\gminia)$-equivariant 
good filtered $\lambda$-flat bundle
$(\vecE'_{\gminia\ast},
 \DD^{\prime\lambda}_{\gminia})$
on $(X',D')$
such that
$(\vecE'_{\gminia\ast},
 \DD^{\prime\lambda}_{\gminia})
_{|\Dhat'}
\simeq
 (\vecEhat'_{\gminia\ast},
 \DDhat^{\prime\lambda}_{\gminia})$.
The restrictions to $X'-D'$ are denoted by $E_{\gminia}$.
We obtain an isomorphism:
\[
 \psihat:
 \bigoplus_{\gminia\in\Irr(\DDlambda)}
 (\vecE'_{\gminia\ast},
 \DD^{\prime\lambda}_{\gminia})_{|
 \Dhat'} 
\simeq
\varphi_e^{\ast}(\vecE_{\ast},\DDlambda)
 _{|\Dhat'}
\]
We may assume to have 
the $\Gal(X'/X)$-action on
$\bigoplus 
 (\vecE'_{\gminia\ast},\,
 \DD^{\prime\lambda}_{\gminia})$
which is equal to the natural 
$\Gal(X'/X)$-action on
$\varphi_e^{\ast}\vecE_{\ast}$
after the completion.
Let $N$ be any large number.
We can take a $\Gal(X'/X)$-equivariant 
holomorphic
(not necessarily flat) isomorphism
$\psi_N:
  \bigoplus \vecE'_{\gminia\ast}
\simeq
 \varphi_e^{\ast}\vecE_{\ast}$
which is equal to $\psihat$
on the $N$-th infinitesimal
neighbourhood $\Dhat^{\prime(N)}$
of $D'$.

We have the induced action of
$\Stab(\gminia)$ on 
$(\vecE'_{\gminia\ast},\DDprimelambda_{\gminia})$.
Let $\DD^{\prime\lambda\,\reg}_{\gminia}:=
 \DDprimelambda_{\gminia}-d\gminia$.
Applying Lemma \ref{lem;08.1.13.50}
to $(\vecE'_{\gminia\ast},
 \DD^{\prime\lambda\,\reg}_{\gminia})$,
we can take a $\Stab(\gminia)$-equivariant
hermitian metric $h'_{\gminia}$ of $E'_{\gminia}$
with the property as in Lemma \ref{lem;08.1.13.50}.
We may assume 
$\omega^{\ast}h'_{\gminia}
=h'_{\omega^{\ast}\gminia}$
on $E'_{\omega^{\ast}\gminia}$.
Hence, the metric $\bigoplus h'_{\gminia}$
of $\bigoplus E'_{\gminia}$ is
$\Gal(X'/X)$-equivariant.
We obtain the induced metrics
$h^{\prime(0)}$ and $h^{(0)}$ of 
$\varphi^{\ast}_eE$ and $E$,
respectively.
The metric $h^{(0)}$ is called
an auxiliary metric.

\begin{lem}
\label{lem;08.1.14.10}
Let $N$ be sufficiently large.
By construction,
we have the following:
\begin{itemize}
\item
 $(E,h^{(0)})$ is acceptable.
\item
 The norm estimate holds for 
 $(\vecE_{\ast},\DDlambda,h^{(0)})$.
 (See Lemma {\rm\ref{lem;08.1.13.50}}.)
\item
 $G(\DDlambda,h^{(0)})$ is bounded
 with respect to $h^{(0)}$ and $g_{\epsilon}$
 for some $\epsilon>0$.
\item
 Let $\theta_{h^{(0)}}$ denote the section
 of $\End(E)\otimes\Omega^{1,0}$
 induced by $h^{(0)}$ and $\DDlambda$.
 Then, $|\theta_{h^{(0)}}|_{h^{(0)}}$ is 
 of polynomial order with respect to $z_1$.
 More precisely,
\[
 \varphi_e^{\ast}\theta_{h^{(0)}}
-\bigoplus 
 (1+|\lambda|^2)^{-1}d\gminia
 \, \id_{\psi(E'_{\gminia})}
 \]
 is $O(1)\, dz_1'/z_1'
 +\sum_{j=2}^{n} O(1)\, dz_j$
 with respect to $h^{\prime(0)}$.
\end{itemize}
\end{lem}
\pf
The claims hold for
$\bigoplus 
 (\vecE'_{\gminia\ast},
 \DDprimelambda_{\gminia})$
with the metric $\bigoplus h'_{\gminia}$
by construction of $h'_{\gminia}$.
Under the identification of
$\bigoplus \vecE'_{\gminia\ast}$
and $\varphi_e^{\ast}\vecE_{\ast}$
via $\psi$,
we have
\[
\varphi_e^{\ast}\DDlambda
-\bigoplus_{\gminia}\DDprimelambda_{\gminia}
=O\bigl(z_1^{N/2}\bigr).
\]
Since $N$ is sufficiently large,
the claim of the lemmas are clear.
\hfill\qed

\subsection{Family version}
\label{subsection;08.1.14.6}

We put $X:=\Delta$ and $D:=\{O\}$.
Let $(\vecE_{\ast},\DDlambda)$
be a good filtered $\lambda$-flat bundle on $(X,D)$.
The restriction to $X-D$
is denoted by $(E,\DDlambda)$.
We take a real number $c\not\in \Par(\vecE_{\ast})$.
We have the induced parabolic structure
$F$ of $\prolongg{c}{E}$.
Take $\eta>0$ such that 
 $10\eta<\gap(\prolongg{c}{E}_{\ast})$.
(See Section 3.1 of \cite{mochi4}
 for $\gap$.)
For any $\epsilon\geq 0$,
let $\omega_{\epsilon}$ denote
the metric of $X-D$ given as follows:
\begin{equation}
 g_{\epsilon}:=
 \bigl(
 \epsilon^2|z|^{2\epsilon}
 +\eta^2|z|^{2\eta}
\bigr)\,\frac{dz\, d\zbar}{|z|^2}
\end{equation}
Let $F^{(\epsilon)}$ denote $\epsilon$-perturbation
of $F$ for any sufficiently small $\epsilon>0$,
as explained in (II) of 
Section \ref{subsection;07.10.15.80} below.

\begin{lem}
\label{lem;08.1.14.7}
If $\epsilon_0>0$ is sufficiently small,
there exists a family of hermitian metrics
$\{h_{0}^{(\epsilon)}\,|\,
 0\leq \epsilon\leq\epsilon_0\}$
of $E$
with the following property:
\begin{itemize}
\item
 $G(\DDlambda,h_{0}^{(\epsilon)})$
 are bounded
 with respect to $h_{0}^{(\epsilon)}$
 and $g_{\epsilon}$.
 The estimate is uniform for $\epsilon$.
\item
 Norm estimate holds for
 $(E,F^{(\epsilon)},h_{0}^{(\epsilon)})$.
\item
 $\{h_0^{(\epsilon)}\,|\,\epsilon>0\}$
 converges to
 $h_0^{(0)}$ in the $C^{\infty}$-sense 
 locally on $X-D$.
\item
 Let $t_{\epsilon}$ be determined by
 $\det(h_0^{(\epsilon)})\big/\det(h_0^{(0)})$.
 Then, $t_{\epsilon}$ and $t_{\epsilon}^{-1}$
 are bounded, uniformly in $\epsilon$.
\end{itemize}
\end{lem}
\pf
In Sections 4.2--4.4 of \cite{mochi5},
we give the construction of such a family of hermitian metric
for a regular $\lambda$-flat bundle.
The result can be extended
to the case in which
a meromorphic $\lambda$-flat bundle
is not necessarily regular but good,
with the same argument in Section
\ref{subsection;08.1.14.2}.
\hfill\qed

\section{Harmonic bundles on curves}
\label{subsection;08.9.29.102}
\subsection{Review of a result due to Biquard-Boalch}
\label{subsection;08.1.14.5}

Let $C$ be a smooth connected complex projective curve,
and $D$ be a finite subset of $C$.
Let $(\vecE_{\ast},\nabla)$ be a good filtered flat bundle
on $(C,D)$.
The restriction to $C-D$ is denoted by
$(E,\nabla)$.
We recall the following result due to Biquard-Boalch
\cite{biquard-boalch}.
(See also \cite{sabbah3}.)

\begin{prop}
\label{prop;07.10.16.10}
$(\vecE_{\ast},\nabla)$ is polystable
with $\pardeg(\vecE_{\ast})=0$,
if and only if
there exists a wild harmonic metric $h$ 
of $(E,\nabla)$ adapted to $\vecE_{\ast}$.
Such a metric is unique up to obvious ambiguity.
\end{prop}
\pf
We give an outline of a proof
based on Simpson's method,
for our later purpose.
The claim is easy in the rank one case.
Hence, we take a harmonic metric $h_{\det(E)}$
of $\bigl(\det(\vecE_{\ast}),\nabla\bigr)$.
Let $g_{\epsilon}$ be a Kahler metric of $C-D$,
which is given by
$|z|^{-2+\epsilon}$ around $P\in D$
for a holomorphic coordinate $z$
such that $z(P)=0$.

\begin{lem}
There exists a hermitian metric $h^{(0)}$
of $E$ with the following property:
\begin{itemize}
\item
 $(E,h^{(0)})$ is acceptable.
\item
 The norm estimate holds for
 $(E,\nabla,h^{(0)})$
 at each $P\in D$.
 (See Lemma {\rm\ref{lem;08.1.14.1}}.)
\item
 $G(\nabla,h^{(0)})$ is 
 bounded with respect to 
 $h^{(0)}$ and $g_{\epsilon}$
 for some $\epsilon>0$.
\item
 Let $\theta^{(0)}$ denote the section
 of $\End(E)\otimes\Omega^{1,0}$
 induced by $h^{(0)}$ and $\nabla$.
 Then, $|\theta^{(0)}|_{h^{(0)}}$ is 
 of polynomial order with respect to $|z|^{-1}$
 around each $P\in D$,
 where $z$ denotes a holomorphic coordinate
 around $P$ such that $z(P)=0$.
\item
 $\det(h^{(0)})=h_{\det(E)}$.
\end{itemize}
\end{lem}
\pf
Applying the construction
in Section \ref{subsection;08.1.14.2}
around each $P\in D$,
we obtain a metric of
$E_{|\nbigu-D}$,
where $\nbigu$ denotes some neighbourhood of $D$.
We extend it to a $C^{\infty}$ hermitian metric
$h^{(0)\star}$ of $E$.
Let $s$ be determined by
$h_{\det(E)}=
 \det(h^{(0)\star})\, s$.
Then, $s$ and $s^{-1}$ are bounded
on $C-D$.
Hence, we obtain the desired metric
$h^{(0)}$ with an obvious modification.
\hfill\qed

\vspace{.1in}

For any $\vecF_{\ast}\subset \vecE_{\ast}$,
we have
$\pardeg(\vecF_{\ast})
=\deg(F,h^{(0)})$ due to a result of Simpson
(Lemma 6.2 of \cite{s2}).
Hence, $(\vecE_{\ast},\nabla)$ is stable
if and only if $(E,\nabla,h)$ is analytically stable.
Then, due to a theorem of Simpson 
in \cite{s1} and \cite{s2},
we obtain the harmonic metric $h=h^{(0)}\, s$
such that (i) $s$ and $s^{-1}$ are bounded
with respect to $h^{(0)}$,
(ii) $\nabla s$ is $L^2$ with respect to $h^{(0)}$,
(iii) $\det(s)=1$.
(See also Proposition 2.49 of \cite{mochi5}.)
Let $\theta$ be the Higgs field
associated to $h$ and $\nabla$.
Let $d''$ denote the $(0,1)$-part of $\nabla$,
and let $\delta'_{h^{(0)}}$ denote the $(1,0)$-operator
induced by $d''$ and $h^{(0)}$.
Because of the $L^2$-property of $d'' s$
and the self-adjointness of $s$,
we obtain that
$\delta'_{h^{(0)}}s$ is $L^2$
with respect to $h^{(0)}$.

\vspace{.1in}

Let us show that
$(E,\nabla,h)$ is a {\em wild} harmonic bundle.
For any  $P\in D$,
we take a coordinate neighbourhood
$(U_P,z)$ with $z(P)=0$.
Let $\psi:\hyperh\lrarr U_P-\{P\}$
be given by $\psi(\zeta)=\exp(2\pi\sqrt{-1}\zeta)$.
We put $K_n:=\{\zeta\,|\,-1<\Re\zeta<1,\,\,
 n-1<\Image\zeta<n+1\}$.
Because 
$\theta=\theta^{(0)}+s^{-1}\delta_{h_0}'s/2$,
there exists a constant $C$,
which is independent of $n$,
such that the following holds:
\[
 \int_{K_n}|\psi^{\ast}\theta|_{h}^2
\leq
2\int_{K_n}|\psi^{\ast}\theta^{(0)}|_{h^{(0)}}^2+C
\]
Since $\int_{K_n}|\psi^{\ast}\theta|_h^2$
is the energy of the harmonic map on $K_n$
corresponding to the harmonic bundle
$\psi^{\ast}(E,\nabla,h)$
up to some positive constant multiplication,
we obtain the following estimate
(see \cite{ES}):
\[
 \sup_{K_n'}|\psi^{\ast}\theta|_h^2\leq
 C_1\,
 \sup_{K_n}|\psi^{\ast}\theta^{(0)}|_{h^{(0)}}^2
+C_2
\]
Here $K_n':=\{-2/3\leq \Re(\zeta)<2/3,\,\,
 n-2/3<\Image(\zeta)<n+2/3\}$.
Since the norm of $\theta^{(0)}$ 
is of polynomial order of $|z|^{-1}$
around $P$,
the norm of $\theta$ 
is also of polynomial order of $|z|^{-1}$.
In particular,
the eigenvalues of $\theta$ are also of polynomial order.
Thus, $(E,\delbar_E,\theta,h)$ is wild.

Conversely,
we can show the ``if'' part
by the same argument as that in Chapter 6 of \cite{s2}.

We know the norm estimate 
for wild harmonic bundles on curves
(Proposition \ref{prop;07.10.6.40}).
Hence, if $h_i$ $(i=1,2)$
are wild harmonic metrics
for $(E,\nabla)$ adapted to $\vecE_{\ast}$,
they are mutually bounded.
Then, 
by the same argument 
in the proof of Proposition 2.6 in \cite{mochi4},
we can show the existence of
a decomposition
$(\vecE_{\ast},\nabla)
=\bigoplus (\vecE_{i\,\ast},\nabla_i)$
such that
(i) it is orthogonal with respect to
 both $h_i$ $(i=1,2)$,
(ii) $h_1=a_i\, h_2$ on $E_i$ for some $a_i>0$.
\hfill\qed

\begin{rem}
We will also prove that
if a harmonic metric $h$ of $(E,\nabla)$
is adapted to $\vecE_{\ast}$,
then $(E,\nabla,h)$ is a wild harmonic bundle
(Proposition {\rm\ref{prop;07.6.20.2}}).
\hfill\qed
\end{rem}

\subsection{Convergence of a sequence of 
harmonic bundles}
\label{subsection;08.1.14.35}

Let $(C,D)$ be as in Section \ref{subsection;08.1.14.5}.
Let $(E,\vecF,\nabla)$ be a stable good parabolic flat bundle
on $(C,D)$ such that $\pardeg(E,\vecF)=0$.
Let $\vecF^{(\epsilon)}$ denote an $\epsilon$-perturbation
of $\vecF$
as in the case (II) of
Section \ref{subsection;07.10.15.80} below.
We have $\det(E,\vecF)=\det(E,\vecF^{(\epsilon)})$.
We take the harmonic metric $h_{\det E}$
of $\det(E,\vecF,\nabla)$.
Let $h^{(\epsilon)}$ be harmonic metrics of $(E,\nabla)$
adapted to $\vecF^{(\epsilon)}$
such that $\det h^{(\epsilon)}=h_{\det E}$.
Let $\theta^{(\epsilon)}$ denote the associated Higgs field.
We have a straightforward generalization
of Proposition 4.1 of \cite{mochi5}.

\begin{prop}
\label{prop;07.6.11.15}
We have the convergences
$h^{(\epsilon)}\lrarr h^{(0)}$
and $\theta^{(\epsilon)}\lrarr \theta^{(0)}$
in the $C^{\infty}$-sense
locally on $C-D$.
\end{prop}
\pf
Let $\eta$ be a small positive number
such that $\eta<\gap(E,\vecF)/10$.
Let $\epsilon_0$ be a small positive number
such that $10\rank(E)\,\epsilon_0<\eta$.
For any $0\leq \epsilon<\epsilon_0$,
let us take Kahler metrics $g_{\epsilon}$ of $C-D$
with the following property:
\begin{itemize}
\item
Let $P\in D$.
Let $(U_P,z)$ be a holomorphic coordinate
neighbourhood of $P$ such that $z(P)=0$.
The restriction of $g_{\epsilon}$
is as in Section \ref{subsection;08.1.14.6}.
\item
$g_{\epsilon}\to g_0$ 
in the $C^{\infty}$-sense
locally on $C-D$.
\end{itemize}

\begin{lem}
We can construct a family of hermitian metrics
$h_0^{(\epsilon)}$ of $E_{|C-D}$
with the following property:
\begin{itemize}
\item 
 $G(\DDlambda,h_{0}^{(\epsilon)})$
 are uniformly bounded
 with respect to $h_{0}^{(\epsilon)}$
 and $g_{\epsilon}$.
\item
 The norm estimate holds for  each
 $(E,\vecF^{(\epsilon)},h_{0}^{(\epsilon)})$.
\item
 $\{h_0^{(\epsilon)}\,|\,\epsilon>0\}$
 converges to
 $h_0^{(0)}$ in the $C^{\infty}$-sense 
 locally on $C-D$.
\item
 $\det(h_0^{(\epsilon)})=h_{\det(E)}$.
\end{itemize}
\end{lem}
\pf
It can be shown by using the argument
in Section 4.5.1 of \cite{mochi5}
together with Lemma \ref{lem;08.1.14.7}
(instead of Lemma 4.11 of \cite{mochi5}).
\hfill\qed

\vspace{.1in}
Then, the claim of Proposition \ref{prop;07.6.11.15}
can be shown
by using the argument in Section 4.5 of \cite{mochi5}.
\hfill\qed

\subsection{Convergence of 
a sequence of hermitian metrics}
\label{subsection;08.1.14.36}

Let $(C,D)$ be as in Section
\ref{subsection;08.1.14.5}.
Let $(E,\vecF,\nabla)$ be a stable good parabolic
flat bundle on $(C,D)$ with $\pardeg(E,\vecF)=0$.
For each $P\in D$,
let $(U_P,z)$ be a holomorphic coordinate
around $P$ such that $z(P)=0$.
Let $\vecF^{(\epsilon)}$ be an $\epsilon$-perturbation
as in the case (II) of Section \ref{subsection;07.10.15.80}.
We have $h_0^{(\epsilon_i)}$ be a harmonic metric
for each $(E,\vecF^{(\epsilon_i)},\nabla)$
for some sequence $\{\epsilon_i\}$
such that $\epsilon_i\lrarr 0$.
For simplicity of the description, we use $\epsilon$
instead of $\epsilon_i$.
We assume $\det h_0^{(\epsilon)}=\det h_0^{(0)}$.
Note that the sequence $h_0^{(\epsilon)}$ $(\epsilon>0)$
converges to $h_0^{(0)}$
in the $C^{\infty}$-sense locally on $C-D$
(Proposition \ref{prop;07.6.11.15}).

Let $N$ be a large positive number,
for example $N>10$.
We use Kahler metrics $g_{\epsilon}$ $(\epsilon\geq\,0)$
of $C-D$ which are as follows on $U_P$ for each $P\in D$:
\[
\bigl(
 \epsilon^{N+2}|z|^{2\epsilon}
+|z|^{2}
\bigr)\frac{dz\, d\zbar}{|z|^2}
\]
We assume that $\{g_{\epsilon}\}$ converges
to $g_{0}$ for $\epsilon\lrarr 0$
in the $C^{\infty}$-sense
locally on $C-D$.
By using the argument in Section 5.1 of \cite{mochi5}
with Proposition \ref{prop;07.6.11.15},
we can show the following lemma,
which is a straightforward generalization of 
Proposition 5.1 of \cite{mochi5}.

\begin{lem}
\label{lem;07.6.12.35}
Let $h^{(\epsilon)}$ $(\epsilon>0)$ be
hermitian metrics of $E_{|C-D}$
with the following properties:
\begin{enumerate}
\item
 Let $s^{(\epsilon)}$ be determined by
 $h^{(\epsilon)}=h_0^{(\epsilon)}\, s^{(\epsilon)}$. 
 Then, $s^{(\epsilon)}$ is bounded with respect to
 $h^{(\epsilon)}_0$,
and we have $\det s^{(\epsilon)}=1$.
 We also have the finiteness of the $L^2$-norm
$\bigl\|\nabla s^{(\epsilon)}\bigr\|
 _{2,h_0^{(\epsilon)},g_{\epsilon}}<\infty$.
(The estimates may depend on $\epsilon$.)
\item
$\|G(h^{(\epsilon)})\|_{2,h^{(\epsilon)},g_{\epsilon}}
 <\infty$
and
$ \lim_{\epsilon\to 0} \|G(h^{(\epsilon)})\|
 _{2,h^{(\epsilon)},g_{\epsilon}}=0$.
\end{enumerate}
Then the following claims hold.
\begin{itemize}
\item
The sequence $\{s^{(\epsilon)}\}$ is
weakly convergent to the identity of $E$
in $L_1^2$ locally on $C-D$.
\item
$|s^{(\epsilon)}|_{h_0^{(\epsilon)}}$
and 
$\bigl|(s^{(\epsilon)})^{-1}
 \bigr|_{h^{(\epsilon)}_0}$ 
are bounded on $C-D$ 
uniformly in $\epsilon$.
\hfill\qed
\end{itemize}
\end{lem}

\begin{cor}
 \label{cor;06.1.18.16}
\mbox{{}}
\begin{itemize}
\item
 The sequence $\bigl\{ h^{(\epsilon)} \bigr\}$
 is convergent to $h^{(0)}_0$ weakly in $L_1^2$
 locally on $C-D$. 
\item
The sequence 
$\bigl\{\nabla s^{(\epsilon)}\bigr\}$ 
is convergent to $0$
weakly in $L^2$ locally on $C-D$.
\item
The sequence
$\{\theta^{(\epsilon)}\}$ 
is convergent to $\theta^{(0)}$
weakly in $L^2$ locally on $C-D$.
\item
In particular,
the sequences are convergent almost everywhere.
\hfill\qed
\end{itemize}
\end{cor}

\subsection{Continuity for a holomorphic family}
\label{subsection;08.9.29.105}

Let $\nbigc\lrarr \Delta$ be a holomorphic family
of smooth projective curve,
and $\nbigd\lrarr\Delta$ be a relative divisor.
Let $(E,\vecF,\nabla)$ be a good filtered flat bundle
on $(\nbigc,\nbigd)$.
Let $t$ be any point of $\Delta$.
We denote the fibers by $\nbigc_t$ and $\nbigd_t$,
and the restriction of $(E,\vecF,\nabla)$
to $(\nbigc_t,\nbigd_t)$ is denoted by
$(E_t,\vecF_t,\nabla_t)$.
We assume 
$\pardeg(E_t,\vecF_t)=0$ 
and that $(E_t,\vecF_t,\nabla_t)$ is stable
for each $t$.
For simplicity, we also assume that
we are given a pluri-harmonic metric $h_{\det (E)}$
of $\det(E,\nabla)_{|\nbigc-\nbigd}$
which is adapted to the induced parabolic structure.

Let $h_{H,t}$ be a harmonic metric
of $(E_t,\vecF_t,\nabla_t)$
such that
$\det(h_{H,t})=h_{\det(E)\,|\,\nbigc_t}$.
They give the metric $h_H$ of $E$.
Let $\theta_{H,t}$ be the Higgs filed
obtained from
$(E_t,\nabla,h^{(\epsilon_t)})$,
which is a section of
$\End(E_t)\otimes\Omega^{1,0}_{\nbigc_t}(\log\nbigd_t)$.
They give the section $\theta_H$
of $\End(E)\otimes\Omega^{1,0}_{\nbigc/\Delta}(\log\nbigd)$,
where $\Omega^{1,0}_{\nbigc/\Delta}(\log \nbigd)$ denotes
the sheaf of the logarithmic relative $(1,0)$-forms.
The following lemma is a straightforward generalization
of Proposition 4.2 of \cite{mochi5},
and it can be proved by using an argument
in the proof of Propositions 4.1 of \cite{mochi5}
with Lemma \ref{lem;08.1.14.10}.

\begin{lem}
\label{lem;07.6.12.36}
$h_H$ and $\theta_H$ are continuous.
Their derivatives of any degree along the fiber directions
are also continuous.
\hfill\qed
\end{lem}

\section{Some characterizations of wildness of
 harmonic bundle}
\label{subsection;08.9.28.150}
We fix a non-zero $\lambda$ in this section.

\subsection{Statements}

Let $X:=\Delta^n$, $D_i:=\{z_i=0\}$
and $D:=\bigcup_{i=1}^{\ell}D_i$.
Let $\pi_i:X\lrarr D_i$ denote the natural projection.
We put $D_i^{\circ}:=
 D_i\setminus \bigcup_{j\neq i,1\leq j\leq \ell}D_j$.
Let $(\vecE_{\ast},\DDlambda)$ be 
a good filtered flat bundle on $(X,D)$.
Let $h$ be a pluri-harmonic metric of
$(E,\DDlambda):=(\vecE_{\ast},\DDlambda)_{|X-D}$,
and let $(E,\delbar_E,\theta)$ be the associated
Higgs bundle.
We will prove the following proposition
in Section \ref{subsection;08.9.29.90},
after showing the special cases
in Sections \ref{subsection;08.9.29.92}
and \ref{subsection;08.9.28.7}.

\begin{prop}
\label{prop;07.10.13.11}
Assume that there exist subsets $Z_i\subset D_i$
$(i=1,\ldots,\ell)$ with the following property:
\begin{itemize}
\item
The Lesbegue measure of $Z_i$ are $0$.
\item $h_{|\pi_i^{-1}(P)}$ is adapted to
$\vecE_{\ast|\pi_i^{-1}(P)}$
for any $P\in D_i^{\circ}\setminus Z_i$.
\end{itemize}
Then, the following holds:
\begin{itemize}
\item
The harmonic bundle
$\harmonicbundle$ is good and wild,
and $h$ is adapted to $\vecE_{\ast}$.
\item
If $(\vecE_{\ast},\DDlambda)$ is unramified,
$\harmonicbundle$ is also unramified,
and the following holds for any $P\in D$:
\begin{equation}
 \label{eq;10.6.16.1}
 \Irr(\theta,P)=
 \bigl\{(1+|\lambda|^2)^{-1}\gminia\,\big|\,
 \gminia\in\Irr(\DDlambda,P) \bigr\}.
\end{equation}
\end{itemize}
\end{prop}

We give a slightly different proposition.
For simplicity, 
we assume that
the determinant
$\det\harmonicbundle$ is good and wild,
and $\det(h)$ is adapted to $\det(\vecE_{\ast})$.
We prove the following proposition
in Section \ref{subsection;08.9.29.91}.

\begin{prop}
\label{prop;07.10.16.12}
If $h$ is adapted to $\vecE_{\ast}$,
then the following holds:
\begin{itemize}
\item
$h_{|\pi_i^{-1}(Q)}$ is adapted to
$\vecE_{\ast\,|\,\pi_i^{-1}(Q)}$
for any $Q\in D_i^{\circ}$ and for any $i=1,\ldots,\ell$.
\item
$\harmonicbundle$ is good and wild.
\item
If $(\vecE_{\ast},\DDlambda)$ is unramified,
$\harmonicbundle$ is also unramified,
and the following holds for any $P\in D$:
\[
 \Irr(\theta,P)=
 \bigl\{(1+|\lambda|^2)^{-1}\gminia\,\big|\,
 \gminia\in\Irr(\DDlambda,P)\bigr\}.
\]
\end{itemize}
\end{prop}

\subsection{A characterization of wildness
 of harmonic bundles on a punctured disc}
\label{subsection;08.9.29.92}

We put $X:=\bigl\{z\in\cnum\,\big|\,|z|<1\bigr\}$
and $D:=\{O\}$.
Let $(\vecE_{\ast},\DDlambda)$
be a good filtered $\lambda$-flat bundle
on $(X,D)$.
The restriction to $X-D$
is denoted by $(E,\DDlambda)$.

\begin{prop}
\label{prop;07.6.20.2}
Let $h$ be a harmonic metric of $(E,\DDlambda)$
adapted to $\vecE_{\ast}$.
Then, $(E,\DDlambda,h)$ is a wild harmonic bundle.
If $(\vecE_{\ast},\DDlambda)$ is regular,
then $h$ is tame.
\end{prop}
\pf
The second claim follows from
the first claim and the comparison of
the irregular values of $\lambda$-connection
and the Higgs field (Theorem \ref{thm;07.11.18.20}).
We can take an auxiliary metric $h_0$ 
for $(\vecE_{\ast},\DDlambda)$
as in Section \ref{subsection;08.1.14.2},
which has the property as in Lemma 
\ref{lem;08.1.14.10}.
Let $\theta_{h_0}$ denote the section of
$\End(E)\otimes\Omega^{1,0}$
associated to $h_0$ and $\DDlambda$.
Let $d''$ denote the $(0,1)$-part of $\DDlambda$,
and let $\delta_{h_0}'$ denote
the $(1,0)$-operator induced by 
$d''$ and $h_0$.
Let $s$ be determined by $h=h_0\, s$,
which is self-adjoint with respect to
both the metrics $h$ and $h_0$.
Because $h$ is adapted to $\vecE_{\ast}$,
we have $|s|_{h_0}=O(|z|^{-\epsilon})$
for any $\epsilon>0$.

\begin{lem}
\label{lem;07.10.16.11}
$|s|_{h_0}$ is bounded.
\end{lem}
\pf
We have the inequality
$ \Delta\log \tr(s)\leq
 C\, \bigl|\Lambda G(h_0)\bigr|_{h_0}$
on $X-D$.
 (See Section 2.2 of \cite{mochi5},
 for example.)
Hence,
we have 
$\Delta(\log \tr(s)+\epsilon\log|z|)
 \leq
 C\, \bigl|\Lambda G(h_0)\bigr|_{h_0}$
on $X-D$ for any $\epsilon>0$.
Since $\log\tr(s)+\epsilon\log|z|$ is bounded above,
the inequality holds on $X$ as distributions.
Hence, the values of
$\log \tr(s)+\epsilon\log|z|$
is dominated by the values on the boundary.
By taking $\epsilon\to 0$,
we obtain the boundedness of $\tr(s)$.
Then, the claim of the lemma follows.
\hfill\qed

\begin{lem}
\label{lem;07.6.20.1}
Let $(L,\DDlambda)$ be a flat $\lambda$-bundle
on $X-D$ of rank one.
Let $h_1$ be a harmonic metric of $(L,\DDlambda)$.
Let $u$ be a holomorphic section of $L$ on $X-D$
with the following property:
\begin{itemize}
\item
 Let $f$ be determined by 
 $\DDlambda u=u\, f\, dz$.
 Then, $f$ is meromorphic.
\item 
 $|u|_h\leq C\, |z|^a$ for some $a\in\real$.
\end{itemize}
Then, there exist $b\in\real$ and $C_i>0$ $(i=1,2)$
such that
$C_1\leq |u|_h\, |z|^b\leq C_2$.
\end{lem}
\pf
Taking the tensor product
with an appropriate wild harmonic bundle on $X-D$
of rank one,
and replacing $u$ with $e^g\, u$ 
for some holomorphic function $g$,
we may assume $\DDlambda u=0$.
Let $d''$ denote the $(0,1)$-part of $\DDlambda$.
We put $A:=|u|_h^2$.
Due to $\rank L=1$,
we have the following:
\begin{equation}
 \label{eq;07.10.12.50}
\delbar\del \log A=
 R(d'',h)=(1+|\lambda|^2)^{-1}\,
 G(\DDlambda,h)=0
\end{equation}
Hence, we have the expression
$\log A=-b\,\log |z|+\Re F(z)$,
where $F$ is a holomorphic function
on $X-D$.
By assumption,
we have $\log A\leq -a\, \log r$
on $X-D$,
and hence $F$ must be holomorphic on $X$.
Then, the claim of the lemma follows.
\hfill\qed

\vspace{.1in}
Let us return to the proof of Proposition
\ref{prop;07.6.20.2}.
We have $\det(h)=\det(h_0)\, \det(s)$.
By Lemma \ref{lem;07.6.20.1},
we have $\det(s)\geq C_1\, |z|^{-N_1}$
for some positive constants $C_1$ and $N_1$.
Hence, we obtain 
$C_2\, |z|^{N_2}\leq |s|_{h_0}$
and 
$|s^{-1}|_{h_0}\leq C_3|z|^{-N_3}$.

Recall the following formula (see \cite{s1} or 
 Section 2.2 of \cite{mochi5}):
\begin{equation}
\label{eq;08.9.29.50}
 (1+|\lambda|^2)\Delta\tr(s)=-\tr\bigl(
 s\sqrt{-1}\Lambda G(h_0)
 \bigr)
-\bigl|\DDlambda(s)\cdot s^{-1/2}\bigr|_{h_0}^2
\end{equation}
\begin{lem}
\label{lem;08.1.13.30}
We have the finiteness
$\int_{X-D}
 \bigl|\DDlambda s\cdot s^{-1/2}\bigr|_{h_0}^2<\infty$.
\end{lem}
\pf
Although it follows from a lemma in \cite{s2},
we give a direct argument.
Let $\rho$ be an $\real_{\geq\,0}$-valued
$C^{\infty}$-function on $\real$ such that
$\rho(t)=1$ for $t\leq 1$ and $\rho(t)=0$ for $t\geq 2$.
We can take such a function
as $\rho'\,\rho^{-1/2}$ is bounded.
We put $\chi_N(z)=\rho\bigl(-N^{-1}\log|z|\bigr)$.
We have the following equality:
\[
 \DDlambda\DDlambdastar(\chi_N\, s)
=\chi_N\,\DDlambda\DDlambdastar s
+\DDlambda\chi_N\cdot \DDlambdastar s
+\DDlambda s\cdot \DDlambdastar \chi_N
+s\,\DDlambda\DDlambdastar\chi_N
\]
(See Section 2.1 of \cite{mochi5}
 for $\DDlambdastar$.)
By using (\ref{eq;08.9.29.50}),
we obtain the following inequality
for some constant $A_1$, which is independent  of $N$:
\[
 \int\chi_N\,
 \bigl|\DDlambda s\cdot s^{-1/2}\bigr|_{h_0}^2
\leq A_1+2\int
 \bigl|\DDlambda s\bigr|_{h_0}\,
 \bigl|\DDlambdastar\chi_N\bigr|
\]
The following inequality holds
for some $A_2>0$:
\begin{multline}
 \int\bigl|\DDlambda s\bigr|_{h_0}
 \bigl|\DDlambdastar \chi_N\bigr|
\leq
 \left(
 \int \chi_N^{-1}
 \bigl|\DDlambdastar\chi_N\bigr|^2
 |s^{1/2}|_{h_0}^2
 \right)^{1/2}
  \left(
 \int\chi_N
 \bigl|\DDlambda s\cdot s^{-1/2}\bigr|_{h_0}^2
 \right)^{1/2} \\
\leq
A_2\,
 \left(
 \int\chi_N\, 
 \bigl|\DDlambda s\cdot s^{-1/2}\bigr|_{h_0}^2
 \right)^{1/2}
\end{multline}
Hence, we obtain
$\int\chi_N\, 
 |\DDlambda s\cdot s^{-1/2}|_{h_0}^2\leq A_3$,
independently of $N$.
Thus, we obtain
$\int\bigl|\DDlambda(s)\cdot s^{-1/2}\bigr|_{h_0}^2
 \leq A_3$.
\hfill\qed

\vspace{.1in}

In particular,
we obtain the finiteness
$\int_{X-D}\bigl|d'' s\bigr|_{h_0}^2<\infty$
from Lemma \ref{lem;08.1.13.30}.
Due to the self-adjointness of $s$,
we also obtain the finiteness
$\int_{X-D}
 \bigl|\delta'_{h_0}s\bigr|^2_{h_0}<\infty$.

We take a universal covering 
$\psi:\hyperh\lrarr\Delta^{\ast}$
given by $\psi(\zeta)=\exp(2\pi\sqrt{-1}\zeta)$.
We put 
$ K_n:=\bigl\{\zeta\in\hyperh\,\big|\,
 -1<\Re\zeta<1,\,\,
 n-1<\Image(\zeta)<n+1\bigr\}$.
Let $\theta$ be the Higgs field corresponding
to the harmonic bundle $(E,\DDlambda,h)$.
We have
$\theta=\theta_{h_0}+s^{-1}\delta'_{h_0}s$.
Due to the $L^2$-property of $\delta'_{h_0}s$
and the estimate $|s^{-1}|_{h_0}\leq C_3|z|^{-N_3}$,
we obtain the following:
\[
 \int_{K_n} 
 \bigl|\psi^{\ast}\theta\bigr|_{h_0}
\leq C_{12}\, e^{C_{13}n}
\]
Since $h$ and $h_0$ are mutually bounded
up to polynomial order of $|z|^{-1}$,
we obtain the following:
\begin{equation}
 \label{eq;08.9.29.51}
 \int_{K_n}\bigl|\psi^{\ast}\theta \bigr|^2_{h}
\leq C_{14}\, e^{C_{15}n}
\end{equation}
Recall that $|\theta|^2$ is the energy function
for a harmonic map up to positive constant
multiplication.
By using a result in \cite{ES},
we obtain $|\theta|_h^2\leq C_{16}|z|^{-C_{17}}$
from (\ref{eq;08.9.29.51}).
In particular, $\theta$ is wild.
\hfill\qed

\begin{cor}
Let $(\vecE_{\ast},\DDlambda,h)$ be as in 
Proposition {\rm\ref{prop;07.6.20.2}}.
Then, the norm estimate holds for $h$.
\hfill\qed
\end{cor}

\subsection{Curve test in the smooth divisor case}
\label{subsection;08.9.28.7}

We put $X:=\Delta^n$ and $D_i:=\{z_i=0\}$.
Let $(\vecE_{\ast},\DDlambda)$ be 
a good filtered $\lambda$-flat bundle on $(X,D_1)$.
The restriction to $X-D_1$ is denoted by
$(E,\DDlambda)$.
Let $h$ be a pluri-harmonic metric of
$(E,\DDlambda)$ on $X-D_1$,
and let $(E,\delbar_E,\theta)$ be the associated
Higgs bundle.
Let $\pi_i:X\lrarr D_i$ denote the projection.
For each $P\in D_1$,
we put $\pi_1^{-1}(P)^{\ast}:=
 \pi_1^{-1}(P)-\{P\}$.

\begin{prop}
\label{prop;07.6.12.37}
Assume that there exists a subset $Z\subset D_1$
with the following property:
\begin{itemize}
\item
The Lesbegue measure of $Z$ is $0$.
\item
For any $P\in D_1-Z$,
the restriction
$h_{|\pi_1^{-1}(P)^{\ast}}$ is adapted to
$\vecE_{\ast|\pi_1^{-1}(P)}$.
\end{itemize}
Then, the following holds:
\begin{itemize}
\item
 $\harmonicbundle$ is good wild harmonic bundle,
 and $h$ is adapted to $\vecE_{\ast}$,
 i.e.,
 $\vecE_{\ast}=\nbigp_{\ast}\nbigelambda$.
\item
 If $(\vecE_{\ast},\DDlambda)$ is unramified,
 $\harmonicbundle$ is also unramified,
 and 
\[
 \Irr(\theta)=
 \bigl\{(1+|\lambda|^2)^{-1}\gminia\,\big|\,
 \gminia\in \Irr(\DDlambda)\bigr\}.
\]
\end{itemize}
\end{prop}
\pf
In the following,
we will shrink $X$ without mention
if it is necessary.
By taking a ramified covering,
we may assume that 
$(\vecE_{\ast},\DDlambda)$ is unramified.
We divide the proof into
several steps.

\subsubsection{Decomposition of $(E,\theta)$}

We have the expression
$\theta=\sum f_i\, dz_i$.
By Proposition \ref{prop;07.6.20.2}
and the assumption,
$(E,\delbar_E,\theta,h)_{|\pi_1^{-1}(P)}$ are
wild harmonic bundles for any $P\in D_1-Z$.
Hence, 
$P(t)=\det\bigl(t-z_1\, f_1\bigr)$
is contained in $M(X,D_1)[t]$.
We put $T:=\bigl\{
 (1+|\lambda|^2)^{-1}z_1\del_1\gminia\,\big|\,
 \gminia\in\Irr(\DDlambda)
 \bigr\}$.
According to Theorem \ref{thm;07.11.18.20},
the following holds:
\begin{itemize}
\item
 Take any $Q\in D_1-Z$.
 Let $\alpha$ be a multi-valued meromorphic function
 on $\pi_1^{-1}(Q)$,
 which is a solution of
 $P(t)_{|\pi_1^{-1}(Q)}=0$.
 Then, there exists $\gminib(\alpha,Q)\in T$
 such that the following holds:
\[
 \bigl|
 \gminib(\alpha,Q)_{|\pi_1^{-1}(Q)}
-\alpha \bigr|=O(1)
\]
\end{itemize}
In other words,
the assumption made in Section
{\rm\ref{subsection;07.10.13.10}} 
(Appendix below) is satisfied
with the good set of irregular values $T$.
Applying Proposition \ref{prop;07.10.13.1} below,
we obtain the decomposition
$P(t)=\prod_{\gminib\in T}
 P_{\gminib}(t-\gminib)$
such that $P_{\gminib}(t)\in H(X)[t]$.
By using it,
we obtain the decomposition
\begin{equation}
 \label{eq;08.1.14.13}
 (E,f_1)=\bigoplus_{\gminia\in \Irr(\DDlambda)} 
 (E'_{\gminia},f_{1,\gminia}),
\end{equation}
such that 
$\beta-(1+|\lambda|^2)^{-1}
 z_1\del_z\gminia$ is bounded
for any solution
$\beta$ of $\det\bigl(t-z_1f_{1,\gminia}\bigr)=0$.
In particular,
the restrictions $(E,\delbar_E,\theta,h)_{|\pi_1^{-1}(P)}$
are wild harmonic bundles
for any $P\in D_1$.
The decomposition (\ref{eq;08.1.14.13})
is preserved by $f_i$ $(i=2,\ldots,n)$
due to the commutativity 
$[f_i,f_1]=0$.
Note the decomposition {\rm(\ref{eq;08.1.14.13})}
is holomorphic with respect to $\delbar_E$,
but not necessarily holomorphic with respect to
the $(0,1)$-part of $\DDlambda$.
Let $\pi_{\gminia}$ denote the projection
onto $E_{\gminia}'$
with respect to {\rm(\ref{eq;08.1.14.13})}.

\begin{lem}
\label{lem;08.1.14.20}
For each $\gminia\in\Irr(\DDlambda)$,
we put 
$f^{\reg}_{1,\gminia}:=
 f_{1,\gminia}-(1+|\lambda|^2)^{-1}
 \del_1\gminia\, \pi_{\gminia}$.
Then,
$\det\bigl(t-z_1f^{\reg}_{1,\gminia}\bigr)$
are contained in $H(X)[t]$,
and the coefficients of
the restriction
$\det\bigl(t-z_1f^{\reg}_{1,\gminia}\bigr)_{|D_1}$
are constant.
\end{lem}
\pf
The first claim follows from the construction.
The second claim follows from 
Proposition \ref{prop;07.7.19.31}
and the assumption
$(\vecE_{\ast},\DDlambda)$ is good.
\hfill\qed

\subsubsection{Norm estimate for
 $(\vecE_{\ast},\DDlambda,h)$}

Recall that Simpson's main estimate
in Section \ref{subsection;07.10.7.1}
depends only on the behaviour of the eigenvalues
of the Higgs field.
Due to Corollary \ref{cor;07.11.22.30},
the curvature
$R(h_{|\pi_1^{-1}(P)})$
of $(E,\delbar_E,h)_{|\pi_1^{-1}(P)}$ is
bounded with respect to
$h_{|\pi_1^{-1}(P)}$ and the Poincar\'e metric
$g_{\poin}$ of $\pi_1^{-1}(P)^{\ast}$,
and the estimate  is uniform in $P\in D_1$.

\vspace{.1in}

As in Section \ref{subsection;08.1.14.2},
we take good filtered flat bundles 
$(\vecE_{\gminia\ast},\DDlambda_{\gminia})$
such that $\DDlambda_{\gminia}$ is $\gminia$-regular
for each $\gminia\in\Irr(\DD)$,
and an isomorphism
$ \psi_N:
 \bigoplus \vecE_{\gminia\ast}
\simeq
  \vecE_{\ast}$
which is an approximation
of the irregular decomposition
in $N$-th order.
Then,
we take an auxiliary metric $h_0$
of $(E,\DDlambda)$,
which has the property in Lemma \ref{lem;08.1.14.10}.
By taking a ramified covering,
we may assume that
$G(\DDlambda,h_0)$ is bounded
with respect to 
$h_0$ and the Euclidean metric of $X$.
Let $s$ be determined by
$h=h_0\, s$.
Let $\Delta_i^{\lambda}$ denote 
$-(1+|\lambda|^2)\del_{z_i}\del_{\zbar_i}$.

\begin{lem}
$s$ and $s^{-1}$ are bounded
with respect to $h_0$.
Namely, $h$ and $h_0$ are mutually bounded.
In particular, 
$h$ is adapted to $\vecE_{\ast}$,
and more strongly,
the norm estimate holds for
$(\vecE_{\ast},\DDlambda,h)$.
\end{lem}
\pf
Due to $G(h,\DDlambda)=0$,
the boundedness of $G(h_0,\DDlambda)$
and an inequality in Subsection 2.2 of \cite{mochi5},
we have a constant $A$
such that the following holds
on $\pi_1^{-1}(P)^{\ast}$,
independently of $P\in D_1$:
\[
 \Delta^{\lambda}_1\log\tr(s_{|\pi_1^{-1}(P)}) \leq A,
\quad\quad
 \Delta^{\lambda}_1\log\tr(s^{-1}_{|\pi_1^{-1}(P)}) \leq A
\]
Since we have already known
the boundedness of $s_{|\pi_1^{-1}(P)^{\ast}}$
and $s^{-1}_{|\pi_1^{-1}(P)^{\ast}}$ for 
any $P\in D_1$,
the inequality holds on $\pi_1^{-1}(P)$
in the sense of distributions for such $P$.
Hence $\bigl|s_{|\pi_1^{-1}(P)^{\ast}}\bigr|_{h_0}$
and $\bigl|s^{-1}_{|\pi_1^{-1}(P)^{\ast}}\bigr|_{h_0}$
are dominated by their values at
$\del\pi_1^{-1}(P)$ for any $P\in D_1$.
Hence, we obtain the uniform boundedness
of $s$ and $s^{-1}$.
\hfill\qed

\subsubsection{Estimate for
 $f_i$ $(i=2,\ldots,n)$}

We put
$F_i':=\sum_{\gminia}
 (1+|\lambda|^2)^{-1}\del_i\gminia\,\pi_{\gminia}$
and $f_i^{\reg}:=f_i-F_i'$
for $i>1$.
We would like to show that
$\det(t-f_i^{\reg})$ is contained in $H(X)[t]$.

\begin{lem}
\label{lem;08.1.14.14}
Let $p_{\gminia,N}$ denote the projection
onto $\psi_N(E_{\gminia})$
with respect to the decomposition
$E=\bigoplus\psi_N(E_{\gminia})$.
Then, we have the estimate
$p_{\gminia,N}-\pi_{\gminia}
=O\bigl(|z_1|^{N/2}\bigr)$
with respect to $h$.
\end{lem}
\pf
We have already known
$\bigl|\pi_{\gminia}-p_{\gminia,N}\bigr|_{|\pi_1^{-1}(P)}=
 O(|z_1|^{N})$ for each $P$
(Proposition \ref{prop;07.10.5.50}).
Let $\delbar_{E,1}$ denote the restriction of
$\delbar_E$ to the $z_1$-direction.
Because the constants in Simpson's Main estimate 
(Section \ref{subsection;07.10.7.1})
depend only on the behaviour of the eigenvalues,
we also have the following estimate on $X-D_1$
with respect to $h$ and 
the Poincar\'e metric of $X-D_1$,
for some $\epsilon>0$:
\[
 \bigl(\delbar_{E,1}
 +\lambda\, f_1^{\dagger}\,
 d\zbar_1\bigr)(\pi_{\gminia}-p_{\gminia,N})
=\lambda\,\bigl[
  f_1^{\dagger}\, d\zbar_1, \,\,
 \pi_{\gminia} \bigr]
=O\bigl(\exp(-\epsilon|z_1|^{-1})\bigr)
\]
Then, the claim follows from 
Lemma \ref{lem;07.7.7.40} below
and the uniform boundedness
$\bigl| R(h_{|\pi_1^{-1}(P)})\bigr|_{h,g_{\poin}}<C$.
\hfill\qed

\begin{cor}
\label{cor;08.1.14.20}
We put
$F_i:=\sum_{\gminia}
 (1+|\lambda|^2)^{-1}
 \del_{i}\gminia\, p_{\gminia,N}$
for $i>1$.
Then,
\[
 F_i-F_i'=O(|z_1|^{N/3}).
\]
\end{cor}
\pf
It follows from Lemma \ref{lem;08.1.14.14}.
\hfill\qed

\vspace{.1in}

Let $\theta_{0}$ denote the section of
$\End(E)\otimes\Omega^{1,0}$
induced by $h_0$ and $\DDlambda$.
We have the expression
$\theta_0=\sum_{i=1}^n f_{0,i}\,dz_i$.
Note the following estimate for $i>1$
with respect to $h$ (Lemma \ref{lem;08.1.14.10}):
\begin{equation}
 \label{eq;08.1.14.16}
 f_{0,i}-F_i
=O(1)
\end{equation}

\begin{lem}
\label{lem;08.1.14.17}
There exists a subset $Z_1\subset D_1$
whose measure is $0$,
such that the following finiteness holds
for any $P\in D_1- Z_1$:
\begin{equation}
 \label{eq;08.1.14.15}
 \int_{\pi_1^{-1}(P)}
 \bigl|f_i-f_{0,i}\bigr|_h^2
 \,\dvol_{g_{\poin}}<\infty
\end{equation}
Here $g_{\poin}$ denote the Poincar\'e metric
of the punctured disc $\pi_1^{-1}(P)^{\ast}$.
\end{lem}
\pf
Let $\DDlambda_i$ denote 
the restriction of $\DDlambda$
to the $z_i$-direction.
We have the equality
\[
  \Delta_i^{\lambda}\tr(s)=
\bigl(s, A_i\bigr)_{h_0}
-\bigl|\DDlambda_i s\cdot s^{-1}\bigr|_{h_0}^2,
\]
for some bounded section $A_i$ 
of $\End(E)$.
(See Section 2.2.5 of \cite{mochi5},
 for example.)
Let $\chi$ be a cut function on a disc.
For any $Q\in D_i\setminus D_1$,
we have the following equality on $\pi_i^{-1}(Q)$:
\[
 \Delta_i^{\lambda} \tr\bigl(\chi \, s\bigr)
=\chi\,(A_i,s)_h
-\chi\,\bigl|\DDlambda_i s\cdot s^{-1/2}\bigr|^2
+2(1+|\lambda|^2)
 \Re\left(\frac{\del \chi}{\del z_i}
 \frac{\del \tr(s)}{\del \zbar_i}\right)
+\bigl(\Delta_i^{\lambda}\chi\bigr) \tr (s)
\]
Therefore, we obtain the following:
\begin{multline*}
 0=\int_{\pi_i^{-1}(Q)}\!\!\chi\,(A_i,s)_h
-\int_{\pi_i^{-1}(Q)}\!\!\!
 \chi\,|\DDlambda_i s\cdot s^{-1/2}|^{2}
 \\
+\int_{\pi_i^{-1}(Q)}\!\!\!\!
 2(1+|\lambda|^2)\Re\left(
 \frac{\del \chi}{\del z_i}\frac{\del \tr s}{\del \zbar_i}
 \right)
+\int_{\pi_i^{-1}(Q)}\!\!(\Delta^{\lambda}_i\chi)\, \tr(s)
\end{multline*}
Thus, we obtain the following inequality
for some constants $B_j$ $(j=1,2,3)$:
\begin{multline*}
 \int_{\pi_i^{-1}(Q)}\!\!
 \chi\,\bigl|\DDlambda_i s\bigr|_h^2
\leq
 B_1\int_{\pi_i^{-1}(Q)}\!\!
 \bigl(\chi+|\Delta_i^{\lambda}\chi|\bigr)
\,|s|
+B_2\int_{\pi_i^{-1}(Q)}
 \!\!|\del_i\chi|\cdot |\DDlambda_i s|_h
 \\
\leq
 B_3+B_2\int_{\pi_i^{-1}(Q)}
 \!\!|\del_i\chi|\cdot |\DDlambda_i s|_h
\end{multline*}
Let $\chi$ be
such that (i) $\chi(r)=1$ for $r\leq r_0/3$,
(ii) $\chi(r)=\exp\bigl(-(r-r_0)^{-1}\bigr)$
for $r_0-\eta<r<r_0$ $(\eta>0)$,
(iii) $\chi(r)=0$ for $r\geq r_0$.
Then, $(\del\chi)\, \chi^{-1/2}$ is bounded.
Hence, we obtain the following inequality
for some constants $B_j$ $(j=4,5)$:
\begin{multline}
 \int_{\pi_i^{-1}(Q)} 
 \chi\,\bigl|\DDlambda_i s\bigr|_h^2
\leq B_3+
B_4\,\left(\int_{\pi_i^{-1}(Q)}|
 (\del_i\chi)\,\chi^{-1/2}|^2\right)^{1/2}
\,
 \left(\int_{\pi_i^{-1}(Q)}
 \chi\,|\DDlambda_i s|_h^2\right)^{1/2} \\
\leq
 B_3+B_5\, \left(\int_{\pi_i^{-1}(Q)}
 \chi\,|\DDlambda_i s|_h^2\right)^{1/2}
\end{multline}
Hence,
we have the uniform boundedness of
the integrals
$\int_{\pi_i^{-1}(Q)}\bigl|\DDlambda_i s\bigr|_h^2$
for $i=2,\ldots,n$
and for $Q\in D_i\setminus D_1$,
when we shrink $X$.
Let $d''$ denote the $(0,1)$-part of $\DDlambda$,
and let $\delta_{h_0}'$ denote the $(1,0)$-operator
induced by $h_0$ and $d''$.
The restriction of $\delta_{h_0}'$ to 
the $z_i$-direction
is denoted by $\delta'_{h_0,i}$.
Since $s$ is self-adjoint with respect to $h_0$,
we obtain the uniform boundedness
of the integrals
$\int_{\pi_i^{-1}(Q)}|\delta'_{h_0,i} s|_h^2$
for $i=2,\ldots,n$
and for $Q\in D_i\setminus D_1$.

Since we have the relation
$\theta-\theta_0=
-(1+|\lambda|^2)^{-1}\, s^{-1}\delta_{h_0}'s$,
we obtain the uniform boundedness of
the integrals
$ \int_{\pi_i^{-1}(Q)}\bigl|f_i-f_{0,i}\bigr|^2 
 \, |dz_i\, d\zbar_i|$
for $i=2,\ldots,n$
and for $Q\in D_i\setminus D_1$.
Then, 
the claim of 
Lemma \ref{lem;08.1.14.17}
follows from Fubini's theorem.
\hfill\qed

\vspace{.1in}

Let us complete the proof of 
Proposition \ref{prop;07.6.12.37}.
We obtain the following finiteness
for any $P\in D_1\setminus Z_1$ and $i>1$,
due to Corollary \ref{cor;08.1.14.20},
the estimate (\ref{eq;08.1.14.16})
and Lemma \ref{lem;08.1.14.17}:
\begin{equation}
 \label{eq;08.1.14.25}
 \int_{\pi_1^{-1}(P)}
 \bigl|f_i^{\reg}\bigr|_h^2
 \,\dvol_{g_{\poin}}<\infty
\end{equation}
Note that
$f_i^{\reg}$ are holomorphic with respect to $\delbar_E$.
Hence, we obtain
the boundedness of
$(f_i^{\reg})_{|\pi_1^{-1}(P)}$ $(i=2,\ldots,n)$
for any $P\in D_1\setminus Z_1$
from (\ref{eq;08.1.14.25}),
by using the norm estimate for wild harmonic bundles
on curves.

By construction,
we have $[\theta,f_i^{\reg}]=0$.
Hence, we have the inequality
$\Delta_1|f_i^{\reg}|_h^2\leq 0$ 
on $\pi_1^{-1}(P)^{\ast}$
for any $P\in D_1$ (Lemma 4.1 of \cite{s2}),
where $\Delta_1=-\del_{z_1}\del_{\zbar_1}$.
Since we have already known the boundedness of
$|f_i^{\reg}|_h^2$ on $\pi_1^{-1}(P)$
for any $P\in D_1-Z_1$,
we obtain $\Delta_1|f_i^{\reg}|_h^2\leq 0$ 
on $\pi_1^{-1}(P)$ as distributions,
for any $P\in D_1-Z_1$.
Hence, the values of $|f_i^{\reg}|_h^2$ 
on $\pi_1^{-1}(P)^{\ast}$
are dominated by the values
on $\del\pi_1^{-1}(P)^{\ast}$
for $P\in D_1-Z_1$.
Then, we obtain the boundedness of
$|f_i^{\reg}|_h^2$ on $X-D$ by using the continuity.
As a result, we obtain
$\det(t-f_i^{\reg})
 \in H(X)[t]$.
Together with Lemma \ref{lem;08.1.14.20},
we can conclude that
$\theta$ is good and wild.
Thus, the proof of Proposition \ref{prop;07.6.12.37}
is finished.
\hfill\qed

\subsection{Proof of 
 Proposition \ref{prop;07.10.13.11}}
\label{subsection;08.9.29.90}

We may and will assume that
$(\vecE_{\ast},\DDlambda)$ is unramified.
We will replace $X$ with a small neighbourhood
of the origin $O$ without mention.
We assume that
the coordinate system is admissible
for $\Irr(\DDlambda):=\Irr(\DDlambda,O)$.
We have the expression
$\theta=\sum_{i=1}^{n} f_i\, dz_i$.
Due to Proposition \ref{prop;07.6.12.37},
we have already known that
$\det(t-z_j\, f_j)\in M(X,D)[t]$ $(j=1,\ldots,\ell)$
and $\det(t-f_j)\in M(X,D)[t]$
$(j=\ell+1,\ldots,n)$.

\begin{lem}
\label{lem;08.9.29.60}
Let $\alpha$ be 
a multi-valued meromorphic function on $X$,
which is a solution of
$\det(t-z_1\, f_1)=0$.
Then, there exists $\gminia\in \Irr(\DDlambda)$
such that 
\[
 \bigl|
 (1+|\lambda|^2)\alpha-z_1\del_1\gminia
\bigr|=O(1).
\]
\end{lem}
\pf
We put $T:=\bigl\{
 (1+|\lambda|^2)^{-1}\, z_1\del_1\gminia\,\big|\,
 \gminia\in \Irr(\DDlambda)
 \bigr\}$.
Due to Proposition \ref{prop;07.6.12.37},
the polynomial $\det(t-z_1\, f_1)$
satisfies the assumption made in 
Section \ref{subsection;07.10.13.10} (Appendix below) 
with the good subset $T$.
Hence, the claim of the lemma follows from
Proposition \ref{prop;07.10.13.1} below.
\hfill\qed

\vspace{.1in}
By Lemma \ref{lem;08.9.29.60},
we obtain the decomposition
$ (E,f_1)=\bigoplus_{\gminia\in \Irr(\DDlambda)}
 (E_{\gminia},f_{1,\gminia})$
such that the eigenvalues of 
\[
 z_1\, f_{1,\gminia}-
 (1+|\lambda|^2)^{-1}z_1\del_1\gminia
 \,\id_{E_{\gminia}}
\]
are bounded.
We have the corresponding decomposition
$f_i=\bigoplus_{\gminia\in\Irr(\DDlambda)} 
 f_{i,\gminia}$.
We set
\[
 f_{i,\gminia}^{\reg}:=
 f_{i,\gminia}
-(1+|\lambda|^2)^{-1}\del_i\gminia
 \,\id_{E_{\gminia}}.
\]

\begin{lem}
\label{lem;08.1.14.30}
$\det(t-z_i\, f^{\reg}_{i,\gminia})\in H(X)[t]$ 
for $i=1,\ldots,\ell$,
and $\det(t-f_{i,\gminia}^{\reg})\in H(X)[t]$ 
for $i=\ell+1,\ldots,n$.
\end{lem}
\pf
Let $1\leq i\leq \ell$
and let $Q\in D_i^{\circ}$.
Let $X_Q$ be a small neighbourhood of
$Q$ such that $D_Q=X_Q\cap D$ is smooth.
Let $\Irr(\DDlambda,i)$  
denote the image of $\Irr(\DDlambda)$
by the map
$M(X,D)/H(X)\lrarr M(X,D)/M(X,D(\neq i))$,
where 
$D(\neq i):=\bigcup_{1\leq j\leq \ell,j\neq i}D_i$.
Due to Proposition \ref{prop;07.6.12.37},
we have the decomposition
$(E,\theta)_{|X_Q}
=\bigoplus_{\gminib\in\Irr(\DDlambda,i)}
 (E_{\gminib,Q},\theta_{\gminib,Q})$
with the following property:
\begin{itemize}
\item
 For the expression 
 $\theta_{\gminib,Q}
 =\sum f_{\gminib,i,Q}\, dz_i$,
 we put 
 $f^{\reg}_{\gminib,i,Q}:=
 f_{\gminib,i,Q}-
 \del_i\gminib\, \id_{E_{\gminib,Q}}$.
 Then, the characteristic polynomials
 of $\det(t-z_if^{\reg}_{\gminib,i,Q})$
 and $\det(t-f^{\reg}_{\gminib,j,Q})$
 $(j\neq i)$
 are contained in $H(X_Q)[t]$.
\end{itemize}
We can observe 
$E_{\gminib,Q}
=\bigoplus_{\gminia\mapsto \gminib}
 E_{\gminia|X_Q-D_Q}$
by the comparison of the eigenvalues
of $z_1 \, f_1$.
Hence, we obtain that
$\det(t-z_i\, f_i^{\reg})$
$(i=1,\ldots,\ell)$
and $\det(t-f_i^{\reg})$
$(i=\ell+1,\ldots,n)$
are contained in $H(X_Q)[t]$.
Then, the claim of Lemma \ref{lem;08.1.14.30}
follows from the Hartogs theorem.
\hfill\qed

\vspace{.1in}

By using Proposition \ref{prop;07.7.19.31},
we can show that the coefficients of
$\det\bigl(t-z_i\, f_{i,\gminia}^{\reg}\bigr)_{|D_i}$
are constants for $i=1,\ldots,\ell$.
Hence, $\harmonicbundle$ is
unramifiedly good wild,
and we have (\ref{eq;10.6.16.1}).
We have the filtered $\lambda$-flat bundles
$\nbigp_{\ast}\nbigelambda$
obtained as in Section \ref{section;07.6.2.1}.
By the assumption,
we have 
$\nbigp_{\ast}
\nbigelambda_{|\pi_i^{-1}(P)}
=\vecE_{\ast|\pi_i^{-1}(P)}$
for $i=1,\ldots,\ell$
and $P\in D_i^{\circ}\setminus Z_i$.
We can deduce 
$\nbigp_{\ast}\nbigelambda
=\vecE_{\ast}$.
Thus, the proof of Proposition \ref{prop;07.10.13.11}
is finished.
\hfill\qed

\subsection{Proof of
 Proposition \ref{prop;07.10.16.12}}
\label{subsection;08.9.29.91}

We may assume $(\vecE_{\ast},\DDlambda)$
is unramified.
Due to Proposition \ref{prop;07.10.13.11},
we have only to show the first claim.
We put $(\vecE_{Q\,\ast},\DDlambda_Q):=
 (\vecE_{\ast},\DDlambda)_{|\pi_i^{-1}(Q)}$
and $E_Q:=E_{|\pi_i^{-1}(Q)}$
for any $Q\in D_i^{\circ}$.
We take an auxiliary metric $h_Q$
for $(\vecE_{Q\,\ast},\DDlambda_Q)$
as in Section \ref{subsection;08.1.14.2}.
We also assume $\det(h_Q)=\det(h_{|\pi_i^{-1}(Q)})$
by an obvious modification.

Let $s_Q$ be determined by
$h_{|\pi_i^{-1}(Q)}=h_Q\, s_Q$.
Because $h$ is adapted to $\vecE_{\ast}$,
we have $|s_Q|_{h_Q}=O\bigl(|z_i|^{-\epsilon}\bigr)$
for any $\epsilon>0$.
Then, we can show $|s_Q|_{h_Q}$ is bounded
by the same argument as that in the proof of
Lemma \ref{lem;07.10.16.11}.
Because $\det(s_Q)=1$,
we also obtain that $s_Q^{-1}$ is bounded
with respect to $h_Q$.
Then, we can conclude that
$h_{|\pi_i^{-1}(Q)}$ is adapted to
$\vecE_{Q,\ast}$,
and the proof of Proposition \ref{prop;07.10.16.12}
is finished.
\hfill\qed

\subsection{Decomposition of polynomials (Appendix)}
\label{subsection;08.9.28.140}

\subsubsection{Statement}

\label{subsection;07.10.13.10}

Let $X:=\Delta^n$, $D_i:=\{z_i=0\}$
and $D:=\bigcup_{i=1}^{\ell}D_i$.
We put $D_i^{\circ}:=D_i\setminus
 \bigcup_{1\leq j\leq \ell,j\neq i}D_j$.
Let $\pi_i$ denote the projection of $X$ to $D_i$.
Let $P(\vecz,t)=\sum_{j=0}^r a_j(\vecz)\, t^j
 \in M(X,D)[t]$ be a monic polynomial 
in the variable $t$ with $M(X,D)$-coefficient.
Assume that we are given a good set of irregular values
$T\subset M(X,D)$ (Definition \ref{df;07.6.1.15})
and subsets
$Z_i\subset D_i$ $(i=1,\ldots,\ell)$
 such that the following holds
for each $i$:
\begin{itemize}
\item
The Lesbegue measure of $Z_i$ is $0$.
\item
 Take any $Q\in D_i^{\circ}\setminus Z_i$.
 Let $\alpha$ be a multi-valued meromorphic function
 on $\pi_i^{-1}(Q)$,
 which is a solution of
 $P(t)_{|\pi_i^{-1}(Q)}=0$.
 Then, there exists $\gminia(\alpha,Q)\in T$
 such that the following holds:
\[
 \bigl|
 \gminia(\alpha,Q)_{|\pi_i^{-1}(Q)}
-\alpha
 \bigr|
=O(1)
\]
In other brief words,
$\alpha$ is decomposed into
a multi-valued holomorphic function
and the polar part
$\gminia(\alpha,Q)_{|\pi_i^{-1}(Q)}$.
\end{itemize}

We will prove the following proposition
in Sections
\ref{subsection;08.9.2.10}--\ref{subsection;08.9.2.11}.

\begin{prop}
\label{prop;07.10.13.1}
We have a splitting
$P(\vecz,t)=\prod_{\gminia\in T}P_{\gminia}(\vecz,t)$
in $M(X,D)[t]$ with the following property:
\begin{itemize}
\item
 Each $P_{\gminia}(\vecz,t)$  is a monic
 in the variable $t$.
\item
Let $\alpha$ be a multi-valued meromorphic section on $X$
which is a solution of $P_{\gminia}(\vecz,t)=0$.
Then, $\alpha-\gminia$ is bounded.
In other words,
$P_{\gminia}(\vecz,t-\gminia)\in H(X)[t]$.
\end{itemize}
\end{prop}

\subsubsection{Preliminary}
\label{subsection;08.9.2.10}

Let $X:=\Delta^n$ and $D:=D_1$.
Let $P\in M(X,D)[t]$ be a monic.
\begin{lem}
\label{lem;08.1.14.12}
Assume that there exists a subset 
$Z\subset D$ whose measure is $0$,
such that the following condition holds:
\begin{itemize}
\item
 Let $Q\in D-Z$.
 For any roots  $\alpha$ of
 $P_{|\pi_1^{-1}(Q)}$,
 we have $|\alpha|=O(1)$.
\end{itemize}
Then, the coefficients of $P$ are holomorphic,
and any roots of $P$ are bounded on $X$.
\end{lem}
\pf
Let $P(\vecz,t)=\sum a_j(\vecz)\, t^j$.
We obtain the boundedness
of $a_{j|\pi_1^{-1}(Q)}$ for any $Q\in D-Z$.
Then, we obtain that $a_j$ are holomorphic on $X$.
\hfill\qed

\vspace{.1in}

Assume that we are given a good set of irregular values
$S=\bigl\{\gminia=\gminia_m\, z_1^m\bigr\}
 \subset M(X,D)$ in the level $m$,
and a subset $Z\subset D$
with the following property:
\begin{itemize}
\item
 The Lesbegue measure of $Z$ is $0$.
\item
 Let $Q\in D-Z$.
 Let $\alpha$ be a root of
 $P_{|\pi_1^{-1}(Q)}$.
 Then, there exists $\gminia(\alpha)\in S$
 such that  the following holds:
\[
  \alpha-\gminia(\alpha)_{|\pi_1^{-1}(Q)}
=O(|z_1|^{m+1})
\]
\end{itemize}
\begin{lem}
\label{lem;07.10.13.2}
We have a decomposition
$P=\prod_{\gminia\in S}P_{\gminia}$
into the monics with the following property:
\begin{itemize}
\item
$\bigl|\alpha-\gminia\bigr|=
 O\bigl(|z_1|^{m+1}\bigr)$
 for any root $\alpha$ of $P_{\gminia}$.
\end{itemize}
\end{lem}
\pf
Let $r:=\deg_tP$.
We set 
$Q(\vecz,t)=z_1^{-r m}
 P(\vecz,z_1^{m}t)$.
By Lemma \ref{lem;08.1.14.12},
we obtain the expansion
$Q(\vecz,t)=
 \sum_{j\geq 0} Q^j(\vecz,t)\, z_1^j$.
We have the decomposition
of $Q^0=\prod_{\gminia\in S}Q^0_{\gminia}$.
It can be lifted to the decomposition
$Q=\prod_{\gminia\in S}Q_{\gminia}$.
By applying Lemma \ref{lem;08.1.14.12}
to $Q_{\gminia}\bigl(z_1^{-1}\, (t-\gminia)\bigr)$,
we obtain that the coefficients of
$Q_{\gminia}\bigl(z_1^{-1}\, (t-\gminia)\bigr)$
are holomorphic on $X$.
Hence, the induced decomposition
$P=\prod_{\gminia\in S}P_{\gminia}$
has the desired property.
\hfill\qed

\subsubsection{Proof of 
Proposition \ref{prop;07.10.13.1}}
\label{subsection;08.9.2.11}

We use an induction on $|T|$.
Let $\vecm(0):=
 \min\bigl\{\ord(\gminia)\,\big|\,\gminia\in T\bigr\}$.
In the case $\bigl|\eta_{\vecm(0)}(T)\bigr|=1$,
we pick an element $\gminia_0\in T$,
and we put $\Ptilde(t)=P(t-\gminia_0)$.
Then, the polynomial $\Ptilde$
satisfies the assumption in 
Section \ref{subsection;07.10.13.10}
with the good set of irregular values
$\Ttilde:=\bigl\{\gminia-\gminia_0\,\big|\,
 \gminia\in T\bigr\}$
and we have
$\min\bigl\{\ord(\gminib)\,\big|\,\gminib\in\Ttilde\bigr\}
\gneq\vecm(0)$.
Hence, we can reduce the problem
to the case $\bigl|\eta_{\vecm}(T)\bigr|\geq 2$,
which we will assume in the following argument.

We put 
$\nbigp:=
 \vecz^{-r\vecm(0)}P(\vecz^{\vecm(0)}t)$,
which is monic in $M(X,D)[t]$.
It can be shown that it is contained
in $H(X)[t]$ by using Lemma \ref{lem;08.1.14.12}
and Hartogs theorem.
Let $k$ be determined by 
$\vecm(0)\in\seisuu_{<0}^k$.
Let $D_{I}:=\bigcap_{i\in I}D_i$
for $I\subset\{1,\ldots,k\}$.
We set
\[
 T_I:=
 \bigl\{(\vecz^{-\vecm(0)}\gminia)_{|D_I}\,\big|\,
 \gminia\in T \bigr\}.
\]
For $I\subset J\subset\{1,\ldots,k\}$,
we have the natural map
$\phi_{I,J}:T_I\lrarr T_J$
by the restriction to $D_{J}$.

Let $Q$ be a point of 
$D_j^{\circ}:=D_j\setminus 
 \bigcup_{i\neq j,1\leq i\leq k} D_i$.
Let $U_Q$ be a small neighbourhood of $Q$
in $X$.
By using Lemma \ref{lem;07.10.13.2},
we obtain the decomposition 
$ \nbigp(t)
=\prod_{\gminib\in T_j}\nbigp_{\gminib,Q}(t)$
on $U_Q$,
such that
$\nbigp_{\gminib,Q}(t)_{|D_j}$ has 
the unique root $\gminib$.
For $\gminic\in T_{\ellsitabar}$,
we put
\[
 \nbigp_{\gminic,Q}(t):=
 \prod_{\phi_{j,\ellsitabar}(\gminib)=\gminic}
 \nbigp_{\gminib,Q}(t).
\]
By varying $Q$ in $D_j^{\circ}$,
we obtain 
the decomposition
$\nbigp=\prod_{\gminic\in T_{\ellsitabar}}
 \nbigp_{\gminic,j}(t)$ around $D_j^{\circ}$.
Because of $\nbigp\in H(X)[t]$,
we have $\nbigp_{\gminic,j}\in H(X)[t]$,
and the decomposition holds on $X$.
Since $\nbigp_{\gminic,j}(t)_{|D_{\ellsitabar}}$
has the unique root $\gminic$,
we obtain that
$\nbigp_{\gminic,j}$ is independent of $j$,
which we denote by $\nbigp_{\gminic}$.

\vspace{.1in}
We put 
$P^{\vecm(0)}_{\gminic}:=
 \vecz^{r\vecm(0)}
 \nbigp_{\gminic}
 \bigl(\vecz^{-\vecm(0)}t\bigr)$
and
$T_{\gminic}:=
 \bigl\{
 \gminia\in T\,\big|\,
 \eta_{\vecm}(\gminia)=\gminic\,\vecz^{\vecm(0)}
 \bigr\}$.
Then, $P^{\vecm(0)}_{\gminic}\in M(X,D)[t]$
has the following property:
\begin{itemize}
\item
Take any $Q\in D_j^{\circ}$.
 Let $\alpha$ be a solution of
 $P^{\vecm(0)}_{\gminic}(t)_{|\pi_j^{-1}(Q)}=0$,
 which is a multi-valued meromorphic function
 on $\pi_j^{-1}(Q)$.
 Then, there exists $\gminia(\alpha,Q)\in T_{\gminic}$
 such that the following holds:
\[
 \bigl|
 \gminia(\alpha,Q)_{|\pi_j^{-1}(Q)}
-\alpha
 \bigr|
=O(1)
\]
\end{itemize}
By applying the hypothesis of the induction,
we obtain the desired decompositions
for $P^{\vecm(0)}_{\gminic}(t)$,
and thus for $P(t)$.
Hence the induction can proceed,
and the proof of Proposition \ref{prop;07.10.13.1}
is finished.
\hfill\qed

\section{The filtered flat bundle associated to
 wild harmonic bundle}
\label{subsection;08.9.29.110}
\subsection{Polystability}

Let $X$ be a connected smooth 
$n$-dimensional projective variety
with an ample line bundle $L$.
Let $D$ be a simple normal crossing hypersurface.
Let $(\vecE_{\ast},\DDlambda)$ be 
a filtered $\lambda$-flat bundle
on $(X,D)$ which is generically good,
i.e., we have a Zariski dense open subset $D'$ of $D$
such that $\vecE_{\ast}$ is good around any point
$P\in D'$.
We say that a pluri-harmonic metric $h$ of $E$ is
generically adapted to $\vecE_{\ast}$,
if it is adapted to $\vecE_{\ast}$
around $P$ for any points $P\in D'$.

\begin{prop}
\label{prop;07.12.6.11}
Let $(\vecE_{\ast},\DDlambda)$
be a filtered $\lambda$-flat sheaf
which is saturated and
generically good.
\begin{itemize}
\item
Assume that we have a pluri-harmonic metric $h$
of $(E,\DDlambda):=(\vecE_{\ast},\DDlambda)_{|X-D}$
which is generically adapted to $\vecE_{\ast}$.
Then, $(\vecE_{\ast},\DDlambda)$ is $\mu_L$-polystable
with $\pardeg_L(\vecE_{\ast})=0$.
The canonical decomposition 
(see Section {\rm\ref{subsection;07.12.6.10}})
is orthogonal with respect to $h$.
The restriction of $h$ to any stable components
of $(\vecE_{\ast},\DDlambda)$ is also pluri-harmonic.
\item
Let $h'$ be another pluri-harmonic metric of $(E,\DDlambda)$
generically adapted to $\vecE_{\ast}$.
Then, we have the decomposition
$(\vecE_{\ast},\DDlambda)
=\bigoplus (\vecE_{i\,\ast},\DDlambda_i)$
such that
(i) it is orthogonal with respect to
both of the metrics $h$ and $h'$,
(ii) $h_{i}=a_i\, h'_i$ for some $a_i>0$,
where $h_i$ and $h'_i$ are the restrictions
of $h$ and $h'$ to $E_i$, respectively.
In particular,
if $(\vecE_{\ast},\DDlambda)$ is $\mu_L$-stable,
$h'=a\, h$ for some $a>0$.
\end{itemize}
\end{prop}
\pf
Let us consider the first claim.
In the one dimensional case,
it can be shown by Simpson's argument
in Section 10 of \cite{s1} and
Section 6 of \cite{s2}.
We give only an outline.
Note that we have already known that
$(E,\DDlambda,h)$ is good wild
according to Proposition \ref{prop;07.6.20.2}.
In particular, $(E,h)$ is acceptable,
which we will implicitly use.
Let $(\vecW_{\ast},\DDlambda)$ be 
a filtered $\lambda$-flat subbundle
of $(\vecE_{\ast},\DDlambda)$.
We put $W:=\vecW_{\ast|X-D}$.
Let $h_W$ be the metric of $W$
induced by $h$.
Let $R(h_W)$ be the curvature 
of $(W,h_W)$.
The analytic degree $\deg^{an}(W)$
is defined to be
$(\sqrt{-1}/2\pi)\int_X\tr R(h_W)$.
According to Lemma 10.5 of \cite{s1}
(see also Lemma 6.2 of \cite{s2}),
it is equal to $\pardeg(\vecW_{\ast})$.
Let $\pi_W$ be the orthogonal projection of $E$
onto $W$.
By the Chern-Weil formula
(see Lemma 2.34 of \cite{mochi5}
for the Chern-Weil formula for $\lambda$-flat bundle),
we have the following formula
\[
\deg^{an}(W)=
 \frac{-1}{2\pi(1+|\lambda|^2)}
 \int_{X-D}\bigl| \DDlambda\pi_W\bigr|_h^2
\leq 0
\]
Hence, we obtain that
$(\vecE_{\ast},\DDlambda)$
is semistable.
Moreover, if $\pardeg(\vecW_{\ast})=0$,
we have $\DDlambda\pi_W=0$.
It implies
$(E,\DDlambda,h)$ is decomposed into
$(W,\DDlambda,h_W)
\oplus (W',\DDlambda,h_{W'})$
as harmonic bundles.
Then, we obtain the decomposition
$(\vecE_{\ast},\DDlambda)
=(\vecW_{\ast},\DDlambda)
\oplus
 (\vecW'_{\ast},\DDlambda)$.
Hence, $(\vecE_{\ast},\DDlambda)$
is polystable.

We can reduce the higher dimensional case
to the one dimensional case
as follows.
By considering the restriction of 
$(\vecE_{\ast},\DDlambda)$
to sufficiently ample and generic curves,
we obtain the $\mu_L$-semistability
and $\pardeg_L(\vecE_{\ast})=0$.
Let $\vecW_{\ast}\subsetneq \vecE_{\ast}$ 
be a saturated filtered subsheaf such that
$\pardeg_L(\vecW_{\ast})=\pardeg_L(\vecE_{\ast})$.
Let $\pi_W$ denote the orthogonal projection
to $W$ which is defined outside of the subset
with codimension two.
By considering the restriction to 
the sufficiently ample general curves,
we obtain $\DDlambda\pi_W=0$.
In particular, $\pi_W$ is holomorphic.
By Hartogs theorem,
$\pi_W$ is defined on whole $X-D$,
and $\pi_W^2=\pi_W$.
Hence, we obtain the decomposition
$E=W\oplus W'$,
where $W'=\Ker\pi_W$.
It is orthogonal and flat.
Hence, we obtain the decomposition
of harmonic bundles
$(W,\DDlambda_W,h_W)\oplus 
 (W',\DDlambda_{W'},h_{W'})$.

Let $P$ be a point of $D$
around which $(E,\DDlambda,h)$ is good.
Then, $(W,\DDlambda_W,h_W)$ and 
$(W',\DDlambda_{W'},h_{W'})$
are also good around $P$,
and the decomposition
is prolonged as
$\vecE_{\ast}=\vecW_{\ast}\oplus\vecW'_{\ast}$
around $P$.
By Hartogs property,
we obtain the decomposition 
$\vecE_{\ast}=\vecW_{\ast}\oplus\vecW'_{\ast}$
on whole $X$.
Hence,
the $\mu_L$-polystability follows.
We also obtain that 
the restriction of $h$ to any $\mu_L$-stable components
are pluri-harmonic.
Let $(\vecE_{\ast},\DDlambda)=
 \bigoplus_{i=1}^{\ell}
 (\vecE_{i\ast},\DDlambda_i)$ be the canonical decomposition.
We also have the decomposition
$(\vecE_{\ast},\DDlambda)
=(\vecE_{1\,\ast},\DDlambda_1)
\oplus
 (\vecE_{1\,\ast}^{\bot},\DD_1^{\lambda\bot})$
whose restriction to $X-D$ is 
orthogonal with respect to $h$.
It is easy to derive
$\vecE_{1\,\ast}^{\bot}
=\bigoplus_{i=2}^{\ell}
 \vecE_{i\ast}$.
Hence the orthogonality of the canonical decomposition
is also obtained.

Let us show the second claim.
In the case $\dim X=1$,
$h$ and $h'$ are mutually bounded
due to the norm estimate.
Therefore, the claim can be shown
using the argument in the proof of 
Proposition 5.2 of \cite{mochi4}.
Let us consider the higher dimensional case.
For any point $P\in X-D$,
we have $s_P$ such that
$h'_{|P}=h_{|P}\, s_P$.
By considering the curve case,
we can show that $s$ is flat with respect to $\DDlambda$.
Then the claim follows.
\hfill\qed

\begin{cor}
\label{cor;10.6.19.2}
Let $(\vecE_{\ast},\DDlambda)$ and $h$
be as in Proposition {\rm\ref{prop;07.12.6.11}}.
Assume that
$(\vecE_{\ast},\DDlambda)$
is the tensor product of
$\mu_L$-stable $(\vecE_{0\ast},\DDlambda_0)$
and a vector space $V$.
Then, $h$ is of the form $h_0\otimes g_V$,
where $h_0$ is a pluri-harmonic metric
for $(\vecE_{0\ast},\DDlambda_0)$
as in Proposition {\rm\ref{prop;07.12.6.11}},
and $g_V$ is a metric of $V$.
\end{cor}
\pf
We take an inclusion
$\vecE_{0\ast}\subset \vecE_{\ast}$.
By restricting $h$,
we obtain a pluri-harmonic metric $h_0$
for $(\vecE_{0\ast},\DDlambda_0)$
as in Proposition \ref{prop;07.12.6.11}.
By using the second claim of Proposition
\ref{prop;07.12.6.11},
we obtain that 
$h$ is isomorphic to a direct sum
of copies of $h_0$.
\hfill\qed

\begin{cor}
\label{cor;10.6.19.1}
Let $(\vecE_{\ast},\DDlambda)$ and $h$
be as in Proposition {\rm\ref{prop;07.12.6.11}}.
Let $\vecE'_{\ast}$ be a filtered
$\lambda$-flat subbundle of $\vecE_{\ast}$
such that $\pardeg_L(\vecE'_{\ast})=0$.
Let $E''$ be the orthogonal complement of
$E':=\vecE'_{\ast|X-D}$ in $E$.
Then, $E''$ is naturally extended to
a filtered $\lambda$-flat subbundle
$\vecE''_{\ast}$ of $\vecE_{\ast}$,
and we have
$\vecE_{\ast}=
 \vecE'_{\ast}\oplus\vecE''_{\ast}$.
\end{cor}
\pf
We can deduce this claim
from the orthogonality of the canonical decomposition
and Corollary \ref{cor;10.6.19.2}.
We can also deduce it 
directly from the proof of Proposition
\ref{prop;07.12.6.11}.
\hfill\qed

\subsection{Vanishing of the characteristic numbers}
\label{subsection;08.9.29.24}

Let $(E,\DDlambda,h)$ be 
a good wild harmonic bundle on $X-D$.
We have the associated good 
filtered $\lambda$-flat bundle
$(\vecE_{\ast},\DDlambda)$ on $(X,D)$.

\begin{prop}
\label{prop;07.12.6.12}
We have the vanishing of the characteristic numbers:
\[
\int_X\parch_{2,L}(\vecE_{\ast})=0,
\quad
\int_X\parchern_{1,L}^2(\vecE_{\ast})=0.
\]
\end{prop}
\pf
We can show it by the essentially same argument
as the proof of Proposition 5.3 of \cite{mochi4}.
We give only an outline
with minor simplification.
(i.e., we may simplify Lemma 5.4 of \cite{mochi4}.)

Let $\pi:\Xtilde\lrarr X$ be the blow up
at cross points of $D$.
Let $\Dtilde$ denote the inverse image of $D$.
Let $(\Etilde,\delbar_{\Etilde},\thetatilde,\htilde)
:=\pi^{-1}\harmonicbundle$.
As remarked in Section \ref{subsection;07.11.5.101},
$\nbigp_{\ast}\nbigelambdatilde$ is obtained 
from $\nbigp_{\ast}\nbigelambda$
by the procedure in Section
{\rm \ref{subsection;07.11.5.60}}.
Hence,
we may have only to consider the integrals
over $\Xtilde$ by Lemma \ref{lem;07.11.5.100}.

Let $\htilde_0$ be an ordinary metric
for the parabolic bundle
$\nbigp_0\nbigelambdatilde$
as given in Chapter 4 of \cite{mochi4},
where 
we considered ordinary metrics
for parabolic Higgs bundles
with possibly irrational parabolic weight.
We apply it in the trivial Higgs field case.
Then, we have only to show the following:
\[
 \int_{\Xtilde} \tr\Bigl(R(\htilde_0)^2\Bigr)=0,
\quad
 \int_{X}\tr\Bigl(R(\htilde_0)\Bigr)^2=0.
\]
Let $\htilde_1$ be the hermitian metric
of $\Etilde$ 
which is as in Section \ref{subsection;07.11.5.101}
around the cross points of $\Dtilde$,
and as in Section 4.2.6 of \cite{mochi4}
around $\Dtilde_i$.
Then, the following equality can be shown
by using Lemmas 4.5 and 4.10 of \cite{mochi4}
and Lemma \ref{lem;07.11.5.110}:
\[
 \int_{\Xtilde} \tr\bigl(R(\htilde_0)^2\bigr)
=\int_{\Xtilde}\tr\bigl(R(\htilde_1)^2\bigr),
\quad
 \int_{X}\tr\bigl(R(\htilde_0)\bigr)^2
=\int_{\Xtilde}\tr\bigl(R(\htilde_1)\bigr)^2
\]
We can show the following equalities
by using the argument in the pages 62--63
of \cite{mochi4} and Lemma 5.2 of \cite{s1}:
\[
 \int_{\Xtilde}\tr\bigl(R(\htilde_1)^2\bigr)
=\int_{\Xtilde}\tr\bigl(R(\htilde)^2\bigr)=0,
\quad
  \int_{\Xtilde}\tr\bigl(R(\htilde_1)\bigr)^2
=\int_{\Xtilde}\tr\bigl(R(\htilde)\bigr)^2=0
\]
Thus, we are done.
\hfill\qed

\section{Perturbation}
\label{subsection;07.10.15.80}
\label{subsection;08.1.13.10}

The construction of this section
will be used in Sections 
\ref{subsection;07.11.9.30}
and \ref{subsection;07.11.9.31}.
(We have already used such a perturbation for 
 good filtered flat bundle on curves 
 in Sections 
 \ref{subsection;08.1.14.35}--\ref{subsection;08.1.14.36}.)
It is essentially the same as that given in 
Section 2.1.6 of \cite{mochi5}.

Let $X$ be a complex {\em surface}
with a simple normal crossing divisor 
$D=\bigcup_{i\in \Lambda}D_i$.
Let $(\vecE_{\ast},\nabla)$ be
a good filtered flat bundle on $(X,D)$.
Let $\vecc\in\real^{\Lambda}$ such that
$c_i\not\in\Par(\vecE_{\ast},i)$,
and let $\prolongg{\vecc}{E}$ denote
the $\vecc$-truncation.
We have the induced filtration $\lefttop{i}F$
of $\lefttop{i}E_{|D_i}$ for each $i\in \Lambda$.
We have
$\pi_a:\lefttop{i}F_a\bigl(\prolongg{\vecc}{E}_{|D_i}\bigr)
 \lrarr 
 \lefttop{i}\Gr^F_a\bigl(\prolongg{\vecc}{E}\bigr)$.
We have the endomorphism
$\Res_i(\nabla)$ on
 $\lefttop{i}\Gr^F_a\bigl(\prolongg{\vecc}{E}\bigr)$.
Since the conjugacy classes of $\Res_i(\nabla)_{|P}$
are independent of the choice of
$P\in D_i^{\circ}:=D_i\setminus \bigcup_{j\neq i}D_j$,
the nilpotent part of $\Res_i(\nabla)$ induces
the weight filtration
$W$ of $\lefttop{i}\Gr^F_{a}
 \bigl(\prolongg{\vecc}{E}\bigr)_{|D_i^{\circ}}$.
It is extended to the filtration
of $\Gr^F_{a}\bigl(\prolongg{\vecc}{E}\bigr)$
in the category of vector bundles.
We put
$\lefttop{i}\Ftilde_{(a,k)}:=
 \pi_a^{-1}(W_k)$,
which is the filtration of $\prolongg{\vecc}{E}_{|D_i}$
indexed by 
$\nbigs_i:=\bigl\{
(a,k)\in
 \openclosed{c_i-1}{c_i}\times\seisuu\,\big|\,
 \Gr^W_k\bigl(\lefttop{i}\Gr^F_a(\prolongg{\vecc}{E})\bigr)
\neq 0\bigr\}$
with the lexicographic order.
For $(\veca,\veck)=
 \bigl(a_i,k_i\,\big|\,i\in \Lambda\bigr)
 \in \prod\nbigs_i$,
we put
\[
 \prolongg{(\veca,\veck)}{\Etilde}
=\Ker\Bigl(
 \prolongg{\vecc}{E}\lrarr
 \bigoplus_{i\in \Lambda}
 \prolongg{\vecc}{E}_{|D_i}\big/
 \lefttop{i}\Ftilde_{(a_i,k_i)}
 \bigl(\prolongg{\vecc}{E}_{|D_i}\bigr)
 \Bigr)
\]
Note that 
$\prolongg{(\veca,\veck)}{\Etilde}$ is
a good lattice.
If $(\vecE_{\ast},\nabla)$ is unramifiedly good,
this claim is obvious.
The general case can be reduced
to the unramified case.

\vspace{.1in}

Let $\epsilon$ be a sufficiently small number.
We take an increasing map
$\varphi_i:\nbigs_i\lrarr \real$
such that $|\varphi_i(a,k)-a|\leq C\, \epsilon$
for some $C>0$.
(Since we are interested in the family 
of the filtrations 
$\vecF^{(\epsilon)}$ $(\epsilon>0)$,
this condition makes sense.)
Then, $\lefttop{i}\widetilde{F}$ and $\varphi_i$ give
the $\vecc$-parabolic filtration 
$\vecF^{(\epsilon)}=\bigl(\lefttop{i}F^{(\epsilon)}\,\big|\,
 i\in \Lambda\bigr)$.
Thus, we obtain
a good $\vecc$-parabolic flat bundle
$\bigl(\prolongg{\vecc}{E},\vecF^{(\epsilon)},\nabla\bigr)$
which is called a $\epsilon$-perturbation
of $(\prolongg{\vecc}{E},\vecF,\nabla)$.
By construction,
we have the following convergence in 
the cohomology group $H^{\ast}(X,\real)$:
\[
 \lim_{\epsilon\to 0}
 \parchern_1(\prolongg{\vecc}{E},\vecF^{(\epsilon)})
=\parchern_1(\prolongg{\vecc}{E},\vecF),
\quad\quad
  \lim_{\epsilon\to 0}
 \parch_2(\prolongg{\vecc}{E},\vecF^{(\epsilon)})
=\parch_2(\prolongg{\vecc}{E},\vecF)
\]

The following proposition is standard.
(See Proposition 3.28 of \cite{mochi4}, for example.)
\begin{prop}
Assume that $\bigl(\prolongg{\vecc}{E},\vecF,\DDlambda\bigr)$
is $\mu_L$-stable.
If $\epsilon$ is sufficiently small,
then the $\epsilon$-perturbation
$\bigl(\prolongg{\vecc}{E},
 \vecF^{(\epsilon)},\DDlambda\bigr)$ is also
$\mu_L$-stable.
\hfill\qed
\end{prop}

We will use two kinds of perturbations 
$\varphi_i$ of parabolic weights.

\begin{description}
\item[(I)]
The image of $\varphi_i$ is contained in $\rnum$
for each $i\in \Lambda$
(Section \ref{subsection;08.2.22.1}).
\item[(II)]
 For simplicity, we assume 
$\epsilon=m^{-1}$
 and $0<10\rank E\, 
 \epsilon<\gap(\prolongg{\vecc}{E},\vecF)$.
(See Section 3.1 of \cite{mochi4}
 for $\gap$.)
Let $i\in \Lambda$.
For each $a\in\Par(\prolongg{\vecc}{E},\vecF)$,
we take $a'(\epsilon,i)\in m^{-1}\,\seisuu$
such that 
$|a'(\epsilon,i)-a|<m^{-1}$.
Let $L(\epsilon,i)\in\real$ be determined by the following:
\[
 L(\epsilon,i)\,\rank(E)
:=\sum (a(\epsilon,i)'-a)\,
 \rank\lefttop{i}\Gr^{F}_a(\prolongg{\vecc}{E})
\]
Then, we put
$a(\epsilon,i):=a'(\epsilon,i)-L(\epsilon,i)$
and 
$\varphi(a,k):=
 a(\epsilon,i)
+k\, \epsilon$.
By construction,
we have the following equality:
\[
 \sum_{a,k} \varphi(a,k)\, 
 \rank\bigl(\lefttop{i}\Gr^{\Ftilde}_{a,k}
 (\prolongg{\vecc}{E})\bigr)
=
 \sum_a a\,
\rank \bigl(
 \lefttop{i}\Gr^{F}_a(\prolongg{\vecc}{E})
\bigr)
\]
Hence, we have
$\parchern_1(\prolongg{\vecc}{E},\vecF)
=\parchern_1(\prolongg{\vecc}{E},\vecF^{(\epsilon)})$.
The parabolic structure satisfies
the SPW-condition in
Definition 2.6 of \cite{mochi5}.
Namely,
for each $i$,
we also have some $-1/m<\gamma_i\leq 0$
such that 
$\Par(\prolongg{\vecc}{E},\vecF^{(\epsilon)},i)$
is contained in
$\bigl\{
 c_i+\gamma_i+p/m\,\big|\,
 p\in\seisuu_{\leq 0},\,
 -1<\gamma_i+p/m\leq 0
 \bigr\}$.
\end{description}

\begin{rem}
The construction given in this section
is valid, when the base manifold $X$ is a curve.
\hfill\qed
\end{rem}

\chapter[Construction of initial metric]
 {Construction of Initial Metric
 and Preliminary Correspondence}
\label{section;08.3.9.3}
In this chapter,
we mainly study 
{\em graded semisimple}
good filtered flat bundle.
Almost all the results are
minor generalization of 
those in Chapter 3 of \cite{mochi5}
for $\lambda=1$.
We will often give only outlines
for the proof.

In Sections
\ref{subsection;07.6.12.21}--\ref{subsection;07.6.12.20},
we explain local constructions of ordinary metrics.
We give in Sections \ref{subsection;08.9.30.1}
and \ref{subsection;08.9.30.2}
the estimates which we will use 
in Sections \ref{subsection;08.2.22.1}
and \ref{subsection;08.3.9.1}.
We also explain in Sections \ref{subsection;08.1.15.2}
and \ref{subsection;08.1.15.3}
the induced metrics on divisors,
which will be used in 
Section \ref{subsection;08.2.22.1}.

In Section \ref{subsection;08.2.22.1},
we explain some formulas for
the parabolic Chern character
of good filtered $\lambda$-flat bundles.
Our main purpose is to show
the vanishing of characteristic numbers
for good Deligne-Malgrange filtered bundle
(Corollary \ref{cor;08.9.30.3}),
which will be significant in the proof of
Theorems \ref{thm;07.10.14.60}
and \ref{thm;07.10.15.1}.
Note that this vanishing holds
even if the good Deligne-Malgrange filtered bundle
is not graded semisimple.

In Section \ref{subsection;08.3.9.1},
we show Kobayashi-Hitchin correspondence
for graded semisimple filtered flat bundles
satisfying the SPW-condition in the surface case.
This result will be used in the proof 
of Kobayashi-Hitchin correspondence
for wild harmonic bundles
(Theorem \ref{thm;06.1.23.100}).

It also implies Bogomolov-Gieseker inequality
for $\mu_L$-stable good filtered flat bundles,
which is explained in Section 
\ref{subsection;07.11.9.30}.
Note that the inequality holds
even if the $\mu_L$-stable good filtered flat bundle
satisfies neither graded semisimplicity
nor the SPW-condition.

\section{Around a cross point}
\label{subsection;07.6.12.21}
\subsection{Taking a ramified covering}

Let $X=\Delta^2$, $D_j=\{z_j=0\}$
and $D=D_1\cup D_2$.
In the following argument,
we will replace $X$
by a small neighbourhood of $O=(0,0)$
without mention.
Let $(\vecE_{\ast},\nabla)$ be a good filtered flat bundle
on $(X,D)$ such that 
$\Res_j(\nabla)$ is graded semisimple,
i.e., the induced endomorphisms
on $\lefttop{j}\Gr^F\bigl(\prolong{E}\bigr)$
are semisimple.
Let $\vecc=(c_1,c_2)\in\real^2$
such that
$c_i\not\in\Par(\vecE_{\ast},i)$.
We assume the SPW-condition
for $(\vecE_{\ast},\nabla)$,
i.e.,
there exist
a positive integer $m$
and $\gamma_i\in\real$ with
$-1/m<\gamma_i\leq 0$,
such that 
\[
\Par(\prolongg{\vecc}{E},i)
\subset 
\bigl\{
 c_i+\gamma_i+p/m\,\big|\,
 p\in\seisuu,\,\,
 -1<\gamma_i+p/m\leq 0
 \bigr\}.
\]

Let $d$ be an integer 
divided by $m\,\rank(E)!^3$.
We put $X':=X$, $D_j':=D_j$, $D':=D$
and $O'=(0,0)$.
Let $\varphi_d:(X',D')\lrarr (X,D)$ be the ramified covering
given by
$\varphi_d(z_1,z_2)=(z_1^d,z_2^d)$.
Let $\Gal(X'/X)$ denote the Galois group of $X'/X$.
We have the filtered flat bundle
$(\vecE'_{\ast},\nabla')$ on $(X',D')$
induced by $(\vecE_{\ast},\nabla)$ and $\varphi_d$
as in Section \ref{subsection;07.11.5.60}.
We put $c'_i:=d(c_i+\gamma_i)$
and $\vecc':=(c_1',c_2')$.
By the assumption,
$\Par(\vecE_{\ast}',i)$
is contained in
$\bigl\{
  c_i'+n\,\big|\,n\in\seisuu
 \bigr\}$.

\subsection{Taking an equivariant decomposition}

Since we have assumed that $(\vecE_{\ast},\nabla)$ is good
and that $d$ is divided by $(\rank E!)^3$,
we have the good subset of irregular values
$\Irr(\nabla',O')\in M(X',D')/H(X')$
and the irregular decomposition 
of the completion of $(\prolongg{\vecc'}{E'},\nabla')$ 
at $O'$.
For simplicity,
we assume that the coordinate system 
$(z_1,z_2)$ is admissible
for $\Irr(\nabla',O')$.
Let $\Irr(\nabla',D_1'):=\Irr(\nabla',O')$,
and let $\Irr(\nabla',D_2')$ denote the image of
$\Irr(\nabla',O')$ via the 
naturally defined map
$M(X',D')/H(X')\lrarr  M(X',D')/M(X',D_1')$.
We have the irregular decompositions 
for $j=1,2$:
\[
 (\prolongg{\vecc'}{E'},\nabla)_{|\Dhat_j'}
=\bigoplus_{\gminia\in \Irr(\nabla',D_j')}
 \bigl(
 \prolongg{\vecc'}{E}'_{\gminia,\Dhat_j'},
 \nablahat_{\gminia,\Dhat_j'}
\bigr).
\]
As remarked in Subsection \ref{subsection;10.5.3.40},
we have the decomposition on $\Dhat'$:
\begin{equation}
\prolongg{\vecc'}{E}'_{|\Dhat'}=
 \bigoplus_{\gminib\in\Irr(\nabla',D_2')}
 \prolongg{\vecc'}{E}'_{\gminib,\Dhat'}
\end{equation}
We have the naturally defined $\Gal(X'/X)$-action
on $\Irr(\nabla',D_i')$.
For any $g\in \Gal(X'/X)$,
we have
$g\cdot \prolongg{\vecc'}{E}'_{\gminib}
 =\prolongg{\vecc'}{E}'_{g\cdot \gminib}$.
Due to the graded semisimplicity assumption,
the endomorphisms $\Res_j(\nabla')$ $(j=1,2)$
are semisimple.
We have the eigen decomposition:
\[
  \prolongg{\vecc'}{E}'_{|D'_j}
=\bigoplus_{\alpha\in\cnum}\lefttop{j}\EE_{\alpha}
\]

\begin{lem}
\label{lem;07.10.15.30}
We can take a decomposition
\[
 \prolongg{\vecc'}{E}'
=\bigoplus_{\gminia\in\Irr(\nabla',O')}
  \bigoplus_{\vecalpha\in \cnum^2}
 U_{\gminia,\vecalpha}
\]
with the following property:
\begin{itemize}
\item
 For $g\in \Gal(X'/X)$,
 we have
 $g\cdot U_{\gminia,\vecalpha}
 =U_{g\cdot\gminia,\vecalpha}$.
\item
 For $\gminia\in\Irr(\nabla',D_1')$,
 we put $U^{(1)}_{\gminia}:=
 \bigoplus_{\vecalpha} U_{\gminia,\vecalpha}$.
 Then,
 $U^{(1)}_{\gminia|\Dhat_1^{\prime(N)}}
 =\prolongg{\vecc'}{E}'_{
 \gminia,\Dhat'_1|\Dhat_1^{\prime(N)}}$.
\item
 For $\gminib\in\Irr(\nabla',D_2')$,
 we put $U^{(2)}_{\gminib}:=
 \bigoplus_{\gminia\in T(\gminib)}
 U^{(1)}_{\gminia}$,
 where $T(\gminib):=\bigl\{\gminia\in
 \Irr(\nabla',O')\,\big|\,\gminia\mapsto\gminib\bigr\}$.
 Then, $U^{(2)}_{\gminib|\Dhat^{\prime(N)}}
 =\prolongg{\vecc'}{E}'_{
 \gminib,\Dhat'|\Dhat^{\prime(N)}}$.
\item
 For $j=1,2$,
 we have
  $\lefttop{j}\EE_{\alpha}=
 \bigoplus_{\gminia\in \Irr(\nabla',O')}
 \bigoplus_{q_i(\vecalpha)=\alpha}
 U_{\gminia,\vecalpha|D_j'}$,
 where $q_j:\cnum^2\lrarr\cnum$
 denotes the projection onto $j$-th component.
\end{itemize}
\end{lem}
\pf
We have the decomposition
$\Irr(\nabla',D'_2)=\coprod \gminio_i$
into the orbits
of the $\Gal(X'/X)$-action.
We take representatives $\gminib_i$
of $\gminio_i$.
For each $\gminib_i$, let $\Stab(\gminib_i)$ denote
the stabilizer of $\gminib_i$ with respect to
the $\Gal(X'/X)$-action.
We can take a $\Stab(\gminio_i)$-invariant subbundle 
$U^{(2)}_{\gminib_i}\subset \prolongg{\vecc'}{E'}$,
such that $U^{(2)}_{\gminib_i|\Dhat^{\prime(N)}}=
 \prolongg{\vecc'}{E}'_{\gminib_i,\Dhat'|\Dhat^{\prime(N)}}$.
For $g\cdot\gminib_i\in \gminio_i$,
we put 
$U^{(2)}_{g\cdot \gminib_i}
:=g\cdot U^{(2)}_{\gminib_i}$.
Thus, we obtain the decomposition
$\prolongg{\vecc'}{E'}
=\bigoplus_{\gminib\in \Irr(\nabla',O')}
 U^{(2)}_{\gminib}$.
Let $\pi^{(2)}_{\gminib}$ denote the projection
onto $U^{(2)}_{\gminib}$ with respect to the decomposition.

We have the decomposition
$\Irr(\nabla',D'_1)=\coprod \gminip_i$
into the orbits of the $\Gal(X'/X)$-action.
We take representatives $\gminia_i$
of $\gminip_i$.
For each $\gminia_i$, let $\Stab(\gminip_i)$ denote
the stabilizer of $\gminia_i$ with respect to
the $\Gal(X'/X)$-action.
By using the equivariant version of
Lemma \ref{lem;10.6.30.3},
we can take a $\Stab(\gminip_i)$-invariant subbundle 
$U^{(1)}_{\gminia_i}\subset \prolongg{\vecc'}{E'}$
with the following property:
\begin{itemize}
\item
$U^{(1)}_{\gminia_i|\Dhat_1^{\prime(N)}}=
 \prolongg{\vecc'}{E}'_{\gminia_i,\Dhat'_1|
 \Dhat_1^{\prime(N)}}$.
\item
$U^{(1)}_{\gminia_i|D_2'}$ is preserved by
$\Res_2(\nabla')$ and
$\pi^{(2)}_{\gminib|D_2'}$
for any $\gminib\in\Irr(\nabla',D_2')$.
\end{itemize}
We put
$U^{(1)}_{g\cdot\gminia_i}:=
 g\cdot U^{(1)}_{\gminia_i}$
for $g\cdot\gminia_i\in\gminip_i$,
and then we obtain the decomposition
$\prolongg{\vecc'}{E}'=
 \bigoplus_{\gminia\in \Irr(\nabla',O')}
 U^{(1)}_{\gminia}$.
By construction,
the following holds:
\begin{description}
\item[(a)]
$U^{(1)}_{\gminia|\Dhat_1^{\prime(N)}}
=\prolongg{\vecc'}{E}'_{
 \gminia,\Dhat_1'|\Dhat_1^{\prime(N)}}$.
\item[(b)]
 $U^{(1)}_{\gminia|D_2'}$ are preserved by
 $\Res_2(\nabla')$.
\item[(c)]
 The decomposition
$\bigoplus_{\gminia\in T(\gminib)}
 U^{(1)}_{\gminia|D_2'}$
 is preserved by 
 $\pi^{(2)}_{\gminib|D_2'}$.
\end{description}
Since we have
$U^{(2)}_{\gminib|O'}
=\bigoplus_{\gminia\in T(\gminib)}
U^{(1)}_{\gminia|O'}$,
we obtain
$U^{(2)}_{\gminib|D_2'}
=\bigoplus_{\gminia\in T(\gminib)}
U^{(1)}_{\gminia|D_2'}$
from (c).
By making the modification to $U^{(1)}_{\gminia}$
which is trivial modulo $z_1^N\,z_2$,
we obtain the decomposition
$U^{(2)}_{\gminib}
=\bigoplus_{\gminia\in T(\gminib)}
 U^{(1)}_{\gminia}$.
The conditions (a) and (b) 
are satisfied for the modified $U^{(1)}_{\gminia}$.
Let $\pi^{(1)}_{\gminia}$ denote the projection
onto $U^{(1)}_{\gminia}$ with respect to the decomposition.

Since $U^{(1)}_{\gminia|D_j'}$ $(j=1,2)$
are preserved by $\Res_j(\nabla')$,
we have the eigen decomposition
$U^{(1)}_{\gminia|D_j'}
=\bigoplus
 \lefttop{j}\EE_{\alpha}
 \bigl(U^{(1)}_{\gminia|D_j'}\bigr)$,
which are $\Stab(\gminip_i)$-equivariant.
Then, we can take a $\Stab(\gminip_i)$-equivariant
decomposition
$U^{(1)}_{\gminia_i}
=\bigoplus_{\vecalpha\in\cnum^2}
 U_{\gminia_i,\vecalpha}$
such that the following holds:
\[
 \bigoplus_{q_j(\vecalpha)=\alpha}
 U_{\gminia_i,\vecalpha|D_j}
=\lefttop{j}\EE_{\alpha}\bigl(
 U_{\gminia_i,\vecalpha|D_j}
 \bigr)
\]
We put 
$U_{g\cdot\gminia_i,\vecalpha}:=
 g\cdot U_{\gminia_i,\vecalpha}$.
Thus, we obtain the decomposition
$\prolongg{\vecc'}{E}'
=\bigoplus_{\gminia}
 \bigoplus_{\vecalpha}
 U_{\gminia,\vecalpha}$,
which has the desired property,
and the proof of Lemma \ref{lem;07.10.15.30}
is finished.
\hfill\qed

\vspace{.1in}
We consider 
$F=F^{\irr}+F^{\reg}$,
where $F^{\irr}$ and $F^{\reg}$
are as follows:
\[
 F^{\irr}=\!\!\!
 \bigoplus_{\gminia\in \Irr(\nabla',O')}
 \!\!\!
 d\gminia\,\id_{U^{(1)}_{\gminia}},
\,\,\,
 F^{\reg}=\bigoplus_{\substack{
 \gminia\in\Irr(\nabla',O'),\\
 \vecalpha\in\cnum^2}}
 \left(\alpha_1\frac{dz_1}{z_1}
 +\alpha_2\frac{dz_2}{z_2}\right)\,
 \id_{U_{\gminia,\vecalpha}}
\]
We put $\nabla^{(0)}:=\nabla-F$,
which gives 
the holomorphic connection of $\prolongg{\vecc'}{E}'$,
although it is not necessarily flat.

\subsection{Frame and metric}

\begin{lem}
We can take frames 
$\vecv_{\gminia,\vecalpha}
 =\bigl(v_{\gminia,\vecalpha,j}\,\big|\,
 j=1,\ldots,\rank U_{\gminia,\vecalpha}
 \bigr)$
of $U_{\gminia,\vecalpha}$
such that
$g\cdot v_{\gminia,\vecalpha,j}
=\omega(g,\gminia,\vecalpha,j)
 \,v_{g\cdot\gminia,\vecalpha,j}$
for some $\omega(g,\gminia,\vecalpha,j)\in\cnum$
with $|\omega(g,\gminia,\vecalpha,j)|=1$.
\end{lem}
\pf
For each $\gminia_i$,
we can take a frame
$\vecv_{\gminia_i,\vecalpha}
=(v_{\gminia_i,\vecalpha,j})$
of $U_{\gminia,\vecalpha}$
such that 
$g\cdot v_{\gminia_i,\vecalpha,j}
=\omega(g,\gminia_i,\vecalpha,j)
 \,v_{\gminia_i,\vecalpha,j}$
for any $g\in \Stab(\gminip_i)$.
For $\gminia\in \gminip_i$,
we take $g(\gminia,\gminia_i)\in G$ such that 
$g(\gminia,\gminia_i)\cdot \gminia_i=\gminia$,
and we put 
$v_{\gminia,\vecalpha,j}:=
 g(\gminia,\gminia_i)\cdot v_{\gminia_i,\vecalpha,j}$.
Thus, we obtain frames
$\vecv_{\gminia,\vecalpha}$
which have the desired property
by construction.
\hfill\qed

\vspace{.1in}

We obtain the frame
$\vecv=\bigl(\vecv_{\gminia,\vecalpha}\bigr)$
of $\prolongg{\vecc'}{E'}$.
Let $h'_0$ denote the metric
given by 
$h'_0(v_i,v_j)=
 |z_1|^{-2c_1'}|z_2|^{-2c_2'}\,\delta_{i,j}$.
Since it is $\Gal(X'/X)$-equivariant,
it induces a hermitian metric $h_0$ of $E$,
which is adapted to $\vecE_{\ast}$.
We have the vanishing of the curvature $R(h_0)=0$.

\subsection{Estimate of the connection one form
 with respect to the frame}

Let $A$ be the connection one form of $\nabla^{(0)}$
with respect to $\vecv$,
i.e., $\nabla^{(0)}\vecv=\vecv\, A$.
Corresponding to the decomposition
$\prolongg{\vecc'}{E}'=
 \bigoplus_{\gminia\in \Irr(\nabla',D'_1)}
 U^{(1)}_{\gminia}$,
we have the decomposition
$\vecv=(\vecv^{(1)}_{\gminia})$.
Correspondingly,
we have the decomposition
$A=\sum A_{\gminia,\gminib}$.
By our choice of the decomposition in Lemma 
\ref{lem;07.10.15.30},
the following holds:
\begin{itemize}
\item
We have
 $A_{\gminia,\gminib}\equiv 0$ modulo 
 $(z_1\,z_2)^N$
 in the case $\ord(\gminia-\gminib)<(0,0)$.
\item
 If $\ord(\gminia-\gminib)=(j,0)$ for some $j<0$,
 $A_{\gminia,\gminib}$ are holomorphic and
 $A_{\gminia,\gminib}\equiv 0$
 modulo $z_1^N$.
\end{itemize}

In the case $\ord(\gminia-\gminib)=(j,0)$,
we have some refined estimate
for the expression $A_{\gminia,\gminib}=
 A_{\gminia,\gminib,1}\, dz_1
+A_{\gminia,\gminib,2}\,dz_2$:
\begin{itemize}
\item
 We may and will assume
 $A_{\gminia,\gminib,2}=O\bigl(z_1^N\,z_2\bigr)$
 after taking some more ramified covering.
\item
 For the decomposition
 $A_{\gminia,\gminib,1}=
 \sum A_{\gminia,\gminib,1,\vecalpha,\vecbeta}$,
 we have
 $A_{\gminia,\gminib,1,\vecalpha,\vecbeta}
 =O\bigl(z_1^N\, z_2\bigr)$
 if $\alpha_2\neq \beta_2$,
 because the eigen decomposition 
 $\prolongg{\vecc'}{E}'_{\gminib,\Dhat'_2|D_2'}
 =\bigoplus \EE_{\alpha}
 \bigl(\prolongg{\vecc'}{E}'_{\gminib,\Dhat'_2|D_2'}
 \bigr)$ is flat with respect to 
 the action of  $\nabla^{(0)}(\del_1)$.
\end{itemize}

Let us look at the term
$A_{\gminia,\gminia}=
 \sum A_{\gminia,\gminia,\vecalpha,\vecbeta}$.
We have the following estimate
in the case $\vecalpha\neq\vecbeta$
for the expression
$A_{\gminia,\gminia,\vecalpha,\vecbeta}
=A_{\gminia,\gminia,\vecalpha,\vecbeta,1}\, dz_1
+A_{\gminia,\gminia,\vecalpha,\vecbeta,2}\, dz_2$:
\begin{itemize}
\item
 We may assume
 $A_{\gminia,\gminia,\vecalpha,\vecbeta,i}
 =O\bigl(z_i\bigr)$
 after taking one more ramified covering.
\item
 If $\alpha_1\neq\beta_1$,
 we have $A_{\gminia,\gminia,\vecalpha,\vecbeta,2}
 =O\bigl(z_1\, z_2\bigr)$.
If $\alpha_2\neq\beta_2$,
 we have $A_{\gminia,\gminia,\vecalpha,\vecbeta,1}
 =O\bigl(z_1\, z_2\bigr)$.
\end{itemize}

\subsection{Estimate of the associated
 $(1,0)$-forms}

Let $\theta_{\nabla'}$ be the $(1,0)$-form
induced by $\nabla'$ and $h'_0$.
(See \cite{s5}
or Section 2.2 of \cite{mochi5}.)
We use the symbol $\theta_{\nabla^{(0)}}$
in a similar meaning.
We put 
$\tau:=2^{-1}\bigl(c'_1\, dz_1/z_1
 +c_2'\,dz_2/z_2\bigr)$.
Then, 
$\theta_{\nabla^{(0)}}-\tau\,\id$
is $C^{\infty}$ on $X'$.
(See Section 2.2.5 of \cite{mochi5}.)
We have the decomposition
$\theta_{\nabla^{(0)}}=
 \sum_{\gminib\neq\gminib'}
 \theta^{(2)}_{\nabla^{(0)},\gminib,\gminib'}
+\sum_{\gminib}\theta^{(2)}_{\nabla^{(0)},\gminib}$
corresponding to
$\prolongg{\vecc'}{E}'
=\bigoplus_{\gminib\in\Irr(\nabla',D_2')}
 U^{(2)}_{\gminib}$.
We also have the decomposition
$\theta_{\nabla^{(0)}}=
 \sum_{\gminia\neq\gminia'}
 \theta^{(1)}_{\nabla^{(0)},\gminia,\gminia'}
+\sum_{\gminia}\theta^{(1)}_{\nabla^{(0)},\gminia}$
corresponding to
$\prolongg{\vecc'}{E}'
=\bigoplus_{\gminia\in\Irr(\nabla',D_1')}
 U^{(1)}_{\gminia}$.
From the above estimate for $A$,
we have the following estimates:
\[
 \theta_{\nabla^{(0)}}
-\bigoplus_{\gminib\in\Irr(\nabla',D_2')}
 \theta^{(2)}_{\nabla^{(0)},\gminib}
=O(|z_1|^N|z_2|^N)
\]
\[
 \theta_{\nabla^{(0)}}
-\bigoplus_{\gminia\in\Irr(\nabla',D_1')}
 \theta^{(1)}_{\nabla^{(0)},\gminia}
=O\bigl(|z_1|^N\bigr)
\]
(See Section 2.2.2 of \cite{mochi5}
for the relation of $A$ and $\theta_{\nabla^{(0)}}$,
for example.)
We have the decomposition
$\theta_{\nabla^{(0)}}=
 \sum 
 \bigl(\theta_{\nabla^{(0)}}\bigr)_{\vecalpha,\vecbeta}$
corresponding to the decomposition
$\prolongg{\vecc'}{E'}
=\bigoplus U_{\vecalpha}$,
where $U_{\vecalpha}=
 \bigoplus_{\gminia} U_{\gminia,\vecalpha}$.
After taking some more ramified covering,
we may and will assume to have the following:
\begin{itemize}
\item
If $\alpha_1\neq\beta_1$,
we have
$(\theta_{\nabla^{(0)}})_{\vecalpha,\vecbeta}
=C^{\infty}\cdot z_1\,dz_1
+C^{\infty}\cdot z_1\, z_2\,dz_2$.
\item
If $\alpha_2\neq\beta_2$,
we have
$(\theta_{\nabla^{(0)}})_{\vecalpha,\vecbeta}
=C^{\infty}\cdot z_1\,z_2\,dz_1
+C^{\infty}\cdot z_2\, dz_2$.
\item
If $\alpha_i\neq\beta_i$ $(i=1,2)$,
we have
$(\theta_{\nabla^{(0)}})_{\vecalpha,\vecbeta}
=C^{\infty}\cdot z_1\,z_2\,dz_1
+C^{\infty}\cdot z_1\, z_2\, dz_2$.
\end{itemize}

We have
$\theta_{\nabla'}=\theta_{\nabla^{(0)}}+F/2$.
\begin{lem}
\label{lem;07.10.15.51}
$\bigl[\theta_{\nabla'},
 \theta_{\nabla'}^{\dagger}\bigr]$
is bounded with respect to the metric $h_0'$
and the Euclidean metric of $X'$.
\end{lem}
\pf
$\bigl[\theta_{\nabla'},
 \theta_{\nabla'}^{\dagger}\bigr]
=\bigl[\theta_{\nabla^{(0)}},
 \theta_{\nabla^{(0)}}^{\dagger}\bigr]
+\frac{1}{2}
 \bigl[F,\theta_{\nabla^{(0)}}^{\dagger}\bigr]
+\frac{1}{2}
 \bigl[\theta_{\nabla^{(0)}},F^{\dagger}\bigr]$.
The first term is $C^{\infty}$ on $X'$.
For the estimate of the other terms,
let us take an auxiliary sequence $\nbigm=(\vecm)$
for the good set $\Irr(\nabla',O')$
as in Section \ref{subsection;07.11.6.1}.
We use the symbols $\etabar_{\vecm}$
and $\zeta_{\vecm}$ in 
Section \ref{subsection;07.11.12.20}.
Let $\Irrbar(\nabla',\vecm)$ denote the image of
$\Irr(\nabla',O')$ via $\etabar_{\vecm}$.
For $\gminic\in\Irrbar(\nabla',\vecm)$,
let $U^{\vecm}_{\gminic}$ denote the sum of
$U^{(1)}_{\gminia}$ such that 
$\etabar_{\vecm}(\gminia)=\gminic$.
Let $p_{\gminic}^{\vecm}$ denote the projection
onto $U_{\gminic}^{\vecm}$
with respect to the decomposition
$\prolongg{\vecc'}{E}'=
\bigoplus_{\gminic\in\Irrbar(\nabla',\vecm)}
 U_{\gminic}^{\vecm}$.
Then, we put as follows:
\[
 F^{\vecm}=\sum_{\gminib\in\Irrbar(\nabla',\vecm)}
 d\zeta_{\vecm}(\gminib)\,
 p^{\vecm}_{\gminib}
\]
Note 
$F^{\irr}=\sum F^{\vecm}$.
We have
$\bigl[F^{\vecm},
 \theta^{\dagger}_{\nabla^{(0)}}\bigr]
=O\bigl(|z_1|^{N/2}\,|z_2|^{N/2}\bigr)$
in the case $\vecm<\veczero$,
and we have
$\bigl[F^{\vecm},
 \theta^{\dagger}_{\nabla^{(0)}}\bigr]
=O(|z_1|^{N/2})$
in the case $\vecm=(j,0)$ for some $j<0$.
Hence, we obtain that
$\bigl[F^{\irr},\theta_{\nabla^{(0)}}^{\dagger}\bigr]$
is bounded.
We also have the following:
\[
 \bigl[F^{\reg},\theta_{\nabla^{(0)}}^{\dagger}\bigr]
\in C^{\infty}+
C^{\infty}\cdot \frac{\zbar_1}{z_1}
+C^{\infty}\cdot\frac{\zbar_2}{z_2}
\]
Thus the proof of Lemma \ref{lem;07.10.15.51}
is finished.
\hfill\qed

\begin{lem}
\label{lem;07.10.15.50}
$\bigl[\theta_{\nabla'},\theta_{\nabla'}
 \bigr]$ is bounded
with respect to
$h_0'$ and the Euclidean metric of $X'$.
After taking a refined ramified covering,
we have
$\bigl[\theta_{\nabla'},\theta_{\nabla'}\bigr]=
O\bigl(|z_1\, z_2|\bigr)\, dz_1\,dz_2$.
\end{lem}
\pf
The second claim immediately
follows from the first one.
We have
$\bigl[\theta_{\nabla'},\theta_{\nabla'}\bigr]
=\bigl[\theta_{\nabla^{(0)}},\theta_{\nabla^{(0)}}\bigr]
+2\bigl[F,\theta_{\nabla^{(0)}}\bigr]$.
The first term is $C^{\infty}$.

Let us look at the second term
$\bigl[F,\theta_{\nabla^{(0)}}\bigr]$.
By construction,
we have 
$\bigl[F^{\vecm},\theta_{\nabla^{(0)}}\bigr]
=O(|z_1|^{N/2}\,|z_2|^{N/2})\,dz_1\,dz_2$
for $\vecm<\veczero$
and 
$\bigl[F^{\vecm},\theta_{\nabla^{(0)}}\bigr]
=O\bigl(|z_1|^{N/2}\bigr)\,dz_1\,dz_2$
for $\vecm=(j,0)$.
We also have the following:
\[
\left[
 F_1^{\reg}\frac{dz_1}{z_1},
 \theta_{\nabla^{(0)}}
\right]
\in
 C^{\infty}\cdot
 \frac{dz_1}{z_1}\, z_1\, z_2\,dz_2
=C^{\infty}\cdot
 z_2\, dz_1\, dz_2
\]
Similarly,
we have
$\left[
 F_2^{\reg}dz_2/z_2,\theta_{\nabla^{(0)}}
 \right]
\in C^{\infty}\cdot z_1\,dz_1\,dz_2$.
Then, the claim of 
Lemma \ref{lem;07.10.15.50}
follows.
\hfill\qed

\subsection{Estimates which will be used later}
\label{subsection;08.9.30.1}

Let $G(h_0')$ denote the pseudo curvature
associated to $\nabla'$ and $h_0'$.
(See \cite{s2}, \cite{s5} or Section 2.2
 of \cite{mochi5}.)

\begin{lem}
\label{lem;07.6.12.6}
$G(h'_0)$ is bounded
with respect to $h_0'$ and the Euclidean metric
of $X'$.
\end{lem}
\pf
We have
$G(h'_0)=2R(h'_0)
-4\bigl(\delbar^2+\theta_{\nabla'}^2-
 [\theta_{\nabla'}^{\dagger},
 \theta_{\nabla'}]\bigr)$.
Recall $R(h'_0)=0$,
and $\delbar^2$ is adjoint of 
$\theta_{\nabla'}^2$ with respect to $h_0'$.
Hence, 
we obtain the boundedness of $G(h'_0)$
from Lemmas \ref{lem;07.10.15.51}
and \ref{lem;07.10.15.50}.
\hfill\qed

\vspace{.1in}

We have
$\tr\bigl(R(h'_0)^2\bigr)
=4^{-1}\tr\bigl(G(h'_0)^2\bigr)
-4\delbar\tr\bigl(\theta_{\nabla'}^2
 \theta_{\nabla'}^{\dagger}\bigr)$
(Lemma 2.31 of \cite{mochi5}).
In Section \ref{subsection;08.3.9.1},
we would like to compare the integrals
of $\tr\bigl(R(h'_0)^2\bigr)=0$
and $\tr\bigl( G(h'_0)^2\bigr)$.
So, let us look at $\tr\bigl(\theta_{\nabla'}^2
 \theta_{\nabla'}^{\dagger}\bigr)$.
We have the following.
\begin{lem}
\label{lem;07.6.12.7}
$\tr\bigl(\theta_{\nabla'}^2
 \,\theta_{\nabla'}^{\dagger}\bigr)$
is bounded
with respect to $h_0'$ and 
the Euclidean metric of $X'$.
\end{lem}
\pf
The boundedness of
$\tr\bigl(\theta_{\nabla'}^2\,
 \theta^{\dagger}_{\nabla^{(0)}}\bigr)$
follows from
Lemma \ref{lem;07.10.15.50}
and the smoothness of 
$\theta^{\dagger}_{\nabla^{(0)}}
 -\overline{\tau}\,\id$.
Let us estimate
$\tr\bigl(\theta_{\nabla'}^2\, F^{\dagger}\bigr)$.
Let 
$\prolongg{\vecc'}{E'}=\bigoplus U^{\vecm}_{\gminib}$
be the decomposition as in the proof of
Lemma \ref{lem;07.10.15.51}.
We have the corresponding decomposition
$\theta_{\nabla'}=
 \sum_{\gminib\neq \gminib'}
 \theta^{\vecm}_{\gminib,\gminib'}
+\sum_{\gminib} \theta^{\vecm}_{\gminib}$.

In the case $\vecm<\veczero$,
we have
$\theta_{\nabla'}^2
-\sum_{\gminib\in\Irrbar(\nabla',\vecm)}
 \bigl(\theta_{\gminib}^{\vecm}\bigr)^2
=O(|z_1\, z_2|^{N/2})$.
Because
$\tr\bigl(
 (\theta^{\vecm}_{\gminib})^2
 \bigr)=0$,
we have the following estimate:
\[
 \tr\bigl(\theta_{\nabla'}^2\,
 F^{\vecm\dagger}\bigr)
=\sum \tr\bigl((\theta^{\vecm}_{\gminib})^2\bigr)\,
 d\overline{\zeta_{\vecm}(\gminib)}
+O\bigl(
 |z_1\, z_2|^{N/3}
 \bigr)
=O\bigl(|z_1\, z_2|^{N/3}\bigr)
\]
In particular,
$\tr\bigl(\theta_{\nabla'}^2\,
 F^{\vecm\dagger}\bigr)$ is bounded
in the case $\vecm<\veczero$.

In the case $\vecm=(j,0)$,
we have 
$ \bigl(
 \theta_{\nabla'}^2
\bigr)
-\sum_{\gminib}
 (\theta^{\vecm}_{\gminib})^2
=O(|z_1|^{N/2})$,
and hence we obtain the following:
\[
 \tr\Bigl(
 \theta_{\nabla'}^2\, F^{\vecm\dagger}
 \Bigr)
=\sum \tr\bigl((\theta^{\vecm}_{\gminib})^2\bigr)
 \, d\overline{\zeta_{\vecm}(\gminib)}
+O\bigl(|z_1|^{N/3}\bigr)
=O\bigl(|z_1|^{N/3}\bigr)
\]
In particular,
$\tr\bigl(\theta_{\nabla'}^2\,
 F^{\vecm\dagger}\bigr)$ is bounded
in the case $\vecm=(j,0)$.

Recall that we have assumed
$\theta_{\nabla'}^2=
 O\bigl(|z_1\, z_2|\bigr)\, dz_1\,dz_2$
after taking some refined ramified covering,
we obviously have the boundedness of
$\tr\bigl(\theta_{\nabla'}^2\, F^{\reg\dagger}\bigr)$.
Thus, we obtain Lemma \ref{lem;07.6.12.7}.
\hfill\qed

\subsection{The induced metric 
of $\lefttop{i}\Gr^{F,\EE}(\prolong{E})$}
\label{subsection;08.1.15.2}

For simplicity,
we assume $c_i=\gamma_i=0$
in this subsection.
Let $(b,\beta)\in \KMS\bigl(\prolong{E}_{\ast},D_i\bigr)$.
We consider the induced metric
on $\lefttop{i}\Gr^{F,\EE}_{b,\beta}(\prolong{E})$.
We restrict ourselves to
the case $i=1$ for simplicity of description.
The other case can be argued similarly.
We put $D_1^{\circ}:=D_1-\{O\}$.
We fix a positive $C^{\infty}$-function $\rho$
on $X$.
Let $\chi:=\rho\, |z_1|^2$.
Then, $h_0$ and $\chi$ naturally
induce the hermitian metric
$h_{b,\beta}$ of 
$\Gr^{F,\EE}_{b,\beta}(\prolong{E})_{|
 D_1^{\circ}}$ as follows.
Let $v_i$ $(i=1,2)$ be sections of 
$\lefttop{1}\Gr^{F,\EE}_{b,\beta}(E)$.
We take lifts $\vtilde_i$ of $v_i$ to $\prolong{E}$,
i.e.,
$\vtilde_{i|D_1}$ is contained 
in $\lefttop{1}F_{b}$,
and mapped to $v_i$
via the projection 
$\lefttop{1}F_b\lrarr \lefttop{1}\Gr^F_b(E)$.
Then,
it can be shown that
$\bigl(
\chi^{b}\, h_0(\vtilde_1,\vtilde_2)
\bigr)_{|D_1^{\circ}}$ is independent
of the choice of lifts $\vtilde_i$,
which is denoted by
$h_{b,\beta}(v_1,v_2)$.

\begin{lem}
\mbox{{}}\label{lem;08.1.15.6}
\begin{itemize}
\item
We have a holomorphic frame $\vecu$ of
$\lefttop{1}\Gr^{F,\EE}_{b,\beta}(\prolong{E})$ 
on $D_1$
which is compatible with the induced parabolic structure
at $O$,
such that 
$h_{b,\beta}(u_i,u_j)=
 \rho^{b}\,|z_1|^{-2\deg^F(u_i)}\,\delta_{i,j}$.
\item
Let $R(h_{b,\beta})$ denote the curvature of
$\bigl(
 \lefttop{1}\Gr^{F,\EE}_{b,\beta}(\prolong{E}),
 h_{b,\beta}
 \bigr)$.
Then, the following holds on $D_1^{\circ}$:
\begin{equation}
 \label{eq;08.1.15.6}
 \tr\bigl(R(h_{b,\beta})\bigr)
-b\,\rank\Gr^{F,\EE}_{b,\beta}(\prolong{E})
 \,\delbar\del\log\rho=0.
\end{equation}
\end{itemize}
\end{lem}
\pf
The claims are clear from the construction.
\hfill\qed

\section{Around smooth point}
\label{subsection;07.6.12.20}
\subsection{Taking a ramified covering}

Let $X:=\Delta^2$ and $D:=D_1$.
Let $(\vecE_{\ast},\nabla)$ be 
a good filtered flat bundle on $(X,D)$,
which is graded semisimple.
Let $c\in\real$ such that
$c\not\in\Par(\vecE_{\ast})$.
We assume the SPW-condition,
i.e.,
there exist a positive integer $m$ 
and $\gamma\in\real$ with
$-1/m<\gamma\leq 0$,
such that 
\[
\Par(\prolongg{c}{E})
\subset
\bigl\{c+\gamma+p/m\,\big|\,
 p\in\seisuu,\,\,
 -1<\gamma+p/m\leq 0 \bigr\}.
\]

Take an integer $d$ 
which is divided by $m\,\rank (E)!^3$.
Let $\varphi_d:X'\lrarr X$
be given by $\varphi_d(z_1,z_2)=(z_1^d,z_2)$.
We have the filtered flat bundle
$(\vecE_{\ast}',\nabla')$ on $(X',D')$
induced by $(\vecE_{\ast},\nabla)$
and $\varphi_d$.
We put $c':=d(c+\gamma)$.
By the assumption,
$\Par(\vecE'_{\ast})$
is contained in
$\bigl\{c'+n\,\big|\,n\in\seisuu\bigr\}$.
Since $d$ is divided by $\rank(E)!^3$,
$(\vecE'_{\ast},\nabla')$ is unramified,
and we have the irregular decomposition:
\[
 \prolongg{c'}{E'}_{|\Dhat'}
=\bigoplus_{\gminia\in\Irr(\nabla')}
 E'_{\gminia}
\]
We have the $\Gal(X'/X)$-action on $\Irr(\nabla')$,
and $g\cdot E'_{\gminia}=E'_{g\cdot \gminia}$.
From the graded semisimplicity assumption,
$\Res(\nabla')$ is semisimple.
We have the eigen decomposition of $\Res(\nabla')$:
\[
 \prolongg{c'}{E}'_{|D'}
=\bigoplus_{\alpha\in\cnum} \EE_{\alpha}
\]

\subsection{Taking a decomposition}

Let $N$ be sufficiently large.
The following lemma can be shown
using an argument similar to 
that employed in the proof of
Lemma \ref{lem;07.10.15.30}.
\begin{lem}
We can take a decomposition
\[
 \prolongg{c'}{E'}
=\bigoplus_{\gminia\in\Irr(\nabla')}
  \bigoplus_{\alpha\in \cnum}
 U_{\gminia,\alpha}
\]
with the following property:
\begin{itemize}
\item
 For $g\in \Gal(X'/X)$,
 we have
 $g\cdot U_{\gminia,\alpha}
 =U_{g\cdot\gminia,\alpha}$.
\item
 We put $U_{\gminia}:=
 \bigoplus_{\alpha} U_{\gminia,\alpha}$.
 Then,
 $U_{\gminia|\Dhat^{\prime(N)}}
 =\prolongg{c'}{E}'_{\gminia|\Dhat^{\prime(N)}}$.
\item
  $\EE_{\alpha}=
 \bigoplus_{\gminia\in \Irr(\nabla')}
 U_{\gminia,\alpha|D'}$.
\hfill\qed
\end{itemize}
\end{lem}

We consider the following:
\begin{equation}
 \label{eq;07.6.12.25}
 F^{\irr}:=\bigoplus_{\gminia} 
 d\gminia\,\id_{U_{\gminia}},
\quad
 F^{\reg}:=\bigoplus_{\gminia,\alpha}
 \left(\alpha\,\frac{dz_1}{z_1} \right)\,
 \id_{U_{\gminia,\alpha}},
\quad
F:=F^{\reg}+F^{\irr}
\end{equation}
We put $\nabla^{(0)}:=\nabla'-F$,
which is a holomorphic connection of
$\prolongg{c'}{E}'$,
although it is not necessarily flat.

\subsection{Metric}

Take a large number $N$.
Let $\htilde'_0$ be a $\Gal(X'/X)$-invariant
$C^{\infty}$-metric of $\prolongg{c'}{E}'$ 
with the following property:
\begin{itemize}
\item 
There exist $C^{\infty}$-metrics 
$\htilde'_{0,\gminia}$ of $U_{\gminia}$
such that 
$\htilde_0'
\equiv\bigoplus \htilde'_{0,\gminia}$
on the $N$-th infinitesimal neighbourhood
of $D$.
\item
$\htilde'_{0|D}$ is compatible with the eigen decomposition,
i.e., $\prolongg{c'}{E}'_{|D'}=
 \bigoplus \EE_{\alpha}$ is orthogonal
with respect to $\htilde'_{0|D'}$
\end{itemize}
We remark that
$R(h'_0)-
 \bigoplus_{\gminia} R(h'_{0,\gminia})$
is $O\bigl(|z_1|^N\bigr)$.
Let $\rho$ be a positive $C^{\infty}$-function on $X$,
and we put $\chi:=\rho\, |z_1|^2$.
Then, we set
\[
 h'_0:=
 \varphi_d^{\ast}(\chi)^{-(c+\gamma)}
 \,\htilde_0'.
\]
The induced metric of $E$
is denoted by $h_0$.

\subsection{Estimate of the associated
  $(1,0)$-form}

Let $\theta_{\nabla'}$ denote the $(1,0)$-form
obtained from $\nabla'$ and $h_0'$.
We use the symbol $\theta_{\nabla^{(0)}}$
in a similar meaning.
Let $\tau:=
 2^{-1}(c+\gamma)\, 
 \del\log(\varphi_d^{\ast}\chi)$,
and then 
$\theta_{\nabla^{(0)}}
-\tau\,\id$ is $C^{\infty}$ 
on $X'$.
We have the decomposition
$\theta_{\nabla^{(0)}}
=\sum_{\gminib}\theta_{\nabla^{(0)},\gminib}
+\sum_{\gminib\neq\gminib'}
 \theta_{\nabla^{(0)},\gminib,\gminib'}$,
corresponding to 
$\prolongg{c'}{E'}=
 \bigoplus_{\gminia}U_{\gminia}$.
By construction, we have the following estimate:
\begin{equation}
 \label{eq;08.2.23.21}
 \theta_{\nabla^{(0)}}
-\bigoplus \theta_{\nabla^{(0)},\gminib}
=O(|z_1|^N)
\end{equation}
We also know that
$\theta_{\nabla^{(0)}|D'}$ is compatible
with the residue $\Res(\nabla')$.
Let $\theta_{\nabla^{(0)}}=
 \sum_{\alpha,\beta}\theta_{\nabla^{(0)},\alpha,\beta}$
be the decomposition corresponding to
$\prolongg{c'}{E'}=\bigoplus_{\alpha}\bigl(
 \bigoplus_{\gminia}U_{\gminia,\alpha}\bigr)$.
After taking a refined ramified covering,
we may assume the following:
\begin{equation}
 \label{eq;07.6.12.5}
 \theta_{\nabla^{(0)}}
-\tau\,\id
=O\bigl(|z_1|)\, dz_1+O(1)\, dz_2,
\quad
 \theta_{\nabla^{(0)},\alpha,\beta}
=O\bigl(|z_1|\bigr)\, dz_1
+O\bigl(|z_1|\bigr)\, dz_2,
\,\,(\alpha\neq \beta)
\end{equation}

\begin{lem}
\label{lem;07.10.15.60}
After taking a refined ramified covering,
we have the boundedness of 
$\bigl[\theta_{\nabla'},
 \theta_{\nabla'}^{\dagger}\bigr]$
with respect to $h_0'$
and the Euclidean metric of $X'$.
\end{lem}
\pf
We have 
$\theta_{\nabla'}=\theta_{\nabla^{(0)}}+F/2$,
and hence
$ \bigl[\theta_{\nabla'},
 \theta_{\nabla'}^{\dagger}\bigr]
=\bigl[\theta_{\nabla^{(0)}},\,
 \theta_{\nabla^{(0)}}^{\dagger}\bigr]
+\frac{1}{2}
 \bigl[F,\,\theta_{\nabla^{(0)}}^{\dagger}\bigr]
+\frac{1}{2}
 \bigl[\theta_{\nabla^{(0)}},\,F^{\dagger}\bigr]
+\frac{1}{4}
 \bigl[F,\,F^{\dagger}\bigr]$.
The first term is $C^{\infty}$ on $X'$.
We have
$\bigl[F^{\irr},\theta^{\dagger}_{\nabla^{(0)}}\bigr]
=O\bigl(|z_1|^N\bigr)$
by (\ref{eq;08.2.23.21}).
We also have the boundedness of
$\bigl[F^{\reg},\theta_{\nabla^{(0)}}^{\dagger}\bigr]$
because of (\ref{eq;07.6.12.5}).
We have 
$\bigl[F^{\irr},(F^{\irr})^{\dagger}\bigr]
=O\bigl(|z_1|^N\bigr)$,
$\bigl[F^{\reg},(F^{\irr})^{\dagger}\bigr]
=O\bigl(|z_1|^N\bigr)$.
We also have
$\bigl[
 F^{\reg},(F^{\reg})^{\dagger}
 \bigr]=O(|z_1|^{-1})\, dz_1\,d\zbar_1$.
Hence, we obtain the desired boundedness of
 $[\theta_{\nabla'},\theta_{\nabla'}^{\dagger}]$
after taking some refined ramified covering.
\hfill\qed

\begin{lem}
\label{lem;07.10.15.61}
$\bigl[\theta_{\nabla'},\theta_{\nabla'}\bigr]$
is bounded with respect to
$h_0'$ and the Euclidean metric of $X'$.
We have
$\bigl[\theta_{\nabla'},\theta_{\nabla'}\bigr]
=O(|z_1|)\, dz_1\,dz_2$
after taking some refined ramified covering.
\end{lem}
\pf
We have
$\bigl[\theta_{\nabla'},\theta_{\nabla'}\bigr]
=\bigl[\theta_{\nabla^{(0)}},\theta_{\nabla^{(0)}}\bigr]
+\bigl[F,\theta_{\nabla^{(0)}}\bigr]$.
The first term is $C^{\infty}$,
and we may assume that it is
$O(|z_1|)\, dz_1\,dz_2$,
after taking some more ramified covering.
We have
$\bigl[F^{\irr},\theta_{\nabla^{(0)}}\bigr]
 =O\bigl(|z_1|^N\bigr)\, dz_1\, dz_2$.
From the compatibility of $\theta_{\nabla^{(0)}|D}$
and the residue,
we obtain
$\bigl[F^{\reg},\theta_{\nabla^{(0)}}\bigr]
=O(1)\, dz_1\,dz_2$.
After taking some more ramified covering,
we may have 
$O(|z_1|)\, dz_1\, dz_2$.
Thus, we are done.
\hfill\qed

\subsection{Estimate which will be used}
\label{subsection;08.9.30.2}

\begin{lem}
\label{lem;07.6.12.8}
$G(h'_0)$ is bounded
with respect to $h_0'$ and
the Euclidean metric of $X'$.
\end{lem}
\pf
It follows from Lemmas \ref{lem;07.10.15.60}
and \ref{lem;07.10.15.61}
with an argument 
similar to that employed in
the proof of Lemma \ref{lem;07.6.12.6}.
\hfill\qed

\vspace{.1in}

We have
$\tr\bigl(R(h'_0)^2\bigr)
=4^{-1}\tr\bigl(G(h'_0)^2\bigr)
 -4\delbar\tr\bigl(\theta_{\nabla'}^2
 \theta_{\nabla'}^{\dagger}\bigr)$.
Let us see the second term.
\begin{lem}
\label{lem;07.6.12.11}
$\tr\bigl(\theta_{\nabla'}^2
 \theta_{\nabla'}^{\dagger}\bigr)$
is bounded
with respect to $h_0'$ and 
the Euclidean metric of $X'$.
\end{lem}
\pf
The boundedness of
$\tr\bigl(\theta_{\nabla'}^2
 \theta_{\nabla^{(0)}}^{\dagger}\bigr)$
follows from Lemma \ref{lem;07.10.15.61}
and the smoothness of 
$\theta_{\nabla^{(0)}}^{\dagger}
-\taubar\,\id$.
Let $\eta_j$ and $\zeta_j$ be the maps
$M(X',D')\lrarr M(X',D')$
given by $\eta_j\bigl(\sum a_k(z_2)z_1^k\bigr)
=\sum_{k\leq j}a_k(z_2)z_1^k$
and $\zeta_j\bigl(\sum a_k(z_2) z_2^k\bigr)
=a_j(z_2) z_1^j$.
Let $\Irr(\nabla',j)$ denote the image of
$\Irr(\nabla')$ by $\eta_{j}:M(X,D)\lrarr M(X,D)$.
For $\gminic\in\Irr(\nabla',j)$,
let $U^{(j)}_{\gminic}$ denote the sum of
$U_{\gminia}$ such that $\eta_{j}(\gminia)=\gminic$.
Let $p_{\gminic}^{(j)}$ denote the projection
onto $U_{\gminic}^{(j)}$
with respect to the decomposition
$\prolongg{c'}{E}'=\bigoplus_{\gminic\in\Irr(\nabla',j)}
 U_{\gminic}^{(j)}$.
We have the corresponding decomposition
$\theta_{\nabla'}
=\sum \theta^{(j)}_{\gminic}$.
We set
\[
 F^{(j)}:=\sum_{\gminib\in\Irr(\nabla',j)}
 d\zeta_{j}(\gminib)\,
 p^{(j)}_{\gminib}
\]
We have the decomposition
$F^{\irr}=\sum F^{(j)}$.
We obtain the following estimate
as in the proof of Lemma \ref{lem;07.6.12.7}:
\[
  \tr\bigl(
 \theta_{\nabla'}^2\, F^{(j)\dagger}
 \bigr)
=
 \sum \tr\bigl(
 (\theta^{(j)}_{\gminib})^2
 \bigr)\, d\overline{\zeta_{j}(\gminib)}
+O(|z_1|^{N/2})
=O(|z_1|^{N/2})
\]
Hence, 
$\tr\bigl(\theta_{\nabla'}^2\, F^{\irr\dagger}\bigr)$
is bounded.
The boundedness of
$\tr\bigl(\theta_{\nabla'}^2\, F^{\reg\dagger}\bigr)$
follows from
$\theta_{\nabla'}^2=O(|z_1|)\, dz_1\, dz_2$.
\hfill\qed

\subsection{The induced metric
of $\Gr^{F,\EE}_{b,\beta}(E)$}
\label{subsection;08.1.15.3}

For simplicity,
we assume $c=\gamma=0$
in this subsection.
We remark 
that the eigenvalues of
$\Res(\nabla')$ are given by
$d\bigl(b+\beta\bigr)$,
where 
 $(b,\beta)\in\KMS\bigl(\prolong{E}_{\ast}\bigr)$.
We have the decomposition
$\EE_{\alpha}
=\bigoplus_{0\leq q\leq d-1} \EE_{\alpha}^{q}$
such that 
$\Gal(X'/X)\ni \omega$ acts on $\EE_{\alpha}^q$
by the multiplication of $\omega^q$.
The decomposition
$\EE_{\alpha}=\bigoplus \EE_{\alpha}^q$
is orthogonal with respect to $\htilde_0'$.

\vspace{.1in}

Let $(b,\beta)\in\KMS(\prolong{E}_{\ast})$.
Then, $h_0$ and $\chi$ naturally
induce a hermitian metric
$h_{b,\beta}$ of 
$\Gr^{F,\EE}_{b,\beta}(E)$ as follows.
Let $v_i$ $(i=1,2)$ be sections of 
$\Gr^{F,\EE}_{b,\beta}(E)$.
We take lifts $\vtilde_i$ of $v_i$
to $\prolong{E}$,
i.e.,
$\vtilde_{i|D}$ is contained in $F_{b}$,
and mapped to $v_i$
via the projection $F_b\lrarr \Gr^F_b(E)$.
Then,
it can be shown that
$\bigl(
\chi^{b}\, h_0(\vtilde_1,\vtilde_2)
\bigr)_{|D}$ is independent
of the choice of $\vtilde_i$,
which is denoted by
$h_{b,\beta}(v_1,v_2)$.

Let $(q,\alpha)$ and $(b,\beta)$
be related by the relation
$b=-q/d$ and $\beta=(\alpha+q)/d$.
Let $h'_{b,\beta}$ denote 
the metric of $\EE_{\alpha}^q$
induced by $\htilde'_0$.

\begin{lem}
\label{lem;08.1.15.5}
Let $R(h_{b,\beta})$
and $R(h'_{b,\beta})$ be the curvatures
of $(\Gr^{F,\EE}_{b,\beta}(E),h_{b,\beta})$
and $(\EE_{\alpha}^q,h'_{b,\beta})$
respectively.
Then, we have the following relation:
\begin{equation}
 \label{eq;08.1.15.1}
 \tr\bigl(R(h'_{b,\beta})\bigr)
=\tr\bigl(R(h_{b,\beta}) \bigr)
-b\,\rank\Gr^{F,\EE}_{b,\beta}(E)\,
 \delbar\del\log\rho
\end{equation}
\end{lem}
\pf
We take the isomorphism
$\Phi:\Gr^{F,\EE}_{b,\beta}(E)
 \simeq
 \EE_{\alpha}^q$
given as follows.
Let $v$ be a section of $\Gr^{F,\EE}_{b,\beta}(E)$.
We take a lift $\vtilde$ of $v$ to $\prolong{E}$.
Then, 
$\Phi(v):=(z_1^{-q}\varphi_d^{\ast}\vtilde)_{|D'}$
is contained in $\EE_{\alpha}^q$,
and it is independent of the choice of $\vtilde$.
Under the isomorphism $\Phi$,
we have
$h'_{b,\beta}=
 h_{b,\beta}\, \rho^{-b}$.
Then, (\ref{eq;08.1.15.1})
follows from a general formula.
\hfill\qed

\section{Some formulas
 for the parabolic Chern character}
\label{subsection;08.2.22.1}
Let $X$ be a smooth projective complex surface,
and $D$ be a simple normal crossing hypersurface
with the irreducible decomposition 
$D=\bigcup_{i\in \Lambda}D_i$.
Let $(\vecE_{\ast},\nabla)$ be 
a good filtered flat bundle on $(X,D)$.
We have the bundles
$\lefttop{i}\Gr^{F,\EE}_{a,\alpha}(\prolongg{\vecc}{E})$
on $D_i$
for each $\vecc\in\real^{\Lambda}$ and
for each $(a,\alpha)\in \KMS\bigl(
 \prolongg{\vecc}{E}_{\ast},D_i\bigr)$.
It is naturally equipped with 
the $\vecc_i$-parabolic structure
at $D_i\cap\bigcup_{j\neq i}D_j$,
where $\vecc_i:=\bigl(c_j\,\big|\,j\in \Lambda,j\neq i\bigr)$.

\begin{prop}
\label{prop;08.2.23.31}
We have the following formula:
\begin{multline}
\int_X
 2\parch_{2}(\vecE_{\ast})=\\
 \sum_{\substack{i\in \Lambda\\
 (a,\alpha)\in \KMS(\prolongg{\vecc}{E}_{\ast},D_i)}}
 \!\!\!\!\!\!\!\!\!\!\!
\bigl(a+\Re(\alpha)\bigr)
\Bigl(
 -\pardeg_{D_i}\bigl(
 \lefttop{i}\Gr_{a,\alpha}^{F,\EE}(
 \prolongg{\vecc}{E})_{\ast}\bigr)
+a\,
 \rank\bigl(\lefttop{i}\Gr_{a,\alpha}^{F,\EE}
 (\prolongg{\vecc}{E})\bigr)\,[D_i]^2
\Bigr)
\end{multline}
Here $[D_i]^2$ denotes the self intersection
number of $D_i$.
\end{prop}
\pf
Note that the right hand side
is independent of the choice of 
$\vecc\in\real^{\Lambda}$.
By using $\epsilon$-perturbations
explained in (I) of Section \ref{subsection;07.10.15.80},
we can reduce the problem to the case
in which the following conditions are satisfied:
\begin{itemize}
\item
$\Par(\vecE_{\ast},i)\subset\rnum$,
$0\not\in \Par(\vecE_{\ast},i)$,
and $(\vecE_{\ast},\nabla)$ is graded semisimple.
\end{itemize}
So, we will assume them in the following argument.
We may assume $\vecc=(0,\ldots,0)$.

We take $d$ such that
(i) $d\, w\in\seisuu$ for any
$w\in\Par(\vecE_{\ast},D_i)$,
(ii) $d$ is divided by $(\rank E !)^3$.
We take a $C^{\infty}$-metric $h_0$
of $E$ on $X-D$ which is
as in Section \ref{subsection;07.6.12.21}
 around cross points $D_i\cap D_j$,
and as in Section \ref{subsection;07.6.12.20}
around smooth points of $D$.

We take a hermitian metrics $g_{i}$ 
of $\nbigo(D_i)$. 
The curvature of $(\nbigo(D_i),g_i)$ is denoted by 
$\omega_i$.
Let $\sigma_i:\nbigo\lrarr\nbigo(D_i)$ 
denote the canonical morphism.
By using the functions $|\sigma_i|_{g_i}^2$,
we obtain the $C^{\infty}$-metrics
$h_{a,\alpha}$ of
$\lefttop{i}\Gr^{F,\EE}_{a,\alpha}
 (\prolong{E})_{|D_i^{\circ}}$
for each $(a,\alpha)\in\KMS\bigl(\prolong{E},D_i\bigr)$,
as explained in Sections \ref{subsection;08.1.15.2}
and \ref{subsection;08.1.15.3}.
It is compatible with the induced parabolic structure
of $\Gr^{F,\EE}_{a,\alpha}(\prolong{E})$.
Hence, we have the following:
\begin{equation}
 \label{eq;08.1.15.7}
 \frac{\sqrt{-1}}{2\pi}\int_{D_i}
 \tr\bigl(R(h_{a,\alpha})\bigr)
=\pardeg_{D_i}\bigl(
 \lefttop{i}\Gr^{F,\EE}_{a,\alpha}(\prolong{E})_{\ast}
 \bigr)
\end{equation}

Let $P$ be a point of smooth part of $D$.
We take a holomorphic coordinate $(U,z_1,z_2)$
around $P$ such that $D_{U}:=D\cap U=\{z_1=0\}$.
Let $\varphi_P:U'\lrarr U$ be a ramified covering
as in Section \ref{subsection;07.6.12.20}.
We have the irregular decomposition
$\prolong{E}'_{|\Dhat'_{U}}
=\bigoplus_{\gminia\in \Irr(\nabla,P)}
  E'_{\gminia}$.
We put
$\Ghat^{\irr}_P:=
 \bigoplus_{\gminia\in\Irr(\nabla)}
 \gminia\,\id_{E'_{\gminia}}$.
We can take a $C^{\infty}$-equivariant section
$G^{\irr}_P$ of $\End(\prolong{E'})$
such that $G^{\irr}_{P|\Dhat^{\prime(N)}_U}
=\Ghat^{\irr}_{P|\Dhat^{\prime(N)}_U}$
for some large $N$.
It induces an endomorphism of $E_{|U\setminus D}$.
For any cross point $P\in D$,
we take a holomorphic coordinate
$(U,z_1,z_2)$ around $P$
such that
$U\cap D=\{z_1=0\}\cup \{z_2=0\}$.
We take  $\varphi_P:U'\lrarr U$
and a decomposition
$\prolong{E}'=\bigoplus U_{\gminia}$
as in Section \ref{subsection;07.6.12.21}.
We put $G^{\irr}_P:=\bigoplus
 \gminia\, \id_{U_{\gminia}}$,
which is $\Gal(U'/U)$-equivariant.
Hence, it induces an endomorphism of 
$E_{|U\setminus D}$.
By varying $P\in D$, gluing $G^{\irr}_P$ in $C^{\infty}$,
and extending it,
we construct an endomorphism $G^{\irr}$ of $E$
on $X-D$.

Let $\nabla^u$ denote the unitary connection
associated to $h_0$ and the $(0,1)$-part $d''$ of $\nabla$.
We obtain the $C^{\infty}$-section
$\nbigf^{\irr}:=\nabla^uG^{\irr}$
of $\End(E)\otimes\Omega^{0,1}$ on $X-D$.
By construction,
we have
$\nbigf^{\irr}_{|\Dhat^{\prime(N)}_U} 
=F^{\irr}_{|\Dhat^{\prime(N)}_U}$ 
around any smooth point $P$ of $D$
(Section \ref{subsection;07.6.12.20}).
We also have
$\nbigf^{\irr}=F^{\irr}$ 
around any cross point of $D$
(Section \ref{subsection;07.6.12.21}).

Let $\sigma_i:\nbigo_X\lrarr\nbigo_X(D_i)$ denote
the canonical section.
The norm of $\sigma_i$ with respect to
a chosen hermitian metric is denoted by 
$|\sigma_i|$.
Let $X_{\delta}:=\bigcap\bigl\{|\sigma_i|\geq\delta\bigr\}$.
We remark $R(h_0)=0$ around the cross point of $D$.
Hence, we have the following:
\[
 \int_{\del X_{\delta}}
 \tr\bigl(\nbigf^{\irr}\,R(h_0)\bigr)
=\int_{\del X_{\delta}}
 \tr\bigl(
 \nabla^u
 G^{\irr}\, R(h_0)
 \bigr)
=\int_{\del X_{\delta}}
 d\tr\bigl(G^{\irr}\, R(h_0)\bigr)
=0.
\]

We put $\nabla^{(1)}:=\nabla-\nbigf^{\irr}$,
and then 
$\theta_{\nabla}=\theta_{\nabla^{(1)}}+\nbigf^{\irr}/2$.
Recall the relation
$R(h_0)=-2d''\theta_{\nabla}$.
Hence we have the following:
\begin{multline}
 \left(\frac{\sqrt{-1}}{2\pi}\right)^2
 \int_{X_{\delta}} \tr\bigl(R(h_0)^2\bigr)
=-2\left(\frac{\sqrt{-1}}{2\pi}\right)^2
 \int_{X_{\delta}} \tr\bigl(
 d''\theta_{\nabla}\, R(h_0)
 \bigr) \\
=-2\left(\frac{\sqrt{-1}}{2\pi}\right)^2
 \int_{X_{\delta}}
 d\tr\bigl(\theta_{\nabla}\, R(h_0)\bigr)
=-2\left(\frac{\sqrt{-1}}{2\pi}\right)^2
 \int_{\del X_{\delta}}
 \tr\bigl(\theta_{\nabla^{(1)}}\, R(h_0)\bigr)
\end{multline}
We put $Y_{\delta,i}:=\bigl\{|\sigma_i|=
 \delta\bigr\}\cap X_{\delta}$.
Let us look at $\int_{Y_{\delta,i}}
 \tr\big(\theta_{\nabla^{(1)}}\, R(h_0)\bigr)$.
The integrand is $0$ around  $D_i\cap D_j$ for any $j\neq i$.
Let $P\in D_i\setminus \bigcup_{j\neq i}D_j$.
Let $(U,z_1,z_2)$, $D_U$
and $\varphi_P:U'\lrarr U$ be as above.
We put $D'_{U}:=\varphi_P^{-1}(D_U)$.
We use the notation in Section \ref{subsection;07.6.12.20}
for this local argument.
On $U'$, we have 
$\nabla^{(1)}=\nabla^{(0)}+F^{\reg}$,
where $F^{\reg}$ is as in (\ref{eq;07.6.12.25}).
Since $\theta_{\nabla^{(0)}}$ is $C^{\infty}$
on $U'$,
it does not contribute to the limit for $\delta\to 0$.
(Note $c=\gamma=0$ in this case.)
Thus we have only to look at the term
$\tr\bigl(F^{\reg}\, R(h'_0)\bigr)/2$.
Then, the limit of the contribution of
$U\cap Y_{\delta,i}$ is as follows:
\begin{multline}
-\frac{1}{d}
 \int_{D'_{U}}
\frac{\sqrt{-1}}{2\pi}
 \sum_{(a,\alpha)\in\KMS(\prolong{E}_{\ast},D_i)}
 (d\alpha+da)\,\tr R(h'_{a,\alpha})
=\\
-\sum_{(a,\alpha)\in\KMS(\prolong{E}_{\ast},D_i)}
 (\alpha+a)
 \frac{\sqrt{-1}}{2\pi}\int_{D_{U}}
\Bigl(
 \tr R(h_{a,\alpha})
-a\,\rank\lefttop{i}\Gr^{F,\EE}_{a,\alpha}(E)
 \,\omega_i
\Bigr)
\end{multline}
Here, $h_{a,\alpha}'$ is 
as in Section \ref{subsection;08.1.15.3},
and we have used Lemma \ref{lem;08.1.15.5}.
We remark $R(h_0)=0$ around any cross points of $D$
and the vanishing (\ref{eq;08.1.15.6}).
Together with (\ref{eq;08.1.15.7}),
we obtain the following:
\begin{multline}
\label{eq;08.2.23.30}
\lim_{\delta\to 0}
\left( \frac{\sqrt{-1}}{2\pi}\right)^2
\int_{X_{\delta}}
 \tr\Bigl(R(h_0)^2\Bigr)=\\
-\sum_i
 \sum_{(a,\alpha)\in\KMS(\prolong{E}_{\ast},D_i)}
 (a+\alpha)\,
 \Bigl(
 \pardeg\bigl(
 \lefttop{i}\Gr^{F,\EE}_{a,\alpha}
 (\prolong{E})_{\ast}\bigr)
-a\, 
 \rank\bigl(\lefttop{i}\Gr^{F,\EE}_{a,\alpha}
 (\prolong{E})\bigr)
\,[D_i]^2
\Bigr)
\end{multline}
As shown in the proof of Proposition 4.18 of \cite{mochi4},
we have the equality:
\[
 \left(\frac{\sqrt{-1}}{2\pi}\right)^2
 \int_{X-D}\tr\bigl(R(h_0)^2\bigr)
=\int_X2\parch_2(\vecE_{\ast})
\]
Taking the real part,
we obtain the desired formula.
\hfill\qed

\begin{prop}
Let $\tau$ be a closed $2$-form on $X$.
Then, we have the following equality:
\[
 \int_X \parchern_1(\vecE_{\ast})\,\tau
=-\sum_i
 \sum_{(a,\alpha)\in
 \KMS(\prolongg{\vecc}{E},D_i)}
 \Re(a+\alpha)\,
 \rank\lefttop{i}\Gr^{F,\EE}_{a,\alpha}
 (\prolongg{\vecc}{E})
 \, \int_{D_i}\tau
\]
\end{prop}
\pf
Recall that
$(\sqrt{-1}/2\pi)\, \tr R(h_0)$ represents
$\parchern_1(\vecE_{\ast})$.
We also have the relation
$\tr R(h_0)=-2\delbar\tr(\theta_{\nabla})$.
Then, it can be shown by an argument
in the proof of Proposition \ref{prop;08.2.23.31}.
\hfill\qed

\begin{rem}
See Section {\rm 3.5} 
of {\rm \cite{mochi5}}
for more formulas
in the regular case.
Similar formulas should hold 
even in the irregular case.
\hfill\qed
\end{rem}

Let $L$ be an ample line bundle on $X$.

\begin{cor}
\label{cor;08.9.30.3}
Let $(\nbige,\nabla)$ be a meromorphic flat connection
on $(X,D)$ such that
the associated 
Deligne-Malgrange filtered flat bundle
$(\vecE_{\ast}^{DM},\nabla)$ is good.
(See Subsection {\rm\ref{subsection;08.2.22.5}}
for Deligne-Malgrange filtered flat bundle.)
Then, 
we have the vanishing of the characteristic numbers
$ \pardeg_L(\vecE_{\ast}^{DM})
=\int_X\parch_{2}(\vecE_{\ast}^{DM})=0$.
\end{cor}
\pf
We have only to remark that
any $(a,\alpha)\in\KMS(\vecE_{\ast}^{DM},i)$
satisfy $a+\Re(\alpha)=0$
by the definition of $\vecE_{\ast}^{DM}$.
\hfill\qed

\subsection{Blow up and 
the parabolic characteristic numbers
(Appendix)}
\label{subsection;08.9.29.22}

Let $X$ be a smooth projective surface
with a simple normal crossing hypersurface $D$.
Let $\pi:\Xtilde\lrarr X$ be a blow up of $X$
at a point $P$.
Let $\vecE_{\ast}$ be a filtered bundle on $(X,D)$.
We have the induced filtered bundle
$\vecEtilde_{\ast}$ on $(\Xtilde,\Dtilde)$.

\begin{lem}
\label{lem;07.11.5.100}
We have the following equality:
\begin{equation}
\label{eq;07.11.5.10}
 \int_{\Xtilde}\parch_2(\vecEtilde_{\ast})
=\int_{X}\parch_2(\vecE_{\ast})
\end{equation}
\begin{equation}
\label{eq;07.11.5.11}
  \parchern_1(\vecEtilde)
=\pi^{\ast}\bigl(
 \parchern_1(\vecE_{\ast})\bigr)
\end{equation}
See {\rm\cite{mochi4}} for $\parchern_1(\vecE_{\ast})$
and $\parch_2(\vecE_{\ast})$,
and see {\rm\cite{i-s}} for a more systematic treatment.
\end{lem}
\pf
The equality (\ref{eq;07.11.5.11})
can be reduced to the rank one case,
which can be checked easily.
Let us show (\ref{eq;07.11.5.10}).
We give only an indication for a direct calculation.
We use the formula (10) in \cite{mochi4}
for $\parch_2(\vecE_{\ast})$.
Let us consider the case $P\in D_1\cap D_2$.
The other case can be shown similarly.
Let $\Dtilde_i$ denote the proper transform of $D_i$,
and let $\Dtilde_P$ denote the exceptional divisor.

Let $\lefttop{i}F$ $(i=1,2)$
denote the parabolic filtration of
$\prolong{E}_{|D_i}$,
and let $\lefttop{i}\Gr^F_a(\prolong{E}):=
 \lefttop{i}F_a(\prolong{E}_{|D_i})\big/
 \lefttop{i}F_{<a}(\prolong{E}_{|D_i})$.
For any $\veca=(a_1,a_2)\in\real^2$,
we put 
\[
\lefttop{\nibar}F_{\veca}(\prolong{E}_{|P}):=
 \lefttop{1}F_{a_1}(\prolong{E}_{|P})
\cap
 \lefttop{2}F_{a_2}(\prolong{E}_{|P}),
\quad
\lefttop{\nibar}
 \Gr^{F}_{\veca}(\prolong{E}):=
 \frac{\lefttop{\nibar}F_{\veca}(\prolong{E}_{|P})}
 {\sum_{\vecb\lneq\veca}
 \lefttop{\nibar}F_{\vecb}(\prolong{E}_{|P})}
\]
We put 
$\Par(i):=\bigl\{
 a\,\big|\,\lefttop{i}\Gr^F_a(\prolong{E})\neq 0 
 \bigr\}$
and
$\Par(P):=\bigl\{
 \veca\in\real^2\,\big|\,
 \lefttop{\nibar}\Gr^F_{\veca}(\prolong{E})\neq 0
 \bigr\}$.

We take a splitting
$\prolong{E}_{|P}
=\bigoplus U_{\veca}$
of the filtrations $\lefttop{i}F$ $(i=1,2)$
indexed by $\Par(P)$,
i.e.,
it is taken as
$\lefttop{\nibar}F_{\veca}
=\bigoplus_{\vecb\leq\veca}U_{\vecb}$
holds for each $\veca=(a_1,a_2)$.
We put
$\lefttop{P}F_{b}(\prolong{E}_{|P}):=
 \bigoplus_{a_1+a_2\leq b}
 U_{\veca}$,
which is independent of the choice of a splitting.
We have
$\rank \lefttop{P}F_{-1}=\sum_{a_1+a_2\leq -1}
 \rank\lefttop{\nibar}\Gr^{F}_{\veca}(\prolong{E})$.
We have the exact sequence:
\[
 0\lrarr\pi^{\ast}(\prolong{E})
\lrarr \prolong{\Etilde}
\lrarr 
 \lefttop{P}F_{-1}\otimes\nbigo_{\Dtilde_P}(\Dtilde_P)
\lrarr 0
\]
Hence, we obtain the following equality,
by using Lemma 3.23 of \cite{mochi4}:
\begin{equation}
\label{eq;07.11.5.15}
  \int_{\Xtilde}\parch_2(\prolong{\Etilde})
-\int_X\parch_2(\prolong{E})
=-\frac{1}{2}\sum_{a_1+a_2\leq -1}
 \rank\lefttop{\nibar}\Gr^{F}_{\veca}(\prolong{E})
\end{equation}
Let $\lefttop{P}\Ftilde$
denote the parabolic filtration
of $\prolong{\Etilde}_{|\Dtilde_P}$,
and we set 
\[
 \Gr^{\lefttop{P}\Ftilde}_a\bigl(
 \prolong{\Etilde}\bigr):=
 \lefttop{P}\Ftilde_a(\prolong{\Etilde}_{|\Dtilde_P})
 \big/\lefttop{P}\Ftilde_{<a}
 (\prolong{\Etilde}_{|\Dtilde_P}).
\]
Let $\Par\bigl(\Dtilde_P\bigr):=
 \bigl\{a\,\big|\,
 \Gr^{\lefttop{P}\Ftilde}_a
 \bigl(\prolong{\Etilde}\bigr)\neq 0
 \bigr\}$.
We have the following equality:
\begin{equation}
\label{eq;07.11.5.16}
-\sum_{a\in \Par\bigl(\Dtilde_P\bigr)}
 a\, \deg_{\Dtilde_P}
 \Gr^{\lefttop{P}\Ftilde}_a(\prolong{\Etilde})
= 
\sum_{\substack{\veca\in\Par(P)\\
 a_1+a_2\leq -1 }}
 (a_1+a_2+1)\, 
 \rank\lefttop{\nibar}\Gr^{F}_{\veca}(\prolong{E})
\end{equation}
Let $\lefttop{i}\Ftilde(\prolong{\Etilde}_{|\Dtilde_i})$
denote the parabolic filtration of
$\prolong{\Etilde}_{|\Dtilde_i}$,
and we set
\[
 \lefttop{i}\Gr^{\Ftilde}_a(\prolong{\Etilde}):=
 \lefttop{i}\Ftilde_a(\prolong{\Etilde}_{|\Dtilde_i})\big/
 \lefttop{i}\Ftilde_{<a}(\prolong{\Etilde}_{|\Dtilde_i}).
\]
We have the following equality:
\begin{multline}
\label{eq;07.11.5.17}
-\sum_{i=1,2}\sum_{a\in\Par(i)}
 a\, \deg_{\Dtilde_i}
 \bigl(\lefttop{i}\Gr^{\Ftilde}_a(\prolong{\Etilde})\bigr)
+\sum_{i=1,2}\sum_{a\in\Par(i)}
 a\, \deg_{D_i}
 \bigl(\lefttop{i}\Gr^{F}_a(\prolong{E})\bigr) \\
=-\sum_{\substack{\veca\in\Par(P)\\a_1+a_2\leq -1}}
 (a_1+a_2)\,
 \rank\lefttop{\nibar}\Gr^F_{\veca}(\prolong{E}).
\end{multline}
We have the following equality:
\begin{multline}
\label{eq;07.11.5.18}
\frac{1}{2}
 \sum_{a\in\Par(\Dtilde_P)}
 a^2\,\rank\bigl(
 \Gr^{\lefttop{P}F}_a(\prolong{\Etilde})
\bigr)\, [\Dtilde_P]^2=\\
-\frac{1}{2}
 \sum_{\substack{\veca\in\Par(P)\\ a_1+a_2\leq -1}}
 (a_1+a_2+1)^2\,
 \rank\lefttop{\nibar}\Gr^F_{\veca}(\prolong{E})
-\frac{1}{2}
  \sum_{\substack{\veca\in\Par(P)\\ a_1+a_2> -1}}
 (a_1+a_2)^2\,
 \rank\lefttop{\nibar}\Gr^F_{\veca}(\prolong{E})
\end{multline}
We have the following:
\begin{multline}
\label{eq;07.11.5.19}
 \frac{1}{2}\sum_{i=1,2}
 \sum_{a\in\Par(i)}
 a^2\,\rank\lefttop{i}\Gr^{\Ftilde}_a(\prolong{\Etilde})
\, [\Dtilde_i]^2
-\frac{1}{2}\sum_{i=1,2}
 \sum_{a\in\Par(i)}
 a^2\,\rank\lefttop{i}\Gr^{F}_a(\prolong{E})
 \,[D_i]^2\\
=
-\frac{1}{2}\sum_{i=1,2}
 \sum_{a\in\Par(i)}
 a^2\,\rank\lefttop{i}\Gr^{\Ftilde}_a(\prolong{\Etilde})
=
-\frac{1}{2}\sum_{\veca\in\Par(P)}
 (a_1^2+a_2^2)\,
 \rank\lefttop{\nibar}\Gr^{F}_{\veca}(\prolong{E})
\end{multline}
We have the following:
\begin{multline}
\label{eq;07.11.5.20}
 \sum_{\substack{\veca\in\Par(P)\\ a_1+a_2>-1}}
 (a_1+a_2)\, a_2\, r_{\veca}
+\sum_{\substack{\veca\in\Par(P)\\ a_1+a_2\leq -1}}
 (a_1+a_2+1)\, a_2\,
 r_{\veca}\\
+\sum_{\substack{\veca\in\Par(P)\\a_1+a_2>-1}}
 a_1 \,(a_1+a_2)\,
 r_{\veca}
+\sum_{\substack{\veca\in\Par(P)\\ a_1+a_2\leq -1}}
 a_1\,(a_1+a_2+1)\,
 r_{\veca}
-\sum_{\veca\in\Par(P)}
 a_1 \, a_2 \,
 r_{\veca}\\
=
 \sum_{\substack{\veca\in\Par(P)\\ a_1+a_2> -1}}
 (a_1^2+a_2^2+a_1 a_2) \,
 r_{\veca}
+ \sum_{\substack{\veca\in\Par(P)\\ a_1+a_2\leq -1}}
 (a_1^2+a_2^2+a_1 a_2+a_1+a_2)\,
 r_{\veca},
\end{multline}
where 
$r_{\veca}:=
 \rank\lefttop{\nibar}\Gr_{\veca}(\prolong{E})$.
Taking summation of
(\ref{eq;07.11.5.15})--(\ref{eq;07.11.5.20})
and using the formula (10) in \cite{mochi4},
we obtain (\ref{eq;07.11.5.10}).
\hfill\qed

\section{Preliminary correspondence}
\label{subsection;08.3.9.1}
Let $X$ be a smooth irreducible projective complex surface,
and $D$ be a simple normal crossing hypersurface
with the irreducible decomposition 
$D=\bigcup_{i\in \Lambda}D_i$.
We also assume that $D$ is ample.
Let $L$ be an ample line bundle on $X$,
and $\omega$ be a Kahler form which represents $c_1(L)$.
We will not distinguish 
$\omega$ and the associated Kahler metric.
Let $(\vecE_{\ast},\nabla)$ be 
a $\mu_L$-stable good filtered flat bundle
on $(X,D)$, which is graded semisimple.
Let $\vecc\in\real^{\Lambda}$
such that 
$c_i\not\in\Par(\vecE_{\ast},i)$.
We assume the SPW-condition for each 
$i\in {\Lambda}$,
i.e.,
there exist
 a positive integer $m$
 and $\gamma_i\in\real$ with $-1/m<\gamma_i\leq 0$,
 such that 
\[
 \Par(\vecE_{\ast})
\subset\bigl\{
 c_i+\gamma_i+p/m\,\big|\,
 p\in\seisuu,\,\,
 -1<\gamma_i+p/m<0
 \bigr\}.
\]

We take a large integer $d$ 
which is divided by $m\, \rank (E)!^3$.
For $\epsilon=1/d$,
we take a Kahler metric $\omega_{\epsilon}$ of $X-D$
as in Section 4.3.1 of \cite{mochi4}.

\begin{prop}
\label{prop;07.6.12.30}
There exists a hermitian metric $h_{HE}$ of $E$ on $X-D$
satisfying the following conditions:
\begin{itemize}
\item
 Hermitian-Einstein condition 
 $2^{-1} 
 \Lambda_{\omega_{\epsilon}} 
 G(h_{HE})=a\, \id_{E}$
 for some constant $a$ determined 
by the following equation:
\begin{equation}
 \label{eq;05.9.17.1}
 a\,\frac{\sqrt{-1}}{2\pi}
 \frac{\rank E}{2}
 \int_{X-D}\omega_{\epsilon}^2
=a\,\frac{\sqrt{-1}}{2\pi}
 \frac{\rank(E)}{2}
 \int_X\omega^2
=\pardeg_{\omega}(\prolongg{\vecc}{E}_{\ast})
\end{equation}
\item
 $h_{HE}$ is adapted to $\vecE_{\ast}$.
\item
 $\deg_{\omega_{\epsilon}}(E,h_{HE})
 =\pardeg_{\omega}(\vecE_{\ast})$.
\item
 We have the following equalities:
\[
 \int_X 8\parch_{2}(\vecE_{\ast})
=\left(\frac{\sqrt{-1}}{2\pi}\right)^2\int_{X-D}
 \tr\Bigl(G(h_{HE})^2\Bigr)
\]
\[
 4 \int_X \parchern_1^2(\prolongg{\vecc}{E}_{\ast})
=\left(\frac{\sqrt{-1}}{2\pi}\right)^2\int_{X-D}
 \tr\bigl(G(h_{HE})\bigr)^2
\]
\end{itemize}
\end{prop}
\pf
We take a $C^{\infty}$-metric $h_0$
of $E_{|X-D}$ which is
as in Section \ref{subsection;07.6.12.21}
 around any cross points of $D$,
and as in Section \ref{subsection;07.6.12.21}
 around any smooth points of $D$.
It satisfies the following property:

\begin{lem}
\mbox{{}}
\begin{itemize}
\item
$G(h_0)$ is bounded with respect to $h_0$
and $\omega_{\epsilon}$.
\item
The following equalities hold:
\begin{equation}
 \label{eq;07.6.12.10}
\frac{1}{4}
\left(\frac{\sqrt{-1}}{2\pi}\right)^2
\int_{X-D}\tr\bigl(G(h_0)^2\bigr)
=
\left(\frac{\sqrt{-1}}{2\pi}\right)^2
\int_{X-D}\tr\bigl(R(h_0)^2\bigr)
=
\int_X2\parch_2(\vecE_{\ast})
\end{equation}
\begin{equation}
\label{eq;07.6.12.12}
 \frac{1}{4}
\left(\frac{\sqrt{-1}}{2\pi}\right)^2
\int_{X-D}\tr\bigl(G(h_0)\bigr)^2
=
\left(\frac{\sqrt{-1}}{2\pi}\right)^2
\int_{X-D}\tr\bigl(R(h_0)\bigr)^2
=
\int_X\parchern_1(\vecE_{\ast})^2
\end{equation}
\end{itemize}
\end{lem}
\pf
The first claim follows from 
Lemmas \ref{lem;07.6.12.6} and
\ref{lem;07.6.12.8}.
The equality (\ref{eq;07.6.12.10})
follows from Lemma \ref{lem;07.6.12.7},
Lemma \ref{lem;07.6.12.11}
and the Chern-Weil formula to express 
the Chern character by the curvature
(See Proposition 4.18 of \cite{mochi4}
 in the trivial Higgs case, for example.)
The equality (\ref{eq;07.6.12.12}) follows from
the general formula $\tr\bigl(G(h_0)\bigr)=2\tr R(h_0)$
and the Chern-Weil formula.
(See Lemma 4.16 of \cite{mochi4},
 for example.)
\hfill\qed

\vspace{.1in}
Make a modification $h_{in}=h_0\, \exp(-g')$
as in the proof of Proposition 6.1 of \cite{mochi4},
such that
$\Delta_{\omega_{\epsilon}}\det(h_{in})$ is constant.
Then, $h_{in}$ satisfies the following conditions:
\begin{itemize}
\item
 $h_{in}$ is adapted to the parabolic structure of
 $\vecE_{\ast}$.
\item
 $G(h_{in})$ is bounded with respect to $h_{in}$
 and $\omega_{\epsilon}$.
\item
 Let $V$ be any saturated coherent subsheaf of $E$,
 and let $\pi_V$ denote the orthogonal projection of $E$
 onto $V$.
  Then $\delbar\pi_V$ is $L^2$ with respect to
 $h_{in}$ and $\omega_{\epsilon}$,
 if and only if 
 there exists a saturated coherent subsheaf
 $\prolongg{\vecc}{V}$
 of $\prolongg{\vecc}{E}$ such that 
 $\prolongg{\vecc}{V}_{|X-D}=V$
 (\cite{li2} and \cite{siu}).
 Moreover we have
 $\pardeg_{\omega}(\prolongg{\vecc}{V}_{\ast})
 =\deg_{\omega_{\epsilon}}(V,h_{in,V})$,
 where $h_{in,V}$ denotes the metric of $V$ induced by $h_{in}$.
\item
 $\tr\Lambda_{\omega_{\epsilon}}G(h_{in})\
  =2\rank(E)\, a$ 
 for the constant $a$
 determined by the equation
 {\rm(\ref{eq;05.9.17.1})}.
\item The following equalities hold:
\[
 \left(\frac{\sqrt{-1}}{2\pi}\right)^2
 \int_{X-D}\tr \Bigl(G(h_{in})^2\Bigr)
=\int_X8\parch_{2}(\prolongg{\vecc}{E}_{\ast}),
\]
\[
 \left(\frac{\sqrt{-1}}{2\pi}\right)^2
 \int_{X-D}\tr\Bigl(G(h_{in})\Bigr)^2
=\int_X4\parchern_{1}^2(\prolongg{\vecc}{E}_{\ast}).
\]
\end{itemize}
Then, the proposition follows from
Simpson's theorem in \cite{s1} and \cite{s2}.
(See also Proposition 2.49 of \cite{mochi5}.)
\hfill\qed

\section{Bogomolov-Gieseker inequality}
\label{subsection;07.11.9.30}
Let $Y$ be an $n$-dimensional smooth connected
projective variety over $\cnum$.
Let $L$ be an ample line bundle on $Y$,
and let $D$ be a simple normal crossing hypersurface 
of $Y$.

\begin{cor}
Let $(\vecE_{\ast},\nabla)$ be a $\mu_L$-stable 
good filtered flat bundle on $(Y,D)$.
Then, Bogomolov-Gieseker inequality
holds for $\vecE_{\ast}$.
Namely, we have the following inequality:
\[
 \int_Y\parch_{2,L}(\vecE_{\ast})
\leq
 \frac{\int_Y\parchern_{1,L}^2(\vecE_{\ast})}{2\rank E}.
\]
\end{cor}
\pf
According to
 Mehta-Ramanathan type theorem
(Proposition \ref{prop;06.8.12.15}),
the problem can be reduced to the case $\dim Y=2$.
By using $\epsilon$-perturbations 
as in the case (II) of
Section \ref{subsection;07.10.15.80},
we can reduce the problem
to the case in which
$(\vecE_{\ast},\nabla)$ is graded semisimple,
and satisfy the SPW-condition.
Then, the claim follows from 
Proposition \ref{prop;07.6.12.30}
and the inequality for the curvature of 
Hermitian-Einstein metric
(see \cite{s1}).
\hfill\qed

\begin{cor}
\label{cor;07.10.16.1}
Let $(\vecE_{\ast},\nabla)$ be a $\mu_L$-stable
good filtered flat bundle on $(Y,D)$
with the trivial characteristic numbers
$ \pardeg_L(\vecE_{\ast})
=\int_Y\parch_{2,L}(\vecE_{\ast})=0$.
Then, we have
$\parchern_{1,L}(\vecE_{\ast})=0$.
\end{cor}
\pf
Due to Proposition \ref{prop;06.8.12.15},
the problem can be reduced to the case $\dim Y=2$.
We have the following: 
\[
 \int_Y\parchern_1(\vecE_{\ast})\,
 c_1(L)=
 \pardeg_L(\vecE_{\ast})=0
\]
\[
 0=2\rank(E)\int_Y\parch_2(\vecE_{\ast})
\leq \int_Y\parchern_{1}^2(\vecE_{\ast})
\]
Then, the claim follows from the Hodge index theorem.
\hfill\qed

\chapter{Preliminary for Resolution}
\label{section;10.6.1.1}
We give a preparation for
resolution of ``turning points'' for Higgs fields
underlying harmonic bundle.
In Section \ref{subsection;07.11.1.1},
we explain a procedure to construct
a projective morphism via which
given ideals are transformed to principal ideals.

When we are given an algebraic equation
on the rational function field of an irreducible variety,
we sometimes hope that
the polar part of a solution may ramify
only along the poles.
We show in Section \ref{subsection;08.9.29.220}
that it can be achieved if we take the pull back
via an appropriate morphism.
It will be used in Section 
\ref{subsection;07.10.14.110}.

It is significant to ask, for a given Higgs field $\theta$,
whether there exists a projective birational map $\varphi$
such that $\varphi^{\ast}\theta$ is good
in the sense of Definition \ref{df;08.9.28.150}.
We show in Section \ref{subsection;07.6.12.51}
the existence of such a map under some condition.
This result will be useful for
the resolution of turning points
(Sections
 \ref{subsection;08.9.28.101}--\ref{subsection;07.10.15.2})
and the correspondence between 
wild harmonic bundles and polarized wild pure twistor
$D$-modules
(Section \ref{subsection;07.10.28.35}).

\section{Resolution for a tuple of ideals}
\label{subsection;07.11.1.1}
Let $G$ be a finite group.
Let $Y$ be a normal complex analytic space
provided with a $G$-action.
Let $I_1,\ldots,I_r$ be ideal sheaves of $\nbigo_Y$.
Assume that the $G$-action induces
the permutation of $\{I_1,\ldots,I_r\}$.
We would like to give a canonical construction
to obtain a normal complex analytic space $Y'$
equipped with a $G$-action and a $G$-equivariant
birational projective morphism 
$\psi:Y'\lrarr Y$
such that the following holds:
\begin{itemize}
\item
 Let $I'_j$ denote the ideal sheaves of $\nbigo_{Y'}$
 generated by $\psi^{-1}(I_j)$.
 Then, $I'_j$ are invertible.
 Moreover,
 for any point $P$ of $Y'$,
 the ideals $I'_j$ are totally ordered with respect to 
 inclusion relation around $P$.
(We may permit that
$I_j'=I_k'$ around $P$.)
\end{itemize}
Note that we have the induced
projective birational morphism
$Y'/G\lrarr X$.

We put $Y^{(0)}:=Y$.
Using an inductive argument,
we will construct $Y^{(k)}$
with the following property:
\begin{itemize}
\item
We have $G$-actions on $Y^{(k)}$
and $G$-equivariant maps
$\psi_k:Y^{(k)}\lrarr Y$.
\item
We put 
$I_j^{(k)}:=\psi_k^{-1}(I_j)\cdot\nbigo_{Y^{(k)}}$.
For any closed point $P\in Y^{(k)}$,
we have an ordered subset
\[
 S(P):=\{i_1(P),\ldots,i_k(P)\}\subset \{1,\ldots,r\}
\]
such that
(i) $I_{i_h(P)}$ $(h=1,\ldots,k)$
are invertible,
(ii) 
$ I^{(k)}_{i_1(P)}
\supset
 I^{(k)}_{i_2(P)}
\supset\cdots\supset
I^{(k)}_{i_k(P)}$,
(iii)
$I_j^{(k)}$ ($j\not\in S(P)$) 
are contained in
$I^{(k)}_{i_k(P)}$,
around $P$.
\end{itemize}
Assume that $Y^{(k)}$ has already been
constructed.
We have the following $G$-equivariant
ideal sheaf on $Y^{(k)}$:
\[
 I_Y^{(k)}:=\sum_{|J|=k+1}
 \bigcap_{i\in J}I_i^{(k)}
\]
For any closed point $P\in Y^{(k)}$
and for any $h\not\in S(P)$,
we have 
$ \bigcap_{j=1}^kI^{(k)}_{i_j(P)}
\cap I^{(k)}_h=I^{(k)}_h$.
Hence, the following holds around $P$:
\[
 I_Y^{(k)}=\sum_{h\not\in S(P)}I_h^{(k)}
\]
Let $\pi_k:Y^{(k+1)}\lrarr Y^{(k)}$ 
be the blow up of $Y^{(k)}$
along $I_Y^{(k)}$.
Let $\psi_{k+1}:=\psi_k\circ\pi_k$.
We put
$I_j^{(k+1)}:=
 \psi_{k+1}^{-1}(I_j)\cdot\nbigo_{Y^{(k+1)}}
=\pi_k^{-1}(I_j^{(k)})\cdot\nbigo_{Y^{(k+1)}}$.
Let $P'$ be any closed point of $Y^{(k+1)}$,
and let $P:=\pi_k(P')$.
Around $P'$,
the ideals $I^{(k+1)}_j$ $(j\in S(P))$
are invertible.
We have
$\sum_{j\not\in S(P)}I^{(k+1)}_j
=\pi_k^{-1}\bigl(I_Y^{(k)}\bigr)\cdot 
\nbigo_{Y^{(k+1)}}$,
which is invertible by definition of blow up.
Hence, one of $I_j^{(k+1)}$ $(j\not\in S(P))$
is invertible,
and we have $I_{h}^{(k+1)}\subset I_j^{(k+1)}$
for any $h\not\in S(P)$.
Thus, the inductive construction can proceed.
Let $Y'$ be the normalization of $Y^{(r+1)}$,
and let $\psi:Y'\lrarr Y$ be 
the naturally induced morphism.
By construction, $(Y',\psi)$ has the desired property.

\vspace{.1in}

Let $\Gtilde$ be a finite group,
and let $\Ytilde$ be a complex analytic space
with a $\Gtilde$-action.
Assume that we are given 
a surjective homomorphism $\Gtilde\lrarr G$,
and an equivariant non-ramified covering
$F:\Ytilde\lrarr Y$
which induces
$\Ytilde/\Gtilde\simeq Y/G$.
We have the induced ideals 
$\Itilde_j$ $(j=1,\ldots,r)$ of $\Ytilde$,
which are the pull back of $I_j$.
The $\Gtilde$-action induces
the permutation of
$\bigl\{\Itilde_1,\ldots,\Itilde_r\bigr\}$.
Applying the above procedure,
we obtain the $\Gtilde$-equivariant map
$\psitilde:\Ytilde'\lrarr \Ytilde$.

\begin{lem}
\label{lem;07.11.1.2}
We have the natural isomorphism
$\Ytilde'\simeq \Ytilde\times_{Y}Y'$.
The induced map
$F':\Ytilde'\lrarr Y'$ is the non-ramified Galois covering,
whose Galois group is 
the kernel $K$ of $\Gtilde\lrarr G$.
In particular,
we have the natural isomorphism
$\Ytilde'/\Gtilde
\simeq Y'/G$.
\end{lem}
\pf
Inductively,
we can check the following:
\begin{itemize}
\item
 We have the natural isomorphism
 $\Ytilde^{(k)}\simeq
 Y^{(k)}\times_Y\Ytilde$,
 and the induced map 
 $F^{(k)}:\Ytilde^{(k)}\lrarr Y^{(k)}$
 is unramified $K$-covering.
\item
 $\Itilde^{(k)}_j=
 (F^{(k)})^{-1}(I^{(k)}_j)\cdot
 \nbigo_{\Ytilde^{(k)}}$
and
 $\Itilde^{(k)}_Y= 
 (F^{(k)})^{-1}(I^{(k)}_Y)\cdot
 \nbigo_{\Ytilde^{(k)}}$.
\end{itemize}
Then, the claim of the lemma follows.
\hfill\qed

\begin{rem}
The construction can also work
in the category of complex algebraic geometry.
\hfill\qed
\end{rem}

\section[Separation of Ramification 
and the Polar part]{Separation of
 the ramification and the polar part
 of a Higgs field}
\label{subsection;08.9.29.220}
\subsection{Statement}

Let $X$ be a smooth projective variety
over $\cnum$.
Let $D$ be a simple normal crossing hypersurface
of $X$.
Let $U=\Spec R$ be 
an affine Zariski open subset of $X$.
For simplicity, we assume that
the ideal sheaf $\nbigo_U(-D\cap U)$ is principal,
and let $g\in R$ be a generator of
$\nbigo_U(-D\cap U)$.
Let $R_g$ denote the localization of $R$
with respect to $g$.
Let $K$ denote the quotient field of $R$,
which is the same as the rational function
field of $X$.
In the following,
an algebraic closure $\Kbar$ of $K$ is fixed.
Let $\nbigp(T)\in R_g[T]$ be a monic
of degree $r$.
We will prove the following proposition
in Sections 
\ref{subsection;08.4.3.25}--\ref{subsection;08.4.3.26}.

\begin{prop}
\label{prop;08.4.3.1}
There exists a birational projective
morphism $F:X_0\lrarr X$
with the following property:
\begin{itemize}
\item
 $D_0:=F^{-1}(D)$ is 
 a simple normal crossing hypersurface,
 and the restriction
 $X_0-D_0\lrarr X-D$ is an isomorphism.
\item
 Let $U_0:=F^{-1}(U)$.
 For any $P\in U_0\cap D_0$,
 we take a holomorphic coordinate neighbourhood
 $(X_P,z_1,\ldots,z_n)$ around $P$
 such that 
 $D_P:=D_0\cap X_P
 =\bigcup_{i=1}^{\ell}\{z_i=0\}$.
 We take a ramified covering
 $\psi_P:\Xtilde_P\lrarr X_P$ given by
 $\psi_P(\zeta_1,\ldots,\zeta_n)
 =(\zeta_1^{M},\ldots,\zeta_{\ell}^M,
 \zeta_{\ell+1},\ldots,\zeta_n)$,
 where $M$ is divisible by $r!$.
 We put $\Dtilde_P:=\psi_P^{-1}(D_P)$.
 Then, 
 there exist a finite subset
 $S\subset M(\Xtilde_P,\Dtilde_P)$
 and monic polynomials
 $\nbigq_{\gminia}\in 
 H(\Xtilde_P)[T]$
 for $\gminia\in S$ such that
 the following holds:
\[
 (F\circ\psi_P)^{\ast}\nbigp(T)
=\prod_{\gminia\in S}
 \nbigq_{\gminia}(T-\gminia)
\]
Here, $M(\Xtilde_P,\Dtilde_P)$
denote the ring of meromorphic functions on $\Xtilde_P$
which admit poles along $\Dtilde_P$,
and $H(\Xtilde_P)$ denote the ring of
holomorphic functions on $\Xtilde_P$.
\end{itemize}
\end{prop}

\begin{rem}
Let $\alpha$ be a (possibly) multi-valued root
of the polynomial $(F\circ\psi_P)^{\ast}\nbigp$.
The claim of Proposition {\rm\ref{prop;08.4.3.1}}
can be reworded that
the polar parts of $\alpha$ can ramify
only along $D_0$.
\hfill\qed
\end{rem}

\subsection{Construction of the space}
\label{subsection;08.4.3.25}

Let $K'/K$ be the Galois extension
associated to $\nbigp(T)$.
Let $G$ denote the Galois group of
$K'$ over $K$.
Let $\pi:X'\lrarr X$ be the normalization of $X$
in $K'$.
We put $D':=\pi^{-1}(D)$.
The pull back $U':=\pi^{-1}(U)$
is the affine scheme $\Spec R'$,
where $R'$ is the normalization of $R$
in $K'$.
We have the roots
$\alpha_i$ $(i=1,\ldots,r)$
of $\nbigp(T)$ in $R'_g$.

We take a large number $N$
such that $g^N\cdot\alpha_i\in R'$ for any $i$.
Let $I_{i,j}$ denote the ideal of $R'$
generated by $g^N$ and 
$g^N\cdot(\alpha_i-\alpha_j)$.
They naturally induce the ideal sheaves
$\nbigi_{i,j}$ of $\nbigo_{X'}$.
Note that the closed subset associated to $\nbigi_{i,j}$
are contained in $D'$.

We apply the construction in 
Section \ref{subsection;07.11.1.1}
to $X'$ with the ideal sheaves
$\nbigi_{i,j}$  and the $G$-action.
Then, we obtain the normal variety  $X_1'$
with the $G$-action
and the $G$-equivariant morphism
$F_1':X_1'\lrarr X'$.
We put $U_1':=(F_1')^{-1}(U')$,
$D_1':=(F_1')^{-1}(D')$
and $\alpha_{1,i}:=(F_1')^{-1}(\alpha_i)$.
The pull back of $g$ is denoted by $g_1$.
Note that $X_1'-D_1'\lrarr X'-D'$ is isomorphic.
Since the ideals sheaves
$(F_1')^{-1}\nbigi_{i,j}\cdot \nbigo_{X_1'}$ 
is invertible,
either one of the following holds,
for each $Q\in U_1'\cap D_1'$
and for each pair $i,j$:
\begin{description}
\item[($A1$)]
 $g_1^N$ generates 
$(F_1')^{-1}\nbigi_{i,j}\cdot \nbigo_{X_1'}$
 around $Q$.
 In this case,
 $\alpha_{1,i}-\alpha_{1,j}\in
 \nbigo_{X_1',Q}$,
 where $\nbigo_{X_1',Q}$ denote the local ring
 at $Q$.
\item[($A2$)]
 $g_1^N\cdot\bigl(
 \alpha_{1,i}-\alpha_{1,j}
 \bigr)$ generates 
$(F_1')^{-1}\nbigi_{i,j}\cdot \nbigo_{X_1'}$
 around $Q$.
 In this case,
 there exists $\chi_{i,j}\in \nbigo_{X_1',Q}$
 such that
 $1=(\alpha_{1,i}-\alpha_{1,j})\cdot \chi_{i,j}$
 in $\nbigo_{X_1',Q}$.
\end{description}
Let $X_1:=X_1'/G$, which is also a normal variety.
We have the induced morphism
$F_1:X_1\lrarr X$.
We put $D_1:=F_1^{-1}(D)$,
which is the same as $D_1'/G$.
Note that the restriction
$X_1-D_1\lrarr X-D$ is isomorphic,
and in particular,
$X_1-D_1$ is smooth.

We take a smooth projective variety
$X_0$ with the birational projective
morphism $F_{1,0}:X_0\lrarr X_1$
such that 
(i) $D_0:=F_{1,0}^{-1}(D_1)$
is a simple normal crossing divisor,
(ii) $X_0-D_0\simeq X_1-D_1$.
Let $\pi_0:X_0'\lrarr X_0$ denote the normalization
of $X_0$ in $K'$,
and we put $D_0':=\pi_0^{-1}(D_0)$.
Thus, we obtain the following commutative diagram:
\[
 \begin{CD}
 X_0' @>{F_{1,0}'}>>  X_1' @>{F_1'}>> X'\\
 @V{\pi_0}VV @V{\pi_1}VV @V{\pi}VV\\
 X_0 @>{F_{1,0}}>> X_1 @>{F_1}>> X
 \end{CD}
\]
We put $F:=F_1\circ F_{1,0}$,
$F':=F_1'\circ F_{1,0}'$,
$U_0:=F^{-1}(U)$
and $U_0':=\pi_0^{-1}(U_0)$.
We also put
$\alpha_{0,i}:=(F')^{-1}\alpha_i$,
which are the algebraic sections of 
$\nbigo_{X_0'}(\ast D_0')$
on $U_0'$.
The pull back of $g$ is denoted by $g_0$.
By construction, either one of the following holds
for any $Q\in U_0'\cap D_0'$
and any $i,j$:
\begin{description}
\item[($A_01$)]
 $\alpha_{0,i}-\alpha_{0,j}
 \in \nbigo_{U_0',Q}$.
\item[($A_02$)]
 There exists 
 $\chi_{i,j}\in\nbigo_{U_0',Q}$
 such that 
 $1=(\alpha_{0,i}-\alpha_{0,j})\cdot
 \chi_{i,j}$.
\end{description}

\subsection{Proof of 
Proposition \ref{prop;08.4.3.1}}
\label{subsection;08.4.3.26}

Let us show that $(X_0,F_0)$
has the desired property.
Let $P$ be any point of $D_0$.
We take an affine Zariski open neighbourhood
$\nbigu_P$ of $P$ in $X_0$
with an etale morphism
$\varphi_P:\nbigu_P\lrarr \cnum^n$
such that 
$D_P:=D_0\cap \nbigu_P=
 \varphi_P^{-1}\bigl(
 \bigcup_{i=1}^{\ell}\{w_i=0\}
 \bigr)$.
Let $g_i:=\varphi_P^{\ast}(w_i)$
$(i=1,\ldots,\ell)$.
Let $\Ktilde_P$ denote the Galois extension 
over $K$ associated to the polynomials
$T^M-g_i$ $(i=1,\ldots,\ell)$,
where $M$ is as in the claim of Proposition
\ref{prop;08.4.3.1}.
Let $\Ktilde_P'$ denote the minimal field
which contains $\Ktilde_P$ and $K'$.
Due to a general theory (\cite{fujisaki}, for example),
$\Ktilde_{P}'$ is the Galois extension over $K$,
and the Galois group $\Gtilde_P$
of $\Ktilde_P'/\Ktilde_P$
is the same as the Galois group of
$K'/K'\cap\Ktilde_P$.

Let $\nbigu_P'$,
$\nbigutilde_P$
and $\nbigutilde_P'$ denote
the normalization of $\nbigu_P$
in $K'$, $\Ktilde_P$ and $\Ktilde_P'$,
respectively.
We have the following commutative diagram:
\[
 \begin{CD}
 \nbigutilde_P' @>{\lambda_P'}>> \nbigu_P'\\
 @V{\pitilde_P}VV @V{\pi_P}VV \\
 \nbigutilde_P @>{\lambda_P}>> \nbigu_P
 \end{CD}
\]
We have the natural isomorphism
$\nbigutilde_P\simeq \nbigutilde_P'/\Gtilde_P$.
Let $\Dtilde_P:=\lambda_P^{-1}(D_P)$.
We remark that $\pitilde_P$ is unramified
at the generic points of $\Dtilde_P$,
because $M$ is divisible by $r!$.

In the following argument,
we will work in the category of 
complex analytic spaces.
Let $\Ptilde=\lambda_P^{-1}(P)$.
We take a point $Q\in \pitilde_P^{-1}(\Ptilde)$.
Let $\Gtilde_Q:=\bigl\{
 \sigma\in\Gtilde_P\,\big|\,
 \sigma(Q)=Q \bigr\}$.
We take a $\Gtilde_Q$-invariant
small open connected neighbourhood $V_Q$
of $Q$ in $\nbigutilde_P'$.
By shrinking $V_Q$,
we may assume that 
$V_{\Ptilde}:=V_Q/\Gtilde_Q$ naturally gives 
a small neighbourhood of $\Ptilde$
in $\nbigutilde_P$.
Let $\pi_Q:V_Q\lrarr V_{\Ptilde}$
denote the induced projection.

We have the minimal reduced closed analytic subspace
$R_{\Ptilde}\subset V_{\Ptilde}$
such that 
$\pi_Q$ is non-ramified covering 
on the non-empty open set
$V_{\Ptilde}-R_{\Ptilde}$.
We fix a base point 
$P_0\in V_{\Ptilde}-R_{\Ptilde}$,
and let $\pi_1(V_{\Ptilde}-R_{\Ptilde},P_0)$
denote the fundamental group.
Since $V_Q$ is normal,
it is smooth in codimension one.
Since $V_Q-\pi^{-1}(R_{\Ptilde})$ is
connected,
we have the naturally induced surjection
$\rho_Q:\pi_1(V_{\Ptilde}-R_{\Ptilde},P_0)
\lrarr \Gtilde_Q$.
For any smooth point $A$ of
$\Rtilde_P$,
let $\gamma_A$ denote 
an element of
$\pi_1(V_{\Ptilde}-R_{\Ptilde},P_0)$
which is induced by a path
connecting $A$ and $P_0$,
and a small loop around $R_{\Ptilde}$ near $A$.
Let $\sigma_A:=\rho_Q(\gamma_A)$.

Let $\alpha_{Q,j}$
denote the pull back of
$\alpha_{0,j}$ by
$V_Q\lrarr \nbigu_P'$,
and let $\Dtilde_Q:=\pi_Q^{-1}(\Dtilde_P)$.
\begin{lem}
\label{lem;08.4.3.2}
$\alpha_{Q,j}$ are $\Gtilde_Q$-invariant
in $\nbigo_{V_Q,Q}(\ast \Dtilde_Q)\big/
 \nbigo_{V_Q,Q}$.
\end{lem}
\pf
Let $A$ be as above.
For each $j$,
either one of the following holds:
\begin{description}
\item[($A_Q1$)]
 $\sigma_A^{\ast}\bigl(\alpha_{Q,j}\bigr)
 -\alpha_{Q,j}
 \in \nbigo_{V_Q,Q}$.
\item[($A_Q2$)]
 There exists $\chi_j\in \nbigo_{V_Q,Q}$
 such that
 $1=\chi_j\cdot
 \bigl(\sigma_A^{\ast}(\alpha_{Q,j})-\alpha_{Q,j}\bigr)$.
\end{description}
We have the closed analytic subset
$Z_A\subset V_Q$ of codimension $1$
such that 
(i) $Q\in Z_A$ and 
 $A\in \pi_Q(Z_A)$,
(ii) $Z_A\not\subset \pi_Q^{-1}(\Dtilde_P)$,
(iii) 
 $\bigl(
 \sigma_{A}^{\ast}(\alpha_{Q,j})
 -\alpha_{Q,j}\bigr)_{Q_1}$ 
 is contained in the maximal ideal
 at $Q_1$
 for any $Q_1\in Z_A$.
Hence, we can conclude that 
$(A_Q2)$ cannot happen.

Since the elements $\gamma_A$ generate
$\pi_1(V_{\Ptilde}-R_{\Ptilde},P_0)$,
the claim of Lemma \ref{lem;08.4.3.2}
follows.
\hfill\qed

\vspace{.1in}
Because the $\Gtilde_Q$-invariant part of
$\nbigo_{V_Q,Q}(\ast \Dtilde_Q)\big/
 \nbigo_{V_Q,Q}$
is naturally isomorphic to
$\nbigo_{V_P,P}(\ast \Dtilde_P)\big/
 \nbigo_{V_P,P}$,
we have
$\beta_{P,j}
 \in \nbigo_{V_P,P}(\ast \Dtilde_P)\big/
 \nbigo_{V_P,P}$
induced by $\alpha_{Q,j}$,
and hence
the finite subset
$S=\bigl\{
 \beta_{P,j}
 \bigr\}$.
We take a lift of $S$
to $\nbigo_{V_P,P}(\ast \Dtilde_P)$,
which is also denoted by $S$.
We have 
$\gminia(j)\in S$ for each $j$
such that
$\alpha_{Q,j}-\pi_Q^{\ast}\bigl(
 \gminia(j)\bigr)
 \in \nbigo_{V_Q}$.
Then, we can conclude that 
$(X_0,D_0)\lrarr (X,D)$ has the desired property,
and the proof of Proposition \ref{prop;08.4.3.1}
is finished.
\hfill\qed

\section{Resolution 
 for generically good Higgs field}
\label{subsection;07.6.12.51}
\subsection{Statement}
\label{subsection;07.11.1.16}

Let $X$ be a complex manifold,
and let $D$ be a normal crossing divisor.
In the following,
Zariski open subset means
a complement of a closed analytic subset.
\index{Zariski open}
Let $\nbiga$ be a $\rnum$-vector subspace of $\cnum$.
Let $(E,\theta)$ be a holomorphic Higgs bundle on $X-D$
of rank $r$.
We state some conditions
on $(E,\theta)$.
\begin{description}
\item[(Generically $\nbiga$-good)]
We have a Zariski open subset $D'$ of $D$
such that
$\theta$ is $\nbiga$-good around any point of $D'$.
(See Definition \ref{df;07.12.28.1} below.)
Note
the condition implies that the eigenvalues of $\theta$
are multi-valued meromorphic $1$-forms
on $(X,D)$,
which may ramify along some hypersurface of $X$.
\end{description}
If $(E,\theta)$ is generically $\nbiga$-good,
the following condition makes sense.
\begin{description}
\item[(RD)]
The ramification of the eigenvalues of
$\theta$ may happen only along $D$.
\end{description}

We have another slightly more complicated
condition.
For any $P\in D$,
 we take a holomorphic coordinate neighbourhood
 $(X_P,z_1,\ldots,z_n)$ around $P$
 such that 
 $D_P:=D\cap X_P
 =\bigcup_{i=1}^{\ell(P)}\{z_i=0\}$.
 We take a ramified covering
 $\psi_P:\Xtilde_P\lrarr X_P$ given by
 $\psi_P(\zeta_1,\ldots,\zeta_n)
 =(\zeta_1^{M},\ldots,\zeta_{\ell(P)}^M,
 \zeta_{\ell(P)+1},\ldots,\zeta_n)$,
 where $M$ is divisible by $r!$.
 We put $\Dtilde_P:=\psi_P^{-1}(D_P)$.
 We have the expression:
\[
 \psi_P^{\ast}(\theta)
=\sum_{i=1}^{\ell(P)} f_i\,\frac{d\zeta_i}{\zeta_i}
+\sum_{i=\ell(P)+1}^{n} f_i\,d\zeta_i
\]
Let $M(\Xtilde_P,\Dtilde_P)$
denote the ring of meromorphic functions on $\Xtilde_P$
which admit poles along $\Dtilde_P$,
and let $H(\Xtilde_P)$ denote
the ring of holomorphic functions on $\Xtilde_P$.
If $(E,\theta)$ is generically $\nbiga$-good,
the characteristic polynomials
$\det(T-f_i)$ are contained in $M(\Xtilde_P,\Dtilde_P)[T]$.
Then, the following condition also makes sense.
\begin{description}
\item[(SRP)]
 For each $P\in D$ and $i=1,\ldots,n$,
 we have a finite subset
 $S_i(P)\subset M(\Xtilde_P,\Dtilde_P)$
 and monic polynomials
 $\nbigq_{P,i,\gminib}\in 
 H(\Xtilde_P)[T]$
 ($\gminib\in S_i(P)$),
 such that the following holds:
\[
 \det(T-f_i)
=\prod_{\gminib\in S_i(P)}
 \nbigq_{P,i,\gminib}(T-\gminib)
\]
\end{description}

We will prove the following proposition
in Section \ref{subsection;08.4.3.20}
after the preparation in
Sections 
\ref{subsection;08.4.4.1}--\ref{subsection;08.4.3.12}.
\begin{prop}
\label{prop;07.11.1.3}
Assume that $(E,\theta)$ is
generically $\nbiga$-good,
and moreover 
it satisfies either (RD) or (SRP).
Then,
there exists a complex manifold $X'$
with a birational projective morphism
$\varphi:X'\lrarr X$ such that
(i) $D':=\varphi^{-1}(D)$ is normal crossing,
(ii) $X'-D'\simeq X-D$,
(iii) $\varphi^{\ast}(\theta)$ is $\nbiga$-good
(Definition {\rm\ref{df;07.12.28.2}}).
\end{prop}

\begin{cor}
Let $(E_i,\theta_i,h_i)$ $(i=1,2)$
be good wild harmonic bundles on $(X,D)$.
There exists a projective birational morphism
$\varphi:(X',D')\lrarr (X,D)$
such that 
$\varphi^{\ast}\bigl(
 (E_1,\theta_1,h_1)\otimes
 (E_2,\theta_2,h_2)
 \bigr)$ is good.
\end{cor}
\pf
It is easy to see that
$\id_{E_1}\otimes\theta_2
+\theta_1\otimes\id_{E_2}$
satisfies the condition (RD).
\hfill\qed

\begin{rem}
If (RD) holds, (SRP) also holds.
Hence, we have only to consider
the case (SRP)
in Proposition {\rm\ref{prop;07.11.1.3}}.
We consider the case (RD) too,
partially because it is easy to state and argue,
and partially because it suffices for 
some of our purposes.
\hfill\qed
\end{rem}

Proposition \ref{prop;07.11.1.3}
is a preparation for resolution of turning points
(Theorem \ref{thm;07.10.14.60})
and the correspondence between
wild harmonic bundles and polarized
wild pure twistor $D$-modules
(Theorem \ref{thm;07.10.28.30}).

\subsection{The irregular values
and the logarithmic one forms}
\label{subsection;08.4.4.1}

Let $X:=\Delta^n$ and
$D:=\bigcup_{i=1}^{\ell}\{z_i=0\}$.
Let $(E,\theta)$ be 
a generically unramifiedly 
$\nbiga$-good Higgs bundle on $X-D$.
We have the expression
\begin{equation}
\label{eq;08.9.2.5}
\theta=
 \sum_{i=1}^{\ell} f_i\,\frac{dz_i}{z_i}
+\sum_{i=\ell+1}^n f_i\, dz_i.
\end{equation}
We shall observe that
a (not necessarily good)
set of irregular values is associated
to $\theta$ in a functorial way,
under some assumptions.

Let $M(X,D)$ denote the ring of
meromorphic functions on $X$
which admit poles along $D$,
and let $H(X)$ denote the ring of
holomorphic functions on $X$.
Let $CL(X,D,\nbiga)$ denote the space of
logarithmic $1$-forms of the form 
$\sum a_i\, dz_i/z_i$
$(a_i\in \nbiga)$.
In the case $\nbiga=\cnum$,
it is simply denoted by $CL(X,D)$.
Recall the following lemma.
(See Proposition 3.13 in Chapter II of \cite{d}.)
\begin{lem}
\label{lem;08.4.4.2}
Let $\omega$ be a meromorphic one form
on $X$ which admits pole along $D$.
\begin{itemize}
\item
Assume that $d\omega$ is logarithmic on $(X,D)$,
then we can take $\gminia\in M(X,D)$
such that $\omega-d\gminia$ is logarithmic on $(X,D)$,
by shrinking $X$ around the origin
$(0,\ldots,0)$ appropriately.
Such $\gminia$ is well defined in $M(X,D)\big/H(X)$
for $\omega$.
\item
If moreover $d\omega$ is holomorphic on $X$,
we have the well defined $\kappa\in CL(X,D)$
such that
$\omega-d\gminia-\kappa$
is holomorphic on $X$.
\hfill\qed
\end{itemize}
\end{lem}

\subsubsection{The case (RD)}
\label{subsection;08.4.3.10}

In this subsection, we assume the following:
\begin{description}
\item[(UR)]
The eigenvalues of $f_i$ are single-valued on $X$.
\end{description}

\begin{rem}
If (RD) is satisfied,
the above condition is satisfied
on some ramified covering of $(X,D)$.
\hfill\qed
\end{rem}

We have the homogeneous
generalized eigen decomposition 
$E=\bigoplus E^{(p)}$
on a Zariski open subset of $X-D$,
and $f_i$ have the single-valued eigenvalues 
$\alpha^{(p)}_i$ on $E^{(p)}$.
We put
$\omega^{(p)}:=
 \sum_{i=1}^{\ell} \alpha^{(p)}_i\, dz_i/z_i
+\sum_{i=\ell+1}^n\alpha^{(p)}_i\, dz_i$.

\begin{lem}
\mbox{{}}
\begin{itemize}
\item
We have 
$\kappa^{(p)}\in CL(X,D,\nbiga)$
and meromorphic functions
$\gminia^{(p)}\in M(X,D)$
such that 
$\omega^{(p)}
 -\bigl(d\gminia^{(p)}+\kappa^{(p)}\bigr)$
are holomorphic on $X$.
\item
$\kappa^{(p)}$ are well defined for
$\omega^{(p)}$,
and 
$\gminia^{(p)}
\in M(X,D)\big/H(X)$ are well defined
for $\omega^{(p)}$.
\end{itemize}
\end{lem}
\pf
Since $\theta$ is generically $\nbiga$-good,
$d\omega^{(p)}$ are holomorphic
around any general point $Q\in D$.
Hence, $d\omega^{(p)}$ are holomorphic $2$-forms
on $X$.
Due to Lemma \ref{lem;08.4.4.2},
we can take
$\gminia^{(p)}\in M(X,D)$
and $\kappa^{(p)}\in CL(X,D)$
such that
$\omega^{(p)}-d\gminia^{(p)}-\kappa^{(p)}$ 
is holomorphic.
Since $\theta$ is generically $\nbiga$-good,
we obtain
$\kappa^{(p)}\in CL(X,D,\nbiga)$.
\hfill\qed

\vspace{.1in}
Let $\veckappa(\theta)$ denote
the tuple $\bigl(\kappa^{(p)}\bigr)$
of logarithmic one forms on $(X,D)$,
and let $\Irr(\theta)$ denote 
the subset $\bigl\{\gminia^{(p)}\bigr\}$
in $M(X,D)\big/H(X)$.
They are well defined for $\theta$.

\vspace{.1in}

Let $\Xtilde:=\Delta_{\zeta}^{n}$ and 
$\Dtilde:=\bigcup_{j=1}^{k}\{\zeta_j=0\}$.
Let $\psi:\Xtilde\lrarr X$ be a morphism
such that $\psi^{-1}(D)=\Dtilde$.
We have the induced tuple
$\psi^{\ast}\veckappa(\theta)$
of logarithmic one forms on $(\Xtilde,\Dtilde)$,
and the induced subset
$\psi^{\ast}\Irr(\theta)
\subset M(\Xtilde,\Dtilde)/H(\Xtilde)$.

\begin{lem}
We have the functoriality
$\Irr(\psi^{\ast}\theta)
=\psi^{\ast}\Irr(\theta)$.
We also have
$\veckappa(\psi^{\ast}\theta)
\equiv\psi^{\ast}\veckappa(\theta)$
modulo the holomorphic exact one forms.
\end{lem}
\pf
Since
$\psi^{\ast}\omega^{(p)}
-d\psi^{\ast}(\gminia^{(p)})
 -\psi^{\ast}\kappa^{(p)}$
are holomorphic,
the claim of the lemma is clear.
\hfill\qed

\subsubsection{The case (SRP)}
\label{subsection;08.4.3.11}

In this subsection,
we assume the following:
\begin{description}
\item[(UR-SRP)]
We have 
the finite subsets $S_i\subset M(X,D)$
and polynomials $\nbigp_{\gminib,i}\in H(X)[T]$
$(\gminib\in S_i)$
such that
\[
 \det(T-f_i)=\prod_{\gminib\in S_i}
 \nbigp_{\gminib,i}(T-\gminib)
\]
We also assume that
the roots of
$\nbigp_{\gminib,i}(T)$ are 
unramified along Zariski dense open subset of $D$.
\end{description}

\begin{rem}
If (SRP) is satisfied,
the above condition is satisfied 
on some ramified covering of $(X,D)$.
\hfill\qed
\end{rem}

We have a normal complex space
$X'$ with a finite morphism $\psi:X'\lrarr X$
such that the roots of $\psi^{\ast}\det(T-f_i)$
are single-valued for any $i$.
We have the homogeneous
generalized eigen decomposition 
$\psi^{\ast}E=\bigoplus E^{(p)}$
on a Zariski open subset of $X'$,
and $\psi^{\ast}f_i$ have 
single-valued eigenvalues 
$\alpha^{(p)}_i$ on $E^{(p)}$.
We put
$\omega^{(p)}:=
 \sum_{i=1}^{\ell}\alpha^{(p)}_i\, 
 \psi^{\ast}dz_i/z_i
+\sum_{i=\ell+1}^n
 \alpha^{(p)}_i\, \psi^{\ast}dz_i$.
We regard them as 
multi-valued meromorphic $1$-forms on $(X,D)$
which may ramify outside of $X-D$.
For each $\alpha^{(p)}_i$,
we have the corresponding
$\gminib^{(p)}_i\in S_i$.
We put
$\omega_0^{(p)}:=
 \sum_{i=1}^{\ell}
 \gminib^{(p)}_i\, dz_i/z_i
+\sum_{i=\ell+1}^n
 \gminib^{(p)}_i\, dz_i$,
which are single-valued meromorphic one forms
on $(X,D)$.

\begin{lem}
\mbox{{}}
\begin{itemize}
\item
 We can take $\gminia^{(p)}\in M(X,D)$ 
 such that
 $\omega_0^{(p)}-d\gminia^{(p)}$ 
 are logarithmic one forms.
 Such $\gminia^{(p)}$ are well defined
 in $M(X,D)/H(X)$ for $\omega^{(p)}$.
\item
 We have the well defined
 $\kappa^{(p)}\in
 CL(X,D,\nbiga)$
 such that
 $\omega^{(p)}-d\gminia^{(p)}-\kappa^{(p)}$
 are multi-valued holomorphic one forms on $X$,
 i.e.,
 it is of the form
 $\sum \beta^{(p)}_i\,\psi^{\ast}dz_i$ 
 where $\beta^{(p)}_i$ are holomorphic 
 on $X'$.
\end{itemize}
\end{lem}
\pf
Let $Q$ be any point of $D$
such that
(i) $(E,\theta)$ is unramifiedly $\nbiga$-good
around $Q$,
(ii) any roots of $\nbigp_{\gminib,i}$
$(i=1,\ldots,\ell,\,\gminib\in S_i)$
are unramified around $Q$.
On a small neighbourhood $X_Q$ of $Q$,
$\omega^{(p)}$ give the single-valued
meromorphic one forms,
and $d\omega^{(p)}$ are holomorphic $2$-forms.
Because 
$\omega^{(p)}-\omega_0^{(p)}$
are logarithmic,
$d\omega^{(p)}_0$ are logarithmic
on $X_Q$.
By varying $Q$ in a Zariski dense subset in $D$,
we obtain that $d\omega_0^{(p)}$ are logarithmic.
Due to Lemma \ref{lem;08.4.4.2},
we can take $\gminia^{(p)}\in M(X,D)$
such that 
$\omega^{(p)}_0-d\gminia^{(p)}$
are logarithmic.
We obtain that
$\omega^{(p)}-d\gminia^{(p)}$
are logarithmic on $X_Q$.
Since $\theta$ is generically $\nbiga$-good,
we can find
$\kappa^{(p)}_Q\in CL(X,D,\nbiga)$
such that
$\tau^{(p)}:=\omega^{(p)}
-\bigl(
 \kappa^{(p)}+d\gminia^{(p)}
 \bigr)$ are holomorphic
on $X_Q$.
Hence, $\tau^{(p)}$ are holomorphic on $X'$.
\hfill\qed

\vspace{.1in}

Let $\veckappa(\theta)$ denote
the tuple of the logarithmic one forms
$\bigl(\kappa^{(p)}\bigr)$,
and let $\Irr(\theta)$ denote the subset 
$\{\gminia^{(p)}\}$ of $M(X,D)/H(X)$.
They are well defined for $\theta$.

\vspace{.1in}

Let $\Xtilde:=\Delta^n$ and 
$\Dtilde:=\bigcup_{j=1}^{m}\{\zeta_j=0\}$.
Let $\psi:\Xtilde\lrarr X$ be a morphism
such that $\psi^{-1}(D)=\Dtilde$.
We have the induced tuple of
logarithmic one forms
$\psi^{\ast}\veckappa(\theta)$,
and the induced subset
$\psi^{\ast}\Irr(\theta)
\subset M(\Xtilde,\Dtilde)/H(\Xtilde)$.

\begin{lem}
$\psi^{\ast}(E,\theta)$ also satisfies
the condition
{\bf (UR-SRP)}.
We have the functoriality
$\Irr(\psi^{\ast}\theta)
=\psi^{\ast}\Irr(\theta)$.
We also have
$\veckappa(\psi^{\ast}\theta)
\equiv\psi^{\ast}\veckappa(\theta)$
modulo the holomorphic exact forms.
\end{lem}
\pf
Since
$\psi^{\ast}\omega^{(p)}
-d\psi^{\ast}(\gminia^{(p)})
-\psi^{\ast}\kappa^{(p)}$
are multi-valued holomorphic one forms on $\Xtilde$,
the claim of the lemma follows.
\hfill\qed

\subsection{The ideal sheaves
 associated to a set of irregular values}
\label{subsection;08.4.3.12}

Let $X:=\Delta^n$
and $D:=\bigcup_{i=1}^{\ell}\{z_i=0\}$.
Let $\nbigi$ be a finite subset of $M(X,D)/H(X)$.
We take a lift of $\nbigi$ to $M(X,D)$,
which is given by
$\bigl\{\gminia^{(p)}\,\big|\,
 p=1,\ldots,s\bigr\}$.
For $j=1,\ldots,\ell$,
we put $m_j(\gminia^{(p)}):=\ord_{z_j}(\gminia^{(p)})$
in the case that 
$\{z_j=0\}$ is contained in the pole of $\gminia^{(p)}$,
and $m_j(\gminia^{(p)}):=0$ otherwise.
We set
$m_j(\nbigi):=\min_p m_j(\gminia^{(p)})$,
which is well defined for $\nbigi$.
We put 
$\xi(\nbigi):=
 \prod_{j=1}^{\ell}z_j^{-m_j(\nbigi)}$,
which is well defined for $\nbigi$
and independent of the coordinate system
up to the multiplication of invertible functions.
Let $I(\gminia^{(p)},\nbigi)$ $(p=1,\ldots,s)$
denote the ideals generated by
$\xi(\nbigi)\cdot \gminia^{(p)}$
and $\xi(\nbigi)$.
The tuple of the ideals
$\vecI(\nbigi)$ is well defined for $\nbigi$.

\begin{rem}
The construction is independent of the choice
of a coordinate system,
and it can be globalized.
Namely,
let $X$ be a complex manifold,
and let $D$ be a normal crossing divisor of $X$.
If we are given 
a finite subset $\nbigi\subset M(X,D)\big/H(X)$,
we obtain the tuple of the ideals 
$\vecI(\nbigi)$,
whose restriction to holomorphic coordinate 
neighbourhoods are given as above.

We also have the purely algebraic construction of
the ideal $\vecI(\nbigi)$,
when $X$ is algebraic.
\hfill\qed
\end{rem}

Let $\Xtilde:=\Delta^n$ and 
$\Dtilde:=\bigcup_{j=1}^{\ell}\{\zeta_j=0\}$.
Let $\psi:(\Xtilde,\Dtilde)\lrarr (X,D)$ 
be a ramified covering.
We have the induced subset
$\psi^{-1}\nbigi\subset M(\Xtilde,\Dtilde)/H(\Xtilde)$
and the induced ideals
$\psi^{-1}\vecI(\nbigi):=\bigl(
\psi^{-1}I(\gminia^{(p)},\nbigi)\cdot
\nbigo_{\Xtilde}\,\big|\,
p=1,\ldots,s\bigr)$.
It is easy to see
$\vecI(\psi^{-1}\nbigi)
=\psi^{-1}\vecI(\nbigi)$.
If $\psi$ is a ramified covering,
we can reconstruct
$\vecI(\nbigi)$ as the Galois descent of
$\vecI(\psi^{-1}\nbigi)$.

\subsection{Construction of a resolution}
\label{subsection;08.4.3.20}

Let us return to the setting in Section 
\ref{subsection;07.11.1.16}.
Let $P$ be any point of $D$.
We take a coordinate neighbourhood $U$ around $P$.
We take a ramified covering
$\varphi:(U^{\circ},D^{\circ})\lrarr (U,D)$
such that either
(UR) or (UR-SRP) is satisfied for 
$\varphi^{\ast}\theta$.
Applying the construction in Sections
\ref{subsection;08.4.3.10}--\ref{subsection;08.4.3.11},
we obtain the subset 
$\nbigi_P:=\Irr(\varphi^{\ast}\theta)
 \subset M(U^{\circ},D^{\circ})/H(U^{\circ})$,
and the set of logarithmic one forms
$\veckappa(\varphi^{\ast}\theta)$.
Applying the procedure of 
Section \ref{subsection;08.4.3.12} to $\nbigi_P$,
we obtain the tuple of ideals $\vecI(\nbigi_P)$
of $\nbigo_{U^{\circ}}$.
Applying the procedure in 
Section \ref{subsection;07.11.1.1} to them,
we obtain the complex analytic space
$U'$ with a birational projective morphism 
$U'\lrarr U$.
According to Lemma \ref{lem;07.11.1.2}
and the remark in the last part of 
Section \ref{subsection;08.4.3.12},
we can globalize the construction,
and thus we obtain the complex analytic space $X'$
with the birational projective morphism
$X'\lrarr X$.
Taking a resolution of singularity,
we obtain a complex manifold $X_1$
with a birational projective morphism 
$F_1:X_1\lrarr X$
such that 
(i) $D_1:=F_1^{-1}(D)$ is normal crossing, 
(ii) $X_1-D_1\simeq X-D$.
We put $(E_1,\theta_1):=F_1^{\ast}(E,\theta)$.

\vspace{.1in}

Let $P_1$ be any point of $D_1$.
We take a small neighbourhood $U$ of $P=F_1(P_1)$.
Let $U_1$ be a coordinate neighbourhood around $P_1$
such that $F_1(U_1)\subset U$.
We take a ramified covering $(U^{\circ},D^{\circ})$
of $U$, as above.
We have $\nbigi_P$ and $\veckappa(\varphi^{\ast}\theta)$
on $(U^{\circ},D^{\circ})$ as above.
We take a lift of $\nbigi_P$ to
$M(U^{\circ},D^{\circ})$
which is given by
$\bigl\{\gminia^{(p)}\bigr\}$.
Let $\xi(\nbigi_P)$ be as in
Section \ref{subsection;08.4.3.12}.
Let $\psi_1:(U_1^{\circ},D_1^{\circ})\lrarr (U_1,D_1)$
be a ramified covering
such that the composite
$U_1^{\circ}\lrarr U$ factors through $U^{\circ}$,
and so we have 
$F^{\circ}_1:(U_1^{\circ},D_1^{\circ})\lrarr
 \bigl(U^{\circ},D^{\circ}\bigr)$.
Then, either (UR) or (UR-SRP)
holds for $\psi_1^{\ast}\theta_1$.
We put 
$\gminia_1^{(p)}:=
 (F_1^{\circ})^{\ast}\gminia^{(p)}$,
and then $\Irr(\psi_1^{\ast}\theta_1) 
\subset M(U_1^{\circ},D_1^{\circ})/H(U_1^{\circ})$
is the same as 
$(F_1^{\circ})^{\ast}\nbigi_P
=\bigl\{ \gminia_1^{(p)}
 \bigr\}$.
We also have
$\veckappa(\psi_1^{\ast}\theta_1)
\equiv
 (F_1^{\circ})^{\ast}\veckappa(\varphi^{\ast}\theta)$
modulo the exact holomorphic one forms.
By construction, either one of the following holds,
for each $p$:
\begin{description}
\item[(Case 1)]
$I(\gminia^{(p)},\nbigi_P)\cdot\nbigo_{U^{\circ}_1}$
is generated by 
$(F_1^{\circ})^{\ast}\xi(\nbigi_P)$.
\item[(Case 2)]
$I(\gminia^{(p)},\nbigi_P)\cdot\nbigo_{U^{\circ}_1}$
is generated by 
$(F_{1}^{\circ})^{\ast}\bigl(
 \xi(\nbigi_P)\cdot \gminia^{(p)}\bigr)$.
\end{description}
In the case 1,
since $(F_1^{\circ})^{\ast}
  \xi(\nbigi_P)\cdot\gminia^{(p)}_1$
is the multiplication of
$(F_{1}^{\circ})^{\ast}\xi(\nbigi_P)$
and a holomorphic function,
$\gminia_1^{(p)}$ is holomorphic
on $U_1^{\circ}$.
Let us consider the case 2.
Let $(z_{1,1},\ldots,z_{1,n})$
be a coordinate system of $U_1^{\circ}$
such that 
$D_1^{\circ}=
 \bigcup_{j=1}^{\ell_1}\{z_{1,j}=0\}$.
Then,
$(F_1^{\circ})^{\ast}\xi(\nbigi_P)$
is the product of an invertible function
and a monomial
$\prod_{j=1}^{\ell_1}z_{1,j}^{q_j}$
for some $q_j\geq 0$.
Since $(F_1^{\circ})^{\ast}\xi(\nbigi_P)
 \cdot\gminia^{(p)}_1$ divides
$(F_1^{\circ})^{\ast}\xi(\nbigi_P)$,
we obtain that
$\gminia_1^{(p)}$ is
the product of an invertible function
and $\prod_{j=1}^{\ell_1}z_{1,j}^{r_j(p)}$
for some $r_j(p)\leq 0$.

\vspace{.1in}

We put $\gminia_1^{(p,q)}:=
 \gminia_1^{(p)}-\gminia_1^{(q)}$.
Let $\nbigitilde_{P_1}\subset 
 M(U_1^{\circ},D_1^{\circ})\big/
 H(U_1^{\circ})$ be given by
$\bigl\{
 \gminia_1^{(p,q)}
 \bigr\}$.
Applying the procedure in Section
\ref{subsection;08.4.3.12} to $\nbigitilde_{P_1}$,
we obtain the tuple of ideals
$I(\nbigitilde_{P_1})$.
Applying the procedure in Section
\ref{subsection;07.11.1.1},
making the globalization,
and taking a resolution of singularity,
we obtain a  complex manifold $X_2$
with a birational projective morphism
$F_{2,1}:X_2\lrarr X_1$
such that 
(i) $D_2:=F_{2,1}^{-1}(D_1)$ is normal crossing,
(ii) $X_2-D_2\simeq X_1-D_1$.
We set
$(E_2,\theta_2):=
  F_{2,1}^{\ast}(E_1,\theta_1)$
and $F_2:=F_1\circ F_{2,1}$.

Let $P_2$ be any point of $D_2$.
Let $P_1:=F_{2,1}(P_2)$ and $P:=F_{2}(P_2)$.
We take neighbourhoods $U_1$ and $U$
of $P_1$ and $P$ as above, respectively.
We also take ramified coverings 
$(U_1^{\circ},D_1^{\circ})\lrarr (U_1,D_1)$
and $(U^{\circ},D^{\circ})\lrarr (U,D)$
as above.
Let $U_2$ be a coordinate neighbourhood around $P_2$
such that $F_{2,1}(U_2)\subset U_1$.
Let $\psi_2:(U_2^{\circ},D_2^{\circ})\lrarr (U_2,D_2)$
be a ramified covering such that the composite
$U_2^{\circ}\lrarr U_1$ factors through $U_1^{\circ}$.
So we have
$F_{2,1}^{\circ}:
 (U_2^{\circ},D_2^{\circ})\lrarr
 (U^{\circ}_1,D_1^{\circ})$.
Then, either (UR) or (UR-SRP) 
holds for $\psi_2^{\ast}\theta_2$.
We put $\gminia_2^{(p)}:=
 (F_{2,1}^{\circ})^{\ast}\gminia_1^{(p)}$,
and then
$\Irr(\psi_2^{\ast}\theta_2)$
is the same as
$\bigl\{
 \gminia_2^{(p)}
 \bigr\}$.
We also have 
$\veckappa(\psi_2^{\ast}\theta_2)
\equiv (F_{2,1}^{\circ})^{\ast}
 \veckappa(\psi_1^{\ast}\theta_1)$
modulo the holomorphic exact one forms.
Let $(w_{1},\ldots,w_{n})$ be
a coordinate system
of $(U^{\circ}_2,D_2^{\circ})$
such that $D_2^{\circ}=
 \bigcup_{k=1}^{\ell_2}\{w_{k}=0\}$.
By construction of $X_2$,
either one of the following holds
for each $\gminia_2^{(p)}-\gminia_2^{(q)}$:
\begin{description}
\item[(Case 1)]
 $\gminia_2^{(p)}-\gminia_2^{(q)}$
 is holomorphic on $U^{\circ}_2$.
\item[(Case 2)]
 If $\gminia_2^{(p)}-\gminia_2^{(q)}$ 
 admits a pole,
 $\gminia_2^{(p)}-\gminia_2^{(q)}$
 is the product of an invertible function 
 and 
 $\prod_{k=1}^{\ell_2}w_{k}^{s_k(p,q)}$
 for some $s_k(p,q)\leq 0$.
In this sense,
$\ord(\gminia_2^{(p)}-\gminia_2^{(q)})$
 exists in $\seisuu_{\leq 0}^{\ell_2}$.
\end{description}
Moreover, 
the ideals
$(F_{2,1}^{\circ})^{-1}
 I\bigl(\gminia^{(p,q)}_1,\nbigitilde_{P_1}\bigr)
\cdot\nbigo_{U^{\circ}_2}$
are totally ordered with respect to inclusion relation.
Hence, the set of
$\ord(\gminia_2^{(p)}-\gminia_2^{(q)})$
for non-holomorphic
$\gminia_2^{(p)}-\gminia_2^{(q)}$
is totally ordered with respect to
$\leq_{\seisuu^{\ell_2}}$.
By construction of $X_1$,
either one of the following holds
for each $p$:
(i) $\gminia_2^{(p)}$ is holomorphic,
or (ii) $\ord(\gminia_2^{(p)})$ exists
in $\seisuu_{\leq0}^{\ell_2}$.
Hence, 
$\Irr(\psi_2^{\ast}\theta_2)$ is 
a good set of irregular values.

By construction,
any eigenvalues of
$\psi_2^{\ast}\theta_2$
are of the form
$d\gminia+\kappa+\tau$,
where
(i) $\gminia$ are contained
 in the good set of irregular values
 $\Irr(\psi_2^{\ast}\theta_2)$,
(ii) $\kappa=\sum_{k=1}^{\ell_2} a_k\cdot dw_k/w_k$
 $(a_k\in\nbiga)$,
(iii) $\tau$ are multi-valued holomorphic one forms
 on $U_2^{\circ}$.
Hence $\theta_2$ is $\nbiga$-good,
and thus the proof of Proposition \ref{prop;07.11.1.3}
is finished.
\hfill\qed

\chapter{Kobayashi-Hitchin Correspondence
and Some Applications}
\label{section;08.10.24.41}
In this chapter, we study existence problem
of pluri-harmonic metrics
for good filtered flat bundles.
As an application,
we show the existence of a resolution
of turning points for a meromorphic flat bundle.
(See {\em Introduction}.)

\vspace{.1in}

In Section \ref{subsection;07.11.9.31},
we establish Kobayashi-Hitchin correspondence
between good wild harmonic bundles
and $\mu_L$-polystable good filtered
flat bundles with trivial characteristic numbers
(Theorem \ref{thm;06.1.23.100}).
This is a rather straightforward generalization
of the correspondence in the tame case.

We study in Sections 
\ref{subsection;08.9.28.101}--\ref{subsection;07.10.15.2}
a characterization of semisimplicity of
meromorphic flat bundles 
by the existence of $\sqrt{-1}\real$-good wild
pluri-harmonic metric
(Theorem \ref{thm;07.10.15.1}).
We also show the resolution of turning points
for a meromorphic flat bundle on a projective variety
(Theorem \ref{thm;07.10.14.60}).
These two results are intimately related.
(They will be slightly refined 
in Section \ref{subsection;08.12.18.3}.)

Contrast to the tame case,
we cannot directly apply Theorem \ref{thm;06.1.23.100}
to show Theorem \ref{thm;07.10.15.1},
because the Deligne-Malgrange filtered sheaf
of a meromorphic flat bundle
is not necessarily {\em good} even in the surface case.
We need the existence of a resolution of turning points
as in Theorem \ref{thm;07.10.14.60}.
Conversely, once we know Theorem \ref{thm;07.10.15.1},
it is rather easy to show Theorem \ref{thm;07.10.14.60}.

Our argument will proceed as follows.
We have already established 
Theorem \ref{thm;07.10.14.60} for the surface case
in \cite{mochi6} and Proposition \ref{prop;10.5.4.3}.
Applying Theorem \ref{thm;06.1.23.100},
we obtain Theorem \ref{thm;07.10.15.1}
for the surface case.
Then, we will show a variant of 
Theorem \ref{thm;07.10.15.1}
in the higher dimensional case,
in which we do not assume
the existence of good Deligne-Malgrange lattice.
By using it,
we show Theorem \ref{thm;07.10.14.60}
in the higher dimensional case.
After that, it is easy to show 
Theorem \ref{thm;07.10.15.1}.
This is a nice interaction between 
the theories of wild harmonic bundles
and meromorphic flat bundles.

In Section \ref{subsection;08.12.18.3},
we refine the results in Section
\ref{subsection;08.9.28.101},
i.e.,
we consider the case that 
$X$ is not necessarily projective
but proper algebraic.

\section[Kobayashi-Hitchin 
correspondence]{Kobayashi-Hitchin correspondence
 for good filtered flat bundles}
\label{subsection;07.11.9.31}
Let $X$ be an $n$-dimensional 
connected smooth projective variety
with an ample line bundle $L$,
and let $D$ be 
a simple normal crossing hypersurface of $X$
with the irreducible decomposition
$D=\bigcup_{i\in\Lambda}D_i$.
Let $(\vecE_{\ast},\nabla)$ be 
a $\mu_L$-stable good filtered flat bundle on $(X,D)$
with trivial characteristic numbers
$\pardeg_{L}(\vecE_{\ast})=
 \int_X\parch_{2,L}(\vecE_{\ast})=0$,
and we put $(E,\nabla):=(\vecE_{\ast},\nabla)_{|X-D}$.
Recall $\parchern_1(\vecE_{\ast})=0$ (Corollary \ref{cor;07.10.16.1}).
For each $\vecc\in\real^{\Lambda}$,
we have the determinant line bundle
$\det\bigl(\prolongg{\vecc}{E}\bigr)$ of the bundle
$\prolongg{\vecc}{E}$,
on which we have the induced parabolic structure
and the induced flat connection.
Thus we obtain the canonically determined
filtered flat bundle
$\bigl(\det\vecE_{\ast},\nabla\bigr)$ on $(X,D)$
of rank one.
We also have
$\parchern_1\bigl(\det\vecE_{\ast}\bigr)
=\parchern_1\bigl(\vecE_{\ast}\bigr)=0$.
Therefore, we can take a pluri-harmonic metric 
$h_{\det E}$ of $(\det(E),\nabla)$
which is adapted to the parabolic structure
of $\det \vecE_{\ast}$.
It is determined up to positive constant multiplication.

\begin{thm}
 \label{thm;06.1.23.100}
There exists a unique pluri-harmonic metric $h$
of $(E,\nabla)$ with the following properties:
\begin{itemize}
\item $\det(h)=h_{\det E}$.
\item $(E,\nabla,h)$ is a good wild harmonic bundle on $X-D$.
\item $h$ is adapted to the parabolic structure of
 $\vecE_{\ast}$.
\end{itemize}
\end{thm}
\pf
It can be proved by the argument
in the tame case \cite{mochi5},
after the preparation in Chapters 
\ref{section;08.3.9.2}--\ref{section;08.3.9.3}.
Hence, we give only an indication.

Let us begin with a remark.
If we have a pluri-harmonic metric $h$
of $(E,\nabla)$ adapted to $\vecE_{\ast}$,
 the corresponding harmonic bundle
$(E,\delbar_E,\theta,h)$ is good and wild,
according to Proposition \ref{prop;07.10.16.12}.
Let us show the uniqueness of such a metric.
Let $h_i$  $(i=1,2)$ be pluri-harmonic metrics
of $(E,\nabla)$ which are adapted to $\vecE_{\ast}$
such that $\det(h_i)=h_{\det(E)}$.
Let $C\subset X$ be 
any sufficiently ample and generic complete intersection curve,
which intersects with $D$ transversally.
The restriction of $(\vecE_{\ast},\nabla)_{|C}$
is stable, according to Proposition \ref{prop;06.8.12.15}.
Because of the weak norm estimate of 
good wild harmonic bundles 
(Theorem \ref{thm;07.12.2.55}),
the metrics $h_{i|C\setminus D}$
are adapted to $\vecE_{\ast|C}$.
Then,
we obtain $h_{1|C\setminus D}=h_{2|C\setminus D}$
due to Proposition \ref{prop;07.10.16.10}.
Since such curves $C$ cover the Zariski open subset
of $X-D$,
we obtain $h_1=h_2$ on $X-D$.
Thus, we obtain the uniqueness.

\vspace{.1in}

Let us show the existence.
In the curve case, it was shown by 
Biquard and Boalch
(Proposition \ref{prop;07.10.16.10}).
Let us consider the surface case.
\begin{lem}
\label{lem;07.10.16.15}
We have the desired pluri-harmonic metric $h$
in the case $\dim X=2$.
\end{lem}
\pf
The argument is essentially the same as that 
in the proof of Theorem 5.4 in \cite{mochi5}.
So, we give only an indication.
Let $(\prolongg{\vecc}{E},\vecF)$
be the $\vecc$-truncation.
For a large integer $d$, we put $\epsilon=1/d$,
and we take an $\epsilon$-perturbation
$\vecF^{(\epsilon)}$ as explained 
in (II) of Section \ref{subsection;07.10.15.80}.
Then, we have the hermitian metric $h_{HE}^{(\epsilon)}$
of $(E,\nabla)$  such that
(i) adapted to $\vecF^{(\epsilon)}$,
(ii) $\det(h_{HE}^{(\epsilon)})=h_{\det}$,
(iii) it satisfies the conditions 
in Proposition \ref{prop;07.6.12.30}.

We would like to take
the limit of 
the sequence $\{h_{HE}^{(\epsilon)}\}$
in $\epsilon\to 0$.
Let $m$ be large.
By using Lemma \ref{lem;07.6.12.35},
Lemma \ref{lem;07.6.12.36}
and the argument in Section 5.2 of \cite{mochi5},
we obtain the pluri-harmonic metric $h$
of $(E,\nabla)$ satisfying
 $\det(h)=h_{\det(E)}$
and the following property:
\begin{itemize}
\item
 $h_{|s^{-1}(0)}$ is adapted to
 $\vecE_{\ast|s^{-1}(0)}$ for
 any sufficiently general $s\in H^0(X,L^m)$.
\end{itemize}
According to Proposition \ref{prop;07.10.13.11},
$(E,\nabla,h)$ is a good wild harmonic bundle,
and  $h$ is adapted to $\vecE_{\ast}$.
Thus, the proof of Lemma \ref{lem;07.10.16.15}
is finished.
\hfill\qed

\vspace{.1in}

Now, we can show the existence in the general case
by using the same argument
as that in the proof of Theorem 5.16 of \cite{mochi5},
and thus the proof of Theorem \ref{thm;06.1.23.100}
is finished.
\hfill\qed

\vspace{.1in}

We say that
a $\mu_L$-polystable good filtered flat bundle 
has trivial characteristic numbers,
if each $\mu_L$-stable component has
trivial characteristic numbers.

\begin{cor}
Let $(\vecE_{\ast},\nabla)$ be a good 
filtered flat bundle on $(X,D)$.
We put $E:=\vecE_{\ast|X-D}$.
It is $\mu_L$-polystable
with trivial characteristic numbers,
if and only if there exists a pluri-harmonic metric
$h$ of $(E,\nabla)$ adapted to 
the prolongment $\vecE_{\ast}$.
\end{cor}
\pf
The ``only if'' part follows from Theorem
\ref{thm;06.1.23.100}.
Let us see the ``if'' part.
Let $(E,\nabla,h)$ be a good wild harmonic bundle,
and let $(\vecE_{\ast},\nabla)$ is
the associated filtered flat bundle.
According to Proposition \ref{prop;07.12.6.11},
$(\vecE_{\ast},\nabla)$ is $\mu_L$-polystable
with $\pardeg_L(\vecE_{\ast})=0$.
Moreover,
if $(\vecE_{0\ast},\nabla_0)$
is a stable component of $(\vecE_{\ast},\nabla)$,
and then the restriction of $h$ to 
$E_0:=\vecE_{0\ast|X-D}$
gives the pluri-harmonic metric of $(E_0,\nabla_0)$
adapted to the prolongment $\vecE_{0\ast}$.
Then, we obtain
$\int_{X}\parch_{2,L}(\vecE_{0\ast})=0$
from Proposition \ref{prop;07.12.6.12}.
\hfill\qed

\vspace{.1in}
We also obtain the following corollary
on the uniqueness of pluri-harmonic metrics,
from Theorem \ref{thm;06.1.23.100}
and Proposition \ref{prop;07.12.6.11}.
\begin{cor}
\label{cor;08.1.22.15}
Let $(\vecE_{\ast},\nabla)$ be
a $\mu_L$-polystable good filtered flat bundle
on $(X,D)$ with the trivial characteristic numbers.
Let $(\vecE_{\ast},\nabla)
=\bigoplus_{i\in \Gamma}
 (\vecE_{i\ast},\nabla_i)\otimes\cnum^{m(i)}$
be the canonical decomposition,
where $(\vecE_{i\ast},\nabla_i)$ are 
$\mu_L$-stable.
We take a pluri-harmonic metric $h_i$
of $(E_i,\nabla_i)$ adapted to the prolongment
$\vecE_{i\ast}$ for each $i\in\Gamma$.
Then, any pluri-harmonic metric $h$
adapted to $\vecE_{\ast}$
is of the form $h=\bigoplus h_i\otimes g_i$,
where $g_i$ are Hermitian metrics of
$\cnum^{m(i)}$.
\hfill\qed
\end{cor}

\begin{rem}
Recall that,
for a given good filtered $\nabla$-flat bundle
$(\vecE_{\ast},\nabla)$,
we have a deformation
$(\vecE_{\ast}^{(2)},\nabla)$
as explained in Section 
{\rm\ref{subsection;10.5.17.51}}.
Although we considered a pluri-harmonic metric
$h$ adapted to $\vecE_{\ast}$,
it might be more natural
to consider a pluri-harmonic metric adapted
to $\vecE_{\ast}^{(2)}$.
See the proof of Theorem {\rm\ref{thm;07.10.14.75}}.
\hfill\qed
\end{rem}

\section{Applications to 
algebraic meromorphic flat bundles}
\label{subsection;08.9.28.101}
\subsection{Existence of a good model}

Let $X$ be an $n$-dimensional smooth 
connected projective variety
with a simple normal crossing hypersurface $D$.
We will prove the following theorem
in Section \ref{subsection;07.10.15.2}.

\begin{thm}
\label{thm;07.10.14.60}
Let $(\nbige,\nabla)$ be a meromorphic flat connection
on $(X,D)$.
Then, there exists a smooth projective variety $X'$
with a birational morphism
$\varphi:X'\lrarr X$ such that
(i) $D':=\varphi^{-1}(D)$ is simply normal crossing,
(ii) $X'-D'\simeq X-D$,
(iii) the Deligne-Malgrange lattice of 
$\varphi^{\ast}(\nbige,\nabla)$ 
is a locally free $\nbigo_X$-module
and good.
\end{thm}

Before going into the proof,
we give a direct corollary.
Let $K$ be a finite extension of $\cnum(x_1,\ldots,x_n)$.
Let $M$ be a finite dimensional $K$-vector space
with $\Der(K/\cnum)$-action.
We have the following corollary.

\begin{cor}
There exist a smooth projective variety $X$
with $K(X)\simeq K$,
and a meromorphic connection 
$(\nbige,\nabla)$ on $X$ with $(\nbige,\nabla)\otimes K\simeq M$,
such that the Deligne-Malgrange lattice of
$(\nbige,\nabla)$ is good.
\hfill\qed
\end{cor}

\begin{rem}
Very recently,
Kedlaya established the existence of
resolution of turning points
in a more general situation
with a completely different method.
See {\rm\cite{kedlaya}} and
{\rm\cite{kedlaya2}}.
\hfill\qed
\end{rem}

\subsection{Characterization of semisimplicity}

Let $X$ be an $n$-dimensional smooth 
connected projective variety,
and let $D$ be a simple normal crossing hypersurface of $X$.
Let $(\nbige,\nabla)$ be a meromorphic flat connection
on $(X,D)$.
We will prove the following theorem in
Section \ref{subsection;07.10.15.2}.
\begin{thm}
\label{thm;07.10.15.1}
$(\nbige,\nabla)$ is semisimple,
if and only if the following holds:
\begin{itemize}
\item
 Let $\varphi:\Xtilde\lrarr X$ be 
 a birational projective morphism
 such that the Deligne-Malgrange filtered bundle
 $\vecEtilde_{\ast}^{DM}$ associated to 
 $\varphi^{\ast}(\nbige,\nabla)$ is good,
 as in Theorem {\rm\ref{thm;07.10.14.60}}.
 Then, there exists a pluri-harmonic metric $h$
 of $(\nbigetilde,\nablatilde):=
 \varphi^{\ast}(\nbige,\nabla)_{|\Xtilde-\varphi^{-1}(D)}$
 which is adapted to $\vecEtilde_{\ast}^{DM}$.
 Note that $(\nbigetilde,\nablatilde,h)$
 is a $\sqrt{-1}\real$-good wild harmonic bundle.
\end{itemize}
\end{thm}

We give a complement on the uniqueness.
Such a metric $h$ is unique up to obvious ambiguity
as in Corollary {\rm\ref{cor;08.1.22.15}},
which follows from Proposition \ref{prop;07.12.6.11}.
\begin{prop}
\label{prop;10.6.18.2}
Let $(\nbige,\nabla)$ be semisimple,
and let $h$ be the metric as in 
Theorem {\rm\ref{thm;07.10.15.1}}.
Then,
 (i) if $(\nbige,\nabla)$ is simple,
 $h$ is uniquely determined up to constant multiplication,
 (ii) the canonical decomposition of $(\nbige,\nabla)$
 is orthogonal  with respect to $h$,
 (iii) if $(\nbige,\nabla)$ is a tensor product of
 a simple meromorphic flat bundle $(\nbige',\nabla')$
 and a vector space $V$,
 then $h$ is the tensor product of
 $h'$ for $(\nbige',\nabla')$
 and the metric of $V$.
\hfill\qed
\end{prop}

\subsection{Good Deligne-Malgrange lattice}

It seems better to give an explicit statement
for the filtered flat bundle
associated to a $\sqrt{-1}\real$-good wild
harmonic bundle,
which is implicitly used in the statement of
Theorem \ref{thm;07.10.15.1}.
Let $X$ be a complex manifold,
and let $D$ be a simple normal crossing hypersurface
of $X$.
Let $\harmonicbundle$ be 
a $\sqrt{-1}\real$-good wild harmonic bundle
on $(X,D)$.
\begin{prop}
\label{prop;08.12.18.10}
The filtered flat bundle
$(\nbigp_{\ast}\nbige^1,\DD^1)$
is the Deligne-Malgrange filtered bundle
associated to the meromorphic flat bundle
$(\nbigp\nbige^1,\DD^1)$.
In particular,
$\nbigp_0\nbige^1$ is the good Deligne-Malgrange
lattice of $(\nbigp\nbige^1,\DD^1)$.
\end{prop}
\pf
In the case $\dim X=1$,
the claim follows from the comparison
of the KMS-spectra in Proposition
\ref{prop;07.7.19.31}.
The general case can be reduced
to the above case.
\hfill\qed

\subsection{Pull back}

We give an application of Theorem
\ref{thm;07.10.15.1}.
Let $X$ and $Y$ be smooth projective varieties.
Let $F:Y\lrarr X$ be a rational morphism.
Let $(\nbigv,\nabla)$ be a semisimple meromorphic
flat connection on $X$.
Assume that $F(Y)$ is not contained
in the pole of $(\nbigv,\nabla)$.
We have the meromorphic flat bundle
$F^{\ast}(\nbigv,\nabla)$ on $Y$.

\begin{thm}
$F^{\ast}(\nbigv,\nabla)$ is also semisimple.
\end{thm}
\pf
By replacing $X$ birationally,
we may assume the following:
\begin{itemize}
\item
 We have a simple normal crossing hypersurface $D_X$
 such that
 $(\nbigv,\nabla)$ is a semisimple meromorphic flat bundle
 on $(X,D_X)$.
\item
 The Deligne-Malgrange filtered flat bundle 
 $\vecV^{DM}_{X\,\ast}$
 associated to $(\nbigv,\nabla)$ is good.
\end{itemize}
We take a good wild pluri-harmonic metric $h_X$
of $(\nbigv,\nabla)_{|X-D_X}$, adapted to
$\vecV^{DM}_{X\,\ast}$.
Let $(E_X,\delbar_{E_X},\theta_X,h_X)$ denote 
the corresponding good wild harmonic bundle
on $X-D_X$.
The associated meromorphic flat bundle
$(\nbigp\nbige_X^1,\DD_X^1)$ is naturally isomorphic to
$(\nbigv,\nabla)$.
By construction, $\theta$ is good and wild.
Moreover, the eigenvalues of the endomorphism
$\Res_{D_{X,i}}(\theta)$ on 
$\lefttop{i}\Gr^F\bigl(
 \prolong{E_X}\bigr)$ on $D_{X,i}$
are purely imaginary.

By replacing $Y$ birationally,
we may assume to have 
a regular morphism $F:Y\lrarr X$
such that $D_Y:=F^{-1}(D_X)$ is 
a simple normal crossing hypersurface.
\begin{lem}
\label{lem;07.10.15.10}
$(E_Y,\delbar_{E_Y},\theta_Y,h_Y):=
 F^{-1}(E_X,\delbar_{E_X},\theta_X,h_X)$
is a $\sqrt{-1}\real$-good wild harmonic bundle on $Y-D_Y$.
\end{lem}
\pf
Note that $F^{\ast}(\prolong{E_X})$ gives a prolongment
of $E_Y$ to a locally free $\nbigo_Y$-module on $Y$.
Let $P$ be any point of $D_Y$.
We take a small holomorphic coordinate neighbourhood 
$V$ around $P$.
We also take a small holomorphic coordinate neighbourhood
$U$ around $F(P)$.
We put $D_U:=D_X\cap U$ and $D_V:=D_Y\cap V$.
We can take a ramified covering
$\varphi_U:(\Utilde,\Dtilde_U)\lrarr (U,D_U)$
such that we have the decomposition:
\begin{equation}
 \label{eq;07.10.15.5}
 \varphi_U^{\ast}(\prolong{E_X},\theta_X)  
 =\bigoplus_{\gminia\in \Irr(\theta_X,F(P))}
 \bigl(
 E_{\gminia},\theta_{\gminia}
\bigr)
\end{equation}
Here, $\theta_{\gminia}-d\gminia\cdot\id_{E_{\gminia}}$
are logarithmic with respect to 
$\varphi_U^{\ast}\bigl(\prolong{E_X}\bigr)$
for any $\gminia\in \Irr(\theta_X,F(P))$.
We can take a ramified covering
$\varphi_V:(\Vtilde,\Dtilde_V)\lrarr (V,D_V)$
such that the composite
$F\circ \varphi_V:\Vtilde\lrarr U$
factors through $\Utilde$,
i.e.,
there exists $\Ftilde:\Vtilde\lrarr \Utilde$
such that
$F\circ\varphi_V=\varphi_U\circ \Ftilde$.
From (\ref{eq;07.10.15.5}),
we obtain the decomposition around $P$:
\[
 \varphi_V^{\ast}(F^{\ast}(\prolong{E_X}),\theta_Y) 
=\bigoplus_{\gminia} 
 \bigl(\Ftilde^{\ast}(E_{\gminia}),\theta_{Y,\gminia}\bigr)
\]
And,
$\theta^{\reg}_{\Vtilde}:=
 \bigoplus\bigl(
 \theta_{Y,\gminia}
-d\Ftilde^{\ast}\gminia
 \cdot\id_{\Ftilde^{\ast}(E_{\gminia})}
 \bigr)$
is logarithmic.

Let us show that the eigenvalues of
the residue of $\theta_Y$ are constant on
each component of $D_Y$, and that 
they are purely imaginary.
Let $P$ be a smooth point of $D_Y$.
Let $V$ and $U$ be as above.
Let $(w_1,\ldots,w_m)$
and $(z_1,\ldots,z_n)$ be holomorphic coordinate
on $\Vtilde$ and $\Utilde$ respectively,
such that $\Dtilde_{V}=\{w_1=0\}$
and $\Dtilde_U=\bigcup_{j=1}^{\ell}\{z_j=0\}$.
Let $d_i$ $(i=1,\ldots,\ell)$ be the order of
$\Ftilde^{\ast}(z_i)$ with respect to $w_1$.
Let $Q$ be any point of $\Dtilde_V$.
Then, 
$\Res_{w_1}(\varphi_V^{\ast}\theta_Y)_{|Q}$ on 
$\varphi_V^{\ast}F^{\ast}(\prolong{E_X})_{|Q}$ 
is given by
 $\sum_{i=1}^{\ell}
 d_i\cdot \Res_{z_i}(\varphi_U^{\ast}\theta_X)_{|F(Q)}$.
Hence, we obtain that
the eigenvalues of
$\Res_{w_1}(\varphi_V^{\ast}\theta_Y)$
are constant on $\Dtilde_V$,
and they are purely imaginary.

We have the expression
$\theta^{\reg}_{\Vtilde}
 =f_1\cdot dw_1/w_1+\sum_{j=2}^m f_j\cdot dw_j$.
Then, we can conclude that 
the coefficients of $P(t):=\det(t-f_1)_{|\Dtilde_{V}}$ 
are constant,
and that the solutions of $P(t)=0$
are purely imaginary.
Thus, the proof of Lemma \ref{lem;07.10.15.10}
is finished.
\hfill\qed

\vspace{.1in}

We have the meromorphic flat bundle
$(\nbigp\nbige_Y^1,\DD^1_Y)$
associated to 
$(E_Y,\delbar_{E_Y},\theta_Y,h_Y)$.
We also have the good lattices
$\nbigp_{\veca}\nbige_Y^1$ of $\nbigp\nbige_Y^1$.
We have the naturally defined  morphism
$F^{\ast}(\nbigp_0\nbige_X^1,\DD_X^1)
\lrarr (\nbigp_0\nbige_Y^1,\DD_Y^1)$,
which induces the morphism
$\iota:F^{\ast}(\nbigp\nbige_X^1)
\lrarr \nbigp\nbige_Y^1$
of locally free $\nbigo_Y(\ast D_Y)$-modules.
Since the restriction of $\iota$ to $Y-D_Y$ is 
an isomorphism,
$\iota$ is an isomorphism on $Y$.
Hence, $F^{\ast}(\nbigv,\nabla)$ has the good lattices
$\nbigp_{\veca}\nbige^1_Y$.
By using Proposition \ref{prop;07.7.19.31},
we can show that 
$\nbigp_{\ast}\nbige^1_Y
 =\bigl(\nbigp_{\veca}\nbige^1_Y\,\big|\,
 \veca\in\real^{\ell}\bigr)$
is the Deligne-Malgrange filtered flat bundle
associated to $F^{\ast}(\nbigv,\nabla)$.
Then, we can conclude that
$F^{\ast}(\nbigv,\nabla)$ is semisimple
due to Theorem \ref{thm;07.10.15.1}.
\hfill\qed

\section{Proof of Theorems \ref{thm;07.10.14.60}
 and \ref{thm;07.10.15.1}}
\label{subsection;07.10.15.2}
\subsection{Preliminary characterization of semisimplicity}
\label{subsection;07.11.9.35}

Let $X$ be an $n$-dimensional 
smooth connected projective variety,
and let $D$ be a simple normal crossing hypersurface 
of $X$.
Let $(\nbige,\nabla)$ be a flat meromorphic connection
on $(X,D)$.

\begin{prop}
\label{prop;07.6.12.50}
If $(\nbige,\nabla)$ is simple,
there exists a pluri-harmonic metric $h$
of $(E,\nabla):=(\nbige,\nabla)_{|X-D}$ 
with the following property,
which is unique up to positive constant multiplication:
\begin{itemize}
\item
 Let $P$ be any smooth point of $D$
 around which $(\nbige,\nabla)$ is good.
 Then, $(E,\nabla,h)$ is 
 a $\sqrt{-1}\real$-good wild harmonic bundle
 around $P$,
 and $h$ is adapted to the Deligne-Malgrange filtered bundle
 $\vecE^{DM}_{\ast}$ around $P$.
\end{itemize}
In particular,
if the Deligne-Malgrange filtered 
bundle $\vecE^{DM}_{\ast}$
associated to $(\nbige,\nabla)$ is good,
then
$h$ is adapted to $\vecE^{DM}_{\ast}$,
and $(E,\nabla,h)$ is 
a $\sqrt{-1}\real$-good wild harmonic bundle
on $(X,D)$.
\end{prop}
\pf
First, let us consider the case $\dim X=2$.
We can take an appropriate 
birational projective morphism
$\varphi:\Xtilde\lrarr X$ such that
$\varphi^{\ast}(\nbige,\nabla)$
has the good formal structure.
(See \cite{mochi6}.)
Note that 
$\varphi^{\ast}(\nbige,\nabla)$ is also simple.
So, we may assume that $(\nbige,\nabla)$
has the good formal structure from the beginning.
Then, the associated Deligne-Malgrange filtered bundle
is good (Proposition \ref{prop;10.5.4.3}).
Because
$\varphi^{\ast}(\nbige,\nabla)$ is simple,
the associated Deligne-Malgrange filtered
bundle is $\mu_L$-stable.
We also have the vanishing of
the parabolic characteristic numbers,
according to Corollary \ref{cor;08.9.30.3}.
Hence, 
we obtain the desired pluri-harmonic metric
due to Theorem \ref{thm;06.1.23.100}
(or Lemma \ref{lem;07.10.16.15}).

The higher dimensional case
can be reduced to the surface case 
by using Corollary \ref{cor;08.9.30.11}
and an argument in the proof of 
Theorem 5.16 of \cite{mochi5}.
\hfill\qed

\subsection{Pull back via
the birational projective morphism}
\label{subsection;08.12.18.1}

Let $\varphi:X'\lrarr X$ be any 
birational projective morphism.
If $(\nbige,\nabla)$ is simple,
then $\varphi^{\ast}(\nbige,\nabla)$ is also simple.
Therefore, we have pluri-harmonic metrics
$h_{\nbige}$ and $h_{\varphi^{\ast}\nbige}$
for $(\nbige,\nabla)$ and $\varphi^{\ast}(\nbige,\nabla)$
as in Proposition \ref{prop;07.6.12.50}.

\begin{lem}
\label{lem;07.10.14.100}
$\varphi^{\ast}h_{\nbige}=
 a\cdot h_{\varphi^{\ast}\nbige}$
for some $a>0$.
\end{lem}
\pf
Let $C$ be a sufficiently generic and ample
curve in $X$.
Then, $h_{\nbige|C}$ and $h_{\varphi^{\ast}\nbige|C}$
are the pluri-harmonic metrics
for the simple meromorphic connection
$(\nbige,\nabla)_{|C}$.
Thus, we have 
$\varphi^{\ast}h_{\nbige|C}=
 a_C\cdot h_{\varphi^{\ast}\nbige|C}$
because of the uniqueness.
Since such curves $C$ cover Zariski open subsets of
$X$ and $X'$,
the claim of the lemma follows.
\hfill\qed

\subsection{Proof of Theorem \ref{thm;07.10.14.60}
 in the case that $(\nbige,\nabla)$ is simple}

\label{subsection;07.10.14.110}

Assume that $(\nbige,\nabla)$ is simple.
We can take a pluri-harmonic metric $h_{\nbige}$
as in Proposition \ref{prop;07.6.12.50}.
Let $(E,\delbar_E,\theta,h_{\nbige})$ be 
the corresponding harmonic bundle on $X-D$.

\begin{lem}
\label{lem;08.4.3.30}
There exists a birational projective morphism
$F:X_1\lrarr X$ with the following property:
\begin{itemize}
\item
 $D_1:=F^{-1}(D)$ is a simple normal crossing hypersurface,
 and $X_1-D_1\simeq X-D$.
\item
 We put $(E_1,\delbar_{E_1},h_1,\theta_1):=
 F^{\ast}(E,\delbar_E,h,\theta)$.
 Then, $\theta_1$ is generically $\sqrt{-1}\real$-good,
 and (SRP) holds.
 (See Section {\rm\ref{subsection;07.11.1.16}}.)
\end{itemize}
\end{lem}
\pf
We take a covering 
$X=\bigcup\nbigu^{(p)}$
by affine Zariski open subsets 
$\nbigu^{(p)}=\Spec R^{(p)}$
with etale morphisms
$\varphi_p:\nbigu^{(p)}\lrarr \cnum^n$
such that
$\varphi_p^{-1}\bigl(
 \bigcup_{i=1}^{\ell(p)}\{z_i=0\}
 \bigr)=\nbigu_p\cap D$.
Let $g^{(p)}=\prod_{i=1}^{\ell(p)}\varphi_p^{-1}(z_i)$
be the defining equation of $D\cap \nbigu^{(p)}$.
Let $R^{(p)}_{g^{(p)}}$
denote the localization of $R^{(p)}$
with respect to $g^{(p)}$.
Let $\tau_i^{(p)}:=\varphi_p^{\ast}dz_i/z_i$ 
for $i=1,\ldots,\ell(p)$,
and $\tau_i^{(p)}:=\varphi_p^{\ast}dz_i$ 
for $i=\ell(p)+1,\ldots,n$.
They give a frame of
$\Omega^1(\log D)$ on $\nbigu^{(p)}$.
We have the expression
$\theta=\sum f_i^{(p)}\tau_i^{(p)}$.
We obtain the characteristic polynomials
$\nbigp^{(p)}_i(T):=
 \det\bigl(T-f_i^{(p)}\bigr)$.
Let $P$ be any point of $D$ around which
$(E,\delbar_E,\theta,h)$ is good wild.
Then, the coefficients of
$\nbigp^{(p)}_i(T)$ are meromorphic around $P$
whose poles are contained in $D$.
Hence, we obtain
$\nbigp^{(p)}_i(T)\in R_{g^{(p)}}^{(p)}[T]$.

Applying Proposition \ref{prop;08.4.3.1}
to each $\nbigp^{(p)}_i$,
we obtain a smooth projective variety $X_{1,i}^{(p)}$
with a regular projective birational morphism
$F_{1,i}^{(p)}:X_{1,i}^{(p)}\lrarr X$
such that $(F_{1,i}^{(p)})^{\ast}\nbigp^{(p)}_i$
satisfies the condition in Proposition \ref{prop;08.4.3.1}.
Taking the fiber products of
$X_{1,i}^{(p)}$ over $X$,
the normalization,
and the resolution of the singularity,
we obtain 
a birational projective morphism
$F_1:X_1\lrarr X$
such that any
$F_1^{\ast}\nbigp^{(p)}_i$
satisfy the condition in Proposition \ref{prop;08.4.3.1}.
The morphism factors through each $X_{1,i}^{(p)}$.

Due to Lemma \ref{lem;07.10.14.100},
$\theta_1$ is generically $\sqrt{-1}\real$-good.
Let us show that (SRP) holds.
We put $\nbigu_1^{(p)}:=F_1^{-1}(\nbigu^{(p)})$.
For any $P\in D_1\cap\nbigu_1^{(p)}$,
we take a holomorphic coordinate
$(U,z_1,\ldots,z_n)$ around $P$
such that
$U\cap D_1=\bigcup_{i=1}^{\ell}\{z_i=0\}$.
We take a ramified covering
$\psi_P:\Utilde\lrarr U$ 
given by
$\psi_P(\zeta_1,\ldots,\zeta_n)
=(\zeta_1^M,\ldots,\zeta_{\ell}^M,
 \zeta_{\ell+1},\ldots,\zeta_n)$,
where $M$ is divisible by $r!$.
We put $\Dtilde:=\psi_P^{-1}D_1$.
Let $\psitilde_P:=F_1\circ\psi_P$.
By construction,
the eigenvalues of $\psitilde_P^{\ast}f_i^{(p)}$
are of the form
$\beta+\gamma$,
where
$\beta\in\nbigo_{\Utilde}(\ast\Dtilde)$
and $\gamma$ are multivalued holomorphic
functions on $\Utilde$.
For the expression
$\psi_P^{\ast}\theta
=\sum_{i=1}^{\ell}f_i\cdot d\zeta_i/\zeta_i
+\sum_{i=\ell+1}^nf_i\cdot d\zeta_i$,
we have
$f_i=\sum_j a_{i,j}\cdot 
 \psitilde_p^{\ast}f_j^{(p)}$,
where $a_{i,j}\in H(\Utilde_P)$.
Note the commutativity of $f_j^{(p)}$.
Then, it is easy to observe that
(SRP) holds for $\theta_1$.
Thus we obtain Lemma \ref{lem;08.4.3.30}.
\hfill\qed

\vspace{.1in}

Hence, we may and will assume that
(SRP) holds from the beginning.
Applying Proposition \ref{prop;07.11.1.3}
to $\theta$,
we can take a birational projective morphism
$\varphi:X'\lrarr X$ such that
$\varphi^{\ast}\theta$ is 
$\sqrt{-1}\real$-good,
i.e.,
$(\Etilde,\nablatilde,\htilde):=
 \varphi^{\ast}(E,\nabla,h)$
is a $\sqrt{-1}\real$-good wild harmonic bundle.
Then, we have the associated prolongment
$(\vecEtilde_{\ast},\nablatilde)$, which is good.
Due to Lemma \ref{lem;07.10.14.100} 
and Proposition \ref{prop;07.6.12.50},
it is the same as the Deligne-Malgrange filtered flat sheaf
of $\varphi^{\ast}\nbige$ 
around a general point $P$ of the smooth part of $D$.
Hence, we can conclude that
$\vecEtilde_{\ast}$ is the same as the Deligne-Malgrange
filtered sheaf associated 
to $(\varphi^{\ast}\nbige,\nablatilde)$,
and it is good.
Thus, we obtain the desired resolution
in the case that $(\nbige,\nabla)$ is simple.

\subsection{End of proof of
 Theorems \ref{thm;07.10.14.60} and
  \ref{thm;07.10.15.1}}
\label{subsection;07.10.14.111}

Let us consider the general case.
A meromorphic flat connection
$(\nbige,\nabla)$ is the extension
of simple meromorphic flat connections
$(\nbige_i,\nabla_i)$ $(i=1,\ldots,N)$.
We can take a birational projective morphism
$\varphi:X'\lrarr X$ such that 
(i) $D':=\varphi^{-1}(D)$ is simple normal crossing,
(ii) $X'-D'\simeq X-D$,
(iii) any $\varphi^{\ast}(\nbige_i,\nabla_i)$
have the good Deligne-Malgrange lattices.
We may also assume that
the union of the set of the irregular values of 
$\varphi^{\ast}(\nbige_i,\nabla_i)$ are also good
around each point of $\varphi^{\ast}D$.
(For example,
we have only to apply Proposition \ref{prop;07.11.1.3},
to the meromorphic Higgs bundle
 corresponding to
 $\bigoplus(\nbige_i,\nabla_i)$.)
Then, 
the Deligne-Malgrange lattice of
$\varphi^{\ast}(\nbige,\nabla)$ is good,
according to Corollary \ref{cor;10.5.4.101}.
Thus, the proof of Theorem \ref{thm;07.10.14.60}
is finished.

\vspace{.1in}
Let us show Theorem \ref{thm;07.10.15.1}.
Let $(\nbige,\nabla)$ be a meromorphic flat connection
on $(X,D)$.
We take a birational projective morphism 
$\varphi:\Xtilde\lrarr X$
 such that the Deligne-Malgrange filtered bundle
 $\vecEtilde_{\ast}^{DM}$ associated to 
 $(\nbigetilde,\nablatilde):=
 \varphi^{\ast}(\nbige,\nabla)$ is good
(Theorem \ref{thm;07.10.14.60}).
Note that $(\nbige,\nabla)$ is semisimple,
if and only if $(\nbigetilde,\nablatilde)$
is semisimple.
If $(\nbige,\nabla)$ is semisimple,
we can take a pluri-harmonic metric $h$
with the desired property
due to Proposition \ref{prop;07.6.12.50}.
Conversely, if we have a pluri-harmonic metric
$h$ of $(\nbigetilde,\nablatilde)$ adapted to
$\vecEtilde_{\ast}^{DM}$,
then $(\vecEtilde_{\ast}^{DM},\nablatilde)$
is $\mu_L$-polystable
due to Proposition \ref{prop;07.12.6.11}.
By Lemma \ref{lem;08.9.30.10},
we obtain that $(\nbigetilde,\nablatilde)$
is semisimple,
and thus $(\nbige,\nabla)$ is semisimple.
\hfill\qed

\section{Minor refinement
 of the result in Section 
 \ref{subsection;08.9.28.101}}
\label{subsection;08.12.18.3}
Let $X$ be a smooth proper algebraic
variety over $\cnum$.
Let $D$ be a simply normal crossing
hypersurface of $X$.
Let $(\nbige,\nabla)$ be a meromorphic flat bundle.
We set $(E,\nabla):=(\nbige,\nabla)_{|X-D}$.

\begin{prop}
\label{prop;08.12.18.2}
\mbox{{}}
$(\nbige,\nabla)$ is semisimple,
if and only if
there exists a pluri-harmonic metric
$h$ of $(E,\nabla)$ such that
the following holds:
\begin{itemize}
\item
 Let $P$ be any smooth point of $D$
 around which $(\nbige,\nabla)$ is good.
 Then, $(E,\nabla,h)$ is 
 a $\sqrt{-1}\real$-good wild harmonic bundle
 around $P$,
 and $h$ is adapted to 
 the Deligne-Malgrange filtered bundle
 $\vecE^{DM}_{\ast}$ around $P$.
\end{itemize}
Such a metric is unique up to obvious ambiguity,
in the following sense:
\begin{itemize}
\item
If $(\nbige,\nabla)$ is simple,
it is unique up to positive constant multiplication.
\item
Let $(\nbige,\nabla)=
\bigoplus_{i}(\nbige_i,\nabla_i)\otimes V_i$
be the canonical decomposition,
i.e., $(\nbige_i,\nabla_i)$ are simple,
$V_i$ are vector spaces,
and $(\nbige_i,\nabla_i)\not\simeq
 (\nbige_j,\nabla_j)$ for $i\neq j$.
Then, the metric $h$ is of the form
$\bigoplus h_i\otimes g_i$,
where $h_i$ is a pluri-harmonic metric
for $(\nbige_i,\nabla_i)$ as above,
and $g_i$ are metrics of $V_i$.
\end{itemize}
\end{prop}
\pf
Assume that $(\nbige,\nabla)$ is semisimple.
By using Chow's lemma and 
Theorem \ref{thm;07.10.14.60},
we can take a projective birational morphism
$\varphi:X'\lrarr X$ such that 
(i) $D':=\varphi^{-1}(D)$ is simply normal crossing,
(ii) the Deligne-Malgrange lattice of 
$(\nbige',\nabla'):=
 \varphi^{\ast}(\nbige,\nabla)$ 
is a locally free $\nbigo_X$-module
and good.
We put
$(E',\nabla'):=(\nbige',\nabla')_{|X'-D'}$.
By using Theorem \ref{thm;06.1.23.100}
or Theorem \ref{thm;07.10.15.1},
we can take a pluri-harmonic metric $h'$
of $(E',\nabla')$ such that
$\nbigp_{\ast}E'$
is the same as the Deligne-Malgrange filtered
bundle of $(\nbige',\nabla')$.
There is a Zariski closed subset $W\subset X$
such that $\varphi$ induces
$X'\setminus \varphi^{-1}(W)\simeq
 X\setminus W$.

\begin{lem}
We have the pluri-harmonic metric $h$
of $(E,\nabla)$ such that
$h_{|X\setminus W}=
h'_{|X'\setminus \varphi^{-1}(W)}$
under the natural identification.
\end{lem}
\pf
According to Theorem 6.1.3 of \cite{sabbah2},
$(E',\nabla',h')$ naturally induces
a pure twistor $D$-module $\nbigt(E')$
of weight $0$
with the polarization $(\id,\id)$ 
on $X'\setminus D'$.
By Hard Lefschetz Theorem
for regular polarized pure twistor $D$-module
(Theorem 6.1.1 of \cite{sabbah2}),
we obtain the pure twistor $D$-module
$\varphi_{\dagger}^0\nbigt(E')$
of weight $0$ on $X\setminus D$.
Let $\nbigt_E$ be the direct summand of
$\varphi_{\dagger}^0\nbigt(E')$
whose strict support is $X\setminus D$.
It is easy to observe that
the underlying $D$-module is
naturally isomorphic to $(E,\nabla)$,
because of the regularity.
Moreover, by the Hard Lefschetz theorem,
$\nbigt_E$ is equipped with the polarization
whose restriction to $X\setminus W$
is the same as the polarization
of $\nbigt(E')$ under the natural identification.
Thus, we obtain a pluri-harmonic metric $h$
of $(E,\nabla)$.
\hfill\qed

\vspace{.1in}

If $P\in D\setminus W$,
$(\nbige,\nabla)$ is good around $P$,
and $h$ is adapted to $\vecE_{\ast}^{DM}$
by construction.
Let $P$ be any point of $D$
around which $(\nbige,\nabla)$ is good.
By applying Proposition \ref{prop;07.10.13.11},
we obtain that $h$ is adapted to $\vecE_{\ast}^{DM}$.
Hence, $h$ has the desired property.

Conversely,
assume that $(\nbige,\nabla)$ is equipped with
a metric $h$ as above.
Let $\nbige_1\subset \nbige$ be 
an $\nbigo_X(\ast D)$-submodule
such that 
$\nabla\nbige_1\subset
 \nbige_1\otimes\Omega_X^1(\ast D)$.
We put $E_1:=\nbige_{1|X\setminus D}$,
which is a flat subbundle of 
$E=\nbige_{|X\setminus D}$.
Let $E_2$ denote the orthogonal complement of
$E_1$ in $E$ with respect to $h$.
\begin{lem}
\label{lem;10.6.19.3}
$E_2$ is a flat subbundle of $E$.
\end{lem}
\pf
Let us consider a morphism $F$ from
smooth projective curve $C$ to $X$ such that
(i) $F(C)$ is not contained in $D$,
(ii) $(\nbige,\nabla)$ is good 
around each point of $F(C)\cap D$.
We have $F^{\ast}(\nbige_1,\nabla)
\subset F^{\ast}(\nbige,\nabla)$,
and $F^{\ast}h$.
Note that semisimplicity is related with
the polystability as in Lemma \ref{lem;08.9.30.10}.
Applying Corollary \ref{cor;08.1.22.15},
we obtain that $F^{\ast}E_2$
is flat with respect to $F^{\ast}\nabla$.
Because a Zariski open set of $X$
is covered by such curves,
we obtain that $E_2$ is also flat.
\hfill\qed

\vspace{.1in}
Let $\vecE_{\ast}^{DM}$
and $\vecE_{1\ast}^{DM}$
be the Deligne-Malgrange filtered sheaf
associated to $(\nbige,\nabla)$
and $(\nbige_1,\nabla)$, respectively.
Let $\pi$ be the orthogonal projection of $E$
onto $E_1$.
Let us observe that it is extended to
a projection $\vecE_{\ast}\lrarr\vecE_{1\ast}$.
Let $P$ be a smooth point of $D$
around which $(\nbige,\nabla)$ is good.
Because
$\vecE^{DM}_{\ast}\simeq
 \nbigp_{\ast}E_1\oplus\nbigp_{\ast}E_2$
around $P$,
$\pi$ has the desired extension
around $P$.
By using the Hartogs property of 
reflexive sheaves,
we obtain the desired extension
on $X$.
It implies that
$E_2$ is extended to a meromorphic 
flat subbundle $\nbige_2$ of $\nbige$,
and we have the flat decomposition
$\nbige=\nbige_1\oplus\nbige_2$.
Moreover, 
$\nbige_i$ are also equipped
with the metrics $h_i$ satisfying the condition
in the statement of Proposition 
\ref{prop;08.12.18.2}.
Hence, by an easy inductive argument,
we can conclude that $\nbige$
is semisimple.
We can also deduce that
the canonical decomposition of
$(\nbige,\nabla)$ is orthogonal
with respect to $h$.

Let $h$ and $h'$ be metrics for $(\nbige,\nabla)$
as in Proposition \ref{prop;08.12.18.2}.
Let $s$ be the endomorphism of $E$
determined by $h=h'\cdot s$,
which is self-adjoint with respect to
both $h$ and $h'$.
Let $C$ and $F$ be as in the proof of
Lemma \ref{lem;10.6.19.3}.
By applying Corollary \ref{cor;08.1.22.15},
we obtain that $F^{\ast}s$ is $\nabla$-flat.
Hence, we obtain that $s$ is $\nabla$-flat.
It induces the $\nabla$-flat decomposition
$E=\bigoplus E_i$ 
such that 
(i) it is orthogonal with respect to
both $h$ and $h'$,
(ii) $h_{E_i}=a_i\,h'_{E_i}$
for some $a_i>0$.
By using the argument in the previous paragraph,
we can show that the decomposition
is extended to that
of meromorphic flat bundles
 $\nbige=\bigoplus\nbige_i$.
In particular,
we obtain the uniqueness
in the case that $(\nbige,\nabla)$ is simple.
By using the argument in the proof of
Corollary \ref{cor;10.6.19.2},
we obtain that the restriction of $h$
to a direct summand of the canonical decomposition
has the form as in the statement of Proposition 
\ref{prop;08.12.18.2}.
Thus, the proof of Proposition \ref{prop;08.12.18.2}
is finished.
\hfill\qed

\begin{cor}
\label{cor;08.12.18.5}
There exists a birational projective morphism
$\varphi:X'\lrarr X$ such that
(i) $D':=\varphi^{-1}(D)$ is simply normal crossing,
(ii) $X'-D'\simeq X-D$,
(iii) the Deligne-Malgrange lattice of 
$\varphi^{\ast}(\nbige,\nabla)$ 
is a locally free $\nbigo_X$-module
and good.
Namely,
Theorem {\rm\ref{thm;07.10.14.60}}
holds without the assumption that
$X$ is projective.
\end{cor}
\pf
We have only to apply the argument in
Subsections 
\ref{subsection;08.12.18.1}--\ref{subsection;07.10.14.111},
by replacing Proposition \ref{prop;07.6.12.50}
with Proposition \ref{prop;08.12.18.2}.
\hfill\qed

\part{Application to Wild Pure Twistor $D$-Modules}
\label{part;08.9.29.130}

\chapter{Wild Pure Twistor $D$-Modules}
In Section \ref{subsection;08.9.8.40},
we recall the notion of wild pure twistor
$D$-modules and their basic property
due to Sabbah (\cite{sabbah2} and \cite{sabbah5}).
Note that we will give a review on 
$\nbigr$-modules, $\nbigr$-triples
and variants in Chapter \ref{section;08.10.20.40}.

We show in Section 
\ref{subsection;08.10.18.120}
the correspondence between
wild pure twistor $D$-modules
and wild harmonic bundles on curves.
Although,
we will study a more general correspondence
in Chapter \ref{section;08.9.29.150},
we need the correspondence in the curve case
for the proof of Hard Lefschetz theorem
in Chapter \ref{section;08.10.18.121}.

In Section \ref{subsection;08.1.18.13},
we study Gysin map for wild pure twistor $D$-modules,
which is a preparation for 
Section \ref{subsection;07.10.23.21}.

\section{Review of  wild pure twistor $D$-modules
 due to Sabbah}
\label{subsection;08.9.8.40}
\subsection{Wild pure twistor $D$-modules
 and some property}

We recall the notions of
pure twistor $D$-module 
and wild pure twistor $D$-modules
defined in \cite{sabbah2} and \cite{sabbah5}.
See those papers for more details and precision.
(See also \cite{saito1} and \cite{saito2}
 for the original work due to Saito
 on pure Hodge modules.)
First, we recall the definition of pure twistor $D$-module
in \cite{sabbah2}.
Let $X$ be a complex manifold.

\begin{df}
\label{df;11.5.10}
\index{pure twistor $D$-module}
Let $w$ be an integer.
The category $\MT_{\leq d}(X,w)$ is defined to be
the full subcategory of the category of $\nbigr_X$-triples
whose objects $\nbigt$ satisfy the following conditions:
\index{category $\MT_{\leq d}(X,w)$}
\begin{description}
\item[($HSD$)]
 The underlying $\nbigr_X$-modules of $\nbigt$
 are holonomic and strictly $S$-decomposable.
 The dimensions of their strict supports are less than $d$.
\item[($MT_{>0}$)]
 Let $U$ be any open subset of $X$,
 and $f$ be any holomorphic function on $U$.
 Let $u$ be any element of
 $\real\times\cnum$,
 and $l$ be any integer.
 Then the induced $\nbigr$-triples
 $\Gr^W_l\psitilde_{f,u}\nbigt_{|U}$
 are objects of $\MT_{\leq d-1}(U,w+l)$.
 Here $W$ denotes the weight filtration
 of the induced nilpotent map on $\psitilde_{f,u}\nbigt_{|U}$.
\item[($MT_{0}$)]
 By (HSD),
 we have the decomposition
 $\nbigt=\bigoplus_Z \nbigt_Z$
 by the strict supports,
 where $Z$ runs through the irreducible closed
 subsets of $X$.
 In the case $\dim Z=0$,
 $\nbigt_{Z}$ is a push forward of a pure twistor
 structure in dimension $0$ via the inclusion
 $Z\lrarr X$.
\end{description}
In the case $d\geq \dim X$,
we denote the category by $\MT(X,w)$.
Any objects of $\MT(X,w)$ are called 
pure twistor $D$-modules with weight $w$ on $X$.
\index{category $\MT(X,w)$}
\hfill\qed
\end{df}
See \cite{sabbah2}
for the fundamental property of pure twistor $D$-modules.
Note that $\nbigt_Z$ with $\dim Z=0$
can be obtained from $\nbigt$
by a successive use of the vanishing cycle functor.

Let $\nbiga$ be a $\rnum$-vector subspace of $\cnum$.
\begin{df}
Let $\nbigt$ be a pure twistor D-module on $X$.
We have the decomposition
$\nbigt=\bigoplus\nbigt_Z$ by the strict supports.
For the definition of $\nbiga$-wild pure twistor $D$-modules,
we impose the following condition
inductively on the dimension of the supports of $Z$:
\index{$\nbiga$-wild pure twistor $D$-module}
\index{wild pure twistor $D$-module}
\index{$\sqrt{-1}\real$-wild pure twistor $D$-module}
\begin{description}
\item[($MT_{>0}(\ramification,\exp,\nbiga)$)]
 Let $U$ be any open subset of $X$,
 and $g$ be a holomorphic function on $U$.
 For any $n\in\seisuu_{>0}$ and $\gminia\in\cnum[t_n^{-1}]$,
$\nbigt_Z$ is strictly specializable along $g$
 with ramification and exponential twist by $\gminia$.
 And
 the induced $\nbigr_X$-triples
 $\Gr^{W(N)}_{\ell}\psitilde_{g,\gminia,u}\bigl(
 \nbigt_Z\bigr)$ 
 are wild pure twistor $D$-modules with weight $w+\ell$ 
 for any $u\in\real\times\cnum$ and $\ell\in\seisuu$.
 Moreover,
 $\psitilde_{g,\gminia,u}(\nbigt_Z)=0$
 unless $u\in\real\times\nbiga$.
\hfill\qed
\end{description}
\end{df}

Let $\MT^{\wild}(X,w,\nbiga)$ 
denote the category of 
$\nbiga$-wild pure twistor $D$-modules
on $X$ with weight $w$.
Let  $\MT_{\leq d}^{\wild}(X,w,\nbiga)$ denote
the category of 
$\nbiga$-wild pure twistor $D$-modules on $X$ 
with weight $w$
such that the dimensions of 
the strict supports are less than $d$.
For a subvariety $Y\subset X$,
let $\MT^{\wild}_Y(X,w,\nbiga)$ denote the category of 
$\nbiga$-wild pure twistor
$D$-modules on $X$ with weight $w$
whose strict supports are contained in $Y$.
\index{category $\MT^{\wild}(X,w,\nbiga)$}
\index{category $\MT_{\leq d}^{\wild}(X,w,\nbiga)$}
\index{category $\MT^{\wild}_Y(X,w,\nbiga)$}

\begin{rem}
In the following argument,
we often omit to distinguish $\nbiga$.
In {\rm\cite{mochi2}},
$\sqrt{-1}\real$-wild pure twistor $D$-module
satisfying the regularity condition is
called a pure imaginary pure twistor $D$-module.
\hfill\qed
\end{rem}

We recall some properties of
wild pure twistor $D$-modules.
We refer to \cite{sabbah2} and \cite{sabbah5}
for more details.

\begin{lem}
 Kashiwara's equivalence holds.
 Namely,
 $\MT^{\wild}_Y(X,w,\nbiga)=\MT^{\wild}(Y,w,\nbiga)$
for a complex submanifold $Y\subset X$.
\end{lem}
\pf
According to \cite{sabbah2},
$\MT_Y(X,w)=\MT(Y,w)$.
It is easy to check that 
the condition ($MT_{>0}(\ramification,\exp,\nbiga)$)
is preserved.
\hfill\qed

\vspace{.1in}
The following lemma is clear from the definition.
\begin{lem}
\label{lem;08.12.4.2}
Let $\nbigt=(\nbigm',\nbigm'',C)$ be
an $\nbigr_X$-triple.
\begin{itemize}
\item
 If $\nbigt$ is a wild pure twistor $D$-module on $X$,
 the restriction $\nbigt_{|U}$ is also 
 a wild pure twistor $D$-module on $U$
 for any open subset $U$ of $X$.
\item
 Assume 
 (i) there exists an open covering
 $X=\bigcup U_i$ such that
 $\nbigt_{|U_i}$ are wild pure twistor $D$-modules
 on $U_i$,
 (ii) $\nbigm'$ and $\nbigm''$ are holonomic
 (see Subsection {\rm\ref{subsection;08.12.4.1}}).
Then, $\nbigt$ is also a wild pure twistor $D$-module
on $X$.
\hfill\qed
\end{itemize}
\end{lem}

\begin{lem}
\label{lem;07.12.6.4}
Let $\nbigt\in \MT^{\wild}(X,w,\nbiga)$
and $\nbigt'\in \MT^{\wild}(X,w',\nbiga)$
for $w>w'$.
Any morphism $\nbigt\lrarr\nbigt'$
of $\nbigr$-triples is trivial.
\end{lem}
\pf
It can be reduced to the claim for 
pure twistor $D$-modules
in the sense of \cite{sabbah2}.
(See Proposition 4.1.8 of \cite{sabbah2}.)
\hfill\qed

\vspace{.1in}

Let $\MTW_{\leq d}^{\wild}(X,\nbiga)$ be the category of
$\nbigr$-triples $\nbigt$ with increasing filtrations $W$
indexed by $\seisuu$
such that
$\Gr^W_{\ell}\nbigt\in \MT_{\leq d}^{\wild}(X,\ell,\nbiga)$.
Let $\MTN^{\wild}(X,w,\nbiga)$ be the category of
$\nbigr$-triples $\nbigt$ with morphisms
$\nbign:\nbigt\lrarr\nbigt\otimes\Tate(-1)$
such that
(i) $\nbign^m=0$ for any sufficiently large $m$,
(ii) $\Gr^{W(\nbign)}_{\ell}(\nbigt)\in
  \MT_{\leq d}^{\wild}(X,w+\ell,\nbiga)$
for the monodromy weight filtration $W(\nbign)$.
A morphism 
$\varphi:(\nbigt_1,\nbign_1)\lrarr(\nbigt_2,\nbign_2)$
is defined to be a morphism of $\nbigr$-triples
such that $\varphi\circ\nbign_1=\nbign_2\circ\varphi$.
\index{$\MTW_{\leq d}^{\wild}(X,\nbiga)$}
\index{category $\MTN^{\wild}(X,w,\nbiga)$}

\begin{lem}
 \mbox{{}}\label{lem;07.12.1.7}
\begin{description}
\item[($a_d$)]
The category
$MT_{\leq d}^{\wild}(X,w,\nbiga)$ is abelian.
All morphisms are strict and strictly specializable.
\item[($b_d$)]
The category $\MTW_{\leq d}^{\wild}(X,\nbiga)$ is abelian.
Any morphism $\varphi:(\nbigt_1,W_1)\lrarr(\nbigt_2,W_2)$
is strict with respect to the filtrations $W_i$.
All morphisms are strict.
\item[($c_d$)]
The category $\MTN_{\leq d}^{\wild}(X,w,\nbiga)$ is abelian.
Any morphism
$\varphi:(\nbigt_1,\nbign_1)\lrarr(\nbigt_2,\nbign_2)$
is strict with respect to $W(\nbign_1)$ and $W(\nbign_2)$.
The filtrations on $\Image(\varphi)$,
$\Ker(\varphi)$ and $\Cok(\varphi)$
induced by $W(\nbign_i)$ $(i=1,2)$
are equal to the monodromy weight filtrations
of the induced morphisms.
All morphisms are strict.
\end{description}
\end{lem}
\pf
The proof of the lemma can be carried out
with the argument in \cite{sabbah2}
(based on \cite{saito1} for pure Hodge modules),
and many claims can be reduced
to the case of pure twistor $D$-modules.
We give only an outline,
and we refer to \cite{sabbah2} and \cite{sabbah5}
for more details.
Let us see $(a_d)\Longrightarrow (c_d)$.
Let $\varphi:(\nbigt_1,\nbign_1)\lrarr (\nbigt_2,\nbign_2)$.
We have the induced morphisms
$\Gr^{W(\nbign)}_{\ell}(\varphi):
 \Gr^{W(\nbign)}_{\ell}(\nbigt_1)
\lrarr
 \Gr^{W(\nbign)}_{\ell}(\nbigt_2)$
of pure twistor $D$-modules.
By using Lemma \ref{lem;07.12.6.4},
it is easy to derive that $\varphi$ is strict
with respect to $W(\nbign_1)$ and $W(\nbign_2)$.
Therefore, $\Ker\varphi$, $\Image(\varphi)$ and $\Cok(\varphi)$
are equipped with the naturally induced filtration $W$,
and we have the natural isomorphisms
$\Gr_{\ell}^W(\Ker(\varphi))\simeq 
 \Ker(\Gr_{\ell}^W(\varphi))$,
$\Gr^W_{\ell}(\Image(\varphi))\simeq
 \Image(\Gr_{\ell}^W(\varphi))$
and 
$\Gr^W_{\ell}(\Cok(\varphi))\simeq
 \Cok(\Gr_{\ell}^W(\varphi))$.
Since we have the isomorphisms
$\nbign_1^{\ell}:
 \Gr^W_{\ell}(\Ker(\varphi))
\simeq
 \Gr^{W}_{-\ell}(\Ker(\varphi))$,
the induced filtration $W$ is the same as 
the weight filtration of the induced nilpotent map $\nbign$
on $\Ker(\varphi)$.
Similar claims hold for $\Image(\varphi)$
and $\Cok(\Image\varphi)$.
We can show $(a_d)\Longrightarrow(b_d)$
in a similar way.

Let us see $(c_{d-1})\Longrightarrow (a_d)$.
Let $\varphi:\nbigt_1\lrarr\nbigt_2$ be 
a morphism of wild pure twistor $D$-modules
with weight $w$.
We may assume that $\nbigt_i$ have 
the same strict support $Z$.
Since, it is a morphism of pure twistor $D$-modules,
we have already known that
(i) $\varphi$ is strict and strictly specializable
along any function $g$, 
(ii) $\Ker(\varphi)$, $\Image(\varphi)$
 and $\Cok(\varphi)$ are pure twistor $D$-modules
 whose strict supports are $Z$.
Let $g$ be any function. 
We have the induced morphism
$\psitilde_{g,\gminia,u}(\varphi):
 (\psitilde_{g,\gminia,u}(\nbigt_1),\nbign_{1})
\lrarr (\psitilde_{g,\gminia,u}(\nbigt_2),\nbign_2)$
to which we can apply $(c_{d-1})$.
We obtain that
$\Gr^{W(N)}_{\ell}\Ker(\psitilde_{g,\gminia,u}(\varphi))$,
$\Gr^{W(N)}_{\ell}\Image(\psitilde_{g,\gminia,u}(\varphi))$
and
$\Gr^{W(N)}_{\ell}\Cok(\psitilde_{g,\gminia,u}(\varphi))$
are wild pure twistor $D$-modules with weight $w+\ell$.
In particular,
we also obtain that 
$\Cok\psitilde_{g,\gminia,u}(\varphi)$
is strict.
Due to Lemma \ref{lem;08.1.16.10},
we obtain the natural isomorphisms
$\Ker(\psitilde_{g,\gminia,u}(\varphi))\simeq
 \psitilde_{g,\gminia,u}(\Ker\varphi)$,
$\Image(\psitilde_{g,\gminia,u}(\varphi))\simeq
 \psitilde_{g,\gminia,u}(\Image\varphi)$
and 
$\Cok(\psitilde_{g,\gminia,u}(\varphi))\simeq
 \psitilde_{g,\gminia,u}(\Cok\varphi)$.
Hence, we obtain that
$\Ker(\varphi)$, $\Image(\varphi)$
and $\Cok(\varphi)$ are also wild pure twistor $D$-modules
with weight $w$.
\hfill\qed

\vspace{.1in}
In the proof,
we have obtained the following.
\begin{cor}
\label{cor;08.1.16.12}
Let $\nbigt_i$ $(i=1,2)$ be
wild pure twistor $D$-modules of weight $w$.
Let $\varphi:\nbigt_1\lrarr\nbigt_2$
be a morphism of wild pure twistor $D$-modules.
For any holomorphic function $g$,
$\gminia\in\cnum[t_n^{-1}]$
and $u\in\real\times\cnum$,
we have the natural isomorphisms
$\psitilde_{g,\gminia,u}(\Ker\varphi)
\simeq
 \Ker\psitilde_{g,\gminia,u}(\varphi)$,
$ \psitilde_{g,\gminia,u}(\Image\varphi)
\simeq
 \Image\psitilde_{g,\gminia,u}(\varphi)$,
and
$\psitilde_{g,\gminia,u}(\Cok\varphi)
\simeq
 \Cok\psitilde_{g,\gminia,u}(\varphi)$.
We also have 
\[
P\Gr^{W(N)}_{\ell}
 \psitilde_{g,\gminia,u}(\Ker\varphi)
\simeq
 \Ker
 P\Gr^{W(N)}_{\ell}
 \psitilde_{g,\gminia,u}(\varphi),
\]
\[
  P\Gr^{W(N)}_{\ell}
 \psitilde_{g,\gminia,u}(\Image\varphi)
\simeq
  \Image
 P\Gr^{W(N)}_{\ell}
 \psitilde_{g,\gminia,u}(\varphi),
\]
\[
P\Gr^{W(N)}_{\ell}
 \psitilde_{g,\gminia,u}(\Cok\varphi)
\simeq
 \Cok
 P\Gr^{W(N)}_{\ell}
  \psitilde_{g,\gminia,u}(\varphi).
\]

In particular,
let $\nbigc$ be a complex
in $\MT^{\wild}(X,w,\nbiga)$.
By taking the functors $\psitilde_{g,\gminia,u}$
and 
$P\Gr^{W(N)}_{\ell}\psitilde_{g,\gminia,u}$,
we have the complex
$\psitilde_{g,\gminia,u}(\nbigc)$
in $\MTN^{\wild}(X,w,\nbiga)$,
and the complex
$P\Gr^{W(N)}_{\ell}\psitilde_{g,\gminia,u}(\nbigc)$
in $\MT^{\wild}(X,w,\nbiga)$.
Then, we have the natural isomorphisms:
\[
 H^i\bigl(\psitilde_{g,\gminia,u}(\nbigc)\bigr)
\simeq 
 \psitilde_{g,\gminia,u}H^i(\nbigc),
\quad
H^i\bigl(
 P\Gr^{W(N)}_{\ell}
 \psitilde_{g,\gminia,u}(\nbigc)\bigr)
\simeq 
 P\Gr^{W(N)}_{\ell}
 \psitilde_{g,\gminia,u}H^i(\nbigc).
\]
\hfill\qed
\end{cor}

\subsection{Polarizable wild pure twistor $D$-module}

We recall the definition of polarization 
of wild pure twistor $D$-modules.
(\cite{sabbah2} and \cite{sabbah5}.
 See also Saito's original work \cite{saito1}.)
\begin{df}
Let $\nbigt\in \MT^{\wild}(X,w,\nbiga)$.
Let $\nbigs:\nbigt\lrarr\nbigt^{\ast}(-w)$
be a Hermitian sesqui-linear duality
with weight $w$.
It is called a polarization,
if the following inductive conditions are satisfied:
\index{polarization}
\begin{description}
\item[($MTP_{>0}(\ramification,\exp)$)]
 Let $U$ be an open subset of $X$,
 and let $f$ be a holomorphic function on $U$.
 For each
 $\gminia\in\cnum[t_n^{-1}]$,
 $u\in\real\times\nbiga$
 and $\ell\in\seisuu_{\geq\,0}$,
 the induced Hermitian sesqui-linear duality
\[
 \nbigs_{f,\gminia,u,{\ell}}:
 P\Gr^{W}_{\ell}\psitilde_{f,\gminia,u}(\nbigt)
 \lrarr P\Gr^{W}_{\ell}
 \psitilde_{f,\gminia,u}(\nbigt)^{\ast}(-w-l) 
\]
 is a polarization of 
 $P\Gr^{W}_{\ell}\psitilde_{f,\gminia,u}(\nbigt)$
 as an $\nbiga$-wild pure twistor $D$-module with
 weight $w+\ell$.
\item[($MTP_0$)]
 We have the decomposition $\nbigs=\bigoplus\nbigs_Z$,
 corresponding to $\nbigt=\bigoplus\nbigt_Z$.
 If $\dim Z=0$,
 $(\nbigt_Z,\nbigs_Z)$ is a push-forward
 of a polarized pure twistor structure with
 weight $w$.
\end{description}
An $\nbiga$-wild pure twistor $D$-module $\nbigt$ is
called polarizable,
if there exists a polarization $\nbigs$ of $\nbigt$.
\hfill\qed
\end{df}
\index{polarized wild pure twistor $D$-module}
\index{polarizable wild pure twistor $D$-module}

Let $\MT^{\wild}(X,w,\nbiga)^{(p)}$ denote
the full subcategory of $\MT^{\wild}(X,w,\nbiga)$
which consists of the polarizable $\nbiga$-wild
pure twistor $D$-modules with
weight $w$.
\index{category $\MT^{\wild}(X,w,\nbiga)^{(p)}$}
When we are given a subvariety $Y\subset X$,
let $\MT^{\wild}_Y(X,w,\nbiga)^{(p)}$
(resp. $\MT^{\wild}_{Y,ss}(X,w,\nbiga)^{(p)}$)
denote the full subcategory of 
$\nbiga$-wild pure twistor
$D$-modules on $X$ with weight $w$
whose supports are contained in $Y$
(resp. whose strict supports are exactly $Y$).
\index{category $\MT^{\wild}_Y(X,w,\nbiga)^{(p)}$}
\index{category $\MT^{\wild}_{Y,ss}(X,w,\nbiga)^{(p)}$}
We say that $\nbigt\in \MT^{\wild}(X,w,\nbiga)$
is simple, if it does not contain any non-trivial subobjects.
We say that $(\nbigt,\nbigs)$ with
weight $w$
is simple, if $\nbigt$ is simple.
\index{simple wild pure twistor $D$-module}

\begin{lem}
\label{lem;10.7.2.1}
 Let $(\nbigt,\nbigs)$ be
 a polarized wild pure twistor $D$-module
 with weight $w$, whose strict support is $Z$.
 We have a Zariski closed subset $Z_0\subset Z$
 such that
 $(\nbigt,\nbigs)_{|X-Z_0}$ is a push-forward of
 a variation of polarized pure twistor structure on $Z-Z_0$.
\end{lem}
\pf
It can be reduced to the claim for 
pure twistor $D$-modules
in the sense of \cite{sabbah2}.
Then, it follows from
Proposition 4.1.9 of \cite{sabbah2}.
\hfill\qed

\begin{prop}
Let $(\nbigt,\nbigs)$ be 
a polarized $\nbiga$-wild pure twistor $D$-module
of weight $w$.
Let $\nbigt'\subset\nbigt$
be wild pure twistor $D$-submodule
of weight $w$.
Then, the composition $\nbigs'$
of the following morphisms
\[
 \nbigt'\lrarr\nbigt\stackrel{\nbigs}\lrarr
 \nbigt^{\ast}(-w)\lrarr \nbigt^{\prime\ast}(-w)
\]
gives a polarization of $\nbigt'$,
and we have the decomposition
$(\nbigt,\nbigs)=(\nbigt',\nbigs')
\oplus(\nbigt'',\nbigs'')$.
\end{prop}
\pf
We have the decomposition of
polarized pure twistor $D$-modules
as above,
which is shown in Proposition 4.2.5 
of \cite{sabbah2}.
It gives the decomposition of
polarized $\nbiga$-wild pure twistor $D$-modules.
\hfill\qed

\begin{cor}
\label{cor;07.12.6.1}
Let $(\nbigt,\nbigs)$ be 
a polarized $\nbiga$-wild pure twistor $D$-module
of weight $w$.
Then, we have a decomposition
$(\nbigt,\nbigs)=\bigoplus (\nbigt_i,\nbigs_i)$,
where $\nbigt_i$ are simple $\nbiga$-wild pure twistor
$D$-modules of weight $w$.
In particular,
the abelian category
$\MT^{\wild}(X,w,\nbiga)^{(p)}$ is semisimple.
\hfill\qed
\end{cor}

\begin{prop}
\label{prop;07.12.6.2}
Let $(\nbigt,\nbigs)$ be a simple polarized
$\nbiga$-wild pure twistor $D$-module with weight $w$.
\begin{itemize}
\item
 Let $(V,S)$ be a polarized pure twistor structure 
 with weight $0$ in dimension $0$.
 Then, the naturally defined Hermitian sesqui-linear duality
 $\nbigt\otimes V\lrarr
 \nbigt^{\ast}\otimes V^{\ast}(-w)$ is a polarization.
\item
 Let $V$ be a pure twistor structure with weight $0$
 in dimension $0$.
 For any polarization $\nbigstilde$ of $\nbigt\otimes V$,
 there exists a polarization $S$ of $V$
 such that $\nbigstilde=\nbigs\otimes S$. 
\end{itemize}
\end{prop}
\pf
The first claim is clear.
Because of Corollary \ref{cor;07.12.6.1},
we have the decomposition 
$V=\bigoplus V_i$ into rank one objects
such that 
$(\nbigt,\nbigstilde)=
 \bigoplus (\nbigt\otimes V_i,\nbigstilde_i)$.
Because $\nbigt$ is simple,
there exists a positive number $a_i>0$
such that $\nbigstilde_i=a_i\cdot \nbigs$,
under some identification $V_i\simeq\Tate^S(0)$.
Then, the second claim follows.
\hfill\qed

\begin{prop}
\label{prop;07.12.6.3}
Let $\nbigt$ be a polarizable $\nbiga$-wild
pure twistor $D$-modules
with weight $w$.
\begin{itemize}
\item
We have the canonical decomposition 
$\nbigt=\bigoplus \nbigt_i$,
where $\nbigt_i$ are the direct sum of 
some simple objects $\nbigt_i'$
such that $\nbigt'_i\not\simeq\nbigt'_j$ $(i\neq j)$.
\item
Any polarization $\nbigs$ of $\nbigt$
is a direct sum of the polarizations $\nbigs_i$
of $\nbigt_i$.
\end{itemize}
\end{prop}
\pf
The first claim follows from 
Corollary \ref{cor;07.12.6.1}.
The second claim is obvious.
\hfill\qed

\begin{cor}
\label{cor;10.6.18.1}
Let $\nbigt$ be an $\nbiga$-wild pure twistor $D$-module
with weight $w$.
Let $\nbigs_i$ $(i=1,2)$ be polarizations of $\nbigt$.
Then, there exists an automorphism 
$\varphi$ of $\nbigt$
such that
$\varphi:(\nbigt,\nbigs_1)\lrarr(\nbigt,\nbigs_2)$
is an isomorphism.
In other words,
a polarization of $\nbigt$ is unique up to
obvious ambiguity.
\end{cor}
\pf
It follows from Proposition \ref{prop;07.12.6.2}
and Proposition \ref{prop;07.12.6.3}.
\hfill\qed

\subsection{Polarized graded 
wild Lefschetz twistor $D$-module}

We recall the definition
and some results for polarized graded
wild Lefschetz twistor $D$-module
due to Sabbah (\cite{sabbah2} and \cite{sabbah5}).
We consider the graded objects with the lower index.
It is easy to give a modification for the upper index case.

Let $\epsilon$ be an integer.
A graded $\nbiga$-wild Lefschetz twistor $D$-module
with weight $w$ of type $\epsilon$
is a pair $(\nbigt,\nbigl)$ as follows:
\index{graded $\nbiga$-wild 
Lefschetz twistor $D$-module}
\begin{itemize}
\item
$\nbigt=\bigoplus_{j\in\seisuu}\nbigt_j$ 
is a graded $\nbigr$-triple.
Each $\nbigt_j$ is an $\nbiga$-wild pure twistor $D$-module
of weight $w-\epsilon\cdot j$.
\item
$\nbigl$ is a graded morphism 
$\nbigt\lrarr\nbigt(\epsilon)[2]$
such that
$\nbigl^j:\nbigt_{j}\lrarr \nbigt_{-j}(\epsilon\cdot j)$
are isomorphisms for any $j\geq 0$.
\end{itemize}

We put $P\nbigt_{j}:=\Ker\nbigl^{j+1}\cap \nbigt_j$
for each $j\geq 0$,
and $P\nbigt_j=0$ for $j<0$.
Then, we obtain 
the naturally defined primitive decomposition
$\nbigt_j=\bigoplus
 \nbigl^k \bigl(P\nbigt_{j+2k}(-\epsilon k)\bigr)$
for any $j$,
where $\nbigl^k \bigl(P\nbigt_{j+2k}(-\epsilon k)\bigr)$
denote the image of
$\nbigl^k:P\nbigt_{j+2k}(-\epsilon k)
\lrarr \nbigt_j$.

Recall that
the Hermitian adjoint $\nbigt^{\ast}$
has the natural grading given by
$(\nbigt^{\ast})_j:=
 (\nbigt_{-j})^{\ast}$.
It is naturally equipped with an induced morphism
$\nbigl^{\ast}:\nbigt^{\ast}\lrarr 
 \nbigt^{\ast}(\epsilon)[2]$,
and $(\nbigt^{\ast},\nbigl^{\ast})$
is a graded $\nbiga$-wild Lefschetz twistor $D$-module
with weight $w$ of type $\epsilon$.
We have the natural isomorphism
$(P\nbigt_j)^{\ast}
 \simeq
 (\nbigl^{\ast})^j
 \bigl(P\nbigt_j^{\ast}(-\epsilon j)\bigr)$.

\vspace{.1in}

Let $\nbigs:\nbigt\lrarr\nbigt^{\ast}(-w)$
be a graded Hermitian sesqui-linear duality satisfying
$\nbigl^{\ast}\circ\nbigs+\nbigs\circ\nbigl=0$.
We have the composite:
\[
 \nbigs_{-j}\circ\nbigl^j:
 P\nbigt_j\lrarr
 \nbigl^j
\bigl(
 P\nbigt_{j}(-\epsilon j)
\bigr)(\epsilon j)
 \lrarr 
 (P\nbigt_j)^{\ast}(-w+\epsilon j)
\]
Then, $\nbigs$ is called a polarization 
of $(\nbigt,\nbigl)$,
if the composite $\nbigs_{-j}\circ\nbigl^j$
gives a polarization of $P\nbigt_j$
for each $j\geq 0$.
\index{polarization}

The following lemma,
due to Sabbah,
is one of the most important properties
of graded $\nbiga$-wild Lefschetz twistor $D$-module.
Saito originally proved the corresponding property
for his Hodge modules.

\begin{lem}
 \label{lem;05.2.14.20}
Assume $\epsilon$ is $1$ or $-1$.
Let $(\nbigt,\nbigl,\nbigs)$ (resp. $(\nbigt',\nbigl',\nbigs')$)
be a polarized graded 
$\nbiga$-wild Lefschetz twistor $D$-module of type $\epsilon$
with weight $w$ (resp. $w-\epsilon$).
Let $c:\nbigt\lrarr\nbigt'[1]$ and
 $v:\nbigt'\lrarr\nbigt[1](\epsilon)$
 be graded morphisms
such that $v\circ c=\nbigl$ and $c\circ v=\nbigl'$.
Assume that $c$ and $v$ are adjoint with respect to $\nbigs$
and $\nbigs'$,
i.e.,
$\nbigs'\circ c=v^{\ast}\circ\nbigs$.
Then we have the decomposition
$\nbigt'=\Image c \oplus\Ker v$.
\end{lem}
\pf
It can be shown by using the argument
in the proof of Proposition 4.2.10 of \cite{sabbah2}.
We give only an outline.
We may assume that
both $\nbigt$ and $\nbigt'$ have
the same strict support $Z$.
By the result for 
smooth polarized graded Lefschetz twistor structures
(Lemma 5.2.15 of \cite{saito1} or
 Proposition 2.19 of \cite{sabbah2}),
we have such decomposition on
the generic part of $Z$.
Because of the strict $S$-decomposability,
we can easily obtain the decomposition
on $Z$.
\hfill\qed

\subsection{Bi-graded wild Lefschetz twistor $D$-module}

We recall the notion of
bi-graded wild Lefschetz twistor $D$-module
due to Sabbah (\cite{sabbah2} and \cite{sabbah5}).
We consider the bi-graded objects with lower indices.
It is easy to give a modification in the other cases.
\index{bi-graded wild Lefschetz twistor $D$-module}

Let $\epsilon_a$ $(a=1,2)$ be integers.
A bi-graded $\nbiga$-wild 
Lefschetz twistor $D$-modules of type
$(\epsilon_1,\epsilon_2)$ with weight $w$
is defined naturally,
i.e.,
it is a tuple $(\nbigt,\nbigl_1,\nbigl_2)$ as follows:
\begin{itemize}
\item
$\nbigt$ is a bi-graded $\nbigr$-triple
$\bigoplus_{i,j\in\seisuu}\nbigt_{i,j}$.
Each $\nbigt_{i,j}$ is a 
$\nbiga$-wild pure twistor $D$-module
of weight $w-i\epsilon_1-j\epsilon_2$.
\item
$\nbigl_1$ is a tuple of morphisms
$\nbigt_{i,j}\lrarr\nbigt_{i-2,j}(\epsilon_1)$
($i,j\in\seisuu$)
such that
\[
 \nbigl_1^{i}:\nbigt_{i,j}\lrarr
 \nbigt_{-i,j}(i\epsilon_1)
\]
is an isomorphism for each $i\geq 0$.
\item
$\nbigl_2$ is a tuple of morphisms
$\nbigt_{i,j}\lrarr\nbigt_{i,j-2}(\epsilon_2)$
($i,j\in\seisuu$)
such that
\[
\nbigl_2^{j}:\nbigt_{i,j}\lrarr
 \nbigt_{i,-j}(j \epsilon_2)
\]
is an isomorphism for each $j\geq 0$.
\item
$\nbigl_1$ and $\nbigl_2$ commute.
\end{itemize}
We put
$P\nbigt_{i,j}=
 \Ker\nbigl_1^{i+1}\cap
 \Ker\nbigl_2^{j+1}\cap \nbigt_{i,j}$.
They obviously induce the primitive decomposition
as in the case of graded $\nbiga$-wild
Lefschetz twistor $D$-modules.
The Hermitian adjoint $\nbigt^{\ast}$
naturally has the grading and the morphisms
$\nbigl_i^{\ast}$ $(i=1,2)$,
for which $(\nbigt^{\ast},\nbigl^{\ast})$
is a bi-graded $\nbiga$-wild
Lefschetz twistor $D$-module.

A bi-graded morphism 
$\nbigs:\nbigt\lrarr\nbigt^{\ast}(-w)$
such that 
$\nbigl^{\ast}_a\circ\nbigs
+\nbigs\circ\nbigl_a=0$ $(a=1,2)$
induces the morphism:
\[
 \nbigs\circ \nbigl_1^i\circ\nbigl_2^j:
 P\nbigt_{i,j}\lrarr
 (P\nbigt_{i,j})^{\ast}(-w+i \epsilon_1+j \epsilon_2).
\]
If $\nbigs\circ\nbigl_1^i\circ\nbigl_2^j$ gives
a polarization of $P\nbigt_{i,j}$,
$\nbigs$ is called a polarization of
$(\nbigt,\nbigl_1,\nbigl_2)$.

\vspace{.1in}
Let $d:\nbigt\lrarr\nbigt(\epsilon_1+\epsilon_2)$ be
a morphism satisfying the following:
\begin{itemize}
\item
 $d:\nbigt_{i,j}\lrarr \nbigt_{i-2,j-2}(\epsilon_1+\epsilon_2)$
 and $d\circ d=0$.
\item
 $d$ is anti-commutative
 with $\nbigl_a$  $(a=1,2)$,
 i.e.,
 $d\circ\nbigl_a+\nbigl_a\circ d=0$.
\item
 $d$ is self-adjoint with respect to $\nbigs$,
 i.e.,
 $d^{\ast}\circ\nbigs=\nbigs\circ d$.
\end{itemize}
We have the naturally defined bi-graded structure
on $\nbigt^{(1)}:=\Ker d/\Image d$.
We also have the naturally defined maps 
$\nbigl^{(1)}_a$ $(a=1,2)$
and the Hermitian sesqui-linear
duality $\nbigs^{(1)}$.

The following lemma,
due to Sabbah,
is one of the most important properties
of bi-graded $\nbiga$-wild Lefschetz twistor $D$-module.
Saito originally proved the corresponding property for
his Hodge modules.
\begin{lem}
\label{lem;05.1.18.25}
$(\nbigt^{(1)},\nbigl^{(1)}_1,\nbigl^{(1)}_2,\nbigs^{(1)})$
is a polarized bi-graded $\nbiga$-wild 
Lefschetz $D$-module with weight $w$
of type $(\epsilon_1,\epsilon_2)$.
\end{lem}
\pf
It can be shown using the same argument
as that in the proof of Proposition 4.2.10 of \cite{sabbah2}.
We give only an outline.
We may assume that $\nbigt$ has the strict support $Z$
which is irreducible.
We use an induction on the dimension of $Z$.
By the result for 
smooth polarized bi-graded Lefschetz twistor structures
(Proposition 4.22 of \cite{saito1},
 Theorem 4.5 of \cite{g-n},
 Lemma 2.1.20 of \cite{sabbah2}),
the claim holds on the generic part of $Z$.
By the strict $S$-decomposability,
it follows that
the morphisms
$(\nbigl^{(1)}_1)^{j_1}:
 \nbigt^{(1)}_{j_1,j_2}\simeq
 \nbigt^{(1)}_{-j_1,j_2}(j_1\cdot\epsilon_1)$ 
are isomorphisms for $j_1\geq 0$,
and that 
the morphisms
$(\nbigl^{(1)}_2)^{j_2}:
 \nbigt^{(1)}_{j_1,j_2}\simeq
 \nbigt^{(1)}_{j_1,-j_2}(j_2\cdot\epsilon_2)$ 
are isomorphisms for $j_2\geq 0$.
We have only to show that
the composite
$\nbigs^{(1)}\circ
 (\nbigl_1^{(1)})^{j_1}\circ(\nbigl_2^{(1)})^{j_2}$
gives a polarization of
$P\nbigt^{(1)}_{j_1,j_2}$ for each 
$(j_1,j_2)\in\seisuu_{\geq 0}^2$.

Let $g$ be any holomorphic function.
For any $\ell\in\seisuu_{\geq 0}$,
$u\in\real\times\cnum$
and $\gminia\in\cnum[t_n^{-1}]$,
we consider the following:
\begin{multline*}
 (\nbigt_{g,\gminia,u,\ell})_{j_1,j_2}:=
 P\Gr^{W}_{\ell}\psitilde_{g,\gminia,u}(\nbigt_{j_1,j_2})
 \\
=\Ker\Bigl(
 \nbign^{\ell+1}:
 \Gr^W_{\ell}\psitilde_{g,\gminia,u}(\nbigt_{j_1,j_2})
\lrarr
 \Gr^{W}_{-\ell-2}\psitilde_{g,\gminia,u}(\nbigt_{j_1,j_2})
 (-\ell-1)
 \Bigr)
\end{multline*}
Then, 
$\nbigt_{g,\gminia,u,\ell}:=
 \bigoplus_{j_1,j_2}
 \bigl(\nbigt_{g,\gminia,u,\ell}\bigr)_{j_1,j_2}$
with the induced morphisms $\nbigl_a$ $(a=1,2)$
is a bi-graded $\nbiga$-wild 
Lefschetz twistor $D$-module
of weight $w+\ell$.
We have the following morphisms
induced by $\nbign^{\ell}$ and $\nbigs$:
\[
  \Gr^W_{\ell}\psitilde_{g,\gminia,u}(\nbigt_{j_1,j_2})
 \lrarr \Gr^{W}_{-\ell}
 \psitilde_{g,\gminia,u}(\nbigt_{j_1,j_2})(-\ell)
\lrarr
 \Gr^{W}_{\ell}\psitilde_{g,\gminia,u}
 (\nbigt_{-j_1,-j_2})^{\ast}(-w-\ell) 
\]
They induce a morphism
$(\nbigt_{g,\gminia,u,\ell})_{j_1,j_2}
\lrarr
 \bigl((\nbigt_{g,\gminia,u,\ell})_{-j_1,-j_2}
 \bigr)^{\ast}(-w-\ell)$,
which gives a  Hermitian sesqui-linear duality
$\nbigs_{g,\gminia,u,\ell}$ of 
the bi-graded $\nbiga$-wild 
Lefschetz twistor $D$-module
$\nbigt_{g,\gminia,u,\ell}$ with weight $w+\ell$.
\begin{lem}
\label{lem;07.10.27.1}
$\nbigs_{g,\gminia,u,\ell}$ is a polarization
of $\nbigt_{g,\gminia,u,\ell}$.
\end{lem}
\pf
We have the natural isomorphism
$P\bigl(\nbigt_{g,\gminia,u,\ell}\bigr)_{j_1,j_2}
\simeq
 P\Gr^{W}_{\ell}\psitilde_{g,\gminia,u}(P\nbigt_{j_1,j_2})$,
via which $\nbigs_{g,\gminia,u,\ell}$
is the same as the polarization of
 $P\Gr^{W}_{\ell}\psitilde_{g,\gminia,u}(P\nbigt_{j_1,j_2})$.
Thus, we obtain Lemma \ref{lem;07.10.27.1}.
\hfill\qed

\vspace{.1in}
We have the induced morphisms
$d:(\nbigt_{g,\gminia,u,\ell})_{j_1,j_2}
 \lrarr
 (\nbigt_{g,\gminia,u,\ell})_{j_1-1,j_2-1}$.
By the hypothesis of the induction,
the induced tuple
$\bigl(\nbigt^{(1)}_{g,\gminia,u,\ell},
 \nbigl_1^{(1)},\nbigl_2^{(1)},
 \nbigs_{g,\gminia,u,\ell}^{(1)}\bigr)$
is a polarized 
bi-graded $\nbiga$-wild Lefschetz twistor $D$-module
of weight $w+\ell$.
In particular,
$\nbigs_{g,\gminia,u,\ell}^{(1)}\circ
 (\nbigl_1^{(1)})^{j_1}\circ
 (\nbigl_2^{(1)})^{j_2}$
gives a polarization of
$P\bigl(\nbigt^{(1)}_{g,\gminia,u,\ell}\bigr)_{j_1,j_2}$.

By Corollary \ref{cor;08.1.16.12},
we have a natural isomorphism
\[
  (\nbigt_{g,\gminia,u,\ell}^{(1)})_{j_1,j_2}
\simeq
 P\Gr^{W}_{\ell}\psitilde_{g,\gminia,u}
 \bigl(\nbigt^{(1)}_{j_1,j_2}\bigr),
\]
compatible with the induced Lefschetz morphisms.
We obtain the isomorphism
$P\bigl(\nbigt^{(1)}_{g,\gminia,u,\ell}\bigr)_{j_1,j_2}
\simeq
 P\Gr^{W}_{\ell}\psitilde_{g,\gminia,u}
 \bigl(P\nbigt^{(1)}_{j_1,j_2}\bigr)$,
and the induced Hermitian sesqui-linear duality
of $P\Gr^{W}_{\ell}\psitilde_{g,\gminia,u}
 \bigl(P\nbigt^{(1)}_{j_1,j_2}\bigr)$ is equal to
the polarization
$\nbigs_{g,\gminia,u,\ell}^{(1)}\circ
 (\nbigl_1^{(1)})^{j_1}\circ
 (\nbigl_2^{(1)})^{j_2}$.
Therefore we can conclude that
$\nbigs^{(1)}\circ(\nbigl^{(1)}_1)^{j_1}
 \circ (\nbigl^{(1)}_2)^{j_2}$
is a polarization of $P\nbigt^{(1)}_{j_1,j_2}$,
and we obtain Lemma \ref{lem;05.1.18.25}.
\hfill\qed

\section[Wild pure twistor D-modules on curves]
{Wild harmonic bundles and wild pure twistor
$D$-modules on curves}
\label{subsection;08.10.18.120}
\subsection{Statement}

Let $X$ be a smooth complex curve.
Let $\nbiga$ be a vector subspace of $\cnum$
over $\rnum$.
Let $D$ be a discrete subset of $X$.
Let $(V,\DD^{\sankaku})$
be a variation of pure twistor structure of weight $w$
on $X-D$ with a pairing.
We say that 
$(V,\DD^{\sankaku},S)$ is 
$\nbiga$-wild on $(X,D)$,
if the corresponding harmonic bundle is wild
on $(X,D)$.
\index{$\nbiga$-wild
variation of pure twistor structure}

Let $\VPTgenwild(X,D,w,\nbiga)$ denote the category of
$\nbiga$-wild variation of polarized 
pure twistor structure of weight $w$ on $(X,D)$.
Let $\MPT^{\wild}_{\strict}(X,D,w,\nbiga)$ denote
the category of polarized $\nbiga$-wild pure 
twistor $D$-modules of weight $w$,
such that 
(i) their strict supports are $X$,
(ii) their restriction to $X-D$ comes from 
a polarized variation of pure twistor structure.
In this subsection, 
for both categories,
morphisms are defined to be isomorphisms.
\index{category $\VPTgenwild(X,D,w,\nbiga)$}
\index{category $\MPT^{\wild}_{\strict}(X,D,w,\nbiga)$}
For a given
$(\nbigt,\nbigs)\in\MPT^{\wild}_{\strict}(X,D,w,\nbiga)$,
its restriction to $X-D$ comes from
a harmonic bundle with the Tate twist.
It is easy to see that the harmonic bundle is wild
on $(X,D)$.
Hence, we have a naturally defined functor
\[
 \Phi:\MPT^{\wild}_{\strict}(X,D,w,\nbiga)
\lrarr \VPTgenwild(X,D,w,\nbiga).
\]

\begin{prop}
\label{prop;07.10.19.50}
$\Phi$ is an equivalence.
\end{prop}
We have only to consider the case $w=0$.
In the following argument,
we omit to distinguish $\nbiga$.

\begin{rem}
We will show the higher dimensional version
(Theorem {\rm\ref{thm;07.10.28.30}})
in Section {\rm\ref{subsection;10.6.6.13}}.
Because we need Proposition 
{\rm\ref{prop;07.10.19.50}}
in the proof of Hard Lefschetz Theorem
(Theorem {\rm\ref{thm;07.10.23.20}}),
we include it here.
See also Corollary {\rm\ref{cor;10.6.6.20}}
for a variant of the statement.
\hfill\qed
\end{rem}

\subsection{From wild harmonic bundle to
 wild pure twistor $D$-module}
\label{subsection;07.10.28.80}

Let $D$ be a discrete subset of $X$.
Let $\harmonicbundle$ be a wild harmonic bundle
defined on $X-D$.
By applying the construction
in Section \ref{subsection;07.10.19.30},
we obtain the $\nbigr_X$-triple
$\gbigt(E):=(\gbige,\gbige,\gbigc)$ 
with the Hermitian sesqui-linear
duality $\gbigs:\gbigt(E)\lrarr\gbigt(E)$.
Let us show that $(\gbigt(E),\gbigs)$ is 
a polarized wild pure twistor $D$-module.
We may and will assume that $X=\Delta$
and $D=\{O\}$.
We have already known that
it is strictly $S$-decomposable
(Lemma \ref{lem;08.12.15.1}).
It is easy to check that $\gbige$ is holonomic.

Let $g(z):=z^n$ and $\gminia\in\cnum[t_m^{-1}]$.
Let us show that
$P\Gr^{W(N)}_l\psitilde_{g,\gminia,u}\gbigt(E)$
with the induced Hermitian sesqui-linear duality
is a polarized pure twistor structure of weight $l$.
Let $\varphi_m:\cnum_{t_m}\lrarr \cnum_t$
given by $\varphi_m(t_m)=t_m^m$.
The induced map
$X\times\cnum_{t_m}\lrarr X\times\cnum_t$
is also denoted by $\varphi_m$.
Let $\pi:\Xtilde\lrarr X$ be given by 
$\pi(\zeta)=\zeta^{me}$
such that $(\Etilde,\delbar_{\Etilde},\thetatilde,\htilde):=
 \pi^{-1}(E,\delbar_E,\theta,h)$
  is unramifiedly good.
The induced morphism
$\Xtilde\times\cnum_{t_m}\lrarr X\times\cnum_{t_m}$
is also denoted by $\pi$.
We put $\pitilde:=\varphi_m\circ\pi$.
Let $i_g:X\lrarr X\times\cnum_t$ denote 
the graph of $g$,
and let $\Gamma$ denote the image.
Let $\omega_m$ be a primitive $m$-th root of $1$.
We have the following:
\[
 \pitilde^{-1}\Gamma
=\bigcup_{p=0}^{m-1}
 \bigl\{(\zeta,t_m)\,\big|\,
 t_m-\omega_m^p\cdot\zeta^{ne}=0
 \bigr\}
\]
Let $j_p:\Xtilde\lrarr\Xtilde\times\cnum_{t_m}$
be the graph of the functions
$\gamma_p:=\omega_m^p\cdot \zeta^{ne}$.
We put $\gminia_p:=
 \gminia\bigl(\omega_m^p\cdot \zeta^{ne}\bigr)$.

We have the $\nbigr_{\Xtilde}$-triple
$\gbigttilde=(\gbigetilde,\gbigetilde,\gbigctilde)$
with the Hermitian sesqui-linear duality 
$\gbigs:=(\id,\id)$
associated to $(\Etilde,\delbar_{\Etilde},\thetatilde,\htilde)$.
We have the unramifiedly good wild harmonic bundles
$(\Etilde_p,\delbar_{\Etilde_p},\thetatilde_p,\htilde_p):=
\bigl(\Etilde,\delbar_{\Etilde},\thetatilde,\htilde\bigr)
 \otimes L(-\gminia_p)$.
We have the associated $\nbigr_{\Xtilde}$-triples
$\gbigt_p:=\bigl(\gbigetilde_p,\gbigetilde_p,\gbigctilde_p\bigr)$
with the Hermitian sesqui-linear duality $\gbigs_p:=(\id,\id)$.
We obtain the
$\nbigr_{\Xtilde\times\cnum_{t_m}}$-triple
$j_{p\dagger}\gbigt_p$
and the localization
$j_{p\dagger}
 \gbigt_p(\ast t_m)$.

\begin{lem}
\label{lem;08.1.16.35}
$P\Gr^{W(N)}_l\psitilde_{t_m,u}(j_{p\dagger}\gbigt_p)$
with the induced Hermitian sesqui-linear dualities
are polarized pure twistor $D$-modules of weight $l$
with $0$-dimensional supports.
\end{lem}
\pf
It follows from Corollary \ref{cor;07.10.19.35}.
\hfill\qed

\begin{lem}
\label{lem;08.1.16.50}
$\pi_{\dagger}j_{p\dagger}\gbigt_p(\ast t_m)$
are strictly specializable along $t_m$,
and 
\[
 P\Gr^{W(N)}_l\psitilde_{t_m,u}(
 \pi_{\dagger}j_{p\dagger}\gbigt_p)
\]
with the induced Hermitian sesqui-linear duality
are polarized pure twistor structures of weight $l$.
\end{lem}
\pf
It can be shown 
by using the argument in \cite{saito1},
\cite{sabbah2} or Sections 14.6 of \cite{mochi2}
(see Lemma \ref{lem;08.1.16.40}
and Proposition \ref{prop;05.1.15.61} below
 in this paper)
with Lemma \ref{lem;08.1.16.35}.
Note Assumption \ref{assumption;05.1.15.31} below
is satisfied in this case.
\hfill\qed

\begin{lem}
\label{lem;08.1.16.30}
$\varphi_m^{\dagger}\bigl(
 i_{g\dagger}\gbigt(E)(\ast t)
 \bigr)
 \otimes \gbigl(-\gminia)$
is a direct summand of
\[
 \bigoplus_{p=0}^{m-1}
 \pi_{\dagger}\bigl(
 j_{p\dagger}\gbigt_p(\ast t_m)
 \bigr).
\]
\end{lem}
\pf
Let $i_{g\dagger}\gbige(\ast t)$
denote 
$i_{g\dagger}\gbige
 \otimes_{\nbigo_{\nbigx\times\cnum_t}}
 \nbigo_{\nbigx\times\cnum_t}(\ast t)$.
Let $j_{p\dagger}\gbigetilde_p(\ast t_m)$
denote $j_{p\dagger}\gbigetilde_p
 \otimes_{\nbigo_{\nbigxtilde\times\cnum_{t_m}}}
 \nbigo_{\nbigxtilde\times\cnum_{t_m}}(\ast t_m)$.
We have the following natural isomorphism:
\[
 \pi^{\dagger}\bigl(
 \varphi_m^{\dagger}
 \bigl(i_{g\dagger}\gbige(\ast t)\bigr)
 \bigr)
\simeq
 \pitilde^{\dagger}
 \bigl(i_{g\dagger}\gbige(\ast t)\bigr)
\simeq
 \bigoplus j_{p\dagger}\gbigetilde(\ast t_m)
\]
Therefore,
we have the following:
\begin{equation}
 \label{eq;08.1.16.26}
 \bigoplus_{p=0}^{m-1}
 j_{p\dagger}\gbigetilde_p(\ast t_m)
\simeq
 \pi^{\dagger}\Bigl(
  \varphi_m^{\dagger}
 \bigl(i_{g\dagger}\gbige(\ast t)\bigr)
\otimes\nbigl(-\gminia)
 \Bigr)
\end{equation}
Note that the multiplication of
$\zeta\cdot t_m$ on 
$j_{p\dagger}\gbigetilde_p(\ast t_m)$
is invertible.
Hence, 
we obtain the following isomorphism
from (\ref{eq;08.1.16.26})
by using Lemma \ref{lem;08.1.16.25}:
\begin{equation}
 \label{eq;08.1.17.1}
 \bigoplus_{p=0}^{m-1} 
 \pi_{\dagger}\bigl(
 j_{p\dagger}\gbigetilde_p(\ast t_m)
 \bigr)
\simeq
  \varphi_m^{\dagger}
 \bigl(i_{g\dagger}\gbige\otimes\nbigo(\ast t)\bigr)
\otimes\nbigl(-\gminia)
 \otimes 
 \pi_{\ast}\nbigo_{\nbigxtilde\times\cnum_{t_m}}(\ast t_m)
\end{equation}
We have the natural $\mu_{em}$-action on
$\pi_{\dagger}\Bigl(
 \bigoplus j_{p\dagger}\gbigetilde_p(\ast t_m)
 \Bigr)$,
which is the same as the one
induced by the $\mu_{em}$-action on
$\pi_{\ast}\nbigo_{\nbigxtilde\times\cnum_{t_m}}(\ast t_m)$
in the right hand side of (\ref{eq;08.1.17.1}).
We have a natural morphism from 
$\varphi_m^{\dagger}
 \bigl(i_{g\dagger}\gbige(\ast t)\bigr)
\otimes\nbigl(-\gminia)$
to the $\mu_{em}$-invariant part.
Its restriction to the outside of $\{t_m=0\}$
is an isomorphism,
and the multiplication of $t_m$
are invertible on both of them.
Hence, 
$\varphi_m^{\dagger}
 \bigl(i_{g\dagger}\gbige(\ast t)\bigr)
\otimes\nbigl(-\gminia)$
is identified with
the $\mu_{em}$-invariant part
of $\pi_{\dagger}\Bigl(
 \bigoplus j_{p\dagger}\gbigetilde_p(\ast t_m)
 \Bigr)$.

To show that the isomorphism is compatible 
with the sesqui-linear pairing,
we have only to compare them
on $(X\times\cnum_m)\setminus(X\times\{0\})$,
where the claim is clear.
\hfill\qed

\begin{lem}
\label{lem;07.12.1.8}
$P\Gr^{W(N)}_l\psitilde_{t,\gminia,u}\bigl(\gbigt(E)\bigr)$
with the naturally induced Hermitian sesqui-linear duality
is a polarized pure twistor structure of weight $l$.
\end{lem}
\pf
Since 
$P\Gr^{W(N)}_l\psitilde_{t,\gminia,u}\bigl(
 \gbigt(E)\bigr)$
is a direct summand of
\[
 \bigoplus_p
P\Gr^{W(N)}_l\psitilde_{t_m,u}\bigl(
 \pi_{\dagger}j_{p\dagger}\gbigt_p
 \bigr),
\]
Lemma \ref{lem;07.12.1.8}
follows from Lemma \ref{lem;08.1.16.50}.
\hfill\qed

\begin{cor}
\label{cor;10.6.6.21}
The above $\nbigr_X$-triple
$\gbigt(E)$ with the Hermitian duality
a polarized wild pure twistor $D$-module
of weight $0$.
\hfill\qed
\end{cor}

\subsection{End of the proof of 
 Proposition \ref{prop;07.10.19.50}}
\label{subsection;10.6.6.11}

By Corollary \ref{cor;10.6.6.21},
we obtain that $\Phi$ is essentially surjective.
Let us show that $\Phi$ is fully faithful.
It is easy to deduce that
$\Phi$ is faithful by Lemma \ref{lem;07.11.2.15}.
Let us show that $\Phi$ is full.

Let $(\nbigt,\nbigs)\in\MPT^{\wild}_{\strict}(X,D,0)$,
where $\nbigt=(\nbigm,\nbigm,C)$
and $\nbigs=(\id,\id)$,
and let $(E,\delbar_E,\theta,h)$ be 
the corresponding wild harmonic bundle on $X-D$.
We have the associated wild pure twistor $D$-module
$\gbigt(E)=(\gbige,\gbige,\gbigc)$ of weight $0$
with the polarization $\gbigs=(\id,\id)$.
We have only to show that
the isomorphism
$\gbigt(E)_{|X-D}\simeq
 \nbigt_{|X-D}$ is extended
to $\gbigt(E)\simeq \nbigt$.
Since the property is local,
we may and will assume $X=\Delta$ and $D=\{O\}$.
We have the natural isomorphism
$\nbigm_{|(X-D)\times\cnum_{\lambda}}
 \simeq
 \gbige_{|(X-D)\times\cnum_{\lambda}}$.
We have only to show that it is
extended to an isomorphism
$\nbigm\simeq\gbige$.
Since both of them are strictly $S$-decomposable,
we have only to show 
$\nbigm(\ast z)\simeq\gbige(\ast z)$.
(See Lemma \ref{lem;07.11.2.15}.)

Let $\varphi:(\Xtilde,\Dtilde)\lrarr (X,D)$ be a ramified covering
such that $\varphi^{-1}\harmonicbundle$ is unramified.
We obtain the $\nbigr_{\Xtilde}(\ast \Dtilde)$-triples
$\varphi^{\dagger}\nbigt(\ast z)$
and $\varphi^{\dagger}\gbigt(E)(\ast z)$.
They are the same on $\Xtilde-\Dtilde$.
By the assumption that $(\nbigt,\nbigs)$ is
a polarized wild pure twistor $D$-module,
$P\Gr^{W(N)}_l\psi_{\ztilde,\gminia,u}\bigl(
 \varphi^{\dagger}\nbigt(\ast z)\bigr)$
with the naturally induced 
Hermitian sesqui-linear duality
are polarized pure twistor structures of weight $l$.
Due to Theorem \ref{thm;10.5.20.4},
we obtain that $\varphi^{\ast}\nbigm(\ast z)$
is naturally isomorphic to 
$\varphi^{\ast}\gbige(\ast z)$,
and hence 
$\nbigm(\ast z)\simeq
 \gbige(\ast z)$.
Thus, we are done.
\hfill\qed

\section{Gysin map for wild pure twistor $D$-modules}
\label{subsection;08.1.18.13}
We consider the Gysin maps for
wild pure twistor $D$-modules,
following M. Saito's argument in the Hodge case
(see \cite{saito1}).

\subsection{$\nbigr$-module $j_{\ast}j^{\ast}\nbigm$}
\label{subsection;08.1.17.25}

Let $X$ be a complex manifold,
and let $i_Y:Y\subset X$ be a smooth hypersurface of $X$.
We put $\nbigx:=\cnum_{\lambda}\times X$
and $\nbigy:=\cnum_{\lambda}\times Y$.
We have the following $\nbigr_X$-submodule of 
$\nbigo_{\nbigx}(\ast \nbigy)$:
\[
 j_{\ast}j^{\ast}\nbigo_{\nbigx}:=
 \nbigr_{X}\cdot\nbigo_{\nbigx}(\nbigy)
\subset
 \nbigo_{\nbigx}(\ast \nbigy)
\]
\index{sheaf $j_{\ast}j^{\ast}\nbigo_{\nbigx}$}
If $Y=\{t=0\}$,
it is equal to
$\sum_m \nbigo_{\nbigx}\cdot 
 (\lambda\cdot t^{-1})^m\cdot t^{-1}$.
Thus, it is holonomic and strict.
For an $\nbigr_X$-module $\nbigm$,
we define
\[
 j_{\ast}j^{\ast}\nbigm:=
 \nbigm\otimes_{\nbigo_{\nbigx}}
j_{\ast}j^{\ast}\nbigo_{\nbigx}
\subset
 \nbigm(\ast \nbigy)
\]
\index{sheaf $j_{\ast}j^{\ast}\nbigm$}
In other words,
$j_{\ast}j^{\ast}\nbigm$
is the submodule of $\nbigm(\ast \nbigy)$
generated by 
$\nbigm(\nbigy):=\nbigm
 \otimes_{\nbigo_{\nbigx}}\nbigo_{\nbigx}(\nbigy)$
over $\nbigr_X$.
If the support of any submodule of $\nbigm$
is not contained in $\nbigy$,
we have the injection 
$\nbigm\lrarr j_{\ast}j^{\ast}\nbigm$.

In the following,
let $\nbigm$ be a strict,
holonomic $\nbigr_X$-module,
and we assume that $Y$ is 
{\em strictly non-characteristic} to $\nbigm$.
(See Section 3.7 of \cite{sabbah2}.)

\begin{lem}
Assume that $X$ has a global coordinate system
$(t,z_1,\ldots,z_{n-1})$ with $Y=\{t=0\}$.
Then, the $\nbigr_X$-modules
$\nbigm$ and
$j_{\ast}j^{\ast}\nbigm$ are strictly specializable
along $t$.
\begin{itemize}
\item
$\KMS(\nbigm,t)$ is contained in
$\seisuu_{\leq -1}\times\{0\}
 \subset\real\times\cnum$,
and the $V$-filtration for $\nbigm$ is given by
$V_{n}(\nbigm)=t^{-1-n}\nbigm$
for $n\in\seisuu_{\leq -1}$
and $V_n(\nbigm)=\nbigm$ for $n\in\seisuu_{\geq 0}$.
\item
$\KMS(j_{\ast}j^{\ast}\nbigm,t)$ is contained in
$\seisuu\times\{0\}\subset\real\times\cnum$.
The $V$-filtration for $j_{\ast}j^{\ast}\nbigm$
is given by
\[
 V_{n}\bigl(
 j_{\ast}j^{\ast} \nbigm
 \bigr)
=\left\{
 \begin{array}{ll}
t^{-n}\nbigm\otimes\nbigo_{\nbigx}(\nbigy)
 &  (n\in\seisuu_{\leq\,0})\\
\sum_{p+q\leq n,\,\,q\leq 0} 
 \deldel_t^p
 V_{q}
 & (n\in\seisuu_{\geq\,1})
 \end{array}
\right.
\]
The induced morphism
$t:\psi_{t,0}(j_{\ast}j^{\ast}\nbigm)
 \lrarr \psi_{t,-\vecdelta_0}(j_{\ast}j^{\ast}\nbigm)$
 is an isomorphism,
and 
$\deldel_t:
\psi_{t,-\vecdelta_0}(j_{\ast}j^{\ast}\nbigm)
 \lrarr 
 \psi_{t,0}(j_{\ast}j^{\ast}\nbigm)$
is trivial.
Note $V_{n}(j_{\ast}j^{\ast}\nbigm)
=V_n(\nbigm)$ for $n\in\seisuu_{\leq -1}$,
and $\vecdelta_0=(1,0)$.
\end{itemize}
Remark that
the filtration $V$ is independent
of the choice of a coordinate system.
\end{lem}
\pf
The first claim follows from
Lemma 3.7.4 of \cite{sabbah2}.
Let us consider the second claim.
It is easy to check 
Condition \ref{condition;08.1.23.10}.
By construction,
$V_n(j_{\ast}j^{\ast}\nbigm)=
 V_n(\nbigm)$ for $n\leq -1$.
Hence,
$\Gr^V_{n}(j_{\ast}j^{\ast}\nbigm)$
are strict for $n\leq -1$.
We can easily check 
(i) $t:\Gr^V_{0}\lrarr \Gr^V_{-1}$ is 
an isomorphism,
(ii) $\deldel_t:\Gr^V_{-1}\lrarr \Gr^V_0$ is $0$.
In particular,
it follows that 
$\Gr^V_{0}(j_{\ast}j^{\ast}\nbigm)$
is strict.
We can also easily check the vanishing
of the action of $-\deldel_tt-n\lambda$ 
on $\Gr_n^V(j_{\ast}j^{\ast}\nbigm)$.
It follows that $t:\Gr^V_{n}\lrarr \Gr^V_{n-1}$
are injective for $n\geq 1$,
and hence
we obtain that
$\Gr^V_{n}(j_{\ast}j^{\ast}\nbigm)$
are strict for $n\geq 1$.
Thus, we are done.
\hfill\qed

\vspace{.1in}

Let $\nbign_{Y/X}$ denote the sheaf of sections of 
the normal bundle of $Y$ in $X$.
We put $\nbign_{\nbigy/\nbigx}:=
 \lambda\cdot p_{\lambda}^{\ast}\nbign_{Y/X}$,
where $p_{\lambda}$ denotes the projection
$\nbigx\lrarr X$.
\index{sheaf $\nbign_{Y/X}$}
\index{sheaf $\nbign_{\nbigy/\nbigx}$}

\begin{lem}
\label{lem;07.12.1.10}
We have a naturally induced isomorphism:
\[
 \Lambda':j_{\ast}j^{\ast}(\nbigm)\big/\nbigm
\simeq
\lambda^{-1}\cdot i_{Y\dagger}i_Y^{\dagger}\nbigm
\]
\end{lem}
\pf
We have naturally induced isomorphisms
\[
j_{\ast}j^{\ast}(\nbigm)\big/\nbigm
\simeq
 \nbigm\otimes_{\nbigo_{\nbigx}}
 j_{\ast}j^{\ast}(\nbigo_{\nbigx})
 \big/\nbigo_{\nbigx},
\quad
i_{Y\dagger}i_Y^{\dagger}\nbigm
\simeq
 \nbigm\otimes_{\nbigo_{\nbigx}} 
i_{Y\dagger}i_Y^{\dagger}\nbigo_{\nbigx}.
\]
Hence,
we have only to consider 
the case $\nbigm=\nbigo_{\nbigx}$.
Recall that
$i_{Y\dagger}i_Y^{\dagger}\nbigo_{\nbigx}$
is isomorphic to
$i_{Y\ast}(\nbign_{\nbigy/\nbigx}\otimes
 \Sym^{\bullet}\nbign_{\nbigy/\nbigx})$.
It is equipped with
the filtration given by
$V_m=i_{Y\ast}(\nbign_{\nbigy/\nbigx}\otimes
 \Sym^{\bullet\leq m}\nbign_{\nbigy/\nbigx})$
for $m\in\seisuu_{\geq\,0}$.
Note that
$V_0=i_{Y\ast}\nbign_{\nbigy/\nbigx}
\subset i_{Y\dagger}i_Y^{\dagger}\nbigo_{\nbigx}$
generates 
$i_{Y\dagger}i_Y^{\dagger}\nbigo_{\nbigx}$
over $\nbigr_X$.

We have the subsheaf
$\nbigo_{\nbigx}(\nbigy)\subset 
j_{\ast}j^{\ast}\nbigo_{\nbigx}$
and a naturally defined surjection
\[
\nbigo_{\nbigx}(\nbigy)\lrarr 
\nbigo_{\nbigx}(\nbigy)_{|\nbigy}
\simeq p_{\lambda}^{\ast}\nbign_{Y/X}
\simeq \lambda^{-1}\cdot V_0.
\]
We claim 
(i) it is uniquely extended to a morphism
$j_{\ast}j^{\ast}\nbigo_{\nbigx}\lrarr 
 \lambda^{-1}\cdot 
 i_{Y\dagger}i_Y^{\dagger}\nbigo_{\nbigx}$
as $\nbigr_X$-modules,
(ii) the kernel is $\nbigo_{\nbigx}$.
Because of the uniqueness,
we have only to check them locally.
Then, the claims can be checked by a direct calculation.
\hfill\qed

\vspace{.1in}
We have the isomorphism
$\Lambda:j_{\ast}j^{\ast}\nbigm\big/\nbigm
 \simeq i_{Y\dagger}i_Y^{\dagger}\nbigm$
given by $\Lambda=-\sqrt{-1}\lambda\cdot \Lambda'$.
In the case
$X=\{(t,z_1,\ldots,z_{n-1})\}$,
$Y=\{t=0\}$ and $\nbigm=\nbigo_{\nbigx}$,
we have
$\Lambda(\sqrt{-1}t^{-1})
=[\deldel_t]$,
where the latter denotes the naturally defined section
of $\nbign_{\nbigy/\nbigx}$.

\subsection{$\nbigr$-module $j_!j^{\ast}\nbigm$}
\label{subsection;08.1.18.1}

Let $\nbigm$ and $Y$ be as above.
Let us construct an $\nbigr_X$-module
$j_!j^{\ast}\nbigm$.
\index{sheaf $\jbikkuri j^{\ast}\nbigm$}
First, let us consider the case $\nbigm=\nbigo_{\nbigx}$.
As an $\nbigo_{\nbigx}$-module,
we define
\[
 j_!j^{\ast}\nbigo_{\nbigx}:=
 \nbigo_{\nbigx}\oplus 
 i_{Y\dagger}i_Y^{\dagger}\nbigo_{\nbigx}
=\nbigo_{\nbigx}\oplus
 i_{Y\ast}\bigl(
 \Sym^{\bullet}\nbign_{\nbigy/\nbigx}
\otimes\nbign_{\nbigy/\nbigx}
 \bigr)
\]
The action of $\Theta_{\nbigx}$ 
on $j_!j^{\ast}\nbigo_{\nbigx}$
is given as follows:
\[
 v\cdot (s_1,s_2)
=\bigl(v\cdot s_1, \,\,
 v\cdot s_2+\pi(v_{|Y})\cdot s_{1|Y}
 \bigr)
\]
Here, ``$|Y$'' denotes the restriction to $Y$,
$\pi$ denotes a projection
of $\Theta_{\nbigx|\nbigy}
 \lrarr \nbign_{\nbigy/\nbigx}$,
and $v\cdot s_2$ is given by
the natural $\nbigr_X$-module structure
on $i_{Y\dagger}i_Y^{\dagger}\nbigo_{\nbigx}$.
It is uniquely extended to an action of $\nbigr_X$
on $j_!j^{\ast}\nbigo_{\nbigx}$.
We have only to check it on 
a coordinate neighbourhood such that
$Y=\{t=0\}$ by a direct calculation.

In the general case,
we put 
$j_!j^{\ast}\nbigm:=
 \nbigm\otimes j_!j^{\ast}\nbigo_{\nbigx}$.
(Note that $Y$ is assumed to be
strictly non-characteristic to $\nbigm$.)
By construction,
we have the following exact sequence
of $\nbigr_X$-modules:
\begin{equation}
 \label{eq;07.12.1.11}
 0\lrarr
 i_{\dagger}i^{\dagger}\nbigm\lrarr
 j_!j^{\ast}\nbigm\lrarr \nbigm\lrarr 0
\end{equation}

\begin{lem}
Assume that $X$ has a global coordinate system
$(t,z_1,\ldots,z_{n-1})$  with $Y=\{t=0\}$.
The $\nbigr_X$-module $j_!j^{\ast}\nbigm$
is strictly specializable along $t$.
The set $\KMS(j_!j^{\ast}\nbigm)$
is contained in $\seisuu\times\{0\}$.
The $V$-filtration is given 
by
$V_n(j_!j^{\ast}\nbigm)
=V_n(\nbigm)\oplus 
 V_n(i_{\dagger}i^{\dagger}\nbigm)$.
More concretely,
\[
 V_n(j_!j^{\ast}\nbigm)
=\left\{
 \begin{array}{ll}
 t^{-n-1}\nbigm & (n\in\seisuu_{\leq -1})\\
 \nbigm\oplus
 \sum_{p\leq n}
 i_{Y\ast}\bigl(
 \Sym^{p}\nbign_{\nbigy/\nbigx}
 \otimes
 \nbign_{\nbigy/\nbigx}
 \otimes\nbigm
 \bigr)
 & (n\in \seisuu_{\geq \,0})
 \end{array}
 \right.
\]
The induced map
$t:\psi_{t,0}(j_!j^{\ast}\nbigm)
\lrarr
 \psi_{t,-\vecdelta_0}(j_!j^{\ast}\nbigm)$
is $0$,
and $\deldel_t:
 \psi_{t,-\vecdelta_0}(j_!j^{\ast}\nbigm)
\lrarr
 \psi_{t,0}(j_!j^{\ast}\nbigm)$
is an isomorphism.

Note that the filtration $V$ is independent
of the choice of a coordinate system.
\hfill\qed
\end{lem}

\subsection{Sesqui-linear pairings
$j_!j^{\ast}C$ and $j_{\ast}j^{\ast}C$}
\index{pairing $\jbikkuri j^{\ast}C$}
\index{pairing $j_{\ast}j^{\ast}C$}

Let $(\nbigm',\nbigm'',C)$ be an $\nbigr_X$-triple
such that $\nbigm'$ and $\nbigm''$ are
strict and holonomic.
Let $Y$ be strictly non-characteristic to 
$\nbigm'$ and $\nbigm''$.
We would like to construct a sesqui-linear pairing
$j_{\ast}j^{\ast}C$ (resp. $j_{!}j^{\ast}C$)
of $j_{!}j^{\ast}\nbigm'$ 
(resp. $j_{\ast}j^{\ast}\nbigm'$)
and $j_{\ast}j^{\ast}\nbigm''$
(resp. $j_!j^{\ast}\nbigm''$).

We explain the construction of $j_{\ast}j^{\ast}C$.
Let us consider the case in which
$X$ has a global coordinate system
with $Y=\{t=0\}$.
Let $\lambda_0\in\vecS$,
and let $U(\lambda_0)$ denote a small neighbourhood
of $\lambda_0$ in $\cnum_{\lambda}^{\ast}$.
We put $\vecI(\lambda_0):=U(\lambda_0)\cap \vecS$.
Let $f$ (resp. $g$) be a local section
of $\nbigm'$ (resp. $\nbigm''$)
on $U(\lambda_0)\times X$
(resp. $\sigma(U(\lambda_0))\times X$).
Let $\phi$ be a $C^{\infty}$-top form on $X$
with compact support.
Then,
$\bigl\langle
 C(f,\sigma^{\ast}g),\,\,
 |t|^{2s}\cdot \tbar^{-1}\phi
 \bigr\rangle$
gives a continuous function on 
$\vecI(\lambda_0)\times
 \bigl\{s\in\cnum\,\big|\,\Re(s)>\!>0\bigr\}$,
which is holomorphic with respect to $s$.

\begin{lem}
\label{lem;08.1.17.10}
$\bigl\langle
 C(f,\sigma^{\ast}g),
 |t|^{2s}\cdot\tbar^{-1}\phi
 \bigr\rangle$
gives a continuous function on
$\vecI(\lambda_0)\times
\bigl\{s\in\cnum\,\big|\,
 \Re(s)>-1\bigr\}$
which is holomorphic with respect to $s$.
\end{lem}
\pf
Since it can be shown by a standard argument,
we give only an outline.
Let $\omega:=dt\cdot d\tbar
 \cdot\prod_{j=1}^{n-1} dz_j\cdot d\zbar_j$.
We take a sufficiently large $N$.
A test function $\phi$ on $X$
can be decomposed as
\[
 \phi=
 \sum_{i+j\leq N} \phi_{i,j}\cdot t^i\cdot \tbar^j
 \cdot\omega
+\varphi\cdot\omega,
\]
where $\phi_{i,j}$ are test functions
which are constant with respect to $t$ around $t=0$,
and $\varphi$ is a test function
such that $\varphi=O\bigl(|t|^N\bigr)$ 
around $t=0$.
If $N$ is sufficiently large,
the contribution of $\varphi$ 
to $\bigl\langle
 C(f,\sigma^{\ast}g),
 |t|^{2s}\cdot\tbar^{-1}\phi
 \bigr\rangle$
is holomorphic around $s=0$.
Let us consider the following:
\[
 F_{i,j}(s):=
 \bigl\langle
 C(f,\sigma^{\ast}g),\,\,|t|^{2s}\cdot
 t^i\cdot \tbar^{j-1}\cdot \phi_{i,j}\cdot\omega
 \bigr\rangle
=\bigl\langle 
 C\bigl(t^{i+1}f,\,\,
 \sigma^{\ast}(t^{j}g)
 \bigr),\,\,
 |t|^{2(s-1)}\phi_{i,j}\cdot\omega
 \bigr\rangle
\]
Note $t^{i+1}f\in V_{-i-2}(\nbigm')$
and $t^jg\in V_{-j-1}(\nbigm'')$.
We put
$b_{M,L}:=
 \prod_{m=0}^L\bigl(
 -\del_tt+i+2+m \bigr)^M$.
Using a standard argument,
we obtain the following equality,
modulo continuous functions on
$\vecI(\lambda_0)\times\cnum_s$
which are holomorphic with respect to $s$:
\begin{multline}
 \label{eq;08.10.13.2}
 \bigl\langle
 C\bigl(b_{M,L}(t^{i+1}f),\,
 \sigma^{\ast}(t^jg)\bigr),\,\,
 |t|^{2(s-1)}\cdot\phi_{i,j}\cdot\omega
 \bigr\rangle
\equiv \\
 \prod_{m=0}^L(s+i+1+m)^M
\bigl\langle
 C\bigl(t^{i+1}f,\,
 \sigma^{\ast}(t^jg)\bigr),\,\,
 |t|^{2(s-1)}\cdot\phi_{i,j}\cdot\omega
 \bigr\rangle
\end{multline}
If $M$ is sufficiently large,
there exists a local section $P\in V_0\nbigr_X$
such that 
$b_{M,L}(t^{i+1}f)
=P\cdot (t^{L+1+i+1}f)$.
Hence, the left hand side of (\ref{eq;08.10.13.2}) is
\begin{equation}
 \label{eq;08.11.25.1}
\bigl\langle
 C\bigl(t^{i+1}f,\,\sigma^{\ast}(t^jg)\bigr),\,\,
 t^{L+1}\cdot|t|^{2(s-1)}
 \cdot Q\phi_{i,j}\cdot\omega
 \bigr\rangle 
\end{equation}
for some $Q\in V_0\nbigr_X$.
If $L$ is sufficiently large,
(\ref{eq;08.11.25.1})
is continuous and holomorphic with respect to $s$
on $\vecI(\lambda_0)\times
 \{\Re(s)>-1\}$.
Thus, the claim of the lemma follows.
\hfill\qed

\vspace{.1in}

Let $\nbigm''(\nbigy):=
 \nbigm''\otimes_{\nbigo_{\nbigx}}\nbigo_{\nbigx}(\nbigy)$,
which is a $V_0\nbigr_X$-submodule
of $j_{\ast}j^{\ast}\nbigm''$.
We obtain the following pairing:
\begin{equation}
\label{eq;08.10.13.3}
  \Ctilde:
 \nbigm'_{|\vecS\times X}\otimes
\sigma^{\ast}\nbigm''(\nbigy)_{|\vecS\times X}
\lrarr
 \distribution_{\vecS\times X/\vecS},
\end{equation}
\begin{equation}
  \bigl\langle
 \Ctilde(f,\sigma^{\ast}(t^{-1}g)),\,
 \phi
 \bigr\rangle
=\bigl\langle
 C(f,\sigma^{\ast}g),\,
 |t|^{2s}\cdot\tbar^{-1}\cdot\phi
 \bigr\rangle_{|s=0}.
\end{equation}
Due to Lemma \ref{lem;07.10.26.25},
the restriction of $\Ctilde$ to
$\nbigm'_{|\vecS\times X}
\otimes
 \sigma^{\ast}\nbigm''_{|\vecS\times X}$
is equal to $C$.

\begin{lem}
The pairing
$\Ctilde$ is a 
$V_0\nbigr_X\otimes 
 \sigma^{\ast}V_0\nbigr_X$-homomorphism.
It is independent of the choice of
a coordinate system..
\end{lem}
\pf
Let us consider
$\bigl\langle
 \Ctilde\bigl(Pf,\sigma^{\ast}(Q(t^{-1}g))\bigr),\phi
 \bigr\rangle$.
We have $Q'\in V_0\nbigr_X$ such that 
$Qt^{-1}=t^{-1}Q'$ in $\nbigr_X$.
By definition,
we have the following:
\begin{multline}
\label{eq;07.8.30.2}
 \Bigl\langle
 \Ctilde\bigl(Pf,\sigma^{\ast}(Q(t^{-1}g))\bigr),\,
 \phi \Bigr\rangle
=\Bigl\langle
 \Ctilde\bigl(Pf,\sigma^{\ast}(t^{-1}Q'g)\bigr),\,
 \phi \Bigr\rangle\\
=\Bigl\langle
 C\bigl(Pf,\sigma^{\ast}(Q'g)\bigr),\,
 |t|^{2s}\tbar^{-1}\cdot\phi
 \Bigr\rangle_{|s=0}
=\Bigl\langle
 C(f,\sigma^{\ast}(g)),\,
 \lefttop{t}P\cdot
 \overline{\sigma^{\ast}(\lefttop{t}Q't^{-1})}
 \bigl(|t|^{2s}\cdot\phi\bigr)
 \Bigr\rangle_{|s=0}
\end{multline}
Note
$\lefttop{t}Q't^{-1}=t^{-1}\lefttop{t}Q$.
There is a $C^{\infty}$-top form $\varphi$
with compact support,
such that the following holds:
\[
 \lefttop{t}P\cdot
 \overline{
 \sigma^{\ast}(\lefttop{t}Q)}
 \bigl(
 |t|^{2s}\cdot\phi
 \bigr)
=|t|^{2s}\cdot
 \bigl(
 \lefttop{t}P\bigr)\cdot
  \overline{\sigma^{\ast}(\lefttop{t}Q)}
 \phi
+s\cdot |t|^{2s}\varphi
\]
Thus, (\ref{eq;07.8.30.2}) can be rewritten as follows:
\[
 \Bigl\langle
 C(f,\sigma^{\ast}(g)),\,\,
 \tbar^{-1}|t|^{2s}\cdot
 \bigl(\lefttop{t}P\bigr)
\cdot \overline{\sigma^{\ast}(\lefttop{t}Q)}\cdot
 \phi
 \Bigr\rangle_{|s=0}
=\Bigl\langle
 \Ctilde(f,\sigma^{\ast}(t^{-1}g)),\,
 \lefttop{t}P\cdot
 \overline{\sigma^{\ast}(\lefttop{t}Q)}
 \phi
 \Bigr\rangle
\]
Thus, we obtain the first claim.
The second claim is clear by construction.
\hfill\qed

\vspace{.1in}

Note that $j_{!}j^{\ast}\nbigm'$ and
$j_{\ast}j^{\ast}\nbigm''$
are generated by $\nbigm'$ and $\nbigm''(\nbigy)$
over $\nbigr_X$,
respectively.
Let $\deldelbar_t:=-\lambda^{-1}\cdot\delbar_t$,
which is identified with
$\sigma^{\ast}(\deldel_t)$.
We would like to extend $\Ctilde$
to the pairing
\begin{equation}
 j_{\ast}j^{\ast}C:
 j_!j^{\ast}\nbigm'\otimes
 \sigma^{\ast}\bigl(j_{\ast}j^{\ast}\nbigm''\bigr)
\lrarr
 \distribution_{\vecS\times X/\vecS}
\end{equation}
by the formula
\begin{equation}
\label{eq;07.8.30.3}
 j_{\ast}j^{\ast}C\Bigl(
 \sum \deldel_t^k a_k,\,
 \sigma^{\ast}\Bigl(
 \sum\deldel_t^l b_l
 \Bigr)
 \Bigr)
=\sum \deldel_t^k\deldelbar_t^l
 \Ctilde(a_k,\sigma^{\ast}b_l)
\end{equation}
for $a_k\in\nbigm'$ and $b_l\in\nbigm''(\nbigy)$.
We have to check the well-definedness.

\begin{lem}
If $\sum_{l}\deldel_t^lb_l=0$ in $j_{\ast}j^{\ast}\nbigm'$,
we have 
$\sum_{l} \deldelbar_t^l \Ctilde(a_k,\sigma^{\ast}b_l)=0$.
\end{lem}
\pf
Since it can be shown by a standard argument
(see \cite{sabbah2} or \cite{mochi2},
 for example),
we give only an outline.
We use the notation in the proof of 
Lemma \ref{lem;08.1.17.10}.
We put 
$\Phi:=\sum_{l}
 \deldelbar_t^l \Ctilde(a_k,\sigma^{\ast}b_l)$.
Since the support of $\Phi$ is contained in $t=0$,
$\Phi$ is of the form
$\sum_{p,q}A_{p,q}\cdot 
 \del_t^p\delbar_t^q\delta_{Y}$,
where $A_{p,q}$ denote distributions
on $Y$, and $\delta_{Y}$
denotes the $\delta$-function at $Y$.
We put 
$b_{M,L}:=
 \prod_{m=0}^L(-\del_t\cdot t+m)^M$.
If $M$ is sufficiently large,
there exists $Q\in V_0\nbigr_X$
such that 
$b_{M,L}a_k=
 Q\cdot t^{L+1}\cdot a_k$.
Hence
$b_{M,L}\Phi=0$ for any large $L$.
We also have the relation
$(-\del_tt+m)\del^p\delta_Y
=(p+m)\del^p\delta_Y$.
Then, it is easy to derive $\Phi=0$.
\hfill\qed

\begin{lem}
If $\sum_{k=0}^N \deldel_t^k a_k=0$ 
 in $j_!j^{\ast}\nbigm'$,
we have
$\sum_{k=0}^N 
 \deldel_t^k \Ctilde(a_k,\sigma^{\ast}b_l)=0$.
\end{lem}
\pf
In the case $N=0$, the claim is clear.
Let us consider the case $N=1$.
Assume
$a_0+\deldel_t a_1=0$
in $j_!j^{\ast}\nbigm'=0$.
Then, we have
$a_{1|Y}=0$ 
and $a_0+\deldel_t a_1=0$ in $\nbigm'$.
In particular, we have $a_1\in V_{-2}\nbigm'$,
and hence $a_1=t\cdot a_1'$
for some $a_1'\in \nbigm'$.
We have the following:
\begin{equation}
 \label{eq;08.10.13.1}
 \deldel_t \Ctilde(a_1,\sigma^{\ast}b_l)
=\deldel_t \Ctilde(t a_1',\sigma^{\ast}b_l)
=\deldel_t t \Ctilde(a_1',\sigma^{\ast}b_l)
=\Ctilde(\deldel_t(t a_1'),\,\sigma^{\ast}b_l)
=\Ctilde(\deldel_ta_1,\sigma^{\ast}b_l)
\end{equation}
In particular, we obtain
$\Ctilde(a_0,\sigma^{\ast}b_l)
+\deldel_t\Ctilde(a_1,\sigma^{\ast}b_l)
=\Ctilde(a_0+\deldel_ta_1,\sigma^{\ast}b_l)
=0$.

Let us show the claim for general $N$,
assuming the claim for $N-1$.
If $\sum_{j=0}^N \deldel_t^ja_j=0$
in $j_!j^{\ast}\nbigm'$,
we have $a_N\in V_{-2}$
and $a_N=t\cdot a_N'$ for some $a_N'\in \nbigm'$.
We have the following vanishing in $j_{!}j^{\ast}\nbigm$:
\[
 \sum_{j=0}^{N-2}\deldel_t^ja_j
+\deldel_t^{N-1}\bigl(a_{N-1}+\deldel_t t a_N'\bigr)=0
\]
By the assumption of the induction on $N$,
the following holds:
\[
 \sum_{j=0}^{N-2}\deldel_t^j
 \Ctilde\bigl(a_j,\,\sigma^{\ast}b_l\bigr)
+\deldel_t^{N-1}
 \Ctilde\bigl(a_{N-1}+\deldel_t (ta_N'),
 \,\,\sigma^{\ast}b_l\bigr)
=0
\]
We have
$\Ctilde(a_{N-1}+\deldel_t (ta_N'),\sigma^{\ast}b_l)
=\Ctilde(a_{N-1},\sigma^{\ast}b_l)
+\deldel_t\Ctilde(t\cdot a_{N}',\sigma^{\ast}b_l)$
as in (\ref{eq;08.10.13.1}).
Then, the induction can proceed.
\hfill\qed

\vspace{.1in}
Thus, we obtain the pairing 
$j_{\ast}j^{\ast}C$ given by 
(\ref{eq;07.8.30.3})
in the case that $X$ is equipped with
a coordinate system such that $Y=\{t=0\}$.
\begin{lem}
$j_{\ast}j^{\ast}C$ is 
an $\nbigr_{X|\vecS\times X}
 \otimes
 \sigma^{\ast}
 \nbigr_{X|\vecS\times X}$-homomorphism.
It is independent of the choice of
a coordinate system.
As a result,
we obtain the globally well defined pairing
$j_{\ast}j^{\ast}C$.
\end{lem}
\pf
The first claim is clear by construction.
The restriction of
$j_{\ast}j^{\ast}C$ to
$\nbigm'
 \otimes
 \nbigm''(\nbigy)$
is equal to $\Ctilde$
given in (\ref{eq;08.10.13.3}),
which is 
is independent of the choice of
a coordinate system.
Since 
$j_{!}j^{\ast}\nbigm'$ and
$j_{\ast}j^{\ast}\nbigm''$
are generated by $\nbigm'$ and $\nbigm''(\nbigy)$
over $\nbigr_X$,
the extension of
$\Ctilde$ to
an $\nbigr_X\otimes
 \sigma^{\ast}\nbigr_X$-homomorphism
of $j_!j^{\ast}\nbigm'$
and $j_{\ast}j^{\ast}\nbigm''$
is unique.
Hence, $j_{\ast}j^{\ast}C$ is independent
of the choice of a coordinate system.
\hfill\qed

\vspace{.1in}

Similarly, we have the globally defined 
sesqui-linear pairing
$j_!j^{\ast}C$ of
$j_{\ast}j^{\ast}\nbigm'$ and
$j_!j^{\ast}\nbigm''$
given by the local formula
\begin{equation}
 j_{!}j^{\ast}C\Bigl(
 \sum\deldel_t^l b_l,\,\,
 \sigma^{\ast}\Bigl(
 \sum \deldel_t^k a_k
 \Bigr)
 \Bigr)
=\sum \deldel_t^l\deldelbar_t^k
 \Ctilde(b_l,\sigma^{\ast}a_k)
\end{equation}
for $b_l\in\nbigm'(\nbigy)$ and $a_k\in\nbigm''$.

\subsection{$\nbigr$-triples $j_{\ast}j^{\ast}\nbigt$
 and $j_{!}j^{\ast}\nbigt$}
\index{$\nbigr$-triple
 $j_{\ast}j^{\ast}\nbigt$}
\index{$\nbigr$-triple
 $\jbikkuri j^{\ast}\nbigt$}

We obtain the following $\nbigr$-triples:
\[
 j_{\ast}j^{\ast}\nbigt:=
 \bigl(j_!j^{\ast}\nbigm',
 j_{\ast}j^{\ast}\nbigm'',
 j_{\ast}j^{\ast}C\bigr),
\quad
 j_{!}j^{\ast}\nbigt:=
 \bigl(j_{\ast}j^{\ast}\nbigm',
 j_{!}j^{\ast}\nbigm'',
 j_{!}j^{\ast}C\bigr)
\]

\begin{lem}
\label{lem;08.1.17.20}
We have the natural identification
$(j_{\ast}j^{\ast}\nbigt)^{\ast}
=j_!j^{\ast}(\nbigt^{\ast})$
and 
$j_{\ast}j^{\ast}(\nbigt^{\ast})
 =j_!j^{\ast}(\nbigt)^{\ast}$.
\end{lem}
\pf
We have only to show one of them.
We have only to consider the case
in which $X$ has a coordinate system
with $Y=\{t=0\}$.
We have 
$ \bigl(
 j_{\ast}j^{\ast}C
\bigr)^{\ast}(t^{-1}g,\sigma^{\ast}(f))
:=
\overline{\sigma^{\ast}
 \bigl(j_{\ast}j^{\ast}C(f,\sigma^{\ast}(t^{-1}g))\bigr)}$
for local sections $f$ and $g$
of $j_!j^{\ast}\nbigm'$
and $j_{\ast}j^{\ast}\nbigm''$,
respectively.
Hence, we obtain the following equality
for local sections $f\in \nbigm'$
and $g\in \nbigm''(\nbigy)$:
\begin{multline}
 \Bigl\langle
 \overline{
 \sigma^{\ast}
 \bigl(j_{\ast}j^{\ast}C(f,\sigma^{\ast}(t^{-1}g))\bigr)},\,
 \phi
 \Bigr\rangle
=\sigma^{\ast}\overline{
 \bigl\langle
 j_{\ast}j^{\ast}C(f,\sigma^{\ast}(t^{-1}g)),\,
 \overline{\sigma^{\ast}\phi}
 \bigr\rangle } =\\
\sigma^{\ast}\overline{
  \bigl\langle
 C(f,\sigma^{\ast}(g)),\,
 |t|^{2s}\tbar^{-1}\overline{\sigma^{\ast}\phi}
 \bigr\rangle_{|s=0}}
=\bigl\langle
 C^{\ast}(g,\sigma^{\ast}f),\,
 |t|^{2s}t^{-1}\phi
 \bigr\rangle_{|s=0} \\
=\bigl\langle
 j_!j^{\ast}C^{\ast}(t^{-1}g,\sigma^{\ast}f),\,
 \phi \bigr\rangle
\end{multline}
Then, the claim of the lemma follows.
\hfill\qed

\vspace{.1in}
It is easy to observe that
$j_{\ast}j^{\ast}$ and $j_!j^{\ast}$
are functorial,
i.e.,
a morphism of $\nbigr$-triples
$\varphi:\nbigt_1\lrarr\nbigt_2$
naturally induces
\[
 j_{\ast}j^{\ast}\varphi:
 j_{\ast}j^{\ast}\nbigt_1
\lrarr
 j_{\ast}j^{\ast}\nbigt_2,
\quad
  j_{!}j^{\ast}\varphi:
 j_{!}j^{\ast}\nbigt_1
\lrarr
 j_{!}j^{\ast}\nbigt_2.
\]

\subsection{$\Cok(\nbigt)$ and $\Ker(\nbigt)$}

\index{$\nbigr$-triple $\Cok(\nbigt)$}
\index{$\nbigr$-triple $\Ker(\nbigt)$}

We have the natural morphism
$\nbigt\lrarr j_{\ast}j^{\ast}\nbigt$.
The cokernel is denoted by $\Cok(\nbigt)$.
We also have the natural morphism
$j_!j^{\ast}\nbigt\lrarr\nbigt$.
The kernel is denoted by $\Ker(\nbigt)$.
It is easy to observe that
$\Ker$ and $\Cok$ are functorial
as in the case of $j_{\ast}j^{\ast}$
and $j_!j^{\ast}$.
The underlying sesqui-linear pairings
of $\Cok(\nbigt)$ and $\Ker(\nbigt)$
are denoted by
$[j_{\ast}j^{\ast}C]$
and $[j_{!}j^{\ast}C]$, respectively.
We also use the symbols
$j_{\ast}j^{\ast}C$
and $j_{!}j^{\ast}C$,
if there are no risk of confusion.

\begin{lem}
\label{lem;08.1.18.15}
We have the natural identification
$\Ker(\nbigt)^{\ast}=\Cok(\nbigt^{\ast})$.
Together with Lemma {\rm\ref{lem;08.1.17.20}},
we obtain the following identification of 
the exact sequences:
\[
 \begin{CD}
 0@>>>\nbigt^{\ast}
 @>>>j_{\ast}j^{\ast}(\nbigt^{\ast})
 @>>> \Cok(\nbigt^{\ast})
 @>>> 0 \\
@. @V{=}VV @V{=}VV @V{=}VV @.\\
 0 @>>>
 \nbigt^{\ast}
 @>>>
 j_!j^{\ast}(\nbigt)^{\ast}
 @>>>
\Ker(\nbigt)^{\ast}
 @>>> 0
 \end{CD}
\]
\end{lem}
\pf
It can be checked directly from the definition.
\hfill\qed

\subsection{Pull back of polarized pure twistor structure
 in the strictly non-characteristic case}
\label{subsection;08.12.16.2}

Let $\nbigt=(\nbigm',\nbigm'',C)$ be
a wild pure twistor $D$-module on $X$ 
of weight $0$
with a polarization 
$\nbigs=(\varphi',\varphi'')$.
Let $i_Y:Y\subset X$ be
a smooth hypersurface
which is strictly non-characteristic
to $\nbigt$.
Recall that we have in this case
the polarized wild pure twistor $D$-module
$i_Y^{\dagger}\nbigt=
 (i_Y^{\dagger}\nbigm',
 i_Y^{\dagger}\nbigm'',
 C_Y)$ on $Y$ of weight $0$
with the induced polarization
$i_Y^{\dagger}\nbigs=
 (i_Y^{\dagger}\varphi',i_Y^{\dagger}\varphi'')$.
(See \cite{sabbah2} for more details.)
It can be seen as follows.
We have only to consider the case
in which $\nbigm'=\nbigm''=:\nbigm$
and $\varphi'=\varphi''=\id$,
to which the general case can be reduced.
\index{$\nbigr$-triple $i_Y^{\dagger}\nbigt$}

Locally $Y$ is defined by a coordinate function $t$.
Note that $i_Y^{\dagger}\nbigm$ is
equal to $\psi_{t,-\vecdelta_0}(\nbigm)$
in this case,
as shown in Lemma 3.7.4 of \cite{sabbah2}.
Moreover, the nilpotent part of $-\deldel_tt$
on $\psi_{t,-\vecdelta_0}(\nbigm)$ is trivial.
We also have the induced sesqui-linear pairing
$C_Y:=\psi_{t,-\vecdelta_0}C$.
Hence, $i_Y^{\dagger}\nbigt$ is given by
$\psi_{t,-\vecdelta_0}\nbigt$ locally.
It is a wild pure twistor $D$-module
with a polarization $(\id,\id)$.

We obtain the strict $S$-decomposability
of $i_Y^{\dagger}\nbigm$ from the local expression
as $\psi_{t,-\vecdelta_0}(\nbigm)$
and the twistor property of $\nbigt$.
We can glue the locally defined
sesqui-linear pairings 
by using the uniqueness of the extension
of the sesqui-linear pairings for
strictly $S$-decomposable $\nbigr$-modules.
(See \cite{sabbah2}, \cite{mochi2}
 or Proposition \ref{prop;a11.23.10} below.)
Thus, we obtain a globally defined
$\nbigr$-triple.
It is a wild pure twistor $D$-module
on $Y$ of weight $0$
with a polarization $(\id,\id)$.

\subsection{Some isomorphisms
 and a polarization of $\Ker(\nbigt)$}
\label{subsection;08.1.18.14}

Let $(\nbigt,\nbigs)$ and $Y$ be as in
Subsection \ref{subsection;08.12.16.2}.
Let $\Lambda$ be the isomorphism
given in Section \ref{subsection;08.1.17.25}.

\begin{lem}
\label{lem;08.1.17.40}
The pair of the morphisms
$\varphi_1:=\bigl(
 \Lambda,
 \id
 \bigr)$
gives an isomorphism
\[
 i_{Y\dagger}i_Y^{\dagger}\nbigt
\simeq
 \Ker(\nbigt)\otimes\Tate^S(-1/2).
\]
\end{lem}
\pf
We have only to consider the case
$\nbigm'=\nbigm''=:\nbigm$.
We have only to check the compatibility
of $\varphi_1$ and the sesqui-linear pairings.
Since $i_{Y\dagger}i_Y^{\dagger}\nbigm$ is strictly
$S$-decomposable,
we have only to check the compatibility
on the generic part.
Hence, we have only to consider the case
in which $\nbigt$ comes from
a harmonic bundle on $X$.

Let $\Tate^S(0)=
(\nbigo_{\nbigx},\nbigo_{\nbigx},C_0)$.
We have the isomorphisms
$i_{Y\dagger}i_Y^{\dagger}\nbigm
\simeq
 \nbigm\otimes 
 i_{Y\dagger}i_Y^{\dagger}\nbigo_{\nbigx}$
and
$j_{\ast}j^{\ast}\nbigm\big/\nbigm
\simeq
 \nbigm\otimes 
 j_{\ast}j^{\ast}\nbigo_{\nbigx}/\nbigo_{\nbigx}$.
For local sections
$f\otimes a\in 
 \nbigm\otimes j_{\ast}j^{\ast}\nbigo_{\nbigx}/\nbigo_{\nbigx}$ 
and $g\otimes b\in
 \nbigm\otimes i_{Y\dagger}i_Y^{\dagger}\nbigo_{\nbigx}$,
we have the following:
\[
 [j_{!}j^{\ast}C]\bigl(
 f\otimes a,\,\sigma^{\ast}(g\otimes b)
 \bigr)
=C(f,\sigma^{\ast}g)
 \cdot [j_{!}j^{\ast}C_0](a,\sigma^{\ast}b)
\]
For local sections
$f\otimes a$ and
$g\otimes b$ of
$\nbigm\otimes i_{Y\dagger}i_Y^{\dagger}\nbigo_{\nbigx}$,
we have the following:
\[
 i_{Y\dagger}i_Y^{\dagger}
 C\bigl(f\otimes a,\,\,
 \sigma^{\ast}(g\otimes b)\bigr)
=C(f,\sigma^{\ast}g)
\cdot
 i_{Y\dagger}i_Y^{\dagger}C_0(a,\sigma^{\ast}b)
\]
Hence, we have only to consider the case
$\nbigt=\Tate^S(0)$,
$X=\Delta$ and $Y=\{0\}$.

Let $\omega=\chi\cdot\sqrt{-1}dt\cdot d\tbar/2\pi$,
where $\chi$ is a test function on $U$
which are constantly $1$ around $0$.
We can check 
$j_!j^{\ast}C_0\bigl(t^{-1},\sigma^{\ast}1\bigr)
=t^{-1}$ directly from the definition.
Then, we have the following:
\[
 \bigl\langle
 j_{!}j^{\ast}C_0\bigl(t^{-1},
 \sigma^{\ast}(\deldel_t\cdot 1)\bigr),\,\,
\omega
 \bigr\rangle
=
 \bigl\langle
 -\lambda^{-1}\delbar_t(t^{-1}),\,
 \omega
 \bigr\rangle
=-\lambda^{-1}
\]
Hence, we obtain the following:
\[
 \Bigl\langle
\sqrt{-1}\lambda
 [j_{!}j^{\ast}C_0]\bigl(
 \sqrt{-1}t^{-1},
 \sigma^{\ast}(\deldel_t\cdot 1)\bigr),\,\,
\omega
 \Bigr\rangle
=1
\]
We also have
$\bigl\langle
 i_{Y\dagger}i_Y^{\dagger}C_0([\deldel_t],
 \sigma^{\ast}[\deldel_t]),\,
 \omega
 \bigr\rangle
=1$. (See Section 1.6.d of \cite{sabbah2}.)
Now, we can check the compatibility of
$\varphi_1$ and the sesqui-linear pairings
easily.
\hfill\qed

\begin{lem}
The pair of morphisms
$\varphi_2:=\bigl(-\id,\Lambda\bigr)$
gives an isomorphism
\[
 \Cok(\nbigt)\otimes\Tate^S(1/2)
\lrarr
 i_{Y\dagger}i_Y^{\dagger}\nbigt.
\]
\end{lem}
\pf
Using an argument in the proof of 
Lemma \ref{lem;08.1.17.40},
we can reduce the issue
to the case $\nbigt=\Tate^S(0)$,
$X=\Delta$ and $Y=\{0\}$.
We have
$ j_{\ast}j^{\ast}C_0\bigl(
 \deldel_t\cdot 1,\,
 \sigma^{\ast}(t^{-1})
 \bigr)
=\lambda\del_t(\tbar^{-1})$
for the sections
$\deldel_t\cdot 1\in j_!j^{\ast}\nbigo_{\nbigx}$
and $t^{-1}\in j_{\ast}j^{\ast}\nbigo_{\nbigx}$.
Let $\omega$ be as in the proof of
Lemma \ref{lem;08.1.17.40}.
Then, we have the following:
\[
 \bigl\langle
 j_{\ast}j^{\ast}C_0\bigl(
 \deldel_t\cdot 1,\,
 \sigma^{\ast}(t^{-1})
 \bigr),\,
\omega
\bigr\rangle
=\bigl\langle
 \lambda\del_t(\tbar^{-1}),\,
\omega
 \bigr\rangle
=\lambda 
\]
Hence, we obtain the following:
\[
 \Bigr\langle
 (\sqrt{-1}\lambda)^{-1}
  [j_{\ast}j^{\ast}C_0]\bigl(
 -\deldel_t\cdot 1,
 \sigma^{\ast}(\sqrt{-1}t^{-1})
 \bigr),\,
\omega
\Bigl\rangle
=1
\]
Then, we can check the desired compatibility
easily.
\hfill\qed

\vspace{.1in}
We obtain an isomorphism
$\varphi_3=(-\Lambda,\Lambda):
 \Cok(\nbigt)\otimes\Tate^S(1/2)
\simeq
 \Ker(\nbigt)\otimes\Tate^S(-1/2)$.
We obtain the following isomorphisms:
\begin{multline}
 \Ker(\nbigt)\otimes\Tate^S(-1/2)
\stackrel{\varphi_3^{-1}}{\lrarr}
 \Cok(\nbigt)\otimes\Tate^S(1/2)
\stackrel{\Cok(\nbigs)}{\lrarr}
 \Cok(\nbigt^{\ast})\otimes\Tate^S(1/2) \\
\lrarr
 \Ker(\nbigt)^{\ast}\otimes\Tate^S(1/2)
\lrarr
 \bigl(
 \Ker(\nbigt)\otimes\Tate^S(-1/2)
\bigr)^{\ast}
\end{multline}
The third map is given in
Lemma \ref{lem;08.1.18.15}.
The composite is denoted by $\varphi_4$.
\begin{lem}
We have the following commutative diagram:
\[
 \begin{CD}
 \Ker(\nbigt)\otimes\Tate^S(-1/2)
 @>{\varphi_4}>>
\bigl(
 \Ker(\nbigt)\otimes\Tate^S(-1/2)
\bigr)^{\ast}\\
 @A{\varphi_1}AA @V{\varphi_1^{\ast}}VV \\
 i_{Y\dagger}i_Y^{\dagger}
 \nbigt @>{i_{Y\dagger}i_Y^{\dagger}\nbigs}>>
 i_{Y\dagger}i_Y^{\dagger}\nbigt^{\ast}
 \end{CD}
\]
Hence, $\varphi_4$ gives a polarization 
of $\Ker(\nbigt)\otimes\Tate^S(-1/2)$.
\end{lem}
\pf
The first component of $\varphi_4$
is the composite of the following:
\[
  j_{\ast}j^{\ast}\nbigm'/\nbigm'
 \stackrel{-\Lambda^{-1}}{\llarr}
 i_{Y\dagger}i_Y^{\dagger}\nbigm'
 \stackrel{i_{Y\dagger}i_Y^{\dagger}\varphi'}{\llarr}
 i_{Y\dagger}i_Y^{\dagger}\nbigm''
 \stackrel{\id}{\llarr}
  i_{Y\dagger}i_Y^{\dagger}\nbigm''
 \stackrel{-1}{\llarr}
   i_{Y\dagger}i_Y^{\dagger}\nbigm''
\]
Hence, the first component of
$\varphi_1^{\ast}\circ
 \varphi_4\circ
 \varphi_1$ is
\[
 \Lambda\circ
 \bigl(
 (-\Lambda)^{-1}
\circ i_{Y\dagger}i_Y^{\dagger}\varphi'
\circ(-1)
 \bigr)
\circ
 \id
=i_{Y\dagger}i_Y^{\dagger}\varphi'
\]
It is equal to the first component of
$i_{Y\dagger}i_Y^{\dagger}\nbigs$.
The equality for the second components
can be checked in a similar way.
\hfill\qed

\vspace{.1in}

Let $\nbigs_{\Ker(\nbigt)}:
 \Ker(\nbigt)\lrarr
 \Ker(\nbigt)^{\ast}\otimes\Tate^S(1)$ denote
the induced polarization of $\Ker(\nbigt)$.

\subsection{Relation with the Lefschetz morphism}
\label{subsection;08.1.19.5}

Let $(\nbigt,\nbigs)$ and $Y$ be as above.
Let $a_X$ denote the obvious map from $X$ to a point.
Let $c$ be the cohomology class $2\pi[Y]$.
When we take a $C^{\infty}$-form $\omega$ 
representing $c$,
we have the induced map
$\nbigl_c=(-L_{\omega},L_{\omega}):
 a^i_{X\dagger}\nbigt
\lrarr a^{i+2}_{X\dagger}\nbigt\otimes\Tate(1)$,
i.e.,
\[
\begin{CD}
 a^{-i}_{X\dagger}\nbigm'
 @<{-L_{\omega}}<<
 a^{-i-2}_{X\dagger}\nbigm',
\end{CD}
\quad\quad
\begin{CD}
 a^i_{X\dagger}\nbigm''
 @>{L_{\omega}}>>
 a^{i+2}_{X\dagger}\nbigm''
\end{CD}
\]
Here, $L_{\omega}$ is given by the multiplication of
$\lambda^{-1}\omega$.
The map $\nbigl_c$ is called the Lefschetz morphism.
(See \cite{sabbah2}.
 See also Subsection \ref{subsection;08.12.16.3}.)
\index{Lefschetz morphism}

\vspace{.1in}

From the exact sequence
$0\lrarr \Ker(\nbigt)\lrarr 
 j_!j^{\ast}\nbigt\lrarr \nbigt\lrarr 0$
of $\nbigr_X$-triples,
we obtain the following morphism:
\[
 i_{Y,\Gys}^{\ast}:
 a_{X\dagger}^i(\nbigt)
\lrarr
 a_{X\dagger}^{i+1}\bigl(\Ker(\nbigt)\bigr)
\]
\index{morphism $i_{Y,\Gys}^{\ast}$}
From the exact sequence
$0\lrarr \nbigt\lrarr j_{\ast}j^{\ast}\nbigt
 \lrarr\Cok(\nbigt)\lrarr 0$,
we obtain the following morphism:
\[
 i_{Y\ast}^{\Gys}:
 a_{X\dagger}^i\Cok(\nbigt)
\lrarr
 a_{X\dagger}^{i+1}(\nbigt)
\]
\index{morphism $i_{Y\ast}^{\Gys}$}
Recall that we have obtained
the isomorphism
\[
 \varphi_3^{-1}:
 \Ker(\nbigt)\otimes\Tate^S(-1/2)
\lrarr
 \Cok(\nbigt)\otimes\Tate^S(1/2).
\]

\begin{lem}
The following diagram is commutative:
\[
 \begin{CD}
 a^i_{X\dagger}(\nbigt)@>{\nbigl_c}>>
 a^{i+2}_{X\dagger}(\nbigt)\otimes\Tate^S(1)\\
 @V{i_{Y,\Gys}^{\ast}}VV @A{i_{Y\,\ast}^{\Gys}}AA \\
 a^{i+1}_{X\dagger}\bigl(\Ker(\nbigt)\bigr)
 @>{a_X^{i+1}(\varphi_3^{-1})}>>
 a_{X\dagger}^{i+1}\bigl(\Cok(\nbigt)\bigr)
 \otimes\Tate^S(1)
 \end{CD}
\]
\end{lem}
\pf
Let $n:=\dim X$.
Note that we recall the construction of
$a_{X\dagger}\nbigm'$
and $a_{X\dagger}\nbigm''$
in Subsection \ref{subsection;08.12.16.3}.
We use the notation there.
For a local defining function $t$ of $Y$,
the section $\lambda^{-1}[\deldel_t] dt$
of $i_{Y\dagger}i_Y^{\dagger}\nbigo_{\nbigx}
 \otimes\Omega^{\bullet}_{\nbigx}$
is independent of the choice of $t$.
Hence, it defines a global section of
$i_{Y\dagger}i_Y^{\dagger}\nbigo_{\nbigx}
 \otimes\Omega^{\bullet}_{\nbigx}$.
Similarly,
we have the well-defined global section
$\lamda^{-1} dt/t$
of $j_{\ast}j^{\ast}
 \nbigo_{\nbigx}/\nbigo_{\nbigx}
 \otimes\Omega^{\bullet}_{\nbigx}$.

Let $s$ be a $C^{\infty}$-section of
$\nbigm''\otimes\Omega_{\nbigx}^{\bullet}$
such that $d s=0$.
The induced map
$a^i_{X\dagger}\nbigm''\lrarr
 a^{i+1}_{X\dagger}\bigl(
 i_{\dagger}i^{\dagger}\nbigm''
 \bigr)$ 
in $i_{Y,\Gys}^{\ast}$
is given by the following correspondence:
\[
 \begin{array}{ccccccc}
 & & 
 j_!j^{\ast}\nbigm''\otimes\Omega_{\nbigx}^{\bullet}
 \lrarr&
 &
 \nbigm''\otimes\Omega_{\nbigx}^{\bullet} & &\\
 & & \darr & &\\
 i_{\dagger}i^{\dagger}
 \nbigm''\otimes\Omega_{\nbigx}^{\bullet +1}
 &\lrarr &
 j_!j^{\ast}\nbigm''
 \otimes\Omega_{\nbigx}^{\bullet+1}
 \end{array}
\]
{\footnotesize
\[
 \begin{array}{ccccc}
 & & (s,0)& \lrarr & s\\
 & & \darr & \\
 (-1)^n\lambda^{-1} [\deldel_t]\cdot dt\cdot s_{|Y}
 & \lrarr & 
 \bigl(0,(-1)^n\lambda^{-1}
 [\deldel_t]\cdot dt\cdot s_{|Y}\bigr)
 \end{array}
\]}
Here, $(-1)^n$ appears because of the shift of the degree.
By the isomorphism
$\Lambda:
 i_{\dagger}i^{\dagger}\nbigm''\simeq
 j_{\ast}j^{\ast}\nbigm''/\nbigm''$ in 
$\varphi_3^{-1}:
 \Ker(\nbigt)\simeq\Cok(\nbigt)\otimes\Tate(1)$,
we have the correspondence:
\[
 \lambda^{-1}[\deldel_t] \cdot dt \cdot s_{|Y}
\longleftrightarrow
\sqrt{-1}\lambda^{-1}dt/t\cdot  s_{|Y} 
\]
We take a hermitian metric $g$
of the line bundle $\nbigo(Y)$.
Let $\sigma:\nbigo\lrarr \nbigo(Y)$ be the canonical section.
We obtain the function
$\log|\sigma|_g^2$.
Recall that the $C^{\infty}$-form
$\sqrt{-1}\delbar\del\log|\sigma|_g^2$
represents the cohomology class $2\pi[Y]$.
Then, 
we obtain the image of 
$(-1)^n\sqrt{-1}\lambda^{-1}dt/t\cdot s_{|Y}$
via the map
$a^{i+1}_{X\dagger}
 \bigl(j_{\ast}j^{\ast}\nbigm''/\nbigm''\bigr)
\lrarr a^{i+2}_{X\dagger}\bigl(\nbigm''\bigr)$ 
in $i_{Y\ast}^{\Gys}$,
by the following correspondence:
\[
 \begin{array}{ccccccc}
 & & 
 j_{\ast}j^{\ast}\nbigm''
 \otimes\Omega_{\nbigx}^{\bullet +1}
 & \lrarr & 
 (j_{\ast}j^{\ast}\nbigm''/\nbigm'')
 \otimes\Omega_{\nbigx}^{\bullet +1}\\
 & & \darr & & \\
 \nbigm''\otimes\Omega_{\nbigx}^{\bullet +2}
 & \lrarr &
 j_{\ast}j^{\ast}\nbigm''
 \otimes\Omega_{\nbigx}^{\bullet +2}
 \end{array}
\]
{\footnotesize
\[
 \begin{array}{ccccc}
 & & (-1)^n\sqrt{-1}\lambda^{-1}
 \del\log|\sigma|_g^2\cdot s
 & \lrarr & (-1)^n \sqrt{-1}\lambda^{-1}dt/t \cdot s_{|Y} \\
 & & \darr \\
 \sqrt{-1}\lambda^{-1}
 \delbar\del\log|\sigma|_g^2\cdot s
 & \lrarr &
 \sqrt{-1}\lambda^{-1}
 \delbar\del\log|\sigma|_g^2\cdot s
 \end{array}
\]}
Thus, 
we can conclude that 
$L_{\omega}$ is equal to
the $\bigl(
 a^i_{X\dagger}\nbigm''
 \lrarr a^{i+2}_{X\dagger}\nbigm''
\bigr)$-component
of $i_{Y\ast}^{\Gys}\circ a_X^{i+1}(\varphi_3^{-1})
 \circ i_{Y,\Gys}^{\ast}$.
Similarly we can check that
$L_{-\omega}$ is equal to
the 
$\bigl(
 a^{-i-2}_{X\dagger}\nbigm'
 \lrarr a^{-i}_{X\dagger}\nbigm'\bigr)$-component
of $i_{Y\ast}^{\Gys}\circ 
 a_X^{i+1}(\varphi_3^{-1})
 \circ i_{Y,\Gys}^{\ast}$.
\hfill\qed

\begin{cor}
\label{cor;07.10.27.20}
We have the following commutative diagram:
\begin{equation}
 \label{eq;08.1.18.21}
{\footnotesize
  \begin{CD}
 a^{-i}_{X\dagger}(\nbigt)
 @>{\nbigl_c^{j-1}}>>
 a^{-i+2j-2}_{X\dagger}(\nbigt)\otimes\Tate^S(j-1)
 @>{\nbigl_c}>>
 a^{-i+2j}_{X\dagger}(\nbigt)\otimes\Tate^S(j)\\
 @V{i_{Y,\Gys}^{\ast}}VV 
 @V{i_{Y,\Gys}^{\ast}}VV
 @A{i_{Y\ast}^{\Gys}}AA \\
 a^{-i+1}_{X\dagger}\bigl(\Ker(\nbigt)\bigr)
 @>{\nbigl_c^{j-1}}>>
 a^{-i+2j-1}_{X\dagger}\bigl(\Ker(\nbigt)\bigr)
 \otimes\Tate^S(j-1)
 @>{\varphi_3^{-1}}>>
 a^{-i+2j-1}_{X\dagger}
 \bigl(\Cok(\nbigt)\bigr)\otimes\Tate^S(j)
 \end{CD}
}
\end{equation}
Here, $a_{X\dagger}^{-i+2j-1}(\varphi_3^{-1})$
is denoted just by $\varphi_3^{-1}$.
\hfill\qed
\end{cor}

\subsection{A lemma}

Let us consider the case in which
$X$ is a projective variety
and $Y$ is a smooth ample hypersurface.
We recall the following lemma.
\begin{lem}
\label{lem;08.1.19.1}
Let $\nbigm$ be an $\nbigr_X$-module
on $\cnum_{\lambda}\times X$,
which is strict, coherent and holonomic.
Then, we have the vanishing
 $a^i_{X\dagger}\bigl(
 j_{\ast}j^{\ast}\nbigm\bigr)=0$
 for $i>0$.

Moreover, the following holds:
 \begin{itemize}
 \item If $i\geq 1$,
 $a^i_{X\dagger}\bigl(
 j_{\ast}j^{\ast}\nbigm/\nbigm
 \bigr)
 \simeq
 a^{i+1}_{X\dagger}\nbigm$.
\item If $i=0$,
 $a^0_{X\dagger}\bigl(
 j_{\ast}j^{\ast}\nbigm/\nbigm
 \bigr)
 \lrarr
 a^1_{X\dagger}\nbigm$ is surjective.
 \end{itemize}
\end{lem}
\pf
Let us show the first claim.
We shall argue such a vanishing
on a small compact
neighbourhood $\nbigk$ of
$\lambda_0\in\cnum_{\lambda}^{\ast}$.
We put $\nbigxzero:=\nbigk\times X$.
For an $\nbigo_{\nbigxzero}$-coherent sheaf 
$\nbigf$, we have
$a^i_{X\ast}\bigl(
 \nbigf\otimes\nbigo_{\nbigx}(\ast \nbigy)
\bigr)=0$ for $i>0$,
because $\nbigo_{\nbigx}(\nbigy)$ is
relatively ample.
If $\nbigk$ is sufficiently small,
$\nbigm_{\nbigxzero}$ has a finite filtration
of $\nbigr$-modules
such that the associated graded $\nbigr$-module
is the limit of $\nbigo_{\nbigxzero}$-coherent
submodules.
Hence, we obtain the vanishing
$a^i_{X\ast}
 \bigl(j_{\ast}j^{\ast}\nbigm
 \otimes p_{\lambda}^{\ast}
 \Omega^{r,0}_X\bigr)=0$
for $i>0$ and any $r\geq 0$,
where $\Omega^{r,0}$ denote
the sheaf of holomorphic $r$-forms on $X$.
We can deduce the first claim
from the vanishing.
The other claim follows from
the exact sequence
$0\lrarr \nbigm\lrarr j_{\ast}j^{\ast}\nbigm
 \lrarr j_{\ast}j^{\ast}\nbigm/\nbigm\lrarr 0$.
\hfill\qed

\chapter{Hard Lefschetz Theorem}
\label{section;08.10.18.121}
We show Hard Lefschetz theorem for
polarized wild pure twistor $D$-modules,
which is essentially due to
M. Saito and C. Sabbah.
This chapter is included for
rather expository purpose.

\section{Statement}
Let $X$ and $Y$ be complex manifolds,
and let $f:X\lrarr Y$ be a projective morphism.
Let $\nbiga$ be a $\rnum$-vector subspace of $\cnum$.
Let $(\nbigt,\nbigs)$ be
a polarized $\nbiga$-wild pure twistor $D$-module 
of weight $w$ on $X$.
We have the graded $\nbigr$-triple
$\bigoplus f^i_{\dagger}\nbigt$ on $Y$.
We have the Lefschetz morphism
$\nbigl_c:f^i_{\dagger}\nbigt\lrarr
 f^{i+2}_{\dagger}\nbigt\otimes\Tate^S(1)$
for the Chern class $c$ of a relatively ample line bundle.
We also have the induced Hermitian sesqui-linear duality
$\bigoplus f^i_{\dagger}\nbigs$
of $\bigoplus f^i_{\dagger}\nbigt$
with weight $w$.
(See \cite{sabbah2}.)
Let us show the following theorem
in this chapter.

\begin{thm}
\label{thm;07.10.23.20}
$\bigl(
 \bigoplus f^i_{\dagger}\nbigt,
 \nbigl_c,\bigoplus f^i_{\dagger}\nbigs
\bigr)$
is a polarized $\nbiga$-wild 
Lefschetz twistor $D$-module
with weight $w$
on $Y$.
\end{thm}

Let $\nbigm$ be the underlying $\nbigr$-module of $\nbigt$.
We obtain the $D$-module
by taking the specialization at $\nbigxlambda$
for $\lambda\neq 0$,
which is denoted by $\nbigm^{\lambda}$.

\begin{cor}
Let $(\nbigt,\nbigs)$ be a polarized wild
pure twistor $D$-module on $X$.
The Hard Lefschetz Theorem holds for 
the $D$-module $\nbigm^{\lambda}$ $(\lambda\neq 0)$,
i.e.,
$L_c^i:f^{-i}_{\dagger}(\nbigm^{\lambda})
\lrarr f^i_{\dagger}(\nbigm^{\lambda})$
is isomorphic.
As a result,
we have an isomorphism
$Rf_{\dagger}(\nbigm^{\lambda})
\simeq
 \bigoplus R^if_{\dagger}(\nbigm^{\lambda})$
in the derived category of 
cohomologically holonomic complexes on $Y$.
(See {\rm \cite{d5}}.)
\hfill\qed
\end{cor}

Hard Lefschetz Theorem was proved
by Beilinson-Bernstein-Deligne-Gabber
\cite{bbd}
for regular holonomic $D$-modules
of geometric origin
by the method of mod $p$-reductions.
Then, it was vastly generalized in
the original work due to M. Saito 
for polarized pure Hodge modules \cite{saito1}.
Using the argument of M. Saito,
C. Sabbah proved it for 
regular polarized pure twistor $D$-modules
\cite{sabbah2}.
We closely follow the argument of 
Saito and Sabbah.
(But, we follow Saito's argument more closely
 in some part.)

\subsection{Plan of the proof}

We use the induction on the dimensions
of the strict supports $\Supp(\nbigt)$
and $f\bigl(\Supp(\nbigt)\bigr)$.
Let us consider the following claim
for any $n\geq m$:
\begin{description}
\item[$P(n,m)$]
The claim of Theorem \ref{thm;07.10.23.20} holds,
in the case that $\dim\Supp(\nbigt)\leq n$
and $\dim f(\Supp(\nbigt))\leq m$ hold.
\end{description}
The claim $P(0,0)$ is obvious.
The proof is divided into the following three steps.
\begin{description}
\item[Step 1]
Prove the Hard Lefschetz Theorem
for a polarized wild pure twistor $D$-module
in the case that $X$ is a smooth projective curve 
and $Y$ is a point
(Proposition \ref{prop;07.10.22.152}).
\item[Step 2]
Give the argument for 
$P(n-1,m-1)\Longrightarrow P(n,m)$.
\item[Step 3]
Give the argument for
$P(n-1,0)\Longrightarrow P(n,0)$.
\end{description}

We will show Step 1
in Section \ref{subsection;07.10.23.30}.
A rather detailed argument for Step 2 is 
explained in Chapter 14.6 of \cite{mochi2}
in the regular case.
Hence, we will indicate how to modify the statements 
in Section \ref{subsection;07.10.23.22}.
We will give a rather detailed argument for Step 3
in Section \ref{subsection;07.10.23.21}
by following M. Saito,
although it is just a translation of his argument
to the case of wild pure twistor $D$-modules.
We also rely on some results of Sabbah in \cite{sabbah2}.
In the following argument,
we omit to distinguish $\nbiga$.

\section{Step 1}
\label{subsection;07.10.23.30}
\subsection{Statement}

Let $X$ be a smooth projective curve,
and let $D$ be a finite subset of $X$.
Let $a_X$ denote the obvious morphism of $X$
to a point.
Let $\omega$ be a Kahler form of $X$.

\begin{prop}
\label{prop;07.10.22.152}
Let $(\nbigt,\nbigs)$ be 
a polarized wild pure twistor $D$-module
on $X$ of weight $w$.
Then, the push-forward
$\bigl(\bigoplus a^i_{X\dagger}\nbigt, \nbigl_{\omega},
 a_{X\dagger}\nbigs\bigr)$
is a polarized graded Lefschetz twistor structure
of weight $w$.
\end{prop}

We have only to consider the case
in which the strict support of $\nbigt$ is $X$.

\subsection{The $L^2$-cohomology
for $(\nbigp_{\ast}\nbigelambda,\DDlambda,h)$}

Let $X$ be a smooth projective curve.
Let $D$ be a finite subset of $X$.
We take a Kahler metric $g_{X-D}$ of $X-D$
which is Poincar\'e like around any point of $D$.
Let $\omega_{1}$ denote the associated 
Kahler form.
We assume 
$\int_{X-D}\omega_{1}=\int_X\omega$.

Let $\harmonicbundle$ be
a wild harmonic bundle on $X-D$.
We have the associated filtered $\lambda$-flat
bundle $(\nbigp_{\ast}\nbigelambda,\DD)$
for $\lambda\in\cnum$.
Recall that the norm estimate holds for
$(\nbigp_{\ast}\nbigelambda,\DDlambda,h)$
(Proposition \ref{prop;07.10.6.40}).
Because of the results
in Section \ref{subsection;07.10.20.1},
we have the following quasi-isomorphisms
of the complexes of sheaves:
\begin{equation}
 \label{eq;07.7.22.50}
 \nbigs\bigl(\nbigp_{\ast}\nbigelambda
 \otimes\Omega^{\bullet,0}_X\bigr)
\lrarr
 \nbigl^{\bullet}_{\poly}\bigl(
 \nbigp_{\ast}\nbigelambda,\DDlambda\bigr)
\lrarr
\nbigl^{\bullet}\bigl(
 \nbigp_{\ast}\nbigelambda,\DDlambda\bigr)
\end{equation}

In particular,
we obtain the following corollary.
\begin{cor}
\label{cor;08.1.19.113}
The cohomology groups 
$H^{\bullet}\Bigl(
 \Gamma\bigl(
 \nbigl^{\bullet}\bigl(\nbigp_{\ast}\nbigelambda,
 \DDlambda\bigr)
 \bigr)
 \Bigr)$
associated to
$\nbigl^{\bullet}\bigl(\nbigp_{\ast}\nbigelambda,
 \DDlambda\bigr)$
is finite dimensional,
where $\Gamma$ denotes the functor taking
global sections.
\end{cor}
\pf
The hyper-cohomology group
associated to 
$ \nbigs\bigl(\nbigp_{\ast}\nbigelambda
 \otimes\Omega^{\bullet,0}_X\bigr)$
is finite dimensional.
Hence, we obtain Corollary \ref{cor;08.1.19.113}
by the quasi isomorphisms
in (\ref{eq;07.7.22.50}).
\hfill\qed

\vspace{.1in}

Let $\Harm^i$ denote the space 
introduced in Section \ref{subsection;08.1.20.13}.

\begin{prop}
\label{prop;08.1.19.112}
The induced map
$\Harm^i\lrarr
 H^i\bigl(
 \Gamma\bigl(\nbigl^{\bullet}(
 \nbigp_{\ast}\nbigelambda,
 \DDlambda)\bigr)
 \bigr)$ is an isomorphism.
\end{prop}
\pf
We have only to apply
Lemma \ref{lem;08.1.19.111}
with Corollary \ref{cor;08.1.19.113}.
\hfill\qed

\vspace{.1in}

We have the factorization
$\Harm^i\subset \Gamma\bigl(
 \nbigl^i_{\poly}(\nbigp_{\ast}\nbigelambda,\DDlambda)
 \bigr)
 \subset
 \Gamma\bigl(
 \nbigl^i(\nbigp_{\ast}\nbigelambda,\DDlambda) \bigr)$.
Hence, we obtain the following factorization:
\[
 \Harm^i\simeq
 H^i\Bigl(
 \Gamma\bigl(
 \nbigl^{\bullet}_{\poly}(\nbigp_{\ast}\nbigelambda,
 \DDlambda)
 \bigr)\Bigr)
 \simeq
 H^i\Bigl(
 \Gamma\bigl(
 \nbigl^{\bullet}(\nbigp_{\ast}\nbigelambda,\DDlambda)
 \bigr)
 \Bigr)
\]

In the case $\lambda\neq 0$,
we also have the factorization
$\Harm^{\bullet}\subset
 \Gamma\bigl(
 \nbiglbar_{\poly}^{\bullet}
 (\nbigp_{\ast}\nbigelambda,\DDlambda)
 \bigr)
 \subset
 \Gamma\bigl(
 \nbigl^{\bullet}_{\poly}
 (\nbigp_{\ast}\nbigelambda,\DDlambda)
 \bigr)
 \bigr) $,
according to Proposition \ref{prop;07.7.19.15}
and Corollary \ref{cor;07.7.7.43}.
It induces the following factorization:
\begin{equation}
 \label{eq;08.1.27.30}
 \Harm^i\simeq
 H^i\Bigl(
 \Gamma\bigl(
 \nbiglbar_{\poly}^{\bullet}(
   \nbigp_{\ast}\nbigelambda,\DDlambda)
 \bigr)
 \Bigr)
\simeq
 H^i\Bigl(
 \Gamma\bigl(
 \nbigl^{\bullet}_{\poly}(\nbigp_{\ast}\nbigelambda,
 \DDlambda)
 \bigr)\Bigr)
\end{equation}

We have the multiplication 
$L_{\omega_{1}}$ of $\omega_{1}$
on $\Harm^{\bullet}$
and $H^i\Bigl(
 \Gamma\bigl(
 \nbiglbar_{\poly}^{\bullet}(
   \nbigp_{\ast}\nbigelambda,\DDlambda)
 \bigr)
 \Bigr)$,
which is compatible with (\ref{eq;08.1.27.30}).
We also have the multiplication $L_{\omega}$ 
of $\omega$ on 
$H^i\Bigl(
 \Gamma\bigl(
 \nbiglbar_{\poly}^{\bullet}(
   \nbigp_{\ast}\nbigelambda,\DDlambda)
 \bigr)
 \Bigr)$.
\begin{lem}
\label{lem;08.1.27.31}
The isomorphism {\rm(\ref{eq;08.1.27.30})}
is compatible with 
$L_{\omega_{1}}$ and $L_{\omega}$.
\end{lem}
\pf
There exists $\tau$ such that
(i) $\omega-\omega_{1}=d\tau$,
(ii) $\tau$ is bounded with respect to $g_{X-D}$.
Hence, it is easy to show
$L_{\omega}=L_{\omega_{1}}$
on $H^i\Bigl(
 \Gamma\bigl(
 \nbiglbar_{\poly}^{\bullet}(
   \nbigp_{\ast}\nbigelambda,\DDlambda)
 \bigr)
 \Bigr)$.
Then, the claim of the lemma follows.
\hfill\qed

\begin{rem}
Recall $\lambda$-connection is equivalent
to ordinary connection in the case $\lambda\neq 0$.
Let $\nbigp_{\min}\nbigelambda$ denote
the $D_X$-submodule of $\nbigp\nbigelambda$
generated by $\nbigp_{<1}\nbigelambda$
over $D_X$.
If $\lambda$ is generic,
it is the same as the standard minimal extension,
and we obtain an isomorphism
between the cohomology group
of the $D$-module
$\nbigp_{\min}\nbigelambda$
and the $L^2$-cohomology group
$H^{\bullet}\Bigl(
 \Gamma\bigl(
 \nbigl^{\bullet}\bigl(\nbigp_{\ast}\nbigelambda,
 \DDlambda\bigr)
 \bigr)
 \Bigr)$
from the quasi isomorphism 
{\rm (\ref{eq;07.7.22.50})}.
For this isomorphism,
we do not need 
the harmonicity of $h$
but the norm estimates.
\hfill\qed
\end{rem}

\subsection{Deformation caused by 
variation of the irregular values}
\label{subsection;10.5.23.40}

Recall
$\nbigq_{\ast}\nbige^0:=
 \nbigp_{\ast}\nbige^0$,
and hence we have the isomorphism
$\Phibar^0: \Harm^i\lrarr 
\hyperh^i\bigl(
 \nbigs(\nbigq_{\ast}\nbige^0
 \otimes\Omega^{\bullet,0}_X)
 \bigr)$.
In the case $\lambda\neq 0$,
let $\nbigq h$ be a hermitian metric of
$\nbigelambda$ on $X-D$
whose restriction to a neighbourhood of $D$
is as in Section \ref{subsection;07.9.11.11}
for the filtered $\lambda$-flat bundle
$(\nbigq\nbigelambda,\DDlambda)$.
We can construct the metric
$\nbigq h^{(1+|\lambda|^2)}$
as in Section \ref{subsection;07.10.20.2}
with $T=1+|\lambda|^2$.
Then, $\nbigq h^{(1+|\lambda|^2)}$
and $h$ are mutually bounded
(Lemma \ref{lem;08.1.19.30}).
Hence,
we obtain the following quasi isomorphisms
because of the results in
Section \ref{subsection;07.10.20.2}:
\begin{multline}
\label{eq;08.1.20.50}
 \nbiglbar_{\poly}\bigl(
 \nbigp_{\ast}\nbigelambda,\DDlambda\bigr)
\stackrel{\simeq}{\llarr}
  \nbiglbar_{\poly}\bigl(\nbigelambda,
 \nbigq h^{(1+|\lambda|^2)},h\bigr)
\stackrel{\simeq}{\lrarr}
  \nbiglbar_{\poly}\bigl(\nbigq_{\ast}\nbigelambda,
 \DDlambda\bigr) \\
\stackrel{\simeq}{\lrarr}
  \nbigl^{\bullet}_{\poly}\bigl(\nbigq_{\ast}\nbigelambda,
 \DDlambda\bigr)
\stackrel{\simeq}{\llarr}
  \nbigs\bigl(\nbigq_{\ast}\nbigelambda
 \otimes\Omega_X^{\bullet,0} \bigr)
\end{multline}
They induce the isomorphisms of
the associated cohomology groups.
Therefore, we obtain a natural isomorphism:
\begin{equation}
 \label{eq;07.7.22.51}
\Phibar^{\lambda}:
\Harm^i
\lrarr
 \hyperh^i\bigl(
 \nbigs\bigl(\nbigq_{\ast}\nbigelambda
 \otimes\Omega_X^{\bullet,0} \bigr)
 \bigr).
\end{equation}

\begin{lem}
\label{lem;08.1.27.32}
We have the compatibility
$\Phibar^{\lambda}\circ
 L_{\omega_{1}}
=L_{\omega}\circ \Phibar^{\lambda}$.
\end{lem}
\pf
The multiplication $L_{\omega}$
is compatible with
the quasi isomorphisms in (\ref{eq;08.1.20.50}).
Then, Lemma \ref{lem;08.1.27.32}
follows from Lemma \ref{lem;08.1.27.31}.
\hfill\qed

\vspace{.1in}

If $|\lambda|$ is sufficiently small,
due to
Proposition \ref{prop;07.7.19.15}
and Proposition \ref{prop;08.9.17.3},
we have the following inclusions 
which are the quasi isomorphisms:
\begin{equation}
 \label{eq;07.9.12.15}
 \Harm^{\bullet}\subset
 \Gamma\bigl(
 \nbiglbar^{\bullet}_{\poly}
 (\nbigelambda,
 \nbigq h^{(1+|\lambda|^2)},h)
 \bigr)
\subset
 \Gamma\bigl(
 \nbiglbar^{\bullet}_{\poly}
 (\nbigq_{\ast}\nbigelambda,\DDlambda)
\bigr)
\subset
 \Gamma\bigl(
\nbigl^{\bullet}_{\poly}
 (\nbigq_{\ast}\nbigelambda,\DDlambda)
 \bigr)
\end{equation}
It is easy to see that
the composite of (\ref{eq;07.9.12.15}) 
is equal to 
 the isomorphism $\Phibar^{\lambda}$,
if $|\lambda|$ is sufficiently small.

\subsection{Quasi-isomorphisms for local families}

Let us consider the family case.
Let $U(\lambda_0)$ denote a small neighbourhood
of $\lambda_0$ in $\cnum_{\lambda}$.
We have the family of flat $\lambda$-connections
$(\nbigq^{(\lambda_0)}_{\ast}\nbige, \DD)$,
for which we take a hermitian metric $h_1$
as in Subsection \ref{subsection;10.5.23.1}.
We have the following quasi-isomorphisms
of the associated complexes of sheaves
(the left hand side is given only 
in the case $\lambda_0\neq 0$):
\begin{equation}
\begin{CD}
\nbiglbar_{\poly}^{\bullet}
 \bigl(\nbigq^{(\lambda_0)}_{\ast}\nbige,\DD\bigr)
@>>>
\nbigl_{\poly}^{\bullet}
 \bigl(\nbigq^{(\lambda_0)}_{\ast}\nbige,\DD\bigr)
@<<<
 \nbigs\bigl(\nbigq^{(\lambda_0)}_{\ast}\nbige
 \otimes\Omega^{\bullet,0}_X\bigr)
\end{CD}
\end{equation}

Let $p_X$ denote the projection
$U(\lambda_0)\times X\lrarr U(\lambda_0)$.
By considering the push-forward for $p_X$,
we obtain isomorphisms:
\[
R^ip_{X\ast}
\nbiglbar_{\poly}^{\bullet}
 \bigl(\nbigq^{(\lambda_0)}_{\ast}\nbige,\DD\bigr)
 \stackrel{\simeq}{\lrarr}
R^ip_{X\ast}
\nbigl_{\poly}^{\bullet}
 \bigl(\nbigq^{(\lambda_0)}_{\ast}\nbige,\DD\bigr)
\stackrel{\simeq}{\llarr}
R^ip_{X\ast}
 \nbigs\bigl(\nbigq^{(\lambda_0)}_{\ast}\nbige
 \otimes\Omega^{\bullet,0}_X\bigr)
\]
They are $\nbigo_{U(\lambda_0)}$-coherent.

Around $\lambda_0=0$,
we have the inclusion
$\Harm^{\bullet}\otimes\nbigo_{U(0)}
\lrarr
 p_{X\ast}\bigl(
 \nbigl_{\poly}^{\bullet}(
 \nbigq^{(0)}_{\ast}\nbige,\DD)
 \bigr)$
due to Proposition \ref{prop;07.7.19.15}
and Proposition \ref{prop;07.10.22.30}.
Therefore, we obtain a morphism
of the coherent $\nbigo_{U(0)}$-modules:
\begin{equation}
\Phibar^{(0)}:
 \bigoplus_i \Harm^{i}\otimes\nbigo_{U(0)}
\lrarr
 \bigoplus_iR^ip_{X\ast}\bigl(
  \nbigs\bigl(\nbigq^{(0)}_{\ast}\nbige
 \otimes\Omega^{\bullet,0}_X\bigr)
 \bigr)
\end{equation}

\begin{lem}
\mbox{{}}\label{lem;08.1.20.51}
\begin{itemize}
\item
$\Phibar^{(0)}$ is an isomorphism.
\item
The specialization of
$R^ip_{X\ast}\bigl(
  \nbigs\bigl(\nbigq^{(0)}_{\ast}\nbige
 \otimes\Omega^{\bullet,0}_X\bigr)
 \bigr)$ at $\lambda\in U(0)$
is naturally isomorphic to
$\hyperh^i\bigl(
 \nbigs\bigl(\nbigq_{\ast}\nbigelambda
 \otimes\Omega^{\bullet}_X\bigr)
 \bigr)$,
and the specialization of
$\Phibar^{(0)}$ is the same as
$\Phibar^{\lambda}$ under the isomorphism.
\end{itemize}
\end{lem}
\pf
For simplicity of the description,
we put $\nbigs^{(0)}:=
 \nbigs(\nbigq^{(0)}_{\ast}\nbige
 \otimes\Omega^{\bullet,0})$.
For any $\lambda_1\in U(0)$,
let $\nbigs^{(0)}_{|\lambda_1}$
denote the specialization of $\nbigs^{(0)}$
at $\{\lambda_1\}\times X$.
We remark that we have the natural 
quasi isomorphism
$\nbigs^{(0)}_{|\lambda_1}
\lrarr
 \nbigs\bigl(
 \nbigq_{\ast}\nbige^{\lambda_1}
 \otimes\Omega^{\bullet,0}
 \bigr)$.
By using the remark in the last paragraph
of Subsection \ref{subsection;10.5.23.40},
we obtain the following commutative diagram
of the $\nbigo_{U(0)}$-coherent sheaves:
\[
 \begin{CD}
 \Harm^i\otimes\nbigo_{U(0)} 
 @>{\lambda-\lambda_1}>>
 \Harm^{i}\otimes \nbigo_{U(0)}
 @>>>
 \Harm^i @>>> 0\\
 @V{\Phibar^{(0)}}VV 
 @V{\Phibar^{(0)}}VV 
 @V{\Phibar^{\lambda_1}}V{\simeq}V 
 @VVV \\
R^ip_{X\ast}\nbigs^{(0)}
 @>{\lambda-\lambda_1}>>
R^ip_{X\ast}\nbigs^{(0)}
 @>>>
R^ip_{X\ast}\nbigs^{(0)}_{|\lambda_1}
 @>>>
 R^{i+1}p_{X\ast}\nbigs^{(0)}
 \end{CD}
\]
Here, $\lambda-\lambda_1$ means the multiplication
of $\lambda-\lambda_1$.
Then, it is easy to show the claims of the lemma.
\hfill\qed

\vspace{.1in}
Let us see around $\lambda_0\neq 0$.
Let $T_1(\lambda):=1+|\lambda|^2$.
We have the metric $h_1^{(T_1(\lambda))}$
as constructed in Section \ref{subsection;07.12.2.25}.
\begin{lem}
\label{lem;08.1.20.41}
If we shrink $U(\lambda_0)$ appropriately,
the metrics $h_1^{(T_1(\lambda))}$ and $h$
are mutually bounded.
\end{lem}
\pf
Let $h_0$ be a hermitian metric
for the family of filtered $\lambda$-flat bundles
$(\nbigpzero_{\ast}\nbige,\DD)$
as in Subsection \ref{subsection;10.5.23.1}.
Let $S$ be a small sector in
$U(\lambda_0)\times (X-D)$.
We can take a $\DD$-flat splitting 
$\nbigp^{(\lambda_0)}_0\nbige_{\Sbar}
=\bigoplus_{\gminia,S}
 \nbigpzero_0\nbige_{\gminia,S}$
of the full Stokes filtration.
Let $p^{(j)}_{\gminia,S}$ denote 
the projection onto 
$\nbigpzero_0\nbige_{\gminia,S}$.
Let $F_S(w)$ be given as in 
(\ref{eq;08.1.20.30}):
\[
  F_S(w):=\exp\bigl(w\,\nbigb_S\bigr),
\quad
  \nbigb_S:=
 \sum_{\gminia\in\Irr(\theta)}
 \gminia\, p_{\gminia,S}
\]
According to Corollary \ref{cor;07.12.2.30},
the metrics
$F_S(-\lambdabar+\lambdabar_0)^{\ast}h_0$
and $h$ are mutually bounded on $S$.
By construction of 
$\nbigqzero\nbige$
(see Subsection \ref{subsection;08.9.28.41}),
the metrics
$F_S(\lambdabar_0)^{\ast}h_0$
and $h_1$ are mutually bounded on $S$.
Hence, 
$h$ and $F_S(-\lambdabar)^{\ast}h_1$
are mutually bounded.
It implies the claim of Lemma \ref{lem;08.1.20.41}.
\hfill\qed

\vspace{.1in}

Then,
we obtain the following inclusion,
due to Proposition \ref{prop;07.7.19.15}
and Proposition \ref{prop;08.9.17.3}:
\[
 \bigoplus_i \Harm^{i}\otimes\nbigo_{U(\lambda_0)}
\lrarr
 p_{X\ast}\Bigl(
 \nbiglbar_{\poly}^{\bullet}\bigl(
 \nbigq_{\ast}^{(\lambda_0)}\nbige,
 h_1^{(T_1(\lambda))}
 \bigr) \Bigr)
\]
We also have the following quasi-isomorphisms
of the complexes of sheaves
due to
the results in Section \ref{subsection;07.10.20.6}:
\begin{multline}
\label{eq;07.10.22.150}
 \nbiglbar_{\poly}^{\bullet}(
 \nbigq_{\ast}^{(\lambda_0)}\nbige,h_1^{(T_1(\lambda))})
\stackrel{\simeq}{\llarr}
 \nbiglbar_{\poly}^{\bullet}
 \bigl(\nbigq^{(\lambda_0)}\nbige,h_1,h_1^{(T_1(\lambda))}
  \bigr)
\stackrel{\simeq}{\lrarr}
 \nbiglbar_{\poly}^{\bullet}
 \bigl(\nbigq_{\ast}^{(\lambda_0)}\nbige, \DD\bigr)
 \\
 \stackrel{\simeq}{\lrarr}
 \nbigl_{\poly}^{\bullet}
 (\nbigq^{(\lambda_0)}_{\ast}\nbige,\DD)
 \stackrel{\simeq}{\llarr}
 \nbigs\bigl(\nbigq^{(\lambda_0)}_{\ast}\nbige
 \otimes\Omega^{\bullet,0}_X\bigr)
\end{multline}
Thus, we obtain a morphism of
coherent $\nbigo_{U(\lambda_0)}$-sheaves:
\[
 \Phibar^{(\lambda_0)}:
\bigoplus_i
 \Harm^{i}\otimes\nbigo_{U(\lambda_0)}
\lrarr
\bigoplus_i
 R^ip_{X\ast}\bigl(
 \nbigs\bigl(\nbigq^{(\lambda_0)}_{\ast}\nbige
 \otimes\Omega^{\bullet,0}_X\bigr)
 \bigr)
\]

\begin{lem}
\mbox{{}}\label{lem;08.1.31.4}
\begin{itemize}
\item
$\Phibar^{(\lambda_0)}$ is an isomorphism.
\item
The specialization of
$R^ip_{X\ast}\bigl(
 \nbigs\bigl(\nbigq^{(\lambda_0)}_{\ast}\nbige
 \otimes\Omega^{\bullet,0}_X\bigr)
 \bigr)$ at $\lambda\in U(\lambda_0)$
is naturally isomorphic to
$\hyperh^i\bigl(
 \nbigs\bigl(\nbigq_{\ast}\nbigelambda
 \otimes\Omega^{\bullet,0}\bigr)
 \bigr)$,
and the specialization of
$\Phibar^{(\lambda_0)}$ at $\lambda$
is the same as $\Phibar^{\lambda}$
under the isomorphism.
\end{itemize}
\end{lem}
\pf
The claims can be shown by 
the arguments in the proof of
Lemma \ref{lem;08.1.20.51}.
For simplicity of the description,
we put $\nbigs^{(\lambda_0)}:=
 \nbigs(\nbigq^{(\lambda_0)}_{\ast}\nbige
 \otimes\Omega^{\bullet,0})$.
For any $\lambda_1\in U(\lambda_0)$,
let $\nbigs^{(\lambda_0)}_{|\lambda_1}$
denote the specialization of $\nbigs^{(\lambda_0)}$
at $\{\lambda_1\}\times X$.
We remark that we have the natural 
quasi isomorphism
$\nbigs^{(\lambda_0)}_{|\lambda_1}
\lrarr
 \nbigs\bigl(
 \nbigq_{\ast}\nbige^{\lambda_1}
 \otimes\Omega^{\bullet,0}
 \bigr)$.
By using the commutative diagrams
(\ref{eq;08.1.31.2}) and (\ref{eq;08.1.31.3}),
we obtain the following commutative diagram
of the $\nbigo_{U(\lambda_0)}$-coherent sheaves:
\[
 \begin{CD}
 \Harm^i\otimes\nbigo_{U(\lambda_0)} 
 @>{\lambda-\lambda_1}>>
 \Harm^{i}\otimes \nbigo_{U(\lambda_0)}
 @>>>
 \Harm^i @>>> 0\\
 @V{\Phibar^{(\lambda_0)}}VV 
 @V{\Phibar^{(\lambda_0)}}VV 
 @V{\Phibar^{\lambda_1}}V{\simeq}V 
 @VVV \\
R^ip_{X\ast}\nbigs^{(\lambda_0)}
 @>{\lambda-\lambda_1}>>
R^ip_{X\ast}\nbigs^{(\lambda_0)}
 @>>>
R^ip_{X\ast}\nbigs^{(\lambda_0)}_{|\lambda_1}
 @>>>
 R^{i+1}p_{X\ast}\nbigs^{(\lambda_0)}
 \end{CD}
\]
Here, $\lambda-\lambda_1$ means the multiplication
of $\lambda-\lambda_1$.
Then, it is easy to show the both claims of the lemma.
\hfill\qed

\subsection{Global isomorphism}

Let $\Omega^{\bullet,0}_X$ denote 
the holomorphic de Rham complex on $X$.
We have the complex 
$\gbige\otimes\Omega^{\bullet,0}_X$
induced by $\DD$ and $\lambda\del_X$.
We would like to compare
the hyper-cohomology of
$\gbige\otimes\Omega^{\bullet,0}_X$
and the harmonic forms.
Let $U(\lambda_0)$ denote a small neighbourhood
of $\lambda_0$ in $\cnum_{\lambda}$.
Let $\gbige^{(\lambda_0)}$ denote 
the restriction of $\gbige$ to $U(\lambda_0)\times X$.

\begin{lem}
The naturally defined morphism
$\nbigs(\nbigq^{(\lambda_0)}\nbige\otimes
 \Omega^{\bullet,0}_X)
\lrarr \gbige^{(\lambda_0)}\otimes\Omega^{\bullet,0}_X$
is a quasi-isomorphism.
Therefore, we obtain an isomorphism:
\[
 \Phi^{(\lambda_0)}:
 \bigoplus_i\Harm^i\otimes\nbigo_{U(\lambda_0)}
\lrarr
 \bigoplus_i R^ip_{X\ast}
 \bigl(\gbige^{(\lambda_0)}\otimes
 \Omega^{\bullet,0}_X\bigr)
\]
\end{lem}
\pf
We have only to show the acyclicity
of the quotient complex
of $\nbigs(\nbigq^{(\lambda_0)}\nbige\otimes
 \Omega^{\bullet,0}_X)
\lrarr \gbige^{(\lambda_0)}\otimes\Omega^{\bullet,0}_X$.
For that purpose,
we may replace them with their completions
at $U(\lambda_0)\times D$.
For the regular part
which comes from a tame harmonic bundle,
the issue was studied in 
\cite{sabbah2} and \cite{mochi2}.
For the irregular parts,
both of them are acyclic.
\hfill\qed

\vspace{.1in}

For $U(\lambda_1)\subset U(\lambda_0)$,
we have 
$\gbige^{(\lambda_1)}
=\gbige^{(\lambda_0)}_{|U(\lambda_1)\times X}$.
According to Lemmas \ref{lem;08.1.20.51}
and \ref{lem;08.1.31.4},
we obtain $\Phi^{(\lambda_0)}_{|U(\lambda_1)}
=\Phi^{(\lambda_1)}$.
Hence, we can glue 
$\{\Phi^{(\lambda_0)}\,|\,\lambda_0\in\cnum\}$,
and we obtain the global isomorphism:
\begin{equation}
\label{eq;07.9.12.30}
 \Phi:
 \bigoplus_i\Harm^{i}\otimes\nbigo_{\cnum_{\lambda}}
\lrarr 
\bigoplus
 R^ip_{X\ast}\bigl(
 \gbige \otimes\Omega^{\bullet,0}_X
 \bigr)
\end{equation}

\begin{lem}
\label{lem;08.1.27.50}
We have the compatibility
$\Phi\circ L_{\omega_{1}}
=L_{\omega}\circ \Phi$.
\end{lem}
\pf
Since the specialization of
$R^ip_{X\ast}\bigl(
 \gbige \otimes\Omega^{\bullet,0}_X
 \bigr)$ at $\lambda$ is 
naturally identified with
$\hyperh^i\bigl(
 \nbigs(
 \nbigq_{\ast}\nbigelambda\otimes\Omega^{\bullet,0})
 \bigr)$,
the claim follows from Lemma \ref{lem;08.1.27.32}.
\hfill\qed

\subsection{Twist}

To consider the push-forward of $\nbigr$-triples,
we twist the de Rham complex
as in \cite{sabbah2}
which we refer to for more details and precision.
Let $\Omega_{\nbigx}^{1,0}:=
 \lambda^{-1}\,
  p_{\lambda}^{\ast}\Omega_{X}^{1,0}
\subset p_{\lambda}^{\ast}\Omega_X^{1,0}
 \bigl(\ast (X\times 0)\bigr)$,
and $\Omega_{\nbigx}^{p,0}:=
 \bigwedge^p
 \Omega_{\nbigx}^{1,0}$.
The derivation $\del_X$ for $\Omega^{\bullet,0}_X$
induces the derivation of $\Omega^{\bullet,0}_{\nbigx}$.
We have the natural isomorphism
$p_{\lambda}^{\ast}\Omega^{\bullet,0}_X
\simeq
 \Omega^{\bullet,0}_{\nbigx}$,
via which $\del_X$ is identified with 
$\lambda\, \del_X$
on $p_{\lambda}^{\ast}\Omega^{\bullet}_X$.

We have the derivation on
$\gbige\otimes\Omega^{\bullet,0}_{\nbigx}$
induced by $\DD^f$ and $\del_X$,
which is also denoted by $\DD^f$.
We have the natural isomorphism
$(\gbige\otimes\Omega^{\bullet,0}_{\nbigx},\DD^f)
\simeq
 (\gbige\otimes \Omega^{\bullet,0}_X,\DD)$
induced by the natural isomorphism
$\Omega^{1,0}_{\nbigx}
\simeq
 p_{\lambda}^{\ast}\Omega^{1,0}_X$
given by the multiplication of $\lambda$.
Shifting the degree,
we consider the complex
$(\gbige\otimes\Omega^{1+\bullet,0}_{\nbigx},\DD^f)$.
We have a similar twist
$\nbigs(
 \nbigq^{(\lambda_0)}\nbige\otimes
 \Omega^{1+\bullet,0}_{\nbigx})$.

We put $\Omega^{p,q}_{\nbigx}:=
 \Omega^{p,0}_{\nbigx}\otimes
 p_{\lambda}^{-1}\Omega^{0,q}_X$.
The derivation
$-\del_X-\delbar_X$ induces 
the derivation on $\Omega^{1+\bullet,\bullet}_{\nbigx}$.
By the natural isomorphism
$\Omega^{1+\bullet,\bullet}_{\nbigx}
\simeq
 p_{\lambda}^{-1}\Omega^{\bullet,\bullet}_X$,
$-\del_X-\delbar_X$
corresponds to $-\lambda\del_X-\delbar_X$.
By the natural isomorphism,
we obtain the complex of sheaves
$\nbigk^{\bullet}_{\poly}
 \bigl(\nbigq^{(\lambda_0)}_{\ast}\nbige,\DD\bigr)$
corresponding to
$\nbigl^{\bullet}_{\poly}
 \bigl(\nbigq^{(\lambda_0)}_{\ast}\nbige,\DD\bigr)$.
Then, we have the following quasi-isomorphisms:
\begin{equation}
 \label{eq;08.1.31.10}
\begin{CD}
 \nbigk^{\bullet}_{\poly}
 \bigl(\nbigq^{(\lambda_0)}_{\ast}\nbige,\DD\bigr)
 @<<<
 \nbigs(\nbigq^{(\lambda_0)}\nbige
\otimes\Omega^{1+\bullet,0}_{\nbigx})
@>>>
 \gbige^{(\lambda_0)}\otimes
 \Omega^{1+\bullet,0}_{\nbigx}
\end{CD}
\end{equation}

If $\lambda_0=0$,
we have the induced inclusion:
\[
 \Harm^{p}\otimes\nbigo_{U(0)}
\lrarr
 \Gamma\Bigl(\nbigl^{p}_{\poly}
 \bigl(\nbigq^{(0)}_{\ast}\nbige,\DD\bigr)
 \Bigr)
\simeq
 \Gamma\Bigl(\nbigk^{p-1}_{\poly}
 \bigl(\nbigq^{(0)}_{\ast}\nbige,\DD\bigr)\Bigr)
\]
Thus, we obtain the isomorphism
$\Phitilde^{(0)}:
 \Harm^{1+\bullet}\simeq
 Rp_{X\ast}\bigl(\gbige\otimes
 \Omega^{1+\bullet,0}_{\nbigx}\bigr)$.

If $\lambda_0\neq 0$,
we have the corresponding twist
for the diagram (\ref{eq;07.10.22.150}):
\begin{multline}
 \nbigkbar_{\poly}^{\bullet}(
 \nbigq_{\ast}^{(\lambda_0)}\nbige,h_1^{(T_1(\lambda))})
\stackrel{\simeq}{\llarr}
 \nbigkbar_{\poly}^{\bullet}
 \bigl(\nbigq^{(\lambda_0)}\nbige,h_1,h_1^{(T_1(\lambda))}
  \bigr)
\stackrel{\simeq}{\lrarr}
 \nbigkbar_{\poly}^{\bullet}
 \bigl(\nbigq_{\ast}^{(\lambda_0)}\nbige,
 \DD\bigr)
 \\
 \stackrel{\simeq}{\lrarr}
 \nbigk_{\poly}^{\bullet}
 (\nbigq^{(\lambda_0)}_{\ast}\nbige,\DD)
 \stackrel{\simeq}{\llarr}
 \nbigs\bigl(\nbigq^{(\lambda_0)}_{\ast}\nbige
 \otimes\Omega^{1+\bullet,0}_{\nbigx}\bigr)
\end{multline}
We also have the induced inclusion:
\[
 \Harm^{p}\otimes\nbigo_{U(\lambda_0)}
\lrarr
 \Gamma\Bigl(
  \nbigkbar_{\poly}^{p-1}(
 \nbigq_{\ast}^{(\lambda_0)}\nbige,h_1^{(T_1(\lambda))})
 \Bigr)
\]
We obtain induced isomorphisms
$ \Phitilde^{(\lambda_0)}:
 \Harm^{1+\bullet}\otimes\nbigo_{U(\lambda_0)}
\lrarr
 Rp_{X\ast}\bigl(\gbige^{(\lambda_0)}\otimes
 \Omega^{1+\bullet,0}_{\nbigx}\bigr)$.
By gluing them,
we obtain
$\Phitilde:
 \Harm^{1+\bullet}\otimes\nbigo_{\cnum_{\lambda}}
\lrarr
 Rp_{X\ast}\bigl(\gbige\otimes
 \Omega^{1+\bullet,0}_{\nbigx}\bigr)$.

We put $\nbigh^p:=
 \Harm^{p+1}\otimes\nbigo_{\cnum_{\lambda}}$.
We obtain an induced isomorphism
\[
 \Phitilde:\nbigh^{\bullet}\lrarr  
 a^{\bullet}_{X\dagger}\gbige.
\]

\subsection{The graded Lefschetz twistor structure
 on the space of harmonic forms}

Let $C^0(\vecS,\cnum)$ denote the sheaf of
continuous functions on 
$\vecS=\bigl\{|\lambda|=1\bigr\}$.
Recall that 
the natural multiplication and the integration
induces the following sesqui-linear pairing:
\[
C_{\gbigh}^p:\nbigh^{-p}_{|\vecS}
\otimes
 \sigma^{\ast}
 \bigl(\nbigh^{p}\bigr)_{|\vecS}
\lrarr
 C^0(\vecS,\cnum)
\]
Thus, we obtain the $\nbigr$-triples
$\gbigh^p:=
 (\nbigh^{-p},\nbigh^p,C^p_{\gbigh})$
($p=-1,0,1$).
We put $\gbigh^{\bullet}:=\bigoplus\gbigh^p$.
We have the Lefschetz morphism
$\nbigl_{\omega_{1}}=
 (-L_{\omega_{1}},L_{\omega_{1}}):
 \gbigh^{-1}\lrarr \gbigh^1\otimes\Tate^S(1)$.
We regard it as the morphism
$\gbigh^{\bullet}\lrarr\gbigh^{\bullet}\otimes\Tate^S(1)$.
We also have the Hermitian sesqui-linear duality
$\nbigs_{\gbigh}:\gbigh^{\bullet}\lrarr (\gbigh^{\bullet})^{\ast}$.
The following lemma can be shown
using the argument for Hodge-Simpson theorem {\rm 2.2.4}
in \cite{sabbah2}.
\begin{lem}
\label{lem;07.10.23.10}
$(\gbigh^{\bullet},\nbigl_{\omega_1},\nbigs_{\gbigh})$
is a polarized graded Lefschetz twistor structure
of weight $0$.
\hfill\qed
\end{lem}

\subsection{Compatibility
 with the sesqui-linear pairing}

Let $C_1^p$ denote the induced sesqui-linear pairing:
\[
C_1^p:
 R^{-p}p_{X\ast}(\gbige\otimes
 \Omega^{1+\bullet,0}_{\nbigx})_{|\vecS}
\otimes
\sigma^{\ast}\bigl(
 R^{p}p_{X\ast}(\gbige\otimes
 \Omega^{1+\bullet,0}_{\nbigx})
 \bigr)_{|\vecS}
\lrarr
 C^0(\vecS,\cnum)
\]
We have 
$a^p_{X\dagger}(\nbigt)
=\bigl(a^{-p}_{X\dagger}(\gbige),
 a^p_{X\dagger}(\gbige),
 C^p_1
 \bigr)$.

\begin{lem}
$C_1^p$ and $C_{\gbigh}^p$ are compatible with
the isomorphism $\Phitilde$.
Namely, the following diagram is commutative:
\begin{equation}
 \label{eq;07.10.22.151}
 \begin{CD}
\Harm^{1-p}\otimes\nbigo_{\cnum_{\lambda}|\vecS}
\otimes
 \sigma^{\ast}
 \bigl(\Harm^{1+p}\otimes
 \nbigo_{\cnum_{\lambda}}\bigr)_{|\vecS}
@>{C_{\gbigh}^p}>>
 C^0(\vecS,\cnum)\\
 @V{\Phitilde\otimes\sigma^{\ast}\Phitilde}V{\simeq}V 
 @V{=}VV \\
 R^{-p}p_{X\ast}(\gbige\otimes
 \Omega^{1+\bullet,0}_{\nbigx})_{|\vecS}
\otimes
\sigma^{\ast}\bigl(
 R^{p}p_{X\ast}(\gbige\otimes
 \Omega^{1+\bullet,0}_{\nbigx})
  \bigr)
 _{|\vecS}
@>{C_1^p}>>
 C^0(\vecS,\cnum)
 \end{CD}
\end{equation}
\end{lem}
\pf
Let $\lambda_0\in\cnum$ such that $|\lambda_0|=1$.
We have the following pairings:
\begin{multline*}
 \nbigkbar^{-p}_{\poly}
 (\nbigq^{(\lambda_0)}_{\ast}\nbige,h_1^{(T_1)})
 _{|\vecI(\lambda_0)\times X}
\otimes
 \sigma^{\ast}
 \bigl(
 \nbigkbar^{p}_{\poly}
 (\nbigq^{(-\lambda_0)}_{\ast}\nbige,h_1^{(T_1)})
 \bigr)_{|\vecI(\lambda_0)\times X}
\lrarr \\
 \distribution_{\vecI(\lambda_0)\times X/\vecI(\lambda_0)}
 \otimes\Omega^{1,1}_{\nbigx}
\end{multline*}
\begin{multline*}
 \nbigkbar^{-p}_{\poly}
 (\nbigq^{(\lambda_0)}_{\ast}\nbige,h_1,h_1^{(T_1)})
 _{|\vecI(\lambda_0)\times X}
\otimes
 \sigma^{\ast}
 \bigl(
 \nbigkbar^{p}_{\poly}
 (\nbigq^{(-\lambda_0)}_{\ast}\nbige,h_1,h_1^{(T_1)})
 \bigr)_{|\vecI(\lambda_0)\times X}
\lrarr \\
 \distribution_{\vecI(\lambda_0)\times X/\vecI(\lambda_0)}
 \otimes\Omega^{1,1}_{\nbigx}
\end{multline*}
\begin{multline*}
 \nbigkbar^{-p}_{\poly}
 (\nbigq^{(\lambda_0)}_{\ast}\nbige,\DD)
 _{|\vecI(\lambda_0)\times X}
\otimes
 \sigma^{\ast}
 \bigl(
 \nbigkbar^{p}_{\poly}
 (\nbigq^{(-\lambda_0)}_{\ast}\nbige,\DD)
 \bigr)_{|\vecI(\lambda_0)\times X}
\lrarr \\
 \distribution_{\vecI(\lambda_0)\times X/\vecI(\lambda_0)}
\otimes\Omega^{1,1}_{\nbigx}
\end{multline*}
\begin{multline*}
 \nbigk^{-p}_{\poly}
 (\nbigq^{(\lambda_0)}_{\ast}\nbige,\DD)
 _{|\vecI(\lambda_0)\times X}
\otimes
 \sigma^{\ast}
 \bigl(
 \nbigk^{p}_{\poly}
 (\nbigq^{(-\lambda_0)}_{\ast}\nbige,\DD)
 \bigr)_{|\vecI(\lambda_0)\times X}
\lrarr \\
 \distribution_{\vecI(\lambda_0)\times X/\vecI(\lambda_0)}
\otimes\Omega^{1,1}_{\nbigx}
\end{multline*}
\begin{multline*}
 \nbigs
 (\nbigq^{(\lambda_0)}_{\ast}\nbige\otimes
 \Omega^{1-p,0}_{\nbigx})
 _{|\vecI(\lambda_0)\times X}
\otimes
 \sigma^{\ast}
 \bigl(
 \nbigs
 \bigl(\nbigq^{(-\lambda_0)}_{\ast}\nbige\otimes
 \Omega^{1+p,0}_{\nbigx}\bigr)
 \bigr)_{|\vecI(\lambda_0)\times X}
\lrarr \\
 \distribution_{\vecI(\lambda_0)\times X/\vecI(\lambda_0)}
 \otimes\Omega^{1,1}_{\nbigx}
\end{multline*}
(The continuity of distributions
 follows from the remark in 
Subsection \ref{subsection;10.5.23.11}.)
They are compatible with the quasi-isomorphisms
in (\ref{eq;07.10.22.150}) and (\ref{eq;08.1.31.10}).
Hence, we obtain the commutativity of the diagram
(\ref{eq;07.10.22.151}).
\hfill\qed

\vspace{.1in}

Then,
we obtain the isomorphism
of $\nbigr$-triples
$\gbigh^{\bullet}
\simeq
 a^{\bullet}_{X\dagger}(\gbigt)$.
It is clearly compatible
with the induced Hermitian sesqui-linear dualities.
We also have the compatibility 
with the Lefschetz morphisms
by Lemma \ref{lem;08.1.27.50}.
Hence, Proposition \ref{prop;07.10.22.152}
follows from Lemma \ref{lem;07.10.23.10}.
\hfill\qed

\section{Step 2}
\label{subsection;07.10.23.22}
In this section,
the most proof is referred 
to our previous paper \cite{mochi2}.
However, we should emphasize
that it is essentially due to Saito and Sabbah,
as remarked in \cite{mochi2}.

\subsection{Preliminary I}
\label{subsection;05.1.15.12}

Let $f:X\lrarr Y$ be a projective morphism 
of complex manifolds.
Let $Z$ be a closed subvariety of $X$.
Let $c$ be the first Chern class of 
relatively ample line bundle on $X$.
In this subsection, the following assumption is imposed.
\begin{assumption}
\label{assumption;05.1.15.11}
The claim of Theorem {\rm \ref{thm;07.10.23.20}}
holds for a morphism $f$ and 
polarized wild pure twistor $D$-modules
whose strict supports are contained in $Z$.
\hfill\qed
\end{assumption}

Let us consider a tuple
$(\nbigt,M,N,\nbigs)$ on $X$
as follows:
\begin{condition}
\label{condition;05.1.18.20}\mbox{{}}
\begin{itemize}
\item
 An $\nbigr_X$-triple $\nbigt$ 
 with any increasing filtration $M$
 such that $\Gr^M_i\nbigt$ are 
 wild pure twistor $D$-modules
 of weight $w+i$.
 The support of $\nbigt$ is contained in $Z$.
\item
 A nilpotent map $N:\nbigt\lrarr\nbigt(-1)$
 whose weight filtration is equal to $M$.
\item
 A Hermitian sesqui-linear duality 
$\nbigs:\nbigt\lrarr \nbigt^{\ast}(-w)$.
 We assume that $N$ is skew adjoint 
 with respect to $\nbigs$,
 i.e., $N^{\ast}\circ\nbigs+\nbigs\circ N=0$.
\item
 $\nbigs\circ N^j$ gives a polarization of
 $P\Gr^M_j\nbigt$ for each $j\geq 0$.
\hfill\qed
\end{itemize}
\end{condition}

We put $\widehat{\nbigt}^i:=\Gr^W_{-i}\nbigt$,
and then we obtain the graded $\nbigr$-triple
$\widehat{\nbigt}=\bigoplus_j\widehat{\nbigt}^j$.
We have the naturally induced nilpotent map
$\widehat{N}:
 \widehat{\nbigt}^{j}\lrarr\widehat{\nbigt}^{j+2}(-1)$
and the hermitian sesqui-linear duality
$\widehat{\nbigs}:
 \widehat{\nbigt}\lrarr\widehat{\nbigt}^{\ast}(-w)$.
By the condition,
the tuple 
$\bigl(\widehat{\nbigt},\widehat{N},\widehat{\nbigs}\bigr)$
is a polarized graded wild Lefschetz 
twistor $D$-module of weight $w$ and type $-1$.

We obtain the complex $f_{\dagger}\nbigt$
of $\nbigr_Y$-triples
with the nilpotent map $f_{\dagger}N$ 
and the Hermitian sesqui-linear duality $f_{\dagger}\nbigs$.
(See \cite{sabbah2}.)
We have the induced filtration $M$ on $f_{\dagger}\nbigt$,
which is the weight filtration of $f_{\dagger}N$.
We also have the Lefschetz map 
$\nbigl_{c}:f_{\dagger}\nbigt\lrarr f_{\dagger}\nbigt[2](1)$.

By taking the cohomology,
we obtain $\bigoplus f^j_{\dagger}\nbigt$
with $\bigoplus f^j_{\dagger}N$,
$\bigoplus f^j_{\dagger}\nbigs$,
and the Lefschetz map $\nbigl_c$.
We also have the induced filtration $M$
on each $f^j_{\dagger}\nbigt$.
We put 
$\widetilde{\nbigt}^{i,j}:=
 \Gr^M_{-i}f^j_{\dagger}\nbigt$.
The induced nilpotent map and sesqui-linear duality
are denoted by $\widetilde{N}$ and $\widetilde{\nbigs}$.
The Lefschetz map is also denoted by $\widetilde{\nbigl}_c$.

\begin{prop}
 \label{prop;05.1.15.1}
\mbox{{}}
\begin{itemize}
\item
Let $(\nbigt,M,N,\nbigs)$ be as in 
Condition {\rm\ref{condition;05.1.18.20}}.
Under Assumption {\rm\ref{assumption;05.1.15.11}},
the tuple
$(\widetilde{\nbigt},
 \widetilde{\nbigs},\widetilde{N},\widetilde{\nbigl}_c)$
is a polarized wild bi-graded Lefschetz twistor $D$-module
of type $(-1,1)$ with weight $w$.
\item
The induced filtration $M$ on
$f^j_{\dagger}\nbigt$ is 
equal to the weight filtration
of $f^j_{\dagger}N$. 
\end{itemize}
\end{prop}
\pf
It can be shown by using 
Lemma \ref{lem;07.12.6.4},
Lemma \ref{lem;05.1.18.25}
and the argument in the proof of
Proposition 14.133 in \cite{mochi2}.
\hfill\qed

\subsection{Preliminary II}
\label{subsection;05.1.15.13}

Let $f:X\lrarr Y$ be a projective morphism of 
a complex manifolds.
Let $g$ be a holomorphic function on $Y$.
Let $c$ be the first Chern class of 
a relatively ample line bundle on $X$.
We put $\widetilde{g}=g\circ f$.
Let $Z$ be a closed subvariety of $X$
such that $\widetilde{g}$ is not constantly $0$ on $Z$.
We impose the following assumption
in this subsection.
\begin{assumption}
 \label{assumption;05.1.15.31}
The claim of Theorem {\rm\ref{thm;07.10.23.20}} holds
for the morphism $f$ and 
polarized wild pure twistor $D$-modules
whose strict supports are contained 
in $\widetilde{g}^{-1}(0)\cap Z$.
\hfill\qed
\end{assumption}

Let us consider a tuple
$(\nbigt,\nbigs,\gminia)$ on $X$ as follows:
\begin{condition}
\label{condition;05.1.15.30}
 \mbox{{}}
\begin{itemize}
\item 
$\nbigt$ is a holonomic $\nbigr_X$-triple,
$\nbigs:\nbigt\lrarr\nbigt^{\ast}(-w)$ is 
a hermitian sesqui-linear duality,
and $\gminia\in\cnum[t_m^{-1}]$.
\item 
$\nbigt$ is strictly $S$-decomposable along $\widetilde{g}$.
 The support of $\nbigt$ is contained in $Z$.
Moreover,
it is strictly specializable along $\widetilde{g}$
with ramification and exponential twist by $\gminia$.
\item
For each $u\in\real\times\cnum$,
we put 
$\nbigt_{\gminia,u}:=
 \widetilde{\psi}_{\widetilde{g},\gminia,u}\nbigt$.
The induced nilpotent map and 
 Hermitian sesqui-linear duality
are denoted by $N_{\gminia,u}$ and $\nbigs_{\gminia,u}$.
Let $M$ denote 
 the monodromy weight filtration of $N_{\gminia,u}$.
Then, 
$(\nbigt_{\gminia,u},M,
 N_{\gminia,u},\nbigs_{\gminia,u})$
satisfies Condition {\rm \ref{condition;05.1.18.20}}.
\hfill\qed
\end{itemize}
\end{condition}

We put $\widehat{\nbigt}_{\gminia,u}
 :=\Gr^M\nbigt_{\gminia,u}$.
The induced nilpotent map and 
the Hermitian sesqui-linear duality
are denoted by $\widehat{N}_{\gminia,u}$
and $\widehat{\nbigs}_{\gminia,u}$.
 The tuple
 $\bigl(\widehat{\nbigt}_{\gminia,u},
 \widehat{N}_{\gminia,u},
 \widehat{\nbigs}_{\gminia,u}\bigr)$
 is a polarized graded wild Lefschetz twistor $D$-module 
 of weight $w$ and type $-1$.

\begin{lem}
\label{lem;08.1.16.40}
Let $(\nbigt,\nbigs,\gminia)$ be as in 
Condition {\rm\ref{condition;05.1.15.30}}.
Then, 
$f^j_{\dagger}\nbigt$
is strictly specializable along $g$
with ramification and exponential twist 
by the $\gminia\in\cnum[t_n^{-1}]$,
and we have the natural isomorphism
$\psitilde_{g,\gminia,u}f^j_{\dagger}\nbigt
\simeq
 f^j_{\dagger}\psitilde_{\gtilde,\gminia,u}\nbigt$.
\end{lem}
\pf
Applying Proposition \ref{prop;05.1.15.1}
to any 
$\bigl(\nbigt_{\gminia,u},M,N_{\gminia,u},
 \nbigs_{\gminia,u}\bigr)$,
we obtain 
$f^j_{\dagger}\nbigt_{\gminia,u}$
are strict.
Hence, the claim follows from 
Proposition \ref{prop;07.10.26.100}
and Proposition \ref{prop;08.1.15.30}.
\hfill\qed

\vspace{.1in}

For any
$u\in\real\times\cnum$,
we put
\[
 \widetilde{\nbigt}_{\gminia,u}^{i,j}:=
 \Gr^M_{-i}\psitilde_{g,\gminia,u}
 f^j_{\dagger}\nbigt.
\]
Here, $M$ denotes the weight filtration
for the induced nilpotent maps for $\tildepsi_{g,\gminia,u}$.
Then, we obtain the bi-graded $\nbigr$-triples
$\widetilde{\nbigt}_{\gminia,u}
=\bigoplus \widetilde{\nbigt}_{\gminia,u}^{i,j}$.
The induced nilpotent map and the sesqui-linear duality
are denoted by $\widetilde{N}_{\gminia,u}$
and $\widetilde{\nbigs}_{\gminia,u}$.
The Lefschetz map $\widetilde{\nbigl}_{\gminia,u}$ 
is also induced.

\begin{prop}
\label{prop;05.1.15.61}
Let $(\nbigt,\nbigs,\gminia)$ be
as in Condition {\rm\ref{condition;05.1.15.30}}.
Under Assumption {\rm\ref{assumption;05.1.15.31}},
$\bigl(\widetilde{\nbigt}_{\gminia,u},\widetilde{\nbigs}_{\gminia,u},
 \widetilde{N}_{\gminia,u},\widetilde{\nbigl}_{\gminia,u}\bigr)$
is a polarized bi-graded wild Lefschetz twistor $D$-module
of type $(-1,1)$ with weight $w$.
\end{prop}
\pf
It can be shown by using Proposition \ref{prop;05.1.15.1}
and the same argument as that in the proof of
Proposition 14.139 of \cite{mochi2}.
\hfill\qed

\vspace{.1in}

Let $(\nbigt,\nbigs)$ be as in
Condition \ref{condition;05.1.15.30}
with $\gminia=0$.
We put
$\nbigt_0:=\phi_{\widetilde{g},0}(\nbigt)(-1/2)$,
and $\widehat{\nbigt}^i_0:=\Gr^M_{-i}\nbigt_0$.
Then, we obtain the graded $\nbigr_X$-triple $\widehat{\nbigt}_0$.
We have the induced nilpotent maps $\widehat{N}_0$
and the sesqui-linear duality $\widehat{\nbigs}_0$.
The tuple 
$\bigl(\widehat{\nbigt}_0,
 \widehat{\nbigs}_0,
 \widehat{N}_0\bigr)$
is a polarized graded wild 
 Lefschetz twistor $D$-module of weight $w+1$.
We have the maps
$\Can:\widehat{\nbigt}_{-\vecdelta_0}
 \lrarr\widehat{\nbigt}_0$
and $\Var:\widehat{\nbigt}_0
 \lrarr\widehat{\nbigt}_{-\vecdelta_0}(-1)$,
such that $\Nhat_{-\vecdelta_0}=\Var\circ\Can$
and $\Nhat_0=\Can\circ\Var$.
The maps $\Can$ and $\Var$ are adjoint with respect to 
the sesqui-linear dualities.
The nilpotent map $N_0$ induces 
the nilpotent map $\widetilde{N}_0$
on $f^j_{\dagger}\nbigt_0$.
Let $M$ denote the weight filtration.
We put
$\widetilde{\nbigt}^{i,j}_0:=
\Gr^M_{-i}f^j_{\dagger}\nbigt_0$.
The induced sesqui-linear duality
is denoted by $\widetilde{\nbigs}_0$.
We also have the induced Lefschetz map 
 $\widetilde{\nbigl}_0$.

\begin{prop}
\label{prop;05.1.15.60}
\mbox{{}}
\begin{itemize}
\item
The underling $\nbigr_Y$-modules
of $f^i_{\dagger}\nbigt$
are strictly $S$-decomposable along $g$.
\item
The tuple
$\bigl(\widetilde{\nbigt}_0,
 \widetilde{\nbigs}_0,\widetilde{N}_0,
 \widetilde{\nbigl}_0\bigr)$
is a bi-graded wild Lefschetz twistor $D$-module 
of weight $w+1$ and type $(-1,1)$.
\end{itemize}
\end{prop}
\pf
It can be shown by using
Proposition \ref{prop;05.1.15.1},
Lemma \ref{lem;05.2.14.20}
and the same argument as that in the proof of
Proposition 14.140 of \cite{mochi2}.
\hfill\qed

\vspace{.1in}

Let $\nbigm_i'$ and $\nbigm_i''$ are the underlying $\nbigr_Y$-modules
of $f^i_{\dagger}\nbigt$.
Due to the $S$-decomposability of $\nbigm_i'$ and $\nbigm_i''$,
we obtain the sesqui-linear pairing
of $\psi_{g,0}(\nbigm'_i)$ and $\psi_{g,0}(\nbigm''_i)$
by the specialization.
As a result, we obtain the 
$\nbigr$-triple
$\phi_{g,0}f^i_{\dagger}\nbigt$.

\begin{lem}
\label{lem;05.1.15.62}
$\phi_{g,0}f^i_{\dagger}\nbigt(-1/2)$
and $f^i_{\dagger}\nbigt_0$ 
are isomorphic.
\end{lem}
\pf
It follows from
Lemma \ref{lem;07.10.27.5}
and Proposition \ref{prop;05.1.15.60}.
\hfill\qed

\subsection{The explanation for Step 2}
\label{subsection;05.1.18.10}

The following assumption for the induction
is imposed in this subsection.
\begin{assumption}
Let $n$ and $m$ be non-negative integers.
Let $f:X\lrarr Y$ be a projective morphism.
Let $(\nbigt,\nbigs)$ be a polarized wild pure twistor $D$-module
of weight $w$ on $X$
whose strict support $\Supp(\nbigt)$ is irreducible.
We put
$n(\nbigt)=\dim\Supp(\nbigt)$
and $m(\nbigt)=\dim f(\Supp(\nbigt))$.
The claim of Theorem {\rm\ref{thm;07.10.23.20}} holds
for $\nbigt$
if $n(\nbigt)\leq n$ and $m(\nbigt)\leq m$ are satisfied.
\hfill\qed
\end{assumption}
Under the assumption, we show that
the claim of Theorem {\rm\ref{thm;07.10.23.20}}
holds for $\nbigt$
such that $n(\nbigt)\leq n+1$ and $m(\nbigt)\leq m+1$ 
are satisfied.

Since the claim is local on $Y$,
we may replace $Y$ with any open subset of $Y$.
Let $g$ be any holomorphic function on $Y$
such that $g$ is not constantly $0$
on $f\bigl(\Supp(\nbigt)\bigr)$.
Due to Proposition \ref{prop;05.1.15.60},
the underlying $\nbigr_Y$-modules of 
$\bigoplus f^i_{\dagger}\nbigt$
are strictly $S$-decomposable along $g$.
Then, we have the decomposition 
$f^j_{\dagger}\nbigt
=\bigoplus_Z(f^j_{\dagger}\nbigt)_Z$
by the strict supports.
The decomposition is compatible with
the Lefschetz map $\nbigl_c$
and the induced Hermitian sesqui-linear duality
$\bigoplus f^j_{\dagger}\nbigs$.
By using
Proposition \ref{prop;05.1.15.60}
and Lemma \ref{lem;05.1.15.62},
we can check that
$\bigoplus_j (f^j_{\dagger}\nbigt)_Z$
with $(\nbigl_c)_Z$ and 
$\bigoplus_j(f^j_{\dagger}\nbigs)_Z$
is a polarized graded Lefschetz twistor structure
if $\dim Z=0$.
Then,
by using Proposition \ref{prop;05.1.15.61},
we can check that
$\bigl(\bigoplus f^j_{\dagger}\nbigt,
 \nbigl_c,\bigoplus_j f^j_{\dagger}\nbigs\bigr)$
is a polarized graded wild Lefschetz twistor $D$-module
of weight $w$.
Thus, we finish Step 2 in the proof of 
Theorem \ref{thm;07.10.23.20}.

\section{Step 3}
\label{subsection;07.10.23.21}
\subsection{Statement}

Let $X$ be a smooth projective variety
with an ample line bundle $\nbigo_{X}(1)$.
Let $a_X$ denote the obvious morphism
of $X$ to a point.
Let $(\nbigt,\nbigs)$ be 
a polarized wild pure twistor $D$-module
on $X$ of weight $w$,
whose strict support is an irreducible closed subset $Z$
of $X$ with $\dim Z=n$.
We will show the following proposition
in this section.
\begin{prop}
\label{prop;07.10.28.10}
Assume $P(n-1,0)$ holds.
Then, 
$\Bigl(\bigoplus_i a^i_{X\dagger}\nbigt,
 \nbigl_c,
 \bigoplus_i a^i_{X\dagger}\nbigs\Bigr)$
is a polarized graded Lefschetz twistor structure 
of weight $w$,
where $\nbigl_c$ denote the Lefschetz morphism
associated to $\nbigo_X(1)$.
\end{prop}

We follow the argument of \cite{saito1} very closely.
We also use the results in \cite{sabbah2}.
We may assume 
$w=0$, $\nbigt=(\nbigm,\nbigm,C)$,
and $\nbigs=(\id,\id)$.

\subsection{Preliminary}
\label{subsection;08.1.18.18}

We fix an embedding $X\subset\proj^N$.
Let $X_0$ be the intersection $X\cap H_1\cap H_2$,
where $H_i$ denote general hyperplanes
such that $X_0$ is 
strictly non-characteristic to $\nbigt$.
Let $\Xtilde$ be the blow up of $X$ along $X_0$.
We obtain the following diagram:
\[
 \begin{CD}
 X @<{\pi}<< \Xtilde @>{p}>> \proj^1
 \end{CD}
\]
Here, $p$ is the Lefschetz pencil for $H_1$ and $H_2$.
We have the $\nbigr$-triple
$(\nbigttilde,\nbigstilde)$ on $\Xtilde$,
as in \cite{sabbah2}.
Let $\nbigmtilde$ denote the underlying 
$\nbigr_{\Xtilde}$-module of $\nbigttilde$.
We have the naturally induced Hermitian sesqui-linear
duality $\nbigstilde$ of $\nbigttilde$.

\begin{lem}
\label{lem;08.1.18.10}
$\Bigl(\bigoplus_{i}
p^i_{\dagger}\nbigttilde,
 \nbigl_c,\bigoplus_i p^i\nbigstilde
\Bigr)$ is a polarized 
graded wild Lefschetz twistor $D$-module
on $\proj^1$,
where $\nbigl_c$ denotes the Lefschetz map
associated to $\pi^{\ast}\nbigo_{X}(1)$.
\end{lem}
\pf
See (1) and (2) in Chapter 6.4 of \cite{sabbah2}.
\hfill\qed

\begin{lem}
\label{lem;08.1.18.11}
We have the vanishing
$\pi^j_{\dagger}\nbigttilde=0$ for $j\neq 0$,
and the canonical decomposition
$(\pi_{\dagger}\nbigttilde,
 \pi_{\dagger}\nbigstilde)
=(\nbigt,\nbigs)\oplus
 (\nbigt_1,\nbigs_1)$
of $\nbigr_X$-triples
with Hermitian sesqui-linear duality.
\end{lem}
\pf
See (3) in Section 6.4 of \cite{sabbah2}.
We just give a remark on 
the vanishing 
$\pi^j_{\dagger}\nbigttilde=0$ for $j\neq 0$.
According to (3) in Section 6.4 of \cite{sabbah2},
it can be shown that the underlying $\nbigr_X$-modules
of $\pi^j_{\dagger}\nbigttilde$
are strict.
Hence, we have only to show 
the vanishing of the specialization 
$\pi_{\dagger}^j\pi^{\dagger}\nbigm^{1}$
at $\lambda=1$,
where $\nbigm^1$ denotes the specialization of
$\nbigm$ along $\lambda=1$.
Hence, we can reduce the issue
to the case of $D_X$-modules.
By the projection formula,
we have the isomorphism
$\pi_{\dagger}
 \pi^{\dagger}\nbigm^{1}
\simeq
 \nbigm^{1}
 \otimes^L\pi_{\dagger}\nbigo_{\Xtilde}$
in the derived category of $D_X$-modules.
We have the vanishing
$\pi^j_{\dagger}\nbigo_{\Xtilde}=0$
for $j\neq 0$,
and $\pi_{\dagger}^0\nbigo_X$
is isomorphic to the direct sum
$\nbigo_X\oplus M_2$,
where $M_2$ is locally isomorphic to
the push-forward of $\nbigo_{X_0}$.
(See Section 5.3.9 of \cite{saito1},
for example.)
Since $X_0$ is non-characteristic to
$\nbigm^{1}$,
we obtain that 
$\nbigm^{1}\otimes^L M_2
\simeq
 \nbigm^{1}\otimes M_2$.
\hfill\qed

\vspace{.1in}
By Lemma \ref{lem;08.1.18.10} and 
Deligne's result \cite{d5},
the spectral sequence
$a^i_{\proj^1\dagger}p^j_{\dagger}\nbigttilde
\Longrightarrow
 a^{i+j}_{\Xtilde\dagger}\nbigttilde$
degenerates at the $E^2$-level.
Hence, due to Proposition \ref{prop;05.1.15.1},
$a_{\Xtilde\dagger}^i\nbigttilde$
is a pure twistor structure of weight $i$.
According to Lemma \ref{lem;08.1.18.11},
$a^i_{X\dagger}\nbigt$ is a direct summand of
$a^{i}_{\Xtilde\dagger}\nbigttilde$.
Hence, $a^i_{X\dagger}\nbigt$
is a pure twistor structure of weight $i$.

\vspace{.1in}

We have the Leray filtration:
\[
  L^{1}a^l_{\Xtilde\dagger}\nbigmtilde
\subset
 L^{0}a^l_{\Xtilde\dagger}\nbigmtilde
\subset
 L^{-1}a^l_{\Xtilde\dagger}\nbigmtilde
\]
Because of the degeneration of the spectral sequence,
we have the isomorphisms
$ \Gr_L^i\bigl(a^l_{\Xtilde\dagger}\nbigmtilde\bigr)
\simeq a^i_{\proj^1\dagger}
 p^{l-i}_{\dagger}\nbigmtilde$.
We also have the Leray filtration 
for pure twistor structures of weight $l$:
\[
 L^{1}a^l_{\Xtilde\dagger}\nbigttilde
\subset
 L^{0}a^l_{\Xtilde\dagger}\nbigttilde
\subset
 L^{-1}a^l_{\Xtilde\dagger}\nbigttilde
\]
We have the isomorphisms
$\Gr_L^i\bigl(a^l_{\Xtilde\dagger}\nbigttilde\bigr)
\simeq a^i_{\proj^1\dagger}
 p^{l-i}_{\dagger}\nbigttilde$.

\vspace{.1in}

We will use the following general lemma
in linear algebra.
\begin{lem}
\label{lem;08.1.19.10}
Let $V$ be a pure twistor structure of
weight $w$ with a polarization $S$.
Let $V_1$ be a pure twistor structure 
of weight $w$
with a monomorphism $\varphi:V_1\lrarr V$.
We have the following maps:
\[
\begin{CD}
 V_1 @>{\varphi}>>
 V @>{S}>> V^{\ast}\otimes\Tate(-w)
 @>{\varphi^{\ast}}>>
 V_1^{\ast}\otimes\Tate(-w)
\end{CD}
\]
Then, the composite $S_1$ gives 
a polarization of $V_1$.
\end{lem}
\pf
We can reduce the problem to the case $w=0$.
Then, the claim is the same as the following:
\begin{itemize}
\item
 Let $V$ be a $\cnum$-vector space,
 and let $V_1$ be a subspace of $V$.
 Let $h$ be a hermitian metric of $V$.
 Let $h_1$ denote the restriction of $h$ to $V_1$.
 We have the induced anti-linear maps
 $\varphi_h:V\lrarr V^{\lor}$
 and $\varphi_{h_1}:V_1\lrarr V_1^{\lor}$.
 Let $i$ denote the inclusion $V_1\lrarr V$.
 Then, the following diagram is commutative:
\[
 \begin{CD}
 V_1 @>{\varphi_{h_1}}>> V_1^{\lor}\\
 @V{i}VV @A{i^{\lor}}AA \\
 V @>{\varphi_{h}}>> V^{\lor}
 \end{CD}
\]
\end{itemize}
The claim can be checked directly.
\hfill\qed

\subsection{Gysin maps}

Let $Y$ denote a general fiber of $p$.
For $Y\subset X$,
we apply the construction in 
Section \ref{subsection;08.1.18.13}.
We have obtained
the polarized wild pure twistor
$D$-module
$(\Ker(\nbigt),\nbigs_{\Ker(\nbigt)})$
with weight $-1$.
(See Section \ref{subsection;08.1.18.14}.)
Since the support of $\Ker(\nbigt)$ is contained in $Y$,
the hypothesis of the induction can be applied.

Recall we have obtained the morphisms
$i_{Y,\Gys}^{\ast}$
and $i_{Y,\ast}^{\Gys}$
in Section \ref{subsection;08.1.19.5}.

\begin{lem}
\label{lem;07.10.28.1}
The following holds:
\begin{itemize}
\item
 $i_{Y,\Gys}^{\ast}:
 a^i_{X\dagger}\nbigt\lrarr
     a^{i+1}_{X\dagger}\Ker(\nbigt)$
 is an isomorphism for $i<-1$
 and a monomorphism for $i=-1$.
\item
 The morphism
$(i_{Y,\Gys}^{\ast})^{\ast}:
 \bigl(a_{X\dagger}^{i+1}\Ker(\nbigt)\bigr)^{\ast}
\lrarr
 a_{X\dagger}^{i}(\nbigt)^{\ast}$
is an isomorphism if $i<-1$,
and an epimorphism if $i=-1$.
In other words,
 $i_{Y\ast}^{\Gys}:
 a^i_{X\dagger}\Cok(\nbigt)\lrarr
 a^{i+1}_{X\dagger}\nbigt$
 is an isomorphism for $i>0$
 and an epimorphism for $i=0$.
\end{itemize}
\end{lem}
\pf
Since we have already known
that
$a^i_{X\dagger}\Ker(\nbigt)$,
$a^i_{X\dagger}\Cok(\nbigt)$
and $a^i_{X\dagger}\nbigt$ 
are pure twistor structures with appropriate weights
(Section \ref{subsection;08.1.18.18}),
the claim follows from Lemma \ref{lem;08.1.19.1}.
\hfill\qed

\begin{lem}
\label{lem;08.1.19.12}
The morphism
$L_c^i:a_{X\dagger}^{-i}\nbigm
\lrarr
 a_{X\dagger}^{i}\nbigm$ is an isomorphism
for $i\geq 2$.
\end{lem}
\pf
We identify 
$j_{\ast}j^{\ast}\nbigm/\nbigm$
and $i_{Y\dagger}i_Y^{\dagger}\nbigm$
by $\Lambda$. (See Section \ref{subsection;08.1.17.25}
for $\Lambda$.)
According to Corollary \ref{cor;07.10.27.20},
we have the factorization of $L_c^i$ 
up to signature,
as follows:
\[
\begin{CD}
a_{X\dagger}^{-i}\nbigm
@>{b_1}>>
 a_{X\dagger}^{-i+1}
 i_{Y\dagger}i_Y^{\dagger}\nbigm
@>{L_c^{i-1}}>>
 a_{X\dagger}^{i-1}
 i_{Y\dagger}i_Y^{\dagger}\nbigm
@>{b_2}>>
 a_{X\dagger}^{i}\nbigm
\end{CD}
\]
Due to Lemma \ref{lem;07.10.28.1},
$b_j$ $(j=1,2)$ are isomorphisms.
Since we can apply the hypothesis of the induction
on $\Ker(\nbigt)$,
we also have the middle arrow is an isomorphism.
Thus, we obtain Lemma \ref{lem;08.1.19.12}.
\hfill\qed

\begin{lem}
The following diagram is commutative:
\begin{equation}
 \label{eq;08.1.18.22}
 \begin{CD}
 a_{X\dagger}^{-i}\nbigt
 @>{\nbigs\circ\nbigl_c^i}>>
 \bigl(
 a_{X\dagger}^{-i}\nbigt
 \bigr)^{\ast}\otimes\Tate^S(i)\\
 @V{i_{Y,\Gys}^{\ast}}VV 
 @A{(i_{Y,\Gys}^{\ast})^{\ast}}AA \\
 a_{X\dagger}^{-i+1}\Ker(\nbigt)
 @>{\nbigs_{\Ker(\nbigt)}\circ \nbigl_c^{i-1}}>>
 \bigl(
 a_{X\dagger}^{-i+1}\Ker(\nbigt)
 \bigr)^{\ast}\otimes\Tate^{S}(i)
 \end{CD}
\end{equation}
Here $(i_{Y,\Gys}^{\ast})^{\ast}$
denotes the Hermitian adjoint of
$i_{Y,\Gys}^{\ast}$,
and $\nbigs_{\Ker(\nbigt)}$
denotes the polarization of
$\Ker(\nbigt)$
given in Section
{\rm \ref{subsection;08.1.18.14}}.
\end{lem}
\pf
We have the commutativity of the diagram
(\ref{eq;08.1.18.21})
in Corollary \ref{cor;07.10.27.20}.
Note that $\nbigs_{\Ker(\nbigt)}$
can be factorized as follows:
\[
 \Ker(\nbigt)
\stackrel{\varphi_3^{-1}}{\lrarr}
 \Cok(\nbigt)\otimes\Tate^S(1)
\stackrel{\varphi_5}{\lrarr}
 \Cok(\nbigt^{\ast})\otimes\Tate^S(1)
\lrarr
 \Ker(\nbigt)^{\ast}\otimes\Tate^S(1)
\]
Here, $\varphi_3$ is given 
in Section \ref{subsection;08.1.18.14},
$\varphi_5$ is induced 
by the polarization $\nbigs$ of $\nbigt$,
and the right arrow is the natural identification
given in Lemma \ref{lem;08.1.18.15}.

We have the following commutative diagram:
\begin{equation}
\label{eq;08.1.18.20}
 \begin{CD}
 a^{i}_{X\dagger}\nbigt
 @>{\nbigs}>> 
 a^i_{X\dagger}(\nbigt^{\ast})
 @>>>
 a^{i}_{X\dagger}(\nbigt^{\ast})
 @>>>
 \bigl(a_{X\dagger}^{-i}\nbigt\bigr)^{\ast}
 \\
 @AAA @AAA @AAA @AAA \\
 a_{X\dagger}^{i-1}\Cok(\nbigt)
 @>{\varphi_5}>>
 a_{X\dagger}^{i-1}\Cok(\nbigt^{\ast})
 @>>>
 a_{X\dagger}^{i-1}\bigl(
 \Ker(\nbigt)^{\ast}\bigr)
 @>>>
 \bigl(a_{X\dagger}^{-i+1}\Ker(\nbigt)\bigr)^{\ast}
 \end{CD}
\end{equation}
The horizontal arrows are the natural identifications.
Composing (\ref{eq;08.1.18.21})
and (\ref{eq;08.1.18.20}),
we obtain the commutativity of the diagram
(\ref{eq;08.1.18.22}).
\hfill\qed

\subsection{The case $i\geq 2$}

Let us consider the case $i\geq 2$.
Let $Pa^{-i}_{X\dagger}\nbigt$ denote the kernel
of $\nbigl_{c}^{i+1}:a^{-i}_{X\dagger}\nbigt
 \lrarr a^{i+2}_{X\dagger}\nbigt\otimes\Tate^S(i+1)$.
\begin{lem}
We have the following commutative diagram:
\begin{equation}
 \label{eq;08.1.19.23}
 \begin{CD}
 Pa^{-i}_{X\dagger}\nbigt 
 @>{\nbigs\circ\nbigl^{i}}>> 
 (Pa^{-i}_{X\dagger}\nbigt)^{\ast}\otimes\Tate^S(i)\\
 @V{\simeq}VV @A{\simeq}AA \\
 Pa^{-i+1}_{X\dagger}\Ker(\nbigt) 
@>{\nbigs_{\Ker(\nbigt)}\circ\nbigl_c^{i-1}}>>
 \bigl(Pa^{-i+1}_{X\dagger}\Ker(\nbigt)\bigr)^{\ast}
 \otimes\Tate^S(i)
 \end{CD}
\end{equation}
\end{lem}
\pf
By Lemma \ref{lem;08.1.19.12},
we have the decompositions
$a_{X\dagger}^{-i}(\nbigm)
=Pa_{X\dagger}^{-i}(\nbigm)
\oplus
 \nbign_1$
and
$a^i_{X\dagger}\nbigm
=Pa^i_{X\dagger}\nbigm
\oplus
 \nbign_2$:
\[
Pa_{X\dagger}^{-i}(\nbigm)
:=\Ker\bigl(
 L_c^{i+1}:a_{X\dagger}^{-i}\nbigm
\lrarr a_{X\dagger}^{i+2}\nbigm
 \bigr),
\quad
\nbign_1:=
 \Image\bigl(
 L_c:a_{X\dagger}^{-i-2}\nbigm
\lrarr a_{X\dagger}^{-i}\nbigm
 \bigr)
\]
\[
 Pa_{X\dagger}^{i}(\nbigm)
:=\Image\bigl(
 L_c^{i}:Pa_{X\dagger}^{-i}\nbigm
\lrarr a_{X\dagger}^{i}\nbigm
 \bigr),
\quad
\nbign_2:=
 \Image\bigl(
 L_c^i: \nbign_1
\lrarr a_{X\dagger}^{i}\nbigm
 \bigr)
\]
Let $\nbigt_1$ denote the image of
$\nbigl_c:a^{-i+2}_{X\dagger}\nbigt\otimes\Tate^S(-1)
\lrarr a^{-i}_{X\dagger}\nbigt$.
Then, the underlying $\nbigr$-modules
of $Pa^{-i}_{X\dagger}\nbigt$
(resp. $\nbigt_1$) are given by
$Pa^{i}_{X\dagger}\nbigm$
and $Pa^{-i}_{X\dagger}\nbigm$
(resp. $\nbign_2$ and $\nbign_1$).
We also have a similar decomposition
for $a_{X\dagger}^{-i+1}(\Ker(\nbigt))$.
Then, we can directly check that 
the morphisms in the diagram (\ref{eq;08.1.18.22})
preserves the decompositions.
Hence, 
we obtain (\ref{eq;08.1.19.23})
from (\ref{eq;08.1.18.22}).
\hfill\qed

\vspace{.1in}
Since the lower horizontal arrow gives a polarization of
$Pa^{-i+1}_{X\dagger}\Ker(\nbigt)$
by the hypothesis of the induction,
the Hermitian sesqui-linear duality 
$\nbigs\circ\nbigl_c^i$ gives a polarization
of $Pa^{-i}_{X\dagger}\nbigt$.

\subsection{The case $i=1$}

In the case $i=1$,
we obtain the following commutative diagram
from (\ref{eq;08.1.18.22}):
\begin{equation}
 \label{eq;07.12.2.35}
 \begin{CD}
 a^{-1}_{X\dagger}\nbigt @>{\nbigs\circ\nbigl_c}>> 
 (a^{-1}_{X\dagger}\nbigt)^{\ast}\otimes\Tate^S(1)\\
 @V{i_{Y,\Gys}^{\ast}}VV 
 @A{(i_{Y,\Gys}^{\ast})^{\ast}}AA \\
 a^0_{X\dagger}\Ker(\nbigt) @>{\nbigs_{\Ker(\nbigt)}}>>
 \bigl(a^0_{X\dagger}\Ker(\nbigt)\bigr)^{\ast}
 \otimes\Tate^S(1)
 \end{CD}
\end{equation}

Since we have already known that
$L_c^3:
 a_{X\dagger}^{-3}\nbigm
\lrarr
 a_{X\dagger}^3\nbigm$ is an isomorphism,
we have the following decompositions:
\[
 a^{-1}_{X\dagger}\nbigm
=L_c\cdot a^{-3}_{X\dagger}\nbigm
 \oplus\Ker\Bigl(
 L_c^2:a^{-1}_{X\dagger}\nbigm\lrarr 
      a^{3}_{X\dagger}\nbigm
 \Bigr)
\]
\[
  a^1_{X\dagger}\nbigm
=L_c^2 \cdot a^{-3}_{X\dagger}\nbigm
 \oplus \Ker \Bigl(
 L_c:a^{1}_{X\dagger}\nbigm\lrarr
 a^3_{X\dagger}\nbigm\Bigr)
\]
This decomposition is compatible 
with the sesqui-linear pairing.
We set
\[
  Pa^{-1}_{X\dagger}\nbigt:=
 \bigl(\Ker L_c,\Ker L_c^2,C_1\bigr),
\quad
 \nbigt_1:=
 \bigl(\Image L_c^2,\Image L_c,C_2\bigr)
\]
Here, $C_i$ denote
the naturally induced sesqui-linear pairings.
We have 
$a_{X\dagger}^{-1}\nbigt
=Pa^{-1}_{X\dagger}\nbigt
 \oplus \nbigt_1$.
We have already known that
$a_{X\dagger}^{-1}\nbigt$
is a pure twistor structure with weight $-1$.
(Section \ref{subsection;08.1.18.18}).
Hence,
$Pa_{X\dagger}^{-1}\nbigt$ and
$\nbigt_1$ are also
pure twistor structures with weight $-1$.
The morphism 
$\nbigs\circ\nbigl_c:
 a_{X\dagger}^{-1}\nbigt
\lrarr
 (a_{X\dagger}^{-1}\nbigt)^{\ast}
 \otimes\Tate^S(1)$ preserves
the decomposition.
It is easy to show that
the induced map
$\nbigt_1\lrarr 
 \nbigt_1^{\ast}\otimes\Tate^S(1)$
is an isomorphism,
by using the fact
that $L_c^3:a_{X\dagger}^{-3}\nbigm
 \lrarr a_{X\dagger}^3\nbigm$
is an isomorphism.

We have the primitive decomposition
$a^0_{X\dagger}\bigl(
 \Ker(\nbigt)\bigr)
=Pa^0_{X\dagger}\bigl(\Ker(\nbigt)\bigr)
\oplus
 \nbigt_2$.
Here, the underlying $\nbigr_X$-modules of
$Pa^0_{X\dagger}\Ker\bigl(\nbigt\bigr)$
are given by the following:
\[
 \Ker\Bigl(
 L_c:a_{X\dagger}^0\bigl(
 j_{\ast}j^{\ast}\nbigm/\nbigm\bigr)
\lrarr
 a_{X\dagger}^2\bigl(
 j_{\ast}j^{\ast}\nbigm/\nbigm\bigr)
 \Bigr),
\]
\[
  \Ker\Bigl(
 L_c:a_{X\dagger}^0\bigl(
 i_{Y\dagger}i_Y^{\dagger}\nbigm
 \bigr)
\lrarr
 a_{X\dagger}^2\bigl(
 i_{Y\dagger}i_Y^{\dagger}\nbigm
 \bigr)
 \Bigr)
\]
The underlying $\nbigr$-modules of $\nbigt_2$
are given by the following:
\[
  \Image\Bigl(
 L_c:a_{X\dagger}^{-2}\bigl(
 j_{\ast}j^{\ast}\nbigm/\nbigm\bigr)
\lrarr
 a_{X\dagger}^0\bigl(
 j_{\ast}j^{\ast}\nbigm/\nbigm\bigr)
 \Bigr),
\]
\[
  \Image\Bigl(
 L_c:a_{X\dagger}^{-2}\bigl(
 i_{Y\dagger}i_Y^{\dagger}\nbigm
 \bigr)
\lrarr
 a_{X\dagger}^0\bigl(
 i_{Y\dagger}i_Y^{\dagger}\nbigm
 \bigr)
 \Bigr)
\]

\begin{lem}
\label{lem;08.1.18.40}
The morphism
$i_{Y,\Gys}^{\ast}:
 a_{X\dagger}^{-1}\nbigt
\lrarr a_{X\dagger}^0\Ker(\nbigt)$
induces the morphisms:
\[
 Pa_{X\dagger}^{-1}\nbigt
\lrarr Pa_{X\dagger}^0\Ker(\nbigt),
\quad\quad
 \nbigt_1\lrarr \nbigt_2
\]
\end{lem}
\pf
We put
$\nbigm_0:=
 i_{Y\dagger}i_Y^{\dagger}\nbigm$.
We will identify $\nbigm_0$
and $j_{\ast}j^{\ast}\nbigm/\nbigm$
via $\Lambda$.
We have the following commutative diagram:
\[
 \begin{CD}
 a_{X\dagger}^{-3}\nbigm
 @>{L_c}>>
 a_{X\dagger}^{-1}\nbigm
 @>{L_c^2}>>
 a_{X\dagger}^3\nbigm\\
 @VVV @V{f}VV @AA{\simeq}A\\
 a^{-2}_{X\dagger}\nbigm_0
 @>{L_c}>>
 a^{0}_{X\dagger}\nbigm_0
 @>{L_c}>>
 a^2_{X\dagger}\nbigm_0
 \end{CD}
\]
Hence, we obtain
$f(\Image L_c)\subset \Image(L_c)$
and 
$f(\Ker L_c^2)\subset \Ker L_c$.

We have the following commutative diagram:
\[
 \begin{CD}
 a^{-2}_{X\dagger}\nbigm_0
 @>{L_c}>>
 a^{0}_{X\dagger}\nbigm_0
 @>{L_c}>>
 a^2_{X\dagger}\nbigm_0\\
 @A{\simeq}AA @V{g}VV @VVV \\
 a_{X\dagger}^{-3}\nbigm
 @>{L_c^2}>>
 a_{X\dagger}^{1}\nbigm
 @>{L_c}>>
 a_{X\dagger}^3\nbigm\\
 \end{CD}
\]
Hence, we obtain
$g\bigl(\Ker(L_c)\bigr)
 \subset \Ker(L_c)$
and $g\bigl(\Image(L_c)\bigr)
 \subset\Image(L_c^2)$.
Thus, Lemma \ref{lem;08.1.18.40} is proved.
\hfill\qed

\vspace{.1in}

We obtain the following commutative diagram
from (\ref{eq;07.12.2.35})
and Lemma \ref{lem;08.1.18.40}:
\[
 \begin{CD}
 Pa^{-1}_{X\dagger}\nbigt
 @>{\nbigs\circ\nbigl_c}>>
 (Pa^{-1}_{X\dagger}\nbigt)^{\ast}
 \otimes\Tate^S(1) \\
 @V{b_1}VV 
 @A{b_1^{\ast}}AA\\
 Pa^0_{X\dagger}\Ker(\nbigt) 
 @>{\nbigs_{\Ker(\nbigt)}}>>
 Pa^0_{X\dagger}\Ker(\nbigt)^{\ast}\otimes\Tate^S(1)
 \end{CD}
\]
We know that $b_1$ is monomorphic,
and $b_1^{\ast}$ is epimorphic
(Lemma \ref{lem;07.10.28.1}).
By the hypothesis of the induction,
$\nbigs_{\Ker(\nbigt)}$ 
gives a polarization of 
$Pa^0_{X\dagger}\Ker(\nbigt)$.
Hence, we can conclude that
$\nbigl_c:
 Pa^{-1}_{X\dagger}\nbigt\lrarr
Pa^1_{X\dagger}\nbigt\otimes\Tate^S(1)$
is an isomorphism,
and 
$(Pa^{-1}_{X\dagger}\nbigt,\nbigs')$ is 
a polarized pure twistor structure of weight $-1$
due to Lemma \ref{lem;08.1.19.10},
where $\nbigs':=\nbigs\circ\nbigl_c$.
Since the induced map
$\nbigt_1\lrarr \nbigt_1^{\ast}\otimes\Tate^S(1)$
is an isomorphism,
we also obtain that
$\nbigl_c:a^{-1}_{X\dagger}\nbigt
\lrarr a^1_{X\dagger}\nbigt\otimes\Tate^S(1)$
is an isomorphism.

\subsection{The case $i=0$}

We have already known
$\nbigl_c^2:a^{-2}_{X\dagger}\nbigt
\lrarr a^2_{X\dagger}\nbigt\otimes\Tate(2)$
is an isomorphism.
In particular,
we obtain the decompositions of $\nbigr$-modules
$a_{X\dagger}^0\nbigm=
 \Image L_c\oplus Pa_{X\dagger}^0\nbigm$
and the $\nbigr$-triples
$a_{X\dagger}^0\nbigt
 =\Image\nbigl_c\oplus Pa_{X\dagger}^0\nbigt$,
where $Pa^0_{X\dagger}\nbigm:=\Ker L_c$
and $Pa^0_{X\dagger}\nbigt:=\Ker\nbigl_c$.
We have already known that
$\Image\nbigl_c$ and $Pa^0_{X\dagger}\nbigt$
are pure twistor structures of weight $0$.
The decomposition is compatible
with the Hermitian sesqui-linear duality.
Let $\nbigs_0$ denote the induced Hermitian
sesqui-linear duality of $Pa^0_{X\dagger}\nbigt$.
We would like to show that
$\nbigs_0$ gives a polarization of
$Pa^0_{X\dagger}\nbigt$.

\vspace{.1in}
We have the following pure twistor structures:
\[
 L^0a^0_{\Xtilde\dagger}\nbigttilde
=\left(
 \frac{a_{\Xtilde\dagger}^0\nbigmtilde}
 {L^1a_{\Xtilde\dagger}^0\nbigmtilde},\,\,
 L^0a_{\Xtilde\dagger}^0\nbigmtilde,\,\,
 C_1
\right),
\quad
\frac{a^0_{\Xtilde\dagger}\nbigttilde^{\ast}}
 {L^1a^0_{\Xtilde\dagger}\nbigttilde^{\ast}}
=(L^0a_{\Xtilde\dagger}^0\nbigttilde)^{\ast}
\]
Here, $C_1$
denotes the naturally induced sesqui-linear pairing.

\begin{prop}
\label{prop;07.10.28.2}
We have the following factorizations:
\[
Pa^0_{X\dagger}\nbigt \lrarr 
L^0a^0_{\Xtilde\dagger}\nbigttilde \lrarr
  a^0_{\Xtilde\dagger}\nbigttilde,
\quad\quad
 a^0_{\Xtilde\dagger}\nbigttilde^{\ast}
\lrarr
 \frac{a^0_{\Xtilde\dagger}\nbigttilde^{\ast}}
 {L^1a^0_{\Xtilde\dagger}\nbigttilde^{\ast}}
\lrarr
 Pa^0_{X\dagger}\nbigt^{\ast}
\]
\end{prop}
\pf
We have only to show the first factorization.
We have already known the twistor property of
$Pa^0_{X\dagger}\nbigt$,
$L^0a^0_{\Xtilde\dagger}\nbigttilde$ and
$a^0_{\Xtilde\dagger}\nbigttilde$.
Hence, we have only to have
such factorization for the specialization at $\lambda=1$.
Namely, we have only to show that
$ Pa_{X\dagger}^0\nbigm^1$
is contained in
$L^0a_{\Xtilde}^0\nbigmtilde^1$.
We put $\Ytilde:=Y$
and let $i_{\Ytilde}:\Ytilde\lrarr \Xtilde$
denote the natural inclusion.

\begin{lem}
We have the following commutative diagram:
\begin{equation}
\label{eq;08.1.18.50}
 \begin{CD}
 0 @>>> \nbigm^1 @>>> 
 j_{\ast}j^{\ast}\nbigm^1
 @>>> j_{\ast}j^{\ast}\nbigm^1/\nbigm^1
 @>>> 0\\
 @. @A{f_1}AA @A{f_2}AA @A{=}AA @. \\
 0 @>>> \pi_{\dagger}\nbigmtilde^1 @>>> 
 \pi_{\dagger}(\jtilde_{\ast}\jtilde^{\ast}\nbigmtilde^1)
 @>>> 
 \pi_{\dagger}
 (\jtilde_{\ast}\jtilde^{\ast}\nbigmtilde^1/\nbigmtilde^1)
 @>>> 0
 \end{CD}
\end{equation}
\begin{equation}
\label{eq;08.1.18.51}
 \begin{CD}
 0 @>>>
 i_{Y\dagger}i_Y^{\dagger}\nbigm^1
 @>>>
 j_!j^{\ast}\nbigm^1
 @>>> \nbigm^1
 @>>> 0\\
 @. @V{=}VV @V{g_2}VV @V{g_1}VV @.\\
 0 @>>>
 \pi_{\dagger}\bigl(
 i_{\Ytilde\dagger}
 i_{\Ytilde}^{\dagger}\nbigmtilde^1
 \bigr)
 @>>>
 \pi_{\dagger}\bigl(
 \jtilde_!\jtilde^{\ast}
 \nbigmtilde^1
 \bigr)
 @>>>
 \pi_{\dagger}\bigl(
 \nbigmtilde^1
 \bigr)
 @>>> 0
 \end{CD}
\end{equation}
Here $f_1$ and $g_1$ are the natural morphisms.
\end{lem}
\pf
We will use commutative diagrams
in Section \ref{subsection;08.9.29.210} below.
Because of the non-characteristic condition,
we obtain (\ref{eq;08.1.18.50})
and (\ref{eq;08.1.18.51})
from (\ref{eq;08.1.18.52})
and (\ref{eq;08.1.18.53}) 
respectively,
by taking the tensor product with $\nbigm^1$.
\hfill\qed

\vspace{.1in}
Let us return to the proof of Proposition 
\ref{prop;07.10.28.2}.
We obtain the following commutative diagrams:
\begin{equation}
 \label{eq;08.1.18.70}
 \begin{CD}
 a_{X\dagger}^i
 \bigl(j_{\ast}j^{\ast}\nbigm^1/\nbigm^1\bigr)
 @>>> a_X^{i+1}(\nbigm^1)\\
@AAA @AAA \\
 a_{\Xtilde\dagger}^i\bigl(
 \jtilde_{\ast}\jtilde^{\ast}\nbigmtilde^1/\nbigmtilde^1
 \bigr)
 @>>>
 a_{\Xtilde\dagger}^{i+1}(\nbigmtilde^1)
 \end{CD}
\quad\quad\quad
 \begin{CD}
 a_{X\dagger}^i(\nbigm^1)
 @>>>
 a_{X\dagger}^{i+1}\bigl(
 i_{Y\dagger}i_{Y}^{\dagger}\nbigm^1
 \bigr) \\
 @VVV @VVV \\
 a_{\Xtilde\dagger}^i(\nbigmtilde^1)
 @>>>
 a_{\Xtilde\dagger}^{i+1}\bigl(
 i_{\Ytilde\dagger}
 i_{\Ytilde}^{\dagger}\nbigmtilde^1
 \bigr) 
 \end{CD}
\end{equation}
We also have the following commutative diagram:
\begin{equation}
 \label{eq;08.1.18.71}
 \begin{CD} 
 a_{X\dagger}^0\nbigm^1
 @>>>
 a_{\Xtilde\dagger}^0\nbigmtilde^1
 @>>>
 a_{\Xtilde\dagger}^1\bigl(
 i_{\Ytilde\dagger}i_{\Ytilde}^{\dagger}
 \nbigmtilde^1\bigr)
 \\
 @V{L_c}VV @V{L_c}VV @V{\simeq}VV \\
 a_{X\dagger}^2\nbigm^1
 @>>>
 a_{\Xtilde\dagger}^2\nbigmtilde^1
 @<<<
 a_{\Xtilde\dagger}^1\bigl(
 \jtilde_{\ast}\jtilde^{\ast}\nbigmtilde^1
 \big/\nbigmtilde^1
 \bigr)
 \end{CD}
\end{equation}
The right arrow is the natural isomorphism.
From (\ref{eq;08.1.18.70})
and (\ref{eq;08.1.18.71}),
the morphism
$L_c:a^0_{X\dagger}\nbigm^1
\lrarr
 a^2_{X\dagger}\nbigm^1$
is factorized as follows:
\[
\begin{CD}
 a^0_{X\dagger}\nbigm^1
 @>>>
 a_{\Xtilde\dagger}^0\nbigmtilde^1
 @>>>
 a^1_{\Xtilde\dagger}
 \bigl(i_{\Ytilde\dagger}i_{\Ytilde}^{\dagger}
 \nbigmtilde^1\bigr)
=
 a^1_{X\dagger}
 \bigl(i_{Y\dagger}i_Y^{\dagger}\nbigm^1\bigr)
@>{\simeq}>>
 a^2_{X\dagger}\nbigm^1
\end{CD}
\]
Hence, we obtain the following:
\[
 Pa^0_{X\dagger}\nbigm^1
\subset
 \Ker\bigl(
 \psi:
a_{\Xtilde\dagger}^0\nbigmtilde^1
\lrarr
a^1_{\Xtilde\dagger}
 \bigl(i_{\Ytilde\dagger}i_{\Ytilde}^{\dagger}
 \nbigmtilde^1\bigr)
 \bigr)
\] 
Note 
we have the decomposition
$p^1_{\dagger}\nbigmtilde^1
=N_1\oplus N_2$,
where the strict support of $N_1$ is $\proj^1$
and the support of $N_2$ is $0$-dimensional.
We have $a^j_{\proj^1\dagger}(N_2)=0$
for $j\neq 0$,
and
$f\in a^{-1}_{\proj^1\dagger}(N_1)$
is $0$ if and only if
the restriction of $f$ to some general point is $0$.
Hence,
we obtain the injectivity of the natural morphism
$\Gr^{-1}_La_{\Xtilde\dagger}^0(\nbigmtilde^1)
\simeq
 a_{\proj^1\dagger}^{-1}
p^1_{\dagger}\nbigmtilde^1
\lrarr 
a^1\bigl(i_{\Ytilde\dagger}i_{\Ytilde}^{\dagger}
 \nbigmtilde^1\bigr)$.
Hence, we obtain
$Pa^0_{X\dagger}\nbigm^1
\subset
 L^0a_{\Xtilde\dagger}^0\nbigmtilde^1$.
Then, the claim of Proposition \ref{prop;07.10.28.2}
follows.
\hfill\qed

\vspace{.1in}
Because $\nbigl_c:L^1a^0_{\Xtilde\dagger}\nbigmtilde\lrarr 
 L^1a^2_{\Xtilde\dagger}\nbigmtilde$ is 
an isomorphism,
the following induced morphism is injective:
\[
 Pa^0_{X\dagger}\nbigm\lrarr 
 \Gr_L^0a^0_{\Xtilde\dagger}\nbigmtilde
=\frac{L^0a^0_{\Xtilde\dagger}\nbigmtilde}
 {L^1a^0_{\Xtilde\dagger}\nbigmtilde}
\subset
\frac{a^0_{\Xtilde\dagger}\nbigmtilde}
 {L^1a^0_{\Xtilde\dagger}\nbigmtilde}
\]

Due to  Proposition \ref{prop;07.10.28.2},
we obtain the following morphisms:
\begin{equation}
\label{eq;07.8.30.4}
 Pa^0_{X\dagger}\nbigt\lrarr
 L^0a^0_{\Xtilde\dagger}\nbigttilde
\lrarr
 \Gr^0_La^0_{\Xtilde\dagger}\nbigttilde
\simeq a^0_{\proj^1\dagger}p^0_{\dagger}\nbigttilde
\end{equation}
\begin{equation}
\label{eq;07.10.28.5}
a^0_{\proj^1\dagger}p^0_{\dagger}\nbigttilde^{\ast}
\simeq
 \Gr^0_La^0_{\Xtilde\dagger}\nbigttilde^{\ast}
\lrarr
 \frac{a^0_{\Xtilde\dagger}\nbigttilde^{\ast}}
 {L^1a^0_{\Xtilde\dagger}\nbigttilde^{\ast}}
\lrarr
Pa^0_{X\dagger}\nbigt^{\ast}
\end{equation}
The composite of (\ref{eq;07.8.30.4}) 
is denoted by $F$.
The composite of (\ref{eq;07.10.28.5}) is
the adjoint $F^{\ast}$.
They are morphisms of 
pure twistor structures of weight $0$.
Since $Pa_{X}^0\nbigm\lrarr 
 a^0_{\proj^1}p^0\nbigm$ is injective,
$F$ is a monomorphism.

\vspace{.1in}
We have the primitive decomposition:
\[
 a^0_{\proj^1\dagger}p^0_{\dagger}\nbigttilde
=\bigoplus_{i\geq 0}
 a^0_{\proj^1\dagger}\Bigl(
 \nbigl_c^i\cdot Pp^{-2i}_{\dagger}\nbigttilde
 \Bigr)
\]
\begin{lem}
\label{lem;07.10.28.6}
We have the following factorization of $F$:
\[
 Pa^0_{X\dagger}\nbigt\lrarr
 a^0_{\proj^1\dagger}Pp^0_{\dagger}\nbigttilde
 \lrarr
 a^0_{\proj^1\dagger}p^0_{\dagger}\nbigttilde
\]
\end{lem}
\pf
Because of the twistor property,
we have only to show that
the image of 
$Pa^0_{X\dagger}\nbigm^1\lrarr 
 a^0_{\proj^1\dagger}p^0_{\dagger}\nbigmtilde^1$
is contained in
$a^0_{\proj^1\dagger}Pp^0_{\dagger}\nbigmtilde^1$.
The composite of the following morphisms is $0$:
\[
\begin{CD}
 Pa_{X\dagger}^0\nbigm^1
@>>>
 L^0a_{\Xtilde\dagger}^0\nbigmtilde^1
@>{L_c}>>
 L^0a_{\Xtilde\dagger}^2\nbigmtilde^1
\end{CD}
\]
Hence, the composite of
$Pa_{X\dagger}^0\nbigm^1
\lrarr
 \Gr^0_La_{\Xtilde\dagger}^0\nbigmtilde^1
\lrarr
 \Gr^0_La_{\Xtilde\dagger}^2\nbigmtilde^1$
is $0$.
It implies the above claim.
\hfill\qed

\begin{lem}
The following diagram is commutative:
\begin{equation}
 \label{eq;08.1.18.80}
 \begin{CD}
 Pa^0_{X\dagger}\nbigt @>{F}>> 
 a^0_{\proj^1\dagger}Pp^0_{\dagger}\nbigttilde \\
 @V{\nbigs_0}VV @V{\nbigs_1}VV \\
 Pa^0_{X\dagger}\nbigt^{\ast} @<{F^{\ast}}<<
 a^0_{\proj^1\dagger}Pp^0_{\dagger}\nbigttilde^{\ast}
 \end{CD}
\end{equation}
Here $\nbigs_1$ denote the naturally induced 
polarization of
$a^0_{\proj^1\dagger}Pp^0_{\dagger}\nbigttilde$.
\end{lem}
\pf
We have started from the following commutative diagram:
\[
 \begin{CD}
 Pa_{X\dagger}^0\nbigt
 @>>>
 a_{X\dagger}^0\nbigt
 @>>>
 a_{\Xtilde\dagger}^0\nbigttilde\\
 @VVV @VVV @VVV\\
 Pa_{X\dagger}^0\nbigt^{\ast}
 @<<<
 a_{X\dagger}^0\nbigt^{\ast}
 @<<<
 a_{\Xtilde\dagger}^0\nbigttilde^{\ast}
 \end{CD}
\]
We obtain the commutativity of
the following diagram,
because it is obtained as the factorization
(Proposition \ref{prop;07.10.28.2}):
\[
 \begin{CD}
 Pa_{X\dagger}^0\nbigt
 @>>>
 L^0a_{\Xtilde\dagger}^0\nbigttilde
 @>>>
 a_{\Xtilde\dagger}^0\nbigttilde\\
 @VVV @VVV @VVV\\
 Pa_{X\dagger}^0\nbigt^{\ast}
 @<<<
 {\displaystyle\frac{a_{\Xtilde\dagger}^0\nbigttilde^{\ast}}
 {L^1a_{\Xtilde\dagger}^0\nbigttilde^{\ast}}}
 @<<<
 a_{\Xtilde\dagger}^0\nbigttilde^{\ast}
 \end{CD}
\]
Then, we obtain the following commutative diagram:
\begin{equation}
 \label{eq;08.1.18.81}
 \begin{CD}
 Pa_{X\dagger}^0\nbigt
 @>>>
 L^0a_{\Xtilde\dagger}^0\nbigttilde
 @>>>
 \Gr_L^0a_{\Xtilde\dagger}^0\nbigttilde
 \\
 @VVV @VVV @VVV\\
 Pa_{X\dagger}^0\nbigt^{\ast}
 @<<<
 {\displaystyle\frac{a_{\Xtilde\dagger}^0\nbigttilde^{\ast}}
 {L^1a_{\Xtilde\dagger}^0\nbigttilde^{\ast}}}
 @<<<
 \Gr_L^0a_{\Xtilde\dagger}^0\nbigttilde^{\ast}
 \end{CD}
\end{equation}
Since (\ref{eq;08.1.18.80})
is obtained as the factorization of (\ref{eq;08.1.18.81})
in Lemma \ref{lem;07.10.28.6},
we obtain the desired commutativity.
\hfill\qed

\vspace{.1in}

Since $F$ is monomorphic,
$(Pa^0_{X\dagger}\nbigt,\nbigs_0)$ is polarized,
and hence Proposition \ref{prop;07.10.28.10}
is proved.
\hfill\qed

\chapter{Correspondences}
\label{section;08.9.29.150}
In this chapter,
we study some correspondences
and their application to a conjecture of Kashiwara.
In Sections 
\ref{subsection;08.10.15.1}--\ref{subsection;07.10.28.35},
we establish the correspondence between
wild harmonic bundles and
polarized wild pure twistor $D$-modules
on complex manifolds
(Theorem \ref{thm;07.10.28.30}).
The argument is essentially the same
as that in the tame case in \cite{mochi2},
although we have some additional difficulties
caused by Stokes structure and ramification.
In Section \ref{subsection;08.10.15.2},
we show the correspondence
between semisimple algebraic holonomic $D$-modules
and polarizable wild pure twistor $D$-modules 
on projective varieties
(Theorem \ref{thm;07.10.14.75}).
As an easy consequence,
we obtain Kashiwara's conjecture
(Theorem \ref{thm;08.3.6.2}).

\section{Wild harmonic bundles and
 wild pure twistor $D$-modules}
\label{subsection;08.10.15.1}
\subsection{Preliminary}

Let $X$ be a complex manifold,
and $Z$ be any closed irreducible subvariety of $X$.
Let $\nbiga$ be a $\rnum$-vector subspace of $\cnum$.
In the following,
``Zariski open'' means
``the complement 
of some closed analytic subset''.
\index{Zariski open}
Let $U$ be a smooth Zariski open subset of $Z$.
Recall that
a harmonic bundle $\harmonicbundle$ on $U$
is called $\nbiga$-wild on $(Z,U)$,
if we have a complex manifold $\Ztilde$
and a birational projective morphism
$\varphi:\Ztilde\lrarr Z$ satisfying
the following property:
\begin{itemize}
\item
$\Dtilde:=\Ztilde\setminus\varphi^{-1}(U)$
is a normal crossing divisor.
\item
$\varphi^{\ast}\harmonicbundle_{|U}$
is an $\nbiga$-wild harmonic bundle on $(\Ztilde,\Dtilde)$.
\end{itemize}

\begin{df}
\label{df;08.1.17.1}
Let $(V,\DD^{\sankaku},S)$ be 
a variation of polarized pure
twistor structure of weight $w$ defined on $U$.
We say that $(V,\DD^{\sankaku},S)$ is 
an $\nbiga$-wild variation of polarized
pure twistor structure of weight $w$ 
on $(Z,U)$,
if the underlying harmonic bundle
is $\nbiga$-wild on $(Z,U)$.
\index{wild variation of polarized
pure twistor structure}
\hfill\qed
\end{df}

Let $(\nbigt,\nbigs)$ be a polarized $\nbiga$-wild
pure twistor $D$-module of weight $w$
whose strict support is $Z$.
As remarked in Lemma \ref{lem;10.7.2.1},
there exists a Zariski open subset $U$ of $Z$ such that
$(\nbigt,\nbigs)_{|X\setminus Y}$ comes from 
a harmonic bundle $\harmonicbundle$ on $U$
with the Tate twist, 
where $Y:=Z\setminus U$.
We will prove  the following lemma
in Section \ref{subsection;07.10.28.35}.
\begin{lem}
\label{lem;10.5.25.1}
$\harmonicbundle$ is 
$\nbiga$-wild on $(Z,U)$.
\end{lem}

\subsection{Statement}
\label{subsection;10.6.6.13}
Let $Z$ be an irreducible closed subset of $X$,
and let $U$ be a Zariski open subset of $Z$,
and $Y:=Z\setminus U$.
Let $\VPTgenwild(Z,U,w,\nbiga)$ denote the category
of $\nbiga$-wild variations
of polarized pure twistor structure of weight $w$ on 
$(Z,U)$.
Let $\MPT^{\wild}_{\strict}(Z,U,w,\nbiga)$ denote 
the category 
of polarized $\nbiga$-wild pure twistor $D$-modules
such that
(i) their strict supports are $Z$,
 i.e., the support of any non-zero direct summand is $Z$,
(ii) their restriction to $X\setminus Y$ come from 
variations of pure polarized twistor structure of weight $w$
on $U$.
In the both categories,
we consider only isomorphisms.
\index{category $\VPTgenwild(Z,U,w,\nbiga)$}
\index{category $\MPT^{\wild}_{\strict}(Z,U,w,\nbiga)$}
By Lemma \ref{lem;10.5.25.1},
we obtain a functor
\[
 \Phi: \MPT^{\wild}_{\strict}(Z,U,w,\nbiga)
 \lrarr \VPTgenwild(Z,U,w,\nbiga).
\]

The following theorem is one of the main results
in this monograph.
\begin{thm} \mbox{{}}
 \label{thm;07.10.28.30}
$\Phi$ is an equivalence.
\end{thm}
Note that this theorem implies the following:
\begin{itemize}
\item
Let $(V,\DD^{\sankaku},S)$
be an $\nbiga$-wild variation of
polarized pure twistor structure of weight $w$
on $(Z,U)$.
Then, 
$(\nbigt,\nbigs)\in
 \MPT^{\wild}_{\strict}(Z,U,w,\nbiga)$
such that
$(\nbigt,\nbigs)_{|X\setminus Y}$
 comes from $(V,\DD^{\sankaku},S)$.
In other words,
$(V,\DD^{\sankaku},S)$
is extended to
a polarized wild pure twistor $D$-module on $X$.
It is called a minimal extension of
$(V,\DD^{\sankaku},S)$.
\item
 If $(\nbigt',\nbigs')$ is another minimal extension
 of $(V,\DD^{\sankaku},S)$,
 we have $(\nbigt,\nbigs)\simeq
 (\nbigt',\nbigs')$,
 where the restriction of the isomorphism to
 $X\setminus Y$ is the identity of
 $(V,\DD^{\sankaku},S)$.
 Namely, a minimal extension is unique.
\item
 In particular,
 if $(V,\DD^{\sankaku},S)=\bigoplus
 (V_i,\DD^{\sankaku}_i,S_i)$,
 we have the corresponding decomposition
 $(\nbigt,\nbigs)=\bigoplus(\nbigt_i,\nbigs_i)$.
\end{itemize}

We will show the essential surjectivity of
$\Phi$ in Section \ref{subsection;07.10.28.31},
and the fully faithfulness in Section
\ref{subsection;07.10.28.35}.

\subsubsection{Variant}
\label{subsection;10.6.6.12}
Let $Z,U,Y$ be as above.
Let $(V,\DD^{\sankaku})$ be a variation of pure
twistor structure of weight $w$ defined on 
a smooth Zariski open subset $U$ of $Z$.
We say that $(V,\DD^{\sankaku})$ is 
a {\em polarizable} $\nbiga$-wild variation of 
pure twistor structure of weight $w$ 
on $(Z,U)$,
if there exists a polarization $S$
of $(V,\DD^{\sankaku})$
such that $(V,\DD^{\sankaku},S)$ 
is an $\nbiga$-wild variation of {\em polarized}
pure twistor structure of weight $w$ 
on $(Z,U)$.
\index{polarizable wild variation of
pure twistor structure}
Let $\VPTgenwild(Z,U,w,\nbiga)'$ denote the category
of the {\em polarizable} $\nbiga$-wild variations
of pure twistor structure of weight $w$ on $(Z,U)$.
Morphisms in this category are defined to be
morphisms for variations of twistor structure
on $U$.

Let $\MPT^{\wild}_{\strict}(Z,U,w,\nbiga)'$ denote
the category of {\em polarizable}
$\nbiga$-wild pure twistor $D$-modules of weight $w$
such that 
(i) their strict supports are $Z$,
(ii) their restriction to 
$X\setminus Y$ come from
variations of twistor structure of weight $w$.
(We consider the full subcategory of
the category of wild pure twistor $D$-modules 
of weight $w$ on $X$.)
\index{category $\VPTgenwild(Z,U,w,\nbiga)'$}
\index{category $\MPT^{\wild}_{\strict}(Z,U,w,\nbiga)'$}
By Lemma \ref{lem;10.5.25.1},
we obtain the following functor:
\[
 \Phi:
\MPT^{\wild}_{\strict}(Z,U,w,\nbiga)'
\lrarr
 \VPTgenwild(Z,U,w,\nbiga)' 
\]
\begin{cor}
\label{cor;10.6.6.20}
$\Phi$ is an equivalence.
\end{cor}
\pf
The essential surjectivity follows from
Theorem \ref{thm;07.10.28.30}.
Let us consider the fully faithfulness.
Note that both the categories
$\MPT^{\wild}_{\strict}(Z,U,w,\nbiga)'$
and $\VPTgen^{\wild}(Z,U,w,\nbiga)$
are semisimple.
Hence, we have only to show the following:
\begin{itemize}
\item
$\Phi(\nbigt)$ is simple if 
$\nbigt\in\MPT^{\wild}_{\strict}(Z,U,w,\nbiga)$ is simple.
\item
Let $\nbigt_i\in \MPT^{\wild}_{\strict}(Z,U,w,\nbiga)$
$(i=1,2)$ be simple.
If $\Phi(\nbigt_1)\simeq\Phi(\nbigt_2)$,
then $\nbigt_1\simeq\nbigt_2$.
\end{itemize}
The both claims follow from 
Theorem \ref{thm;07.10.28.30}.
\hfill\qed

\section{Prolongment to 
 polarized wild pure twistor $D$-modules}
\label{subsection;07.10.28.31}
We shall prove the essential surjectivity of $\Phi$
in Theorem \ref{thm;07.10.28.30}.
The argument is essentially the same
as that in the tame case \cite{mochi2}.

In Subsection \ref{subsection;10.5.25.2},
we will reduce the issue to the local and
unramified case,
in which we have already constructed
an $\nbigr$-triple in Chapter
\ref{section;07.10.12.10}.
We have only to show that it is
a wild polarized pure twistor $D$-module
(Proposition \ref{prop;07.8.30.5}).
In Subsection \ref{subsection;08.1.21.2},
we give a preparation for the functoriality
of the family of 
filtered $\lambda$-flat bundles
associated to 
an unramifiedly good wild harmonic bundle.
In Subsection \ref{subsection;07.12.3.15},
we shall show that the $\nbigr$-triple
in Chapter \ref{section;07.10.12.10}
is strictly $S$-decomposable
along any function.
In Subsection \ref{subsection;10.5.25.3},
we study the specialization along a monomial
function with an exponential twist.
Then, in Subsection \ref{subsection;07.10.28.42},
we argue the specialization along any function
with exponential twist.

Recall that we have already studied
the specialization of the $\nbigr$-triple
along monomial functions
without exponential twist,
in Sections \ref{subsection;08.2.2.1}
and \ref{subsection;08.9.28.10}.
Namely, we have already known that
it is strictly $S$-decomposable
along any monomial function
(Proposition \ref{prop;07.10.28.50}),
and we have also known
that we obtain polarized graded wild Lefschetz
twistor $D$-module as the specialization
along a monomial function
(Proposition \ref{prop;07.9.8.26}.
More precisely,
we use the hypothesis of the induction
on the dimension of the strict support.)
The specialization along any functions
can be reduced to that along monomial functions
by Hironaka's resolution and the propositions
in Subsection \ref{subsection;05.1.15.13}.

\subsection{Reduction to the 
local and unramified case}
\label{subsection;10.5.25.2}

Let $X=\Delta^n_z$ and $D=\bigcup_{i=1}^{\ell} D_i$.
Let $\harmonicbundle$ be an unramifiedly $\nbiga$-good 
wild harmonic bundle on $X-D$.
We have constructed the $\nbigr_X$-module $\gbige$
with the sesqui-linear pairing
$\gbigc:\gbige_{|\vecS\times X}
 \otimes\sigma^{\ast}\gbige_{|\vecS\times X}
\lrarr\distribution_{\vecS\times X/\vecS}$.
(See Sections \ref{subsection;08.9.28.1}
and \ref{subsection;08.9.28.11}.)
We will prove the following proposition in Sections
\ref{subsection;08.1.21.2}--\ref{subsection;07.10.28.42}.

\begin{prop}
\label{prop;07.8.30.5}
$\gbigt(E):=(\gbige,\gbige,\gbigc)$ is 
an $\nbiga$-wild pure twistor $D$-module of weight $0$.
The natural Hermitian sesqui-linear duality
$\nbigs=(\id,\id)$ gives a polarization.
\end{prop}

Let $X$ and $D$ be as above.
Let $\harmonicbundle$ be an $\nbiga$-good wild harmonic bundle,
which is not necessarily unramified.
We show the next lemma 
by assuming Proposition \ref{prop;07.8.30.5}.
\begin{lem}
\label{lem;07.10.28.40}
There exists 
an $\nbiga$-wild pure twistor $D$-module 
$(\gbige_1,\gbige_1,\gbigc_1)$ of weight $0$
with the natural polarization
$(\id,\id)$ satisfying the following:
\begin{itemize}
\item
The restriction to $X-D$ is isomorphic to 
the polarized pure twistor $D$-module
associated to $\harmonicbundle$.
\item
Let $\nbigx:=\cnum_{\lambda}\times X$
and $\nbigd:=\cnum_{\lambda}\times D$.
Then,
$\gbige_1\otimes\nbigo_{\nbigx}(\ast \nbigd)$
is naturally isomorphic to $\nbigq\nbige$.
In particular,
we naturally have
$\gbige_1\subset\nbigq\nbige$.
\end{itemize}
\end{lem}
\pf
Take a ramified covering 
$\varphi:(\Xtilde,\Dtilde)\lrarr (X,D)$
such that $(\Etilde,\delbar_{\Etilde},\thetatilde,\htilde):=
 \varphi^{-1}\harmonicbundle$ is unramified.
We obtain the associated polarized 
$\nbiga$-wild pure twistor $D$-module 
$\gbigt(\Etilde)=\bigl(\gbigetilde,\gbigetilde,\gbigctilde\bigr)$
with $\gbigstilde=(\id,\id)$
as in Proposition \ref{prop;07.8.30.5},
which is $\Gal(\Xtilde/X)$-equivariant.
Hence, we obtain a polarized $\nbiga$-wild
pure twistor $D$-module
$P\varphi^0_{\dagger}\gbigt(\Etilde)$ on $X$
with the $\Gal(\Xtilde/X)$-action.
The invariant part of
$P\varphi^0_{\dagger}\gbigt(\Etilde)_{|X-D}$
is isomorphic to the polarized pure twistor $D$-module
associated to $\harmonicbundle$.
Note that a wild pure twistor $D$-module
has the decomposition by strict supports,
because it is assumed to be strict $S$-decomposable
in definition.
We have the direct summand 
$\gbigt_1=(\gbige_1,\gbige_1,\gbigc_1)$ of
the invariant part
$P\varphi_{\dagger}\gbigt(\Etilde)^{\Gal(\Xtilde/X)}$
such that the strict support of $\gbigt_1$ is $X$,
with the naturally induced polarization $\gbigs_1=(\id,\id)$.
Then, $(\gbigt_1,\gbigs_1)$ gives a prolongment
of $(E,\delbar_E,\theta,h)$.

We have $\varphi^{\ast}\nbigq\nbige
=\nbigq\nbigetilde$,
and hence
$\nbigq\nbige$
is naturally isomorphic to
the $\Gal(\Xtilde/X)$-invariant part
of $\varphi_{\dagger}\nbigq\nbigetilde$.
We also have
$\gbigetilde\otimes\nbigo_{\nbigxtilde}(\ast\nbigdtilde)
=\nbigq\nbigetilde$.
Hence,
we obtain a natural morphism
$\gbige_1\lrarr \nbigq\nbige$
satisfying 
$\gbige_1\otimes\nbigo_{\nbigx}(\ast\nbigd)
\simeq\nbigq\nbige$.
\hfill\qed

\vspace{.1in}
Let $X$ be a complex manifold
with a normal crossing hypersurface $D$.
Let $\harmonicbundle$ be 
an $\nbiga$-good wild harmonic bundle
on $(X,D)$.
We show the next lemma
by assuming Proposition \ref{prop;07.8.30.5}
and hence Lemma \ref{lem;07.10.28.40}.
\begin{lem}
\label{lem;07.12.3.10}
There exists 
an $\nbiga$-wild pure twistor $D$-module 
$(\gbige_1,\gbige_1,\gbigc_1)$ of weight $0$
with the polarization $(\id,\id)$,
whose restriction to $X-D$ is 
isomorphic to the polarized pure twistor $D$-module
associated to $\harmonicbundle$.
\end{lem}
\pf
Let $U$ be an open subset of $X$
with a holomorphic coordinate system 
$(z_1,\ldots,z_n)$
such that $D\cap U=\bigcup_{j=1}^{\ell}\{z_j=0\}$.
By applying Lemma \ref{lem;07.10.28.40},
we obtain the wild pure twistor $D$-module
$\gbigt_U:=
 (\gbige_U,\gbige_U,\gbigc_U)$ with the polarization
$(\id,\id)$.
Let $U_i$ $(i=1,2)$ be such two open sets.
Then, the restriction of
$\gbige_{U_i}\otimes \nbigo(\ast\nbigd)$
to $\cnum_{\lambda}\times (U_1\cap U_2)$
are naturally isomorphic to 
$\nbigq\nbige_{
 |\cnum_{\lambda}\times (U_1\cap U_2)}$.
Let $g$ be a holomorphic function on $U_1\cap U_2$
such that $U_1\cap U_2\cap D=\{g=0\}$.
Since both $\gbige_{U_i|U_1\cap U_2}$ are strictly 
$S$-decomposable along $g$,
the isomorphism on $(U_1\cap U_2)\setminus D$
is extended to the isomorphism on $U_1\cap U_2$.
(See Lemma \ref{lem;07.11.2.15}.)
Hence, 
by varying $U$ and gluing $\gbigt_U$,
we obtain the globally defined
$\nbigr_X$-triple 
$\gbigt(E):=(\gbige_1,\gbige_1,\gbigc_1)$
with the Hermitian sesqui-linear duality
$\gbigs(E):=(\id,\id)$.
Let us check that $\bigl(\gbigt(E),\gbigs(E)\bigr)$
is a polarized wild pure twistor $D$-module 
of weight $0$.
According to Lemma \ref{lem;08.12.4.2},
we have only to show that
$\gbige_1$ is a good $\nbigr_X$-module.
(See Subsection \ref{subsection;08.12.4.1}.)

Note that $\gbige_1$ is naturally contained
in $\nbigq\nbige$.
Let $K$ be any compact subset of $X$.
For any $\lambda_0\in\cnum_{\lambda}$,
there exist a large number $N$
and a neighbourhood $\nbigu$ of
$\{\lambda_0\}\times K$ in $\nbigx$,
such that an $\nbigo_{\nbigu}$-module
$\gbige_1\cap \nbigqzero_{N\vecdelta}\nbige$
generates $\gbige_1$.
Hence, $\gbige_1$ is good.
Thus, the proof of Lemma \ref{lem;07.12.3.10}
is completed.
\hfill\qed

\vspace{.1in}

Let us show the essential surjectivity of $\Phi$
assuming Proposition \ref{prop;07.8.30.5},
and hence Lemma \ref{lem;07.12.3.10}.
Let $\bigl(\harmonicbundle\bigr)$ be 
an $\nbiga$-wild harmonic bundle on $(Z,U)$.
We take a complex manifold $\Ztilde$
and a birational projective morphism
$\varphi:\Ztilde\lrarr Z$
such that (i) the complement of $\varphi^{-1}(U)$ 
is normal crossing,
(ii) $\varphi^{-1}\harmonicbundle$
is a good $\nbiga$-wild harmonic bundle.
By applying Lemma \ref{lem;07.12.3.10},
we take a polarized $\nbiga$-wild pure twistor $D$-module
$(\gbigt_1,\gbigs_1)$ on $\Ztilde$
whose restriction to $\varphi^{-1}(U)$
is isomorphic to the polarized pure twistor $D$-module
associated to $\varphi^{-1}\harmonicbundle$.
According to Theorem \ref{thm;07.10.23.20},
$\bigl(
 \bigoplus_i\varphi^i_{\dagger}\gbigt_1,
 \nbigl_c,\bigoplus_i\varphi^i_{\dagger}\gbigs_1
\bigr)$
is a polarized graded $\nbiga$-wild 
Lefschetz twistor $D$-module
on $X$,
where $\nbigl_c$ denotes a Lefschetz map
associated to a line bundle relatively ample with respect 
to $\varphi$.
Let $P\varphi^0_{\dagger}\gbigt_1$
be the primitive part of $\varphi^0_{\dagger}\gbigt_1$.
We have the decomposition
$P\varphi^0_{\dagger}\gbigt_1
=\gbigt_2\oplus \gbigt_2'$,
where the strict support of $\gbigt_2$ is $Z$,
and the support of $\gbigt_2'$ is strictly smaller than $Z$.
We have the naturally induced polarization $\gbigs_2$
on $\gbigt_2$.
Then, $(\gbigt_2,\gbigs_2)$
gives the desired prolongment of $\harmonicbundle$.
Thus, the proof of 
the essential surjectivity of $\Phi$
is reduced to that of Proposition \ref{prop;07.8.30.5}.

\vspace{.1in}

In the rest of this section,
we will prove Proposition \ref{prop;07.8.30.5}.
We use an induction on $\dim X$.

\subsection{Pull back of the associated family of
 meromorphic $\lambda$-flat bundles}
\label{subsection;08.1.21.2}

Let $X:=\Delta^n_z$ and 
$D:=\bigcup_{i=1}^{\ell} D_i$.
Let $\harmonicbundle$ be 
an unramifiedly $\nbiga$-good wild harmonic bundle
on $X-D$.
Let $\varphi:X_1\lrarr X$ be 
a birational projective morphism
such that $D_1:=\varphi^{-1}(D)$ is 
a normal crossing divisor.
\begin{lem}
$(E_1,\delbar_{E_1},\theta_1,h_1):=
\varphi^{-1}\harmonicbundle$
is an unramifiedly 
$\nbiga$-good wild harmonic bundle on $X_1-D_1$.
\end{lem}
\pf
Because the claim is local on $X$,
we may assume to have the decomposition
$(E,\delbar_E,\theta)=\bigoplus
 (E_{\gminia},\delbar_{E_{\gminia}},\theta_{\gminia})
 \otimes L(\gminia)$.
Then, the claim is clear.
\hfill\qed

\vspace{.1in}
Let $\nbigp_{<0}\nbigelambda$
be the union of
$\nbigp_{\vecb}\nbigelambda$
for $\vecb\in\real^{\ell}$
such that $b_i<0$ $(i=1,\ldots,\ell)$.

\begin{lem}
We have 
$\varphi^{\ast}\nbigp\nbigelambda
 =\nbigp\nbigelambda_1$
and
$\varphi^{\ast}\nbigp_{<0}\nbigelambda
\subset
 \nbigp_{<0}\nbigelambda_1$.
\end{lem}
\pf
Because 
$\varphi^{\ast}\nbigp_{<0}\nbigelambda
=\nbigo_{X_1}\otimes_{\varphi^{-1}\nbigo_X}
 \varphi^{-1}\nbigp_{<0}\nbigelambda$,
we have only to show
$\varphi^{-1}\nbigp_{<0}\nbigelambda
\subset
 \nbigp_{<0}\nbigelambda_1$,
which follows from the definition.
(See Subsection
 \ref{subsection;07.11.24.20}.)
\hfill\qed

\vspace{.1in}
Let $\nbigq_{<0}\nbigelambda$
be the union of
$\nbigq_{\vecb}\nbigelambda$
for $\vecb\in\real^{\ell}$
such that $b_i<0$ $(i=1,\ldots,\ell)$.
We obtain the following lemma
from Lemma \ref{lem;10.5.25.31}.
\begin{lem}
\label{lem;07.10.28.45}
We have 
$\varphi^{\ast}\nbigq_{<0}\nbigelambda
\subset \nbigq_{<0}\nbigelambda_1$
and 
$\varphi^{\ast}\nbigq\nbigelambda
= \nbigq\nbigelambda_1$.
\hfill\qed
\end{lem}

It is easy to derive the following lemma
from Lemma \ref{lem;07.10.28.45}
by using the local freeness
of $\nbigqzero\nbige_1$
and $\nbigqzero_{<0}\nbige_1$.
\begin{lem}
\label{lem;07.9.6.1}
For any $\lambda_0$,
we have
$\varphi^{\ast}\nbigq^{(\lambda_0)}_{<0}\nbige
\subset
 \nbigq^{(\lambda_0)}_{<0}\nbige_1$
and 
$\varphi^{\ast}\nbigq^{(\lambda_0)}\nbige
=\nbigq^{(\lambda_0)}\nbige_1$.
\hfill\qed
\end{lem}

\subsection{Strict $S$-decomposability}
\label{subsection;07.12.3.15}

Let $X$, $D$ and $\harmonicbundle$
be as above.
Let $g$ be any holomorphic function on $X$.
Let us show that
$\gbige$ is strictly $S$-decomposable along $g$.
Let $\Xtilde$ be a complex manifold
with a projective birational morphism
$\varphi:\Xtilde\lrarr X$ such that
(i) $\Xtilde-\varphi^{-1}\bigl(g^{-1}(0)\cup D\bigr)
\simeq
 X-\bigl(g^{-1}(0)\cup D\bigr)$,
(ii) $\varphi^{-1}\bigl(g^{-1}(0)\cup D\bigr)$ is 
a normal crossing divisor.
We obtain the $\nbiga$-good wild harmonic bundle
$(\Etilde,\delbar_{\Etilde},\thetatilde,\htilde):=
 \varphi^{-1}(E,\delbar_E,\theta,h)$
and the associated $\nbigr_{\Xtilde}$-triple
$\gbigt(\Etilde):=(\gbigetilde,\gbigetilde,\gbigctilde)$.
Let $g_1:=g^{a}\cdot \prod_{j=1}^{\ell} z_j^{n_j}$
for some $a,n_j\in\seisuu_{\geq \,0}$.
We put $\gtilde_1:=g_1\circ\varphi$.
\begin{itemize}
\item
By Proposition \ref{prop;07.10.28.50},
$\gbigt(\Etilde)$ is strictly $S$-decomposable
along $\gtilde_1$.
\item
For any $u\in\real\times\cnum$,
according to Proposition \ref{prop;07.9.8.26}
and the hypothesis of the induction on $\dim X$,
$\bigoplus_{l}\Gr^W_l\psitilde_{\gtilde_1,u}(\gbigt(\Etilde))$
with the induced morphism $\nbign$
and Hermitian sesqui-linear duality
is a polarized graded $\nbiga$-wild 
Lefschetz twistor $D$-module of weight $0$.
\end{itemize}
Namely, $\gbigt(\Etilde)$ and
the Hermitian sesqui-linear duality $(\id,\id)$
satisfy Condition \ref{condition;05.1.15.30}
with $\gminia=0$.
According to Proposition \ref{prop;05.1.15.60},
$\varphi^i_{\dagger}\gbigt(\Etilde)$ is
strictly $S$-decomposable along $g_1$.
Hence, we have the decomposition
$\varphi^0_{\dagger}\gbigetilde
=\nbigm_1\oplus\nbigm_2$ such that 
(i) $\nbigm_1$ contains
no non-trivial $\nbigr$-submodules whose supports
are contained in $D\cup g^{-1}(0)$,
(ii) the support of $\nbigm_2$ is contained in 
$g^{-1}(0)\cup D$.
By construction,
we have the natural identification 
of the restrictions of
$\gbige$ and $\nbigm_1$
to $\cnum_{\lambda}\times
 \bigl(X-\bigl(D\cup\{g=0\}\bigr)\bigr)$.

\begin{lem}
\label{lem;08.10.14.1}
The above identification
is naturally extended to a morphism 
$\kappa:\gbige\lrarr \nbigm_1$ on 
 $\cnum_{\lambda}\times X$.
\end{lem}
\pf
It can be shown using an argument
in Section 19.3.1 of \cite{mochi2}.
We give only an outline.
We consider the corresponding right $\nbigr_X$-modules,
$\gbige\otimes\omega_{\nbigx}$
and $\gbigetilde\otimes\omega_{\nbigxtilde}$.
By Lemma \ref{lem;07.9.6.1},
any section of
$\nbigq^{(\lambda_0)}_{<\vecdelta}\nbige
 \otimes\omega_{\nbigx}$
naturally induces the section of
$\nbigq^{(\lambda_0)}_{<\vecdelta}\nbigetilde
\otimes\omega_{\nbigxtilde}\subset
\gbigetilde\otimes\omega_{\nbigxtilde}$.
Hence, it naturally induces the section of
$\varphi^0_{\dagger}
 \bigl(\gbigetilde\otimes\omega_{\nbigxtilde}\bigr)$.
(See Lemma 14.11 of \cite{mochi2}.)
Hence, we obtain 
$F:\nbigqzero_{<\vecdelta}(\nbige)
\otimes\omega_{\nbigx}
\lrarr  \nbigm_1\otimes\omega_{\nbigx}$.

Let $f_i$ $(i=1,2)$ be local sections of 
$\nbigqzero_{<\vecdelta}(\nbige)
 \otimes\omega_{\nbigx}$,
and let $P_i$ $(i=1,2)$ be local sections of
$\nbigr_X$.
If $f_1\cdot P_1=f_2\cdot P_2$ in
$\gbige\otimes\omega_{\nbigx}$,
the restrictions of
$F(f_i)\cdot P_i$ to 
$\cnum_{\lambda}\times\bigl(
 X-\bigl(D\cup g^{-1}(0)\bigr)\bigr)$
are the same.
By using the strict $S$-decomposability of $\nbigm_1$
along $g\cdot\prod_{j=1}^{\ell}z_j$,
we obtain
$F(f_1)\cdot P_1=F(f_2)\cdot P_2$ on $X$.

Since $\gbige\otimes\omega_{\nbigx}$
is generated by 
$\nbigq^{(\lambda_0)}_{<\vecdelta}\nbige
 \otimes\omega_{\nbigx}$
over $\nbigr_X$,
we obtain the desired map.
\hfill\qed

\vspace{.1in}

Let us show that
$\kappa$ is an isomorphism.
Since $\gbige_{|\cnum_{\lambda}\times(X-D)}$ 
comes from a harmonic bundle,
we know that $\gbige_{|\cnum_{\lambda}\times (X-D)}$
is strictly $S$-decomposable along $g_{|X-D}$.
(See \cite{sabbah2} or \cite{mochi2}.)
Hence, $\kappa_{|X-D}$ is an isomorphism,
due to Lemma \ref{lem;07.11.2.15}.
We have already known that
$\gbige$ is strictly $S$-decomposable
$\prod_{j=1}^{\ell} z_j$.
Hence, $\kappa$ is an isomorphism,
due to Lemma \ref{lem;07.11.2.15} again.
In particular,
we obtain that $\gbige$ is strictly $S$-decomposable
along $g$.

\subsection{Specialization along monomial functions
with exponential twist}
\label{subsection;10.5.25.3}

We study the specialization along
monomial functions
before considering the general case.

\begin{prop}
\label{prop;07.10.28.72}
\mbox{{}}
Let $g:=\prod_{i=1}^{k}z_i^{m_i}$.
\begin{itemize}
\item
$\gbigt(E)$ is strictly specializable along $g$
with ramification and exponential twist by 
any $\gminia\in\cnum[t_m^{-1}]$
\item
$P\Gr^{W(N)}_l\psitilde_{g,\gminia,u}\gbigt(E)$
are $\nbiga$-wild pure twistor $D$-modules of weight $l$
with the naturally induced polarization
for any $u\in\real\times\cnum$
and $\gminia\in\cnum[t_m^{-1}]$.
\item
Moreover,
$P\Gr^{W(N)}_l\psitilde_{g,\gminia,u}\gbigt(E)=0$
unless $u\in\real\times\nbiga$.
\end{itemize}
\end{prop}
\pf
We first consider the unramified case,
i.e., $\gminia\in\cnum[t^{-1}]$,
and then we argue the general case.

\subsubsection{Unramified case}

Note that 
$(E,\delbar_E,\theta,h)\otimes
 L(-g^{\ast}\gminia)$
is not necessarily good.
According to Proposition \ref{prop;07.11.1.3},
we can take a birational projective morphism
$(\Xtilde,\Dtilde)\lrarr (X,D)$
such that 
(i) $\Xtilde-(\varphi^{\ast}g)^{-1}(0)
 \simeq X-g^{-1}(0)$,
(ii) $\varphi^{-1}\bigl(
 (E,\delbar_E,\theta,h)
 \otimes L(-g^{\ast}\gminia)\bigr)$
is an unramifiedly $\nbiga$-good wild harmonic bundle
on $(\Xtilde,\Dtilde)$.
We set $\gtilde:=\varphi^{\ast}g$.
We set
\[
 (\Etilde,\delbar_{\Etilde},\thetatilde,\htilde):=
 \varphi^{\ast}(E,\delbar_E,\theta,h),
\quad
 (\Etilde',\delbar_{\Etilde'},\thetatilde',\htilde'):=
 (\Etilde,\delbar_{\Etilde},\thetatilde,\htilde)
 \otimes
 L\bigl(-\gtilde^{\ast}\gminia\bigr).
\]
We have the associated $\nbigr_{\Xtilde}$-triples
$\nbigt(\Etilde):=
 (\gbigetilde,\gbigetilde,\gbigctilde)$
and
$\nbigt(\Etilde'):=
 (\gbigetilde',\gbigetilde',\gbigctilde')$
with the Hermitian sesqui-linear dualities $(\id,\id)$.
Let $i_g$ (resp. $i_{\gtilde}$) denote 
the graph $X\lrarr X\times\cnum_t$
(resp. $\Xtilde\lrarr \Xtilde\times\cnum_t$)
for the function $g$ (resp. $\gtilde$).
Since $(\Etilde',\delbar_{\Etilde'},\thetatilde',\htilde')$
is also unramifiedly $\nbiga$-good wild,
the following holds:
\begin{itemize}
\item
The $\nbigr_{\Xtilde\times\cnum_t}$-triple
$i_{\gtilde\dagger}\gbigt(\Etilde')$
is strictly $S$-decomposable along $t$.
\item
According to Proposition \ref{prop;07.9.8.26}
and the hypothesis of the induction on $\dim X$,
$\bigoplus_{l}
 \Gr^W_l\psitilde_{\gtilde,u}
 \bigl(\gbigt(\Etilde')\bigr)$
with the induced morphism $\nbign$
and Hermitian sesqui-linear duality
is a polarized graded $\nbiga$-wild 
Lefschetz twistor $D$-module of weight $0$
for each $u\in\real\times\cnum$.
Moreover, 
$\psitilde_{\gtilde,u}
 \bigl(\gbigt(\Etilde')\bigr)=0$
unless $u\in\real\times\nbiga$,
by Corollary \ref{cor;07.11.2.40}.
\end{itemize}
Because
$i_{\gtilde\dagger}\gbigt(\Etilde')(\ast t)
\simeq
 i_{\gtilde\dagger}\gbigt(\Etilde)(\ast t)
\otimes\gbigl(-\gtilde^{\ast}\gminia)$,
we obtain the following:
\begin{itemize}
\item
The $\nbigr_{\Xtilde\times\cnum_t}(\ast t)$-triple
$i_{\gtilde\dagger}\gbigt(\Etilde)(\ast t)
\otimes\gbigl(-\gtilde^{\ast}\gminia)$
is strictly specializable along $t$.
\item
$\bigoplus_{l}
 \Gr^W_l\psitilde_{\gtilde,u,\gminia}
 (\gbigt(\Etilde))$
with the induced morphism $\nbign$
and Hermitian sesqui-linear duality
is a polarized graded $\nbiga$-wild 
Lefschetz twistor $D$-module of weight $0$
for each $u\in\real\times\cnum$.
Moreover,
$\psitilde_{\gtilde,u,\gminia}
 (\gbigt(\Etilde))=0$
unless $u\in\real\times\nbiga$.
\end{itemize}
Namely,
$\gbigt(\Etilde)$ 
and the Hermitian sesqui-linear duality $(\id,\id)$
satisfy Condition \ref{condition;05.1.15.30}
with $\gminia$.
We also have
$i_{g\dagger}
 \varphi_{\dagger}^i\gbigt(\Etilde)(\ast t)
 \otimes \gbigl(-g^{\ast}\gminia)=0$
unless $i=0$.
According to Lemma \ref{lem;08.1.16.40}
and Proposition \ref{prop;05.1.15.61},
the following holds:
\begin{itemize}
\item
The $\nbigr_{X\times\cnum_t}(\ast t)$-triple
$i_{g\dagger}\varphi_{\dagger}^0
 \gbigt(\Etilde)(\ast t)
\otimes\gbigl(-g^{\ast}\gminia)$
is strictly specializable along $t$.
\item
$\bigoplus_l\Gr^W_{l}\psitilde_{g,u,\gminia}
 \varphi_{\dagger}^0\bigl(\gbigt(\Etilde)\bigr)$
with the induced morphism $\nbign$
and Hermitian sesqui-linear duality
is a polarized graded wild Lefschetz
twistor $D$-module of weight $0$
for each $u\in\real\times\cnum$.
Moreover,
$\psitilde_{g,u,\gminia}
 \varphi_{\dagger}^0\bigl(\gbigt(\Etilde)\bigr)=0$
unless $u\in\real\times\nbiga$.
\end{itemize}

Let us show that
$i_{g\dagger}\varphi_{\dagger}^0
 \gbigt(\Etilde)(\ast t)$
and 
$i_{g\dagger}\gbigt(E)(\ast t)$
are isomorphic.
We have only to show that
$\varphi_{\dagger}^0
 \gbigt(\Etilde)(\ast g)$
and $\gbigt(E)(\ast g)$
are isomorphic.
By construction,
their restrictions to $X-\{g=0\}$ 
are naturally isomorphic.
As in Section \ref{subsection;07.12.3.15},
the isomorphism of the underlying $\nbigr$-modules
is extended to 
$\gbige(\ast g)\simeq
 \varphi_{\dagger}^0\gbigetilde(\ast g)$,
and hence 
$\varphi_{\dagger}^0
 \gbigt(\Etilde)(\ast g)\simeq\gbigt(E)(\ast g)$.
Hence, we are done
in the unramified case.

\subsubsection{General case}
Let $\varphi_m:X\times\cnum_{t_m}\lrarr X\times\cnum_t$
be induced by  $\varphi_m(t_m)=t_m^m$.
Let $\pi_m:\Xtilde\lrarr X$ be the ramified covering given by
$\pi_m(\zeta_1,\ldots,\zeta_n)
=(\zeta_1^m,\ldots,\zeta_k^m,\zeta_{k+1},\ldots,\zeta_n)$.
The induced morphism
$\Xtilde\times\cnum_{t_m}\lrarr X\times\cnum_{t_m}$
is also denoted by $\pi_m$.
We put $\pitilde_m:=\varphi_m\circ\pi_m$.
Let $\Gamma_g:=\{t-g=0\}\subset X\times\cnum_{t}$.
Let $\omega_m$ be a primitive $m$-th root of $1$.
We have the following decomposition:
\[
 \pitilde_m^{-1}\Gamma_g
=\bigcup_{p=0}^{m-1}
 \bigl\{
 t_m-\omega_m^p\cdot g(\zeta_1,\ldots,\zeta_n)
 =0
\bigr\}
\]
Let $j_p:\Xtilde\lrarr\Xtilde\times\cnum_{t_m}$
be the graph of $\omega_m^p\cdot g(\zeta_1,\ldots,\zeta_n)$.

We set
$(\Etilde,\delbar_{\Etilde},\thetatilde,\htilde):=
 \pi_m^{-1}(E,\delbar_E,\theta,h)$,
which is unramifiedly $\nbiga$-good wild.
We have the associated $\nbigr$-module $\gbigetilde$
and the associated $\nbigr$-triple $\gbigt(\Etilde)$.
The direct sum
$\bigoplus_p j_{p\dagger}
 \gbigt(\Etilde)(\ast t_m)
 \otimes\gbigl(-\gminia)$
is $\Gal(\Xtilde/X)$-equivariant.

\begin{lem}
\label{lem;07.10.28.71}
$\varphi_m^{\ast}i_{g\dagger}\gbigt(E)(\ast t)
  \otimes\gbigl(-\gminia)$
is identified with
the $\Gal(\Xtilde/X)$-invariant part of
\[
 \pi_{m\dagger} \Bigl(
 \bigoplus_pj_{p\dagger}\gbigt(\Etilde)(\ast t_m)
 \otimes\gbigl(-\gminia)
 \Bigr).
\]
\end{lem}
\pf
We have the following natural isomorphisms:
\[
 \pi_m^{\ast}\bigl(
 \varphi_m^{\ast}
 (i_{g\dagger}\gbige\otimes\nbigo(\ast t))
 \bigr)
\simeq
 \pitilde_m^{\ast}(i_{g\dagger}\gbige\otimes\nbigo(\ast t))
\simeq
 \bigoplus j_{p\dagger}\gbigetilde\otimes\nbigo(\ast t_m)
\]
Therefore,
we have the following:
\[
 \bigoplus j_{p\dagger}\gbigetilde
 \otimes\nbigo(\ast t_m)
 \otimes\nbigl(-\gminia)
\simeq
 \pi_m^{\ast}\Bigl(
  \varphi_m^{\ast}
 \bigl(i_{g\dagger}\gbige
 \otimes\nbigo(\ast t)\bigr)
\otimes\nbigl(-\gminia)
 \Bigr)
\]
By using Lemma \ref{lem;08.1.16.25},
we obtain the following:
\[
 \pi_{m\dagger}\Bigl(
 \bigoplus j_{p\dagger}\gbigetilde\otimes\nbigo(\ast t_m)
 \otimes\nbigl(-\gminia)
 \Bigr)
\simeq
  \varphi_m^{\ast}
 \bigl(i_{g\dagger}\gbige\otimes\nbigo(\ast t)\bigr)
\otimes\nbigl(-\gminia)
 \otimes \pi_{m\dagger}\pi_m^{\ast}\nbigo(\ast t_m)
\]
Hence,
$ \varphi_m^{\ast}
 \bigl(i_{g\dagger}\gbige\otimes\nbigo(\ast t)\bigr)
\otimes\nbigl(-\gminia)$
is identified with the $\Gal(\Xtilde/X)$-invariant part of 
$\pi_{m\dagger}\Bigl(
 \bigoplus j_{p\dagger}
 \gbigetilde \otimes\nbigl(-\gminia)
 \Bigr)$.

We have the sesqui-linear pairing of
$ \varphi_m^{\ast}
 \bigl(i_{g\dagger}\gbige\otimes\nbigo(\ast t)\bigr)
\otimes\nbigl(-\gminia)$
induced by $\varphi_m^{\ast}C$.
The sesqui-linear pairing of
$\bigoplus j_{p\dagger}\gbige\otimes\nbigl(-\gminia)$
also induces a sesqui-linear pairing of
$ \varphi_m^{\ast}
 \bigl(i_{g\dagger}\gbige\otimes\nbigo(\ast t)\bigr)
\otimes\nbigl(-\gminia)$.
Since the restriction of them to
$X\times\cnum_{t_m}\setminus 
\{t_m= 0\}$ are the same,
they are the same on $X\times\cnum_{t_m}$.
(Note the sesqui-linear pairing 
of $\nbigr(\ast t_m)$-triples 
has the values in the moderate distributions.)
Thus, we obtain Lemma \ref{lem;07.10.28.71}.
\hfill\qed

\vspace{.1in}
Let us return to the proof of
Proposition \ref{prop;07.10.28.72}.
By the previous result in the unramified case,
$j_{p\dagger}\gbigt(\Etilde)(\ast t_m)
 \otimes\gbigl(-\gminia)$
are strictly specializable along $t_m$,
and 
$\Bigl(
 \Gr^W\psitilde_{t_m,u,\gminia}\bigl(
 j_{p\dagger}\gbigt(\Etilde) \bigr),
\nbign_p,\nbigs_p
\Bigr)$
are  polarized graded $\nbiga$-wild Lefschetz twistor $D$-modules
of weight $0$,
where $\nbign_p$ and $\nbigs_p$
are the naturally induced nilpotent maps
and Hermitian sesqui-linear duality.
Moreover,
$\psitilde_{t_m,u,\gminia}\bigl(
 j_{p\dagger}\gbigt(\Etilde)
 \bigr)=0$
unless $u\in\real\times\nbiga$.
By Proposition \ref{prop;05.1.15.61}
and Lemma \ref{lem;08.1.16.40},
the following holds:
\begin{itemize}
\item
 $\pi_{m\dagger}\bigl(
 j_{p\dagger}\gbigt(\Etilde)
 \otimes\nbigl(-\gminia)
 \bigr)$ are strictly specializable along $t_m$.
\item
$\bigl(
\Gr^W\psitilde_{t_m,u,\gminia}
 \pi_{m\dagger}\bigl(
 j_{p\dagger}\gbigt(\Etilde)
 \bigr),\nbign_p',\nbigs_p'
\bigr)$ 
are polarized graded 
$\nbiga$-wild Lefschetz twistor $D$-modules
of weight $0$,
where $\nbign_p'$  and $\nbigs_p'$
denote the induced nilpotent maps
and Hermitian sesqui-linear duality.
Moreover,
$\psitilde_{t_m,u,\gminia}
 \pi_{m\dagger}\bigl(
 j_{p\dagger}\gbigt(\Etilde)
 \bigr)=0$
unless $u\in\real\times\nbiga$.
\end{itemize}
Then, the first claim of 
Proposition \ref{prop;07.10.28.72}
follows from Lemma \ref{lem;07.10.28.71}.
Since 
$P\Gr^{W(N)}_l\psitilde_{g,\gminia,u}\gbigt(E)$
is a direct summand of
$\bigoplus P\Gr^W_l\psitilde_{t_m,u,\gminia}
 \pi_{m\dagger}\bigl(
 j_{p\dagger}\gbigt(\Etilde)
 \bigr)$,
the second and third claims 
of Proposition \ref{prop;07.10.28.72}
follows.
\hfill\qed

\subsection{End of Proof of
Proposition \ref{prop;07.8.30.5}}
\label{subsection;07.10.28.42}

Let $g$ be any function on $X$.
Let $\gminia\in \cnum[t_m^{-1}]$.
We take a complex manifold $X_1$ and
a birational projective morphism
$\varphi:X_1\lrarr X$ such that
(i) $X_1-\varphi^{-1}(g^{-1}(0)) \simeq X-g^{-1}(0)$,
(ii) $\varphi^{-1}(g^{-1}(0)\cup D)$ is normal crossing.
We have the unramifiedly $\nbiga$-good wild harmonic bundle
$(\Etilde,\delbar_{\Etilde},\thetatilde,\htilde):=
 \varphi^{-1}(E,\delbar_E,\theta,h)$.
We have the associated $\nbigr$-triple
$\gbigt(\Etilde)=(\gbigetilde,\gbigetilde,\gbigctilde)$.
Let $\gtilde:=g\circ\varphi$.
Due to Proposition \ref{prop;07.10.28.72},
the following holds:
\begin{itemize}
\item
$\gbigt(\Etilde)$ is strictly specializable
along $\gtilde$
with ramification and exponential twist by $\gminia$.
\item
$P\Gr^{W(N)}_l\psitilde_{\gtilde,u,\gminia}\gbigt(\Etilde)$
is an $\nbiga$-wild pure twistor $D$-module of weight $l$,
and the naturally induced Hermitian sesqui-linear duality
gives a polarization.
Moreover,
$\psitilde_{\gtilde,u,\gminia}\gbigt(\Etilde)=0$
unless $u\in\real\times\nbiga$.
\end{itemize}

By construction,
we have
$\varphi^i_{\dagger}\gbigt(\Etilde)(\ast g)=0$
unless $i=0$.
According to Proposition \ref{prop;05.1.15.61},
the following holds:
\begin{itemize}
\item
 $\varphi^0_{\dagger}\gbigt(\Etilde)$
 is strictly specializable along $g$
 with ramification and exponential twist by $\gminia$.
\item
$P\Gr^{W(N)}_l\psitilde_{g,u,\gminia}
 \varphi^0_{\dagger}\gbigt(\Etilde)$
is an $\nbiga$-wild pure twistor $D$-module of weight $l$,
and the naturally induced Hermitian sesqui-linear duality
gives a polarization.
Moreover,
$\psitilde_{g,u,\gminia}
 \varphi^0_{\dagger}\gbigt(\Etilde)=0$
unless $u\in\real\times\nbiga$.
\end{itemize}
As in Section \ref{subsection;07.12.3.15},
we have the natural isomorphism
$\gbigt(E)(\ast g)
\simeq
 \varphi_{\dagger}^0\gbigt(\Etilde)(\ast g)$.
Thus, the proof of Proposition \ref{prop;07.8.30.5}
is finished.
\hfill\qed

\section{Wildness and uniqueness}
\label{subsection;07.10.28.35}
We shall show Lemma \ref{lem;10.5.25.1}
and the fully faithfulness of $\Phi$
in Theorem \ref{thm;07.10.28.30}.
The faithfulness is clear 
by strict S-decomposability 
of wild pure twistor $D$-modules.
Indeed, let $f:\nbigm_1\lrarr\nbigm_2$
be a morphism of strictly $S$-decomposable
$\nbigr$-modules whose strict supports are $Z$.
If its restriction to a Zariski open subset of $Z$
is $0$, it is $0$ on $Z$.
Hence, we have only to show
Lemma \ref{lem;10.5.25.1}
and that the induced functor $\Phi$ is full.
Let us rewrite the claims.
Let $Z$ be an irreducible closed subset of $X$.
Let $(\nbigt,\nbigs)$ be a polarized $\nbiga$-wild
pure twistor $D$-module of weight $w$
whose strict support is $Z$.
There exists a Zariski open subset $U$ of $Z$ such that
$(\nbigt,\nbigs)_{|X-Y}$ comes from 
a harmonic bundle $\harmonicbundle$ on $U$
with the Tate twist, 
where $Y:=Z-U$.

\begin{prop}
\label{prop;07.10.28.36}
 \mbox{{}}
\begin{itemize}
\item
$\harmonicbundle$ is $\nbiga$-wild on $(Z,U)$.
In particular, we obtain a functor
\[
 \Phi: \MPT^{\wild}_{\strict}(Z,U,w,\nbiga)
 \lrarr \VPTgenwild(Z,U,w,\nbiga).
\]
\item
Recall that we have constructed
the polarized $\nbiga$-wild
pure twistor $D$-module
$(\gbigt,\gbigs)$ of weight $0$
associated to $\harmonicbundle$ 
which is $\nbiga$-wild on $(Z,U)$.
Then, we have a natural isomorphism
$(\nbigt,\nbigs)\simeq
 (\gbigt,\gbigs)\otimes\Tate^S(-w/2)$.
\end{itemize}
\end{prop}
The first claim is Lemma \ref{lem;10.5.25.1},
and the second claim implies that
$\Phi$ is full.

In Subsection \ref{subsection;10.5.25.10},
we consider the case $\dim Z=1$.
To argue the higher dimensional case,
we reduce the issue to the local case
in Subsection \ref{subsection;10.5.25.11}.
We give a preparation in
Subsection \ref{subsection;10.5.25.15}.
We study the restriction to curves
in Subsection \ref{subsection;07.11.2.30},
and the behaviour around a good point 
in Subsection \ref{subsection;08.9.3.1}.
Then, we finish the proof
in Subsection \ref{subsection;10.5.25.20}.

Although much part of the argument
is essentially the same
as that in the tame case,
we have some additional difficulties.
In the tame case,
it is rather easy to show that
the corresponding harmonic bundle is tame
by using a convenient curve test,
i.e., for a given harmonic bundle,
it is tame if its restrictions to curves
are tame.
Although we have a curve test
in the wild case
(Proposition \ref{prop;07.10.13.11}),
we need preliminary to apply it,
and additional argument
to use Proposition \ref{prop;07.11.1.3}.

Another difficulty is caused 
for the second claim by Stokes structure.
Recall that we have the uniqueness of
an extension of a flat bundle
to regular singular meromorphic flat bundle.
Hence,
if we are given regular polarized pure twistor
$D$-modules $\nbigt_i$
$(i=1,2)$ whose restriction to a Zariski open subset
is isomorphic, it is rather easy to obtain
that the underlying $\nbigr$-modules of $\nbigt_i$
are isomorphic.
Because we do not have such uniqueness
in the irregular case,
we need some arguments
to obtain the desired isomorphisms.
Hence, we will study the comparison
of the associated meromorphic objects.

Recall that we have already established
the correspondence
in the case that $Z$ is smooth and one dimensional
(Section \ref{subsection;08.10.18.120}),
in which the wildness is easy to show,
but the uniqueness is not so easy.
Note that we have used the uniqueness result
(Theorem \ref{thm;10.5.20.4})
in an essential way.

\subsection{The case $\dim Z=1$}
\label{subsection;10.5.25.10}

Assume $\dim Z=1$.
Let $P\in Z\setminus U$.
We take a small neighbourhood $X_P$ of $P$.
We have the decomposition
$\nbigt_{|X_P}=\bigoplus_i\nbigt_i$
by strict supports,
corresponding to the irreducible decomposition
$Z\cap X_P=\bigcup_{i\in\Lambda}Z_{i}$.
We have only to show the claims
of Proposition \ref{prop;07.10.28.36}
for each $\nbigt_i$.
Hence, we may assume (i) $X=\Delta^n$,
(ii) $Z$ is not contained in $\{z_n=0\}$.
Let $(\nbigt,\nbigs)$ be 
a polarized $\nbiga$-wild pure twistor $D$-module
of weight $0$, whose strict support is $Z$.
We assume that
the restriction of $(\nbigt,\nbigs)$
to $X-\{z_n=0\}$ comes from
a harmonic bundle $\harmonicbundle$ on
$U=Z\setminus\{z_n=0\}$.
We may assume
$\nbigt=(\nbigm,\nbigm,C)$
and $\nbigs=(\id,\id)$.
Let $\pi:X\lrarr \Delta$ be the projection
onto the $n$-th component.
Let $\Ztilde:=\Delta_z$, and
let $\varphi:\Ztilde\lrarr Z$ be a normalization.
The composite $\pi\circ\varphi$ is denoted by $\pitilde$.
We may assume $\pitilde(z)=z^m$ for some $m>0$.

Let $\nbigt'=(\nbigm',\nbigm',\nbigc')$ denote
the direct summand of $P\pi^0_{\dagger}\nbigt$
whose strict support is $\Delta$.
It is an $\nbiga$-wild pure twistor $D$-module of weight $0$
with the naturally induced polarization $\nbigs'=(\id,\id)$.
We also have the corresponding $\nbiga$-wild harmonic bundle
$(E',\delbar_{E'},\theta',h'):=
 \pi_{\ast}(E,\delbar_E,\theta,h)$
on $\Delta^{\ast}$.
According to Proposition \ref{prop;07.10.19.50},
$(\nbigt',\nbigs')$ is naturally isomorphic 
to the $\nbiga$-wild pure twistor $D$-module
associated to $(E',\delbar_{E'},\theta',h')$
by the construction
in Section \ref{subsection;07.10.28.80}
or Section \ref{subsection;07.10.28.31}.

\begin{lem}
$(\Etilde,\delbar_{\Etilde},\thetatilde,\htilde)
:=\varphi^{\ast}(E,\delbar_E,\theta,h)$
is an $\nbiga$-wild harmonic bundle.
Namely, the first claim of Proposition 
{\rm\ref{prop;07.10.28.36}}
holds in the case $\dim Z=1$.
\end{lem}
\pf
Since $(\Etilde,\delbar_{\Etilde},\thetatilde,\htilde)$
is a direct summand of
$\pitilde^{\ast}(E',\delbar_{E'},\theta',h')$,
the claim is clear.
\hfill\qed

\vspace{.1in}
Let us show the second claim
of Proposition \ref{prop;07.10.28.36}
in the case $\dim Z=1$.
We have the wild pure twistor $D$-module
$\gbigt(\Etilde)=(\gbigetilde,\gbigetilde,\gbigctilde)$
with the polarization $\gbigstilde=(\id,\id)$ on $\Ztilde$,
associated to $(\Etilde,\delbar_{\Etilde},\thetatilde,\htilde)$.
Let $\nbigl_c$ be the Lefschetz map of
$\varphi^{\bullet}_{\dagger}\gbigt(\Etilde)$
for a relatively ample line bundle
with respect to $\varphi$.
Let $P\varphi^0_{\dagger}\gbigt(\Etilde)$
be the kernel of
$\nbigl_c:\varphi^0_{\dagger}\gbigt(\Etilde)
 \lrarr \varphi^2_{\dagger}\gbigt(\Etilde)$.
Let $\nbigt_1=(\nbigm_1,\nbigm_1,C_1)$ be 
the direct summand of 
$P\varphi^0_{\dagger}\gbigt(\Etilde)$
whose strict support is $Z$,
which is equipped
with the naturally induced polarization $\nbigs_1$.
Let $\nbigm_1$ denote the underlying $\nbigr_X$-module
of $\nbigt_1$.
We would like to show
$(\nbigt_1,\nbigs_1)\simeq(\nbigt,\nbigs)$.
For that purpose,
we have only to show that
the natural isomorphism
$\nbigm_{1|X-\{z_n=0\}}\simeq\nbigm_{|X-\{z_n=0\}}$
is extended to an isomorphism on $X$.

\begin{lem}
 \label{lem;07.11.2.3}
We have a natural isomorphism
$\pi_{\dagger}\nbigm(\ast z_n)
\simeq
 \nbigm^{\prime}(\ast z_n)$.
\end{lem}
\pf
By construction of $\nbigm^{\prime}(\ast z_n)$,
we have a natural morphism
$\iota:\nbigm^{\prime}
\lrarr \pi_{\dagger}\nbigm$,
and the support of $\Ker(\iota)$
and $\Cok(\iota)$ are contained in $\{z_n=0\}$.
Because of the coherence,
we obtain the vanishing of $\Cok(\iota)$
and $\Ker(\iota)$ after the localization.
\hfill\qed

\vspace{.1in}
In general,
for a given $\nbigr_X$-module $\nbign$ on 
$\cnum_{\lambda}\times X$,
let $\nbign^{\lambda}$ denote the specialization
of $\nbign$ to $\{\lambda\}\times X$.
Let $\lambda\neq 0$.
We set $\nbigv^{\lambda}:=
L^{-n+1}\varphi^{\dagger}\nbigm^{\lambda}(\ast z_n)$.
Note $L^j\varphi^{\dagger}\nbigm^{\lambda}(\ast z_n)=0$
for any $j\neq -n+1$.
As remarked in Lemma \ref{lem;07.11.2.1},
it gives a meromorphic $\lambda$-flat bundle
on $(\Ztilde,\Dtilde)$.

\begin{lem}
\label{lem;07.11.2.5}
We have a natural isomorphism
$\tr:\varphi^0_{\dagger}\nbigv^{\lambda}
\lrarr 
 \nbigm^{\lambda}(\ast z_n)$.
\end{lem}
\pf
As remarked in Section \ref{subsection;08.1.21.1},
we have the trace map
\[
 \tr:\varphi_{\dagger}L\varphi^{\ast}
 \nbigm^{\lambda}[1-n]
\lrarr
 \nbigm^{\lambda}.
\]
If we take the localization with respect to $z_n$,
we have 
$\varphi_{\dagger}L\varphi^{\ast}
 \nbigm^{\lambda}[1-n] (\ast z_n)
\simeq
 \varphi^0_{\dagger}\nbigv^{\lambda}$.
Thus, we obtain the desired morphism.
To check that it is an isomorphism,
we have only to see that the restriction
is an isomorphism on $X-\{z_n=0\}$,
and it is clear.
\hfill\qed

\begin{lem}
\label{lem;07.11.2.6}
We have a natural isomorphism
$\nbigv^{\lambda}\simeq
 \nbigq\nbigelambdatilde=
 \gbigetilde^{\lambda}(\ast z)$.
\end{lem}
\pf
Because
$\nbigm^{\prime\lambda}(\ast z_n)
=\nbigq\nbige^{\prime\lambda}$
as in Proposition \ref{prop;07.10.19.50},
we have the natural inclusion
$\nbigq\nbigelambdatilde
 \subset \pitilde^{\ast}\nbigm^{\prime\lambda}(\ast z_n)$.
Note $\nbigq\nbigetilde^{\lambda}$,
$\nbigv^{\lambda}$ and
$\pitilde^{\ast}\nbigm^{\prime\lambda}(\ast z_n)$
are locally free $\nbigo_{\Ztilde}(\ast\Dtilde)$-modules.
Because the restriction of 
$\nbigv^{\lambda}$ and
$\nbigq\nbigelambdatilde$ to $\varphi^{-1}(U)$
are the same in 
$\pitilde^{\ast}\nbigm^{\prime\lambda}(\ast z_n)
 _{|\varphi^{-1}(U)}$,
we have only to show that
$\nbigv^{\lambda}$ is also contained in
$\pitilde^{\ast}\nbigm^{\prime\lambda}(\ast z_n)$.
We obtain the following natural isomorphisms
from Lemma \ref{lem;07.11.2.3} and 
Lemma \ref{lem;07.11.2.5}:
\[
 \pitilde^{\ast}\nbigm^{\prime\lambda}(\ast z_n)
\simeq
 \pitilde^{\ast}\bigl(
 \pi_{\dagger}^0\nbigm^{\lambda}(\ast z_n)
 \bigr)
\simeq
 \pitilde^{\ast}
 \pi_{\dagger}^0\varphi_{\dagger}^0
 \bigl(\nbigv^{\lambda}\bigr) 
\simeq
 \pitilde^{\ast}\pitilde_{\dagger}
 \bigl(\nbigv^{\lambda}\bigr)
\]
Since $\pitilde$ is a ramified covering,
$\pitilde_{\dagger}$ is the same as
the push-forward for $\nbigo$-modules,
as explained in Lemma \ref{lem;08.1.16.25}.
Hence, 
$\nbigv^{\lambda}$ is naturally identified with
a direct summand of
$\pitilde^{\ast}\pitilde_{\dagger}(\nbigv^{\lambda})$.
Thus, we are done.
\hfill\qed

\begin{lem}
We have a natural isomorphism
$\nbigm^{\lambda}(\ast z_n)\simeq
 \nbigm^{\lambda}_1(\ast z_n)$
for $\lambda\neq 0$.
\end{lem}
\pf
It follows from Lemma \ref{lem;07.11.2.5}
and Lemma \ref{lem;07.11.2.6}.
\hfill\qed

\vspace{.1in}

The rest is essentially the same as 
the argument given in Section 19.4.2
in \cite{mochi2}.
We give an outline with minor simplification.
Let $\lambda_0\in\cnum_{\lambda}$.
Let $U(\lambda_0)$ denote a small neighbourhood 
of $\lambda_0$.
Let $\vecv=(v_i)$ be a frame of $\nbigqzero_{<0}\nbigetilde$.
Each $v_i$ naturally induces a
section $\varphi_{\dagger}(v_i)$
of $\nbigm_1(\ast z_n)$.
The tuple of $\varphi_{\dagger}(v_i)$
is denoted by $\varphi_{\dagger}(\vecv)$.
Let $\varphi_{\dagger}(\vecv)\cdot \nbigr_X(\ast z_n)$ 
denote the submodule of $\nbigm_1(\ast z_n)$
generated by $\varphi_{\dagger}(\vecv)$
over $\nbigr_X(\ast z_n)$.
\begin{lem}
\label{lem;07.11.2.10}
We have 
$\varphi_{\dagger}(\vecv)\cdot\nbigr_X(\ast z_n)
=\nbigm_1(\ast z_n)$.
In other words,
$\varphi_{\dagger}(\vecv)$
generates $\nbigm_1(\ast z_n)$
over $\nbigr_X(\ast z_n)$.

In particular,
$\nbigm_1^{\lambda}(\ast z_n)$
is generated by the restriction of
$\varphi_{\dagger}(\vecv)$
to $\{\lambda\}\times X$.
\end{lem}
\pf
It is easy to show the coincidence of
the restriction of them to $X-\{z_n=0\}$.
Then, the claim follows from Lemma \ref{lem;08.1.21.11},
for example.
\hfill\qed

\vspace{.1in}

Let $\nbige$ denote the $\nbigr_U$-module 
on $\cnum_{\lambda}\times U$ 
associated to $\harmonicbundle$.
Let $\iota_U:U\lrarr X$.
Then, we obtain the $\nbigr_X$-module
$\iota_{U\dagger}\nbige$.
As remarked in Lemma 19.25 of \cite{mochi2},
any section $f$ of $\iota_{U\dagger}\nbige$
(resp. $\iota_{U\dagger}\nbige^{\lambda}$)
on $U(\lambda_0)\times X$
(resp. $X$)
is uniquely expressed as follows:
\begin{equation}
 \label{eq;07.11.2.12}
 f=\sum_{\vecn\in\seisuu^{n-1}_{\geq 0}}
 \varphi_{\dagger}(v_i)\cdot
 f_{\vecn,i}\cdot \prod_{j=1}^{n-1}\deldel_j^{n_j}
\end{equation}
Here, $f_{\vecn,i}$ are holomorphic functions on 
$U(\lambda_0)\times U$ (resp. $U$).

\begin{lem}
\label{lem;07.11.2.11}
Let $f$ be a section of $\nbigm_1^{\lambda}(\ast z_n)$ 
on $X$.
Then, each $f_{\vecn,i}$ is meromorphic
on $\Ztilde$.
\end{lem}
\pf
Let $z$ be a holomorphic coordinate of $\Ztilde$.
Then, $\varphi^{-1}(\del/\del z_n)$ 
is expressed as the linear combination of
$\del/\del z$ and $\varphi^{-1}(\del/\del z_i)$
$(i=1,\ldots,n-1)$
with meromorphic coefficients.
Then, the claim of Lemma \ref{lem;07.11.2.11}
follows from Lemma \ref{lem;07.11.2.10}.
\hfill\qed

\vspace{.1in}
\begin{lem}
\label{lem;07.11.2.13}
We have a natural isomorphism
$\nbigm(\ast z_n)\simeq\nbigm_1(\ast z_n)$.
\end{lem}
\pf
We have only to show
$\nbigm(\ast z_n)\subset\nbigm_1(\ast z_n)$
in $\iota_{U\dagger}\nbige$.
Let $f$ be a section of
$\nbigm(\ast z_n)$ on $U(\lambda_0)\times X$.
We have the expression of $f$
as in (\ref{eq;07.11.2.12}).
The restrictions $f_{\vecn,i|\{\lambda\}\times X}$
are meromorphic on $X$
for each $\lambda\neq 0$ contained in $U(\lambda_0)$.
Hence, as remarked in Lemma 19.23 of \cite{mochi2},
we obtain that $f_{\vecn,i}$ is meromorphic 
on $U(\lambda_0)\times X$.
Then, it is easy to see that
$f$ is contained in $\nbigm_1(\ast z_n)$.
Thus, we obtain Lemma \ref{lem;07.11.2.13}.
\hfill\qed

\vspace{.1in}

Then, we obtain $\nbigm=\nbigm_1$
due to Lemma \ref{lem;07.11.2.15}.
Thus, the second claim of Proposition 
\ref{prop;07.10.28.36} is finished
in the case $\dim Z=1$.

\subsection{Preliminary for the general case}
\label{subsection;10.5.25.11}

Before going to the case $\dim Z>1$,
we give some lemmas to show
that we can shrink $X$ and $U$
arbitrarily in the proof.
Let $Z$ be a closed analytic subset of $X$.
Let $(E,\delbar_E,\theta,h)$ be a harmonic bundle
defined on a Zariski open subset $U$ of $Z$.

\begin{lem}
\label{lem;08.1.21.30}
Assume the following:
\begin{description}
\item[(Locally wild)]
For any point $P$ of $Z$,
there exists a neighbourhood $\nbigu_P$ of $P$
in $X$ such that 
$\harmonicbundle_{|\nbigu_P\cap U}$
is wild on $(Z\cap \nbigu_P,U\cap\nbigu_P)$.
\end{description}
Then, $\harmonicbundle$ is
wild on $(Z,U)$.
\end{lem}
\pf
We take a birational projective
morphism $\varphi:Z_1\lrarr Z$
such that $D_1:=\varphi^{-1}(Z-U)$ is 
a normal crossing hypersurface of $Z_1$.
By the assumption,
the eigenvalues of
$\varphi^{-1}(\theta)$ are 
(possibly multi-valued) meromorphic $1$-forms.
By replacing $\varphi$ and $U$ appropriately,
we can assume that the ramification of
the $1$-forms may happen only at $D_1$.
For any point of $P_1\in Z_1$,
the local wildness condition is satisfied.
Hence, $\varphi^{-1}(\theta)$ is generically good.
Due to Proposition \ref{prop;07.11.1.3},
we can take an appropriate refinement of $\varphi$
such that $\varphi^{\ast}\theta$ is good.
Then, $\varphi^{\ast}\harmonicbundle$ is good wild.
\hfill\qed

\vspace{.1in}

Let $(E,\delbar_E,\theta,h)$ be 
a wild harmonic bundle on $(Z,U)$.
We have 
the polarized wild pure twistor $D$-module
$(\gbigt,\gbigs)$ of weight $0$,
associated to $(E,\delbar_E,\theta,h)$
as in Section \ref{subsection;07.10.28.31}.
Assume that we are given
another polarized wild pure twistor $D$-module 
 $(\nbigt,\nbigs)$ of weight $0$
whose strict support is $Z$,
such that the restriction to $X-Y$
comes from $(E,\delbar_E,\theta,h)$,
where $Y:=Z-U$.
Then, we have the natural isomorphism
$F:(\gbigt,\gbigs)_{|X-Y}
\simeq
 (\nbigt,\nbigs)_{|X-Y}$

\begin{lem}
\label{lem;08.1.21.31}
Assume the following:
\begin{description}
\item[(Locally extendable)]
 For any point $P\in Z$,
 there exists a neighbourhood $\nbigu_P$ of $P$
 in $X$
 such that the isomorphism
 $F_{|\nbigu_P\setminus Y}$
 can be extended to an isomorphism
 $F_{\nbigu_P}:
 (\gbigt,\gbigs)_{|\nbigu_P} 
\lrarr(\nbigt,\nbigs)_{|\nbigu_P}$.
\end{description}
Then,
$F$ is extended to an isomorphism
$\Ftilde:(\gbigt,\gbigs)\lrarr (\nbigt,\nbigs)$.
\end{lem}
\pf
By the strict $S$-decomposability of
the underlying $\nbigr$-modules of pure twistor $D$-modules,
we have the following general and easy fact:
\begin{itemize}
\item 
 Let $\nbigt$ be a pure twistor $D$-module 
 whose strict support is $Z_1$.
 Let $\varphi$ be an automorphism of $\nbigt$.
 If the restriction of $\varphi$ to some Zariski open subset
 of $Z_1$ is equal to the identity,
 then $\varphi$ is equal to the identity.
\end{itemize}
Then, we obtain
$F_{\nbigu_P|\nbigu_P\cap\nbigu_Q}
=F_{\nbigu_Q|\nbigu_P\cap\nbigu_Q}$,
and thus we obtain the global isomorphism.
\hfill\qed

\vspace{.1in}

By Lemmas \ref{lem;08.1.21.30}
and \ref{lem;08.1.21.31},
we may shrink $X$ arbitrarily to show 
Proposition \ref{prop;07.10.28.36}
in the following argument
without mention.

Let us discuss shrinking $U$.
Let $U'\subset U$ be a Zariski open subset.
We put $Y:=Z\setminus U$
and $Y':=Z\setminus U'$.

\begin{lem}
Let $\harmonicbundle$ be 
a harmonic bundle on $U$.
Then,
$\harmonicbundle$ is wild on $(Z,U)$
if and only if
$\harmonicbundle_{|U'}$ is wild on $(Z,U')$.
\end{lem}
\pf
We argue only the only if part.
The other case can be discussed similarly.
We can take a projective birational morphism
$\varphi:\Ztilde\lrarr Z$
such that
(i) $\Ztilde$ is smooth,
(ii) $\varphi^{-1}(Y)$ and 
$\varphi^{-1}(Y')$ are normal crossing,
(iii) $\varphi^{-1}\harmonicbundle_{U'}$
is good wild on $(\Ztilde,\varphi^{-1}(Y'))$.
Then, it is easy to check that
$\varphi^{-1}\harmonicbundle_{U}$
is good wild on $(\Ztilde,\varphi^{-1}(Y))$.
\hfill\qed

\vspace{.1in}
Let $\harmonicbundle$ be wild on $(Z,U)$.
We have constructed
polarized wild pure twistor $D$-modules
$(\gbigt,\gbigs)$ 
and $(\gbigt',\gbigs')$
associated to
$(E,\delbar_E,\theta,h)$ on $(Z,U)$, and 
$(E,\delbar_E,\theta,h)_{|U'}$ on $(Z,U')$
respectively.
By construction,
we have a natural isomorphism
$(\gbigt,\gbigs)\simeq(\gbigt',\gbigs')$
whose restriction to $U'$ is the identity.
Hence, we may and will shrink $U$ arbitrarily
for the proof of Proposition \ref{prop;07.10.28.36}
in the following argument.

\subsection{Meromorphic flat connections}
\label{subsection;10.5.25.15}

In the following argument,
we will not distinguish meromorphic $\lambda$-flat bundles
and meromorphic flat bundles
in the case $\lambda\neq 0$.
Let us consider the case $\dim Z=k+1>1$.
By shrinking $U$,
we may assume $Z-U=g^{-1}(0)\cap Z$
for some function $g$.
Let $\varphi':\Ztilde'\lrarr Z$ be 
a birational projective morphism
such that $\Ztilde'-(\varphi^{\prime})^{-1}(U)$ is
a normal crossing hypersurface.
Note that the Higgs bundle
$(\varphi')^{-1}(E,\delbar_{E},\theta)$ can be extended
to a meromorphic Higgs sheaf on $\Ztilde'$
given by $(\varphi')^{\ast}\nbigm^0$.
By shrinking $U$, if necessary,
we can take a birational projective morphism
$\varphi:\Ztilde\lrarr Z$, which factors through $\Ztilde'$,
such that the following holds:
\begin{itemize}
\item
 $\Ztilde-\varphi^{-1}(U)$ is a normal crossing
 hypersurface.
\item
 The ramification of the eigenvalues of
 $\varphi^{-1}(\theta)$
 may happen only at $\Dtilde:=\Ztilde-\varphi^{-1}(U)$.
\end{itemize}
We put $\nbigv^{\lambda}:=
 L^{-n+k+1}\varphi^{\ast}\nbigm^{\lambda}
 \otimes\nbigo(\ast \Dtilde)$ for $\lambda\neq 0$.
It is a holonomic $D$-module
whose characteristic variety is contained in
$\Ztilde\cup \pi^{-1}(\Dtilde)$,
where $\pi:T^{\ast}\Ztilde\lrarr\Ztilde$
denotes the natural projection
of the cotangent bundle.
Hence $\nbigv^{\lambda}$ gives 
a meromorphic $\lambda$-flat connection.
It is an $\nbigo_{\Ztilde}(\ast \Dtilde)$-reflexive module
\cite{malgrange}.

\begin{lem}
\label{lem;07.11.1.5}
$\varphi_{\dagger}^0\nbigv^{\lambda}$
and $ \nbigm^{\lambda}(\ast g)$
are naturally isomorphic.
\end{lem}
\pf
As remarked in Section \ref{subsection;08.1.21.1},
we have the trace map
$\tr:\varphi_{\dagger}
 L\varphi^{\ast}\nbigm^{\lambda}[k+1-n]
\lrarr \nbigmlambda$.
After localization with respect to $g$,
we have 
$\varphi_{\dagger}
 L\varphi^{\ast}\nbigm^{\lambda}(\ast g)[k+1-n]
\simeq
 \varphi^0_{\dagger}\nbigv^{\lambda}$.
Then,
we obtain a naturally defined morphism
$ \tr:
 \varphi^0_{\dagger}\nbigv^{\lambda}
\lrarr
 \nbigm^{\lambda}(\ast g)$.
Since the restriction $\tr_{|U}$ is an isomorphism,
the support of $\Ker(\tr)$ and $\Cok(\tr)$ are 
contained in $g^{-1}(0)$.
Since 
the action of $g$ on 
$\varphi_{\dagger}^0\nbigv^{\lambda}$
and $\nbigm^{\lambda}(\ast g)$ are invertible,
we obtain $\tr$ induces an isomorphism on $Z$.
Thus, we obtain Lemma \ref{lem;07.11.1.5}.
\hfill\qed

\subsection{Restriction to a curve}
\label{subsection;07.11.2.30}

Let $\Ctilde$ be a smooth curve in $\Ztilde$
transversal with the smooth part of $\Dtilde$.
We have the wild harmonic bundles 
$(E_{\Ctilde},
 \delbar_{E_{\Ctilde}},\theta_{\Ctilde},
 h_{\Ctilde}):=
 \varphi^{-1}\harmonicbundle_{|\Ctilde\setminus \Dtilde}$.
We have the associated meromorphic $\lambda$-flat bundle
$(\nbigq\nbigelambda_{\Ctilde},\DDlambda_{\Ctilde})$ on 
$\Ctilde$.
We put $C:=\varphi(\Ctilde)$.

\begin{prop}
\label{prop;07.11.1.6}
We have a natural isomorphism
$\nbigq\nbigelambda_{\Ctilde}\simeq 
  \nbigv^{\lambda}_{|\Ctilde}$.
The harmonic bundle 
$(E_{\Ctilde},\delbar_{E_{\Ctilde}},
 \theta_{\Ctilde},h_{\Ctilde})$
is $\nbiga$-wild.
\end{prop}
\pf
In the proof of the proposition,
we may and will shrink $X$ and $Z$ without mention.

\begin{lem}
\label{lem;07.11.1.7}
We can take a function $f$ on $X$
with the following property:
\begin{itemize}
\item
$C\cap U$ is contained in
the smooth part of $f^{-1}(0)\cap Z$.
\end{itemize}
\end{lem}
\pf
Let $I_C$ denote the ideal sheaf of $\nbigo_X$
corresponding to $C$.
Let $P\in C\cap \{g=0\}$.
We have some functions $f_1,\ldots,f_N$
which generates the ideal sheaf $I_C$ around $P$.
If $Q\in C\cap U$ is sufficiently close to $P$,
$f_1,\ldots,f_N$ generate $I_C$ at $Q$.
In particular,
one of $df_{i|Q}$ is not $0$ on $T_QZ$,
which gives a desired function.
\hfill\qed

\begin{lem}
We can take functions
$f_1,\ldots,f_k$ with the following property:
\begin{itemize}
\item
 $C\subset Z\cap\bigcap_{i=1}^kf_i^{-1}(0)$.
\item
 $C\cap U$ is contained in the smooth part of
 $Z\cap \bigcap_{i=1}^jf_i^{-1}(0)$.
\item
 $f_{j+1}^{-1}(0)$ and 
 $Z\cap \bigcap_{i=1}^jf_i^{-1}(0)$
 are transversal at $C\cap U$.
\end{itemize}
\end{lem}
\pf
We have only to apply Lemma \ref{lem;07.11.1.7}
inductively.
\hfill\qed

\vspace{.1in}

We put $Z^{(0)}:=Z$,
$\Ztilde^{(0)}:=\Ztilde$,
$\Ctilde^{(0)}:=\Ctilde$ and $g^{(0)}:=g$.
Let $Z^{(j)}$ denote the irreducible component of
$Z\cap\bigcap_{i\leq j}f_i^{-1}(0)$,
which contains $C$.
We take functions $g^{(j)}=a^{(j)}g^{(j-1)}$ 
such that the singular part of
$Z\cap\bigcap_{i\leq j}f_i^{-1}(0)$
is contained in $(g^{(j)})^{-1}(0)$.
We take $\Ztilde^{(j)}$ inductively as follows:
Assume that we are given 
$\Ztilde^{(j)}$
with a birational projective morphism
$\varphi^{(j)}:\Ztilde^{(j)}\lrarr Z$
and the curve $\Ctilde^{(j)}\subset \Ztilde^{(j)}$
such that $\varphi^{(j)}(\Ctilde^{(j)})=C$.
We take a birational projective morphism
$\mu^{(j)}:\Ztilde^{(j)}_1\lrarr \Ztilde^{(j)}$
with the following property:
\begin{itemize}
\item
 It is bi-holomorphic on 
 $\Ztilde^{(j)}\setminus (g^{(j)})^{-1}(0)$.
\item
 We put $F_j:=
 (\varphi^{(j)}\circ \mu^{(j)})^{\ast}(f_{j+1})$
 and $G_j:=
  (\varphi^{(j)}\circ\mu^{(j)})^{\ast}(g^{(j+1)})$.
 Then, $F_j^{-1}(0)\cup G_j^{-1}(0)$
 is normal crossing.
\end{itemize}
Let $\Ctilde^{(j+1)}$ be the proper transform
of $\Ctilde^{(j)}$.
Let $\Ztilde^{(j+1)}$ denote the irreducible component
of $F_j^{-1}(0)$ which contains $\Ctilde^{(j+1)}$.
We have the induced map
$\varphi^{(j+1)}:\Ztilde^{(j+1)}\lrarr Z^{(j+1)}$.
Note that the intersection of
$\Ztilde^{(j+1)}$ and another irreducible component
of $F_j^{-1}(0)$ is contained in $G_j^{-1}(0)$.
Let $\nbigv^{\lambda(j)}$ be the pull back of
$\nbigvlambda$ via the induced morphism
$\Ztilde^{(j)}\lrarr \Ztilde$.

We obtain the polarized $\nbiga$-wild 
pure twistor $D$-modules
$\nbigt^{(j)}$ of weight $0$
whose strict supports are $Z^{(j)}$,
inductively as follows:
\begin{itemize}
\item
We put $\nbigt^{(0)}:=\nbigt$.
Assume that 
we are given $\nbigt^{(j)}$ on $Z^{(j)}$.
We have the polarized $\nbiga$-wild pure twistor $D$-module
$P\Gr^{W}_0\psi_{f_{j+1},-\vecdelta_0}
 \nbigt^{(j)}$,
 which is decomposed by the strict supports.
Let $\nbigt^{(j+1)}$ denote the 
 $Z^{(j+1)}$-component.
\end{itemize}
Let $\nbigm^{(j)}$ denote the underlying $\nbigr$-module
of $\nbigt^{(j)}$.

\begin{lem}
\label{lem;07.11.1.11}
We have a natural isomorphism
$\varphi^{(j)}_{\dagger}\nbigv^{\lambda(j)}
\simeq
 \nbigm^{(j)\lambda}(\ast g^{(j)})$.
\end{lem}
\pf
We use an induction.
The case $j=0$ is the claim of 
Lemma \ref{lem;07.11.1.5}.
Assume that the claim holds for $j$,
and we will derive the isomorphism for $j+1$.
We put 
$\nbigv_1^{\lambda(j)}:=
 \mu^{(j)\ast}\nbigv^{\lambda(j)}
\otimes\nbigo_{\Ztilde_1^{(j)}}(\ast G_{j})$.
We put $\kappa^{(j)}:=\varphi^{(j)}\circ\mu^{(j)}$.
Then, we have the following natural isomorphisms:
\[
 \kappa^{(j)}_{\dagger}
 \nbigv_1^{\lambda(j)}
\simeq
 \varphi^{(j)}_{\dagger}\mu^{(j)}_{\dagger}
 \mu^{(j)\ast}\nbigv^{\lambda(j)}
 (\ast g^{(j+1)})
\simeq
 \varphi^{(j)}_{\dagger}\nbigv^{\lambda(j)}
 (\ast g^{(j+1)})
\simeq
 \nbigm^{(j)\lambda}(\ast g^{(j+1)})
\]
Here, we have used the isomorphism
$\mu_{\dagger}^{(j)}\mu^{(j)\ast}
 \nbigv^{\lambda(j)}(\ast g^{(j+1)})
\simeq \nbigv^{\lambda(j)}(\ast g^{(j+1)})$
induced by the trace map.
Then, we obtain the following isomorphisms:
\begin{equation}
 \label{eq;07.11.1.10}
 \kappa^{(j)}_{\dagger}\bigl(
 \psi_{F_j,-1}\nbigv_1^{\lambda(j)}
 \bigr)
\simeq
 \psi_{f_{j+1},-1}\bigl(
 \kappa^{(j)}_{\dagger}\nbigv_1^{\lambda(j)}
 \bigr)
\simeq
 \psi_{f_{j+1},-1}\nbigm^{(j)\lambda}(\ast g^{(j+1)})
\end{equation}
Let $\Ztilde^{(j+1)}_i$ denote the irreducible components
of $F_j^{-1}(0)$ which are not contained in
$G_j^{-1}(0)$.
(We put $\Ztilde^{(j+1)}_1:=\Ztilde^{(j+1)}$.)
We have the decomposition by the supports:
\[
 \psi_{F_j,-1}\nbigv_1^{\lambda(j)}
=\bigoplus_i \nbign_{\Ztilde^{(j+1)}_i}
\]
And, $\nbign_{\Ztilde^{(j+1)}_1}
\simeq
 \iota_{\dagger}
 \nbigv^{\lambda(j+1)}$,
where $\iota$ denotes the inclusion
$\Ztilde^{(j+1)}\lrarr \Ztilde^{(j)}_1$.

We have the following natural isomorphisms:
\[
 \nbigm^{(j+1)}(\ast g^{(j+1)})
\simeq
 P\Gr^{W}_0\psi^{(\lambda_0)}_{f_{j+1},-\vecdelta_0}
 \nbigm^{(j)}(\ast g^{(j+1)})
\simeq
 \psi^{(\lambda_0)}_{f_{j+1},-\vecdelta_0}
 \nbigm^{(j)}(\ast g^{(j+1)})
\]
We also have the following isomorphism:
\[
 \Bigl(
\psi^{(\lambda_0)}_{f_{j+1},-\vecdelta_0}
 \nbigm^{(j)}(\ast g^{(j+1)}) 
\Bigr)^{\lambda}
\simeq
 \psi_{f_{j+1},-1}
 \nbigm^{(j)\lambda}(\ast g^{(j+1)})
\]
Hence, the restriction of (\ref{eq;07.11.1.10})
to the $Z^{(j+1)}$-component
gives the desired isomorphism.
Thus, we obtain Lemma \ref{lem;07.11.1.11}.
\hfill\qed

\vspace{.1in}

Let us return to the proof of Proposition \ref{prop;07.11.1.6}.
By construction, we have
$\Ztilde^{(k)}=\Ctilde^{(k)}\simeq \Ctilde$.
Under the isomorphism,
we have 
$\nbigv^{\lambda(k)}\simeq 
 \nbigv^{\lambda}_{|\Ctilde}$
and 
$\varphi^{(k)}=\varphi_{|\Ctilde}$.
We also have $C=Z^{(k)}$,
and $\nbigt^{(k)}$ gives 
the polarized $\nbiga$-wild pure twistor $D$-module
of weight $0$ whose strict supports is $C$.
We have obtained the following isomorphism:
\begin{equation}
 \label{eq;07.11.2.20}
 \varphi^{(k)}_{\dagger}\nbigv^{\lambda(k)}
\simeq
 \nbigm^{(k)\lambda}(\ast g^{(k)})
\end{equation}
The rest of the argument
is essentially the same as the proof of 
Lemma \ref{lem;07.11.2.6}.
We give only an outline.
We may assume to have a coordinate $w$
on $\Ctilde$ such that
$\Ctilde\cap \Dtilde=\{w=0\}$.
We take a projection
$\pi:X\lrarr \Delta$
such that the composite
$\Psi:=\pi\circ\varphi^{(k)}:\Ztilde^{(k)}\lrarr \Delta$
is given by $\Psi(w)=w^{\ell}$ for some $\ell>0$.
Let $\nbigt'=(\nbigm',\nbigm',\nbigc')$ denote
the direct summand of 
the polarized $\nbiga$-wild pure twistor $D$-module
$\pi^0_{\dagger}\nbigt^{(k)}$
whose strict support is $\Delta$.
We have the corresponding $\nbiga$-wild harmonic bundle
$(E',\delbar_{E'},\theta',h')$ on $\Delta^{\ast}$,
which is obtained as the push-forward of
 $(E,\delbar_E,\theta,h)_{|C}$.
As in Lemma \ref{lem;07.11.2.3},
we have the natural isomorphisms:
\[
 \Psi_{\dagger}\nbigv^{\lambda(k)}
\simeq 
 \pi_{\dagger}\varphi^{(k)}_{\dagger}\nbigv^{\lambda(k)}
\simeq
 \pi_{\dagger}\nbigm^{(k)\lambda}(\ast g^{(k)})
\simeq
 \nbigm^{\prime\lambda}(\ast w)
\]
As in the proof of Lemma \ref{lem;07.11.2.6},
we have natural inclusions
$\nbigv^{\lambda(k)}\lrarr
 \Psi^{\ast}\nbigm^{\prime\lambda}(\ast w)$
and 
$\nbigq\nbigelambda_{\Ctilde}\lrarr
 \Psi^{\ast}\nbigm^{\prime\lambda}(\ast w)$.
Since the restriction of
$\nbigv^{\lambda(k)}$
and $\nbigq\nbigelambda_{\Ctilde}$
to $\Ztilde^{(k)}-\{G^{(k)}=0\}$ are the same,
we obtain 
$\nbigv^{\lambda(k)}\simeq
\nbigq\nbigelambda_{\Ctilde}$,
and thus the first claim of
Proposition \ref{prop;07.11.1.6} is proved.
Since 
$(E_{\Ctilde},\delbar_{E_{\Ctilde}},
\theta_{\Ctilde},h_{\Ctilde})$
is a direct summand of
$(\varphi^{(k)})^{-1}(E',\delbar_{E'},\theta',h')$,
the harmonic bundle
$(E_{\Ctilde},\delbar_{E_{\Ctilde}},
\theta_{\Ctilde},h_{\Ctilde})$
is also $\nbiga$-wild.
Thus, we obtain the second claim of
Proposition \ref{prop;07.11.1.6}.
\hfill\qed

\subsection{Around a good point}
\label{subsection;08.9.3.1}

For distinction,
we set $(\Etilde,\delbar_{\Etilde},\thetatilde,\htilde):=
 \varphi^{-1}(E,\delbar_E,\theta,h)$.
Let us fix $\lambda_1\neq 0$.
We continue to use the notation in
Subsection \ref{subsection;07.11.2.30}.
Let $P$ be a smooth point of $\Dtilde$
around which $\nbigv^{\lambda_1}$ is good.
(Recall that the set of good points 
is non-empty and Zariski open.
See \cite{malgrange}.)
\begin{lem}
\label{lem;07.11.1.15}
On an appropriate neighbourhood $\nbigu$ of $P$,
$\bigl(\Etilde,\delbar_{\Etilde},\thetatilde,
 \htilde\bigr)_{|\nbigu\setminus \Dtilde}$ is 
$\nbiga$-good wild on
$(\nbigu,\Dtilde\cap\nbigu)$,
and $\nbigq\nbigelambda_{|\nbigu}$ 
is naturally isomorphic to
$\nbigvlambda_{|\nbigu}$
for any $\lambda\neq 0$.
\end{lem}
\pf
Take a ramified covering
$\eta:(\nbigu',D')\lrarr 
 (\nbigu,\Dtilde\cap \nbigu)$ such that
$\eta^{\ast}
 (\nbigv^{\lambda_1},\DD^{\lambda_1})$
is unramified.
Let $\nbigi$ denote the set of the irregular values of
$\eta^{\ast}(\nbigv^{\lambda_1},
 \DD^{\lambda_1})$.
We may assume that $\nbigu$ is equipped with
a coordinate system $(z_1,\ldots,z_n)$ such that
$\Dtilde\cap\nbigu=\{z_1=0\}$.
Let $\pi:\nbigu\lrarr \Dtilde\cap \nbigu$ 
denote the projection given by the coordinate.

We apply Proposition \ref{prop;07.11.1.6}
to the harmonic bundle
\[
\bigl(
 E_{\pi^{-1}(Q)},
 \delbar_{E_{\pi^{-1}(Q)}},
 \theta_{\pi^{-1}(Q)},
 h_{\pi^{-1}(Q)}
 \bigr) =
 (\Etilde,\delbar_{\Etilde},\thetatilde,\htilde)_{|
 \pi^{-1}(Q)}.
\]
Since 
the meromorphic $\lambda_1$-flat bundle
$\eta^{\ast}\nbigq\nbige^{\lambda_1}_{\pi^{-1}(Q)}$
is unramified,
the meromorphic $\lambda$-flat bundle
$\eta^{\ast}\nbigq\nbige^{\lambda}_{\pi^{-1}(Q)}$
are also unramified for any $\lambda$,
according to Theorem \ref{thm;07.11.18.20}.
Moreover,
the set of irregular values are given as follows:
\[
 \Irr\bigl(
 \eta^{\ast}\nbigq\nbige^{\lambda}_{\pi^{-1}(Q)},
 \DDlambda\bigr)
=\bigl\{
 \gminia_{|\eta^{-1}\pi^{-1}(Q)}\,\big|\,
 \gminia\in \nbigi
 \bigr\}
\]
According to
Lemma \ref{lem;08.1.21.50},
$\eta^{\ast}\nbigvlambda$ also has
an unramifiedly good lattice $V^{\lambda}$.

We have the generalized eigen decomposition
$V^{\lambda}_{|Q}=
 \bigoplus_{\alpha} \EE_{\alpha}
 \bigl(V^{\lambda}_{|Q}\bigr)$
with respect to $\Res(\DDlambda)$.
Let us consider the following set:
\[
 \Sp(Q,\lambda):=
 \Bigl\{
 \alpha+n\cdot\lambda\,\Big|\,
 \EE_{\alpha}\bigl(V^{\lambda}_{|Q}\bigr)\neq 0,
\,\, n\in\seisuu
 \Bigr\}
\]
It is well known and easy to see that the set
$\Sp(Q,\lambda)$
is independent of the choice of 
lattices $V^{\lambda}$.
It is also independent of $Q$.

We have the set 
$\KMS\bigl(\eta^{-1}\Etilde_{|\pi^{-1}(Q)}\bigr)$
of the KMS-spectra of the unramifiedly good
 wild harmonic bundle
$\eta^{-1}(\Etilde,
 \delbar_{\Etilde},\thetatilde,\htilde)_{|\pi^{-1}(Q)}$
at $\lambda=0$
for each $Q$.
We have a discrete subset
$A(Q)\subset\cnum^{\ast}$ such that
the map $\eigenmap(\lambda)$
gives the bijection between
$\KMS\bigl(\eta^{-1}\Etilde_{|\pi^{-1}(Q)}\bigr)$
and $\Sp(Q,\lambda)$
for any $\lambda\in\cnum^{\ast}-A(Q)$.
(See Proposition \ref{prop;07.7.19.31}.)
Then, we can conclude that
the sets 
$\KMS\bigl(\eta^{-1}\Etilde_{|\pi^{-1}(Q)}\bigr)$
are independent of the choice of $Q$,
which is denoted by 
$\KMS\bigl(\eta^{-1}\Etilde,D'\bigr)$.

Let $\lambda$ be generic with respect to
the set $\KMS\bigl(\eta^{-1}\Etilde,D'\bigr)$,
namely, the map
$\eigenmap(\lambda):
 \KMS\bigl(\eta^{-1}\Etilde,D'\bigr)
\lrarr \cnum$ is injective.
For any $c\in\real$,
we can take a good lattice 
$\prolongg{c}{V}'$ of $\eta^{\ast}\nbigvlambda$
with the following property:
\begin{itemize}
\item
 The set of the residue
 $\Res_{D'}(\DDlambda)
 \in \End(\prolongg{c}{V'}_{|D'})$
 is given by the following set:
\[
 \Bigl\{
 \eigenmap(\lambda,u)\,\Big|\,
 u\in \KMS(\eta^{-1}\Etilde,D'),\,\,
 c-1<\paramap(\lambda,u)\leq c
 \Bigr\}
\]
\end{itemize}
Thus, we obtain 
an unramifiedly good filtered $\lambda$-flat bundle
$(V_{\ast}',\DDlambda)$
on $\bigl(\nbigu',\,D'\bigr)$.
By taking the descent,
we obtain a good filtered $\lambda$-flat bundles
$(V_{\ast},\DDlambda)$
on $\bigl(\nbigu,\Dtilde\cap\nbigu\bigr)$.
By construction,
we have 
$\nbigvlambda=
 \bigcup \prolongg{c}V$.

Let $T:=(1+|\lambda|^2)$,
and let us consider the deformations
$V_{\ast}^{\prime(T)}$
and $V_{\ast}^{(T)}$
(Section \ref{subsection;10.5.17.51}).
We have the natural isomorphism
$(\eta^{\ast}\nbigvlambda)^{(T)}_{|\pi^{-1}(Q)}
\simeq 
 \eta^{\ast}\nbigp\nbige_{\pi^{-1}(Q)}$
(Section \ref{subsection;08.9.28.1})
and hence
$(\nbigvlambda)^{(T)}_{|\pi^{-1}(Q)}
\simeq 
 \nbigp\nbige_{\pi^{-1}(Q)}$.
By using the characterization 
of the lattices (Proposition \ref{prop;08.1.22.2}),
we obtain 
$(V_{\ast}^{(T)})_{|\pi^{-1}(Q)}
\simeq 
 \nbigp_{\ast}\nbigelambda_{\pi^{-1}(Q)}$.
Hence, we obtain that
$(\Etilde,\delbar_{\Etilde},
 \thetatilde,\htilde)_{|\nbigu\setminus\Dtilde}$
is wild and good,
due to Proposition \ref{prop;07.10.13.11}.
Since the eigenvalues of
$\Res(\DD^0)$  on $\nbigq\nbigetilde^0$
are the same as those of
the residue on $\nbigq\nbigetilde^0_{|\Ctilde}$,
the harmonic bundle is $\nbiga$-good wild.
For any $\lambda\neq 0$,
we have the isomorphism
$\nbigq\nbige^{\lambda}_{|\pi^{-1}(Q)}
\simeq
 \nbigvlambda_{|\pi^{-1}(Q)}$,
and hence
we obtain the isomorphism
$\nbigq\nbige^{\lambda}
\simeq
 \nbigvlambda$.
Thus, we obtain Lemma \ref{lem;07.11.1.15}.
\hfill\qed

\subsection{End of Proof of Proposition 
 \ref{prop;07.10.28.36}}
\label{subsection;10.5.25.20}

Let us show the first claim of Proposition 
\ref{prop;07.10.28.36}.
Due to Lemma \ref{lem;07.11.1.15},
$(\Etilde,\delbar_{\Etilde},\thetatilde,\htilde)$
satisfies the assumption made in
Section \ref{subsection;07.11.1.16}.
Hence,
we can take a birational projective morphism
$\varphi':\Ztilde'\lrarr \Ztilde$
such that 
$\varphi'^{\ast}\varphi^{\ast}\theta$ is good
due to Proposition \ref{prop;07.11.1.3},
and so we can assume
$(\Etilde,\delbar_{\Etilde},\thetatilde,\htilde)$
is a good wild harmonic bundle
on $(\Ztilde,\Dtilde)$
from the beginning.
\begin{lem}
\mbox{{}}\label{lem;08.1.22.3}
\begin{itemize}
\item
For any $\lambda\neq 0$,
we have a natural isomorphism
$\nbigq\nbigelambda
\simeq
 \nbigv^{\lambda}$.
\item
$(\Etilde,\delbar_{\Etilde},\thetatilde,\htilde)$
is an $\nbiga$-good wild harmonic bundle.
\end{itemize}
\end{lem}
\pf
There exists a closed analytic subset
$W\subset D$ such that
$\nbigv^{\lambda}$ is good
on $\Ztilde-W$.
Then, we have the isomorphism
$\nbigq\nbigelambda_{|\Ztilde-W}
 \simeq
 \nbigv^{\lambda}_{|\Ztilde-W}$
because of Lemma \ref{lem;07.11.1.15}.
Then, it is extended to an isomorphism on $\Ztilde$
by the reflexivity of
$\nbigq\nbigelambda$ and
$\nbigv^{\lambda}$.
Thus, we obtain the first claim of
Lemma \ref{lem;08.1.22.3}.
The second claim also follows from
Lemma \ref{lem;07.11.1.15}.
\hfill\qed

\vspace{.1in}

Hence, 
we can conclude that
$\harmonicbundle$
$\nbiga$-wild harmonic bundle on $(Z,U)$.
Thus, the first claim of Proposition
\ref{prop;07.10.28.36} is proved.

\vspace{.1in}

In the above argument,
we obtain the $\nbiga$-good wild harmonic bundle
$(\Etilde,\delbar_{\Etilde},\thetatilde,\htilde)
=\varphi^{-1}(E,\delbar_E,\theta,h)$
on $(\Ztilde,\Dtilde)$.
We have the associated polarized $\nbiga$-wild
pure twistor $D$-module 
$(\gbigt(\Etilde),\gbigstilde)$ on $\Ztilde$
whose underlying $\nbigr_{\Ztilde}$-module 
is denoted by $\gbigetilde$.
Let $\nbigt_1$ denote the direct summand of
$\varphi^0_{\dagger}\gbigt(\Etilde)$
whose strict support is $Z$.
We have the naturally induced polarization
$\nbigs_1$ of $\nbigt_1$.
We would like to show
$(\nbigt_1,\nbigs_1)\simeq
 (\nbigt,\nbigs)$,
which is the second claim of Proposition
\ref{prop;07.10.28.36}.
Let $\nbigm_1$ denote the underlying $\nbigr$-module
of $\nbigt_1$.
We have only to show $\nbigm_1\simeq\nbigm$.

The rest is essentially the same as
the argument in Section 19.4.3 of \cite{mochi2}.
We give only an outline.
We put $r:=n-k-1$.

\begin{lem}
\label{lem;07.11.2.35}
We have a natural isomorphism
$\varphi_{\dagger}\gbigetilde(\ast g)
\simeq
 \nbigm(\ast g)$.
\end{lem}
\pf
By shrinking $U$,
we may assume to have holomorphic functions
$a_1,\ldots,a_{r}$ with the following property:
\begin{itemize}
\item
$Z$ is one of the irreducible components
of $\bigcap_{j=1}^{r} a_j^{-1}(0)$.
\item
We put $w_j:=\sum_{i=1}^n\del_{z_i}a_j\cdot \del_{z_i}$.
Then, $w_1,\ldots,w_{r}$ give
a frame of the normal bundle of $U$ in $X$.
\end{itemize}
We have the $\nbigr_U$-module $\nbige$
on $\cnum_{\lambda}\times U$,
associated to $(E,\delbar_E,\theta,h)$.
Let $\iota$ denote the immersion $U\lrarr X$.
We have only to show that
$\nbigm(\ast g)$ and 
$\varphi_{\dagger}\gbigetilde(\ast g)$
are the same in $\iota_{\dagger}\nbige$.

Let $s$ be any section of $\nbigm(\ast g)$
on $U(\lambda_0)\times X$.
When we regard $s$ as the section of
$\iota_{\dagger}\nbige$,
it is expressed as a finite sum
of the following form:
\[
 s=\sum_{\vecn\in\seisuu_{\geq 0}^{r}}
 u_{\vecn}(s)\cdot \prod_{j=1}^rw_j^{n_j}
\]
Here $u_{\vecn}(s)$ are holomorphic sections
of $\nbige$ on $U(\lambda_0)\times U$.
We can naturally regard them
as the sections of $\nbigetilde$ on 
$U(\lambda_0)\times(\Ztilde-\Dtilde)$.

\begin{lem}
\label{lem;07.11.2.31}
$u_{\vecn}(s)$ are meromorphic sections of
$\nbigq\nbigetilde$
on $U(\lambda_0)\times\Ztilde$
\end{lem}
\pf
Let $\Ctilde$ be a curve in $\Ztilde$
as in Section \ref{subsection;07.11.2.30},
and let $C=\varphi(\Ctilde)$.
We also use the other objects and the notation
in Section \ref{subsection;07.11.2.30}.
Let us see that $s$ induces
sections $s^{(j)}$
of $\nbigm^{(j)}(\ast g^{(j)})$
inductively.
We may assume that
$f_j$ are coordinate functions on $X$,
by replacing $X$ with $X\times\cnum^{k}$.
Let us consider the specialization along $f_1$.
Since $s_{|X-\{g=0\}}$ is contained in
$\Vzero_{-1}\bigl(\nbigm(\ast g)\bigr)_{|X-\{g=0\}}$,
where $\Vzero$ is the $V$-filtration along $f_1$,
we obtain $s\in \Vzero_{-1}\bigl(\nbigm(\ast g)\bigr)$.
Hence, it induces the section $s^{(1)}$ of
$\Gr^{\Vzero}_{-1}\bigl(\nbigm(\ast g)\bigr)$.
Since the support of 
$\psitilde_{f_1,u}(\nbigm)$ is contained in
$\{g=0\}$ unless $u\in\seisuu_{<0}\times\{0\}$,
we have 
$\Gr^{\Vzero}_{-1}\bigl(\nbigm(\ast g)\bigr)
=\psizero_{f_1,-\vecdelta_0}
 \bigl(\nbigm(\ast g)\bigr)$.
For a similar reason,
we obtain the natural isomorphism
$\psizero_{f_1,-\vecdelta_0}\nbigm(\ast g)
\simeq
 \nbigm^{(1)}(\ast g)$.
Thus, $s^{(1)}$ is induced.
By the same procedure,
we obtain $s^{(j+1)}$ from $s^{(j)}$.
By applying Lemma \ref{lem;07.11.2.11},
we obtain that the restrictions of
$u_{\vecn}(s)$ to $\Ctilde$
are meromorphic sections of
$\nbigq\nbigetilde_{|U(\lambda_0)\times\Ctilde}$.
Then, we can conclude that
$u_{\vecn}(s)$ are meromorphic sections
of $\nbigq\nbigetilde$ on $U(\lambda_0)\times\Ztilde$.
Thus, we obtain Lemma \ref{lem;07.11.2.31}.
\hfill\qed

\vspace{.1in}

Let us return to the proof of 
Lemma \ref{lem;07.11.2.35}.
Let $\Gamma_{\varphi}:\Ztilde\lrarr\Ztilde\times X$
denote the graph.
Each $u_{\vecn}(s)\cdot \prod_{j=1}^rw_j^{n_j}$
gives a section of
$\Gamma_{\varphi\dagger}\nbigq\nbigetilde$,
and hence it induces a section of
$\varphi_{\dagger}\nbigq\nbigetilde$
(See Lemma 14.11 of \cite{mochi2}, for example.)
Therefore, $s$ is naturally contained in
$\varphi_{\dagger}\gbigetilde(\ast g)
\simeq \varphi_{\dagger}\nbigq\nbigetilde$,
and we obtain
$\nbigm(\ast g)
 \subset\varphi_{\dagger}\gbigetilde(\ast g)$.
By using the coherence property of
$\nbigm$ and $\gbigetilde$
and Lemma \ref{lem;08.1.21.11},
we obtain $\nbigm(\ast g)=
 \varphi_{\dagger}\gbigetilde(\ast g)$.
Thus, we obtain Lemma \ref{lem;07.11.2.35}.
\hfill\qed

\vspace{.1in}

We have the natural isomorphism
$\nbigm_1(\ast g)
\simeq
 \varphi^0_{\dagger}\gbigetilde(\ast g)$,
and hence
$\nbigm_1(\ast g)
\simeq
 \nbigm(\ast g)$.
Since both $\nbigm_1(\ast g)$
and $\nbigm$ are strictly $S$-decomposable
along $g$,
we obtain that $\nbigm_1$
and $\nbigm$ are isomorphic
due to Lemma \ref{lem;07.11.2.15}.
Thus, the proof of Proposition \ref{prop;07.10.28.36}
is accomplished
\hfill\qed

\section{Application to
 algebraic semisimple holonomic $D$-modules}
\label{subsection;08.10.15.2}
\subsection{Main theorems}

Let $X$ be a smooth proper 
complex algebraic variety.
Let $Z$ be  an irreducible closed subvariety of $X$.
Let $\Hol_{Z}(X)$
denote the category of holonomic $D$-modules
whose supports are contained in $Z$.
\index{category $\Hol_{Z}(X)$}
(We consider the full subcategory
of $D$-modules,
i.e., morphisms are not necessarily isomorphisms.)
Let $\Hol_{Z,ss}(X)$ 
denote the category of semisimple holonomic $D$-modules
whose strict supports are exactly $Z$,
and let $\Hol_{Z,s}(X)$ denote the category
of simple holonomic $D$-modules
whose strict supports are $Z$.
\index{category $\Hol_{Z,s}(X)$}

For an integer $w$,
let $\MT^{\wild}_{Z,ss}(X,w,\sqrt{-1}\real)^{(p)}$
denote the category of polarizable 
$\sqrt{-1}\real$-wild pure twistor $D$-modules
of weight $w$,
whose strict supports are exactly $Z$,
i.e.,
the support of any non-zero direct summand
is $Z$.
(We consider the full subcategory of
$\MT^{\wild}(X,w)$,
i.e.,
morphisms are not necessarily isomorphisms.)
Let $\MT^{\wild}_{Z,s}(X,w,\sqrt{-1}\real)^{(p)}$
denote the category of simple polarizable
$\sqrt{-1}\real$-wild pure twistor $D$-modules
whose strict supports are $Z$.
\index{category 
 $\MT^{\wild}_{Z,ss}(X,w,\sqrt{-1}\real)^{(p)}$}
\index{category 
 $\MT^{\wild}_{Z,s}(X,w,\sqrt{-1}\real)^{(p)}$}

We have the naturally defined functor
$\Xi_{DR}:\MT^{\wild}_{Z,ss}(X,w,\sqrt{-1}\real)^{(p)}
\lrarr \Hol_{Z}(X)$ 
by taking the specialization of
the underlying $\nbigr$-modules
at $\lambda=1$.
\index{functor $\Xi_{DR}$}
We will prove the following theorem
in Sections
\ref{subsection;08.1.22.10}--\ref{subsection;08.1.22.11}.
\begin{thm}
\label{thm;07.10.14.75}
Fix an integer $w$.
The functor $\Xi_{DR}$ naturally induces
an equivalence of the categories
$\MT^{\wild}_{Z,ss}(X,w,\sqrt{-1}\real)^{(p)}$
and $\Hol_{Z,ss}(X)$.
It also induces an equivalence of the categories
$\MT^{\wild}_{Z,s}(X,w,\sqrt{-1}\real)^{(p)}$
and $\Hol_{Z,s}(X)$.
\end{thm}

Before going to the proof,
we give the following consequence,
which is one of the main motivations
for this study.

\begin{thm}[Kashiwara's conjecture]
\label{thm;08.3.6.2}
Let $X$ be a smooth complex algebraic variety,
and let $\nbigf$ be an algebraic semisimple holonomic 
$D$-module.
\begin{itemize}
\item
 Let $F:X\lrarr Y$ be a projective morphism.
 Let $L$ denote the Lefschetz morphism
 for some line bundle which is relatively ample
 with respect to $F$.
 Then, $F_{\dagger}^j(\nbigf)$ are also semisimple
 for any $j$,
 and the induced morphisms
 $L^j:F_{\dagger}^{-j}\nbigf\lrarr F_{\dagger}^j\nbigf$
 $(j\geq 0)$
 are isomorphic for $j\geq 0$.
 In particular, 
 $F_{\dagger}(\nbigf)$
 is isomorphic to
 $\bigoplus F^i_{\dagger}(\nbigf)[-i]$
 in the derived category of
 cohomologically holonomic $D_Y$-modules.
\item
 Let $g$ be an algebraic function on $X$,
 and let $\gminia\in\cnum[t_n^{-1}]$.
 Then,
 $\Gr^W\psi_{g,\gminia}(\nbigf)$
 and 
 $\Gr^W\phi_{g}(\nbigf)$ 
 are also semisimple, 
 where  
 (i) $\psi_{g,\gminia}$ denotes 
 the nearby cycle functor
 with ramification and exponential twist by $\gminia$
 (see Section {\rm\ref{subsection;08.10.25.1}}),
 (ii) $\phi_g$ denotes the vanishing cycle functor,
 (iii) $\Gr^W$ is taken with respect to
 the weight filtration of the naturally induced
 nilpotent maps.
\end{itemize}
\end{thm}
\pf
We take a smooth proper variety $\Xbar$
which contains $X$ as a Zariski open subset
such that $\Xbar-X$ is a normal crossing hypersurface,
by using Nagata's embedding and
Hironaka's resolution.
We take the minimal extension
$\overline{\nbigf}$ of $\nbigf$ on $\Xbar$.
By Theorem \ref{thm;07.10.14.75},
we can take a polarizable $\sqrt{-1}\real$-wild 
pure twistor $D$-module 
$\overline{\nbigt}$ on $\Xbar$
such that $\Xi_{DR}(\overline{\nbigt})
 \simeq \overline{\nbigf}$.
Let $\nbigt$ denote the restriction of
$\overline{\nbigt}$ to $X$.
We have 
$\Xi_{DR}(\nbigt)=\nbigf$.
Because
$F^j_{\dagger}(\nbigf)
\simeq
 \Xi_{DR}(F^j_{\dagger}(\nbigt))$,
the first claim follows from
Theorem \ref{thm;07.10.23.20}.

We also have the following:
\[
 \Gr^W\psi_{g,\gminia}(\nbigf)
\simeq\bigoplus_{\substack{u=(a,\alpha)\\-1\leq a\leq 0}}
 \Xi_{DR}\bigl(
 \Gr^W\psitilde_{g,\gminia,u}(\nbigt)
 \bigr),
\quad
 \Gr^W\phi_{g}(\nbigf)
\simeq
\Xi_{DR}\bigl(
 \Gr^W\phi_{g}(\nbigt)
\bigr)
\]
Hence, the second claim is also clear.
\hfill\qed

\subsubsection{}

For the proof of Theorem \ref{thm;07.10.14.75},
we have only to show the following
two claims:
\begin{description}
\item[(A)]
 $\Xi_{DR}(\nbigt)$ is simple
 for any $\nbigt\in 
 \MT^{\wild}_{Z,s}\bigl(X,w,\sqrt{-1}\real\bigr)^{(p)}$.
\item[(B)]
 $\Xi_{DR}$ induces an equivalence
 between the categories
 $\MT^{\wild}_{Z,s}
 \bigl(X,w,\sqrt{-1}\real\bigr)^{(p)}$
and $\Hol_{Z,s}(X)$.
\end{description}

\subsection{Construction of 
 wild pure twistor $D$-module}
\label{subsection;08.1.22.10}

Let $M$ be a simple holonomic $D$-module
whose strict support is $Z$.
There exists a smooth Zariski open subset $U\subset Z$
such that $M_{|X-(Z-U)}$ comes from
a flat bundle on $U$.
We can take a smooth projective variety
$\Ztilde$ and
a birational projective morphism
$\pi:\Ztilde\lrarr Z$
such that
 $\Mtilde:=L^{-\dim X+\dim Z}
 \pi^{\ast}M\otimes\nbigo(\ast\Dtilde)$
 is a meromorphic flat connection
on $(\Ztilde,\Dtilde)$,
where $\Dtilde:=\Ztilde-\pi^{-1}(U)$.
(See Proposition \ref{prop;07.11.2.50} 
and Lemma \ref{lem;07.11.2.1}.)
According to Theorem \ref{thm;07.10.14.60},
we may and will assume that
it is a meromorphic flat bundle,
and that the Deligne-Malgrange lattice
associated to $\Mtilde$ is good.

We have the meromorphic flat bundle 
$V=\Mtilde^{(2)}$,
obtained as the deformation of $\Mtilde$ 
by the procedure explained
in Section \ref{subsection;10.5.17.51}
with $T=2$.
(Recall that we used the deformation
in the construction of $\nbigq\nbige$.)

\begin{lem}
\label{lem;07.10.14.70}
$V$ is also simple.
\end{lem}
\pf
We have only to consider the case $\dim Z=1$
due to Mehta-Ramanathan type theorem
(Proposition \ref{prop;06.8.12.15}).
Let $V'\subset V$ be a flat subbundle
such that $0<\rank(V')<\rank(V)$.
Let $P$ be a point of $\Ztilde$,
where the connection of $V$ has the pole.
Let $U_P$ be a small neighbourhood of $P$,
and let $\varphi_P:\Utilde_P\lrarr U_P$ 
be a ramified covering such that
both $\varphi_P^{\ast}V$ and $\varphi_P^{\ast}V'$
have the unramifiedly good Deligne-Malgrange lattice.
Then, it can be easily checked that
the irregular decompositions are compatible.
Hence, we have the naturally defined morphism
$V^{\prime(1/2)}\lrarr V^{(1/2)}=\Mtilde$
with $0<\rank V^{\prime(1/2)}<\rank \Mtilde$,
which contradicts with the simplicity of $\Mtilde$.
\hfill\qed

\vspace{.1in}

Due to Theorem \ref{thm;06.1.23.100},
we can take a pluri-harmonic metric $h$
of $V_{|\Ztilde-\Dtilde}$,
which is adapted to the Deligne-Malgrange 
filtered bundle $\vecV^{DM}_{\ast}$
associated to $V$.
Let $(E,\delbar_E,h,\theta)$ be the corresponding
$\sqrt{-1}\real$-good wild harmonic 
bundle on $\Ztilde-\Dtilde$.
We have $\nbigp\nbige^1=V$
and $\nbigq\nbige^1=\Mtilde$.
Let $(\gbigt(E),\gbigs)$ be
the associated $\sqrt{-1}\real$-wild 
polarized pure twistor $D$-module
of weight $0$ on $\Ztilde$.
Let $\nbigt=(\nbigm,\nbigm,C)$ be
the $\sqrt{-1}\real$-wild polarized 
pure twistor $D$-module of weight $0$,
which is a direct summand of 
$P\pi^0_{\dagger}\gbigt(E)$
whose strict support is $Z$.
We have the induced polarization $\nbigs$
of $\nbigt$.

We would like to show
$\nbigm_{|\lambda=1}\simeq M$.
We can take a divisor $D$ of $X$
such that $(Z-U)\subset D$.
We may assume $U=Z\setminus D$
from the beginning.
We have the trace map
$\tr:\pi^0_{\dagger}\Mtilde\lrarr M(\ast D)$.
The restriction of them to $X-D$ are isomorphic.
Hence the support of
$\Ker(\tr)$ and $\Cok(\tr)$ are contained in $D$.
Both $\pi^0_{\dagger}\Mtilde$ and $M(\ast D)$
are algebraically localized on $X-D$,
and hence
we can conclude that they are isomorphic on $X$.
Thus, we obtain the isomorphism
$\nbigm_{|\lambda=1}\otimes\nbigo(\ast D)
\simeq M\otimes\nbigo(\ast D)$.
Since $M$ is simple,
it is strictly $S$-decomposable along $D$.
We know that
$\nbigm_{|\lambda=1}$ is also strictly $S$-decomposable
along $D$.
Hence, we can conclude 
$\nbigm_{|\lambda=1}\simeq M$
by the same argument as the proof of 
Lemma \ref{lem;07.11.2.15}.

\subsection{End of the proof of 
Theorem \ref{thm;07.10.14.75}}
\label{subsection;08.1.22.11}

Let $(\nbigt,\nbigs)$ be a simple polarized
$\sqrt{-1}\real$-wild pure twistor $D$-module 
of weight $0$
whose strict support is $Z$.
Let $M:=\Xi_{DR}(\nbigt)$.
We have a Zariski open subset $U\subset Z$
such that $(\nbigt,\nbigs)\in 
\MPT_{\strict}^{\wild}(Z,U,0,\sqrt{-1}\real)$.
(See Subsection \ref{subsection;10.6.6.13}
for $\MPT_{\strict}^{\wild}$.)
Let $(E,\delbar_E,\theta,h)$ be
the corresponding harmonic bundle on $U$,
which is $\sqrt{-1}\real$-wild on $(Z,U)$.
By shrinking $U$,
we may assume to have a divisor $D$ of $X$
such that $Z-U=Z\cap D$.
We can take a smooth projective variety $\Ztilde$
with a projective birational morphism
$\varphi:\Ztilde\lrarr Z$
such that 
(i) $\Dtilde:=\varphi^{-1}(Z-U)$ is 
 simply normal crossing,
(ii) $(\Etilde,\delbar_{\Etilde},\thetatilde,\htilde):=
 \varphi^{-1}(E,\delbar_E,\theta,h)$ 
is a $\sqrt{-1}\real$-good wild harmonic bundle
 on $(\Ztilde,\Dtilde)$.
We have the associated
polarized $\sqrt{-1}\real$-wild pure twistor $D$-module
$(\gbigt(\Etilde),\gbigstilde)$ on $\Ztilde$.
We know that
$\nbigt$ is isomorphic to the direct summand of
$P\varphi_{\dagger}^0\gbigt(\Etilde)$
by the correspondence 
in Theorem \ref{thm;07.10.28.30}.

Due to Theorem \ref{thm;07.10.15.1},
$(\nbigp\nbigetilde^1,\DD^1)$ is semisimple.
By using Proposition \ref{prop;07.12.6.11},
we can take an orthogonal decomposition
$(\Etilde,\delbar_{\Etilde},\thetatilde,\htilde)
=\bigoplus 
 (\Etilde_i,\delbar_{\Etilde_i},
 \thetatilde_i,\htilde_i)$
such that 
$(\nbigp\nbige_i,\DD_i^1)$ associated to
$(\Etilde_i,\delbar_{\Etilde_i},
 \thetatilde_i,\htilde_i)$
are simple.
We have the decomposition
$\gbigt(\Etilde)=\bigoplus_i \gbigt(\Etilde_i)$,
and 
$P\varphi_{\dagger}^0\gbigt(\Etilde)
=\bigoplus P\varphi_{\dagger}^0\gbigt(\Etilde_i)$,
which induces a decomposition of 
$(\nbigt,\nbigs)$.
Since we have assumed $(\nbigt,\nbigs)$ is simple,
we obtain that
$(\nbigp\nbigetilde^1,\DD^1)$ is also simple.
(Note that the restriction of $\varphi$ 
is isomorphic on a generic part of $Z$.)
Using the argument in the proof of Lemma 
\ref{lem;07.10.14.70},
we obtain that $\nbigq\nbige^1$ is simple.
By construction of $M$,
we have the isomorphism
$M\otimes\nbigo(\ast D)
\simeq
 \varphi^{0}_{\dagger}
 \nbigq\nbigetilde^1\otimes\nbigo(\ast D)$
induced by the trace morphism.
Hence, the support of any non-trivial submodule of $M$
is contained in $D$.
Because $M$ is strictly $S$-decomposable 
along any function $g$ on $X$ such that
$Z\not\subset g^{-1}(0)$,
we obtain that $M$ is also simple.
Thus, the claim (A) is proved.

\vspace{.1in}
By the construction in 
Section \ref{subsection;08.1.22.10},
we obtain the essential surjectivity in the claim of (B).
Let us show the fully faithfulness.
Let $(\nbigt_i,\nbigs_i)$ $(i=1,2)$ be
simple polarized 
$\sqrt{-1}\real$-wild pure twistor $D$-modules
of weight $0$.
We put $M_i:=\Xi_{DR}(\nbigt_i)$.
Let $(E_i,\delbar_{E_i},\theta_i,h_i)$ be
$\sqrt{-1}\real$-wild harmonic bundles
on $(Z,U)$ underlying $(\nbigt_i,\nbigs_i)$.
Let $\varphi:\Ztilde\lrarr Z$
be a birational projective morphism
such that
$(\Etilde_i,\delbar_{\Etilde_i},
 \thetatilde_i,\htilde_i):=
 \varphi^{\ast}(E_i,\delbar_{E_i},\theta_i,h_i)$
are $\sqrt{-1}\real$-good wild harmonic bundle
on $(\Ztilde,\Dtilde)$,
where $\Dtilde$ denotes some normal crossing divisor 
of $\Ztilde$.
Let $(\nbigq\nbigetilde^1_i,\DDtilde_i^1)$
denote the meromorphic flat bundles
associated to $(\Etilde_i,\delbar_{\Etilde_i},
 \thetatilde_i,\htilde_i)$.

Let $f:M_1\simeq M_2$ be a morphism.
By Lemma \ref{lem;08.1.22.3},
we obtain the induced isomorphism
$(\nbigq\nbigetilde^1_1,\DDtilde_1^1)
\simeq
 (\nbigq\nbigetilde_2^1,\DDtilde_2^1)$
as meromorphic flat bundles.
It induces an isomorphism
of the associated Deligne-Malgrange
filtered flat bundles
$\nbigq_{\ast}\nbigetilde^1_{1}
\simeq
 \nbigq_{\ast}\nbigetilde^1_{2}$.
It induces
$\nbigp_{\ast}\nbigetilde^1_{1}
\simeq
 \nbigp_{\ast}\nbigetilde^1_{2}$.
Hence, we obtain
$\ftilde:
 (\Etilde_1,\delbar_{\Etilde_1},\thetatilde_1,\htilde_1)
\simeq
 (\Etilde_2,\delbar_{\Etilde_2},\thetatilde_2,\htilde_2)$.
(See Corollary \ref{cor;08.1.22.15}
or Proposition \ref{prop;07.12.6.11}.)
Because of the correspondence in Theorem 
\ref{thm;07.10.28.30},
we obtain an isomorphism
$F:(\nbigt_1,\nbigs_1)\simeq
 (\nbigt_2,\nbigs_2)$.
By construction,
we have the coincidence of
the restrictions of $\Xi_{DR}(F)$
and $f$ to a Zariski open subset of $Z$.
Because $M_i$ are simple, 
we obtain $\Xi_{DR}(F)=f$,
i.e.,
$\Xi_{DR}$ is full.

Let $F:\nbigt_1\lrarr\nbigt_2$
be a morphism such that
$\Xi_{DR}(F)=0$.
We obtain the vanishing of the induced morphism
of variations of pure twistor structure on $U$.
Then, we obtain $F=0$.
Hence, $\Xi_{DR}$ is fully faithful.
Thus, the proof of Theorem \ref{thm;07.10.14.75}
is finished.
\hfill\qed

\part{Appendix}

\chapter{Preliminary from Analysis on Multi-sectors}
In Section \ref{subsection;08.9.28.123},
we give a generalization
of Hukuhara-Malmquist type theorem.
Namely, we study a lifting of a formal solution
of some non-linear differential equation
to a solution on small open multi-sectors.
It is one of the fundamental tools
in the study of Stokes structure in Chapters 
\ref{section;07.12.26.6}--\ref{chapter;10.5.26.11}.
For the proof,
we use a classical argument
explained in \cite{wasow}.

In Section \ref{subsection;08.9.28.130},
we give some estimates of sections on a sector.
Lemma \ref{lem;07.7.19.102} will be 
used for comparison of irregular decompositions
(Sections \ref{subsection;07.10.7.21}
and \ref{subsection;08.9.15.41}).
Lemmas \ref{lem;07.7.22.20}
and \ref{lem;07.7.22.21}
will be used in the study of $L^2$-cohomology
associated to wild harmonic bundle on curves
(Chapter \ref{section;07.10.7.1}).

In Section \ref{subsection;08.9.28.131},
we give some estimates of the growth order of solutions
of some differential equations on multi-sectors.
In Section \ref{subsection;07.6.4.6},
we give an estimate of the growth order
of a section whose derivative rapidly decays.
In Section \ref{subsection;08.9.28.132},
we give an estimate of the growth order
of a flat section.
We reformulate it in Sections 
\ref{subsection;08.9.28.133}--\ref{subsection;08.9.28.134}.
These results will be used implicitly
in many places.
For example, it will be used
to give a characterization of Stokes filtrations.

\section{Hukuhara-Malmquist type theorem}
\label{subsection;08.9.28.123}
\subsection{Statement}
\label{subsection;07.6.4.5}

We study a lifting of a formal solution
of some non-linear differential equation
to a solution on small open multi-sectors.
See \cite{majima} and \cite{wasow}
for the history and the classical arguments
in the one dimensional case.
For the higher dimensional case,
it was studied in \cite{majima} more generally
but in a different manner.
We give some statements and 
the outline of a proof
in our convenient way.
Although we will use it for analysis
on a multi-sector around $0$
(Chapters
\ref{section;07.12.26.6}--\ref{chapter;10.5.26.11}),
it is more conventional and easier to consider it
on a multi-sector around $\infty$.
The translation can be done 
in a straightforward way.

We set $S^1:=\bigl\{w\in\cnum\,\big|\,|w|=1\bigr\}$
for which we use a polar coordinate.
Let $(\theta^{(0)},\theta^{(0)}_1,\ldots,\theta^{(0)}_n)$
be a point of $(S^1)^{n+1}$.
Let $W$ be a compact region in $\cnum^M$,
where $M$ denotes some non-negative integer.
For given 
\[
 R\in\real_{>0},
\quad
\vecR_y=(R_1,\ldots,R_n)\in\real_{>0}^n,
\quad
\epsilon\in\openopen{0}{2\pi},
\quad
\vecepsilon_y
=(\epsilon_1,\ldots,\epsilon_n)
 \in \openopen{0}{2\pi}^{n},
\]
let $S(R,\epsilon,\vecR_y,\vecepsilon_y)$ denote
the multi-sector in 
$\cnum^{\ast}\times
 (\cnum^{\ast})^n\times W$:
\begin{multline*}
 S(R,\epsilon,\vecR_y,\vecepsilon_y):=
 \bigl\{x\in\cnum\,\big|\,
 |x|>R,|\arg(x)-\theta^{(0)}|<\epsilon\bigr\}
\times \\
 \prod_{i=1}^n
 \bigl\{y\in\cnum\,\big|\,
 |y_i|>R_i,|\arg(y_i)-\theta^{(0)}_i|<\epsilon_i\bigr\}
\times W
\end{multline*}
We fix a multi-sector
$S\bigl(R^{(0)},\epsilon^{(0)},
 \vecR^{(0)}_y,\vecepsilon_y^{(0)}\bigr)$
for some given $R^{(0)}$, $\epsilon^{(0)}$,
$\vecR^{(0)}$ and $\vecepsilon^{(0)}_y$.
Let $a$ be a non-negative integer
such that $2\epsilon^{(0)}\, (a+1)<\pi$,
and let $\vecb=(b_1,\ldots,b_n)\in\seisuu_{>0}^n$.
Let  $\lambda_i(w)$ $(i=1,\ldots,d)$ 
be holomorphic functions on $W$
satisfying the following condition:
\begin{itemize}
\item
 Let us consider the region
\[
 H(\epsilon^{(0)},\vecepsilon_y^{(0)}):=
 \bigl\{(\theta,\theta_1,\ldots,\theta_n)
 \in (S^1)^{n+1}
 \,\big|\,|\theta-\theta^{(0)}|\leq \epsilon^{(0)},
 |\theta_i-\theta_i^{(0)}|\leq\epsilon_i^{(0)}
 \bigr\}
\]
and the functions
\[
 F_j(\theta,\theta_1,\ldots,\theta_n,w):=
 \Re\Bigl(\lambda_j(w)\,
 \exp\bigl(\sqrt{-1}
 \bigl(\sum b_i\theta_i+(a+1)\theta\bigr)\bigr)
 \Bigr)
\]
for $j=1,\ldots,d$.
Then, the following holds for each $j$:
 \begin{itemize}
\item
 If $\{F_j=0\}\cap
 \bigl(H(\epsilon^{(0)},\vecepsilon^{(0)}_y)\times W\bigr)$
 is not empty,
 it is connected and contained in
 $H(\epsilonbar_j^{(0)},\vecepsilon^{(0)}_y)\times W$
 for some $\epsilonbar_j^{(0)}<\epsilon^{(0)}$.
 \end{itemize}

\end{itemize}
Let us consider the following differential equation
for $\cnum^d$-valued functions $u$
on $S(R^{(0)},\epsilon^{(0)},
 \vecR^{(0)}_y,\vecepsilon_y^{(0)})$:
\begin{equation}
 \label{eq;07.6.14.2}
 x^{-a}\frac{\del u}{\del x}=\Lambda\cdot u
+p\bigl(x,\vecy,w,u(x,\vecy,w)\bigr)
\end{equation}

\begin{itemize}
\item
$\Lambda$ is a diagonal $d$-th square matrix
whose $(j,j)$-th entries are of the form
$\lambda_j(\vecy,w)\vecy^{\vecb}$.
\begin{itemize}
\item
 $\vecy^{\vecb}:=\prod_{i=1}^n y_i^{b_i}$.
\item
$\lambda_j(\vecy,w)$ are holomorphic functions
independent of the variable $x$,
and satisfy
\[
 \bigl|
 \lambda_j(\vecy,w)-\lambda_j(w)
\bigr|
=O\Bigl(\sum_{i=1}^n|y_i|^{-1}\Bigr).
\]
\end{itemize}
\item
 $p(x,\vecy,w,u)=
 \bigl(p_j(x,\vecy,w,u)\,\big|\,j=1,\ldots,d\bigr)$
 is a $\cnum^d$-valued
 holomorphic function
 on $S(R^{(0)},\epsilon^{(0)},
 \vecR^{(0)}_y,\vecepsilon_y^{(0)})\times\cnum^d$,
 and  it is decomposed into
\[
 p^{(0)}(x,\vecy,w)+p^{(1)}(x,\vecy,w)\cdot u
+p^{(2)}(x,\vecy,w,u)
\]
with the following property:
\begin{itemize}
\item
 $p^{(0)}(x,\vecy,w)=O\bigl(|x|^{-m}|\vecy|^{-m}\bigr)$ 
 for any $m$.
 Here, $|\vecy|$ denotes $\prod_{i=1}^n|y_i|$.
\item
 $p^{(1)}(x,\vecy,w)\cdot u$ is linear with respect to 
 $u=(u_1,\ldots,u_d)$.
 We have 
 $|p^{(1)}(x,\vecy,w)|\leq
 g(x)\, |\vecy^{\vecb}|$,
 where $g(x)\to 0$  in $|x|\to \infty$.
\item
 $p^{(2)}(x,\vecy,w,u)$ is polynomial with respect to $u$,
 and the coefficients are
 $O\bigl(|x|^{-m}\, |\vecy|^{-m}\bigr)$
 for any $m$.
 We also assume that it does not contain
 the constant and linear terms with respect to $u$.
 \end{itemize}
\end{itemize}

Take a positive number $\epsilon'<\epsilon^{(0)}$.
We assume $\epsilon'>\epsilonbar_j^{(0)}$,
if $\{F_j=0\}\cap
 \bigl(H(\epsilon^{(0)},\vecepsilon^{(0)}_y)\times W\bigr)$
 is not empty.
We will show the following proposition
in Sections 
\ref{subsection;08.9.2.20}--\ref{subsection;07.9.23.15}.

\begin{prop}
\label{prop;07.6.15.3}
If $R$ and $\vecR_y^{(0)}$ are sufficiently large,
we have a solution $u$ of {\rm(\ref{eq;07.6.14.2})}
on $S(R,\epsilon',\vecR_y^{(0)},\vecepsilon_y^{(0)})$,
such that $|u|=O\bigl(|x|^{-m}|\vecy|^{-m}\bigr)$
for any $m>0$.
\end{prop}

We follow the argument in Chapter 14 of 
the standard text book \cite{wasow}.
We recommend the reader to read it
to understand the idea.
In the following argument,
we may and will assume
$\theta^{(0)}=\theta^{(0)}_i=0$.
We will also consider only the case $a=0$,
because
the general case can easily be reduced to the case $a=0$
by the change of the variables $\xi=x^{a+1}$.
(See {\rm \cite{wasow}}.)

\subsubsection{Remark for 
 simplification of the proof}

Assume the following:
\begin{description}
\item[($\ast$)]
$\{F_j=0\}\cap
 \bigl(H(\epsilon^{(0)},\vecepsilon^{(0)}_y)\times W\bigr)$
 is not empty for each $j=1,\ldots,d$
\end{description}
Then, the proof below can be simplified,
i.e.,
the arguments in Sections
\ref{subsection;07.9.23.10}--\ref{subsection;07.9.23.11}
are not necessary:
Since we have only to consider the case {\bf (B)}
in Section \ref{subsection;08.9.2.20}
under the assumption $(\ast)$,
the integral equation 
$\nbigp_{x_1}(u)=u$ is independent 
of the choice of $x_1$,
which implies that
$U(m_0)$ in Section \ref{subsection;07.9.23.5}
is independent of $m_0$.
Hence, we obtain 
$U(m_0)=O\bigl(
 |x|^{-m}\,|\vecy|^{-N}
 \bigr)$ for any $m$
in the end of Section \ref{subsection;07.9.23.5}.
Actually, we have the uniqueness in 
Proposition \ref{prop;07.6.15.3}
if the condition $(\ast)$ is satisfied.

\vspace{.1in}
We will apply Proposition \ref{prop;07.6.15.3}
for lifting of formal irregular decomposition
to a decomposition on multi-sectors,
for example,
in Sections \ref{subsection;08.9.2.25},
and \ref{subsection;10.5.26.12}.
For this application,
we have only to consider the case in which
the condition $(\ast)$ is satisfied.
Although we keep the proof for generality,
the reader can skip Sections
\ref{subsection;07.9.23.10}--\ref{subsection;07.9.23.11}.

\subsection{Integral transforms}
\label{subsection;08.9.2.20}

We put 
$T:=\bigl\{z\,\big|\, |\arg(z)|\leq \epsilon^{(0)} \bigr\}$.
We take a positive number $x_1$.
Let $x_1+T$ denote the subset $\bigl\{
 x_1+z\,\big|\,z\in T  \bigr\}$ of $\cnum$.
Note that the set
$\bigl\{x\,\big|\,
 |\arg(x)|<\epsilon',\,|x|>R\bigr\}$
is contained in $x_1+T$,
if $R$ is sufficiently large.
For each $x\in x_1+T$,
we take a tuple of paths
$\Gamma(x,x_1)=
 \bigl(\gamma_{1}(x,x_1),\ldots,\gamma_{d}(x,x_1)\bigr)$
contained in $x_1+T$
as in Chapter 14.3 of \cite{wasow}:
\begin{description}
\item[(A)]
If $F_j<0$ on
 $H(\epsilon^{(0)},\vecepsilon^{(0)}_y)$,
let $\gamma_{j}(x,x_1)$ be the segment
connecting $x_1$ to $x$.
\item[(B)]
Otherwise, 
we can take some $\theta^{(1)}$
such that 
 $F_j(\theta^{(1)},\theta_1,\ldots,\theta_n,w)>0$
 for any $|\theta_i|<\epsilon_i^{(0)}$
 and $w\in W$.
Let $\gamma_{j}(x,x_1)$ denote the path connecting 
$\infty$ and $x$ on the line
$\bigl\{x+s\,\exp(\sqrt{-1}\theta^{(1)})\,\big|\,
 s>0\bigr\}$.
 In this case, $\gamma_{j}(x,x_1)$ is independent of $x_1$.
\end{description}

We put 
$ S(x_1,\vecR_y):=
 (x_1+T)
\times
 \prod_{i=1}^n
 \bigl\{y\in\cnum\,\big|\,
 |y_i|>R_i,|\arg(y_i)|<\epsilon^{(0)}_i\bigr\}
\times W$.
For a $\cnum^d$-valued holomorphic function
$v=(v_1,\ldots,v_d)$ on $S(x_1,\vecR^{(0)}_y)$,
we consider the following integral,
if it is convergent:
\begin{multline*}
 \int_{\Gamma(x,x_1)}
 \exp\bigl((x-t)\Lambda\bigr)\,
 v(t,\vecy,w)\, dt:= \\
 \left(
 \int_{\gamma_{j}(x,x_1)}\!\!\!
  \exp\bigl((x-t)\,\lambda_j(\vecy,w)
 \,\vecy^{\vecb}\bigr)
 \, v(t,\vecy,w)\, dt\,\,\Big|\,\,
 j=1,\ldots,d
\right)
\end{multline*}
Then, we consider the following integral transform:
\[
 \nbigp_{x_1}(v)(x,\vecy,w):=
 \int_{\Gamma(x,x_1)}
 \exp\bigl((x-t)\Lambda\bigr)
 \,p\bigl(t,\vecy,w,v(t,\vecy,w)\bigr)\, dt
\]
We will construct the solution
of the integral equation $\nbigp_{x_1}(u)=u$
in Sections
\ref{subsection;07.9.23.20}--\ref{subsection;07.9.23.21}.

\subsection{Preliminary estimate}
\label{subsection;07.9.23.20}

We put $\lambda_0:=\min_{j,\vecw}|\lambda_j(\vecw)|$.
If $R^{(0)}_i$ are sufficiently large,
there is a positive constant $\mu_0$
such that the following holds for any
$(t,\vecy,w)\in
 \gamma_j(x,x_1)
 \times
 \prod_{i=1}^n\bigl\{y\in\cnum\,\big|\,
 |y_i|>R_i^{(0)},\,|\arg(y_i)|<\epsilon_i^{(0)}\bigr\}
 \times W$.
(See Lemma 14.1 in \cite{wasow}):
\begin{equation}
 \label{eq;07.9.22.20}
 \Re\Bigl(
 (x-t)\,\lambda_j(\vecw)\, \vecy^{\vecb}
 \Bigr)
\leq
 -|x-t|\,
 \lambda_0\,\mu_0\,\bigl|\vecy^{\vecb}\bigr|
\end{equation}
Then, we can show the following lemma 
by using the argument 
in the proof of Lemma 14.2 in \cite{wasow}.
\begin{lem}
 \label{lem;07.6.14.1}
There exists a constant $K_m>0$,
which is taken independently from 
large $\vecR_y$ and $x_1$,
with the following property:
\begin{itemize}
\item
 Let $\chi(x,\vecy,w)$ be a $\cnum^d$-valued
 holomorphic function on 
 $S(x_1,\vecR_y)$ such that
\[
 \bigl|\chi(x,\vecy,w)\bigr|\leq 
 c\,|x|^{-m}|\vecy|^{-N}
\]
for some $c>0$, $m>0$ and $N>0$.
We put
\[
 \psi(x,\vecy,w):=\int_{\Gamma(x,x_1)}\!\!\!
 \exp\bigl(
 (x-t)\, \Lambda \bigr)\,\chi(t,\vecy,w)\, dt.
\]
Then, the following inequality holds:
\[
 \bigl|\psi(x,\vecy,w)\bigr|\leq
 \left(
 \frac{K_m}{|\vecy^{\vecb}|}
\right)\, c\,|x|^{-m}|\vecy|^{-N}
\]
\hfill\qed
\end{itemize}
\end{lem}

We may assume $K_m\leq K_{m+1}$ holds.
For each $m$,
we take a small $\gamma_m>0$ such that
$\gamma_m\, K_m<1/2$
and $\gamma_m\geq \gamma_{m+1}$.
We can take $x_1^{(1)}(m)$ such that
the following holds for any 
$(x,\vecy,w)\in S\bigl(x_1^{(1)}(m),\vecR_y^{(0)}\bigr)$:
\[
 \bigl|p_1(x,\vecy,w)\bigr|
\leq
 \frac{1}{2}\gamma_m\, |\vecy^{\vecb}|
\]
Take a constant $C_0$ which is independent of $m$ and $N$.
We can take $x_1^{(2)}(m)>x_1^{(1)}(m)$
such that the following holds
for any $|u^{(i)}|\leq C_0$ $(i=1,2)$
and $(x,\vecy,w)\in S\bigl(x_1^{(2)}(m),\vecR_y^{(0)}\bigr)$:
\[
 \Bigl|
 p_2(x,\vecy,w,u^{(1)})
-p_2(x,\vecy,w,u^{(2)})
 \Bigr|
\leq\frac{1}{2}\gamma_m\,
 \bigl|u^{(2)}-u^{(1)}\bigr|
\]
Thus, 
we obtain the following
for any $|u^{(i)}|\leq C_0$ $(i=1,2)$
and $(x,\vecy,w)\in S(x^{(2)}_1(m),\vecR^{(0)}_y)$:
\[
 \Bigl|
 p(x,\vecy,w,u^{(1)})
-p(x,\vecy,w,u^{(2)})
 \Bigr|
\leq\gamma_m\,
 \bigl|u^{(2)}-u^{(1)}\bigr|
\]

\subsection{Construction of limit for each $(m,N)$}
\label{subsection;07.9.23.21}

For each $m>0$ and $N>0$,
we can take $C_{m,N}>0$ such that
$\bigl|p_0(x,\vecy,w)\bigr|<
 C_{m,N}\, |x|^{-m}\, |\vecy|^{-N}$.
We can take $x^{(3)}_1(m,N)>x^{(2)}_1(m)$
such that the following inequality holds on
$S\bigl(x_1^{(3)}(m,N),\vecR^{(0)}_y\bigr)$:
\begin{equation}
 \label{eq;07.9.22.30}
 \frac{1}{1-\gamma_mK_m}
 \frac{K_m}{|\vecy^{\vecb}|}
 \, C_{m,N}
 \, |x|^{-m}\, |\vecy|^{-N}
\leq C_0,
\quad\quad
C_{m,N}
 \, |x|^{-m}\, |\vecy|^{-N}
\leq C_0.
\end{equation}

We put $v_0:=0$
and $v_i:=\nbigp_{x_1^{(3)}(m,N)}(v_{i-1})$
for $i\geq 1$, inductively.

\begin{lem}
\label{lem;07.6.18.3}
The limit 
$v_{\infty}(m,N):=\lim_{i\to\infty}v_i$
exists
on $S\bigl(x^{(3)}(m,N),\vecR^{(0)}_y\bigr)$,
which satisfies the following estimate
on $S\bigl(x^{(3)}(m,N),\vecR^{(0)}_y\bigr)$:
\[
 \bigl|v_{\infty}(m,N)\bigr|
\leq
 \frac{1}{1-\gamma_mK_m}
 C_{m,N}
 \frac{K_m}{|\vecy^{\vecb}|}
\, |x|^{-m}\,|\vecy|^{-N}
\]
\end{lem}
\pf
Using an inductive argument,
we can show the following:
\[
 |v_{i+1}-v_i|
\leq
 \bigl(\gamma_m K_m \bigr)^{i}
 \, C_{m,N}\,
 \frac{K_m}{|\vecy^{\vecb}|}
 \, |x|^{-m}\, |\vecy|^{-N}
\]
\[
 |v_{i+1}|
\leq 
 \frac{1}{1-\gamma_m K_m}
 \, C_{m,N}
 \, \frac{K_m}{|\vecy^{\vecb}|}
 \, |x|^{-m}\, |\vecy|^{-N}
\]
Then, the claim follows.
\hfill\qed

\subsection{Estimate for any $N$ 
and some fixed $m_0$}
\label{subsection;07.9.23.5}

We fix $m_0$ and $N_0$,
and we put $U(m_0):=v_{\infty}(m_0,N_0)$
which is given on 
$S\bigl(x^{(3)}(m_0,N_0),\vecR^{(0)}_y\bigr)$.
\begin{lem}
$\bigl|U(m_0)\bigr|\leq
 C'_{N}\,|x|^{-m_0}\, |\vecy|^{-N}$ 
on $S\bigl(x^{(3)}(m_0,N_0),\vecR^{(0)}_y\bigr)$
for any $N$.
\end{lem}
\pf
Let $C_{m_0,N}$ be the constant as in Section
\ref{subsection;07.9.23.21}.
For each $N$,
we take large $R^{(3)}_i(m_0,N)$ $(i=1,\ldots,n)$
such that the following holds on the region
$S(m_0,N):=
 \bigcup_{i=1}^n
 S\bigl(x_1^{(3)}(m_0,N_0),\vecR^{(3)}_{y,i}(m_0,N)\bigr)$:
\[
 \frac{1}{1-\gamma_{m_0}K_{m_0}}
 C_{m_0,N}\,
 \frac{K_{m_0}}{|\vecy^{\vecb}|}
 \, |x|^{-m_0}\, |\vecy|^{-N}
\leq C_0,
\quad
 C_{m_0,N}\, |x|^{-m_0}\, |\vecy|^{-N}
\leq C_0
\]
Here, $\vecR^{(3)}_{y,i}(m_0,N)$ denote the tuple
such that $j$-th entry is $R^{(0)}_j$ for $(j\neq i)$
and the $i$-th entry is $R^{(3)}_i(m_0,N)$.

We put $v_0=0$ and
$v_i:=\nbigp_{x^{(3)}_1(m_0,N_0)}(v_{i-1})$
for $i\geq 1$, inductively.
Using the same argument as that 
in the proof of Lemma \ref{lem;07.6.18.3},
we can show the existence of the limit
$U(m_0)^{(N)}:=\lim_{i\to\infty}v_i$
which satisfies the following:
\[
 \bigl|U(m_0)^{(N)}\bigr|
\leq
 \frac{1}{1-\gamma_{m_0}K_{m_0}}
 C_{m_0,N}\,
 \frac{K_{m_0}}{|\vecy^{\vecb}|}
 \,|x|^{-m_0}\,|\vecy|^{-N}
\]
By construction,
the restriction of $U(m_0)$ to 
$S\bigl(x_1^{(3)}(m_0,N_0),
 \vecR^{(3)}_{y,i}(m_0,N)\bigr)$
is equal to $U(m_0)^{(N)}$.
Hence,  we obtain
$\bigl|U(m_0)\bigr|\leq
 C''_{m_0,N}|x|^{-m_0}|\vecy|^{-N}$
on $S(m_0,N)$.
Thus, we are done.
\hfill\qed

\vspace{.1in}

To obtain the estimate
$\bigl|U(m_0)\bigr|\leq
 C_{1,m,N}|x|^{-m}|\vecy|^{-N}$
for any $m$ and $N$,
we would like to compare 
the two solutions $U(m_0)$ and $U(m)$
of the differential equation (\ref{eq;07.6.14.2})
in the case $a=0$.
Note that we do not have the uniqueness of the solutions
of (\ref{eq;07.6.14.2}) in general,
and that 
the integral equation $\nbigp_{x_1}(u)=u$
depends on the choice of $x_1$.
(As remarked in the end of 
Subsection \ref{subsection;07.6.4.5},
if the condition $(\ast)$ is satisfied,
 the proof is completed in this stage.)
We give some preparations
in Sections
\ref{subsection;07.9.23.10}--\ref{subsection;07.9.23.11}
to estimate the ambiguity of the solutions
of (\ref{eq;07.6.14.2}),
and we will finish the proof of Proposition
\ref{prop;07.6.15.3}
in Section \ref{subsection;07.9.23.15}.

In the following argument,
we may and will assume
$x_1^{(i)}(m)\leq x_1^{(i)}(m+1)$
for $i=1,2,3$.

\subsection{Estimate of a solution
 of the linear differential equation}
\label{subsection;07.9.23.10}

Take any $x_2(m)>x_1^{(1)}(m)$.
Let us consider the following linear differential equation
for $\cnum^d$-valued holomorphic functions $u$
on $S(x_2(m),\vecR^{(0)}_y)$:
\begin{equation}
 \label{eq;07.9.23.1}
 \frac{\del u}{\del x}
=\Lambda\cdot u+p_1(x,\vecy,w)\cdot u
\end{equation}
We also consider the following integral transform:
\[
 \nbigr_{x_2(m)}u
=\int_{\Gamma(x,x_2(m))}
 \exp\Bigl((x-t)\Lambda\Bigr)\,
 p_1(t,\vecy,w)\cdot u(t,\vecy,w)\, dt
\]
The following lemma is easy to show.
\begin{lem}
\label{lem;07.9.23.1}
Let $m_1\leq m$.
If $|u|\leq C\, |x|^{-m_1}|\vecy|^{-N}$,
we have the following:
\[
 \bigl|\nbigr_{x_2(m)}(u)\bigr|
\leq
\frac{1}{2}
 (K_{m_1}\gamma_{m_1})\,
 C\,|x|^{-m_1}\, |\vecy|^{-N}
\]
Hence, if $u$ satisfies
$|u|\leq C\,|x|^{-m_1}\, |\vecy|^{-N}$
and $\nbigr_{x_2(m)}(u)=u$,
we obtain $u=0$.
\end{lem}
\pf
Recall that we have assumed
$x_1^{(1)}(m_1)\leq x_1^{(1)}(m)$.
The first claim immediately follows from
our choice of $K_{m_1}$ and $\gamma_{m_1}$.
If $u$ satisfies 
$|u|\leq C\,|x|^{-m_1}\, |\vecy|^{-N}$
and $\nbigr_{x_2(m)}(u)=u$,
we obtain
\[
 |u|\leq C(K_{m_1}\gamma_{m_1}/2)^{j}
 \,|x|^{-m_1}\, |\vecy|^{-N}
\]
for any $j$ using an inductive argument.
Hence, we obtain $u=0$.
\hfill\qed

\vspace{.1in}

Let $u$ be a solution of the equation (\ref{eq;07.9.23.1})
on $S\bigl(x_2(m),\vecR^{(0)}_y\bigr)$
satisfying
$|u|\leq C_N\, |x|^{-m_1}\, |\vecy|^{-N}$
for some $m_1$ and $N$.

\begin{lem}
\label{lem;07.9.23.4}
We have the estimate 
$|u|\leq C_{N}'\, |x|^{-m}\, |\vecy|^{-N}$.
\end{lem}
\pf
We have nothing to prove
in the case $m_1\geq m$,
and so we assume $m_1<m$.
We have the estimate
$\bigl| \nbigr_{x_2(m)}u\bigr|
\leq
 \gamma_{m_1}\, K_{m_1}\,
 C_N\, |x|^{-m_1}\, |\vecy|^{-N}$.
It is easy to show the following equality:
\[
 \frac{\del}{\del x}
 \bigl(\nbigr_{x_2(m)}(u)-u\bigr)
=\Lambda\bigl(
 \nbigr_{x_2(m)}(u)-u
 \bigr)
\]
Then, we have the estimate
$\nbigr_{x_2(m)}u-u
=O\bigl(
 \exp(-\eta\,|x|\, |\vecy^{\vecb}|)
 \bigr)$ for some $\eta>0$.
In particular,
we have the following:
\[
 \bigl|
 \nbigr_{x_2(m)}u-u
\bigr|
 \leq C_{10}\,
 |x|^{-m_1}\,|\vecy|^{-N}.
\]
By an easy induction,
we can show the following:
\[
\bigl|
 \nbigr_{x_2(m)}^{i+1}(u)
-\nbigr_{x_2(m)}^i(u)
\bigr|
 \leq (\gamma_{m_1}\, K_{m_1})^i\, 
 C_{10}\, |x|^{-m_1}\, |\vecy|^{-N}.
\]
Hence, we have the limit
$\nbigr^{\infty}_{x_2(m)}(u):=
 \lim_{i\to\infty}\nbigr_{x_2(m)}^i(u)$
which satisfies the following:
\[
 \nbigr_{x_2(m)}\bigl(
 \nbigr^{\infty}_{x_2(m)}(u)
 \bigr)
=\nbigr^{\infty}_{x_2(m)}(u),
\]
\[
  \bigl|
  \nbigr^{\infty}_{x_2(m)}(u)
 \bigr|
\leq
\left(
 C_N+
 \frac{C_{10}}{1-\gamma_{m_1}K_{m_1}}
\right)
 |x|^{-m_1}\, |\vecy|^{-N}
\]
Hence, we can conclude
$\nbigr^{\infty}_{x_2(m)}(u)=0$
because of Lemma \ref{lem;07.9.23.1}.

Because
$\nbigr_{x_2(m)}u-u
=O\bigl(
 \exp(-\eta\,|x|\, |\vecy^{\vecb}|)
 \bigr)$,
we also have 
\[
 \bigl|
 \nbigr_{x_2(m)}u-u
\bigr|
 \leq C_{11}\,
 |x|^{-m}\,|\vecy|^{-N}.
\]
We can show 
$\bigl|
 \nbigr_{x_2(m)}^{i+1}(u)
-\nbigr_{x_2(m)}^i(u)
\bigr|
 \leq (\gamma_{m}\, K_{m})^i\, 
 C_{11}\, |x|^{-m}\, |\vecy|^{-N}$
by an easy induction.
Hence, we have the following:
\[
 |u|=
 \bigl| u-\nbigr^{\infty}_{x_2(m)}(u)
 \bigr|
\leq 
 \frac{C_{11}}{1-\gamma_mK_m}
 |x|^{-m}|\vecy|^{-N}
\]
Thus, the proof of Lemma \ref{lem;07.9.23.4}
is finished.
\hfill\qed

\subsection{Existence of a small solution
 of some equation}

Let us consider the following differential equation
on $S\bigl(x_2(m),\vecR^{(0)}_y\bigr)$:
\begin{equation}
\label{eq;07.9.23.2}
 \frac{\del u}{\del x}
=\Lambda \cdot u
+p_1(x,\vecy,w)\cdot u
+q(x,\vecy,w)
\end{equation}
Here, $q(x,\vecy,w)$ satisfy
$\bigl|q(x,\vecy,w)\bigr|
\leq C\, |x|^{-m}\, |\vecy|^{-N}$.

\begin{lem}
\label{lem;07.9.23.3}
We have a solution $v$ of
the equation {\rm(\ref{eq;07.9.23.2})}
satisfying $|v|\leq C'|x|^{-m}\,|\vecy|^{-N}$.
\end{lem}
\pf
We consider the following integral transform:
\[
 \nbigq_{x_2(m)}(u)
=\int_{\Gamma(x,x_2(m))}
 \exp\bigl((x-t)\,\Lambda\bigr)
 \Bigl(
 p_1(t,\vecy,w)\cdot u(t,\vecy,w)+q(t,\vecy,w)
 \Bigr)\, dt
\]
If $|u_i|\leq C\, |x|^{-m}|\vecy|^{-N}$ $(i=1,2)$,
we have the following:
\[
 \bigl| \nbigq_{x_2(m)}(u_1)-\nbigq_{x_2(m)}(u_2)
 \bigr|
=\bigl|\nbigr_{x_2(m)}(u_1-u_2)\bigr|
\leq
 \gamma_m\, K_m\, C\,
 |x|^{-m}\, |\vecy|^{-N}
\]
We put $v_0:=0$
and $v_i=\nbigq_{x_2(m)}(v_{i-1})$ for $i\geq 1$
inductively.
Then, we have 
$|v_{i+1}-v_i|\leq
 (\gamma_mK_m)^i\, C\, |x|^{-m}|\vecy|^{-N}$.
Hence, we have the limit
$v_{\infty}:=\lim_{i\to\infty}v_i$
which satisfies the following:
\[
 \nbigq_{x_2(m)}v_{\infty}=v_{\infty},
\quad
 |v_{\infty}|
 \leq
 \frac{1}{1-\gamma_mK_m}
 C\, |x|^{-m}\,|\vecy|^{-N}
\]
Thus, the proof of Lemma \ref{lem;07.9.23.3}
is finished.
\hfill\qed

\subsection{Estimate of 
 the difference of two solutions
 of the equation (\ref{eq;07.6.14.2})}
\label{subsection;07.9.23.11}

Let us consider the equation (\ref{eq;07.6.14.2})
in the case $a=0$.

\begin{lem}
\label{lem;07.9.23.6}
Let $u_i$ $(i=1,2)$ be solutions
of {\rm (\ref{eq;07.6.14.2})}
on $S\bigl(x_2(m),\vecR^{(0)}_y\bigr)$,
which satisfy
$|u_i|\leq C_N\,|x|^{-m_0}\, |\vecy|^{-N}$
for any $N$ and some $m_0>0$.
Then, we have 
$|u_1-u_2|\leq
 C_{N}'\, |x|^{-m}|\vecy|^{-N}$
for any $N$.
\end{lem}
\pf
We may and will assume $m_0<m$.
We put $V:=u_1-u_2$.
We have the following equality:
\[
 \frac{\del V}{\del x}
=\Lambda\cdot V
+p_1(x,\vecy,w)\cdot V
+q(x,\vecy,w)
\]
Here, 
$q(x,\vecy,w):=
 p_2\bigl(x,\vecy,w,u_1(x,\vecy,w)
 \bigr)
-p_2\bigl(x,\vecy,w,u_2(x,\vecy,w)\bigr)$.
By assumption,
we have 
$\bigl| q(x,\vecy,w)\bigr|
\leq B_{N}\,|x|^{-m}\,|\vecy|^{-N}$
for any $N$.
Due to Lemma \ref{lem;07.9.23.3},
we can take $v(m)$ on 
$S(x_2(m),\vecR^{(0)}_y)$
satisfying the following:
\[
 \frac{\del v(m)}{\del x}
=\Lambda\cdot v(m)
+p_1(x,\vecy,w)\cdot v(m)
+q(x,\vecy,w),
\quad
 \bigl|v(m)\bigr|
 \leq
 C_1\,|x|^{-m}|\vecy|^{-N}
\]
We have 
$\bigl| V-v(m)\bigr|
\leq C_{1,N}\, |x|^{-m_0}|\vecy|^{-N}$
and the following equality:
\[
 \frac{\del (V-v(m))}{\del x}
=\Lambda\bigl(V-v(m)\bigr)
+p_1(x,\vecy,w)\bigl(V-v(m)\bigr)
\]
Due to Lemma \ref{lem;07.9.23.4},
we obtain 
$\bigl|V-v(m)\bigr|
\leq C_{2,N}|x|^{-m}\,
 |\vecy|^{-N}$.
Hence, we obtain
$\bigl|V\bigr|
\leq C_{3,N}|x|^{-m}\,|\vecy|^{-N}$.
Thus, we are done.
\hfill\qed

\subsection{End of the proof of
 Proposition \ref{prop;07.6.15.3}}
\label{subsection;07.9.23.15}

Let us show 
$|U(m_0)|\leq C_{m,N}|x|^{-m}|\vecy|^{-N}$
for any $m$ and $N$,
where $U(m_0)$ was constructed in 
Section \ref{subsection;07.9.23.5}.
Take $m>m_0$,
and let us compare $U(m_0)$ and $U(m)$
on $S\bigl(x_1^{(3)}(m,N_0),\vecR^{(0)}_y\bigr)$.
Both of them satisfy the differential equation
(\ref{eq;07.6.14.2}),
and both of them are dominated by
$C_N\,|x|^{-m_0}\, |\vecy|^{-N}$
for any $N$.
Because of Lemma \ref{lem;07.9.23.6},
we obtain the estimate
$\bigl|U(m_0)-U(m)\bigr|
\leq A_{m,N}|x|^{-m}\,|\vecy|^{-N}$
for any $N$.
Thus, we obtain the desired estimate
$|U(m_0)|
\leq A'_{m,N}|x|^{-m}\,|\vecy|^{-N}$,
and the proof of Proposition \ref{prop;07.6.15.3}
is finished.
\hfill\qed

\section{Estimates of some integrals on a sector}
\label{subsection;08.9.28.130}
\subsection{Exponential decay}
\label{subsection;07.12.20.30}

Let $S_{\infty}(r_0,\theta_0,\theta_1)$
denote a sector
$\bigl\{x=r\, e^{\sqrt{-1}\theta}\,\big|\,
 \theta_0<\theta<\theta_1,r_0<r \bigr\}$
around $\infty$
for $\theta_0,\theta_1\in\real$ and $r_0>0$.
Fix any $\theta^{(0)}$ such that
$\theta_0<\theta^{(0)}<\theta_1$.
In the following argument, 
$r_0'$ and $\theta_i'$ $(i=0,1)$
denote the numbers
such that 
(i) $r_0'>r_0$ is appropriately large,
(ii) $\theta_0<\theta_0'<\theta^{(0)}<\theta_1'<\theta_1$,
(iii) $|\theta_i'-\theta^{(0)}|$ are sufficiently small.
Let $Y$ be a complex manifold.
Let $P$ be any point of $Y$.
In the following argument,
$U_P$ denotes a small compact neighbourhood
of $P$.

Let $(V,\nabla)$ be a flat bundle on 
$\nbigy:=S_{\infty}(r_0,\theta_0,\theta_1)\times Y$ 
relative to $Y$
(i.e., we consider derivations
only along $S_{\infty}(r_0,\theta_0,\theta_1)$-direction).
Assume
we are given a frame $\vecv$ such that
$\nabla\vecv=\vecv\, (d\gminia+A\, dx/x)$,
where $A$ is a constant matrix,
and $\gminia$ is a polynomial
$\sum_{j=1}^k\gminia_j\, x^j$
whose coefficients $\gminia_j$ are holomorphic 
functions on $Y$.
We assume $\gminia_k$ is nowhere vanishing.
Let $\omega=\sum \omega_i\, v_i$ be
a $\nabla$-closed section of
$V\otimes\Omega^{1}_{\nbigy/Y}$ on $\nbigy$
such that
$\bigl|\omega_i\bigr|=
   O\bigl(\exp(-\epsilon|x|^L)\bigr)$
for some $\epsilon>0$,
where 
we use the Euclidean metric
$dr\, dr+r^2\, d\theta \, d\theta$
on $S_{\infty}(r_0,\theta_0,\theta_1)$.

\begin{lem}
\label{lem;07.7.19.102}
There exists a section
 $\tau=\sum\tau_i\, v_i$ of $V$ on 
 some small $S_{\infty}(r_0',\theta_0',\theta_1') \times U_P$
satisfying
$\nabla\tau=\omega$
and $O(|\tau_i|)=O\bigl(\exp(-2^{-1}\epsilon|x|^L)\bigr)$.
\end{lem}
\pf
We may replace  $\vecv$
with $\vecv\,\exp(-A\log x)$,
and so we may assume $\rank(V)=1$ and $A=0$.
Formally,
we put 
\begin{equation}
 \label{eq;08.9.2.30}
 \tau(r,\theta,Q):=\exp(\gminia)
 \int_{\gamma(r,\theta)}\exp(-\gminia)\,\omega
\end{equation}
for some path $\gamma(r,\theta)$
connecting $(r,\theta,Q)$ and $(r_2,\theta_2,Q)$,
where $(r_2,\theta_2)$ is a base point of
$S_{\infty}(r_0,\theta_0,\theta_1)\cup\{\infty\}$.
Then, it satisfies $\nabla\tau=\omega$,
if the integral (\ref{eq;08.9.2.30}) converges.
The problem is how to choose the path
$\gamma(r,\theta)$
so that the desired estimate holds.

In the case $L>k$,
we choose the path $\gamma(r,\theta)$
connecting $(r,\theta)$ and $\infty$ on the ray
$\bigl\{(s,\theta)\,\big|\,s\geq r\bigr\}$.
Then, $\tau$ 
satisfies the desired estimate.
So we have only to consider the case $L\leq k$.

We have 
$\Re(\gminia)=
 r^k\,
 \Re\bigl(\gminia_k\, e^{\sqrt{-1}k\theta}\bigr)
+O\bigl(r^{k-1}\bigr)$.
We also have the following equality:
\begin{multline}
\label{eq;07.7.19.101}
 \frac{\del}{\del r}
 \Bigl( r^{-k+1}
 \exp\bigl(-\Re(\gminia)
 -\epsilon \, r^L\bigr)
 \Bigr) 
=(-k+1)\, r^{-k}
 \exp\bigl(-\Re(\gminia)
   -\epsilon\, r^L\bigr)\\
-\left(
 r^{-k+1}\frac{\del \Re(\gminia)}{\del r}
+\epsilon\, L\, r^{-k+L}
\right)\,
 \exp\bigl(-\Re(\gminia)-\epsilon \, r^L\bigr)
\end{multline}
Note we have
$r^{-k+1}\del_r\Re(\gminia)
=k\, 
 \Re\bigl(
 \gminia_k\, e^{\sqrt{-1}k\theta}\bigr)
+O(r^{-1})$.
In the case $k=L$,
we will have to be concerned with the $0$-set of
the function
$r^{-k+1}\del_r\Re(\gminia)
+\epsilon\, k 
=k\Bigl(
 \Re\bigl(
 \gminia_k\, e^{\sqrt{-1}k\theta}\bigr)
+\epsilon\Bigr)+O(r^{-1})$
to use (\ref{eq;07.7.19.101})
for the estimate of the integral of
$\exp\bigl(-\Re(\gminia)-\epsilon\, r^k\bigr)$.

\vspace{.1in}
In the following,
$S_{\infty}(r_0',\theta_0',\theta_1')$ 
will be denoted by $S$,
for simplicity of the description.
If $S\times U_P$ is sufficiently small,
we may assume that one of the following holds:
\begin{description}
\item[(A)]
 $\del_{\theta}\Re(\gminia)\neq 0$
 on $S\times U_P$.
\item[(B)]
 $\del_{r}\Re(\gminia)\neq 0$
 on $S\times U_P$.
\end{description}

Let us consider the case (A).
We may assume that $-\Re(\gminia)$ is increasing
with respect to $\theta$.
We may also assume 
that one of the following holds
for sufficiently small $U_P$,
if we choose $\theta_0'$ appropriately:
\begin{description}
\item[(A1)] 
 $L=k$ and
 $-\Re\bigl(
 \gminia_k\, e^{\sqrt{-1}k\theta_0'}
 \bigr)-\epsilon<0$ 
 on $U_P$.
\item[(A2)] $L=k$ and 
 $-\Re\bigl(
 \gminia_k\, e^{\sqrt{-1}k\theta_0'}
 \bigr)-\epsilon>0$ 
 on $U_P$.
\item[(A3)] $L<k$ and
 $-\Re\bigl(\gminia_k\, e^{\sqrt{-1}k\theta_0'}
 \bigr)<0$ on $U_P$.
\item[(A4)]  $L<k$ and 
 $-\Re\bigl(\gminia_k\, e^{\sqrt{-1}k\theta_0'}
 \bigr)>0$ on $U_P$.
\end{description}
In the case (A1),
let $\gamma_1$ be the path connecting
$\infty$ and $(r,\theta_0')$ on the ray 
$\bigl\{(s,\theta_0')\,\big|\,s\geq r\bigr\}$,
and let $\gamma_2$ be the path connecting
$(r,\theta_0')$ and $(r,\theta)$ on the arc
$\bigl\{re^{\sqrt{-1}\varphi}
 \,\big|\,\theta_0'\leq\varphi\leq\theta\bigr\}$.
The contribution of $\gamma_2$ to the integral
(\ref{eq;08.9.2.30}) is dominated by
$\exp(-\epsilon \, r^L)$.
If $S\times U_P$ is sufficiently small,
the contribution of
$\gamma_1$ is dominated by the following
for some $C>0$:
\[
  \exp\bigl(\Re(\gminia)_{\theta_0',r,Q}\bigr)
 \int_{r}^{\infty}
 \exp\bigl(
 -\Re(\gminia)_{\theta_0',s,Q}
 -\epsilon\, s^{k}\bigr)\, ds
\leq 
 C\,\exp\bigl(-\epsilon\, r^k\bigr)
\]
Here, we have used (\ref{eq;07.7.19.101}).
Thus we are done in the case (A1).

In the case (A2),
let $\gamma_1$ be the path connecting
$(r_0',\theta_0')$ and $(r_0',\theta)$ on the arc
$\{r_0'\, e^{\sqrt{-1}\varphi}
 \,|\,\theta_0'\leq \varphi\leq\theta\}$,
and let $\gamma_2$ be the path connecting
$(r_0',\theta)$ and $(r,\theta)$ on the ray
$\bigl\{s\, e^{\sqrt{-1}\theta}\,\big|\,
 r_0\leq s\leq r \bigr\}$.
We have 
$\Re(\gminia_k\, e^{\sqrt{-1}k\theta})<-\epsilon$
on $S\times U_P$.
Hence, if $S\times U_P$ is sufficiently small,
we may have $\Re(\gminia)<-(2/3)\,\epsilon\, r^k$ 
on $S\times U_P$.
Then, the contribution of $\gamma_1$ is dominated
by $\exp\bigl(-(2/3)\,\epsilon\, r^k\bigr)$.
If $r_0'$ is sufficiently large,
the contribution of $\gamma_2$ 
can be estimated as follows for some $C_1>0$,
by using (\ref{eq;07.7.19.101}):
\[
 \exp\bigl(\Re(\gminia)_{\theta,r,Q}\bigr)
\int_{r_0'}^r
 \exp\bigl(
 -\Re(\gminia)_{\theta,s,Q}
-\epsilon\, s^k
 \bigr)\, ds
\leq
 C_1\, \Bigl(
 \exp\bigl(\Re(\gminia)_{\theta,r,Q}\bigr)
+\exp\bigl(-\epsilon\, r^k\bigr)
 \Bigr)
\]
Thus we are done in the case (A2).

In the case (A3),
we take the paths $\gamma_j$ $(j=1,2)$
as in the case (A1).
The contribution of $\gamma_2$ is dominated by
$\exp\bigl(-\epsilon\,|x|^L\bigr)$.
We have the following estimate,
by using (\ref{eq;07.7.19.101}):
\[
 \exp\bigl(\Re(\gminia)_{\theta_0',r,Q}\bigr)
\int_r^{\infty}
 \exp\bigl(-\Re(\gminia)_{\theta_0',s,Q} 
 -\epsilon \, s^L\bigr)\, ds
\sim r^{-k+1}\exp(-\epsilon \, r^L)
\]
Thus, we are done in the case (A3).

In the case (A4),
we take the paths $\gamma_j$ $(j=1,2)$
as in the case (A2).
The contribution of $\gamma_1$ is dominated by
$\exp(\Re(\gminia)_{\theta,r,Q})$.
If $r_0'$ is sufficiently large,
we have
$\Re(\gminia)_{\theta,r,Q}\leq
 \Re(\gminia)_{\theta_0',r,Q}
<2^{-1}
 \Re(\gminia_k(Q)e^{\sqrt{-1}\theta_0'})\, r^k
\leq -\epsilon \, r^L$ for any $r>r_0'$.
Hence, the contribution of $\gamma_1$ is
dominated as desired.
By using (\ref{eq;07.7.19.101}),
we obtain the following estimate:
\begin{multline*}
  \exp\bigl(\Re(\gminia)_{\theta,r,Q}\bigr)
\int_{r_0'}^{r}
 \exp\bigl(-\Re(\gminia)_{\theta,s,Q} 
 -\epsilon \, s^L\bigr)\, ds
\leq \\
 C_1\, r^{-k+1}\exp(-\epsilon \, r^L)
+C_2\, 
 \exp\bigl(\Re(\gminia)_{\theta,r,Q}\bigr)
\end{multline*}
Then, we can obtain the 
desired estimate for the contribution of $\gamma_2$.

\vspace{.1in}

Let us consider the case (B).
If $S\times U_P$ is sufficiently small,
we may assume that one of the following holds:
\begin{description}
\item[(B1)]
 $L=k$ and
 $-\Re\bigl(
 \gminia_k\, e^{\sqrt{-1}k\theta}\bigr)
 -\epsilon 
 <-\epsilon/10$
on $S\times U_P$.
\item[(B2)]
 $L=k$ and
 $-\Re\bigl(\gminia_k\, e^{\sqrt{-1}k\theta}\bigr)
-2\epsilon/3>\epsilon/10$
on $S\times U_P$.
\item[(B3)]
 $L<k$ and 
$-\Re\bigl(\gminia_k\, e^{\sqrt{-1}k\theta}\bigr)
 <0$
on $S\times U_P$.
\item[(B4)] 
 $L<k$ and 
$-\Re\bigl(\gminia_k\, e^{\sqrt{-1}k\theta}\bigr)
 >0$
on $S\times U_P$.
\end{description}
In the cases (B1) and (B3),
let $\gamma$ be the path connecting
$\infty$ and $(r,\theta)$ on the ray 
$\bigl\{s\, e^{\sqrt{-1}\theta}\,\big|\,s\geq r\bigr\}$.
In the case (B1),
we obtain the following estimate,
by using (\ref{eq;07.7.19.101}):
\[
  \exp\bigl(
 \Re(\gminia)_{\theta,r,Q}
 \bigr)
\,
 \int_{r}^{\infty}
 \exp\bigl(
 -\Re(\gminia)_{\theta,s,Q}
 -\epsilon\, s^k
 \bigr) 
\leq
 C_1\,
 \exp(-\epsilon\, r^k)
\]
The case (B3) can be estimated similarly 
and easily.

In the cases (B2) and (B4),
let us take the paths $\gamma_j$ $(j=1,2)$
as in the case (A2).
In the case (B2),
the contribution of $\gamma_1$
is dominated by 
$\exp\bigl(\Re(\gminia)_{\theta,r,Q}\bigr)$.
If $r_0'$ is sufficiently large,
$\exp\bigl(\Re(\gminia)_{\theta,r,Q}\bigr)$
is dominated by 
$\exp\bigl(-2^{-1}\epsilon\, r^k\bigr)$.
The contribution of $\gamma_2$ can be
dominated as follows,
by using (\ref{eq;07.7.19.101}):
\begin{multline}
\exp\bigl(
 \Re(\gminia)_{\theta,r,Q}
 \bigr)\,
 \int_{r_0'}^{r}
 \exp\bigl(
 -\Re(\gminia)_{\theta,s,Q}
 -\epsilon\, s^k
 \bigr)\, ds \\
\leq
  \exp\bigl(
 \Re(\gminia)_{\theta,r,Q}
 \bigr)\,
 \int_{r_0'}^{r}
 \exp\bigl(
 -\Re(\gminia)_{\theta,s,Q}
 -(2/3)\,\epsilon\, s^k
 \bigr)\, ds \\
\leq
 C_1\Bigl(
 \exp\big(-(2/3)\,\epsilon\, r^k)
+\exp\bigl(
 \Re(\gminia)_{\theta,r,Q}
 \bigr)
 \Bigr)
\end{multline}
Note that the integral
$ \int_{r_0'}^{r}
 \exp\bigl(
 -\Re(\gminia)_{\theta,s,Q}
 -(2/3)\,\epsilon\, s^k
 \bigr)\, ds$
is estimated more easily
than $ \int_{r_0'}^{r}
 \exp\bigl(
 -\Re(\gminia)_{\theta,s,Q}
 -\epsilon\, s^k
 \bigr)\, ds$,
because 
$r^{-k+1}\del_r \Re(\gminia)+2\epsilon\, k/3
<-\eta<0$ for some $\eta>0$ on $S\times U_P$.
The case (B4) can be estimated
as in the case (A4).
Thus the proof of Lemma \ref{lem;07.7.19.102} is finished.
\hfill\qed

\subsection{Polynomial order}
\label{subsection;07.11.16.100}

We continue to use the setting in 
Section \ref{subsection;07.12.20.30}.
Let $\omega=\sum\omega_i\, v_i$ be
a $\nabla$-closed $C^{\infty}$-section of
$V\otimes\Omega^1$ 
on $\nbigy$ such that $|\omega_i|=O\bigl(|x|^N\bigr)$
for some $N$.
\begin{lem}
\label{lem;07.7.22.20}
There exists a $C^{\infty}$-section 
 $\tau=\sum \tau_i\, v_i$ of $V$
on $S=S_{\infty}(r_0',\theta_0',\theta_1')$
such that (i) $\nabla\tau=\omega$,
(ii) $|\tau_i|=O\bigl(|x|^{N+l}\bigr)$,
where $l\geq 0$ is independent of $\omega$.
\end{lem}
\pf
As in the proof of Lemma \ref{lem;07.7.19.102}
we may assume $\rank(V)=1$ and $A=0$
by making a meromorphic transform.
We remark that
the additional growth order $l$ may appear,
but it is independent of $\omega$.
We may assume that one of (A) or (B) holds
as in the proof of Lemma \ref{lem;07.7.19.102},
if $S\times U_P$ is sufficiently small.

Let us consider the case (A).
We may assume that $-\Re(\gminia)$ is increasing
with respect to $\theta$.
If $U_P$ is sufficiently small,
and if we choose $\theta_0'$ appropriately,
one of the following holds:
\begin{description}
\item[(A1')] 
 $-\Re(\gminia_k\, e^{\sqrt{-1}k\theta_0'})<0$
 on $U_P$.
\item[(A2')] 
 $-\Re(\gminia_k\, e^{\sqrt{-1}k\theta_0'})>0$
 on $U_P$.
\end{description}
Let us consider the case (A1').
We take the paths $\gamma_1$ and $\gamma_2$
as in the case (A1) in the proof of Lemma \ref{lem;07.7.19.102}.
Then, $|\tau(r,\theta,Q)|$ is dominated by the following:
\begin{multline*}
 \exp\bigl(\Re(\gminia)_{\theta,r,Q}\bigr)
\int_{\gamma_1}
 \exp\bigl(-\Re(\gminia)_{\theta_0',s,Q}\bigr)
    \, s^N\, ds \\
+\exp\bigl(\Re(\gminia)_{\theta,r,Q}\bigr)
\int_{\gamma_2}
 \exp\bigl(-\Re(\gminia)_{\varphi,r,Q}\bigr)
   \, r^N\, d\varphi
\end{multline*}
The contribution of $\gamma_2$ is dominated by $r^N$.
We have the following equality:
\begin{multline}
\label{eq;07.12.20.35}
 \frac{\del}{\del s}
\Bigl(
 s^{-k+1+N}
 \exp\bigl(-\Re(\gminia)_{\theta,s,Q}\bigr)
\Bigr)\\
=(N-k+1)\, s^{-k+N}\,
 \exp\bigl(-\Re(\gminia)_{\theta,s,Q}\bigr)
-s^{-k+1+N}
 \,\frac{\del \Re(\gminia)_{\theta,s,Q}}{\del s}\,
\exp\bigl(-\Re(\gminia)_{\theta,s,Q}\bigr)
\end{multline}
Note $s^{-k+1}\, \del_s\Re(\gminia)_{\theta,s,Q}
=k\,\Re(\gminia_k)_{\theta,s,Q}+O(s^{-1})$.
Hence, we obtain the following:
\[
 \int_{\gamma_1}
 \exp\bigl(-\Re(\gminia)_{\theta_0',s,Q}\bigr)
 \, s^N\, ds
\leq
 C_1\, r^{N-k+1}
 \exp\bigl(-\Re(\gminia)_{\theta_0',r,Q}\bigr)
\]
Therefore, 
the contribution of $\gamma_1$ is also dominated by
$r^N$.

For the case (A2'),
we take the paths $\gamma_i$ ($i=1,2$)
as in the case (A2)
in the proof of Lemma \ref{lem;07.7.19.102}.
The contribution of $\gamma_1$
is dominated by
$\exp\bigl(\Re(\gminia)\bigr)$.
It can be shown 
as in the proof of
Lemma \ref{lem;07.7.19.102}
by using (\ref{eq;07.12.20.35}),
that the contribution of $\gamma_2$
is dominated by $r^N$.

Let us consider the case (B).
If $S\times U_P$ is sufficiently small,
one of the following holds
for some $\delta>0$:
\begin{description}
\item[(B1')] $-\Re(\gminia)<-\delta$ on $S\times U_P$.
\item[(B2')] $-\Re(\gminia)>\delta$ on $S\times U_P$.
\end{description}
We take the paths $\gamma$
as in the cases (B1) or (B2)
in the proof of Lemma \ref{lem;07.7.19.102},
respectively.
We can obtain the desired estimate as in the cases
(A1') and (A2') above,
by using (\ref{eq;07.12.20.35}).
\hfill\qed

\vspace{.1in}
Let $\omega=\sum\omega_i\, v_i$
be a $C^{\infty}$-section of 
$V\otimes\Omega^2$ on $\nbigy$ such that
$|\omega_i|=O(|x|^N)$.

\begin{lem}
\label{lem;07.7.22.21}
There exists
a $C^{\infty}$-section $\tau=\sum\tau_i\, v_i$
of $V\otimes\Omega^1_{\nbigy/Y}$ on 
$S_{\infty}(r_0',\theta_0',\theta_1')\times U_P$
such that
(i) $\nabla\tau=\omega$,
(ii) $\tau_i=O(|x|^{N+l})$,
where $l\geq 0$  is independent of $\omega$.
\end{lem}
\pf
As in the proof of Lemma \ref{lem;07.7.22.20},
we may assume $\rank(V)=1$ and $A=0$,
and we may assume that one of (A), (B1') or (B2') holds.
We have the expression 
$\omega=f(r,\theta,Q)\, dr\, d\theta$.
In the case (A), we put
\[
 \tau(r,\theta,Q):=
-\Bigl(
 \exp\bigl(\gminia_{\theta,r,Q}\bigr)
 \int_{\theta_0'}^{\theta}
 \exp\bigl(-\gminia_{\varphi,r,Q}\bigr)
 f(r,\varphi,Q)\, d\varphi
\Bigr)\, dr.
\]
In the case (B1') or (B2'),
we take the paths $\gamma$ as in the proof of
Lemma \ref{lem;07.7.22.20},
and we put
\[
 \tau(r,\theta,Q):=
\left(
 \exp(\gminia_{r,\theta,Q})
 \int_{\gamma}
 \exp(-\gminia_{s,\theta,Q})
 f(s,\theta,Q)\, ds
\right)\, d\theta.
\]
(We regard that the contribution of
 $\gamma_1$ is $0$ 
 in the case (B2').)
Then, $\tau$ has the desired property.
\hfill\qed

\section{Some Estimates on a multi-sector}
\label{subsection;08.9.28.131}
\subsection{Estimate of growth order of some integrals 
 on multi-sectors}
\label{subsection;07.6.4.6}

Let $S_x$ and $S_w$
be small multi-sectors around $\infty$:
\[
S_x=
 \bigl\{
 (x_1,\ldots,x_{\ell})\in\cnum^{\ell}\,\big|\,
 R_{x,i}<|x_i|,\,\,
 |\arg(x_i)-\theta_{x,i}|<\delta_{x,i},\,\,
 (i=1,\ldots,\ell)
 \bigr\}
\]
\[
S_w=
 \bigl\{
 (w_1,\ldots,w_p)\in\cnum^p\,\big|\,
 R_{w,i}<|w_i|,\,\,
 |\arg(w_i)-\theta_{w,i}|<\delta_{w,i},\,\,
 (i=1,\ldots,p)
 \bigr\}
\]
Let $Y$ be a compact region in $\cnum^m$.
Let $\lambda$ be a holomorphic function on $Y$.
Let $\vecm\in\seisuu^{\ell}_{>0}$
and $\vecn\in\seisuu^p_{>0}$.
Let $\vecf$ be a $\cnum^d$-valued holomorphic function
on $S_x\times S_w\times Y$
satisfying the following differential equation:
\[
 d_x\vecf=\left(
 d_x\bigl(\lambda(y)\,
 \vecx^{\vecm}\,\vecw^{\vecn}\bigr)
+\sum_{i=1}^{\ell} A_i(\vecx,\vecw,y)\frac{dx_i}{x_i}
 \right)\vecf+\vecomega
\]
Here, $A_i$ denotes $M_d(\cnum)$-valued holomorphic functions,
and we consider the derivations
only with respect to the variables $x_1,\ldots,x_{\ell}$.
Assume the following for some $\epsilon>0$:
\begin{itemize}
\item
$\vecomega=O\bigl(|\vecx|^{-N}|\vecw|^{-N}\bigr)$
holds
for any $N>0$,
where we put $|\vecx|:=\prod_{i=1}^{\ell}|x_i|$
and $|\vecw|:=\prod_{j=1}^p|w_j|$.
\item
$\Re\bigl(\lambda(y)\,
 \vecx^{\vecm}\,\vecw^{\vecn}\bigr)
 <-\epsilon\,
 \bigl|\vecx^{\vecm}\,\vecw^{\vecn}\bigr|$ 
 on $S$.
\item
$|A_i|$ and $|\vecomega|$
 are sufficiently smaller than
$\epsilon\,
 \bigl|\vecx^{\vecm}\,\vecw^{\vecn}\bigr|$.
\item
Take some point $\vecx_0\in S_{x}$.
Then, 
$ \bigl|
 \vecf_{|\{\vecx_0\}\times S_{w}\times Y}
 \bigr|
=O\bigl(|\vecw|^{-N}\bigr)$ holds
for any $N>0$.
\end{itemize}

Let $T_x:=\bigl\{
 (x_1,\ldots,x_{\ell})\in\cnum^{\ell}\,\big|\,
 \bigl|\arg(x_i)-\theta_{x,i}\bigr|<\delta_{x,i}
 \bigr\}$.
\begin{lem}
\label{lem;07.6.14.15}
\label{lem;07.9.30.50}
Under the conditions,
we have $\vecf=O\bigl(|\vecx|^{-N}|\vecw|^{-N}\bigr)$ 
for any $N$
and for any 
$\vecx\in (\vecx_0+T_x)\times S_w\times Y$.
\end{lem}
\pf
The following is minor generalization
of the argument in Lemma 14.2 of \cite{wasow}.
Let $\xi_i:=x_i^{m_i}$.
Let $\vecxi_0=\bigl(\xi_1^{(0)},\ldots,\xi_{\ell}^{(0)}\bigr)$
correspond to $\vecx_0$.
Let $T_{\xi}:=\bigl\{
 (\xi_1,\ldots,\xi_{\ell})\,\big|\,
 |\arg(\xi_i)-m_i\, \theta_{x,i}|<m_i\,\delta_{x,i},\,
 (i=1,\ldots,\ell) \bigr\}$.
By the above correspondence, 
$\vecx_0+T_x$ is contained in
$\vecxi_0+T_{\xi}$.
Hence, we have only to obtain the estimate
on $(\vecxi_0+T_{\xi})\times S_w\times Y$.

We put
$\vecxi^{(i)}:=
 \bigl(\xi_1,\ldots,\xi_i,\xi_{i+1}^{(0)},
 \ldots,\xi_{\ell}^{(0)};\vecw,y\bigr)$
for $i=0,\ldots, \ell$.
We have
$\vecxi^{(0)}=(\vecxi_0,\vecw,y)$
and 
$\vecxi^{(\ell)}=
 \bigl(\xi_1,\ldots,\xi_{\ell};\vecw,y\bigr)$.
Let $\tau_i(s)$ be the segment
connecting $\vecxi^{(i-1)}$ and $\vecxi^{(i)}$:
\[
 \tau_i(s)=
 \bigl(
 \xi_1,\ldots,\xi_{i-1},
 \xi_{i}^{(0)}+s\,(\xi_i-\xi_i^{(0)}),\,
 \xi^{(0)}_{i+1},\ldots,\xi^{(0)}_{\ell};
 \vecw,y
 \bigr),
\quad
(0\leq s\leq 1)
\]
On the path $\tau_i(s)$,
we have the following
inequalities for some constants
$C_i$:
\begin{multline}
\label{eq;07.6.4.1}
 \frac{d |\vecf|^2}{ds}
\leq
 2\Re\left(
 \lambda\prod_{j<i}\xi_j
 \,(\xi_i-\xi_i^{(0)})\,
 \prod_{j>i}\xi_j^{(0)}
\,\vecw^{\vecn}
 \right)\,|\vecf|^2 \\
+C_0\,\Bigl(
|\xi_i|^{-1}|\xi_i-\xi_i^{(0)}|\,
 |A_i|\, |\vecf|^2
+|\xi_i-\xi_i^{(0)}|\,
 \bigl|\iota_{\del_{\xi_i}}\vecomega\bigr|
\, \bigl|\vecf\bigr|
 \Bigr) \\
 \leq
 -C_1\prod_{j<i}\bigl|\xi_j\bigr|
 \,\bigl|\xi_i-\xi_i^{(0)}\bigr|\,
 \prod_{j>i}\bigl|\xi_j^{(0)}\bigr|
 \,|\vecw^{\vecn}|
\,
 \bigl|\vecf\bigr|^2
+C_2\,|\xi_i-\xi_i^{(0)}|\,
 \bigl|\iota_{\del_{\xi_i}}\vecomega\bigr|
\end{multline}
Here, we have used
$|\vecf|\leq 2(1+|\vecf|^2)$.

Let us show the following estimate
for some constant $B_N^{(i)}$ 
independent of $\vecxi$,
by an induction on $i$:
\[
 \bigl|\vecf(\vecxi^{(i)})\bigr|^2
\leq
 B_N^{(i)}\,\prod_{j\leq i}\bigl|\xi_j\bigr|^{-N}
\,\prod_{j>i}\bigl|\xi_j^{(0)}\bigr|^{-N}
\,|\vecw|^{-N}
\]
From (\ref{eq;07.6.4.1}),
we obtain the following:
\begin{multline}
 \label{eq;07.9.22.1}
 \bigl|\vecf(\vecxi^{(i)})\bigr|^2
\leq
 \exp\Bigl(
 -C_1\prod_{j<i}\bigl|\xi_j\bigr|\,
 \bigl|\xi_i-\xi_i^{(0)}\bigr|\,
 \prod_{j>i}\bigl|\xi_j^{(0)}\bigr| 
\,|\vecw^{\vecn}|
 \Bigr)\,\bigl|\vecf(\vecxi^{(i-1)})\bigr|^2
 \\
+\int_{0}^1
 \exp\Bigl(-C_1(1-s)\prod_{j<i}\bigl|\xi_j\bigr|
\,\bigl|\xi_i-\xi_i^{(0)}\bigr|\,
\prod_{j>i}\bigl|\xi_j^{(0)}\bigr|
 \,|\vecw^{\vecn}|
 \Bigr)\,C_2
 \,|\xi_i-\xi_i^{(0)}|\,
 \bigl|\iota_{\del_{\xi_i}}\vecomega\bigr|\, ds
\end{multline}
The first term can be dominated by the following:
\[
 B_N^{(i-1)}\exp\Bigl(
 -C_1\prod_{j<i}\bigl|\xi_j\bigr|
 \, \bigl|\xi_i-\xi_i^{(0)}\bigr|
 \,\prod_{j>i}\bigl|\xi_j^{(0)}\bigr|
\,|\vecw^{\vecn}|
 \Bigr)
\,\prod_{j<i}\bigl|\xi_j\bigr|^{-N}
\,\prod_{j\geq i}\bigl|\xi_j^{(0)}\bigr|^{-N}
\,|\vecw|^{-N}
\]
In the case $\bigl|\xi_i^{(0)}\bigr|\leq \bigl|\xi_i\bigr|/2$,
we have
 $\bigl|\xi_i-\xi_i^{(0)}\bigr|\geq\bigl|\xi_i\bigr|/2$,
and hence the exponential term is very small.
In the case $\bigl|\xi_i^{(0)}\bigr|>\bigl|\xi_i\bigr|/2$,
we have $\bigl|\xi_i^{(0)}\bigr|^{-N}\leq
 C_4\, \bigl|\xi_i\bigr|^{-N}$,
and hence the first term can be estimated appropriately.

Let us look at the second term of (\ref{eq;07.9.22.1}).
For $s>1/2$, we have
$\bigl|
 \xi_i^{(0)}+s(\xi_i-\xi_i^{(0)})
 \bigr|
\geq \bigl|\xi_i\bigr|/2$.
Hence, the integrand can be dominated by the following:
\begin{multline}
 C_{4,N}\,
 |\xi_i-\xi_i^{(0)}|\,
 \prod_{j<i}\bigl|\xi_i\bigr|^{-N}
\,\bigl|\xi_i^{(0)}+s(\xi_i-\xi_i^{(0)})\bigr|^{-N}
\,\prod_{j>i}\bigl|\xi_j^{(0)}\bigr|^{-N}
\,|\vecw|^{-N}\\
\leq
 C_{5,N}\,\prod_{j\leq i}\bigl|\xi_j\bigr|^{-N+1}
 \,\prod_{j>i}\bigl|\xi_j^{(0)}\bigr|^{-N}
 \,|\vecw|^{-N}
\end{multline}
Thus, $\int_{1/2}^1$ can be estimated appropriately.
For $s<1/2$, the exponential term in the integrand
is dominated by
\[
 \exp\Bigl(
 \frac{-C_1}{2}\,
 \prod_{j<i}\bigl|\xi_j\bigr|
\,\bigl|\xi_i-\xi_i^{(0)}\bigr|
\,\prod_{j>i}\bigl|\xi_j^{(0)}\bigr|
\,|\vecw^{\vecn}|
 \Bigr)
\]
The term $\bigl|\iota_{\del_{\xi_i}}\vecomega\bigr|$
is dominated by
$C_{6,N}\prod_{j<i}\bigl|\xi_j\bigr|^{-N}\,
 \prod_{j\geq i}\bigl|\xi_j^{(0)}\bigr|^{-N}
 \,|\vecw|^{-N}$,
and $|\xi_i-\xi_i^{(0)}|$ is dominated by $|\xi_i|$.
Hence, using the same argument for the estimate
of the first term,
we can estimate $\int_{0}^{1/2}$.
Thus, the induction can proceed.
The estimate for 
$\bigl|\vecf(\vecxi^{(\ell)})\bigr|
=\bigl|\vecf(\vecxi)\bigr|$
gives the claim of the lemma.
\hfill\qed

\subsection{Estimate of flat sections}
\label{subsection;08.9.28.132}

Let $S_x$, $S_w$, $T_x$ and $Y$ be as in Section
\ref{subsection;07.6.4.6}.
Let $1\leq k\leq \ell$,
$\vecm\in\seisuu^{k}_{>0}$
and $\vecn\in\seisuu^p_{>0}$.
Let $R$ be a $C^{\infty}$-section of
$M_d(\cnum)\otimes
 \Omega^1_{S_x}$ on $S_x\times S_w\times Y$:
\[
 R=
 \sum_{i=1}^{\ell}R_i^{(1,0)}\, 
 \frac{dx_i}{x_i}
+\sum_{i=1}^{\ell}R_i^{(0,1)}
 \, \frac{d\xbar_i}{\xbar_i}
\]
We assume the estimates
$\bigl|R_i^{(p,q)}\bigr|\leq 
 B_0\bigl|\vecx^{\vecm}\,\vecw^{\vecn}\bigr|$
for $(p,q)=(1,0),(0,1)$
and $i=1,\ldots,\ell$.
Let $\vecf$ be a $\cnum^d$-valued
$C^{\infty}$-function on $S_x\times S_w\times Y$
satisfying the following differential equation:
\[
 d_x\vecf=R\,\vecf
\]
We assume $\vecf(\vecx_0,\vecw,y)\neq 0$
for some $\vecx_0\in S_x$.

\begin{lem}
\label{lem;07.7.10.1}
\label{lem;07.11.26.10}
We have the following 
estimate on $(\vecx_0+T_x)\times S_w\times Y$
for some $C_0>0$:
\[
 \left|
\log\Bigl(
 \frac{\bigl|\vecf(\vecx,\vecw,y)\bigr|^2}
 {\bigl|\vecf(\vecx_0,\vecw,y)\bigr|^2}
 \Bigr)
\right|
\leq
 B_0\, C_0\,
 \bigl|\vecx^{\vecm}\,\vecw^{\vecn}\bigr|
\,\log\Bigl(
 \prod_{i=k+1}^{\ell}|x_i|
 \Bigr)
\]
The constant $C_0$ is independent of $\vecf$.
\end{lem}
\pf
Let $\xi_i=x_i^{m_i}$ for $i=1,\ldots,k$,
and $\xi_i=x_i$ for $i=k+1,\ldots,\ell$.
We use the notation 
in the proof of Lemma \ref{lem;07.6.14.15}.
We have only to obtain the estimate on
$(\vecxi_0+T_{\xi})\times S_w\times Y$.
On the paths
$\tau_i(s)$ $(i=1,\ldots,k)$, we have the following:
\[
\left|
 \frac{d}{ds}\bigl|\vecf\bigr|^2
\right|
\leq
 C_1\, B_0\prod_{j=1}^{i-1}|\xi_j|
 \,\bigl|\xi_i-\xi_i^{(0)}\bigr|
 \,\prod_{j=i+1}^{k}|\xi_j^{(0)}|
 \,|\vecw^{\vecn}|
 \,\bigl|\vecf\bigr|^2
\]
Hence, we obtain the following:
\[
 \left|
 \log\Bigl(
\frac{\bigl|\vecf(\vecxi^{(i)})\bigr|^2}
  {\bigl|\vecf(\vecxi^{(i-1)})\bigr|^2}
 \Bigr)\right| \\
\leq
 C_1\, B_0\prod_{j=1}^{i-1}|\xi_j|\,|\xi_i-\xi_i^{(0)}|
\,\prod_{j=i+1}^{k}|\xi_j^{(0)}|
\,|\vecw^{\vecn}|
\]
Note we have
$\prod_{j=1}^{i-1}|\xi_j|\,|\xi_i-\xi_i^{(0)}|
\,\prod_{j=i+1}^{k}|\xi_j^{(0)}|
\leq C_2\,
 \prod_{j=1}^k|\xi_j|$.
Hence,
we obtain the following:
\begin{equation}
 \label{eq;07.12.26.20}
 \left|
\log\Bigl(
 \frac{\bigl|\vecf(\vecxi^{(k)})\bigr|^2}
  {\bigl|\vecf(\vecxi^{(0)})\bigr|^2}
\Bigr)
\right|
\leq
 C_3\, B_0\prod_{j=1}^{k}|\xi_j|
 \,|\vecw^{\vecn}|
\end{equation}

On the paths $\tau_i(s)$ $(i=k+1,\ldots,\ell)$,
we have the following:
\[
 \left|
 \frac{d}{ds}|\vecf|^2
\right|
\leq
 B_0\,
 \frac{\bigl|\xi_i-\xi_i^{(0)}\bigr|}
 {\bigl|\xi_i^{(0)}+s\,(\xi_i-\xi_i^{(0)})\bigr|}
 \prod_{j=1}^k|\xi_j|\,
 |\vecw^{\vecn}|\,
 |\vecf|^2
\]
By an elementary geometric arguments,
we can show the following inequality on $\tau_i(s)$
for some $0<C_4<1$:
\begin{equation}
 \label{eq;07.12.26.11}
\frac{\bigl|\xi_i-\xi_i^{(0)}\bigr|}
 {\bigl|\xi_i^{(0)}+s\,(\xi_i-\xi_i^{(0)})\bigr|}
\leq
 \frac{C_{4}\, \bigl|\xi_i-\xi_i^{(0)}\bigr|}
 {|\xi_i^{(0)}|+s\,\bigl|\xi_i-\xi_i^{(0)}\bigr|}
\end{equation}
Hence, we obtain 
\[
 \left|
 \log\Bigl(
\frac{\bigl|\vecf(\vecxi^{(i)})\bigr|}
 {\bigl|\vecf(\vecxi^{(i-1)})\bigr|}
\Bigr) \right|
\leq
 C_5\, B_0\,\log\Bigl(
 1+\bigl|\xi_i^{(0)}\bigr|^{-1}
 \, \bigl|\xi_i-\xi_i^{(0)}\bigr|
 \Bigr)\,
 \prod_{j=1}^{k}\bigl|\xi_j\bigr|
\,\bigl|\vecw^{\vecn}\bigr|
\]
Note 
$1+|\xi_i^{(0)}|^{-1}\,|\xi_i-\xi_i^{(0)}|
\leq C_6|\xi_i|$.
Thus we are done.
\hfill\qed

\begin{lem}
\label{lem;07.12.27.5}
If moreover 
$|R_i^{(p,q)}|\leq 
 B_0\,\bigl|x_i^{-1}\vecx^{\vecm}\vecw^{\vecn}\bigr|
+B_0$ are satisfied for 
$(p,q)=(1,0),(0,1)$ and $i=k+1,\ldots,\ell$,
the following holds for some $C_0'$:
\[
 \left|
\log\Bigl(
 \frac{|\vecf(\vecx,\vecw,y)|^2}{|\vecf(\vecx_0,\vecw,y)|^2}
 \Bigr)
\right|
\leq
 B_0\, C_0'\,\Bigl(
 |\vecx^{\vecm}\,\vecw^{\vecn}|
+\log\Bigl(
 \prod_{i=k+1}^{\ell}|x_i|
 \Bigr)
 \Bigr)
\]
\end{lem}
\pf
We obtain (\ref{eq;07.12.26.20})
without any change.
On $\tau_i(s)$ $(i=k+1,\ldots,\ell)$,
we have the following:
\[
 \frac{d}{ds}\log|\vecf|^2
\leq
 B_0\,
 \frac{\bigl|\xi_i-\xi_i^{(0)}\bigr|}
 {\bigl|\xi_i^{(0)}+s\, (\xi_i-\xi_i^{(0)})\bigr|}
+
 B_0\,
 \frac{\bigl|\xi_i-\xi_i^{(0)}\bigr|}
 {\bigl|\xi_i^{(0)}+s\, (\xi_i-\xi_i^{(0)})\bigr|^2}
 \prod_{j=1}^k|\xi_j|\,
 |\vecw^{\vecn}|
\]
The contribution of the first term can be 
dominated by
$\log\bigl(1+|\xi_i^{(0)}|^{-1}|\xi_i-\xi_i^{(0)}|\bigr)$.
By using (\ref{eq;07.12.26.11}),
it is easy to show that
the contribution of the second term
can be dominated by
$\prod_{i=1}^k|\xi_i|\,
 |\vecw^{\vecn}|$.
\hfill\qed

\begin{rem}
\label{rem;10.7.2.3}
We have a variant of Lemma 
{\rm\ref{lem;07.12.27.1}} corresponding to
Lemma {\rm\ref{lem;07.12.27.5}}.
Namely,
if $\ord(R)\geq \vecm(1)$,
then we have an estimate 
\[
  \left|
 \log\Bigl(
 \frac{\bigl|g\,
 \exp\bigl(\lambda^{-1}\gminia\bigr)
 \bigr|(\lambda,\vecz,\vecy)}
 {\bigl|g\,
 \exp\bigl(\lambda^{-1}\gminia\bigr)
 \bigr|(\lambda,\vecz_0,\vecy)}
 \Bigr) \right|
\leq
 C|\vecz^{\vecm(1)}|
+C\log\Bigl(
 |z_{i(0)}|^{-1}
\prod_{i=k+1}^{\ell}|z_i|^{-1}
\Bigr)
\]
If $m_{i(0)}<-1$,
the term $|z_{i(0)}|^{-1}$ in the logarithm
is not necessary.
See Subsection {\rm\ref{subsection;10.7.2.2}}
for the order $\ord R$.
\hfill\qed
\end{rem}

\subsection{Estimate for a flat section of
 a family of $\lambda$-flat bundles for $\lambda\neq 0$}
\label{subsection;08.9.28.133}

We give a special version of
Lemma \ref{lem;07.7.10.1},
which will often be used in our subsequent arguments.
Let $X$ be a product of
$\Delta_z^{\ell}$ and a complex manifold $Y$.
Let $\nbigk$ denote a point or a compact region
of $\cnum_{\lambda}^{\ast}$.
Let $\nbigx:=X\times\nbigk$
and $\nbigd:=D\times \nbigk$.
Let $\pi:\nbigxtilde(\nbigd)\lrarr\nbigx$ denote
the real blow up of $\nbigx$ along $\nbigd$.
Let $S_{\vecz}$ be a multi-sector of 
$(\Delta^{\ast})^{\ell}$.
Let $S$ be a product of 
$S_{\vecz}$
and a compact region $U$ of $Y\times\nbigk$.
Let $\Sbar$ denote the closure of $S$
in $\nbigxtilde(\nbigd)$.

Let $\vecm\in\seisuu_{<0}^k$ for some
$k\in[1,\ell]$,
and let $\vecm(1)=\vecm+\vecdelta_{i(0)}$
for some $i(0)\in [1,k]$.
Let $\gminia$ be the product of
$\vecz^{\vecm}$
and a nowhere vanishing
holomorphic function $\gminia_{\vecm}$
on $\nbigx$.
Let $R$ be a holomorphic section of
$M_d(\cnum)\otimes 
 p_{\lambda}^{\ast}\bigl(
 \vecz^{\vecm(1)}\Omega^{1,0}_{X}(\log D)
 \bigr)$
on $\Sbar$.
\begin{lem}
\label{lem;07.12.27.1}
Let $g$ be a $\cnum^d$-valued holomorphic
function on $S$ such that
$\lambda\, d_{\vecz}g=(d\gminia+R)\, g$.
Let $\vecz_0$ be any point of $S_{\vecz}$.
Then,
we have the following estimate for some $C>0$:
\begin{equation}
\label{eq;07.12.27.2}
 \left|
 \log\Bigl(
 \frac{\bigl|g\,
 \exp\bigl(\lambda^{-1}\gminia\bigr)
 \bigr|(\lambda,\vecz,\vecy)}
 {\bigl|g\,
 \exp\bigl(\lambda^{-1}\gminia\bigr)
 \bigr|(\lambda,\vecz_0,\vecy)}
 \Bigr) \right|
\leq
 C|\vecz^{\vecm(1)}|\,
\log\Bigl(
 |z_{i(0)}|^{-1}
\prod_{i=k+1}^{\ell}|z_i|^{-1}
\Bigr)
\end{equation}
If $m_{i(0)}<-1$,
the term $|z_{i(0)}|^{-1}$ in the logarithm
is not necessary.
(If $g_{|\vecz_0\times U}=0$,
 we use the convention $\log(0/0)=0$.)
\end{lem}
\pf
We put 
$ f=g\, 
 \exp\bigl(\lambda^{-1}
 (\gminia-\gminia_{|\vecz_0\times U})
 \bigr)$.
By changing the variables
$x_i=z_i^{-1}$,
we apply Lemma \ref{lem;07.7.10.1}
to $f$.
Then, we obtain the desired estimate.
(Note we may assume that
 $R$ and $g$ are defined on 
 $S'$ which is slightly larger than $S$.)
\hfill\qed

\begin{cor}
\label{cor;07.12.27.4}
Let $g$ be as in Lemma {\rm\ref{lem;07.12.27.1}}.
Assume the following:
\begin{itemize}
\item
$k=\ell$ and 
$\Re(\lambda^{-1}\gminia)
 <-\delta|\vecz^{\vecm}|$
 for some $\delta>0$.
\end{itemize}
If we shrink $S$ appropriately in the radius direction,
we have
$|g|=O\bigl(
 \exp\bigl(-\epsilon|\vecz^{\vecm}|\bigr)
 \bigr)$
for some $\epsilon>0$.
\hfill\qed
\end{cor}

\begin{cor}
\label{cor;07.12.27.3}
Let $g$ be as in Lemma {\rm\ref{lem;07.12.27.1}}.
Assume $g_{|\vecz_0\times U}$ is nowhere vanishing.
We also assume the following:
\begin{itemize}
\item
$k=\ell$ and 
$\Re(\lambda^{-1}\gminia)
 >\delta|\vecz^{\vecm}|$
 for some $\delta>0$.
\end{itemize}
If we shrink $S$ appropriately in the radius direction,
we have
$|g|\geq C_1
 \exp\bigl(\epsilon|\vecz^{\vecm}|\bigr)$
for some $C_1>0$ and $\epsilon>0$.
\hfill\qed
\end{cor}

\subsection{Estimate for a flat section of
 a family of $\lambda$-flat bundles around $\lambda=0$}
\label{subsection;08.9.28.134}

We give another special version of
Lemma \ref{lem;07.7.10.1}.
Let $X$ be a product of
$\Delta_z^{\ell}$ and a complex manifold $Y$.
Let $\nbigk$ denote a neighbourhood of 
$0_{\lambda}$ in $\cnum_{\lambda}$.
Let $\nbigx:=X\times\nbigk$
and $\nbigd:=D\times \nbigk$.
Let $\nbigx^0:=X\times\{0_{\lambda}\}$.
We put $W:=\nbigx^0\cup\nbigd$.
Let $\nbigxtilde(W)$ denote the real blow up
of $\nbigx$ along $W$.
Let $S_{\lambda}$ and $S_{\vecz}$
be multi-sectors of $\nbigk-\{0_{\lambda}\}$ and
$(\Delta^{\ast})^{\ell}$,
respectively.
Let $S$ be a product of 
$S_{\lambda}$, $S_{\vecz}$
and a compact region $U$ of $Y$,
and let $\Sbar$ denote the closure of $S$
in $\nbigxtilde(W)$.

Let $\vecm\in\seisuu_{<0}^k$ for some
$k\in[1,\ell]$,
and let $\vecm(1)=\vecm+\vecdelta_{i(0)}$
for some $i(0)\in [1,k]$.
Let $\gminia$ be the product of
$\vecz^{\vecm}$
and a nowhere vanishing
holomorphic function $\gminia_{\vecm}$
on $\nbigx$.
Let $R$ be a holomorphic section of
$M_d(\cnum)\otimes 
 p_{\lambda}^{\ast}\bigl(
 \vecz^{\vecm(1)}\Omega_{X}(\log D)
 \bigr)$
on $\Sbar$.

\begin{lem}
\label{lem;07.12.17.12}
Let $g$ be a $\cnum^d$-valued holomorphic
function on $S$ such that
\[
 \lambda\, d_{\vecz}g=(d\gminia+R)\, g.
\]
Let $\vecz_0$ be any point of $S_{\vecz}$.
Then, we have the following estimate for some $C>0$:
\begin{equation}
\label{eq;07.12.17.1}
 \left|
 \log\Bigl(
 \frac{\bigl|g\,
 \exp\bigl(\lambda^{-1}\gminia\bigr)
 \bigr|(\lambda,\vecz,\vecy)}
 {\bigl|g\,
 \exp\bigl(\lambda^{-1}\gminia\bigr)
 \bigr|(\lambda,\vecz_0,\vecy)}
 \Bigr) 
\right|
\leq
 C|\vecz^{\vecm(1)}\,\lambda^{-1}|\,
\log\Bigl(
 |z_{i(0)}|^{-1}
\prod_{i=k+1}^{\ell}|z_i|^{-1}
\Bigr)
\end{equation}
In the case $m_{i(0)}<-1$,
the term $|z_{i(0)}|^{-1}$ in the logarithm
is not necessary.
(We use the convention
 $\log(0/0)=0$
 in the case
 $g_{|S_{\lambda}\times\vecz_0\times U}=0$.)
\end{lem}
\pf
We put 
$ f=g\, 
 \exp\bigl(\lambda^{-1}
 (\gminia-\gminia_{|S_{\lambda}\times\vecz_0\times U})
 \bigr)$.
By changing the variables
$w=\lambda^{-1}$
and $x_i=z_i^{-1}$,
we apply Lemma \ref{lem;07.7.10.1}
to $f$.
Then, we obtain the desired estimate.
\hfill\qed

\begin{cor}
\label{cor;07.12.17.10}
Let $g$ be as in Lemma 
{\rm\ref{lem;07.12.17.12}}.
Assume the following:
\begin{itemize}
\item
$|g|$ is bounded
on $S_{\lambda}\times\vecz_0\times U$.
\item
$k=\ell$ and 
$\Re(\lambda^{-1}\gminia)
 <-\delta|\lambda^{-1}\vecz^{\vecm}|$
 for some $\delta>0$.
\end{itemize}
If we shrink $S$ appropriately in the radius direction,
we have
$|g|=O\bigl(
 \exp\bigl(-\epsilon|\lambda^{-1}\vecz^{\vecm}|\bigr)
 \bigr)$
for some $\epsilon>0$.
\hfill\qed
\end{cor}

\begin{rem}
We have a variant of Lemma
{\rm\ref{lem;07.12.27.1}}
corresponding to 
Lemma {\rm\ref{lem;07.12.27.5}}.
See Remark {\rm\ref{rem;10.7.2.3}}.
\hfill\qed
\end{rem}

\chapter{Acceptable Bundle}
\label{section;08.10.23.2}

We studied acceptable bundles
in \cite{mochi}, \cite{mochi4}
and \cite{mochi2}
after \cite{cg}, \cite{s1} and \cite{s2}.
In this chapter,
we give a review and some complements.
The most important one is 
local freeness (Theorem \ref{thm;07.3.16.1})
in Section \ref{subsection;07.9.22.10}.
Although there is a considerable overlap
with \cite{mochi} and \cite{mochi2},
the author thinks it appropriate
to give the details
in view of significance of the theory of
acceptable bundles for us.

In Section \ref{subsection;08.10.18.30},
we recall basic facts on holomorphic vector bundles
on Kahler manifolds.
We review in Section \ref{subsection;d11.14.35}
the fundamental property of
the twisted metrics for acceptable bundles,
with minor refinement.
Roughly speaking,
if $N$ is sufficiently negative,
we have the vanishing of higher $L^2$-cohomology,
due to which we can extend holomorphic sections
on a hyperplane to those on the whole space.
And,  if $N$ is sufficiently positive,
the norms of holomorphic sections
are pluri-subharmonic,
due to which we do not have to care
the distinctions between 
$L^2$-estimates and growth estimates,
or curve-wise estimates and global estimates.

\vspace{.1in}
In Section \ref{subsection;07.9.22.10},
we give the statement of the main theorem
(Theorem \ref{thm;07.3.16.1}) in this chapter.
It briefly means the local freeness
of the sheaves of holomorphic sections
satisfying some growth estimate.
It was shown in \cite{mochi2}
for acceptable bundles underlying
tame harmonic bundles.
The main body of the proof is almost the same.
We need only some easy changes for the general case.
However, it is used to obtain a filtered $\lambda$-flat bundle
from an unramifiedly good wild harmonic bundle,
which is one of the most fundamental results
for us. 
Hence, we give a rather detailed outline of the proof
in Sections 
\ref{subsection;08.10.18.20}--\ref{subsection;08.10.18.21}.

The theory of acceptable bundles on curves
was well established by Simpson
in \cite{s1} and \cite{s2}.
We explain some complements
in Section \ref{subsection;08.10.18.20}
for our use,
which are more or less well known.
In Section \ref{subsection;08.10.18.40},
we explain that the parabolic filtration
has the splitting by the weight decomposition
with respect to the Galois group
after taking an appropriate ramified covering.
In Section \ref{subsection;08.10.18.41},
we review that the parabolic weights
are essentially the logarithm of 
the limit of the monodromy.
Then, we show in Section \ref{subsection;08.10.18.42}
that the parabolic weights
are well defined for the irreducible components of 
the divisors in the higher dimensional case.
(We have used the property of tame harmonic bundles
 for the control of the parabolic weights 
 in \cite{mochi2}.)

To show an $\nbigo$-module is locally free,
it is always essential to show that holomorphic sections
on a hypersurface can be extended to
those on the whole space,
in an appropriate sense.
We study such problems in Sections 
\ref{subsection;08.10.18.45}--\ref{subsection;10.11.80}
by the method in \cite{mochi} and \cite{mochi2}.
We change the way of construction of a cocycle
(Section \ref{subsection;08.10.18.50}).
Then, we show Theorem \ref{thm;07.3.16.1}
in Section \ref{subsection;08.10.18.21}
by using the argument 
in \cite{mochi} and \cite{mochi2}.

\vspace{.1in}

In Section \ref{subsection;08.10.18.51},
we show that the small deformation
of an acceptable bundle is also naturally extended 
to a filtered bundle,
which is useful to show Theorem \ref{thm;07.11.20.10}.

In Section \ref{subsection;08.10.18.52}
we give some complements.
In Section \ref{subsection;08.9.29.5},
we study estimate of sup norms
for $L^2$-solutions of the $\delbar$-equation.
Lemma \ref{lem;07.9.24.5}
will be used in Sections
\ref{subsection;08.9.29.3},
\ref{subsection;07.7.7.23}
and \ref{subsection;08.9.29.4}.
We also give a variant of such an estimate
in Section \ref{subsection;08.9.29.3},
which is used in Sections
\ref{subsection;08.9.28.6}
and \ref{subsection;08.9.28.7}.
We give 
in Section \ref{subsection;08.9.28.8}
an estimate of a unitary connection form
with respect to a holomorphic frame 
compatible with the parabolic structure.
We show that it is bounded up to log order.
We give in Section \ref{subsection;07.10.5.51}
a refinement if the acceptable
bundle comes from a good wild harmonic bundle.
(See also Section \ref{subsection;08.9.29.9}.)
Such an estimate will be used 
to obtain an estimate of harmonic forms
on a punctured disc
(Section \ref{subsection;08.9.29.10}).

\section{Some general results
 on vector bundles on Kahler manifolds}
\label{subsection;08.10.18.30}
\subsection{Vanishing of $L^2$-cohomology
and some estimates}

We recall some results of Andreotti-Vesentini in \cite{av}.
Let $(Y,g)$ be a complete Kahler manifold, not necessarily
compact.
The volume form of $Y$ is denoted by $\vol$.
Let $(E,\delbar_E,h)$ be a hermitian holomorphic bundle
over $Y$.
Let $(\cdot,\cdot)_{h,g}$ denote
the induced fiber-wise hermitian metric
of $E\otimes\Omega^{p,q}_Y$.
The space of $C^{\infty}$ $(p,q)$-forms 
with compact support
is denoted by $A_c^{p,q}(E)$.
\index{space $A_c^{p,q}(E)$}
For any $\eta_i\in A_c^{p,q}(E)$ $(i=1,2)$,
we define
\[
 \langle
 \eta_1,\eta_2
 \rangle_h:=
 \int (\eta_1,\eta_2)_{h,g}\cdot\vol,
\quad
 ||\eta||^2_h:=\langle \eta,\eta\rangle_h.
\]
The completion of $A_c^{p,q}(E)$
with respect to the norm $||\cdot||_h$
is denoted by $A_h^{p,q}(E)$.
\index{space $A_h^{p,q}(E)$}
\index{$L^2$-norm $\|\cdot\|_h$}
\index{$L^2$-inner product 
 $\langle\cdot,\cdot\rangle_h$}

Let 
$\delbar_E^{\ast}:A_c^{p,q}(E)\lrarr A_c^{p,q-1}(E)$
denote the formal adjoint of
$\delbar_E:A_c^{p,q}(E)\lrarr A_c^{p,q+1}(E)$.
We set
$\laplacian=
\delbar_E^{\ast}\delbar_E+\delbar_E\delbar_E^{\ast}$.
We have the maximal closed extensions
$\delbar_E:A_h^{p,q}(E)\lrarr A_h^{p,q+1}(E)$
and 
$\delbar_E^{\ast}:A_h^{p,q}(E)\lrarr A_h^{p,q-1}(E)$.
We denote the domains of $\delbar_E$
and $\delbar_E^{\ast}$ by $Dom(\delbar_E)$
and $Dom(\delbar_E^{\ast})$ respectively.

\begin{prop}[Proposition 5 of 
 \cite{av}]
 \label{prop;11.28.6}
$A_c^{p,q}(E)$ is dense
in $W^{p,q}:=
 Dom(\delbar_E)\cap Dom(\delbar_E^{\ast})$
with respect to the the graph norm
$||\eta||_h^2+||\delbar_E\eta||^2_h
+||\delbar_E^{\ast}\eta||^2_h$.
(See also {\rm\cite{cg}}).
\hfill\qed
\end{prop}

\begin{prop}[Theorem 21 
in \cite{av}] 
\label{prop;12.11.1}
Assume that 
there exists a positive number $c>0$ satisfying the following:
\begin{quote}
For any $\eta\in W^{p,q}$,
we have $||\delbar_E\eta||_h^2+||\delbar_E^{\ast}\eta||_h^2
 \geq c\cdot ||\eta||_{h}^2$.
\end{quote}
Then, for any $C^{\infty}$-element $\eta\in A_h^{p,q}(E)$
such that $\delbar_E(\eta)=0$,
we have a $C^{\infty}$-solution $\rho\in A_h^{p,q-1}(E)$
satisfying the equation $\delbar_E(\rho)=\eta$.
\hfill\qed
\end{prop}

\subsection{Kodaira inequality}

We recall some inequality due to
K. Kodaira \cite{ko1} and 
M. Cornalba--P. Griffiths \cite{cg}.
(See also \cite{mochi} and \cite{mochi2}.)
For a Kahler manifold $Y$,
we have the operator 
$\Lambda:\Omega_Y^{p,q}\lrarr \Omega_Y^{p-1,q-1}$
which is the adjoint of the multiplication of
the Kahler form
(see  Page 62 of \cite{koba}).
Let $E$ be a holomorphic vector bundle
with a hermitian metric $h$
over $Y$.
We have the metric connection of $E$
induced by the holomorphic structure $\delbar_E$
and the hermitian metric $h$.
We denote the curvature by $R(h)$.
Let $R(\Omega^{0,q}_Y)$ denote the curvature
of the Levi-Civita connection on $\Omega^{0,q}_Y$.
We have the following inequality:
\begin{equation}
\label{eq;08.10.16.2}
 \bigl\|
 \delbar_E\eta
 \bigr\|_h^2
+\bigl\|
 \delbar^{\ast}_E\eta
 \bigr\|_h^2
  \geq \\
\sqrt{-1}\big\langle
 \Lambda(R(\Omega_Y^{0,q}))\cdot\eta,\,\eta
 \big\rangle_h
-\sqrt{-1}
 \big\langle
  \Lambda\bigl(R(h) \eta\bigr)
 -\Lambda R(h)\cdot \eta,\,
 \eta
 \big\rangle_h
\end{equation}
(See Section 2.8.2 of \cite{mochi2}.)

\vspace{.1in}
We give a more specific formula for $(0,1)$-forms.
We denote the Ricci curvature of the Kahler metric $g$
by $\Ric(g)$.
We can naturally
regard $\Ric(g)$ as a section of 
$End(E)\otimes \Omega_Y^{1,1}$.
Let $f$ be a local section of 
$End(E)\otimes\Omega^{1,1}_Y$,
and $\eta$ be a local element of $A_c^{0,1}(E)$.
We put
\begin{equation} 
\label{eq;10.11.60}
 \doublelangle
 f,\eta
 \doublerangle_h:=
 -\sqrt{-1}
\Bigl(
 \Lambda
 \bigl(
 f\cdot\eta
 \bigr)
-\Lambda(f)\cdot\eta,\,\,
 \eta
\Bigr)_h.
\end{equation}
Let $\varphi_1,\ldots,\varphi_d$ be
a local orthonormal frame of $\Omega^{1,0}$,
and let $e_1,\ldots,e_r$ be a local orthonormal
frame of $E$.
For local expressions
$ f=\sum f_{\mu,\nu,i,\bar{j}}
 \cdot e_{\mu}^{\lor}\otimes e_{\nu}
 \otimes(\varphi_i\wedge \overline{\varphi}_j),
\quad
 \eta=\sum \eta_{\mu,i}\cdot
 e_{\mu}\otimes\overline{\varphi}_i$,
we can rewrite (\ref{eq;10.11.60}) as
\begin{equation}
 \label{eq;08.10.19.30}
 \sum f_{\mu,\nu,i,\bar{j}}\cdot
 \eta_{\mu,i}\cdot \etabar_{\nu,j},
\end{equation}
where $e^{\lor}_1,\ldots,e^{\lor}_r$ denote
the dual frame.

For any $\eta\in
 Dom(\delbar_E)\cap Dom(\delbar_E^{\ast})
\subset A^{0,1}_h(E)$,
we have the following inequality:
\[
 ||\delbar_E(\eta)||^2_h
+||\delbar_E^{\ast}(\eta)||_h^2
\geq
\int\!\!
 \doublelangle
 R(h)+\Ric(g),\eta
 \doublerangle_h
\vol.
\]
(See \cite{ko1} and \cite{cg}.
See also Proposition 4.5 of \cite{mochi}.)

\section{Twist of the metric 
 of an acceptable bundle}
\label{subsection;d11.14.35}
Let $X$ be a complex manifold,
and $D$ be a normal crossing divisor of $X$.
Let $g_{\poin}$ be a Poincar\'e like metric
of $X-D$,
i.e., for any point $P\in D$,
we can take a neighbourhood $U$ of $P$ in $X$
such that $U\setminus D$
is bi-holomorphic to
the product of some discs and punctured discs,
and then $g_{\poin}$ and
the Poincar\'e metric of $U\setminus D$
are mutually bounded around $P$.
Recall that
$(E,\delbar_E,h)$ is called acceptable,
if the curvature $R(h)$ is bounded
with respect to $h$ and $g_{\poin}$.
(See also Definition 2.44 of \cite{mochi2}.)
\index{acceptable bundle}

\vspace{.1in}

Let us consider the case $X:=\Delta^n$ 
and $D:=\bigcup_{i=1}^{\ell} D_i$,
where we put $D_i:=\{z_i=0\}$.
Let $(E,\delbar_E,h)$ be an acceptable bundle over $X-D$.
For $\veca\in\real^{\ell}$
and $N\in\real$,
we set
\begin{equation}
\label{eq;a11.9.7}
 \tau(\veca,N):=
 -\sum_{i=1}^la_i\cdot \log|z_i|^2
 +N\cdot
 \Bigl(
 \sum_{i=1}^{\ell}\log(-\log|z_i|^2)
+\sum_{i=\ell+1}^n\log(1-|z_i|^2)
 \Bigr).
\end{equation}
\index{function $\tau(\veca,N)$}
We consider the following metric
for any $\veca$ and $N$:
\begin{equation} 
\label{eq;04.12.27.100}
 h_{\veca,N}:=h\cdot e^{-\tau(\veca,N)}
=h\cdot \prod_{i=1}^{\ell}
 |z_i|^{2a_i}\cdot\bigl(-\log|z_i|^2\bigr)^{-N}
\cdot 
 \prod_{i=\ell+1}^{n}
 \bigl(1-|z_i|^2\bigr)^{-N}
\end{equation}
We use the symbols
$|\cdot|_{\veca,N}$,
$||\cdot||_{\veca,N}$,
$(\cdot,\cdot)_{\veca,N}$ and
$\doublelangle\cdot,\cdot\doublerangle_{\veca,N}$
instead of
$|\cdot|_{h_{\veca,N}}$,
$||\cdot||_{h_{\veca,N}}$,
$(\cdot,\cdot)_{h_{\veca,N}}$ and
$\langle\cdot,\cdot\rangle_{h_{\veca,N}}$.
We also use the symbols $A_{\veca,N}^{p,q}(E)$
instead of $A_{h_{\veca,N}}^{p,q}(E)$.
\index{metric $h_{\veca,N}$}
\index{norm $| \cdot |_{\veca,N}$}
\index{$L^2$-norm $||\cdot||_{\veca,N}$}
\index{inner product $(\cdot,\cdot)_{\veca,N}$}
\index{$L^2$-inner product
 $langle\cdot,\cdot\rangle_{\veca,N}$}
\index{space $A_{\veca,N}^{p,q}(E)$}

\subsection{The case in which $N$ is sufficiently negative}
\label{subsection;04.12.23.50}

We have the following equality 
for any section $\eta\in A^{0,q}_c(E)$
(See the proof of Lemma 2.45 of \cite{mochi2}):
\footnote{The signature in Lemma 2.45 of
\cite{mochi2} is opposite.}
\begin{equation}
 \label{eq;08.10.16.1}
 \sqrt{-1}
 \Bigl(
 \Lambda \bigl(
 \delbar\del\tau(\veca,N)
 \eta
   \bigr)
-\Lambda\bigl(
 \delbar\del\tau(\veca,N)\bigr)
 \cdot\eta
 \Bigr)
 =
-N\cdot q\cdot \eta
\end{equation}
We obtain the following inequality
from (\ref{eq;08.10.16.2}) and (\ref{eq;08.10.16.1}):
\begin{multline}
\label{eq;11.9.3}
 \bigl\|
 \delbar_E\eta
 \bigr\|_{\veca,N}^2
+\bigl\|
 \delbar^{\ast}_E\eta
 \bigr\|_{\veca,N}^2 \geq
 \\
 \sqrt{-1}
 \bigl\langle
  \Lambda(\Omega_{X-D}^{0,q})\cdot\eta,\,\,\eta
 \bigr\rangle_{\veca,N}
-\sqrt{-1}
 \Bigl\langle
 \Lambda\bigl(R(h)\cdot\eta\bigr)
-\Lambda R(h)\cdot\eta,\,\,
\eta
 \Bigr\rangle_{\veca,N}
-N\cdot q\cdot||\eta||_{\veca,N}^2
\end{multline}
There exists a positive constant $C>0$,
depending only on the bound of the curvature $R(h)$,
such that the following holds
on $X-D$ for any $q=1,\ldots,n$
and for any $\eta\in A_c^{0,q}(E)$:
\[
 \Bigl|
 \Bigl\langle
 \Lambda(\Omega_{X-D}^{0,q})\cdot
 \eta-\Lambda\bigl(R(h)\cdot\eta\bigr)
 +\Lambda(R(h))\cdot\eta,\,\,
\eta
 \Bigr\rangle_{h}
 \Bigr|
\leq
 C\cdot|\eta|_{h}^2
\]
If we take a negative integer $N$
such that $N<-C-1$ for the above constant $C$,
we obtain the following inequalities for any $q\geq 1$
and for any $\eta\in A_{c}^{0,q}(E)$,
due to the inequality (\ref{eq;11.9.3}):
\begin{equation}
 \label{eq;11.9.5}
 \bigl\|
  \delbar_E\eta
 \bigr\|^2_{\veca,N}
+\bigl\|
 \delbar_E^{\ast}\eta
 \bigr\|_{\veca,N}^2
\geq ||\eta||_{\veca,N}^2
\end{equation}

\begin{lem} 
\label{lem;9.10.72}
Let $C$ be a positive constant as above.
If $N<-C-1$,
we have the vanishing of
any higher cohomology group
$H^i\bigl(A^{0,\cdot}_{\veca,N}\bigl(
 E\bigr),\delbar_{E}\bigr)$ $(i>0)$.
\end{lem}
\pf
It follows from Proposition \ref{prop;11.28.6}
and Proposition \ref{prop;12.11.1}
and (\ref{eq;11.9.5}).
\hfill\qed

\subsubsection{Refinement}

Let $Z_{\veca,N}^{0,q}(E)$ denote
the kernel of 
the natural morphism
$\delbar_E:
 A_{\veca,N}^{0,q}(E)\lrarr
 A_{\veca,N}^{0,q+1}(E)$.
Let $A_{\veca,N}^{\ast,0,q}(E)$ be the space of
the sections $\tau$ of $E\otimes\Omega^{0,q}$
such that $\tau$ and $\delbar_E^{\ast}\tau$ are $L^2$
with respect to $h_{\veca,N}$ and $g_{\poin}$.
It is also obtained as the completion of
$A_c^{0,q}(E)$ with respect to the norm
$\|\tau\|_{\veca,N}+\|\delbar_E^{\ast}\tau\|_{\veca,N}$.
The kernel of the morphism
$\delbar_E^{\ast}:
 A_{\veca,N}^{\ast,0,q}(E)\lrarr 
 A_{\veca,N}^{\ast,0,q-1}(E)$
is denoted by $Z^{\ast,0,q}_{\veca,N}(E)$.

Let $\Atilde^{0,q}_{\veca,N}(E)$ denote the space
of the sections of $E\otimes\Omega^{0,q}$
such that $\tau$, $\delbar_E\tau$
and $\delbar_E^{\ast}\tau$ are $L^2$.
It is also obtained as the completion of
$A^{0,q}_c(E)$ with respect to the norm
$\|\tau\|_{\Atilde,\veca,N}':=\|\tau\|_{\veca,N}
+\|\delbar_E\tau\|_{\veca,N}
+\|\delbar_E^{\ast}\tau\|_{\veca,N}$.
We set
\[
\Ztilde^{0,q}_{\veca,N}(E):=
 \Ker\Bigl(
\delbar_E:\Atilde_{\veca,N}^{0,q}(E)
 \lrarr A^{0,q+1}_{\veca,N}(E) \Bigr),
\]
\[
 \Ztilde^{\ast,0,q}_{\veca,N}(E):=
 \Ker\Bigl(
\delbar_E^{\ast}:\Atilde_{\veca,N}^{0,q}(E)
 \lrarr A^{\ast,0,q-1}_{\veca,N}(E) \Bigr).
\]
If $N$ is as in Lemma \ref{lem;9.10.72},
for any $q\geq 1$,
$\|\tau\|_{\Atilde,\veca,N}:=
 \| \delbar_E\tau\|_{\veca,N}
+\|\delbar_E^{\ast}\tau\|_{\veca,N}$
determines the norm equivalent to
$\|\tau\|_{\Atilde,\veca,N}'$
due to (\ref{eq;11.9.5}).
The hermitian pairing
$\langle\tau_1,\tau_2\rangle_{\Atilde,\veca,N}
:=(\delbar_E\tau_1,\delbar_E\tau_2)_{\veca,N}
 +(\delbar_E^{\ast}\tau_1,
     \delbar_E^{\ast}\tau_2)_{\veca,N}$
induces the norm $\|\cdot\|_{\Atilde,\veca,N}$.

\begin{lem}
\label{lem;07.9.24.2}
If $N$ is as in Lemma {\rm\ref{lem;9.10.72}},
the following holds for any $q\geq 1$
and any $\veca\in\real^{\ell}$.
\begin{itemize}
\item
We have the decomposition
$\Atilde_{\veca,N}^{0,q}(E)
=\Ztilde_{\veca,N}^{0,q}(E)\oplus
 \Ztilde_{\veca,N}^{\ast,0,q}(E)$.
\item
 $\delbar_E: \bigl(
 \Ztilde_{\veca,N}^{\ast,0,q}(E),
 \|\cdot\|_{\Atilde,\veca,N} \bigr)
 \lrarr
 \bigl( Z_{\veca,N}^{0,q+1}(E),\|\cdot\|_{\veca,N}\bigr)$
is an isomorphism.
\item
For any $\rho\in\Ztilde^{\ast,0,q}_{\veca,N}(E)$,
there exists $\mu\in \Ztilde^{0,q+1}_{\veca,N}(E)$
such that $\rho=\delbar_E^{\ast}\mu$.
\end{itemize}
\end{lem}
\pf
We use a descending induction on $q$.
The claim is trivial for any sufficiently large $q$.
Assume that we have already obtained
the decomposition
$\Atilde^{0,q+1}_{\veca,N}(E)
=\Ztilde^{0,q+1}_{\veca,N}(E)
\oplus
 \Ztilde^{\ast,0,q+1}_{\veca,N}(E)$,
and $\Ztilde^{\ast,0,q+1}_{\veca,N}(E)
\subset\delbar_E^{\ast}\bigl(
 \Ztilde^{0,q+2}_{\veca,N}(E)
 \bigr)$.
For any $\kappa \in Z^{0,q+1}_{\veca,N}(E)$,
let us consider the linear map
$F_{\kappa}:\Ztilde^{0,q+1}_{\veca,N}(E)
\lrarr \cnum$
given by
$F_{\kappa}(\tau):=(\tau,\kappa)_{\veca,N}$.
It is continuous with respect to $\|\cdot\|_{\Atilde,\veca,N}$.
By Riesz representation theorem,
we have $\nu\in \Ztilde^{0,q+1}_{\veca,N}(E)$
such that
$(\tau,\kappa)_{\veca,N}
=\langle\tau,\nu\rangle_{\Atilde,\veca,N}
=(\delbar_E^{\ast}\tau,\delbar_E^{\ast}\nu)_{\veca,N}$
for any $\tau\in\Ztilde^{0,q+1}_{\veca,N}(E)$.
For any
$\tau\in \Ztilde^{\ast,0,q+1}_{\veca,N}(E)$,
we obviously have
$(\delbar_E^{\ast}\tau,
 \delbar_E^{\ast}\nu)_{\veca,N}=0$,
and we also have 
$(\tau,\kappa)_{\veca,N}=0$
because $\tau=\delbar_E^{\ast}\tau'$
for some $\tau'$.
Hence, 
$(\tau,\kappa)_{\veca,N}
=(\delbar_E^{\ast}\tau,\delbar_E^{\ast}\nu)_{\veca,N}$
holds for any $\tau\in \Atilde^{0,q+1}_{\veca,N}(E)$,
which implies
$\kappa=\delbar_E\bigl(\delbar_E^{\ast}\nu\bigr)$.
In particular,
$\kappa$ is contained in the image of
$\delbar_E:
 \Ztilde^{\ast,0,q}_{\veca,N}(E)
\lrarr Z^{0,q+1}_{\veca,N}(E)$.
Then,
$\delbar_E: \Ztilde^{\ast,0,q}_{\veca,N}(E)
\lrarr Z^{0,q+1}_{\veca,N}(E)$ is 
an isomorphism,
by definition of $\|\cdot\|_{\Atilde,\veca,N}$
and $\|\cdot\|_{\veca,N}$.

Due to (\ref{eq;11.9.5}),
we have 
$\Ztilde^{0,q}_{\veca,N}(E)\cap
 \Ztilde^{\ast,0,q}_{\veca,N}(E)=0$.
Let $\omega\in \Atilde^{0,q}_{\veca,N}(E)$.
Applying the above result to $\delbar_E\omega$,
we can take $\omega_1\in \Ztilde^{\ast,0,q}_{\veca,N}(E)$
such that $\delbar_E\omega_1=\delbar_E\omega$.
Then, $\omega-\omega_1\in \Ztilde^{0,q}_{\veca,N}(E)$.
Thus, we obtain the decomposition.
$\Atilde^{0,q}_{\veca,N}(E)
=\Ztilde^{0,q}_{\veca,N}(E)
\oplus
 \Ztilde^{\ast,0,q}_{\veca,N}(E)$.

Let $\rho\in\Ztilde^{\ast,0,q}_{\veca,N}(E)$.
Applying the above to 
$\delbar_E\rho\in Z^{0,q+1}_{\veca,N}(E)$,
we can find $\nu\in \Ztilde^{0,q+1}_{\veca,N}(E)$
such that $\delbar_E\delbar_E^{\ast}\nu=\delbar_E\rho$.
Then, $\delbar_E^{\ast}\nu\in
 \Ztilde^{\ast,0,q}_{\veca,N}(E)$
and we have $\delbar_E(\delbar_E^{\ast}\nu-\rho)=0$.
Then, we obtain $\rho=\delbar_E^{\ast}\nu$.
Thus, the induction can proceed.
\hfill\qed

\begin{lem}
\label{lem;07.9.23.36}
Let $N$ be as in Lemma {\rm\ref{lem;9.10.72}}.
For any  $\omega\in Z^{0,1}_{\veca,N}(E)$,
we have a section $f$ of $E$
such that 
$\delbar f=\omega$
and $\|f\|_{\veca,N}\leq \|\omega\|_{\veca,N}$.
\end{lem}
\pf
Let us consider the linear map
$F_{\omega}:\Ztilde^{0,1}_{\veca,N}(E)\lrarr \cnum$
as in the proof of Lemma \ref{lem;07.9.24.2},
which is continuous with respect to
$\|\cdot\|_{\Atilde,\veca,N}$.
By Riesz representation theorem,
we have $\rho\in \Ztilde^{0,1}_{\veca,N}(E)$
with $\|\rho\|_{\Atilde,\veca,N}\leq \|\omega\|_{\veca,N}$,
such that
$(\tau,\omega)_{\veca,N}
=(\delbar_E^{\ast}\tau,
 \delbar_E^{\ast}\rho)_{\veca,N}$.
As in the proof of Lemma \ref{lem;07.9.24.2},
we have $\omega=\delbar_E(\delbar_E^{\ast}\rho)$.
Since we have
$\|\delbar_E^{\ast}\rho\|_{\veca,N}
=\|\rho\|_{\Atilde,\veca,N}\leq \|\omega\|_{\veca,N}$,
we are done.
\hfill\qed

\subsection{The case in which $N$ is sufficiently positive}

Let $\pi_i:X-D\lrarr D_i$ denote the natural projection
for $i=1,\ldots,n$.
We put $D_i^{\circ}:=
 D_i\setminus \bigcup_{j\neq i,j\leq \ell}D_j$.
Let $P$ be a point of $D_i^{\circ}$,
then we obtain the curve $\pi_i^{-1}(P)$
which is isomorphic to $\Delta^{\ast}$ $(i\leq \ell)$
or $\Delta$ $(\ell<i\leq n)$.
Let $\Delta''$ denote the Laplacian
on $\pi_i^{-1}(P)$
with respect to the Euclidean metric.
The restriction of the metric $h_{\veca,N}$
to $E_{|\pi^{-1}_i(P)}$
is also denoted by $h_{\veca,N}$.
Similarly, the restriction of 
$|\cdot|_{\veca,N}$
to $\pi_i^{-1}(P)$ is denoted by the same symbol.
We can easily show the following lemma.
(See Lemma 2.49 and Corollary 2.50 of \cite{mochi2}.)
\begin{lem} 
\label{lem;04.12.23.1}
There exists a positive constant $N_0$,
depending only on the estimate of the curvature $R(h)$,
such that the following holds:
\begin{itemize}
\item
Let $P$ be a point of $D_i^{\circ}$,
and $U$ be an open subset of the curve $\pi_i^{-1}(P)$.
Let $F$ be a holomorphic section of $E_{|U}$.
Then the inequalities 
$\Delta''|F|_{\veca,N}^2\leq 0$
and 
$\Delta''\log |F|_{\veca,N}\leq 0$
hold on $U$ for any $N\geq N_0$.
In particular, $|F|^2_{\veca,N}$ and 
$\log |F|_{\veca,N}$ are subharmonic on $U$.
\end{itemize}
As a result,
for any holomorphic section $F$ of $E$,
the functions
$\log|F|_{\veca,N}$ and $|F|_{\veca,N}$
are pluri-subharmonic, if $N\geq N_0$.
\hfill\qed
\end{lem}

We explain various consequences of 
this pluri-subharmonicity.
Since we will be interested 
in the behaviour of holomorphic sections
around the origin $O$,
we put 
$X(C):=\bigl\{(z_1,\ldots,z_n)\in X\,\big|\,|z_i|\leq C\bigr\}$,
$D_i(C):=D_i\cap X(C)$ and $D(C)=D\cap X(C)$
for any $0<C<1$,
and we will consider the restrictions of holomorphic sections 
to $X(C)^{\ast}:=X(C)-D(C)$.

\subsubsection{$L^2$-estimates on curves and $X-D$}

The following corollary makes us
to derive the $L^2$-property on curves
from the $L^2$-property on $X-D$,
which easily follows from the pluri-subharmonicity
in Lemma \ref{lem;04.12.23.1}.
(See Corollary 2.51 of \cite{mochi2}.)
\begin{cor} 
\label{cor;04.12.23.30}
Let $F$ be a holomorphic section of $E$ on $X-D$
such that $F\in A^{0,0}_{\veca,N}(E)$.
Let $N_0$ be as in Lemma {\rm\ref{lem;04.12.23.1}},
and let $M>\max\{ N_0,N\}$.
For any $1\leq j\leq \ell$ and $P\in D_j^{\circ}$,
we have the following finiteness:
\[
 \int_{\pi_j^{-1}(P)\cap X^{\ast}(C)}
 \bigl|
  F_{|\pi_j^{-1}(P)\cap X^{\ast}(C)}
 \bigr|^2_{\veca,M}
 \dvol_{\pi_j^{-1}(P)}<\infty
\]
Here $\dvol_{\pi_j^{-1}(P)}$ denotes
the volume form of $\pi_j^{-1}(P)$
with respect to the restriction $g_{\poin\,|\,\pi_j^{-1}(P)}$.
\hfill\qed
\end{cor}

\subsubsection{From $L^2$-estimate to growth estimate
 on a curve}

Let us consider the case $X=\Delta$ and $D=\{O\}$.
We can show the following lemma
by the argument in Lemma 7.12 of \cite{mochi4},
although it is stated for acceptable bundles
induced by tame harmonic bundles.
\begin{lem}
\label{lem;08.10.17.11}
Let $f$ be a holomorphic section of $E$ on $X-D$
such that $\|f\|_{b,N}<\infty$ for some 
$b,N\in\real$.
Let $N_0$ be as in Lemma {\rm\ref{lem;04.12.23.1}},
and $M>\max\{N_0,|N|+2\}$.
Then, the following holds on $X(1/4)$:
\[
 |f(z)|_{h}^2
\leq
 B \|f\|_{b,N}^2
 |z|^{-2b} (-\log|z|)^{M}
\]
Here, $B>0$ is independent of $f$.
\end{lem}
\pf
Let $\dvol$ (resp. $\dvol_{g_{\poin}}$)
denote the volume form
associated to Euclidean metric
(resp. Poincar\'e metric).
From the subharmonicity,
we obtain the following inequalities
for any $0<|z|\leq 1/4$:
\begin{multline}
 \log|f(z)|_{b,M}^2
\leq
 \frac{4}{\pi|z|^2}
\int_{|w-z|\leq |z|/2}
 \log|f(w)|^2_{b,M}\dvol  \\
\leq
 \log\Bigl(
  \frac{4}{\pi|z|^2}
\int_{|w-z|\leq |z|/2}
 |f(w)|^2_{b,M}\dvol
 \Bigr) 
\leq
\log\Bigl(
  \frac{20}{\pi}
\int_{|w-z|\leq |z|/2}
 |f(w)|^2_{b,N}\dvol_{g_{\poin}}
 \Bigr) \\
\leq \log\Bigl(
  \frac{20}{\pi}
 \|f\|_{b,N}
 \Bigr)
\end{multline}
Thus, we are done.
\hfill\qed

\subsubsection{Refinement of growth estimates
 on a curve}

Let us consider the case $X=\Delta$ and $D=\{O\}$.
We can show the following lemma 
by the argument in Lemma 7.13 of \cite{mochi4},
although it is stated for acceptable bundles
induced by tame harmonic bundles.
\begin{lem}
\label{lem;08.10.17.10}
Let $f$ be a holomorphic section 
such that
$|f|_h=O\bigl(
 |z|^{-a-\epsilon}
 \bigr)$ for any $\epsilon>0$ on $X(C)$.
Let $N_0$ be as in Lemma {\rm\ref{lem;04.12.23.1}},
and let $M>N_0$.
We set 
$H(z):=|f(z)|_h^2\cdot
 |z|^{2a}\cdot(-\log|z|)^{-M}$.
Then, 
$H(z)$ is bounded.
More strongly,
$\max_{|z|\leq C'}|H(z)|
=\max_{|z|=C'}H(z)$
for any $0<C'<C$.
\end{lem}
\pf
We put 
$H_{\epsilon}:=
 H(z)\cdot |z|^{2\epsilon}$
for any $\epsilon>0$.
We have the subharmonicity of
$\log H_{\epsilon}$ on $X^{\ast}(C)$.
By the assumption,
we have $\lim_{|z|\to 0}\log H_{\epsilon}(z)=-\infty$.
Hence, $\max_{|z|\leq C'}|H_{\epsilon}(z)|
=\max_{|z|=C'}H_{\epsilon}(z)$ for any $0<C'<C$.
By taking limit $\epsilon\to 0$,
we obtain the desired estimate.
\hfill\qed

\subsubsection{Growth estimates on curves and $X-D$}

Let us return to the case,
$X=\Delta^n$ and $D=\bigcup_{i=1}^{\ell} D_i$.
We would like to remark that
the growth estimates on curves imply
the growth estimate on $X-D$.
The following lemma is a refinement of
Corollary 2.53 of \cite{mochi2}.

\begin{prop}
\label{prop;08.10.17.12}
\label{prop;07.7.7.3}
Let $F$ be a holomorphic section of $E$ on $X-D$.
Assume that we are given numbers
$a_i\in\real$ $(i=1,\ldots,\ell)$
such that the following holds:
\begin{itemize}
\item
$\bigl|F_{|\pi_i^{-1}(P)}\bigr|_h
=O\bigl(
 |z_i|^{-a_i-\epsilon}
 \bigr)$ for any $\epsilon>0$
 and for any $i=1,\ldots,\ell$
 and $P\in D_i^{\circ}$.
\end{itemize}
Let $N_0$ be as in Lemma {\rm\ref{lem;04.12.23.1}},
and let $M>N_0$.
Take $0<C<1$.
Then, there exists a constant $B$,
which is independent of $F$,
such that the following holds on $X^{\ast}(C)$:
\[
 \bigl|F\bigr|_h^2\leq
 B\cdot
 \prod_{j=1}^{\ell}
 |z_j|^{-2a_j}\cdot
 \bigl(-\log|z_j|\bigr)^{M}
\cdot
 \max_{\substack{|z_j|=C\\ 1\leq j\leq \ell}}
 \bigl|F\bigr|_h^2
\]
\end{prop}
\pf
We have only to use Lemma \ref{lem;08.10.17.10}
inductively.
\hfill\qed

\begin{cor}
Let $F$ be a holomorphic section of $E$ over $X-D$
such that $\|F\|_{\veca,N}<\infty$.
Let $N_0$ be as in Lemma {\rm\ref{lem;04.12.23.1}},
and let $M>N_0$.
Then, the following holds
on $X^{\ast}(C)$ for some $0<C<1$:
\[
 \bigl|F\bigr|_{h}^2
\leq
 B\cdot
 \prod_{j=1}^{\ell}
 |z_j|^{-2a_j}\cdot
 \bigl(-\log|z_j|\bigr)^{M}\cdot
 \max_{\substack{|z_j|=C\\ 1\leq j\leq \ell}}
 \bigl|F\bigr|_h^2
\]
In particular,
$F\in \prolongg{\veca}{E}$.
(See Section {\rm\ref{subsection;07.9.22.10}}
for $\prolongg{\veca}{E}$.)
\end{cor}
\pf
We obtain the $L^2$-estimates
for the restrictions of $F$ to curves
transversal to the smooth part of $D$,
due to Corollary {\rm\ref{cor;04.12.23.30}}.
It implies the growth estimates
of the restriction to curves
due to Lemma {\rm\ref{lem;08.10.17.11}}.
Then, we obtain the growth estimate
of $|F|_h$ on $X^{\ast}(C)$
by Proposition {\rm\ref{prop;08.10.17.12}}.
\hfill\qed

\section{Prolongation to filtered bundle}
\label{subsection;07.9.22.10}
Let $X:=\Delta^n$,
and $D:=\bigcup_{i=1}^{\ell}\{z_i=0\}$.
Let $(E,\delbar_E,h)$ be an acceptable bundle on $X-D$
of rank $r$.
We naturally identify $E$
with the sheaf of holomorphic sections.
Let $\veca\in\real^{\ell}$.
The $i$-th components of $\veca$ are 
denoted by $a_i$.
For any open subset $U\subset X$,
we define
\[
 \prolongg{\veca}{E}(U)
:=\Bigl\{
 f\in E(U\setminus D)\,\Big|\,
 |f|_h=O\Bigl(\prod_{i=1}^{\ell}|z|^{-a_i-\epsilon}\Bigr)\,\,
 \forall \epsilon>0
 \Bigr\}.
\]
By taking sheafification,
we obtain an $\nbigo_X$-module
$\prolongg{\veca}{E}$.
In the case $\veca=(0,\ldots,0)$,
we prefer the symbol $\prolong{E}$.
\index{sheaf $\prolongg{\veca}{E}$}
\index{sheaf $\prolong{E}$}

\begin{thm}
\label{thm;07.3.16.1}
Assume that 
$\det(E,\delbar_E,h)$ is flat,
for simplicity.
Then,
$\prolongg{\veca}{E}$ is locally free
for any $\veca\in\real^{\ell}$.
Moreover,
the family $\bigl\{\prolongg{\veca}{E}\,\big|\,
 \veca\in\real^{\ell}\bigr\}$ forms a filtered bundle.
\end{thm}
It implies that
(i) the images
$\lefttop{i}F_b(\prolongg{\veca}{E}_{|D_i})$ of
$\prolongg{\veca(i,b)}{E}_{|D_i}
\lrarr
 \prolongg{\veca}{E}_{|D_i}$
are subbundles for any $b\in \openclosed{a_i-1}{a_i}$,
where the $j$-the components of
$\veca(i,b)$ and $\veca$ are equal for $j\neq i$,
and the $i$-th component of $\veca(i,b)$ is $b$,
(ii) the induced filtrations
$\lefttop{i}F$ $(i=1,\ldots,\ell)$ of 
$\prolongg{\veca}{E}_{|D_i}$ are 
compatible in the sense of 
Definition {\rm 4.37} of {\rm \cite{mochi2}}.
(See Section \ref{subsection;10.5.28.1}.)
\index{filtered bundle}
\index{filtration $\lefttop{i}F$}

The tuple of filtrations
$\bigl(\lefttop{i}F\,\big|\,i=1,\ldots,\ell\bigr)$
is often denoted by $\vecF$ in this situation.
Let $\vecv=(v_1,\ldots,v_r)$ be a frame of 
$\prolongg{\veca}{E}$ compatible with $\vecF$.
Namely, the numbers $a_i(v_k)$ are attached to 
any $v_k$ and $i=1,\ldots,\ell$,
such that 
$\bigl\{
 v_{k|D_i}\,\big|\,
 a_i(v_k)\leq b
 \bigr\}$ gives a frame of
$\lefttop{i}F_b(\prolongg{\veca}{E}_{|D_i})$
for each pair of $i=1,\ldots,\ell$ and
$b\in\openclosed{a_i-1}{a_i}$.
\index{frame compatible with filtrations}
The numbers $a_j(v_i)$
are often denoted by 
$\lefttop{j}\deg^{F}(v_i)$ in this situation,
and called the degree
with respect to $\lefttop{j}F$.
\index{degree $\lefttop{j}\deg^{F}(v_i)$}
We set 
\begin{equation}
 \label{eq;08.10.16.10}
 v_k':=v_k\cdot \prod_{i=1}^{\ell}|z_i|^{a_i(v_k)}.
\end{equation}
Let $H(h,\vecv')$ denote the hermitian-matrix valued
function whose $(p,q)$-entries are given by
$h(v_p',v_q')$.

\begin{thm}
\label{thm;07.10.9.1}
There exist positive constants $C$
and $N$ such that the following holds
around $O$:
\begin{equation}
\label{eq;08.10.17.13}
 C^{-1}\cdot
\left(-\sum_{i=1}^{\ell}\log|z_i|\right)^{-N}
\leq
 H(h,\vecv')
\leq 
 C\cdot
\left(-\sum_{i=1}^{\ell}\log|z_i|\right)^{N}
\end{equation}
In other words,
$\vecv'$ is adapted to $h$
up to log order.
\end{thm}
\index{adapted up to log order}

The holomorphic bundle $\End(E)$
with the induced hermitian metric $h$
is also acceptable.
We can show the following proposition
by using the weak norm estimate
(Theorem \ref{thm;07.10.9.1}).
\begin{prop}
\label{prop;07.11.18.5}
$\prolong{\End(E)}$
is naturally isomorphic to the sheaf
of the endomorphisms $f$
of $\prolongg{\veca}{E}$
for any $\veca\in\real^{\ell}$,
such that $f_{|D_i}$ preserves
the filtration $\lefttop{i}F$
for each $i=1,\ldots,\ell$.
\hfill\qed
\end{prop}

We will prove Theorems \ref{thm;07.3.16.1}
and \ref{thm;07.10.9.1}
in Sections 
\ref{subsection;08.10.18.20}--\ref{subsection;08.10.18.21}.

\section{One dimensional case}
\label{subsection;08.10.18.20}
\subsection{Weak norm estimate}
\label{subsection;08.10.19.120}

In the one dimensional case,
Theorem \ref{thm;07.3.16.1}
was proved by Simpson in \cite{s1}.
(See also \cite{s2}).
For any $a\in\real$,
we have the associated vector bundle
$\prolongg{a}{E}$
with the parabolic filtration $F$.
We set 
\[
 \Par(\prolongg{a}{E}):=
 \bigl\{
 b\in \openclosed{a-1}{a}\,\big|\,
 \Gr^F_b(\prolongg{a}{E})\neq 0
 \bigr\}
\]
\index{set $\Par(\prolongg{a}{E})$}
Moreover, he showed the functoriality
with respect to dual, tensor and direct sum.
In particular,
if $\vecv$ is a frame of $\prolongg{a}{E}$ 
compatible with the parabolic filtration,
the dual frame $\vecv^{\lor}$ of $E^{\lor}$
naturally gives a frame of
$\prolongg{-a-1+\epsilon}{E^{\lor}}$
for some sufficiently small $\epsilon>0$.
Let us observe
Theorem \ref{thm;07.10.9.1}
in the one dimensional case, by using this fact.
Let $\vecv'=(v_i')$ be as in (\ref{eq;08.10.16.10}).
By using Lemma \ref{lem;08.10.17.10},
we obtain
$ H(h,\vecv')\leq
 C_1\cdot\bigl(-\log|z|\bigr)^{M_1}$
for some $C_1,M_1>0$.
We also obtain
$H(h^{\lor},\vecv^{\lor\prime})
\leq C_2\cdot \bigl(-\log|z|\bigr)^{M_2}$,
where $h^{\lor}$ denotes the induced metric of $E^{\lor}$.
Then, we obtain that $\vecv'$  is adapted.

\subsection{Pull back and descent}
\label{subsection;08.10.18.40}

For any positive integer $c$,
let $\psi_c:X\lrarr X$ be given by $\psi_c(z)=z^c$.
We put $(\Etilde,\delbar_{\Etilde},\htilde):=
 \psi_c^{-1}(E,\delbar_E,h)$,
which is also acceptable.
For any $a,b\in\real$,
let $\nu(a,b)\in\seisuu$ be determined by
$b-1<\nu(a,b)+a\leq b$.
For any $b\in\real$,
it is easy to check
\[
 \prolongg{b}{\Etilde}
=\sum_{a}
 z^{-\nu(ca,b)}\cdot
 \psi_c^{-1}(\prolongg{a}{E}).
\]
More concretely,
let $\vecv$ be a frame of
$\prolong{E}$ compatible
with the parabolic filtration $F$.
Let $a_i:=\deg^F(v_i)$.
We set
\[
 \vtilde_i:=
 z^{-\nu(ca_i,b)}\cdot \psi_c^{-1}(v_i).
\]
It is easy to check that
$\vecvtilde=(\vtilde_i)$ is a frame of
$\prolongg{b}{\Etilde}$
compatible with the parabolic filtration,
by using the weak norm estimate.
Hence, we have
\[
 \Par\bigl(\prolongg{b}{\Etilde}\bigr)
=\bigl\{
 ca+\nu(ca,b)\,\big|\,
 a\in\Par\bigl(\prolong{E}\bigr)
 \bigr\}.
\]

Let $\mu_c:=\bigl\{
 w\in\cnum\,\big|\,w^c=1
 \bigr\}$.
Let $\omega$ be a generator.
We have the natural $\mu_c$-action on $X$
given by the multiplication,
which is lifted to the action on $\prolongg{b}{\Etilde}$.
Note $(\omega^{-1})^{\ast}(\vtilde_i)=
 \omega^{\nu(ca_i,b)}\vtilde_i$.

Assume that $0\leq b<1/2$.
If $c$ is a sufficiently large integer,
we have (i) $0\leq \nu(ca,b)\leq c-1$,
(ii) $\Par(\prolong{E})\lrarr \seisuu$
 given by $\nu(ca,b)$ is injective.
We have the decomposition
\begin{equation}
 \label{eq;08.10.17.15}
 \prolongg{b}{\Etilde}_{|O}
=\bigoplus_{0\leq p\leq c-1} V_p,
\quad
 V_p=\bigl\langle
 \vtilde_{j|O}\,\big|\,
 \nu(c a_j,b)=p
 \bigr\rangle,
\end{equation}
where $\omega$ acts on
$V_p$ as the multiplication of $\omega^p$.

For each $p\in 
\bigl\{
 0\leq p\leq c-1\,\big|\,
 V_p\neq 0
 \bigr\}$,
there exists $\chi(p)\in \Par(\prolong{E})$
such that $p=\nu\bigl(\chi(p)a, b\bigr)$.
Thus, we obtain the map
\[
 \chi:
  \bigl\{
 0\leq p\leq c-1\,\big|\, V_p\neq 0
 \bigr\}
\lrarr
 \Par\bigl(\prolong{E}\bigr).
\]
By setting
$\varphi(p):=
 \nu\bigl(\chi(p)a, b\bigr)+\chi(p)a\in
 \Par(\prolongg{b}{E})$,
we obtain the map
\[
 \varphi:
 \bigl\{
 0\leq p\leq c-1\,\big|\, V_p\neq 0
 \bigr\}
\lrarr
 \Par(\prolongg{b}{\Etilde}).
\]
The decomposition (\ref{eq;08.10.17.15})
gives a splitting of the parabolic filtration
$F$ of $\prolongg{b}{\Etilde}$
in the following sense:
\[
 F_d(\prolongg{b}{\Etilde}_{|O})= 
 \bigoplus_{\varphi(p)\leq d}
 V_p
\]
Conversely,
let $\vecutilde$ be a frame of 
$\prolongg{b}\Etilde$.
Assume that it is $\mu_c$-equivariant
in the sense that
$\omega^{\ast}\utilde_j=
 \omega^{-p_j}\cdot \utilde_j$
for some $0\leq p_j\leq c-1$.
Note that $\utilde_{j|O}\in V_{p_j}$,
and hence the frame $\vecutilde$ is 
compatible with the filtration $F$.
We set $u_j:=z^{p_j}\cdot \utilde$.
Since they are $\mu_c$-invariant,
they induce sections of $E$,
which is also denoted by $u_j$.
It is easy to check that 
$\vecu=(u_j)$ is a frame of $\prolong{E}$
compatible with the parabolic filtration
by using the weak norm estimate.

\subsection{Parabolic structure and monodromy}
\label{subsection;08.10.18.41}

Let $X:=\Delta$ and $D:=\{O\}$.
Let $(E,\delbar_E,h)$ be an acceptable
bundle of rank $R$ on $X-D$.
\begin{lem}
\label{lem;08.10.3.100}
There exists a $C^1$-orthonormal frame
$\vece$ of $E$
with the following property.
\begin{itemize}
\item
There exists a diagonal matrix $\Gamma$
whose $(i,i)$-entries $\alpha_i$
satisfy $0\leq \alpha_i<1$.
\item
 Let $A$ be determined by the following:
\[
 \delbar_E\vece
=\vece\cdot \Bigl(
 -\frac{\Gamma}{2}+A
 \Bigr)\cdot \frac{d\zbar}{\zbar}
\]
Then $A=O\bigl((-\log|z|)^{-1}\bigr)$.
\end{itemize}
\end{lem}
\pf
See Lemma 7.10 of \cite{mochi4}
with the simplified proof due to the referee
of the paper.
Although the lemma is stated for tame harmonic bundles,
it holds for any acceptable bundles.
\hfill\qed

\vspace{.1in}
We set $S(\Gamma):=\bigl\{
 -\alpha_1,\ldots,-\alpha_R \bigr\}$.
\footnote{$S(\Gamma)$ in Page {\rm 75} 
of {\rm\cite{mochi4}}
 should be corrected.}
We recall the following lemma.
Although it is also stated
for tame harmonic bundles,
it holds for any acceptable bundles.
\begin{lem}[Lemma 7.17,
\cite{mochi4}]
\label{lem;08.10.3.1}
$\Par(\prolong{E})=S(\Gamma)$
holds.
The multiplicity of 
an eigenvalue $\alpha_i$
is equal to 
$\dim \Gr^F_{-\alpha_i}(E)$.
\hfill\qed
\end{lem}

For $0<r<1$,
let $P(r)$ denote the characteristic polynomial
of the monodromy along the loop
$r\cdot e^{\sqrt{-1}\theta}$
($0\leq \theta\leq 2\pi$)
of the unitary connection $\del_E+\delbar_E$.
We have the limit $\lim_{r\to 0}P(r)$
whose roots are $\alpha_1,\ldots,\alpha_R$.
By Lemma \ref{lem;08.10.3.1},
$\lim_{r\to 0}P(r)$ determines
$\Par(\prolong{E})$.

\subsection{Control of the parabolic weights
 in the higher dimensional case}
\label{subsection;08.10.18.42}

Let $X:=\Delta^n$ and $D:=\{z_1=0\}$.
Let $\pi:X\lrarr D$ be the projection.
Let $(E,\delbar_E,h)$ be an acceptable
bundle on $X-D$.
For any point $Q\in D$,
we have the acceptable bundle
$(E_Q,\delbar_Q,h_Q):=
 (E,\delbar_E,h)_{|\pi^{-1}(Q)\setminus D}$.

\begin{lem}
\label{lem;08.10.3.5}
The following holds for any $Q_i\in D$ $(i=1,2)$:
\[
\Par(\prolong{E_{Q_1}}):=
 \Par(\prolong{E}_{Q_2}),
\quad
\dim \Gr^F_a(\prolong{E}_{Q_1})
=\dim \Gr^F_a(\prolong{E}_{Q_2})
\]
\end{lem}
\pf
Let $\nabla:=\delbar_E+\del_E$
denote the unitary connection
associated to $(E,\delbar_E,h)$.
We use the natural identification
$X=\Delta_{z_1}\times D$.
For $0<r<1$,
let $P(r,Q_i)$ denote the characteristic polynomials
of the monodromy along the loops
$ \bigl(r\cdot e^{\sqrt{-1}\theta},Q_i\bigr)$
($ 0\leq\theta\leq 2\pi$).
We have only to show
$\lim_{r\to 0}P(r,Q_1)=
 \lim_{r\to 0}P(r,Q_2)$.

Let $\gamma(t):[0,1]\lrarr D$
be a $C^{\infty}$-map
such that $\gamma(0)=Q_1$
and $\gamma(1)=Q_2$.
For any $0<r<1$,
we have the map
\[
 \Phi_{r}:[0,1]\times [0,2\pi]
\lrarr
 X-D,
\quad
 \Phi_{r}(t,\theta)
=\bigl(r\cdot e^{\sqrt{-1}\theta},
 \gamma(t)\bigr)
\]
Take any $v\in E_{|(r,Q)}$ such that $|v|=1$.
We have the section $V(r)$ 
of $\Phi_{r}^{\ast}E$
determined by the following conditions:
\[
 V(r)_{|(r,Q)}=v,
\quad
 \nabla(\del_{\theta})V(r)=0,
\quad
\bigl(
 \nabla(\del_t)V(r)
\bigr)_{|\theta=0}=0
\]
Let $R$ denote the curvature of $\nabla$.
Then, we have the following:
\[
 \nabla(\del_{\theta})
 \nabla(\del_t)V(r)
=\Phi_{r}^{\ast}R(\del_t,\del_{\theta})V(r)
=O\bigl(
 (-\log r)^{-1}
 \bigr)
\]
Hence, we obtain the following
estimate with respect to $h$:
\[
  \nabla(\del_t)V(r)_{|\theta=2\pi}
=O\bigl(
 (-\log r)^{-1}
 \bigr)
\]
Then, we obtain
$\lim_{r\to 0}P(r,Q_1)=
 \lim_{r\to 0}P(r,Q_2)$.
\hfill\qed

\section{Extension of holomorphic sections I}
\label{subsection;08.10.18.45}
\subsection{Statement}

We set $X_0:=\Delta_z^{n-1}=
 \bigl\{(z_1,\ldots,z_{n-1})\,\big|\,
 |z_i|<1\bigr\}$,
and $D_0:=\bigcup_{i=1}^{\ell}\{z_i=0\}$
for some $\ell\leq n-1$.
Let $\Delta_w:=\bigl\{w\in\cnum\,\big|\,|w|<1\bigr\}$.
We put $X:=X_0\times\Delta_w$
and $D:=D_0\times \Delta_w$.
We also set
$X^{(0)}:=X_0\times\{0\}$
and $D^{(0)}:=D_0\times\{0\}$.

Let $(E,\delbar_E,h)$ be an acceptable bundle
on $X-D$.
Let $f$ be a holomorphic section of
$\prolong{(E_{|X^{(0)}-D^{(0)}})}$
on $X^{(0)}-D^{(0)}$.
Take $0<R<1$,
and we set $X(R):=\bigl\{|z_i|<R, |w|<R \bigr\}$.
The following lemma is the counterpart
of Lemma 8.51 of \cite{mochi2}.

\begin{lem}
\label{lem;08.10.3.10}
There exists a holomorphic section $F$ of 
$\prolong{E}$ on $X(R)$
such that $F_{|X(R)\cap(X^{(0)}-D^{(0)})}=f$.
\end{lem}

\subsection{Construction of a cocycle}
\label{subsection;08.10.18.50}

The argument is essentially the same
as that in \cite{mochi2}.
We have to modify
the construction of an appropriate cocycle 
(Section 8.4.2 of \cite{mochi2}).
We can argue as follows.

We use the polar coordinate 
$w=r\cdot e^{\sqrt{-1}\theta}$
for $\Delta_w$.
Let $\del_r$ denote the vector field $\del/\del r$.
Let $f$ be a holomorphic section of
$\prolong(E_{|X^{(0)}-D^{(0)}})$.
By using Proposition \ref{prop;08.10.17.12},
we obtain that $f$ is bounded up to log order,
i.e., the following holds for some $C,N>0$:
\[
 |f|_h\leq C\cdot\left(-\sum_{i=1}^{\ell}\log|z_i|\right)^N
\]
By using the unitary connection
$\nabla:=\delbar_E+\del_E$,
we extend $f$ to a continuous section $F$ of $E$ on $X-D$
by the condition $\nabla_{\del_r}F=0$ and $F_{|X_0}=f$.
We have $|F(z,w)|_h\leq |f(z)|_h$ 
by the construction.

\begin{lem}
We have the following estimate
on $X-\bigl(D\cup (X^{(0)})\bigr)$
for some $C>0$
with respect to $h$ and $g_{\poin}$:
\[
 \bigl|\nabla F(z,w)\bigr|\leq 
 C\cdot \bigl|f(z)\bigr|\cdot |w|
\]
\end{lem}
\pf
In the following, $C_i>0$ denote positive constants.
We put $\vbar_i:=(-\log|z_i|)\cdot\zbar_i\del/\del \zbar_i$
$(i=1,\ldots,\ell)$ 
and $\vbar_i:=\del/\del\zbar_i$
$(i=\ell+1,\ldots,n-1)$
for simplicity of description.
Because $\nabla_{\del_r}F=0$ and $[\del_r,\vbar_i]=0$,
we have $\nabla_{\del_r}\nabla_{\vbar_i}F=R(h)(\del_r,\vbar_i)F$.
Hence, we have the following:
\begin{multline*}
 \frac{\del}{\del r}\bigl|\nabla_{\vbar_i}F\bigr|^2_h(z,w)
=2\Re h\bigl(\nabla_{\del_r}\nabla_{\vbar_i}F,\nabla_{\vbar_i}F\bigr)
=2\Re h\Bigl(
 R(h)(\del_r,\vbar_i)F,\nabla_{\vbar_i}F
 \Bigr) \\
\leq
 C_1\cdot\bigl|f(z)\bigr|_h
 \cdot \bigl|\nabla_{\vbar_i}F\bigr|_h(z,w)
\end{multline*}
Note $\nabla_{\vbar_i}F(z,0)=0$.
Hence, we obtain
$\bigl|\nabla_{\vbar_i}F\bigr|_h(z,w)
 \leq
 C_2\cdot \bigl|f(z)\bigr|_h\cdot r$.
We also have the following:
\[
 \frac{\del} {\del r}
 \bigl|\nabla_{\del_{\theta}}F\bigr|_h^2(z,w)
=2\Re h\Bigl( R(h)(\del_r,\del_{\theta})F,
 \nabla_{\del_{\theta}}F \Bigr)
\leq
 C_3\cdot r\cdot\bigl|f(z)\bigr|_h\cdot
 \bigl|\nabla_{\del_{\theta}}F\bigr|_h(z,w)
\]
Note $\nabla_{\del_{\theta}}F(z,0)=0$
on the real blow up of $X-D$ along 
$X^{(0)}-D^{(0)}$.
Hence, we obtain
$\bigl|\nabla_{r^{-1}\del_{\theta}}F\bigr|_h
 \leq C_4\cdot r\cdot |f(z)|_h$.
Thus we are done.
\hfill\qed

\vspace{.1in}
We take a $C^{\infty}$-metric 
$h_{\nbigo(-X^{(0)})}$ of
the line bundle $\nbigo\bigl(-X^{(0)}\bigr)$ 
on $X$.
We naturally regard $\delbar F(z,w)$ as a section of
$E\otimes\nbigo\bigl(-X^{(0)}\bigr)
 \otimes\Omega^{0,1}_{X-D}$
on $X-D$.
Then, the norm of $\delbar F(z,w)$ is bounded
with respect to $h'=h\otimes h_{\nbigo(-X^{(0)})}$ 
and $g_{\poin}$ up to 
polynomial orders
in $-\log|z_i|$ $(i=1,\ldots,\ell)$.

Let $\chi$ be a non-negative $C^{\infty}$-function
on $\Delta_w$
such that $\chi(w)=1$ for $|w|\leq 1/2$ and 
$\chi(w)=0$ for $|w|\geq 2/3$.
We obtain the following $C^1$-section of
$E$ on $X-D$:
\[
 \rho:=\chi\cdot F(z,w)
\]
Then, $\rho$ satisfies the following property
in Proposition 8.40 in \cite{mochi2}.
\begin{itemize}
\item
 $\rho$ is bounded
 up to polynomial order 
 in $-\log|z_i|$ $(i=1,\ldots,\ell)$
 with respect to $h$.
\item
 We regard $\delbar\rho$
 as a section of 
 $E\otimes\nbigo_X(-X^{(0)})
 \otimes\Omega^{0,1}$
 on $X-D$.
 Then, it is bounded
 with respect to 
 $h\otimes h_{\nbigo(-X^{(0)})}$
 and $g_{\poin}$
 up to polynomial order in
 $-\log|z_i|$ $(i=1,\ldots,\ell)$.
\end{itemize}

\subsection{Proof of Lemma \ref{lem;08.10.3.10}}

By using $\rho$ as above,
Lemma \ref{lem;08.10.3.5},
and the argument in Sections
 8.4.3--8.4.4 of \cite{mochi2}
we can show the extension property of holomorphic 
sections as in Proposition 8.46 and Lemma 8.51
of \cite{mochi2}.
To explain how to use
Lemma \ref{lem;08.10.3.5}
and the above cocycle,
we give an outline.

Let $g_{\poin}$ be the Poincar\'e metric
of $X-D$.
Let $\pi_i$ denote the projection of $X-D$
to $\{z_i=0\}\times\Delta_w$.
Let $D_i^{\circ}$ be the image.
For any point $Q\in D_i^{\circ}$,
let $(E_Q,\delbar_Q,h_Q)$ denote
the restriction of
$(E,\delbar_E,h)$ to $\pi_i^{-1}(Q)\setminus D$.
By Lemma \ref{lem;08.10.3.5},
the set
$\Par(\prolong{E},i):=
 \Par(\prolong{E_Q})$ is independent
of the choice of $Q\in D_i^{\circ}$.

For $\veca\in\real^{\ell}$
and $N\in\real$,
we consider the following metrics:
\[
 h_{\veca,N}=h\cdot 
 \prod_{i=1}^{\ell}
 |z_i|^{2a_i}
 \bigl(-\log|z_i|^2\bigr)^{-N}
\cdot
 \prod_{i=\ell+1}^{n-1}
 (1-|z_i|^2)^{-N}
\cdot 
 (1-|w|^2)^{-N}
\]
\[
  \htilde_{\veca,N}:=
 h_{\veca,N}\otimes h_{\nbigo(-X^{(0)})}
\]
Let $\|\cdot \|_{\htilde_{\veca,N}}$
denote the $L^2$-norm with respect to
$\htilde_{\veca,N}$ and $g_{\poin}$.
As in Section \ref{subsection;04.12.23.50},
due to the argument
in Section 2.8.6 of \cite{mochi2},
there exists a negative number $N_0$ such that
the following holds for any $N<N_0$
and for any $C^{\infty}$-section $\eta$ of
$E\otimes\nbigo(-X^{(0)})\otimes\Omega^{0,1}$
with compact support:
\[
 \bigl\|\delbar_E \eta\bigr\|_{\htilde_{\veca,N}}^2
+\bigl\|\delbar_E^{\ast}\eta\bigr\|
 _{\htilde_{\veca,N}}^2
\geq
 \bigl\|\eta\bigr\|_{\htilde_{\veca,N}}^2
\]
According to Andreotti-Vesentini
(Propositions \ref{prop;11.28.6}
and \ref{prop;12.11.1}),
it implies the following:
\begin{itemize}
\item
Let $N<N_0$.
Let $\veca$ be any element of $\real^{\ell}$.
Let $\tau$ be a $\delbar_E$-closed section of 
$E\otimes\nbigo(-X^{(0)})\otimes\Omega^{0,1}$
which is $L^2$ with respect to
$\htilde_{\veca,N}$ and $g_{\poin}$.
Then, there exists a section
$\omega$ of $E\otimes\nbigo(-X^{(0)})$
such that 
(i) it is $L^2$ with respect to
 $\htilde_{\veca,N}$ and $g_{\poin}$,
(ii) $\delbar_E\omega=\tau$.
\end{itemize}

We take $N<N_0$.
We take a small number $\epsilon>0$  such that
$\openclosed{-1}{-1+\epsilon}
\cap \Par(\prolong{E},i)=\emptyset$
for any $i=1,\ldots,\ell$.
Let $\vecdelta=(1,\ldots,1)\in\real^{\ell}$.
We regard $\delbar_E\rho$
as an $L^2$-section of
$E\otimes\nbigo(-X^{(0)})
 \otimes\Omega^{0,1}$
with respect to 
$\htilde_{\epsilon\vecdelta,N}$.
Then, we can find a section $\omega$
of $E\otimes\nbigo(-X^{(0)})$
such that 
(i) it is $L^2$ with respect to
$\htilde_{\epsilon\vecdelta,N}$
and $g_{\poin}$,
(ii) $\delbar_E\omega=\delbar_E\rho$.

We obtain a section $F:=\rho-\omega$ of $E$
on $X-D$,
which is $L^2$ with respect to
$h_{\epsilon\vecdelta,N}$ and $g_{\poin}$.
By our construction, $\delbar_EF=0$.
By using the $L^2$-property of $\omega$
as a section of $E\otimes\nbigo(-X_0)$,
we obtain that $\omega_{|X^{(0)}-D^{(0)}}=0$.
Therefore, $F_{|X^{(0)}-D^{(0)}}=f$.

By using Corollary 2.51 of \cite{mochi2},
we obtain that the restrictions
$F_{|\pi_i^{-1}(Q)\setminus D}$
are $L^2$ with respect to
$h_{\epsilon\vecdelta,N}$
for any $Q\in D_i^{\circ}$.
By Lemma \ref{lem;08.10.17.11},
we obtain 
\[
 F_{|\pi_i^{-1}(Q)}
\in \prolongg{\epsilon}{(E_{|\pi_i^{-1}(Q)})}.
\]
Because $\openclosed{-1}{-1+\epsilon}
\cap \Par(\prolong{E},i)=\emptyset$
for any $i=1,\ldots,\ell$,
we obtain 
\[
 F_{|\pi_i^{-1}(Q)}
\in \prolong{(E_{|\pi_i^{-1}(Q)})}.
\]
Hence, $F$ is a section of $\prolong{E}$
due to Proposition \ref{prop;08.10.17.12}.
\hfill\qed

\section{Extension of holomorphic sections II}
\label{subsection;10.11.80}
\subsection{Statement}

We put $X:=\Delta^n$ and
$D:=\bigcup_{i=1}^{\ell} D_i$.
We also set
\[
\begin{array}{l}
 X^{(1)}:=\bigl\{(z_1,\ldots,z_n)\in X\,\big|\,z_1=z_2\bigr\},
\quad
 D^{(1)}:=X^{{(1)}}\cap D,
 \\
 \mbox{{}}
 \\
 X_0:=\bigl\{(z_1,\ldots,z_n)\in X\,\big|\,z_1=z_2=0\bigr\}.
\end{array}
\]

\begin{condition} 
\label{condition;10.11.45}
Let $(E,\delbar_E,\theta)$ be 
an acceptable bundle over $X-D$.
Let $\epsilon_1$ and $\epsilon_2$ be positive numbers
such that $\epsilon_1+\epsilon_2<1$.
Assume that $\Par\bigl(\prolong{E},i\bigr)$
is contained in $\openclosed{-\epsilon_i}{0}$
for $i=1,2$.
\hfill\qed
\end{condition}
We remark 
$\Par\bigl(E,i\bigr)\cap
 \openclosed{0}{1-\epsilon_i}=\emptyset$
under the condition.
The following proposition
is the counterpart of Proposition 8.46
of \cite{mochi2}.

\begin{prop} 
\label{prop;9.9.81}
Let $(E,\delbar_E,h)$ be an acceptable bundle
satisfying Condition {\rm\ref{condition;10.11.45}}.
Let $f$ be a holomorphic section of
$\prolong{(E_{|X^{(1)}-D^{(1)}})}$
on $X^{(1)}$.
Then there exists a neighbourhood $U$ of $X_0$
in $X$,
and there exists a holomorphic section
$\widetilde{f}\in \Gamma(U,\prolong{E})$,
such that 
$\widetilde{f}_{|X^{(1)}\cap U}=f_{|X^{(1)}\cap U}$.
\end{prop}

\subsection{Preliminary from
 geometry on the blow up}
\label{subsection;08.10.16.10}

\subsubsection{Metrics and 
the curvatures of $\nbigo_{\proj^1}(i)$}

We recall the contents in Section 4.7.3 in \cite{mochi}
with minor corrections.
Let $\proj^1$ denote 
the one dimensional complex projective space.
We use the homogeneous coordinate $[t_0:t_1]$.
The points $[0:1]$ and $[1:0]$ are denoted by $0$ and $\infty$
respectively.
We use the coordinates $t=t_0/t_1$ and $s=t_1/t_0$.
We have the line bundle $\nbigo_{\proj^1}(i)$ 
over $\proj^1$.
It is the gluing of
$\cnum^2=\bigl\{(t,\zeta_1)\bigr\}
=\nbigo_{\proj^1}(i)_{|\proj^1-\{\infty\}}$
and $\cnum^2=\bigl\{
 (s,\zeta_2) \bigr\}
=\nbigo_{\proj^1}(i)_{|\proj^1-\{0\}}$.
The relations are given by
$s=t^{-1}$ and $t^{-i}\cdot \zeta_1=\zeta_2$.

For $\xi=(t,\zeta_1)=(s,\zeta_2)\in\nbigo_{\proj^1}(i)$,
we define
\[
 h_i(\xi,\xi):=
 |\zeta_1|^2 \bigl(1+|t|^2\bigr)^{-i}
=|\zeta_2|^2 \bigl( 1+|s|^2\bigr)^{-i}.
\]
Then, $h_i$ is a smooth hermitian metric of
the line bundle $\nbigo_{\proj^1}(i)$.
For any $a,b\in\real$,
we have the possibly singular metrics $h_{i,(a,b)}$
of $\nbigo_{\proj^1}(i)$
given by
\[
 h_{i,(a,b)}(\xi,\xi):=
 h_i(\xi,\xi) \bigl(1+|t|^{-2}\bigr)^{-a}
     \bigl(1+|t|^2\bigr)^{-b}
=h_i(\xi,\xi) \bigl(1+|s|^{2}\bigr)^{-a}
     \bigl(1+|s|^{-2}\bigr)^{-b},
\]
for $\xi=(t,\zeta_1)=(s,\zeta_2)\in\nbigo_{\proj^1}(i)$.
Around $t=0$ (resp. $s=0$),
$|t|^{-2a} h_{i,(a,b)}$ 
(resp. $|s|^{-2b}h_{i,(a,b)}$)
is a $C^{\infty}$-metric.
The curvature $R(h_{i,(a,b)})$ is as follows:
\begin{equation}
 \label{eq;a11.9.6}
 R(h_{i,(a,b)})=
 (a+b+i)
 \frac{dt\cdot d\tbar}{\bigl(1+|t|^2\bigr)^2}
\end{equation}

\subsubsection{An open subset of
 the line bundle $\nbigo_{\proj^1}(-1)$
 with a complete Kahler metric}

We are mainly interested in the case $i=-1$.
We regard $\nbigo_{\proj^1}(-1)$ as a complex manifold.
We set 
\begin{equation}
 \label{eq;08.10.17.1}
Y:=\bigl\{\xi\in\nbigo_{\proj^1}(-1)\,
 \big|\,h_{-1,(0,0)}(\xi,\xi)<1\bigr\}.
\end{equation}
Let $\pi$ denote the projection of $Y$ onto $\proj^1$.
The image of the $0$-section $\proj^1\lrarr Y$
is also denoted by $\proj^1$.
We have the normal crossing divisor
$D'=\proj^1\cup \pi^{-1}(0)\cup \pi^{-1}(\infty)$
of $Y$.
The manifold $Y-D'$ is isomorphic to
$\bigl\{(t,x)\in \cnum^{\ast\,2}\,\big|\,
 |x|^2\bigl(1+|t|^2\bigr)<1\bigr\}$.

We have a complete Poincar\'e like
Kahler metric of $Y-D'$.
For example, we can construct it as follows.
As a contribution of the $0$-section $\proj^1$,
we put $\tau_1:=
 -\log\Bigl[\bigl(1+|t|^2\bigr)\cdot |x|^2\Bigr]$
and
\[
  g_1:=\frac{1}{\tau_1^2}
 \left(
 \frac{\tbar\cdot dt}{1+|t|^2}
+\frac{dx}{x}
 \right)
\cdot
 \left(
 \frac{t\cdot d\tbar}{1+|t|^2}
+\frac{d\xbar}{\xbar}
 \right)
+\frac{1}{\tau_1}
 \frac{dt\cdot d\tbar}{\bigl(1+|t|^2\bigr)^2}
\]
As a contribution of $\pi^{-1}(\infty)$,
we put $\tau_2:=\log\bigl(1+|t|^2\bigr)$ and
\[
 g_2:=\frac{1}{\tau_2}
 \Bigl(
 -1+\frac{|t|^2}{\tau_2}
 \Bigr)
 \cdot
 \frac{dt\cdot d\tbar}{\bigl(1+|t|^2\bigr)^2}.
\]
As the contribution of the divisor $\pi^{-1}(0)$,
we put $\tau_3:=\log(1+|t|^2)-\log|t|^2=\log(1+|s|^2)$,
where we use $s=t^{-1}$ and
\[
 g_3:=\frac{1}{\tau_3}
 \cdot
 \Bigl(
 -1+\frac{|s|^2}{\tau_3}
 \Bigr)
 \frac{ds\cdot d\sbar}{\bigl(1+|s|^2\bigr)^2}.
\]
Then, we set $g:=g_1+g_2+g_3$.
It is easy to check that
$g$ is a complete Kahler metric.
Note the following formulas:
\[
 \begin{array}{l}
 {\displaystyle
 \delbar\del\log \tau_1=
 \frac{1}{\tau_1^2}
 \Bigl(
  \frac{\tbar\cdot dt}{1+|t|^2}+\frac{dx}{x}
 \Bigr)
 \wedge
 \Bigl(
  \frac{t\cdot d\tbar}{1+|t|^2}+\frac{d\xbar}{\xbar}
 \Bigr)
 +\frac{1}{\tau_1}
 \frac{dt\wedge d\tbar}{\bigl(1+|t|^2\bigr)^2}=
 :\omega_1}\\
\mbox{{}}\\
{\displaystyle
 \delbar\del\log \tau_2=
 \frac{1}{\tau_2}
 \Bigl(
 -1+ \frac{|t|^2}{\tau_2}
 \Bigr)
\cdot
 \frac{dt\wedge d\tbar}
 {\bigl(1+|t|^2\bigr)^2}=:\omega_2}\\
{\displaystyle
 \delbar\del\log \tau_3=
 \frac{1}{\tau_3}
 \Bigl(
 -1+ \frac{|s|^2}{\tau_3}
 \Bigr)
\cdot
 \frac{ds\wedge d\sbar}
 {\bigl(1+|s|^2\bigr)^2}=:\omega_3}
 \end{array}
\]
We put $\omega:=\omega_1+\omega_2+\omega_3$.
Then, $\sqrt{-1}\omega$
 is the Kahler form corresponding to $g$.
We set
\[
 H_0:=
 \frac{1}{\tau_1}
+\frac{1}{\tau_2}
 \Bigl(
   -1+\frac{|t|^2}{\tau_2}
 \Bigr)
+\frac{1}{\tau_3}
 \Bigl(
  -1+\frac{|s|^2}{\tau_3}
 \Bigr)>0.
\]
Then, the following holds:
\[
  \omega^2=
\det(g)
dt\wedge d\tbar
 \wedge dx\wedge d\xbar
= \Bigl(
 \frac{2}{\tau_1^2
 |x|^2\bigl(1+|t|^2\bigr)^2}
 \times H_0
 \Bigr)
 \cdot dt\wedge d\tbar
 \wedge 
 dx\wedge d\xbar
\]
We put 
$H_1:=H_0\cdot 
 \bigl(1+|t|^2\bigr)^{-1}
 \cdot \bigl(1+|s|^2\bigr)^{-1}$.
Recall $\Ric(g)=\delbar\del(\det(g))$.
\begin{lem}\mbox{{}}
\begin{itemize}
\item
Let $C$ be a number such that $0<C<1$.
We have
$H_{1}\sim \bigl(\bigl|\log|t|\bigr|+1\bigr)^{-2}$
on the domain
$\bigl\{\xi\in\nbigo_{\proj^1}(-1)\,\big|\,
 h_{-1,(0,0)}(\xi,\xi) \leq C \bigr\}$.
\item
We have the equality
$Ric(g)-\delbar\del\log(H_1)=-\delbar\del\log\tau_1^2$.
\hfill\qed
\end{itemize}
\end{lem}

\subsubsection{Inequality and vanishing}

We put $X:=\Delta^n$, $D_i:=\{z_i=0\}$ 
and $D:=\bigcup_{i=1}^{\ell}D_i$.
We also put $\Delta_z^2:=\{(z_1,z_2)\,|\,|z_i|<1\}$
and $D_i':=\{z_i=0\}\subset \Delta_z^2$ $(i=1,2)$.
Let $\varphi:\blowup{\Delta_z^2}\lrarr \Delta_z^2$
denote the blow up of $\Delta_z^2$ 
at the origin $O=(0,0)$.
We have the exceptional divisor $\varphi^{-1}(O)$
and the proper transforms 
$\blowup{D'}_i$ of $D'_i$ $(i=1,2)$.

We put $\blowup{X}:=
\blowup{\Delta_z^2}\times\Delta_w^{n-2}$.
Then we have the composite 
$\psi$ of the natural morphisms:
\[
\begin{CD}
 \blowup{X}@>{\varphi\times id}>>
 \Delta^2_z\times\Delta^{n-2}_w
 @>>>\Delta_z^{n}
\end{CD}
\]
Here the latter morphism is the natural isomorphism
given by $w_{i}=z_{i+2}$ $(i=1,\ldots,n-2)$.
We put $\blowup{D}:=\psi^{-1}(D)$, 
which is the same as the following:
\[
 \left[
 \bigl(
 \varphi^{-1}(O)\cup \blowup{D_1'}\cup\blowup{D_2'}
 \bigr)\times \Delta_w^{n-2}
\right]
\cup
\left[
 \widetilde{\Delta}^2_z\times
 \Bigl(
 \bigcup_{i=1}^{\ell-2}\{w_i=0\}
 \Bigr)\right]
\]
The restriction of $\psi$ to $\blowup{X}-\blowup{D}$
gives an isomorphism $\blowup{X}-\blowup{D}\simeq X-D$.

We can take a holomorphic embedding $\iota$
of $Y$ in (\ref{eq;08.10.17.1})
to $\widetilde{\Delta}^2$
satisfying the following:
\begin{itemize}
\item The image of the $0$-section $\proj^1$ is the exceptional divisor
 $\varphi^{-1}(O)$.
\item
We have $\iota^{-1}(\Dtilde_1')=\pi^{-1}(\infty)$
and $\iota^{-1}(\Dtilde_2')=\pi^{-1}(0)$.
\end{itemize}

We put $\overline{X}:=Y\times\Delta_w^{n-2}$.
Then we have the naturally induced morphism
$\overline{X}\lrarr\Xtilde$,
which is also denoted by $\iota$.
We have the point $[1:1]\in \proj^1$,
and we put 
\[
 \overline{D}:=\iota^{-1}(\blowup{D}),\quad
 \overline{X}^{(1)}:=\pi^{-1}([1:1])\times
 \Delta_w^{n-2},\quad
 \overline{D}^{(1)}:=\overline{X}^{(1)}\cap \overline{D}.
\]
The composite $\psi\circ\iota$ is denoted by $\psi_1$.
The metric $g_{\overline{X}-\overline{D}}$
is induced from 
the metric $g$ of $Y-D'$ 
and the Poincar\'{e} metric of
$\Delta_w^{n-2}\setminus 
 \bigcup_{i=1}^{\ell-2}\{w_i=0\}$.

Let $(E,\delbar_E)$ be a holomorphic bundle
with a hermitian metric $h$ over $X-D$.
We assume that $(E,\delbar_E,h)$ is acceptable.
We denote the curvature of $\psi_1^{-1}(E,\delbar_E,h)$
by $\psi_1^{-1}R(h)$.

Let $\epsilon_i$ $(i=1,2)$ be 
as in Condition \ref{condition;10.11.45}.
We can pick positive numbers $\epsilon$, $a$ and $b$
satisfying the following:
\begin{equation}
 \label{eq;08.10.16.11}
 a+b=1,
\quad
 0<a+2\epsilon<1-\epsilon_1,
\quad
 0<b+2\epsilon<1-\epsilon_2
\end{equation}
We set $\vecdelta:=(1,\ldots,1)\in \real^{\ell}$.
Let $h_{\epsilon\vecdelta,N}$
be as in (\ref{eq;04.12.27.100}).
Let $\htilde_{N,\epsilon,a,b}$ denote 
the metric of the bundle
$\psi_1^{-1}(E)(-\overline{X}^{(1)}):=
\psi_1^{-1}(E)\otimes\nbigo_{\proj^1}(-1)$
on $\overline{X}-\overline{D}$
 given as follows:
\begin{equation}
\label{eq;10.11.75}
 \htilde_{N,\epsilon,a,b}:=
 \psi_1^{-1}h_{\epsilon\vecdelta,N}\cdot
 h_{-1,a,b}\cdot
 H_1^{-1}\cdot
  \tau_1^{2+\epsilon}
 \bigl(
 \tau_2\cdot\tau_3
 \bigr)^{\epsilon}
\end{equation}
For simplicity, we use the symbol $\tilde{h}$
instead of $\tilde{h}_{N,\epsilon,a,b}$.

\begin{lem} 
\label{lem;5.22.30}
If $N$ is sufficiently negative,
the following inequality holds for any
$\eta\in 
 A^{0,1}_c(
 \psi_1^{-1}E(-\overline{X}^{(1)}))$:
\[
 \doublelangle
 R(\tilde{h})+\Ric(g),\eta
 \doublerangle_{\tilde{h}}
\geq \epsilon||\eta||_{\tilde{h}}^2
\]
(See {\rm(\ref{eq;10.11.60})}
and {\rm (\ref{eq;08.10.19.30})}
for the definition of
$\doublelangle\cdot,
 \cdot\doublerangle_{\tilde{h}}$.)
\end{lem}
\pf
Recall that $\overline{X}-\overline{D}$ is isomorphic to
the product of $Y-D'$ and 
$\Delta_w^{n-2}-\bigcup_{i=1}^{\ell-2}\{w_i=0\}$,
which is compatible with the metrics.
We can apply the inequality (\ref{eq;11.9.5})
to the
$\Bigl(\Delta_w^{n-2}
 -\bigcup_{i=1}^{\ell-2}\{w_i=0\}\Bigr)$-direction.
Hence we have only to consider 
$(Y-D')$-direction,
i.e.,
we may and will assume $n=2$.
We have the following equality:
\begin{multline}
 R(\tilde{h})+\Ric(g)
=R(\psi_1^{-1}h_{\epsilon\vecdelta,N})+R(h_{-1,a,b})
-\delbar\del\log H_1 \\
+(2+\epsilon)\cdot\delbar\del\log \tau_1
+\epsilon \cdot\delbar\del(\log \tau_2+\log\tau_3)
+\Ric(g)\\
= R(\psi_1^{-1}h_{\epsilon\vecdelta,N})
+ \epsilon\cdot(\omega_1+\omega_2+\omega_3)
\end{multline}
Here we have used 
$R(h_{-1,a,b})=0$ due to (\ref{eq;a11.9.6})
and our choice of $a$ and $b$.
By taking sufficiently negative $N$,
we can assume the following inequality
for any $\eta\in A_c^{0,1}(E)$ on $X-D$:
\begin{equation} 
 \label{eq;04.12.23.51}
 \doublelangle
 R(h_{\epsilon\vecdelta,N}),\eta
 \doublerangle_{\epsilon\vecdelta,N}\geq\,0.
\end{equation}
Then, by a fiber-wise linear algebraic argument
using the expression (\ref{eq;08.10.19.30}),
it is easy to see that the following inequality holds
for any $\eta\in A_c^{0,1}(\psi_1^{-1}(E))$:
\[
 \doublelangle
 \psi^{-1}R(h_{\epsilon\vecdelta,N}),\eta
 \doublerangle_{\tilde{h}}\geq\,0.
\]
We also obtain
$\epsilon
 \big\langle\big\langle
 \omega_1+\omega_2+\omega_3,\eta
 \big\rangle\big\rangle_{\tilde{h}}
\geq \epsilon\cdot||\eta||_{\tilde{h}}$,
by using (\ref{eq;08.10.19.30}).
Thus we are done.
\hfill\qed

\begin{cor} 
\label{cor;10.11.56}
If $N$ is sufficiently negative,
the first cohomology group
$H^1\bigl(
 A^{0,\cdot}_{\tilde{h}}\bigl(\psi_1^{-1}E(-\overline{X}^{(1)})
 \bigr)\bigr)$ vanishes.
\end{cor}
\pf
It immediately follows from Lemma \ref{lem;5.22.30},
Proposition \ref{prop;11.28.6}
and Proposition \ref{prop;12.11.1}.
\hfill\qed

\vspace{.1in}

The contribution of
$h_{-1,a,b}\cdot H_1^{-1}\cdot\tau_1^{2+\epsilon}
 \cdot \bigl(\tau_2\cdot\tau_3\bigr)^{\epsilon}$
to the metric $\htilde$
is equivalent to the following, 
on a curve transversal with $\{t=0\}$:
\[
 |t|^{2a}\cdot\bigl(-\log |t|\bigr)^{2}\cdot |t|^{2\epsilon}
 \bigl(-\log|t|\bigr)^{2\epsilon}
=|t|^{2(a+\epsilon)}\cdot \bigl(-\log|t|\bigr)^{2+2\epsilon}.
\]
We have a similar estimate 
on a curve transversal with $\{s=0\}$.
The contribution is equivalent to 
$(-\log |x|)^{2+\epsilon}$ 
on a curve transversal with $\{x=0\}$.

\subsection{Proof of Proposition \ref{prop;9.9.81}}

We impose the additional condition
$\psi_1\bigl(\overline{X}^{(1)}\bigr)=X^{(1)}$.
We remark that $\psi_1(\overline{X})$ 
gives a neighbourhood of $X_0$.
We take numbers $\epsilon$, $a$ and $b$
as in (\ref{eq;08.10.16.11}).
Note
\[
 \Par\bigl(E,1\bigr)
\cap\openclosed{0}{a+2\epsilon}
\subset
 \Par\bigl(E,1\bigr)
\cap \openclosed{0}{1-\epsilon_1}=\emptyset.
\]
Similarly,
we have $\Par\bigl(E,2\bigr)
\cap\openclosed{0}{b+2\epsilon}=\emptyset$.
Moreover, we may assume that
$\Par\bigl(E,i\bigr)\cap 
 \openclosed{0}{\epsilon}=\emptyset$
for $i=3,\ldots,\ell$,
if $\epsilon$ is sufficiently small.

Let us take a sufficiently negative number $N$
such that {\rm(\ref{eq;11.9.5})} holds.
We take the metric 
$\widetilde{h}=\htilde_{N,\epsilon,a,b}$ of
the bundle $\psi_1^{-1}E(-\overline{X}^{(1)})$
as in (\ref{eq;10.11.75}).
We also put
$\widetilde{h}_0:=
 \psi_1^{-1}h_{\veczero,N}\cdot
 h_{-1,a,b}\cdot
 H_1^{-1}\cdot
  \tau_1^{2+\epsilon}
 \bigl(
 \tau_2\cdot\tau_3
 \bigr)^{\epsilon}$.
We remark that we can use Corollary \ref{cor;10.11.56}
in this setting.
We take the metrics $\widehat{h}$ and $\widehat{h}_0$
of $\psi_1^{-1}E$:
\[
\begin{array}{l}
 \widehat{h}:=
 \psi_1^{-1}h_{\epsilon\vecdelta,N}\cdot
 h_{0,a,b}\cdot
 H_1^{-1}\cdot
  \tau_1^{2+\epsilon}
 \bigl(
 \tau_2\cdot\tau_3
 \bigr)^{\epsilon}\\
 \mbox{{}}
 \widehat{h}_0:=
 \psi_1^{-1}h_{\veczero,N}\cdot
 h_{0,a,b}\cdot
 H_1^{-1}\cdot
  \tau_1^{2+\epsilon}
 \bigl(
 \tau_2\cdot\tau_3
 \bigr)^{\epsilon}\\
\end{array}
\]
Take an embedding 
$\kappa:\Delta_{\zeta}\lrarr \proj^1-\{0,\infty\}$
such that $\kappa(0)=[1:1]\in\proj^1$.
We take a holomorphic function $\eta$
on $\pi^{-1}(\kappa(\Delta_{\zeta}))$
such that $\eta^{-1}(0)$ is the intersection of
the exceptional divisor $\varphi^{-1}(O)$
and $\pi^{-1}(\kappa(\Delta_{\zeta}))$.

The section $f$ induces a holomorphic section of
$\bigl(\psi_1^{-1}E\bigr)_{|\overline{X}^{(1)}}$
over $\overline{X}^{(1)}\setminus \overline{D}^{(1)}$,
which is denoted by $f$.
By using the argument in Subsection
\ref{subsection;08.10.18.50},
we can take a $C^{1}$-section $\rho$ 
of $\psi_1^{-1}E$ over $\overline{X}-\overline{D}$
satisfying the following:
\begin{itemize}
\item
 The support of $\rho$ is contained in
 $\pi^{-1}(\kappa(\Delta_{\zeta}))\times\Delta_w^{n-2}$.
\item
 $\rho$ is bounded as a section of $\psi_1^{-1}E$
 up to polynomial order of $-\log|z_i|$ $(i=3,\ldots,\ell)$
 and $-\log|\eta|$,
 with respect to $\widehat{h}_0$.
 In particular, $\rho$ is an $L^2$-section of
 $\psi_1^{-1}E$
 with respect to $\widehat{h}$.
\item
 $\delbar \rho$ is bounded as a section of
 $\psi_1^{-1}E(-\overline{X}^{(1)})
 \otimes\Omega^{1,0}_{\overline{X}-\overline{D}}$
 up to polynomial order of $-\log|z_i|$ $(i=3,\ldots,\ell)$
 and $-\log|\eta|$,
 with respect to the metric $\widetilde{h}_0$ of
 $\psi_1^{-1}E(-\overline{X}^{(1)})\otimes\Omega^{1,0}$
 and the metric $g_{\overline{X}-\overline{D}}$
 of $\Omega^{0,1}_{\overline{X}-\overline{D}}$.
 The restriction of $\delbar\rho$ to 
 $\overline{X}-(\overline{D}\cup \overline{X}^{(1)})$
 are $C^{\infty}$.
 In particular,
 $\delbar\rho$ is contained in
 $A^{0,1}_{\widetilde{h}}
 \bigl(\psi_1^{-1}E(-\overline{X}^{(1)})\bigr)$.
\item
 We have $\rho_{|\overline{X}^{(1)}-\overline{D}^{(1)}}=f$.
\end{itemize}

Due to Corollary \ref{cor;10.11.56},
we can pick $G\in
 A^{0,0}_{\widetilde{h}}
 \bigl(\psi_1^{-1}E(-\overline{X}^{(1)})\bigr)$
such that $\delbar G=\delbar\rho$.
We put $\ftilde:=\rho-G$,
which satisfies 
$\delbar\ftilde=0$ and
$\ftilde\in A_{\widehat{h}}^{0,0}(\psi_1^{-1}E)$.
By using the $L^2$-property of $G$
as a section of 
$\psi_1^{-1}E\otimes\nbigo(-\Xbar^{(1)})$,
we obtain that $G_{|\Xbar^{(1)}-\Dbar^{(1)}}=0$.
Therefore, $\ftilde_{|\Xbar^{(1)}-\Dbar^{(1)}}=f$.

\vspace{.1in}

We have an open subset $V$ of $X$
such that
$\psi_1\bigl(\overline{X}-\overline{D}\bigr)
=V\setminus D$.
By the identification,
we can regard $\widetilde{f}$
as a holomorphic section of $E_{|V\setminus D}$.
We would like to show that
$\widetilde{f}$ gives a section of
$\prolong{E}$ on $V$.
We put $D_i^{\circ}:=
 D_i\setminus \bigcup_{j\leq \ell,j\neq i}D_j$.
\begin{lem}
Let $P$ be any point of $V\cap D_1^{\circ}$.
We have 
$\bigl|
 \widetilde{f}_{|\pi_1^{-1}(P)\cap V}\bigr|
=O\bigl(|z_1|^{-\kappa}\bigr)$
for any $\kappa>0$.
\end{lem}
\pf
We put $C_P:=\pi_1^{-1}(P)\cap V$.
We may assume that the closure of $C_P$ in $\cnum^n$
is contained in $X$, by shrinking $V$.
Let us take a small neighbourhood $U$ of $P$
in $D_1^{\circ}$.
We may assume $C_P\times U\subset V$.
The metrics $\widehat{h}_{|C_P\times U}$
and $h\cdot |z_1|^{2(a+\epsilon)}$
are mutually bounded up to polynomial order
of $-\log|z_1|$.
Therefore $\widetilde{f}_{|C_P\times U}$
is $L^2$, with respect to the metric
$h\cdot |z_1|^{2(a+\epsilon)}\cdot
 \bigl(-\log|z_1|\bigr)^{-M}$
for $M>0$ and the metric $g_{X-D\,|\,C_P\times U}$.
Due to Corollary \ref{cor;04.12.23.30},
the restriction
$\widetilde{f}_{|C_P}$ is $L^2$
with respect to
$h\cdot |z_1|^{2(a+\epsilon)}\cdot 
\bigl(-\log|z_1|\bigr)^{-M}$
and the metric $g_{X-D\,|\,C_P}$.
We also remark Lemma \ref{lem;08.10.17.11}.
Because
$\Par\bigl(\prolong{E},1\bigr)
\cap \openclosed{0}{a+2\epsilon}=\emptyset$,
we obtain the desired estimate for
$\widetilde{f}_{|C_P}$.
\hfill\qed

\vspace{.1in}

Similarly, we can show the following lemma.
\begin{lem}
 \label{lem;04.12.25.200}
We have 
$\bigl|
 \widetilde{f}_{|\pi_i^{-1}(P)\cap V}
 \bigr|_h=O\bigl(|z_i|^{-\kappa}\bigr)$
for any $\kappa>0$,
any $P\in D_i^{\circ}\cap V$,
and $i=1,2,\ldots,\ell$.
\hfill\qed
\end{lem}

\begin{lem}
$\widetilde{f}$ is a section of
$\prolong{E}$ over $V$.
\end{lem}
\pf
It follows from Lemma \ref{lem;04.12.25.200}
and Proposition \ref{prop;08.10.17.12}.
\hfill\qed

\vspace{.1in}
As a result,
we obtain a holomorphic section $\widetilde{f}$
of $\prolong{E}$ over $V$
such that $\widetilde{f}_{|X^{(1)}}=f$.
Thus the proof of Proposition \ref{prop;9.9.81}
is accomplished.
\hfill\qed

\section{Proof of Theorems \ref{thm;07.3.16.1}
 and \ref{thm;07.10.9.1}}
\label{subsection;08.10.18.21}
We use an induction on the dimension of $X$.
As we have already remarked,
the one dimensional case of
Theorem \ref{thm;07.3.16.1}
was proved by Simpson in \cite{s1} and \cite{s2}.
We remark that the claim for
$\det(E,\delbar_E,h)$ is easy
because it is assumed to be flat.

\subsection{Preliminary prolongation}
\label{subsection;08.10.19.100}

Let $(E,\delbar_E,h)$ be an acceptable bundle
on $X-D$.
We impose the following assumptions
in this section.
\begin{itemize}
\item
Let $\epsilon_i$ $(i=1,2)$ be small positive numbers
such that $\rank(E)\cdot(\epsilon_1+\epsilon_2)<1/2$.
Then,
$\Par(\prolong{E},i)$ is contained in
$\openopen{-\epsilon_i}{0}$
for $i=1,2$.
\item
Pick a point
$(\lefttop{0}z_3,\ldots,\lefttop{0}z_{n})
 \in(\Delta^{\ast})^{n-2}$,
and we set
$C_0:=\{(z,z,\lefttop{0}z_3,\ldots,\lefttop{0}z_n)\in X\}$.
Then,
$\Par(E_{|C_0})$ is contained in
 $\openopen{-\epsilon_1-\epsilon_2}{0}$.
Note that this condition is independent of the choice of
$(\lefttop{0}z_3,\ldots,\lefttop{0}z_{n})$.
\end{itemize}

We set $\btilde_i:=
 \sum_{b\in\Par(\prolong{E},i)} 
 b\cdot \rank\lefttop{i}\Gr^F_b(\prolong{E})$,
and 
$\vecbtilde:=(\btilde_i)$.

\begin{lem}
\label{lem;08.10.16.20}
$\prolongg{\vecbtilde}{\det(E)}_{|C_0}$
is isomorphic to
$\det\bigl(\prolong{(E_{|C_0})}\bigr)$.
\end{lem}
\pf
Take any $P_i\in D_i^{\circ}$ $(i=1,2)$.
We have the following:
\[
 \sum_{b\in\Par(\prolong{E_{|C_0}})}
 b\cdot\dim\Gr^F_b(\prolong{E}_{|C_0})
-\sum_{i=1,2}
 \sum_{b\in\Par(\prolong{E},i)}
 b\cdot \dim \lefttop{i}\Gr^F_b
 (\prolong{E}_{|\pi_i^{-1}(P_i)})
\in\seisuu
\]
By the assumption,
it has to be $0$.
Then, the claim follows.
\hfill\qed

\begin{lem} 
\label{lem;9.9.95}
Under the assumption,
the $\nbigo_X$-sheaf $\prolong{E}$ is locally free.
\end{lem}
\pf
We use the notation in Subsection \ref{subsection;10.11.80}.
Note that we are also using the induction on $\dim X$
to show the theorems.
Due to the hypothesis of the induction on $\dim X$,
we have the local freeness of
$\prolong{(E_{|X^{(1)}-D^{(1)}})}$.
Pick a frame $\vecv=(v_1,\ldots,v_r)$ of 
$\prolong{(E_{|X^{(1)}-D^{(1)}})}$
over $X^{(1)}$.
For each $v_i$,
we pick a section $\vtilde_i$ of 
$\prolong{E}$ over 
a neighbourhood $U$ of $X_0$
such that $\vtilde_{i\,|\,U\cap X^{(1)}}=v_{i\,|\,U\cap X^{(1)}}$,
by using Proposition \ref{prop;9.9.81}.
We may assume that $v_i$ are defined on $X$
by shrinking $X$.

Let us show that 
$\vecvtilde$ gives a frame of $\prolong{E}$
around $X_0$.
We set
$\Omega(\vecvtilde):=
 \vtilde_1\wedge\cdots\wedge \vtilde_r$.
The restriction $\vecvtilde_{|\pi_i^{-1}(P)}$
gives a tuple of holomorphic sections of
$\prolong{(E_{|\pi_i^{-1}(P)})}$
for any $P\in D_i^{\circ}$.
Hence 
$\Omega(\vecvtilde)_{|\pi_i^{-1}(P)}$ 
is a holomorphic section
of $\det\bigl(\prolong{(E_{|\pi_i^{-1}(P)})} \bigr)=
 \prolongg{\btilde_i}{\bigl(\det(E)_{|\pi_i^{-1}(P)}\bigr)}$.
It implies $\Omega(\vecvtilde)$ is a holomorphic section
of $\prolongg{\vecbtilde}{\det(E)}$.

We have the natural isomorphism
$(\prolongg{\vecbtilde}{\det(E)})_{|X^{(1)}}
\simeq
 \det\bigl(\prolong{(E_{|X^{(1)}})}\bigr)$
due to Lemma \ref{lem;08.10.16.20}.
Since $\Omega(\vecvtilde)_{|X^{(1)}}$
gives a frame of
$ \det\bigl(\prolong{(E_{|X^{(1)}})}\bigr)$,
we obtain $\Omega(\vecvtilde)_{|O}\neq 0$.
By shrinking $X$,
we may assume that 
$\Omega(\vecvtilde)$ gives a frame of
$\prolongg{\vecb}{\det E}$.

Let $f$ be any section of $\prolong{E}$.
We have an expression 
$f=\sum f_i\cdot \vtilde_i$,
where $f_i$ are holomorphic on $X-D$.
It is easy to observe that
$f\wedge \vtilde_2\wedge
 \cdots
 \wedge \vtilde_r
=f_1\cdot\Omega(\vecvtilde)$
is a section of
$\prolongg{\vecbtilde}{\det(E)}$.
Since $\Omega(\vecvtilde)$ is a frame,
$f_1$ is holomorphic on $X$.
Similarly, we obtain that
$f_i$ $(i=2,\ldots,r)$ are also holomorphic on $X$.
Hence, 
$\vecvtilde$ is a frame of 
$\prolong{E}$ on $X$.
\hfill\qed

\subsection{Proof of Theorem \ref{thm;07.3.16.1}}

\subsubsection{Step 1}

For a real number $a$,
let $\kappa_{1/2}(a)\in\openclosed{-1/2}{1/2}$ 
and $\nu_{1/2}(a)\in\seisuu$ be determined by
$a=\kappa_{1/2}(a)+\nu_{1/2}(a)$.
We set $\eta:=(10\rank E)^{-1}$.
We can take a positive integer $c$
such that the following  conditions are satisfied:
\begin{enumerate}
\item
 The maps $\Par(\prolong{E},i)\lrarr \seisuu$
 given by 
 $a\longmapsto \nu_{1/2}(c\cdot a)$
 are injective for $i=1,\ldots,\ell$.
\item
$ \{\kappa_{1/2}(c\cdot a_i)\,|\,
  a_i\in\Par(E,i) \}\subset
  \openopen{-\eta}{\eta}$
 for $i=1,2,\ldots,\ell$.
\item
$ \{\kappa_{1/2}(c\cdot a)\,|\,
  a\in\Par(E_{|C_0}) \}\subset
  \openopen{-\eta}{\eta}$,
where $C_0$ is a curve as in 
 Section \ref{subsection;08.10.19.100}.
\end{enumerate}
We set $\veceta:=(\eta,\ldots,\eta)\in\real^{\ell}$.
For a positive integer $c$,
let $\psi_c:X-D\lrarr X-D$
be given by
\[
 \psi_c(z_1,\ldots,z_n)
=(z_1^c,\ldots,z_{\ell}^c,z_{\ell+1},\ldots,z_n).
\]
We put
$(E_1,\delbar_{E_1},h_1):=
 \psi_{c}^{-1}(E,\delbar_E,h)
 \otimes L(-\veceta)$.
We have the natural isomorphism
$\prolongg{\veceta}{\bigl( \psi_{c}^{-1}E\bigr)}
\simeq
 \prolong{E_1}$.
We also have the following,
due to the result in the one dimensional case
(Section \ref{subsection;08.10.18.40}):
\[
  \Par(\prolong{E_1},i)
=\bigl\{
 -\eta+\kappa_{1/2}(c\cdot a_i)\,|\,
  a_i\in\Par(E,i)
 \bigr\}
\subset\openopen{-(5\rank E)^{-1}}{0}
\]
\[
 \Par\bigl(
 \prolong{E_{1|C_0}}
 \bigr)
=\bigl\{
 -2\eta+\kappa_{1/2}(c\cdot a)
 \,\big|\,
 a\in\Par(E_{|C_0})
 \bigr\}
\subset
 \openopen{-2\cdot(5\rank E)^{-1}}{0}
\]
By Lemma \ref{lem;9.9.95},
$\prolong{E_1}$ is locally free.
Hence,
the sheaf 
$\prolongg{\veceta}{\bigl(\psi_{c}^{-1}E\bigr)}$
is a locally free $\nbigo_X$-module.

\vspace{.1in}

Let $\mu_c:=\bigl\{z\in\cnum\,\big|\,
 z^c=1 \bigr\}$.
We fix a generator $\omega$.
We have the natural $\mu_c^{\ell}$-action on $X-D$
given by
\[
 (\omega_1,\ldots,\omega_{\ell})\cdot
 (z_1,\ldots,z_{n})
=\bigl(
 \omega_1\cdot z_1,\ldots,
 \omega_{\ell}\cdot z_{\ell},\,
 z_{\ell+1},\ldots,z_n
 \bigr)
\]
It is lifted to the action
on $\prolongg{\veceta}{\bigl(\psi_{c}^{-1}E\bigr)}$.
The $i$-th component of $\mu_c^{\ell}$
is denoted by $\mu_c^{(i)}$,
which acts on 
$\prolongg{\veceta}
 {\bigl(\psi_{c}^{-1}E\bigr)}_{|D_i}$.
We have the decomposition:
\[
 \prolongg{\veceta}
 {(\psi_{c}^{-1}E)}_{|\,D_i}
=\bigoplus_{0\leq p\leq c-1}
 \lefttop{i}V_p.
\]
Here the generator $\omega$ of $\mu_{c}^{(i)}$
acts as the multiplication of $\omega^p$ on 
$\lefttop{i}V_p$.
We have the following morphism 
given as in Section \ref{subsection;08.10.18.40}:
\[
 \varphi_i:
 \bigl\{p\,\big|\,0\leq p\leq c-1,\,\,
 \lefttop{i}V_p\neq 0\bigr\}
\lrarr
 \Par\bigl(
 \prolongg{\veceta}{(\psi_c^{-1}E)},i
 \bigr)
\]
We consider the filtration $\lefttop{i}F'$
of $\prolongg{\veceta}{\bigl(
 \psi_{\vecc}^{-1}E\bigr)}_{|\,D_i}$
in the category of vector bundles on $D_i$,
given as follows for any $\eta-1<b<\eta$:
\begin{equation}
\label{eq;08.10.19.101}
 \lefttop{i}F'_{b}:=
 \bigoplus_{\substack{\varphi_i(p)\leq b}}
 \lefttop{i}V_p
\end{equation}
By the splittings (\ref{eq;08.10.19.101}),
it is easy to see that the tuple of filtrations
$(\lefttop{i}F'\,|\,i=1,\ldots,\ell)$
are compatible.

We set
$\vecdelta_i:=
 (\overbrace{0,\ldots,0}^{i-1},1,0,\ldots,0)
 \in\real^{\ell}$.
For any $-1<b<0$,
we consider the subsheaf
$\prolongg{\veceta+b\vecdelta_i}
 {\bigl(\psi_{\vecc}^{-1}E\bigr)}'$
of $\prolongg{\veceta}
 {\psi_{\vecc}^{-1}\bigl(E\bigr)}$
given as follows:
\[
 \prolongg{\veceta+b\vecdelta_i}{\bigl(
 \psi_{\vecc}^{-1}E\bigr)}'
=\Ker\Bigl(\pi:
 \prolongg{\veceta}{\bigl(\psi_{\vecc}^{-1}E\bigr)}
\lrarr
 \frac{\prolongg{\veceta}{\bigl(
 \psi_{\vecc}^{-1}E\bigr)}_{|D_i}}
  {\lefttop{i}F'_{\eta+b}}
 \Bigr)
\]
Here $\pi$ denotes the naturally defined morphism
of $\nbigo_X$-modules.

\begin{lem}
 \label{lem;9.10.13}
We have the following, for any $-1<b<0$:
\[
 \prolongg{\veceta+b\vecdelta_i}
 {\bigl(\psi_{\vecc}^{-1}E\bigr)}'
=\prolongg{\veceta+b\vecdelta_i}
 {\bigl(\psi_{\vecc}^{-1}
 E\bigr)}
\]
In other words,
the parabolic filtration $\lefttop{i}F$
is equal to $ \lefttop{i}F'$.
\end{lem}
\pf
Let $f$ be a holomorphic section
of $\prolongg{\veceta+b\vecdelta_i}
 {\bigl(\psi_{\vecc}^{-1}E\bigr)}$.
It can be also regarded as a section of
$\prolongg{\veceta}
 {\bigl(\psi_{\vecc}^{-1}E\bigr)}$.
Let $P$ be a point of $D_i^{\circ}$.
We have the element $f(P)$ of
$\prolongg{\veceta}{\bigl(
 \psi_{\vecc}^{-1}E\bigr)}_{|P}
 =\prolongg{\eta}{
 \bigl({\psi_{\vecc}^{-1}E}_{|\pi_i^{-1}(P)}\bigr)}_{|P}$.
By using the result in the curve case,
we obtain that $f(P)$ is contained in
$\lefttop{i}F'_{b+\eta\,|\,P}$.
(See Section \ref{subsection;08.10.18.40}.)
Let $\overline{f}$ denote the image of $f$
via the projection $\pi$.
Then $\overline{f}(P)=0$ for any $P\in D_i^{\circ}$.
It implies $\overline{f}=0$ on $D_i$.
Hence, we obtain
$f\in
 \prolongg{\veceta+b\vecdelta_i}
 {\bigl(\psi_{\vecc}^{-1}E\bigr)}'$.

Conversely, pick a section 
$f\in
 \prolongg{\veceta+b\vecdelta_i}
 {\bigl(\psi_{\vecc}^{-1}E\bigr)}'$.
By using the result in the curve case
(Section \ref{subsection;08.10.18.40}),
we obtain the following inequality 
for any $P\in D_i^{\circ}$:
\[
\bigl|f_{|\pi_i^{-1}(P)}\bigr|_h
=O\bigl(
 |z_i|^{-b-\eta-\epsilon}
 \bigr)
\quad (\forall \epsilon>0)
\]
Then we obtain $f\in
 \prolongg{\veceta+b\vecdelta_i}
 {\bigl(\psi_{\vecc}^{-1}E\bigr)}$
due to Proposition \ref{prop;08.10.17.12}.
In all, we obtain
$\prolongg{\veceta+b\vecdelta_i}
 {\bigl(\psi_{\vecc}^{-1}E\bigr)}
=\prolongg{\veceta+b\vecdelta_i}
 {\bigl(\psi_{\vecc}^{-1}E\bigr)}'$.
\hfill\qed

\subsubsection{Step 2}

By shrinking $X$,
let us take a $\mu_c^{\ell}$-equivariant frame 
$\vecv=(v_i)$ of
$\prolongg{\veceta}{\psi_{\vecc}^{-1}E}$
in the sense
\[
 (\omega_1,\ldots,\omega_{\ell})^{\ast}
 v_i=
 \prod_{j=1}^{\ell}
 \omega_j^{-p_j(v_i)}\cdot v_i
\]
for some $0\leq p_j(v_i)\leq c-1$.
Note that $\vecv$ is compatible with
the parabolic filtrations 
$\lefttop{j}F$ $(j=1,\ldots,\ell)$.
We set 
\[
 \vbar_i:=
 \prod z_j^{p_j(v_i)}\cdot v_i.
\]
Since they are $\mu_c^{\ell}$-invariant,
we obtain a tuple 
$\overline{\vecv}=(\overline{v}_1,\ldots,\overline{v}_r)$ 
of sections of $E$.
By using the result in the curve case,
we obtain that they are sections of
$\prolong{E}$.
(See Section \ref{subsection;08.10.18.40}).
Moreover, the restrictions
$\vecvbar_{|\pi_i^{-1}(P)}$
for $P\in D_i^{\circ}$
are frames of
$\prolong{\bigl(E_{|\pi_i^{-1}(P)}\bigr)}$.
Hence, we can conclude that
$\vecvbar$ gives a frame of
$\prolong{E}$ on $X$.
Therefore, we obtain that
$\prolong{E}$ is a locally free $\nbigo_X$-module.
Then, Theorem \ref{thm;07.3.16.1} follows.

\subsubsection{Complement}

Let us show directly that
the induced filtrations are compatible.
We have the map
$\chi_i:\bigl\{p\,\big|\,0\leq p\leq c-1,\,\,
 \lefttop{i}V_p\neq 0\bigr\}
\lrarr
 \Par\bigl( \prolong{E},i \bigr)$
as in Section \ref{subsection;08.10.18.40}.
We set $a_i(v_j):=\chi_i\bigl(p_i(v_j)\bigr)$.
We consider the filtration 
$\lefttop{i}F'_b$
of $\prolong{E}_{|D_i}$
by vector subbundles over $D_i$,
given as follows:
\[
 \lefttop{i}F'_{b}:=
 \big\langle
 \bar{v}_{j\,|\,D_i}
\,\big|\,a_i(v_j)\leq b
 \big\rangle
\]
For any $-1<b\leq 0$,
we consider the subsheaf 
$\prolongg{b\cdot\vecdelta_i}{(E)}'$
of $\prolong{E}$ given as follows:
\[
 \prolongg{b\cdot\vecdelta_i}{\bigl(E\bigr)}'
:=
 \ker\Bigl(\pi:
 \prolong{E}\lrarr
 \frac{\prolong{E}_{|D_i}}{\lefttop{i}F'_{b}}
 \Bigr).
\]
Here $\pi$ denotes the naturally defined morphism.
Then $\prolongg{b\cdot\vecdelta_i}{\bigl(E\bigr)}'$
is locally free.

\begin{lem} \label{lem;9.10.16}
$\prolongg{b\vecdelta_i}{E}=
 \prolongg{b\vecdelta_i}{\bigl(E\bigr)}'$
and
$\lefttop{i}F'_b=\lefttop{i}F_b$.
\end{lem}
\pf
Let $f$ be a holomorphic section of
$\prolongg{b\vecdelta_i}{E}$.
We can also regard it as a section of 
$\prolong{E}$.
By applying the result in the curve case
(Section \ref{subsection;08.10.18.40})
to $f_{|\pi_i^{-1}(P)}\in
   \prolong{\bigl(E_{|\pi_i^{-1}(P)}\bigr)}$,
we obtain that $f(P)\in \lefttop{i}F'_{|P}$
for any $P\in D_i^{\circ}$.
Then it is easy to derive
that $f$ is contained in
$\prolongg{b\vecdelta_i}{\bigl(E\bigr)}'$.

Conversely,
let $f$ be a holomorphic section of 
$\prolongg{b\vecdelta_i}{\bigl(E\bigr)}'$.
Applying the result in the curve case
(Section \ref{subsection;08.10.18.40})
to $f_{|\pi_i^{-1}(P)}$,
we obtain 
$f_{|\pi_i^{-1}(P)}
 \in \prolongg{b}{(E_{|\pi_i^{-1}(P)})}$
for any $P\in D_i^{\circ}$.
Then we obtain $f\in\prolongg{b\vecdelta_i}{E}$
by Proposition \ref{prop;08.10.17.12}.
Therefore  we obtain
$\prolongg{b\vecdelta_i}{E}
=\prolongg{b\vecdelta_i}{\bigl(E\bigr)}'$,
and thus $\lefttop{i}F_b=\lefttop{i}F'_b$.
\hfill\qed

\vspace{.1in}

By our construction,
$\lefttop{i}F'$ is the filtration in the category of
the vector bundles over $D_i$,
and the tuple $(\lefttop{i}F'\,|\,i=1,\ldots,\ell)$ is compatible.
Hence, 
the filtration $\lefttop{i}F$ is a filtration
in the category of the vector bundles over $D_i$,
and the tuple $(\lefttop{i}F\,|\,i=1,\ldots,\ell)$ is compatible.

\subsection{Proof of Theorem \ref{thm;07.10.9.1}}

Let $\vecv$ be a frame of $\prolongg{\vecb}{E}$
compatible with the parabolic filtrations
$(\lefttop{i}F\,|\,i=1,\ldots,\ell)$.
We obtain the numbers
$\lefttop{i}b(v_j):=\lefttop{i}\deg^F(v_j)$.
We put
\[
 v_j':=v_j\cdot\prod_{i=1}^{\ell}|z_i|^{\lefttop{i}b(v_j)},
\quad\quad
 \vecv'=(v_j').
\]
Let us show that
$\vecv'$ is adapted up to log order.
By our construction of $\vecv'$
and Proposition \ref{prop;08.10.17.12},
there exists $C_1>0$ such that the following holds:
\[
 H(h,\vecv')\leq
 C_1\cdot\bigl(-\sum \log|z_i|\bigr)^M
\]
Let $\vecv^{\lor}$ denote the dual frame of $\vecv$.
Let $P$ be a point of $D_i^{\circ}$.
According to the functoriality in the curve case
(Section \ref{subsection;08.10.19.120}),
$\vecv^{\lor}_{|\pi_i^{-1}(P)}$ gives a frame of
$\prolongg{-b_i+(1-\epsilon)}
   {E^{\lor}}_{|\pi_i^{-1}(P)}$,
which is compatible with the parabolic filtration.
Hence, $\vecv^{\lor}$ gives a frame
of $\prolongg{-\vecb+(1-\epsilon)\vecdelta}{E^{\lor}}$
for some $\epsilon>0$.
We have
$\lefttop{i}\deg^F(v^{\lor}_j)
=\deg^F\bigl(v^{\lor}_{j\,|\,\pi_i^{-1}(P)}\bigr)
=-\lefttop{i}b(v_j)$
for any point $P\in D^{\circ}_i$.
We put
\[
 \vecv^{\lor\,\prime}=(v_j^{\lor\,\prime}),
\quad\quad
 v_j^{\lor\,\prime}
:=v_j^{\lor}\cdot\prod_{i=1}^{\ell}|z_i|^{-\lefttop{i}b(v_j)}.
\]
Due to Proposition \ref{prop;08.10.17.12},
we obtain the following:
\[
 H(h^{\lor},\vecv^{\lor\,\prime})
\leq
 C_2\cdot
 \Bigl(
 -\sum_{i=1}^{\ell}\log |z_i|
 \Bigr)^M
\]
Here, $h^{\lor}$ denote the induced metric of $E^{\lor}$.
It implies the following:
\[
 C_3\cdot \Bigl(
 -\sum_{i=1}^{\ell}\log |z_i|
 \Bigr)^{-M}
\leq H(h,\vecv')
\]
Thus, we obtain Theorem \ref{thm;07.10.9.1}.
\hfill\qed

\section{Small deformation}
\label{subsection;08.10.18.51}
\subsection{Statement}

Let $X:=\Delta^{\ell}$,
$D_i:=\{z_i=0\}$
and $D:=\bigcup_{i=1}^{\ell}D_i$.
Let $g_{\poin}$ denote 
the Poincar\'e metric of $X-D$.
Let $(E,\delbar_E,h)$ be an acceptable bundle
on $X-D$ of rank $r$.
For any $0<C<1$,
we set $X(C):=\bigl\{
 (z_1,\ldots,z_n)\in X\,\big|\,
 |z_i|<C \bigr\}$,
$D(C):=D\cap X(C)$
and $X^{\ast}(C):=X(C)\setminus D(C)$.

Let $Y$ be an open subset of $\cnum_{\lambda}$
with a base point $\vecy_0$,
and let $p$ denote the projection
$Y\times (X-D)\lrarr X-D$.
Let $\nbiga$ be a section of
$p^{\ast}\bigl(
 \End(E)\otimes\Omega^{0,1}_{X-D}
 \bigr)$ with the following property:
\begin{itemize}
\item
$\bigl(\delbar_E+\delbar_{\lambda}+\nbiga\bigr)^2=0$.
The holomorphic bundle
$\bigl(p^{-1}E,
 \delbar_E+\delbar_{\lambda}+\nbiga
 \bigr)$ is denoted by $(\nbige,\delbar_{\nbige})$.
\item
$\bigl|\nbiga\bigr|_{h,g_{\poin}}$ is bounded,
and 
$\bigl|\nbiga_{\vecy}
-\nbiga_{\vecy_0}\bigr|_{h,g_{\poin}}
\leq B\cdot |\vecy-\vecy_0|$ for some $B>0$,
where $\nbiga_{\vecy}:=
 \nbiga_{|\vecy\times (X-D)}$.
\item
$\nbiga_{\vecy_0}=0$
and $\tr\nbiga=0$.
\end{itemize}
We have the associated filtered sheaf
$\nbige_{\ast}$
as in Subsection \ref{subsection;07.9.22.10}.

\begin{thm}
\label{thm;08.10.18.11}
There exist a positive number $C>0$
and a neighbourhood $\nbigu_0$ of $\vecy_0$
such that the following holds:
\begin{itemize}
\item
The restriction of $\nbige_{\ast}$ 
to $\nbigu_0\times X(C)$ is a filtered bundle.
\item
Let $\vecv$ be a frame of
$\prolongg{\veca}\nbige$
compatible with the parabolic filtrations
$\bigl(\lefttop{i}F\,\big|\,i=1,\ldots,\ell \bigr)$.
We set $a_i(v_j):=\lefttop{i}\deg^F(v_j)$
and 
\[
 v_j':=v_j\cdot\prod_{i=1}^{\ell}|z_i|^{a_i(v_j)}.
\]
Then, for any $\epsilon>0$, there exists
 $C_{\epsilon}>0$ such that the following holds:
\[
 C_{\epsilon}^{-1}\cdot
 \prod_{i=1}^{\ell}|z_i|^{\epsilon}
\leq
 H(h,\vecv')
\leq
 C_{\epsilon}\cdot
 \prod_{i=1}^{\ell}|z_i|^{-\epsilon}
\]
\end{itemize}
\end{thm}

In the following,
$\nbigu_i$ will denote a neighbourhood of $\vecy_0$,
and $C_i$ will denote a constant such that $0<C_i<1$.

\subsection{Extension of holomorphic sections}

Let us argue the extension problem
of holomorphic sections on
$\{\vecy_0\}\times X$
to those on $\nbigu\times X$.
First, we remark that 
an $L^2$-estimate implies a growth estimate.

\begin{lem}
\label{lem;08.10.18.1}
Let $F$ be a holomorphic section of $\nbige$
on $\nbigu\times (X-D)$
such that
$\bigl\| F \bigr\|_{\veca,N}<\infty$.
Let $\nbigu'$ be relatively compact in $\nbigu$.
Then, there exist $M>0$  and $0<C<1$
such that the following holds
on $\nbigu'\times X^{\ast}(C)$:
\[
 |F|_h=O\Bigl(
 \prod_{i=1}^{\ell}|z_i|^{-a_i}
 \bigl(-\log|z_i|\bigr)^M
 \Bigr)
\]
\end{lem}
\pf
Let $\hyperh$ denote the upper half plane.
We use the coordinate
$\zeta=\xi+\sqrt{-1}\eta$.
Let $\varphi_C:\hyperh^{\ell}\times\Delta^{n-\ell}
\lrarr (\Delta^{\ast})^{\ell}
 \times\Delta^{n-\ell}=X-D$
be the map given by
\[
 \varphi_C(\zeta_1,\ldots,\zeta_{\ell},
 z_{\ell+1},\ldots,z_n)
=\bigl(
 e^{2\pi\sqrt{-1}\zeta_1},\ldots,
 e^{2\pi\sqrt{-1}\zeta_{\ell}},
 Cz_{\ell+1},\ldots,Cz_n
 \bigr)
\]
Let $g$ be the Euclidean metric of
$\hyperh^{\ell}\times\Delta^{n-\ell}$.
For any $\kappa>0$,
there exists $C(\kappa)>0$ such that
$\bigl|\varphi_{C(\kappa)}^{\ast}R(h)
 \bigr|_{\varphi_{C(\kappa)}^{\ast}h,g}
 \leq \kappa$
on 
$\bigl\{
 \eta_j\geq C(\kappa),\,j=1,\ldots,\ell
 \bigr\}\subset\hyperh^{\ell}\times\Delta^{n-\ell}$.

For $\vecn=(n_1,\ldots,n_{\ell})$,
we put
\[
K_{\vecn}:=
 \prod_{i=1}^{\ell}\bigl\{(\xi_i,\eta_i)\,\big|\,
 -1<\xi_i<1,n_i-1<\eta_i<n_i+1\bigr\}
\times \Delta^{n-\ell}.
\]
Let $\kappa$ be sufficiently small.
If $n_i>C(\kappa)$,
due to a theorem of Uhlenbeck,
we can take an orthonormal frame
$\vece_{\vecn}$
of $\varphi_{C(\kappa)}^{\ast}(E,\delbar_E,h)$
on $K_{\vecn}$ such that the connection form
$A_{\vecn}$ is sufficiently small
with respect to $h$ and $g$,
where $A_{\vecn}$ is determined by
$\delbar_E\vece_{\vecn}
=\vece_{\vecn}\cdot A_{\vecn}$.

The induced map
$\nbigu\times \hyperh^{\ell}\times\Delta^{n-\ell}
\lrarr\nbigu\times (X-D)$ is also denoted by 
$\varphi_{C(\kappa)}$.
We have the orthonormal frame of
$(\varphi_{C(\kappa)}^{\ast}\nbige)
_{|\nbigu\times K_{\vecn}}$
induced by $\vece_{\vecn}$,
which is denoted by $\vecetilde_{\vecn}$.
Let $\Atilde_{\vecn}$ be given by
$\delbar_{\nbige}\vecetilde_{\vecn}
 =\vecetilde_{\vecn}\cdot \Atilde_{\vecn}$.
By the assumption,
$\Atilde_{\vecn}$ are
uniformly bounded with respect to $g$.
We have the expression
$F_{|K_{\vecn}}=
 \sum F_{\vecn,i}\cdot \etilde_{\vecn,i}$
on each $\nbigu\times K_{\vecn}$.
We put $\vecF_{\vecn}
     =\bigl(F_{\vecn,i}\bigr)$
which satisfies
$\delbar\vecF_{\vecn}
+\Atilde_{\vecn}\vecF_{\vecn}=0$.
We also have the following:
\[
 \int_{K_{\vecn}}
 \bigl|
 \vecF_{\vecn}
 \bigr|^2\cdot \dvol_g
\leq
 B\cdot
 \bigl\| F \bigr\|_{\veca,N}^2\cdot
  \prod_{i=1}^{\ell}
 e^{2a_i\cdot n_i}
 \cdot n_i^{|N|+2}
\]
Here, $B$ is a positive constant
depending on $C(\kappa)$.
If we take a relatively compact subset
$K'_{\vecn}\subset K_{\vecn}$,
we can obtain the estimate of the sup norm of 
$\vecF_{\vecn}$
on $\nbigu'\times K'_{\vecn}$
by a standard boot strapping argument.
Thus, the claim of Lemma \ref{lem;08.10.18.1} follows.
\hfill\qed

\vspace{.1in}

We set $\vecdelta:=(1,\ldots,1)\in\real^{\ell}$.
\begin{lem}
There exist $\nbigu_1$ and $C_1>0$
satisfying the following:
\begin{itemize}
\item
 Let $f$ be any holomorphic section of 
 $\prolongg{\veca}{E}$ on $X-D$.
 For any $\epsilon>0$,
 there exists a section
 $F^{(\epsilon)}$ of 
 $\prolongg{\veca+\epsilon\vecdelta}{\nbige}$
 on $\nbigu_1\times X(C_1)$
 such that
 $F^{(\epsilon)}_{|\vecy_0\times X(C)}
 =f_{|X(C)}$.
\end{itemize}
\end{lem}
\pf
Let $N<0$ be as in Lemma \ref{lem;07.9.24.2}.
For any $\kappa>0$,
there exists a small neighbourhood $\nbigu_2$
such that 
$\bigl\|\nbiga_{\vecy}\bigr\|_h
 \leq \kappa$
for any $\vecy\in \nbigu_2$.
Note that the norm of
the morphism
\[
 \nbiga_{\vecy}:
 A^{0,q}_{\vecb,N}(E)
\lrarr
 A^{0,q+1}_{\vecb,N}(E)
\]
is dominated by $\bigl|\nbiga_{\vecy}\bigr|_h$
for any $\vecb\in\real^{\ell}$.
We can regard $f$ as a section of
$A^{0,0}_{\veca+\epsilon\vecdelta,N}(E)$.
Then, by using a standard argument
in Section 2.9.1 of \cite{mochi2},
we can find a holomorphic section $F^{(\epsilon)}$
of $\nbige$ on $\nbigu_2\times (X-D)$
such that
(i) $F^{(\epsilon)}_{|\vecy_0\times (X-D)}=f$,
(ii) $\bigl\|F^{(\epsilon)}
 \bigr\|_{\veca+\epsilon\vecdelta,N}<\infty$.
By using Lemma \ref{lem;08.10.18.1},
we obtain 
$F^{(\epsilon)}
 \in \prolongg{\veca+\epsilon\vecdelta}{\nbige}$
on $\nbigu_1\times X(C_1)$
for some appropriate $\nbigu_1\subset \nbigu_2$
and $0<C_1<1$.
\hfill\qed

\subsection{Construction of local frames}

By shrinking $X$,
we take a frame $\vecu=(u_i)$ of
$\prolongg{\veca}{E}$ on $X-D$
compatible with the parabolic structure.
Let $a_j(u_i)$ denote the parabolic degree of $u_i$
with respect to the filtration $\lefttop{j}F$,
and we put 
$\veca(u_i)=\bigl(a_j(u_i)\,|\,j=1,\ldots,\ell\bigr)$.

By the assumption $\tr\nbiga=0$,
we have $\det\nbige=p^{\ast}\det(E)$.
We put
\[
 \atilde_j:=
  \sum_{b\in\Par(\prolongg{\veca}{E},j)} 
 b\cdot \rank\lefttop{j}\Gr^F_b\bigl(
 \prolongg{\veca}{E} \bigr)
=\sum_{i=1}^r a_j(u_i)
\]
The tuple $(\atilde_j\,|\,1\leq j\leq \ell)$
is denoted by $\vecatilde$.
Note $\det(\prolongg{\veca}{E})
\simeq \prolongg{\vecatilde}{\det E}$.

Take $\epsilon>0$
such that 
$10 r^2\cdot \epsilon<|b-b'|$
for any distinct $b,b'\in \Par(E,i)$
and any $i=1,\ldots,\ell$.
We extend $u_i$ to the section $u^{(\epsilon)}_i$
of $\prolongg{\veca(u_i)+\epsilon\vecdelta}\nbige$
on $\nbigu_1\times X(C_1)$.
Since $\epsilon$ is sufficiently small,
$\Omega(\vecu^{(\epsilon)}):=
 u_1^{(\epsilon)}\wedge\cdots
 \wedge u_r^{(\epsilon)}$
gives a holomorphic section of
$\prolongg{\vecatilde}{\det\nbige}$.
Let $Z^{(\epsilon)}$ denote 
the $0$-set of $\Omega(\vecu^{(\epsilon)})$.
Since
$\Omega(\vecu^{(\epsilon)})_{|(\vecy_0,O)}\neq 0$
in 
$\prolongg{\vecatilde}{\det\nbige}_{|(\lambda_0,O)}$,
there exist 
$\nbigu_3(\epsilon)$ and $0<C_3(\epsilon)<C_1$
such that
$Z^{(\epsilon)}\cap \bigl(
\nbigu_3(\epsilon)\times X(C_3(\epsilon))
 \bigr)
=\emptyset$.

\begin{lem}
$\vecu^{(\epsilon)}$ gives a frame of
$\prolongg{\veca+\epsilon\vecdelta}{\nbige}$
on $(\nbigu_1\times X(C))\setminus Z^{(\epsilon)}$,
in particular
on $\nbigu_3(\epsilon)\times X(C_3(\epsilon))$.
\end{lem}
\pf
Let $P$ be any point of
$\nbigu_1\times D(C)\setminus Z^{(\epsilon)}$.
Let $\nbigv$ be a small neighbourhood of $P$.
The restriction of
$\vecu^{(\epsilon)}$ to
$\nbigv\setminus (\nbigu_1\times D)$ give a frame.
Let $f$ be a section of 
$\prolongg{\veca+\epsilon\vecdelta}{\nbige}$
on $\nbigv$.
We have the expression
$f=\sum f_i\cdot u_i^{(\epsilon)}$,
where $f_i$ are holomorphic on 
$\nbigv\setminus (\nbigu_1\times D)$.
Since $\epsilon$ is sufficiently small,
it is easy to observe that
$f\wedge u_2^{(\epsilon)}
 \wedge\cdots\wedge u_r^{(\epsilon)}$
is a holomorphic section of 
$\prolongg{\vecatilde}{\det\nbige}$
on $\nbigv$.
It is equal to
$f_1\cdot \Omega(\vecu^{(\epsilon)})$,
and hence $f_1$ are holomorphic on $\nbigv$.
We obtain that the other $f_i$
are also holomorphic on $\nbigv$
in the same way.
Hence,
$\vecu^{(\epsilon)}$ gives a frame of
$\prolongg{\veca+\epsilon\vecdelta}\nbige$
on $\nbigv$.
\hfill\qed

\vspace{.1in}

We set 
$u^{(\epsilon)\prime}_j:=
 \prod_{i=1}^{\ell}
 |z_i|^{a_i(u_j)}\cdot
 u^{(\epsilon)}_j$.
\begin{lem}
\label{lem;08.10.18.10}
Let $P$ be any point of
$(\nbigu_1\times D(C_1))\setminus 
 \nbigz^{(\epsilon)}$,
and let $\nbigv_P$ be any neighbourhood of $P$.
Then, there exists a constant $B_P>0$ 
such that
the following holds on $\nbigv_P$:
\[
 B_P^{-1}\prod_{i=1}^{\ell}
 |z_i|^{2r^2\epsilon}
\leq
 H(h,\vecu^{(\epsilon)\prime})
\leq
 B_P\prod_{i=1}^{\ell}
 |z_i|^{-2r^2\epsilon}
\]
\end{lem}
\pf
The right inequality is clear.
Let $\vecu^{(\epsilon)\lor}$ be the dual frame 
of $\vecu^{(\epsilon)}$
on $(\nbigu_1\times X(C_1))\setminus Z^{(\epsilon)}$.
There exists a constant $B_{P\,i}>0$ $(i=1,2)$
such that the following holds:
\[
 \bigl|u_j^{(\epsilon)\lor}\bigr|_{h}
\leq
B_{P\,1}
 \frac{\bigl|
 u^{(\epsilon)}_1
 \wedge
 \cdots\wedge
 \check{u}_j^{(\epsilon)}
 \cdots\wedge
 u_r^{(\epsilon)}
 \bigr|}{|\Omega(\vecu^{(\epsilon)})|_h}
\leq
B_{P\,2}
\prod_{i=1}^{\ell}
 \bigl|z_i\bigr|^{a(u_i)-r\epsilon}
\]
Hence, we obtain the following for some $B_{3\,i}>0$:
\[
 H(h^{\lor},\vecu^{(\epsilon)\prime\lor})
\leq
 B_{P\,3} \prod_{i=1}^{\ell}
 |z_i|^{-2r^2\epsilon}
\]
Thus we obtain Lemma \ref{lem;08.10.18.10}.
\hfill\qed

\subsection{Proof of Theorem
 \ref{thm;08.10.18.11}}

Let us fix small $\epsilon_0>0$ as above,
and we take a frame $\vecu^{(\epsilon_0)}$
of $\prolongg{\veca+\epsilon_0\vecdelta}{\nbige}$
on $\nbigu_3(\epsilon_0)\times 
 X\bigl(C_3(\epsilon_0)\bigr)$
as above.
Let us show that $u_j^{(\epsilon_0)}$
are actually sections of
$\prolongg{\veca}{\nbige}$ on
$\nbigu_3(\epsilon_0)\times X\big(C_3(\epsilon_0)\bigr)$,
and then 
the first claim of Theorem \ref{thm;08.10.18.11} follows.

Let $0<\epsilon\leq \epsilon_0$.
Let $P$ be a point of
$(\nbigu_3(\epsilon_0)\times D)
 \setminus Z^{(\epsilon)}$,
and let $\nbigv$ be a small neighbourhood
of $P$ such that 
$\vecu^{(\epsilon)}$ gives a frame
of $\prolongg{\veca+\epsilon\vecdelta}{\nbige}$
on $\nbigv$.
We have the expression
$u^{(\epsilon_0)}_i
=\sum b_{j,i}\cdot u^{(\epsilon)}_j$.
By using Lemma \ref{lem;08.10.18.10},
we obtain that $b_{j,i}$ are holomorphic
on $\nbigv$,
and thus
$\prolongg{\veca+\epsilon_0\vecdelta}
 {\nbige}
=\prolongg{\veca+\epsilon\vecdelta}
 {\nbige}$ on $\nbigv$.

For each $I\subset\ellsitabar$,
we put $I^c:=\ellsitabar-I$,
and 
\[
D_I^{\circ}:=
 \bigcap_{i\in I}D_i
\setminus
 \Bigl(
 \bigcup_{j\in I^c}D_j
 \Bigr)
\]
We have the stratification
$Z^{(\epsilon)}=
 \coprod_I
 \coprod_j Z_{I,j}^{(\epsilon)}$
satisfying the following:
\begin{itemize}
\item
$Z_{I,j}^{(\epsilon)}
 \subset \nbigu_1\times D_I^{\circ}$
and $\dim Z_{I,j}=j$.
\item
$Z_{I,j}^{(\epsilon)}$ are 
smooth and locally closed in
$\nbigu_1\times X$.
Moreover,
the induced morphism
$Z_{I,j}^{(\epsilon)}\lrarr
 \nbigu_1\times X\lrarr X$
is immersion.
\end{itemize}
Let $P$ be any point of 
$Z_{I,j}$.
We show that
$\prolongg{\veca+\epsilon\vecdelta}{\nbige}
=\prolongg{\veca+\epsilon_0\vecdelta}{\nbige}$
around $P$,
by using the ascending induction on $|I|$
and the the descending induction on $j$.

Recall $\dim Y=1$.
We take a small neighbourhood $T_P$ of $P$ 
in $Z_{I,j}$.
We take a small tubular neighbourhood 
$T^{(1)}_P$ of $T_P$
in $\nbigu_1\times D_I^{\circ}$.
We have identifications
$T^{(1)}_P\simeq T_P\times N_P$
and $T_P\simeq T_P\times\{(0,\ldots,0)\}$,
where $N_P$ denotes
a $(1+n-|I|-j)$-dimensional multi-disc.
We take a small tubular neighbourhood
$T^{(2)}_P$ of $T^{(1)}_P$
in $\nbigu_1\times X(C_1)$.
We have identifications
$T^{(2)}_P\simeq
 T^{(1)}_P\times \Delta_w^{|I|}$
and $T^{(1)}_P\simeq
 T^{(1)}_P\times \{(0,\ldots,0)\}$,
where $\Delta_w^{|I|}$ is
an $|I|$-dimensional multi-disc.
Moreover, 
$T^{(2)}_P\cap 
 (\nbigu_1\times D)
=T^{(1)}_P\times
 \bigcup_{i\in I}\{z_i=0\}$.
By the assumption of the induction,
we have the following estimate 
on $T_P\times \del N_P\times (\Delta^{\ast})^I$
for any $\epsilon'>\epsilon$:
\begin{equation}
 \label{eq;08.10.18.20}
 \bigl|u_i^{(\epsilon_0)}\bigr|_{h}
=O\Bigl(
 \prod_{k\in I}
 |z_k|^{-a_k-\epsilon'}
 \Bigr)
\end{equation}
Since $\bigl|u_i^{(\epsilon_0)}\bigr|_{h}$
is pluri-subharmonic 
in the $Y$-direction,
(\ref{eq;08.10.18.20}) holds
on $T_P\times N_P\times (\Delta^{\ast})^I$
for any $\epsilon'>\epsilon$.
It means
$u_i^{(\epsilon_0)}
\in \prolongg{\veca+\epsilon\vecdelta}{\nbige}$
around $P$.
Since $\epsilon>0$ is arbitrary,
we obtain that 
$u_i^{(\epsilon_0)}
\in \prolongg{\veca}{\nbige}$.
Thus, 
the proof of the first claim of
Theorem \ref{thm;08.10.18.11}
is completed.

\vspace{.1in}

Let us show the second claim of the theorem.
Because we have also obtained
$\prolongg{\veca(u_i)+\epsilon_0\vecdelta}
 {\nbige}
=\prolongg{\veca(u_i)}{\nbige}$,
the following holds for any $\epsilon>0$:
\[
 \bigl|
 u^{(\epsilon_0)}_{i}
 \bigr|_h=O\Bigl(
\prod_{j=1}^{\ell}|z_j|^{-a_j(u_i)-\epsilon}
 \Bigr)
\]
Hence, there exists $C_{\epsilon}>0$
such that the following holds:
\[
 H(h,\vecu^{(\epsilon_0)\prime})
\leq
 C_{\epsilon}\prod_{j=1}^{\ell}
 |z_j|^{-\epsilon}
\]
By the argument in the proof of Lemma \ref{lem;08.10.18.10},
we also obtain that there exists $C_{\epsilon}'$
such that
\[
 H(h^{\lor},\vecu^{(\epsilon_0)\prime\lor})
\leq
 C_{\epsilon}'\prod_{j=1}^{\ell}
 |z_j|^{-\epsilon}.
\]
Thus, the second claim of Theorem
\ref{thm;08.10.18.11}
is proved.
\hfill\qed

\section{Complement}
\label{subsection;08.10.18.52}
\subsection{Preliminary estimate for sup norms}
\label{subsection;08.9.29.5}

Let $X:=\Delta^n$,
$D_i:=\{z_i=0\}$ and 
$D:=\bigcup_{i=1}^{k} D_i$
for some $k\leq n$.
We use the Poincar\'e metric 
$g_{\poin}$ of $X-D$.
We put
$X(C):=\bigl\{(z_1,\ldots,z_n)\in X\,\big|\,
 |z_i|\leq C\bigr\}$,
$D_i(C):=D_i\cap X(C)$
and $D(C):=D\cap X(C)$
for $0<C<1$.
We set
\begin{equation}
 \label{eq;07.9.24.10}
 Z(C):=
 \bigl\{(z_1,\ldots,z_n)\in X(C)-D(C)\,\big|\,
 |z_i|=C\,\,(i=1,\ldots,k)
 \bigr\}.
\end{equation}

Let $(E,\delbar_E,h)$ be an acceptable bundle on $X-D$.
We use the notation in
Section \ref{subsection;d11.14.35}.
We  assume $\bigl|R(h)\bigr|_{h,g_{\poin}}\leq C_0$
for some given constant $C_0$.
Let $\omega$ be a $C^{\infty}$-section
of $E\otimes\Omega^{0,1}$ on $X-D$.
Assume
$\sup_{X-D}|\omega|_{\veca,N}<\infty$
for simplicity.
Let $f$ be a $C^{\infty}$-section of
$E\otimes\Omega^{0,1}$ on $X-D$
such that
$\delbar f=\omega$
and 
$\|f\|_{\veca,N}\leq\|\omega\|_{\veca,N}$.
Note 
$\|\omega\|_{\veca,N}
\leq
 \eta\cdot\sup|\omega|_{\veca,N}$,
where $\eta>0$ is independent of $\omega$.

\begin{lem}
\label{lem;07.9.24.5}
There exists a constant $C_1$,
depending on $C_0$, $C$ and $N$,
such that the following holds:
\[
 \sup_{X(C)-D(C)}|f|_{\veca,N+k}
\leq C_1\cdot \sup_{X-D}|\omega|_{\veca,N}.
\]
\end{lem}
\pf
We give only an outline.
Let $\hyperh$ denote the upper half plane
$\bigl\{\zeta\in\cnum\,\big|\,
 \Image(\zeta)>0\bigr\}$.
We have the universal covering map
$\varphi:\hyperh\lrarr \Delta^{\ast}$
given by $\varphi(\zeta)=\exp(\sqrt{-1}\zeta)$,
which induces
$\varphi:\hyperh^{k}\times\Delta^{n-k}\lrarr X-D$.
We put $\nbigu(C):=\varphi^{-1}\bigl(X(C)-D(C)\bigr)$.
Let $g$ and $\dvol_g$
be the Euclidean metric of
$\hyperh^k\times \Delta^{n-k}$
and the associated volume form.
Note that the fiber-wise norm of
$\varphi^{\ast}\omega$ with respect to
$\varphi^{\ast}h$ and $g$
is dominated by $\varphi^{\ast}|\omega|_{\veca,N}$.
Let $B_{\epsilon}$ be an $\epsilon$-multi-ball
contained in $\nbigu(C)$.
We have 
$ \int_{B_{\epsilon}}
 \bigl|f\bigr|_{\veca,N+k}^2 \dvol_g
\leq
 \bigl\|f\bigr\|_{\veca,N}^2
\leq \|\omega\|_{\veca,N}^2$.
We take the natural diffeomorphism
$\psi:B_1\simeq B_{\epsilon}\subset 
 \hyperh^k\times\Delta^{n-k}$,
where $B_1:=\bigl\{(w_1,\ldots,w_n)\,\big|\,|w_i|<1\bigr\}$.
Let $g_{1}$ and $\dvol_{g_1}$ denote
the Euclidean metric and the associated volume form
of $B_1$.
We fix a sufficiently small $\epsilon$
depending only on $C$, $C_0$ and $N$,
such that
the curvature of 
$\psi^{\ast}\varphi^{\ast}
 (E,\delbar_E,h_{\veca,N+k})$
is sufficiently small
with respect to 
$\psi^{\ast}\varphi^{\ast}h_{\veca,N+k}$
and $g_{1}$ to which
a theorem of Uhlenbeck can be applied.
Then, we can take an orthogonal frame $\vece$
of $\psi^{\ast}\varphi^{\ast}(E,h_{\veca,N+k})$
for which the $(0,1)$-form $A$ is sufficiently small
with respect to $g_1$,
where $A$ is given by $\delbar_E\vece=\vece\cdot A$.
We have the expressions
$\psi^{\ast}\varphi^{\ast}f
=\sum F_i\cdot e_i$
and 
$\psi^{\ast}\varphi^{\ast}\omega
=\sum G_i\cdot e_i$.
The $L^2$-norm of $\vecF=(F_i)$
with respect to $\dvol_{g_1}$
is dominated by 
$C_1'(C_0,C,N)\cdot \|\omega\|_{\veca,N}$,
where $C_1'(C_0,C,N)$ denotes a constant
depending only on 
$C_0$, $C$ and $N$.
The norm of $\vecG=(G_i)$ with respect to $g_1$
is dominated by $|\omega|_{\veca,N}$.
We have the relation
$\delbar\vecF+[A,\vecF]=\vecG$.
By using the standard bootstrapping argument,
we obtain the desired estimate for the sup-norm of $\vecF$
in terms of $\sup_{X-D}|\omega|_{\veca,N}$.
\hfill\qed

\subsection{A priori refinement of 
estimate for the sup norm}
\label{subsection;08.9.29.3}

Let $\rho$ be a $C^{\infty}$-section of $E$ on $X-D$
such that the following holds for some large $M$:
\[
|\rho|_h=O\left(\prod_{i=1}^k|z_i|^M\right),
\quad
 \bigl|\delbar_E\rho\bigr|_{h,g_{\poin}}
=O\left(\prod_{i=1}^k|z_i|^M\right)
\]
Let $Z(C)$ be as in (\ref{eq;07.9.24.10}).
We define
\[
 \bigl|\rho\bigr|_{\infty,Z(C)}:=
 \sup_{Z(C)}\bigl|\rho\bigr|_h,
\quad
 \|\delbar_E\rho\|_{\infty}:=
 \sup_{X-D}\left|\delbar_E\rho\cdot\prod_{i=1}^k|z_i|^{-M}
 \right|_{h,g_{\poin}}
\]
\begin{lem}
\label{lem;07.7.7.40}
Let $M'<M$.
There exists a constant $C_2$,
depending only on $C_0$, $C$ and $M'$,
such that the following holds:
\[
\sup_{X(C)-D(C)}\Bigl(
 |\rho|_h\cdot
\prod_{i=1}^k|z_i|^{-M'}
\Bigr)
 \leq
 C_2
 \Bigl(
 |\rho|_{\infty,Z(C)}
+\|\delbar_E\rho\|_{\infty}
 \Bigr)
\]
\end{lem}
\pf
In the following, $C_i>0$ denote
constants depending only on 
$C_0$, $C$ and $M'$.
Take $M'<M_1<M_2<M$.
According to Lemma \ref{lem;07.9.23.36},
we can take a $C^{\infty}$-section $g$ of $E$
on $X-D$
satisfying the equality
$\delbar_Eg=\delbar_E\rho$ and the following estimate:
\[
 \int_{X-D}|g|_h^2\cdot\prod_{i=1}^k|z_i|^{-2M_2}
 \left(-\sum_{i=1}^k\log|z_i|\right)^{C_{11}}\dvol
 \leq
 C_{10}\cdot\|\delbar_E\rho\|_{\infty}^2
\]
According to Lemma \ref{lem;07.9.24.5},
we obtain
$|g|_h\leq C_{15}\|\delbar_E\rho\|_{\infty}\cdot
 \prod_{i=1}^k|z_i|^{M_1}$
on $X(C)-D(C)$.
Then, we apply Proposition \ref{prop;07.7.7.3}
to $\rho-g$.
\hfill\qed

\subsection{Connection form of an acceptable bundle}
\label{subsection;08.9.28.8}

Let $X,D$ and $(E,\delbar_E,h)$ be as in
Section \ref{subsection;07.9.22.10}.
As remarked in Theorem \ref{thm;07.3.16.1},
we obtain the prolongment $\prolongg{\veca}{E}$
for each $\veca\in\real^{\ell}$
from $(E,\delbar_E,h)$.
Let $\vecv$ be a holomorphic frame of
$\prolongg{\veca}{E}$
compatible with the parabolic structures.
Then, the $C^{\infty}$-section $F$
of $\End(E)\otimes\Omega^{1,0}$
is determined by $F(\vecv)=\del_E\vecv$.

\begin{lem}
\label{lem;07.6.2.20}
$F$ is bounded up to log order,
with respect to $h$ and the Poincar\'e metric 
$g_{\poin}$ of $X-D$.
\end{lem}
\pf
Let $\vecv'$ be as in Section \ref{subsection;07.9.22.10}.
Let $h_0$ be the hermitian metric of $E$
determined by $h_0(v_i',v_j')=\delta_{i,j}$.
We have $R(h_0)=0$.
Let $s$ be the endomorphism of $E$
determined by $h=h_0\cdot s$,
which is self-adjoint with respect to
both $h$ and $h_0$.
We have $\del_h-\del_{h_0}=s^{-1}\del_{h_0}s$
and $R(h)=\delbar_E(s^{-1}\del_{h_0}s)$.

Let $\pi_1$ be the natural projection of  
$X-D$ to $D_1$.
We have only to show the following
inequality for some positive constants $C$ and $N$
which are independent of
$P\in D_1\setminus \bigcup_{1<j\leq \ell}D_j$:
\[
 \bigl|s^{-1}\del_{h_0}s_{|\pi_1^{-1}(P)}
 \bigr|_{h_0,g_{\poin}}
\leq
 C\cdot\left(-\sum_{j=1}^{\ell} \log|z_j|\right)^N
\]

Let $\kappa_1$ be 
an $\real_{\geq \,0}$-valued $C^{\infty}$-function
on $\real$ such that
$\kappa_1(t)=1$ for $t\leq 1/2$
and $\kappa_1(t)=0$ for $t\geq 2/3$.
Let $\kappa_2$ be 
an $\real_{\geq\,0}$-valued $C^{\infty}$-function
on $\real$ such that 
$\kappa_2(t)=0$ for $t\leq 1/3$
and $\kappa_2(t)=1$ for $t\geq 1$.
We put $\chi_L(z):=
\kappa_1\bigl(-L^{-1}\log|z|\bigr)
\cdot\kappa_2(-\log|z|)$ for any positive number $L$.
Let $\rho$ be any $C^{\infty}$-function on
$\pi_1^{-1}(P)$.
In the following,
we will consider integrals over $\pi_1^{-1}(P)$.
We often use
$\int$ instead of $\int_{\pi_1^{-1}(P)}$
to save the space.
We have the following:
\begin{multline}
 \label{eq;07.5.19.1}
 \int_{\pi_1^{-1}(P)}
 \bigl(s^{-1}\del_{h_0}(\chi_L s),\,
 \del_{h_0}(\chi_L s)\bigr)_{h_0}\!\cdot\!
 \rho
=\\
 \int \bigl(\delbar(s^{-1}\del_{h_0}(\chi_L s)),\,
 \chi_L s\bigr)_{h_0}\!\cdot\!
 \rho
-\!\int \bigl(s^{-1}\del_{h_0}(\chi_L s),\,
 \chi_L s\bigr)_{h_0}\!\cdot\!
 \delbar\rho
\end{multline}
We have the following equalities:
\begin{multline}
\label{eq;07.5.19.2}
 \int_{\pi_1^{-1}(P)}
 \bigl(\delbar (s^{-1}\del_{h_0}(\chi_L \cdot s)),\,
 \chi_L s\bigr)_{h_0}
 \rho=
 \\
\int\bigl(R(h)\chi_L,\chi_L s\bigr)_{h_0}\rho
+\int\tr\bigl(\delbar\chi_L\cdot\del_{h_0}s\cdot\chi_L\bigr)\rho
+\int\bigl(\delbar\del\chi_L,\,\chi_L s\bigr)_{h_0}\rho \\
=\int\bigl(R(h)\chi_L,\,\chi_L s\bigr)_{h_0}\rho
-\int\tr\bigl(\delbar\chi_L \cdot s \cdot\del\chi_L\bigr)\rho
-\int \tr\bigl(\chi_L\cdot\delbar\chi_L \cdot s\bigr)\del \rho
\end{multline}
\begin{multline}
 \label{eq;07.5.19.3}
 \int_{\pi_1^{-1}(P)}
 \bigl(s^{-1}\del_{h_0}(\chi_L s),\,\chi_L s\bigr)_{h_0}
 \delbar\rho
=\int
 \tr\bigl(\del_{h_0}(\chi_L s)\cdot\chi_L \cdot\delbar\rho\bigr) \\
=\int\del\bigl(\tr(\chi_L^2s)\cdot\delbar\rho\bigr)
-\int\tr\bigl(\chi_L s\cdot\del\chi_L\bigr)\delbar\rho
-\int \tr\bigl(\chi_L^2 s\bigr)\cdot\del\delbar\rho
\end{multline}
Let $C_2(P):=-\sum_{2\leq j\leq \ell}\log|z_j(P)|$.
We have the following estimates,
for some sufficiently large $N$:
\begin{equation}
 \label{eq;07.5.19.4}
 \bigl|s_{|\pi_1^{-1}(P)}\bigr|_{h_0}\leq
 \bigl(-\log|z_1|+C_2(P)\bigr)^N,
\quad
\bigl|R(h)_{|\pi_1^{-1}(P)}\bigr|_{h_0,g_{\poin}}\leq
  \bigl(-\log|z_1|+C_2(P)\bigr)^N
\end{equation}
We put $\rho_M:=(-\log|z_1|)^{-M}$
for some sufficiently large $M$.
By using (\ref{eq;07.5.19.1}),
(\ref{eq;07.5.19.2}),
(\ref{eq;07.5.19.3}) and (\ref{eq;07.5.19.4}),
we obtain the following for some constants $C_3$
and $M'$:
\[
  \int_{\pi_1^{-1}(P)}
 \bigl(s^{-1}\del_{h_0}(\chi_L s),\,\del_{h_0}(\chi_L s)\bigr)
 \cdot\rho_M
\leq
 C_3\cdot C_2(P)^{M'}
\]
By making $L\to \infty$,
we obtain the following:
\[
 \int_{\substack{\pi_1^{-1}(P)\\|z_1|\leq e^{-1}}}
 \bigl(s^{-1}\del_{h_0}(s),\,\del_{h_0}(s)\bigr)
 \cdot\rho_M
\leq
 C_3\cdot C_2(P)^{M'}
\]
By making $M$ larger,
we obtain the following estimate
for some constants $C_4$ and $M''$:
\[
 \int_{\substack{\pi_1^{-1}(P)\\|z_1|\leq e^{-1}}}
 \bigl(s^{-1}\del_{h_0}(s),
 s^{-1}\del_{h_0}(s)\bigr)
 \cdot\rho_M
\leq
 C_4\cdot C_2(P)^{M''}
\]

We can take $G_P\cdot dz_1/z_1$
for each $P$
with the following property 
for some $N_1$:
\[
 \delbar_{z_1}\left(G_P\frac{dz_1}{z_1}\right)
=R(h)_{|\pi_1^{-1}(P)},
\quad
 \bigl|G_P\bigr|_{h_0}\leq
 C_{5}\cdot \bigl(-\log|z_1|+C_2(P)\bigr)^{N_1}
\]
(We can use Lemma \ref{lem;07.2.6.11},
for example.)
We put
$\psi_P:=s^{-1}\del_{h_0} s_{|\pi_1^{-1}(P)}
 -G_P\cdot dz_1/z_1$.
Then, we obtain $\delbar_{z_1}\psi_P=0$
and the following inequality
for some $C$, $M_1$ and $N_3$,
which are independent of $P$:
\begin{equation}
 \label{eq;08.1.22.50}
 \int_{\pi_1^{-1}(P)}
 \bigl|\psi_P\bigr|_{h_0,g_{\poin}}^2
 \bigl(-\log|z_1|+C_2(P)\bigr)^{-M_1}
\dvol_{g_{\poin}}
 <C\cdot C_2(P)^{N_3}
\end{equation}
By the holomorphic property,
we obtain
$\bigl|\psi_P\bigr|_{h_0,g_{\poin}}\leq
 C_{10}\cdot\bigl(-\log|z_1|+C_2(P)\bigr)^{N_4}$
for some sufficiently large $N_4$ and $C_{10}$,
which are independent of the choice of $P$,
by the following standard argument.
Let $\vecu=(u_j)$ be a frame of
$\prolong{\End(E)}$ compatible with
the parabolic structure.
We have the expression
$\psi_P=\sum \psi_{P,j}\cdot u_{j|\pi_1^{-1}(P)}$.
We obtain the estimate for
the integrals of $|\psi_{P,j}|^2$
with appropriate weight from (\ref{eq;08.1.22.50}).
Since $\psi_{P,j}$ are holomorphic,
we obtain the estimate for the sup norms
of $|\psi_{P,j}|$.

Thus, we obtain the desired estimate for
$s^{-1}\del_{h_0}s$,
and the proof of Lemma \ref{lem;07.6.2.20}
is finished.
\hfill\qed

\chapter{Review on $\nbigr$-Modules,
 $\nbigr$-Triples and Variants}
\label{section;08.10.20.40}
In Section \ref{subsection;08.10.20.2},
we recall the basic facts on finiteness of
filtered rings and filtered modules
from \cite{kashiwara_text}
for reader's convenience.

In Sections 
\ref{subsection;08.10.18.100}--\ref{subsection;08.10.18.101},
we recall some basic property of
$\nbigr_X$-modules and their specializations.
In Section \ref{subsection;08.10.18.102},
we recall the notion of $\nbigr_X(\ast t)$-module.
We refer to \cite{sabbah2} for more details
and precision.
($\nbigr(\ast t)$-modules are called
$\widetilde{\nbigr}$-modules there.)
See \cite{saito1} for the original work due to M. Saito
on filtered $D$-modules.
See also \cite{mochi2}.

We would like to reduce the study
of wild objects to the study of tame objects.
One of the important tools is formal completion,
for which we give a preparation in
Section \ref{subsection;08.9.29.122}.
In Subsection \ref{subsection;07.11.6.10},
we review formal complex spaces.
We refer to \cite{banica},
 \cite{bingener} and \cite{krasnov}
for more details and precision.
Then, we recall basic facts on 
formal $D$-modules and formal $\nbigr$-modules
in Subsections \ref{subsection;07.7.30.4}
and \ref{subsection;08.9.29.211}.
We can obtain the claims in these subsections
by directly applying a general theory
explained in the appendix of \cite{kashiwara_text}.
We explain in Subsection
\ref{subsection;08.9.29.212}
how to use completions
in showing the strict $S$-decomposability
which is the main motivation for us
to consider formal $\nbigr$-modules.

In Section \ref{subsection;08.10.18.105},
we give some miscellaneous preliminary 
on $D$-modules for references in our argument.
We just indicate where they are used
in this monograph.
The contents in Subsections
\ref{subsection;08.1.21.1} and
\ref{subsection;08.9.29.200}
are implicitly used for the proof of theorems
in Chapter \ref{section;08.9.29.150}.
We prepare
the nearby cycle functor with ramification
and exponential twist in Subsection
\ref{subsection;08.10.25.1}
to state Theorem \ref{thm;08.3.6.2}.
Subsections \ref{subsection;08.1.19.20}
and \ref{subsection;08.9.29.210}
are preliminary for Section 
\ref{subsection;07.10.23.21}.

Section \ref{subsection;08.10.18.110}
is also complements.
We give in Subsection \ref{subsection;08.1.20.100}
a remark on push-forward of
$\nbigr$-modules via a ramified covering,
which will be implicitly used in
Part \ref{part;08.9.29.130}.
Subsection \ref{subsection;08.9.26.1}
is a preparation for 
the argument in Sections
\ref{subsection;08.9.28.2}--\ref{subsection;08.2.2.1}.

We review in Section \ref{subsection;08.1.15.20}
some basic property of
the sheaves of distributions with moderate growth
satisfying some continuity,
introduced by Sabbah in \cite{sabbah2},
which is one of the basic ingredients
for the definition of wild pure twistor $D$-modules.
It is used in Chapter \ref{section;07.10.12.10}
and Part \ref{part;08.9.29.130}.

In Sections 
\ref{subsection;08.10.18.112}--\ref{subsection;08.10.18.111},
we recall the notion of $\nbigr$-triples
and their specialization.
In Section \ref{subsection;08.10.18.113},
we recall the notion of $\nbigr(\ast t)$-triples.
We refer to \cite{sabbah2}
for more details and precise on their property.
($\nbigr(\ast t)$-triples are called
$\widetilde{\nbigr}$-triples there.)
In Section \ref{section;10.6.6.10},
we compare pure twistor structure
and $\nbigr$-triple in dimension $0$.
We also give, for a regular singular $\nbigr$-triple 
on a disc, a comparison of specializations
as $\nbigr$-triple and variation of twistor structure.

\section{Filtered rings}
\label{subsection;08.10.20.2}
We recall some general results
on filtered rings and their modules
given in Appendix of 
a nice textbook
\cite{kashiwara_text} due to Kashiwara
for reader's convenience,
which we refer to for more details and precision.

\subsection{Coherence}

Let $\nbiga$ be a sheaf of rings
on topological space $X$.
We have the standard notions of
{\em  locally finitely generated $\nbiga$-modules}
and {\em locally finitely presented $\nbiga$-modules}.
\begin{prop}[Proposition A.4 of
\cite{kashiwara_text}]
\label{prop;08.10.20.3}
Let $\nbigf$ be a locally finitely generated
$\nbiga$-module.
Let $x\in X$, and 
suppose that the stalk $\nbigf_x$ vanishes.
Then, there exists a neighbourhood $U$ of $x$
such that $\nbigf_{|U}=0$.
\hfill\qed
\end{prop}

An $\nbiga$-module $\nbigf$ is called
{\em pseudo-coherent} if the following holds:
\index{pseudo-coherent}
\begin{itemize}
\item
 Let $U$ be any open subset of $X$.
 Let $\nbigg$ be a locally finitely generated
 $\nbiga_{|U}$-submodule  of $\nbigf_{|U}$.
 Then, $\nbigg$ is locally finitely presented.
\end{itemize}
Any $\nbiga$-submodule of 
a pseudo-coherent $\nbiga$-module
is also pseudo-coherent.

An $\nbiga$-module $\nbigf$ is called
{\em coherent},
if it is pseudo-coherent and
locally finitely generated.
\index{coherent}
The category of coherent $\nbiga$-modules
is abelian.
It is stable under extensions.
(See Proposition A.6 of \cite{kashiwara_text}
for more details.)
The sheaf of rings $\nbiga$ is called {\em coherent}
if it is coherent as an $\nbiga$-module.
If $\nbiga$ is coherent,
an $\nbiga$-module is coherent
if and only if it is locally finitely presented.

An $\nbiga$-module $\nbigf$ is called
{\em Noetherian},
if the following holds:
\index{Noetherian}
\begin{itemize}
\item
 $\nbigf$ is a coherent $\nbiga$-module.
\item
 For any $x\in X$,
 the stalk $\nbigf_x$ is 
 a Noetherian $\nbiga_x$-module.
\item
 Let $U$ be any open subset of $X$.
 Let $\{\nbigg_i\}$ be any family of
 coherent $\nbiga_{|U}$-submodules
 of $\nbigf_{|U}$.
 Then, $\sum \nbigg_i$ is 
 a coherent $\nbiga_{|U}$-module.
\end{itemize}
A sheaf of rings $\nbiga$ is called
Noetherian, if it is Noetherian
as an $\nbiga$-module.

\subsection{Filtered rings and filtered modules}

A sheaf of rings $\nbiga$ is called filtered,
if it is equipped with an increasing sequence of
submodules
$\{F_n(\nbiga)\,|\,n=0,1,2,\ldots\}$
such that
(i) $\nbiga=\bigcup F_n(\nbiga)$,
(ii) $1\in F_0(\nbiga)$,
(iii) $F_n(\nbiga)\cdot F_m(\nbiga)
 \subset F_{n+m}(\nbiga)$.
If an $\nbiga$-module $\nbigm$
is equipped with an increasing sequence
of submodules
$\bigl\{F_n(\nbigm)\,\big|\,n=0,1,2,\ldots\bigr\}$
such that
$F_n(\nbiga)F_m(\nbigm)\subset
 F_{n+m}(\nbigm)$,
the pair $(\nbigm,F)$ is called
a filtered $\nbiga$-module.
We say that
$(\nbigm,F)$ is a {\em coherent} $\nbiga$-module,
if $\bigoplus F_n(\nbigm)$ is a coherent
$\bigoplus F_n(\nbiga)$-module.
In that case,
$F$ is called a {\em coherent filtration}.
\index{coherent filtration}
Similarly,
$F$ is called a {\em locally finitely generated},
if $\bigoplus F_n(\nbigm)$ is locally
finitely generated over $\bigoplus F_n(\nbiga)$.

The following proposition is useful
to check the coherence of a filtration.
\begin{prop}[Lemma A.27
 of \cite{kashiwara_text}]
\label{prop;08.10.20.21}
Let $\nbigm$ be a coherent $\nbiga$-module,
and $\nbign$ be an $\nbiga$-submodule.
A filtration $F(\nbign)$ of $\nbign$
is coherent, if and only if
it is locally finitely generated.
\hfill\qed
\end{prop}

Recall the following general and
fundamental theorem.

\begin{prop}[Theorem A.20 of 
\cite{kashiwara_text}]
\label{prop;08.10.20.1}
Let $(\nbiga,F)$ be a filtered ring.
Assume the following:
\begin{itemize}
\item
 $F_0(\nbiga)$ is Noetherian.
\item
 $\Gr_k(\nbiga)$ are coherent
 $F_0(\nbiga)$-module.
\item
 Let $m\in\seisuu_{>0}$,
 and let $U$ be an open subset of $X$.
 Let $\nbign$ be an $\nbiga$-submodule of
 $\nbiga^{\oplus\,m}_{|U}$
 such that
 $\nbign\cap F_k(\nbiga)^{\oplus\,m}$
 are coherent $F_0(\nbiga)_{|U}$-modules
 for all $k$.
 Then, $\nbign$ is locally finitely generated
 over $\nbiga$.
\end{itemize}
Then, $\nbiga$ is Noetherian.
\hfill\qed
\end{prop}

The following proposition is useful to check
coherence of an $\nbiga$-module.

\begin{prop}[Lemma A.22 of 
 \cite{kashiwara_text}]
\label{prop;08.10.20.31}
Let $(\nbiga,F)$ be a filtered ring
satisfying the conditions in 
Proposition {\rm\ref{prop;08.10.20.1}}.
Let $\nbigm$ be an $\nbiga$-module
such that
(i) locally finitely generated over $\nbiga$,
(ii) pseudo-coherent as an $F_0(\nbiga)$-module.
Then, $\nbigm$ is a coherent $\nbiga$-module.
\hfill\qed
\end{prop}

The following proposition is useful
to show a Noetherian property of
a sheaf of rings $\nbiga$.
\begin{prop}[Theorem A.29, Theorem A.31
 of \cite{kashiwara_text}]
\label{prop;08.10.20.30}
Let $(\nbiga,F)$ be a filtered ring on $X$.
Assume that $F_0(\nbiga)$
and $\Gr^F(\nbiga)$ are Noetherian rings,
and $\Gr^F_k(\nbiga)$ are locally 
finitely generated $F_0(\nbiga)$-modules
for any $k$.
Then, $\nbiga$ is Noetherian ring,
and it satisfies the conditions in 
Proposition {\rm\ref{prop;08.10.20.1}}.
Moreover,
the Rees ring $\bigoplus F_n(\nbiga)$
is also Noetherian.
\hfill\qed
\end{prop}

\section{$\nbigr$-modules}
\label{subsection;08.10.18.100}
\subsection{The sheaf of algebras $\nbigr_X$}
Let us recall the notion of $\nbigr_X$-modules
introduced in \cite{sabbah2}.
Let $X$ be a complex manifold.
We set $\nbigx:=\cnum_{\lambda}\times X$,
where $\cnum_{\lambda}$ denotes the complex line
with the coordinate $\lambda$.
Let $p:\nbigx\lrarr X$ denote the natural projection.
Let $\Theta_X$ denote the tangent sheaf of $X$.
\index{sheaf $\Theta_X$}
We have the subsheaf $\lambda\cdot p^{\ast}\Theta_X$
of $p^{\ast}\Theta_X$.
We denote $\lambda\cdot p^{\ast}\Theta_X$ 
by $\Theta_{\nbigx}$.
\index{sheaf $\Theta_{\nbigx}$}
Since we do not consider the tangent sheaf of $\nbigx$,
there is no risk of confusion.
Let $\nbigd_X$ denote the sheaf of 
differential operators on $X$.
\index{sheaf $\nbigd_X$}
We have the sheaf $p^{\ast}\nbigd_X=
\nbigo_{\nbigx}\otimes_{p^{-1}\nbigo_X}p^{-1}\nbigd_X$,
which is the sheaf of relative differential operators
of $\nbigx$ over $\cnum_{\lambda}$.
Let $\nbigr_X$ denote the sheaf of 
subalgebras of $p^{\ast}\nbigd_{X}$
generated by $\nbigo_{\nbigx}$ and $\Theta_{\nbigx}$.
\index{sheaf $\nbigr_X$}
The notions of left $\nbigr$-modules 
and right $\nbigr$-modules
are naturally defined
as in the case of $\nbigd_X$-modules.
{\em In this paper,
we usually consider left $\nbigr_X$-modules.}

It is equipped with the increasing filtration $F$
by the order of differential operators,
via which $\nbigr_X$ is a sheaf of filtered algebras.
The associated graded sheaf is
isomorphic to $p^{\ast}\Sym^{\bullet}\Theta_X$.
According to Propositions
\ref{prop;08.10.20.31}
and \ref{prop;08.10.20.30},
we obtain the following.
\begin{prop}
\mbox{{}}
\begin{itemize}
\item
The sheaf of algebras $\nbigr_X$ is Noetherian,
and it satisfies the conditions in
Proposition {\rm\ref{prop;08.10.20.1}}.
The Rees ring is also Noetherian.
\item
Let $M$ be an $\nbigr_X$-module
such that (i) pseudo-coherent as an $\nbigo_X$-module,
(ii) locally finitely generated as an $\nbigr_X$-module.
Then, $M$ is a coherent $\nbigr_X$-module.
\hfill\qed
\end{itemize}
\end{prop}

 Let $\lambda_0$  be any point of $\cnum$.
 We say that an $\nbigr_X$-module 
 $\nbigm$ is {\em strict} at $\lambda_0$,
 if the multiplication of $\lambda-\lambda_0$ 
 on $\nbigm$ is injective.
 If $\nbigm$ is strict at any $\lambda_0\in\cnum$,
 then it is called {\em strict}.
\index{strict}

\subsection{Left and right $\nbigr$-module}

Let $\omega_X$ denote the canonical line bundle of $X$,
i.e., it is the sheaf of holomorphic $(n,0)$-forms,
where $n=\dim X$.
\index{sheaf $\omega_X$}
Then, $p^{\ast}\omega_X$ is a right $\nbigr_X$-module,
on which the right $\nbigr_X$-action is given 
as the restriction of
the natural right $p^{\ast}\nbigd_X$-action.
Let 
$\omega_{\nbigx}:=\lambda^{-n}\cdot p^{\ast}\omega_X$
as the subsheaf of 
$p^{\ast}\omega_X\otimes
 \nbigo_{\nbigx}\bigl(\ast(\{0\}\times X)\bigr)$.
It is naturally a right $\nbigr_X$-module.
\index{sheaf $\omega_{\nbigx}$}

As in the case of $D_X$-modules,
for any left $\nbigr_X$-module $\nbigm$,
we have the natural right $\nbigr_X$-module structure on
$\nbigm^r:=\omega_{\nbigx}\otimes\nbigm$.
By this functor, the categories of
left $\nbigr_X$-modules and right $\nbigr_X$-modules 
are naturally equivalent.

\subsection{Holonomic $\nbigr$-module}
\label{subsection;08.12.4.1}

 Let $\nbigm$ be a coherent $\nbigr_X$-module
 on $U\subset\nbigx$.
 As in Section \ref{subsection;08.10.20.2},
 we have the notion of coherent filtration for $\nbigm$.
 If $\nbigm$ has a coherent filtration $F$,
we obtain the coherent 
$p^{\ast}\Sym \Theta_X$-module $\Gr^F\nbigm$.
It naturally induces the coherent sheaf
on $\cnum_{\lambda}\times T^{\ast}X$.
The support $Ch(\nbigm)$ is called 
by the characteristic variety.
We can show in a standard way
that $Ch(\nbigm)$ is independent of the choice of
a coherent filtration.
\index{characteristic variety $Ch(\nbigm)$}

Let $\nbigm$ be an $\nbigr_X$-module on 
 $U\subset\nbigx$.
We say that $\nbigm$ is {\em good},
if the following condition is satisfied:
\begin{itemize}
\item
Let $K$ be any compact subset of $X$
such that $\{\lambda_0\}\times K\subset U$.
Then, there exist a neighbourhood $U'$ of 
$\{\lambda_0\}\times K$ in $U$,
a finite filtration $F$ of $\nbigm_{|U'}$
such that
$\Gr^F(\nbigm_{|U'})$ has a coherent filtration.
\end{itemize}
If $\nbigm$ is good,
we have the globally well-defined characteristic variety
$Ch(\nbigm)\subset \cnum_{\lambda}\times T^{\ast}X$.
If there is a Lagrangian subvariety $L$ of $T^{\ast}X$
such that $Ch(\nbigm)\subset 
\cnum_{\lambda}\times L$,
then $\nbigm$ is called holonomic.
\index{holonomic $\nbigr_X$-module}
As a similar but weaker condition,
we say that
$\nbigm$ is locally good,
if there exists a covering $U=\bigcup U_j$
such that $\nbigm_{|U_j}$ are good.
The characteristic variety makes sense
for any locally good $\nbigr_X$-module.
Goodness is required for well behaviour
with respect to the push-forward
via a proper morphism.

\subsection{Push forward and pull back}
\label{subsection;08.12.16.3}
\index{pull back}
\index{push forward}

The push forward and pull back
of $\nbigr$-modules are defined 
as in the case of $D$-modules.
Their property is also similar.
Let $f:X\lrarr Y$ be a holomorphic map of complex manifolds.
The left and right $(\nbigr_X,f^{-1}\nbigr_Y)$-module
$\nbigr_{X\to Y}$ is given as follows:
\[
 \nbigr_{X\to Y}:=
 \nbigo_{\nbigx}\otimes_{f^{-1}\nbigo_{\nbigy}}
 f^{-1}\nbigr_Y 
\]
The pull back of a {\em left} $\nbigr_Y$-module $\nbign$
is given as follows:
\[
 f^{\dagger}\nbign:=
 \nbigr_{X\to Y}
 \otimes^L_{f^{-1}\nbigr_Y}
 f^{-1}\nbign
\in D^b(\nbigr_X)
\]
The push forward of 
a {\em right} $\nbigr_X$-module $\nbigm$ is given
as follows:
\[
 f_{\dagger}\nbigm
:=Rf_{!}\bigl(
 \nbigm\otimes^L_{\nbigr_X}
 \nbigr_{X\to Y}\bigr)
\in D^b(\nbigr_Y)
\]
Here $\otimes^L$ denotes the derived tensor product,
and $Rf_{!}$ denotes the ordinary push forward with compact support.
(See \cite{sabbah2} for more precision.)
The push forward of a left module is defined 
by using the equivalence of 
the categories of left and right $\nbigr$-modules
as mentioned above.
The $i$-th cohomology sheaves
of $f_{\dagger}\nbigm$ are denoted by
$f^i_{\dagger}\nbigm$.
We refer to Section 1.4 of \cite{sabbah2}
for more details and precision.
We mention one result:
if $\nbigm$ is good (holonomic),
then $f^i_{\dagger}\nbigm$ are also
good (holonomic),
which can be shown in a standard way
as in \cite{kashiwara_text}.

\vspace{.1in}
We recall the construction of
the push-forward for a left $\nbigr_X$-module
in the case that $Y$ is a point
for explanation in 
Subsection \ref{subsection;08.1.19.5}.
Let $\Omega^{r,s}_X$ denote 
the $C^{\infty}$-vector bundle
of $(r,s)$-forms on $X$.
Let $\nbigx^0:=\{0\}\times X$.
We obtain the $\cnum_{\lambda}$-holomorphic bundle
$\Omega^{r,s}_{\nbigx}:=
 p^{\ast}\Omega^{r,s}_X
 \otimes
 \nbigo_{\nbigx}(r\nbigx^0)$,
and the double complex
$\bigl(
 \Omega^{r,s}_{\nbigx},
 \del_X,\delbar_X \bigr)$.
The associated total complex is
denoted by
$(\Omega^{\bullet}_{\nbigx},d)$.

The left $\nbigr_X$-module structure on $\nbigm$
and the exterior derivative $d$
naturally induce the complex
$(\nbigm\otimes\Omega^{\bullet}_{\nbigx},d)$.
Then, $f_{\dagger}\nbigm$ is given
as follows in this case:
\[
 f_{\dagger}\nbigm=
 F_{!}\bigl(
 \nbigm\otimes
 _{\nbigo_{\nbigx}}
 \Omega^{\bullet}_{\nbigx}[\dim X]
 \bigr)
\]
Here, $F_!$ denotes the sheaf theoretic
push-forward with compact support
for the induced map
$F:\nbigx\lrarr\cnum_{\lambda}$.
A closed $q$-form $\omega$ on $X$ induces
a morphism of complexes
$L_{\omega}:
 \nbigm\otimes\Omega^{\bullet}_{\nbigx}
\lrarr
 \nbigm\otimes\Omega^{\bullet+q}_{\nbigx}$,
which induces
$L_{\omega}:
 f_{\dagger}^i\nbigm
\lrarr f_{\dagger}^{i+q}\nbigm$.

\section{Specialization of $\nbigr$-modules}
\label{subsection;08.10.18.101}
Sabbah introduced 
the notions of $V$-filtration,
nearby cycle functor and vanishing cycle functor
for $\nbigr$-modules in \cite{sabbah2},
which are natural generalization
of those for filtered $D$-modules 
in the original work by M. Saito
\cite{saito1}.
(See also \cite{mochi2}.)
\index{nearby cycle functor}
\index{vanishing cycle functor}

\subsection{Strict specializability}

Let $\cnum_t$ be a complex line 
with a coordinate $t$.
Let $X_0$ be a complex manifold.
We put $X:=X_0\times \cnum_t$.
We identify $X_0$ and $X_0\times\{0\}$.
Recall $\nbigx:=\cnum_{\lambda}\times X$
and $\nbigx_0:=\cnum_{\lambda}\times X_0$.
For any point $\lambda_0\in\cnum_{\lambda}$,
let $\nbigxzero$ denote a small neighbourhood 
of $\{\lambda_0\}\times X$.
We use the symbol $\nbigxzero_0$ in a similar meaning.
Let $p_{X_0}:X\lrarr X_0$ denote the projection.
Let $V_0\nbigr_X$ denote the sheaf of subalgebras
of $\nbigr_X$ generated by
$p_{X_0}^{\ast}\nbigr_{X_0}$ and $\deldel_t t$.

Let $\nbigm$ be a coherent $\nbigr_X$-module
on $\nbigx^{(\lambda_0)}$.
Let $\Vzero$ be an increasing filtration of
coherent $V_0\nbigr_X$-submodules 
indexed by $\real$.
We recall some conditions for $\Vzero$.

\begin{condition}
\mbox{{}}\label{condition;08.1.23.10}
\begin{itemize}
\item
 For any $a\in\real$, and for any $P\in\nbigxzero_0$,
 there exists $\epsilon>0$ such that
 $\Vzero_a(\nbigm)=\Vzero_{a+\epsilon}(\nbigm)$
 around $P$.
 Moreover,
 $\bigcup_{a\in\real}\Vzero_a(\nbigm)=\nbigm$.
\item
$t\cdot \Vzero_a(\nbigm)\subset \Vzero_{a-1}(\nbigm)$ 
 for any $a\in\real$,
and $t\cdot \Vzero_a(\nbigm)=\Vzero_{a-1}(\nbigm)$ 
 for any $a<0$.
\item
 $\deldel_t\cdot \Vzero_a(\nbigm)
 \subset \Vzero_{a+1}(\nbigm)$ 
 for any $a\in\real$,
and the induced morphisms
 $\deldel_t:\Gr^{\Vzero}_a(\nbigm)
 \lrarr \Gr^{\Vzero}_{a+1}(\nbigm)$
 are surjective for any $a>-1$.
\hfill\qed
\end{itemize}
\end{condition}

\begin{condition}
\mbox{{}} \label{condition;08.1.23.11}
Let $(\nbigm,\Vzero)$ be as in
Condition {\rm\ref{condition;08.1.23.10}}.
Recall that $\nbigm$ is called 
strictly specializable along $t$ at $\lambda_0$
with respect to $\Vzero$
if the following holds:
\index{strictly specializable}
\begin{itemize}
\item
 $\Gr^{\Vzero}_a(\nbigm)$
 is strict for each $a\in\real$,
 i.e.,
 the multiplication of $\lambda-\lambda_1$
 on $\Gr^{\Vzero}_a(\nbigm)$ is
 injective for any $\lambda_1\in\cnum$.
\item
 For any $P\in\nbigxzero_0$,
 there exists a finite subset 
 $\nbigk(a,\lambda_0,P)\subset \real\times\cnum$
 such that the action of 
 $\prod_{u\in\nbigk(a,\lambda_0)}
 \bigl(-\deldel_tt+\eigenmap(\lambda,u)\bigr)$
 on $\Gr^{\Vzero}_a(\nbigm)$ is nilpotent around $P$.
\end{itemize}
The union $\bigcup_{a\in\real}\nbigk(a,\lambda_0,P)$
is denoted by $\KMS(\nbigm,t,P)$.
We also put
$\nbigk(a,\lambda_0):=
 \bigcup_{P\in\nbigx_0}
 \nbigk(a,\lambda_0,P)$
and $\KMS(\nbigm,t):=
 \bigcup_{P\in\nbigx_0}\KMS(\nbigm,t,P)$.
\hfill\qed
\end{condition}

If $(\nbigm,\Vzero)$ is strictly specializable along $t$,
for $u\in\KMS(\nbigm,t)$,
we put 
\[
 \psizero_{t,u}(\nbigm):=
 \bigcup_N\Ker\Bigl(
 \bigl(-\deldel_t t+\eigenmap(\lambda,u)\bigr)^N
 \Bigr).
\]
Thus, we obtain an $\nbigr_{X_0}$-module
$\psizero_{t,u}(\nbigm)$,
and the decomposition
\begin{equation}
 \label{eq;08.10.15.15}
 \Gr^{\Vzero}_a(\nbigm)
=\bigoplus_{u\in\nbigk(a,\lambda_0)}
 \psizero_{t,u}(\nbigm).
\end{equation}
\index{functor $\psizero_{t,u}\nbigm$}

\begin{rem}
If we consider local issues,
we may assume that 
the set $\nbigk(a,\lambda_0)$ is finite.
We will often assume it implicitly.
\hfill\qed
\end{rem}

\begin{lem}
\label{lem;08.10.15.10}
Such a filtration $\Vzero$ is unique,
if it exists.
The index set
$\KMS(\nbigm,t)$ is also uniquely determined
if $\psizero_{t,u}(\nbigm)\neq 0$
for any $u\in\KMS(\nbigm,t)$.
\end{lem}
\pf
Although it is proved in \cite{sabbah2} and \cite{mochi2},
we give a sketch of the proof.
We have only to argue the issue
locally around any point of $\nbigxzero_0$.
Let $\Vzero$ and $\Vtilde^{(\lambda_0)}$
be filtrations satisfying the above conditions,
with the index sets
$\nbigk(a,\lambda_0)$
and $\nbigktilde(a,\lambda_0)$.
Let $a<0$.
We have $\Vzero_a\subset
 \Vtilde^{(\lambda_0)}_{b}$ for some $b$.
Note $\Vzero_{a-1}\subset \Vtilde^{(\lambda_0)}_{b-1}$.
Hence, we obtain the induced map
$ \Vzero_a\big/\Vzero_{a-1} 
\lrarr
 \Vtilde_{b}^{(\lambda_0)}
 \big/\Vtilde_{b-1}^{(\lambda_0)}$.
For any $a-1<e\leq a$,
the action of
$ \prod_{a-1<c\leq e}
 \prod_{u\in\nbigk(c,\lambda_0)}
 \bigl(-\deldel_tt+\eigenmap(\lambda,u)\bigr)$
is locally nilpotent
on $\Vzero_e\bigl/\Vzero_{a-1}$.
For any $b-1<e\leq b$,
the action of
$ \prod_{b-1<c\leq e}
 \prod_{u\in\nbigktilde(c,\lambda_0)}
 \bigl(-\deldel_tt+\eigenmap(\lambda,u)\bigr)$
is locally nilpotent
on $\Vtilde^{(\lambda_0)}_e
 \bigl/\Vtilde^{(\lambda_0)}_{b-1}$.
For any $a-1<c\leq a$
and $b-1<e\leq b$,
we have the induced map
\[
 \Phi_{c,e}:
 \Vzero_{c}/\Vzero_{a-1}
\lrarr
 \Vtilde^{(\lambda_0)}_b/
 \Vtilde^{(\lambda_0)}_{e}.
\]
We obtain the vanishing of
the restriction of $\Phi_{c,e}$ to
$\nbigx_0^{(\lambda_0)}-(\{\lambda_0\}\times X_0)$ 
if $c<e$.
By using the strictness of 
$\Gr^{\Vzero}$ and $\Gr^{\Vtilde^{(\lambda_0)}}$,
we obtain the vanishing of $\Phi_{c,e}$ if $c<e$.
Hence, we obtain that 
$\Vzero_a\subset \Vtilde^{(\lambda_0)}_a$.
Moreover, we obtain that
$\nbigk(a,\lambda_0)
=\nbigktilde(a,\lambda_0)$
if $-\deldel_tt+\eigenmap(\lambda_0,u)$ has 
non-trivial kernel for each
$u\in \nbigk(a,\lambda_0)
 \cup\nbigktilde(a,\lambda_0)$.
\hfill\qed

\begin{lem} 
\label{lem;9.17.10}
Assume that $\nbigm$ on $\nbigxzero$
is strictly specializable along $t$ at $\lambda_0$
with the index set $\KMS(\nbigm,t)$.
For simplicity, $\nbigxzero$ is assumed to be
the product of $X$ and a neighbourhood of $\lambda_0$,
and $\KMS(\nbigm)$ is assumed to be finite.

Assume that $|\lambda_1-\lambda_0|$ is sufficiently small.
Then $\nbigm$ is strictly specializable 
along $t$ at $\lambda_1$
with the index set $\KMS(\nbigm,t)$.
Moreover,
on a neighbourhood
 $\nbigx^{(\lambda_1)}\subset \nbigxzero$
 of $\{\lambda_1\}\times X$,
we have natural isomorphisms
\[
 \psizero_{t,u}(\nbigm)_{|\nbigx^{(\lambda_1)}}
\simeq
 \psi^{(\lambda_1)}_{t,u}(\nbigm).
\]
\end{lem}
\pf
We will construct the filtration $V^{(\lambda_1)}$
of $\nbigm$ on $\nbigx^{(\lambda_1)}$.
In the following argument,
we restrict $\nbigm$ and the associated modules
to $\nbigx^{(\lambda_1)}$.
For any real number $d\in\real$,
we consider 
$ S(d):=
 \bigl\{c\in\real\,\big|\,
 \exists u\in \nbigk(c,\lambdazero),\,\,
 \paramap(\lambda_1,u)=d \bigr\}$.
Since $|\lambda_1-\lambda_0|$ is small,
we have $\big|S(d)\big|\leq 1$.
First, let us consider the case
$S(d)=\{c\}$.
Let $\pi_c:\Vzero_c(\nbigm)\lrarr \Gr^{\Vzero}_c(\nbigm)$
be the projection,
and we define
\[
 V^{(\lamda_1)}_d(\nbigm):=
 \pi_c^{-1}\Bigl(
 \bigoplus_{\substack{
 u\in\nbigk(c,\lamda_0)\\
 \paramap(\lambda_1,u)\leq d
 }}
 \psizero_{t,u}(\nbigm)
 \Bigr).
\]
Let us consider the case $S(d)=\emptyset$.
In that case,
we put $d_0:=\max\{d'\leq d\,|\,S(d')\neq \emptyset\}$,
and we define
$V_{d}^{(\lambda_1)}(\nbigm):=
 V_{d_0}^{(\lambda_1)}(\nbigm)$.
Then it is easy to check that
$V^{(\lambda_1)}$ is the desired filtration.
\hfill\qed

\begin{lem}
\label{lem;08.10.15.20}
Assume that
$\nbigm$ is strictly specializable along $t$
at $\lambda_0$ with a filtration $\Vzero$.
\begin{itemize}
\item
The multiplication
$t:\Vzero_{<0}(\nbigm)\lrarr \Vzero_{<0}(\nbigm)$
is injective.
\item
The induced maps
$\deldel_t:\Gr^{\Vzero}_a(\nbigm)
\lrarr \Gr^{\Vzero}_{a+1}(\nbigm)$
are isomorphisms
for any $a>-1$.
\end{itemize}
\end{lem}
\pf
Let us show the first claim.
Let $f\in \Vzero_{<0}$ such that $t\cdot f=0$.
Assume $f\neq 0$,
and we will deduce a contradiction.
Note that $\deldel_t^Nf\neq 0$ for any $N$,
and that
$\Vzero_{<0}\cap\nbigr_X\cdot f$
is $V_0\nbigr_X$-coherent.
Hence, there exists an $N$ such that
$\deldel_t^Nf\not\in\Vzero_{<0}$.
It implies the existence of $a\in\real_{<0}$
such that $f\in \Vzero_a\setminus \Vzero_{<a}$.
Hence, we obtain a non-zero element
$[f]\in \Gr^{\Vzero}_a(\nbigm)$.
Because $t\cdot f=0$,
we have $\deldel_tt[f]=0$.
Hence, $[f]\not\in \psizero_{t,u}(\nbigm)$
for any $u\in\nbigk(a,\lambda_0)$,
which is a contradiction.
Therefore, we can conclude that $f=0$.

As for the second claim,
it is easy to check that
the endomorphisms $t\cdot \deldel_t$
of $\Gr^{\Vzero}_a(\nbigm)$
are injective for any $a>-1$.
Thus, we are done.
\hfill\qed

\subsection{Compatibility of 
a morphism and $V$-filtrations}

\begin{prop}
 \label{prop;9.17.30}
Let $\nbigm$ and $\nbign$ be $\nbigr_X$-modules on 
$\nbigxzero$,
which are strictly specializable
along $t$ at $\lambda_0$.
Let $\phi$ be a morphism of $\nbigm$ to $\nbign$.
\begin{enumerate}
\item
 $\phi\bigl(\Vzero_c(\nbigm)\bigr)$
 is contained in $\Vzero_c(\nbign)$.
\item
 Assume that
 $\Gr^{\Vzero}_c(\phi):
 \Gr^{\Vzero}_c(\nbigm)\lrarr \Gr^{\Vzero}_c(\nbign)$
 is strict, \index{strict}
 namely,
 $\Cok(\Gr^{\Vzero}_c(\phi))$ is strict.
Then $\phi$ is strict with respect to the filtrations
$\Vzero(\nbigm)$ and $\Vzero(\nbign)$,
i.e.,
$\Vzero_c(\nbigm)\cap \Image(\phi)=
\phi\bigl(\Vzero_c(\nbign)\bigr)$.
\item
Under the above assumption,
$\Ker(\phi)$, $\Cok(\phi)$ and $\Image(\phi)$
are strictly specializable along $t$ at $\lambda_0$
with the induced filtrations,
and we have the natural isomorphisms
\[
 \Ker(\Gr^{\Vzero}_a(\phi))
\simeq
 \Gr^{\Vzero}_a(\Ker(\phi)),
\quad
 \Cok(\Gr^{\Vzero}_a(\phi))
\simeq
 \Gr^{\Vzero}_a(\Cok(\phi))
\]
\[
 \Image(\Gr^{\Vzero}_a(\phi))
\simeq
 \Gr^{\Vzero}_a(\Image(\phi)).
\]
\end{enumerate}
\end{prop}
\pf
The first claim can be shown 
by the argument in Lemma \ref{lem;08.10.15.10}.
To show the second claim,
let us consider the induced morphism
$\phi':
 \Vzero_d/\Vzero_{<c}(\nbigm)
 \lrarr\Vzero_d/\Vzero_{<c}(\nbign)$.
\begin{lem}
\label{lem;08.10.15.12}
$\phi'$ is strict with respect to the induced filtration
$\Vzero$,
i.e., the following holds for any $c\leq d'\leq d$:
\[
 \Image(\phi')\cap 
\Bigl(
\Vzero_{d'}/\Vzero_{<c}(\nbign)\Bigr)
=\phi'\Bigl(\Vzero_{d'}/\Vzero_{<c}(\nbigm)
 \Bigr).
\]
\end{lem}
\pf
We have only to consider the case in which
$\nbigxzero$ is the product of $X$
and a neighbourhood $U(\lambda_0)$ of $\lambda_0$
in $\cnum_{\lambda}$.
There exists a dense subset 
$U^{\ast}(\lambda_0)\subset U(\lambda_0)$
such that
$\eigenmap(\lambda):
 \KMS(\nbigm,t)\lrarr\cnum$
is injective for each $\lambda\in U^{\ast}(\lambda_0)$.
We set 
$\nbigx_0^{(\lambda_0)\ast}:=
 U^{\ast}(\lambda_0)\times X_0$.
We have the decomposition:
\[
 \Vzero_d/\Vzero_{<c}(\nbigm)_{|\nbigx_0^{(\lambda_0)\ast}}
=\bigoplus_{
 \substack{u\in \KMS(\nbigm,t)\\
 c\leq \paramap(\lambda_0,u)\leq d
 }}
 \psizero_{u}(\nbigm)_{|\nbigx_0^{(\lambda_0)\ast}}
\]
It gives a splitting of the induced filtration $\Vzero$
on $\Vzero_d/\Vzero_{<c}(\nbigm)_{|\nbigx_0^{(\lambda_0)\ast}}$.
We have similar decomposition for
$ \Vzero_d/\Vzero_{<c}(\nbign)_{|\nbigx_0^{(\lambda_0)\ast}}$.
Since $\phi'_{|\nbigx_0^{(\lambda_0)\ast}}$
preserves the decompositions,
it is strict with respect to $\Vzero$.
By using the strictness of
$\Gr^{\Vzero}(\nbigm)$,
$\Gr^{\Vzero}(\nbign)$
and $\Cok\Gr^{\Vzero}(\phi)$,
we obtain that $\phi'$ is strict
with respect to the filtrations $\Vzero$.
\hfill\qed

\begin{lem}
\label{lem;08.10.15.11}
Let $(\lambda_0,P)\in\nbigxzero$.
For any $c$,
there exists $e<0$
such that 
$\Vzero_e(\nbign)\cap \Image\phi
\subset
\phi\bigl(\Vzero_{c}(\nbigm)\bigr)$
around $(\lambda_0,P)$.
\end{lem}
\pf
Let $\nbigq_N$ denote the $V_0\nbigr_X$-submodule
of $\Vzero_{-1}(\nbign)$ 
which consists of the sections $a$
such that $t^Na\in
 \Image(\phi)\cap\Vzero_{-N-1}(\nbign)$.
There exists a large $N_1$ such that
$\nbigq_{N}=\nbigq_{N+1}$ 
for any $N\geq N_1$
around $(\lambda_0,P)$.
We set $\nbigq:=\nbigq_{N_1}$.
We have
$t^N\nbigq=\Image\phi\cap\Vzero_{-N-1}(\nbign)$
for any $N\geq N_1$ by our choice.
There exists a large number $f$ such that
$t^{N_1}\nbigq\subset
\phi\bigl(\Vzero_{f}(\nbigm)\bigr)$.
Take $N_2$ such that $f-N_2\leq c$ and $N_2\geq N_1$.
Then, we obtain
$\Vzero_{-1-N_1-N_2}(\nbign)\cap\Image\phi
\subset \phi\bigl(\Vzero_{c}(\nbigm)\bigr)$.
\hfill\qed

\vspace{.1in}

Let $h$ be an element of $\Vzero_d(\nbigm)$
such that $\phi(h)\in\Vzero_c(\nbign)$
for some $c<d$.
Let $e$ be as in Lemma \ref{lem;08.10.15.11}
for the $c$.
By using Lemma \ref{lem;08.10.15.12} inductively,
there exists an element 
$h_1\in\Vzero_c(\nbigm)$
such that
$\phi(h-h_1)\in \Vzero_{e}(\nbign)$.
By the choice of $e$,
there exists $h_2\in \Vzero_c(\nbigm)$
such that 
$\phi(h_2)=\phi(h-h_1)$,
i.e.,
$\phi(h)=\phi(h_1+h_2)$,
which means the strictness of $\phi$
with respect to the filtrations $\Vzero$.
Thus, we obtain the second claim of
Proposition \ref{prop;9.17.30}.

The third claim easily follows from the second claim.
Note we use Lemma \ref{lem;08.10.15.20}
to show $V_{a-1}(\Ker(\phi))=t\cdot V_a(\Ker(\phi))$
for $a<0$.
\hfill\qed

\subsection{The functor $\psi_{t,u}$}

Let $\nbigm$ be an $\nbigr_X$-module on $\nbigx$.
We say that $\nbigm$ is strictly specializable along $t$
if it is strictly specializable along $t$ at any $\lambda_0$.
If $\nbigm$ is strictly specializable along $t$,
according to Lemma \ref{lem;9.17.10},
$\bigl\{\psi^{(\lambda_0)}_{t,u}\nbigm\,
 \big|\,\lambda_0\in\cnum\bigr\}$
determines the globally defined $\nbigr_{X_0}$-module
$\psi_{t,u}(\nbigm)$
for each $u\in\KMS(\nbigm,t)$.
\index{functor $\psi_{t,u}\nbigm$}

Let $\nbigm$ and $\nbign$ be 
$\nbigr_X$-modules on $\nbigx$,
which are strictly specializable along $t$.
Let $\phi:\nbigm\lrarr\nbign$ be a morphism of $\nbigr_X$-modules.
According to Proposition \ref{prop;9.17.30},
the morphism $\phi$ preserves the $V$-filtrations 
$\Vzero$
at any $\lambda_0\in\cnum_{\lambda}$.
We obtain induced morphisms
$\Gr^{\Vzero}(\phi):
 \Gr^{\Vzero}(\nbigm)\lrarr\Gr^{\Vzero}(\nbign)$.
Since the decomposition (\ref{eq;08.10.15.15})
is obtained as a generalized eigen decomposition,
we obtain induced morphisms
$\psi_{t,u}^{(\lambda_0)}(\phi):
 \psi^{(\lambda_0)}_{t,u}(\nbigm)
 \lrarr\psi^{(\lambda_0)}_{t,u}(\nbign)$.
We can glue them to obtain a morphism
$\psi_{t,u}(\phi):\psi_{t,u}(\nbigm)
 \lrarr\psi_{t,u}(\nbign)$
of $\nbigr_{X_0}$-modules on $\nbigx_0$,
which is clear from the construction
of $V^{(\lambda_1)}$ from $V^{(\lambda_0)}$
given in Lemma \ref{lem;9.17.10}.

\vspace{.1in}

Let $\vecdelta_0$ denote 
the element $(1,0)\in\real\times\cnum$.
If an $\nbigr_X$-module $\nbigm$ is 
strictly specializable along $t$,
we have the naturally induced morphisms 
of $\nbigr_{X_0}$-modules
\[
 t:\psi_{t,u}\nbigm
   \lrarr
   \psi_{t,u-\vecdelta_0}\nbigm,
\quad\quad
 \deldel_t:
 \psi_{t,u}\nbigm
 \lrarr
 \psi_{t,u+\vecdelta_0}\nbigm.
\]
In particular,
we put 
\[
 \can=-\deldel_t:\psi_{t,-\vecdelta_0}\nbigm
 \lrarr\psi_{t,0}\nbigm,\\
 \quad\quad
 \var=t:\psi_{t,0}\nbigm\lrarr
 \psi_{t,-\vecdelta_0}\nbigm.
\]
\index{morphism $\can$}
\index{morphism $\var$}

\begin{rem}
If we consider right $\nbigr$-modules,
$\can$ is given by $\deldel_t$.
\hfill\qed
\end{rem}

\subsection{Strictly specializable morphisms}

Let $\nbigm$ and $\nbign$ be 
$\nbigr_X$-modules on $\nbigx$
which are strictly specializable along $t$.
A morphism $\phi:\nbigm\lrarr\nbign$ is called
strictly specializable,
if the induced morphisms
$\psi_{t,u}(\phi):
 \psi_{t,u}\nbigm\lrarr\psi_{t,u}\nbign$ are strict
for any $u\in\real\times\cnum$,
i.e.,
the cokernel $\Cok\bigl(\psi_{t,u}(\phi)\bigr)$ is strict.
\index{strictly specializable morphism}
We obtain the category 
$\nbigs^2(X,t)$ with
strictly specializable $\nbigr_X$-modules along $t$
and strictly specializable morphisms.
Note that it is abelian
according to Proposition \ref{prop;9.17.30}.
\index{category $\nbigs^2(X,t)$}

\begin{lem}
Assume that an $\nbigr_X$-module
$\nbigm$ is strictly specializable along $t$.
\begin{itemize}
\item
 If we have a decomposition 
 $\nbigm=\nbigm_1\oplus\nbigm_2$,
 then $\nbigm_i$ $(i=1,2)$ are strictly specializable
 along $t$.
\item
 Assume that $\nbigm$ is supported in $\nbigx_0$.
 Then we have $V^{(\lambda_0)}_{<0}\nbigm=0$
 for any $\lambda_0\in\cnum_{\lambda}$.
 We also have $\psi_{t,u}\nbigm=0$,
if $u$ does not contained in $\seisuu_{\geq\,0}\times\{0\}$.
\end{itemize}
\end{lem}
\pf
Let us show the first claim.
Since $\id\oplus 0$ preserves the filtration $\Vzero$ 
on $\nbigxzero$, the decomposition and 
$\Vzero$ are compatible.
Then, the first claim is clear.
See \cite{sabbah2} or Proposition 14.42 of \cite{mochi2}
for the second claim.
\hfill\qed

\vspace{.1in}

Let $\nbigs^2_{X_0}(X,t)$ denote the subcategory
of $\nbigs^2(X,t)$,
whose objects have the supports contained in $\nbigx_0$.
\begin{cor}
We have the equivalence
$ \nbigs^{2}_{X_0}(X,t)
\simeq
 \bigl(\mbox{\rm strict}
 \nbigr_{X_0}\mbox{\rm -modules}\bigr)$.
\hfill\qed
\end{cor}

\subsection{Strict $S$-decomposability}

The following condition essentially appeared
in \cite{saito1}.
\begin{condition}
\label{condition;08.1.23.12}
Let $\nbigm$ on $\nbigxzero$ be 
strictly specializable along $t$
at $\lambda_0$ with the filtration $\Vzero$.
Recall that $\nbigm$ is 
called strictly $S$-decomposable along $t$ at $\lambda_0$,
if moreover 
 $\psizero_{t,0}(\nbigm)
 =\Image (\can)\oplus \Ker(\var)$ holds.

We say that $\nbigm$ on $\nbigx$ is 
strictly $S$-decomposable along $t$,
if it is strictly $S$-decomposable along $t$
at any $\lambda_0$.
\hfill\qed
\index{strictly $S$-decomposable}
\end{condition}

\begin{prop} 
\label{prop;a11.23.1}
Let $\nbigm$ be an $\nbigr_X$-module
which is strictly specializable along $t$.
\begin{enumerate}
\item \label{number;10.2.1}
 The following conditions are equivalent:
   \begin{itemize}
  \item
    $\var:\psi_{t,0}\nbigm\lrarr\psi_{t,-\vecdelta_0}\nbigm$
    is injective.
  \item
    Let $\nbigm'$ be a submodule of $\nbigm$ such that
    the support of $\nbigm'$ is contained in $\nbigx_0$.
    Then $\nbigm'=0$.
  \item
    Let $\nbigm'$ be a submodule of $\nbigm$ such that
    the support of $\nbigm'$ is contained in $\nbigx_0$.
    Assume that $\nbigm'\in \nbigs^2(X,t)$.
    Then $\nbigm'=0$.
   \end{itemize}
 \item \label{number;10.2.2}
  Assume $\can:\psi_{t,-\vecdelta_0}\nbigm\lrarr\psi_{t,0}\nbigm$ is
      surjective.
 Let $\nbigm''\in \nbigs^2(X,t)$ be a quotient of $\nbigm$ such that
 the support of $\nbigm''$ is contained in $\nbigx_0$.
 Then $\nbigm''=0$.
 \item \label{number;05.1.15.100}
$\nbigm$ is strictly $S$-decomposable along $t$
if and only if the following holds:
 \begin{itemize}
  \item We have the decomposition $\nbigm=\nbigm'\oplus\nbigm''$.
 Here the support of $\nbigm''$ is contained in $X_0$,
 and $\nbigm'$ has neither submodules or quotients
 contained in $\nbigs^2(X,t)$
 whose support is contained in $\nbigx_0$.
 \end{itemize}
 \end{enumerate}
\end{prop}
\pf
See \cite{sabbah2}.
(See also \cite{saito1}.)
\hfill\qed

\subsection{The functors
 $\tildepsizero_{t,u}$ and $\tildepsi_{t,u}$}

Let $\nbigm$ be an $\nbigr_{X}$-module,
which is strictly specializable along $t$.
For any $u=(a,\alpha)\in\KMS(\nbigm,t)$
such that $u\not\in\seisuu_{\geq\,0}\times \{0\}$,
the $\nbigr_{X_0}$-module $\tildepsizero_{t,u}(\nbigm)$
\index{functor $\tildepsizero_{t,u}(\nbigm)$}
is defined in \cite{mochi2}
as the inductive limit of
$t:\psizero_{t,u-N\cdot\vecdelta_0}(\nbigm)
\lrarr
 \psizero_{t,u-(N+1)\cdot\vecdelta_0}(\nbigm)$.
It is equal to 
$\psizero_{t,u}\bigl(\nbigm(\ast t)\bigr)$.
(See \cite{sabbah2} or 
Section \ref{subsection;08.9.4.11}.)
It is easy to observe that
$\psitilde^{(\lambda_0)}_{t,u}
 (\nbigm)_{|\nbigx^{(\lambda_1)}}
\simeq
 \psitilde^{(\lambda_1)}_{t,u}(\nbigm)$
in the situation of Lemma \ref{lem;9.17.10}.
Hence, we obtain 
the $\nbigr_{X_0}$-module 
$\tildepsi_{t,u}(\nbigm)$
on $\nbigx_0$
as the gluing of 
$\tildepsizero_{t,u}(\nbigm)$.
\index{functor $\tildepsi_{t,u}(\nbigm)$}
We have a natural isomorphism
$\tildepsi_{t,u}(\nbigm)\simeq
 \lim_{N\to\infty}\psi_{t,u-N\cdot\vecdelta_0}(\nbigm)$.
We have the canonical inclusions
$\psizero_{t,u}(\nbigm)\lrarr \tildepsizero_{t,u}(\nbigm)$
and
$\psi_{t,u}(\nbigm)\lrarr \tildepsi_{t,u}(\nbigm)$
(Lemma 14.52 of \cite{mochi2}).

\begin{rem}
If $u$ is contained in 
$\real_{<0}\times\{0\}\subset\real\times\cnum$,
we have the canonical morphism
$\tildepsi_{t,u}(\nbigm)\simeq \psi_{t,u}(\nbigm)$.
We will not distinguish them.
\hfill\qed
\end{rem}

\begin{rem}
In the case $u=(0,0)\in\real\times\cnum$
(it is also denoted by $0$),
we also use the notation $\phi_{t,0}$ 
instead of $\psi_{t,0}$.
It is called the vanishing cycle functor.
\hfill\qed
\end{rem}

\begin{rem}
Let $u$ be an element of $\KMS(\nbigm,t)$.
Let us consider the following set:
\[
 S(u):=
 \bigl\{
 \lambda\in\cnum^{\ast}\,\big|\,
 \exists b\in \seisuu_{\geq 0},\,\,\mbox{\rm s.t.}\,\,
 \eigenmap\bigl(\lambda,u-b\cdot\vecdelta_0\bigr)=0,\,\,
 \paramap\bigl(\lambda_0,u-b\cdot\vecdelta_0\bigr)\geq 0
 \bigr\}
\]
Then, $S(u)$ is discrete in $\cnum_{\lambda}$.
For any $\lambda_0\in\cnum^{\ast}-S(u)$,
we have the canonical isomorphism
$\psizero_{t,u}(\nbigm)\simeq \tildepsizero_{t,u}(\nbigm)$.
\hfill\qed
\end{rem}

\subsection{The general case}

Let $Y$ be a complex manifold,
and let $\nbigm$ be an $\nbigr_{Y}$-module
on $\cnum_{\lambda}\times Y$.
Let $U$ be an open subset of $Y$,
and let $f$ be a holomorphic function on $U$.
Then we obtain the $\nbigr_{U\times\cnum_t}$-module
$\iota_{\dagger} \nbigm_{|\cnum_{\lambda}\times U}$
on $\cnum_{\lambda}\times U\times\cnum_t$,
where $\iota:U\lrarr U\times\cnum_t$ 
denotes the graph embedding.

\begin{df} \mbox{{}}
\begin{itemize}
\item
Let $f$ and $U$ be as above.
An $\nbigr_Y$-module $\nbigm$ is called 
strictly specializable (resp. $S$-decomposable) along $f$,
if $i_{f\,\dagger}\nbigm_{|U\times\cnum_t}$
is strictly specializable (resp. $S$-decomposable) along $t$.
\item 
An $\nbigr_Y$-module $\nbigm$ is called strictly specializable
(resp. $S$-decomposable),
if it is strictly specializable (resp. $S$-decomposable)
along any holomorphic function
defined on any open subset of $Y$.
\index{strictly specializable}
\index{strictly $S$-decomposable}
\hfill\qed
\end{itemize}
\end{df}

If $\nbigm$ is strictly specializable along 
a holomorphic function $f$ on $Y$,
we define
\[
\begin{array}{l}
\psi_{f,u}(\nbigm):=
\psi_{t,u}\bigl(\iota_{\dagger}\nbigm\bigr),
\quad\quad
\psitilde_{f,u}(\nbigm):=
\psitilde_{t,u}\bigl(\iota_{\dagger}\nbigm\bigr).
\end{array}
\]
\index{functor $\psi_{f,u}(\nbigm)$}
\index{functor $\psitilde_{f,u}(\nbigm)$}

Recall the compatibility of
push-forward and specialization due to Sabbah 
(Theorem 3.3.15 in \cite{sabbah2}).
See also Saito's work \cite{saito1}.
\begin{lem}
Let $F:Y\lrarr Z$ be a proper morphism.
Let $g$ be a holomorphic function on $Z$,
and we set $\gtilde:=F^{\ast}g$.
Let $\nbigm$ be an $\nbigr_Y$-module 
which is strictly specializable along $\gtilde$.
If $F^i_{0\,\dagger}\psitilde_{\gtilde,u}(\nbigm)$ 
are strict
for any $u\in \real\times\cnum$ and any $i$,
then $F^i_{\dagger}(\nbigm)$ are 
also strictly specializable along $g$,
and $\psitilde_{g,u}F^i_{\dagger}(\nbigm)
 \simeq F^i_{0\dagger}\psitilde_{\gtilde,u}(\nbigm)$.
\hfill\qed
\end{lem}

\subsection{Strict support}

Let $Z$ be an irreducible subvariety of $X$.
Let $\nbigm$ be a strictly $S$-decomposable
$\nbigr_X$-module,
whose support is contained in $Z$.
If there are no $\nbigr_X$-submodule
whose support is strictly smaller than $Z$,
then $Z$ is called the strict support of $\nbigm$.

We obtain the following lemma by
using the decomposition given in 
Proposition \ref{prop;a11.23.1}
and an easy inductive argument.
(See \cite{saito1} and \cite{sabbah2}.)
\begin{lem}
Let $\nbigm$ be 
a holonomic strictly $S$-decomposable $\nbigr_X$-module.
Then we have the decomposition 
$\nbigm=\bigoplus_Z\nbigm_Z$ satisfying the following:
\begin{itemize}
\item
$\nbigm_Z$ are
strictly $S$-decomposable $\nbigr_X$-module.
\item
The strict support of $\nbigm_Z$ is $Z$.
\end{itemize}
It is called the decomposition by the strict supports.
\hfill\qed
\end{lem}

\section{$\nbigr_X(\ast t)$-modules}
\label{subsection;08.10.18.102}

\subsection{Strict specializable $\nbigr_X(\ast t)$-modules}
\label{subsection;08.9.4.11}

Let $X=X_0\times\cnum_t$.
We identify $X_0$ with the $\{t=0\}$ in $X$.
We put $\nbigx:=\cnum_{\lambda}\times X$
and $\nbigx_0:=\cnum_{\lambda}\times X_0$.
Let $\nbigx^{(\lambda_0)}$ denote
a small neighbourhood of $\{\lambda_0\}\times X$
in $\cnum_{\lambda}\times X$.
We use the symbol $\nbigxzero_0$
in a similar meaning.
For simplicity,
we assume that $\nbigx^{(\lambda_0)}$
is the product of a neighbourhood of $\lambda_0$
and $X$.
We use the symbols
$\nbigx_0$ and $\nbigx_0^{(\lambda_0)}$
in similar meanings.

Let $\nbigr_X(\ast t)$ denote 
$\nbigr_X\otimes\nbigo_{\nbigx}(\ast t)$,
which is naturally the sheaf of algebras.
\index{sheaf $\nbigr_X(\ast t)$}
We refer to \cite{sabbah2} and \cite{sabbah5}
for general and foundational properties
of $\nbigr_X(\ast t)$-modules.
We have the standard correspondence between
the left $\nbigr_X(\ast t)$-modules
and the right $\nbigr_X(\ast t)$-modules
given by $\nbigm$ and 
$\nbigm\otimes_{\nbigo_{\nbigx}(\ast t)}
 \omega_{\nbigx}(\ast t)$.
In this paper,
we consider {\em left} $\nbigr_X(\ast t)$-modules.
Coherence, holonomicity and strictness for
$\nbigr_X(\ast t)$-modules
are defined as in the case of $\nbigr_X$-modules.

\vspace{.1in}
We recall the notion of strict specializability in this situation.

\begin{df}
\label{df;07.10.24.1}
Let $\nbigm$ be a left $\nbigr_X(\ast t)$-module,
which is coherent, holonomic and strict.
It is called strictly specializable along $t$ at $\lambda_0$,
if we have an increasing filtration $\Vzero$ of
$\nbigm_{|\nbigx^{(\lambda_0)}}$
by coherent $V_0\nbigr_X$-modules
indexed by $\real$, such that 
the following holds:
\begin{itemize}
\item
For any $a\in\real$
and $P\in\nbigx_0$,
there exists $\epsilon>0$
such that 
$\Vzero_a(\nbigm)=\Vzero_{a+\epsilon}(\nbigm)$
around $P$.
Moreover,
$\bigcup_{a}\Vzero_a(\nbigm)=\nbigm$.
\item
 $\Gr^{\Vzero}(\nbigm)$ is a strict
 $\nbigr_{X_0}$-module,
 i.e., the multiplication of
 $\lambda-\lambda_1$ is injective
 for any $\lambda_1\in\cnum$.
\item
 $t\cdot \Vzero_a(\nbigm)=\Vzero_{a-1}(\nbigm)$
and
 $\deldel_t\cdot \Vzero_{a}\subset \Vzero_{a+1}(\nbigm)$
 for any $a\in\real$.
\item
 For each $a\in\real$ and $P\in\nbigx_0$,
 there is a finite subset
 $\nbigk(a,\lambda_0,P)\subset\real\times\cnum$
 such that
 the action of
$ \prod_{u\in\nbigk(a,\lambda_0)}
 \bigl(-\deldel_t t+\eigenmap(\lambda,u)\bigr)$
 on $\Gr^{\Vzero}_a(\nbigm)$ is nilpotent
on a neighbourhood of $P$.
\end{itemize}
We say that $\nbigm$ is strictly specializable along $t$,
if it is strictly specializable along $t$ at any $\lambda_0$.
\index{strictly specializable}
\hfill\qed
\end{df}
The union
$\bigcup_a\nbigk(a,\lambda_0,P)$
is denoted by $\KMS(\nbigm,t,P)$.
We also put
$\nbigk(a,\lambda_0):=
 \bigcup_{P\in\nbigxzero_0}\nbigk(a,\lambda_0,P)$
and
$\KMS(\nbigm,t):=
 \bigcup_{P\in\nbigxzero_0}
 \KMS(\nbigm,t,P)$.

\begin{rem}
If we consider local issues,
we implicitly assume that
the set $\nbigk(a,\lambda_0)$ is finite.
\hfill\qed
\end{rem}

Let $\nbigm$ and $\nbigk(a,\lambda_0)$
be as in Definition \ref{df;07.10.24.1}.
For $u\in\nbigk(a,\lambda_0)$,
we put on $\nbigx^{(\lambda_0)}_0$
\[
 \tildepsizero_{t,u}(\nbigm):=
\bigcup_{N}
 \Ker\Bigl(
 (-\deldel_tt+\eigenmap(\lambda,u))^N:
 \Gr^{\Vzero}_a(\nbigm)\lrarr
 \Gr^{\Vzero}_a(\nbigm)
 \Bigr).
\]
\index{functor $\tildepsizero_{t,u}(\nbigm)$}
By the strictness,
we have the decomposition:
\[
 \Gr^{\Vzero}_a(\nbigm)
=\bigoplus_{u\in \nbigk(a,\lambda_0)}
\tildepsizero_{t,u}(\nbigm)
\]
We implicitly assume 
$\tildepsizero_u(\nbigm)\neq 0$ for each
$u\in\nbigk(a,\lambda_0)$.
Because the multiplication of $t$
induces the isomorphism
$\Gr^{\Vzero}_a(\nbigm)\simeq
 \Gr^{\Vzero}_{a-1}(\nbigm)$,
we may have the bijection:
\[
 \nbigk(a,\lambda_0)
\simeq
 \nbigk(a-1,\lambda_0),
\quad
 u\longleftrightarrow
 u-\vecdelta_0
\]
Here $\vecdelta_0=(1,0)\in\real\times\cnum$.

\begin{lem} \mbox{{}}
\label{lem;07.10.29.1}
\begin{itemize}
\item
The filtration $\Vzero$ as in Definition 
{\rm\ref{df;07.10.24.1}} is unique,
if it exists.
\item
Assume that $\nbigm$ is strictly specializable 
along $t$ at $\lambda_0$.
For simplicity,
$\nbigxzero$ is assumed to be the product
of $X$ and a neighbourhood of $\lambda_0$,
and $\KMS(\nbigm,t)$ is assumed to be finite.
Then, $\nbigm$ is also strictly specializable 
along $t$ at any $\lambda_1$
such that $|\lambda_1-\lambda_0|$ is sufficiently small.
The filtration $V^{(\lambda_1)}$ 
of $\nbigm_{|\nbigx^{(\lambda_1)}}$
can be constructed from $\Vzero$,
and we have
$\tildepsizero_{t,u|\nbigx^{(\lambda_1)}_0}
=\tildepsi^{(\lambda_1)}_{t,u}$
on $\nbigx^{(\lambda_1)}_0
 \subset\nbigx^{(\lambda_0)}_0$.
\end{itemize}
\end{lem}
\pf
See \cite{sabbah2}
or the arguments in Lemma \ref{lem;08.10.15.10}
and Lemma \ref{lem;9.17.10}.
\hfill\qed

\vspace{.1in}
As a result,
if $\nbigm$ on $\nbigx$ is strictly specializable along $t$,
the set $\KMS(\nbigm,t)$ is independent of
the choice of $\lambda_0$.
By gluing
$\psitilde^{(\lambda_0)}_{t,u}(\nbigm)$,
we obtain
the $\nbigr_{X_0}$-module $\psitilde_{t,u}(\nbigm)$
on $\nbigx_0$ for each $u\in \KMS(\nbigm,t)$.
Formally, we put
$\psitilde_{t,u}(\nbigm):=0$
for $u\not\in\KMS(\nbigm,t)$.
The multiplication of $t$ induces the isomorphism
$\psitilde_{t,u}(\nbigm)
 \simeq\psitilde_{t,u-\vecdelta_0}(\nbigm)$.
We put $N:=-\deldel_tt+\eigenmap(\lambda,u)$
on $\psitilde_{t,u}(\nbigm)$,
which is the nilpotent part of
$-\deldel_tt$.

The following lemma can be shown
by an argument in the proof of 
Lemma \ref{lem;08.10.15.10}.
\begin{lem}
Let $\nbigm_i$ $(i=1,2)$ be strictly specializable
along $t$.
Any morphism $\varphi:\nbigm_1\lrarr\nbigm_2$
preserves the filtrations $\Vzero$
for any $\lambda_0$.
In particular,
we have the induced morphism
$\psitilde_{t,u}(\varphi):
 \psitilde_{t,u}(\nbigm_1)\lrarr\psitilde_{t,u}(\nbigm_2)$
for each $u\in\real\times\cnum$.
\hfill\qed
\end{lem}

Let $\nbigm_i$ $(i=1,2)$ be $\nbigr_X(\ast t)$-modules
which are strictly specializable along $t$.
A morphism $\varphi:\nbigm_1\lrarr\nbigm_2$ 
is called strictly specializable along $t$
if $\Cok\bigl(\psitilde_{t,u}(\varphi)\bigr)$ is strict.
\index{strictly specializable morphism}
The following lemma can be shown
by using an argument of the proof of 
Proposition \ref{prop;9.17.30}.
\begin{lem}
\label{lem;08.1.16.11}
Let $\nbigm_i$ be strictly specializable along $t$.
Let $\varphi:\nbigm_1\lrarr\nbigm_2$
be a morphism which is strictly specializable along $t$.
Then, $\varphi$ is strict
with respect to the filtrations $\Vzero$,
i.e.,
$\Vzero_a(\nbigm_2)\cap \Image(\varphi)
=\varphi\bigl(\Vzero_a(\nbigm_1)\bigr)$.

As a result,
$\Ker(\varphi)$,
$\Image(\varphi)$
and $\Cok(\varphi)$ are also
strictly specializable along $t$,
and we have the natural isomorphisms
$\psitilde_{t,u}\Ker(\varphi)
\simeq
 \Ker\psitilde_{t,u}(\varphi)$,
$\psitilde_{t,u}\Image(\varphi)
\simeq
 \Image\psitilde_{t,u}(\varphi)$
and 
$\psitilde_{t,u}\Cok(\varphi)
\simeq
 \Cok\psitilde_{t,u}(\varphi)$.
\hfill\qed
\end{lem}

We refer to Section 3.4 of \cite{sabbah2}
for the following lemma.
\begin{lem}
Let $\nbigm$ be a coherent $\nbigr_X$-module
which is strictly specializable along $t$ at $\lambda_0$
with the filtration $\Vzero$.
\begin{itemize}
\item
$\nbigm(\ast t):=
 \nbigr_X(\ast t)\otimes_{\nbigr_X}\nbigm$
is a coherent $\nbigr_X(\ast t)$-module
which is strictly specializable along $t$ at $\lambda_0$.
\item
The filtration $\Vzero(\nbigm(\ast t))$
is given by
$\Vzero_a(\nbigm(\ast t))
=\Vzero_a(\nbigm)$ for $a<0$,
and $\Vzero_a(\nbigm(\ast t))=
 t^{-n}\Vzero_{a-n}(\nbigm)$
for $a\geq 0$, where $n$ is chosen as $a-n<0$.
\item
We have the natural isomorphism
$\psitilde_{t,u}(\nbigm)
\simeq
 \psitilde_{t,u}(\nbigm(\ast t))$
for any 
 $u\in\real\times\cnum$.
\hfill\qed
\end{itemize}
\end{lem}

\subsection{Pull back via $n$-th ramified covering}
\label{subsection;08.10.4.1}
\index{pull back via ramified covering}

We put $X^{(n)}:=X_0\times\cnum_{t_n}$.
Let $\varphi_n:X^{(n)}\lrarr X$
be the morphism induced by $\varphi_n^{\ast}(t)=t_n^n$.
Let $\nbigm$ be an $\nbigr_X(\ast t)$-module.
Since $\nbigr_{X^{(n)}}(\ast t_n)
 =\nbigo_{X^{(n)}}
 \otimes_{\varphi_n^{-1}(\nbigo_X)}
 \varphi_n^{-1}\nbigr_X(t)$ is flat 
over $\varphi_n^{-1}\nbigr_X(\ast t)$,
we have
$\varphi_n^{\dagger}\nbigm \simeq
 \varphi_n^{\ast}\nbigm$.
If $\nbigm$ is $\nbigr_{X}(\ast t)$-coherent,
$\varphi_n^{\dagger}\nbigm$ is 
$\nbigr_{X^{(n)}}(\ast t_n)$-coherent.
If $\nbigm$ is a holonomic $\nbigr_X(\ast t)$-module,
$\varphi_n^{\dagger}\nbigm$ is 
a holonomic $\nbigr_{X^{(n)}}(\ast t_n)$-module.
The following lemma can be checked directly.
\begin{lem}
If $\nbigm$ is strictly specializable along $t$,
then so is $\varphi_n^{\dagger}\nbigm$.
The $V$-filtration $\Vzero$ of $\varphi_n^{\dagger}\nbigm$
is given as follows:
\[
 \Vzero_a(\varphi_n^{\dagger}\nbigm)
=\sum_{nb-c\leq a}
 \varphi_n^{-1}(\Vzero_b\nbigm)\cdot t_n^c
\]
In particular, we have the following natural isomorphism
of $\nbigr_{X_0}$-modules:
\[
  \psitilde_{t_n,u}(\varphi_n^{\dagger}\nbigm)
\simeq
 \bigoplus_{(u',c)\in \nbigs}
 \psitilde_{t,u'}(\nbigm),
\]
\[
 \nbigs=\bigl\{(u',c)\,\big|\,
 u'\in \KMS(\nbigm,t),\,\,
 c\in\seisuu,\,\,0\leq c\leq n-1,\,\,\,
n \cdot u'-c\cdot\vecdelta_0=u
 \bigr\}
\]
Here, $\vecdelta_0:=(1,0)\in\real\times\cnum$.
\hfill\qed
\end{lem}

\subsection{Exponential twist}
\label{subsection;08.9.4.12}
\index{exponential twist}

Let $\gminia\in \cnum[t_n^{-1}]$.
We have the $\nbigr_{\cnum_{t_n}}(\ast t_n)$-module
$\nbigl(\gminia)$
given as follows:
\[
 \nbigl(\gminia)=\nbigo_{\cnum_{t_n}}(\ast t_n)\cdot e,
\quad\quad
 \deldel_{t_n}e=\del_{t_n}\gminia \cdot e
\]
The pull back via the projection $X^{(n)}\lrarr \cnum_{t_n}$
is also denoted by $\nbigl(\gminia)$.

For an $\nbigr_X(\ast t)$-module $\nbigm$,
we have the $\nbigr_{X^{(n)}}(\ast t_n)$-module
$\varphi_n^{\dagger}\nbigm\otimes_{\nbigo_{X^{(n)}}}
 \nbigl(-\gminia)$.
If it is strictly specializable along $t_n$,
we say that $\nbigm$ is strictly specializable
along $t$ with ramification and 
exponential twist by $\gminia$,
and we define for any $u\in\real\times\cnum$
\[
 \psitilde_{t,\gminia,u}(\nbigm):=
 \psitilde_{t_n,u}\bigl(
 \varphi_n^{\dagger}\nbigm\otimes\nbigl(-\gminia)
 \bigr).
\]
\index{functor $\psitilde_{t,\gminia,u}$}
\index{strictly specializable with ramification
 and exponential twist}

For an $\nbigr_X$-module $\nbigm$,
we have the induced $\nbigr_X(\ast t)$-module
$\nbigm(\ast t)$.
We define 
\[
 \psitilde_{t,\gminia,u}(\nbigm):=
 \psitilde_{t,\gminia,u}\bigl(\nbigm(\ast t)\bigr),
\]
if the right hand side can be defined.

\subsection{Comparison of 
 strictly $S$-decomposable $\nbigr$-modules}

Let $X=\cnum_t\times X_0$ 
for some complex manifold $X_0$.
We identify $X_0$ with $\{t=0\}$.

\begin{lem}
\label{lem;08.1.21.10}
Let $\nbigm$ be a coherent $\nbigr_X$-module
whose support is contained in $X_0$.
Then, $\nbigm(\ast t)=0$.

A similar claim holds for coherent $D$-modules.
\end{lem}
\pf
We have only to show the claim locally.
Let $\pi:T^{\ast}X\lrarr X$ denote the projection
of the cotangent bundle.
We take a coherent filtration $F$ of $\nbigm$.
The associated graded module $\Gr^F(\nbigm)$ 
induces a coherent $\nbigo_{T^{\ast}X}$-module,
and the support is contained in $\pi^{-1}(X_0)$.
Hence, the action of $t$ on $\nbigm$ is locally nilpotent,
and we obtain $\nbigm(\ast t)=0$.
\hfill\qed

\vspace{.1in}

Let $\nbigm_i$ $(i=1,2)$ be coherent $\nbigr_X$-modules,
and let $f:\nbigm_1\lrarr\nbigm_2$ be a morphism.
\begin{lem}
\label{lem;08.1.21.11}
If the restriction
$f_{|X-\{t=0\}}$ is an isomorphism,
the induced morphism
$\nbigm_1(\ast t)\lrarr \nbigm_2(\ast t)$
is an isomorphism.
A similar claim holds for $D$-modules.
\end{lem}
\pf
We obtain $\Ker(f)(\ast t)=\Cok(f)(\ast t)=0$
due to Lemma \ref{lem;08.1.21.10}.
\hfill\qed

\vspace{.1in}

Let $X$ be a general complex manifold,
and let $g$ be a holomorphic function on $X$.
Let $\nbigm_i$ $(i=1,2)$ be  $\nbigr_X$-modules
such that 
(i) they are strictly $S$-decomposable along $g$,
(ii) they do not have any $\nbigr_X$-submodules
 whose supports are contained in $\{g=0\}$.
By the condition (ii),
we have the inclusions
$\nbigm_i\subset\nbigm_i(\ast g)$.

\begin{lem}
\label{lem;07.11.2.15}
If we have an isomorphism
$\nbigm_1(\ast g)\simeq
 \nbigm_2(\ast g)$,
we obtain an isomorphism
$\nbigm_1\simeq\nbigm_2$
as the restriction.
A similar claim holds for $D$-modules.
\end{lem}
\pf
Let $\nbigm_i^{(\lambda_0)}$
denote the restriction of $\nbigm_i$
to a small neighbourhood of $\{\lambda_0\}\times X$.
We may assume that $X=X_0\times\cnum_t$
and $g=t$.
The $V$-filtrations of 
$\nbigm_i^{(\lambda_0)}$ $(i=1,2)$
induce the $V$-filtrations of
$\nbigm_i^{(\lambda_0)}(\ast t)$
via which $\nbigm_i^{(\lambda_0)}(\ast t)$
are strictly specializable along $t$ at $\lambda_0$.
By the uniqueness of such $V$-filtrations,
they are the same under the identification
$\nbigv^{(\lambda_0)}
:=\nbigm_1^{(\lambda_0)}(\ast t)
=\nbigm_2^{(\lambda_0)}(\ast t)$.
In particular,
we have $A^{(\lambda_0)}:=
 \Vzero_{<0}(\nbigm_1^{(\lambda_0)})
=\Vzero_{<0}(\nbigm_2^{(\lambda_0)})$
in $\nbigv^{(\lambda_0)}$.
Since both $\nbigm^{(\lambda_0)}_i$ $(i=1,2)$
are generated by $A^{(\lambda_0)}$ in 
$\nbigv^{(\lambda_0)}$,
we obtain
$\nbigm_1^{(\lambda_0)}
=\nbigm_2^{(\lambda_0)}$
for any $\lambda_0$,
and thus $\nbigm_1=\nbigm_2$.
\hfill\qed

\section{Formal $\nbigr$-modules}
\label{subsection;08.9.29.122}
\subsection{Formal complex spaces}
\label{subsection;07.11.6.10}

We recall some basic facts on
formal complex spaces.
We refer to \cite{banica},
 \cite{bingener} and \cite{krasnov}
for more details and precision.
\index{formal complex space}
Let $X$ be a complex space.
Let $T$ be an analytic subvariety of $X$,
which consists of the underlying
topological space $|T|$
and the structure sheaf $\nbigo_T$.
Let $\nbigi_T$ denote the sheaf of ideals of
$\nbigo_X$, corresponding to $T$.
We have the analytic subspace $T^{(n)}$
corresponding to $\nbigi_T^n$,
i.e.,
$T^{(n)}=\bigl(|T|,\nbigo_X/\nbigi_T^n\bigr)$.
As the limit,
we obtain the ringed space
$\That:=\bigl(|T|,\varprojlim\nbigo_X/\nbigi_T^n\bigr)$,
which is called the completion of $X$
along $T$.
\index{completion}

\begin{prop}
\label{prop;07.7.30.1}
\mbox{{}}
\begin{enumerate}
\item
$\That$ is a formal complex space
in the sense of {\rm\cite{bingener}},
and the morphism $\iota_T:\That\lrarr X$ is flat.
\item
The sheaf of the algebras
$\nbigo_{\That}$ is coherent and Noetherian.
\item
Coherent sheaves $\nbigf$ on $\That$
are equivalent to 
systems of coherent sheaves
$\nbigf^{(n)}$ $(n=1,2,\ldots,)$ on $T^{(n)}$
such that 
$\nbigf^{(n)}\otimes\nbigo_{T^{(m)}}=\nbigf^{(m)}$
for $n\geq m$.
\end{enumerate}
\end{prop}
\pf
The first claim is Lemma 1.6 of \cite{bingener}.
The second claim is
Lemma 1.1 and Corollary 1.5 of \cite{bingener}.
(See Section \ref{subsection;08.10.20.2}
 for Noetherian property in this situation.)
The third claim is Lemma 1.2 of \cite{bingener}.
\hfill\qed

\begin{lem}
\label{lem;07.7.30.2}
The extension $\iota_T^{-1}\nbigo_X\lrarr \nbigo_{\That}$ is
faithfully flat.
We have $\Supp(\nbigf)\cap |T|=\emptyset$
for a coherent $\nbigo_X$-module $\nbigf$ 
with $\iota_T^{\ast}\nbigf=0$.
\end{lem}
\pf
Let $P$ be any point of $|T|$.
We have the completion $\Phat$.
Let $\iota_P:\Phat\lrarr X$ an
$\iota_{P,T}:\Phat\lrarr \That$ 
denote the canonical morphisms.
Since $\nbigo_{\Phat}$ is the completion of
the local ring $\nbigo_{X,P}$ at $P$,
the morphism $\iota_P$ is faithfully flat.
In the case $\iota_T^{\ast}\nbigf=0$,
we have $\iota_P^{\ast}\nbigf=0$,
and hence $\iota_P^{-1}\nbigf=0$.
Since $P$ can be any point of $|T|$,
we obtain $\iota_T^{-1}\nbigf=0$.
The second claim follows from 
the first claim and 
Proposition \ref{prop;08.10.20.3}.
\hfill\qed

\subsection{Formal $D$-modules}
\label{subsection;07.7.30.4}
\index{formal $D$-modules}

Let $X$ be a complex manifold.
Let $Z$ be an analytic subvariety of $X$.
Let $\Zhat$ be the completion of $X$ along $Z$.
Let $\iota:\Zhat\lrarr X$ denote the canonical morphism.
Let $\Theta_X$ denote the tangent sheaf of $X$.
We put $\Theta_{\Zhat}:=\iota^{\ast}\Theta_X$,
which acts on $\nbigo_{\Zhat}$ as differential operators.
Let $\nbigd_{\Zhat}$ denote the sheaf of differential operators
of $\nbigo_{\Zhat}$,
which is generated by $\Theta_{\Zhat}$ and $\nbigo_{\Zhat}$
with the standard relations.
If we are given a coordinate 
system $(z_1,\ldots,z_n)$ of $X$,
$\Theta_{\Zhat}$ is the free $\nbigo_{\Zhat}$-module
with the base $\del_i:=\del/\del z_i$,
and $\nbigd_{\Zhat}$ is the sheaf of algebras
generated by $\nbigo_{\Zhat}$ and $\Theta_{\Zhat}$
with the relation $[\del_i,z_j]=\delta_{i,j}$.
By applying the argument
in Appendix A.1 of \cite{kashiwara_text},
we obtain the following standard proposition.
\begin{prop}
\label{prop;07.7.30.3} \mbox{{}}
\begin{itemize}
\item
The sheaf of algebras $\nbigd_{\Zhat}$ is Noetherian,
and it has the property in 
Proposition {\rm\ref{prop;08.10.20.1}}.
The Rees ring is also Noetherian.
\item
Let $M$ be a $\nbigd_{\Zhat}$-module
such that (i) pseudo-coherent as an $\nbigo_{\Zhat}$-module,
(ii) locally finitely generated 
 as a $\nbigd_{\Zhat}$-module.
Then, $M$ is a coherent $\nbigd_{\Zhat}$-module.
\end{itemize}
\end{prop}
\pf
We give only an outline.
We have the standard filtration of $\nbigd_{\Zhat}$
by the order of the differential operators,
i.e., $F_m\nbigd_{\Zhat}=
 \bigl\{\sum_{|J|\leq m} a_J\cdot\del^J\bigr\}$
on a coordinate neighbourhood,
where $\del^J=\prod \del_i^{j_i}$
and $|J|=\sum j_i$ for $J=(j_1,\ldots,j_n)$.
Then, $F_0D_{\Zhat}=\nbigo_{\Zhat}$ is Noetherian
(Proposition \ref{prop;07.7.30.1}).
Since $\Gr^FD_{\Zhat}$ is locally a polynomial algebra
over $F_0D_{\Zhat}$, it is also Noetherian.
Then, the first claim of the proposition follows
from  Proposition \ref{prop;08.10.20.30}.
The second claim follows from 
Proposition \ref{prop;08.10.20.31}.
\hfill\qed

\vspace{.1in}

Let $\nbigm$ be a $\nbigd_X$-module.
By the standard formalism,
$\iota^{\ast}\nbigm$ has the $\nbigd_{\Zhat}$-module structure.
Thus we obtain the functor
$\iota^{\ast}$ of the category of $\nbigd_X$-modules
to the category of $\nbigd_{\Zhat}$-modules.

\begin{lem}
\label{lem;07.7.30.5}
The functor $\iota^{\ast}$ is exact.
If $\iota^{\ast}\nbigm=0$ for a coherent
$\nbigd_{\Zhat}$-module,
then $\Supp(\nbigm)\cap |Z|=\emptyset$.
\end{lem}
\pf
The first claim follows from 
Proposition \ref{prop;07.7.30.1}.
If $\iota^{\ast}\nbigm=0$,
we obtain $\iota^{-1}\nbigm=0$
because of Lemma \ref{lem;07.7.30.2}.
Then, the second claim follows from
Proposition \ref{prop;08.10.20.3}.
\hfill\qed

\begin{lem}
\label{lem;07.7.30.6}
If $\nbigm$ is a coherent $\nbigd_X$-module,
then $\iota^{\ast}\nbigm$ is a coherent
$\nbigd_{\Zhat}$-module.
\end{lem}
\pf
Because of the Noetherian property 
of $\nbigd_X$ and $\nbigd_{\Zhat}$,
the coherence is equivalent to 
the locally finitely presentedness.
Then, the claim is clear.
\hfill\qed

\subsection{Formal $\nbigr$-modules}
\label{subsection;08.9.29.211}
\index{formal $\nbigr$-module}

Let $X$ be a complex manifold,
and let $Z$ be an analytic subvariety of $X$.
Let $\cnum_{\lambda}$ denote the complex line
with the coordinate $\lambda$.
We put $\nbigx:=\cnum_{\lambda}\times X$
and $\nbigz:=\cnum_{\lambda}\times Z$.
The morphisms induced by $Z\lrarr X$
are denoted by $\iota$.
The completion of $\nbigx$ along $\nbigz$
is denoted by $\nbigzhat$.
Let $p_{\lambda}$ denote the projection of
$\nbigx$ or $\nbigzhat$ onto $X$ or $\Zhat$.
Let $\nbigr_{\Zhat}$ be the sheaf of subalgebras
of $p_{\lambda}^{\ast}\nbigd_{\Zhat}$
generated by $\nbigo_{\nbigzhat}$
and $\lambda\cdot p_{\lambda}^{\ast}\Theta_{\Zhat}$.
If a holomorphic coordinate 
system $(z_1,\ldots,z_n)$ of $X$ is given,
let $\deldel_i$ denote $\lambda\cdot\del_i$,
as usual.
The following lemma can be shown
using the argument for
Proposition \ref{prop;07.7.30.3}.
\begin{prop}
\label{prop;07.9.7.1}
 \mbox{{}}
\begin{itemize}
\item
The sheaf of algebras $\nbigr_{\Zhat}$ is Noetherian,
and it has the property in
Proposition {\rm\ref{prop;08.10.20.1}}.
The Rees algebra is also Noetherian.
\item
Let $\nbigm$ be an $\nbigr_{\Zhat}$-module
such that 
(i) pseudo-coherent as an $\nbigo_{\nbigzhat}$-module,
(ii) locally finitely generated 
 as an $\nbigr_{\Zhat}$-module.
Then, $\nbigm$ is a coherent $\nbigr_{\Zhat}$-module.
\hfill\qed
\end{itemize}
\end{prop}

Let $\nbigm$ be an $\nbigr_X$-module on $\nbigx$.
By the standard formalism, 
$\iota^{\ast}\nbigm$ has the $\nbigr_{\Zhat}$-module structure
on $\nbigzhat$.
Thus we obtain the functor $\iota^{\ast}$
of the category of $\nbigr_X$-modules on $\nbigx$
to the category of $\nbigr_{\Zhat}$-modules on $\nbigzhat$.
The following claims can be shown
by the arguments in Section \ref{subsection;07.7.30.4}.
\begin{lem}
\mbox{{}}\label{lem;08.1.23.15}
\begin{itemize}
\item
 The functor $\iota^{\ast}$ is exact.
 If $\iota^{\ast}\nbigm=0$ for a
 coherent $\nbigr_X$-module,
 then $\Supp(\nbigm)\cap|\nbigz|=\emptyset$.
\item
If $\nbigm$ is a coherent $\nbigr_X$-module,
then $\iota^{\ast}\nbigm$ is a coherent
$\nbigr_{\Zhat}$-module.
\hfill\qed
\end{itemize}
\end{lem}

\subsection{A criterion for 
 strict $S$-decomposability}
\label{subsection;08.9.29.212}

Let $X_0$ be a complex manifold,
and let $Z_0$ be a closed submanifold of $X_0$.
Let $X:=X_0\times\cnum_t$
and $Z:=Z_0\times\cnum_t$.
Let $\iota:\Zhat_0\lrarr X_0$ denote the natural morphism.
The induced morphisms are also denoted by $\iota$.
Let $\nbigx^{(\lambda_0)}$ denote $U(\lambda_0)\times X$,
where $U(\lambda_0)$ is 
some neighbourhood of $\lambda_0$ in $\cnum_{\lambda}$.
We use the symbols like $\nbigz^{(\lambda_0)}$,
$\nbigx_0^{(\lambda_0)}$
and $\nbigzhat^{(\lambda_0)}_0$
in similar meanings.

Let $\nbigm_i$ $(i=1,2)$ be coherent
$\nbigr_X$-modules equipped with
filtrations $\Vzero(\nbigm_i)$
satisfying Condition \ref{condition;08.1.23.10}.
We obtain the $\nbigr_{\Zhat}$-modules
$\iota^{\ast}\nbigm_i$ $(i=1,2)$
on $\nbigzhat^{(\lambda_0)}$
with the induced filtrations
$\iota^{\ast}\Vzero$ indexed by $\real$.
Assume that 
we are given an isomorphism
$\varphi:\iota^{\ast}\nbigm_1
\lrarr \iota^{\ast}\nbigm_2$
of $\nbigr_{\Zhat}$-modules
preserving the filtrations.

\begin{prop}
 \label{prop;08.1.23.16}
 If $\nbigm_1$
 is strictly specializable 
 ($S$-decomposable)
 along $t$ at $\lambda_0$
 with respect to
 $\Vzero(\nbigm_1)$,
 then $\bigl(\nbigm_2,\Vzero(\nbigm_2)\bigr)$
 is also strictly specializable 
 ($S$-decomposable)
 along $t$ at $\lambda_0$
 with respect to $\Vzero(\nbigm_2)$,
 after shrinking $X_0$ around $Z_0$.
(See Condition {\rm\ref{condition;08.1.23.11}}
and Condition {\rm\ref{condition;08.1.23.12}}
 for strictly specializability and
 strict $S$-decomposability.)
\end{prop}
\pf
We have an induced isomorphism
$\iota^{\ast}
 \Gr^{\Vzero}_{a}(\nbigm_1)
\simeq
 \iota^{\ast}
 \Gr^{\Vzero}_a(\nbigm_2)$
as in Proposition \ref{prop;9.17.30}.
Let $N_2$ be the kernel of
$\lambda-\lambda_0$
on $\Gr^{\Vzero}_a(\nbigm_2)$.
Since $\Gr^{\Vzero}_a(\nbigm_1)$ is strict,
$\iota^{\ast}\Gr^{\Vzero}_a(\nbigm_2)$ is also strict
due to Lemma \ref{lem;08.1.23.15}.
Hence, we obtain $\iota^{\ast}N_2=0$.
Then, we obtain $N_2=0$ by shrinking $X_0$,
due to Lemma \ref{lem;08.1.23.15}.
By using coherence,
we obtain that $\Gr^{\Vzero}_a(\nbigm_2)$
is strict, after shrinking $U(\lambda_0)$ and $X_0$
appropriately.
We can check the other conditions
by using Lemma \ref{lem;08.1.23.15}.
\hfill\qed

\section{Preliminary for $D$-modules}
\label{subsection;08.10.18.105}
We recall some basic facts for $D$-modules
and $\nbigr$-modules.
(See \cite{bjork}, \cite{kashiwara_text},
\cite{hotta-tanisaki}, \cite{sabbah2},
 \cite{mochi2}, for example.)

\subsection{Pull back of $D$-modules}
\label{subsection;08.1.21.1}

Let $f:X\lrarr Y$ be a morphism of complex manifolds,
and let $\nbigf$ be a $D_Y$-module.
We put $D_{X\rarr Y}:=
 \nbigo_X\otimes_{f^{-1}\nbigo_Y}f^{-1}D_Y$,
which is equipped with the left $D_X$-action
and the right $f^{-1}D_Y$-action.
Recall that the pull back
$f^{\dagger}\nbigf$ is defined to be
$f^{\dagger}\nbigf:=
 D_{X\rarr Y}\otimes_{f^{-1}D_Y}f^{-1}\nbigf$.
Let $L^{i}f^{\dagger}$ denote 
the $i$-th derived functor of $f^{\dagger}$.
Let $\forget$ denote the natural
functor from the category of $D$-modules 
to the category of $\nbigo$-modules.
Let $L^if^{\ast}$ denote the $i$-th derived
functor of the pull back $f^{\ast}$ for $\nbigo$-modules.
According to the following lemma,
we do not have to distinguish 
the sheaves $L^if^{\dagger}\nbigf$
and $L^if^{\ast}\nbigf$.
\begin{lem}
$\forget(L^if^{\dagger}\nbigm)
\simeq
 L^if^{\ast}(\forget(\nbigm))$
\end{lem}
\pf
See Proposition 2.3.8 of \cite{bjork},
for example.
\hfill\qed

\vspace{.1in}
Recall the following well known results.

\begin{prop}
\label{prop;07.11.2.50}
\mbox{{}}
\begin{itemize}
\item
If $\nbigf$ is a holonomic $D_Y$-module,
then $L^if^{\dagger}\nbigf$ are also holonomic.
\item
Assume that $f$ is a closed embedding of complex manifolds,
and $\nbigf=f_{\dagger}\nbigf_0$ for some
$D_X$-module $\nbigf_0$.
Let $k:=\dim X-\dim Y$.
Then, we have $L^if^{\ast}\nbigf=0$ for $i\neq k$
and $L^{k}f^{\ast}\nbigf\simeq\nbigf_0$.
\end{itemize}
\end{prop}
\pf
See Theorem 3.2.13 of \cite{bjork} for the first claim,
for example.
See Proposition 4.32 of \cite{kashiwara_text}
for the second claim,
for example.
\hfill\qed

\vspace{.1in}
Recall that we have the trace map
whose construction is explained in Chapter 4.9
of \cite{kashiwara_text}, for example:
\[
 \tr_f:
 f_{\dagger}\bigl(Lf^{\dagger}\nbigf\bigr)[\dim X]
\lrarr \nbigf[\dim Y]
\]

\subsection{Holonomic $D$-modules
and meromorphic connections}
\label{subsection;08.9.29.200}

We recall the following general lemma.
Let $X$ be a complex manifold,
and let $D$ be a divisor of $X$.
Let $\pi:T^{\ast}X\lrarr X$ denote the natural projection
of the cotangent bundle.
\begin{lem}
\label{lem;07.11.2.1}
Let $\nbigf$ be a holonomic $D_X$-module
such that 
(i) the characteristic variety of $\nbigf$ is contained in
$\pi^{-1}(D)\cup T_X^{\ast}X$,
(ii) $\nbigf$ is also an $\nbigo_X(\ast D)$-module.
Then, there is an $\nbigo_X$-coherent subsheaf 
$\nbigf_0\subset\nbigf$
such that $\nbigf=\nbigo_{X}(\ast D)\cdot\nbigf_0$.
In particular,
$\nbigf$ gives a meromorphic flat connection
on $(X,D)$.
\end{lem}
\pf
Since the claim is local,
we may assume that $D$ is given as $g^{-1}(0)$.
We may also take a coherent filtration $F$
of $\nbigf$.
Then, $\Gr^F(\nbigf)$ is 
a coherent $\nbigo_{T^{\ast}X}$-module,
and the support of $\Gr^F(\nbigf)$ is contained in
$T_X^{\ast}X\cup \pi^{-1}(D)$.
Hence, there exists a large number $i_0$
such that the support of $\Gr^F_{i}(\nbigf)$
is contained in $\pi^{-1}(D)$
for any $i\geq i_0$.
For such an $i$, 
there exists a large number $N$
such that $g^N\cdot F_i(\nbigf)\subset F_{i_0}(\nbigf)$.
Hence, $F_{i_0}(\nbigf)$ generates $\nbigf$
over $\nbigo_X(\ast D)$.
\hfill\qed

\subsection{Nearby cycle functor
 with ramification and exponential twist}
\label{subsection;08.10.25.1}

Recall the notion of
nearby cycle functor with ramification
and exponential twist,
by following Deligne and Sabbah
(\cite{Deligne-Malgrange-Ramis},
 \cite{sabbah5}).
Let $X_0$ be a complex manifold,
and let $X:=X_0\times\cnum_t$.
We set $X^{(n)}:=X_0\times\cnum_{t_n}$.
Let $\varphi_n:X^{(n)}\lrarr X$ be a morphism
induced by $\varphi_n(t_n)=t_n^n$.
For a given $\gminia\in\cnum[t_n^{-1}]$,
we set
$L(\gminia):=\nbigo_{X^{(n)}}(\ast t_n)\cdot e$
with the meromorphic flat connection
$\nabla e=e\cdot d\gminia$.
It is naturally a holonomic $D_{X^{(n)}}$-module.

Let $\nbigf$ be a holonomic $D_X$-module.
By taking the pulling back 
and the tensor product with $L(-\gminia)$,
we obtain a holonomic $D_{X^{(n)}}$-module
$\varphi_n^{\dagger}\nbigf
 \otimes L(-\gminia)$.
Applying the nearby cycle functor,
we obtain a holonomic $D_{X_0}$-module
$\psi_{t_n}\bigl(
 \varphi_n^{\dagger}\nbigf
 \otimes L(-\gminia)
 \bigr)$.

\vspace{.1in}
Let $Y$ be a complex manifold,
$g$ be a holomorphic function on $Y$,
and $\nbigf$ be a holonomic $D_Y$-module.
Then, for any $\gminia\in\cnum[t_n^{-1}]$,
we define
\[
 \psi_{g,\gminia}(\nbigf):=
 \psi_{t_n}\Bigl(
 \varphi_n^{\ast}\bigl(
 i_{g\dagger}\nbigf
 \bigr)\otimes\nbigl(-\gminia)
 \Bigr)
\]
The functor $\psi_{g,\gminia}$
is called the nearby cycle functor
with ramification and exponential twist by $\gminia$.

\subsection{The functors $j_!j^{\ast}$ 
and $j_{\ast}j^{\ast}$}
\label{subsection;08.1.19.20}
\index{functor $j_{\ast}j^{\ast}$}
\index{functor $\jbikkuri j^{\ast}$}

Let $X$ be a complex manifold,
and let $i_Y:Y\subset X$ be a smooth hypersurface
of $X$.
Let $\Theta_X$ denote the tangent sheaf of $X$.
Let $\nbign_{Y/X}$ denote the sheaf of
the sections of the normal bundle of $Y$ in $X$.
We put $j_{\ast}j^{\ast}\nbigo_X:=\nbigo_X(\ast Y)$.

We also have the following $D_{X}$-module
$j_{!}j^{\ast}\nbigo_X$.
As an $\nbigo_{X}$-module,
we set
\[
  j_!j^{\ast}\nbigo_{X}:=
 \nbigo_{\nbigx}\oplus 
 i_{Y\dagger}i_Y^{\dagger}\nbigo_{X}
=\nbigo_{X}\oplus
 i_{Y\ast}\bigl(
 \Sym^{\cdot}\nbign_{Y/X}
\otimes\nbign_{Y/X}
 \bigr).
\]
The action of $\Theta_{X}$ is given by
$ v\cdot (s,t)
=\bigl(v\cdot s, \,\,
 v\cdot t+\pi(v_{|Y})\cdot s_{|Y}
 \bigr)$,
where ``$|Y$'' denotes the restriction to $Y$,
$\pi$ denotes the projection
of $\Theta_{X|Y}\lrarr \nbign_{Y/X}$,
and $v\cdot t$ is given by the natural
$D_X$-module structure on
$i_{Y\dagger}i_Y^{\dagger}\nbigo_X$.
It is uniquely extended to the action of $D_X$.

Let $\Omega_X$ denote the canonical line bundle of $X$.
For a $D_X$-module $M$,
we have the derived dual module
$\DDD_X(M):=
 \nrhom_{D_X}(M,D_X)
 \otimes\Omega_X^{-1}[\dim X]$
in the derived category of $D_X$-modules.

\begin{lem}
\label{lem;08.1.19.3}
We have a natural isomorphism
$\DDD_X\bigl(j_{\ast}j^{\ast}\nbigo_X\bigr)
\simeq
 j_!j^{\ast}\nbigo_X$
\end{lem}
\pf
Note that we have the canonical isomorphisms
$(j_!j^{\ast}\nbigo_X)_{|X-Y}
\simeq
 \nbigo_{X|X-Y}
\simeq
 \DDD_X\bigl(j_{\ast}j^{\ast}\nbigo_X\bigr)
 _{|X-Y}$.
We have the $\nbigo_X$-submodule
$\nbigo_X$ of $j_!j^{\ast}\nbigo_X$,
which generates $j_!j^{\ast}\nbigo_X$
over $D_X$.
Hence,
an automorphism
$\varphi$ of $j_!j^{\ast}\nbigo_X$
is the identity,
if the restriction $\varphi_{|X-D}$ is the identity.
Therefore,
we have only to show that 
the canonical isomorphism 
can be extended locally.

Then, we have only to consider the case 
$X=\Delta^n$ and $Y=\{z_1=0\}$.
Moreover, it can be reduced to the case $n=1$.
We consider the following exact sequence
on the left $D_{\Delta}$-modules on a disc $\Delta$:
\[
 0\lrarr D_{\Delta}\stackrel{\varphi}{\lrarr}
 D_{\Delta}\stackrel{\psi}{\lrarr}
 \nbigo_{\Delta}(\ast O)\lrarr 0,
\quad
\varphi(P_1)=P_1\cdot \del_t\cdot t,
\quad
\psi(P_2)=P_2\cdot t^{-1}
\]
Hence, as the right $D_{\Delta}$-module,
$\nrhom_{D_{\Delta}}
 (\nbigo_{\Delta}(\ast O),D_{\Delta})[1]$
is the cokernel of the following morphism:
\[
 D_{\Delta} \lrarr D_{\Delta},
\quad
 Q\longmapsto \del_t\cdot t\cdot Q
\]
We have the following exact sequence 
of the right $D$-modules:
\[
\begin{CD}
 0@>>> \del_t D_{\Delta}/\del_t t D_{\Delta}
 @>>> D_{\Delta}/\del_t t D_{\Delta}
 @>>> D_{\Delta}/\del_t D_{\Delta}
 @>>> 0\\
 @.@V{\simeq}VV @V{\simeq}VV @V{\simeq}VV @.\\
0@>>>\cnum[\del_t] @>>>
 D_{\Delta}/\del_t t D_{\Delta} @>>>
 \nbigo_{\Delta} @>>> 0
\end{CD}
\]
Then, it is easy to check
$D_{\Delta}/\del_t t D_{\Delta}\simeq
 j_!j^{\ast}\nbigo_{\Delta}$.
\hfill\qed

\vspace{.1in}

Let $M$ be a $D_X$-module.
We put 
$j_{\ast}j^{\ast}M:=
M\otimes^L\nbigo_X(\ast Y)$
and
$j_!j^{\ast}M:=
 M\otimes^Lj_!j^{\ast}\nbigo_{X}$,
where ``$\otimes^L$'' denotes 
the derived tensor product for $D$-modules.
\begin{lem}
\label{lem;08.1.19.2}
We have a natural isomorphism
$j_!j^{\ast}M
\simeq
 \DDD_X\bigl(j_{\ast}j^{\ast}\DDD_X(M)\bigr)$.
\end{lem}
\pf
Because
$j_{\ast}j^{\ast}\DDD_XM
\simeq\DDD_X(M)\otimes^L\nbigo_X(\ast Y)$,
we have
$\DDD_X\bigl(j_{\ast}j^{\ast}\DDD_X(M)\bigr)
 \simeq
 M\otimes^L \DDD_X\bigl(\nbigo_X(\ast Y)\bigr)$.
We also have the natural isomorphism
$j_!j^{\ast}M\simeq
 M\otimes^L j_!j^{\ast}\nbigo_X$.
Then, Lemma \ref{lem;08.1.19.2}
follows from Lemma \ref{lem;08.1.19.3}.
\hfill\qed

\vspace{.1in}

If $Y$ is non-characteristic to $M$,
we have $j_!j^{\ast}M=M\otimes j_!j^{\ast}\nbigo_X$.
In that case,
we have the following exact sequences:
\[
 0\lrarr M\lrarr j_{\ast}j^{\ast}M\lrarr 
 j_{\ast}j^{\ast}M\big/M\lrarr 0
\]
\[
 0\lrarr
 i_{\dagger}i^{\dagger}M\lrarr
 j_!j^{\ast}M\lrarr M\lrarr 0
\]

\subsection{Some commutative diagrams}
\label{subsection;08.9.29.210}

Let $X$ be a projective variety with 
an embedding $X\subset\proj^N$.
Let $X_0$ be the intersection $X\cap H_1\cap H_2$,
where $H_i$ denote general hyperplanes.
Let $\pi:\Xtilde\lrarr X$ be the blow up of $X$ 
along $X_0$.
We put $Y:=H_1\cap X$.
We have the natural embedding
$Y\subset \Xtilde$
whose image is denoted by $\Ytilde$.

We have the exact sequences
for $D_{\Xtilde}$-modules
as in Section \ref{subsection;08.1.19.20}:
\[
 0\lrarr \nbigo_{\Xtilde}
 \lrarr 
 \jtilde_{\ast}\jtilde^{\ast}\nbigo_{\Xtilde}
 \lrarr \jtilde_{\ast}\jtilde^{\ast}
 \nbigo_{\Xtilde}\big/\nbigo_{\Xtilde}
 \lrarr 0
\]
\[
 0\lrarr
 i_{\Ytilde\dagger}i_{\Ytilde}^{\dagger}\nbigo_{\Xtilde} 
 \lrarr \jtilde_!\jtilde^{\ast}\nbigo_{\Xtilde}
 \lrarr \nbigo_{\Xtilde}\lrarr 0
\]
Because
$\pi^{p}_{\dagger}\nbigo_{\Xtilde}=0$
and
$\pi^p_{\dagger}\bigl(
 i_{\Ytilde}\nbigo_{\Ytilde}
 \bigr)=0$ for $p\neq 0$,
we obtain 
$\pi^p_{\dagger}\bigl(
 \jtilde_{\ast}\jtilde^{\ast}\nbigo_{\Xtilde}
 \bigr)=0$
and $\pi^p_{\dagger}\bigl(
 \jtilde_{!}\jtilde^{\ast}\nbigo_{\Xtilde}
 \bigr)=0$ for $p\neq 0$.
We also recall the following lemma.
\begin{lem}
We have the following commutative diagram
of $D_X$-modules:
\begin{equation}
 \label{eq;08.1.18.52}
 \begin{CD}
 0 @>>> \nbigo_X @>>>
 j_{\ast}j^{\ast}\nbigo_X
 @>>> j_{\ast}j^{\ast}\nbigo_X/\nbigo_X
 @>>> 0\\
 @. @A{f_1}AA @A{f_2}AA @A{=}AA @.\\
 0 @>>> \pi_{\dagger}(\nbigo_{\Xtilde}) @>>> 
 \pi_{\dagger}
 (\jtilde_{\ast}\jtilde^{\ast}\nbigo_{\Xtilde})
 @>>>
 \pi_{\dagger}\bigl(
 \jtilde_{\ast}\jtilde^{\ast}
 \nbigo_{\Xtilde}\big/\nbigo_{\Xtilde}
 \bigr)
 @>>> 0
 \end{CD}
\end{equation}
Here, $f_1$ is the natural morphism.

By taking the dual,
we also have the following:
\begin{equation}
 \label{eq;08.1.18.53}
 \begin{CD}
 0 @>>>
 i_{Y\dagger}i_Y^{\dagger}\nbigo_X
 @>>>
 j_!j^{\ast}\nbigo_X
 @>>>
 \nbigo_X
 @>>> 0\\
 @. @V{=}VV @V{g_2}VV @V{g_1}VV @.\\
 0 @>>>
 \pi_{\dagger}\bigl(
 i_{\Ytilde\dagger}i_{\Ytilde}^{\dagger}\nbigo_{\Xtilde}
 \bigr)
 @>>>
 \pi_{\dagger}\bigl(
 \jtilde_!\jtilde^{\ast}
 \nbigo_{\Xtilde}
 \bigr)
 @>>>
 \pi_{\dagger}\bigl(
 \nbigo_{\Xtilde}
 \bigr)
 @>>> 0
 \end{CD}
\end{equation}
Here, $g_1$ is the natural morphism.
\end{lem}
\pf
We have the following commutative diagram:
\[
 \begin{CD}
 \pi_{\dagger}\nbigo_{\Xtilde}
 @>>>
 \pi_{\dagger}\nbigo_{\Xtilde}
 @>>>
 \nbigo_X\\
 @VVV @VVV @VVV \\
 \pi_{\dagger}\bigl(
 \jtilde_{\ast}\jtilde^{\ast}\nbigo_{\Xtilde}\bigr)
 @>>>
 \pi_{\dagger}\pi^{\dagger}
 (j_{\ast}j^{\ast}\nbigo_X)
 @>>>
 j_{\ast}j^{\ast}\nbigo_X
 \end{CD}
\]
Let $f_2$ be the composite of
the lower horizontal arrows.
Then, we have only to show 
that the induced map
$f_3:\pi_{\dagger}\bigl(
 \jtilde_{\ast}\jtilde^{\ast}\nbigo_{\Xtilde}/\nbigo_{\Xtilde}
 \bigr)
\lrarr j_{\ast}j^{\ast}\nbigo_X/\nbigo_X$
is the identity.
We have only to compare 
the restrictions of $f_3$ and the identity
to $X-X_0$.
Then, the claim is clear.
\hfill\qed

\section{Complement}
\label{subsection;08.10.18.110}
\subsection{Push-forward of $\nbigr$-modules
for ramified covering}
\label{subsection;08.1.20.100}

Let $X:=\Delta_z^n$,
$\Xtilde:=\Delta_{\zeta}^n$,
$D:=\bigcup_{i=1}^{\ell}\{z_i=0\}$,
and $\Dtilde:=\bigcup_{i=1}^{\ell}\{\zeta_i=0\}$.
Let $f:(\Xtilde,\Dtilde)\lrarr (X,D)$
be a ramified covering
given by $f(\zeta_1,\ldots,\zeta_n)
=(\zeta_1^{m_1},\ldots,\zeta_{\ell}^{m_{\ell}},
 \zeta_{\ell+1},\ldots,\zeta_n)$.
Let $\nbigm$ be an $\nbigr_{\Xtilde}$-module.
We obtain the $\nbigr_{\Xtilde}$-module
$\nbigm':=\nbigm\otimes_{\nbigo_{\nbigxtilde}}
 \nbigo_{\nbigxtilde}(\ast \nbigdtilde)$,
on which $\zeta_i$ $(i=1,\ldots,\ell)$ are invertible.
Let us consider the push forward via $f$.
\begin{lem}
\label{lem;08.1.16.25}
The push forward of $\nbigm'$
as $\nbigr$-modules
is isomorphic to 
the push forward
as $\nbigo$-modules.
Similar claim holds also for $D$-modules.
\end{lem}
\pf
We have only to consider
the case of right $\nbigr$-modules.
Let $\nbigc$ be an $\nbigr_{\Xtilde}$-free resolution of
$\nbigo_{\nbigxtilde}
 \otimes_{f^{-1}\nbigo_{\nbigx}}f^{-1}\nbigr_{X}$.
The push forward of $\nbigm'$ as $\nbigr$-modules
is given by the following:
\begin{equation}
 \label{eq;08.1.16.20}
 f_{\ast}\bigl(
 \nbigm'\otimes_{\nbigr_{\Xtilde}}\nbigc
 \bigr)
\end{equation}
Here $f_{\ast}$ denotes the push forward of
$\nbigo$-modules.
Let $\nbigr_{\Xtilde}(\ast \Dtilde):=
 \nbigr_{\Xtilde}\otimes_{\nbigo_{\nbigxtilde}}
 \nbigo_{\nbigxtilde}(\ast\nbigdtilde)$
and $\nbigr_X(\ast D):=
 \nbigr_{X}\otimes_{\nbigo_{\nbigx}}
 \nbigo_{\nbigx}(\ast\nbigd)$.
Because 
$\nbigm':=\nbigm\otimes_{\nbigo_{\nbigxtilde}}
 \nbigo_{\nbigxtilde}(\ast\nbigdtilde)$,
(\ref{eq;08.1.16.20}) is equal to
the following:
\begin{equation}
\label{eq;08.1.16.21}
 f_{\ast}\Bigl(
 \nbigm'\otimes_{\nbigr_{\Xtilde}(\ast\Dtilde)}
 \bigl(
 \nbigr_{\Xtilde}(\ast\Dtilde)
 \otimes_{\nbigr_{\Xtilde}}\nbigc
 \bigr)
\Bigr)
\end{equation}
Note $\nbigr_{\Xtilde}(\ast \Dtilde)
\otimes_{\nbigr_{\Xtilde}}\nbigc$
is quasi-isomorphic to
$\nbigr_{\Xtilde}(\ast \Dtilde)$.
Hence, (\ref{eq;08.1.16.21})
is naturally isomorphic to 
$f_{\ast}(\nbigm')$.
\hfill\qed

\subsection{Some sheaves of algebras}
\label{subsection;08.9.26.1}
Let $X:=\Delta^n$.
For any subset $J\subset\nbar$,
let $\lefttop{J}V_0\nbigr_X$ denote 
the sheaf of subalgebras of
$\nbigr_X$ generated by
$\nbigo_{\nbigx}$, $\deldel_jz_j$ $(j\in J)$
and $\deldel_j$ $(j\in \nbar-J)$.
\index{sheaf $\lefttop{J}V_0\nbigr_X$}
It is equipped with the filtration
by orders of differential operators.
The following lemma can be shown 
by the argument in the proof of
Proposition \ref{prop;07.7.30.3}.
\begin{prop}
\mbox{{}}\label{prop;08.9.26.10}
\begin{itemize}
\item
 The sheaf of algebras
 $\lefttop{I}V_0\nbigr_X$ is Noetherian,
 and it has the property in 
 Proposition {\rm\ref{prop;08.10.20.1}}.
 The Rees ring is also Noetherian.
\item 
 Let $\nbigm$ be an $\lefttop{I}V_0\nbigr_X$-module
 such that 
 (i) pseudo-coherent as an $\nbigo_{\nbigx}$-module,
 (ii) locally finitely generated as a
 $\lefttop{I}V_0\nbigr_X$-module.
Then, it is a coherent $\lefttop{I}V_0\nbigr_X$-module. 
\hfill\qed
\end{itemize}
\end{prop}

\section{The sheaves
 $\distribution_{X\times T/T}$
 and $\distribution^{\moderate}_{X\times T/T}$}
\label{subsection;08.1.15.20}
We recall some sheaves, following Sabbah
in \cite{sabbah4} and \cite{sabbah5},
to which we refer for more details and precise.
For simplicity, we consider the case in which 
$X$ is an open subset of $\cnum^n$.
For $J=(j_1,\ldots,j_n)$,
we put $\del^J:=\prod_{i=1}^n\del_i^{j_i}$
and $|J|:=\sum_{i=1}^n j_i$.
Let $V$ be an open subset of $X\times T$.
Let $\nbige^{(n,n)}_{X\times T/T,c}(V)$
denote the space of $C^{\infty}$-sections
of $\Omega^{n,n}_{X\times T/T}$ on $V$
with compact supports.
\index{sheaf $\nbige^{(n,n)}_{X\times T/T,c}$}
For any compact subset $K\subset V$
and $m\in\seisuu_{>0}$,
we have the semi-norm
$\bigl\|f\bigr\|_{m,K}=\sup_{|J|\leq m}\sup_K|\del_Jf|$.
For any closed subset $Z\subset X$,
let $\nbige^{<Z,(n,n)}_{X\times T/T,c}(V)$
denote the subspace of
$\nbige^{(n,n)}_{X\times T/T,c}(V)$,
which consists of the sections $f$
such that $(\del^Jf)_{|Z}=0$ for any $J$.
\index{sheaf $\nbige^{<Z,(n,n)}_{X\times T/T,c}$}
We have the induced semi-norms
$\|\cdot\|_{m,K}$ on the space 
$\nbige^{<Z,(n,n)}_{X\times T/T,c}(V)$.
By the semi-norms,
the spaces 
$\nbige^{(n,n)}_{X\times T/T,c}(V)$
and $\nbige^{<Z,(n,n)}_{X\times T/T,c}(V)$
are locally convex topological spaces.
Let $C^0_c(T)$ denote the space of 
continuous functions on $T$
with compact supports.
It is a normed vector space
with the sup norms.

Let $\distribution_{X\times T/T}(V)$ denote 
the space of continuous $C^{\infty}(T)$-linear maps 
from $\nbige^{(n,n)}_{X\times T/T,c}(V)$ to
$C^0_c(T)$.
\index{sheaf $\distribution_{X\times T/T}$}
In the case
$X=X_0\times\cnum_t$,
let $\distribution^{\moderate}_{X\times T/T}(V)$
denote the space of
continuous $C^{\infty}(T)$-linear maps
from $\nbige^{<X_0,(n,n)}_{X\times T/T,c}(V)$
to $C^0_c(T)$.
\index{sheaf $\distribution^{\moderate}_{X\times T/T}$}
Any elements of
$\distribution^{\moderate}_{X\times T/T}(V)$
are called distributions with moderate growth.
Let $\distribution_{X\times T/T,X_0}(V)$ be 
the space of continuous $C^{\infty}(T)$-linear maps
$\Phi:\nbige^{(n,n)}_{X\times T/T,c}(V)\lrarr
 C^0_c(T)$
whose supports are contained in $X_0\times T$,
i.e., 
$\Phi(f)=0$ if 
$f=0$ on some neighbourhood of $X_0\times T$.
\index{sheaf $\distribution_{X\times T/T,X_0}$}
They give the sheaves
$\distribution_{X\times T/T}$,
$\distribution^{\moderate}_{X\times T/T}$
and
$\distribution_{X\times T/T,X_0}$ on $X\times T$.
We have the following lemma
as in \cite{malgrange2}.

\begin{lem}
\label{lem;07.10.25.1}
The natural sequence
\[
 0\lrarr \distribution_{X\times T/T,X_0}(V)
\stackrel{\varphi}{\lrarr}
 \distribution_{X\times T/T}(V)
\stackrel{\psi}{\lrarr}
 \distribution^{\moderate}_{X\times T/T}(V)
\lrarr 0
\] is exact.
\end{lem}
\pf
We obtain the surjectivity of $\psi$
by using the Hahn-Banach theorem.
(See \cite{rudin}.)
The injectivity of $\varphi$ is clear.
Let $H$ be the space of functions
$f\in\nbige^{(n,n)}_{X\times T/T,c}(V)$
such that $f=0$ on some neighbourhood of $X_0\times T$.
We can show that $H$ is dense in
$\nbige^{<X_0,(n,n)}_{X\times T/T,c}(V)$
using an argument of Lemma 4.3 in \cite{malgrange2}.
Then, $\Image(\varphi)=\Ker(\psi)$ follows.
\hfill\qed

\vspace{.1in}
Note that $\distribution_{X\times T/T}(V)$
are naturally $\cnum[t]$-modules.
Let $\distribution_{X\times T/T}(\ast t)$
denote the sheafification of
presheaves
$V\longmapsto 
 \distribution_{X\times T/T}(V)
 \otimes_{\cnum[t]}\cnum[t,t^{-1}]$.

\begin{lem}
$\distribution^{\moderate}_{X\times T/T}(V)$
is naturally isomorphic to 
$\distribution_{X\times T/T}(\ast t)(V)$.
It is also isomorphic to the image of
$\distribution_{X\times T/T}(V)
\lrarr \distribution_{X\times T/T}(V\setminus X_0)$.
\end{lem}
\pf
Due to Lemma \ref{lem;07.10.25.1},
we have the natural injection
$\iota:\distribution^{\moderate}_{X\times T/T}(V)
\lrarr
 \distribution_{X\times T/T}(\ast t)(V)$.
It is easy to observe that
any element of $\distribution_{X\times T/T}(V)(\ast t)(V)$
naturally determines
a continuous $C^{\infty}(T)$-linear function
$\nbige^{<X_0(n,n)}_{X\times T/T,c}(V)
\lrarr C^0_c(T)$.
Hence, $\iota$ is surjective.
The second claim immediately follows
from Lemma \ref{lem;07.10.25.1}.
\hfill\qed

\vspace{.1in}
The following lemma can be shown
by using a standard argument.
\begin{lem}
Let $\Phi$ be an element of
$\distribution^{\moderate}_{X\times T/T}(V)$.
Let $K$ be any compact region.
We have some $m\in\seisuu_{>0}$
and $C>0$ such that
$ \sup_{T}\bigl|\Phi(f)\bigr|
\leq
 C\cdot\|f\|_{m,K}$
for any 
$f\in 
 \nbige^{<X_0\,(n,n)}_{X\times T/T,c}(V)$
with $\Supp(f)\subset K$.
\end{lem}
\pf
Assume that the claim does not hold.
We have a sequence $C_m\to\infty$
and $f_m\in \nbige^{<X_0\,(n,n)}_{X\times T/T,c}(V)$
with $\Supp(f_m)\subset K$
such that
$\sup_{T}\bigl|\Phi(f_m)\bigr|
>C_m\cdot\|f_m\|_{m,K}$.
For $g_m:=f_m\cdot \|f_m\|^{-1}\cdot C_m^{-1}$,
we have $\sup_T\bigl|\Phi(g_m)\bigr|>1$
and $\|g_m\|_{m,K}\to 0$.
Thus, we have arrived at the contradiction.
\hfill\qed

\begin{cor}
\label{cor;08.12.24.1}
Let $\Phi\in \distribution^{\moderate}_{X\times T/T}(V)$.
We have the well defined pairing
$\Phi\bigl(|t|^{s}\cdot f\bigr)\in C^0_c(T)$
for $f\in \nbige^{(n,n)}_{X\times T/T,c}(V)$
with $\Supp(f)\subset K$,
if $\Re(s)$ is sufficiently large.
\hfill\qed
\end{cor}

\subsection{Pull back via a ramified covering}

Let $\varphi_n:\cnum_{t_n}\lrarr\cnum_{t}$
be given by $\varphi_n^{\ast}(t)=t_n^n$.
Let $X^{(n)}:=X_0\times \cnum_{t_n}$,
and let $\varphi_n:X^{(n)}\times T\lrarr X\times T$
be the induced map.
Let $V^{(n)}:=\varphi_n^{-1}(V)$.
We have the well defined continuous map
\[
\varphi_{n\,\ast}:
 \nbige^{<X_0\,(n,n)}_{X^{(n)}\times T/T,c}(V^{(n)})
\lrarr
 \nbige^{<X_0\,(n,n)}_{X\times T/T,c}(V). 
\]
Hence, we have the induced map
\[
\varphi_n^{\ast}:
 \distribution^{\moderate}_{X\times T/T}(V)
\lrarr
 \distribution^{\moderate}
 _{X^{(n)}\times T/T}(V^{(n)}). 
\]
We have the following general lemma.
\begin{lem}
\label{lem;07.10.26.1}
Let $\rho$ be a test function on $X_0$,
and let $\tau\in
 \distribution^{\moderate}_{X\times T/T}(V)$.
If $m\not\equiv 0$ modulo $n$,
we have the vanishing
$\big\langle
 \varphi_n^{\ast}\tau,\,
 |t_n|^{2s}t_n^m
 \varphi_n^{\ast}(\chi)\cdot\rho
 \big\rangle=0$
for any $s$ such that $\Re(s)$ is sufficiently large.
\end{lem}
\pf
We have only to consider 
$t_n\longmapsto t_n\cdot a$
for primitive $n$-th root $a$.
\hfill\qed

\subsection{Exponential twist}
Sabbah observed the following in \cite{sabbah5}.
\begin{lem}
\label{lem;07.10.26.3}
Let
$\Phi\in\distribution^{\moderate}_{X\times T/T}(V)$
and $\gminia\in\cnum[t^{-1}]$.
Then,
\[
 \exp\bigl(2\sqrt{-1}\Image(\lambda\gminia)\bigr)
 \Phi
 \in \distribution_{X\times T/T}(V\setminus X_0)
\]
is contained in 
$\distribution^{\moderate}_{X\times T/T}(V)$.
\end{lem}
\pf
For any $f\in \nbige^{<X_0,(n,n)}_{X\times T/T,c}(V)$,
we have
$\exp\bigl(2\sqrt{-1}\Image(\lambda\gminiabar)\bigr)f
 \in \nbige^{<X_0,(n,n)}_{X\times T/T,c}(V)$.
Hence, the claim of the lemma follows.
\hfill\qed

\section{$\nbigr$-triples}
\label{subsection;08.10.18.112}
\subsection{Hermitian sesqui-linear pairing}

We recall the notion of $\nbigr$-triples
in \cite{sabbah2} for the convenience of readers.
We refer to it for more details and precision.
Let $X$ be a complex manifold.
We put $\nbigx:=\cnum_{\lambda}\times X$.
We set 
$\vecS:=\bigl\{
 \lambda\in\cnum\,\big|\,|\lambda|=1
 \bigr\}$.
\index{set $\vecS$}
Let $\distribution_{\vecS\times X/\vecS}$
be as in Section \ref{subsection;08.1.15.20}.
Let $\sigma:
 \cnum_{\lambda}^{\ast}\lrarr
 \cnum_{\lambda}^{\ast}$
given by $\sigma(\lambda)=-\lambdabar^{-1}$.
\index{map $\sigma$}
The induced map
$\vecS\times X\lrarr\vecS\times X$ is also denoted by
$\sigma$.
We have the natural 
$(\nbigr_{X})_{|\vecS\times X}$-action
on $\distribution_{\vecS\times X/\vecS}$.
We also have the 
$(\sigma^{\ast}\nbigr_{X})_{|\vecS}$-action
on $\distribution_{\vecS\times X/\vecS}$
locally given by
$\sigma^{\ast}(\deldel_i)=
 -\lambda^{-1}\delbar_i=:\deldelbar_i$.
Thus, $\distribution_{\vecS\times X/\vecS}$
is a left $(\nbigr_X\otimes
 \sigma^{\ast}\nbigr_{X})_{|\vecS\times X}$-module.

Let $\nbigm'$ and $\nbigm''$ be left $\nbigr$-modules
on $\nbigx$.
Recall that a Hermitian sesqui-linear pairing of 
$\nbigm'$ and $\nbigm''$ is defined to be
a morphism of left
$(\nbigr_X\otimes
 \sigma^{\ast}\nbigr_X)_{|\vecS\times X}$-modules:
\index{Hermitian sesqui-linear pairing}
\index{sesqui-linear pairing}
\[
 C:\nbigm'_{|\vecS\times X}
 \otimes\sigma^{\ast}\nbigm''_{|\vecS\times X}
\lrarr \distribution_{\vecS\times X/\vecS}.
\]
It is also called a sesqui-linear pairing.
We often denote the pairing $C(x,\sigma^{\ast}y)$
by $C(x,\ybar)$ as in \cite{sabbah2}.
See \cite{sabbah2}
for a sesqui-linear pairing of right $\nbigr$-modules.

\vspace{.1in}

Let $C$ be a sesqui-linear pairing of
$\nbigm'$ and $\nbigm''$.
Then a sesqui-linear pairing
$C^{\ast}$ of $\nbigm''$ and $\nbigm'$ is 
given by the following formula:
\begin{equation} 
\label{eq;05.1.14.1}
 C^{\ast}(x,\sigma^{\ast}y)=
 \overline{ \sigma^{\ast}C(y,\sigma^{\ast}x)}.
\end{equation}
Here $x$ and $y$ denote local sections
of $\nbigm''_{|\vecS\times X}$ and 
$\nbigm'_{|\vecS\times X}$ respectively.

The following lemma is proved in \cite{sabbah2}.
\begin{lem}
\label{lem;05.1.16.400}
Let $\nbigm_i$ $(i=1,2)$
be strictly $S$-decomposable $\nbigr_X$-modules
whose strict supports are $Z_i$.
If $Z_1\neq Z_2$,
then there does not exist 
any non-trivial sesqui-linear pairing
of $\nbigm_1$ and $\nbigm_2$.
\hfill\qed
\end{lem}

\subsection{$\nbigr_X$-triple}

A left $\nbigr_X$-triple is a tuple $(\nbigm',\nbigm'',C)$
of $\nbigr_X$-modules  $\nbigm'$, $\nbigm''$
and a sesqui-linear pairing $C$ of $\nbigm'$ and $\nbigm''$.
Let $\nbigt_i=(\nbigm_i',\nbigm_i'',C_i)$ $(i=1,2)$ be
$\nbigr_X$-triples.
A morphism $\varphi:\nbigt_1\lrarr \nbigt_2$ is a tuple
of morphisms
$\varphi':\nbigm_2'\lrarr\nbigm_1'$ and
$\varphi'':\nbigm_1''\lrarr\nbigm_2''$
satisfying
\[
 C_1\bigl(\varphi'(x),\sigma^{\ast}(y)\bigr)
=C_2\bigl(x,\sigma^{\ast}(\varphi''(y))\bigr).
\]
Here $x$ and $y$ denote local sections
of $\nbigm'_{2|\vecS\times X}$
and $\nbigm''_{1|\vecS\times X}$.
The category 
$\nbigr_{X}{\textrm -Triples}$
of $\nbigr_X$-triples is abelian.
\index{category $\nbigr_{X}{\textrm -Triples}$}

\subsection{Tate twist and Hermitian adjoint}
\index{Tate twist}
\index{Hermitian adjoint}
\index{Tate object $\Tate^S(k)$}

For any half integer $k\in\frac{1}{2}\seisuu$,
the $k$-th Tate object $\Tate^S(k)$ is defined to be
the tuple $\bigl(\nbigo_{\nbigx},\nbigo_{\nbigx},C_k\bigr)$,
where $C_k$ is given as follows:
\[
 C_k(f,\sigma^{\ast}g)=
\bigl(\sqrt{-1}\lambda\bigr)^{-2k}\cdot
 f\cdot \overline{\sigma^{\ast}g}
\]
The $k$-th Tate twist $\nbigt(k)$
of $\nbigt=(\nbigm',\nbigm'',C)$ is
the tuple
\[
 \bigl(\nbigm',\nbigm'',(\sqrt{-1}\lambda)^{-2k}C \bigr)
\]
In other words, 
it is the tensor product
$\nbigt\otimes\Tate^S(k)$.
We will not distinguish $\nbigt(k)$
and $\nbigt\otimes\Tate^S(k)$.
A morphism $\varphi=(\varphi',\varphi''):\nbigt\lrarr\nbigt'$
naturally induces the morphism
$\varphi=(\varphi',\varphi''):
 \nbigt(k)\lrarr\nbigt'(k)$ for any $k\in\frac{1}{2}\seisuu$.

\vspace{.1in}

For a left $\nbigr$-triple $\nbigt=(\nbigm',\nbigm'',C)$,
the Hermitian adjoint $\nbigt^{\ast}$ is defined to be
the tuple $(\nbigm'',\nbigm',C^{\ast})$,
where $C^{\ast}$ is given as in (\ref{eq;05.1.14.1}).
For any morphism 
$\varphi=(\varphi',\varphi''):
 \nbigt_1\lrarr \nbigt_2$,
we obtain the morphism
$\varphi^{\ast}=(\varphi'',\varphi'):\nbigt_2\lrarr\nbigt_1$.

The Hermitian adjoint of the Tate object $\Tate^S(-k)$ is 
$\bigl(\nbigo_{\nbigx},\nbigo_{\nbigx},(-1)^{-2k}C_{k}\bigr)$.
We fix the isomorphism 
\[
\bigl((-1)^{2k}\id,\id\bigr):
\Tate^S(k)\simeq \Tate^S(-k)^{\ast}.
\]
It induces the isomorphism
$\nbigt(k)\simeq \bigl(\nbigt^{\ast}(-k)\bigr)^{\ast}$.

\subsection{Hermitian sesqui-linear duality 
of $\nbigr$-triples}
\index{Hermitian sesqui-linear duality}

Let $\nbigt$ be an $\nbigr_X$-triple,
and let $w$ be an integer.
Recall that
a morphism $\nbigs:\nbigt\lrarr \nbigt^{\ast}(-w)$
is called a Hermitian sesqui-linear duality of weight $w$,
if the equality $\nbigs=(-1)^w\nbigs^{\ast}$ holds.
In the case $\nbigs=(S',S'')$,
the condition is equivalent to $S'=(-1)^w S''$.

\vspace{.1in}
Let $\nbigs=(S',S'')$ be a morphism
$\nbigt\lrarr \nbigt^{\ast}(-w)$.
Via the canonical isomorphism
$\Tate^S(k)\simeq \Tate^S(-k)^{\ast}$,
we obtain the map
$\nbigt(k)\lrarr \nbigt(k)^{\ast}(-w+2k)$,
which is given by $\nbigs(k):=((-1)^{2k}S',S'')$.
Then, $\nbigs$ is a Hermitian sesqui-linear duality 
of weight $w$,
if and only if $\nbigs(k)$ is 
a Hermitian sesqui-linear duality
of weight $w+k$.

\vspace{.1in}

Let $\nbigt=(\nbigm',\nbigm'',C)$ be an $\nbigr_X$-triple
with a Hermitian sesqui-linear duality 
$\nbigs=(S',S'')$ of weight $w$.
Later, we will mainly interested in the case 
that $\nbigs$ is an isomorphism,
i.e., $S'$ and $S''$ are isomorphisms.
Then we have the $\nbigr_X$-triple
$\widetilde{\nbigt}=(\nbigm'',\nbigm'',\widetilde{C})$
and the canonical isomorphism
$(S',\id):\nbigt\lrarr\widetilde{\nbigt}$,
where $\widetilde{C}$ is given by
$\widetilde{C}(x,\sigma^{\ast}y)
=C(S'x,\sigma^{\ast}y)$ for local sections
$x$ and $y$ of $\nbigm''$.

For an $\nbigr_X$-triple $\nbigt$
of the form $(\nbigm',\nbigm',C)$,
an isomorphism $(\id,\id):\nbigt\lrarr\nbigt^{\ast}$
is a Hermitian sesqui-linear duality of weight $0$,
if and only if the following equality holds
for local sections $x$ and $y$ on appropriate open subsets:
\[
 C(x,\sigma^{\ast}y)
=C^{\ast}(x,\sigma^{\ast}y)
:=\sigma^{\ast}\overline{C(y,\sigma^{\ast}x)}.
\]

\subsection{Push forward}

Let $f:X\lrarr Y$ be a morphism of complex manifolds.
For an $\nbigr_X$-triple
$\nbigt=(\nbigm',\nbigm'',C)$,
the push-forward
$f_{\dagger}\nbigt\in 
 D\bigl(\nbigr_{Y}{\textrm -Triples}\bigr)$ 
is defined.
More precisely, they are given in the level of
the complexes of $\nbigr_Y$-modules.
We omit the details,
and refer to \cite{sabbah2} 
for more details.

\section{Specialization of $\nbigr$-triples}
\label{subsection;08.10.18.111}
We recall the specialization of sesqui-linear pairings
introduced by Sabbah \cite{sabbah2},
which we refer to for more details and precision.

\subsection{Specialization along
a coordinate function}
\label{subsection;08.12.16.1}

\subsubsection{Preliminary I}

Let $\cnum_t$ be a complex line 
with a coordinate $t$.
Let $X_0$ be an $(n-1)$-dimensional complex manifold.
We put $X:=X_0\times \cnum_t$.
We identify $X_0$ and $X_0\times\{0\}$.
Let $\nbigm'$ and $\nbigm''$ be objects of $\nbigs^2(X,t)$.
Let $C$ be a sesqui-linear pairing of
$\nbigm'$ and $\nbigm''$.
Let $u$ be any element of $\real\times\cnum$.
We recall the construction due to Sabbah to obtain
the specialization along $t$:
\[
 \psitilde_{t,u}C:
 \bigl(
 \psitilde_{t,u}\nbigm'
 \otimes
 \sigma^{\ast}{\psitilde_{t,u}\nbigm''}
\bigr)_{|\vecS\times X_0}
\lrarr
 \distribution_{\vecS\times X_0/\vecS}
\]

For $\lambda_0\in\vecS$,
let $U$ be a small neighbourhood of $\lambda_0$,
and $\vecI:=U\cap\vecS$.
Let $W_0$ be an open subset of $X_0$.
We put $W:=W_0\times \cnum_t$.
Let $m$ be a section of $\nbigm'$ 
on $U\times W$,
and $\mu$ be a section of $\nbigm''$
on $\sigma(U)\times W$.
Let us take any $C^{\infty}$ $(n-1,n-1)$-form 
$\phi$ on $W_0$
whose support is compact.
We also take a $C^{\infty}$-function $\chi$ 
on $\cnum_t$
with the compact support
such that $\chi=1$ around the origin $O\in\cnum_t$.
By considering the push forward for
$\vecS\times X\lrarr\vecS$,
we obtain the following continuous function 
on $\vecI\times\bigl\{
 s\in\cnum\,\big|\,2\Re(s)+|k|>R_0 \bigr\}$,
which is holomorphic with respect to $s$:
\begin{equation} 
 \label{eq;05.1.14.10}
 \nbigi^{(k)}_{C(m,\overline{\mu}),\phi}(s)
:=\left\{
 \begin{array}{ll}
 \bigl\langle C(m,\overline{\mu}),
 \,\,|t|^{2s} t^k\cdot\chi(t)\cdot\phi\wedge
  \frac{\sqrt{-1}}{2\pi}dt\wedge d\overline{t}\bigr\rangle
 & (k\geq 0)\\
\mbox{{}}\\
 \bigl\langle C(m,\overline{\mu}),
 \,\,|t|^{2s}\overline{t}^{-k}\cdot\chi(t)\cdot\phi\wedge
  \frac{\sqrt{-1}}{2\pi}dt\wedge d\overline{t}\bigr\rangle
 & (k\leq 0)
\end{array}
\right.
\end{equation}
Here, $R_0$ denotes some real number
depending only on $m$ and $\mu$.
(See Corollary \ref{cor;08.12.24.1}.)
The following lemma can be shown
by the argument in \cite{sabbah2}
or the proof of
Lemma 14.72 and Lemma 14.73
of \cite{mochi2}.
\begin{lem}
 \label{lem;9.17.50}
Let $k\geq 0$.
We have the following formula:
\begin{equation}\label{eq;9.17.42}
 \lambda\cdot
 \bigl(
 s+k+\lambda^{-1}\eigenmap(\lambda,u)
 \bigr)\cdot
 \nbigi^{(k)}_{C(m,\overline{\mu}),\phi}(s)
=\nbigi^{(k)}_{C(m',\overline{\mu}),\phi}(s)
+F
\end{equation}
Here we put $m':=(-\deldel_tt+\eigenmap(\lambda,u))\cdot m$,
and $F$ is a continuous function on
$\vecI\times \cnum$
which is holomorphic with respect to the variable $s$.

We also have the following equality:
\begin{equation} 
 \label{eq;9.17.43}
 -\lambda^{-1}\cdot
\Bigl(s+k+\frac{\eigenmap(\lambda,u)}{\lambda}\Bigr)\cdot
\nbigi^{(-k)}_{C(m,\overline{\mu}),\phi}(s)
=\nbigi^{(-k)}_{C(m,\overline{\mu}'),\phi}
+F
\end{equation}
Here we put 
 $\mu'=\bigl(-\deldel_tt
 +\eigenmap(\lambda,u)\bigr)\cdot\mu$,
and $F$ is as above.
\hfill\qed
\end{lem}

\subsubsection{Preliminary II}
Let $m$ be a section of 
$V^{(\lambda_0)}_c\nbigm'$
such that 
$0\neq \pi_c(m)\in 
 \psi^{(\lambda_0)}_{t,u_0}\nbigm'$
via the projection
$\pi_c:\Vzero_c\nbigm'\lrarr \Gr^{\Vzero}_c\nbigm'$
for some $u_0\in\real\times\cnum$ such that 
$\paramap(\lambda_0,u_0)=c$.
Let $b_m$ be a polynomial
such that
(i) $b_m(-\deldel_t t)m\in V_{-1}\nbigr\cdot m$,
(ii) it is of the following form:
\[
 b_m(x)=
 \bigl(x+\eigenmap(\lambda,u_0)\bigr)^{\nu(u_0)}
\cdot
 \prod_{u\in S_0}
 \bigl(x+\eigenmap(\lambda,u)\bigr)^{\nu(u)}.
\]
Here  $S_0$ denotes a finite subset of $\real\times\cnum$
such that $\paramap(\lambda_0,u)<c$ for any $u\in S_0$.
For any positive integer $M$, we put
\[
 B_m^{(M)}(x):=
 \prod_{\nu=0}^{M-1} b_m(x+\nu\lambda).
\]
By construction,
there exists a finite subset $S_1(M)\subset \real\times\cnum$
such that the following holds:
\begin{itemize}
\item
 $B_m^{(M)}(x)=\prod_{u\in S_1(M)}
  \bigl(x+\eigenmap(\lambda,u)\bigr)^{\nu'(u)}$.
\item
 For any $u\in S_1(M)$, 
 we have $\paramap(\lambda_0,u)\leq c$.
 If $\paramap(\lambda_0,u)=c$, then $u=u_0$.
\end{itemize}
Moreover,
$\nu'$ is bounded on $\bigcup_{M}S_1(M)$.

\begin{lem} 
\label{lem;9.17.51}
We have the following equality:
\begin{equation}
\label{eq;9.19.3}
 \Bigl(
 \prod_{u\in S_1(M)}
 \lambda^{\nu'(u)}
 \bigl(s+k+\lambda^{-1}
 \eigenmap(\lambda,u)\bigr)^{\nu'(u)}
 \Bigr)
 \cdot \nbigi^{(k)}_{C(m,\overline{\mu}),\phi}(s) 
=\nbigi^{(k)}_{C(m',\overline{\mu}),\phi}
+F
\end{equation}
Here $m'= B_m^{(M)}(-\deldel_tt)\cdot m$,
and $F$ has the property 
as in Lemma {\rm\ref{lem;9.17.50}}.
The first term in the right hand side is
holomorphic with respect to the variable $s$
on the half plane
$\bigl\{s\in\cnum\,\big|\,2\Re(s)+k>R_0-M\bigr\}$.
\end{lem}
\pf
We obtain the equality (\ref{eq;9.19.3}) 
by using Lemma \ref{lem;9.17.50} inductively.
By the construction of $B_m^{(M)}$,
we have 
$ B_m^{(M)}(-\deldel_tt)\cdot m
=t^{M}P\cdot m$
for some $P\in V_0\nbigr_{X}$.
Hence, $\nbigi^{(k)}_{C(m',\overline{\mu}),\phi}$
in (\ref{eq;9.19.3}) is given on
$\vecI\times
 \bigl\{s\in\cnum\,\big|\,2\Re(s)+k>R_0-M\bigr\}$.
\hfill\qed

\vspace{.1in}
Then, we obtain the following lemma.

\begin{lem}
\mbox{{}} \label{lem;a11.18.1}
There exist a discrete subset $S_2$ of $\real\times\cnum$ 
and a number $\nu$ such that
the following holds:
\begin{itemize}
\item
$\nbigi^{(k)}_{C(m,\overline{\mu}),\phi}\cdot
 \prod_{u\in S_2}
 \bigl(s+k+\lambda^{-1}\cdot\eigenmap(\lambda,u)\bigr)^{\nu}$ 
is continuous on $\vecI\times \cnum$,
and it is holomorphic with respect to $s$.
\item
For any element $u\in S_2$,
we have $\paramap(\lambda_0,u)\leq c$.
If $\paramap(\lambda_0,u)=c$,
then $u=u_0$.
\hfill\qed
\end{itemize}
\end{lem}

Note $\sigma(\lambda_0)=-\lambda_0$
for $\lambda_0\in\vecS$.
Let $\mu$ be an element of 
$V^{(-\lambda_0)}_d\nbigm''$
such that
$0\neq \pi_d(\mu)\in 
 \psi_{t,u_1}^{(-\lambda_0)}\nbigm''$
via the projection
$V^{(-\lambda_0)}_d\nbigm''
\lrarr
\Gr^{V^{(-\lambda_0)}}_d\nbigm''$,
where $d=\paramap(-\lambda_0,u_1)$.
We obtain the following lemma
in a similar way.
\begin{lem}
\label{lem;04.2.17.2}
There exist a discrete subset $S_3\subset\real\times\cnum$
and a number  $\nu$ with the following property:
\begin{itemize}
\item
$\nbigi^{(-k)}_{C(m,\overline{\mu}),\phi}\cdot
 \prod_{u\in S_3}
 \bigl(s+k+\lambda^{-1}
 \eigenmap(\lambda,u)\bigr)^{\nu}$ 
 is given on $\vecI\times \cnum$,
and it is holomorphic with respect to $s$.
\item
For any $u\in S_3$, we have
$\paramap(-\lambda_0,u)\leq d$.
If $\paramap(-\lambda_0,u)=d$,
then $u=u_0$.
\hfill\qed
\end{itemize}
\end{lem}

\subsubsection{Construction of 
the specialization $\psi_{t,u}C$}

Let $[m]$ be a section of 
$\psizero_{t,u_0}\nbigm'$ on 
$U\times W_0$,
and $m$ be a section of $\nbigm'$ on 
$U\times W_0\times\cnum_t$
such that $\pi_c(m)=[m]$,
where we put $c:=\paramap(\lambda_0,u_0)$,
and $\pi_c$ denotes the projection 
$\Vzero_c\nbigm'\lrarr \Gr^{\Vzero}_c\nbigm'$.
Let $[\mu]$ be a section of 
$\psi^{(-\lambda_0)}_{t,u_0}\nbigm''$
on $\sigma(U)\times W_0$, 
and $\mu$ be a section of $\nbigm''$
on $\sigma(U)\times W_0\times\cnum_t$
such that $\pi_d(\mu)=[\mu]$,
where we put 
$d=\paramap\bigl(-\lambda_0,u_0\bigr)$,
and $\pi_d$ denotes the projection
$V^{(-\lambda_0)}_d\nbigm''
 \lrarr \Gr^{V^{(-\lambda_0)}}_d\nbigm''$.

We put
$\nbigi_{C(m,\overline{\mu}),\phi}(s)
:=\nbigi^{(0)}_{C(m,\overline{\mu}),\phi}(s)$.
Recall that we have a discrete subset $S$ 
of $\real\times\cnum$
such that
$\nbigi_{C(m,\overline{\mu}),\phi}(s)
 \prod_{u\in S}
\bigl(s+\lambda^{-1}\,\eigenmap(\lambda,u)\bigr)^{\nu}$
is given on $\vecI\times \cnum$.
For any $u\in S$,
we have 
$\paramap(\lambda_0,u)\leq 
 \paramap(\lambda_0,u_0)$
and $\paramap(-\lambda_0,u)
\leq
 \paramap(-\lambda_0,u_0)$.
If one of the equalities holds, $u=u_0$.
Hence, we have 
$\eigenmap(\lambda_0,u)\neq \eigenmap(\lambda_0,u_0)$
for any $u\in S\setminus \{u_0\}$.
(See Lemma 14.80 in \cite{mochi2},
 for example.)

\vspace{.1in}

We would like to define
\index{pairing $\psizero_{t,u_0}C$}
\begin{equation}
 \label{eq;9.19.5}
\big\langle
 \psi^{(\lambda_0)}_{t,u_0}
 C\bigl([m],[\overline{\mu}]\bigr),\,
 \phi
\big\rangle
:=
\underset{s+\lambda^{-1}\eigenmap(\lambda,u_0)}{\Res}
 \bigl(
 \nbigi_{C(m,\overline{\mu}),\phi}(s)
 \bigr).
\end{equation}
The right hand side of {\rm(\ref{eq;9.19.5})} gives
a continuous function on $\vecI$.
The well definedness of
{\rm(\ref{eq;9.19.5})} can be checked
as in \cite{sabbah2} or 
Lemma 14.82 of \cite{mochi2}.
Thus, we obtain
the local section 
$\psi^{(\lambda_0)}_{t,u_0}
C\bigl([m],[\overline{\mu}]\bigr)$
of $\distribution_{\vecI\times X_0/\vecI}$.
By varying $\lambda_0\in\vecS$ and gluing them,
we obtain 
\[
 \psi_{t,u_0}C:
 \psi_{t,u_0}\nbigm'
\otimes
 \sigma^{\ast}{\psi_{t,u_0}\nbigm''}
\lrarr\distribution_{\vecS\times X_0/\vecS}.
\]
We can easily check that
$\psi_{t,u_0}C$ gives a Hermitian sesqui-linear pairing
of 
$\psi_{t,u_0}\nbigm'$
and $\psi_{t,u_0}\nbigm''$.

We put $N:=-\deldel_tt+\eigenmap(\lambda,u_0)$,
which induces the nilpotent map on
$\psi_{t,u_0}\nbigm'$
and $\psi_{t,u_0}\nbigm''$.
As in \cite{sabbah2} or Lemma 14.84 of \cite{mochi2},
we have the following equality:
\[
 \psi_{t,u_0}C\bigl(N[m],[\overline{\mu}]\bigr)
=(\sqrt{-1}\lambda)^2\cdot
  \psi_{t,u_0}C\bigl([m],[\overline{N\mu}]\bigr)
\]

\subsubsection{The induced pairing $\psitilde_{t,u}C$}

We construct the induced pairing:
\[
 \tildepsi_{t,u}C:
 \tildepsi_{t,u}\nbigm'
 \otimes
 \sigma^{\ast}
 \tildepsi_{t,u}\nbigm''
\lrarr\distribution_{\vecS\times X_0/\vecS}
\]
For local sections $m$ and $\mu$ 
of $\psizero_{t,u}(\nbigm')$
and $\psi^{(-\lambda_0)}_{t,u}(\nbigm'')$,
we have the following equality
as in Lemma 14.85:
\begin{equation}
 \label{eq;9.25.12}
 \psizero_{t,u-\vecdelta_0}C\bigl(
 [t\cdot m],\overline{[t\cdot\mu]}
 \bigr)
=\psizero_{t,u}C\bigl(
 [m],\overline{[\mu]}
 \bigr)
\end{equation}
Here, we put $\vecdelta_0=(1,0)\in\real\times\cnum$.

\vspace{.1in}

Let $m$ be a section of 
$\tildepsizero_{t,u}(\nbigm')$,
and $\mu$ be a section of 
$\tildepsi^{(-\lambda_0)}_{t,u}(\nbigm'')$.
We pick a sufficiently large integer $N$,
$m_1\in\psizero_{t,u-N\vecdelta_0}(\nbigm')$
and $\mu_1\in 
 \psi^{(-\lambda_0)}_{t,u-N\vecdelta_0}(\nbigm'')$
corresponding to $m$ and $\mu$ respectively.
Then, the pairing 
$\tildepsizero_{t,u} C(m,\overline{\mu})$
is defined to be
$\psizero_{t,u-N\vecdelta_0} C(m_1,\overline{\mu_1})$.
By varying $\lambda_0$ and gluing 
$\tildepsizero_{t,u}C$,
we obtain $\psitilde_{t,u}C$.
\index{pairing $\tildepsizero_{t,u} C$}

\subsubsection{The nearby cycle functor of
 a strictly specializable $\nbigr$-triple}

\index{pairing $\psi_{t,u}C$}
\index{pairing $\tildepsi_{t,u}C$}
\index{functor $\psi_{t,u}$}
\index{functor $\psizero_{t,u}$}

Let $\nbigt=(\nbigm',\nbigm'',C)$ be an $\nbigr_X$-triple.
We say that $\nbigt$ is strictly specializable
along $t$,
if $\nbigm'$ and $\nbigm''$ are strictly specializable along $t$.

Let $\nbigt$ be an $\nbigr_X$-triple
which is strictly specializable along $t$.
For any $u\in\real\times\cnum$,
we obtain the following induced sesqui-linear pairings:
\[
 \psi_{t,u}C:\bigl(
 \psi_{t,u}(\nbigm')
\otimes
 \sigma^{\ast}{\psi_{t,u}(\nbigm'')}
 \bigr)_{|\vecS\times X_0}
\lrarr
 \distribution_{\vecS\times X_0/\vecS}
\]
Thus, we obtain the $\nbigr_{X_0}$-triple
$\psi_{t,u}(\nbigt)=\bigl(
 \psi_{t,u}\nbigm',
 \psi_{t,u}\nbigm'',
 \psi_{t,u}C
 \bigr)$.
We also have the following modified pairing:
\[
 \tildepsi_{t,u} C:
 \bigl(
 \tildepsi_{t,u}(\nbigm')
\otimes
 \sigma^{\ast}{\tildepsi_{t,u}(\nbigm'')}
\bigr)_{|\vecS\times X_0}
\lrarr
 \distribution_{\vecS\times X_0/\vecS}
\]
Thus, we obtain the $\nbigr_{X_0}$-triple
$ \psitilde_{t,u}(\nbigt):=
 \bigl(\tildepsi_{t,u}\nbigm',\tildepsi_{t,u}\nbigm'',
 \tildepsi_{t,u}C\bigr)$.

We have the following relations:
\[
\begin{array}{l}
 \tildepsi_{t,u}C\bigl(N[m],[\overline{\mu}]\bigr)
=(\sqrt{-1}\lambda)^2\cdot
  \tildepsi_{t,u}C\bigl([m],\overline{N}[\overline{\mu}]\bigr)
 \\ \mbox{{}} \\
\psi_{t,u}C\bigl(N[m],[\overline{\mu}]\bigr)
=(\sqrt{-1}\lambda)^2\cdot
  \psi_{t,u}C\bigl([m],\overline{N}[\overline{\mu}]\bigr)
\end{array}
\]
Here $N:=-\deldel_tt+\eigenmap(\lambda,u)$.
Hence, we have  the following morphism
\[
 \nbign:=(-\sqrt{-1}N,\sqrt{-1}N):
\psitilde_{t,u}\nbigt
 \lrarr\psitilde_{t,u}\nbigt\otimes\Tate^S(-1).
\]
It induces the weight filtrations $W$ on 
$\tildepsi_{t,u}\nbigt$ or $\psi_{t,u}\nbigt$.
The primitive part of Gr is denoted by 
$P\Gr^W_h\tildepsi_{t,u}\nbigt$.
\index{$\nbigr$-triple
 $P\Gr^W_h\tildepsi_{t,u}\nbigt$}

\subsubsection{Vanishing cycle functor}

The specialization $\psi_{t,0}C=\psitilde_{t,0}C$ 
is not appropriate for our use.
We recall the construction of 
the Hermitian sesqui-linear pairing
$\phi_{t,0}C$ of
$\psi_{t,0}\nbigm'$ and
$\psi_{t,0}\nbigm''$
only in the simple case 
that $\nbigm'$ and $\nbigm''$ are strictly
$S$-decomposable along $t$,
by following the older version of \cite{sabbah2}.
See \cite{sabbah2} for a more general case.

Since we have assumed that 
$\nbigm'$ is strictly $S$-decomposable along $t$,
we have the decomposition 
$\nbigm'=\nbigm'_1\oplus \nbigm'_2$
as in the claim (\ref{number;05.1.15.100}) of
Proposition \ref{prop;a11.23.1}.
Similarly we have the decomposition 
$\nbigm''=\nbigm''_1\oplus\nbigm''_2$.
We also have the decomposition
of the Hermitian sesqui-linear pairing $C=C_1\oplus C_2$,
where $C_i$ are pairings
of $\nbigm_i'$ and $\nbigm_i''$.
(See Lemma \ref{lem;05.1.16.400}.)
Since the supports of $\nbigm'_2$ and $\nbigm''_2$ 
are contained in $X_0$,
we have $\psi_{t,0}\nbigm'_2=\nbigm'_2$
and $\psi_{t,0}\nbigm''_2=\nbigm''_2$.
Therefore, we put $\phi_{t,0}C_2:=C_2$.
We have only to define $\phi_{t,0}C_1$.
Hence, we may assume
$\psi_{t,0}=\Image \can$
for  $\nbigm'$ and $\nbigm''$
from the beginning.
Recall that $\can$ and $\var$ are induced by 
the left action of $-\deldel_t$ and $t$.
For $x\in\psi_{t,0}\nbigm'$ 
and $y\in \psi_{t,0}\nbigm''$,
we define
\[
 \phi_{t,0}C(x,y)
:=(\sqrt{-1}\lambda)^{-1}
 \psi_{t,-\vecdelta_0}C(\var x,y_{-1})
=(\sqrt{-1}\lambda)
 \cdot\psi_{t,-\vecdelta_0}C(x_{-1},\var y),
\]
where $x_{-1}\in\psi_{t,-\vecdelta_0}\nbigm'$
and $y_{-1}\in\psi_{t,-\vecdelta_0}\nbigm''$
are sections such that
$x=-\sqrt{-1}\can x_{-1}$ and 
$y=\sqrt{-1}\can y_{-1}$.

Let $\nbigt=(\nbigm',\nbigm'',C)$ be an $\nbigr_X$-triple.
We say that it is strictly $S$-decomposable along $t$,
if $\nbigm'$ and $\nbigm''$ are strictly $S$-decomposable
along $t$.
If it is strictly $S$-decomposable,
the vanishing cycle $\phi_{t,0}\nbigt$
is defined to be
$(\psi_{t,0}\nbigm',\psi_{t,0}\nbigm'',\phi_{t,0}C)$.
By construction,
we have the following morphisms
of $\nbigr_{X_0}$-triples:
\[
\begin{array}{l}
 \Can:=(\var,\sqrt{-1}\can):
 \psi_{t,-\vecdelta_0}\nbigt
\lrarr \phi_{t,0}\nbigt(-1/2)\\
 \mbox{{}}\\
 \Var:=(-\sqrt{-1}\can,\var):
 \phi_{t,0}\nbigt(1/2)
\lrarr \psi_{t,-\vecdelta_0}\nbigt
\end{array}
\]

\subsection{The general case}

Let $Y$ be a complex manifold,
and let $\nbigt$ be an $\nbigr_{Y}$-triple.
Let $W$ be an open subset of $Y$,
and let $f$ be a holomorphic function on $W$.
Let $\iota:W\lrarr W\times\cnum_t$ denote 
the graph embedding.
We obtain the $\nbigr_{W\times\cnum_t}$-triple
$\iota_{\dagger}\bigl(
 \nbigt_{|\cnum_{\lambda}\times W}\bigr)$
on $\cnum_{\lambda}\times(W\times\cnum_t)$.
We say that
$\nbigt$ is strictly specializable ($S$-decomposable)
along $f$,
if $\iota_{\dagger}\bigl(
 \nbigt_{|\cnum_{\lambda}\times U}\bigr)$
is strictly specializable ($S$-decomposable)
along $t$.
If $\nbigt$ is strictly specializable,
 $\psitilde_{f,u}\nbigt$
 (resp. $\psi_{f,u}\nbigt$) are defined to be
 $\psitilde_{t,u}\bigl(
 \iota_{\dagger}\nbigt\bigr)$
 (resp. $\psi_{t,u}\iota_{\dagger}\nbigt$)
 for any $u\in\real\times\cnum$.
 If $\nbigt$ is strictly $S$-decomposable along $f$,
 $\phi_{f,0}(\nbigt)$ is defined to be
 $\phi_{t,0}\bigl(
 \iota_{\dagger}\nbigt\bigr)$.

We can show the following lemma
by using the argument due to Sabbah.
(See Theorem 3.3.15 
 and Corollary 3.6.35 in \cite{sabbah2}.
 See also \cite{saito1} for the original work 
 due to Saito.)
\begin{lem}
 \label{lem;05.1.25.150}
Let $F:X\lrarr Y$ be a proper morphism 
of complex manifolds.
Let $g$ be any holomorphic function on $Y$.
We put $\widetilde{g}:=g\circ F$.
Let $\nbigt$ be an $\nbigr_X$-triple
such that
(i) $\nbigt$ is strictly specializable along $\gtilde$,
(ii) $F^j_{\dagger}\psi_{\gtilde,u}(\nbigt)$
 are strict.
Then, $F^j_{\dagger}\nbigt$
are also strictly specializable along $g$.
Moreover, we have natural isomorphisms
\[
\psitilde_{g,u}F^j_{\dagger}\nbigt
\simeq
 F^j_{\dagger}\psitilde_{\widetilde{g},u}\nbigt,
\quad
 \psi_{g,u}F^j_{\dagger}\nbigt
\simeq
 F^j_{\dagger}\psi_{\widetilde{g},u}\nbigt
\]
for any $u\in\real\times\cnum$.
\hfill\qed
\end{lem}

\begin{lem}
\label{lem;07.10.27.5}
Let $X,Y,F,g$, and $\nbigt$ be 
as in Lemma {\rm \ref{lem;05.1.25.150}}.
Assume moreover that 
$\nbigt$ and $F^j_{\dagger}\nbigt$ $(j\in\seisuu)$
are strictly $S$-decomposable along $\gtilde$
and $g$, respectively.
Then, we have a natural isomorphism
$F^j_{\dagger}\phi_{\gtilde,0}\nbigt
\simeq
 \phi_{g,0}F^j_{\dagger}\nbigt$
for any $j$.
\hfill\qed
\end{lem}
The lemma was proved by Sabbah
without the additional assumption
of strict $S$-decomposability.
(Theorem 3.3.15
and Corollary 3.6.35 of \cite{sabbah2}.)
We will later give a direct argument 
in this restricted case
(Section \ref{subsection;08.10.18.150}).

\subsection{Uniqueness}

We recall a rather general remark
on the uniqueness of
the prolongation of Hermitian sesqui-linear pairings
(Proposition 14.97 of \cite{mochi2}).
Let $X$ be a complex manifold,
and let $f$ be a holomorphic function on $X$.
We put $X':=X-f^{-1}(0)$.
Let $\nbigm'$ and $\nbigm''$ be $\nbigr_X$-modules
such that 
(i) they are strictly specializable along $f$,
(ii) the morphism
$\can:\psi_{f,-\vecdelta_0}(\nbigm')\lrarr
 \psi_{f,0}(\nbigm')$ is surjective.
Let $C_a$ $(a=1,2)$ be sesqui-linear pairings
of $\nbigm'$ and $\nbigm''$.

\begin{prop}
 \label{prop;a11.23.10}
If $C_1=C_2$ on $\vecS\times X'$,
we have $C_1=C_2$ on $\vecS\times X$.
\end{prop}
\pf
Let $i_f:X\lrarr X\times\cnum$ denote the graph embedding.
We have only to show 
$i_{f\,\dagger}C_1=i_{f\dagger}C_2$.
Thus, we may assume that 
$X$ is of the form $X_0\times \cnum$,
and that
$f=t$ is the coordinate of $\cnum$ from the beginning.

For $\lambda_0\in\vecS$,
let $U$ be a small neighbourhood of $\lambda_0$
in $\cnum_{\lambda}$,
and $\vecI:=\vecS\cap U$.
Let $W_0$ be an open subset of $X_0$,
and $W$ be the product of $W_0$ 
and an open subset $W_1$ of $\cnum_t$.
Let $m'$ and $m''$ be sections of
$\nbigm'$ and $\nbigm''$
respectively on $U\times W$ 
and $\sigma(U)\times W$ respectively.
We put
$A(m',\overline{m''}):=
 C_1(m',\overline{m''})-C_2(m',\overline{m''})$.

As a preparation,
let us consider the case 
in which $\lambda_0$ is generic.
\begin{lem}\label{lem;9.29.1}
Assume that $\lambda_0$ is generic with respect to
the set 
$\KMS(\nbigm',t)
 \cup
 \KMS(\nbigm'',t)\cup\{(0,0)\}$.
Then we have 
$A(m',\overline{m''})=0$.
\end{lem}
\pf
We use an argument due to Sabbah 
(Proposition 3.7.6 in \cite{sabbah2}).
Since the support of $A(m',\overline{m''})$ is contained in 
$\vecI\times W_0$,
we have the following expression:
\[
 A(m',\overline{m''})=
 \sum_{a+b\leq p}
 \eta_{a,b}\cdot
 \del_t^a\cdot\delbar^b_t
 \cdot\delta_{\vecI\times W_0}
\]
Here $\eta_{a,b}$ are sections of
$\distribution_{\vecI\times W_0/\vecI}$,
and $\delta_{\vecI\times W_0}$
denote the delta distribution
for $\vecI\times W_0$ in $\vecI\times W$.
We have only to show $\eta_{a,b}=0$
for any $a$ and $b$.

Let us consider the case 
in which $m'\in \Vzero_{<0}(\nbigm')$.
For any $p\in\seisuu_{\geq 0}$,
we have a finite subset
$S(p)\subset \KMS(\nbigm',t)$ 
such that the following holds:
\begin{itemize}
\item
$\paramap(\lambda_0,u)<0$ 
for any element $u\in S(p)$.
\item
We put 
$B_p(x):=
 \prod_{u\in S(p)}
 \bigl(x+\eigenmap(\lambda,u)\bigr)$.
Then there is a section $P_p$ of $V_0\nbigr_X$
such that 
$B_p(-\deldel_tt)\cdot m'=P_p\cdot t^{p+1}\cdot m'$.
\end{itemize}
If $p$ is sufficiently large,
we have the vanishing
\[
B_p(-\deldel_tt)
 A(m',m'')
=P_p t^{p+1} A(m',m'')=0.
\]
Note the following equality:
\[
 \bigl(-\deldel_tt+\eigenmap(\lambda,u)\bigr)
 \cdot \deldel_t^a\cdot
 \delta_{W_0}
=\bigl(
 a\lambda+\eigenmap(\lambda,u)
 \bigr)\cdot
 \deldel_t^a\cdot
 \delta_{W_0}
=\eigenmap(\lambda,u-a\cdot\vecdelta_{0})
 \cdot\deldel_t^a\cdot\delta_{W_0}.
\]
Note
$\paramap(\lambda_0,u-a\cdot\vecdelta_0)<0$.
Since $\lambda_0$ is assumed to be generic,
we have 
$a\lambda_0+\eigenmap(\lambda_0,u)\neq 0$.
Thus we obtain $\eta_{a,b}=0$ 
in the case $m'\in \Vzero_{<0}(\nbigm')$.
Since $\nbigm'$ is generated by 
$\Vzero_{<0}\nbigm'$ around $\lambda_0$,
the general case can be reduced to the case $m'\in \Vzero_{<0}(\nbigm')$.
\hfill\qed

\vspace{.1in}

Let us return to the proof of Proposition \ref{prop;a11.23.10}.
Let $\phi$ be any test function on $W$.
By taking the push forward via 
the projection $p:\vecI\times W\lrarr \vecI$,
we obtain the distribution
$F:=p_{\ast}\bigl( \phi\cdot A(m',\overline{m''})\bigr)$,
which gives a continuous function on $\vecI$.
Due to Lemma \ref{lem;9.29.1},
$F$ vanishes on neighbourhoods of 
any generic $\lambda\in \vecI$.
Therefore we obtain the vanishing of $F$ on $\vecI$.
It means $C_1=C_2$ on $\vecI\times W$.
\hfill\qed

\subsection{The evaluation of
distributions and the value of holomorphic functions}
\label{subsection;08.9.2.11.4}

Let $X$ and $X_0$ be as in 
Subsection \ref{subsection;08.12.16.1}.
For any $M\in\real$,
we set  $\nbigu(M):=
 \bigl\{s\in\cnum\,\big|\,\Re(s)>M\bigr\}$.
Let $\nbigm'$ and $\nbigm''$ be
strictly specializable $\nbigr_X$-modules.
Let $C$ be a Hermitian sesqui-linear pairing
of $\nbigm'$ and $\nbigm''$.
Let $\lambda_0\in\vecS$,
and let $U(\lambda_0)$ be a small
neighbourhood of $\lambda_0$ in $\cnum$.
We put $\vecI(\lambda_0):=
 \vecS\cap U(\lambda_0)$.
Let $m\in V^{(\lambda_0)}_{<0}
 (\nbigm')\otimes C^{\infty}$
and $\mu\in V^{(-\lambda_0)}_{0}
 (\nbigm'')\otimes C^{\infty}$
be $C^{\infty}$-sections
on $U(\lambda_0)\times X$
and $\sigma(U(\lambda_0))\times X$,
respectively.
There exists a large $M$,
depending on $m$ and $\mu$
such that 
(i) the following pairing makes sense
for any $s\in \nbigu(M-k)$, $k\in\seisuu_{\geq 0}$
and 
a $C^{\infty}$ $(n,n)$-form $\varphi$
on $X$ with compact support
\[
 \nbigi^{(k)}_{C(m,\sigma^{\ast}\mu),\varphi}(s)
 :=\Bigl\langle
 C(m,\sigma^{\ast}\mu),\,|t|^{2s}\cdot t^k\cdot\varphi
 \Bigr\rangle,
\]
(ii) it gives a continuous function 
on $\vecI(\lambda_0)\times \nbigu(M-k)$,
which is holomorphic with respect to $s$.
We put
$\nbigi_{C(m,\sigma^{\ast}\mu),\varphi}
:=\nbigi^{(0)}_{C(m,\sigma^{\ast}\mu),\varphi}$.

\begin{lem}
\label{lem;07.10.26.25}
$\nbigi_{C(m,\sigma^{\ast}\mu),\varphi}(s)$ is 
naturally extended to a continuous function
on $\vecI(\lambda_0)\times \nbigu(-\delta)$
for some $\delta>0$,
which are holomorphic with respect to $s$.
And, we have 
$\nbigi_{C(m,\sigma^{\ast}\mu),\varphi}(0)
 =\Bigl\langle
 C\bigl(m,\sigma^{\ast}\mu\bigr),\,
 \varphi\Bigr\rangle$.
\end{lem}
\pf
We may and will assume that 
$m$ and $\mu$ are holomorphic sections.
Fix $R>0$.
Using an argument in \cite{sabbah2}
(or the proof of Lemma 14.76 
and Lemma 14.77 in \cite{mochi2}),
we can show that there exist
a large number $N>0$ and
a finite subset $S\subset\real\times\cnum$
with the following property:
\begin{itemize}
\item
 $\prod_{u\in S}
 \bigl(s+\lambda^{-1}\eigenmap(\lambda,u)
 \bigr)^{N}
 \cdot 
 \nbigi_{C(m,\sigma^{\ast}\mu),\varphi}(s)$
 is naturally extended to a continuous function
 on $\vecI(\lambda_0)\times\nbigu(-R)$,
 which is holomorphic with respect to $s$.
\item
 $\paramap(\lambda_0,u)<0$
 and $\paramap(-\lambda_0,u)\leq 0$
 for any $u\in S$.
\end{itemize}
Note the general formula
$\paramap(\lambda,u)+
 \paramap(\sigma(\lambda),u)
=-2\Re\bigl(\eigenmap(\lambda,u)/\lambda\bigr)$.
Then, the first claim follows.

Let us show the second claim.
By the above argument,
we have known that 
$\nbigi^{(k)}_{C(m,\sigma^{\ast}\mu),\varphi}(s)$
is a continuous function on
$\vecI(\lambda_0)\times\nbigu(-\delta)$,
which is holomorphic with respect to $s$.
We prepare the following lemma.
\begin{lem}
\label{lem;07.10.26.10}
There exists a $k_0\in\seisuu_{>0}$
such that the following holds for any 
$m\in V^{(\lambda_0)}_{<0}$,
$\mu\in V^{(-\lambda_0)}_0$,
any $C^{\infty}$ $(n,n)$-form $\varphi$
with compact support,
and any $k\geq k_0$:
\[
 \nbigi^{(k)}_{C(m,\sigma^{\ast}\mu),\varphi}(0)
=\bigl\langle
 C(m,\sigma^{\ast}\mu),\,\varphi t^k
 \bigr\rangle
\]
\end{lem}
\pf
Let $m_i$ and $\mu_j$ be finite generators
of $V^{(\lambda_0)}_{<0}\nbigm'$
and $V^{(-\lambda_0)}_0\nbigm''$
over $V_0\nbigr$.
We can take some $k_0$ such that
the following holds for any $m_i$, $\mu_j$,
$\varphi$,
and $k\geq k_0$:
\[
 \nbigi^{(k)}_{C(m_i,\sigma^{\ast}\mu_j),\varphi}(0)
=\bigl\langle
 C(m_i,\sigma^{\ast}\mu_j),\,\varphi t^k
 \bigr\rangle
\]
For any $m=\sum a_i\cdot m_i$
and $\mu=\sum b_j\cdot \mu_j$,
we have the following:
\[
 \bigl\langle
 C\bigl(m,\sigma^{\ast}\mu\bigr),
 \,\varphi |t|^{2s} t^k
 \bigr\rangle
=\sum
 \Bigl\langle
 C\bigl(m_i,\sigma^{\ast}\mu_j\bigr),\,
 a_i^{\ast}\overline{\sigma^{\ast}(b_j^{\ast})}
 \bigl(
 \varphi |t|^{2s} t^k\bigr)
 \Bigr\rangle
\]
Here, $a^{\ast}=\sum(-1)^{|J|}\deldel_J\cdot a_J $
for $a=\sum a_J\cdot \deldel^J$.
Note $a_i^{\ast}
 \overline{\sigma^{\ast}(b_j^{\ast})}
 (\varphi |t|^{2s} t^k)$
is of the form $|t|^{2s} t^k G_{i,j}$,
where $G_{i,j}$ are $C^{\infty}$ $(n,n)$-forms
with compact supports.
Then, the claim of Lemma \ref{lem;07.10.26.10}
follows.
\hfill\qed

\vspace{.1in}

By a descending induction on $k\geq 0$, 
let us show
the following equality for any 
$m\in \Vzero_{<0}$, $\mu\in V^{(-\lambda_0)}_0$
and any $C^{\infty}$ $(n,n)$-form $\varphi$
with compact supports:
\[
 \nbigi^{(k)}_{C(m,\sigma^{\ast}\mu),\varphi}(0)
=\bigl\langle
 C(m,\sigma^{\ast}\mu),\,\varphi t^k
 \bigr\rangle
\]
Assume that it holds for $k+1$.
Let $m':=\bigl(-\deldel_t t+\eigenmap(\lambda,u)\bigr)m$.
Then, we have the following:
\begin{multline}
\label{eq;07.10.26.20}
 \bigl\langle
 C(m',\sigma^{\ast}\mu),\,\,\varphi|t|^{2s}t^k
 \bigr\rangle
=\bigl\langle
 C(m,\sigma^{\ast}\mu),\,\,
 (t\deldel_t+
 \eigenmap(\lambda,u))(\varphi |t|^{2s} t^k)
 \bigr\rangle\\
=\bigl(\lambda(s+k)+\eigenmap(\lambda,u)\bigr)\cdot
 \bigl\langle
 C(m,\sigma^{\ast}\mu),\,
 \varphi |t|^{2s} t^k
 \bigr\rangle
+\bigl\langle
 C(m,\sigma^{\ast}\mu),\,
 \deldel_t\varphi \cdot |t|^{2s} t^{k+1}
 \bigr\rangle
\end{multline}
By a similar calculation,
we obtain the following:
\begin{equation}
\label{eq;07.10.26.21}
\bigl\langle
 C(m',\sigma^{\ast}\mu),\varphi t^k
\bigr\rangle
=\bigl(\lambda k+\eigenmap(\lambda,u)\bigr)
 \bigl\langle
 C(m,\sigma^{\ast}\mu),\varphi t^k
 \bigr\rangle
+\bigl\langle
 C(m,\sigma^{\ast}\mu),\deldel_t\varphi\cdot t^{k+1}
 \bigr\rangle
\end{equation}

We take a polynomial $b_m$
such that
$b_m(-\deldel_tt)m
=t\cdot P(-\deldel_t t)m
\in V_{-1}\nbigr\cdot m$.
It is of the following form:
\[
 b_m(x)=
 \prod_{i=1}^M\bigl(x+\eigenmap(\lambda,u_i)\bigr)
\]
Here we have $\paramap(\lambda_0,u_i)<0$
for each $u_i$.
We put $e_i:=\eigenmap(\lambda,u_i)$,
for simplicity of description.
By using (\ref{eq;07.10.26.20}) inductively,
we obtain the following formula:
\begin{multline}
\label{eq;07.10.26.22}
 \Bigl\langle
 C\bigl(P(-\deldel_tt)m,\sigma^{\ast}\mu\bigr),\,
 \varphi |t|^{2s} t^{k+1}
 \Bigr\rangle
=\prod_{i=1}^M\bigl(\lambda(s+k)+e_i\bigr)
\cdot
 \Bigl\langle
 C(m,\sigma^{\ast}\mu),\,\varphi |t|^{2s} t^k
 \Bigr\rangle
 \\
+\sum_{j=0}^{M-1}
 \prod_{i=j+1}^{M}
 \bigl(\lambda(s+k)+e_i\bigr)
\cdot\Bigl\langle
 C\Bigl(\prod_{i=1}^j(-\deldel_tt+e_i)m,
 \sigma^{\ast}\mu\Bigr),\,\,
 \deldel_t\varphi\cdot |t|^{2s} t^{k+1}
 \Bigr\rangle
\end{multline}
We also have the following formula
from (\ref{eq;07.10.26.21}):
\begin{multline}
\label{eq;07.10.26.23}
 \Bigl\langle
 C\bigl(P(-\deldel_tt)m,\sigma^{\ast}\mu\bigr),\,
 \varphi\cdot t^{k+1}
 \Bigr\rangle
=\Bigl\langle
 C(m,\sigma^{\ast}\mu),\varphi t^k
 \Bigr\rangle
\cdot\prod_{i=1}^M\bigl(\lambda k+e_i\bigr)
 \\
+\sum_{j=0}^{M-1}
 \prod_{i=j+1}^{M}
 \bigl(\lambda k+e_i\bigr)
\cdot\Bigl\langle
 C\Bigl(
 \prod_{i=1}^j(-\deldel_tt+e_i)m,
 \sigma^{\ast}\mu\Bigr),\,
 \deldel_t\varphi\cdot t^{k+1}
 \Bigr\rangle
\end{multline}
Then, we obtain
$\nbigi^{(k)}_{C(m,\sigma^{\ast}\mu),\varphi}(0)
=\bigl\langle
 C(m,\sigma^{\ast}\mu),\,\varphi t^k
 \bigr\rangle$
from (\ref{eq;07.10.26.22}) and
(\ref{eq;07.10.26.23}).
Thus, the induction can proceed,
and the proof of Lemma \ref{lem;07.10.26.25}
is finished.
\hfill\qed

\subsection{Proof of Lemma \ref{lem;07.10.27.5}}
\label{subsection;08.10.18.150}

Under the isomorphism
$\psi_{t,0}F^j_{\dagger}\nbigm_i\simeq
 F^j_{\dagger}\psi_{t,0}(\nbigm_i)$
for any $j$,
we would like to show the equality
$\phi_{t,0}(F^j_{\dagger}C)
=F^j_{\dagger}\phi_{t,0}(C)$.
We have the morphisms
$\can:
 \psi_{t,-\vecdelta_0}(\nbigm')
\lrarr \psi_{t,0}(\nbigm')$
and $\var:
 \psi_{t,0}(\nbigm')\lrarr\psi_{t,-\vecdelta_0}(\nbigm')$.
We also have similar morphisms for $\nbigm''$,
which are denoted by the same symbols.
They induce the morphisms
$\can:
 F^{-i}_{0\dagger}\psi_{t,-\vecdelta_0}(\nbigm')
\lrarr
 F^{-i}_{0\dagger}\psi_{t,0}(\nbigm')$
and 
$\var:
 F^{-i}_{0\dagger}\psi_{t,0}(\nbigm')\lrarr
 F^{-i}_{0\dagger}\psi_{t,-\vecdelta_0}(\nbigm')$.
We have similar morphisms for $\nbigm''$.
Recall the following relations:
\begin{equation}
\label{eq;07.12.1.1}
 \phi_{t,0}\bigl(F_{\dagger}^iC\bigr)
 \bigl(-\sqrt{-1}\can (a),\,\sigma^{\ast}b\bigr)
=(\sqrt{-1}\lambda)
 \psi_{t,-\vecdelta_0}\bigl(F_{\dagger}^iC\bigr)
 \bigl(a,\,\sigma^{\ast}\var (b)\bigr),
\end{equation}
\begin{equation}
\label{eq;07.12.1.2}
 \phi_{t,0}\bigl(F_{\dagger}^iC\bigr)
 \bigl(a,\,\sigma^{\ast}\sqrt{-1}\can(b)\bigr) \\
=(\sqrt{-1}\lambda)
 \psi_{t,-\vecdelta_0}\bigl(
 F_{\dagger}^iC\bigr)\bigl(\var(a),\,\sigma^{\ast}b\bigr)
\end{equation}
The pairing $F_{\dagger}^i\phi_{t,0}(C)$
satisfy similar relation with 
$\psi_{t,-\vecdelta_0}F_{\dagger}^iC
=F_{\dagger}^i\psi_{t,-\vecdelta_0}C$.

By the assumption,
we have the decomposition
$F^{-i}_{\dagger}(\nbigm')=M_1'\oplus M_2'$
such that 
(i) $M_1'$ has no non-zero $\nbigr_Y$-submodule
 whose support is contained in $Y_0$,
(ii) the support of $M_2'$ is contained in $Y_0$.
We have a similar decomposition 
$F^i_{\dagger}\nbigm''=M_1''\oplus M_2''$.
By using Proposition \ref{prop;08.1.15.30}
and the relations between the pairings
$\phi_{t,0}(F_{\dagger}^iC)$,
$F^i_{\dagger}\phi_{t,0}(C)$
and $\psi_{t,-\vecdelta_0}(F_{\dagger}^iC)$
like (\ref{eq;07.12.1.1}) and (\ref{eq;07.12.1.2}),
we can show the following:
\begin{itemize}
\item
 $\psi_{t,0}(M_1')$ and $\psi_{t,0}(M_2'')$ 
 are orthogonal with respect to
 both of the pairings
 $\phi_{t,0}(F_{\dagger}^iC)$ and
 $F^i_{\dagger}\phi_{t,0}(C)$.
 Similarly,
 $\psi_{t,0}(M_1'')$ and $\psi_{t,0}(M_2')$ 
 are orthogonal with respect to
 the both pairings.
\item
$\phi_{t,0}(F_{\dagger}^iC)$ and
 $F^i_{\dagger}\phi_{t,0}(C)$
are the same on 
 $\psi_{t,0}(M_1')_{|\vecS\times Y}\otimes
 \sigma^{\ast}\psi_{t,0}(M_1'')_{|\vecS\times Y}$.
\end{itemize}
Hence, we have only to compare 
the induced sesqui-linear pairings of
$\psi_{t,0}(M_2')$ and 
$\sigma^{\ast}\psi_{t,0}(M_2'')$.

\vspace{.1in}
We may assume that $F_0$ is the projection,
i.e., $X_0=Y_0\times Z$.
Let $n=\dim Z$.
Then, $F_{\dagger}\nbigm'$
can be expressed by 
$F_{\ast}(\nbigm'\otimes \Omega^{n+\bullet}_{\nbigz})$.
(See Subsection \ref{subsection;08.12.16.3}
for $\Omega_{\nbigz}^{\bullet}$.)
Similarly,
we have
$F_{\dagger}\nbigm''
=F_{\ast}(\nbigm''\otimes
 \Omega^{n+\bullet}_{\nbigz})$.
Let $\lambda_0\in\vecS$.
Let $u$ and $v$ be local sections of $\Vzero_0M_2'$
and $V^{(-\lambda_0)}_0M_2''$ respectively.
We can regard them as sections of
$F^{-i}_{\dagger}\nbigm'$
and $F^i_{\dagger}(\nbigm'')$
respectively,
and we have 
$\phi_{t,0}F^i_{\dagger}(C)(u,v)=
 F^i_{\dagger}C(u,v)$.
We take lifts $\utilde$ and $\vtilde$
of $u$ and $v$
to $F_{\ast}(\nbigm'\otimes
 \Omega^{n-i}_{\nbigz})$
and $F_{\ast}(\nbigm''\otimes
 \Omega^{n+i}_{\nbigz})$,
respectively.
Let $\varphi$ be a test function on $Y_0$.
Let $\rho$ be a test function on $\cnum_t$
which is constantly $1$ around $t=0$,
and let $\omega:=\rho\cdot \sqrt{-1}dt\cdot d\tbar/2\pi$.
Then, we have the following equalities:
\begin{equation}
\label{eq;07.12.1.6}
 \bigl\langle
 \phi_{t,0}(F^i_{\dagger}C)(u,\sigma^{\ast}v),\,
 \varphi
 \bigr\rangle
=\Bigl\langle
 F^i_{\dagger}C(u,\sigma^{\ast}v),\,
 \varphi\cdot \omega
 \Bigr\rangle
=\Bigl\langle
 C(\utilde,\sigma^{\ast}\vtilde),\,
 F_0^{\ast}\varphi\cdot\omega
 \Bigr\rangle
\end{equation}

Let $[\utilde]$ and $[\vtilde]$
denote the lift of $u$ and $v$
to $F_{0\ast}\bigl(
 \nbigm'\otimes\Omega^{n-i}_{\nbigz}\bigr)$
and $F_{0\ast}\bigl(
 \nbigm''\otimes\Omega^{n+i}_{\nbigz}\bigr)$
respectively,
which are induced by $\utilde$ and $\vtilde$.
We have the following equality:
\begin{equation}
\label{eq;07.12.1.3}
\bigl\langle
 F^i_{0\dagger}\phi_{t,0}(C)(u,\sigma^{\ast}v),\,
 \varphi
\bigr\rangle
=\bigl\langle
 \phi_{t,0}(C)([\utilde],\sigma^{\ast}[\vtilde]),\,
 F_0^{\ast}\varphi
 \bigr\rangle
\end{equation}
We take an appropriate open locally finite covering
$Z=\bigcup_{j\in S} U_j$
and the partition of unity $\{\chi_j\,|\,j\in S\}$
which is subordinate to the covering.
We may assume to have
$u_j^{(1)}\in 
V^{(\lambda_0)}_{-1}
 (\nbigm')\otimes\Omega_{\nbigz}^{n-i}$
and 
$u_j^{(2)}\in V^{(\lambda_0)}_{<0}
 (\nbigm')\otimes\Omega_{\nbigz}^{n-i}$
on $U_j$ 
such that 
$\utilde_{|U_j}=\deldel_tu_j^{(1)}+u_j^{(2)}$.
Then, we have the following equality:
\begin{multline}
\label{eq;07.12.1.4}
\Bigl\langle
 \phi_{t,0}(C)\bigl([\utilde],
 \sigma^{\ast}[\vtilde]\bigr),\,
 F_0^{\ast}\varphi
\Bigr\rangle
=\sum_{j\in S}
 \Bigl\langle
 \phi_{t,0}(C)\bigl([\utilde],
 \sigma^{\ast}[\vtilde]\bigr),\,
 \chi_j F_0^{\ast}\varphi
 \Bigr\rangle  \\
=\sum_{j\in S}
 \Bigl\langle
 \phi_{t,0}(C)\bigl([\deldel_tu_j^{(1)}],
 \sigma^{\ast}[\vtilde]\bigr),\,
 \chi_j F_0^{\ast}\varphi
 \Bigr\rangle 
=\sum_{j\in S}
 \lambda\Bigl\langle
 \psi_{t,-\vecdelta_0}(C)\bigl([u_j^{(1)}],\,
 \sigma^{\ast}[t\vtilde]\bigr),\,\,
 \chi_j F_0^{\ast}\varphi
 \Bigr\rangle  \\
=\sum_{j\in S}
 \underset{s=-1}\Res\,
 \lambda\Bigl\langle
 C\bigl(u_j^{(1)},\sigma^{\ast}\vtilde\bigr),
 \,\,\chi_j F_0^{\ast}\varphi
 \cdot
 |t|^{2s}\tbar
 \!\cdot\!\omega
 \Bigr\rangle\\
=
 \underset{s=-1}\Res\,
 \frac{1}{s+1}
 \Bigl\langle
 C\Bigl(\sum_{j\in S}\chi_j\!\cdot\! u_j^{(1)},\,
 \sigma^{\ast}\vtilde\Bigr),\,\,
 F_0^{\ast}\varphi
 \cdot\deldel_t(|t|^{2(s+1)})
 \cdot\omega
 \Bigr\rangle
\end{multline}
It can be rewritten as follows:
\begin{equation}
\begin{split}
-\underset{s=-1}\Res
\frac{1}{s+1}
 \Bigl\langle
 C\Bigl(\deldel_t 
 \Bigl(\sum_{j\in S}
 \chi_j\cdot u_j^{(1)}\Bigr),\,\sigma^{\ast}\vtilde
\Bigr),\,\,
 F_0^{\ast}\varphi\cdot |t|^{2(s+1)}
 \cdot\omega
 \Bigr\rangle \\
-\underset{s=-1}\Res
 \frac{1}{s+1}
 \Bigl\langle
 C\Bigl(\sum_{j\in S}\chi_j\!\cdot\! u_j^{(1)},\,
 \sigma^{\ast}(\deldelbar_t\vtilde)\Bigr),\,\,
 F_0^{\ast}\varphi\cdot
 |t|^{2(s+1)}\cdot \omega
 \Bigr\rangle\\
-\underset{s=-1}\Res
 \frac{1}{s+1}
 \Bigl\langle
 C\Bigl(\sum_{j\in S}
 \chi_j\!\cdot\! u_j^{(1)},\,\sigma^{\ast}\vtilde\Bigr),\,\,
 F_0^{\ast}\varphi\cdot
 |t|^{2(s+1)}\cdot \deldel_t\omega
 \Bigr\rangle\\
\end{split}
\end{equation}
It is equal to the following:
\begin{equation}
 \label{eq;07.10.26.6}
\begin{split}
 \underset{s=-1}\Res\frac{1}{s+1}
 \Bigl\langle
  C(\utilde,\sigma^{\ast}\vtilde),\,\,
 F_0^{\ast}\varphi\cdot
 |t|^{2(s+1)}\cdot\omega
 \Bigr\rangle \\
-\underset{s=-1}\Res\frac{1}{s+1}
 \Bigl\langle
 C\Bigl(\sum_{j\in S} \chi_j\!\cdot\! u^{(2)}_j,\, 
 \sigma^{\ast}(\vtilde)\Bigr),\,\,
 F_0^{\ast}\varphi\cdot |t|^{2s+2}\cdot\omega
 \Bigr\rangle\\
-\underset{s=-1}\Res
 \frac{1}{s+1}
 \Bigl\langle
 C\Bigl(\sum_{j\in S}\chi_j\!\cdot\! u_j^{(1)},
  \,\sigma^{\ast}(\deldelbar_t\vtilde)\Bigr),\,\,
 F_0^{\ast}\varphi\cdot
 |t|^{2(s+1)}\cdot \omega
 \Bigr\rangle\\
-\underset{s=-1}\Res
 \frac{1}{s+1}
 \Bigl\langle
 C\Bigl(\sum_{j\in S}\chi_j\!\cdot\! u_j^{(1)},\,
 \sigma^{\ast}(\vtilde)\Bigr),\,\,\,
 F_0^{\ast}\varphi\cdot
  |t|^{2(s+1)}\cdot
 \deldel_t\omega
 \Bigr\rangle
\end{split}
\end{equation}
By the assumption,
the support of 
$F_{0\ast}C(\utilde,\sigma^{\ast}\vtilde)$
is contained in $t=0$.
Hence, the first term 
of (\ref{eq;07.10.26.6}) is $0$.
Let us look at the second term.
We have 
$u^{(2)}_j\in V^{(\lambda_0)}_{<0}$
and $\vtilde\in V^{(-\lambda_0)}_{0}$.
Due to Lemma \ref{lem;07.10.26.25},
\[
 G(s):=\Bigl\langle
 C\bigl(\sum \chi_j\!\cdot\! u^{(2)}_j,\,
 \sigma^{\ast}(\vtilde)\bigr),\,
 F_0^{\ast}\varphi\!\cdot\!
 |t|^{2s+2}\!\cdot\!\omega
 \Bigr\rangle
\]
is holomorphic on the region
$\bigl\{\Re(s)>-1-\delta\bigr\}$ for some $\delta>0$,
and the following holds:
\begin{equation}
 \label{eq;07.10.26.30}
 G(-1)=
 \Bigl\langle
 C\Bigl(\sum_{j\in S} \chi_j\cdot u^{(2)}_j,\, 
 \sigma^{\ast}\vtilde\Bigr),\,\,
 F_0^{\ast}\varphi
  \cdot\omega
 \Bigr\rangle
\end{equation}
Hence,
the second term in (\ref{eq;07.10.26.6}) 
is given by (\ref{eq;07.10.26.30}).
For the same reason,
the third term is as follows:
\begin{equation}
-\Bigl\langle
 C\Bigl(\sum_{j\in S}\chi_j\!\cdot\! u_j^{(1)},\,
 \sigma^{\ast}(\deldelbar_t\vtilde)\Bigr),\,\,
 F_0^{\ast}\varphi\cdot\omega
 \Bigr\rangle
\end{equation}
We have a similar equality for the fourth term.
Hence, the sum of
the third and the forth terms in (\ref{eq;07.10.26.6})
can be rewritten as follows:
\begin{multline}
-\Bigl\langle
 C\Bigl(\sum\chi_j\!\cdot\! u_j^{(1)},\,
 \sigma^{\ast}(\deldelbar_t\vtilde)\Bigr),\,
 F_0^{\ast}\varphi\cdot\omega
 \Bigr\rangle 
 -\Bigl\langle
 C\Bigl(\sum\chi_j\!\cdot\! u_j^{(1)},\,\sigma^{\ast}\vtilde\Bigr),\,
 F_0^{\ast}\varphi\cdot\deldel_t\omega
 \Bigr\rangle  \\
=
\Bigl\langle
 C\Bigl(\sum\chi_j\cdot
 \deldel_tu_j^{(1)},\,\sigma^{\ast}\vtilde\Bigr),\,\,
 F_0^{\ast}\varphi\cdot\omega
 \Bigr\rangle
\end{multline}
Hence (\ref{eq;07.10.26.6})
can be written as follows:
\begin{multline}
 \label{eq;07.12.1.5}
  \Bigl\langle
 C\Bigl(\sum \chi_j\cdot u^{(2)}_j,\, 
 \sigma^{\ast}\vtilde\Bigr),\,\,
 F_0^{\ast}\varphi\cdot\omega
 \Bigr\rangle
+\Bigl\langle
 C\Bigl(\sum\chi_j\cdot \deldel_tu_j^{(1)},\,
 \sigma^{\ast}\vtilde\Bigr),\,\,
 F_0^{\ast}\varphi\cdot\omega
 \Bigr\rangle \\
=\bigl\langle
 C\bigl(\utilde,\sigma^{\ast}\vtilde\bigr),\,
 F_0^{\ast}\varphi\cdot\omega
 \bigr\rangle
\end{multline}
From (\ref{eq;07.12.1.3}), (\ref{eq;07.12.1.4})
(\ref{eq;07.10.26.6}) and (\ref{eq;07.12.1.5}),
we obtain 
\[
\bigl\langle
 F^i_{0\dagger}\phi_{t,0}(C)(u,\sigma^{\ast}v),\,
 \varphi
\bigr\rangle
=\bigl\langle
 C\bigl(\utilde,\sigma^{\ast}\vtilde\bigr),\,
 F_0^{\ast}\varphi\cdot\omega
 \bigr\rangle. 
\]
Together with (\ref{eq;07.12.1.6}),
we obtain 
$\phi_{t,0}(F^i_{\dagger}C)(u,\sigma^{\ast}v)
=F^i_{0\dagger}\phi_{t,0}(C)(u,\sigma^{\ast}v)$.
Thus, Lemma \ref{lem;07.10.27.5} is proved.
\hfill\qed

\section{$\nbigr(\ast t)$-triples}
\label{subsection;08.10.18.113}
\label{subsection;08.9.14.120}
We recall the notion of $\nbigr_X(\ast t)$-triple
due to Sabbah in \cite{sabbah2},
where it is called $\nbigrtilde_X$-triples.
We use the notation in 
Section \ref{subsection;08.12.16.1}.
Let $\nbigm'$  and $\nbigm''$ be
$\nbigr_X(\ast t)$-modules.
An $\nbigr_{X}(\ast t)\otimes
 \sigma^{\ast}\nbigr_{X}(\ast t)$-homomorphism
$C:\nbigm'_{|\vecS\times X}\otimes
 \sigma^{\ast}\nbigm''_{|\vecS\times X}
\lrarr
 \distribution^{\moderate}_{\vecS\times X/\vecS}$
is called a Hermitian sesqui-linear pairing of 
$\nbigm'$ and $\nbigm''$.
It is also called a sesqui-linear pairing.
Such $\nbigt=(\nbigm',\nbigm'',C)$ is called
an $\nbigr_X(\ast t)$-triple.
\index{Hermitian sesqui-linear pairing}
\index{sesqui-linear pairing}
\index{$\nbigr_X(\ast)$-triple}

Let $\nbigt_i=(\nbigm_i',\nbigm_i'',C_i)$ $(i=1,2)$
be $\nbigr_X(\ast t)$-triples.
A morphism of $\nbigr_X(\ast t)$-triples
$\varphi:\nbigt_1\lrarr\nbigt_2$
is defined to be a pair 
of $\varphi':\nbigm_2'\lrarr \nbigm_1'$
and $\varphi'':\nbigm_1''\lrarr\nbigm_2''$
satisfying 
$C_1\bigl(
 \varphi'(x),\sigma^{\ast}(y)
 \bigr)
=C_2\bigl(x,\sigma^{\ast}(\varphi''(y))\bigr)$.
The category of 
$\nbigr_X(\ast t)$-triples is abelian.

For an $\nbigr_X(\ast t)$-triple
 $\nbigt=(\nbigm',\nbigm'',C)$,
we have the induced sesqui-linear pairing
\[
 C^{\ast}:\nbigm''_{|\vecS\times X}
 \otimes\sigma^{\ast}\nbigm'_{|\vecS\times X}
\lrarr
 \distribution^{\moderate}_{\vecS\times X/\vecS}
\]
given by
$C^{\ast}(x,\sigma^{\ast}y)
=\overline{\sigma^{\ast}C(y,\sigma^{\ast}x)}$.
The tuple
$\nbigt^{\ast}=(\nbigm'',\nbigm',C^{\ast})$
is called the Hermitian adjoint of $\nbigt$.
A morphism 
$\varphi=(\varphi',\varphi''):
 \nbigt_1\lrarr\nbigt_2$ naturally induces
$\varphi^{\ast}=(\varphi'',\varphi'):
 \nbigt_2^{\ast}\lrarr\nbigt_1^{\ast}$.

Let $k$ be any half integer.
For an $\nbigr_X(\ast t)$-triple 
$\nbigt=(\nbigm',\nbigm'',C)$,
we put $\nbigt(k):=
 \bigl(\nbigm',\nbigm'',
 (\sqrt{-1}\lambda)^{-2k}C\bigr)$,
which is called the $k$-th Tate twist
of $\nbigt$.
It can be regarded as $\nbigt\otimes\Tate^S(k)$.
We have the isomorphism
$\nbigt(k)\simeq \bigl(\nbigt^{\ast}(-k)\bigr)^{\ast}$
as in the case of $\nbigr_X$-triples.

For an $\nbigr_X(\ast t)$-triple $\nbigt$,
a morphism $\nbigs:\nbigt\lrarr\nbigt^{\ast}(-w)$
is called a Hermitian sesqui-linear duality of weight $w$,
if the equality $\nbigs=(-1)^w\nbigs^{\ast}$
holds.

\subsection{Specialization}

Let $\nbigt=(\nbigm',\nbigm'',C)$ be
an $\nbigr_X(\ast t)$-triple.
We say that
$\nbigt$ is strictly specializable,
if $\nbigm'$ and $\nbigm''$ are strictly specializable
along $t=0$.
In that case,
we have the induced pairing
for any $u\in\real\times\cnum$
(\cite{sabbah2})
as in the case of $\nbigr_X$-triples:
\[
  \psitilde_{t,u}C:
 \psitilde_{t,u}(\nbigm')_{|\vecS\times X_0}
\otimes
 \sigma^{\ast}\psitilde_{t,u}(\nbigm'')_{|\vecS\times X_0}
\lrarr
 \distribution_{\vecS\times X_{0}/\vecS}
\]
Thus, we obtain the $\nbigr_{X_0}$-triple
$\psitilde_{t,u}(\nbigt)=
\bigl(\psitilde_{t,u}(\nbigm'),\psitilde_{t,u}(\nbigm''),
 \psitilde_{t,u}(C) \bigr)$.

\begin{lem}
\label{lem;08.1.16.10}
Let $\nbigt_i$ $(i=1,2)$ be 
strictly specializable $\nbigr_X(\ast t)$-triples.
Let $\varphi:\nbigt_1\lrarr\nbigt_2$
be a morphism which is strictly specializable
along $t$.
Then, $\Ker(\varphi)$,
$\Image(\varphi)$ and $\Cok(\varphi)$
are also strictly specializable along $t$,
and we have the natural isomorphisms
$\psitilde_{t,u}\Ker(\varphi)
\simeq
 \Ker\psitilde_{t,u}(\varphi)$,
$\psitilde_{t,u}\Image(\varphi)
\simeq
 \Image\psitilde_{t,u}(\varphi)$
and 
$\psitilde_{t,u}\Cok(\varphi)
\simeq
 \Cok\psitilde_{t,u}(\varphi)$.
\end{lem}
\pf
Because Lemma \ref{lem;08.1.16.11},
we have only to compare 
the induced sesqui-linear pairings,
which can be checked easily.
\hfill\qed

\subsection{Pull back via unramified covering}
\index{pull back via unramified covering}

Let $\varphi_n$ be as in 
Section \ref{subsection;08.10.4.1}.
Let $\nbigt=(\nbigm',\nbigm'',C)$
be an $\nbigr_X(\ast t)$-triple.
We have the induced pairing
as explained in \cite{sabbah5}
(see also Section \ref{subsection;08.1.15.20}):
\[
 \varphi_n^{\dagger}C:
 \varphi_n^{\dagger}\nbigm'_{|\vecS\times X^{(n)}}
\otimes
 \sigma^{\ast}
 \varphi_n^{\dagger}\nbigm''_{|\vecS\times X^{(n)}}
\lrarr
 \distribution^{\moderate}_{\vecS\times X^{(n)}/\vecS}
\]
And thus, we obtain the $\nbigr_{X^{(n)}}(\ast t_n)$-triple
$\varphi_n^{\dagger}(\nbigt)
=\bigl(\varphi_n^{\dagger}\nbigm',
 \varphi_n^{\dagger}\nbigm'',
 \varphi_n^{\dagger}C\bigr)$.

Let $\vecI(\lambda_0)$ denote a small neighbourhood
of $\lambda_0$ in $\vecS$,
and let $\vecI(-\lambda_0)=\sigma(\vecI(\lambda_0))$.
\begin{lem}
Assume that $\nbigt$ is strictly specializable along $t$.
Let 
$[\tau_1]\in
 \psitilde_{t,(u+m_1\vecdelta)/n}
 (\nbigm')_{|\vecI(\lambda_0)\times X_0}$
and 
 $[\tau_2]\in
 \psitilde_{t,(u+m_2\vecdelta)/n}
 (\nbigm'')_{|\vecI(-\lambda_0)\times X_0}$.
We have the following relation:
\begin{multline*}
 \psitilde_{t_n,u}\bigl(\varphi_n^{\dagger}C\bigr)
 \Bigl(
 [t_n^{m_1}\varphi_n^{\dagger}\tau_1], 
 \sigma^{\ast}\bigl(
 [ t_n^{m_2}\varphi_n^{\dagger}\tau_2]
 \bigr)
 \Bigr)
=\\
\left\{
 \begin{array}{ll}
 \psitilde_{t,(u+m\vecdelta_0)/n}(C)
 \bigl([\tau_1],\sigma^{\ast}[\tau_2]\bigr)
 & m_1=m_2=m\\
 0 & m_1\neq m_2
 \end{array}
 \right.
\end{multline*}
As a result,
we have the following decomposition
of $\nbigr_{X_0}$-triples:
\[
 \psitilde_{t_n,u}\bigl(
 \varphi_n^{\dagger}\nbigt\bigr)
\simeq
 \bigoplus_{\substack{n u'-c\vecdelta_0=u\\ 0\leq c\leq n-1}}
 \psitilde_{t,u'}(\nbigt)
\]
\end{lem}
\pf
Let $\rho$ be a $C^{\infty}$-top form on $X_0$
with compact support.
Let $\chi_1=\chi_1(|t_n|)$ be 
a test function on $\cnum_{t_n}$
which are constantly $1$ around $0$,
and let 
$\chi:=\chi_1\cdot\sqrt{-1}dt\cdot d\tbar/2\pi$.
We have the following:
\begin{multline}
\label{eq;07.10.26.2}
\Bigl\langle
 \psitilde_{t_n,u}(\varphi_n^{\dagger}C)
 \bigl([t_n^{m_1}\varphi_n^{\dagger}\tau_1],\,
 \sigma^{\ast}[t_n^{m_2}\varphi_n^{\dagger}\tau_2]\bigr),\,\,
 \rho
\Bigr\rangle=\\
\underset{s+\eigenmap(\lambda,u)}
\Res\,
 \Bigl\langle
 \varphi_n^{\dagger}C\bigl(t_n^{m_1}\varphi_n^{\dagger}\tau_1,\,
 \sigma^{\ast}(t_n^{m_2}\varphi_n^{\dagger}\tau_2)\bigr),\,\,
 |t_n|^{2s}\chi\cdot \rho
 \Bigr\rangle
=\\
 \underset{s+\eigenmap(\lambda,u)}\Res\,
 \Bigl\langle
 \varphi_n^{\dagger}
 C\bigl(\tau_1,
 \sigma^{\ast}\tau_2\bigr),\,
 |t_n|^{2s}t_n^{m_1}\tbar_n^{m_2}\cdot\chi\cdot\rho
 \Bigr\rangle
\end{multline}
If $m_1\neq m_2$,
(\ref{eq;07.10.26.2}) is $0$
by Lemma \ref{lem;07.10.26.1}.
In the case $m_1=m_2=m$,
(\ref{eq;07.10.26.2}) is equal to the following:
\begin{multline}
\underset{s+\eigenmap(\lambda,u)}\Res\,
 \Bigl\langle
 C(\tau_1,\sigma^{\ast}\tau_2),\,
 |t|^{2(s+m)/n}\chi\cdot\rho
 \Bigr\rangle\cdot n
=\underset{n\sbar-m+\eigenmap(\lambda,u)}\Res\,
 \Bigl\langle
 C(\tau_1,\sigma^{\ast}\tau_2),\,
 |t|^{2\sbar}\chi\cdot\rho
 \Bigr\rangle\cdot n \\
=\underset{\sbar+\eigenmap(\lambda,(u+m\vecdelta_0)/n)}
 \Res\,
 \Bigl\langle
 C(\tau_1,\sigma^{\ast}\tau_2),\,|t|^{2\sbar}\chi\cdot\rho
 \Bigr\rangle
=\Bigl\langle
 \psitilde_{t,(u+m\vecdelta_0)/n}C
 \bigl([\tau_1],\sigma^{\ast}[\tau_2]\bigr),\,
 \rho\Bigr\rangle
\end{multline}
Thus we are done.
\hfill\qed

\subsection{Exponential twist}
\index{exponential twist}

Let $\nbigl(\gminia)$ be as in Section
\ref{subsection;08.9.4.12}.
For any $\gminia\in\cnum[t_n^{-1}]$,
we have the naturally defined sesqui-linear pairing
\[
 C_{\gminia}:\nbigl(\gminia)_{|\vecS\times X^{(n)}}\otimes
 \sigma^{\ast}\nbigl(\gminia)_{|\vecS\times X^{(n)}}
\lrarr \distribution_{\vecS\times X^{(n)}/\vecS},
\]
given by
$C_{\gminia}(e,\sigma^{\ast}e)
=\exp\bigl(-2\sqrt{-1}\Image(\lambda\gminiabar)\bigr)$.
We put $\gbigl(\gminia):=
 \bigl(\nbigl(\gminia),\nbigl(\gminia),C_{\gminia}\bigr)$.

Let $\nbigt$ be an $\nbigr_{X}(\ast t)$-triple.
Sabbah observed that
$\varphi_n^{\dagger}(\nbigt)\otimes \gbigl(-\gminia)$
is well defined
as the $\nbigr_{X^{(n)}}(\ast t_n)$-triple
(see Lemma \ref{lem;07.10.26.3}):
\[
 \varphi_n^{\dagger}(\nbigt)\otimes \gbigl(-\gminia)
:=
 \Bigl(\varphi_n^{\dagger}\nbigm'\otimes\nbigl(-\gminia),\,\,
 \varphi_n^{\dagger}\nbigm''\otimes\nbigl(-\gminia),\,\,
\varphi_n^{\dagger}C\cdot
 C_{-\gminia}
 \Bigr)
\]
It is naturally identified with
$ \Bigl(\varphi_n^{\dagger}\nbigm',\,\,
 \varphi_n^{\dagger}\nbigm'',\,\,
 \exp\bigl(2\sqrt{-1}\Image(\lambda\gminiabar)\bigr)
 \varphi_n^{\dagger}C
 \Bigr)$.
If the underlying $\nbigr_{X^{(n)}}(\ast t_n)$-modules of
$\varphi_n^{\dagger}(\nbigt)\otimes \gbigl(-\gminia)$
are strictly specializable,
we obtain the following $\nbigr_{X_0}$-triples
for any $u\in\real\times\cnum$:
\[
 \psitilde_{t,\gminia,u}(\nbigt):=
\psitilde_{t_n,u}\bigl(
 \varphi_n^{\dagger}(\nbigt)\otimes
 \gbigl(-\gminia)\bigr)
\]
\index{functor $\psitilde_{t,\gminia,u}$}
We have the well defined map:
\[
 \nbign=\bigl(-\sqrt{-1}N,\sqrt{-1}N\bigr)
 :\psitilde_{t,u,\gminia}(\nbigt)\lrarr
 \psitilde_{t,u,\gminia}(\nbigt)
 \otimes\Tate^S(-1)
\]
We have the monodromy weight filtration $W(\nbign)$
of $\nbign$ in the category of 
$\nbigr_{X_0}$-triples.

\vspace{.1in}

Note that we have the natural isomorphism
$(\id,\id):\gbigl(\gminia)^{\ast}\simeq\gbigl(\gminia)$.
A Hermitian sesqui-linear duality 
$\nbigs:\nbigt\lrarr\nbigt^{\ast}(-w)$
of weight $w$
naturally induces
$\varphi_n^{\dagger}\nbigt\otimes \gbigl(-\gminia)
\lrarr
\bigl(
 \varphi_n^{\dagger}
 \nbigt\otimes \gbigl(-\gminia)
 \bigr)^{\ast}(-w)$,
which is also a Hermitian sesqui-linear duality
of weight $w$.
If $\varphi_n^{\dagger}(\nbigt)\otimes \gbigl(-\gminia)$ is
strictly specializable along $t_n$,
we have the induced morphisms:
\[
 \psitilde_{t,u,\gminia}(\nbigs):
 \psitilde_{t,u,\gminia}(\nbigt)\lrarr
 \psitilde_{t,u,\gminia}(\nbigt)^{\ast}(-w)
\]
\[
 \psitilde_{t,u,\gminia}(\nbigs):
 \Gr^{W(\nbign)}_{\ell}\psitilde_{t,u,\gminia}(\nbigt)
\lrarr
 \Gr^{W(\nbign)}_{-\ell}\psitilde_{t,u,\gminia}(\nbigt)^{\ast}(-w)
\]
Therefore, we have the induced Hermitian sesqui-linear duality:
\[
\nbigs_{t,u,\gminia,\ell}:=
 \psitilde_{t,u,\gminia}(\nbigs)\circ
 \nbign^{\ell}:
 P\Gr_{\ell}^{W(\nbign)}\psitilde_{t,u,\gminia}(\nbigt)
\lrarr
 P\Gr_{\ell}^{W(\nbign)}\psitilde_{t,u,\gminia}(\nbigt)^{\ast}(-w-\ell)
\]

\subsection{Push-forward}

We put $Y:=Y_0\times\cnum_t$.
Let $F_0:X_0\lrarr Y_0$ be a morphism,
and let $F:X\lrarr Y$ denote the induced morphism.
Let $\nbigm$ be an $\nbigr_X(\ast t)$-module.
We assume that $F$ is proper on the support of $\nbigm$.
We have the well defined push-forward
$F_{\dagger}\nbigm$ in
the derive category of $\nbigr_Y(\ast t)$-modules.
The $i$-th cohomology sheaves are denoted by
$F^i_{\dagger}\nbigm$.
Recall the following compatibility
in \cite{sabbah2} and \cite{sabbah5}.
(See also \cite{saito1}.)

\begin{prop}
\label{prop;07.10.26.100}
Assume that $\nbigm$ is strictly specializable
along $t$.
If $F^i_{0\,\dagger}\psitilde_{t,u}(\nbigm)$ are strict
for any $u\in \real\times\cnum$ and any $i$,
then $F^i_{\dagger}(\nbigm)$ are also strictly specializable
along $t=0$,
and $\psitilde_{t,u}F^i_{\dagger}(\nbigm)
 \simeq F^i_{0\dagger}\psitilde_{t,u}(\nbigm)$.
\hfill\qed
\end{prop}

Let $Y^{(n)}:=Y_0\times\cnum_{t_n}$.
The induced morphism $Y^{(n)}\lrarr Y$ is
also denoted by $\varphi_n$.
Recall the following compatibility 
in \cite{sabbah5},
which can be shown directly.
\begin{prop}
For any $\gminia\in\cnum[t_n^{-1}]$
and for any holonomic
$\nbigr_X(\ast t)$-module $\nbigm$,
we have
$F_{\dagger}\bigl(
 \varphi_n^{\dagger}(\nbigm)
 \otimes\nbigl(-\gminia)
\bigr)
\simeq
 \varphi_n^{\dagger}F_{\dagger}(\nbigm)
 \otimes\nbigl(-\gminia)$.
\hfill\qed
\end{prop}

Let $\nbigt=(\nbigm',\nbigm'',C)$ be an $\nbigr_X$-triple
(resp. $\nbigr_X(\ast t)$-triple).
We have the well defined push-forward
$F_{\dagger}(\nbigt)=
 \bigl(F_{\dagger}(\nbigm')^{\op},F_{\dagger}(\nbigm''),
 F_{\dagger}C\bigr)$
in the derived category of $\nbigr_Y$-triples
(resp. $\nbigr_Y(\ast t)$-triples).
The $i$-th cohomology triples are denoted by
$F^i_{\dagger}(\nbigt)
=\bigl(F^{-i}_{\dagger}\nbigm',
 F^i_{\dagger}\nbigm'',
 F^i_{\dagger}C\bigr)$.

Recall the compatibility of 
the specialization and the push-forward 
for the sesqui-linear pairing
in \cite{sabbah2}.
\begin{prop}
\label{prop;08.1.15.30}
Assume that both $\nbigm'$ and $\nbigm''$
satisfy the assumption in Proposition 
{\rm\ref{prop;07.10.26.100}}.
Then, we have the compatibility
$\psitilde_{t,u}(F^i_{\dagger}C)
=F^i_{0\dagger}(\psitilde_{t,u}C)$.
Namely,
$\psitilde_{t,u}F^i_{\dagger}\nbigt
\simeq
 F^i_{0\dagger}(\psitilde_{t,u}\nbigt)$.

We also have 
$F^i_{\dagger}\bigl(
 \varphi_n^{\dagger}\nbigt\otimes\gbigl(-\gminia)\bigr)
\simeq
 \varphi_n^{\dagger}F_{\dagger}^i(\nbigm)
\otimes\gbigl(-\gminia)$.
\hfill\qed
\end{prop}

\subsection{Specialization 
 of $\nbigr_X$-triples along a function $g$
 with ramification and exponential twist}

Let $g$ be a non-trivial function on $X$.
Let $i_g:X\lrarr X\times\cnum$ denote the graph. 
Let $\nbigt=(\nbigm,\nbigm,C)$ be an $\nbigr_X$-triple.
We have the $\nbigr_{X\times\cnum}$-triple 
$i_{g\,\dagger}\nbigt$,
and the induced $\nbigr_{X\times\cnum_t}(\ast t)$-triple
$(i_{g\,\dagger}\nbigt)(\ast t) $.
Let $\gminia\in\cnum[t_n^{-1}]$.
We say that $\nbigt$ is 
strictly specializable along $g$
with ramification and exponential twist by $\gminia$,
if $\varphi_n^{\dagger}\bigl(i_{g\,\dagger}\nbigt(\ast t) \bigr)
 \otimes\nbigl(-\gminia)$ is strictly specializable
along $t_n$.
In that case,
we define
\[
 \psitilde_{g,\gminia,u}(\nbigt):=
 \psitilde_{t,\gminia,u}
 \bigl(i_{g\dagger}\nbigt(\ast t)\bigr)
\]
for any  $u\in\real\times\cnum$.
\index{functor $\psitilde_{g,\gminia,u}$}
If a Hermitian sesqui-linear duality
$\nbigs:\nbigt\lrarr\nbigt^{\ast}(-w)$ is given,
we have the induced Hermitian sesqui-linear duality
\[
\nbigs_{g,\gminia,u,\ell}:
 P\Gr^{W(\nbign)}_{\ell}
 \bigl(\psitilde_{g,\gminia,u}(\nbigt)\bigr)
\lrarr 
 P\Gr^{W(\nbign)}_{\ell}
 \bigl(\psitilde_{g,\gminia,u}(\nbigt)\bigr)^{\ast}(-w-\ell).
\]

\section{Comparison}
\label{section;10.6.6.10}
\subsection{Twistor structure and 
$\nbigr$-module in dimension $0$}
\label{subsection;10.5.20.21}

Recall that an $\nbigr$-triple in dimension $0$
is a tuple $\nbigt=(\nbigm',\nbigm'',C)$ 
of coherent $\nbigo_{\cnum_{\lambda}}$-modules
$\nbigm'$ and $\nbigm''$ 
with an 
$\nbigo_{\cnum_{\lambda}^{\ast}}$-homomorphism
\[
 C:\nbigm'_{|\cnum_{\lambda}^{\ast}}
\otimes
 \sigma^{\ast}\nbigm''_{|\cnum_{\lambda}^{\ast}}
\lrarr\nbigo_{\cnum_{\lambda}^{\ast}}.
\]
If $\nbigm'$ and $\nbigm''$ are locally free,
$\nbigt$ is called strict.
If the pairing is perfect, $\nbigt$ is called perfect.

Let $V$ be a twistor structure on $\proj^1$.
We have the associated $\nbigr$-triple
\[
 \Theta(V)=\bigl(
 V_0^{\lor},\sigma^{\ast}(V_{\infty}),C_V
 \bigr),
\]
where $V_0:=V_{|\cnum_{\lambda}}$,
$V_{\infty}:=V_{|\cnum_{\mu}}$
and $C_V$ is the naturally defined 
Hermitian sesqui-linear pairing.
It is easy to observe that 
$\Theta$ gives an equivalence 
of the categories of
algebraic vector bundles on $\proj^1$
and strict perfect $\nbigr$-triples
in dimension $0$,
which preserves tensor products
and direct sums.
We also have
$\Theta(\sigma^{\ast}(V)^{\lor})
=\Theta(V)^{\ast}$.
A perfect strict $\nbigr$-triple
is called pure twistor structure of weight $n$,
if it is isomorphic to $\Theta(V)$
for some pure twistor structure $V$ of weight $n$.

\vspace{.1in}

Recall $\Tate^S(k):=
 \bigl(
 \nbigo_{\cnum_{\lamda}}\,e^{(2k)}_{0},
 \nbigo_{\cnum_{\lamda}}\, e^{(2k)}_{\infty},
 (\sqrt{-1}\lambda)^{-2k}\bigr)$
for a half-integer $k$.
The Tate twist in the category of
$\nbigr$-triples in dimension $0$
is given by the tensor product with $\Tate^S(k)$.
The canonical isomorphism
$\Tate^{S}(k)\lrarr
 \Tate^S(-k)^{\ast}$
is given by $\bigl((-1)^{2k},1\bigr)$.
We take a complex number $a$
such that $a^2=-\sqrt{-1}$.
We fix the isomorphism
\[
 \Phi^{(n)}=\bigl(\Phi^{(n)\prime},
 \Phi^{(n)\prime\prime}\bigr):
 \Theta(\nbigo_{\proj^1}(n))
\lrarr \Tate^S(-n/2),
\]
given by 
$\Phi^{(n)\prime}(e_0^{(-n)})
=a^n\,f_0^{(-n)}$
and 
$\Phi^{(n)\prime\prime}
 \bigl(\sigma^{\ast}(f_{\infty}^{(n)})\bigr)
=a^{-n}\, e^{(-n)}_{\infty}$.
Note the commutativity 
of the following induced diagram
(Lemma 3.149 of \cite{mochi2}):
\begin{equation}
 \label{eq;08.1.29.5}
 \begin{CD}
 \Theta\bigl(\nbigo_{\proj^1}(n)\bigr)^{\ast}
@>>>
 \Theta\bigl(\nbigo_{\proj^1}(-n)\bigr) \\
 @AAA @VVV \\
 \Tate^S(-n/2)^{\ast}
 @>>>
 \Tate^S(n/2)
 \end{CD}
\end{equation}
The functor $\Theta$ preserves the Tate twists
in the both categories,
and it also preserves the compatibility 
between the Tate twist and the adjunction
because of the commutativity of (\ref{eq;08.1.29.5}).

For an integer $n$, a morphism
$\nbigs:\nbigt\lrarr\nbigt^{\ast}\otimes\Tate^S(-n)$
is called a Hermitian sesqui-linear duality of weight $n$,
if $(-1)^n\nbigs^{\ast}=\nbigs$ is satisfied.
Note that a morphism
$S:V\otimes\sigma^{\ast}(V)\lrarr \Tate(-n)$
induces 
$\Psi_S:V\lrarr \sigma^{\ast}(V)^{\lor}
 \otimes\Tate(-n)$,
and hence
$\Theta(\Psi_S):
 \Theta(V)\lrarr \Theta(V)^{\ast}\otimes\Tate^S(-n)$.

\begin{lem}
$S$ is a pairing of weight $n$,
if and only if 
$\Theta(\Psi_S)$ 
is a Hermitian sesqui-linear duality
of weight $n$.
\end{lem}
\pf
$S$ is a pairing of weight $n$,
if and only if
$(1\otimes\iota)\circ\sigma^{\ast}\Psi_S^{\lor}
=(-1)^n\Psi_S$.
By the functor $\Theta$,
it is transferred to the condition
that $\Theta(\Psi_S)$ is 
a Hermitian sesqui-linear duality 
of weight $n$.
\hfill\qed

\begin{rem}
Let $(V^{\sankaku},\DD^{\sankaku})$
be a variation of twistor structure
on $\proj^1\times X$,
obtained as the gluing of
$(V_0,\DD_0)$ and 
$(V_{\infty},\DD_{\infty}^{\dagger})$.
Then, we obtain an $\nbigr_X$-triple
$(V_0^{\lor},\sigma^{\ast}V_{\infty},C_V)$.
\hfill\qed
\end{rem}

Let $\nbigt$ be a pure twistor structure 
of weight $n$.
A Hermitian-sesqui-linear duality
$\nbigs:\nbigt\lrarr\nbigt^{\ast}\otimes\Tate^S(-n)$
of weight $n$ is a polarization,
if $(\nbigt,\nbigs)$ is isomorphic to
$\bigl(\Theta(V),\Theta(\Psi_S)\bigr)$
for some polarized pure twistor structure
$(V,S)$ of weight $n$
(Corollary 3.151 of \cite{mochi2}).

\subsection{Comparison of specializations}
\label{subsection;08.9.4.21}

Let $X=\Delta^n$ and $D=\bigl\{O\bigr\}$.
We put $\nbigx:=\cnum_{\lambda}\times X$
and $\nbigd:=\cnum_{\lambda}\times D$.
Let $\nbigt=(\nbigm_1,\nbigm_2,C)$ be 
a regular coherent strictly specializable
$\nbigr_X(\ast z)$-triple on $X$ 
such that $\nbigs=(\id,\id)$ 
is a Hermitian sesqui-linear duality
of weight $0$.
Note that $\nbigm_i$ has the KMS-structure
at each $\lambda_0\in\cnum_{\lambda}$.
Assume that
the restriction $(\nbigt,\nbigs)_{|X\setminus D}
=\Theta(V^{\lor},\DD^{\sankaku},\nbigs_1)$
for a variation of twistor structure
$(V,\DD^{\sankaku},\nbigs)$
with a pairing of weight $0$
on $\proj^1\times(X\setminus D)$.
We use the identifications
$V_0=\nbigm_1$
and $\sigma^{\ast}V^{\lor}_{\infty}=\nbigm_2$.

As explained in Section \ref{subsection;08.1.28.10},
for any $u\in\real\times\cnum$,
we obtain the vector bundle
$S_u(V)$ with the nilpotent map
and the symmetric pairing 
\[
 \nbign^{\sankaku}_{u}:
 S_u(V)\lrarr S_u(V)\otimes\Tate(-1),
\quad
 \nbigs_u:S_u(V)\otimes\sigma^{\ast}S_u(V)
 \lrarr \Tate(0).
\]
We have the specialization
$\psitilde_{z,u}(\nbigt)
=\bigl(\psitilde_{z,u}\nbigm,
 \psitilde_{z,u}\nbigm,
 \psitilde_{z,u}C\bigr)$ 
 for any $u\in T$
with the induced nilpotent map
\[
 \nbign_u=\bigl(-\sqrt{-1}N,\sqrt{-1}N\bigr):
 \psitilde_{z,u}(\nbigt)\lrarr
 \psitilde_{z,u}(\nbigt)
 \otimes\Tate^S(-1),
\]
where $N=-\deldel_tt+\eigenmap(\lambda,u)$
is the nilpotent part of $-\deldel_tt$.
We also have the induced
Hermitian sesqui-linear duality
of weight $0$
\[
 \psitilde_{z,u}(\nbigs):
 \psitilde_{z,u}(\nbigt)\lrarr
 \psitilde_{z,u}(\nbigt)^{\ast}.
\]
\begin{prop}
\label{prop;10.5.20.10}
We have a natural isomorphism
\[
 \Theta\bigl(S_{u}(V),\nbign_{u}^{\sankaku},
 \nbigs_{u}\bigr)
\simeq 
 \bigl(
 \psitilde_{z,u}(\nbigt),
 \nbign_u,
 \psitilde_{z,u}(\nbigs)
 \bigr).
\]
\end{prop}

We give some consequences.
By taking the primitive part 
with respect to the weight filtration of $\nbign_u$,
we obtain the $\nbigr$-triples
$P\Gr^{W}_{\ell}\psitilde_{z,u}(\nbigt)$ 
for any $u\in T$ and $\ell\in\seisuu$
with the Hermitian sesqui-linear duality
\[
 \nbigs_{u,\ell}=
 \psitilde_{z,u}(\nbigs)\circ \nbign_u^{\ell}:
 P\Gr^{W}_{\ell}\psitilde_{z,u}(\nbigt)
\lrarr
 P\Gr^W_{\ell}\psitilde_{z,u}(\nbigt)^{\ast}
 \otimes\Tate^S(-\ell).
\]
The following fact implicitly appeared
in \cite{mochi2}.

\begin{prop}
\label{prop;08.1.28.120}
$\bigl(S_{u}(V),\nbign_{u}^{\sankaku},
 \nbigs_{u}\bigr)$
is a polarized mixed twistor structure of weight $0$,
if and only if
$\bigl(P\Gr^W_{\ell}\psitilde_{z,u}(\nbigt),
\nbigs_{u,\ell}\bigr)$ 
 are polarized pure twistor structure of weight $\ell$
for any $\ell\in\seisuu_{\geq \,0}$.
\end{prop}
\pf
Recall that 
$\bigl(S_{u}(V),\nbign_u^{\sankaku},
\nbigs_u\bigr)$
is a polarized mixed twistor structure,
if (i) $\Gr^{W}_{\ell}S_u(V)$ are 
 pure twistor structure of weight $\ell$
 for any $\ell\in\seisuu$,
(ii) the induced pairings
$\nbigs_u\bigl(
 \nbign_u^{\ell}\otimes\id
 \bigr)$ on 
$P\Gr^{W}_{\ell}S_u(V)$
$(\ell\geq 0)$
give polarizations.
Here, $W$ denote the weight filtration
associated to $\nbign$.
Then, Proposition \ref{prop;08.1.28.120}
follows from Proposition \ref{prop;10.5.20.10}.
\hfill\qed

\begin{cor}
\label{cor;08.1.29.2}
$(V,\DD^{\sankaku},\nbigs)$ is 
a variation of polarized pure twistor structure,
if $\bigl(
 P\Gr^W_{\ell}\psitilde_{z,u}(\nbigt),\nbigs_{u,\ell}
\bigr)$ are pure twistor structures
for any $u\in\real\times\cnum$ and $\ell\in\seisuu_{\geq\,0}$.
It comes from a tame harmonic bundle.
\end{cor}
\pf
The first claim follows from Proposition
\ref{prop;08.1.28.120}
and Lemma \ref{lem;07.7.11.12}.
We obtain the tameness
from the regularity of 
the meromorphic Higgs field.
\hfill\qed

\subsubsection{Proof of Proposition
\ref{prop;10.5.20.10}}

Let $\uhat=u+(1,0)\in\real\times\cnum$.
We have the natural identification
$\nbigg_{\uhat}(V_0)
\simeq
 \psitilde_{z,u}(\nbigm_1)$
by construction.
We also have the natural identification
$\nbigg_{\uhat}(V_0)
\simeq
 \nbigg_{u}(V_0)$
given by the multiplication of $z$.
We will not distinguish them in the following argument.
We also have natural identification
$\psitilde_{z,u}(\nbigm_2)\simeq
 \nbigg_u\bigl(\sigma^{\ast}V_{\infty}^{\lor}\bigr)
\simeq
 \sigma^{\ast}\nbigg_{-u}\bigl(
 V_{\infty}^{\lor}\bigr)
\simeq
 \sigma^{\ast}\nbigg_u(V_{\infty})$.
Note Lemma \ref{lem;08.1.29.1}
for the signature of the KMS-spectrum.

By construction,
$\Theta(S_u(V))$ is a tuple of
$\nbigo_{\cnum_{\lambda}}$-modules
$\nbigg_u(V_0)$ and 
$\sigma^{\ast}\nbigg_u(V_{\infty})^{\lor}
\simeq
 \nbigg_{u}(\sigma^{\ast}(V_{\infty}^{\lor}))$,
and an induced pairing $C_{0,u}$.
We shall show that $C_{0,u}=\psi_uC$
under the above identification,
which implies
$\Theta(S_u(V))\simeq\psitilde_{z,u}(\nbigt)$.

We have only to compare them
around generic
$\lambda_0\in\cnum_{\lambda}^{\ast}$.
In the following,
$U(\lambda_0)$ denotes a neighbourhood of $\lambda_0$,
and $\nbigx^{(\lambda_0)}:=U(\lambda_0)\times X$
and $\nbigd^{(\lambda_0)}:=U(\lambda_0)\times D$.
We put
$\nbigx^{(\sigma(\lambda_0))}
:=\sigma(U(\lambda_0))\times X$
and $\nbigd^{(\sigma(\lambda_0))}:=
 \sigma(U(\lambda_0))\times D$.
We have the generalized eigen decomposition
of the monodromy 
along the loop with the counter-clockwise direction:
\[
 \nbigp^{(\lambda_0)}_{a}(V_0)
=\bigoplus_{a-1<\paramap(\lambda_0,u)\leq a}
 \nbigp^{(\lambda_0)}_{a}
 (V_0)_{\eigenmap^f(\lambda,u)}
\]
Here, eigenvalues of the monodromy
on $\nbigp^{(\lambda_0)}_a
 (V_0)_{\eigenmap^f(\lambda,u)}$
are $\eigenmap^f(\lambda,u)$.
Note that
$\nbigp^{(\lambda_0)}_{\paramap(\lambda,u)}
 (V_0)_{\eigenmap^f(\lambda,u)}$
is naturally isomorphic to
$\nbigg_u\nbigv_{0|\nbigx^{(\lambda_0)}}$
in Section \ref{subsection;08.1.28.110}.
We will use the identification in the following.
We also have the generalized eigen decomposition
of $\nbigp^{(\sigma(\lambda_0))}_aV_0$
of the monodromy along the loop.
It is easy to observe that
the sesqui-linear pairing is compatible 
with the decompositions.
Hence,
we obtain the sesqui-linear pairing
$C_{u}$ of
$\nbigg_u\nbigv_{0|
 \nbigx^{(\lambda_0)}\setminus\nbigd^{(\lambda_0)}}$
and
$\sigma^{\ast}
 \nbigg_u\nbigv_{0|
 \nbigx^{(\sigma(\lambda_0))}
\setminus\nbigd^{(\sigma(\lambda_0))}}$
to the sheaf of $\lambda$-holomorphic
$C^{\infty}$-functions on 
$\nbigx^{(\lambda_0)}\setminus
 \nbigd^{(\lambda_0)}$.

\vspace{.1in}
Let $f\in\nbigg_{\uhat}(V_0)$ 
on $U(\lambda_0)$
and $g\in\sigma^{\ast}
 \nbigg_{\uhat}(V_{\infty})^{\lor}$
on $\sigma(U(\lambda_0))$.
Let $F_0$ be the section of
$\nbigg_{\uhat}\nbigv_0$ on $\nbigx^{(\lambda_0)}$
such that $F_{0|\nbigd^{(\lambda_0)}}=f$
and $(-\deldel_zz+\eigenmap(\lambda,u))^LF_0=0$
for some large $L$.
Such $F_0$ is called the lift of $f$
to $\nbigg_{\uhat}\nbigv_0$.
Similarly,
we take the lift $G_0$ of $g$ to
$\sigma^{\ast}
 \nbigg_{\uhat}\nbigv_{\infty}^{\lor}$ on
$\nbigx^{(\sigma(\lambda_0))}$.
Let $\omega_0=\sqrt{-1}(2\pi)^{-1}\, dz\, d\zbar$.
Let $\chi$ be a test function on $X$
which is constantly $1$ around $0$.
Hence, we obtain the following:
\begin{equation}
 \label{eq;08.1.28.20}
 \psitilde_{z,u}C_{u}\bigl(f,\sigma^{\ast}g\bigr)
=
\underset{s+\lambda^{-1}\eigenmap(\lambda,u)}{\Res}
 \Bigl\langle
C_u\bigl(F_0,
 \sigma^{\ast}G_0\bigr),\,
|z|^{2s}\,\chi\,\omega_0
 \Bigr\rangle
\end{equation}
By the elementary argument in the proof of 
Lemma 18.11 of \cite{mochi2},
we can show that the function
$C_u\bigl(F_0,\sigma^{\ast}G_0\bigr)$
is of the following form:
\begin{equation}
 \label{eq;10.5.20.20}
 |z|^{2\lambda^{-1}\eigenmap(\lambda,u)-2}
 \sum_{k=0}^M a_k(\lambda)\, (\log|z|^2)^{k}
\end{equation}
Here, $a_k$ are holomorphic functions of $\lambda$,
and $M$ denotes some positive integer.
And, we can show the following equality
by a direct calculation
(see Lemma 2.34 of \cite{mochi2} for more general case):
\[
\underset{s+\lambda^{-1}\eigenmap(\lambda,u)}{\Res}
 \Bigl\langle
C_u\bigl(F_0,\sigma^{\ast}G_0\bigr),\,
|z|^{2s}\,\chi\,\omega_0\Bigr\rangle
=a_0(\lambda)
\]

Let $F$ be the section of 
$\nbigg_{\uhat}\nbigh(V_0)$
on $U(\lambda_0)$
such that $\Phi^{\can}(F)=f$.
If we regard $F$ as a multi-valued
flat section of $\nbigg_{\uhat}\nbigv_0$,
we have an expression
\[
 F=\exp\bigl(-\lambda^{-1}
 \eigenmap(\lambda,\uhat)\,
 \log z\bigr)\,
\Bigl(
 F_0+\sum_{j=1}^NF_{j}(\log z)^{j},
\Bigr)
\]
such that
 $(-\deldel_zz+\eigenmap(\lambda,u))^L
 F_{j}=0$ for some large $L$.
Similarly,
let $G$ be the section of
$\sigma^{\ast}\nbigg_{\uhat}
 \nbigh(V_{\infty})^{\lor}$
on $\sigma(U(\lambda_0))$
such that $\Phi^{\can}(G)=f$.
We have similar expression for $G$,
i.e.,
\[
 G=\exp\bigl(-\lambda^{-1}
 \eigenmap(\lambda,\uhat)\log z\bigr)
 \Bigl(
 G_0+\sum_{j=1}^{N}
 G_{j}\,(\log z)^j
\Bigr)
\]
such that
$(-\deldel_zz+\eigenmap(\lambda,u))^L
 G_{j}=0$ for some large $L$.
Because $C_u\bigl(F_{i},\sigma^{\ast}G_j\bigr)$
has an expression as in (\ref{eq;10.5.20.20}),
$C_u\bigl(
 F,\sigma^{\ast}G \bigr)$
is described as a polynomial in $(\log |z|^2)$
whose coefficients are holomorphic functions
in $\lambda$.
Because $C_u\bigl(
 F,\sigma^{\ast}G \bigr)$
is constant in the $X$-direction,
we obtain that
\[
 C_{0,u}(f,\sigma^{\ast}g)=C_u\bigl(
 F,\sigma^{\ast}G \bigr)=a_0(\lambda). 
\]
Hence,
we obtain $\psitilde_{z,u}C=C_{0,u}$.

It is easy to compare the Hermitian sesqui-linear
dualities under the identification.
Let us compare the nilpotent maps.
Let $\Res^{\nil}$ denote the nilpotent part
of the residue.
We have
\[
 \Theta(\nbign_u^{\sankaku})
=\bigl(
 \Res^{\nil}(\DD_0^{\lor})^{\lor}\otimes t_0^{(-1)},\,
 \sigma^{\ast}\Res^{\nil}(\DD_{V_{\infty}^{\lor}})
 \otimes\sigma^{\ast}t_{\infty}^{(-1)}
 \bigr)
\]
We have 
$\Res^{\nil}(\DD_0^{\lor})^{\lor}\otimes t_0^{(-1)}
=-\Res^{\nil}(\DD_0)\otimes t_0^{(-1)}$.
Under the identification of
$\Theta^S(n)$ and $\Theta(n)$
in Subsection \ref{subsection;10.5.20.21},
it is identified with
$-\sqrt{-1}N\otimes e_0^{(2)}$,
where $N$ is the nilpotent part of $-\deldel_tt$.
We have 
$\sigma^{\ast}\Res^{\nil}(\DD_{V_{\infty}^{\lor}})
\otimes\sigma^{\ast}t_{\infty}^{(-1)}
=-\Res^{\nil}(\DD_{\sigma^{\ast}V^{\lor}_{\infty}})
 \otimes\sigma^{\ast}t_{\infty}^{(-1)}$,
which is identified with
$\sqrt{-1}N\otimes e_0^{(-2)}$.
Hence, we have
$\nbign_u^{\sankaku}=\nbign_u$.
Thus, the proof of Proposition \ref{prop;10.5.20.10}
is finished.
\hfill\qed

\backmatter

\printindex

\end{document}